\documentclass{book}
\usepackage{amsmath,amssymb,amscd,amsthm}
\usepackage{amsfonts}

\usepackage{tabularx}
\usepackage{graphicx}
\usepackage{multicol}
\usepackage{enumerate}

\usepackage{tikz}
\usepackage[all]{xy}
\xyoption{poly}
\xyoption{arc}
\xyoption{curve}

\CompileMatrices

\usepackage{imakeidx}

\begin{document}

\newcommand{\nc}{\newcommand}

\nc{\pr}{\noindent{\em Proof. }} 
\nc{\g}{\mathfrak g}
\nc{\n}{\mathfrak n} 
\nc{\opn}{\overline{\n}}
\nc{\h}{\mathfrak h}
\renewcommand{\b}{\mathfrak b}
\nc{\Ug}{U(\g)} 
\nc{\Uh}{U(\h)} 
\nc{\Un}{U(\n)}
\nc{\Uopn}{U(\opn)}
\nc{\Ub}{U(\b)} 
\nc{\p}{\mathfrak p}
\renewcommand{\l}{\mathfrak l}
\nc{\z}{\mathfrak z} 
\renewcommand{\h}{\mathfrak h}
\nc{\m}{\mathfrak m}
\renewcommand{\k}{\mathfrak k}
\nc{\opk}{\overline{\k}}
\nc{\opb}{\overline{\b}}
\nc{\e}{{\epsilon}}
\nc{\ke}{{\bf k}_\e}
\nc{\Hk}{{\rm Hk}^{\gr}(A,A_0,\e )}
\nc{\gr}{\bullet}
\nc{\ra}{\rightarrow}
\nc{\Alm}{A-{\rm mod}}
\nc{\DAl}{{D}^-(A)}
\nc{\HA}{{\rm Hom}_A}

\nc{\ha}{\mathfrak h}
\nc{\ba}{\mathfrak b}

\nc{\la}{\mathfrak l}
\nc{\ta}{\mathfrak t}
\nc{\ka}{\mathfrak k}

\newtheorem{theorem}{Theorem}{}
\newtheorem{lemma}[theorem]{Lemma}{}
\newtheorem{corollary}[theorem]{Corollary}{}
\newtheorem{conjecture}[theorem]{Conjecture}{}
\newtheorem{proposition}[theorem]{Proposition}{}
\newtheorem{axiom}{Axiom}{}
\newtheorem{remark}[theorem]{Remark}{}
\newtheorem{example}{Example}{}
\newtheorem{exercise}{Exercise}{}
\newtheorem{definition}[theorem]{Definition}{}

\title{Q-W--algebras, Zhelobenko operators and a proof of  \\ De Concini--Kac--Procesi conjecture}

\author{A. Sevostyanov}




\maketitle


\chapter*{Introduction}
\addcontentsline{toc}{chapter}{Introduction}

\setcounter{equation}{0}


\subsection*{Schur--Weyl duality philosophy in the Hecke algebra setting}\label{phil}

In this monograph we consider a phenomenon which occurs in the study of certain classes and categories of representations of semisimple Lie algebras, groups of Lie type, and the related quantum groups. This phenomenon is similar to the classical Schur--Weyl duality. However, the relevant classes of representations are quite different from the finite-dimensional irreducible representations of the general linear group or, more generally, of complex semisimple Lie groups which appear in the Schur--Weyl setting. 

All examples of the above mentioned type are realizations of the following quite general construction a homological version of which was suggested in \cite{S1} (in fact, in the quantum group case the construction presented below requires some technical modifications; we shall not discuss them in the introduction). Let $A$ be an associative algebra over a unital ring $\bf k$, $B\subset A$ a subalgebra with a character $\chi:B\rightarrow {\bf k}$. Denote by ${\bf k}_\chi$ the corresponding rank one representation of $B$. Let $Q_\chi=A\otimes_B {\bf k}_\chi$ be the induced representation of $A$. 

Let ${\rm Hk}(A,B,\chi)={\rm End}_A(Q_\chi)^{opp}$ be the algebra of $A$--endomorphisms of $Q_\chi$ with the opposite multiplication. One says that the algebra ${\rm Hk}(A,B, \chi)$ is obtained from $A$ by a quantum constrained reduction with respect to the subalgebra $B$. ${\rm Hk}(A,B, \chi)$ is an algebra of Hecke type. Indeed, if $A$ is the group algebra of a Chevalley group over a finite field, $B$ the group algebra  of a Borel subgroup in it, and $\chi$ is the trivial complex representation of the Borel subgroup one obtains the Iwahori--Hecke algebra this way (see \cite{I}). 

For any representation $V$ of $A$ the algebra ${\rm Hk}(A,B,\chi)$ naturally acts in the space 
$$
V_\chi:={\rm Hom}_{A}(Q_\chi,V)\simeq {\rm Hom}_{B}({\bf k}_\chi,V)
$$ 
by compositions of homomorphisms and for any ${\rm Hk}(A,B, \chi)$--module $W$ the induced representation $Q_\chi\otimes_{{\rm Hk}(A,B,\chi)}W$ is a left $A$--module. Let ${\rm Hk}(A,B,\chi)-{\rm mod}$ be the category of left ${\rm Hk}(A,B,\chi)$--modules and $A-{\rm  mod}_B^\chi$ the category of left $A$--modules of the form $Q_\chi\otimes_{{\rm Hk}(A,B,\chi)}W$, where $W\in {\rm Hk}(A,B,\chi)-{\rm mod}$, with morphisms induced by morphisms of left ${\rm Hk}(A,B,\chi)$--modules. 
The point is that in many important examples $A-{\rm  mod}_B^\chi$ is a full subcategory in the category of left $A$--modules, and the  functors ${\rm Hom}_{A}(Q_\chi,\cdot): A-{\rm  mod}_B^\chi\to {\rm Hk}(A,B,\chi)-{\rm mod}$ and $Q_\chi\otimes_{{\rm Hk}(A,B,\chi)}\cdot: {\rm Hk}(A,B,\chi)-{\rm mod} \to A-{\rm  mod}_B^\chi$ yield mutually inverse Schur--Weyl type equivalences of the categories,
\begin{equation}\label{catequiv}
A-{\rm  mod}_B^\chi\simeq {\rm Hk}(A,B,\chi)-{\rm mod}.
\end{equation}

In these cases for $W\in {\rm Hk}(A,B,\chi)-{\rm mod}$ one has 
$$
(Q_\chi\otimes_{{\rm Hk}(A,B,\chi)}W)_\chi= {\rm Hom}_A(Q_\chi,Q_\chi\otimes_{{\rm Hk}(A,B,\chi)}W)\simeq {\rm Hom}_B({\bf k}_\chi,Q_\chi\otimes_{{\rm Hk}(A,B,\chi)}W)\simeq W,
$$
and for  $W,W'\in {\rm Hk}(A,B,\chi)-{\rm mod}$ by the formula above 
$$
{\rm Hom}_A(Q_\chi\otimes_{{\rm Hk}(A,B,\chi)}W',Q_\chi\otimes_{{\rm Hk}(A,B,\chi)}W)\simeq {\rm Hom}_{{\rm Hk}(A,B,\chi)}(W',{\rm Hom}_A(Q_\chi,Q_\chi\otimes_{{\rm Hk}(A,B,\chi)}W))\simeq
$$
$$
\simeq {\rm Hom}_{Hk(A,B,\chi)}(W',W).
$$

At first sight the category $A-{\rm  mod}_B^\chi$ looks a bit exotic. But it turns out that in many situations it has alternative descriptions in terms of the algebra $A$, its subalgebra $B$ and the character $\chi$ only, and actually such categories and the related algebras of Hecke type played a very important, if not central, role in representation theory for at least last sixty years.


An important example of equivalences of type (\ref{catequiv}) was considered in \cite{K}. In that paper Kostant showed that algebraic analogues of the principal series representations, called irreducible Whittaker modules, for a complex semisimple Lie algebra $\g$ are in one--to--one correspondence with the one--dimensional representations of the center $Z(U(\g))$ of the enveloping algebra $U(\g)$

In the situation considered in \cite{K} the algebra $Z(U(\g))$ is isomorphic to ${\rm Hk}(A,B,\chi)$ with $A=U(\g)$, $B=U(\n_-)$, where $\n_-$ is a nilradical of $\g$, and $\chi$ being a non--singular character of $U(\n_-)$, i.e. it does not vanish on all simple root vectors in $\n_-$. The category $A-{\rm  mod}_B^\chi$ in this case can be described as the category of $\g$--modules on which $x-\chi(x)$ acts locally nilpotently for any $x\in \n_-$. 

This correspondence was generalized in \cite{Ly} and in the Appendix to \cite{Pr} to a more general categorical setting and the categorical equivalence established in the Appendix to \cite{Pr} is called the Skryabin equivalence. 

Similar equivalences were obtained in \cite{Pr} in the case of semisimple Lie algebras over fields of prime characteristic and in \cite{SDM} in the case of quantum groups associated to complex semisimple Lie algebras  for generic values of the deformation parameter. Various approaches to the proofs of the above mentioned statements have been developed in \cite{GG,Sk}.

An analogous construction appears also in the case of finite groups of Lie type (see \cite{Car1}, Chapter 10) and of finite Chevalley groups (see \cite{Ka1,Ka2}). In the latter case the corresponding modules $Q_\chi$ are called generalized Gelfand–Graev representations.

Note that in general the problem of classification of ${\rm Hk}(A,B,\chi)$--modules is usually very difficult. Sometimes it is easier to classify irreducible objects in the category $A-{\rm  mod}_B^\chi$ and then to translate the result to the category ${\rm Hk}(A,B,\chi)-{\rm mod}$ (see \cite{Los1,Los2} for the case of algebras ${\rm Hk}(A,B,\chi)$ considered in \cite{Pr}). 


\subsection*{A strategy for establishing Schur--Weyl type equivalences in the Hecke algebra setting for Lie algebras and quantum groups, and its relation to Zhelobenko operators}\label{strategy}

In all cases considered in \cite{K,Ly,Pr,SDM} the algebras $A$ and $B$, the characters $\chi$ and the appropriate categories $A-{\rm  mod}_B^\chi$ and ${\rm Hk}(A,B,\chi)-{\rm mod}$ of representations of $A$ and of ${\rm Hk}(A,B,\chi)$ are relatively easy to define. It is much more difficult to obtain alternative descriptions of the category $A-{\rm  mod}_B^\chi$. However, one should note that the approach to this problem in all papers mentioned above is slightly different: all those papers start with the description of a category of $A$--modules in intrinsic terms using the algebra $A$, its subalgebra $B$ and the character $\chi$. And then one proves that this category is equivalent to the category ${\rm Hk}(A,B,\chi)-{\rm mod}$, the equivalence being established using the functors ${\rm Hom}_{A}(Q_\chi,\cdot)$ and $Q_\chi\otimes_{{\rm Hk}(A,B,\chi)}\cdot$. Finally one deduces that this category actually coincides with $A-{\rm  mod}_B^\chi$.

 In the Lie algebra case the most simple proofs of statements of this kind were proposed in the Appendix to \cite{Pr} in the zero characteristic case and in \cite{Sk} in the prime characteristic case. But the phenomenon behind these proofs is already manifest in \cite{K}. Namely, in the case of Lie algebras over fields of zero characteristic one always has $A=U(\g)$ and $B=U(\m)$ for some reductive Lie algebra $\g$ and a nilpotent Lie subalgebra $\m\subset \g$, and the above mentioned phenomenon amounts to introducing a second $U(\m)$--module structure on $Q_\chi$ by tensoring with the one--dimensional representation ${\bf k}_{-\chi}$ and to demonstrating that for ${\bf k}=\mathbb{C}$ a certain ``classical limit'' of the $U(\m)$--module  $Q_\chi\otimes {\bf k}_{-\chi}$ is isomorphic to the algebra of regular functions $\mathbb{C}[\mathcal{C}]$ on a closed algebraic variety $\mathcal{C}$, and the ``classical limit'' of the $U(\m)$--action on $Q_\chi\otimes {\bf k}_{-\chi}$ is induced by a free action of the complex unipotent algebraic group $M$ corresponding to the Lie algebra $\m$ on $\mathcal{C}$. The ``classical limits'' here are understood in the sense of taking associate graded objects with respect to suitable filtrations.

The action $$M\times \mathcal{C}\rightarrow \mathcal{C}$$ has a global cross-section $\Sigma \subset \mathcal{C}$, called a Slodowy slice, so that the action map 
\begin{equation}\label{isom}
M\times \Sigma \rightarrow \mathcal{C}	
\end{equation}
is an isomorphism of varieties, and 
\begin{equation}\label{isof}
\mathbb{C}[\mathcal{C}]\simeq \mathbb{C}[\Sigma]\otimes \mathbb{C}[M].	
\end{equation}
The space $W_0:=\mathbb{C}[\Sigma]\simeq\mathbb{C}[\mathcal{C}]^M$ can be regarded as a ``classical limit''  of ${\rm Hk}(A,B,\chi)$ which is called a W--algebra in this case. We can also write $\mathbb{C}[\mathcal{C}]\simeq W_0[M]$, where $W_0[M]$ is the algebra of regular functions on $M$ with values in $W_0$. In fact $W_0$ carries the natural structure of a Poisson algebra. It is called a Poisson W--algebra.  

Let $A-{\rm  mod}_B^\chi$ be the category of left $A$--modules $V$ for which the $U(\m)$--action on $V\otimes {\bf k}_{-\chi}$ is locally nilpotent.
In the Appendix to \cite{Pr} it is shown that if one equips $V\in A-{\rm  mod}_B^\chi$ with the second $U(\m)$--module structure by tensoring with ${\bf k}_{-\chi}$ then, as a $U(\m)$--module, $V\otimes {\bf k}_{-\chi}$ is isomorphic to ${\rm hom}_{\bf k}(U(\m),V_\chi)$, 
\begin{equation}\label{isochi}
V\otimes {\bf k}_{-\chi}\simeq {\rm hom}_{\bf k}(U(\m),V_\chi)\simeq V_\chi[M],
\end{equation}
where ${\rm hom}_{\bf k}$ stands for the space of homomorphisms vanishing on some power of the augmentation ideal of $U(\m)$, $V_\chi={\rm Hom}_{U(\m)}({\bf k}_\chi,V)$ is called the space of Whittaker vectors in $V$, and, as above, the latter isomorphism holds if ${\bf k}=\mathbb{C}$. In the Appendix to \cite{Pr} it is shown that isomorphisms (\ref{isochi}) directly imply an equivalence between the category of left $A$--modules $V$ for which the $U(\m)$--action on $V\otimes {\bf k}_{-\chi}$ is locally nilpotent and the corresponding category $A-{\rm  mod}_B^\chi$ introduced before formula (\ref{catequiv}).

Isomorphisms of type (\ref{isom}) occur in the quantum group setting as well (see \cite{S6,S10,SDM}), and the same idea is applied in \cite{SDM} to establish similar categorical equivalences in the quantum group case for generic values of the deformation parameter.

In \cite{SZ,S13} it was observed that an isomorphism of type (\ref{isom}) gives rise to a natural projection operator $\Pi: \mathbb{C}[\Sigma]\rightarrow \mathbb{C}[\mathcal{C}]^M\simeq \mathbb{C}[\Sigma]=W_0$.  
Namely, according to (\ref{isom}) any $x\in \mathcal{C}$ can be uniquely represented in the form 
\begin{equation}\label{gslices}
x=n(x)\circ \sigma(x), n(x)\in M, \sigma(x)\in \Sigma.
\end{equation}
If for $f\in \mathbb{C}[\mathcal{C}]$ we define $\Pi f\in \mathbb{C}[\mathcal{C}]$ by
\begin{equation}\label{Pidefs}
(\Pi f)(x)=f(n^{-1}(x)\circ x)=f(\sigma(x))
\end{equation}
then $\Pi f$ is an $M$--invariant function, and any $M$--invariant regular function on $\mathcal{C}$ can be obtained this way. Moreover, by the definition $\Pi^2=\Pi$, i.e. $\Pi$ is a projection onto $\mathbb{C}[\mathcal{C}]^{M}$.
 
In the quantum group setting considered in \cite{S6,S10,SDM} the ``classical limiting'' variety $\mathcal{C}$ is always a closed subvariety in a complex semisimple algebraic Lie group $G$, $\Sigma$ is an analogue of a Slodowy slice for $G$ introduced in \cite{S6,S12}, and $M$ is a unipotent subgroup of $G$, where the ``classical limit'' simply corresponds now to the $q=1$ specialization of the deformation parameter $q$. The peculiarity of the quantum group case is that every element of $M$ can be uniquely represented as an ordered product of elements of some one--parameter subgroups $M_i\subset G$, $i=1,\ldots ,c$ corresponding to roots, i.e. $M=M_1\ldots M_c$. If we denote by $u_i$ the parameter in $M_i$ and  by $X_i(u_i)$ the element of $M_i$ corresponding to the value $u_i\in \mathbb{C}$ of the parameter then factorizing $n(x)$ in (\ref{gslices}) as follows
\begin{equation}\label{factoriz}
n(x)=X_1(u_1(x))\ldots X_c(u_c(x))
\end{equation}
one can express the operator $\Pi$ as a composition of operators $\Pi_i$,
\begin{equation}\label{Pis}
(\Pi_i f)(x)=f(X_i(-u_i(x))\circ x),
\end{equation}
\begin{equation}\label{PCs}
\Pi f=\Pi_1\ldots \Pi_c f .
\end{equation} 
$u_i(x)$ here can be regarded as regular functions on $\mathcal{C}\subset G$.

The first miracle of the quantum group case is that there are explicit formulas for the functions $u_i(x)$ in (\ref{factoriz}) expressing them in terms of matrix elements of finite-dimensional irreducible representations of $G$. These formulas were obtained in \cite{S13}. In fact in this book we use these formulas in a  modified form suitable for quantization.

The main objective of this book is to obtain quantum group counterparts of these formulas. This provides a description of quantum group analogues of W--algebras, called q-W--algebras, as images of operators $\Pi_c^q$ which are quantum analogues of $\Pi$, or, more precisely, of closely related operators $\Pi_c$ introduced by formula (\ref{pic}). This description implies that q-W--algebras belong to the class of the so-called Mickelsson algebras (see e.g. \cite{Z1,Z5,Z6,Z7,Z9} and \cite{Z2}, Ch. 4).

Magically, the appropriate modifications of the classical formulas for $u_i(x)$ and of formulas (\ref{Pis}) can be directly extrapolated to the quantum case, so the operator $\Pi_c^q$ is given in a factorized form similar to (\ref{PCs}). Note that no operators similar to $\Pi_c^q$ can be defined in the Lie algebra setting discussed above.

Using the quantum group analogues $B_{jk}$ of the functions $u_i(x)$ one can also construct natural bases in modules $V$ from the corresponding category $A-{\rm  mod}_B^\chi$ and establish isomorphisms similar to (\ref{isochi}) in the case when the deformation parameter is not a root of unity. Recall that in the Lie algebra case with ${\bf k}=\mathbb{C}$ for any $V\in A-{\rm  mod}_B^\chi$ the Skryabin equivalence provides an isomorphism $V\simeq Q_\chi\otimes_{{\rm Hk}(A,B,\chi)}V_\chi$. If we denote by $V_\chi^0$ the ``classical limit'', i.e. the $q=1$ specialization, of $V_\chi$ then recalling that the  ``classical limit'' of $Q_\chi$ is $\mathbb{C}[\mathcal{C}]$ and the  ``classical limit'' of ${\rm Hk}(A,B,\chi)$ is $\mathbb{C}[\Sigma]$ we infer from (\ref{isof})  that the ``classical limit'' of $Q_\chi\otimes_{{\rm Hk}(A,B,\chi)}V_\chi\simeq V$ is 
$$
\mathbb{C}[\mathcal{C}]\otimes_{\mathbb{C}[\Sigma]}V_\chi^0\simeq (\mathbb{C}[\Sigma]\otimes \mathbb{C}[M])\otimes_{\mathbb{C}[\Sigma]}V_\chi^0\simeq \mathbb{C}[M]\otimes V_\chi^0.
$$
These isomorphisms together with (\ref{gslices}) and (\ref{factoriz}) give a hint how to construct natural bases in modules from the category $V\in A-{\rm  mod}_B^\chi$ in the quantum group case. 
Namely, if $V$ is such a module it is natural to expect that if one picks up a linear basis $v_p$, $p\in \mathbb{N}$ in the space of Whittaker vectors $V_\chi$ then the elements of $V$ given by properly defined ordered monomials in $B_{jk}$ applied to $v_p$, $p\in \mathbb{N}$ form a linear basis in $V$. We show that this is indeed the case. These bases are key ingredients for an alternative proof of the Skryabin equivalence for quantum groups. 

Operators conceptually similar to $\Pi_c^q$ appeared in the literature a long time ago as the projection operators onto subspaces of singular vectors in some modules over a complex finite-dimensional semisimple Lie algebra $\g$ the action of a nilradical $\n_-\subset \g$ on which is locally nilpotent. The first example of such operators, called extremal projection operators, for $\g=\mathfrak{sl}_2$ was explicitly constructed in \cite{L}. In papers \cite{AST1,AST2,AST3} the results of \cite{L} were generalized to the case of arbitrary complex semisimple Lie algebras, and explicit formulas for extremal projection operators were obtained. A summary of these results can be found in \cite{T}. Later, using a certain completion of an extension of the universal enveloping algebra of $\g$, Zhelobenko observed in \cite{Z1} that the existence of extremal projection operators is an almost trivial fact. In \cite{Z1} he also introduced a family of operators which are analogous to our operators $P_i$. These operators are called now Zhelobenko operators. Properties of extremal projection operators and of the Zhelobenko operators have been extensively studied in \cite{Z1}--\cite{Z9}, and the results obtained in these papers were summarized in book \cite{Z2}. 

In our terminology the situation considered in these works corresponds to the case when $A=U(\g)$, $B=U(\n_-)$ and $\chi$ is the trivial character of $U(\n_-)$. As observed in \cite{SZ}, in this case $\mathcal{C}=\b_-$, the Borel subalgebra $\b_-\subset \g$ containing $\n_-$, and the action of the unipotent group $N_-$ corresponding to $\n_-$ on $\mathcal{C}$ is induced by the adjoint action of a Lie group $G$ with the Lie algebra $\g$ on $\g$. This action is not free but is gives rise to a birational equivalence
$$
N_-\times \h \rightarrow \b_-,
$$
where $\h=\b_-/\n_-$ is a Cartan subalgebra. In \cite{SZ} it is shown that using this birational equivalence one can still define operators similar to $\Pi_i$ and $\Pi$ acting on a certain localization of the algebra of regular functions $\mathbb{C}[\b_-]$ and these operators are ``classical limits'' of the Zhelobenko and of the extremal projection operators, respectively. 


\subsection*{Kac--Weisfeiler and De Concini--Kac--Procesi conjectures}\label{conjKWDKP}

Remarkably, as observed in \cite{Sk}, the arguments from the Appendix to \cite{Pr} are applicable to obtain alternative descriptions of the corresponding categories $A-{\rm  mod}_B^\chi$ from \cite{Pr} for Lie algebras over fields of prime characteristic. 

Along the same line, formulas for the quantum group analogues $B_{jk}$ of the functions $u_i(x)$ and for the operator $\Pi_c^q$ can be specialized to the case when $q$ is a primitive odd $\bar{m}$-th root of unity $\varepsilon$ subject to a few other conditions depending on the Cartan matrix of the corresponding semisimple Lie algebra $\g$. This provides technical tools for the proof of a root of unity version of the Skryabin equivalence for quantum groups. Similarly to the case of generic $q$ one can construct bases in modules $V$ from the corresponding category $A-{\rm  mod}_B^\chi$. In case when $q$ is a root of unity all such irreducible modules are finite-dimensional, and if one picks up a linear basis $v_p$, $p=1,\ldots, n$ in the space of Whittaker vectors $V_\chi$ then the elements of $V$ given by applied to $v_p$, $p=1,\ldots, n$ properly defined ordered monomials in $c$ variables $B_{jk}$ powers of which are truncated at the degree $\bar{m}$ form a linear basis in $V$. In particular, the dimension of $V$ is divisible by $\bar{m}^c$. It turns out that any finite-dimensional module over the standard quantum group $U_\varepsilon(\g)$, where $\varepsilon$ is a primitive odd $\bar{m}$-th root of unity subject to the extra conditions mentioned in the beginning of this paragraph, belongs to one of the categories $A-{\rm  mod}_B^\chi$ with appropriate $A$, $B$ and $\chi$, so its dimension is divisible by $b:=\bar{m}^c$. Moreover, the number $b$ is equal to the number from the De Concini--Kac--Procesi conjecture on dimensions of irreducible modules over quantum groups at roots of unity suggested in \cite{DKP1}. Thus our result confirms this conjecture. Due to its importance we are going to discuss the De Concini--Kac--Procesi conjecture in more detail.

It is very well known that the number of simple modules for a finite-dimensional algebra  is finite. However, often it is very difficult to classify such representations. In some important particular examples even dimensions of simple modules over finite-dimensional algebras are not known.

One of the important examples of that kind is representation theory of semisimple Lie algebras over algebraically closed fields of prime characteristic.
Let $\frak g'$ be the Lie algebra of a semisimple algebraic group $G'$ over an algebraically closed field $\bf k$ of characteristic $p>0$. Let $x\mapsto x^{[p]}$ be the $p$-th power map of $\frak g'$ into itself. The structure of the enveloping algebra of $\frak g'$ is quite different from the zero characteristic case. Namely, the elements $x^{p}-x^{[p]},~x\in \frak g'$ are central.
For any linear form $\theta$ on $\frak g'$, let $U_{\theta}$ be the
quotient of
the enveloping algebra of $\frak g'$ by the ideal generated by the central elements
$x^{p}-x^{[p]}-\theta (x)^{p}$ with $x\in \frak g'$.  Then $U_{\theta }$ is a finite-dimensional algebra. Kac and Weisfeiler proved that any simple $\frak g'$-module can be regarded as a module over $U_{\theta }$ for a unique $\theta$ as above (this explains why all simple $\frak g'$--modules are finite-dimensional). The Kac--Weisfeiler conjecture formulated in \cite{KW} and proved in \cite{Pr1} says that if the $G'$--coadjoint orbit of $\theta$ has dimension ${\rm dim}~{\mathcal{O}}_\theta$ and $p$ is good for the root system of $G'$ then $p^{\frac {{\rm dim}~{\mathcal{O}}_\theta}{2}}$ divides the dimension of every finite-dimensional $U_{\theta }$--module.

One can identify $\theta$ with an element of $\frak g'$ via the Killing form and
reduce the proof of the Kac--Weisfeiler conjecture to the case of nilpotent $\theta$. In that case Premet defines in \cite{Pr1} a subalgebra $U_\theta(\m_\theta)\subset U_{\theta }$ generated by a Lie subalgebra $\m_\theta\subset \g'$ such that $U_\theta(\m_\theta)$ has dimension $p^{\frac {{\rm dim}~{\mathcal{O}}_\theta}{2}}$ and every finite-dimensional $U_{\theta }$--module is $U_\theta(\m_\theta)$--free. Verification of the latter fact uses the theory of support varieties (see \cite{FP1,FP2,FP, Pr3}). Namely, according to the theory of support varieties, in order to prove that a $U_{\theta }$--module is $U_\theta(\m_\theta)$--free one should check that it is free over every subalgebra $U_\theta(x)$ generated in $U_\theta(\m_\theta)$ by a single element $x\in \m_\theta$.

There is a more elementary and straightforward proof of the Kac--Weisfeiler conjecture given in \cite{PS}. The simplest proof of this conjecture follows from the results of \cite{Sk} on a prime characteristic version of the Skryabin equivalence which we already discussed above. A proof of the conjecture for $p>h$, where $h$ is the Coxeter number of the corresponding root system, using localization of $\mathcal{D}$--modules was presented in \cite{BMR}.

Another important example of finite-dimensional algebras is related to the theory of quantum groups at roots of unity. Let $\mathfrak g$ be a complex finite-dimensional semisimple Lie algebra.
A remarkable property of the standard Drinfeld-Jimbo quantum group $U_\varepsilon(\mathfrak g)$ associated to $\mathfrak g$, where $\varepsilon$ is a primitive  $\bar{m}$-th root of unity, is that its center  contains a large commutative subalgebra. In this book we consider the simply connected version of $U_\varepsilon(\mathfrak g)$ and the case when $\bar{m}$ is odd. In this case the large commutative subalgebra is isomorphic to the algebra $Z_G$ of regular functions on (a finite covering of a big Bruhat cell in) the connected, simply connected complex algebraic group corresponding to $\g$.

Consider finite-dimensional representations of $U_\varepsilon(\mathfrak g)$, on which $Z_G$ acts according to non--trivial characters $\eta_g$ given by evaluation of regular functions at various points $g\in G$. Note that all irreducible representations of $U_\varepsilon(\mathfrak g)$ are of that kind, and every such representation is a representation of the algebra $U_{\eta_g}=U_\varepsilon(\mathfrak g)/U_\varepsilon(\mathfrak g){\rm Ker}~\eta_g$ for some $\eta_g$. In \cite{DKP1} De Concini, Kac and Procesi showed that if $g_1$ and $g_2$ are two conjugate elements of $G$ then the algebras $U_{\eta_{g_1}}$ and $U_{\eta_{g_2}}$ are isomorphic. Moreover in \cite{DKP1} De Concini, Kac and Procesi formulated the following conjecture.

\vskip 0.3cm

{\bf De Concini--Kac--Procesi conjecture.}
{\em The dimension of any finite-dimensional representation of the algebra $U_{\eta_g}$ is divisible by $b_g:=\bar{m}^{\frac{1}{2}{\rm dim}~\mathcal{O}_g}$, where $\mathcal{O}_g$ is the conjugacy class of $g$.}

\vskip 0.3cm

This conjecture is the quantum group counterpart of the Kac--Weisfeiler conjecture for semisimple Lie algebras over fields of prime characteristic.

As it is shown in \cite{DK1} it suffices to verify the De Concini--Kac--Procesi conjecture in case of exceptional elements $g\in G$ (an element $g\in G$ is called exceptional if the centralizer in $G$ of its semisimple part has a finite center). However, the De Concini--Kac--Procesi conjecture is related to the geometry of the group $G$ which is more complicated than the geometry of the vector space $\g'$ in the case of the Kac--Weisfeiler conjecture.

The De Concini--Kac--Procesi conjecture is known to be true for the conjugacy classes of regular elements (see \cite{DKP2}), for the subregular
unipotent conjugacy classes in type $A_n$ when $\bar{m}$ is a power of a prime number (see \cite{10}), for all conjugacy classes in $A_n$
when $\bar{m}$ is a prime number (see \cite{8}), for the conjugacy classes $\mathcal{O}_g$ of $g\in  SL_n$ when the conjugacy class of the unipotent part of $g$ is spherical (see \cite{9}), and for spherical conjugacy classes (see \cite{CCC}).
In \cite{KR} a proof of the De Concini--Kac--Procesi conjecture using localization of quantum $\mathcal{D}$--modules was outlined in the case of unipotent conjugacy classes. In contract to many papers quoted above the strategy of the proof of the De Concini--Kac--Procesi conjecture developed in this book does not use the reduction to the case of exceptional elements, and all conjugacy classes are treated uniformly.

Namely, following Premet's philosophy we use certain subalgebras $U_{\eta_g}(\m_-)\subset U_{\eta_g}$ introduced in \cite{S11}. These subalgebras have non--trivial characters $\chi: U_{\eta_g}(\m_-) \rightarrow \mathbb{C}$. In terms of the previously introduced notation, we show that for $A=U_{\eta_g}$, $B=U_{\eta_g}(\m_-)$ and an appropriate $\chi$ the category $A-{\rm  mod}_B^\chi$ can be identified with the category of finite-dimensional representations of $U_{\eta_g}$ and equivalence (\ref{catequiv}) holds if ${\rm Hk}(A,B,\chi)-{\rm mod}$ is the category of finite-dimensional representations of the corresponding algebra ${\rm Hk}(A,B,\chi)$. 

As observed in \cite{S11} every finite-dimensional $U_{\eta_g}$--module is also equipped with an action of the algebra $U_{\eta_1}(\m_-)$ corresponding to the trivial character $\eta_1$ of $Z_G$ given by the evaluation at the identity element of $G$. In the setting of quantum groups at roots of unity this action is a counterpart of the second $U(\m)$--module structure on objects $V$ of the category $A-{\rm  mod}_B^\chi$  which appeared in (\ref{isochi}) in the case of Lie algebras over fields of zero characteristic.

Since the De Concini--Kac--Procesi conjecture is related to the structure of the set of conjugacy classes in $G$ it is natural to look at transversal slices to the set of conjugacy classes. It turns out that the definition of the subalgebras $U_{\eta_g}(\m_-)$ is related to the existence of some special transversal slices $\Sigma_s$ to conjugacy classes in $G$. These slices $\Sigma_s$ associated to (conjugacy classes of) elements $s$ in the Weyl group of $\g$ were introduced by the author in \cite{S6}. The slices $\Sigma_s$ play the role of Slodowy slices in algebraic group theory. In the particular case of elliptic Weyl group elements these slices were also introduced later by He and Lusztig in paper \cite{XL} within a different framework.

A remarkable property of a slice $\Sigma_s$ observed in \cite{S11} is that if $g$ is conjugate to an element in $\Sigma_s$ then the dimension of the corresponding subalgebra  $U_{\eta_g}(\m_-)\subset U_{\eta_g}$ is equal to $\bar{m}^{\frac{1}{2}{\rm codim}~\Sigma_s}$. The dimension of the algebra  $U_{\eta_1}(\m_-)$ is also equal to $\bar{m}^{\frac{1}{2}{\rm codim}~\Sigma_s}$. If $g\in \Sigma_s$ (in fact $g$ may belong to a larger variety) then the corresponding subalgebras $U_{\eta_g}(\m_-)$ and $U_{\eta_1}(\m_-)$ can be explicitly described in terms of quantum group analogues of root vectors. Note that one can also define  analogues $U_h^s(\m_-)$ of subalgebras $U_{\eta_g}(\m_-)$ in the standard Drinfeld--Jimbo quantum group $U_h(\mathfrak g)$ over the ring of formal power series $\mathbb{C}[[h]]$ (see \cite{S10}).

In \cite{S10}, Theorem 5.2 it is shown that for every conjugacy class $\mathcal{O}$ in $G$ one can find a transversal slice $\Sigma_s$ such that $\mathcal{O}$ intersects $\Sigma_s$ and ${\rm dim}~\mathcal{O}={\rm codim}~\Sigma_s$. Using this result we showed in \cite{S11} that for every element $g\in G$ one can find a subalgebra $U_{\eta_g}(\m_-)$ in $U_{\eta_g}$ of dimension $\bar{m}^{\frac{1}{2}{\rm dim}~\mathcal{O}_g}$ with a non--trivial character $\chi$. The dimension of the corresponding algebra $U_{\eta_1}(\m_-)$ is also equal to $\bar{m}^{\frac{1}{2}{\rm dim}~\mathcal{O}_g}$.

Following the strategy outlined in the beginning of this section we show that if $\bar{m}$ satisfies a certain condition then every finite-dimensional $U_{\eta_g}$--module is free over $U_{\eta_1}(\m_-)$. 
Thus the dimension of every such module is divisible by $\bar{m}^{\frac{1}{2}{\rm dim}~\mathcal{O}_g}$. This establishes the De Concini--Kac--Procesi conjecture.

Note that in the case of restricted representations of a small quantum group similar results were obtained in \cite{Dr}. The situation in \cite{Dr} is rather similar to the case of the trivial character $\eta=\eta_1$ in our setting. 

We also show that the rank of every finite-dimensional $U_{\eta_g}$--module $V$ over ${U}_{\eta_1}(\m_-)$ is equal to the dimension of the space $V_\chi$ and that $U_{\eta_g}$ is the algebra of matrices of size $\bar{m}^{\frac{1}{2}{\rm dim}~\mathcal{O}_g}$ over the corresponding q-W--algebra ${\rm Hk}(A,B,\chi)={\rm Hk}(U_{\eta_g},U_{\eta_g}(\m_-),\chi)$ which has dimension $\bar{m}^{{\rm dim}~\Sigma_s}$. In the  case of Lie algebras over fields of prime characteristic similar results were obtained in \cite{Pr}.

Note that the support variety technique used in \cite{Pr1} to prove the Kac--Weisfeiler conjecture can not be transferred to the case of quantum groups straightforwardly. The notion of the support variety is still available in case of quantum groups (see \cite{Dr,GK,O}). But in practical applications it is much less efficient since in the case of quantum groups there is no any underlying vector space.


\subsection*{The structure of the book}\label{structure}

In conclusion we would like to make a few remarks on the structure of the book. It consists of six chapters. In this introduction we have given a very superficial and incomplete review of the content of the book which rather aims to provide the reader with a general guide outlining the main ideas and the strategy of the main proofs. More technical comments are given in the beginning of each chapter. 

In Chapters \ref{part1} and \ref{part2} we summarize results from \cite{S6,S10,SDM,S12,S13} on the algebraic group analogues $\Sigma_s$ of the Slodowy slices and the related results on quantum groups and on the subalgebras $U_h^s(\m_-)\subset U_h(\mathfrak g)$. Chapter \ref{part1} also contains some results on combinatorics of Weyl groups and on root systems required for the definition of the slices $\Sigma_s$, and Chapter \ref{part2} contains some advanced results on quantum groups required later for the study of q-W--algebras.

In Chapter \ref{part3}, following \cite{S10,SDM}, we recall the definition of q-W--algebras and the description of their classical Poisson counterparts given in \cite{S13} in terms of the Zhelobenko type operators $\Pi_{jk}$ and $\Pi$. The main purpose of this chapter is to bring this description to a form suitable for quantization. Formulas (\ref{Apc}), (\ref{Pip}) and (\ref{pic}) obtained in this chapter for $\Pi_{jk}$ and $\Pi_c$ have direct quantum analogues (\ref{Bpdef}), (\ref{Ppk}) and (\ref{Piqdef}) obtained in Chapter \ref{part4} for $P_{jk}$ and $\Pi_c^q$. The main result of Chapter \ref{part4} (Theorem \ref{Piqmain}) is the description of the q-W--algebra as the image of the operator $\Pi_c^q$.

In Chapter \ref{part5} we prove a version of the Skryabin equivalence of type (\ref{catequiv}) for equivariant modules over quantum groups established in \cite{SDM}. The new proof of this equivalence in Theorem \ref{sqeq} is based on Corollary \ref{Bpbas} which allows to construct some nice bases in modules from the category $A-{\rm  mod}_B^\chi$ (see the discussion in the introduction above). Theorem \ref{sqeq} also gives precise values of $\varepsilon$ of the deformation parameter $q$ for which the categorical equivalence holds while in \cite{SDM} it was established for generic $\varepsilon$ only.

Finally in Chapter \ref{part6} we apply the results of Chapter \ref{part4} to the study of representations of quantum groups at roots of unity and prove the De Concini--Kac--Procesi conjecture. The strategy of this proof has already been discussed above.

Some historic remarks are given in the bibliographic comments after each chapter.

\vskip 0.3cm 
{\bf Acknowledgements}
\vskip 0.3cm

The author is greatly indebted to Giovanna Carnovale, Iulian Ion Simion, Lewis Topley and to the members of the representation theory seminars at the Universities of Bologna and Padua for careful reading of some parts of the manuscript. 

The results presented in this book have been partially obtained during several research stays at Institut des Haut \'{E}tudes Scientifiques, Paris, Max--Planck--Instut f\"{u}r Mathematik, Bonn and at the Representation Theory thematic program held at Institut Henri Poincar\'{e}, Paris, January 6--April 3, 2020. The author is grateful to these institutions for hospitality. 

The research on this project received funding from the European Research Council (ERC) under the European Union's Horizon 2020 research and innovation program (QUASIFT grant agreement 677368) during the visit of the author to Institut des Haut \'{E}tudes Scientifiques, Paris.


\newpage

\tableofcontents

\renewcommand{\thetheorem}{\thesection.\arabic{theorem}}

\renewcommand{\thelemma}{\thesection.\arabic{lemma}}

\renewcommand{\theproposition}{\thesection.\arabic{proposition}}

\renewcommand{\thecorollary}{\thesection.\arabic{corollary}}

\renewcommand{\theremark}{\thesection.\arabic{remark}}

\renewcommand{\thedefinition}{\thesection.\arabic{definition}}

\renewcommand{\theequation}{\thesection.\arabic{equation}}


\chapter{Algebraic group analogues of Slodowy slices}\label{part1}

The q-W--algebras are non-commutative deformations of algebras of regular functions on certain algebraic varieties in algebraic groups transversal to conjugacy classes. In this book these varieties play a role similar to that of the Slodowy slices in the theory of W--algebras and of generalized Gelfand-Geraev representations of semisimple Lie algebras. In this chapter we define these varieties and study their properties. We also develop the relevant Weyl group combinatorics.


\section{Notation}\label{notation}

\setcounter{equation}{0}
\setcounter{theorem}{0}

Fix the notation used throughout the book.
In this book we denote by $\mathbb{N}$ \index[not]{n@$\mathbb{N}$} the set of non-negative integer numbers, $\mathbb{N}=\{0,1,\ldots \}$.

Let $G_{\bf k}$ \index[not]{G@$G_{\bf k}$} be a connected \index{group!connected} finite--dimensional semisimple \index{group!semisimple} algebraic \index{group!algebraic} group over an algebraically closed \index{field!algebraically closed} field $\bf k$. \index[not]{k@$\bf k$} Denote by $\g_{\bf k}$ \index[not]{g@$\g_{\bf k}$} the Lie algebra \index{Lie!algebra} of $G_{\bf k}$. Let $H_{\bf k}\subset G_{\bf k}$ \index[not]{H@$H_{\bf k}$} be a maximal \index{torus, maximal} torus in $G_{\bf k}$, $\h_{\bf k}\subset \g_{\bf k}$ \index[not]{h@$\h_{\bf k}$} the corresponding Cartan subalgebra. \index{Lie!subalgebra! Cartan} If the characteristic exponent \index{field!characteristic exponent} of $\bf k$ is $p$, i.e. $p={\rm char}~{\bf k}$ if ${\rm char}~{\bf k}>0$ and $p=1$ if ${\rm char}~{\bf k}=0$, \index[not]{c@${\rm char}~{\bf k}$} we also write $G_{\bf k}=G_p$, \index[not]{G@$G_p$} $H_{\bf k}=H_p$, \index[not]{H@$H_p$} and if ${\bf k}=\mathbb{C}$ we write $G_{\mathbb{C}}=G$, \index[not]{G@$G$} $H_{\mathbb{C}}=H$, \index[not]{H@$H$} $\g_{\mathbb{C}}=\g$, \index[not]{g@$\g$} $\h_{\mathbb{C}}=\h$. \index[not]{h@$\h$} Note that $G$ is also a connected finite-dimensional complex \index{group!complex} semisimple Lie group \index{group!Lie}.

Let $\Delta=\Delta(\g, \h)$ \index[not]{D@$\Delta$} \index[not]{D@$\Delta(\g, \h)$} be the root system \index{root!system} of the pair $\left( {\frak g},{\frak h}\right)$, $Q$ \index[not]{Q@$Q$} the corresponding root lattice, \index{lattice!root} and $P$ \index[not]{P@$P$} the weight lattice. \index{lattice!weight}  Let $\Gamma=\{\alpha_i|~i=1,\ldots, l\},~~l={\rm rank}({\frak g})$ \index[not]{G@$\Gamma$} \index[not]{l@$l={\rm rank}({\frak g})$} \index[not]{a@$\alpha_i$} be a system of simple roots, \index{root!simple} $\Delta_+=\{ \beta_1, \ldots ,\beta_D \}$ \index[not]{D@$\Delta_+$} the corresponding set of positive roots, \index{roots!positive set of}  $Q_+=\mathbb{N}\Delta_+$, $P_+$ \index[not]{Q@$Q_+$} \index[not]{P@$P_+$} the set of the corresponding integral dominant weights, \index{weight!integral dominant} $\omega_1,\ldots ,\omega_l$ the fundamental weights. \index{weight!fundamental} Let also $H_1,\ldots ,H_l$ be the set of simple root generators of $\frak h$. \index{generators, simple root of a Cartan subalgebra}

Denote by  $a_{ij}$ \index[not]{a@$a_{ij}$} the corresponding Cartan matrix, \index{Cartan matrix}
and let $d_1,\ldots , d_l$, $d_i\in \{1,2,3\}$, \index[not]{d@$d_i$} $i=1,\ldots ,l$ be coprime positive integers such that the matrix $b_{ij}=d_ia_{ij}$ \index[not]{b@$b_{ij}$} is symmetric. There exists a unique canonical non--degenerate invariant
symmetric bilinear form \index{form!non--degenerate symmetric invariant bilinear!on a complex semisimple Lie algebra} $\left\langle ~\cdot~ ,~\cdot~ \right\rangle$ \index[not]{ZZ@$\left\langle ~\cdot~ ,~\cdot~ \right\rangle$} on ${\frak g}$ such that
$\left\langle H_i , H_j\right\rangle=d_j^{-1}a_{ij}$. It induces an isomorphism of vector spaces
${\frak h}\simeq {\frak h}^*$ \index[not]{h@$\h^*$} under which $\alpha_i \in {\frak h}^*$ corresponds
to $d_iH_i \in {\frak h}$. We denote by $h^\vee$ \index[not]{h@$h^\vee$} the element of $\frak h$ that
corresponds to $h \in {\frak h}^*$ under this isomorphism. For a root $\alpha \in \Delta$ the element $\alpha^\vee \in \h$ \index[not]{a@$\alpha^\vee$} is called the corresponding coroot.
The induced bilinear form on ${\frak h}^*$ is given by $\left\langle \alpha_i , \alpha_j\right\rangle=b_{ij}$. Given the restriction $d_i\in \{1,2,3\}$ the bilinear form is normalized in such a way that $\left\langle \alpha,\alpha\right\rangle=2$ for short roots $\alpha\in \Delta$.

Let $W$ \index[not]{W@$W$} be the Weyl group \index{group!Weyl} of the root system $\Delta$. $W$ can also be defined as the Weyl group of the pair $(\g,\h)$, $W=W(\g,\h)$. \index[not]{W@$W(\g,\h)$} $W$ is the subgroup of $GL({\frak h})$
generated by the fundamental reflections \index{reflection!fundamental} $s_1,\ldots ,s_l$, \index[not]{s@$s_i$}
$$
s_i(h)=h-\alpha_i(h)H_i,~~h\in{\frak h}.
$$
The action of $W$ on $\h$ \index{action!of a Weyl group on a Cartan subalgebra} preserves the restriction of the bilinear form $\left\langle ~\cdot~ ,~\cdot~ \right\rangle$ to $\frak h$.

For any root $\alpha\in \Delta$ we also denote by $s_\alpha$ \index[not]{s@$s_\alpha$} the corresponding reflection with respect to $\alpha$. \index{reflection!with respect to a root}

For every element $w\in W$ one can introduce the set $\Delta_w=\{\alpha \in \Delta_+: w(\alpha)\in -\Delta_+\}$, \index[not]{D@$\Delta_w$} and the number of the elements in the set $\Delta_w$ is equal to the length $l(w)$ \index[not]{l@$l(w)$} of the element $w$ with respect to the system $\Gamma$ of simple roots in $\Delta_+$. We also write $\Delta_-=-\Delta_+$. \index[not]{D@$\Delta_-$}

Let $\h_{\mathbb{R}}$ \index[not]{h@$\h_{\mathbb{R}}$} be the real form of $\h$, the real linear span of simple coroots in $\h$. The set of roots $\Delta$ is a subset of the dual space $\h_\mathbb{R}^*$. \index[not]{h@$\h_{\mathbb{R}}^*$}

For any $w\in W$ we denote by $\h^w$ \index[not]{h@$\h^w$} the fixed point space of $w$ in $\h$, $\h^w=\{x\in \h| wx=x\}$, and by $\h^w_{\mathbb{R}}$ the fixed point space of $w$ in $\h_{\mathbb{R}}$, $\h^w_{\mathbb{R}}=\{x\in \h_{\mathbb{R}}| wx=x\}=\h^w\cap \h_{\mathbb{R}}$. \index[not]{h@$\h^w_{\mathbb{R}}$}

One can define $\Delta$ as the root system of the pair $(G,H)$, $\Delta=\Delta(G,H)$, \index[not]{D@$\Delta(G, H)$} and if, for some algebraically closed field $\bf k$, $G$ and $G_{\bf k}$ have the same root system then one can also define $\Delta$ as the root system of the pair $(G_{\bf k},H_{\bf k})$, $\Delta=\Delta(G_{\bf k},H_{\bf k})$. \index[not]{D@$\Delta(G_{\bf k},H_{\bf k})$}

Similarly, one can define $W$ as the Weyl group of the pair $(G,H)$, $W=W(G,H)$, \index[not]{W@$W(G,H)$} and if, for some algebraically closed field $\bf k$, $G$ and $G_{\bf k}$ have the same root system then one can also define $W$ as the Weyl group of the pair $(G_{\bf k},H_{\bf k})$. We shall also use the symbol $W(R_{\bf k},H_{\bf k})$ for the Weyl group of the pair $(R_{\bf k},H_{\bf k})$ in the case when $R_{\bf k}$ \index[not]{R@$R_{\bf k}$} is a reductive algebraic group over $\bf k$, and $H_{\bf k}$ is a maximal torus in it. The Weyl group $W(R_{\bf k},H_{\bf k})$ can be defined by the following formula
\begin{equation}\label{Wgrdef}
W(R_{\bf k},H_{\bf k})=N_{R_{\bf k}}(H_{\bf k})/Z_{R_{\bf k}}(H_{\bf k}),
\end{equation}
where $N_{R_{\bf k}}(H_{\bf k})$ and $Z_{R_{\bf k}}(H_{\bf k})$ are the normalizer and the centralizer of $H_{\bf k}$ in $R_{\bf k}$, respectively, with respect to the conjugation action. \index{action!conjugation of an algebraic group} \index[not]{W@$W(R_{\bf k},H_{\bf k})$} \index[not]{N@$N_{R_{\bf k}}(H_{\bf k})$} \index[not]{Z@$Z_{R_{\bf k}}(H_{\bf k})$}

For an arbitrary algebraically closed field $\bf k$, we denote by $\dot{w}$ \index[not]{w@$\dot{w}$} a representative of $w\in W(R_{\bf k},H_{\bf k})=N_{R_{\bf k}}(H_{\bf k})/Z_{R_{\bf k}}(H_{\bf k})$ in $N_{R_{\bf k}}(H_{\bf k})$. 
If ${\bf k}=\mathbb{C}$ we simply denote this representative by $w$. 

For $w\in W, g\in G$ we write $w(g)=wgw^{-1}$. \index[not]{w@$w(g)$}
A representative of $w\in W=W(G,H)$ in $N_G(H)$ can also be regarded as an element of the normalizer $N_G(\h)$ \index[not]{N@$N_G(\h)$} of $\h$ in $G$ with respect to the adjoint action.

Let $B_+$ be the Borel subgroup \index{subgroup!Borel} of $G$ corresponding to $\Delta_+$ and $B_-$ \index[not]{B@$B_\pm$}
the opposite Borel subgroup of $G$, $N_\pm$ \index[not]{N@$N_\pm$} their unipotent radicals, \index{radical!unipotent} respectively. 

We denote by ${{\frak b}_\pm}$ and ${\frak n}_\pm$ \index[not]{b@$\b_\pm$} \index[not]{n@$\n_\pm$} the Lie subalgebras \index{Lie!subalgebra} of $\g$ corresponding to $B_\pm$ and $N_\pm$, respectively.

We identify $\frak g$ and its dual by means of the canonical bilinear form $\left\langle ~\cdot~ ,~\cdot~ \right\rangle$.
Then the coadjoint action of $G$ on ${\frak g}^*$ \index[not]{g@${\frak g}^*$} \index{action!Lie group!adjoint} \index{action!Lie group!coadjoint} is naturally identified with the adjoint one. Using the canonical bilinear form we shall also identify ${{\frak n}_+}^*\simeq {\frak n}_-,~{{\frak b}_+}^*\simeq {\frak b}_-$, $\h\simeq \h^*$. \index[not]{n@${\frak n}_\pm^*$} \index[not]{b@${\frak b}_\pm^*$}

Let ${\frak g}_\beta$ \index[not]{g@${\frak g}_\beta$} be the root subspace \index{root!subspace} corresponding to a root $\beta \in \Delta$,
${\frak g}_\beta=\{ x\in {\frak g}| [h,x]=\beta(h)x \mbox{ for every }h\in {\frak h}\}$.
${\frak g}_\beta\subset {\frak g}$ is a one--dimensional subspace.
The following result is well known. 
\begin{lemma}\label{orthogr}
For $\alpha\neq -\beta$ the root subspaces ${\frak g}_\alpha$ and ${\frak g}_\beta$ are orthogonal with respect
to the canonical invariant bilinear form $\left\langle ~\cdot~ ,~\cdot~ \right\rangle$. Moreover ${\frak g}_\alpha$ and ${\frak g}_{-\alpha}$
are non--degenerately paired by this form.
\end{lemma}

Let $X_\alpha\in \g$ \index[not]{X@$X_\alpha$} be a non--zero root vector \index{root!vector} corresponding to a root $\alpha\in \Delta$.
Root vectors $X_{\alpha}\in {\frak g}_\alpha$ satisfy the following relations:
$$
[X_\alpha,X_{-\alpha}]=\left\langle X_\alpha,X_{-\alpha}\right\rangle\alpha^\vee.
$$


\section{Systems of positive roots associated to Weyl group elements}\label{background}

\setcounter{equation}{0}
\setcounter{theorem}{0}

Algebraic group analogues of the Slodowy slices are associated to (conjugacy classes) in the Weyl group. \index{conjugacy class!in a Weyl group} In this section we introduce the relevant combinatorics of the Weyl group and of root systems. We start by defining systems of positive roots associated to Weyl group elements which play the key role in the definition of the algebraic group analogues of the Slodowy slices. 

Let $s$ be an element of the Weyl group $W$ and denote by $\h'$ \index[not]{h@$\h'$} the orthogonal complement in $\h$, with respect to the canonical bilinear form on $\g$, to the subspace $\h^s=\{h\in \h| sh=h\}$ fixed by the natural action of $s$ on $\h$, $\h'=(\h^s)^\perp$, so that $\h'^\perp=\h^s$. Let $\h'^*$ \index[not]{h@$\h'^*$} be the image of $\h'$ in $\h^*$ under the identification $\h^*\simeq \h$ induced by the canonical bilinear form on $\g$. Thus $\h'^*$ embeds into $\h$ thanks to the direct vector space decomposition $\h = \h^s + \h'$.
By Theorem C in \cite{C} $s$ can be represented as a product of two involutions,
\begin{equation}\label{inv}
s=s^1s^2,
\end{equation}
where $s^1=s_{\gamma_1}\ldots s_{\gamma_{\widetilde{l}}}$, \index[not]{l@$\widetilde{l}$} $s^2=s_{\gamma_{\widetilde{l}+1}}\ldots s_{\gamma_{l'}}$, \index[not]{s@$s^1$} \index[not]{s@$s^2$} \index[not]{l@$l'$} the roots in each of the sets $\gamma_1, \ldots, \gamma_{\widetilde{l}}$ and $\gamma_{\widetilde{l}+1}, \ldots, \gamma_{l'}$ belong to some set of positive roots, and mutually orthogonal, and
the roots $\gamma_1, \ldots, \gamma_{l'}$ \index[not]{g@$\gamma_i$} form a linear basis of $\h'^*$.

Weyl group elements naturally act on $\h_{\mathbb{R}}$ as orthogonal transformations with respect to the scalar product induced by the symmetric bilinear form on $\g$, and one can define the real forms $\h_{\mathbb{R}}^s=\h^s\cap \h_{\mathbb{R}}, \h_{\mathbb{R}}'=\h'\cap \h_{\mathbb{R}}$. \index[not]{h@$\h_{\mathbb{R}}'$}

Let $\widehat{\gamma}_1, \ldots \widehat{\gamma}_{l'}$ be the vectors of unit length in the directions of $\gamma_1, \ldots \gamma_{l'}$, and $\widehat{\gamma}^*_1,\ldots,\widehat{\gamma}^*_{l'}$ the basis of $\h_{\mathbb{R}}'$ dual to $\widehat{\gamma}_1, \ldots \widehat{\gamma}_{l'}$. Let $O$ be the $l'\times l'$ symmetric matrix with real entries $O_{ij}=\left\langle \widehat{\gamma}_i,\widehat{\gamma}_j\right\rangle$. $I-O$ is also a symmetric real matrix, and hence it is diagonalizable and has real eigenvalues.

The following proposition gives a recipe for constructing a spectral decomposition for the action of the orthogonal transformation $s$ on $\h_{\mathbb{R}}$.

\begin{proposition}\label{IM}
Let $\lambda$ be a (real) eigenvalue of the symmetric matrix $I-O$, and $u\in \mathbb{R}^{l'}$ a corresponding non--zero real eigenvector with components $u_i$, $i=1,\ldots ,l'$. Let $a_u,b_u\in \h_{\mathbb{R}}$ be defined by
\begin{equation}\label{ab}
a_u = \sum_{i=1}^{\widetilde{l}} u_i\widehat{\gamma}^*_i,~b_u = \sum_{i=\widetilde{l}+1}^{l'} u_i\widehat{\gamma}^*_i.
\end{equation}

(i) If $\lambda$ is an eigenvalue of $I-O$ then $-\lambda$ is also an eigenvalue of $I-O$, and $\lambda=\pm 1$ are not eigenvalues of $I-O$.

(ii) If $\lambda \neq 0$ then the angle $\theta$ between $a_u$ and $b_u$ satisfies $\cos \theta=\lambda$, the plane $\h_\lambda \subset \h_{\mathbb{R}}$ \index[not]{h@$h_\lambda$} spanned by $a_u$ and $b_u$ is invariant with respect to the involutions $s^{i}$, $i=1,2$, $s^1$ acts on $\h_\lambda$ as the reflection in the line spanned by $b_u$, and $s^2$ acts on $\h_\lambda$ as the reflection in the line spanned by $a_u$. If $\lambda>0$ the orthogonal transformation $s=s^1s^2$ acts on $\h_\lambda$ as a rotation through the angle $2\theta$.

(iii) If $\lambda\neq \mu$ are two positive eigenvalues of $I-O$ then the planes $\h_\lambda$ and $\h_\mu$ are mutually orthogonal.

(iv) Let $\lambda\neq 0$ be an eigenvalue of $I-O$ of multiplicity greater than $1$, and $u^k \in \mathbb{R}^{l'}$, $k=1,\ldots, ~{\rm mult}~\lambda$ a basis of the eigenspace corresponding to $\lambda$. If the basis $u^k$ is orthonormal with respect to the standard scalar product on $\mathbb{R}^{l'}$ then the corresponding planes $\h_\lambda^k$ \index[not]{h@$\h_\lambda^k$} defined with the help of $u^k$, $k=1,\ldots, {\rm mult}~\lambda$ are mutually orthogonal.

(v) If $\lambda= 0$ is an eigenvalue of $I-O$, then there is a basis $u^k \in \mathbb{R}^{l'}$, $k=1,\ldots, ~{\rm mult}~0$ of the eigenspace corresponding to $0$ orthonormal with respect to the standard scalar product on $\mathbb{R}^{l'}$ and such that the corresponding non--zero elements $a_{u^k}$, $b_{u^k}$ are all mutually orthogonal. Moreover, $s^1a_{u^k}=-a_{u^k}$, $s^2a_{u^k}=a_{u^k}$, $s^1b_{u^k}=b_{u^k}$, $s^2b_{u^k}=-b_{u^k}$ for non--zero elements $a_{u^k}$, $b_{u^k}$. In particular, for  non--zero elements $a_{u^k}$, $b_{u^k}$ we have $sa_{u^k}=-a_{u^k}$, $sb_{u^k}=-b_{u^k}$, and non--zero elements $a_{u^k}$, $b_{u^k}$ is a basis of the subspace of $\h_{\mathbb{R}}$ on which $s$ acts by multiplication by $-1$.
\end{proposition}

\begin{proof}
Firstly, we study some general properties of the eigenvalues and eigenvectors of the matrix $O$. 
By definition the matrix $O$ can be written in a block form,
\begin{equation}\label{m}
O=\left(
  \begin{array}{cc}
    I_{\widetilde{l}} & A \\
    A^\top & I_{l'-\widetilde{l}} \\
  \end{array}
\right),
\end{equation}
where $A$ is an $n\times (l'-\widetilde{l})$ matrix, $A^\top$ is the transpose to $A$, $I_{\widetilde{l}}$ and $I_{l'-\widetilde{l}}$ are the unit matrices of sizes $\widetilde{l}$ and $l'-\widetilde{l}$. $O^{-1}$ is also symmetric and has a similar block form,
\begin{equation}\label{m-1}
O^{-1}=\left(
  \begin{array}{cc}
    B & C \\
    C^\top & D \\
  \end{array}
\right),~B=B^\top,~D=D^\top,
\end{equation}
with the entries $O^{-1}_{ij}=\left\langle \widehat{\gamma}^*_i,\widehat{\gamma}^*_j\right\rangle$.

For any vector $u\in \mathbb{R}^{l'}$ we introduce its $\mathbb{R}^{\widetilde{l}}$ and $\mathbb{R}^{l'-\widetilde{l}}$ components $\widetilde{u}$ and $\widetilde{\widetilde{u}}$ in a similar way,
\begin{equation}\label{udec}
u=\left(
    \begin{array}{c}
      \widetilde{u} \\
      \widetilde{\widetilde{u}} \\
    \end{array}
  \right).
\end{equation}
We shall consider both $\widetilde{u}$ and $\widetilde{\widetilde{u}}$ as elements of $\mathbb{R}^{l'}$ using natural embeddings $\mathbb{R}^{\widetilde{l}}, \mathbb{R}^{l'-\widetilde{l}} \subset \mathbb{R}^{l'}$ associated to decomposition (\ref{udec}).

If $u$ is a non--zero eigenvector of $I-O$ corresponding to an eigenvalue $\lambda$ then the equation $(I-O)u=\lambda u$ gives
\begin{equation}\label{ueq}
-A\widetilde{\widetilde{u}}=\lambda\widetilde{u},~-A^\top\widetilde{u}=\lambda\widetilde{\widetilde{u}}.
\end{equation}
From these equations we deduce that 
$$
\left(
    \begin{array}{c}
      -\widetilde{u} \\
      \widetilde{\widetilde{u}} \\
    \end{array}
  \right)
$$
is a non--zero eigenvector of $I-O$ corresponding to the eigenvalue $-\lambda$ This proves the first claim in (i). 

$\lambda=1$ is not an eigenvalue of $I-O$ since the matrix $O$ is invertible. Therefore $\lambda=-1$ is also not an eigenvalue of $I-O$ by the first part of part (i) which is already proved. This justifies (i). 

Since $O^{-1}O=I$ one has
\begin{equation}\label{mmi}
BA+C=0,~C^\top+DA^\top=0.
\end{equation}
Multiplying the first and the second equations in (\ref{ueq}) from the left by $B$ and $D$, respectively, and using (\ref{mmi}) we obtain that
\begin{equation}\label{ident}
C\widetilde{\widetilde{u}}=\lambda B\widetilde{u},~C^\top\widetilde{u}=\lambda D\widetilde{\widetilde{u}}.
\end{equation}

Now if $u^{1}$ and $u^{2}$ are two non--zero eigenvectors of $I-O$ corresponding to an eigenvalue $\lambda$ then by (\ref{m-1}) we have
\begin{equation}\label{1}
(a_{u^1},a_{u^2})=\sum_{i,j=1}^{\widetilde{l}}u_i^1u_j^2\left\langle \widehat{\gamma}^*_i,\widehat{\gamma}^*_j\right\rangle=\sum_{i,j=1}^{\widetilde{l}}u_i^1u_j^2B_{ij}=
\widetilde{u}^1\cdot B\widetilde{u}^2,
\end{equation}
where $\cdot$ stands for the standard scalar product in $\mathbb{R}^{l'}$.

Similarly,
\begin{equation}\label{others}
\left\langle b_{u^1},b_{u^2}\right\rangle=D\widetilde{\widetilde{u}}^1\cdot \widetilde{\widetilde{u}}^2, \left\langle a_{u^1},b_{u^2}\right\rangle=\widetilde{u}^1\cdot C\widetilde{\widetilde{u}}^2.
\end{equation}

From (\ref{ident}), (\ref{1}) and the first identity in (\ref{others}) we also obtain that if $\lambda \neq 0$ then
\begin{equation}\label{2}
\left\langle a_{u^1},a_{u^2}\right\rangle=\widetilde{u}^1\cdot B\widetilde{u}^2=\frac{1}{\lambda}\widetilde{u}^1\cdot C\widetilde{\widetilde{u}}^2= \frac{1}{\lambda}C^\top\widetilde{u}^1\cdot \widetilde{\widetilde{u}}^2=D\widetilde{\widetilde{u}}^1\cdot \widetilde{\widetilde{u}}^2=\left\langle b_{u^1},b_{u^2}\right\rangle.
\end{equation}

Similarly, for any real eigenvalue $\lambda$ we have
\begin{equation}\label{3}
\left\langle a_{u^1},b_{u^2}\right\rangle=\widetilde{u}^1\cdot C\widetilde{\widetilde{u}}^2=\lambda\left\langle a_{u^1},a_{u^2}\right\rangle,~\left\langle b_{u^1},a_{u^2}\right\rangle=\widetilde{\widetilde{u}}^1\cdot C^\top{\widetilde{u}}^2=\lambda\left\langle a_{u^1},a_{u^2}\right\rangle.
\end{equation}

Therefore if $\lambda\neq 0$, taking into account (\ref{2}), we obtain for $u^1=u^2=u$
$$
\lambda=\frac{\left\langle a_{u},b_{u}\right\rangle}{\left\langle a_{u},a_{u}\right\rangle}=\frac{\left\langle a_{u},b_{u}\right\rangle}{\sqrt{\left\langle b_{u},b_{u}\right\rangle}\sqrt{\left\langle a_{u},a_{u}\right\rangle}}=\cos \theta,
$$
which justifies the first claim in (ii). Note that $\left\langle b_{u},b_{u}\right\rangle=\left\langle a_{u},a_{u}\right\rangle\neq 0$ for otherwise $a_u=b_u=0$, and hence $u=0$ which contradicts the choice of $u$.

Let again  $u$ be a non--zero eigenvector of $I-O$ corresponding to an eigenvalue $\lambda$.
For $i=1,\ldots, \widetilde{l}$ by the definition of the matrices $B$ and $C$ we have
$$
\left\langle \widehat{\gamma}^*_i,\lambda a_u-b_u\right\rangle=\lambda (B\widetilde{u})_i-(C\widetilde{\widetilde{u}})_i=0,
$$
where at the last step we used the first identity in (\ref{ident}). From the last identity we deduce that $\lambda a_u-b_u$ is a linear combination of $\widehat{\gamma}^*_{\widetilde{l}+1},\ldots, \widehat{\gamma}^*_{l'}$, and hence
$$
s^2(\lambda a_u-b_u)=-(\lambda a_u-b_u).
$$

However, by the definition of $a_u$, $s^2a_u=a_u$. Therefore
\begin{equation}\label{aubu}
s^2b_u=2\lambda a_u-b_u.
\end{equation}
Let $\lambda \neq 0$. Then recalling that by (\ref{2}) $\left\langle a_{u},a_{u}\right\rangle=\left\langle b_{u},b_{u}\right\rangle$ we conclude that $\lambda a_u=\cos (\theta) a_u$ is the orthogonal projection of $b_u$ onto the line spanned by $a_u$ and that $s^2b_u$ is obtained from $b_u$ by the reflection in the line spanned by $a_u$ as shown at Figure 1. 
\vskip -0.3cm
\begin{center}
$$
\xy/r10pc/: ="A","A", {\ar+(0.53,0.45)},"A",+(0.6,0.52)*{b_u}, "A", *\xycircle(0.7,0.7){}, "A", {\ar+(0.7,0)},"A",+(0.78,0.05)*{a_u},"A", {\ar+(0.53,-0.45)},"A",+(0.6,-0.52)*{s^2b_u}, "A", +(0.53,0.45)="E", "A", +(0.53,-0.45)="F", "E";"F"**@{--}, "A",+(0.42,0.06)*{\lambda a_u}, "A", {\ar+(0.53,0)} 
\endxy
$$
\vskip 0.3cm
Fig. 1
\end{center}

Similarly, $s^1b_u=b_u$, $s^1a_u$ is obtained from $a_u$ by the reflection in the line spanned by $b_u$.

Thus the plane $\h_\lambda \subset \h_{\mathbb{R}}$ spanned by $a_u$ and $b_u$ is invariant with respect to the involutions $s^{i}, ~i=1,2$, $s^1$ acts on $\h_\lambda$ as the reflection in the line spanned by $b_u$, and $s^2$ acts on $\h_\lambda$ as the reflection in the line spanned by $a_u$. Since the angle between $a_u$ and $b_u$ is $\theta$, for $\lambda>0$ the orthogonal transformation $s=s^1s^2$ acts on $\h_\lambda$ as a rotation through the angle $2\theta$ which completes the proof of (ii).

From the general theory of orthogonal transformations it follows that if $\lambda\neq \mu$ are two positive eigenvalues of $I-O$, $\lambda,\mu \neq 1$ then the planes $\h_\lambda$ and $\h_\mu$ are mutually orthogonal which confirms part (iii).

Now if $u^{1}$ and $u^{2}$ are two non--zero eigenvectors of $I-O$ corresponding to an eigenvalue $\lambda \neq 0$ then by part (i) $\lambda \neq \pm 1$, and (\ref{1}), (\ref{2}), (\ref{3}) and the identity $O^{-1}u^2=\frac{1}{1-\lambda}u^2$ yield
$$
\left\langle a_{u^1}+b_{u^1},a_{u^2}+b_{u^2}\right\rangle=2\left\langle a_{u^1},a_{u^2}\right\rangle(\lambda + 1)=u^1\cdot O^{-1}u^2=\frac{1}{1-\lambda}u^1 \cdot u^2.
$$
Thus if $\lambda \neq 0$ and $u^{1,2}$ are mutually orthogonal, $a_{u^1}$, $a_{u^2}$ are also mutually orthogonal, and from (\ref{2}) and (\ref{3}) we obtain that $b_{u^1}$ and $b_{u^2}$, $a_{u^1}$ and  $b_{u^2}$, $a_{u^2}$ and $b_{u^1}$  are mutually orthogonal. Therefore the planes spanned by $a_{u^1},b_{u^1}$ and by $a_{u^2},b_{u^2}$ are mutually orthogonal. Part (iv) immediately follows from this property.

It remains to prove part (v). If $\lambda= 0$ is an eigenvalue of $I-O$ then $\widetilde{u}$ and $\widetilde{\widetilde{u}}$ are the components of an eigenvector $u$ of $I-O$ with eigenvalue $0$ if and only if $A\widetilde{\widetilde{u}}=0$ and $A^\top\widetilde{u}=0$. Therefore using the usual orthogonalization procedure one can construct a basis $u^k \in \mathbb{R}^{l'}$, $k=1,\ldots, ~{\rm mult}~0$ of the eigenspace corresponding to $0$ orthonormal with respect to the standard scalar product on $\mathbb{R}^{l'}$ and such that the components $\widetilde{u}^k$ and $\widetilde{\widetilde{u}}^k$ $k=1,\ldots, ~{\rm mult}~0$ are all mutually orthogonal.

By (\ref{aubu}) $s^2b_{u^k}=-b_{u^k}$. Also by the definition of $a_{u^k}$ $s^2a_{u^k}=a_{u^k}$.  Similarly, $s^1a_{u^k}=-a_{u^k}$ and $s^1b_{u^k}=b_{u^k}$. 

Now using the definition of eigenvectors in the form $Ou^k=u^k$ and (\ref{ident}) with $u=u^k$
we deduce that for the basis $u^k$ the following relations hold: $B\widetilde{u}^k=\widetilde{u}^k$, $D\widetilde{\widetilde{u}}^k=\widetilde{\widetilde{u}}^k$. 

From these relations we obtain for $k\neq k'$ by (\ref{1}) 
$$
\left\langle a_{u^k},a_{u^{k'}}\right\rangle=\widetilde{u}^k\cdot B\widetilde{u}^{k'}=\widetilde{u}^k\cdot \widetilde{u}^{k'}=0
$$
and by (\ref{others}) 
$$
\left\langle b_{u^k},b_{u^{k'}}\right\rangle=D\widetilde{\widetilde{u}}^k\cdot \widetilde{\widetilde{u}}^{k'}=\widetilde{\widetilde{u}}^k\cdot \widetilde{\widetilde{u}}^{k'}=0.
$$

By (\ref{3}) we always have
$$
\left\langle a_{u^k},b_{u^{k'}}\right\rangle=\lambda\left\langle a_{u^k},a_{u^{k'}}\right\rangle=0.
$$
This completes the proof of part (v).

\end{proof}


Using the previous proposition we can decompose $\h_{\mathbb{R}}$ into a direct orthogonal sum of $s$--invariant subspaces,
\begin{equation}\label{hdec}
\h_\mathbb{R}=\bigoplus_{i=0}^{M(s)} \h_i, \index[not]{h@$\h_i$} \index[not]{M@$M(s)$}
\end{equation}
where  each of the subspaces $\h_i\subset \h_\mathbb{R}$, $i=1,\ldots, M(s)$ is invariant with respect to both involutions $s^{i}, ~i=1,2$ in the decomposition $s=s^1s^2$, and there are the following three possibilities for each $\h_i$: $\h_i$ is two--dimensional ($\h_i=\h_\lambda^k$ for an eigenvalue $0<\lambda<1$ of the matrix $I-O$, and $k=1,\ldots, {\rm mult}~\lambda$) and the Weyl group element $s$ acts on it as rotation with angle $\theta_i$, $0<\theta_i<\pi$ \index[not]{t@$\theta_i$} or $\h_i=\h_{\lambda}^k$, $\lambda=0$, $k=1,\ldots, {\rm mult}~\lambda$ has dimension $1$ and $s$ acts on it by multiplication by $-1$ or $\h_i$ coincides with the vector subspace of $\h_{\mathbb{R}}$ fixed by the action of $s$. Note that since $s$ has finite order, we have $\theta_i=\frac{2\pi n_i}{m_i}$, $n_i,m_i\in \{1,2,\ldots \}$.

Since the number of roots in the root system $\Delta$ is finite one can always choose elements $h_i\in \h_i$, $i=0,\ldots, M(s)$, \index[not]{h@$h_i$} such that $h_i(\alpha)\neq 0$ for any root $\alpha \in \Delta$ which is not orthogonal to the $s$--invariant subspace $\h_i$ with respect to the natural pairing between $\h_{\mathbb{R}}$ and $\h_{\mathbb{R}}^*$.

Now we consider certain $s$--invariant subsets of roots ${\Delta}_i$, $i=0,\ldots, M(s)$, \index[not]{D@$\Delta_i$} defined as follows
\begin{equation}\label{di}
\Delta_i=\{ \alpha\in \Delta: h_j(\alpha)=0, j>i,~h_i(\alpha)\neq 0 \},
\end{equation}
where we formally assume that $h_{M(s)+1}=0$.
Note that for some indexes $i$ the subsets $\Delta_i$ are empty, and that the definition of these subsets depends on the order of the terms in direct sum (\ref{hdec}).

 Now consider the nonempty $s$--invariant subsets of roots $\Delta_{i_k}$, $k=0,\ldots, M$. \index[not]{D@$\Delta_{i_k}$} \index[not]{M@$M$} For convenience we assume that indexes $i_k$ are labeled in such a way that $i_j<i_k$ if and only if $j<k$.

Observe also that the root system $\Delta$ is the disjoint union of the subsets ${\Delta}_{i_k}$,
\begin{equation}\label{decD}
\Delta=\bigcup_{k=0}^{M}{\Delta}_{i_k}.
\end{equation}

Now assume that
\begin{equation}\label{cond}
|h_{i_k}(\alpha)|>|\sum_{k'\leq j<k}h_{i_j}(\alpha)|, ~{\rm for~any}~\alpha\in {\Delta}_{i_k},~k,k'=0,\ldots, M,~k'<k.
\end{equation}
Condition (\ref{cond}) can be always fulfilled by suitable rescalings of the elements $h_{i_k}$.

Consider the element
\begin{equation}\label{hwb}
\bar{h}=\sum_{k=0}^{M}h_{i_k}\in \h_\mathbb{R}. \index[not]{h@$\bar{h}$}
\end{equation}

From definition (\ref{di}) of the sets ${\Delta}_i$ we obtain that for $\alpha \in {\Delta}_{i_k}$
\begin{equation}\label{dech}
\bar{h}(\alpha)=\sum_{j\leq k}h_{i_j}(\alpha)=h_{i_k}(\alpha)+\sum_{j< k}h_{i_j}(\alpha)
\end{equation}
Now condition (\ref{cond}), the previous identity and the inequality $|x+y|\geq ||x|-|y||$ imply that for $\alpha \in {\Delta}_{i_k}$ we have
$$
|\bar{h}(\alpha)|\geq ||h_{i_k}(\alpha)|-|\sum_{j< k}h_{i_j}(\alpha)||>0.
$$
Since $\Delta$ is the disjoint union of the subsets ${\Delta}_{i_k}$, $\Delta=\bigcup_{k=0}^{M}{\Delta}_{i_k}$, the last inequality ensures that  $\bar{h}$ belongs to a Weyl chamber of the root system $\Delta$. 

Denote by $\Delta_+^s$ \index[not]{D@$\Delta_\pm^s$} the subset of positive roots with respect to the Weyl chamber \index{Weyl chamber} containing $\bar{h}$,
\begin{equation}\label{D+sdef}
\Delta_+^s=\{\alpha\in \Delta| \alpha(\bar{h})>0\}.
\end{equation} 
Let $\Delta_-^s=-\Delta_+^s$.

From condition (\ref{cond}) and formula (\ref{dech}) we also obtain that a root $\alpha \in {\Delta}_{i_k}$ is positive if and only if
\begin{equation}\label{wc}
h_{i_k}(\alpha)>0.
\end{equation}
We denote by $({\Delta}_{i_k})_+$ the set of positive roots contained in ${\Delta}_{i_k}$, $({\Delta}_{i_k})_+=\Delta_+^s\cap {\Delta}_{i_k}$. \index[not]{D@$({\Delta}_{i_k})_+$}

We also define other $s$--invariant subsets of roots $\overline{\Delta}_{i_k}$, $k=0,\ldots, M$,
\begin{equation}\label{dik}
\overline{\Delta}_{i_k}=\bigcup_{i_j\leq i_k}{\Delta}_{i_j}. \index[not]{D@$\overline{\Delta}_{i_k}$}
\end{equation}
According to this definition we have a chain of strict inclusions
\begin{equation}\label{inc}
\overline{\Delta}_{i_M}\supset\overline{\Delta}_{i_{M-1}}\supset\ldots\supset\overline{\Delta}_{i_0},
\end{equation}
such that $\overline{\Delta}_{i_M}=\Delta$, $\overline{\Delta}_{i_0}=\Delta_{i_0}$, and $\overline{\Delta}_{i_k}\setminus \overline{\Delta}_{i_{k-1}}={\Delta}_{i_k}$.

The following lemma shows that the subsets of roots $\overline{\Delta}_{i_k}\subset \Delta$ are root systems of some standard Levi subalgebras in $\g$.
\begin{lemma}\label{parab}
Let $\Gamma^s$ \index[not]{G@$\Gamma^s$} be the set of simple roots in $\Delta_+^s$. Then
$\Gamma^s\cap \overline{\Delta}_{i_k}$ is a set of simple roots in $\overline{\Delta}_{i_k}$.
\end{lemma}

\begin{proof}
Indeed, let $\alpha\in \overline{\Delta}_{i_k}\cap \Delta_+^s$, $\alpha=\sum_{i=1}^l n_i\alpha_i$, where $n_i\in \{0,1,2,\ldots \}$ and $\Gamma^s=\{\alpha_1, \ldots ,\alpha_l\}$. Assume that $\alpha$ does not belong to the linear span of roots from $\Gamma^s\cap \overline{\Delta}_{i_k}$ and $t>i_k$ is maximal possible such that for some $\alpha_q\in {\Delta}_{t}$ one has $n_q>0$. Then by (\ref{di}) and (\ref{wc}) $h_t(\alpha)=\sum_{i=1}^l n_ih_t(\alpha_i)=\sum_{\alpha_i\in {\Delta}_{t}}  n_ih_t(\alpha_i)>0$, and by the choice of $t$ $h_{r}(\alpha)=0$ for $r>t$. Therefore $\alpha\in {\Delta}_{t}$, and hence $\alpha \not \in \overline{\Delta}_{i_k}$. Thus we arrive at a contradiction.

 \end{proof}


\pagestyle{myheadings}
\markboth{CHAPTER \thechapter.~ALGEBRAIC GROUP ANALOGUES OF SLODOWY SLICES}{\thesection.~NORMAL ORDERINGS OF POSITIVE ROOT SYSTEMS}

\section{Some normal orderings of positive root systems associated to Weyl group elements}\label{wgrord}

\setcounter{equation}{0}
\setcounter{theorem}{0}

For the purpose of quantization we shall need a certain normal ordering on the root system $\Delta_+^s$. Some properties of this ordering will also be used in the proof of the properties of transversal slices to conjugacy classes in algebraic groups \index{conjugacy class!in an algebraic group} in the next section.

An ordering of a set of positive roots $\Delta_+$ is called {\it normal} \index{ordering!normal of a set of positive roots} if it satisfies the following property.
\vskip 0.3cm
{\bf Property N.} ~~{\it
For any three roots $\alpha,~\beta,~\gamma$ such that $\gamma=\alpha+\beta$ we have either $\alpha<\gamma<\beta$ or $\beta<\gamma<\alpha$. 
}
\vskip 0.3cm

If 
$$
\beta_1, \ldots , \beta_D 
$$
is a normal ordering of $\Delta_+$, any part of it of the form
$$
\beta_h, \beta_{h+1}, \ldots , \beta_{h+b}
$$
is called a {\it segment} \index{segment!in a normally ordered set of positive roots} and denoted by $[\beta_h,\beta_{h+b}]$.

If $A, B\subset \Delta_+$ are two subsets, we write $A< B$ if for any $\theta\in A$, and $\vartheta\in B$ one has $\theta< \vartheta$ with respect to a normal ordering $<$ on $\Delta_+$. 

\begin{proposition}{\bf (\cite{Z1}, \S 3,  Proposition 2)}
Let $\alpha_1,\ldots, \alpha_l$ be the simple roots in $\Delta_+$, $s_1, \ldots ,s_l$ the corresponding simple reflections. Let $\overline{w}$ \index[not]{w@$\overline{w}$} be the element of $W$ of maximal length with respect to the system $s_1, \ldots ,s_l$ of simple reflections. For any reduced decomposition $\overline{w}=s_{i_1}\ldots s_{i_D}$ \index{Weyl group element!reduced decomposition} of $\overline{w}$ the ordering
\begin{equation}\label{normorddef}
\beta_1=\alpha_{i_1},\beta_2=s_{i_1}\alpha_{i_2},\ldots,\beta_D=s_{i_1}\ldots s_{i_{D-1}}\alpha_{i_D}
\end{equation}
is a normal ordering in $\Delta_+$, and there is a one--to--one correspondence between normal orderings of $\Delta_+$ and reduced decompositions of $\overline{w}$.
\end{proposition}

From this proposition and from properties of Coxeter groups it follows that any two normal orderings in $\Delta_+$ can be reduced to each other by the so--called {\it elementary transpositions}. \index{root!system!elementary transposition} The elementary transpositions for rank 2 root systems are inversions of the following normal orderings (or the inverse normal orderings):
\begin{equation}\label{rank2}
\begin{array}{lr}
\alpha,~\beta & A_1+A_1 \\
\\
\alpha,~\alpha+\beta,~\beta & A_2 \\
\\
\alpha,~\alpha+\beta,~\alpha+2\beta,~\beta & B_2 \\
\\
\alpha,~\alpha+\beta,~2\alpha+3\beta,~\alpha+2\beta,~\alpha+3\beta,~\beta & G_2
\end{array}
\end{equation}
where it is assumed that $(\alpha,\alpha)\geq (\beta,\beta)$. Moreover, for any rank 2 root system there exist two normal orderings one of which is contained in the list (\ref{rank2}) and the other is the inverse ordering.

In general, an elementary transposition of a normal ordering of a set of positive roots $\Delta_+$  is the inversion of an ordered segment of form (\ref{rank2}) (or of a segment with the inverse ordering) in the ordered set $\Delta_+$, where $\alpha-\beta\not\in \Delta$.

We recall some results of \S 3 and \S 4 in \cite{Z1} on normal orderings and reduced decompositions. 
For $w\in W$ a decomposition $w=uv$, $u,v \in W$, is called {\it reduced} \index{Weyl group element!reduced decomposition} if $l(w)=l(u)+l(v)$, where $l(~\cdot~)$ is the length function with respect to the system of simple roots in $\Delta_+$. In this case $\Delta_w=\Delta_v\cup v^{-1}\Delta_u$. This definition has an obvious generalization for products of several Weyl group elements.

Recall that for $w\in W$, a decomposition $w=s_{i_1}\ldots s_{i_m}$ is called {\it reduced} \index{Weyl group element!reduced decomposition} if it has minimal possible length equal to $l(w)$. Obviously, if $w=uv$, $u,v,w \in W$ is a reduced decomposition then the product of reduced decompositions of $u$ and $v$ is a reduced decomposition of $w$.


An ordering of the set $\Delta_w=\{\alpha\in \Delta_+ | w\alpha\in -\Delta_+\}$ is called {\it normal} \index{ordering!normal of a set $\Delta_w$} if it coincides with an initial segment of some normal ordering of $\Delta_+$. 
\begin{lemma}{\bf (\cite{Z1}, \S 4, Theorem 1)}\label{wN}
Let $w\in W$, $\Delta_+$ a system of positive roots in $\Delta$, $\alpha_1,\ldots, \alpha_l$ and $s_1,\ldots, s_l$ the corresponding simple roots and simple reflections, respectively. Then the following statements are true.

(i) Any simple root in $\Delta_w$ can be moved to the first position in a normal ordering of $\Delta_w$ by a composition of elementary transpositions. 

(ii) Any two  normal orderings of the set $\Delta_w$ can be obtained from each other by elementary transpositions within $\Delta_w$. 

(iii) Any two reduced decompositions of $w$ can be obtained from each other using braid group relations in $W$.

(iv) If $w=s_{j_{b}}\ldots s_{j_{1}}$ is a reduced decomposition of $w$ then $\beta_{1}=\alpha_{j_1},\beta_{2}=s_{j_1}\alpha_{j_2}\ldots, \beta_{b}=s_{j_1}\ldots s_{j_{b-1}}\alpha_{j_b}$ is a normal order of the set $\Delta_w$, and 
\begin{equation}\label{wbeta}
w=s_{\beta_1}\ldots s_{\beta_b}.
\end{equation} 

Moreover, there is a one-to-one correspondence between the reduced decompositions of $w$, the presentations of $w$ of the form (\ref{wbeta}), and the normal orderings of the set $\Delta_w$. In particular, any reduced decomposition of $w$ is an initial part of a reduced decomposition of the longest element $\overline{w}\in W$.
\end{lemma}

We shall also need the following generalization of the previous lemma.
\begin{lemma}\label{segmord}
Let $\alpha_1,\ldots, \alpha_l$ be the simple roots in $\Delta_+$, $s_1, \ldots ,s_l$ the corresponding simple reflections, $\overline{w}=s_{i_1}\ldots s_{i_D}$ a reduced decomposition of the longest element $\overline{w}$ of $W$ with respect to the system $s_1, \ldots ,s_l$ of its generators, and
\begin{equation}\label{initdec}
\beta_1=\alpha_{i_1},\beta_2=s_{i_1}\alpha_{i_2},\ldots,\beta_D=s_{i_1}\ldots s_{i_{D-1}}\alpha_{i_D}
\end{equation}
the corresponding normal ordering in $\Delta_+$.

Let $\beta_h, \beta_{h+1}, \ldots , \beta_{h+b}$ be a segment of this normal ordering, $w=s_{i_h}\ldots s_{i_{h+b}}$ the corresponding part of the reduced decomposition of $\overline{w}$. 

Let $w=s_{j_{h}}\ldots s_{j_{h+b}}$ be another reduced decomposition of $w$. Then the following statements are true.

(i) \begin{equation}\label{wreord}
\overline{w}=s_{i_1}\ldots s_{i_{h-1}}s_{j_{h}}\ldots s_{j_{h+b}}s_{i_{h+b+1}}\ldots s_{i_D}
\end{equation} 
is a reduced expression of $\overline{w}$.

(ii) The normal ordering of $\Delta_+$ corresponding to reduced decomposition (\ref{wreord}) has the form
\begin{equation}\label{rootreord}
\beta_1, \beta_2, \beta_{h-1}, \beta_{h_1},\ldots, \beta_{h_b}, \beta_{h+b+1}\ldots,\beta_D,
\end{equation}
where the segment $\beta_{h_1},\ldots, \beta_{h_b}$ is obtained from $\beta_h, \beta_{h+1}, \ldots , \beta_{h+b}$ by applying elementary transpositions to normal ordering (\ref{initdec}) such that each such transposition does not change the positions of all positive roots, except for $\beta_h, \beta_{h+1}, \ldots , \beta_{h+b}$.

(iii) Let $u=s_{i_1}\ldots s_{i_{h-1}}$. Then $\Delta_{uw^{-1}u^{-1}}=u\Delta_{w^{-1}}=\{\beta_h, \beta_{h+1}, \ldots , \beta_{h+b}\}$, and the segments $\beta_h, \beta_{h+1}, \ldots , \beta_{h+b}$ and $\beta_{h_1},\ldots, \beta_{h_b}$ correspond to the following two decompositions of the element $uw^{-1}u^{-1}$,
\begin{equation}\label{wbeta1}
uw^{-1}u^{-1}=s_{\beta_h} \ldots  s_{\beta_{h+b}}=s_{\beta_{h_1}}\ldots s_{\beta_{h_b}}.
\end{equation}
Moreover, if a segment $\beta_{h_1},\ldots, \beta_{h_b}$ is obtained from $\beta_h, \beta_{h+1}, \ldots , \beta_{h+b}$ by applying elementary transpositions to normal ordering (\ref{initdec}) such that each such transposition does not change the positions of all positive roots, except for $\beta_h, \beta_{h+1}, \ldots , \beta_{h+b}$, then the corresponding ordering (\ref{rootreord}) is normal, and there is a one-to-one correspondence between the reduced decompositions of $w$, the decompositions of the element $uw^{-1}u^{-1}$ of the form (\ref{wbeta1}), the segments obtained from $\beta_h, \beta_{h+1}, \ldots , \beta_{h+b}$ by applying elementary transpositions within the segment, and the normal orderings of $\Delta_+$ of type (\ref{rootreord}).
\end{lemma}

\begin{proof}
For part (i) we observe that the decompositions $u=s_{i_1}\ldots s_{i_{h-1}}$, $v=s_{i_{h+1}}\ldots s_{i_D}$ are reduced being parts of a reduced decomposition of $\overline{w}$. The decomposition $\overline{w}=uwv$ given by (\ref{wreord}) must be reduced as its length is equal to the length of $\overline{w}$, i.e. to the  cardinality of the set $\Delta_+$.

To justify (ii) we observe that, according to the definition of the normal ordering given in (\ref{initdec}), after replacing the part $w=s_{i_{h}}\ldots s_{i_{h+b}}$ of the corresponding reduced decomposition $\overline{w}=s_{i_1}\ldots s_{i_D}$ with $w=s_{j_{h}}\ldots s_{j_{h+b}}$, the roots at the first $h-1$ and the last $D-(h+b)$ positions will remain at the same places in new ordering (\ref{rootreord}) as for the roots at these positions we still have the same expressions,
$$
\beta_k=s_{i_1}\ldots s_{i_{k-1}}\alpha_{i_k},~1\leq k< h,
$$
$$
\beta_k=uws_{i_{h+1}}\ldots s_{i_{k-1}}\alpha_{i_k},~h+b< k\leq D.
$$

Let $\beta'_{1}=\alpha_{i_h},\beta'_{2}=s_{i_h}\alpha_{i_{h+1}},\ldots, \beta'_{b}=s_{i_h}\ldots s_{i_{h+b-1}}\alpha_{i_{h+b}}$ be the normal ordering of $\Delta_{w^{-1}}$ associated to the reduced decomposition $w=s_{i_{h}}\ldots s_{i_{h+b}}$ according to Lemma \ref{wN} (iv). Then the segment $\beta_h, \beta_{h+1}, \ldots , \beta_{h+b}$ coincides with $u\beta'_{1}=u\alpha_{i_h},u\beta'_{2}=us_{i_h}\alpha_{i_{h+1}},\ldots, u\beta'_{b}=us_{i_h}\ldots s_{i_{h+b-1}}\alpha_{i_{h+b}}$. As a set this segment coincides with $u\Delta_{w^{-1}}$. 

According to definition (\ref{normorddef}) in the new normal ordering (\ref{wreord}) the roots from the set $u\Delta_{w^{-1}}$ forming the segment $\beta_{h_1},\ldots, \beta_{h_b}$ are placed as follows $u\alpha_{j_h}, us_{j_h}\alpha_{j_{h+1}},\ldots, us_{j_h}\ldots s_{j_{h+b-1}}\alpha_{j_{h+b}}$. Note that according to Lemma \ref{wN} (iv) the segment $\alpha_{j_h}, s_{j_h}\alpha_{j_{h+1}},\ldots, s_{j_h}\ldots s_{j_{h+b-1}}\alpha_{j_{h+b}}$ is the normal ordering of the set $\Delta_{w^{-1}}$ which is obtained from the segment $\beta'_{1}=\alpha_{i_h},\beta'_{2}=s_{i_h}\alpha_{i_{h+1}},\ldots, \beta'_{b}=s_{i_h}\ldots s_{i_{h+b-1}}\alpha_{i_{h+b}}$ by applying elementary transpositions by Lemma \ref{wN} (ii). Therefore the segment 
$$
\beta_{h_1},\ldots, \beta_{h_b}=u\alpha_{j_h}, us_{j_h}\alpha_{j_{h+1}},\ldots, us_{j_h}\ldots s_{j_{h+b-1}}\alpha_{j_{h+b}}
$$
is obtained from the segment $\beta_h, \beta_{h+1}, \ldots , \beta_{h+b}=u\beta'_{1},u\beta'_{2},\ldots, u\beta'_{b}$ by applying elementary transpositions. At the same time part (iv) of Lemma \ref{wN} and the fact that the decomposition $uw$ is reduced imply statement (iii) of this lemma. This completes the proof.

\end{proof}

We shall also need the notion circular normal orderings of the root system ${\Delta}$.
Let $\beta_{1}, \beta_{2}, \ldots, \beta_{D}$ be a normal ordering
of a positive root system $\Delta_+\subset \Delta$. Then one can introduce the corresponding {\it circular normal ordering} of the root system ${\Delta}$ \index{ordering!circular normal of a root system} where
the roots in ${\Delta}$ are located on a circle in
the following way
\begin{center}
\setlength{\unitlength}{0.6mm}
     \begin{picture}(180,120)(-40,0)
     \put(0,50){\makebox(0,0){$\beta_1$}}
     \put(4,69){\makebox(0,0){$\beta_2$}}
     \put(15,85){\circle*{1.5}}
     \put(31,96){\circle*{1.5}}
     \put(50,100){\circle*{1.5}}
     \put(69,96){\circle*{1.5}}
     \put(85,85){\circle*{1.5}}
     \put(96,69){\makebox(0,0){$\beta_D$}}
     \put(100,50){\makebox(0,0){$-\beta_1$}}
     \put(96,31){\makebox(0,0){$-\beta_2$}}
     \put(85,15){\circle*{1.5}}
     \put(69,4){\circle*{1.5}}
     \put(50,0){\circle*{1.5}}
     \put(31,4){\circle*{1.5}}
     \put(15,15){\circle*{1.5}}
     \put(4,31) {\makebox(0,0){-$\beta_D$}}
     \put(64,88){\vector (3,-2){10}}
     \put(36,12){\vector (-3,2){10}}
     \end{picture}
\end{center}
\begin{center}
 Fig. 2
\end{center}

Let $\alpha,\beta\in \Delta$. One says that the segment $[\alpha, \beta]$ of the circle
is {\it minimal} \index{segment!minimal in a circularly ordered root system} if it does not contain the opposite roots $-\alpha$ and $-\beta$ and the root $\beta$ follows after $\alpha$ on the circle above, the circle being oriented clockwise.
In that case one also says that $\alpha < \beta$ in the sense of the circular normal ordering,
\begin{equation}\label{noc}
\alpha < \beta \Leftrightarrow {\rm the ~segment}~ [\alpha, \beta]~{\rm  of~ the ~circle~
is~ minimal}.
\end{equation}

Later we shall need the following property of minimal segments which is a direct consequence of Proposition 3.3 in \cite{kh-t}.
\begin{lemma}\label{minsegm}
Let $[\alpha, \beta]$ be a minimal segment in a circular normal ordering of a root system $\Delta$. Then if $\alpha+\beta$ is a root we have
$$
\alpha<\alpha+\beta<\beta.
$$
\end{lemma}

This lemma immediately implies the following property of minimal segments.
\begin{lemma}\label{segmsub}
Let $G_{\bf k}$ be a semisimple connected algebraic group, $\g_{\bf k}$ its Lie algebra, $H_{\bf k}\subset G_{\bf k}$ a maximal torus, $\Delta$ the root system of the pair $(G_{\bf k},H_{\bf k})$. Then for any minimal segment $[\alpha,\beta]\subset \Delta$ with respect to any circular normal ordering of $\Delta$ the vector subspace of $\g_{\bf k}$ spanned by the root vectors corresponding to the roots from $[\alpha,\beta]$ is an algebraic Lie subalgebra $\n_{[\alpha,\beta]}$ of $\g_{\bf k}$. We denote the subgroup of $G_{\bf k}$ corresponding to this subalgebra by $N_{[\alpha,\beta]}\subset G_{\bf k}$. 
\end{lemma}

The following simple consequences of the definition of the circular ordering will be frequently used in the future.
\begin{lemma}\label{circ+}
(i) Any minimal segment in a circular normal ordering of a root system $\Delta$ of length equal to the number of roots in a positive roots system in $\Delta$ is a system of positive roots.

(ii) Let $\Delta_+$ be a system of positive roots in a root system $\Delta$. Then for any normal ordering of $\Delta_+$ and any $\alpha\in \Delta_+$ the root $\alpha$ cannot be represented as a linear combination of roots from $\Delta_+$ strictly greater xor strictly less than $\alpha$ with real coefficients of the same sign.
\end{lemma}

\begin{proof}
(i) follows from the fact that by Lemma \ref{minsegm} any minimal segment is an additively closed set of roots which does not contain opposite roots. If its length is equal to the number of positive roots then this segment must be a system of positive roots.

(ii) Equip $\Delta$ with the circular ordering associated to the given normal ordering of $\Delta_+$. Let $\Delta_\alpha$ be the minimal segment of length equal to the number of positive roots for which $\alpha$ is the first root. By part (i) $\Delta_\alpha$ is a system of positive roots in $\Delta$. Let $h\in \h$ be an element of the corresponding Weyl chabmer, so that $h(\beta)>0$ for all $\beta\in \Delta_\alpha$. 

Suppose that 
\begin{equation}\label{apr}
\alpha=\sum_{\beta<\alpha, \beta\in \Delta_+}c_\beta \beta,
\end{equation} 
where $c_\beta\geq 0$. Then 
\begin{equation}\label{hab}
h(\alpha)=\sum_{\beta<\alpha, \beta\in \Delta_+}c_\beta h(\beta).
\end{equation}
 
Note that $\{\beta\in \Delta: \beta<\alpha, \beta\in \Delta_+\}\subset -\Delta_\alpha$ by the choice of $\Delta_\alpha$, and $\alpha\in \Delta_\alpha$. Thus $h(\alpha)>0$ and $h(\beta)\leq 0$, and hence (\ref{hab}) leads to a contradiction. We deduce that presentation (\ref{apr}) is impossible. 

The other case is considered in a similar way. This completes the proof.

\end{proof}

Unless otherwise explicitly stated, we shall assume from now on that in sum (\ref{hdec}) $\h_0=\h_{\mathbb{R}}^s=\h^s\cap \h_{\mathbb{R}}$ \index[not]{h@$\h_0$} is the vector subspace of $\h_{\mathbb{R}}$ fixed by the action of $s$ pointwise. If $\h_{\mathbb{R}}^s$ is trivial it will be still convenient to keep the same notation with $\h_0=0$ and include it into sum (\ref{hdec}). 

According to this convention we always have that $\Delta_0=\{\alpha \in \Delta: s\alpha=\alpha\}$ \index[not]{D@$\Delta_0$}is the set of roots fixed by the action of $s$ elementwise, and hence ${\Delta}_{0}$ may be empty. In this case it will be also convenient to add the empty set ${\Delta}_{0}$ to union (\ref{decD}), so that we shall always have $i_0=0$ in (\ref{decD}) and $\Delta_{i_0}=\Delta_0$.

To formulate the main statement of this section we shall consider a special type of factorizations of Weyl group elements introduced in (\ref{inv}). This type is described in the following lemma.
\begin{lemma}
Let $s\in W$ be an element of the Weyl group $W$. There there is a decomposition $s=s^1s^2$ of the form (\ref{inv}) such that any root $\alpha\in \Delta$ fixed by the action of $s^2$ is also fixed by the action of $s$, i.e.
\begin{equation}\label{cond2}
s^2\alpha=\alpha \Rightarrow \alpha \in \Delta_0.
\end{equation}
\end{lemma}

\begin{proof}
Let $s=s^1s^2$ be any presentation of the form (\ref{inv}) of a Weyl group element $s\in W$. Consider a direct vector space decomposition of the form (\ref{hdec}) defined with the help of the factorization $s=s^1s^2$ such that the one--dimensional subspaces $\h_i$ on which $s^1$ acts by multiplication by $-1$ and $s^2$ acts trivially, if there are any non--trivial subspaces of this type, are placed immediately after $\h_0=\h_{\mathbb{R}}^s$, i.e. they are labeled by indexes forming a set of the form $\{1,2,\ldots ,t\}$. 


Suppose that the direct sum $\bigoplus_{k=1}^u\h_{i_k}$ of  the subspaces $\h_{i_k}$, which correspond to the non--empty sets $\Delta_{i_k}$, $k=1,\ldots , M$ and on which $s^1$ acts by multiplication by $-1$ and $s^2$ acts trivially, is not trivial.
Since  the one--dimensional subspaces $\h_i$ on which $s^1$ acts by multiplication by $-1$ and $s^2$ acts trivially are placed immediately after $\h_0$ in sum (\ref{hdec}),
 the roots from the union $\bigcup_{k=1}^u{\Delta}_{i_k}$ must be orthogonal to all subspaces $\h_{i_k}$, $k> 0$ on which $s^1$ does not act by multiplication by $-1$ and to all roots from the set $\gamma_{\widetilde{l}+1}, \ldots \gamma_{l'}$ as $s^2$ acts trivially on $\bigoplus_{k=0}^u\h_{i_k}$. Pick up a root $\gamma \in \bigcup_{k=1}^u{\Delta}_{i_k}$ which belongs to the same set of positive roots as $\gamma_1,\ldots, \gamma_{l'}$. Then $\gamma$ is orthogonal to the roots $\gamma_{\widetilde{l}+1}, \ldots \gamma_{l'}$. Therefore, by the choice of $\gamma$, $s^1_0:=s^1s_\gamma$ is an involution the dimension of the fixed point space of which is equal to the dimension of the fixed point space of the involution $s^1$ plus one, and $s^2_0:=s_\gamma s^2$ is another involution the dimension of the fixed point space of which is equal to the dimension of the fixed point space of the involution $s^2$ minus one,
$$
{\rm dim}~\h^{s^1_0} ={\rm dim}~\h^{s^1}+1, {\rm dim}~\h^{s^2_0} ={\rm dim}~\h^{s^2}-1.
$$ 
By the construction we also have a decomposition $s=s^1_0s^2_0$.

Now we can construct decomposition (\ref{hdec}) using the new factorization $s=s^1_0s^2_0$ of $s$. After that we can also apply the algorithm described in the previous paragraph. Iterating these two steps we shall eventually arrive at the situation when the direct sum $\bigoplus_{k=1}^u\h_{i_k}$ of  the subspaces $\h_{i_k}$ on which $s^1$ acts by multiplication by $-1$ and $s^2$ acts trivially is trivial. 

This property also implies that any root fixed by the action of $s^2$ is fixed by the action of $s$ as well,
$$
s^2\alpha=\alpha \Rightarrow \alpha \in \Delta_0,
$$
which is property (\ref{cond2}).

Indeed, if $s^2\alpha=\alpha$ then by the definition of the subspaces $\h_{i_k}\subset \h_{\mathbb{R}}$, $\alpha$ annihilates any $s$-invariant subspace $\h_{i_k}\subset \h_{\mathbb{R}}$ on which $s^2$ acts in a non--trivial way. Therefore $\alpha$ may not vanish identically only on the direct sum $\bigoplus_{k=1}^u\h_{i_k}$ of  the subspaces $\h_{i_k}$ on which $s^1$ acts by multiplication by $-1$ and $s^2$ acts trivially or on $\h_0$. But the sum $\bigoplus_{k=1}^u\h_{i_k}$ is trivial by the construction. Therefore $\alpha$ does not vanish identically only on $\h_0$ (if $\h_0$ is not trivial) which implies $\alpha \in \Delta_0$.

Note that condition (\ref{cond2}) is a property of the corresponding decomposition (\ref{inv}) and that this condition does not depend on decomposition (\ref{hdec}) used in the arguments in the previous paragraph. This completes the proof.

\end{proof}

From now on, unless otherwise explicitly stated, we shall make the following assumptions.

\begin{enumerate}
\item
We shall only consider decompositions (\ref{inv}) which satisfy property (\ref{cond2}).
\item
As in the proof of the previous lemma, we shall also always consider decompositions of the type (\ref{hdec}) for $s=s^1s^2$ such that the one--dimensional subspaces $\h_i$ in sum (\ref{hdec}) on which $s^1$ acts by multiplication by $-1$ and $s^2$ acts trivially, if there are any non--trivial subspaces of this type, are placed immediately after $\h_0$ in (\ref{hdec}), i.e. they are labeled by indexes forming a set of the form $\{1,2,\ldots ,t\}$.
\item
According to the arguments in the proof of the previous lemma, under assumptions 1 and 2 we can suppose that in a decomposition (\ref{hdec}) associated to $s=s^1s^2$ the direct sum $\bigoplus_{k=1}^u\h_{i_k}$ of the subspaces $\h_{i_k}$ on which $s^1$ acts by multiplication by $-1$, $s^2$ acts trivially, and which correspond to the non--empty sets in the collection $\Delta_{i}$, $i=1,\ldots , M(s)$, is trivial.
\end{enumerate}

Normal orderings of sets of positive roots introduced in the following proposition play a key role in many statements in this book.
\begin{proposition}\label{pord}
Let $s\in W$ be an element of the Weyl group $W$ of the pair $(\g,\ha)$, $\Delta$ the root system of the pair $(\g,\ha)$, $s=s^1s^2$ a factorization (\ref{inv}) for $s$ satisfying condition (\ref{cond2}). 

Let $\Delta_+^s$ be a system of positive roots defined by (\ref{D+sdef}) using the factorization $s=s^1s^2$, and a decomposition (\ref{hdec}) associated to it and satisfying assumptions 1, 2 and 3. Changing the signs of the roots $\gamma_1, \ldots, {\gamma_{l'}}$, which are used in the definition of the factorization (\ref{inv}), we can assume that the roots in each of the sets $\gamma_1, \ldots, \gamma_{\widetilde{l}}$ and ${\gamma_{\widetilde{l}+1}},\ldots, {\gamma_{l'}}$ belong to $\Delta_+^s$ and are mutually orthogonal.

For any $w\in W$ denote $\Delta_w^s=\{\alpha \in \Delta_+^s:w\alpha \in -\Delta_+^s\}$. \index[not]{D@$\Delta_w^s$} Then the following statements are true.

\vskip 0.3cm
\noindent
(i) The decomposition $s=s^1s^2$ is reduced in the sense that ${l}(s)={l}(s^1)+{l}(s^2)$, where ${l}(~\cdot~)$ is the length function in $W$ with respect to the system of simple roots in $\Delta_+^s$, and $\Delta_{s}^s=\Delta_{s^{2}}^s\cup s^2(\Delta_{s^{1}}^s)$, $\Delta_{s^{-1}}^s=\Delta_{s^{1}}^s\cup s^1(\Delta_{s^{2}}^s)$ (disjoint unions). 

\vskip 0.3cm
\noindent
(ii) If $\alpha\in \Delta_{s^{1}}^s$ (resp. $\alpha\in \Delta_{s^{2}}^s$, $\alpha\in \Delta_{s}^s$, or $\alpha\in \Delta_{s^{-1}}^s$), $\beta\in \Delta_0$, and $\alpha+\beta\in \Delta$ then $\alpha+\beta\in \Delta_{s^{1}}^s$ (resp. $\alpha+\beta\in \Delta_{s^{2}}^s$, $\alpha+\beta\in \Delta_{s}^s$, or $\alpha+\beta\in \Delta_{s^{-1}}^s$).  

\vskip 0.3cm
\noindent
(iii) $s^2(\Delta_{s^{1}}^s) \subset \Delta_+^s\setminus \left(\Delta_{s^{1}}^s \cup \Delta_{s^{2}}^s\cup \Delta_0\right)$.

\vskip 0.3cm
\noindent
(iv) There is a normal ordering of the root system $\Delta_+^s$ of the following form
\begin{equation}\label{NO}
\beta_1^1,\ldots, \beta_t^1,\beta_{t+1}^1, \ldots, \beta_{t+\frac{p-\widetilde{l}}{2}}^1, \gamma_1, \ldots , \gamma_2, \ldots,  \gamma_3,\ldots, \gamma_{\widetilde{l}}, \beta_{t+p+1}^1,\ldots, \beta_{l(s^1)}^1,\ldots, 
\end{equation}
$$ 
\beta_1^2,\ldots, \beta_{t'}^2, \gamma_{\widetilde{l}+1}, \ldots, \gamma_{\widetilde{l}+2}, \ldots, \gamma_{\widetilde{l}+3},\ldots, \gamma_{l'}, \beta_{t'+\frac{p'+l'-\widetilde{l}}{2}+1}^2, \ldots,\beta_{t'+p'}^2, \beta_{t'+p'+1}^2,\ldots, \beta_{l(s^2)}^2, \beta_1^0, \ldots, \beta_{D_0}^0, 
$$
where \index[not]{b@$\beta^1_i,\beta^2_i$}
$$
\Delta_{s^1}^s=\{\beta_1^1,\ldots,\beta_{l(s^1)}^1\}
=\{\beta_1^1,\ldots, \beta_t^1,\beta_{t+1}^1, \ldots,\beta_{t+\frac{p-\widetilde{l}}{2}}^1, \gamma_1, \ldots, \gamma_2, \ldots, \gamma_3,\ldots, \gamma_{\widetilde{l}}, \beta_{t+p+1}^1,\ldots, \beta_{l(s^1)}^1\},
$$
$$
\Delta_{s^1}^{-1}:=\{\alpha\in \Delta_+^s|s^1\alpha=-\alpha\}=\{\beta_{t+1}^1, \ldots,\beta_{t+p}^1\}
=\{\beta_{t+1}^1, \ldots, \beta_{t+\frac{p-\widetilde{l}}{2}}^1, \gamma_1, \ldots, \gamma_2, \ldots, \gamma_3,\ldots, \gamma_{\widetilde{l}}\},
$$
\begin{equation}\label{tnum}
p=|\Delta_{s^1}^{-1}|, \index[not]{D@$\Delta_{s^1}^{-1}$} t=|[\beta_1^1,\beta_t^1]|=|[\beta_{t+p+1}^1,\beta_{l(s^1)}^1]|=\frac{l(s^1)-p}{2},
\end{equation} 
so that
$$
2t=|\Delta_{s^1}^s\setminus \Delta_{s^1}^{-1}|=l(s^1)-|\Delta_{s^1}^{-1}|,
$$
\begin{equation}\label{pnum}
|[\beta_{t+1}^1,\beta_{t+\frac{p-\widetilde{l}}{2}}^1]|=\frac{p-\widetilde{l}}{2}, ~|[\gamma_1,\gamma_{\widetilde{l}}]|=\frac{p+\widetilde{l}}{2};
\end{equation}

$$
\Delta_{s^2}^s=\{\beta_1^2,\ldots, \beta_{l(s^2)}^2\}
=\{\beta_1^2,\ldots, \beta_{t'}^2, \gamma_{\widetilde{l}+1}, \ldots , \gamma_{\widetilde{l}+2},\ldots, \gamma_{\widetilde{l}+3},\ldots, 
\gamma_{l'},\beta_{t'+\frac{p'+l'-\widetilde{l}}{2}+1}^2, \ldots,\beta_{t'+p'}^2, \beta_{t'+p'+1}^2,\ldots, \beta_{l(s^2)}^2\},
$$
$$
\Delta_{s^2}^{-1}:=\{\alpha\in \Delta_+^s|s^2\alpha=-\alpha\}=\{\beta_{t'+1}^2,\ldots , \beta_{t'+p'}^2\}
=\{\gamma_{\widetilde{l}+1}, \ldots, \gamma_{\widetilde{l}+2}, \ldots , \gamma_{\widetilde{l}+3},\ldots,
\gamma_{l'},\beta_{t'+\frac{p'+l'-\widetilde{l}}{2}+1}^2, \ldots,\beta_{t'+p'}^2\},
$$
\begin{equation}\label{t'num}
p'=|\Delta_{s^2}^{-1}|, \index[not]{D@$\Delta_{s^2}^{-1}$} t'=|[\beta_1^2,\beta_{t'}^2]|=|[\beta_{t'+p'+1}^2,\beta_{l(s^2)}^2]|=\frac{l(s^2)-p'}{2},
\end{equation} 
so that
$$
2t'=|\Delta_{s^2}^s\setminus \Delta_{s^2}^{-1}|=l(s^2)-|\Delta_{s^2}^{-1}|,
$$
\begin{equation}\label{qnum}
|[\beta_{t'+\frac{p'+l'-\widetilde{l}}{2}+1}^2, \ldots,\beta_{t'+p'}^2]|=\frac{p'-l'+\widetilde{l}}{2}, ~|[\gamma_{\widetilde{l}+1},\gamma_{l'}]|=\frac{p'+l'-\widetilde{l}}{2};
\end{equation}
$$
(\Delta_0)_+:=\{\beta_1^0, \ldots, \beta_{D_0}^0\}=\{\alpha\in \Delta_+^s|s(\alpha)=\alpha\}=\Delta_+^s\cap \Delta_0,  
$$
and $D_0$ is the number of positive roots in $\Delta_+^s$ fixed by the action of $s$. \index[not]{b@$\beta_i^0$} \index[not]{D@$D_0$}

\vskip 0.3cm
Moreover, normal ordering (\ref{NO}) has the following properties.
\vskip 0.3cm
\noindent
(v) The length of the ordered segment $\Delta_{\m_+}\subset \Delta$ in normal ordering (\ref{NO}),
\begin{equation}\label{dn}
\Delta_{\m_+}:=\{\gamma_1, \ldots , \gamma_2, \ldots,  \gamma_3,\ldots, \gamma_{\widetilde{l}}, \beta_{t+p+1}^1,\ldots, \beta_{l(s^1)}^1,\ldots, 
\beta_1^2,\ldots, \beta_{t'}^2, \gamma_{\widetilde{l}+1}, \ldots, \gamma_{\widetilde{l}+2}, \ldots, \gamma_{\widetilde{l}+3},\ldots, \gamma_{l'}\}, 
\end{equation}
is equal to \index[not]{D@$\Delta_{\m_+}$}
\begin{equation}\label{dimm}
|\Delta_{\m_+}|=D-\left(\frac{l(s)-l'}{2}+D_0\right),
\end{equation}
where $D$ is the number of roots in $\Delta_+^s$. \index[not]{D@$D$}

\vskip 0.3cm
\noindent
(vi) For any two roots $\alpha, \beta\in \Delta_{\m_+}$ such that $\alpha<\beta$ the sum $\alpha+\beta$ cannot be represented as a linear combination $\sum_{k=1}^jm_k\gamma_{i_k}$, where $m_k\in \mathbb{N}$ and $\alpha<\gamma_{i_1}<\ldots <\gamma_{i_j}<\beta$.

\vskip 0.3cm
\noindent
(vii) The roots from the set $s^2(\Delta_{s^1}^s)$ form a segment in normal ordering (\ref{NO}) preceding the segment formed by the roots from the set $\Delta_{s^2}^s$ which is, in turn, followed by the final segment $(\Delta_0)_+=\Delta_0\cap \Delta_+^s$. Thus the roots from the set $\Delta_s^s=\Delta_{s^2}^s\cup s^2(\Delta_{s^1}^s)$  form a segment in $\Delta_+^s$ which contains $\Delta_{s^2}^s$.

The roots from the set $\Delta_{s^{1}}^s$ form an initial segment which does not intersect the final segment $\Delta_s^s\cup (\Delta_0)_+$.

\vskip 0.3cm
\noindent
(viii) For any $\alpha\in ({\Delta}_{i_k})_+$, $i_k>0$ such that $s\alpha \in ({\Delta}_{i_k})_+$ one has $s\alpha>\alpha$, and if $\beta, \gamma\in {\Delta}_{i_{j}}\cup \{0\}$, $j<k$ and $s\alpha+\beta, \alpha+\gamma\in \Delta$ then $s\alpha+\beta,\alpha+\gamma\in \Delta_+^s$ and $s\alpha+\beta>\alpha+\gamma$.

In particular, for any $\alpha \in \Delta_+^s$, $\alpha\not\in \Delta_0$ and any $\alpha_0\in \Delta_0$ such that $s\alpha \in \Delta_+^s$ one has $s\alpha>\alpha$ and if $s\alpha +\alpha_0\in \Delta$ then $s\alpha +\alpha_0\in \Delta_+^s$ and $s\alpha +\alpha_0>\alpha$.

\vskip 0.3cm
\noindent
(ix) If $\alpha, \beta\in \Delta_{i_k}\cap \Delta_+^s$, $i_k>0$, $\alpha\leq \beta$ and $s\beta\in \Delta_+^s$ then $\h_{i_k}$ is two--dimensional, and the orthogonal projection of $s\beta$ onto $\h_{i_k}$ is obtained by a clockwise rotation with a non--zero angle and by a rescaling with a positive coefficient from the orthogonal projection of $\alpha$ onto $\h_{i_k}$.
\end{proposition}

\begin{remark}
According to the conditions 1-3 imposed before Proposition \ref{pord} in the decomposition $s=s^1s^2$ we always have $s^2\neq id$ provided $s\neq id$ but in the case when $s$ is an involution we have $s=s^2$ and $s^1=id$. Hence in this case $\Delta_{s^1}^s=\emptyset$, conditions 1-3 imply $\Delta_+^s=\Delta_{s^2}^s \cup (\Delta_0)_+$ (disjoint union), so that normal ordering (\ref{NO}) has the form
$$
\beta_1^2,\ldots, \beta_{t'}^2, \gamma_{\widetilde{l}+1}, \ldots, \gamma_{\widetilde{l}+2}, \ldots, \gamma_{\widetilde{l}+3},\ldots, \gamma_{l'}, \beta_{t'+\frac{p'+l'-\widetilde{l}}{2}+1}^2, \ldots,\beta_{t'+p'}^2, \beta_{t'+p'+1}^2,\ldots, \beta_{l(s^2)}^2, \beta_1^0, \ldots, \beta_{D_0}^0,
$$
and
$$
\Delta_{\m_+}=\{\beta_1^2,\ldots, \beta_{t'}^2, \gamma_{\widetilde{l}+1}, \ldots, \gamma_{\widetilde{l}+2}, \ldots, \gamma_{\widetilde{l}+3},\ldots, \gamma_{l'}\}. 
$$
\end{remark}

\begin{remark}
The most important properties of ordering (\ref{NO}), which we shall need later, are (v) and (vi). As $\Delta_{\m_+}$ is a segment, it is additively closed, and hence the roots subspaces in $\g$ corresponding to the roots from it generate a Lie subalgebra $\m_+$ \index[not]{m@$\m_+$} of dimension $D-\left(\frac{l(s)-l'}{2}+D_0\right)$. We shall define a quantum group counterpart of its enveloping algebra $U(\m_+)$, \index[not]{U@$U(\m_+)$} and condition (vi) will ensure that the quantum group counterpart of $U(\m_+)$ has  a non--trivial character which does not vanish on suitably chosen quantum root vectors corresponding to the roots $\gamma_1,\ldots, \gamma_{l'}$. This subalgebra of the quantum group and the non--trivial character will be used to define q-W--algebras in the framework of the Hecke algebra philosophy outlined in the Introduction. Formula (\ref{dimm}) for the dimension of $\m_+$ will play a crucial role in the proof of the De Concini-Kac-Procesi conjecture. More precisely, we shall obtain and use the following relation
$$
2{\rm dim}~\m_++{\rm dim}~\Sigma_s={\rm dim}~G,
$$
where $\Sigma_s\subset G$ is the transversal slice to the set of conjugacy classes \index{conjugacy class!in an algebraic group} in $G$ that will be defined in part (ii) of Proposition \ref{crosssect} with the help of the same system of positive roots $\Delta_+^s$.
\end{remark}

\begin{proof}

Firstly we describe the sets $({\Delta}_{i_k})_+={\Delta}_{i_k}\cap \Delta_+^s$ for $i_k>0$ and some their subsets. We consider in detail the case when the corresponding $s$--invariant subspaces $\ha_{i_k}$ are two--dimensional planes. The case when $\ha_{i_k}$ are invariant lines on which $s^2$ acts by reflection and $s^1$ acts trivially can be treated in a similar way. 

Assume that $k>0$ is such that $\ha_{i_k}$ is two--dimensional. The plane $\ha_{i_k}$ is shown at Figure 3.

$$
\xy/r10pc/: ="A",-(1,0)="B", "A",+(1,0)="C","A",-(0,1)="D","A",+(0.1,0.57)*{h_{i_k}},"A", {\ar+(0,+0.5)},"A", {\ar@{-}+(0,+1)}, "A";"B"**@{-},"A";"C"**@{-},"A";"D"**@{-},"A", {\ar+(0.6,0.13)},"A",+(0.65,0.18)*{v^2_k},"A",+(0.87,0.18)*{{\Delta}_{i_k}^{2}},"A", {\ar@{-}+(0.9,0.41)}, "A",+(0.32,0)="D", +(-0.038,0.134)="E","D";"E" **\crv{(1.33,0.07)},"A",+(0.45,0.15)*{\psi_k},"A",+(0.49,0.05)*{\psi_k},"A", {\ar+(-0.6,0.23)},"A",+(-0.65,0.28)*{v^1_k},"A",+(-0.87,0.28)*{{\Delta}_{i_k}^{1}},"A", {\ar@{-}+(-0.85,0.72)}, "A",+(-0.32,0)="F", +(0.066,0.21)="G","F";"G" **\crv{(0.67,0.07)},"A",+(-0.45,0.27)*{\varphi_k},"A",+(-0.49,0.08)*{\varphi_k},"A", {\ar@{--}+(0.85,-0.72)},"A",+(0.87,-0.28)*{s^1{\Delta}_{i_k}^{1}},"A", {\ar@{--}+(0.4,0.9)},"A",+(0.80,0.73)*{s^2{\Delta}_{i_k}^{1}},"A", {\ar@{--}+(-0.4,0.9)},"A",+(-0.57,0.79)*{s^1{\Delta}_{i_k}^{2}}
\endxy
$$
\begin{center}
 Fig. 3
\end{center}

The vector $h_{i_k}$ is directed upwards at the picture. By (\ref{wc}) a root $\alpha\in {\Delta}_{i_k}$ belongs to the set $({\Delta}_{i_k})_+$ if and only if $h_{i_k}(\alpha)>0$. Identifying $\ha_\mathbb{R}$ and $\ha_\mathbb{R}^*$ with the help of the bilinear form one can deduce that $\alpha\in {\Delta}_{i_k}$ is in $\Delta_+^s$ if and only if its orthogonal projection onto $\ha_{i_k}$, with respect to the bilinear form, is contained in the upper--half plane shown at Figure 3.

The involutions $s^1$ and $s^2$ act on $\h_{i_k}$ as reflections with respect to the lines orthogonal to the vectors labeled by $v^1_k$ and $v^2_k$, respectively, at Figure 3, the angle between $v^1_k$ and $v^2_k$ being equal to $\pi-\theta_{i_k}/2$. The nonzero projections of the roots from the set $\{\gamma_1, \ldots \gamma_{\widetilde{l}}\}\cap {\Delta}_{i_k}$ onto the plane $\h_{i_k}$ are proportional to the vector $v^1_k$, and the nonzero projections of the roots from the set $\{\gamma_{\widetilde{l}+1}, \ldots , \gamma_{l'}\}\cap {\Delta}_{i_k}$ onto the plane $\h_{i_k}$ are proportional to the vector $v^2_k$. The element $s$ acts on $\ha_{i_k}$ by clockwise rotation with the angle $\theta_{i_k}=2(\varphi_k+\psi_k)$.

For $k=1,\ldots, M$ we denote by ${\Delta}_{i_k}^{(r)}$ \index[not]{D@${\Delta}_{i_k}^{(r)}$} the subset of roots in $({\Delta}_{i_k})_+$ orthogonal projections of which onto $\ha_{i_k}$ are directed along a ray $r\subset \ha_{i_k}$ starting at the origin. We call ${\Delta}_{i_k}^{(r)}$ {\it the family corresponding to the ray} $r$. \index{roots!family corresponding to a ray} Below we shall only consider rays $r$ which correspond to nonempty sets of the form ${\Delta}_{i_k}^{(r)}$, so that $({\Delta}_{i_{k}})_+=\bigcup_{j=1}^{M_k}{\Delta}_{i_k}^{(r_j)}$ (disjoint union of non-empty sets) for some $M_k>0$, $M_k\in \mathbb{N}$. \index[not]{M@$M_k$}

\begin{lemma}\label{lem}
For $k=1,\ldots, M$ and for each ray $r$ the following statements are true.  

\noindent
(i) Each ${\Delta}_{i_k}^{(r)}$ is an additively closed set of roots.

\noindent
(ii) Let ${\Delta}_{i_k}^{(r_1)}$ and  ${\Delta}_{i_k}^{(r_2)}$ be two families corresponding to rays $r_1$ and $r_2$, and $\zeta_1\in {\Delta}_{i_k}^{(r_1)}$, $\zeta_2 \in {\Delta}_{i_k}^{(r_2)}$ two roots such that $\zeta_1+\zeta_2 =\zeta_3\in \Delta$. Then the rays $r_1$ and $r_2$ form an angle strictly less than $\pi$, and  $\zeta_3 \in {\Delta}_{i_k}^{(r_3)}$, where ${\Delta}_{i_k}^{(r_3)}$ is the family corresponding to a ray $r_3$ such that $r_3$ lies inside of the angle formed by $r_1$ and $r_2$.

(iii) Let $0\leq j<k$, and $\zeta_1\in {\Delta}_{i_j}$, $\zeta_2 \in {\Delta}_{i_k}^{(r)}$ two roots such that  $\zeta_1+\zeta_2 =\zeta_3\in \Delta$. Then $\zeta_3\in {\Delta}_{i_k}^{(r)}$.
\end{lemma}

\begin{proof}
All statements are simple consequences of the fact that the sum of the orthogonal projections of any two roots onto $\ha_{i_k}$ is equal to the orthogonal projection of the sum.

For part (i) we observe that the orthogonal projections of any two roots $\alpha, \beta$ from ${\Delta}_{i_k}^{(r)}$ onto $\ha_{i_k}$ have the same direction therefore the orthogonal projection of the sum $\alpha+\beta$ onto $\ha_{i_k}$ has the same direction as the orthogonal projections of $\alpha$ and $\beta$, and hence $\alpha+\beta \in {\Delta}_{i_k}^{(r)}$.

For (ii) it suffices to observe that the sum of the orthogonal projections of $\zeta_1$ and $\zeta_2$ onto $\ha_{i_k}$ is equal to the orthogonal projection of the sum, and the sum of the orthogonal projections of $\zeta_1$ and $\zeta_2$ onto $\ha_{i_k}$ lies inside of the angle formed by $r_1$ and $r_2$. The rays $r_1$ and $r_2$ form an angle strictly less than $\pi$ since by (\ref{wc}) the orthogonal projections of all roots from $({\Delta}_{i_k})_+$ onto $\ha_{i_k}$ belong to an open half plane if $\ha_{i_k}$ is two--dimensional or a half space if $\ha_{i_k}$ is one--dimensional.

Part (iii) follows from the fact that $\zeta_1$ has zero orthogonal projection onto $\ha_{i_k}$.

\end{proof}

Now we prove properties (i), (ii) and (iii) in the statement of the proposition.
Recall that by Theorem C in \cite{C} the roots $\gamma_1,\ldots ,\gamma_{l'}$ form a linear basis in the annihilator ${\ha_{\mathbb{R}}'}^*$ of $\ha_0$ with respect to the pairing between $\ha_{\mathbb{R}}$ and $\ha_{\mathbb{R}}^*$. Therefore $s^1=s_{\gamma_1}\ldots s_{\gamma_{\widetilde{l}}}$ and $s^2=s_{\gamma_{\widetilde{l}+1}}\ldots s_{\gamma_{l'}}$ fix all roots from $\Delta_0\subset \ha_0$.

Thus, from (\ref{decD}) with the convention $\h_0=\h^s_{\mathbb{R}}$, we have disjoint union decompositions 
\begin{equation}\label{decds}
\Delta_{s^{1}}^s=\bigcup_{k=1}^M{{\Delta}_{i_k}^{1}}, \Delta_{s^{2}}^s=\bigcup_{k=1}^M{{\Delta}_{i_k}^{2}}, \Delta_{s}^s=\bigcup_{k=1}^M{{\Delta}_{i_k}^{s}},
\end{equation} 
where $\Delta_{i_k}^1={{\Delta}_{i_k}}\cap \Delta_{s^{1}}^s$, ${{\Delta}_{i_k}^{2}}={{\Delta}_{i_k}}\cap \Delta_{s^{2}}^s$, ${{\Delta}_{i_k}^{s}}={{\Delta}_{i_k}}\cap \Delta_{s}^s$. \index[not]{D@$\Delta_{i_k}^1$} \index[not]{D@$\Delta_{i_k}^2$} \index[not]{D@$\Delta_{i_k}^s$} 

\begin{lemma}\label{lemik}
If $\ha_{i_k}$ is an $s$--invariant plane then the sets ${{\Delta}_{i_k}^{1}}$ and ${{\Delta}_{i_k}^{2}}$ have empty intersection and are the unions of the sets ${\Delta}_{i_k}^{(r)}$ with $r$ belonging to the sectors ${{\Delta}_{i_k}^{1}}$ (resp. ${{\Delta}_{i_k}^{2}}$) at Figure 3. The set ${\Delta}_{i_k}^s$ is the union of the sets ${\Delta}_{i_k}^{(r)}$ with $r$ belonging to the union of the non-overlapping sectors $s^2{\Delta}_{i_k}^{1}$ and ${\Delta}_{i_k}^{2}$ at Figure 3. 

If $\ha_{i_k}$ is an $s$--invariant line on which $s$ acts by multiplication by $-1$ the set ${{\Delta}_{i_k}^{1}}$ is empty and ${{\Delta}_{i_k}^{s}}={{\Delta}_{i_k}^{2}}$. This set is the set ${\Delta}_{i_k}^{(r)}=({\Delta}_{i_k})_+$, where $r$ is the semi-axis in $\ha_{i_k}$ directed along $h_{i_k}$.

\end{lemma}

\begin{proof}
Consider the case when $\ha_{i_k}$ is a plane.
At Figure 3 the elements from the sets ${{\Delta}_{i_k}^{1}}$ and ${{\Delta}_{i_k}^{2}}$ are orthogonally projected onto the interiors of the sectors in the plane $\ha_{i_k}$ labeled by ${{\Delta}_{i_k}^{1}}$ and ${{\Delta}_{i_k}^{2}}$, respectively. Therefore the sets ${{\Delta}_{i_k}^{1}}$ and ${{\Delta}_{i_k}^{2}}$ have empty intersection and are the unions of the sets ${\Delta}_{i_k}^{(r)}$ with $r$ belonging to the sectors ${{\Delta}_{i_k}^{1}}$ (resp. ${{\Delta}_{i_k}^{2}}$).

In the case when $\ha_{i_k}$ is a plane the element $s$ acts on $\ha_{i_k}$ by clockwise rotation with the angle $\theta_{i_k}=2(\varphi_k+\psi_k)$ (see Figure 3). Therefore the set ${\Delta}_{i_k}^s$ consists of the roots the orthogonal projections of which onto $\ha_{i_k}$ belong to the union of the non-overlapping sectors labeled $s^2{\Delta}_{i_k}^{1}$ and ${\Delta}_{i_k}^{2}$ at Figure 3.

If $\ha_{i_k}$ is one--dimensional, we recall that by the assumption 3 made before this proposition there are no one--dimensional subspaces $\ha_{i_k}$ on which $s^1$ acts by multiplication by $-1$. Thus in this case the set ${{\Delta}_{i_k}^{1}}$ is empty, and hence ${{\Delta}_{i_k}^{s}}={{\Delta}_{i_k}^{2}}$. This set is the set ${\Delta}_{i_k}^{(r)}=({\Delta}_{i_k})_+$, where $r$ is the semi-axis in $\ha_{i_k}$ directed along $h_{i_k}$.

\end{proof}
 
By the previous lemma ${{\Delta}_{i_k}^{s}}={{\Delta}_{i_k}^{2}}\cup s^2{{\Delta}_{i_k}^{1}}$ (disjoint union), and hence by (\ref{decds}) 
$$
\Delta_{s}^s=\bigcup_{k=1}^M{{\Delta}_{i_k}^{s}}=\bigcup_{k=1}^M\left({{\Delta}_{i_k}^{2}}\cup s^2{{\Delta}_{i_k}^{1}}\right)=\left(\bigcup_{k=1}^M{{\Delta}_{i_k}^{2}}\right)\bigcup s^2\left(\bigcup_{k=1}^M{{\Delta}_{i_k}^{1}}\right)=\Delta_{s^{2}}^s\cup s^2\Delta_{s^{1}}^s ~\text{(disjoint union)}.
$$ 
In particular, by the results of \S 3 in \cite{Z1} the decomposition $s=s^1s^2$ is reduced in the sense that $l(s)=l(s^1)+l(s^2)$, as $l(s)$, $l(s^1)$ and $l(s^2)$ are equal to the cardinalities of the sets $\Delta_{s}^s$, $\Delta_{s^1}^s$ and $\Delta_{s^2}^s$, respectively. 
Similarly, $\Delta_{s^{-1}}^s=\Delta_{s^{1}}^s\cup s^1(\Delta_{s^{2}}^s)$ (disjoint union). This proves (i).

Part (ii) immediately follows from Lemma \ref{lem} (iii) with $j=0$, Lemma \ref {lemik}, decompositions (\ref{decds}), and similar results for $s^{-1}$.

To justify (iii) we recall that by (\ref{decds}) for any root $\alpha\in \Delta_{s^{1}}^s$ one has $\alpha \not \in \Delta_0$. In fact, in this case $\alpha \in {{\Delta}_{i_k}^{1}}$, where $\ha_{i_k}$ is a two--dimensional plane, as by assumption 3 there are no one--dimensional subspaces $\ha_{i_k}$ on which $s^1$ acts by multiplication by $-1$. Hence using Figure 3 we deduce that the orthogonal projection of $s^2\alpha$ onto $\ha_{i_k}$ belongs to the interior of the sector labeled $s^2{{\Delta}_{i_k}^{1}}$ which is a subset of the upper half plane and does not overlap with the sectors ${\Delta}_{i_k}^{1}$ and ${\Delta}_{i_k}^{2}$. Obviously $s^2\alpha \not \in \Delta_0$. This implies $s^2\alpha \in \Delta_+^s\setminus \left(\Delta_{s^{1}}^s \cup \Delta_{s^{2}}^s\cup \Delta_0\right)$ which proves (iii).

Next, we construct in several steps normal ordering (\ref{NO}) satisfying properties (iv)-(ix).
\vskip 0.3cm
\noindent
{\em Step 1.}
\vskip 0.3cm
First we construct an auxiliary normal ordering on $\Delta_+^s$ satisfying properties summarized in Lemmas \ref{lem0} and \ref{lem1} below. We do it by induction over the sets $(\overline{\Delta}_{i_k})_+:=\overline{\Delta}_{i_k} \cap \Delta_+^s$, $k=0,\ldots , M$, where the sets $\overline{\Delta}_{i_k}$ are defined by (\ref{dik}). \index[not]{D@$(\overline{\Delta}_{i_k})_+$}
Note that by Lemma \ref{parab} each $\overline{\Delta}_{i_k}$ is the root system of a standard Levi subalgebra \index{Lie!subalgebra!Levi standard} in $\g$ with respect to the positive root system $\Delta_+^s$.

For the base of the induction, consider the set $(\overline{\Delta}_{0})_+=({\Delta}_{0})_+$. If it is not empty we fix an arbitrary normal ordering on $(\overline{\Delta}_{0})_+=({\Delta}_{0})_+$. 

If $({\Delta}_{0})_+$ is empty and $\ha_{i_1}$ is one--dimensional then we fix an arbitrary normal ordering on $(\overline{\Delta}_{i_1})_+=({\Delta}_{i_1})_+$.

If $({\Delta}_{0})_+$ is empty and $\ha_{i_1}$ is two--dimensional then we choose a normal ordering in $({\Delta}_{i_1})_+$ in the following way.
First fix an initial arbitrary normal ordering on $({\Delta}_{i_1})_+$. Since by Lemma \ref{lem} each set ${\Delta}_{i_1}^{(r)}$ is additively closed we obtain an induced ordering for ${\Delta}_{i_1}^{(r)}$ which satisfies the defining property for the normal ordering.

Now using these induced orderings on the sets ${\Delta}_{i_1}^{(r)}$ we define an auxiliary normal ordering on $({\Delta}_{i_1})_+$ such that on the sets ${\Delta}_{i_1}^{(r)}$ it coincides with the induced normal ordering defined in the previous paragraph, and
if ${\Delta}_{i_1}^{(r_1)}$ and  ${\Delta}_{i_1}^{(r_2)}$ are two families corresponding to rays $r_1$ and $r_2$ such that $r_2$ lies on the right from $r_1$ in $\ha_{i_1}$ then for any $\alpha\in {\Delta}_{i_1}^{(r_1)}$ and  $\beta \in {\Delta}_{i_1}^{(r_2)}$ one has $\alpha<\beta$. By Lemma \ref{lem} the two conditions imposed on the auxiliary normal ordering in $({\Delta}_{i_1})_+$ are compatible and define it in a unique way for the given initial normal ordering on $({\Delta}_{i_1})_+$. Since $s$ acts by a clockwise rotation on $\ha_{i_1}$ we have $s({\Delta}_{i_1}^{(r)})={\Delta}_{i_1}^{(s(r))}$ for $s(r)$ in the upper--half plane, and hence the new normal ordering satisfies the condition that for any $\alpha\in ({\Delta}_{i_1})_+$ such that $s\alpha \in ({\Delta}_{i_1})_+$ one has $s\alpha>\alpha$.

Now assume that an auxiliary normal ordering has already been constructed for the set $\overline{\Delta}_{i_{k-1}}$ and define it for the set $\overline{\Delta}_{i_k}$.

By Lemma \ref{parab} $\overline{\Delta}_{i_{k-1}}$ is generated by some subset of simple roots of the set of simple roots of $(\overline{\Delta}_{i_{k}})_+$. Therefore there exists an initial normal ordering on $(\overline{\Delta}_{i_{k}})_+$ in which the roots from the set $(\overline{\Delta}_{i_k})_+\setminus (\overline{\Delta}_{i_{k-1}})_+=({\Delta}_{i_k})_+$ form an initial segment and the remaining roots from $(\overline{\Delta}_{i_{k-1}})_+$ are ordered according to the previously defined auxiliary normal ordering. As in the case of the induction base this initial normal ordering gives rise to an induced ordering on each set ${\Delta}_{i_k}^{(r)}$.

Now using these induced orderings on the sets ${\Delta}_{i_k}^{(r)}$ we define an auxiliary normal ordering on $(\overline{\Delta}_{i_{k}})_+$. We impose the following conditions on it. Firstly we require that the roots from the set $({\Delta}_{i_k})_+$ form an initial segment and the remaining roots from $(\overline{\Delta}_{i_{k-1}})_+$ are ordered according to the previously defined auxiliary normal ordering.
Secondly, on the sets ${\Delta}_{i_k}^{(r)}$ the auxiliary normal ordering coincides with the induced normal ordering which is already defined, and if ${\Delta}_{i_k}^{(r_1)}$ and  ${\Delta}_{i_k}^{(r_2)}$ are two families corresponding to rays $r_1$ and $r_2$ such that $r_2$ lies on the right from $r_1$ in $\ha_{i_k}$ then for any $\alpha\in {\Delta}_{i_k}^{(r_1)}$ and  $\beta \in {\Delta}_{i_k}^{(r_2)}$ one has $\alpha<\beta$. By Lemma \ref{lem} the conditions imposed on the auxiliary normal ordering in $({\Delta}_{i_k})_+$ are compatible and define it in a unique way. 

Now we can proceed in a similar way by induction which yields an auxiliary normal ordering on $\Delta_+^s$ of the following form:
\begin{equation}\label{auxord1}
({\Delta}_{i_M})_+,\ldots, ({\Delta}_{i_0})_+,
\end{equation}
where $({\Delta}_{i_k})_+$, $k=0,\ldots, M$ are disjoint segments placed in the auxiliary normal ordering as in (\ref{auxord1}), and for each $k>0$ the segment $({\Delta}_{i_k})_+$ is the disjoint union of the disjoint segments ${\Delta}_{i_k}^{(r_j)}$, $j=1,\ldots, M_k$ which are placed in the auxiliary normal ordering in the following way
\begin{equation}\label{auxord2}
{\Delta}_{i_k}^{(r_1)},\ldots, {\Delta}_{i_k}^{(r_{M_k})},
\end{equation}
where for $i<j$ the ray $r_j$ lies on the right from $r_i$ in $\ha_{i_k}$.

Since $s$ acts by a clockwise rotation on two-dimensional $\ha_{i_k}$ we have $s({\Delta}_{i_k}^{(r)})={\Delta}_{i_k}^{(s(r))}$ for $s(r)$ in the upper--half plane in $\ha_{i_k}$, and hence the new normal ordering satisfies the condition that for any $\alpha\in ({\Delta}_{i_k})_+$ such that $s\alpha \in ({\Delta}_{i_k})_+$ one has $s\alpha>\alpha$.

Note also that the roots from $\overline{\Delta}_{i_{k-1}}$ have zero orthogonal projections onto $\ha_{i_k}$. Therefore if $\alpha\in {\Delta}_{i_k}^{(r)}\subset ({\Delta}_{i_k})_+$, $\beta, \gamma\in \overline{\Delta}_{i_{k-1}}$ are such that $s\alpha \in {\Delta}_{i_k}^{(s(r))}\subset ({\Delta}_{i_k})_+$, $s\alpha+\beta, \alpha+\gamma\in \Delta$ then by (\ref{di}) and (\ref{wc})  $s\alpha+\beta\in {\Delta}_{i_k}^{(s(r))} \subset \Delta_+^s\cap ({\Delta}_{i_k})_+,\alpha+\gamma\in {\Delta}_{i_k}^{(r)}\subset \Delta_+^s\cap ({\Delta}_{i_k})_+$ and $s\alpha+\beta>\alpha+\gamma$ as $s({\Delta}_{i_k}^{(r)})={\Delta}_{i_k}^{(s(r))}$ for $s(r)$ in the upper--half plane.

These properties of the new normal ordering are summarized in the following lemma.
\begin{lemma}\label{lem0}
(i) For $k=1,\ldots, M$ such that $\ha_{i_k}$ is two-dimensional one has $s({\Delta}_{i_k}^{(r)})={\Delta}_{i_k}^{(s(r))}$ if $s(r)$ in the upper--half plane, and ${\Delta}_{i_k}^{(r)}<s({\Delta}_{i_k}^{(r)})$. 

(ii) In particular, for $k=1,\ldots, M$ and for any $\alpha\in {\Delta}_{i_k}^{(r)}$ such that $s\alpha \in ({\Delta}_{i_k})_+$ one has $s\alpha>\alpha$ and if $\beta, \gamma\in \overline{\Delta}_{i_{k-1}}\cup \{0\}$, $s\alpha+\beta, \alpha+\gamma\in \Delta$ then $\alpha+\gamma\in {\Delta}_{i_k}^{(r)}$, $s\alpha+\beta\in {\Delta}_{i_k}^{(s(r))}$ and $s\alpha+\beta>\alpha+\gamma$.
\end{lemma}

Observe that, according to the definition of the auxiliary normal ordering of $\Delta_+^s$ we have the following properties of this normal ordering which can be seen from (\ref{auxord1}) and (\ref{auxord2}).
\begin{lemma}\label{lem1}
The auxiliary normal ordering of $\Delta_+^s=\bigcup_{k=0}^M ({\Delta}_{i_{k}})_+$ (disjoint union) has the following properties.

(i) For any $k=0,\ldots, M$ the roots from the disjoint sets $({\Delta}_{i_{k}})_+$ form segments, and for any $k=0,\ldots, M-1$ one has $({\Delta}_{i_{k}})_+>({\Delta}_{i_{k+1}})_+$; 

(ii) For any $k=1,\ldots, M$ the roots from the sets ${\Delta}_{i_k}^{(r_j)}$, $j=1,\ldots, M_k$ form segments, $({\Delta}_{i_{k}})_+=\bigcup_{j=1}^{M_k}{\Delta}_{i_k}^{(r_j)}$ (disjoint union), and the roots from the set $(\Delta_0)_+=\Delta_0\cap \Delta_+^s$, if it is non-empty, form a final segment;

(iii) For any $k=0,\ldots, M$ such that the corresponding $\ha_{i_k}$ is two--dimensional and for any two segments ${\Delta}_{i_k}^{(r_1)}$ and  ${\Delta}_{i_k}^{(r_2)}$ corresponding to rays $r_1$ and $r_2$ such that $r_2$ lies on the right from $r_1$ in $\ha_{i_k}$ one has ${\Delta}_{i_k}^{(r_1)}<{\Delta}_{i_k}^{(r_2)}$.
\end{lemma}

\vskip 0.3cm
\noindent
{\em Step 2.}
\vskip 0.3cm
Now we modify the auxiliary normal ordering in $\Delta_+^s$ in such a way that the roots from the set $\Delta_{s^{1}}^s$ will form an initial segment, the roots from the set $s^2(\Delta_{s^1}^s)$ will form a segment preceding the segment formed by the roots from the set $\Delta_{s^2}^s$ which will, in turn, be followed by the final segment $(\Delta_0)_+=\Delta_0\cap \Delta_+^s$.

We shall do it with the help of the following lemma.
\begin{lemma}\label{lem2}
Assume that $\Delta_+^s$ is equipped with an arbitrary normal ordering such that the roots from a set ${\Delta}_{i_k}^{(r)}=\{\zeta_1,\ldots, \zeta_a\}$ for some ray $r\subset \ha_{i_k}$  form a segment $\zeta_1,\ldots, \zeta_a$. Suppose also that for some natural $p$ such that $0< p<k$ the roots from a set ${\Delta}_{i_p}^{(t)}=\{\xi_1,\ldots, \xi_b\}$ for some ray $t\subset \ha_{i_p}$ form a segment $\xi_1,\ldots, \xi_b$ and that $\zeta_1,\ldots, \zeta_a,\xi_1,\ldots, \xi_b$ is also a segment. Then applying elementary transpositions within the segment $\zeta_1,\ldots, \zeta_a,\xi_1,\ldots, \xi_b$ one can reduce it to the form $\xi_{i_1},\ldots, \xi_{i_b},\zeta_{j_1},\ldots, \zeta_{j_a}$.
\end{lemma}

\begin{proof}
The proof is by induction. First consider the segment $\zeta_1,\ldots, \zeta_a,\xi_1$.

Since the orthogonal projections of the roots from the set ${\Delta}_{i_p}$ onto $\ha_{i_k}$ are equal to zero, for any $\varsigma\in {\Delta}_{i_p}^{(t)}$ and $\upsilon \in {\Delta}_{i_k}^{(r)}$ such that $\upsilon\pm\varsigma\in \Delta$ we have $\upsilon\pm\varsigma\in {\Delta}_{i_k}^{(r)}$. Assume now that $\varsigma$ and $\upsilon$ are contained in an ordered segment of form (\ref{rank2}) or in a segment with the inverse ordering. By the previous observation this segment contains no other roots from ${\Delta}_{i_p}^{(t)}$, and $\varsigma$ is the first or the last element in that segment. For the same reason the other roots in that segment must also belong to ${\Delta}_{i_k}^{(r)}$. Therefore applying an elementary transposition, if necessarily, one can move $\varsigma$ to the first position in that segment.

Applying this procedure iteratively to the segment $\zeta_1,\ldots, \zeta_a,\xi_1$ we can reduce it to the form $\xi_1,\zeta_{m_1},\ldots, \zeta_{m_a}$ by applying elementary transpositions within the segment $\zeta_1,\ldots, \zeta_a,\xi_1$.

Now we can apply the same procedure to the segment $\zeta_{m_1},\ldots, \zeta_{m_a},\xi_2$ to reduce the segment $\xi_1,\zeta_{m_1},\ldots, \zeta_{m_a},\xi_2$ to the form $\xi_1,\xi_2,\zeta_{l_1},\ldots, \zeta_{l_a}$ by applying elementary transpositions within the segment $\zeta_{m_1},\ldots, \zeta_{m_a},\xi_2$.

Iterating this procedure we obtain the statement of the lemma.

\end{proof}

Now observe that according to Lemma \ref{lem1} (i) for $k=0,\ldots, M$ the roots from each of the sets $({\Delta}_{i_{k}})_+$ form a segment in the auxiliary normal ordering of $\Delta_+^s$, and by Lemma \ref{lem1} (ii) for $k=1,\ldots, M$ the roots from the sets ${\Delta}_{i_k}^{(r)}$ form segments inside $({\Delta}_{i_{k}})_+$. As we observed in Lemma \ref{lemik}, for $k=1,\ldots, M$ the sets ${{\Delta}_{i_k}^{1}}$ (resp. ${{\Delta}_{i_k}^{2}}$)  are the unions of the sets ${\Delta}_{i_k}^{(r)}$ with $r$ belonging to the sectors ${{\Delta}_{i_k}^{1}}$ (resp. ${{\Delta}_{i_k}^{2}}$) at Figure 3, if $\h_{i_k}$ is two--dimensional, and if it is one--dimensional then the set ${{\Delta}_{i_k}^{1}}$ is empty and ${{\Delta}_{i_k}^{s}}={{\Delta}_{i_k}^{2}}$. This set is the set ${\Delta}_{i_k}^{(r)}=({\Delta}_{i_k})_+$, where $r$ is the semi-axis in $\ha_{i_k}$ directed along $h_{i_k}$.
Hence by Lemma \ref{lem1} (i) and (ii) (see also (\ref{auxord1}), (\ref{auxord2}))  the roots from the sets ${{\Delta}_{i_k}^{1}}$ and ${{\Delta}_{i_k}^{2}}$ form an initial and a final segment, respectively, inside $({\Delta}_{i_{k}})_+$.

Therefore we can apply Lemma \ref{lem2} to move all roots from the segments ${{\Delta}_{i_k}^{1}}$, $k=1,\ldots , M$ to the left and to move all roots from the segments ${{\Delta}_{i_k}^{2}}$, $k=1,\ldots , M$ to the right to positions preceding the final segment formed by the roots from $(\Delta_0)_+$. After this modification the roots from the set $\Delta_{s^{1}}^s$ form an initial segment.

Now using similar arguments the roots from the sets $s^2{{\Delta}_{i_k}^{1}}$, $k=1,\ldots , M$ forming segments by Lemma \ref{lemik} and by Lemma \ref{lem1} (i) and (ii) (see also (\ref{auxord1}), (\ref{auxord2})) as well can be moved to the right to positions preceding the final segment formed by the roots from the set $\Delta_{s^{2}}^s\cup(\Delta_0)_+$. Note that by parts (i) and (iii) $\Delta_{s^{1}}^s\cap s^2(\Delta_{s^1}^s)=\emptyset$ and $\Delta_{s^{2}}^s\cap s^2(\Delta_{s^1}^s)=\emptyset$. So the last modification does not affect the positions of the roots in $\Delta_{s^{1}}^s$ and in $\Delta_{s^{2}}^s$ and after applying it the roots from the set $s^2(\Delta_{s^1}^s)$ will form a segment preceding the segment formed by the roots from the set $\Delta_{s^2}^s$ which is, in turn, followed by the final segment $(\Delta_0)_+=\Delta_0\cap \Delta_+^s$. Thus the roots from the set $\Delta_s^s=\Delta_{s^2}^s\cup s^2(\Delta_{s^1}^s)$ form a segment in $\Delta_+^s$ which contains $\Delta_{s^2}^s$. By construction the segment $\Delta_{s^{1}}^s$ does not intersect the final segment $\Delta_s^s\cup (\Delta_0)_+$. This proves that property (vii) holds for the ordering constructed at Step 2.

Thus at Step 2 we obtain a normal ordering on $\Delta_+^s$ of the following form:
\begin{equation}\label{auxord1'}
\Delta_{s^1}^s,({\Delta}_{i_M})_+',\ldots, ({\Delta}_{i_1})_+', s^2(\Delta_{s^1}^s), \Delta_{s^2}^s, (\Delta_0)_+, 
\end{equation}
where $({\Delta}_{i_k})_+'=({\Delta}_{i_k})_+\setminus (\Delta_s^s\cup \Delta_{s^1}^s \cup (\Delta_0)_+)$, \index[not]{D@$({\Delta}_{i_k})_+'$} $k=1,\ldots, M$ are disjoint segments placed in the normal ordering as in (\ref{auxord1'}), and for each $k>0$ the segment $({\Delta}_{i_k})_+'$ is the disjoint union of the disjoint segments ${\Delta}_{i_k}^{'(r_j)}:={\Delta}_{i_k}^{(r_j)}\setminus (\Delta_s^s \cup \Delta_{s^1}^s)\cup (\Delta_0)_+)$, \index[not]{D@${\Delta}_{i_k}^{'(r_j)}$} $j=1,\ldots, M_k$,
\begin{equation}\label{dik'}
({\Delta}_{i_k})_+'=\bigcup_{j=1}^{M_k} {\Delta}_{i_k}^{'(r_j)}, k=1,\ldots, M,
\end{equation}
which are placed in the normal ordering in the following way
\begin{equation}\label{auxord2'}
{\Delta}_{i_k}^{'(r_1)},\ldots, {\Delta}_{i_k}^{'(r_{M_k})},
\end{equation}
where for $i<j$ the ray $r_j$ lies on the right from $r_i$ in $\ha_{i_k}$.

The segment $\Delta_{s^1}^s$ (resp. $\Delta_{s^2}^s$) is the disjoint union of the segments ${\Delta}_{i_k}^{1}$ (resp. ${\Delta}_{i_k}^{2}$), $k=1,\ldots, M$,
\begin{equation}\label{auxord3'}
\Delta_{s^1}^s=\bigcup_{k=1}^M {\Delta}_{i_k}^{1}, \Delta_{s^2}^s=\bigcup_{k=1}^M {\Delta}_{i_k}^{2}
\end{equation}
which are placed in the normal ordering obtained at Step 2 in the following way
\begin{equation}\label{except}
{\Delta}_{i_M}^{1}, \ldots, {\Delta}_{i_1}^{1},~(\text{resp.}~{\Delta}_{i_M}^{2}, \ldots, {\Delta}_{i_1}^{2}),
\end{equation}
and the segment $s^2(\Delta_{s^1}^s)$ is the disjoint union of the segments $s^2({\Delta}_{i_k}^{1})$, $k=1,\ldots, M$
\begin{equation}\label{auxord4'}
s^2(\Delta_{s^1}^s)=\bigcup_{k=1}^M s^2({\Delta}_{i_k}^{1})
\end{equation}
which are placed in the normal ordering obtained at Step 2 in the following way
\begin{equation}\label{auxord5'}
s^2({\Delta}_{i_M}^{1}), \ldots, s^2({\Delta}_{i_1}^{1}).
\end{equation}

Note that ${{\Delta}_{i_k}^{1}}$ (resp. ${{\Delta}_{i_k}^{2}}$) is the disjoint union of the segments ${\Delta}_{i_k}^{(r)}$ with $r$ belonging to the sectors ${{\Delta}_{i_k}^{1}}$ (resp. ${{\Delta}_{i_k}^{2}}$) at Figure 3 if $\h_{i_k}$ is two--dimensional. If $\h_{i_k}$ is one--dimensional then the set ${{\Delta}_{i_k}^{2}}$ is identified with ${\Delta}_{i_k}^{(r)}$, where $r$ is the positive semiaxis directed along $h_{i_k}$. For ${{\Delta}_{i_k}^{1}}$ this situation does not occur. Thus 
\begin{equation}\label{auxord6'}
{{\Delta}_{i_k}^{1}}=\bigcup_{r\in {{\Delta}_{i_k}^{1}}~\text{at Figure 3}}{\Delta}_{i_k}^{(r)}, 
\end{equation}
\begin{equation}\label{auxord6''}
{{\Delta}_{i_k}^{2}}=\left\{\begin{array}{ll} \bigcup_{r\in {{\Delta}_{i_k}^{2}}~\text{at Figure 3}}{\Delta}_{i_k}^{(r)}& {\rm if}~{\rm dim}~\h_{i_k}=2 \\ {\Delta}_{i_k}^{(r)},~\text{where $r$ is the positive semiaxis directed along $h_{i_k}$} & {\rm if}~{\rm dim}~\h_{i_k}=1 \end{array} \right. ,
\end{equation}
so that for $k=1,\ldots, M$
\begin{equation}\label{auxord7'}
{\Delta}_{i_k}^{'(r)}={\Delta}_{i_k}^{(r)}\setminus (\Delta_s^s\cup \Delta_{s^1}^s \cup (\Delta_0)_+)=\left\{\begin{array}{ll} {\Delta}_{i_k}^{(r)} & \text{if} ~r \not\in {\Delta}_{i_k}^{1}\cup s^2({\Delta}_{i_k}^{1})\cup {\Delta}_{i_k}^{2}~\text{at Figure 3} \\ \emptyset & \text{otherwise} \end{array} \right. ,
\end{equation}

The segments $s^2({\Delta}_{i_k}^{1})$, $k=1,\ldots, M$ are non--empty for two--dimensional $\h_{i_k}$ only. They can be represented as the following disjoint unions of segments
\begin{equation}\label{auxord8'}
s^2({\Delta}_{i_k}^{1})=\bigcup_{r\in s^2({\Delta}_{i_k}^{1})~\text{at Figure 3}}{\Delta}_{i_k}^{(r)},
\end{equation}
which are placed in the normal ordering in such a way that
\begin{equation}\label{auxord9'}
{\Delta}_{i_k}^{(r)} < {\Delta}_{i_k}^{(r')},
\end{equation}
where for any $r,r'\in s^2({\Delta}_{i_k}^{1})$ at Figure 3, the ray $r'$ lies on the right from $r$ in $\ha_{i_k}$.
 
The empty sets which may formally appear in the description of the ordering in formulas (\ref{auxord1'}), (\ref{auxord2'}), (\ref{except}), (\ref{auxord5'}), (\ref{auxord8'}), (\ref{auxord9'}) should be omitted.

Note that these formulas imply that the ordering defined at this step has the form
\begin{equation}\label{auxord10'}
\Delta^1,\Delta^2,\ldots, \Delta^{R-1},\Delta^R,
\end{equation}
where $\Delta^j$, \index[not]{D@$\Delta^j$} $j=1,\ldots, R$ \index[not]{R@$R$} are disjoint segments, $\Delta^1=\Delta_{s^1}^s$, $\Delta^{R-1}=\Delta_{s^2}^s$, $\Delta^R=(\Delta_0)_+$, and for $j=2,\ldots, R-2$ one has $\Delta^j={\Delta}_{i_k}^{(r)}$ for some $k$ and $r$.

In formula (\ref{auxord10'}) we could partition further the segments $\Delta^1=\Delta_{s^1}^s$ and $\Delta^{R-1}=\Delta_{s^2}^s$ as in (\ref{auxord6'}) and (\ref{auxord6''}). We shall not do that since at the next step these segments will be reordered.

Note that according to the algorithm given in Lemma \ref{lem2} for each fixed $k$ the
order of the segments formed by the roots from the sets ${\Delta}_{i_k}^{(r)}$ is preserved after applying that lemma. Therefore the new normal ordering obtained this way still satisfies the properties stated in Lemmas \ref{lem0} and \ref{lem1} (ii), (iii). In particular, by Lemma \ref{lem0} (i) for any $j$ such that $\Delta^j\not\subset \Delta_s^s\cup \Delta_0$ one has 
\begin{equation}\label{auxord11'}
\Delta^j<s(\Delta^j).
\end{equation}

\vskip 0.3cm
\noindent
{\em Step 3.}
\noindent
\vskip 0.3cm
Now we can apply elementary transpositions to bring the initial segment formed by the roots from $\Delta_{s^{1}}^s=\{\beta_1^1, \ldots, \beta_{l(s^1)}^1\}$ and the segment formed by the roots from $\Delta_{s^{2}}^s=\{\beta_{1}^2,\ldots, \beta_{l(s^2)}^2\}$ and preceding the final segment $(\Delta_0)_+$ to the form described in (\ref{NO}). 

For this purpose we shall use the following lemma.

\begin{lemma}{\bf (\cite{Ric}, Theorem A; \cite{C}, Lemma 5; \cite{Sp1}, Proposition 3.3)}\label{invdec} 
Let $\Delta_+$ be a system of positive roots in $\Delta=\Delta(\g,\h)$, and $s_1,\ldots, s_l$ be the corresponding simple reflections in the Weyl group $W=W(\g,\h)$ of $\Delta$. Let $z\in W$ be an involution. Then the following statements are true.

(i) There is a Levi subalagbra \index{Lie!subalgebra!Levi} of $\g$ with Cartan subalgebra $\h$ and semisimple part \index{semisimple part!of a Levi subalgebra} $\m_z$, \index[not]{m@$\m_z$} which has Cartan subalgebra $\h_z\subset \h$, \index[not]{h@$\h_z$} such that $z$ is the longest element of the Weyl group $W(\m_z,\h_z)$ \index[not]{W@$W(\m_z,\h_z)$} with respect to the system of simple roots in $\Delta_+(\m_z,\h_z):=\Delta_+\cap \Delta(\m_z,\h_z)$, \index[not]{D@$\Delta(\m_z,\h_z)$} and $z$ acts by multiplication by $-1$ on the Cartan subalgebra $\h_z$. \index[not]{D@$\Delta_+(\m_z,\h_z)$}

(ii) $z$ can be expressed as a product of ${\rm dim}~\h_z$ reflections from the Weyl group $W(\m_z,\h_z)$ of the pair $(\m_z,\h_z)$, with respect to mutually orthogonal roots.

(iii) There is a reduced decomposition of $z$, with respect to the system of simple roots in $\Delta_+$, of the form
\begin{equation}\label{invred}
z=xyx^{-1},
\end{equation}
where $x\in W$, and $y\in W'$ is the longest element in a standard parabolic subgroup $W'$ of $W$ such that $W(\m_z,\h_z)=xW'x^{-1}$, and if $y=s_{j_1}\ldots s_{j_q}$ is a reduced decomposition of $y$ in $W'$ and $x=s_{i_1}\ldots s_{i_t}$ is a reduced decomposition of $x$ then 
\begin{equation}\label{eqroot}
z=s_{i_1}\ldots s_{i_t}s_{j_1}\ldots s_{j_p}s_{i_t} \ldots s_{i_1}
\end{equation}
is a reduced decomposition of $z$ in $W$ with respect to the set of simple roots in $\Delta_+$, so that the corresponding normal ordering of the set $\Delta_{z^{-1}}=\Delta_z$ has the form
\begin{equation}\label{zsegform}
\beta_1,\ldots, \beta_t, \beta_{t+1}, \ldots, \beta_{t+p},\beta_{t+p+1}, \ldots, \beta_{2t+p},
\end{equation}
where $z\beta_{j}=-\beta_{j}$, $j=t+1, \ldots, t+p$, and $\{\beta_{t+1}, \ldots, \beta_{t+p}\}=\Delta_+(\m_z,\h_z)=\Delta(\m_z,\h_z)\cap \Delta_+$. 
\end{lemma}

Let $s_{i_1}\ldots s_{i_{l(s^1)}}$ be the initial part of the reduced decomposition of $\overline {w}$ corresponding to the normal ordering of $\Delta_+^s$ obtained at Step 2, so that $(s^1)^{-1}=s^1=s_{i_1}\ldots s_{i_{l(s^1)}}$.
By parts (i) and (iii) of Lemma \ref{invdec} and by parts (ii) and (iii) of Lemma \ref{segmord} for $z=w=(s^1)^{-1}=s^1=s_{i_1}\ldots s_{i_{l(s^1)}}$ and $u=1$, using elementary transpositions within the initial segment $\Delta_{s^1}^s$ one can reduce its ordering to form (\ref{zsegform}),
\begin{equation}\label{o5}
\beta_1^1,\ldots, \beta_t^1,\beta_{t+1}^1, \ldots,\beta_{t+p}^1,\beta_{t+p+1}^1,\ldots, \beta_{l(s^1)}^1,
\end{equation}
where $\beta_{t+1}^1, \ldots,\beta_{t+p}^1$ is a normal ordering of the system $\Delta_+^s(\m_{s^1},\h_{s^1}):=\Delta(\m_{s^1},\h_{s^1})\cap \Delta_+^s=\Delta_{s^1}^{-1}$ \index[not]{D@$\Delta_+^s(\m_{s^1},\h_{s^1})$} of positive roots in the root system $\Delta(\m_{s^1},\h_{s^1})$ \index[not]{D@$\Delta(\m_{s^1},\h_{s^1})$} of the pair $(\m_{s^1},\h_{s^1})$, \index[not]{m@$\m_{s^1}$} so that \index[not]{h@$\h_{s^1}$} $p$ is the number of positive roots in $\Delta_+^s(\m_{s^1},\h_{s^1})=\Delta_{s^1}^{-1}$. 

Note that by (\ref{eqroot}) one has $t=l(s^1)-(t+p)$, i. e. there are equal numbers of roots on the left and on the right from the segment $\beta_{t+1}^1, \ldots,\beta_{t+p}^1$ in the segment $$\beta_1^1,\ldots, \beta_t^1,\beta_{t+1}^1, \ldots,\beta_{t+p}^1,\beta_{t+p+1}^1,\ldots, \beta_{l(s^1)}^1,$$ and 
\begin{equation}\label{l1*}
t=\frac{l(s^1)-p}{2}.
\end{equation}
This proves (\ref{tnum}).

Observe that $s^1=s_{\gamma_1}\ldots s_{\gamma_{\widetilde{l}}}$ is the expression mentioned in part (ii) of Lemma \ref{invdec}, and the roots ${\gamma_1},\ldots ,{\gamma_{\widetilde{l}}}$ form a linear basis of the Cartan subalgebra $\h_{s^1}$.

Now according to the results of Appendix 1, applying elementary transpositions we can reduce the ordering $\beta_{t+1}^1, \ldots,\beta_{t+p}^1$ to the form compatible with the decomposition $s^1=s_{\gamma_1}\ldots s_{\gamma_{\widetilde{l}}}$, i.e., we can bring it to the form
\begin{equation}\label{norm1}
\beta_{t+1}^1, \ldots, \beta_{t+\frac{p-\widetilde{l}}{2}}^1, \gamma_1, \ldots, \gamma_2, \ldots, \gamma_3,\ldots, \gamma_{\widetilde{l}},
\end{equation}
where in the last formula we relabeled the roots $\beta_i^1$ in such a way that after reordering $\beta_i^1<\beta_j^1$ if and only if $i<j$, and, according to the definition of normal orderings compatible with decompositions of Weyl group involutions given in Appendix 1, the new ordering has the property that for any two positive roots $\alpha, \beta\in \Delta_+^s(\m_{s^1},\h_{s^1})$ such that $\gamma_1\leq \alpha<\beta$ the sum $\alpha+\beta$ cannot be represented as a linear combination $\sum_{k=1}^jc_k\gamma_{i_k}$, where $c_k\in \mathbb{N}$ and $\alpha<\gamma_{i_1}<\ldots <\gamma_{i_k}<\beta$.
Formula (\ref{pnum}) follows from the definition given in Appendix 1 (see formula (\ref{pnuma})).

Note that by Lemma \ref{segmord} the elementary transpositions applied to obtain ordering (\ref{norm1}) do not affect the positions of the roots which do not belong to the set $\Delta_+^s(\m_{s^1},\h_{s^1})=\Delta_{s^1}^{-1}$.

Now let $\overline{w}=s_{i_1}\ldots s_{i_D}$ be the reduced decomposition $\overline{w}$ corresponding to the normal ordering of $\Delta_+^s$ constructed so far. Denote by $s_{i_h}$ the simple reflection in this decomposition the position of which corresponds to the position of the last root preceding the segment $\Delta_{s^2}^s$ in the normal ordering of $\Delta_+^s$ constructed so far. Let $w=s_{i_{h+1}}\ldots s_{i_{h+l(s^2)}}$, $u=s_{i_1}\ldots s_{i_{h}}$. Then $s^2=(s^2)^{-1}=uw^{-1}u^{-1}=uwu^{-1}$ by Lemma \ref{segmord} (iii).


Now, similarly to the case of the segment $\Delta_{s^{1}}^s$, we can apply parts (i) and (iii) of Lemma \ref{invdec} and parts (ii) and (iii) of Lemma \ref{segmord} with $z=w=s_{i_{h+1}}\ldots s_{i_{h+l(s^2)}}$, $u=s_{i_1}\ldots s_{i_{h}}$ to the ordered segment $\Delta_{s^2}^s$ to reduce its ordering to form (\ref{zsegform}). 

Also, as Appendix 1, applying elementary transpositions to the segment $\Delta_+^s(\m_{s^2},\h_{s^2}):=\Delta_{s^2}^{-1}=\Delta(\m_{s^2},\h_{s^2})\cap \Delta_+^s\subset \Delta_{s^2}^s$, \index[not]{D@$\Delta_+^s(\m_{s^2},\h_{s^2})$} which is a system of positive roots in the root system $\Delta(\m_{s^2},\h_{s^2})$ \index[not]{D@$\Delta(\m_{s^2},\h_{s^2})$} of the \index[not]{h@$\h_{s^2}$} pair $(\m_{s^2},\h_{s^2})$, \index[not]{m@$\m_{s^2}$} we can reduce its normal ordering to the form inverse to that compatible with the decomposition $s^2=s_{\gamma_{\widetilde{l}+1}}\ldots s_{\gamma_{l'}}$. So that finally we obtain the following normal ordering of the set $\Delta_+^s$
\begin{equation}\label{no}
\beta_1^1,\ldots, \beta_t^1,\beta_{t+1}^1, \ldots, \beta_{t+\frac{p-\widetilde{l}}{2}}^1, \gamma_1, \ldots , \gamma_2, \ldots,  \gamma_3,\ldots, \gamma_{\widetilde{l}}, \beta_{t+p+1}^1,\ldots, \beta_{l(s^1)}^1,\ldots, 
\end{equation}
$$ 
\beta_1^2,\ldots, \beta_{t'}^2, \gamma_{\widetilde{l}+1}, \ldots, \gamma_{\widetilde{l}+2}, \ldots, \gamma_{\widetilde{l}+3},\ldots, \gamma_{l'}, \beta_{t'+\frac{p'+l'-\widetilde{l}}{2}+1}^2, \ldots,\beta_{t'+p'}^2, \beta_{t'+p'+1}^2,\ldots, \beta_{l(s^2)}^2, \beta_1^0, \ldots, \beta_{D_0}^0,
$$
where
$$
\gamma_{\widetilde{l}+1}, \ldots, \gamma_{\widetilde{l}+2}, \ldots , \gamma_{\widetilde{l}+3},\ldots,
\gamma_{l'},\beta_{t'+\frac{p'+l'-\widetilde{l}}{2}+1}^2, \ldots,\beta_{t'+p'}^2
$$
is the normal ordering of the system of positive roots $\Delta_+^s(\m_{s^2},\h_{s^2})$ inverse to that compatible with the decomposition $s^2=s_{\gamma_{\widetilde{l}+1}}\ldots s_{\gamma_{l'}}$, and we use the notation from part (iv) of Proposition \ref{pord}. 

By construction normal ordering (\ref{no}) has the required form (\ref{NO}). Formulas (\ref{qnum}) and (\ref{t'num}) are established similarly to (\ref{pnum}) and (\ref{tnum}), respectively. This completes the proof of part (iv).

\vskip 0.3cm

We claim that normal ordering (\ref{no}) has properties (v)-(ix) listed in the statement of this proposition. 

Firstly we recall that according to Lemma \ref{segmord} the elementary transpositions used at Step 3 do not affect the positions of the roots which do not belong to $\Delta_{s^{1}}^s$ and $\Delta_{s^{2}}^s$. We state this property for future references as a lemma.
\begin{lemma}\label{s-1}
The elementary transpositions used at Step 3 to bring the segments $\Delta_{s^{1}}^s$ and $\Delta_{s^{2}}^s$ to the form required in (\ref{NO}) do not affect the positions of the roots which do not belong to $\Delta_{s^{1}}^s$ and $\Delta_{s^{2}}^s$. Thus the normal ordering obtained at Step 3 has properties described in (\ref{auxord1'})-(\ref{auxord11'}), except for (\ref{except}). 
\end{lemma}  

From this lemma we deduce that the normal ordering constructed at Step 3 still has form (\ref{auxord1'}), and hence property (vii) holds for normal ordering (\ref{no}) as it holds for the ordering constructed at Step 2.

For for $\alpha\in {\Delta}_{i_k}^{(r)}$, $i_k>0$ we still have $\alpha\in {\Delta}_{i_k}^{(s(r))}$ if $s\alpha\in \Delta_+^s$ by the definition of the sets ${\Delta}_{i_k}^{(r)}$.
We claim that Lemma \ref{lem0} still holds for the normal ordering obtained at Step 3. We shall need a slightly more detailed version of it which implies property (viii).  
\begin{lemma}\label{lem0'}
(i) Suppose that $k>0$. Then for any two-dimensional $\ha_{i_k}$ one has $s({\Delta}_{i_k}^{(r)})={\Delta}_{i_k}^{(s(r))}$ if $s(r)$ in the upper--half plane, and ${\Delta}_{i_k}^{(r)}<s({\Delta}_{i_k}^{(r)})$. 

(ii) In particular, for any $\alpha\in {\Delta}_{i_k}^{(r)}$ such that $s\alpha \in ({\Delta}_{i_k})_+$ one has $s\alpha>\alpha$ and if $\beta, \gamma\in \overline{\Delta}_{i_{k-1}}\cup \{0\}$, $s\alpha+\beta, \alpha+\gamma\in \Delta$ then $\alpha+\gamma\in {\Delta}_{i_k}^{(r)}$, $s\alpha+\beta\in {\Delta}_{i_k}^{(s(r))}$ and $s\alpha+\beta>\alpha+\gamma$.

(iii)
\begin{equation}\label{<<}
s(\Delta_{s^1}^s)\subset \Delta_+^s,~s(\Delta_{s^1}^s)>\Delta_{s^1}^s.
\end{equation}
\end{lemma}
\begin{proof}
First observe that by Proposition \ref{pord} (iii), which has been already proved, $s^2(\Delta_{s^1}^s)\subset \Delta_+^s\setminus(\Delta_{s^1}^s\cup \Delta_{s^2}^s\cup \Delta_0)$. Since $\Delta_{s^1}^s=s^1(-\Delta_{s^1}^s)$ we deduce that $s(\Delta_{s^1}^s)=s^1s^2(\Delta_{s^1}^s)\subset \Delta_+^s$, and $s(\Delta_{s^1}^s)=s^1s^2(\Delta_{s^1}^s)>\Delta_{s^1}^s$.
This proves (\ref{<<}). 

Note that by Lemma \ref{s-1} the normal ordering obtained at Step 3 still has properties described in (\ref{auxord1'})-(\ref{auxord7'}), except for (\ref{except}). By (\ref{auxord3'}), (\ref{auxord4'}) and (\ref{auxord6'}) formula (\ref{<<}) also implies that part (i), and hence (ii), hold for $r$ which belong to the sector ${{\Delta}_{i_k}^{1}}$ at Figure 3. 

If $r$ belongs to the union of the sectors ${{\Delta}_{i_k}^{2}}$ and $s^2({{\Delta}_{i_k}^{1}})$ at Figure 3 then ${\Delta}_{i_k}^{(r)}\subset \Delta_s^s$, so $s({\Delta}_{i_k}^{(r)})\subset \Delta_-^s$, and properties (i) and (ii) are void.

If $r \not\in {\Delta}_{i_k}^{1}\cup {\Delta}_{i_k}^{2}\cup s^1({\Delta}_{i_k}^{2})$ at Figure 3 then ${\Delta}_{i_k}^{(r)}={\Delta}_{i_k}^{'(r)}\subset (\Delta_{i_k})_+'$, and properties (i) and (ii) hold for the normal ordering defined at Step 2 which may only differ by positions of roots in $\Delta_{s^{1}}^s$ and $\Delta_{s^{2}}^s$ from the normal ordering defined at Step 3 by Lemma \ref{s-1}. In particular, $\Delta_{s^{2}}^s>{\Delta}_{i_k}^{(r)}$ with respect to both the normal ordering defined at Step 2 and the normal ordering defined at Step 3. Therefore (i) and (ii) hold for the normal ordering defined at Step 3 as well. This completes the proof.  

\end{proof}

Parts (i) and (ii) of the previous lemma imply property (viii) since when the corresponding $\h_{i_k}$ is one--dimensional this property is void.

If $\alpha, \beta\in \Delta_{i_k}\cap \Delta_+^s$, $\alpha\leq\beta$ and $s\beta\in \Delta_+^s$ then $\h_{i_k}$ must be two dimensional, as otherwise $s\beta\in \Delta_-^s$. We show that the orthogonal projection of $s\beta$ onto $\h_{i_k}$ is obtained by a clockwise rotation with a non--zero angle and by a rescaling with a positive coefficient from the orthogonal projection of $\alpha$ onto $\h_{i_k}$.

Indeed, observe that by Lemma \ref{s-1} the elementary transpositions which we used at Step 3 to bring the segments formed by the roots from the sets $\Delta_{s^1}^s$ and $\Delta_{s^2}^s$ to the form required in normal ordering (\ref{NO}) do not affect the positions of the other roots and after this rearrangement the orthogonal projections of the roots from $\Delta_{i_k}^{1}$ (resp. $\Delta_{i_k}^{2}$) onto $\h_{i_k}$ still belong to the sectors labeled $\Delta_{i_k}^{1}$ (resp. $\Delta_{i_k}^{2}$) at Figure 3.

Therefore by Lemma \ref{lem1} (iii) if $\alpha, \beta\in \Delta_{i_k}\cap \Delta_+^s$, $\beta \in \Delta_{i_k}^1$ and $\alpha\leq \beta$ then the orthogonal projection of $\alpha$ onto $\h_{i_k}$ belongs to the sector labeled $\Delta_{i_k}^1$ at Figure 3.
On the other hand since $s\beta\in \Delta_+^s$ the orthogonal  projection of $s\beta$ onto $\h_{i_k}$ belongs to the upper half plane and does not belong to the sector labeled $\Delta_{i_k}^1$ at Figure 3 as $s$ acts on $\h_{i_k}$ by clockwise rotation by the angle $\theta_{i_k}=2(\varphi_k+\psi_k)>2\varphi_k$. Thus the orthogonal projection of $s\beta$ onto $\h_{i_k}$ is obtained by a clockwise rotation with a non--zero angle and by a rescaling with a positive coefficient from the orthogonal projection of $\alpha$ onto $\h_{i_k}$.

The case when $\beta \in ({\Delta}_{i_k})_+$, $\beta \not \in {{\Delta}_{i_k}^{1}}$, $s\beta\in ({{\Delta}_{i_k}})_+$ is treated in a similar way with the help of Lemma \ref{lem1} (iii). This proves (ix).

By the definition of normal ordering (\ref{NO}) the number of roots in the segment $\Delta_{\m_+}$,
$$
\Delta_{\m_+}=\{\gamma_1, \ldots , \gamma_2, \ldots,  \gamma_3,\ldots, \gamma_{\widetilde{l}}, \beta_{t+p+1}^1,\ldots, \beta_{l(s^1)}^1,\ldots, 
\beta_1^2,\ldots, \beta_{t'}^2, \gamma_{\widetilde{l}+1}, \ldots, \gamma_{\widetilde{l}+2}, \ldots, \gamma_{\widetilde{l}+3},\ldots, \gamma_{l'}\}
$$
is  
$$
|\Delta_{\m_+}|=D-|[\beta_1^1,\beta_{t+\frac{p-\widetilde{l}}{2}}^1]|-|[\beta_{t'+\frac{p'+l'-\widetilde{l}}{2}+1}^2,\beta_{l(s^2)}^2]|-|[\beta_1^0,\beta_{D_0}^0]|=
$$
$$
=D-D_0-|[\beta_1^1,\beta_t^1]|-|[\beta_{t+1}^1,\beta_{t+\frac{p-\widetilde{l}}{2}}^1]|-|[\beta_{t'+\frac{p'+l'-\widetilde{l}}{2}+1}^2, \beta_{t'+p'}^2]|-|[\beta_{t'+p'+1}^2,\beta_{l(s^2)}^2]|,
$$
where $D_0=|[\beta_1^0,\beta_{D_0}^0]|$ is the number of positive roots fixed by the action of $s$.
By formulas (\ref{tnum}), (\ref{pnum}), (\ref{t'num}) and (\ref{qnum}) the last expression takes the form 
$$
|\Delta_{\m_+}|=D-D_0-t-\frac{p-\widetilde{l}}{2}-t'-\frac{p'-(l'-\widetilde{l})}{2}=
$$
$$
=D-D_0-\frac{l(s^1)-p}{2}-\frac{l(s^2)-p'}{2}-\frac{p-\widetilde{l}}{2}-\frac{p'-(l'-\widetilde{l})}{2}=D-\left(\frac{l(s)-l'}{2}+D_0\right) ,
$$
where $l(s)=l(s^1)+l(s^2)$ is the length of $s$. This establishes (v).

Now let $\alpha, \beta\in \Delta_{\m_+}$, be any two roots such that $\alpha<\beta$. We shall show that the sum $\alpha+\beta$ cannot be represented as a linear combination $\sum_{k=1}^jm_k\gamma_{i_k}$, where $m_k\in \mathbb{N}$ and $\alpha<\gamma_{i_1}<\ldots <\gamma_{i_j}<\beta$.

Suppose that such a decomposition exists, $\alpha+\beta=\sum_{k=1}^jm_k\gamma_{i_k}$. Obviously at least one of the roots $\alpha, \beta$ must belong to the set $\Delta_+^s(\m_{s^1},\h_{s^1})\cap \Delta_{\m_+}$ or to the set $\Delta_+^s(\m_{s^2},\h_{s^2})\cap \Delta_{\m_+}$ for otherwise the set of roots $\gamma_{i_k}$ such that $\alpha<\gamma_{i_k}<\beta$ is empty because by the definition of $\Delta_{\m_+}$, $(\Delta_{\m_+}\setminus (\Delta_+^s(\m_{s^1},\h_{s^1})\cup\Delta_+^s(\m_{s^2},\h_{s^2})))\cap \{\gamma_1,\ldots , \gamma_{l'}\}=\emptyset$.

Suppose that $\alpha \in \Delta_+^s(\m_{s^1},\h_{s^1})\cap \Delta_{\m_+}$. The other cases are considered in a similar way.

If $\beta \not \in \Delta_+^s(\m_{s^2},\h_{s^2})\cap \Delta_{\m_+}$ then $\alpha+\beta=\sum_{k=1}^jm_k\gamma_{i_k}$, and $\gamma_{i_k}\leq \gamma_{\widetilde{l}}$ for $k=1,\ldots, j$. In particular, since $\alpha \in \h_{s^1}$ and $\gamma_{i_k}\in \h_{s^1}$ if $\gamma_{i_k}\leq \gamma_{\widetilde{l}}$, we have $\beta=\sum_{k=1}^jm_k\gamma_{i_k}-\alpha \in \h_{s^1}$. This is impossible by the definition of the ordering of the set $\Delta_+^s(\m_{s^1},\h_{s^1})$ compatible with the decomposition $s^1=s_{\gamma_1}\ldots s_{\gamma_{\widetilde{l}}}$.

If $\beta \in \Delta_+(\m_{s^2},\h_{s^2})\cap \Delta_{\m_+}$ then $\alpha+\beta=\sum_{k=1}^jm_k\gamma_{i_k}=\sum_{i_k\leq \widetilde{l}}m_k\gamma_{i_k}+\sum_{i_k>\widetilde{l}}m_k\gamma_{i_k}$. This implies
$$
\alpha-\sum_{i_k\leq \widetilde{l}}m_k\gamma_{i_k}=\sum_{i_k>\widetilde{l}}m_k\gamma_{i_k}-\beta.
$$
The left hand side of the last formula is an element of $\h_{s^1}$ and the right hand side is an element $\h_{s^2}$. Since $\h'=\h_{s^1}+\h_{s^2}$ is a direct vector space decomposition we infer that
$$
\alpha=\sum_{i_k\leq \widetilde{l},\alpha<\gamma_{i_k}}m_k\gamma_{i_k}
$$
and
$$
\beta=\sum_{i_k>\widetilde{l},\gamma_{i_k}<\beta}m_k\gamma_{i_k}.
$$
But this is impossible by Lemma \ref{circ+} (i). Therefore the sum $\alpha+\beta$, $\alpha<\beta$, $\alpha,\beta\in \Delta_{\m_+}$ cannot be represented as a linear combination $\sum_{k=1}^jm_k\gamma_{i_k}$, where $m_k\in \mathbb{N}$ and $\alpha<\gamma_{i_1}<\ldots <\gamma_{i_j}<\beta$. This confirms (vi) and completes the proof of the proposition.



\end{proof}

We shall also need another system of positive roots associated to (the conjugacy class of) the Weyl group element $s$. \index{conjugacy class!in a Weyl group} Consider the circular normal ordering of $\Delta$ corresponding to a system of positive roots $\Delta_+^s$ and to its normal ordering introduced in Proposition \ref{pord}. The minimal segment which consists of the roots $\alpha$ satisfying $\gamma_1\leq \alpha < -\gamma_1$ is a system of positive roots by Lemma \ref{circ+} (i).

\begin{definition}\label{circorddef}
Let $\Delta_+=\{\alpha\in \Delta: \gamma_1\leq \alpha < -\gamma_1\}$, \index[not]{D@$\Delta_+$} where the inequalities are with respect to the circular normal ordering of $\Delta$ corresponding to a system of positive roots $\Delta_+^s$ and to its normal ordering (\ref{NO}). The system of positive roots $\Delta_+$ equipped with the normal ordering induced by the circular normal ordering is called a normally ordered system of positive roots associated to the (conjugacy class of) the Weyl group element $s\in W$. \index{ordering!circular normal of a root system!associated to (the conjugacy class of) a Weyl group element}
\end{definition}

Note that for the root system $\Delta_+$ introduced in Definition \ref{circorddef} $\Delta^s_+\cap \Delta_+ =\Delta^s_+\setminus \Delta_{s^1}^-$, where $\Delta_{s^1}^-=\{\beta_1^1,\ldots, \beta_{t+\frac{p-\widetilde{l}}{2}}^1\}=\Delta^s_+\cap \Delta_-$. \index[not]{D@$\Delta_{s^1}^-$} Therefore
$$
\Delta_+=(\Delta^s_+\cap \Delta_+)\cup \{-\beta_1^1,\ldots, -\beta_{t+\frac{p-\widetilde{l}}{2}}^1\}=(\Delta^s_+\cap \Delta_+)\cup (-\Delta_{s^1}^-).
$$
By Lemma \ref{wN} (iv) this implies that $w_s\Delta_+=\Delta^s_+$, where $w_s=s_{\beta_1^1}\ldots, s_{\beta_{t+\frac{p-\widetilde{l}}{2}}^1}$, \index[not]{w@$w_s$} and if $\overline{w}=s_{i_1}\ldots s_{i_{t+\frac{p-\widetilde{l}}{2}}}\ldots s_{i_D}$ is the reduced decomposition corresponding to normal ordering (\ref{NO}) then $w_s=s_{i_{t+\frac{p-\widetilde{l}}{2}}}\ldots s_{i_1}$.


We have the following property of the length of $s$ with respect to the sets of simple roots in $\Delta_+^s$.
\begin{proposition}
For all systems of positive roots $\Delta_+^s$ the lengths $l(s)$ of $s$ with respect to the sets of simple roots in $\Delta_+^s$ are the same, and they are equal to the length of $s$ with respect to the set of simple roots in any system of positive roots $\Delta_+$ associated to $s$.
\end{proposition}

\begin{proof}
The first statement is a consequence of the definition of $\Delta_+^s$.  

To prove the second assertion we recall that by (\ref{wc}) a root $\alpha\in {\Delta}_{i_k}$, $k>0$ belongs to the set $({\Delta}_{i_k})_+=({\Delta}_{i_k})_+\cap \Delta_+^s$ if and only if $h_{i_k}(\alpha)>0$. Identifying $\h_\mathbb{R}$ and $\h_\mathbb{R}^*$ with the help of the bilinear form one can deduce that $\alpha\in {\Delta}_{i_k}$ is in $\Delta_+^s$ if and only if its orthogonal projection onto $\h_{i_k}$ is contained in the upper--half plane shown at Figure 3 for two--dimensional $\h_{i_k}$, or belongs to the positive semiaxis in $\h_{i_k}$ directed along $h_{i_k}$ if $\h_{i_k}$ is one--dimensional.

According to the definition of the set $\Delta_+$,
\begin{equation}\label{dik+}
{\Delta}_{i_k}\cap \Delta_+=\left(({\Delta}_{i_k})_+\setminus \{\alpha\in ({\Delta}_{i_k})_+\cap \Delta_{s^1}^s: \alpha <\gamma_1\}\right)\cup \{-\alpha: \alpha\in ({\Delta}_{i_k})_+\cap \Delta_{s^1}^s, \alpha <\gamma_1\},
\end{equation}
i.e ${\Delta}_{i_k}\cap \Delta_+$ is obtained from $({\Delta}_{i_k})_+$ by removing some roots which belong to the set ${\Delta}_{i_k}^{1}$ and the orthogonal projections of which onto $\h_{i_k}$ belong to the sector labeled ${\Delta}_{i_k}^{1}$ at Figure 3, and by adding the opposite negative roots the orthogonal projections of which onto $\h_{i_k}$ belong to the sector labeled $s^1{\Delta}_{i_k}^{1}$ at Figure 3. Note that if $\h_{i_k}$ is one--dimensional the set ${\Delta}_{i_k}^{1}$ is empty.

Recall that the involutions $s^1$ and $s^2$ act in $\h_{i_k}$ as reflections with respect to the lines orthogonal to the vectors labeled by $v^1_k$ and $v^2_k$, respectively, at Figure 3, the angle between $v^1_k$ and $v^2_k$ being equal to $\pi-\theta_{i_k}/2$. Therefore the element $s$ acts on $\h_{i_k}$ by clockwise rotation with the angle $\theta_{i_k}=2(\varphi_k+\psi_k)$, and hence the set $\Delta_s^s\cap {\Delta}_{i_k}$ consists of the roots the orthogonal projections of which onto $\h_{i_k}$ belong to the union of the sectors labeled $s^2{\Delta}_{i_k}^{1}$ and ${\Delta}_{i_k}^{2}$ at Figure 3. Together with the description of the set ${\Delta}_{i_k}\cap \Delta_+$ given in (\ref{dik+}) this implies that the number of the roots in the set $\Delta_s\cap {\Delta}_{i_k}=\{\alpha\in {\Delta}_{i_k}\cap \Delta_+: s\alpha \in \Delta_-\}$, where $\Delta_s=\{\alpha\in \Delta_+: s\alpha \in \Delta_-\}$, is equal to the number of roots in the set $\Delta_s^s\cap {\Delta}_{i_k}$. From this observation we deduce that the length $l(s)$ of $s$ with respect to the system of simple roots in $\Delta_+^s$ is the same as the length of $s$ with respect to the system of simple roots in $\Delta_+$, as both of them are equal to the cardinality to the set $\bigcup_{k=0}^M\Delta_s^s\cap {\Delta}_{i_k}$ (disjoint union) which is the same as the cardinality of the set $\bigcup_{k=0}^M\Delta_s\cap {\Delta}_{i_k}$ (disjoint union).

\end{proof}

We shall also need a family of systems of positive roots in $\Delta$ related to the circular normal ordering associated to (\ref{NO}). According to Lemma \ref{s-1} (see, in particular, formula (\ref{auxord10'})) for $j=1,\ldots R$ the minimal segment 
\begin{equation}
\Delta^{j},\Delta^{j+1},\ldots, -\Delta^{j-1}
\end{equation}
of the circular normal ordering associated to (\ref{NO}) is a system of positive roots in $\Delta$ as its length is equal to the number of positive roots and it is closed under addition of roots by Lemma \ref{minsegm}. Denote this system of positive roots by $\Delta_+^j$. \index[not]{D@$\Delta_+^j$} Note that by this definition $\Delta_+^1=\Delta_+^s$.
\begin{lemma}\label{kparab}
Let $\Gamma^k\subset \Delta_+^k$, \index[not]{G@$\Gamma^k$} $k=1,\ldots, R$ be the set of simple roots. Then $\Gamma^k\cap \Delta_0$ is a set of simple roots in $(\Delta_0)_+\subset \Delta_+^k$, and hence the set of roots $\Delta_+^k\cup \Delta_0$ is parabolic. \index{roots!parabolic set of}
\end{lemma}
\begin{proof}
We show first that for any $\alpha,\beta\in \Delta_+^k$ such that $\alpha$ or $\beta$ is not an element of $\Delta_0$ one has $\alpha+\beta\not\in \h_0$.

Indeed if $\alpha,\beta\in \Delta_+^k$ are such that $\alpha$ or $\beta$ does not belong to $\Delta_0$ then $\alpha \in \Delta_{i_p}$, $\beta\in \Delta_{i_q}$ for some $i_p,i_q$ and $p>0$ or $q>0$. If $\alpha+\beta\in \h_0$ then the orthogonal projections of $\alpha+\beta$ onto $\ha_{i_v}$ with any $v>0$ must be equal to zero. The definition of the sets $\Delta_{i_k}$ implies now that $p=q$, $i_p=i_q>0$ and either $\alpha \in \Delta_{i_p}^{(r)}$, $\beta\in -\Delta_{i_p}^{(r)}$ or  $\alpha \in -\Delta_{i_p}^{(r)}$, $\beta\in \Delta_{i_p}^{(r)}$ for some ray $r$. But this is impossible as by the definition of $\Delta_+^k$ and by the definition of the sets $\Delta^j$ in (\ref{auxord10'})  the set $\Delta_+^k$ only contains elements of one of the sets $\Delta_{i_p}^{(r)}$ xor $-\Delta_{i_p}^{(r)}$. Thus $\alpha+\beta\not\in \h_0$. In particular, $\alpha+\beta\not\in \Delta_0\subset \h_0$.

Now let $\Gamma^k=\{\alpha_1,\ldots, \alpha_l\}$ and $\Gamma^k\cap \Delta_0=\{\alpha_1,\ldots, \alpha_u\}$ with $u\leq l$. In order to show that $\Gamma^k\cap \Delta_0$ is a set of simple roots in $(\Delta_0)_+\subset \Delta_+^k$ we shall verify that any element of $(\Delta_0)_+\subset \Delta_+^k$ is a linear combination of the roots from the set $\Gamma^k\cap \Delta_0$ with non-negative integer coefficients.

Let $\delta=\sum_{i=1}^{l}n_i\alpha_i\in \Delta_+^k$ with $n_j>0$ for some $j>u$. If $\delta$ is simple then $\delta=\alpha_j\not \in \Delta_0\cap \Delta_+^k=(\Delta_0)_+$.

If $\delta$ is not simple then $\delta=\alpha+\beta$, where $\alpha,\beta \in \Delta_+^k$ and the heights of both $\alpha$ and $\beta$ are strictly less than the height of $\delta$. Assume that $\alpha$ or $\beta$ is not an element of $\Delta_0$. Then by the first part of the proof $\delta=\alpha+\beta\not\in \h_0$. In particular, $\delta\not \in \Delta_0\cap \Delta_+^k=(\Delta_0)_+\subset \h_0$. 

Now simple induction over the height of $\delta$ implies that any $\delta=\sum_{i=1}^{l}n_i\alpha_i\in \Delta_+^k$ with $n_j>0$ for some $j>u$ is not an element of $(\Delta_0)_+$. Thus any element of $(\Delta_0)_+\subset \Delta_+^k$ is a linear combination of the roots from the set $\Gamma^k\cap \Delta_0$ with non-negative integer coefficients, i.e. $\Gamma^k\cap \Delta_0$ is a set of simple roots in $(\Delta_0)_+$. This completes the proof.

\end{proof}

The following lemma will be used for defining Zhelobenko type operators for q-W--algebras both in the classical and in the quantum case.
\begin{lemma}\label{jsegm}
For $j=1,\ldots, R-1$ the following statements are true.

(i) $s(\Delta^j)\subset \Delta_+^{j+1}$.

(ii) $\Delta^j\cap \mathbb{N}s(\Delta_+^j\cup \Delta_0)=\emptyset$.

(iii) $\Delta^j\cap \mathbb{N}s(\Delta_s^s\cup \Delta_0)=\emptyset$.

(iv) For any $\alpha\in \Delta^j$, $\alpha_0\in \Delta_0$ such that $\alpha+\alpha_0\in \Delta$ one has $\alpha+\alpha_0\in \Delta^j$.
\end{lemma}

\begin{proof}
(i) For $j$ such that $\Delta^j\not\subset \Delta_s^s\cup \Delta_0$ one has by (\ref{auxord11'}) $\Delta^j<s(\Delta^j)\subset \Delta_+^s$, and hence (i) holds in this case by the definition of $\Delta_+^{j+1}$.

If $j$ such that $\Delta^j\subset s^2\Delta_{s^1}^s$ then 
$$
s(\Delta^j)\subset s^1s^2s^2\Delta_{s^1}^s=s^1\Delta_{s^1}^s=-\Delta_{s^1}^s.
$$
Observing that $\Delta_{s^1}^s<s^2\Delta_{s^1}^s$ by Proposition \ref{pord} (vii) we deduce again that (i) holds by the definition of $\Delta_+^{j+1}$.

Finally, if $j=R-1$ then $\Delta^{R-1}=\Delta_{s^2}^s$, and
$$
s(\Delta^j)= s^1s^2\Delta_{s^2}^s=-s^1\Delta_{s^2}^s\subset -(\Delta_-^s\setminus \Delta_0)\subset -(\Delta_-^s\setminus \Delta_0)\cup (\Delta_0)_+=\Delta_+^R.
$$
This completes the proof of part (i).

(ii) First observe that by Lemma \ref{kparab} the set of roots $\Delta_+^j\cup \Delta_0$ is parabolic. Therefore to establish (ii) it suffices to show that $\Delta^j\cap s(\Delta_+^j\cup \Delta_0)=\emptyset$.

Indeed, if $\alpha\in (\Delta_+^j\cup \Delta_0)\cap \Delta_+^s$ is such that $s\alpha \in \Delta_+^s$ then $\alpha \in \Delta^k\not\subset \Delta_s^s$, $k\geq j$, so $s\Delta^k\subset \Delta_+^s$, and by (\ref{auxord11'}) $s\alpha \in s\Delta^k>\Delta^k\geq \Delta^j$ for $k\neq R$, $s\alpha\in \Delta^R>\Delta^j$ for $k=R$, and hence $s\alpha\not\in \Delta^j$ in this case. 

If $\alpha\in \Delta_+^j\cup \Delta_0$ is such that $s\alpha \in \Delta_-^s$ then $s\alpha\not\in \Delta^j$ as $\Delta^j\subset \Delta_+^s$.

If $\alpha\in (\Delta_+^j\cup \Delta_0)\cap \Delta_-^s$ is such that $s\alpha \in \Delta_+^s$ then $\alpha \in  -\Delta_s^s=-(\Delta_{s^2}^s\cup s^2\Delta_{s^1}^s)$. 

Assume that $s\alpha\in \Delta^j$. If $\alpha \in -s^2\Delta_{s^1}^s$ then $s\alpha \in \Delta_{s^1}^s=\Delta^1$, so $j=1$, $\Delta_+^j=\Delta_+^s$, and hence $\alpha\in (\Delta_+^1\cup \Delta_0)\cap \Delta_-^s=(\Delta_+^s\cup \Delta_0)\cap \Delta_-^s=-(\Delta_0)_+$. Thus we arrive at a contradiction as by the assumption $\alpha \in -s^2\Delta_{s^1}^s$, and $s^2\Delta_{s^1}^s\cap (\Delta_0)_+=\emptyset$. 

If $\alpha \in -\Delta_{s^2}^s$ then $\alpha\in \Delta_+^R$ by the definition of $\Delta_+^R$, and $\alpha\not \in \Delta_+^j$ for any $j<R$. Thus we again arrive at a contradiction as by the assumption $\alpha \in \Delta_+^j$ for $j<R$.
This completes the proof. 

(iii) is obvious as $s(\Delta_s^s\cup \Delta_0)\subset \Delta_-^s\cup \Delta_0$, the set of roots $\Delta_-^s\cup \Delta_0$ is parabolic, and $\Delta^j\subset \Delta_+^s\setminus \Delta_0$.

(iv)  The statement of part (iv) for $j=1,R-1$ follows from Proposition \ref{pord} (ii).

For $j\neq 1,R-1$ one has $\Delta^j=\Delta_{i_k}^{(r)}$ for some $k$ and $r$ with $k>0$. So the same property follows from the definitions of the sets $\Delta_{i_k}^{(r)}$.

\end{proof}

The relative positions of the systems of positive rots $\Delta_+^s$, $\Delta_+$ and of  the minimal segments introduced in Proposition \ref{pord} are shown at the following picture where all the segments are placed on a circle according to the circular normal ordering of roots corresponding to normal ordering (\ref{NO}) of $\Delta_+^s$.
\begin{center}
\begin{tikzpicture}
\draw [lightgray, ultra thin, dashed] (0,0) circle [radius=6cm];
\draw [olive, thin] (xyz polar cs:angle=0, radius=6.2) arc [start angle=0, end angle=180, radius=6.2];
\node [above, olive] at (xyz polar cs:angle=90, radius=6.2) {$\Delta_+^s$};
\draw [teal, thin] (xyz polar cs:angle=168, radius=5.9) arc [start angle=168, end angle=-11, radius=5.9];
\draw [brown, ultra thick, dotted] (xyz polar cs:angle=168, radius=5.83) arc [start angle=167, end angle=22, radius=5.83];
\node [teal] at (xyz polar cs:angle=80, radius=5.4) {$\Delta_+$};
\node [red] at (xyz polar cs:angle=165, radius=5.4) {$\Delta_{s^1}^s$};
\node [blue] at (xyz polar cs:angle=165, radius=6.7) {$\Delta_{s^1}^{-1}$};
\node [red] at (xyz polar cs:angle=-15, radius=5.4) {$-\Delta_{s^1}^s$};
\node [blue] at (xyz polar cs:angle=-15, radius=6.7) {$-\Delta_{s^1}^{-1}$};
\node [green] at (xyz polar cs:angle=25, radius=5.4) {$\Delta_{s^2}^s$};
\node [cyan] at (xyz polar cs:angle=25, radius=6.7) {$\Delta_{s^2}^{-1}$};
\node [green] at (xyz polar cs:angle=-155, radius=5.4) {$-\Delta_{s^2}^s$};
\node [cyan] at (xyz polar cs:angle=-155, radius=6.7) {$-\Delta_{s^2}^{-1}$};
\node [violet] at (xyz polar cs:angle=5, radius=6.7) {$(\Delta_0)_+$};
\node [violet] at (xyz polar cs:angle=-175, radius=6.7) {$-(\Delta_0)_+$};
\node [brown] at (xyz polar cs:angle=100, radius=5.4) {$\bf \Delta_{\bf \m_+}$};
\draw [red, ultra thick] (-6,0) arc [start angle=180, end angle=150, radius=6];
\draw [blue, ultra thick] (xyz polar cs:angle=175, radius=6.1) arc [start angle=175, end angle=155, radius=6.1];
\draw [green, ultra thick] (xyz polar cs:angle=40, radius=6) arc [start angle=40, end angle=10, radius=6];
\draw [cyan, ultra thick] (xyz polar cs:angle=35, radius=6.1) arc [start angle=35, end angle=15, radius=6.1];
\draw [violet, ultra thick] (xyz polar cs:angle=10, radius=6) arc [start angle=10, end angle=0, radius=6];
\draw [red, ultra thick] (6,0) arc [start angle=0, end angle=-30, radius=6];
\draw [blue, ultra thick] (xyz polar cs:angle=-5, radius=6.1) arc [start angle=-5, end angle=-25, radius=6.1];
\draw [green, ultra thick] (xyz polar cs:angle=-140, radius=6) arc [start angle=-140, end angle=-170, radius=6];
\draw [cyan, ultra thick] (xyz polar cs:angle=-145, radius=6.1) arc [start angle=-145, end angle=-165, radius=6.1];
\draw [violet, ultra thick] (xyz polar cs:angle=-170, radius=6) arc [start angle=-170, end angle=-180, radius=6];
\draw[fill] (xyz polar cs:angle=168, radius=6) circle [radius=2pt];
\node at (xyz polar cs:angle=170, radius=5.7) {$\bf \gamma_1$};
\draw[fill] (xyz polar cs:angle=22, radius=6) circle [radius=2pt];
\node at (xyz polar cs:angle=20, radius=5.6) {$\bf \gamma_{l'}$};
\draw[fill] (xyz polar cs:angle=155, radius=6) circle [radius=2pt];
\node at (xyz polar cs:angle=155, radius=5.5) {$\bf \gamma_{\widetilde{l}}$};
\draw[fill] (xyz polar cs:angle=35, radius=6) circle [radius=2pt];
\node at (xyz polar cs:angle=35, radius=5.3) {$\bf \gamma_{\widetilde{l}+1}$};
\draw[fill] (xyz polar cs:angle=-12, radius=6) circle [radius=2pt];
\node at (xyz polar cs:angle=-11, radius=6.5) {$\bf -\gamma_1$};
\end{tikzpicture}

Fig. 4
\end{center}

The reader may find this picture useful in combination with Lemma \ref{minsegm} when adding roots or commuting roots vectors. This picture can be also useful for deriving some formulas containing q-commutators of quantum root vectors as explained in the next chapter.


\section{Algebraic group analogues of Slodowy slices}\label{slodowy}

\setcounter{equation}{0}
\setcounter{theorem}{0}

In this section we define analogues of the Slodowy slices \index{algebraic group analogues of Slodowy slices} for algebraic groups. As this construction is important for semisimple algebraic groups over arbitrary algebraically closed fields, especially for the study of the Lusztig partition in Sections \ref{luspart} and \ref{stt}, we work over an arbitrary algebraically closed field $\bf k$ in this and in the next two sections. Applications to quantum groups only require ${\bf k}=\mathbb{C}$.

Let $s\in W=W(G_{\bf k}, H_{\bf k})$ be a Weyl group element, $\Delta_+^s$ a system of positive roots defined in (\ref{D+sdef}) for $s$, $\Gamma^s$ the set of simple roots in $\Delta_+^s$.

Denote by $P_{\bf k}^s$ \index[not]{P@$P_{\bf k}^s$} the parabolic subgroup \index{subgroup!parabolic} of $G_{\bf k}$ containing the Borel subgroup $B_{{\bf k},-}^s$ \index[not]{B@$B_{{\bf k},-}^s$} corresponding to $-\Delta_+^s$ and associated to the subset $-\Gamma_0^s$ of the set of simple roots in $-\Gamma^s$, where $\Gamma_0^s=\Gamma^s\cap {\Delta}_{0}$. \index[not]{G@$\Gamma_0^s$}

Let $N_{\bf k}^s$ \index[not]{N@$N_{\bf k}^s$} and $L_{\bf k}^s$ \index[not]{L@$L_{\bf k}^s$} the unipotent radical \index{radical!unipotent} and the Levi factor \index{subgroup!parabolic!Levi factor of} of $P_{\bf k}^s$, respectively, and $\overline{N}_{\bf k}^s$ \index[not]{N@$\overline{N}_{\bf k}^s$} the opposite unipotent radical.  

Note that we have a natural inclusion $P_{\bf k}^s\supset N_{\bf k}^s$, and by Lemma \ref{parab} ${\Delta}_{0}$ is the root system of the reductive algebraic group $L_{\bf k}^s$, while, by the definition of $P_{\bf k}^s$, its unipotent radical $N_{\bf k}^s$ is generated by the one-parameter subgroups corresponding to the roots from the set $(-\Delta_+^s)\setminus \Delta_0$.  

As in Section \ref{notation}, denote a representative for the Weyl group element $s$ in $N_{G_{\bf k}}(H_{\bf k})$ by $\dot{s}$.
Let $Z_{\bf k}^s$ \index[not]{Z@$Z_{\bf k}^s$} be the connected subgroup of $G_{\bf k}$ generated by the semisimple part of the standard Levi subgroup \index{semisimple part!of a Levi subgroup} $L_{\bf k}^s$ \index{subgroup!Levi standard} and by the identity component \index{component!identity of an algebraic group} $H^0_{\bf k}$ \index[not]{H@$H^0_{\bf k}$} of centralizer of $\dot{s}$ in $H_{\bf k}$. 

Let $N_{{\bf k},s}=\{ v \in N_{\bf k}^s|\dot{s}v\dot{s}^{-1}\in \overline{N}_{\bf k}^s \}$. \index[not]{N@$N_{{\bf k},s}$} Observe that $N_{{\bf k},s}\subset N_{\bf k}^s$ is the algebraic subgroup generated by the one--parameter subgroups \index{subgroup!one--parameter corresponding to a root} corresponding to the roots from the set $-\Delta_s^s$, where $\Delta_s^s=\{\alpha \in  \Delta_+^s: s\alpha\in -\Delta_+^s\}$, and the cardinality of the set $-\Delta_s^s$ is equal to $l(s)$, where $l(s)$ is the length of the Weyl group element $s\in W$ with respect to the system of simple roots in $\Delta_+^s$(see e.g. \cite{Car}, \S 2.2, 8.4). Therefore ${\rm dim}~N_{{\bf k},s}=l(s)$.

Denote by $\p_{\bf k}^s$, \index[not]{p@$\p_{\bf k}^s$} $\n_{\bf k}^s$, \index[not]{n@$\n_{\bf k}^s$} $\opn_{\bf k}^s$, \index[not]{n@$\opn_{\bf k}^s$} $\l_{\bf k}^s$, \index[not]{l@$\l_{\bf k}^s$} and $\z_{\bf k}^s$ \index[not]{z@$\z_{\bf k}^s$} the Lie subalgebras \index{Lie!subalgebra} of $\g_{\bf k}$ corresponding to $P_{\bf k}^s$, $N_{\bf k}^s$, $\overline{N}_{\bf k}^s$, $L_{\bf k}^s$, and $Z_{\bf k}^s$, respectively. Note that the Lie subalgebra corresponding to $H^0_{\bf k}$ is the fixed point subspace $\h^s_{\bf k}$ \index[not]{h@$\h^s_{\bf k}$} for the action of $s$ on $\h_{\bf k}$.

As in Section \ref{notation}, when ${\bf k}=\mathbb{C}$ we drop the subscript $\bf k$ in the symbols above, so that $P_{\mathbb{C}}^s=P^s$, \index[not]{P@$P^s$} $\p_{\mathbb{C}}^s=\p^s$, \index[not]{p@$\p^s$} $N_{\mathbb{C},s}=N_s$, \index[not]{N@$N_s$} etc. Recall that in this case, according to our convention, we also write $\dot{s}=s$ when ${\bf k}=\mathbb{C}$. 

In the proofs below we shall frequently use the following two lemmas. The first one is a direct consequence of Proposition 8.1.1, Corollary 8.1.2 and Lemma 8.2.2 in \cite{Sp}.
\begin{lemma}\label{decNr}
Let $\Psi, \Psi_i\subset \Delta$, $i=1,\ldots ,m$ be additively closed subsets of roots such that
$$
\Psi=\bigcup_{i=1}^m\Psi_i
$$
is a disjoint union and $\Psi$ does not contain opposite roots. For any additively closed subset $\Xi\subset \Delta$ \index{roots!additively closed subset of} which does not contain opposite roots, denote by $N_\Xi\subset G_{\bf k}$ \index[not]{N@$N_\Xi$} the algebraic subgroup generated by the one--parameter subgroups \index{subgroup!algebraic, generated by one--parameter subgroups} of $G_{\bf k}$ corresponding to the roots from $\Xi$. Then multiplication in $G_{\bf k}$ yields an isomorphism of algebraic varieties \index{variety!algebraic} \index{varieties!isomorphism of}
$$
N_\Psi= N_{\Psi_1}\ldots N_{\Psi_m}\simeq N_{\Psi_1}\times\ldots \times N_{\Psi_m}.
$$
\end{lemma}

Part (i) of the second lemma below is Corollary in Section 7.4 of \cite{Hu}, part (ii) is Lemma 2.3.3 and 1.6.10(4) in \cite{Sp}, and part (iii) is Theorem 5.3.2(iii) in \cite{Sp}.
\begin{lemma}\label{vart}
(i) Let $G_{\bf k}^1$ and $G_{\bf k}^2$ be algebraic subgroups in an algebraic group $A_{\bf k}$ over an algebraically closed field $\bf k$, and $G_{\bf k}^1$ normalizes $G_{\bf k}^2$. Then $G_{\bf k}^1G_{\bf k}^2$ is an algebraic subgroup of $A_{\bf k}$ and multiplication in $A_{\bf k}$ defines an isomorphism of varieties, $G_{\bf k}^1\times G_{\bf k}^2\to G_{\bf k}^1G_{\bf k}^2= G_{\bf k}^2G_{\bf k}^1$.

(ii) Let $A_{\bf k}$ be an algebraic group over an algebraically closed field $\bf k$. Then any orbit $\mathcal{O}$ \index[not]{O@$\mathcal{O}$} \index{action!of an algebraic group!orbit} for a regular algebraic group action \index{action!of an algebraic group!regular} of $A_{\bf k}$ on a variety $V$ over $\bf k$ is open in its closure. Thus $\mathcal{O}$ has the natural structure of an algebraic variety. 

(iii) Let $\phi:X\to Y$ be an equivariant morphism of homogeneous spaces \index{space!homogeneous} for an algebraic group $A_{\bf k}$ over an algebraically closed field $\bf k$. Then $\phi$ is an isomorphism if and only if it is bijective and induces a bijection of the tangent spaces at the some points $x\in X$ and $y=\phi(x)\in Y$. 
\end{lemma}

To formulate our main statement in this section we need the following lemma. 
\begin{lemma}\label{NZsN}
Let $N_{\bf k}^sZ_{\bf k}^s\dot{s}N_{\bf k}^s\subset G_{\bf k}$ be the image in $G_{\bf k}$ of $N_{\bf k}^s\times Z_{\bf k}^s\times N_{\bf k}^s$ under the map 
\begin{equation}\label{defNZsN}
N_{\bf k}^s\times Z_{\bf k}^s\times N_{\bf k}^s \to G_{\bf k}, (n,z,n')\mapsto nz\dot{s}n'
\end{equation}
induced by the group multiplication in $G_{\bf k}$. Then 

(i) $Z_{\bf k}^sN_{\bf k}^s$ is an algebraic subgroup of $G_{\bf k}$ and multiplication defines an isomorphism of varieties, $Z_{\bf k}^s\times N_{\bf k}^s\to Z_{\bf k}^sN_{\bf k}^s= N_{\bf k}^sZ_{\bf k}^s$;

(ii) $N_{\bf k}^sZ_{\bf k}^s\dot{s}N_{\bf k}^s$ is an algebraic subvariety \index{subvariety} of $G_{\bf k}$, $N_{\bf k}^sZ_{\bf k}^s\dot{s}N_{\bf k}^s=N_{\bf k}^sZ_{\bf k}^s\dot{s}N_{{\bf k},s}$, and the map 
\begin{equation}\label{descrNZsN}
N_{\bf k}^s\times Z_{\bf k}^s\times N_{{\bf k},s} \to N_{\bf k}^sZ_{\bf k}^s\dot{s}N_{{\bf k},s}=N_{\bf k}^sZ_{\bf k}^s\dot{s}N_{\bf k}^s, (n,z,n')\mapsto nz\dot{s}n'
\end{equation}
induced by the group multiplication in $G_{\bf k}$ is an isomorphism of algebraic varieties;

(iii) Moreover, 
\begin{equation}\label{isodescr}
N_{\bf k}^sZ_{\bf k}^s\dot{s}N_{\bf k}^s= N_{\bf k}^sZ_{\bf k}^s\dot{s}N_{{\bf k},s}= N_{\bf k}^s\dot{s}Z_{\bf k}^sN_{{\bf k},s}.
\end{equation}
\end{lemma}

\begin{proof}
The proof of part (i) follows from Lemma \ref{vart} (i) with $A_{\bf k}=G_{\bf k}$, $G_{\bf k}^1=Z_{\bf k}^s$ and $G_{\bf k}^2=N_{\bf k}^s$. 

(ii) The proof is parallel to that of Lemma 8.3.6 in \cite{Sp} where a similar statement is justified in the case of Bruhat cells. 

Using the isomorphism of part (i) and definition (\ref{defNZsN}) one can describe the set $N_{\bf k}^sZ_{\bf k}^s\dot{s}N_{\bf k}^s$ as the orbit of the element $\dot{s}$ in $G_{\bf k}$ for the following regular action of $N_{\bf k}^sZ_{\bf k}^s\times N_{\bf k}^s$ on $G_{\bf k}$:
$$
(N_{\bf k}^sZ_{\bf k}^s\times N_{\bf k}^s)\times G_{\bf k} \to G_{\bf k}, ((nz,n'), g)\mapsto nzg{n'}^{-1}, n,n'\in N_{\bf k}^s, z\in Z_{\bf k}^s, g\in G_{\bf k}.
$$
Now the first claim in part (ii) follows from Lemma \ref{vart} (ii) with $A_{\bf k}=N_{\bf k}^sZ_{\bf k}^s\times N_{\bf k}^s$, $V=G_{\bf k}$, $\mathcal{O}=N_{\bf k}^sZ_{\bf k}^s\dot{s}N_{\bf k}^s$.

To establish isomorphism (\ref{descrNZsN}) we introduce the algebraic subgroup $N_{{\bf k},s}':=N_{\bf k}^s\cap \dot{s}^{-1}N_{\bf k}^s\dot{s}\subset N_{\bf k}^s$. \index[not]{N@$N_{{\bf k},s}'$} Applying Lemma \ref{decNr} for $\Psi_1=-\Delta_s^s$ and $\Psi_2=-(\Delta_+^s\setminus (\Delta_s^s\cup \Delta_0))$, so that $N_{\Psi_1}=N_{{\bf k},s}$ and $N_{\Psi_2}=N_{{\bf k},s}'$, one has an isomorphism of varieties 
\begin{equation}\label{sdecN}
N_{\bf k}^s= N_{{\bf k},s}'N_{{\bf k},s}\simeq N_{{\bf k},s}'\times N_{{\bf k},s},
\end{equation} 
induced by the group multiplication in $G_{\bf k}$. Hence
\begin{equation}\label{NsZN}
N_{\bf k}^sZ_{\bf k}^s\dot{s}N_{\bf k}^s=N_{\bf k}^sZ_{\bf k}^s\dot{s}N_{{\bf k},s}'N_{{\bf k},s}= N_{\bf k}^sZ_{\bf k}^s\dot{s}N_{{\bf k},s},
\end{equation}
as $Z_{\bf k}^s$ normalizes $N_{\bf k}^s$. This implies that morphism of varieties (\ref{descrNZsN}) is bijective. \index{varieties!morphism of}
So in fact one can view $N_{\bf k}^sZ_{\bf k}^s\dot{s}N_{\bf k}^s$ as the orbit of the element $\dot{s}$ in $G_{\bf k}$ for the following regular action of the algebraic group $N_{\bf k}^sZ_{\bf k}^s\times N_{{\bf k},s}$ on $G_{\bf k}$:
$$
(N_{\bf k}^sZ_{\bf k}^s\times N_{{\bf k},s})\times G_{\bf k} \to G_{\bf k}, ((nz,n'), g)\mapsto nzg{n'}^{-1}, n\in N_{\bf k}^s, n'\in N_{{\bf k},s}, z\in Z_{\bf k}^s, g\in G_{\bf k},
$$
and the action map gives rise to a bijective equivariant morphism \index{varieties!equivariant morphism of} 
\begin{equation}\label{morph}
N_{\bf k}^sZ_{\bf k}^s\times N_{{\bf k},s}\to N_{\bf k}^sZ_{\bf k}^s\dot{s}N_{{\bf k},s}, (nz,n')\mapsto nz\dot{s}{n'}^{-1}
\end{equation}
of $N_{\bf k}^sZ_{\bf k}^s\times N_{{\bf k},s}$ and $N_{\bf k}^sZ_{\bf k}^s\dot{s}N_{{\bf k},s}$ viewed as homogeneous spaces for $N_{\bf k}^sZ_{\bf k}^s\times N_{{\bf k},s}$.

One verifies straightforwardly that this morphism induces an isomorphism of the tangent spaces at the points $(1,1)$ and $\dot{s}$. Therefore by Lemma \ref{vart} (iii) with $X=N_{\bf k}^sZ_{\bf k}^s\times N_{{\bf k},s}$, $Y=N_{\bf k}^sZ_{\bf k}^s\dot{s}N_{{\bf k},s}$, $A_{\bf k}=N_{\bf k}^sZ_{\bf k}^s\times N_{{\bf k},s}$ morphism (\ref{morph}) is an isomorphism of varieties. Composing it with the algebraic group isomorphism $N_{\bf k}^sZ_{\bf k}^s\simeq N_{\bf k}^s\times Z_{\bf k}^s$ established in part (i) and with the map taking inverse on $N_{{\bf k},s}$ we deduce that (\ref{descrNZsN}) is an isomorphism of varieties. 

(iii) The first identity in (\ref{isodescr}) was established in part (ii), the second one follows from the fact that $s$ fixes the root system $\Delta_0$ of $Z_{\bf k}^s$, and hence $\dot{s}Z_{\bf k}^s\dot{s}^{-1}=Z_{\bf k}^s$.
\end{proof}



Now we can state the main proposition in this section in which we define transversal slices to conjugacy classes in algebraic groups.
\begin{proposition}\label{crosssect}
(i) $Z_{\bf k}^sN_{{\bf k},s}$ is an algebraic subgroup in $G_{\bf k}$, with the variety structure inherited from $Z_{\bf k}^s\times N_{{\bf k},s}$ via the bijection $Z_{\bf k}^sN_{{\bf k}, s}\simeq Z_{\bf k}^s\times N_{{\bf k},s}$ induced by the multiplication map in $G_{\bf k}$. The set $\dot{s}Z_{\bf k}^sN_{{\bf k}, s}= Z_{\bf k}^s\dot{s}N_{{\bf k}, s}$ with the variety structure induced from $Z_{\bf k}^sN_{{\bf k},s}$ is a closed subvariety \index{subvariety!closed} of $N_{\bf k}^sZ_{\bf k}^s\dot{s}N_{\bf k}^s$ and of $G_{\bf k}$, and the conjugation map
\begin{equation}\label{cross}
N_{\bf k}^s\times \dot{s}Z_{\bf k}^sN_{{\bf k}, s}\rightarrow N_{\bf k}^sZ_{\bf k}^s\dot{s}N_{\bf k}^s, (n,g)\mapsto ngn^{-1},
\end{equation}
is an isomorphism of varieties;

(ii) The variety $\Sigma_{{\bf k},s}=\dot{s}Z_{\bf k}^sN_{{\bf k},s}=Z_{\bf k}^s\dot{s}N_{{\bf k}, s}$ \index[not]{S@$\Sigma_{{\bf k},s}$} is a transversal slice to the set of conjugacy classes in $G_{\bf k}$; \index{conjugacy class!in an algebraic group}

(iii) Assume that ${\bf k}=\mathbb{C}$ and that $G_{\bf k}=G$ is simply connected. \index{group!algebraic!simply connected} Then $N^ssZ^sN^s$ is a closed subvariety of $G$. \index[not]{N@$N^s$} \index[not]{Z@$Z^s$} 
\end{proposition}

\begin{proof}

(i) Note first that $Z_{\bf k}^sN_{{\bf k}, s}$ is an algebraic subgroup in $G_{\bf k}$. Indeed, by the definition $s$ fixes the root system $\Delta_0$ of $Z_{\bf k}^s$, so $\dot{s}Z_{\bf k}^s\dot{s}^{-1}=Z_{\bf k}^s$, and also $Z_{\bf k}^s$ normalizes $N_{\bf k}^s$ and $\overline{N}_{\bf k}^s$. Hence for any $n_s\in N_{{\bf k},s}$ and any $z\in Z_{\bf k}^s$ one has $zn_sz^{-1}\in N_{\bf k}^s$ and $\dot{s}zn_sz^{-1}\dot{s}^{-1}=(\dot{s}z\dot{s}^{-1})(\dot{s}n_s\dot{s}^{-1})(\dot{s}z\dot{s}^{-1})^{-1}\in \overline{N}_{\bf k}^s$ as $\dot{s}z\dot{s}^{-1}\in Z_{\bf k}^s$ and $\dot{s}n_s\dot{s}^{-1}\in \overline{N}_{\bf k}^s$ by the definition of $N_{{\bf k}, s}$. We deduce that $Z_{\bf k}^s$ normalizes $N_{{\bf k}, s}$ in $G_{\bf k}$. Thus by Lemma \ref{vart} (i) with $A_{\bf k}=G_{\bf k}$, $G_{\bf k}^1=Z_{\bf k}^s$ and $G_{\bf k}^2=N_{{\bf k},s}$, $Z_{\bf k}^sN_{{\bf k},s}$ is an algebraic subgroup in $G_{\bf k}$, with the variety structure inherited from $Z_{\bf k}^s\times N_{{\bf k},s}$ via the bijection $Z_{\bf k}^sN_{{\bf k}, s}\simeq Z_{\bf k}^s\times N_{{\bf k},s}$ induced by the multiplication map in $G_{\bf k}$. 

$\dot{s}Z_{\bf k}^sN_{{\bf k},s}\simeq Z_{\bf k}^sN_{{\bf k}, s}$ is a closed subvariety of $G_{\bf k}$ as the left multiplication by $\dot{s}$ in $G_{\bf k}$ is a morphism of varieties.


$\dot{s}Z_{\bf k}^sN_{{\bf k}, s}\simeq Z_{\bf k}^s\dot{s}N_{{\bf k}, s}$ is a closed subvariety of $N_{\bf k}^sZ_{\bf k}^s\dot{s}N_{\bf k}^s$ due to the closed embedding
$$
Z_{\bf k}^s\times N_{{\bf k}, s}\hookrightarrow N_{\bf k}^s\times Z_{\bf k}^s\times N_{{\bf k},s}, (z,n)\mapsto (1,z,n)
$$
and isomorphism of varieties (\ref{descrNZsN}). 

We show that map (\ref{cross}) is an isomorphism of varieties. By the definition this map is a morphism of varieties. We shall define the inverse map and show that it is also a morphism of varieties.

Observe that map (\ref{cross}) is bijective if and only if for any given $k_s\in N_{{\bf k},s}, u\in N_{\bf k}^s$ and $z\in Z_{\bf k}^s$ the equation
\begin{equation}\label{tpr'}
uz\dot{s}k_s=nz'\dot{s}n_sn^{-1}
\end{equation}
has a unique solution $n\in N_{\bf k}^s,n_s\in N_{{\bf k},s},z'\in Z_{\bf k}^s$. In this case the solution defines the inverse map to morphism of varieties (\ref{cross}). We shall construct this map explicitly. From this construction it will be clear that the inverse map is also a morphism of varieties.

First observe that any element $uz\dot{s}k_s$ is uniquely conjugated by $k_s\in N_{{\bf k},s}$ to $vz\dot{s}\in N_{\bf k}^sZ_{\bf k}^s\dot{s}, v=k_su$, and hence we can assume that $k_s=1$ in (\ref{tpr'}),
\begin{equation}\label{tpr}
vz\dot{s}=nz'\dot{s}n_sn^{-1}.
\end{equation}
Using isomorphism (\ref{descrNZsN}) we deduce that the corresponding map
\begin{equation}\label{map1}
N_{\bf k}^sZ_{\bf k}^s\dot{s}N_{{\bf k},s}\to N_{{\bf k},s}\times N_{\bf k}^sZ_{\bf k}^s\dot{s}, uz\dot{s}k_s\mapsto (k_s,k_suz\dot{s})
\end{equation}
is an isomorphism of varieties.

Now we show that for any given $v\in N_{\bf k}^s$ and $z\in Z_{\bf k}^s$ equation (\ref{tpr}) has a unique solution $n\in N_{\bf k}^s,n_s\in N_{{\bf k},s},z'\in Z_{\bf k}^s$ which is expressed in terms of $vz\in N_{\bf k}^sZ_{\bf k}^s$ using only a composition of morphisms of algebraic varieties: the isomorphism of varieties $N_{\bf k}^sZ_{\bf k}^s\simeq N_{\bf k}^s\times Z_{\bf k}^s$, algebraic factorization of elements of $N_{\bf k}^s$ as products of elements from one--parameter subgroups, Chevalley commutation relations  \index{commutation relations!Chevalley} between one--parameter subgroups of $G_{\bf k}$, conjugation by elements of $Z_{\bf k}^s$ and by the element $\dot{s}\in G_{\bf k}$. This implies that the corresponding map $N_{\bf k}^sZ_{\bf k}^s\dot{s}\to N_{\bf k}^s\times \dot{s}Z_{\bf k}^sN_{{\bf k},s}$, $vz\dot{s}\mapsto (n,z'\dot{s}n_s)$ is an injective morphism of varieties, and hence the composition of this map with isomorphism (\ref{map1}) is the inverse to (\ref{cross}) morphism of varieties. 



First we rewrite equation (\ref{tpr}) in a slightly different form,
\begin{equation}\label{tpr1}
vz\dot{s}n\dot{s}^{-1}=nz'\dot{s}n_s\dot{s}^{-1}.
\end{equation}

Observe that both the left hand side and the right hand side of equation (\ref{tpr1}) belong to the subvariety $N_{\bf k}^sL_{\bf k}^s\overline{N}_{\bf k}^s\subset G_{\bf k}$, and that the $L_{\bf k}^s$--component of equation (\ref{tpr1}) in $N_{\bf k}^sL_{\bf k}^s\overline{N}_{\bf k}^s$ is reduced to
\begin{equation}\label{eqz}
z=z'.
\end{equation}

Indeed, using (\ref{sdecN}) we can write
\begin{equation}\label{sN}
\dot{s}N_{\bf k}^s\dot{s}^{-1}=\dot{s}N_{{\bf k},s}'\dot{s}^{-1}\dot{s}N_{{\bf k},s}\dot{s}^{-1}\subset N_{\bf k}^s\overline{N}_{\bf k}^s.
\end{equation}
If $n=mm_s$ is the factorization of $n$ corresponding to the factorization $N_{\bf k}^s=N'_{{\bf k},s}N_{{\bf k},s}$  then recalling that $Z_{\bf k}^s$ normalizes both $N_{\bf k}^s$ and $\overline{N}_{\bf k}^s$ we deduce that the factorizations of the right hand side and of the left hand side of equation (\ref{tpr1}) corresponding to the unique factorization $N_{\bf k}^sL_{\bf k}^s\overline{N}_{\bf k}^s$ take the form
\begin{equation}\label{tpr4}
(vz\dot{s}m\dot{s}^{-1}z^{-1})z(\dot{s}m_s\dot{s}^{-1})=nz'(\dot{s}n_s\dot{s}^{-1}),
\end{equation}
where $vz\dot{s}m\dot{s}^{-1}z^{-1},n\in N_{\bf k}^s$ and $\dot{s}m_s\dot{s}^{-1}, \dot{s}n_s\dot{s}^{-1}\in \overline{N}_{\bf k}^s$, $z,z'\in Z_{\bf k}^s\subset L_{\bf k}^s$. This implies (\ref{eqz}).

Note also that the map $N_{\bf k}^sZ_{\bf k}^s\dot{s}\to Z_{\bf k}^s$, $vz\dot{s}\mapsto z$ is a morphism of varieties. 

Now equation (\ref{tpr4}) can be rewritten as follows
\begin{equation}\label{tpr5}
vz\dot{s}m\dot{s}^{-1}\dot{s}m_s\dot{s}^{-1}=
nz\dot{s}n_s\dot{s}^{-1},
\end{equation}
and the $\overline{N}_{\bf k}^s$--component of the last equation in $N_{\bf k}^sL_{\bf k}^s\overline{N}_{\bf k}^s$ is
$$
\dot{s}m_s\dot{s}^{-1}=\dot{s}n_s\dot{s}^{-1}.
$$
From this relation we obtain that
\begin{equation}\label{nss}
m_s=n_s,
\end{equation}
and hence (\ref{tpr5}) yields
\begin{equation}\label{tpr6}
vz\dot{s}m\dot{s}^{-1}z^{-1}=
n.
\end{equation}
Now we show that the last equation defines $n$ in a unique way.

First observe that by (\ref{auxord10'}) we have a disjoint union of minimal segments
\begin{equation}\label{dik-dec}
-(\Delta_+^s\setminus (\Delta_0)_+)=\bigcup_{p=1}^{R-1} (-\Delta^p),
\end{equation}
and by (\ref{auxord1'}) 
\begin{equation}\label{dik-dec1}
-(\Delta_+^s\setminus (\Delta_s^s\cup(\Delta_0)_+))=\bigcup_{p=1}^{q} (-\Delta^p)
\end{equation}
for some $0<q\leq R-1$.

By Lemma \ref{decNr} applied to disjoint unions (\ref{dik-dec}), (\ref{dik-dec1}) we have isomorphisms of varieties
\begin{equation}\label{decnk}
N_{\bf k}^s= N_{-\Delta^1}N_{-\Delta^2}\ldots N_{-\Delta^{R-1}}\simeq N_{-\Delta^1}\times N_{-\Delta^2}\times \ldots \times N_{-\Delta^{R-1}}, 
\end{equation}
$$
N'_{{\bf k},s}=N_{-\Delta^1}N_{-\Delta^2}\ldots N_{-\Delta^q}\simeq N_{-\Delta^1}\times N_{-\Delta^2}\times \ldots \times N_{-\Delta^q}.
$$
Let
\begin{equation}\label{decnu}
n=n_1\ldots n_{R-1},~m=n_1\ldots n_q,~v=v_1 v_2\ldots v_{R-1}
\end{equation}
be the corresponding factorizations of the elements $n\in N_{\bf k}^s$, $m\in N'_{{\bf k},s}$, and $v\in N_{\bf k}^s$, respectively.

We claim that the factors $n_p$, $p=1,\ldots, q$ can be uniquely found by induction starting with $n_1$. Indeed, substituting factorizations (\ref{decnu}) into (\ref{tpr6}) we obtain
$$
v_1v_2\ldots v_{R-1}z\dot{s}n_1\ldots n_{q}\dot{s}^{-1}z^{-1}=n_1\ldots n_{R-1}.
$$
Now equating the $N_{-\Delta^p}$--components of the left hand side and of the right hand side of the last equation, with respect to the first factorization in (\ref{decnk}), and using the fact that by Lemma \ref{jsegm} (i) for $p=1,\ldots, q$ one has $\dot{s}N_{-\Delta^p}\dot{s}^{-1}\subset N_{-\Delta^{p+1}}\ldots N_{-\Delta^{R-1}}$, and by by Lemma \ref{jsegm} (iv) $Z_{\bf k}^s$ normalizes $N_{-\Delta^p}$ for $p=1,\ldots, R-1$, we obtain
\begin{equation}\label{indf1}
n_1=v_1, n_p=(v_p\ldots v_{R-1}z\dot{s}n_1\ldots n_{p-1}\dot{s}^{-1}z^{-1})_p, p=2,\ldots ,q,
\end{equation}
\begin{equation}\label{indf2}
n_s=m_s=(vz\dot{s}m\dot{s}^{-1}z^{-1})_s,
\end{equation}
where the subscript $(\ldots)_p$ stands for the $N_{-\Delta^p}$--component with respect to the first factorization in (\ref{decnk}), and $(vz\dot{s}m\dot{s}^{-1}z^{-1})_s$ is the $N_{{\bf k},s}$-component of $vz\dot{s}m\dot{s}^{-1}z^{-1}\in N_{\bf k}^s$ with respect to the factorization $N_{\bf k}^s=N'_{{\bf k},s}N_{{\bf k},s}$. From formulas (\ref{indf1}) one can recursively find the components $n_p$, $p=1,\ldots ,q$ starting from $n_1=v_1$, and finally one can find $n_s=m_s$ using (\ref{indf2}) with $m=n_1\ldots n_q$. 

Note that by the construction the maps $N_{\bf k}^sZ_{\bf k}^s\dot{s}\to N_{\bf k}^s$, $vz\dot{s}\mapsto n=mm_s$ and $N_{\bf k}^sZ_{\bf k}^s\dot{s}\to N_{{\bf k},s}$, $vz\dot{s}\mapsto n_s$ are injective morphisms of varieties. 

Finally we conclude that
the map $N_{\bf k}^sZ_{\bf k}^s\dot{s}\to N_{\bf k}^s\times \dot{s}Z_{\bf k}^sN_{{\bf k},s}$, $vz\dot{s}\mapsto (n,z\dot{s}n_s)$ is an injective morphism of varieties. Thus its composition with isomorphism of varieties (\ref{map1}) is the injective morphism inverse to morphism (\ref{cross}). This establishes isomorphism (\ref{cross}) and completes the proof of part (i).

(ii) Next we have to show that the variety $\dot{s}Z_{\bf k}^sN_{{\bf k}, s}\subset G_{\bf k}$ is a transversal slice to the set of conjugacy classes in $G_{\bf k}$, i.e. that the differential of the conjugation map
\begin{equation}\label{map}
\gamma:G_{\bf k}\times \dot{s}Z_{\bf k}^sN_{{\bf k}, s}\rightarrow G_{\bf k}
\end{equation}
is surjective.

Note that the set of smooth points of map (\ref{map}) is stable
under the $G_{\bf k}$--action by left translations \index{action!of an algebraic group!by left translations} on the first factor of
$G_{\bf k}\times \dot{s}Z_{\bf k}^sN_{{\bf k}, s}$. Therefore it suffices to show that the
differential of map (\ref{map}) is surjective at points
$(1,szn_s)$, $n_s\in N_{{\bf k}, s}$, $z\in Z_{\bf k}^s$.

In terms of the left trivialization of the tangent bundle $TG_{\bf k}$
and of the induced trivialization of $T(\dot{s}Z_{\bf k}^sN_{{\bf k}, s})$ the differential
of map (\ref{map}) at points $(1,\dot{s}zn_s)$ takes the form
\begin{eqnarray}\label{diff8*}
d\gamma_{(1,\dot{s}zn_s)}:(x,(n,w))\rightarrow -(Id-{\rm
Ad}(\dot{s}zn_s)^{-1})x+n+w, \\x\in\g_{\bf k}\simeq T_1(G_{\bf k}), (n,w)\in
\n_{{\bf k},s}+\z_{\bf k}^s\simeq T_{\dot{s}zn_s}(\dot{s}Z_{\bf k}^sN_{{\bf k}, s}), \nonumber
\end{eqnarray}
where $\n_{{\bf k}, s}\subset \g_{\bf k}$ is the Lie subalgebra corresponding to the subgroup $N_{{\bf k}, s}\subset G_{\bf k}$.

In order to show that the image of map (\ref{diff8*}) coincides
with $T_{\dot{s}zn_s}G_{\bf k}\simeq \g_{\bf k}$ we shall need a direct sum decomposition of the Lie algebra $\g_{\bf k}$ as a vector space,
\begin{equation}\label{dec}
\g_{\bf k}=\n_{\bf k}^s+ \z_{\bf k}^s + \opn_{\bf k}^s+ \h'_{\bf k},
\end{equation}
where $\h'_{\bf k}$ \index[not]{h@$\h'_{\bf k}$} is a complementary $s$--invariant subspace to $\h_{\bf k}^s$ in $\h_{\bf k}$.

We shall use isomorphism (\ref{cross}), $\alpha: N_{\bf k}^s\times \dot{s}Z_{\bf k}^sN_{{\bf k}, s}\rightarrow N_{\bf k}^s\dot{s}Z_{\bf k}^sN_{{\bf k}, s}=N_{\bf k}^s\dot{s}Z_{\bf k}^sN_{\bf k}^s$. By the definition $\alpha$ is the restriction of the map $\gamma$ to the subvariety $N_{\bf k}^s\times \dot{s}Z_{\bf k}^sN_{{\bf k}, s}\subset G_{\bf k}\times \dot{s}Z_{\bf k}^sN_{{\bf k}, s}$. Observe that in terms of the left trivialization of the tangent bundle $TG_{\bf k}$ the differential of the map $\alpha$ at points $(1,\dot{s}zn_s)\in N_{\bf k}^s\times \dot{s}Z_{\bf k}^sN_{{\bf k}, s}$, $n_s\in N_{{\bf k}, s},z\in Z_{\bf k}^s$ is given by
\begin{eqnarray}\label{diff8}
d\alpha_{(1,szn_s)}:(x,(n,w))\rightarrow -(Id-{\rm
Ad}(\dot{s}zn_s)^{-1})x+n+w, \\x\in\n\simeq T_1(N_{\bf k}^s), (n,w)\in
\n_{{\bf k}, s}+\z_{\bf k}^s\simeq T_{\dot{s}zn_s}(\dot{s}Z_{\bf k}^sN_{{\bf k},s}). \nonumber
\end{eqnarray}

Recall that the conjugation map $\alpha: N_{\bf k}^s\times \dot{s}Z_{\bf k}^sN_{{\bf k},s}\rightarrow N_{\bf k}^s\dot{s}Z_{\bf k}^sN_{{\bf k}, s}$ is an isomorphism of varieties, and hence its differential is a vector space isomorphism of the corresponding tangent spaces at all points. Using the left trivialization of the tangent bundle $TG_{\bf k}$ the tangent space $T_{\dot{s}zn_s}(N_{\bf k}^s\dot{s}Z_{\bf k}^sN_{{\bf k}, s})$ can be identified with $\n_{{\bf k}, s}+\z_{\bf k}^s+{\rm Ad}(\dot{s}zn_s)^{-1}\n_{\bf k}^s$, $T_{\dot{s}zn_sz}(N_{\bf k}^s\dot{s}Z_{\bf k}^sN_{{\bf k}, s})\simeq \n_{{\bf k}, s}+\z_{\bf k}^s+{\rm Ad}(\dot{s}zn_s)^{-1}\n_{\bf k}^s$. Therefore using (\ref{diff8}) and the fact that $d\alpha_{(1,\dot{s}zn_s)}$ is a vector space isomorphism we deduce that 
\begin{equation}\label{inc*}
(Id-{\rm
Ad}(\dot{s}zn_s)^{-1})\n_{\bf k}^s+\n_{{\bf k}, s}+\z_{\bf k}^s={\rm Ad}(\dot{s}zn_s)^{-1}\n_{\bf k}^s +\n_{{\bf k}, s}+\z_{\bf k}^s.
\end{equation}

Now observe that by the definition the vector subspace $(Id-{\rm Ad}(\dot{s}zn_s)^{-1})\n_{\bf k}^s\subset \g_{\bf k}$ is contained in the image of $d\alpha_{(1,\dot{s}zn_s)}$, and by (\ref{inc*}) the vector subspace ${\rm Ad}(\dot{s}zn_s)^{-1}\n_{\bf k}^s\subset \g_{\bf k}$ is also contained in the image of $d\alpha_{(1,\dot{s}zn_s)}$.
Since $\n_{\bf k}^s=(Id-{\rm Ad}(\dot{s}zn_s)^{-1})\n_{\bf k}^s+{\rm Ad}(\dot{s}zn_s)^{-1}\n_{\bf k}^s$, we deduce that $\n_{\bf k}^s$ is contained in the image of $d\alpha_{(1,\dot{s}zn_s)}$, and hence in the image of $d\gamma_{(1,\dot{s}zn_s)}$,
\begin{equation}\label{inc1}
\n_{\bf k}^s \subset {\rm Im}~d\gamma_{(1,\dot{s}zn_s)}.
\end{equation}

Next observe that similarly to (\ref{cross}) one can show that the conjugation map
\begin{equation}\label{ao}
N_{\bf k}^s\times N_{{\bf k}, s^{-1}}Z_{\bf k}^s\dot{s}\rightarrow N_{\bf k}^sZ_{\bf k}^s\dot{s}N_{\bf k}^s=N_{{\bf k}, s^{-1}}Z_{\bf k}^s\dot{s}N_{\bf k}^s, N_{{\bf k}, s^{-1}}:=\{n\in N_{\bf k}^s:\dot{s}^{-1}n\dot{s}\in \overline{N}_{\bf k}^s\}
\end{equation}
is an isomorphism of varieties.

Interchanging the roles of $N_{\bf k}^s$ and $\overline{N}_{\bf k}^s$ in (\ref{ao}) we immediately obtain that the conjugation map $\overline{\alpha}: \overline{N}_{\bf k}^s\times \overline{N}_{{\bf k}, s^{-1}}Z_{\bf k}^s\dot{s}\rightarrow \overline{N}_{{\bf k}, s^{-1}}Z_{\bf k}^s\dot{s}\overline{N}_{\bf k}^s$, where $\overline{N}_{{\bf k}, s^{-1}}=\dot{s}N_{{\bf k}, s}\dot{s}^{-1}$, is an isomorphism of varieties. Observe also that by the definition the map $\overline{\alpha}$ is the restriction of $\gamma$ to the subvariety $\overline{N}_{\bf k}^s\times \overline{N}_{{\bf k}, s^{-1}}Z_{\bf k}^s\dot{s}=\overline{N}_{\bf k}^s\times \dot{s}Z_{\bf k}^sN_{{\bf k}, s}\subset G_{\bf k} \times \dot{s}Z_{\bf k}^sN_{{\bf k}, s}$. 

Using the differential of the map $\overline{\alpha}$ we immediately infer, similarly to inclusion (\ref{inc1}), that
\begin{equation}\label{inc2}
\opn_{\bf k}^s \subset {\rm Im}~d\gamma_{(1,\dot{s}zn_s)}.
\end{equation}

Now observe that $\h_{\bf k}$ normalizes $\n_{\bf k}^s$, and $Z_{\bf k}^s$ is a subgroup of $L_{\bf k}^s$ the adjoint action of which has the property that for any $x'\in \h'_{\bf k}$, $l\in L_{\bf k}^s$ one has $({\rm Ad}(l)x')_{\h'_{\bf k}}=x'$, where for any element $y\in \g_{\bf k}$ we denote by $(y)_{\h'_{\bf k}}$ the $\h'_{\bf k}$--component of $y$ with respect to decomposition (\ref{dec}). Therefore we have for any $x'\in \h'_{\bf k}$
\begin{equation}\label{s22}
\left((Id-{\rm
Ad}(\dot{s}zn_s)^{-1})x'\right)_{\h'_{\bf k}}=(Id-{\rm
Ad}\dot{s}^{-1})x'.
\end{equation}

Since by the definition of $\h'_{\bf k}$ the operator ${\rm Ad}\dot{s}^{-1}$ has no fixed points in it, from formula (\ref{s22}) it follows that the operator $\h'_{\bf k}\to \h'_{\bf k}$, $x'\mapsto  (Id-{\rm Ad}\dot{s}^{-1})x'$ is invertible, so that $\h'_{\bf k}$ is contained in the image of $(Id-{\rm Ad}(\dot{s}zn_s)^{-1})$, and hence, by formula (\ref{diff8*}), in the image of $d\gamma_{(1,\dot{s}zn_s)}$. Recalling also inclusions (\ref{inc1}) and (\ref{inc2}) and taking into account the obvious inclusion $\z_{\bf k}^s\subset {\rm Im}~d\gamma_{(1,\dot{s}zn_s)}$ and decomposition (\ref{dec}) we deduce that the image of the map $d\gamma_{(1,\dot{s}zn_s)}$ coincides with $\g_{\bf k} \simeq T_{\dot{s}zn_s}G_{\bf k}$. Therefore the differential of the map $\gamma$ is surjective at all points. This completes the proof of part (ii).

(iii) Finally we show that $N^ssZ^sN^s$ is a closed subvariety in $G$ provided that $G$ is simply connected.

Recall that by Lemma \ref{NZsN} (ii) $N^ssZ^sN^s= N^ssZ^sN_s= N^sZ^ssN_s\simeq N^s\times Z^s\times N_s$ is a subvariety of $G$, and observe that multiplication by $s^{-1}$ from the right induces an isomorphiam of varieties, $N^sZ^ssN_s\simeq N^sZ^ssN_ss^{-1}$. By the definition of $N_s$, $N^sZ^ssN_ss^{-1}$ is also a subset of $N^sZ^s\overline{N}^s$. 

\begin{lemma}\label{NNclosed}
(i) $N^sZ^s\overline{N}^s$ \index[not]{N@$\overline{N}^s$} and $\overline{N}^sZ^sN^s$ are subvarieties of $G$ with the variety structure induced from $N^s\times Z^s\times \overline{N}^s$ using the bijective maps 
\begin{equation}\label{isoNZON}
N^s\times Z^s\times \overline{N}^s\to N^sZ^s\overline{N}^s, (n,z,\bar{n})\mapsto nz\bar{n}, 
\end{equation}
$$
\overline{N}^s\times Z^s\times N^s\to \overline{N}^sZ^sN^s, (\bar{n},z,n)\mapsto \bar{n}zn,
$$
induced by the multiplication in $G$.

(ii) The varieties $N^sZ^s\overline{N}^s$ and $\overline{N}^sZ^sN^s$ are closed in $G$.
\end{lemma}

\begin{proof}
(i) We prove the statement for $N^sZ^s\overline{N}^s$. The other case is treated in a similar way.

By Lemma \ref{vart} (i) with $A_{\bf k}=G$, $G_{\bf k}^1=Z^s$ and $G_{\bf k}^2=N^s$, $N^sZ^s\subset G$ is an algebraic subgroup with the variety structure induced from $N^s\times Z^s$ by the multiplication map $N^s\times Z^s\to N^sZ^s$ in $G$.

One can describe the set $N^sZ^s\overline{N}^s$ as the orbit of the element $1$ in $G$ for the following regular action of $N^sZ^s\times \overline{N}^s$ on $G$:
\begin{equation}\label{NZON}
(N^sZ^s\times \overline{N}^s)\times G \to G, ((nz,\bar{n}), g)\mapsto nzg\bar{n}^{-1}, n\in N^s,z\in Z^s, \bar{n}\in \overline{N}^s, g\in G.
\end{equation}
Thus $N^sZ^s\overline{N}^s$ has the natural structure of an algebraic variety by part (ii) of Lemma \ref{vart} with $A_{\bf k}=N^sZ^s\times \overline{N}^s$, $V=G_{\bf k}$, $\mathcal{O}=N^sZ^s\overline{N}^s$.

By definition (\ref{NZON}) with $g=1$ gives rise to a bijective $N^sZ^s\times \overline{N}^s$-equivariant morphism of $N^sZ^s\times \overline{N}^s$ and $N^sZ^s\overline{N}^s$ viewed as homogeneous spaces for $N^sZ^s\times \overline{N}^s$,
\begin{equation}\label{emor}
N^sZ^s\times \overline{N}^s\to N^sZ^s\overline{N}^s, (nz,\bar{n})\mapsto nz\bar{n}^{-1}, n\in N^s,z\in Z^s, \bar{n}\in \overline{N}^s.
\end{equation}

One verifies straightforwardly that this morphism induces an isomorphism of the tangent spaces at the points $(1,1)$ and $1$. Therefore by Lemma \ref{vart} (iii) with $X=N^sZ^s\times \overline{N}^s$, $Y=N^sZ^s\overline{N}^s$, $A_{\bf k}=N^sZ^s\times \overline{N}^s$ morphism (\ref{emor}) is an isomorphism of varieties. Composing it with the algebraic group isomorphism $N^sZ^s\simeq N^s\times Z^s$ and with the map taking the inverse in $\overline{N}^s$ we deduce that the first map in (\ref{isoNZON}) is an isomorphism of varieties.

(ii) We shall prove the statement for $N^sZ^s\overline{N}^s$. The other case is treated in a similar way. We shall consider the case when $\Delta_0$ is not empty. The other case can be considered in a similar way.

Recall that by the definition $\h'_{\mathbb{R}}$ and $\h_0$ are annihilators of each other with respect to the restriction of the bilinear form on $\g$ to $\h_{\mathbb{R}}$. As before, let $\h'^*_{\mathbb{R}}$ and $\h_0^*$ be the images of $\h'_{\mathbb{R}}$ and $\h_0$, respectively, under the isomorphism $\h_{\mathbb{R}}\simeq \h_{\mathbb{R}}^*$ induced by the bilinear form on $\g$.

Introduce the element 
\begin{equation}\label{barh0}
\bar{h}_0=\sum_{k=1}^{M}h_{i_k}\in \ha_\mathbb{R}'. 
\end{equation}
\index[not]{h@$\bar{h}_0$}
By the definition of $\Delta_+^s$ for any $x\in \ha_0^*$ one has $\bar{h}_0(x)=0$ and a root $\alpha\in \Delta\setminus \Delta_0$ belongs to $\Delta_+^s$ if and only if $\bar{h}_0(\alpha)>0$,
\begin{equation}\label{D+sdef0}
\Delta_+^s\setminus \Delta_0=\{\alpha\in \Delta| \bar{h}_0(\alpha)>0\},~ \bar{h}_0(x)=0,~x\in \ha_0^*.
\end{equation} 

Let $\bar{h}_0^*\in \h_{\mathbb{R}}^*$ \index[not]{h@$\bar{h}_0^*$} be the image in $\h_{\mathbb{R}}$ of the element $\bar{h}_0\in \h'_{\mathbb{R}}$ under the isomorphism $\h_{\mathbb{R}}\simeq \h_{\mathbb{R}}^*$ induced by the bilinear form on $\g$. Since $\bar{h}_0\in \h'_{\mathbb{R}}$ we actually have $\bar{h}_0^*\in \h'^*_{\mathbb{R}}$. 

Let $\alpha_1,\ldots ,\alpha_p$ be the simple roots in $\Gamma^s$ which do not belong to $\Delta_0$, $\omega_1,\ldots, \omega_p$ the corresponding fundamental weights. $\h'^*_{\mathbb{R}}$ is a vector subspace in the real linear span $\Pi$ of $\omega_1,\ldots, \omega_p$ as $\Pi$ is the annihilator of the subspace of $\h_{\mathbb{R}}^*$ spanned by the roots from $\Delta_0$ which is contained in $\h_0^*$. The subset $\Pi_+$ of $\Pi$ which consists of $x$ satisfying the condition $\left\langle x,\alpha\right\rangle>0$, $\alpha\in \Delta_+^s \setminus \Delta_0$ is open in $\Pi$ and by the definition $\bar{h}_0^*\in \Pi_+\cap \h'^*_{\mathbb{R}}$. Therefore the intersection $\Pi_+\cap \h'^*_{\mathbb{R}}$ is not empty and open in $\h'^*_{\mathbb{R}}$.

The roots $\gamma_1, \ldots, \gamma_{l'}$ form a linear basis of $\h'^*_{\mathbb{R}}$. They also span a $\mathbb{Z}$--sublattice $Q'$ in the $\mathbb{Z}$--lattice generated by $\omega_1,\ldots, \omega_p$ as every root is a linear combination of fundamental weights with integer coefficients and $\gamma_1, \ldots, \gamma_{l'}$ form a linear basis of $\h'^*_{\mathbb{R}}\subset \Pi$. Linear combinations of elements of $Q'$ with rational coefficients are dense in $\h'^*_{\mathbb{R}}$, and, in particular, in the open set $\Pi_+\cap \h'^*_{\mathbb{R}}$. Since the subset $\Pi_+$ of $\Pi$ consists of $x$ satisfying the condition $\left\langle x,\alpha\right\rangle>0$, $\alpha\in \Delta_+^s \setminus \Delta_0$, there is a linear basis of $\h'^*_{\mathbb{R}}$ which consists of linear combinations of $\omega_1,\ldots, \omega_p$ with positive rational coefficients. Multiplying the elements of this basis by appropriate positive integer numbers we obtain a linear basis $\Omega_i$, $i=1,\ldots, l'$ of $\h'^*_{\mathbb{R}}$ which consists of integral dominant weights of the form $\Omega_i=\sum_{j=1}^pg_{ij}\omega_j$, $g_{ij}\in \mathbb{Z},g_{ij}>0$. 

Recall that by the definition $\h^s=\h'^\perp$ is the orthogonal compliment to $\h'$ in $\h$ with respect to the restriction of the symmetric bilinear form on $\g$ to $\h$. $\h^s$ is the complexification of $\h_0$, and hence we deduce that an element $x\in \h$ belongs to $\h^s=\h'^\perp$ if and only if $\Omega_i(x)=0$, $i=1,\ldots, l'$,
\begin{equation}\label{xh'}
x\in \h^s=\h'^\perp \Longleftrightarrow \Omega_i(x)=0, i=1,\ldots, l'.
\end{equation}

Let $B_+^s$ \index[not]{B@$B_\pm^s$} be the Borel subgroup of $G$ corresponding to the system $\Delta_+^s$ of positive roots, $B_-^s$ the opposite Borel subgroup, $N_\pm^s$ \index[not]{N@$N_\pm^s$} their unipotent radicals, respectively.

Next observe that since $G$ is simply connected, for each integral dominant weight $\lambda$ with respect to the system $\Delta_+^s$ of positive roots there exists an irreducible finite--dimensional representation of $G$ with highest weight $\lambda$.

Let $V_{\Omega_i}$, $i=1,\ldots ,l'$ be the irreducible finite--dimensional representation of $G$ with highest weight $\Omega_i$ \index{representation!of a semisimple algebraic group!irreducible highest weight} with respect to the system $\Delta_+^s$ of positive roots. Denote by $v_{\Omega_i}$ a nonzero highest weight vector \index{representation!of a semisimple algebraic group!highest weight vector} in $V_{\Omega_i}$ and by $(~\cdot~ ,~\cdot~ )$ \index[not]{ZZ@$(~\cdot~ ,~\cdot~ )$} the non--degenerate contravariant form \index{form!contravariant non--degenerate!on a representation of a semisimple algebraic group} on $V_{\Omega_i}$ normalized in such a way that $( v_{\Omega_i},v_{\Omega_i}) =1$. The matrix element $g\mapsto ( v_{\Omega_i},g v_{\Omega_i})$, $g\in G$ can be regarded as a regular function on $G$ whose restriction to the big dense Bruhat cell $N_-^sHN_+^s$ is given by the character $\Omega_i$ of $H$, 
\begin{equation}\label{oi}
(v_{\Omega_i},n_-hn_+v_{\Omega_i}) =( v_{\Omega_i},hv_{\Omega_i}) =\Omega_i(h), n_-\in N_-^s,h\in H, n_+\in N_+^s. 
\end{equation}

Each fundamental weight $\omega_i$ can be regarded as a regular function $g\mapsto ( v_{\omega_j},g v_{\omega_j})$ on $G$ defined similarly to the regular function $g\mapsto ( v_{\Omega_i},g v_{\Omega_i})$ with $V_{\Omega_i}$ replaced by the irreducible finite--dimensional representation $V_{\omega_i}$ with highest weight $\omega_i$. By the definition of $\Omega_i$ the function $( v_{\Omega_i}, g v_{\Omega_i})$ can be expressed as a product of the functions $( v_{\omega_j}, g v_{\omega_j}) $, 
\begin{equation}\label{prodf}
( v_{\Omega_i},g v_{\Omega_i}) =\prod_{j=1}^{p}( v_{\omega_j},g v_{\omega_j}) ^{g_{ij}}, g\in G.
\end{equation}

Consider the closed subvariety in $G$ defined by the equations $( v_{\Omega_i},gv_{\Omega_i}) =1$, $i=1,\ldots ,l'$, $g\in G$. According to the Bruhat decomposition \index{Bruhat!decomposition} every element $g\in G$ belongs to $g\in B_-^swB_+^s$ for some $w\in W$. In this case $g$ can be written in the form $g=n_-whn_+$ for some $n_\pm \in N_\pm^s$, $h\in H$. Now by (\ref{oi}) and (\ref{prodf}) we have $( v_{\Omega_i},gv_{\Omega_i}) =\Omega_i(h)( v_{\Omega_i},wv_{\Omega_i}) =\Omega_i(h)\prod_{j=1}^{p}( v_{\omega_j},w v_{\omega_j}) ^{g_{ij}}$. As different weight spaces of $V_{\omega_j}$ are orthogonal with respect to the contravariant form, the right hand side of the last identity is not zero for all $i=1,\ldots ,l'$ if and only if $w$ fixes all weights $\omega_i$, $i=1,\ldots ,p$, i.e. if and only if $w$ belongs to the Weyl group of the root subsystem ${\Delta}_0$. Since ${\Delta}_0$ is the root system of the Levi factor $L^s$, and $( v_{\omega_i},v_{\omega_i}) =1$, one has $( v_{\Omega_i},wv_{\Omega_i}) \neq 0$, $i=1,\ldots ,l'$ if and only if $g\in N^sL^s\overline{N}^s$, and in that case $( v_{\Omega_i},gv_{\Omega_i}) =\Omega_i(h)$, where $g=n_-whn_+$ for some $n_\pm \in N_\pm^s$, $h\in H$, and $w$ is an element of the Weyl group of the root subsystem ${\Delta}_0$. 

As we already proved in (\ref{xh'}), an element $x\in \h$ belongs to $\h^s$ if and only if $\Omega_i(x)=0$, $i=1,\ldots, l'$.
Therefore the conditions $(v_{\Omega_i},gv_{\Omega_i}) =\Omega_i(h)=1$, $i=1,\ldots ,l'$ are equivalent to the fact that $h$ belongs to a subgroup ${H^0}^\prime$ of $H$ with the Lie algebra $\h^s$. Hence the equations $(v_{\Omega_i},gv_{\Omega_i}) =1$, $i=1,\ldots ,l'$ hold if and only if $g\in N^s{Z^s}^\prime\overline{N}^s$, where ${Z^s}^\prime\subset L^s$ is a subgroup of $L^s$ with the same Lie algebra as $Z^s$.
Thus $N^s{Z^s}^\prime\overline{N}^s$ is a closed subvariety of $G$. Its connected component \index{component!connected of a variety} containing the identity element of $G$ is obviously $N^sZ^s \overline{N}^s$. Thus $N^sZ^s \overline{N}^s$ is a closed subvariety of $G$.

\end{proof}

Now recall that $N^sZ^ssN_s\simeq N^sZ^ssN_ss^{-1}$ is a subset of $N^sZ^s\overline{N}^s$ as by the definition of $N_s$, $sN_ss^{-1}\subset \overline{N}^s$. In fact $sN_ss^{-1}$ is the algebraic subgroup in $\overline{N}^s$ generated by the one--parameter subgroups corresponding to the roots from the set $\{\alpha \in  \Delta_+^s: s^{-1}\alpha\in -\Delta_+^s\}$. So by isomorphisms (\ref{descrNZsN}) and (\ref{isoNZON}), and by Lemma \ref{NNclosed} (i), $N^sZ^ssN_ss^{-1}\simeq N^sZ^ssN_s$ is a closed subvariety of $N^sZ^s\overline{N}^s$ due to the closed embedding
$$
N^s\times Z^s\times sN_ss^{-1}\hookrightarrow N^s\times Z^s\times \overline{N}^s.
$$

Also Lemma \ref{NNclosed} (ii) immediately implies that the closed subvariety $N^sZ^ssN_ss^{-1}\subset N^sZ^s\overline{N}^s$ is also closed in $G$, and hence the variety $N^ssZ^sN^s= N^sZ^ssN_s=  N^sZ^ssN_ss^{-1}s\simeq N^sZ^ssN_ss^{-1}$ is closed in $G$.

\end{proof}

The subvarieties $\Sigma_{{\bf k}, s}\subset G_{\bf k}$ are analogues of the Slodowy slices in algebraic group theory.

\begin{remark}\label{crosssectrem}
In fact, in the construction of the inverse map to morphism (\ref{cross}) suggested in the proof of the previous proposition we only used relations in $G_{\bf k}$ arising from the corresponding Chevalley group over $\mathbb Z$. Therefore isomorphisms similar to (\ref{cross}) hold in the case when $G_{\bf k}$ is replaced with the Chevalley group (or even group scheme) over an arbitrary ring. \index{group!Chevalley}
\end{remark}


\setcounter{equation}{0}
\setcounter{theorem}{0}

\section{The Lusztig partition}\label{luspart}

In this section, as before, $G_{\bf k}$ is a connected finite--dimensional semisimple algebraic group over an algebraically closed field $\bf k$. At the same time we shall also consider connected finite--dimensional semisimple algebraic groups $G_p$ of the same type as $G_{\bf k}$, i.e. with the same root system, over algebraically closed fields for all characteristic exponents $p$.  

In the next section we shall show that for every conjugacy class $\mathcal{O}$ \index[not]{O@$\mathcal{O}$} \index{conjugacy class!in an algebraic group} in $G_{\bf k}$ one can find a subvariety $\Sigma_{s,{\bf k}}\subset G_{\bf k}$ such that $\mathcal{O}$ intersects $\Sigma_{s,{\bf k}}$ and ${\rm dim}~\mathcal{O}={\rm codim}~\Sigma_{s,{\bf k}}$. It turns out that there is a remarkable partition of the group $G_{\bf k}$ introduced in \cite{L2} the strata of which are unions of conjugacy classes of the same dimension. For each stratum of this partition there is a Weyl group element $s$ such that all conjugacy classes $\mathcal{O}$ from that stratum intersect $\Sigma_{s,{\bf k}}$, and ${\rm dim}~\mathcal{O}={\rm codim}~\Sigma_{s,{\bf k}}$. This, in particular, determines the dimensions of the conjugacy classes in each stratum. In this section, which is rather descriptive, we recall the definition of this partition called the Lusztig partition. The main property of this partition gives an affirmative answer to an old question about intersection of the conjugacy classes in $G_{\bf k}$ with Bruhat cells. The exposition in this section mainly follows paper \cite{L2} to which we refer the reader for technical details.

Note that all objects introduced in this section, except for the map $\phi_{G_p}$ and for the Lusztig partition itself, only depend on the Weyl group of $G_{\bf k}$ or of $G_p$, on the characteristic of $\bf k$ and on $p$. The fundamental reason behind this phenomenon is that in the definitions of these objects only representation theory of the Weyl group and the sets of unipotent conjugacy classes in $G_{\bf k}$ or in $G_p$ are used. These sets only depend on the characteristic of $\bf k$, on $p$ and on the Weyl group. The reader may always assume that $G_{\bf k}$ and all $G_p$ have the same root datum. \index{root!datum}

For any Weyl group $W$ let $\widehat{W}$ \index[not]{W@$\widehat{W}$} be the set of isomorphism classes of irreducible representations of $W$ over $\mathbb{Q}$. For any $E\in \widehat{W}$ let $b_E$ \index[not]{b@$b_E$} be the smallest nonnegative integer such that $E$ appears with non--zero multiplicity in the $b_E$-th symmetric power of the reflection representation of $W$. \index{representation!of a Weyl group!reflection} If this multiplicity is equal to $1$ then one says that $E$ is {\it good}. \index{representation!of a Weyl group!good} If $W'\subset W$ are two Weyl groups, and $E\in \widehat{W}'$ is good then there is a unique $\widetilde{E}\in \widehat{W}$ \index[not]{E@$\widetilde{E}$} such that $\widetilde{E}$ appears in the decomposition of the induced representation \index{representation!of a Weyl group!induced} ${\rm Ind}_{W'}^WE$, $b_{\widetilde{E}}=b_E$, and $\widetilde{E}$ is good. The representation $\widetilde{E}$ is called {\it j-induced} from $E$, and we write in this case $\widetilde{E}=j_{W'}^WE$. \index{representation!of a Weyl group!j-induced}

Let $g\in G_p$, and $g=g_sg_u$ its decomposition as a product of the semisimple part $g_s$ \index{semisimple part!of an algebraic group element} and of the unipotent part $g_u$. \index{unipotent part of an algebraic group element} Let $C=Z_{G_p}(g_s)^\circ$ \index[not]{C@$C$} be the identity component of the centralizer of $g_s$ in $G_p$. $C$ is a reductive subgroup of $G_p$ of the same rank as $G_p$. Let $H_p$ be a maximal torus in $C$. $H_p$ is also a maximal torus in $G_p$, and hence one has a natural imbedding
$$
W(C,H_p)=N_C(H_p)/H_p \hookrightarrow N_{G_p}(H_p)/H_p=W(G_p,H_p)=W,
$$
where $N_C(H_p), N_{G_p}(H_p)$ stand for the normalizers of $H_p$ in $C$ and in $G_p$, respectively, $W(C,H_p)$ is the Weyl group of the pair $(C,H_p)$ and $W(G_p,H_p)$ is the Weyl group of the pair $(G_p,H_p)$.

Let $E$ be the irreducible representation of $W(C,H_p)$ associated with the help of the Springer correspondence \index{Springer correspondence} to the conjugacy class of $g_u$ and the trivial local system \index{local system, trivial} on it. Then $E$ is good, and let $\widetilde{E}$ be the j-induced representation of $W$. This gives a well-defined map $\phi_{G_p}:G_p\rightarrow \widehat{W}$, $\phi_{G_p}(g)=\widetilde{E}$. \index[not]{f@$\phi_{G_p}$} The fibers of this map are called {\it the strata} \index{partition!Lusztig!strata} of $G_p$. By the definition the map $\phi_{G_p}$ is constant on each conjugacy class in $G_p$. Therefore the strata are unions of conjugacy classes.

Moreover, by 2.4 in \cite{L2} we have the following formula for the dimension of the centralizer $Z_{G_p}(g)$ \index{centralizer of an element in an algebraic group} \index[not]{Z@$Z_{G_p}(g)$} of any element $g\in G_p$ in $G_p$:
\begin{equation}\label{dimcon}
{\rm dim}~Z_{G_p}(g)={\rm rank}~G_p+2b_{\phi_{G_p}(g)},
\end{equation}
where ${\rm rank}~G_p$ is the rank of $G_p$.

It turns out that the image $\mathcal{R}(W)$ \index[not]{R@$\mathcal{R}(W)$} of $\phi_{G_p}$ only depends on $W$ and not on $p$ or the underlying root datum of $G_p$. It can be described as follows. Let $\mathcal{N}(G_p)$ \index[not]{N@$\mathcal{N}(G_p)$} be the unipotent variety of $G_p$ and $\underline{\mathcal{N}}(G_p)$ \index[not]{N@$\underline{\mathcal{N}}(G_p)$} the set of unipotent conjugacy classes in $G_p$. Let $\mathcal{X}^p(W)$ \index[not]{X@$\mathcal{X}^p(W)$} be the set of irreducible representations of $W$ associated by the Springer correspondence to unipotent classes in $\underline{\mathcal{N}}(G_p)$ and the trivial local systems on them. We shall identify $\mathcal{X}^p(W)$ and $\underline{\mathcal{N}}(G_p)$. Let $\mathcal{I}_p:\underline{\mathcal{N}}(G_p) \rightarrow \mathcal{X}^p(W)$ \index[not]{I@$\mathcal{I}_p$} be the corresponding bijective map.

\begin{proposition}{\bf (\cite{L2}, Sections 1 and 2)}\label{RWdescr}
We have
$$
\mathcal{R}(W)=\mathcal{X}^1(W)\cup_{r~{\rm prime}}\mathcal{X}^r(W).
$$

If $G_p$ is of type $A_n$, $(n\geq 1)$ or $E_6$ then $\mathcal{R}(W)=\mathcal{X}^1(W)$.

If $G_p$ is of type $B_n$ $(n\geq 2)$, $C_n$ $(n\geq 3)$, $D_n$ $(n\geq 4)$, $F_4$ or $E_7$ then $\mathcal{R}(W)=\mathcal{X}^2(W)$.

If $G_p$ is of type $G_2$ then $\mathcal{R}(W)=\mathcal{X}^3(W)$.

If $G_p$ is of type $E_8$ then $\mathcal{R}(W)=\mathcal{X}^2(W)\cup\mathcal{X}^3(W)$, and $\mathcal{X}^2(W)\cap\mathcal{X}^3(W)=\mathcal{X}^1(W)$.
\end{proposition}

The description of the set $\mathcal{R}(W)$ given in Proposition \ref{RWdescr} and the bijections $\underline{\mathcal{N}}(G_p) \rightarrow \mathcal{X}^p(W)$ yield certain maps between the sets $\underline{\mathcal{N}}(G_p)$ which preserve dimensions of conjugacy classes by (\ref{dimcon}). For instance, one always has an inclusion $\mathcal{X}^1(W) \subset \mathcal{X}^r(W)$ for any $r\geq 2$. The corresponding inclusion $\underline{\mathcal{N}}(G_1)\subset \underline{\mathcal{N}}(G_p)$ coincides with the Spaltenstein map \index{Spaltenstein map} $\pi^{G_{1}}_p: \underline{\mathcal{N}}(G_1)\rightarrow \underline{\mathcal{N}}(G_p)$ \index[not]{p@$\pi^{G_{1}}_p$} which is a bijection for good $p$ (see \cite{Spal}, Th\'{e}or\`{e}me III.5.2).

Now we discuss an alternative description of the strata in terms of intersections of conjugacy classes with Bruhat cells. \index{Bruhat!cell}
Fix a system of positive roots in $\Delta$. Note that $\Delta$ can be regarded as the root system of the pair $(G_p,H_p)$, $\Delta=\Delta(G_p,H_p)$. Let $B_{p,-}$ \index[not]{B@$B_{p,-}$} be the Borel subgroup in $G_p$ associated to the corresponding system of negative roots, $H_p\subset B_{p,-}$ the maximal torus contained in $B_{p,-}$, and $l(~\cdot~)$ the corresponding length function on $W$. Denote by $\underline{W}$ \index[not]{W@$\underline{W}$} the set of conjugacy classes in $W$. For each $w\in W=N_{G_p}(H_p)/H_p$ one can pick up a representative $\dot{w}\in N_{G_p}(H_p)$.
If $p$ is the characteristic exponent of $\bf k$, we write as before $B_{p,-}=B_{{\bf k},-}$, $\underline{\mathcal{N}}(G_p)=\underline{\mathcal{N}}(G_{\bf k})$, etc. \index[not]{B@$B_{{\bf k},-}$} \index[not]{N@$\underline{\mathcal{N}}(G_{\bf k})$}

Let $\mathcal{C}$ \index[not]{C@$\mathcal{C}$} be a conjugacy class in $W$. \index{conjugacy class!in a Weyl group} Pick up a representative $w\in \mathcal{C}$ of minimal possible length with respect to $l(~\cdot~)$. By Theorem 0.4 in \cite{L4'} there is a unique conjugacy class $\mathcal{O}\in \underline{\mathcal{N}}(G_{1})$ of minimal possible dimension which intersects the Bruhat cell $B_{1,-}\dot{w}B_{1,-}$ and does not depend on the choice of the minimal length representative $w$ in $\mathcal{C}$. We denote this class by $\Phi_1^{G_{1}}(\mathcal{C})$. \index[not]{F@$\Phi_1^{G_{1}}$}

As shown in Section 1.1 in \cite{L4'}, one can always find a representative $w\in \mathcal{C}$ of minimal possible length \index{Weyl group element!of minimal possible length in its conjugacy class} with respect to $l(~\cdot~)$ which is elliptic in a parabolic Weyl subgroup $W'\subset W$, i.e. $w$ acts without fixed points in the reflection representation of $W'$. \index{Weyl group element!elliptic} Indeed, by Theorem 3.2.12 in \cite{Gck} there is a parabolic subgroup $W'\subset W$ such that $\mathcal{C}\cap W'$ is an elliptic conjugacy class in $W'$, i.e. every element in it is elliptic in $W'$. By Lemma 3.1.14 in \cite{Gck} if $w\in \mathcal{C}\cap W'$ is of minimal possible length in its conjugacy class in $W'$ with respect to the restriction of $l(~\cdot~)$ to $W'$ then it is also of minimal possible length in $\mathcal{C}$ with respect to $l(~\cdot~)$.

Let $P'_{p}\subset G_{p}$ \index[not]{P@$P'_{p}$} be the parabolic subgroup which contains $B_{p,-}$ and corresponds to $W'$, and $M'_{p}$ \index[not]{M@$M'_{p}$} the semi-simple part of the Levi factor of $P'_{p}$, so that $W'$ is the Weyl group of $M'_{p}$. Let $\Phi_p^{G_{1}}(\mathcal{C})$ \index[not]{F@$\Phi_p^{G_{1}}$} be the unipotent class in $G_p$ containing the class $\pi_p^{M'_{1}}\Phi_1^{M'_{1}}(\mathcal{C})$. This class only depends on the conjugacy class $\mathcal{C}$, and hence one has a map $\Phi_p^{G_{1}}: \underline{W}\rightarrow \underline{\mathcal{N}}(G_p)$ which is in fact surjective by 4.5(a) in \cite{L4'}.

Let $\mathcal{C}\in \underline{W}$, and $m_\mathcal{C}$ \index[not]{m@$m_\mathcal{C}$} the dimension of the fixed point space for the action of any $w\in \mathcal{C}$ in the reflection representation.
Then by Theorem 0.2 in \cite{L3'} for any $\gamma \in \underline{\mathcal{N}}(G_p)$ there is a unique $\mathcal{C}_0\in (\Phi_p^{G_{1}})^{-1}(\gamma)$ \index[not]{C@$\mathcal{C}_0$} such that the function $m_\mathcal{C}: (\Phi_p^{G_{1}})^{-1}(\gamma) \rightarrow \mathbb{N}$ reaches its minimum at $\mathcal{C}_0$. We denote $\mathcal{C}_0$ by $\Psi_p^{G_{1}}(\gamma)$. \index[not]{P@$\Psi_p^{G_{1}}$} Thus one obtains an injective map $\Psi_p^{G_{1}}: \underline{\mathcal{N}}(G_p)\rightarrow \underline{W}$.

Now using the bijections $\mathcal{I}_p:\underline{\mathcal{N}}(G_p) \rightarrow \mathcal{X}^p(W)$ one can define the union $\widehat{\underline{\mathcal{N}}}(W)$ \index[not]{N@$\widehat{\underline{\mathcal{N}}}(W)$} of the sets $\underline{\mathcal{N}}(G_1)$ and of $\underline{\mathcal{N}}(G_r)$ over all prime $r$ as the union $\mathcal{X}^1(W)\cup_{r~{\rm prime}}\mathcal{X}^r(W)=\mathcal{R}(W)$.
Thus we have a bijection
\begin{equation}\label{bij}
F:\widehat{\underline{\mathcal{N}}}(W)=\underline{\mathcal{N}}(G_1)\cup_{r~{\rm prime}}\underline{\mathcal{N}}(G_r)\rightarrow \mathcal{X}^1(W)\cup_{r~{\rm prime}}\mathcal{X}^r(W)=\mathcal{R}(W).
\end{equation}
\index[not]{F@$F$}

Using the maps $\Phi_p^{G_{1}}$ one can also define a surjective map $\Phi^W: \underline{W}\rightarrow \widehat{\underline{\mathcal{N}}}(W)$ \index[not]{F@$\Phi^W$} as described in Section 4.1 in \cite{L2}. As it is observed in Section 4.1 of \cite{L2} (see also 0.4 in \cite{L3'}), if $\Phi_r^{G_1}(\mathcal{C})\in \underline{\mathcal{N}}(G_1)$ for all $r>1$ then $\Phi_r^{G_{1}}(\mathcal{C})$ is independent of $r$, and one puts $\Phi^W(\mathcal{C})=\Phi_r^{G_{1}}(\mathcal{C})$ for any $r>1$, and if $\Phi_r^{G_{1}}(\mathcal{C}) \not\in \underline{\mathcal{N}}(G_1)$ for some $r>1$ then $r$ is unique, and one defines $\Phi^W(\mathcal{C})=\Phi_r^{G_{1}}(\mathcal{C})$. By Proposition \ref{RWdescr} and by formula (\ref{bij}) the map $\Phi^W$ introduced this way is well defined.

By the definition there is a right-sided injective inverse $\Psi^W$ \index[not]{P@$\Psi^W$} to $\Phi^W$ such that if $\gamma \in \underline{\mathcal{N}}(G_1)$ then $\Psi^W(\gamma)=\Psi_1^{G_{1}}(\gamma)$, and if $\gamma \not \in \underline{\mathcal{N}}(G_1)$, and $\gamma \in \underline{\mathcal{N}}(G_r)$ then $\Psi^W(\gamma)=\Psi_r^{G_{1}}(\gamma)$.

Denote by $C(W)$ the image of $\widehat{\underline{\mathcal{N}}}(W)$ in $\underline{W}$ under the map $\Psi^W$, $C(W)=\Psi^W(\widehat{\underline{\mathcal{N}}}(W))$. \index[not]{C@$C(W)$} We shall identify $C(W)$, $\widehat{\underline{\mathcal{N}}}(W)$ and  $\mathcal{R}(W)$.

Now the strata of the Lusztig partition \index{partition!Lusztig} can be described geometrically as follows. Let $\mathcal{C} \in C(W)$. Pick up a representative $w\in \mathcal{C}$ of minimal possible length with respect to $l(~\cdot~)$. Denote by $\underline{G}_p$ \index[not]{G@$\underline{G}_p$} the set of conjugacy classes in $G_p$, and by $G_\mathcal{C}'$ \index[not]{G@$G_\mathcal{C}'$} the set of all conjugacy classes in $G_p$ which intersect the Bruhat cell $B_{p,-}\dot{w}B_{p,-}$. \index{Bruhat!cell} This definition does not depend on the choice of the minimal possible length representative $w$. Let
$$
d_\mathcal{C}={\hbox{\raise-1.5mm\hbox{${\textstyle \rm min}\atop {\scriptstyle \gamma \in G_\mathcal{C}'}$}}}~{\rm dim}~\gamma.
$$
\index[not]{d@$d_\mathcal{C}$}
Then the stratum $G_\mathcal{C}=\phi_{G_p}^{-1}(F(\Phi^W(\mathcal{C})))$ \index[not]{G@$G_\mathcal{C}$} \index{partition!Lusztig!strata} can be described as follows (see Theorem 5.2, \cite{L2}),
\begin{equation}\label{charstr}
G_\mathcal{C}=\bigcup_{\gamma \in G_\mathcal{C}', ~{\rm dim}\gamma=d_\mathcal{C}}\gamma.
\end{equation}

Thus we have a disjoint union
$$
G_p=\bigcup_{\mathcal{C}\in C(W)}G_\mathcal{C}.
$$
Note that by the definition of the stratum, for good $p$, if $\mathcal{C}\in {\rm Im}(\Psi_1^{G_{1}})$ then $G_\mathcal{C}$ contains a unique unipotent class, and if $\mathcal{C}\not \in {\rm Im}(\Psi_1^{G_{1}})$ then $G_\mathcal{C}$ does not contain unipotent classes.

The maps introduced above are summarized in the following diagram
\begin{equation}\label{diag}
\begin{array}{ccccccc}
   &  & \mathcal{X}^1(W) & \stackrel{\mathcal{I}_1}{\longleftarrow} & \underline{\mathcal{N}}(G_{1}) &  & \\
   &  & \downarrow ~{\iota} &  & \downarrow ~\pi^{G_{1}} &  & \\
  G_{\bf k} & \stackrel{\phi_{G_{\bf k}}}{\longrightarrow} & \mathcal{R}(W) & \stackrel{F}{\longleftarrow} & \widehat{\underline{\mathcal{N}}}(W) & {\hbox{\raise-1mm\hbox{${\textstyle \stackrel{\Phi^W}{\longleftarrow}}\atop {\hbox{${\textstyle \longrightarrow}\atop {\scriptstyle \Psi^W}$}}$}}} & \underline{W},
\end{array}
\end{equation}
where $\iota$ is an inclusion, bijections $\mathcal{I}_1$ and $F$ are induced by the Springer correspondence with the trivial local data, and the inclusion $\pi^{G_{1}}$ is induced by the Spaltenstein map.

For exceptional groups the maps $\mathcal{I}_1$ and $F$ can be described explicitly using tables in \cite{Spal1}, the maps $\Phi^W$ and $\Psi^W$ can be described using the tables in Section 2 in \cite{L3'}, and the maps $\iota$ and $\pi^{G_1}$ can be described explicitly using the tables of unipotent classes in \cite{Li}, Chapter 22 or \cite{Spal1} (note that the labeling for unipotent classes in bad characteristics in \cite{Li} differs from that in \cite{Spal1}). The dimensions of the conjugacy classes in the strata in $G_{\bf k}$ can be obtained using the dimension tables of the centralizers of unipotent elements in the case when a stratum contains a unipotent class (see \cite{Car1,Li}), the tables for the dimensions of the centralizers of the unipotent elements in bad characteristic when a stratum does not contain a unipotent class (see \cite{Li}) or formula (\ref{dimcon}) and the tables of the values of the $b$--invariant $b_E$ for representations of Weyl groups (see \cite{Car1,Gck}). Note that formula (\ref{dimcon}) implies that
if $\mathcal{O}$ is any conjugacy class in $G_\mathcal{C}$, $\mathcal{O}\in G_\mathcal{C}$ then
\begin{equation}\label{od}
{\rm dim}~\mathcal{O}={\rm dim}~\Phi^W(\mathcal{C}).
\end{equation}

In the case of classical groups all those maps and dimensions are described in terms of partitions (see \cite{Car1,Gk1,Li,L2',L3',L4',Spal}). In the case of classical matrix groups the strata can also be described explicitly (see \cite{L2}).
We recall this description below. By (\ref{od}) the dimensions of the conjugacy classes in every stratum of $G_{\bf k}$ are equal to the dimension of the corresponding conjugacy class in $\widehat{\underline{\mathcal{N}}}(W)$. The dimensions of centralizers of unipotent elements in arbitrary characteristic can be found in \cite{Hess,Li}.

If $\bar{\lambda} =(\bar{\lambda}_1\geq \bar{\lambda}_2\geq \ldots\geq \bar{\lambda}_m)$ \index[not]{l@$\bar{\lambda}$} is a partition \index{partition} we denote by $\bar{\lambda}^* =(\bar{\lambda}_1^*\geq \bar{\lambda}_2^*\geq \ldots\geq \bar{\lambda}_m^*)$ \index[not]{l@$\bar{\lambda}^*$} the corresponding dual partition. \index{partition!dual} It is defined by the property that $\bar{\lambda}_1^*=m$ and $\bar{\lambda}_i^*-\bar{\lambda}_{i+1}^*=l_i(\bar{\lambda})$, where $l_i(\bar{\lambda})$ \index[not]{l@$l_i(\bar{\lambda})$} is the number of times $i$ appears in the partition $\bar{\lambda}$. We also denote by $\bar{\tau}(\bar{\lambda})$ the length \index{partition!length} of $\bar{\lambda}$, $\bar{\tau}(\bar{\lambda})=m$. \index[not]{t@$\bar{\tau}(\bar{\lambda})$} If a partition $\bar{\mu}$ is obtained from $\bar{\lambda}$ by adding a number of zeroes, we shall identify $\bar{\lambda}$ and $\bar{\mu}$.


\subsection*{$\bf A_n$} \index[not]{A@$A_n$}

$G_{\bf k}$ is of type ${\rm SL}(V)$, \index[not]{S@${\rm SL}(V)$} where $V$ is a vector space of dimension $n+1 \geq 1$ over an algebraically closed field ${\bf k}$ of characteristic exponent $p\geq 1$. $W$ is the group of permutations \index{group!permutation} of $n+1$ elements.
All sets in (\ref{diag}), except for $G_{\bf k}$, are identified with the set of partitions of $n+1$, and under this identification all the maps, except for $\phi_{G_{\bf k}}$, are the identity maps.

To describe $\phi_{G_{\bf k}}$ for $G_{\bf k}={\rm SL}(V)$ we choose a sufficiently large $m \in \mathbb{N}$. Let $g \in G_{\bf k}$. For any $x \in {\bf k}^*$ let $V_x$ be the
generalized $x$--eigenspace of $g:V \rightarrow V$ and let $\lambda^x_1
\geq \lambda^x_2
\geq \ldots \geq \lambda^x_m$ be the sequence
in $\mathbb{N}$ whose terms are the sizes of the Jordan blocks \index{blocks!Jordan} of $x^{-1}g: V_x \rightarrow V_x$. Then $\phi_{G_{\bf k}}(g)$ is
the partition $\bar{\lambda}(g):=(\bar{\lambda}_1(g) \geq \bar{\lambda}_2(g) \geq \ldots \geq \bar{\lambda}_m(g))$ given by $\bar{\lambda}_j(g) = \sum_{x\in {\bf k}^*} \lambda^x_j$. \index[not]{l@$\bar{\lambda}(g)$}

If $g$ is any element in the stratum $G_{\bar{\lambda}}$ \index[not]{G@$G_{\bar{\lambda}}$} corresponding to a partition $\bar{\lambda}=(\bar{\lambda}_1
\geq \bar{\lambda}_2\geq \ldots \geq \bar{\lambda}_m)$, $\bar{\lambda}_m\geq 1$, then
\begin{equation}\label{dimsl}
{\rm dim}~Z_{G_{\bf k}}(g)=n+2\sum_{i=1}^m(i-1)\bar{\lambda}_i.
\end{equation}

The element of $\underline{W}$ which corresponds to $\bar{\lambda}$ is the Coxeter class in the Weyl subgroup \index{conjugacy class!in a Weyl group!Coxeter} of the type
\begin{equation}\label{san}
A_{\bar{\lambda}_1-1}+A_{\bar{\lambda}_2-1}+\ldots +A_{\bar{\lambda}_m-1}.
\end{equation}
The summands in diagram (\ref{san}) are called {\it blocks} of $W$. \index{blocks!in a Weyl group} Blocks of type $A_0$ are called trivial.


\subsection*{$\bf C_n$} \index[not]{C@$C_n$}

$G_{\bf k}$ is of type ${\rm Sp}(V)$, \index[not]{S@${\rm Sp}(V)$} where $V$ is a symplectic space \index{space!symplectic} of dimension $2n$, $n \geq 2$ over an algebraically closed field ${\bf k}$ of characteristic exponent $p$.
$W$ is the group of permutations of the set $E=\{\varepsilon_1,\ldots, \varepsilon_n,-\varepsilon_1,\ldots,-\varepsilon_n\}$ \index[not]{e@$\varepsilon_i$} which also commute with the involution  $\varepsilon_i\mapsto -\varepsilon_i$.
Each element $s\in W$ can be expressed as a product of disjoint cycles \index{cycle!in a permutation group} of the form
$$
\varepsilon_{m_1}\rightarrow \pm \varepsilon_{m_2}\rightarrow \pm \varepsilon_{m_3}\rightarrow \ldots \rightarrow  \pm \varepsilon_{m_r}\rightarrow \pm \varepsilon_{m_1}.
$$
This cycle is of length $r$; it is called positive \index{cycle!in a permutation group!positive} if $s^r(\varepsilon_{m_1})=\varepsilon_{m_1}$ and negative \index{cycle!in a permutation group!negative} if $s^r(\varepsilon_{m_1})=-\varepsilon_{m_1}$. The lengths of the cycles together with their signs give a set of positive or negative integers called the signed cycle-type of $s$. To each positive cycle of $s$ of length $r$ there corresponds a pair of orbits $X,-X$, $X\neq -X$, $|X|=r$, for the action of the group $\left< s \right>$ generated by $s$ on the set $E=\{\varepsilon_1,\ldots, \varepsilon_n,-\varepsilon_1,\ldots,-\varepsilon_n\}$, and to each negative cycle of $s$ of length $r$ there corresponds an orbit $X$, $X=-X$, $|X|=2r$, for the action of $\left< s \right>$ on $E$. We call orbits of the former type positive \index{group!permutation!positive orbit} and orbits of the latter type negative. \index{group!permutation!negative orbit} A positive cycle of length $1$ is called trivial. \index{cycle!in a permutation group!trivial} It corresponds to a pair of fixed points for the action of $\left< s \right>$ on $E$.

Elements of $\underline{W}$ are parametrized by pairs of partitions $(\bar{\lambda},\bar{\mu})$  satisfying the following conditions. 
\begin{itemize}
\item
The parts of $\bar{\lambda}$ are even (for any $w\in \mathcal{C}\in \underline{W}$ they are the numbers of elements in the negative orbits $X$, $X=-X$,  in $E$ for the action of the group $\left< w \right>$ generated by $w$);
\item 
$\bar{\mu}$ consists of pairs of equal parts (they are the numbers of elements in the positive $\left< w \right>$--orbits $X$ in $E$; these orbits appear in pairs $X,-X$, $X\neq -X$);  
\item
$\sum \bar{\lambda}_i +\sum \bar{\mu}_j=2n$. 
\end{itemize}
We denote this set of pairs of partitions by $\mathcal{A}^1_{2n}$. \index[not]{A@$\mathcal{A}^1_{2n}$}

An element of $\underline{W}$ which corresponds to a pair $(\bar{\lambda},\bar{\mu})$, $\bar{\lambda}=(\bar{\lambda}_1
\leq \bar{\lambda}_2\leq \ldots \leq \bar{\lambda}_m)$ and $\bar{\mu}=(\bar{\mu}_1=\bar{\mu}_2\leq \ldots \leq \bar{\mu}_{2k-1}=\bar{\mu}_{2k})$ is the Coxeter class in the Weyl subgroup of the type
\begin{equation}\label{scn}
C_{\frac{\bar{\lambda}_1}{2}}+C_{\frac{\bar{\lambda}_2}{2}}+\ldots +C_{\frac{\bar{\lambda}_m}{2}}+A_{\bar{\mu}_1-1}+A_{\bar{\mu}_3-1}+\ldots +A_{\bar{\mu}_{2k-1}-1}.
\end{equation}

If the characteristic exponent of ${\bf k}$ is not equal to 2, the elements of $\underline{\mathcal{N}}(G_{\bf k})$ are parametrized by the partitions $\bar{\lambda}$ of $2n$ for which $l_j(\bar{\lambda})$ is even for odd $j$. We denote this set of partitions by $\mathcal{T}_{2n}$. \index[not]{T@$\mathcal{T}_{2n}$} In the case of $G_{\bf k}={\rm Sp}(V)$ the parts of $\bar{\lambda}$ are just the sizes of the Jordan blocks in $V$ of the unipotent elements from the conjugacy class corresponding to $\bar{\lambda}$.

In this case $\widehat{\underline{\mathcal{N}}}(W)=\underline{\mathcal{N}}(G_2)$, and $G_2$ is of type ${\rm Sp}(V)$ where $V$ is a symplectic space of dimension $2n$ over an algebraically closed field of characteristic $2$. Elements of $\underline{\mathcal{N}}(G_2)$ are parametrized by pairs $(\bar{\lambda},\varepsilon)$, where $\bar{\lambda} = (\bar{\lambda}_1
\leq \bar{\lambda}_2\leq \ldots \leq \bar{\lambda}_m)\in \mathcal{T}_{2n}$, and $\varepsilon:\{\bar{\lambda}_1, \bar{\lambda}_2, \ldots , \bar{\lambda}_m\}\rightarrow \{0,1,\omega\}$ is a function such that
\begin{equation}\label{eps}
\varepsilon(k)=\left\{
                 \begin{array}{ll}
                 \omega  & \hbox{if $k$ is odd;} \\
                   1  & \hbox{if $k=0$;} \\
                   1  & \hbox{if $k>0$ is even, $l_k(\bar{\lambda})$ is odd;} \\
                   0~{\rm or}~1  & \hbox{if $k>0$ is even, $l_k(\bar{\lambda})$ is even.}
                 \end{array}
               \right.
\end{equation}
We denote the set of such pairs $(\bar{\lambda},\varepsilon)$ by $\mathcal{T}_{2n}^2$. \index[not]{T@$\mathcal{T}_{2n}^2$}

Elements of $\widehat{W}$ are parametrized by pairs of partitions $(\bar{\alpha},\bar{\beta})$ written in non--decreasing order, $\bar{\alpha}_1\leq \bar{\alpha}_2\leq \ldots \leq \bar{\alpha}_{\bar{\tau}(\bar{\alpha})}$, $\bar{\beta}_1\leq \bar{\beta}_2 \leq \ldots \leq \bar{\beta}_{\bar{\tau}(\bar{\beta})}$, and such that $\sum \bar{\alpha}_i+\sum \bar{\beta}_i=n$. By adding zeroes we can assume that the length $\bar{\tau}(\bar{\alpha})$ of $\bar{\alpha}$ is related to the length of $\bar{\beta}$ by $\bar{\tau}(\bar{\alpha})=\bar{\tau}(\bar{\beta})+1$. The set of such pairs is denoted by $X_{n,1}$. \index[not]{X@$X_{n,1}$}

The maps $\mathcal{I}_1,F$ can be described as follows. Let $\bar{\lambda} = (\bar{\lambda}_1
\leq \bar{\lambda}_2\leq \ldots \leq \bar{\lambda}_{2m+1})\in \mathcal{T}_{2n}$, and assume that $\bar{\lambda}_1=0$. If $\mathcal{I}_1(\bar{\lambda})=((c'_1,c'_3,\ldots,c'_{2m+1}),(c'_2,c'_4,\ldots,c'_{2m}))$ then the parts $c'_i$ are defined by induction starting from $c'_1=0$,
$$
  \begin{array}{ll}
    c'_i=\frac{\bar{\lambda}_i}{2} & \hbox{if $\bar{\lambda}_i$ is even and $c'_{i-1}$ is already defined;} \\
    c'_i=\frac{\bar{\lambda}_i+1}{2} & \hbox{if $\bar{\lambda}_i=\bar{\lambda}_{i+1}$ is odd and $c'_{i-1}$ is already defined;} \\
    c'_{i+1}=\frac{\bar{\lambda}_i-1}{2} & \hbox{if $\bar{\lambda}_i=\bar{\lambda}_{i+1}$ is odd and $c'_{i}$ is already defined.}
  \end{array}
$$

The image of $\mathcal{I}_1$ consists of all pairs $((c'_1,c'_3,\ldots,c'_{2m+1}),(c'_2,c'_4,\ldots,c'_{2m}))\in X_{n,1}$ such that $c'_i\leq c'_{i+1}+1$ for all $i$.

If $F(\bar{\lambda},\varepsilon)=((c_1,c_3,\ldots,c_{2m+1}),(c_2,c_4,\ldots,c_{2m}))$ then the parts $c_i$ are defined by induction starting from $c_1=0$,
$$
  \begin{array}{ll}
    c_i=\frac{\bar{\lambda}_i}{2} & \hbox{if $\bar{\lambda}_i$ is even, $\varepsilon(\bar{\lambda}_i)=1$ and $c_{i-1}$ is already defined;} \\
    c_i=\frac{\bar{\lambda}_i+1}{2} & \hbox{if $\bar{\lambda}_i=\bar{\lambda}_{i+1}$ is odd and $c_{i-1}$ is already defined;} \\
    c_{i+1}=\frac{\bar{\lambda}_i-1}{2} & \hbox{if $\bar{\lambda}_i=\bar{\lambda}_{i+1}$ is odd and $c_{i}$ is already defined;} \\
    c_i=\frac{\bar{\lambda}_i+2}{2} & \hbox{if $\bar{\lambda}_i=\bar{\lambda}_{i+1}$ is even, $\varepsilon(\bar{\lambda}_i)=\varepsilon(\bar{\lambda}_{i+1})=0$ and $c_{i-1}$ is already defined;} \\
    c_{i+1}=\frac{\bar{\lambda}_i-2}{2} & \hbox{if $\bar{\lambda}_i=\bar{\lambda}_{i+1}$ is even, $\varepsilon(\bar{\lambda}_i)=\varepsilon(\bar{\lambda}_{i+1})=0$ and $c_{i}$ is already defined.}
  \end{array}
$$

The image $\mathcal{R}(W)$ of $F$ consists of all pairs $((c_1,c_3,\ldots,c_{2m+1}),(c_2,c_4,\ldots,c_{2m}))\in X_{n,1}$ such that $c_i\leq c_{i+1}+2$ for all $i$.

The map $\Phi^W$ is defined by $\Phi^W(\bar{\lambda},\bar{\mu})=(\bar{\nu},\varepsilon)$, where the set of parts of $\bar{\nu}$ is just the union of the sets of parts of $\bar{\lambda}$ and $\bar{\mu}$, and
$$
\varepsilon(k)=\left\{
                 \begin{array}{ll}
                   1 & \hbox{if $k\in 2\mathbb{N}$ is a part of $\bar{\lambda}$;} \\
                   0 & \hbox{if $k\in 2\mathbb{N}$ is not a part of $\bar{\lambda}$;} \\
                   \omega & \hbox{if $k$ is odd.}
                 \end{array}
               \right.
$$

The map $\Psi^W$ associates to each pair $(\bar{\nu}, \varepsilon)$ a unique point $(\bar{\lambda},\bar{\mu})$ in the preimage $(\Phi^W)^{-1}(\bar{\nu}, \varepsilon)$ such that the number of parts of $\bar{\mu}$ is minimal possible. This point is defined by the conditions
$$
l_k(\bar{\lambda})=\left\{
               \begin{array}{ll}
                 0 & \hbox{if $k$ is odd or $k$ is even, $l_k(\bar{\nu})\geq 2$ is even and $\varepsilon(k)=0$;} \\
                 l_k(\bar{\nu}) & \hbox{otherwise,}
               \end{array}
             \right.
$$
$$
l_k(\bar{\mu})=\left\{
               \begin{array}{ll}
                 l_k(\bar{\nu}) & \hbox{if $k$ is odd or $k$ is even, $l_k(\bar{\nu})\geq 2$ is even and $\varepsilon(k)=0$;} \\
                 0 & \hbox{otherwise.}
               \end{array}
             \right.
$$

The map $\pi^{G_{1}}$ is given by $\pi^{G_{1}}(\bar{\lambda})=(\bar{\lambda}, \varepsilon')$, where
\begin{equation}\label{econd}
\varepsilon'(k)=\left\{
                 \begin{array}{ll}
                 \omega  & \hbox{if $k$ is odd;} \\
                   1  & \hbox{if $k$ is even.}
                 \end{array}
               \right.
\end{equation}
The map ${{\pi}}^{G_{1}}$ is injective and its image consists of pairs $(\bar{\lambda}, \varepsilon')\in {\mathcal{T}}_{2n}^2$, where $\varepsilon'$ satisfies conditions (\ref{econd}).

To describe $\phi_{G_{\bf k}}$ for $G_{\bf k}={\rm Sp}(V)$ we choose a sufficiently large $m \in \mathbb{N}$. Let $g \in G_{\bf k}$. For any $x \in {\bf k}^*$ let $V_x$ be the
generalized $x$--eigenspace of $g:V \rightarrow V$. For any $x \in {\bf k}^*$ such that $x^2\neq 1$ let $\lambda^x_1
\geq \lambda^x_2
\geq \ldots \geq \lambda^x_{2m+1}$ be the sequence
in $\mathbb{N}$ whose terms are the sizes of the Jordan blocks of $x^{-1}g: V_x \rightarrow V_x$.

For any $x \in {\bf k}^*$ with $x^2= 1$ let $\lambda^x_1
\geq \lambda^x_2\geq \ldots \geq \lambda^x_{2m+1}$ be the sequence in $\mathbb{N}$,
where $((\lambda^x_1
\geq \lambda^x_3\geq \ldots \geq \lambda^x_{2m+1}),(\lambda^x_2
\geq \lambda^x_4\geq \ldots \geq \lambda^x_{2m}))$ is the pair of partitions
such that the corresponding irreducible representation of the Weyl
group of type $B_{{\rm dim}~V_x/2}$ is the Springer representation \index{representation!of a Weyl group!Springer} attached to the unipotent
element $x^{-1}g \in {\rm Sp}(V_x)$ and to the trivial local data.

Let $\bar{\lambda}(g)$ be
the partition $\bar{\lambda}_1(g) \geq \bar{\lambda}_2(g) \geq \ldots \geq \bar{\lambda}_{2m+1}(g)$ \index[not]{l@$\bar{\lambda}(g)$} given by $\bar{\lambda}_j(g) = \sum_{x} \lambda^x_j$, where $x$ runs over a set of representatives for the orbits of the
involution $a\mapsto a^{-1}$ of ${\bf k}^*$.
Now $\phi_{G_{\bf k}}(g)$ is the pair of partitions $((\bar{\lambda}_1(g)
\geq \bar{\lambda}_3(g)\geq \ldots \geq \bar{\lambda}_{2m+1}(g)),(\bar{\lambda}_2(g)
\geq \bar{\lambda}_4(g)\geq \ldots \geq \bar{\lambda}_{2m}(g)))$.

If $g$ is any element in the stratum $G_{(\bar{\lambda},\varepsilon)}$ \index[not]{G@$G_{(\bar{\lambda},\varepsilon)}$} corresponding to a pair $(\bar{\lambda},\varepsilon)\in \mathcal{T}_{2n}^2$, $\bar{\lambda}=(\bar{\lambda}_1
\geq \bar{\lambda}_2\geq \ldots \geq \bar{\lambda}_m)$ then
\begin{equation}\label{dimsp}
{\rm dim}~Z_{G_{\bf k}}(g)=n+\sum_{i=1}^m(i-1)\bar{\lambda}_i+\frac{1}{2}|\{i:\bar{\lambda}_i ~\hbox{is odd}\}|+|\{i:\bar{\lambda}_i ~\hbox{is even and}~\varepsilon(\bar{\lambda}_i)=0\}|.
\end{equation}


\subsection*{$\bf B_n$} \index[not]{B@$B_n$}

$G_{\bf k}$ is of type ${\rm SO}(V)$, \index[not]{S@${\rm SO}(V)$} where $V$ is a vector space of dimension $2n+1$, $n \geq 2$ over an algebraically closed field ${\bf k}$ of characteristic exponent $p$. If $p\neq 2$ then $V$ is equipped with a non--degenerate symmetric bilinear form. \index{form!non--degenerate symmetric bilinear!on a vector space} If $p=2$ then $V$ is equipped with a bilinear form $B(~\cdot~,~\cdot~)$ \index[not]{B@$B(~\cdot~,~\cdot~)$} and a non--zero quadratic form $Q(~\cdot~)$ \index[not]{Q@$Q(~\cdot~)$} \index{form!quadratic on a vector space} such that
$$
B(x,y)=Q(x+y)-Q(x)-Q(y),~x,y\in V,
$$
and the restriction of $Q(~\cdot~)$ to the null space $V^\perp=\{x\in V:B(x,y)=0~\forall~y\in V\}$ of $B(~\cdot~,~\cdot~)$ has zero kernel.

$W$ is the same as in the case of $C_n$. Therefore we can use the description of the set $\underline{W}$ in terms of pairs of partitions introduced in the case of $C_n$. 

An element of $\underline{W}$ which corresponds to a pair $(\bar{\lambda},\bar{\mu})$, $\bar{\lambda}=(\bar{\lambda}_1\geq \bar{\lambda}_2\geq \ldots \geq \bar{\lambda}_m)$ and $\bar{\mu}=(\bar{\mu}_1=\bar{\mu}_2\geq \ldots \geq \bar{\mu}_{2k-1}=\bar{\mu}_{2k})$ is the class represented by the sum of the blocks in the following diagram (we use the notation of \cite{C}, Section 7 for the conjugacy classes in $W$) \index{blocks!in a Weyl group}
\begin{eqnarray}
A_{\bar{\mu}_1-1}+A_{\bar{\mu}_3-1}+\ldots +A_{\bar{\mu}_{2k-1}-1}+ \qquad \qquad \qquad \qquad \qquad \qquad \qquad \qquad \qquad \qquad \qquad \qquad \nonumber \\
+D_{\frac{\bar{\lambda}_1+\bar{\lambda}_2}{2}}(a_{\frac{\bar{\lambda}_2}{2}-1})+
D_{\frac{\bar{\lambda}_3+\bar{\lambda}_4}{2}}(a_{\frac{\bar{\lambda}_4}{2}-1})+\ldots +D_{\frac{\bar{\lambda}_{m-2}+\bar{\lambda}_{m-1}}{2}}(a_{\frac{\bar{\lambda}_{m-1}}{2}-1})+B_{\frac{\bar{\lambda}_m}{2}} ~(\mbox{$m$ is odd}),\label{sbn} \\
A_{\bar{\mu}_1-1}+A_{\bar{\mu}_3-1}+\ldots +A_{\bar{\mu}_{2k-1}-1}+ \qquad \qquad \qquad \qquad \qquad \qquad \qquad \qquad \qquad \qquad \qquad \qquad \nonumber \\
+D_{\frac{\bar{\lambda}_1+\bar{\lambda}_2}{2}}(a_{\frac{\bar{\lambda}_2}{2}-1})+D_{\frac{\bar{\lambda}_3+\bar{\lambda}_4}{2}}(a_{\frac{\bar{\lambda}_4}{2}-1})+\ldots +D_{\frac{\bar{\lambda}_{m-1}+\bar{\lambda}_{m}}{2}}(a_{\frac{\bar{\lambda}_{m}}{2}-1})~(\mbox{$m$ is even}),\nonumber
\end{eqnarray}
where it is assumed that $D_k(a_0)=D_k$. \index[not]{D@$D_i(a_j)$}

If the characteristic of $\bf k$ is not equal to 2, the elements of $\underline{\mathcal{N}}(G_{\bf k})$ are parametrized by the partitions $\bar{\lambda}$ of $2n+1$ for which $l_j(\bar{\lambda})$ is even for even $j$. We denote this set of partitions by $\mathcal{Q}_{2n+1}$. \index[not]{Q@$\mathcal{Q}_{2n+1}$} In the case of $G_{\bf k}={\rm SO}(V)$ the parts of $\bar{\lambda}$ are just the sizes of the Jordan blocks in $V$ of the unipotent elements from the conjugacy class corresponding to $\bar{\lambda}$.

In this case $\widehat{\underline{\mathcal{N}}}(W)=\underline{\mathcal{N}}(G_2)$, and $G_2$ is of type ${\rm SO}(V)$. In fact $G_2$ is isomorphic to a group of type ${\rm Sp}(V')$, ${\rm dim}~V'=2n$ (see e.g. Section 8.1 in \cite{Spal}), and hence
$\underline{\mathcal{N}}(G_2) \simeq \mathcal{T}_{2n}^2$.

We also have $\widehat{W}\simeq X_{n,1}$, and the map $F$ is the same as in case of $C_n$.

The map $\mathcal{I}_1$ can be described as follows. Let $\bar{\lambda} = (\bar{\lambda}_1
\leq \bar{\lambda}_2\leq \ldots \leq \bar{\lambda}_{2m+1})\in \mathcal{Q}_{2n+1}$. If $$\mathcal{I}_1(\bar{\lambda})=((c'_1,c'_3,\ldots,c'_{2m+1}),(c'_2,c'_4,\ldots,c'_{2m}))$$ then the parts $c'_i$ are defined by induction starting from $c'_1$,
\begin{equation}\label{algf}
  \begin{array}{ll}
    c'_i=\frac{\bar{\lambda}_i-1}{2}+i-1-2\left[\frac{i-1}{2}\right] & \hbox{if $\bar{\lambda}_i$ is odd and $c'_{i-1}$ is already defined;} \\
    c'_i=\frac{\bar{\lambda}_i}{2} & \hbox{if $\bar{\lambda}_i=\bar{\lambda}_{i+1}$ is even and $c'_{i-1}$ is already defined;} \\
    c'_{i+1}=\frac{\bar{\lambda}_i}{2} & \hbox{if $\bar{\lambda}_i=\bar{\lambda}_{i+1}$ is even and $c'_{i}$ is already defined.}
  \end{array}
\end{equation}

The image of $\mathcal{I}_1$ consists of all pairs $((c'_1,c'_3,\ldots,c'_{2m+1}),(c'_2,c'_4,\ldots,c'_{2m}))\in X_{n,1}$ such that $c'_i\leq c'_{i+1}$ for all odd $i$ and $c'_i\leq c'_{i+1}+2$ for all even $i$.

The image $\mathcal{R}(W)$ of $F$ consists of all pairs $((c_1,c_3,\ldots,c_{2m+1}),(c_2,c_4,\ldots,c_{2m}))\in X_{n,1}$ such that $c_i\leq c_{i+1}+2$ for all $i$.

The maps $\Phi^W$ and $\Psi^W$ are the same as in case of $C_n$.

The map $\pi^{G_1}$ is given by $\pi^{G_1}(\bar{\lambda})=(\bar{\nu}, \varepsilon')$, $\bar{\lambda} = (\bar{\lambda}_1
\leq \bar{\lambda}_2\leq \ldots \leq \bar{\lambda}_{2m+1})\in \mathcal{Q}_{2n+1}$, where
$$
\bar{\nu}_i=\left\{
        \begin{array}{ll}
          \bar{\lambda}_i-1 & \hbox{if $\bar{\lambda}_i$ and $i$ are odd and $\bar{\lambda}_{i-1}<\bar{\lambda}_i$;} \\
          \bar{\lambda}_i+1 & \hbox{if $\bar{\lambda}_i$ is odd, $i$ is even and $\bar{\lambda}_{i}<\bar{\lambda}_{i+1}$;}  \\
          \bar{\lambda}_i & \hbox{otherwise},
        \end{array}
      \right.
$$
and
$$
\varepsilon'(k)=\left\{
                 \begin{array}{ll}
                 \omega  & \hbox{if $k$ is odd;} \\
                   0  & \hbox{if $k$ is even, there exists even $\bar{\lambda}_i=k$ with even $i$ such that $\bar{\lambda}_{i-1}<\bar{\lambda}_i$;} \\
1 & \hbox{otherwise.}
                 \end{array}
               \right.
$$

The map ${\pi}^{G_{1}}$ is injective and its image consists of pairs $(\bar{\nu}, \varepsilon)\in {\mathcal{T}}_{2n}^2$ such that $\varepsilon(k)\neq 0$ if $\bar{\nu}_k^*$ is odd and for each even $i$ such that $\bar{\nu}_i^*$ is even we have $\bar{\nu}_{i-1}^*=\bar{\nu}_i^*$, i.e. $i-1$ does not appear in the partition $\bar{\nu}$. Here $\bar{\nu}_1^*\geq\bar{\nu}_2^*\geq \ldots \geq \bar{\nu}_{m}^*$ is the partition dual to $\bar{\nu}$.

To describe $\phi_{G_{\bf k}}$ for $G_{\bf k}={\rm SO}(V)$ we choose a sufficiently large $m \in \mathbb{N}$. Let $g \in G_{\bf k}$. For any $x \in {\bf k}^*$ let $V_x$ be the
generalized $x$--eigenspace of $g:V \rightarrow V$. For any $x \in {\bf k}^*$ such that $x^2\neq 1$ let $\lambda^x_1
\geq \lambda^x_2
\geq \ldots \geq \lambda^x_{2m+1}$ be the sequence
in $\mathbb{N}$ whose terms are the sizes of the Jordan blocks of $x^{-1}g: V_x \rightarrow V_x$.

For any $x \in {\bf k}^*$ with $x^2= 1$ let $\lambda^x_1
\geq \lambda^x_2\geq \ldots \geq \lambda^x_{2m+1}$ be the sequence in $\mathbb{N}$, where
$((\lambda^x_1
\geq \lambda^x_3\geq \ldots \geq \lambda^x_{2m+1}),(\lambda^x_2
\geq \lambda^x_4\geq \ldots \geq \lambda^x_{2m}))$ is the pair of partitions
such that the corresponding irreducible representation of the Weyl
group of type $B_{({\rm dim}~V_x-1)/2}$ (if $x\neq -1$ or $p=2$) or $D_{{\rm dim}~V_x/2}$ (if $x=-1$ or $p\neq 2$) is the Springer representation attached to the unipotent
element $x^{-1}g \in {\rm SO}(V_x)$ and to the trivial local data.

Let $\bar{\lambda}(g)$ be
the partition $\bar{\lambda}_1(g) \geq \bar{\lambda}_2(g) \geq \ldots \geq \bar{\lambda}_{2m+1}(g)$ \index[not]{l@$\bar{\lambda}(g)$} given by $\bar{\lambda}_j(g) = \sum_{x} \lambda^x_j$, where $x$ runs over a set of representatives for the orbits of the
involution $a\mapsto a^{-1}$ of ${\bf k}^*$.
Now $\phi_{G_{\bf k}}(g)$ is the pair of partitions $((\bar{\lambda}_1(g)
\geq \bar{\lambda}_3(g)\geq \ldots \geq \bar{\lambda}_{2m+1}(g)),(\bar{\lambda}_2(g)
\geq \bar{\lambda}_4(g)\geq \ldots \geq \bar{\lambda}_{2m}(g)))$.

If $g$ is any element in the stratum $G_{(\bar{\lambda},\varepsilon)}$ \index[not]{G@$G_{(\bar{\lambda},\varepsilon)}$} corresponding to a pair $(\bar{\lambda},\varepsilon)\in \mathcal{T}_{2n}^2$, $\bar{\lambda}=(\bar{\lambda}_1
\geq \bar{\lambda}_2\geq \ldots \geq \bar{\lambda}_m)$ then the dimension of the centralizer of $g$ in $G_{\bf k}$ is given by formula (\ref{dimsp}),
\begin{equation}\label{dimsoodd}
{\rm dim}~Z_{G_{\bf k}}(g)=n+\sum_{i=1}^m(i-1)\bar{\lambda}_i+\frac{1}{2}|\{i:\bar{\lambda}_i ~\hbox{is odd}\}|+|\{i:\bar{\lambda}_i ~\hbox{is even and}~\varepsilon(\bar{\lambda}_i)=0\}|.
\end{equation}


\subsection*{$\bf D_n$} \index[not]{D@$D_n$}

$G_{\bf k}$ is of type ${\rm SO}(V)$, \index[not]{S@${\rm SO}(V)$} where $V$ is a vector space of dimension $2n$, $n \geq 3$ over an algebraically closed field ${\bf k}$ of characteristic exponent $p$. If $p\neq 2$ $V$ is equipped with a non--degenerate symmetric bilinear form. If $p=2$ $V$ is equipped with a non--degenerate bilinear form \index{form!non--degenerate symmetric bilinear on a vector space} $B(~\cdot~,~\cdot~)$ \index[not]{B@$B(~\cdot~,~\cdot~)$} and a non--zero quadratic form \index{form!quadratic on a vector space} $Q(~\cdot~)$ \index[not]{Q@$Q(~\cdot~)$} such that
$$
B(x,y)=Q(x+y)-Q(x)-Q(y),~x,y\in V.
$$
We remind that in this case ${\rm SO}(V)$ is the connected component containing the identity of the group of linear automorphisms of $V$ preserving the quadratic, and hence the bilinear, form.

$W$ is the group of even permutations of the set $E=\{\varepsilon_1,\ldots, \varepsilon_n,-\varepsilon_1,\ldots,-\varepsilon_n\}$ which also commute with the involution  $\varepsilon_i\mapsto -\varepsilon_i$. $W$ can be regarded as a subgroup in the Weyl group $W'$ of type $C_n$.

Let $\underline{\widetilde{W}}$ \index[not]{W@$\underline{\widetilde{W}}$} be the set of $W'$--conjugacy classes in $W$.
The elements of $\underline{\widetilde{W}}$ are parametrized by the pairs of partitions $(\bar{\lambda},\bar{\mu})$ satisfying the following conditions.
\begin{itemize}
\item
The parts of $\bar{\lambda}$ are even (for any $w\in \mathcal{C}\in \underline{\widetilde{W}}$ they are the numbers of the elements in the negative orbits $X$, $X=-X$,  in $E$ for the action of the group $\left< w \right>$ generated by $w$);
\item
The number of parts of $\bar{\lambda}$ is even;
\item
$\bar{\mu}$ consists of pairs of equal parts (they are the numbers of the elements in the positive $\left< w \right>$--orbits $X$ in $E$; these orbits appear in pairs $X,-X$, $X\neq -X$);
\item 
$\sum \bar{\lambda}_i +\sum \bar{\mu}_j=2n$.
\end{itemize}
We denote this set of pairs of partitions by $\mathcal{A}^0_{2n}$. \index[not]{A@$\mathcal{A}^0_{2n}$}

To each pair $(-,\bar{\mu})$, where all parts of $\bar{\mu}$ are even, there correspond two conjugacy classes in $W$. To any other element of $\mathcal{A}^0_{2n}$ there corresponds a unique conjugacy class in $W$.

An element of $\underline{\widetilde{W}}$ which corresponds to a pair $(\bar{\lambda},\bar{\mu})$, $\bar{\lambda}=(\bar{\lambda}_1
\geq \bar{\lambda}_2\geq \ldots \geq \bar{\lambda}_m)$ and $\bar{\mu}=(\bar{\mu}_1=\bar{\mu}_2\geq \ldots \geq \bar{\mu}_{2k-1}=\bar{\mu}_{2k})$ is the class represented by the sum of the blocks in the following diagram (we use the notation of \cite{C}, Section 7) \index{blocks!in a Weyl group}
\begin{equation}\label{sdn}
A_{\bar{\mu}_1-1}+A_{\bar{\mu}_3-1}+\ldots +A_{\bar{\mu}_{2k-1}-1}+D_{\frac{\bar{\lambda}_1+\bar{\lambda}_2}{2}}(a_{\frac{\bar{\lambda}_2}{2}-1})+
D_{\frac{\bar{\lambda}_3+\bar{\lambda}_4}{2}}(a_{\frac{\bar{\lambda}_4}{2}-1})+\ldots +D_{\frac{\bar{\lambda}_{m-1}+\bar{\lambda}_{m}}{2}}(a_{\frac{\bar{\lambda}_{m}}{2}-1}).
\end{equation}

Consider now the case $p\neq 2$. Let $G'_{\bf k}$ \index[not]{G@$G'_{\bf k}$} be the extension of $G_{\bf k}$ by the Dynkin graph automorphism of order $2$. Then $G'_{\bf k}$ is of type ${\rm O}(V)$. \index[not]{O@${\rm O}(V)$} Denote by $\underline{\widetilde{\mathcal{N}}}(G_{\bf k})$ \index[not]{N@$\underline{\widetilde{\mathcal{N}}}(G_{\bf k})$} the set of the unipotent classes of $G'_{\bf k}$. Note that they are all contained in $G_{\bf k}$.
The elements of $\underline{\widetilde{\mathcal{N}}}(G_{\bf k})$ are parametrized by partitions $\bar{\lambda}$ of $2n$ for which $l_j(\bar{\lambda})$ is even for even $j$. Note that the number of parts of such partitions is even. We denote this set of partitions by $\mathcal{Q}_{2n}$. \index[not]{Q@$\mathcal{Q}_{2n}$} In the case when $G_{\bf k}={\rm SO}(V)$ the parts of $\bar{\lambda}$ are just the sizes of the Jordan blocks in $V$ of the unipotent elements from the conjugacy class corresponding to $\bar{\lambda}$. If $\bar{\lambda}$ has only even parts then $\bar{\lambda}$ corresponds to two unipotent classes in $G_{\bf k}$ of the same dimension. In all other cases there is a unique unipotent class in $G_{\bf k}$ which corresponds to $\bar{\lambda}$.

One has $\widehat{\underline{\mathcal{N}}}(W)=\underline{\mathcal{N}}(G_2)$, and $G_2$ is of type ${\rm SO}(V)$.

Let $G'_2$ \index[not]{G@$G'_2$} be the extension of $G_2$ by the Dynkin graph automorphism of order $2$. Then $G'_2$ is of type ${\rm O}(V)$. Denote by $\underline{\widetilde{\mathcal{N}}}(G_2)$ the set of the unipotent classes of $G'_2$ contained in $G_2$. Since the bilinear form $B(~\cdot~,~\cdot~)$ is also alternating \index{form!alternating, on a vector space} in characteristic 2 there is a natural injective homomorphism from ${\rm O}(V)$ to ${\rm Sp}(V)$, ${\rm dim}~V=2n$, and $\underline{\widetilde{\mathcal{N}}}(G_2) \simeq \widetilde{\mathcal{T}}_{2n}^2$, where $\widetilde{\mathcal{T}}_{2n}^2$ \index[not]{T@$\widetilde{\mathcal{T}}_{2n}^2$} is the set of the elements $(\bar{\lambda},\varepsilon)\in \mathcal{T}_{2n}^2$ such that $\bar{\lambda}$ has an even number of parts (see I.2.6 in \cite{Spal}).

Let $\widehat{\widetilde{W}}$ \index[not]{W@$\widehat{\widetilde{W}}$} be the set of the orbits of irreducible characters of $W$ under the action of $W'$.
Elements of $\widehat{\widetilde{W}}$ are parametrized by unordered pairs of partitions $(\bar{\alpha},\bar{\beta})$ written in non--decreasing order, $\bar{\alpha}_1\leq \bar{\alpha}_2\leq \ldots \leq \bar{\alpha}_{\bar{\tau}(\bar{\alpha})}$, $\bar{\beta}_1\leq \bar{\beta}_2 \leq \ldots \leq \bar{\beta}_{\bar{\tau}(\bar{\beta})}$, and such that $\sum \bar{\alpha}_i+\sum \bar{\beta}_i=n$. By adding zeroes we can assume that the length of $\bar{\alpha}$ is equal to the length of $\bar{\beta}$. The set of such pairs is denoted by $Y_{n,0}$. \index[not]{Y@$Y_{n,0}$}

Instead of the maps in (\ref{diag}) we shall describe the following maps
\begin{equation}\label{diagso2n}
\begin{array}{ccccccc}
   &  & \widetilde{\mathcal{X}}^1(W) & \stackrel{\widetilde{\mathcal{I}}_1}{\longleftarrow} & \underline{\widetilde{\mathcal{N}}}(G_{1}) &  & \\
   &  & \downarrow ~{\iota} &  & \downarrow ~\widetilde{\pi}^{G_{1}} &  & \\
  G_{\bf k} & \stackrel{{\widetilde{\phi}_{G_{\bf k}}}}{\longrightarrow} & {\widetilde{\mathcal{R}}(W)} & \stackrel{\widetilde{F}}{\longleftarrow} & {\underline{\widetilde{\mathcal{N}}}}(G_2) & {\hbox{\raise-1mm\hbox{${\textstyle \stackrel{\widetilde{\Phi}^W}{\longleftarrow}}\atop {\hbox{${\textstyle \longrightarrow}\atop {\scriptstyle \widetilde{\Psi}^W}$}}$}}} & \underline{\widetilde{W}},
\end{array}
\end{equation}
where $\widetilde{\mathcal{I}}_1$ \index[not]{I@$\widetilde{\mathcal{I}}_1$} and $\widetilde{F}$ \index[not]{F@$\widetilde{F}$} are induced by the restrictions of the maps $\mathcal{I}_1$ and $F$ for $G'_1$, $G'_2$ to $\underline{\widetilde{\mathcal{N}}}(G_1)$, ${\underline{\widetilde{\mathcal{N}}}}(G_2)$, respectively, $\widetilde{\mathcal{X}}^1(W)$ \index[not]{X@$\widetilde{\mathcal{X}}^1(W)$} and ${\widetilde{\mathcal{R}}(W)}$ \index[not]{R@${\widetilde{\mathcal{R}}(W)}$} are their images, ${\widetilde{\phi}_{G_{\bf k}}}, \widetilde{\Psi}^W,\widetilde{\Phi}^W$ \index[not]{f@${\widetilde{\phi}_{G_{\bf k}}}$} \index[not]{P@$\widetilde{\Psi}^W$} \index[not]{F@$\widetilde{\Phi}^W$} and  $\widetilde{\pi}^{G_1}$ \index[not]{p@$\widetilde{\pi}^{G_1}$} are also induced by the corresponding maps for $G'_1$, $G'_2$ and $W'$.

The map $\widetilde{\mathcal{I}}_1$ is defined by the same algorithm as in the case of $B_n$ (see (\ref{algf})). The image of $\widetilde{\mathcal{I}}_1$ consists of all pairs $((c'_1,c'_3,\ldots,c'_{2m+1}),(c'_2,c'_4,\ldots,c'_{2m}))\in Y_{n,0}$ such that $c'_i\leq c'_{i+1}$ for all odd $i$ and $c'_i\leq c'_{i+1}+2$ for all even $i$.

If $(\bar{\lambda},\varepsilon)\in \widetilde{\mathcal{T}}_{2n}^2$, $\bar{\lambda} = (\bar{\lambda}_1
\leq \bar{\lambda}_2\leq \ldots \leq \bar{\lambda}_{2m})$ and $\widetilde{F}(\bar{\lambda},\varepsilon)=((c_1,c_3,\ldots,c_{2m-1}),(c_2,c_4,\ldots,c_{2m}))$ then the parts $c_i$ are defined by induction starting from $c_1$,
$$
  \begin{array}{ll}
    c_i=\frac{\bar{\lambda}_i-2}{2}+2(i-1)-4\left[\frac{i-1}{2}\right] & \hbox{if $\bar{\lambda}_i$ is even, $\varepsilon(\bar{\lambda}_i)=1$  and $c_{i-1}$ is already defined;} \\
    c_i=\frac{\bar{\lambda}_i-1}{2} +2(i-1)-4\left[\frac{i-1}{2}\right] & \hbox{if $\bar{\lambda}_i=\bar{\lambda}_{i+1}$ is odd and $c_{i-1}$ is already defined;} \\
    c_{i+1}=\frac{\bar{\lambda}_i-3}{2}+2i-4\left[\frac{i}{2}\right]  & \hbox{if $\bar{\lambda}_i=\bar{\lambda}_{i+1}$ is odd and $c_{i}$ is already defined;}
\\
    c_i=\frac{\bar{\lambda}_i}{2} +2(i-1)-4\left[\frac{i-1}{2}\right] & \hbox{if $\bar{\lambda}_i=\bar{\lambda}_{i+1}$ is even, $\varepsilon(\bar{\lambda}_i)=0$ and $c_{i-1}$ is already defined;} \\
    c_{i+1}=\frac{\bar{\lambda}_i}{2} +2(i-1)-4\left[\frac{i-1}{2}\right] & \hbox{if $\bar{\lambda}_i=\bar{\lambda}_{i+1}$ is even, $\varepsilon(\bar{\lambda}_i)=0$ and $c_{i}$ is already defined.}
  \end{array}
$$

The image ${\widetilde{\mathcal{R}}(W)}$ of $\widetilde{F}$ consists of all pairs $((c_1,c_3,\ldots,c_{2m+1}),(c_2,c_4,\ldots,c_{2m}))\in Y_{n,0}$ such that $c_i\leq c_{i+1}$ for all odd $i$ and $c_i\leq c_{i+1}+4$ for all even $i$.

The maps $\widetilde{\Phi}^W$ and $\widetilde{\Psi}^W$ are defined by the same formulas as in the case of $C_n$.

The map ${\widetilde{\pi}}^{G_{1}}$ is given by $\widetilde{\pi}^{G_{1}}(\bar{\lambda})=(\bar{\nu}, \varepsilon')$, $\bar{\lambda} = (\bar{\lambda}_1
\leq \bar{\lambda}_2\leq \ldots \leq \bar{\lambda}_{2m})\in \mathcal{Q}_{2n}$, where
$$
\bar{\nu}_i=\left\{
        \begin{array}{ll}
          \bar{\lambda}_i-1 & \hbox{if $\bar{\lambda}_i$ is odd, $i$ is even and $\bar{\lambda}_{i-1}<\bar{\lambda}_i$;} \\
          \bar{\lambda}_i+1 & \hbox{if $\bar{\lambda}_i$ and $i$ are odd, and $\bar{\lambda}_{i}<\bar{\lambda}_{i+1}$;}  \\
          \bar{\lambda}_i & \hbox{otherwise},
        \end{array}
      \right.
$$
and
$$
\varepsilon'(k)=\left\{
                 \begin{array}{ll}
                 \omega  & \hbox{if $k$ is odd;} \\
                   0  & \hbox{if $k$ is even, there exists even $\bar{\lambda}_i=k$ with odd $i$ such that $\bar{\lambda}_{i-1}<\bar{\lambda}_i$;} \\
1 & \hbox{otherwise.}
                 \end{array}
               \right.
$$

The map ${\widetilde{\pi}}^{G_{1}}$ is injective and its image consists of the pairs $(\bar{\nu}, \varepsilon)\in \widetilde{\mathcal{T}}_{2n}^2$ such that $\varepsilon(k)\neq 0$ if $\bar{\nu}_k^*$ is odd and for each even $i$ such that $\bar{\nu}_i^*$ is even we have $\bar{\nu}_{i-1}^*=\bar{\nu}_i^*$, i.e. $i-1$ does not appear in the partition $\bar{\nu}$. Here $\bar{\nu}_1^*\geq\bar{\nu}_2^*\geq \ldots \geq \bar{\nu}_{m}^*$ is the partition dual to $\bar{\nu}$.

To describe ${\widetilde{\phi}_{G_{\bf k}}}$ for $G_{\bf k}={\rm SO}(V)$ we choose a sufficiently large $m \in \mathbb{N}$. Let $g \in G_{\bf k}$. For any $x \in {\bf k}^*$ let $V_x$ be the
generalized $x$--eigenspace of $g:V \rightarrow V$. For any $x \in {\bf k}^*$ such that $x^2\neq 1$ let $\lambda^x_1
\geq \lambda^x_2
\geq \ldots \geq \lambda^x_{2m}$ be the sequence
in $\mathbb{N}$ whose terms are the sizes of the Jordan blocks of $x^{-1}g: V_x \rightarrow V_x$.

For any $x \in {\bf k}^*$ with $x^2= 1$ let $\lambda^x_1
\geq \lambda^x_2\geq \ldots \geq \lambda^x_{2m}$ be the sequence in $\mathbb{N}$, where $((\lambda^x_1
\geq \lambda^x_3\geq \ldots \geq \lambda^x_{2m-1}),(\lambda^x_2
\geq \lambda^x_4\geq \ldots \geq \lambda^x_{2m}))$ is the pair of partitions
such that the corresponding irreducible representation of the Weyl
group of type $D_{{\rm dim}~V_x/2}$ is the Springer representation attached to the unipotent
element $x^{-1}g \in {\rm SO}(V_x)$ and to the trivial local data.

Let $\bar{\lambda}(g)$ be
the partition $\bar{\lambda}_1(g) \geq \bar{\lambda}_2(g) \geq \ldots \geq \bar{\lambda}_{2m+1}(g)$ \index[not]{l@$\bar{\lambda}(g)$} given by $\bar{\lambda}_j(g) = \sum_{x} \lambda^x_j$, where $x$ runs over a set of representatives for the orbits of the
involution $a\mapsto a^{-1}$ of ${\bf k}^*$.
Now $\widetilde{\phi}_{G_{\bf k}}(g)$ is the pair of partitions $((\bar{\lambda}_1(g)
\geq \bar{\lambda}_3(g)\geq \ldots \geq \bar{\lambda}_{2m-1}(g)),(\bar{\lambda}_2(g)
\geq \bar{\lambda}_4(g)\geq \ldots \geq \bar{\lambda}_{2m}(g)))$.

The preimage $\widetilde{\phi}_{G_{\bf k}}^{-1}(\bar{\lambda},\bar{\mu})$ is a stratum in $G_{\bf k}$ in all cases except for the one when the pair $(\bar{\lambda},\bar{\mu})$ is of the form $((\bar{\lambda}_1
\geq \bar{\lambda}_3\geq \ldots \geq \bar{\lambda}_{2m-1}),(\bar{\lambda}_1
\geq \bar{\lambda}_3\geq \ldots \geq \bar{\lambda}_{2m-1}))$. In that case $\widetilde{\phi}_{G_{\bf k}}^{-1}(\bar{\lambda},\bar{\mu})$ is a union of two strata, and the conjugacy classes in each of them have the same dimension.

If $g$ is any element in the stratum $G_{(\bar{\lambda},\varepsilon)}$ \index[not]{G@$G_{(\bar{\lambda},\varepsilon)}$} corresponding to a pair $(\bar{\lambda},\varepsilon)\in \widetilde{\mathcal{T}}_{2n}^2$, $\bar{\lambda}=(\bar{\lambda}_1
\geq \bar{\lambda}_2\geq \ldots \geq \bar{\lambda}_m)$ then the dimension of the centralizer of $g$ in $G_{\bf k}$ is given by the following formula
\begin{equation}\label{dimsoev}
{\rm dim}~Z_{G_{\bf k}}(g)=n+\sum_{i=1}^m(i-1)\bar{\lambda}_i-\frac{1}{2}|\{i:\bar{\lambda}_i ~\hbox{is odd}\}|-|\{i:\bar{\lambda}_i ~\hbox{is even and}~\varepsilon(\bar{\lambda}_i)=1\}|.
\end{equation}


\setcounter{equation}{0}
\setcounter{theorem}{0}

\section{The strict transversality condition}\label{stt}

Recall that the definition of $\Delta_+^s$, and hence of $\Sigma_{{\bf k}, s}$, depends on the choice of the ordering of the terms in decomposition (\ref{hdec}).
In this section for every conjugacy class $\mathcal{C} \in C(W)$ we define a variety $\Sigma_{{\bf k}, s}$, $s\in \mathcal{C}$ such that every conjugacy class $\mathcal{O}\in G_\mathcal{C}$ intersects $\Sigma_{{\bf k}, s}$ and
\begin{equation}\label{str1}
{\rm dim}~\mathcal{O}={\rm codim}~\Sigma_{{\bf k}, s}.
\end{equation}
 
It turns out that in order to fulfill condition (\ref{str1}) the subspaces $\h_i$ in (\ref{hdec}) should be ordered in such a way that $\h_0\subset \h_{\mathbb{R}}$ is the subspace fixed by the action of $s$, and if $\h_i=\h_\lambda^k$, $\h_j=\h_\mu^l$ and $0\leq\lambda <\mu< 1$ then $i<j$, where $\lambda$ and $\mu$ are eigenvalues of the corresponding matrix $I-O$ for $s$. In the case of exceptional root systems this is verified using a computer code, and in the case of classical root systems this is confirmed by explicit computation based on a technical lemma. In order to formulate this lemma we recall realizations of classical irreducible root systems.

Let $V$ be a real Euclidean $n$--dimensional vector space with an orthonormal basis $\varepsilon_1,\ldots ,\varepsilon_n$. \index[not]{e@$\varepsilon_i$} The root systems of types $A_{n-1},B_n,C_n$ and $D_n$ can be realized in $V$ as follows.

\subsection*{$\bf A_{n-1}$} The roots are $\varepsilon_i-\varepsilon_j$, $1\leq i,j\leq n$, $i\neq j$, $\h_\mathbb{R}$ is the hyperplane in $V$ consisting of the points the sum of the  coordinates of which is zero.

\subsection*{$\bf B_n$} The roots are $\pm\varepsilon_i\pm\varepsilon_j$, $1\leq i<j\leq n$, $\pm \varepsilon_i$, $1\leq i\leq n$, $\h_\mathbb{R}=V$.

\subsection*{$\bf C_n$} The roots are $\pm\varepsilon_i\pm\varepsilon_j$, $1\leq i<j\leq n$, $\pm 2\varepsilon_i$, $1\leq i\leq n$, $\h_\mathbb{R}=V$.

\subsection*{$\bf D_n$} The roots are $\pm\varepsilon_i\pm\varepsilon_j$, $1\leq i<j\leq n$, $\h_\mathbb{R}=V$.

In all these cases the corresponding Weyl group $W$ is a subgroup of the Weyl group of type $C_n$ acting on the elements of the basis $\varepsilon_1,\ldots ,\varepsilon_n$ by permuting them and by changing the sign of an arbitrary subset of them.

Now we formulate the main lemma.
\begin{lemma}\label{mainl*}
Let $s$ be an element of the Weyl group of type $C_n$ operating on the set $E=\{\varepsilon_1,\ldots, \varepsilon_n,-\varepsilon_1,\ldots,-\varepsilon_n\}$ as indicated in Section \ref{luspart}, where $\varepsilon_1,\ldots ,\varepsilon_n$ is the basis of $V$ introduced above. Assume that $s$ has either only one nontrivial cycle of length $k/2$ ($k$ is even), which is negative, or only one nontrivial cycle of length $k$, which is positive, $1<k\leq n$. Let $\Delta$ be a root system of type $A_{n-1},B_n,C_n$ or $D_n$ realized in $V$ as above.

(i) If $s$ has only one nontrivial cycle of length $k/2$, which is negative, then $k$ is even, the spectrum of $s$ in the complexification $V_\mathbb{C}$ \index[not]{V@$V_\mathbb{C}$} of $V$ is $\epsilon_r=\exp(\frac{2\pi i (k-2r+1)}{k})$, $r=1,\ldots,k/2$, \index[not]{e@$\epsilon_r$} and possibly $\epsilon_0=1$, all eigenvalues are simple except for possibly $1$.

(ii) If $s$ only has one nontrivial cycle of length $k$, which is positive, then the spectrum of $s$ in the complexification of $V$ is $\epsilon_r=\exp(\frac{2\pi i (k-r)}{k})$, $r=1,\ldots,k-1$, and $\epsilon_0=1$, all eigenvalues are simple except for possibly $1$.

In both cases we denote by $V_r$ \index[not]{V@$V_r$} the invariant subspace in $V$ which corresponds to $\epsilon_r=\exp(\frac{2\pi i ([k/2]+1-r)}{k})$, $r=1,\ldots,\left[\frac{k}{2} \right]$ or $\epsilon_0=1$ in the case of a positive nontrivial cycle and to $\epsilon_r=\exp(\frac{2\pi i (2\left[\frac{k/2+1}{2} \right]+1-2r)}{k})$, $r=1,\ldots,\left[\frac{k/2+1}{2} \right]$ or $\epsilon_0=1$ in the case of a negative cycle. For $r\neq 0$ the space $V_r$ is spanned by the real and the imaginary parts of a nonzero eigenvector of $s$ in $V_\mathbb{C}$ corresponding to $\epsilon_r$, and $V_0$ is the subspace of fixed points of $s$ in $V$.

$V_r$ is two--dimensional if $\epsilon_r\neq \pm 1$, one--dimensional if $\epsilon_r=-1$ or may have an arbitrary dimension if $\epsilon_r=1$.

Let $\Delta_+^s$ be a system of positive roots associated to $s$ and defined as in Section \ref{background}, where we use the decomposition
\begin{equation}\label{decv}
V=\bigoplus_{r} V_r
\end{equation}
as (\ref{hdec}) in the definition of $\Delta_+^s$. Denote by $\Delta_{r}\subset \Delta$ the corresponding subsets of roots defined as in (\ref{di}).

Let $\Delta_0^s\subset \Delta$  be the root subsystem fixed by the action of $s$ and $\underline{l}(s)$ \index[not]{l@$\underline{l}(s)$} the number of positive roots which become negative under the action of $s$.

(iii) If $s$ has only one nontrivial cycle of length $k$, which is positive, we have
\begin{enumerate}
\item if $\Delta=A_{n-1}$ then $\Delta_0^s=A_{n-k-1}$, $\underline{l}(s)=2n-k-1$;

\item if $\Delta=B_{n}$ (resp. $C_n$) then $\Delta_0^s=B_{n-k}$ (resp. $C_{n-k}$), $\underline{l}(s)=4n-2k$ for odd $k$ and $\underline{l}(s)=4n-2k+1$ for even $k$;

\item if $\Delta=D_{n}$ then $\Delta_0^s=D_{n-k}$, $\underline{l}(s)=4n-2k-2$ for odd $k$ and $\underline{l}(s)=4n-2k-1$ for even $k$.
\end{enumerate}

(iv) If $s$ has only one nontrivial cycle of length $\frac{k}{2}$,  which is negative, we have

\begin{enumerate}
\item if $\Delta=B_{n}$ (resp. $C_n$) then $\Delta_0^s=B_{n-k/2}$ (resp. $C_{n-k/2}$), $\underline{l}(s)=2n-k/2$;

\item if $\Delta=D_{n}$ then $\Delta_0^s=D_{n-k/2}$, $\underline{l}(s)=2n-k/2-1$.
\end{enumerate}

(v) If $s$ has only one nontrivial cycle of length $k$, which is positive, $\Delta$ is of type $B_n,C_n$ or $D_n$, and $k$ is even then
$\Delta=\Delta_{k/2}\cup \Delta_{k/2-1}\cup \Delta_0^s$ (disjoint union),
and all roots in $\Delta_{k/2-1}$ are orthogonal to the fixed point subspace for the action of $s$ on $V$.

(vi) In all the other cases $\Delta=\Delta_{i_{\rm max}}\cup \Delta_0^s$ (disjoint union), where $i_{\rm max}$ is the maximal possible index $i$ which appears in decomposition (\ref{decv}).

\end{lemma}

\begin{proof}
The proof is similar in all cases. We only give details in the most complicated case when $s$ has only one nontrivial cycle, which is positive, $\Delta$ is of type $B_n$(resp. $C_n$), and $k$ is even. Without loss of generality one can assume that $s$ corresponds to the cycle of the form
$$
\varepsilon_1\rightarrow\varepsilon_2\rightarrow\varepsilon_4\rightarrow\varepsilon_6\rightarrow \cdots \rightarrow\varepsilon_{k-2}\rightarrow\varepsilon_k\rightarrow\varepsilon_{k-1}\rightarrow\varepsilon_{k-3}\rightarrow \cdots \rightarrow\varepsilon_3\rightarrow\varepsilon_1~(k>2),~ \varepsilon_1\rightarrow\varepsilon_2\rightarrow\varepsilon_1~(k=2).
$$

From this definition one easily sees that $\Delta_0^s=B_{n-k}$(resp. $C_{n-k}$)$=\Delta\cap V'$, where $V'\subset V$ is the subspace spanned by $\varepsilon_{k+1},\ldots,\varepsilon_n$. Computing the eigenvalues of $s$ in $V_\mathbb{C}$ is a standard exercise in linear algebra. The eigenvalues are expressed in terms of the exponents of the root system of type $A_{k-1}$ (see \cite{Car}, Ch. 10).

The invariant subspace $V_r$ is spanned by the real and the imaginary parts of a nonzero eigenvector of $s$ in $V_\mathbb{C}$ corresponding to the eigenvalue $\epsilon_r$.
If $\epsilon_r\neq \pm 1$ then $V_r$ is two--dimensional, and for $\epsilon_r=-1$ $V_r$ is one--dimensional. In the former case $V_r$ will be regarded as the real form of a complex plane with the orthonormal basis $1,i$. Under this identification the orthogonal projection operator onto $V_r$ acts on the basic vectors $\varepsilon_j$ as follows
\begin{equation}\label{pr}
\varepsilon_{2j+1}\mapsto \sqrt{\frac{2}{k}}\epsilon_r^j,~j=0,\ldots,\frac{k}{2}-1,\varepsilon_{2j}\mapsto \sqrt{\frac{2}{k}}\epsilon_r^{-j},~j=1,\ldots,\frac{k}{2}.
\end{equation}
Consider the case when $k>2$; the case $k=2$ can be analyzed in a similar way.

To compute $\underline{l}(s)$ using the definition of $\Delta_+^s$ given in Section \ref{background} one should first look at all roots which have nonzero projections onto $V_{k/2}$ on which $s$ acts by rotation by the angle $\frac{2\pi }{k}$.

From (\ref{pr}) we deduce that the roots which are not fixed by $s$ and have zero orthogonal projections onto $V_{k/2}$ are $\pm(\varepsilon_j+\varepsilon_{k-j+1})$, $j=1,\ldots \frac{k}{2}$. The number of these roots is equal to $k$, and they all have nonzero orthogonal projections onto $V_{k/2-1}$. From (\ref{pr}) we also obtain that all the other roots which are not fixed by $s$ have nonzero orthogonal projections onto $V_{k/2}$, hence $|\Delta_{k/2-1}|=k$. The number of roots fixed by $s$ is $2(n-k)^2$ since it is equal to the number of roots in $\Delta_0=\Delta_0^s=B_{n-k}$(resp. $C_{n-k}$). Hence $\Delta=\Delta_{k/2}\cup \Delta_{k/2-1}\cup \Delta_0$ (disjoint union), the number of roots in $\Delta_{k/2}$ is   $|\Delta|-|\Delta_0|-|\Delta_{k/2-1}|=2n^2-2(n-k)^2-k=4nk-2k^2-k$, $|\Delta_{k/2}|=4nk-2k^2-k$.

Now using the symmetry of the root system $\Delta$ as a subset of $V$ and the fact that $s$ acts as rotation by the angles $\frac{2\pi }{k}$ and $\frac{4\pi }{k}$ in $V_{k/2}$ and $V_{k/2-1}$, respectively, we deduce that the number of positive roots in $\Delta_{k/2}$ (resp. $\Delta_{k/2-1}$) which become negative under the action of $s$ is equal to the number of roots in $\Delta_{k/2}$ (resp. $\Delta_{k/2-1}$) divided by the order of $s$ in $V_{k/2}$ (resp. $V_{k/2-1}$). Therefore
$$
\underline{l}(s)=\frac{|\Delta_{k/2}|}{k}+\frac{|\Delta_{k/2-1}|}{k/2}=\frac{4nk-2k^2-k}{k}+\frac{k}{k/2}=
4n-2k+1.
$$
This completes the proof in the considered case.
\end{proof}

Now we are in a position to prove the main statement of this section.
\begin{theorem}\label{mainth}
Let $G_{\bf k}$ be a connected semisimple algebraic group over an algebraically closed filed $\bf k$, and $\mathcal{O}\in \widehat{\underline{\mathcal{N}}}(W)$. Let $H_{\bf k}$ be a maximal torus of $G_{\bf k}$, $W$ the Weyl group of the pair $(G_{\bf k},H_{\bf k})$, and $s\in W$ an element from the conjugacy class $\Psi^W(\mathcal{O})$. Let $\Delta$ be the root system of the pair $(G_{\bf k},H_{\bf k})$ and $\Delta_+^s$ a system of positive roots in $\Delta$ associated to $s$ and defined in Section \ref{background} with the help of decomposition (\ref{hdec}), where the subspaces $\h_i$ are ordered in such a way that $\h_0$ is the vector subspace of $\h_{\mathbb{R}}$ fixed by the action of $s$, and if $\h_i=\h_\lambda^k$, $\h_j=\h_\mu^l$ and $0\leq\lambda <\mu< 1$ then $i<j$. In the case of exceptional root systems we assume, in addition, that $\Delta_+^s$ is chosen as in the tables in Appendix 2, so that $s=s^1s^2$ is defined by the data from columns three and four in the tables in Appendix 2. Then all conjugacy classes in the stratum $G_\mathcal{O}:=\phi_{G_{\bf k}}^{-1}(F(\mathcal{O}))$ \index[not]{G@$G_\mathcal{O}$} intersect the corresponding variety $\Sigma_{{\bf k}, s}$ at some points of the subvariety $\dot{s}H_{\bf k}^0N_{{\bf k}, s}$, where $H_{\bf k}^0\subset H_{\bf k}$ is the identity component of the centralizer of $\dot{s}$  in $H_{\bf k}$. Moreover, if $\mathcal{O}\in \underline{\mathcal{N}}(G_p)\subset \widehat{\underline{\mathcal{N}}}(W)$ for some $p$, then for any $g\in G_\mathcal{O}$
\begin{equation}\label{dimid}
{\rm dim}~Z_{G_{\bf k}}(g)={\rm codim}_{G_p}~\mathcal{O}={\rm dim}~\Sigma_{{\bf k}, s}.
\end{equation}
\end{theorem}

\begin{proof}
Firstly we establish identity (\ref{dimid}). Note that the first identity in (\ref{dimid}) is equivalent to identity (\ref{od}). 
We shall split the proof of the second identity into several lemmas.

First we compute the dimension of the slice $\Sigma_{{\bf k}, s}$, $s\in \Psi^W(\mathcal{O})$ and justify that for any $g\in G_\mathcal{O}$ the second equality in (\ref{dimid}) holds.
\begin{lemma}\label{dimsect1}
Assume that the conditions of Theorem \ref{mainth} are satisfied.  Then
for any $g\in G_\mathcal{O}$, where $\mathcal{O}\in \underline{\mathcal{N}}(G_p)\subset \widehat{\underline{\mathcal{N}}}(W)$ for some $p$, the equalities in (\ref{dimid}) hold, i.e.
$$
{\rm dim}~Z_{G_{\bf k}}(g)={\rm dim}~\Sigma_{{\bf k}, s}={\rm codim}_{G_p}~\mathcal{O}.
$$
\end{lemma}

\begin{proof}
Observe that by the definition of the slice $\Sigma_{{\bf k}, s}$
$$
{\rm dim}~\Sigma_{{\bf k}, s}=l(s)+|\Delta_0|+{\rm dim}~\h_0,
$$
where $l(s)$ is the length of $s$ with respect to the system of simple roots in $\Delta_+^s$.
Hence to compute ${\rm dim}~\Sigma_{{\bf k}, s}$ we have to find all numbers in the right hand side of the last identity.

Consider the case of classical groups when each Weyl group element is a product of cycles in a permutation group. In this case identity (\ref{dimid}) is proved by a straightforward calculation using Lemma \ref{mainl*}.

Let $G_{\bf k}$ be of type $A_n$ and $s$ a representative in the conjugacy class of the Weyl group which corresponds to a partition $\bar{\lambda}=(\bar{\lambda}_1\geq \bar{\lambda}_2\geq \ldots\geq \bar{\lambda}_m)$. Recall that $s$ is the product of the cycles which correspond to the parts of $\bar{\lambda}$. The particular ordering of the invariant subspaces $\h_i$ in the statement of Theorem \ref{mainth} implies that the length $l(s)$ equal to the number $|\Delta_s^s|$ of positive roots in $\Delta_+^s$ which become negative under the action of $s$ should be computed by successive application of Lemma \ref{mainl*} to the cycles $\hat{s}_i$ \index[not]{s@$\hat{s}_i$} of $s$, which correspond to $\bar{\lambda}_i$ placed in a non--increasing order. We claim that according to this observation one has
\begin{equation}\label{ls}
l(s)= \sum_{\scriptsize \begin{array}{l} k=1 \\ \bar{\lambda}_k>1\end{array}}^m \underline{l}(\hat{s}_k),~\underline{l}(\hat{s}_k)=2(n-\sum_{i=1}^{k-1}\bar{\lambda}_i)-\bar{\lambda}_k+1,
\end{equation}
where the first sum in (\ref{ls}) is taken over $k$ for which $\bar{\lambda}_k>1$, and we keep the notation of Lemma \ref{mainl*}.

Indeed, recall that according to formula (\ref{decD}) we have the following disjoint union decomposition
$$
\Delta=\bigcup_{k=0}^{M}{\Delta}_{i_k}
$$
which implies another disjoint union decomposition
$$
\Delta^s_s=\bigcup_{k=1}^{M}({\Delta}_{i_k}\cap \Delta^s_s).
$$
Thus the length $l(s)$ equal to the cardinality of the set $\Delta^s_s$ can be found as the sum of the cardinalities of the sets ${\Delta}_{i_k}\cap \Delta^s_s$,
\begin{equation}\label{ls1}
l(s)=\left| \Delta^s_s \right|=\bigcup_{k=1}^{M}\left| {\Delta}_{i_k}\cap \Delta^s_s \right|.
\end{equation}

We find the cardinalities of the sets ${\Delta}_{i_k}\cap \Delta^s_s$  by successive application of Lemma \ref{mainl*} to the cycles $\hat{s}_i$ of $s$, which correspond to $\bar{\lambda}_i$ placed in a non--increasing order.  

According to our convention for the ordering of the terms in sum (\ref{hdec}) used to define $\Delta_+^s$, the cycle $\hat{s}_1$ of $s$ which corresponds to the maximal part $\bar{\lambda}_1$ is the only cycle of $s$ non-trivially acting on $\h_{i_M}$ by Lemma \ref{mainl*} (vi). Applying part (iii) 1. of Lemma \ref{mainl*} to the cycle $\hat{s}_1$ we obtain $\left|{\Delta}_{i_M}\cap \Delta^s_s\right|=\underline{l}(\hat{s}_1)=2n-\bar{\lambda}_1+1$ and $\Delta_0^{\hat{s}_1}=A_{n-\bar{\lambda}_1}=\Delta\setminus \Delta_{i_M}$ . 

The remaining cycles $\hat{s}_2,\ldots ,\hat{s}_m$ of $s$ corresponding to $\bar{\lambda}_2\geq \ldots\geq \bar{\lambda}_m$ act on $\Delta_0^{\hat{s}_1}$ and leave the set $\Delta_{i_M}$ and its subset $\Delta_s^s\cap \Delta_{i_M}$ invariant elementwise, so we can apply Lemma \ref{mainl*} to $\hat{s}_2$ acting on $\Delta_0^{\hat{s}_1}$ to get $\left|{\Delta}_{i_{M-1}}\cap \Delta^s_s\right|=\underline{l}(\hat{s}_2)=2(n-\bar{\lambda}_1)-\bar{\lambda}_2+1$ and $\Delta_0^{\hat{s}_2}=A_{n-\bar{\lambda}_1-\bar{\lambda}_2}=\Delta\setminus \left(\Delta_{i_M}\cup \Delta_{i_{M-1}}\right)$. Iterating this procedure and using (\ref{ls1}) we obtain (\ref{ls}).

The number of roots fixed by $s$ can be represented in a similar form,
\begin{equation}\label{fx}
 |\Delta_0|= \sum_{\scriptsize \begin{array}{l} k=1 \\ \bar{\lambda}_k=1\end{array}}^m \underline{l}(\hat{s}_k),~\underline{l}(\hat{s}_k)=2(n-\sum_{i=1}^{k-1}\bar{\lambda}_i)-\bar{\lambda}_k+1,
\end{equation}
where the sum in (\ref{fx}) is taken over $k$ for which $\bar{\lambda}_k=1$.

Finally the dimension of the fixed point space $\h_0$ of $s$ in $\h$ is $m-1$, ${\rm dim}~\h_0=m-1$.

Recall now that
\begin{equation}\label{dimsect}
 {\rm dim}~\Sigma_{{\bf k}, s}=l(s)+|\Delta_0|+{\rm dim}~\h_0,
\end{equation}
and hence
$$
{\rm dim}~\Sigma_{{\bf k}, s}=\sum_{k=1}^m \underline{l}(\hat{s}_k)+m-1=\sum_{k=1}^m\left(2(n-\sum_{i=1}^{k-1}\bar{\lambda}_i)-\bar{\lambda}_k+1\right)+m-1.
$$

Exchanging the order of summation and simplifying this expression we obtain that
$$
{\rm dim}~\Sigma_{{\bf k}, s}=n+2\sum_{i=1}^m(i-1)\bar{\lambda}_i
$$
which coincides with (\ref{dimsl}).

The computations of ${\rm dim}~\Sigma_{{\bf k}, s}$ in the cases of $B_n$ and of $C_n$ are similar. If $(\bar{\nu}, \varepsilon)\in \mathcal{T}_{2n}^2$, $\bar{\nu}=(\bar{\nu}_1\geq \bar{\nu}_2\geq \ldots \geq \bar{\nu}_m)$, corresponds to $\mathcal{O}\in \widehat{\underline{\mathcal{N}}}(W)=\underline{\mathcal{N}}(G_2)$ then $\Psi^W(\bar{\nu}, \varepsilon)=(\bar{\lambda},\bar{\mu})\in \mathcal{A}^1_{2n}\simeq\underline{W}$ is defined in Section \ref{luspart}, part $\bf C_n$. $\bar{\lambda}$ consists of even parts $\bar{\nu}_i$ of $\bar{\nu}$ for which $\varepsilon(\bar{\nu}_i)=1$, and $\bar{\mu}$ consists of all odd parts of $\bar{\nu}$ and of even parts $\bar{\nu}_i$ of $\bar{\nu}$ for which $\varepsilon(\bar{\nu}_i)=0$, the last two types of parts appear in pairs of equal parts. Let $s\in W$ be an element from the conjugacy class $\Psi^W(\bar{\nu}, \varepsilon)$. Then  each part $\bar{\lambda}_i$ corresponds to a negative cycle of $s$ of length $\frac{\bar{\lambda}_i}{2}$, and each pair $\bar{\mu}_i=\bar{\mu}_{i+1}$ of equal parts of $\bar{\mu}$ corresponds to a positive cycle of $s$ of length $\bar{\mu}_i$. We order the cycles $\hat{s}_k$ of $s$ associated to the (pairs of equal) parts of the partition $\bar{\nu}$ in a way compatible with a non--increasing ordering of the parts of the partition $\bar{\nu}=(\bar{\nu}_1\geq \bar{\nu}_2\geq \ldots \geq \bar{\nu}_m)$, i.e. if we denote by $\hat{s}_k$ the cycle that corresponds to an even part $\bar{\nu}_k$ of $\bar{\nu}$ for which $\varepsilon(\bar{\nu}_k)=1$ or to a pair $\bar{\nu}_k=\bar{\nu}_{k+1}$ of odd parts of $\bar{\nu}$ or of even parts of $\bar{\nu}$ for which $\varepsilon(\bar{\nu}_k)=0$ then $\hat{s}_k\geq \hat{s}_l$ if $\bar{\nu}_k\geq \bar{\nu}_l$.

Similarly to the case of $A_n$, by the definition of $\Delta_+^s$ and by Lemma \ref{mainl*} applied iteratively to the cycles $\hat{s}_k$ in the order defined in the previous paragraph, the length $l(s)$ of $s$ is the sum of the following terms $\underline{l}(\hat{s}_k)$.

To each even part $\bar{\nu}_k$ of $\bar{\nu}$ for which $\varepsilon(\bar{\nu}_k)=1$ we associate the term $$\underline{l}(\hat{s}_k)=2(n-\sum_{i=1}^{k-1}\frac{\bar{\nu}_i}{2})-\frac{\bar{\nu}_k}{2};$$

to each pair of odd parts $\bar{\nu}_k=\bar{\nu}_{k+1}>1$ we associate the term
$$\underline{l}(\hat{s}_k)=4(n-\sum_{i=1}^{k-1}\frac{\bar{\nu}_i}{2})-2\bar{\nu}_k=
\left(2(n-\sum_{i=1}^{k-1}\frac{\bar{\nu}_i}{2})-\frac{\bar{\nu}_k}{2}\right)+\left(2(n-\sum_{i=1}^{k}\frac{\bar{\nu}_i}{2})-
\frac{\bar{\nu}_{k+1}}{2}\right);$$
note that the sum of these terms over all pairs $\bar{\nu}_k=\bar{\nu}_{k+1}=1$ gives the number $|\Delta_0|$ of the roots fixed by $s$;

to each pair of even parts $\bar{\nu}_k=\bar{\nu}_{k+1}$ for which $\varepsilon(\bar{\nu}_k)=0$ we associate the term
$$\underline{l}(s_k)=4(n-\sum_{i=1}^{k-1}\frac{\bar{\nu}_i}{2})-2\bar{\nu}_k+1=
\left(2(n-\sum_{i=1}^{k-1}\frac{\bar{\nu}_i}{2})-\frac{\bar{\nu}_k}{2}+\frac{1}{2}\right)+\left(2(n-\sum_{i=1}^{k}\frac{\bar{\nu}_i}{2})-
\frac{\bar{\nu}_{k+1}}{2}+\frac{1}{2}\right).$$

The dimension of the fixed point space $\h_0$ of $s$ in $\h_{\mathbb{R}}$ is equal to a half of the sum of the number of all even parts $\bar{\nu}_k$ for which $\varepsilon(\bar{\nu}_k)=0$ and of the number of all odd parts $\bar{\nu}_k$,

\begin{equation}\label{hofix}
{\rm dim}~\h_0= \frac{1}{2}|\{i:\bar{\nu}_i ~\hbox{is odd}\}|+\frac{1}{2}|\{i:\bar{\nu}_i ~\hbox{is even and}~\varepsilon(\bar{\nu}_i)=0\}|.
\end{equation}

Finally substituting all the computed contributions into formula (\ref{dimsect}) we obtain
\begin{eqnarray*}
{\rm dim}~\Sigma_{{\bf k}, s} = \sum_{k=1}^m\left(2(n-\sum_{i=1}^{k-1}\frac{\bar{\nu}_i}{2})-\frac{\bar{\nu}_k}{2}\right)+\frac{1}{2}|\{i:\bar{\nu}_i ~\hbox{is even and}~\varepsilon(\bar{\nu}_i)=0\}|+ \\
    +\frac{1}{2}|\{i:\bar{\nu}_i ~\hbox{is odd}\}|+\frac{1}{2}|\{i:\bar{\nu}_i ~\hbox{is even and}~\varepsilon(\bar{\nu}_i)=0\}|.
\end{eqnarray*}

Exchanging the order of summation and simplifying this expression we obtain that
\begin{equation}\label{dimsp1}
{\rm dim}~\Sigma_{{\bf k}, s}= n+\sum_{i=1}^m(i-1)\bar{\nu}_i+\frac{1}{2}|\{i:\bar{\nu}_i ~\hbox{is odd}\}|+|\{i:\bar{\nu}_i ~\hbox{is even and}~\varepsilon(\bar{\nu}_i)=0\}|
\end{equation}
which coincides with (\ref{dimsp}) or (\ref{dimsoodd}).

In the case of $D_n$ the number ${\rm dim}~\Sigma_{{\bf k}, s} $ can be easily obtained if we observe that the map $\widetilde{\Psi}^W$ is defined by the same formula as $\Psi^W$ in case of $C_n$. In case when $\widetilde{\Psi}^W(\bar{\nu},\varepsilon)=(-,\bar{\mu})$, where all parts of $\bar{\mu}$ are even, there are two conjugacy classes in $W$ which correspond to $\widetilde{\Psi}^W(\bar{\nu},\varepsilon)$. However, the numbers $l(s)$, $|\Delta_0|$ and ${\rm dim}~\h_0$ are the same in both cases. They only depend on $\widetilde{\Psi}^W(\bar{\nu},\varepsilon)$ in all cases. Let $s\in W$ be a representative from the conjugacy class $\widetilde{\Psi}^W(\bar{\nu},\varepsilon)$, $\bar{\nu}=(\bar{\nu}_1\geq \bar{\nu}_2\geq \ldots \geq \bar{\nu}_m)$.

From Lemma \ref{mainl*} we deduce that in the case of $D_n$ the contributions of the cycles $\hat{s}_k$ of $s$ to the formula for ${\rm dim}~\Sigma_{{\bf k}, s}$ can be obtained from the corresponding contributions in case of $C_n$ in the following way:
for each pair of odd parts $\bar{\nu}_k=\bar{\nu}_{k+1}$ and for each pair of even parts $\bar{\nu}_k=\bar{\nu}_{k+1}$ with $\varepsilon(\bar{\nu}_k)=0$ the corresponding contribution $\underline{l}(\hat{s}_k)$ to $l(s)$ should be reduced by 2 and for each even part $\bar{\nu}_k$ of $\bar{\nu}$ with $\varepsilon(\bar{\nu}_k)=1$ the corresponding contribution $\underline{l}(\hat{s}_k)$ to $l(s)$ should be reduced by 1. This observation and formula (\ref{dimsp1}) yield
\begin{eqnarray*}
{\rm dim}~\Sigma_{{\bf k}, s}= n+\sum_{i=1}^m(i-1)\bar{\nu}_i+\frac{1}{2}|\{i:\bar{\nu}_i ~\hbox{is odd}\}|+|\{i:\bar{\nu}_i ~\hbox{is even and}~\varepsilon(\bar{\nu}_i)=0\}|- \\
-|\{i:\bar{\nu}_i ~\hbox{is odd}\}|-|\{i:\bar{\nu}_i ~\hbox{is even}\}|= \\
= n+\sum_{i=1}^m(i-1)\bar{\nu}_i-\frac{1}{2}|\{i:\bar{\nu}_i ~\hbox{is odd}\}|-|\{i:\bar{\nu}_i ~\hbox{is even and}~\varepsilon(\bar{\nu}_i)=1\}|
\end{eqnarray*}
which coincides with (\ref{dimsoev}).

In the case of root systems of exceptional types ${\rm dim}~\Sigma_{{\bf k}, s}$ can be found in the tables in Appendix 2. According to those tables equality (\ref{dimid}) holds in all cases.

\end{proof}

Now we show that all conjugacy classes in the stratum $G_\mathcal{O}=\phi_{G_{\bf k}}^{-1}(F(\mathcal{O}))$ intersect the corresponding variety $\Sigma_{{\bf k}, s}$, $s\in \Psi^W(\mathcal{O})$. The strategy of the proof is as follows. We are going to use characterization (\ref{charstr}) of the stratum $G_\mathcal{O}$.   

Let $\Delta_+^s$ be the system of the positive roots introduced in the statement of Theorem \ref{mainth},
\begin{equation}\label{dsdef1}
\h_\mathbb{R}=\bigoplus_{i=0}^{M(s)} \h_i
\end{equation}
the corresponding decomposition of $\h_\mathbb{R}$ and $h_i\in \h_i$ the corresponding elements of the subspaces $\h_i$.

We are going to rearrange the terms in the sum (\ref{dsdef1}) and define a system of positive roots $\Delta_+^{s'}$, using the rearranged sum as described in Section \ref{background}, in such a way that $s$ is elliptic in a standard parabolic subgroup $W^{s'}\subset W$ with respect to the system of simple roots in $\Delta_+^{s'}$. 

Let $w$ be a minimal length representative in the conjugacy class of $s$ in $W^{s'}$ with respect to the system of simple reflections in $W^{s'}$.

Let $B_{{\bf k},-}^{s'}$ be the  Borel subgroup in $G_{\bf k}$ corresponding to $-\Delta_+^{s'}$, $L_{\bf k}^{s'}$ the standard Levi subgroup with the Weyl group $W^{s'}$.
Denote by $B_{{\bf k},-}^{s''}=B_{{\bf k},-}^{s'}\cap L_{\bf k}^{s'}$ the Borel subgroup in $L_{\bf k}^{s'}$. One can always find a representative $\dot{w}\in N_{L_{\bf k}^{s'}}(H_{\bf k})$ of $w$. 

By characterization (\ref{charstr}) any conjugacy class in $G_\mathcal{O}$ intersects $B_{{\bf k},-}^{s'}\dot{w}B_{{\bf k},-}^{s'}$. In Lemma \ref{inters} we show that in fact
any conjugacy class in $G_\mathcal{O}$ intersects $B_{{\bf k},-}^{s''}\dot{w}B_{{\bf k},-}^{s''}\subset B_{{\bf k},-}^{s'}\dot{w}B_{{\bf k},-}^{s'}$.

Next, in Lemmas \ref{goodk} and \ref{minl} we prove that in fact one can take $w=s$.
In Lemma \ref{lemmaend} we verify that $B_{{\bf k},-}^{s''}\dot{s}B_{{\bf k},-}^{s''}\subset N_{\bf k}^s\dot{s}H_{\bf k}N_{\bf k}^s$, and hence any conjugacy class in $G_\mathcal{O}$ intersects $N_{\bf k}^s\dot{s}H_{\bf k}N_{\bf k}^s\supset B_{{\bf k},-}^{s''}\dot{s}B_{{\bf k},-}^{s''}$. 

We prove that any element of $N_{\bf k}^s\dot{s}H_{\bf k}N_{\bf k}^s$ can be conjugated to an element of $N_{\bf k}^s\dot{s}H_{\bf k}^0N_{\bf k}^s\subset N_{\bf k}^s\dot{s}Z_{\bf k}^sN_{\bf k}^s$. By Proposition \ref{crosssect} (i) this implies that any conjugacy class in $G_\mathcal{O}$ also intersects $\Sigma_{{\bf k}, s}$. Finally we verify that in fact any conjugacy class $G_\mathcal{O}$ intersects  $\Sigma_{{\bf k}, s}$ at some point of $\dot{s}H_{\bf k}^0N_{{\bf k}, s}$.

We start realizing this program by defining $\Delta_+^{s'}$. Recall that $\h_0$ is the subspace of $\h_{\mathbb{R}}$ fixed by the action of $s$.
If $\h_0=0$, let $\h_i'=\h_i$ and $h_i'=h_i$, $i=0,\ldots, M(s)$. Otherwise let $\h_{M(s)}'=\h_0$, $\h_i'=\h_{i+1}$, $h_i'=h_{i+1}$, $i=0,\ldots ,M(s)-1$ and choose an element $h_{M(s)}'\in \h_{M(s)}'$ such that $h_{M(s)}'(\alpha)\neq 0$ for any root $\alpha \in \Delta$ which is not orthogonal to the $s$--invariant subspace $\h_{M(s)}'$ with respect to the natural pairing between $\h_{\mathbb{R}}$ and $\h_{\mathbb{R}}^*$. \index[not]{h@$\h_i'$}  \index[not]{h@$h_i'$}

By a suitable rescaling of $h_{M(s)}'$ we can assume that conditions (\ref{cond}) are satisfied for the elements $h_i'$ and roots $\alpha$ from the sets $\Delta_{i}'$ defined as in (\ref{dik}) with $h_j,h_i$ replaced by $h_j',h_i'$. Indeed, observe that
$$
\Delta_{i}'=\{ \alpha\in \Delta: h_j'(\alpha)=0, j>i,~h_i'(\alpha)\neq 0 \}\subset \{ \alpha\in \Delta: h_j(\alpha)=0, j>{i+1},~h_{i+1}(\alpha)\neq 0 \}= \Delta_{i+1}, i=0, \ldots ,M(s)-1
$$
by the definition of the elements $h_i'$. Thus, since (\ref{cond}) is satisfied for $h_i$, $i=0,\ldots ,M(s)$, it is also satisfied for $h_i'=h_i$ if $i<M(s)$. By a suitable rescaling of $h_{M(s)}'$ we can assume that (\ref{cond}) is satisfied for $h_{M(s)}'$ as well.

Let $\Delta_+^{s'}$ \index[not]{D@$\Delta_+^{s'}$} be the system of positive roots in $\Delta=\Delta(G_{\bf k},H_{\bf k})$ which corresponds to the Weyl chamber containing the element $\bar{h}':=\sum_{i=0}^{M(s)}h_i'$.
By Lemma \ref{parab} the set of roots annihilating $\h_0$ is the root system of a standard Levi subgroup $L_{\bf k}^{s'}\subset G_{\bf k}$ \index[not]{L@$L_{\bf k}^{s'}$} with respect to the system of simple roots in $\Delta_+^{s'}$. As before, we denote the root system of the pair $(L_{\bf k}^{s'},H_{\bf k})$ by $\Delta(L_{\bf k}^{s'},H_{\bf k})$. 

Using formula (\ref{inv}) and recalling that the roots $\gamma_1, \ldots \gamma_{l'}$ form a linear basis of $\h'^*$, so that for $i=1,\ldots , l'$ $\gamma_i$ annihilate $\h_0$, we deduce that for $i=1,\ldots , l'$ $\gamma_i\in \Delta(L_{\bf k}^{s'},H_{\bf k})$, and hence $s$ belongs to the Weyl group $W^{s'}\subset W$ \index[not]{W@$W^{s'}$} of the root system $\Delta(L_{\bf k}^{s'},H_{\bf k})$. Note that, as $L_{\bf k}^{s'}$ is a standard Levi subgroup in $G_{\bf k}$ with respect to the system of simple roots in $\Delta_+^{s'}$, $W^{s'}$ is a parabolic subgroup in $W$ with respect to the system of simple roots in $\Delta_+^{s'}$. Since $\gamma_1, \ldots \gamma_{l'}$ form a linear basis of $\h'^*$, the linear span of the roots from $\Delta(L_{\bf k}^{s'},H_{\bf k})$ coincides with $\h'^*$, and hence the element $s$ is elliptic in $W^{s'}$ as $s$ acts without fixed points on $\h'^*$.

Let $w$ be a minimal length representative in the conjugacy class of $s$ in $W^{s'}$ with respect to the system of simple roots in $\Delta(L_{\bf k}^{s'},H_{\bf k})_+:=\Delta_+^{s'}\cap\Delta(L_{\bf k}^{s'},H_{\bf k})$. By Lemma 3.1.14 in \cite{Gck} if $w\in \Psi^W(\mathcal{O})\cap W^{s'}$ is of minimal possible length with respect to the system of simple reflections in $W^{s'}$ then it is also of minimal possible length with respect to the system of simple reflections in $W$, where in both cases the simple reflections are the reflections with respect to the simple roots in $\Delta_+^{s'}$. Note that $w$ is elliptic in $W^{s'}$ as well.

Let $B_{{\bf k},-}^{s'}$ \index[not]{B@$B_{{\bf k},-}^{s'}$} be the Borel subgroup in $G_{\bf k}$ corresponding to $-\Delta_+^{s'}$, $P_{\bf k}^{s'}\supset B_{{\bf k},-}^{s'}$ \index[not]{P@$P_{\bf k}^{s'}$} the parabolic subgroup of $G_{\bf k}$ corresponding to $W^{s'}$. Thus $L_{\bf k}^{s'}$ is the Levi factor of $P_{\bf k}^{s'}$.

Denote by $B_{{\bf k},-}^{s''}=B_{{\bf k},-}^{s'}\cap L_{\bf k}^{s'}$ \index[not]{B@$B_{{\bf k},-}^{s''}$} the Borel subgroup in $L_{\bf k}^{s'}$. One can always find a representative $\dot{w}\in N_{L_{\bf k}^{s'}}(H_{\bf k})$ of $w$.

\begin{lemma}\label{inters}
Any conjugacy class in $G_\mathcal{O}$ intersects $B_{{\bf k},-}^{s''}\dot{w}B_{{\bf k},-}^{s''}\subset B_{{\bf k},-}^{s'}\dot{w}B_{{\bf k},-}^{s'}$.
\end{lemma}

\begin{proof}
By characterization (\ref{charstr}) the stratum $G_\mathcal{O}$ consists of all conjugacy classes of minimal possible dimension which intersect the Bruhat cell $B_{{\bf k},-}^{s'}\dot{w}B_{{\bf k},-}^{s'}$. 

Denote by $U_{\bf k}^{s'}$ \index[not]{U@$U_{\bf k}^{s'}$} the unipotent radical \index{radical!unipotent} of $P_{\bf k}^{s'}$, and by $Z_{G_{\bf k}}(L_{\bf k}^{s'})$ \index[not]{Z@$Z_{G_{\bf k}}(L_{\bf k}^{s'})$} the centralizer of $L_{\bf k}^{s'}$ in $G_{\bf k}$. Then by the definition of parabolic subgroups one can always find a one parameter subgroup $\rho:{\bf k}^*\rightarrow Z_{G_{\bf k}}(L_{\bf k}^{s'})$ \index[not]{r@$\rho(t)$} such that
\begin{equation}\label{lm*}
\lim_{t\rightarrow 0}\rho(t)n\rho(t^{-1})=1
\end{equation}
for any $n\in U_{\bf k}^{s'}$.

Let $\gamma \in G_\mathcal{O}$ be a conjugacy class which intersects $B_{{\bf k},-}^{s'}\dot{w}B_{{\bf k},-}^{s'}$ at point $b\dot{w}b'$, $b,b'\in B_{{\bf k},-}^{s'}$ such that $b\dot{w}b'\not \in B_{{\bf k},-}^{s''}\dot{w}B_{{\bf k},-}^{s''}$. Since by the  definitions of $B_{{\bf k},-}^{s''}$ and $U_{\bf k}^{s'}$ we have a unique factorization $B_{{\bf k},-}^{s'}=B_{{\bf k},-}^{s''}U_{\bf k}^{s'}$, there are unique factorizations $b=un$, $b'=u'n'$, $u,u'\in B_{{\bf k},-}^{s''}$, $n,n'\in U_{\bf k}^{s'}$.
By (\ref{lm*}) we have
$$
\lim_{t\rightarrow 0}\rho(t)b\dot{w}b'\rho(t^{-1})=\lim_{t\rightarrow 0}u\rho(t)n\rho(t^{-1})\dot{w}u'\rho(t)n'\rho(t^{-1})=u\dot{w}u'\in B_{{\bf k},-}^{s''}\dot{w}B_{{\bf k},-}^{s''},
$$
and hence the closure of $\gamma$ contains a conjugacy class $\gamma'$ which intersects $B_{{\bf k},-}^{s'}\dot{w}B_{{\bf k},-}^{s'}$ at some point of $B_{{\bf k},-}^{s''}\dot{w}B_{{\bf k},-}^{s''}\subset B_{{\bf k},-}^{s'}\dot{w}B_{{\bf k},-}^{s'}$. In particular, ${\rm dim}~\gamma > {\rm dim}~\gamma'$. This is impossible by characterization (\ref{charstr}) of $G_\mathcal{O}$ as $\gamma$ has minimal possible dimension among the conjugacy classes intersecting $B_{{\bf k},-}^{s'}\dot{w}B_{{\bf k},-}^{s'}$.  Hence $\gamma$ intersects $B_{{\bf k},-}^{s'}\dot{w}B_{{\bf k},-}^{s'}$ at some point of $B_{{\bf k},-}^{s''}\dot{w}B_{{\bf k},-}^{s''}\subset B_{{\bf k},-}^{s'}\dot{w}B_{{\bf k},-}^{s'}$.

\end{proof}

\begin{lemma}\label{goodk}
Let $G_{\bf k}$ be a connected semisimple algebraic group over an algebraically closed field $\bf k$ of characteristic good for $G_{\bf k}$. Let $H_{\bf k}$ be a maximal torus of $G_{\bf k}$, $W$ the Weyl group of the pair $(G_{\bf k},H_{\bf k})$, and $s\in W$ an elliptic element. Denote by $\mathcal{C}_s$ \index[not]{C@$\mathcal{C}_s$} the conjugacy class of $s$ in $W$. Then $\Phi^{W}(\mathcal{C}_s)\subset {\underline{{\mathcal{N}}}}(G_{\bf k})$.
\end{lemma}

\begin{proof}
The statement of this lemma is a consequence of the fact that $s$ is elliptic. Indeed, it suffices to consider the case when $G_{\bf k}$ is simple. 

In the case when $G_{\bf k}$ is of type $A_n$ this is obvious since ${\underline{\widehat{\mathcal{N}}}}(W)$ consists of the unipotent classes of $G_{\bf k}$, ${\underline{\widehat{\mathcal{N}}}}(W)={\underline{{\mathcal{N}}}}(G_{\bf k})$. In fact in this case $\mathcal{C}_s$ is the Coxeter class, and $\Phi^{W}(\mathcal{C}_s)$ is the class of regular unipotent elements. 

If $G_{\bf k}$ is of type $B_n$, $C_n$ or $D_n$, formula (\ref{hofix}) implies that if $\Phi^{W}(\mathcal{C}_s)$ corresponds to $(\bar{\nu},\varepsilon)\in \mathcal{T}^2_{2n}$ (or $(\bar{\nu},\varepsilon)\in\widetilde{\mathcal{T}}^2_{2n}$) then $\bar{\nu}$ has no odd parts and no even parts $\bar{\nu}_i$ with $\varepsilon(\bar{\nu}_i)=0$. According to the description given in the previous section the map ${{\pi}}^{G_1}$ (resp. ${\widetilde{\pi}}^{G_1}$) is injective and its image consists of pairs $(\bar{\nu}, \varepsilon)\in {\mathcal{T}}_{2n}^2$ (resp. $\widetilde{\mathcal{T}}_{2n}^2$) such that $\varepsilon(k)\neq 0$ if $\bar{\nu}_k^*$ is odd and for each even $i$ such that $\bar{\nu}_i^*$ is even we have $\bar{\nu}_{i-1}^*=\bar{\nu}_i^*$, i.e. $i-1$ does not appear in the partition $\bar{\nu}$. We deduce that $\Phi^{W}(\mathcal{C}_s)$ is contained in the image of ${{\pi}}^{G_1}$ (resp. ${\widetilde{\pi}}^{G_1}$), i.e. $\Phi^{W}(\mathcal{C}_s)\in {\underline{{\mathcal{N}}}}(G_{\bf k})$ is a unipotent class in $G_{\bf k}$. 

In the case when $G_{\bf k}$ is of exceptional type this can be checked by examining the tables in Appendix 2.

\end{proof}

Now we show that in fact one can always take $w=s$.

\begin{lemma}\label{minl}
The element $s$ is of minimal length in its conjugacy class in $W^{s'}$ with respect to the system of simple roots in $\Delta(L_{\bf k}^{s'},H_{\bf k})_+=\Delta_+^{s'}\cap\Delta(L_{\bf k}^{s'},H_{\bf k})$.
\end{lemma}

\begin{proof}
First observe that the formulation of this lemma only uses root systems which can be  described independently of $\bf k$. Therefore we can assume in the proof that the characteristic of $\bf k$ is good for $G_{\bf k}$.

Let $M_{\bf k}^{s'}$ \index[not]{M@$M_{\bf k}^{s'}$} be the semisimple part of $L_{\bf k}^{s'}$ and $\mathcal{O}_n:=\Phi^{W^{s'}}(\mathcal{O}_w)\subset M_{\bf k}^{s'}$, \index[not]{O@$\mathcal{O}_n$} where $\mathcal{O}_w$ is the conjugacy class of $w$ in the Weyl group $W^{s'}=W(L_{\bf k}^{s'},H_{\bf k})$. By the previous lemma applied to the group $M_{\bf k}^{s'}$ and the elliptic element $w\in W^{s'}$ we have $\mathcal{O}_n\in {\underline{{\mathcal{N}}}}(M_{\bf k}^{s'})$.

Therefore $\mathcal{O}_n$ is the unipotent class of minimal possible dimension which intersects $B_{{\bf k},-}^{s''}\dot{w}B_{{\bf k},-}^{s''}$.
By Theorem 0.7 in \cite{L4'} the codimension of $\mathcal{O}_n$ in $M_{\bf k}^{s'}$ is equal to $l_1(w)$, where $l_1$ \index[not]{l@$l_1$} is the length function in $W^{s'}$ with respect to the system of simple roots in $\Delta(L_{\bf k}^{s'},H_{\bf k})_+=\Delta_+^{s'}\cap\Delta(L_{\bf k}^{s'},H_{\bf k})$,
\begin{equation}\label{cdm}
{\rm codim}_{M_{\bf k}^{s'}}~\mathcal{O}_n=l_1(w).
\end{equation}

Now we show that $s$ has minimal length in the Weyl group $W^{s'}$ with respect to the system of simple roots in the set of positive roots $\Delta(L_{\bf k}^{s'},H_{\bf k})_+$.

Indeed, let $\Sigma_{{\bf k}, s}'$ be the variety in $M_{\bf k}^{s'}$ associated to $s\in W^{s'}$ and defined similarly to $\Sigma_{{\bf k}, s}\subset G_{\bf k}$, where we use $\Delta(L_{\bf k}^{s'},H_{\bf k})_+$ as the system of positive roots in the definition of $\Sigma_{{\bf k}, s}'$.

Formula (\ref{dimid}) confirmed in Lemma \ref{dimsect1} is applicable to the slice $\Sigma_{{\bf k}, s}'$ and yields
$$
{\rm codim}_{M_{\bf k}^{s'}}~\mathcal{O}_n={\rm dim}~\Sigma_{{\bf k}, s}'.
$$

Formula (\ref{dimsect}) and the fact that $s$ is elliptic in $W^{s'}$ imply that
$$
{\rm dim}~\Sigma_{{\bf k}, s}'=l_1(s),
$$

From the last two formulas we infer
$$
{\rm codim}_{M_{\bf k}^{s'}}~\mathcal{O}_n={\rm dim}~\Sigma_{{\bf k}, s}'=l_1(s).
$$

The last formula and (\ref{cdm}) yield $l_1(w)=l_1(s)$, and hence $s$ has minimal possible length in its conjugacy class in $W^{s'}$ with respect to the system of simple roots in $\Delta(L_{\bf k}^{s'},H_{\bf k})_+$.

\end{proof}

Now we can assume that $s=w$.

\begin{lemma}\label{lemmaend}
Any conjugacy class $\gamma \in G_\mathcal{O}$ intersects $\Sigma_{{\bf k}, s}$ at some point of $\dot{s}H_{\bf k}^0N_{{\bf k}, s}$.
\end{lemma}

\begin{proof}
Let $\alpha\in \Delta\cap\h'^*$. Then $\alpha\in \Delta_{i}\cap\h'^*$ for some $i>0$. Observe that by the definition of the subspaces $\h_k'$ and by the choice of the elements $h_k$, $k=0,\ldots, M(s)$ one has
$$
\Delta_{i}\cap\h'^*=\{ \beta \in \Delta  : h_j(\beta)=0, j>i, h_i(\beta)\neq 0, h_0(\beta)=0\}=\{ \beta \in \Delta  : h_j'(\beta)=0, j>i-1, h_{i-1}'(\beta)\neq 0, \}=\Delta_{i-1}',
$$
and hence by (\ref{wc}) $\alpha\in (\Delta_{i})_+\cap\h'^*$ if and only if $h_{i-1}'(\alpha)=h_i(\alpha)>0$, i.e. $\alpha\in (\Delta_{i})_+\cap\h'^*$ if and only if $\alpha \in \Delta_+^{s'}\cap \Delta_{i}'$.
Therefore if we denote $B_{{\bf k},-}^{s'''}=B_{{\bf k},-}^s\cap L_{\bf k}^{s'}$, \index[not]{B@$B_{{\bf k},-}^{s'''}$} where $B_{{\bf k},-}^s$ is the Borel subgroup corresponding to $-\Delta_+^s$, then $B_{{\bf k},-}^{s'''}=B_{{\bf k},-}^{s'}\cap L_{\bf k}^{s'}=B_{{\bf k},-}^{s''}$.

By Lemma \ref{inters} applied to the element $s\in W$ any conjugacy class $\gamma \in G_\mathcal{O}$ intersects $B_{{\bf k},-}^{s''}\dot{s}B_{{\bf k},-}^{s''}$, and hence it also intersects $B_{{\bf k},-}^{s'''}\dot{s}B_{{\bf k},-}^{s'''}\subset B_{{\bf k},-}^s\dot{s}B_{{\bf k},-}^s$ as $B_{{\bf k},-}^{s'''}=B_{{\bf k},-}^{s''}$. But by the definition of $L_{\bf k}^{s'}$ $s$ acts on the root system of the pair $(L_{\bf k}^{s'},H_{\bf k})$ without fixed points. Since $B_{{\bf k},-}^{s'''}=B_{{\bf k},-}^s\cap L_{\bf k}^{s'}$ and $s$ fixes all the roots of the pair $(Z_{\bf k}^sH_{\bf k},H_{\bf k})$ we have an inclusion $B_{{\bf k},-}^{s'''}\subset H_{\bf k}N_{\bf k}^s$, where $N_{\bf k}^s$ is the unipotent radical of the parabolic subgroup $P_{\bf k}^s$ associated to $s$. Hence  $B_{{\bf k},-}^{s'''}\dot{s}B_{{\bf k},-}^{s'''}\subset N_{\bf k}^s\dot{s}H_{\bf k}N_{\bf k}^s$.

Recall that we denote by $H_{\bf k}^0$ the identity component of the centralizer $H_{\bf k}^s$ \index[not]{H@$H_{\bf k}^s$} of $\dot{s}$ in $H_{\bf k}$. Next we show that all elements of $N_{\bf k}^s\dot{s}H_{\bf k}N_{\bf k}^s$ are conjugate to elements of $N_{\bf k}^s\dot{s}H_{\bf k}^0N_{\bf k}^s\subset N_{\bf k}^s\dot{s}Z_{\bf k}^sN_{\bf k}^s$ by elements of $H_{\bf k}$. 

Observe that by parts (2) and (3) of the proof of the Theorem in Section 18.3 of \cite{Hu} the set of elements $h\dot{s}h^{-1}\dot{s}^{-1}$, $h\in H_{\bf k}$ is a closed subgroup $H_{\bf k}'$ \index[not]{H@$H_{\bf k}'$} of $H_{\bf k}$, and the map
$$
H_{\bf k}\to H'_{\bf k}, h\mapsto h\dot{s}h^{-1}\dot{s}^{-1}
$$
is a surjective group homomorphism with the kernel being $H_{\bf k}^s$. Thus ${\rm dim}~H_{\bf k}={\rm dim}~H_{\bf k}^s+{\rm dim}~H_{\bf k}'$.

Now consider the group homomorphism
$$
\phi: H_{\bf k}^s\times H_{\bf k}'\to H_{\bf k}
$$
induced by the group multiplication in $H_{\bf k}$. By Corollary B in Section 18.1 of \cite{Hu} its differential at the identity element of $H_{\bf k}^s\times H_{\bf k}'$ is a linear isomorphism of the corresponding tangent spaces, i.e. of the Lie algebras of $H_{\bf k}^s\times H_{\bf k}'$ and of $H_{\bf k}$. Since $\phi$ is a group homomorphism, this differential is an isomorphism of the corresponding tangent spaces at all points of $H_{\bf k}^s\times H_{\bf k}'$. Thus its image must have the dimension equal to ${\rm dim}~H_{\bf k}^s+{\rm dim}~H_{\bf k}'={\rm dim}~H_{\bf k}$.  

By 4.4, Proposition B in \cite{Hu} the image of the homomorphism $\phi$ is a closed subgroup in $H_{\bf k}$. The identity ${\rm dim}~H_{\bf k}^s+{\rm dim}~H_{\bf k}'={\rm dim}~H_{\bf k}$ implies that this image contains the identity component of $H_{\bf k}$. But $H_{\bf k}$ is irreducible, and hence the image of $\phi$ coincides with $H_{\bf k}$, i.e. $\phi$ is surjective.

For the same reason the image of $H_{\bf k}^0\times H_{\bf k}'$ in $H_{\bf k}$ under $\phi$ coincides with $H_{\bf k}$, i.e. the restriction of $\phi$ to $H_{\bf k}^0\times H_{\bf k}'$ is surjective. 


This implies that for any $h\in H_{\bf k}$ there exist $h^0\in H_{\bf k}^0$ and $h'\dot{s}{h'}^{-1}\dot{s}^{-1}\in H_{\bf k}'$ ($h'\in H_{\bf k}$) such that
$$
h=h^0h'\dot{s}{h'}^{-1}\dot{s}^{-1}=h'h^0\dot{s}{h'}^{-1}\dot{s}^{-1},
$$
or
\begin{equation}\label{h'def}
{h'}^{-1}h\dot{s}h'=h^0\dot{s}.
\end{equation}
Since $H_{\bf k}$ normalizes $N_{\bf k}^s$, we deduce from the last identity that any element $nh\dot{s}n'\in N_{\bf k}^sH_{\bf k}\dot{s}N_{\bf k}^s=N_{\bf k}^s\dot{s}H_{\bf k}N_{\bf k}^s$, $n,n'\in N_{\bf k}^s$, $h\in H_{\bf k}$, can be conjugated by an element $h'\in H_{\bf k}$ to the element ${h'}^{-1}nh'h^0\dot{s}{h'}^{-1}n'h'\in N_{\bf k}^sH_{\bf k}^0\dot{s}N_{\bf k}^s=N_{\bf k}^s\dot{s}H_{\bf k}^0N_{\bf k}^s$, where $h^0\in H_{\bf k}^0$ and $h'$ are related to $h$ by (\ref{h'def}).

Finally observe that $N_{\bf k}^s\dot{s}H_{\bf k}^0N_{\bf k}^s\subset N_{\bf k}^s\dot{s}Z_{\bf k}^sN_{\bf k}^s$, and hence any conjugacy class $\gamma \in G_\mathcal{O}$ intersects $N_{\bf k}^s\dot{s}H_{\bf k}^0N_{\bf k}^s\subset N_{\bf k}^s\dot{s}Z_{\bf k}^sN_{\bf k}^s$. By Proposition \ref{crosssect} (i)  $\gamma$ also intersects $\Sigma_{{\bf k}, s}$. Formula (\ref{eqz}) implies that the $Z_{\bf k}^s$--component of any element from $N_{\bf k}^s\dot{s}Z_{\bf k}^sN_{\bf k}^s$ is equal to the $Z_{\bf k}^s$--component in $\Sigma_{{\bf k}, s}=\dot{s}Z_{\bf k}^sN_{{\bf k}, s}$ of its image under the isomorphism $N_{\bf k}^s\dot{s}Z_{\bf k}^sN_{\bf k}^s\simeq N_{\bf k}^s\times \dot{s}Z_{\bf k}^sN_{{\bf k}, s}$. Therefore  any conjugacy class $\gamma \in G_\mathcal{O}$ intersects  $\Sigma_{{\bf k}, s}$ at some point of $\dot{s}H_{\bf k}^0N_{{\bf k}, s}$.
This completes the proof.

\end{proof}

By the previous lemma the first claim in Theorem \ref{mainth} holds. This completes the proof.

\end{proof}


\section{Bibliographic comments}

\pagestyle{myheadings}
\markboth{CHAPTER \thechapter.~ALGEBRAIC GROUP ANALOGUES OF SLODOWY SLICES}{\thesection.~BIBLIOGRAPHIC COMMENTS}

\setcounter{equation}{0}
\setcounter{theorem}{0}

A uniform classification of conjugacy classes of Weyl group elements from which one can obtain presentation (\ref{inv}) was obtained in \cite{C}. 

The definition of systems of positive roots $\Delta_+^s$ associated to (conjugacy classes of) Weyl group elements was introduced in \cite{S6}. It is based on a deep generalization of the results by Coxeter and Steinberg on the properties of the Coxeter elements. In our notation this corresponds to the case when $\gamma_1, \ldots, \gamma_{l'}$ is a set of simple roots in $\Delta$, so that according to (\ref{inv}) $s$ is a product of simple reflections, i.e. a Coxeter element. In this case there is a unique plane in $\h_{\mathbb{R}}$, called a Coxeter plane, on which $s$ acts by rotation by the angle $2\pi/h$, where $h$ is the Coxeter number of $\g$. This plane was introduced by Coxeter in book \cite{Cox}, and the pictures of root systems of Lie algebras of high ranks which one can find in many textbooks are obtained using orthogonal projections of roots onto these planes. The key observation is that all these projections are non--zero. Coxeter originally applied this property to study regular polytopes.

Later in paper \cite{St1} Steinberg proved interesting properties of Coxeter elements using the properties of the action of Coxeter elements on Coxeter planes.

The construction of the spectral decomposition for Weyl group elements in Proposition \ref{IM} suggested in \cite{S12} is a generalization of similar results on the properties of the Coxeter plane which can be found in \cite{Car}, Section 10.4.

Normal orderings of positive root systems of the form $\Delta_+^s$ described in  Proposition \ref{pord} were firstly introduced in \cite{S10} where one can also find the construction of normal orderings of positive root systems compatible with Weyl group involutions from Appendix 1. Later the original definition was refined in \cite{S13}. Proposition \ref{pord} is a modified version of Proposition 5.1 in \cite{S10} and of Proposition 2.2 in \cite{S13}.

Circular orderings of root systems were defined in \cite{KT3} to describe commutation relations between quantum group analogues of root vectors. In \cite{SDM} this construction was used to modify the positive root systems $\Delta_+^s$ in order to construct positive root systems associated to (conjugacy classes of) Weyl group elements which appear in the end of Section \ref{wgrord}.

The slices $\Sigma_{{\bf k}, s}$ introduced in \cite{S6} in the case ${\bf k}=\mathbb{C}$ are generalizations of the Steinberg cross-sections to the set of conjugacy classes of regular elements in $G_{\bf k}$ suggested in \cite{St2}. $\Sigma_{{\bf k}, s}$ reduces to a Steinberg cross-section when $\gamma_1, \ldots, \gamma_{l'}$ is a set of simple roots, i.e. when $s$ is a Coxeter element. In this case isomorphism (\ref{cross}) is stated in \cite{St2} without proof. The first proof of this result appeared in \cite{SemSev}. However, that proof is not applicable for the root system of type $E_6$. The proof of isomorphism (\ref{cross}) in Proposition \ref{crosssect} based on the properties of normal ordering (\ref{NO}) is significantly simpler than the first complete proof of this result presented in \cite{S6}. 

Another construction of the slices $\Sigma_{{\bf k}, s}$ in the case when $s$ are elliptic can be found in \cite{XL}.

A more general approach to the definition of the transverse slices to conjugacy classes in algebraic groups generalizing the definition given in this chapter and the results of \cite{XL} was given in \cite{Malt2}, and the relevant Weyl combinatorics was developed in \cite{Malt1}. 
 
The closedness of the varieties $N^sZ^ssN^s$ was justified in \cite{SDM}, Proposition 6.2. 

In book \cite{SL} Slodowy proved Brieskorn's conjectures announced in \cite{Br} on the realization of simple singularities using the adjoint quotient of complex semisimple Lie algebras. Although a significant part of Slodowy's book is devoted to the study of the conjugation quotient for semisimple algebraic groups and to constructing some its resolutions, he ended up with a Lie algebra version of the construction of simple singularities and introduced transverse slices for the adjoint action for this purpose. These slices are called now the Slodowy slices. The slices $\Sigma_{{\bf k}, s}$ can be regarded as algebraic group analogues of the Slodowy slices.

The Lusztig partition was introduced in \cite{L2}. Its definition is related to the study of the properties of intersections of conjugacy classes in $G_{\bf k}$ with Bruhat cells established in \cite{L3',L4'} where the map $\Phi^W$ from the set of Weyl group conjugacy classes to the set of unipotent classes and its one sided inverse $\Psi^W$, which we use in Section \ref{luspart}, are defined using these properties. These properties are also related to the generalized Springer correspondence. We only briefly discussed the relevant results in this book. 

The study of intersections of conjugacy classes in $G_{\bf k}$ with Bruhat cells was initiated in \cite{St2}. Some results on these intersections were obtained in \cite{EG}, and another map from nilpotent orbits in a complex semisimple Lie algebra to conjugacy classes in the Weyl group was defined in \cite{KL}. In \cite{L3'} it is mentioned that this map is likely to coincide with the map $\Psi^W$ introduced in \cite{L3'}. The results of \cite{YU} imply that the restriction of $\Psi^W$ to $\underline{\mathcal{N}}(G_1)$ coincides with the Kazhdan-Lusztig map defined in \cite{KL}.

The main result of Theorem \ref{mainth} on the dimensions of the slices $\Sigma_{{\bf k}, s}$ is an experimental observation made in \cite{S12}, Theorem 5.2. Other results of Section \ref{stt} can also be found in \cite{S12}.

Note that the slices $\Sigma_{{\bf k}, s}$ listed in Appendix 2 are slightly different from those from Appendix B to \cite{S12}. The corresponding slices in both sets have the same dimensions. But in this book the roots $\gamma_1,\ldots, \gamma_{l'}$ in the tables in Appendix 2 are chosen in such a way that the corresponding root systems $\Delta_+^s$ satisfy condition (\ref{cond2}). The algorithm for constructing the slices $\Sigma_s$ listed in the tables in Appendix B to \cite{S12} was modified to fulfill this condition. The description of the original algorithm can be found in \cite{S12}. 

The ordering of the $s$--invariant planes in $\h_{\mathbb{R}}$ according to the angles of rotations by which $s$ acts in the planes as in Theorem \ref{mainth} was used in \cite{XN} to prove properties of minimal length elements in finite Coxeter groups.


\chapter{Quantum groups}\label{part2}

\pagestyle{myheadings}
\markboth{CHAPTER \thechapter.~QUANTUM GROUPS}{\thesection.~THE DEFINITION OF QUANTUM GROUPS}

In this chapter we recall some definitions and results on quantum groups required for the study of q-W--algebras. Besides the standard definitions and results related to quantum groups we shall need some rather non--standard realizations of the Drinfeld--Jimbo quantum group in terms of which q-W--algebras are defined. These realizations are related to the definition of the algebraic group analogues of the Slodowy slices in the previous chapter. We shall consider the Drinfeld--Jimbo quantum group $U_h({\frak g})$ defined over the ring of formal power series ${\Bbb C}[[h]]$, \index[not]{C@${\Bbb C}[[h]]$} where $h$ is an indeterminate, and some its specializations defined over smaller rings. In this chapter we shall also present some new results on the adjoint action of quantum groups required for the definition of Zhelobenko type operators for q-W--algebras in Chapter \ref{part4}.


\section{The definition of quantum groups}\label{QGdef}

\setcounter{equation}{0}
\setcounter{theorem}{0}

In this section we remind the definition of the standard Drinfeld-Jimbo quantum group $U_h({\frak g})$. We mainly follow the notation of \cite{ChP}. 

Let $V$ be a ${\Bbb C}[[h]]$--module equipped with the $h$--adic
topology. \index{topology!$h$--adic} This topology is characterized by requiring that
$\{ h^nV ~|~n\geq 0\}$ is a base of the neighborhoods of $0$ in $V$, and that translations
in $V$ are continuous. In this book all ${\Bbb C}[[h]]$--modules are supposed to be complete with respect to this topology.

{\it A topological Hopf algebra} over ${\Bbb C}[[h]]$ is a complete ${\Bbb C}[[h]]$--module 
equipped with a structure of ${\Bbb C}[[h]]$--Hopf algebra,
the algebraic tensor products entering the axioms of the Hopf algebra are replaced by their
completions in the $h$--adic topology. \index{Hopf algebra!topological}

{\it The standard quantum group} $U_h({\frak g})$ \index{quantum group!standard} associated to a complex finite-dimensional semisimple Lie algebra
$\frak g$ is a topological Hopf algebra over ${\Bbb C}[[h]]$ topologically generated by elements
$H_i,~X_i^+,~X_i^-,~i=1,\ldots ,l$, \index[not]{H@$H_i$} \index[not]{X@$X_i^\pm$} subject to the following defining relations:
\begin{equation}\label{defQGrel}
\begin{array}{l}
[H_i,H_j]=0,~~ [H_i,X_j^\pm]=\pm a_{ij}X_j^\pm, ~~X_i^+X_j^- -X_j^-X_i^+ = \delta _{i,j}{K_i -K_i^{-1} \over q_i -q_i^{-1}},\\ \index[not]{d@$\delta_{i,j}$}
\\
\sum_{r=0}^{1-a_{ij}}(-1)^r
\left[ \begin{array}{c} 1-a_{ij} \\ r \end{array} \right]_{q_i}
(X_i^\pm )^{1-a_{ij}-r}X_j^\pm(X_i^\pm)^r =0 ,~ i \neq j , \index{commutation relations!defining for $U_h(\g)$}
\end{array}
\end{equation}
where $\delta_{i,j}$ is the Kronecker delta, \index{Kronecker delta}
$$
K_i=e^{d_ihH_i},~~e^h=q,~~q_i=q^{d_i}=e^{d_ih}, \index[not]{K@$K_i$} \index[not]{q@$q_i$}
$$
$$
\left[ \begin{array}{c} m \\ n \end{array} \right]_q={[m]_q! \over [n]_q![n-m]_q!} ,~\index[not]{ZZZ@$\left[ \begin{array}{c} m \\ n \end{array} \right]_q$}
[n]_q!=[n]_q\ldots [1]_q ,~ [n]_q={q^n - q^{-n} \over q-q^{-1} }, \index[not]{n@$[n]_q$}  
\index[not]{n@$[n]_q"!$}
$$
with comultiplication defined by \index{quantum group!standard!comultiplication}
$$
\Delta_h(H_i)=H_i\otimes 1+1\otimes H_i,~~ 
\Delta_h(X_i^+)=X_i^+\otimes K_i^{-1}+1\otimes X_i^+,~~\Delta_h(X_i^-)=X_i^-\otimes 1 +K_i\otimes X_i^-, \index[not]{D@$\Delta_h(~\cdot~)$}
$$
antipode defined by
$$
S_h(H_i)=-H_i,~~S_h(X_i^+)=-X_i^+K_i,~~S_h(X_i^-)=-K_i^{-1}X_i^-, \index[not]{S@$S_h(~\cdot~)$} \index{quantum group!standard!antipode}
$$
and counit defined by
$$
\varepsilon_h(H_i)=\varepsilon_h(X_i^\pm)=0. \index[not]{e@$\varepsilon_h(~\cdot~)$} \index{quantum group!standard!counit}
$$

We shall also use the weight--type generators
$$
Y_i=\sum_{j=1}^l d_i a^{-1}_{ij}H_j, \index[not]{Y@$Y_i$}
$$
where $a^{-1}_{ij}$ are the entries of the matrix inverse to the Cartan matrix $a_{ij}$.

Let $L_i^{\pm 1}=e^{\pm hY_i}$. \index[not]{L@$L_i$} These elements commute with the quantum simple root vectors $X_i^\pm$ as follows:
\begin{equation}\label{weight-root}
L_iX_j^\pm L_i^{-1}=q_i^{\pm \delta_{i,j}}X_j^\pm.
\end{equation}
We also obviously have
\begin{equation}\label{comml}
L_iL_j=L_jL_i.
\end{equation}

The Hopf algebra $U_h({\frak g})$ is a quantization of the standard Lie bialgebra \index{Lie!bialgebra} structure on $\frak g$ in the sense that $U_h({\frak g})/hU_h({\frak g})\simeq U({\frak g})$, \index[not]{U@$U(\g)$} $\Delta_h=\Delta$ (mod $h$), where $\Delta$ \index[not]{D@$\Delta(~\cdot~)$} is
the standard comultiplication on $U({\frak g})$, and
$$
{\Delta_h -\Delta_h^{opp} \over h}~(\mbox{mod }h)=-\delta.
$$
Here
$\delta: {\frak g}\rightarrow {\frak g}\otimes {\frak g}$ \index[not]{d@$\delta$} is the standard cocycle \index{cocycle!standard} on $\frak g$, and $\Delta^{opp}_h=\sigma \Delta_h$, \index[not]{D@$\Delta^{opp}_h(~\cdot~)$} $\sigma$ \index[not]{s@$\sigma$} is the permutation in $U_h({\frak g})^{\otimes 2}$,
$\sigma (x\otimes y)=y\otimes x$.
Recall that
$$
\delta (x)=({\rm ad}x\otimes 1+1\otimes {\rm ad}x)2r_\pm,~~ r_\pm\in {\frak g}\otimes {\frak g}, \index[not]{a@${\rm ad}$}
$$
\begin{equation}\label{rcl}
r_\pm=\pm\frac 12 \sum_{i=1}^lY_i \otimes H_i \pm \sum_{\beta \in \Delta_+}\left\langle X_{\beta},X_{-\beta}\right\rangle^{-1} X_{\pm\beta}\otimes X_{\mp\beta}. \index[not]{r@$r_\pm$}
\end{equation}
Here $X_{\pm \beta}\in {\frak g}_{\pm \beta}$ are non--zero root vectors of $\frak g$.
The element $r=r_++r_-\in {\frak g}\otimes {\frak g}$ \index[not]{r@$r$} is called {\it the classical standard r--matrix}. \index{r--matrix!classical!standard}


\section{The braid group action}\label{BGact}

\pagestyle{myheadings}
\markboth{CHAPTER \thechapter.~QUANTUM GROUPS}{\thesection.~THE BRAID GROUP ACTION}

\setcounter{equation}{0}
\setcounter{theorem}{0}

One can define a quantum group analogue of the braid group action on $\g$. The material covered in this section can be found, e.g., in Section 8.1 of \cite{ChP}.
Let $r_{ij}$, $i\neq j$ be equal to $2,3,4,6$ if $a_{ij}a_{ji}$ is equal to $0,1,2,3$, respectively. The {\it braid group} $\mathcal{B}_\g$ \index{group!braid} \index[not]{B@$\mathcal{B}_\g$} associated to $\g$ has generators $T_i$, $i=1,\ldots, l$, \index[not]{T@$T_i$} and defining relations
$$
T_iT_jT_iT_j\ldots=T_jT_iT_jT_i\ldots
$$
for all $i\neq j$, where there are $r_{ij}$ $T$'s on each side of the equation.

Recall that if $X_{\pm \alpha_i}$ are non--zero simple root vectors of $\g$ then one can introduce an action of the braid group $\mathcal{B}_\g$ by algebra automorphisms of ${\frak g}$ defined on the standard generators as follows: \index{action!braid group!on $\g$}
\begin{eqnarray}
T_i(X_{\pm \alpha_i})=-X_{\mp \alpha_i},~T_i(H_j)=H_j-a_{ji}H_i, \nonumber \\
\nonumber \\
T_i(X_{\alpha_j})=\frac{1}{(-a_{ij})!}
{\rm ad}_{X_{\alpha_i} }^{-a_{ij}}X_{\alpha_j},~i\neq j, \label{BGactg} \\
\nonumber \\
T_i(X_{-\alpha_j})=\frac{(-1)^{a_{ij}}}{(-a_{ij})!}
{\rm ad}_{X_{-\alpha_i} }^{-a_{ij}}X_{-\alpha_j},~i\neq j. \nonumber
\end{eqnarray}

Similarly, $\mathcal{B}_\g$ acts by algebra automorphisms of $U_h({\frak g})$ \index{action!braid group!on $U_h({\frak g})$} as follows:
\begin{eqnarray}
T_i(X_i^+)=-X_i^-e^{hd_iH_i},~T_i(X_i^-)=-e^{-hd_iH_i}X_i^+,~T_i(H_j)=H_j-a_{ji}H_i, \nonumber \\
\nonumber \\
T_i(X_j^+)=\sum_{r=0}^{-a_{ij}}(-1)^{r-a_{ij}}q_i^{-r}
(X_i^+ )^{(-a_{ij}-r)}X_j^+(X_i^+)^{(r)},~i\neq j, \label{BGactqg} \\
\nonumber \\
T_i(X_j^-)=\sum_{r=0}^{-a_{ij}}(-1)^{r-a_{ij}}q_i^{r}
(X_i^-)^{(r)}X_j^-(X_i^-)^{(-a_{ij}-r)},~i\neq j, \nonumber
\end{eqnarray}
where
$$
(X_i^+)^{(r)}=\frac{(X_i^+)^{r}}{[r]_{q_i}!},~(X_i^-)^{(r)}=\frac{(X_i^-)^{r}}{[r]_{q_i}!},~r\geq 0,~i=1,\ldots,l. \index[not]{X@$(X_i^\pm)^{(r)}$}
$$

Recall that action (\ref{BGactg}) of the generators $T_i$ is induced by the adjoint action of certain representatives of the Weyl group elements $s_i$ in $G$. Similarly, action (\ref{BGactqg}) is induced by conjugation by certain elements of a completion of $U_h({\frak g})$.

To define these elements consider the {\it restricted dual} \index{Hopf algebra!restricted dual} $\mathbb{C}_h[G]$ \index[not]{C@$\mathbb{C}_h[G]$} of $U_h({\frak g})$, i.e. the algebra topologically generated by the matrix elements of finite rank representations of $U_h({\frak g})$, \index{representation!of a quantum group!finite rank} with the multiplication induced by the comultiplication on $U_h({\frak g})$. $\mathbb{C}_h[G]$ is naturally a Hopf algebra with the comultiplication induced by the multiplication on $U_h({\frak g})$ (see e.g. \cite{ChP}, p. 113). 

The topological dual $\mathbb{C}_h[G]^*$ \index[not]{C@$\mathbb{C}_h[G]^*$} of the Hopf algebra $\mathbb{C}_h[G]$ \index{Hopf algebra!dual} is just an algebra, and there is a natural embedding of algebras $U_h({\frak g})\hookrightarrow \mathbb{C}_h[G]^*$. 

Now define a q-exponential by \index{q-exponential}
\begin{equation}\label{qexp'}
{\exp}_q'(x)=\sum_{k=0}^\infty q^{\frac{1}{2}k(k-1)}{x^k \over [k]_q!}. \index[not]{e@${\exp}_q'(~\cdot~)$}
\end{equation}

Then the automorphism $T_i$ in (\ref{BGactqg}) is given by the conjugation of elements of $U_h({\frak g})\hookrightarrow \mathbb{C}_h[G]^*$ in $\mathbb{C}_h[G]^*$ by the invertible element $T_i\in \mathbb{C}_h[G]^*$ (see \cite{Saito})
\begin{align}
T_i={\exp}_{q_i^{-1}}'(-q_i^{-1}X_i^-K_i){\exp}_{q_i^{-1}}'(X_i^+){\exp}_{q_i^{-1}}'(-q_iX_i^-K_i^{-1})q_i^{\frac{H_i(H_i+1)}{2}}= \label{T1}
\\
={\exp}_{q_i^{-1}}'(q_i^{-1}X_i^+K_i^{-1}){\exp}_{q_i^{-1}}'(-X_i^-){\exp}_{q_i^{-1}}'(q_iX_i^+K_i)q_i^{\frac{H_i(H_i+1)}{2}}. \nonumber \index[not]{T@$T_i$}
\end{align}
Note that the right hand side of the previous formula defines an endomorphism for every finite rank representation of $U_h({\frak g})$ as the elements $X_i^\pm$ act nilpotently on every such representation, and the action of the elements $H_i$ is semisimple.

The inverse of $T_i$ in (\ref{T1}) can be found using the identity
$$
{\exp}_q'(x){\exp}_{q^{-1}}'(-x)=1
$$
which implies
\begin{align} \label{T2}
T_i^{-1}=q_i^{-\frac{H_i(H_i+1)}{2}}{\exp}_{q_i}'(q_iX_i^-K_i^{-1}){\exp}_{q_i}'(-X_i^+){\exp}_{q_i}'(q_i^{-1}X_i^-K_i)=
\\
=q_i^{-\frac{H_i(H_i+1)}{2}}{\exp}_{q_i}'(-q_iX_i^+K_i){\exp}_{q_i}'(X_i^-){\exp}_{q_i}'(-q_i^{-1}X_i^+K_i^{-1}).
\nonumber
\end{align}

From formula (\ref{T2}) we obtain the following relations in $\mathbb{C}_h[G]^*$
\begin{align} \label{T3}
{\exp}_{q_i}'(-X_i^+)=
{\exp}_{q_i^{-1}}'(-q_iX_i^-K_i^{-1})q_i^{\frac{H_i(H_i+1)}{2}}T_i^{-1}{\exp}_{q_i^{-1}}'(-q_i^{-1}X_i^-K_i)=
\\
={\exp}_{q_i^{-1}}'(-q_iX_i^-K_i^{-1})q_i^{\frac{H_i(H_i+1)}{2}}{\exp}_{q_i^{-1}}'(q_i^{-1}X_i^+)T_i^{-1}.
\nonumber
\end{align}

Let
\begin{equation}\label{expq}
{\exp}_q(x)=exp_q'(qx)=\sum_{k=0}^\infty q^{\frac{1}{2}k(k+1)}{x^k \over [k]_q!}. \index[not]{e@${\exp}_q(~\cdot~)$}
\end{equation}

The multiplication in $\mathbb{C}_h[G]$ induces a map of $\mathbb{C}[[h]]$--modules $\Delta_h:\mathbb{C}_h[G]^*\to (\mathbb{C}_h[G]\otimes \mathbb{C}_h[G])^*$, $(\Delta_h(f))(x\otimes y)=f(xy)$, $f\in\mathbb{C}_h[G]^*$, $x,y \in \mathbb{C}_h[G]$, the restriction of which to $U_h({\frak g})$ coincides with the comultiplication on $U_h({\frak g})$. With respect to this map we have (see \cite{ChP}, Proposition 8.2.6)

\begin{equation}\label{DT}
\Delta_h(T_i)=\theta_iT_i\otimes T_i=T_i\otimes T_i\overline{\theta}_i,
\end{equation}
$$
\theta_i={\exp}_{q_{i}}[(1-q_{i}^{-2})X_{i}^+\otimes X_{i}^-], \overline{\theta}_i={\exp}_{q_{i}}[(1-q_{i}^{-2})K_i^{-1}X_{i}^-\otimes X_{i}^+K_i], \index[not]{t@$\theta_i$} \index[not]{t@$\overline{\theta}_i$}
$$
\begin{equation}\label{DT-1}
\Delta_h(T_i^{-1})={\overline{\theta}_i}^{-1}T_i^{-1}\otimes T_i^{-1}=T_i^{-1}\otimes T_i^{-1}\theta_i^{-1},
\end{equation}
\begin{equation}\label{ti-1}
\theta_i^{-1}={\exp}_{q_{i}^{-1}}[(1-q_{i}^{2})X_{i}^+\otimes X_{i}^-], {\overline{\theta}_i}^{-1}={\exp}_{q_{i}^{-1}}[(1-q_{i}^{2})K_i^{-1}X_{i}^-\otimes X_{i}^+K_i],
\end{equation}
where the right hand sides of these identities are well--defined automorphisms of tensor products of finite rank representations of $U_h({\frak g})$, so that they can be evaluated on products in $\mathbb{C}_h[G]$ of matrix elements of such representations.

Denote $q_\alpha =q^{d_i}$ \index[not]{q@$q_\alpha$} 
if the positive root $\alpha$ is Weyl group conjugate to the simple root $\alpha_i$. By Proposition 8.1.3 in \cite{ChP}, for a reduced decomposition $w=s_{i_1}\ldots s_{i_k}$, $T_w=T_{i_1}\ldots T_{i_k}$ \index[not]{T@$T_w$} only depends on $w$ and (\ref{DT}) implies
\begin{equation}\label{dhtw}
\Delta_h(T_w)=\prod^{k}_{p=1}\theta_{\beta_p}T_w\otimes T_w=T_w\otimes T_w\prod_{p=1}^k\overline{\theta}_{\beta_p'},
\end{equation}
where in the products $\theta_{\beta_p}$ (resp. $\overline{\theta}_{\beta_p'}$) appears on the left from $\theta_{\beta_q}$ (resp. $\overline{\theta}_{\beta_q'}$) if $p<q$, and for $p=1, \ldots, k$
$$
\beta_p=s_{i_1}\ldots s_{i_{p-1}}\alpha_{i_p},~\beta_p'=s_{i_k}\ldots s_{i_{p+1}}\alpha_{i_p},
$$
$$
X_{\beta_p}^\pm=T_{i_1}\ldots T_{i_{p-1}}X_{i_p}^\pm,~{\overline{X}_{\beta_p'}^\pm}=T_{i_k}^{-1}\ldots T_{i_{p+1}}^{-1}X_{i_p}^\pm,~ K_{\beta_p'}=T_{i_k}^{-1}\ldots T_{i_{p+1}}^{-1}K_{i_p}, 
$$
$$
\theta_{\beta_p}={\exp}_{q_{\beta_p}}[(1-q_{\beta_p}^{-2})X_{\beta_p}^+\otimes X_{\beta_p}^-], ~\overline{\theta}_{\beta_p'}={\exp}_{q_{\beta_p'}}[(1-q_{\beta_p'}^{-2})K_{\beta_p'}^{-1}{\overline{X}_{\beta_p'}^-}\otimes {\overline{X}_{\beta_p'}^+}K_{\beta_p'}]. \index[not]{t@$\theta_{\beta_p}$} \index[not]{t@$\overline{\theta}_{\beta_p'}$}  
$$

Note that for a reduced decomposition $w=s_{i_1}\ldots s_{i_k}$ one has $T_{i_1}^{-1}\ldots T_{i_k}^{-1}=(T_{i_k}\ldots T_{i_1})^{-1}=T_{w^{-1}}^{-1}$, and $T_{w^{-1}}$ only depends on $w$. Therefore $\overline{T}_w:=T_{i_1}^{-1}\ldots T_{i_k}^{-1}=T_{w^{-1}}^{-1}$ \index[not]{T@$\overline{T}_w$} only depends on $w$ and (\ref{DT-1}) yields
\begin{equation}\label{dhotw}
\Delta_h(\overline{T}_w)=\prod^{k}_{p=1}\overline{\theta}_{\beta_p}'\overline{T}_w\otimes \overline{T}_w=\overline{T}_w\otimes \overline{T}_w\prod_{p=1}^k{\theta}_{\beta_p'}',
\end{equation}
where in the products $\overline{\theta}_{\beta_p}'$ (resp. ${\theta}_{\beta_p'}'$) appears on the left from $\overline{\theta}_{\beta_q}'$ (resp. ${\theta}_{\beta_q'}'$) if $p<q$, and for $p=1,\ldots, k$
$$
{\overline{X}_{\beta_p}^\pm}=T_{i_1}^{-1}\ldots T_{i_{p-1}}^{-1}X_{i_p}^\pm,~{{X}_{\beta_p'}^\pm}'=T_{i_k}\ldots T_{i_{p+1}}X_{i_p}^\pm,~ \overline{K}_{\beta_p}=T_{i_1}^{-1}\ldots T_{i_{p-1}}^{-1}K_{i_p}, 
$$
$$
{\theta}_{\beta_p'}'={\exp}_{q_{\beta_p'}^{-1}}[(1-q_{\beta_p'}^{2}){{X}_{\beta_p'}^+}'\otimes {{X}_{\beta_p'}^-}'], ~\overline{\theta}_{\beta_p}'={\exp}_{q_{\beta_p}^{-1}}[(1-q_{\beta_p}^{2})\overline{K}_{\beta_p}^{-1}{\overline{X}_{\beta_p}^-}\otimes {\overline{X}_{\beta_p}^+}\overline{K}_{\beta_p}]. \index[not]{t@${\theta}_{\beta_p'}'$} \index[not]{t@$\overline{\theta}_{\beta_p}'$} 
$$

If $w\alpha_i=\alpha_j$ for some $i$ and $j$ then by Proposition 8.1.6. in \cite{ChP}
\begin{equation}\label{treda}
T_wX_i^\pm=X_j^\pm,~\overline{T}_wX_i^\pm=X_j^\pm.
\end{equation}


\section{Quantum root vectors}\label{QRvect}

\pagestyle{myheadings}
\markboth{CHAPTER \thechapter.~QUANTUM GROUPS}{\thesection.~QUANTUM ROOT VECTORS}

\setcounter{equation}{0}
\setcounter{theorem}{0}

In this section we recall the construction of analogues of root vectors for $U_h({\frak g})$  in terms of the braid group action on $U_h({\frak g})$.
Recall that for any reduced decomposition $\overline{w}=s_{i_1}\ldots s_{i_D}$ of the longest element $\overline{w}$ of the Weyl group $W$ of $\g$ the ordering
$$
\beta_1=\alpha_{i_1},\beta_2=s_{i_1}\alpha_{i_2},\ldots,\beta_D=s_{i_1}\ldots s_{i_{D-1}}\alpha_{i_D}
$$
is a normal ordering in $\Delta_+$, and there is a one--to--one correspondence between normal orderings of $\Delta_+$ and reduced decompositions of $\overline{w}$.

Fix a reduced decomposition $\overline{w}=s_{i_1}\ldots s_{i_D}$ of $\overline{w}$ and define the corresponding {\it quantum root vectors} \index{root!vector!quantum} in $U_h({\frak g})$ by
\begin{equation}\label{rootvect}
X_{\beta_k}^\pm=T_{i_1}\ldots T_{i_{k-1}}X_{i_k}^\pm. \index[not]{X@$X_{\beta_k}^\pm$}
\end{equation}

We also define
$$
K_{\beta_k}=T_{i_1}\ldots T_{i_{k-1}}K_{i_k}. \index[not]{K@$K_{\beta_k}$}
$$

Note that one can construct root vectors in the Lie algebra $\g$ in a similar way. 
Namely, the root vectors $X_{\pm \beta_k}\in \g_{\pm \beta_k}$ of $\g$ can be defined by
\begin{equation}\label{rootvectg}
X_{\pm \beta_k}=T_{i_1}\ldots T_{i_{k-1}}X_{\pm \alpha_{i_k}},
\end{equation}
where $X_{\pm \alpha_{i_k}}$ are as in (\ref{BGactg}).

The quantum root vectors $X_{\beta}^\pm$ satisfy the following relations:
\begin{eqnarray}\label{qcom}
X_{\alpha}^\pm X_{\beta}^\pm - q^{\left\langle \alpha,\beta\right\rangle}X_{\beta}^\pm X_{\alpha}^\pm= \sum_{m_1,\ldots, m_k\in \mathbb{N}}C(m_1,\ldots,m_k)
{(X_{\zeta_1}^\pm)}^{m_1}{(X_{\zeta_2}^\pm)}^{m_2}\ldots {(X_{\zeta_k}^\pm)}^{m_k}= \nonumber \\
=\sum_{m_1,\ldots, m_k\in \mathbb{N}}C'(m_1,\ldots,m_k)
{(X_{\zeta_1}^\pm)}^{(m_1)}{(X_{\zeta_2}^\pm)}^{(m_2)}\ldots {(X_{\zeta_k}^\pm)}^{(m_k)},~\alpha<\beta, \index{commutation relations!for quantum root vectors}
\end{eqnarray}
where $\alpha<\zeta_1<\ldots<\zeta_k<\beta$, $[\alpha,\beta]=\{\alpha,\zeta_1,\ldots,\zeta_k,\beta\}$ as a set, for any $\alpha \in \Delta_+$ we put ${(X_{\alpha}^\pm)}^{(k)}=\frac{(X_\alpha^\pm)^{k}}{[k]_{q_\alpha}!}$, $k\geq 0$, \index[not]{X@${(X_{\alpha}^\pm)}^{(k)}$} $q_\alpha =q^{d_i}$ if the positive root $\alpha$ is Weyl group conjugate to the simple root $\alpha_i$, $C'(m_1,\ldots,m_k)\in {\Bbb C}[q,q^{-1}]$, $C(m_1,\ldots,m_k)\in {\mathcal{P}}$, only finitely many of these coefficients are non--zero, and ${\mathcal{P}}$ is the algebra defined by
\begin{equation}
{\mathcal{P}}=\left\{ \begin{array}{ll} {\Bbb C}[q,q^{-1}] & \text{ if $\g$ is simply-laced} \\
{\Bbb C}[q,q^{-1}, \frac{1}{[2]_q}] & \text{ if $\g$ is of type $B_l$, $C_l$ or $F_4$} \\ {\Bbb C}[q,q^{-1}, \frac{1}{[2]_q}, \frac{1}{[3]_q}] & \text{ if $\g$ is of type $G_2$}
\end{array}\right. . \index[not]{P@$\mathcal{P}$}
\end{equation}
\index{Lie!algebra!simply-laced}

The fact that the coefficients $C(m_1,\ldots,m_k)$ in the right hand sides of formulas (\ref{qcom}) belong to the algebra ${\mathcal{P}}$ was noted in \cite{Dr}, Lemma 1.1.1.

If $X=c_1H_1+\ldots +c_lH_l$ for some $c_1,\ldots, c_l\in \mathbb{R}$ then from commutation relations (\ref{defQGrel}) and from the fact that the action of the braid group on the elements $H_i$ coincides with the action of the Weyl group on the corresponding simple root generators of $\h$ we obtain 
\begin{equation}\label{Ccomm}
e^{hX}X_\alpha^\pm e^{-hX}=q^{\pm \alpha(X)}X_\alpha^\pm,
\end{equation}
where in the expression $\alpha(X)$ $X$ is regarded as an element of $\h$ under the natural identification of the elements $H_1,\ldots, H_l$ with the simple root generators of $\h$.

In particular, (\ref{Ccomm}) implies that for any $k\in \mathbb{N}$
\begin{equation}\label{Ccomm1}
(X_i^\pm K_i)^k =q_i^{\pm k(k-1)}(X_i^\pm)^k K_i^k. 
\end{equation}

Note that by construction
\begin{equation}\label{qrootclass}
\begin{array}{l}
X_\beta^+~(\mbox{mod }h)=X_\beta \in {\frak g}_\beta,\\
\\
X_\beta^-~(\mbox{mod }h)=X_{-\beta} \in {\frak g}_{-\beta}
\end{array}
\end{equation}
are root vectors of $\frak g$.

Define an algebra involution $\tau_0$ \index[not]{t@$\tau_0$} of $U_h({\frak g})$ by 
$$
\tau_0(X_i^\pm)=X_i^\mp, \tau_0(H_i)=-H_i, \tau_0(h)=-h.
$$
It satisfies the relations $T_i^{-1}=\tau_0 T_i\tau_0$ and hence for any $\alpha\in \Delta_+$
$$
\tau_0(X_{\alpha}^\pm)=\overline{X}_{\alpha}^\mp, 
$$
where
$$
\overline{X}_{\beta_k}^\pm=T_{i_1}^{-1}\ldots T_{i_{k-1}}^{-1}X_{i_k}^\pm. \index[not]{X@$\overline{X}_{\beta_k}^\pm$}
$$

By applying $\tau_0$ to relations (\ref{qcom}) one can obtain the following relations for the quantum root vectors $\overline{X}_{\beta}^\pm$
\begin{eqnarray}\label{qcom1}
\overline{X}_{\alpha}^\pm \overline{X}_{\beta}^\pm - q^{-\left\langle \alpha,\beta\right\rangle}\overline{X}_{\beta}^\pm \overline{X}_{\alpha}^\pm= \sum_{m_1,\ldots, m_k\in \mathbb{N}}D(m_1,\ldots,m_k)
{(\overline{X}_{\zeta_1}^\pm)}^{m_1}{(\overline{X}_{\zeta_2}^\pm)}^{m_2}\ldots {(\overline{X}_{\zeta_k}^\pm)}^{m_k}= \nonumber \\
=\sum_{m_1,\ldots, m_k\in \mathbb{N}}D'(m_1,\ldots,m_k)
{(\overline{X}_{\zeta_1}^\pm)}^{(m_1)}{(\overline{X}_{\zeta_2}^\pm)}^{(m_2)}\ldots {(\overline{X}_{\zeta_k}^\pm)}^{(m_k)},~\alpha<\beta, \index{commutation relations!for quantum root vectors}
\end{eqnarray}
where $\alpha<\zeta_1<\ldots<\zeta_k<\beta$, $[\alpha,\beta]=\{\alpha,\zeta_1,\ldots,\zeta_k,\beta\}$ as a set, for any $\alpha \in \Delta_+$ we put ${(\overline{X}_{\alpha}^\pm)}^{(k)}=\frac{(\overline{X}_\alpha^\pm)^{k}}{[k]_{q_\alpha}!}$, \index[not]{X@${(\overline{X}_{\alpha}^\pm)}^{(k)}$} $k\geq 0$, $D(m_1,\ldots,m_k)\in {\mathcal{P}}$, $D'(m_1,\ldots,m_k)\in {\Bbb C}[q,q^{-1}]$, and only finitely many of these coefficients are non--zero.

One can also obtain commutation relations between positive and negative quantum root vectors. These relations are known in some form. For completeness we give a proof of them using (\ref{qcom}) and (\ref{qcom1}) only.
\begin{lemma}
Let $[-\beta,\alpha]$, $\alpha,\beta\in \Delta_+$ be a minimal segment with respect to the circular normal ordering of $\Delta$ corresponding to a normal ordering $\beta_1,\ldots, \beta_D$ of $\Delta_+$, so that $[-\beta,\alpha]=\{-\beta,\zeta_1,\ldots,\zeta_k,\alpha\}$ as a set, where $-\beta<\zeta_1<\ldots<\zeta_k<\alpha$ and the inequalities for the roots are with respect to the circular normal ordering in $\Delta$ corresponding to the normal ordering $\beta_1,\ldots, \beta_D$ of $\Delta_+$. Then the corresponding quantum root vectors satisfy the following relations
\begin{eqnarray}\label{qcom2}
X_{\alpha}^+ X_{\beta}^- - X_{\beta}^- X_{\alpha}^+= \sum_{m_1,\ldots, m_k\in \mathbb{N}}\bar{C}(m_1,\ldots,m_k)
{(X_{\zeta_1})}^{m_1}{(X_{\zeta_2})}^{m_2}\ldots {(X_{\zeta_k})}^{m_k}= \nonumber \\
=\sum_{m_1,\ldots, m_k\in \mathbb{N}}\bar{C}'(m_1,\ldots,m_k)
{(X_{\zeta_1})}^{(m_1)}{(X_{\zeta_2})}^{(m_2)}\ldots {(X_{\zeta_k})}^{(m_k)}, 
\end{eqnarray}
where in the right hand side $X_{\zeta}=X_{\zeta}^+$ for $\zeta\in \Delta_+$ and  $X_{\zeta}=X_{\zeta}^-$ for $\zeta\in \Delta_-$, $\bar{C}'(m_1,\ldots,m_k)\in U_q(H)$, $U_q(H)$ \index[not]{U@$U_q(H)$} is the ${\Bbb C}[q,q^{-1}]$--subalgebra of $U_h(\g)$ generated by $K_i^{\pm 1}$, $i=1,\ldots ,l$, $\bar{C}(m_1,\ldots,m_k)\in {\mathcal{P}'}$, where ${\mathcal{P}'}$ \index[not]{P@${\mathcal{P}'}$} is the ${\mathcal{P}}$--algebra generated by $U_q(H)$, and and only finitely many of the coefficients $\bar{C}'(m_1,\ldots,m_k)$ and $\bar{C}(m_1,\ldots,m_k)$ are non--zero.

Also
\begin{eqnarray}\label{qcom3}
\overline{X}_{\alpha}^+ \overline{X}_{\beta}^- - \overline{X}_{\beta}^- \overline{X}_{\alpha}^+= \sum_{m_1,\ldots, m_k\in \mathbb{N}}\bar{D}(m_1,\ldots,m_k)
{(\overline{X}_{\zeta_1})}^{m_1}{(\overline{X}_{\zeta_2})}^{m_2}\ldots {(\overline{X}_{\zeta_k})}^{m_k}= \nonumber \\
=\sum_{m_1,\ldots, m_k\in \mathbb{N}}\bar{D}'(m_1,\ldots,m_k)
{(\overline{X}_{\zeta_1})}^{(m_1)}{(\overline{X}_{\zeta_2})}^{(m_2)}\ldots {(\overline{X}_{\zeta_k})}^{(m_k)},~\alpha<\beta, 
\end{eqnarray}
where $\overline{X}_{\zeta}=\overline{X}_{\zeta}^+$ for $\zeta\in \Delta_+$ and  $\overline{X}_{\zeta}=\overline{X}_{\zeta}^-$ for $\zeta\in \Delta_-$, $\bar{D}'(m_1,\ldots,m_k)\in U_q(H)$, $\bar{D}(m_1,\ldots,m_k)\in {\mathcal{P}}'$, and only finitely many of these coefficients are non--zero.
\end{lemma}

\begin{proof}
The proof is by induction over the length of the positive part of the segment $[-\beta, \alpha]$, i.e. over the cardinality of the set $[-\beta, \alpha]\cap \Delta_+$. We shall consider the first identity in the case when $\alpha<\beta$. The others are proved in a similar way. 

Let $\overline{w}=s_{i_1}\ldots s_{i_D}$ be the reduced decomposition of the longest element of the Weyl group corresponding to the normal ordering $\beta_1,\ldots, \beta_D$ of $\Delta_+$. 

First assume that $\alpha=\beta_1=\alpha_{i_1}$. Let $\beta=\beta_m=s_{i_1}\ldots s_{i_{m-1}}\alpha_{i_m}$. 
Then $s_{i_1}^{-1}\beta_2,\ldots, s_{i_1}^{-1}\beta_D, \alpha_{i_1}$ is another normal ordering of $\Delta_+$, and by (\ref{qcom}) for this normal ordering
\begin{align*}
T_{i_1}^{-1}(X_{\alpha}^+ X_{\beta}^- - X_{\beta}^- X_{\alpha}^+)=T_{i_1}^{-1}(X_{i_1}^+ X_{\beta_m}^- - X_{\beta_m}^- X_{i_1}^+)= \\
=K_{i_1}^{-1}(-X_{i_1}^- X_{s_{i_1}^{-1}\beta_m}^- +q^{-\left\langle \alpha_{i_1}, s_{i_1}^{-1}\beta_m\right\rangle} X_{s_{i_1}^{-1}\beta_m}^- X_{i_1}^-)= \\
=\sum_{m_1,\ldots, m_k\in \mathbb{N}}K_{i_1}^{-1}C(m_1,\ldots,m_k)
{(X_{\zeta_1})}^{m_1}{(X_{\zeta_2})}^{m_2}\ldots {(X_{\zeta_k})}^{m_k},
\end{align*}
where $[-s_{i_1}^{-1}\beta_m,-\alpha_{i_1}]=\{-s_{i_1}^{-1}\beta_m,\zeta_1, \ldots, \zeta_k,-\alpha_{i_1}\}$ as a set, $-s_{i_1}^{-1}\beta_m<\zeta_1<\ldots<\zeta_k<-\alpha_{i_1}$, the inequalities for the roots are with respect to the circular normal ordering of $\Delta$ associated to the ordering $s_{i_1}^{-1}\beta_2,\ldots, s_{i_1}^{-1}\beta_D, \alpha_{i_1}$ of $\Delta_+$, $C(m_1,\ldots,m_k)\in {\mathcal{P}}$, only finitely many of these coefficients are non--zero, and the quantum root vectors are defined using the ordering $s_{i_1}^{-1}\beta_2,\ldots, s_{i_1}^{-1}\beta_D, \alpha_{i_1}$ of $\Delta_+$.

Now applying $T_{i_1}$ to the last identity we get
$$
X_{\alpha}^+ X_{\beta}^- - X_{\beta}^- X_{\alpha}^+= 
=\sum_{m_1,\ldots, m_k\in \mathbb{N}}K_{i_1}C(m_1,\ldots,m_k)
{(X_{\zeta_1})}^{m_1}{(X_{\zeta_2})}^{m_2}\ldots {(X_{\zeta_k})}^{m_k},
$$
where $[-\beta,\alpha]=\{-\beta,\zeta_1,\ldots,\zeta_k,\alpha\}$ as a set, $-\beta<\zeta_1<\ldots<\zeta_k<\alpha$, the inequalities for the roots are with respect to the circular normal ordering in $\Delta$ corresponding to the normal ordering $\beta_1,\ldots, \beta_D$ of $\Delta_+$, and the quantum root vectors are defined using the ordering $\beta_1,\ldots, \beta_D$ of $\Delta_+$. This establishes the base of the induction.

Now assume that the identity in question is proved for all normal orderings of $\Delta_+$ and for all $\alpha=\beta_k$ with $k<n$ for some $n>0$ and for all possible $\beta$ such that $[-\beta,\alpha]$, $\alpha,\beta\in \Delta_+$ is a minimal segment.

Let $\alpha=\beta_n=s_{i_1}\ldots s_{i_{n-1}}\alpha_{i_n}$, $\beta=\beta_m=s_{i_1}\ldots s_{i_{m-1}}\alpha_{i_m}$, $n<m$. 
Then $s_{i_1}^{-1}\beta_2,\ldots, s_{i_1}^{-1}\beta_D, \alpha_{i_1}$ is another normal ordering of $\Delta_+$, and by the induction hypothesis for this normal ordering with $s_{i_1}^{-1}\alpha=s_{i_2}\ldots s_{i_{n-1}}\alpha_{i_n}$, $s_{i_1}^{-1}\beta=\beta_m=s_{i_2}\ldots s_{i_{m-1}}\alpha_{i_m}$ we have
\begin{align*}
X_{s_{i_1}^{-1}\alpha}^+ X_{s_{i_1}^{-1}\beta}^- - X_{s_{i_1}^{-1}\beta}^- X_{s_{i_1}^{-1}\alpha}^+= \\
=\sum_{m_1',\ldots, m_k'\in \mathbb{N}}\bar{C}(m_1',\ldots,m_k')
{(X_{\zeta_1'})}^{m_1'}{(X_{\zeta_2'})}^{m_2'}\ldots {(X_{\zeta_k'})}^{m_k'},
\end{align*}
where $[-s_{i_1}^{-1}\beta_m,s_{i_1}^{-1}\alpha]=\{-s_{i_1}^{-1}\beta_m,\zeta_1',\ldots,\zeta_k',s_{i_1}^{-1}\alpha\}$ as a set, $-s_{i_1}^{-1}\beta_m<\zeta_1'<\ldots<\zeta_k'<s_{i_1}^{-1}\alpha$, the inequalities for the roots are with respect to the circular normal ordering in $\Delta$ associated to the ordering $s_{i_1}^{-1}\beta_2,\ldots, s_{i_1}^{-1}\beta_D, \alpha_{i_1}$ of $\Delta_+$, $\bar{C}(m_1',\ldots, m_k')\in {\mathcal{P}}'$, only finitely many of these coefficients are non--zero, and the quantum root vectors are defined using the ordering $s_{i_1}^{-1}\beta_2,\ldots, s_{i_1}^{-1}\beta_D, \alpha_{i_1}$ of $\Delta_+$.

Now applying $T_{i_1}$ to the last identity we get
\begin{align*}
X_{\alpha}^+ X_{\beta}^- - X_{\beta}^- X_{\alpha}^+=\sum_{m_1,\ldots, m_k\in \mathbb{N}}\bar{C}''(m_1,\ldots,m_k)
{(X_{\zeta_1})}^{m_1}{(X_{\zeta_2})}^{m_2}\ldots {(X_{\zeta_k})}^{m_k},
\end{align*}
where $[-\beta,\alpha]=\{-\beta,\zeta_1,\ldots,\zeta_k,\alpha\}$ as a set, where $-\beta<\zeta_1<\ldots<\zeta_k<\alpha$ and the inequalities for the roots are with respect to the circular normal ordering in $\Delta$ corresponding to the normal ordering $\beta_1,\ldots, \beta_D$ of $\Delta_+$, the quantum root vectors are defined using the ordering $\beta_1,\ldots, \beta_D$ of $\Delta_+$, $\bar{C}''(m_1,\ldots,m_k)\in {\mathcal{P}}'$, and only finitely many of these coefficients are non--zero. This establishes the induction step and completes the proof.

\end{proof}

Define an algebra anti-involution $\omega$ \index[not]{o@$\omega$} of $U_h({\frak g})$ by
\begin{equation}\label{omega}
\omega(X_i^\pm)=X_i^\mp, \omega(H_i)=H_i, \omega(h)=-h.
\end{equation}
It commutes with the braid group action and for any $\alpha\in \Delta_+$ satisfies
\begin{equation}\label{omegaX}
\omega(X_{\alpha}^+)=X_{\alpha}^-.
\end{equation}

Define also an algebra anti-involution $\omega_0'$ \index[not]{o@$\omega_0'$} of $U_h({\frak g})$ by
\begin{equation}\label{omega0'}
\omega_0'(X_i^\pm)=X_i^\pm, \omega_0'(H_i)=-H_i, \omega_0'(h)=-h.
\end{equation}
It satisfies
$$
\omega\omega_0'=\omega_0'\omega,
$$
$$
\omega_0'(T_iX_j^\pm)=(-1)^{a_{ij}}q_i^{\mp a_{ij}}T_i(\omega_0' X_j^\pm)=(-1)^{a_{ij}}q_i^{\mp a_{ij}}T_i(X_j^\pm), i\neq j,
$$
$$
\omega_0'(T_iX_i^+)=q_i^{-2}T_i(\omega_0' X_i^+)=q_i^{-2}T_i(X_i^+),\omega_0'(T_iX_i^-)=q_i^{2}T_i(\omega_0' X_i^-)=q_i^{2}T_i(X_i^-).
$$
As a consequence we obtain that if $X$ is a homogeneous polynomial in quantum simple root vectors then
$$
T_i(\omega_0'X)=c_X\omega_0'(T_iX), 
$$
where $c_X=\epsilon p$, $\epsilon= \pm 1$, $p\in q^{\mathbb{Z}}$, and hence
\begin{equation}\label{Xo0}
\omega_0'(X_\alpha^+)=c_\alpha X_\alpha^+,	
\end{equation}
where $c_{\alpha}=\epsilon_\alpha p_\alpha$, \index[not]{c@$c_\alpha$} $\epsilon_\alpha= \pm 1$, $p_\alpha\in q^{\mathbb{Z}}$. 
We also have
\begin{equation}\label{X-o0}
\omega_0'(X_\alpha^-)=\omega_0'\omega(X_\alpha^+)=\omega\omega_0'(X_\alpha^+)=\omega(c_\alpha X_\alpha^+)=c_\alpha^{-1}\omega(X_\alpha^+)=c_\alpha^{-1}X_\alpha^-	
\end{equation}

We shall also need the algebra involution $\tau$ \index[not]{t@$\tau$} of $U_h({\frak g})$defined by 
$$
\tau(X_i^\pm)=X_i^\pm, \tau(H_i)=H_i, \tau(h)=-h.
$$

Then the anti-involution $\omega_0:=\tau\omega_0'$ \index[not]{o@$\omega_0$} satisfies 
\begin{equation}\label{omega0}
\omega_0(X_i^\pm)=X_i^\pm, \omega_0(H_i)=-H_i, \omega_0(h)=h.
\end{equation}

Note that the anti-involutions $\omega$, $\omega_0'$, $\omega_0$ and the involution $\tau$ commute with each other.

For any braid group automorphism $T$ of $U_h(\g)$ we denote $T^\tau=\tau T \tau$. \index[not]{T@$T^\tau$}
Denote the natural extension of $\tau$ to $\mathbb{C}_h[G]^*$ by the same letter. Then the  automorphism $T_i^\tau$, $i=1,\ldots, l$ is given by conjugation of elements of $U_h(\g)\subset \mathbb{C}_h[G]^*$ by the invertible element $\tau(T_i)\in \mathbb{C}_h[G]^*$, where $T_i$ is given by (\ref{T1}). 

\begin{remark}
In \cite{Lus} the automorphisms $T_i$, $T_i^{-1}$, $T_i^\tau$, ${T_i^\tau}^{-1}$ are denoted $T''_{i,1}$, $T'_{i,-1}$, $T''_{i,-1}$, $T'_{i,1}$, respectively.
\end{remark}


\section{Some subalgebras in quantum groups and their Poincar\'{e}--Birkhoff--Witt bases}\label{PBWbases}

\pagestyle{myheadings}
\markboth{CHAPTER \thechapter.~QUANTUM GROUPS}{\thesection.~ POINCAR\'{E}--BIRKHOFF-WITT BASES}

\setcounter{equation}{0}
\setcounter{theorem}{0}

Now we shall explicitly describe a topological  ${\Bbb C}[[h]]$--basis for $U_h({\frak g})$. We shall also recall the definition of some rational forms of $U_h({\frak g})$ and of their bases.

Denote by $U_h({\frak n}_+)$ \index[not]{U@$U_h(\n_+)$}, $U_h({\frak n}_-)$ \index[not]{U@$U_h(\n_-)$} and $U_h(\h)$ \index[not]{U@$U_h(\h)$} the ${\Bbb C}[[h]]$--subalgebras of $U_h({\frak g})$ topologically generated by the
$X_i^+$, by the $X_i^-$ and by the $H_i$, respectively. For any $\alpha\in \Delta_+$ one has $X_{\alpha}^\pm\in U_h({\frak n}_\pm)$.
From the definition of the quantum root vectors it also follows that $[H_i,X_{\alpha}^\pm]=\pm \alpha(H_i)X_{\alpha}^\pm$, $i=1,\ldots, l$. Therefore using the uniqueness of the presentation of any positive root as a sum of simple roots we immediately deduce the following property of the quantum root vectors.
\begin{proposition}\label{rootprop}
For $\beta =\sum_{i=1}^lm_i\alpha_i,~m_i\in {\Bbb N}$  $X_{\beta}^\pm $ is a polynomial in
the noncommutative variables $X_i^\pm$ homogeneous in each $X_i^\pm$ of degree $m_i$.
\end{proposition}

Denote by $U_q^{res}(\g)$ \index[not]{U@$U_q^{res}(\g)$} the subalgebra in $U_h(\g)$ generated over ${\mathbb{C}[q,q^{-1}]}$ by the elements 
$$
K_i^{\pm 1},~(X_i^\pm)^{(k)},~i=1,\ldots ,l, k\geq 1.
$$ 

The elements
\begin{equation}\label{Kic}
\left[ \begin{array}{l}
K_i;c \\
r
\end{array} \right]_{q_i}=\prod_{s=1}^r \frac{K_i q_i^{c+1-s}-K_i^{-1}q_i^{s-1-c}}{q_i^s-q_i^{-s}}~,~i=1,\ldots,l,~c\in \mathbb{Z},~r\in \mathbb{N} \index[not]{ZZZ@$\left[ \begin{array}{l}
K_i;c \\
r
\end{array} \right]_{q_i}$}
\end{equation}
belong to $U_q^{res}(\g)$. Denote by $U_q^{res}(H)$ \index[not]{U@$U_q^{res}(H)$} the subalgebra of $U_q^{res}(\g)$ generated by these elements and by $K_i^{\pm 1}$, $i=1,\ldots,l$.

The subalgebras $U_q(H)$, $U_q^{res}(\g)$ and $U_q^{res}(H)$ of $U_h(\g)$ are invariant under the braid group action.

Let $U_\mathcal{P}({\frak n}_+)$, \index[not]{U@$U_\mathcal{P}(\n_+)$} $U_\mathcal{P}({\frak n}_-)$, \index[not]{U@$U_\mathcal{P}(\n_-)$} (resp. $U_q^{res}({\frak n}_+)$, \index[not]{U@$U_q^{res}(\n_+)$} and $U_q^{res}({\frak n}_-)$) \index[not]{U@$U_q^{res}(\n_-)$} be the $\mathcal{P}$ (resp. $\mathbb{C}[q,q^{-1}]$)--subalgebras of $U_h({\frak g})$ (resp. $U_q^{res}(\g)$) generated by the
$X_i^+$, by the $X_i^-$, $i=1,\ldots,l$, (resp. by the $(X_i^+)^{(r)}$, and by the $(X_i^-)^{(r)}$, $i=1,\ldots,l$, $r\geq 0$), respectively.
Denote also by $U_h({\frak b}_\pm)$ \index[not]{U@$U_h(\b_\pm)$} (resp. $U_q^{res}({\frak b}_\pm)$) \index[not]{U@$U_q^{res}(\b_\pm)$} the ${\Bbb C}[[h]]$ (resp. $\mathbb{C}[q,q^{-1}]$)--subalgebras of $U_h(\g)$ (resp. $U_q^{res}(\g)$) generated by $U_h(\n_\pm)$ and by $U_h(\h)$ (resp. by $U_q^{res}(\n_\pm)$ and by $U_q^{res}(H)$).

Using the root vectors $X_{\beta}^\pm$ and the elements $(X_{\beta}^\pm)^{(r)}$  we can construct bases for these subalgebras.
Namely, let $(X^\pm)^{\bf r}:=(X^\pm_{\beta_1})^{r_1}\ldots (X^\pm_{\beta_D})^{r_D}$, \index[not]{X@$(X^\pm)^{\bf r}$} $(X^\pm)^{(\bf r)}:=(X^\pm_{\beta_1})^{(r_1)}\ldots (X^\pm_{\beta_D})^{(r_D)}$, \index[not]{X@$(X^\pm)^{(\bf r)}$} ${\bf r}:=(r_1,\ldots, r_D)\in {\Bbb N}^D$, $H^{\bf m}:=H_1^{m_1}\ldots H_l^{m_l}$, \index[not]{H@$H^{\bf m}$} ${\bf m}:=(m_1,\ldots , m_l)\in {\Bbb N}^l$. 

Commutation relations (\ref{qcom}), (\ref{qcom1}), (\ref{qcom2}) and (\ref{qcom3}) between quantum root vectors imply the following lemma.

\begin{lemma}\label{segmPBW}
(i) The elements $(X^+)^{\bf r}$, $(X^-)^{\bf r}$, $(X^+)^{\bf (r)}$, and $(X^-)^{\bf (r)}$ for ${\bf r}\in {\Bbb N}^D$
form bases of $U_\mathcal{P}({\frak n}_+)$, $U_\mathcal{P}({\frak n}_-)$, $U_q^{res}({\frak n}_+)$, and $U_q^{res}({\frak n}_-)$, \index{quantum group!basis}
respectively.

(ii) The elements $(X^+)^{\bf r}$, $(X^-)^{\bf r}$ and $H^{\bf m}$ form topological bases of $U_h({\frak n}_+), U_h({\frak n}_-)$ and $U_h(\h)$, respectively.

(iii) The multiplication defines isomorphisms of ${\mathbb{C}[q,q^{-1}]}$--modules:
\begin{equation}\label{resPBW}
U_q^{res}({\frak n}_+)\otimes U_q^{res}(H) \otimes U_q^{res}({\frak n}_-)\rightarrow U_q^{res}({\frak g}),	
\end{equation}
$$
U_q^{res}({\frak n}_-)\otimes U_q^{res}(H) \otimes U_q^{res}({\frak n}_+)\rightarrow U_q^{res}({\frak g}),
$$
and of complete ${\Bbb C}[[h]]$--modules
$$
U_h({\frak n}_+)\otimes U_h(\h) \otimes U_h({\frak n}_-)\rightarrow U_h({\frak g}),
$$
$$
U_h({\frak n}_-)\otimes U_h(\h) \otimes U_h({\frak n}_+)\rightarrow U_h({\frak g}),
$$
where the tensor products in the left hand side are completed in the $h$--adic topology.

(iv) Let $[\alpha,\beta]$=$\{\beta_p,\ldots , \beta_q\}$ be a minimal segment in $\Delta_+$,
$U_\mathcal{P}([\alpha,\beta])$, \index[not]{U@$U_\mathcal{P}([\pm\alpha,\pm\beta])$} $U_\mathcal{P}([-\alpha,-\beta])$ (resp. $U_q^{res}([\alpha,\beta])$, \index[not]{U@$U_q^{res}([\pm\alpha,\pm\beta])$} $U_q^{res}([-\alpha,-\beta])$) the $\mathcal{P}$ (resp. $\mathbb{C}[q,q^{-1}]$)--subalgebras of $U_h({\frak g})$
generated by the
$X_\gamma^+$ and by the $X_\gamma^-$, $\gamma \in [\alpha, \beta]$ (resp. by the $(X_\gamma^+)^{(r)}$ and by the $(X_\gamma^-)^{(r)}$, $\gamma \in [\alpha, \beta]$, $r\geq 0$), respectively.
Then $U_\mathcal{P}([\alpha,\beta])\subset U_\mathcal{P}(\n_+)$, $U_\mathcal{P}([-\alpha,-\beta])\subset U_\mathcal{P}(\n_-)$ (resp. $U_q^{res}([\alpha,\beta])\subset U_q^{res}(\n_+)$, $U_q^{res}([-\alpha,-\beta])\subset U_q^{res}(\n_-)$), and the elements $(X^\pm_{\beta_p})^{r_p}\ldots (X^\pm_{\beta_q})^{r_q}$ (resp. $(X^\pm_{\beta_p})^{(r_p)}\ldots (X^\pm_{\beta_q})^{(r_q)}$), ${ r_i}\in {\Bbb N}$
form bases of $U_\mathcal{P}([\alpha,\beta])$, $U_\mathcal{P}([-\alpha,-\beta])$ (resp. $U_q^{res}([\alpha,\beta])$, $U_q^{res}([-\alpha,-\beta])$).

(v) Let 
$\overline{U}_\mathcal{P}([\alpha,\beta])$, \index[not]{U@$\overline{U}_\mathcal{P}([\pm\alpha,\pm\beta])$} $\overline{U}_\mathcal{P}([-\alpha,-\beta])$ (resp. $\overline{U}_q^{res}([\alpha,\beta])$, \index[not]{U@$\overline{U}_q^{res}([\pm\alpha,\pm\beta])$} $\overline{U}_q^{res}([-\alpha,-\beta])$) be the $\mathcal{P}$ (resp. $\mathbb{C}[q,q^{-1}]$)--subalgebras of $U_h({\frak g})$
generated by the $\overline{X}_\gamma^+$ and by the $\overline{X}_\gamma^-$, $\gamma \in [\alpha, \beta]$ (resp. by the $(\overline{X}_\gamma^+)^{(r)}$ and by the $(\overline{X}_\gamma^-)^{(r)}$, $\gamma \in [\alpha, \beta]$, $r\geq 0$), respectively.
Then $\overline{U}_\mathcal{P}([\alpha,\beta])\subset U_\mathcal{P}(\n_+)$, $\overline{U}_\mathcal{P}([-\alpha,-\beta])\subset U_\mathcal{P}(\n_-)$ (resp. $\overline{U}_q^{res}([\alpha,\beta])\subset U_q^{res}(\n_+)$, $\overline{U}_q^{res}([-\alpha,-\beta])\subset U_q^{res}(\n_-)$), and the elements $(\overline{X}^\pm_{\beta_p})^{r_p}\ldots (\overline{X}^\pm_{\beta_q})^{r_q}$ (resp. $(\overline{X}^\pm_{\beta_p})^{(r_p)}\ldots (\overline{X}^\pm_{\beta_q})^{(r_q)}$), ${ r_i}\in {\Bbb N}$
form bases of $\overline{U}_\mathcal{P}([\alpha,\beta])$, $\overline{U}_\mathcal{P}([-\alpha,-\beta])$ (resp. $\overline{U}_q^{res}([\alpha,\beta])$, $\overline{U}_q^{res}([-\alpha,-\beta])$).

(vi) Let $[\alpha,-\beta]$=$\{\beta_p,\ldots , -\beta_q\}$, $\alpha, \beta \in \Delta_+$ be a minimal segment in $\Delta$,
$U_{\mathcal{P}'}([\alpha,-\beta])$, \index[not]{U@$U_{\mathcal{P}'}([\alpha,-\beta])$} $U_{\mathcal{P}'}([-\alpha,\beta])$ (resp. $U_{U_q^{res}(H)}^{res}([\alpha,-\beta])$, \index[not]{U@$U_{U_q^{res}(H)}^{res}([\alpha,-\beta])$} $U_{U_q^{res}(H)}^{res}([-\alpha,\beta])$) the $\mathcal{P}'$ (resp. $U_q^{res}(H)$)--subalgebras of $U_h({\frak g})$ generated by the $X_\gamma$ (resp. $(X_\gamma)^{(r)}$), where $\gamma \in [\alpha,-\beta]$ or $\gamma \in [-\alpha,\beta]$, respectively, and $X_\gamma=X_\gamma^+$ if $\gamma \in \Delta_+$, $X_\gamma = X_\gamma^-$ if  $\gamma \in \Delta_-$.
Then the elements $(X_{\beta_p})^{r_p}\ldots (X_{-\beta_q})^{r_q}$ (resp. $(X_{\beta_p})^{(r_p)}\ldots (X_{-\beta_q})^{(r_q)}$), ${ r_i}\in {\Bbb N}$
form bases of $U_{\mathcal{P}'}([\alpha,-\beta])$ (resp. $U_{U_q^{res}(H)}^{res}([\alpha,-\beta])$), and the elements $(X_{-\beta_p})^{r_p}\ldots (X_{\beta_q})^{r_q}$ (resp. $(X_{-\beta_p})^{(r_p)}\ldots (X_{\beta_q})^{(r_q)}$), ${ r_i}\in {\Bbb N}$
form bases of $U_{\mathcal{P}'}([-\alpha,\beta])$ (resp. $U_{U_q^{res}(H)}^{res}([-\alpha,\beta])$).

(vii) Let 
$\overline{U}_{\mathcal{P}'}([\alpha,-\beta])$, \index[not]{U@$\overline{U}_{\mathcal{P}'}([\alpha,-\beta])$} $\overline{U}_{\mathcal{P}'}([-\alpha,\beta])$ (resp. $\overline{U}_{U_q^{res}(H)}^{res}([\alpha,-\beta])$, \index[not]{U@$\overline{U}_{U_q^{res}(H)}^{res}([\alpha,-\beta])$} $\overline{U}_{U_q^{res}(H)}^{res}([-\alpha,\beta])$) be the $\mathcal{P}'$ (resp. $U_q^{res}(H)$)--subalgebras of $U_h({\frak g})$ generated by the $\overline{X}_\gamma$ (resp. by the $(\overline{X}_\gamma)^{(r)}$), where $\gamma \in [\alpha,-\beta]$ or $\gamma \in [-\alpha,\beta]$, respectively, and $\overline{X}_\gamma=\overline{X}_\gamma^+$ if $\gamma \in \Delta_+$, $\overline{X}_\gamma = \overline{X}_\gamma^-$ if  $\gamma \in \Delta_-$.
Then the elements $(\overline{X}_{\beta_p})^{r_p}\ldots (\overline{X}_{-\beta_q})^{r_q}$ (resp. $(\overline{X}_{\beta_p})^{(r_p)}\ldots (\overline{X}_{-\beta_q})^{(r_q)}$), ${r_i}\in {\Bbb N}$
form bases of $\overline{U}_{\mathcal{P}'}([\alpha,-\beta])$ (resp. $\overline{U}_{U_q^{res}(H)}^{res}([\alpha,-\beta])$), and the elements $(\overline{X}_{-\beta_p})^{r_p}\ldots (\overline{X}_{\beta_q})^{r_q}$ (resp. $(\overline{X}_{-\beta_p})^{(r_p)}\ldots (\overline{X}_{\beta_q})^{(r_q)}$), ${r_i}\in {\Bbb N}$
form bases of $\overline{U}_{\mathcal{P}'}([-\alpha,\beta])$ (resp. $\overline{U}_{U_q^{res}(H)}^{res}([-\alpha,\beta])$).
\end{lemma}

\begin{proof}
The first four statements of this lemma are just Propositions 8.1.7, 9.1.3 and 9.3.3 in \cite{ChP}, and Proposition 40.2.1, Corollary 40.2.2, and Proposition 41.1.4 in \cite{Lus} from which statement (v) also follows.
The proofs of the other claims are similar to each other. Consider, for instance, the case of the algebra $U_{U_q^{res}(H)}^{res}([\alpha,-\beta])$. 

Let $[\alpha,-\beta]$=$\{\beta_p,\ldots , -\beta_q\}$, $\alpha, \beta \in \Delta_+$ be a minimal segment in $\Delta$, $U_{U_q^{res}(H)}^{res}([\alpha,-\beta])$ the $U_q^{res}(H)$--subalgebra of $U_h({\frak g})$ generated by the $(X_\gamma)^{(r)}$, where $\gamma \in [\alpha,-\beta]$, and $X_\gamma=X_\gamma^+$ if $\gamma \in \Delta_+$, $X_\gamma = X_\gamma^-$ if  $\gamma \in \Delta_-$.
We show that the elements $(X_{\beta_p})^{(r_p)}\ldots (X_{-\beta_q})^{(r_q)}$, ${r_i}\in {\Bbb N}$
form a basis of $U_{U_q^{res}(H)}^{res}([\alpha,-\beta])$.

Consider the algebra $U_q^{res}(\g)\otimes_{\mathbb{C}[q,q^{-1}]}\mathbb{C}(q)$. If $x\in U_{U_q^{res}(H)}^{res}([\alpha,-\beta])\subset U_q(\g)$ then using commutation relations (\ref{qcom}) and (\ref{qcom2}) one can represent $x$ as a $U_q^{res}(H)\otimes_{\mathbb{C}[q,q^{-1}]}\mathbb{C}(q)$--linear combination of the elements $(X_{\beta_p})^{(r_p)}\ldots (X_{\beta_q})^{(r_q)}$, ${r_i}\in {\Bbb N}$. We can also consider $x$ as an element of $U_q^{res}(\g)\otimes_{\mathbb{C}[q,q^{-1}]}\mathbb{C}(q)$ and by the Poincar\'{e}--Birkhoff--Witt theorem \index{Poincar\'{e}--Birkhoff--Witt theorem} for $U_q^{res}(\g)\otimes_{\mathbb{C}[q,q^{-1}]}\mathbb{C}(q)$ (see Proposition 9.1.3 in \cite{ChP}) the above mentioned presentation of $x$ is unique. Now by the Poincar\'{e}--Birkhoff--Witt theorem for $U_q^{res}(\g)$ (see Proposition 9.3.3 in \cite{ChP} or parts (iii) and (iv) of this lemma) the coefficients in this presentation must belong to $U_q^{res}(H)$. This completes the proof in the considered case. 

\end{proof}


A basis for $U_q^{res}(H)$ is a little bit more difficult to describe. We do not need its explicit description.

\begin{remark}\label{segmPBWrev}
The anti-involutions $\omega$, $\omega_0'$ and $\omega_0$ (resp. the involution $\tau$) give rise to ant-involutions (resp. to an involution) of $U_q({\frak g})$ and $U_q^{res}(\g)$,  which we denote by the same letters.
Applying the anti-involution $\omega_0'$ to the elements of the bases constructed in Lemma \ref{segmPBW} and using (\ref{Xo0}), (\ref{X-o0}) we obtain other bases of similar types where the order of the quantum root vectors in the products defining the elements of the bases is reversed.
\end{remark}

For any minimal segment $[\alpha,\beta]\subset \Delta_+$, let $U_{U_q^{res}(H)}^{res}([\pm\alpha,\pm\beta])$ \index[not]{U@$U_{U_q^{res}(H)}^{res}([\pm\alpha,\pm\beta])$} (resp. $\overline{U}_{U_q^{res}(H)}^{res}([\pm\alpha,\pm\beta])$) \index[not]{U@$\overline{U}_{U_q^{res}(H)}^{res}([\pm\alpha,\pm\beta])$} be the subalgebra in $U_q^{res}(\g)$ generated by $U_q^{res}([\pm\alpha,\pm\beta])$ (resp. by $\overline{U}_q^{res}([\pm\alpha,\pm\beta])$) and by $U_q^{res}(H)$. 
Note that by this definition
\begin{equation}\label{bpmres}
U_q^{res}({\frak b}_\pm)=U_{U_q^{res}(H)}^{res}([\pm\beta_1,\pm\beta_D])=\overline{U}_{U_q^{res}(H)}^{res}([\pm\beta_1,\pm\beta_D]),
\end{equation}
and that 
\begin{equation}\label{npmres}
U_q^{res}(\n_\pm)=U_q^{res}([\pm\beta_1,\pm\beta_D])=\overline{U}_q^{res}([\pm\beta_1,\pm\beta_D]).
\end{equation}

Using Remark \ref{segmPBWrev} we obtain from parts (iii)-(vii) of Lemma \ref{segmPBW} the following corollary which is a quantum group counterpart of the properties of algebraic groups stated in Lemma \ref{comm} below.
\begin{corollary}\label{segmq}
For any two subalgebras $A,B\subset C$ of an algebra $C$ denote by $AB\subset C$ the image of the map $A\otimes B\to C$ induced by the multiplication in $C$. 

Let $[\alpha,\beta]\subset \Delta$ be any minimal segment,
such that $[\alpha, \beta]=[\alpha,\gamma]\cup [\delta,\beta]$ (disjoint union of minimal segments). Then the following statements are true.

(i) 
$$
U_{U_q^{res}(H)}^{res}([\alpha,\beta])=U_{U_q^{res}(H)}^{res}([\alpha,\gamma])U_{U_q^{res}(H)}^{res}([\delta,\beta])=U_{U_q^{res}(H)}^{res}([\delta,\beta])U_{U_q^{res}(H)}^{res}([\alpha,\gamma]),
$$
and
$$
\overline{U}_{U_q^{res}(H)}^{res}([\alpha,\beta])=\overline{U}_{U_q^{res}(H)}^{res}([\alpha,\gamma])\overline{U}_{U_q^{res}(H)}^{res}([\delta,\beta])=\overline{U}_{U_q^{res}(H)}^{res}([\delta,\beta])\overline{U}_{U_q^{res}(H)}^{res}([\alpha,\gamma]).
$$

(ii)
If $[\alpha,\beta]\subset \Delta_+$ or $[\alpha,\beta]\subset \Delta_-$, the multiplication in $U_q^{res}(\g)$ defines isomorphisms of $\mathbb{C}[q,q^{-1}]$--modules
$$
U_q^{res}([\alpha,\beta])\otimes U_q^{res}(H)\to U_q^{res}([\alpha,\beta])U_q^{res}(H)=U_{U_q^{res}(H)}^{res}([\alpha,\beta]),
$$
$$
U_q^{res}(H)\otimes U_q^{res}([\alpha,\beta])\to U_q^{res}(H)U_q^{res}([\alpha,\beta])=U_{U_q^{res}(H)}^{res}([\alpha,\beta]),
$$
$$
\overline{U}_q^{res}([\alpha,\beta])\otimes U_q^{res}(H)\to \overline{U}_q^{res}([\alpha,\beta])U_q^{res}(H)=\overline{U}_{U_q^{res}(H)}^{res}([\alpha,\beta]),
$$
$$
U_q^{res}(H)\otimes \overline{U}_q^{res}([\alpha,\beta])\to U_q^{res}(H)\overline{U}_q^{res}([\alpha,\beta])=U_{U_q^{res}(H)}^{res}([\alpha,\beta]).
$$

(iii)
If $[\alpha,\beta]\subset \Delta_+$ or $[\alpha,\beta]\subset \Delta_-$, the multiplication in $U_q^{res}(\g)$ defines isomorphisms of $\mathbb{C}[q,q^{-1}]$--modules
$$
U_q^{res}([\alpha,\gamma])\otimes U_q^{res}([\delta,\beta])\to U_q^{res}([\alpha,\gamma])U_q^{res}([\delta,\beta])=U_q^{res}([\alpha,\beta]),
$$
$$
U_q^{res}([\delta,\beta])\otimes U_q^{res}([\alpha,\gamma])\to U_q^{res}([\delta,\beta])U_q^{res}([\alpha,\gamma])=U_q^{res}([\alpha,\beta]),
$$
and
$$
\overline{U}_q^{res}([\alpha,\gamma])\otimes \overline{U}_q^{res}([\delta,\beta])\to \overline{U}_q^{res}([\alpha,\gamma])\overline{U}_q^{res}([\delta,\beta])=\overline{U}_q^{res}([\alpha,\beta]),
$$
$$
\overline{U}_q^{res}([\delta,\beta])\otimes \overline{U}_q^{res}([\alpha,\gamma])\to \overline{U}_q^{res}([\delta,\beta])\overline{U}_q^{res}([\alpha,\gamma])=\overline{U}_q^{res}([\alpha,\beta]).
$$
\end{corollary}

\begin{remark}
By part (a) of the Proposition in Section 2.2 of \cite{DKP3}, for any minimal segment $[\alpha,\beta]\subset \Delta$ the algebras $U_{U_q^{res}(H)}^{res}([\alpha,\beta])$, $\overline{U}_{U_q^{res}(H)}^{res}([\alpha,\beta])$, and other algebras of similar type which are defined in parts (iv)-(vii) of Lemma \ref{segmPBW}, only depend on the Weyl group element $s_{\alpha}\ldots s_\beta$, where the product of reflections is taken over all roots contained in the segment $[\alpha,\beta]$ in the order on $[\alpha,\beta]$ induced by the circular normal ordering of $\Delta$. We shall not need this result in this book. 
\end{remark}

We shall also need the following simple lemma.
\begin{lemma}\label{wlong}
Let $\overline{w}\in W$ be the longest element of the Weyl group. Then $T_{\overline{w}}U_q^{res}(\n_\pm)\subset U_q^{res}(\b_\mp)$.
\end{lemma}

\begin{proof}
Let $\overline{w}=s_{i_1}\ldots s_{i_D}$ be a reduced decomposition of the longest element $\overline{w}$ of the Weyl group $W$ of $\g$,
$$
\beta_1=\alpha_{i_1},\beta_2=s_{i_1}\alpha_{i_2},\ldots,\beta_D=s_{i_1}\ldots s_{i_{D-1}}\alpha_{i_D}
$$
the corresponding normal ordering in $\Delta_+$, and
$$
\overline{X}_{\beta_k}^\pm=T_{i_1}^{-1}\ldots T_{i_{k-1}}^{-1}X_{i_k}^\pm,~k=1,\ldots,D
$$
the quantum root vectors. 

Consider the inverse reduced decomposition of $\overline{w}$, $\overline{w}=s_{i_D}\ldots s_{i_1}$. Since $T_{\overline{w}}$ only depends on $\overline{w}$, we have for $k=1,\ldots,D$, $r\in \mathbb{N}$, using commutation relations (\ref{Ccomm1}),
$$
T_{\overline{w}}(\overline{X}_{\beta_k}^+)^{(r)}=T_{i_D}\ldots T_{i_1}T_{i_1}^{-1}\ldots T_{i_{k-1}}^{-1}(X_{i_k}^+)^{(r)}=T_{i_D}\ldots T_{i_k}(X_{i_k}^+)^{(r)}=
$$
$$
=\frac{1}{[r]_{q_{i_k}}!}(-T_{i_D}\ldots T_{i_{k-1}}(X_{i_k}^-K_{i_k}))^r=\frac{(-1)^r}{[r]_{q_{i_k}}!}T_{i_D}\ldots T_{i_{k-1}}((X_{i_k}^-)^rK_{i_k}^r)=
$$
$$
=(-1)^rq_{i_k}^{-r(r-1)}T_{i_D}\ldots T_{i_{k-1}}(X_{i_k}^-)^{(r)}T_{i_D}\ldots T_{i_{k-1}}K_{i_k}^r\in U_q^{res}(\b_-),
$$
where to justify the last inclusion we used the fact that $T_{i_D}\ldots T_{i_{k-1}}(X_{i_k}^-)^{(r)}\in U_q^{res}(\n_-)$ by part (i) of Lemma \ref{segmPBW}. 

Since the elements $(\overline{X}_{\beta_k}^+)^{(r)}$, $k=1,\ldots,D$, $r\in \mathbb{N}$ generate the algebra $U_q^{res}(\n_+)$ by Lemma \ref{segmPBW} (v) and by (\ref{npmres}), we deduce that $T_{\overline{w}}U_q^{res}(\n_+)\subset U_q^{res}(\b_-)$. The other inclusion is established in a similar way. This completes the proof.

\end{proof}

In conclusion we note that the algebras $U_q^{res}(\g)\otimes_{\mathbb{C}[q,q^{-1}]}\mathbb{C}(q)$ and $U_q^{res}(\g)$ are graded by the elements of the root lattice $Q$. Indeed, one can assign the weight \index{weight!of a quantum group element} $-r_1\beta_1-\ldots-r_D\beta_D+m_1\beta_1+\ldots+m_D\beta_D$ to each element
$$
(X^-_{\beta_1})^{(r_1)}\ldots (X^-_{\beta_D})^{(r_D)}X(X^+_{\beta_1})^{(m_1)}\ldots (X^+_{\beta_D})^{(m_D)},
$$
$$
{\rm wt}((X^-_{\beta_1})^{(r_1)}\ldots (X^-_{\beta_D})^{(r_D)}X(X^+_{\beta_1})^{(m_1)}\ldots (X^+_{\beta_D})^{(m_D)})=-r_1\beta_1-\ldots-r_D\beta_D+m_1\beta_1+\ldots+m_D\beta_D, \index[not]{w@${\rm wt}$}
$$
where $X\in U_q^{res}(H)$, use parts (i) and (ii) of Lemma \ref{segmPBW} to extend the definition of the weight to $U_q^{res}(\g)$ and observe that the relations in $U_q^{res}(\g)$ following from (\ref{defQGrel}) are homogeneous with respect to this grading. Then one can naturally extend this grading to $U_q^{res}(\g)\otimes_{\mathbb{C}[q,q^{-1}]}\mathbb{C}(q)$. 

For $\mu \in Q$ we introduce the $\mathbb{C}[q,q^{-1}]$--submodule $(U_q^{res}(\g))_\mu\subset U_q^{res}(\g)$, \index[not]{U@$(U_q^{res}(\g))_\mu$} $(U_q^{res}(\g))_\mu=\{x\in U_q^{res}(\g): {\rm wt}(x)=\mu \}$, and call it {\it the weight subspace} of weight $\mu$ of $U_q^{res}(\g)$. \index{weight!subspace of a quantum group} Alternatively, one has 
\begin{equation}\label{Qgr*}
(U_q^{res}(\g))_\mu=\{x\in U_q^{res}(\g): [h,x]=\mu(h)x \text{ for all } h\in \h \subset U_h(\g) \},
\end{equation} 
where the commutator is defined in $U_h(\g)$, and $\h$ is regarded as a complex vector subspace of $U_h(\g)$ spanned by the elements $H_i$, $i=1,\ldots,l$.

By parts (i) and (iii) of  Lemma \ref{segmPBW} we have the following direct sum of $\mathbb{C}[q,q^{-1}]$--modules
$$
U_q^{res}(\g)=\bigoplus_{\mu\in Q}(U_q^{res}(\g))_\mu.
$$ 

From the definition of the braid group action on $U_h(\g)$ and from (\ref{Qgr*}) it follows that for all $w\in W$ one has $T_w((U_q^{res}(\g))_\mu)=(U_q^{res}(\g))_{w\mu}$.

The anti-involutions $\omega_0'$, $\omega$ and the involution $\tau$ preserve weights of elements while the anti-involution $\omega$ inverses the signs of their weights. As a consequence, one has $T_w^\tau((U_q^{res}(\g))_\mu)=(U_q^{res}(\g))_{w\mu}$.

The algebra $U_h(\g)$ is, of course, not graded by the elements of the root lattice. But one can still define weights of the elements of its topological basis introduced in parts (ii) and (iii) of  Lemma \ref{segmPBW},
$$
{\rm wt}((X^-_{\beta_1})^{r_1}\ldots (X^-_{\beta_D})^{r_D}X(X^+_{\beta_1})^{m_1}\ldots (X^+_{\beta_D})^{m_D})=-r_1\beta_1-\ldots-r_D\beta_D+m_1\beta_1+\ldots+m_D\beta_D, \index[not]{w@${\rm wt}$}
$$
where $X\in U_h(\h)$.

\section{The universal R--matrix}\label{URmatr}

\pagestyle{myheadings}
\markboth{CHAPTER~\thechapter.~QUANTUM GROUPS}{\thesection.~THE UNIVERSAL R--MATRIX}

\setcounter{equation}{0}
\setcounter{theorem}{0}

$U_h({\frak g})$ is a {\it quasitriangular} Hopf algebra, \index{Hopf algebra!quasitriangular} i.e. there exists an invertible element
${\mathcal R}\in U_h({\frak g})\otimes U_h({\frak g})$ \index[not]{R@${\mathcal R}$} (completed tensor product), called {\it a quantum universal R--matrix}, \index{r--matrix!quantum universal} such that
\begin{equation}\label{quasitr}
\Delta^{opp}_h(a)={\mathcal R}\Delta_h(a){\mathcal R}^{-1}\mbox{ for all } a\in U_h({\frak g}),
\end{equation}
and
\begin{equation}\label{rmprop}
\begin{array}{l}
(\Delta_h \otimes id){\mathcal R}={\mathcal R}_{13}{\mathcal R}_{23},\\
\\
(id \otimes \Delta_h){\mathcal R}={\mathcal R}_{13}{\mathcal R}_{12},
\end{array}
\end{equation}
where ${\mathcal R}_{12}={\mathcal R}\otimes 1,~{\mathcal R}_{23}=1\otimes {\mathcal R},
~{\mathcal R}_{13}=(\sigma \otimes id){\mathcal R}_{23}$, and $\sigma (x\otimes y)=y\otimes x$, $x,y\in U_h({\frak g})$.

More generally, if $X=\sum_kx_k\otimes y_k\in A\otimes A$ is an element of the (completed) tensor product $A\otimes A$ of a (topological) Hopf algebra $A$ we shall write $X_{ij}$ for the element 
$$
X_{ij}=\sum_k 1\otimes \ldots \otimes 1\otimes x_k\otimes 1\otimes \ldots \otimes 1\otimes y_k\otimes 1\otimes \ldots \otimes 1\in \underbrace{A\otimes \ldots \otimes A}_{\text{n factors}},
$$
where in the right hand side $x_k$ and $y_k$ belong to the i-th and the j-th factors, respectively, of the n-fold (completed) tensor product $A\otimes \ldots \otimes A$.  

From (\ref{quasitr}) and (\ref{rmprop}) it follows that $\mathcal R$ satisfies {\it the quantum Yang--Baxter equation} \index{Yang-Baxter equation!quantum}
\begin{equation}\label{YB}
{\mathcal R}_{12}{\mathcal R}_{13}{\mathcal R}_{23}={\mathcal R}_{23}{\mathcal R}_{13}{\mathcal R}_{12}.
\end{equation}

For every quasitriangular Hopf algebra with a universal r--matrix $\mathcal R$ we also have 
\begin{equation}\label{SR}
(S\otimes id){\mathcal R}=(id \otimes S^{-1}){\mathcal R}={\mathcal R}^{-1},
\end{equation}
and
\begin{equation}\label{S}
(S\otimes S){\mathcal R}={\mathcal R},
\end{equation}
where $S$ \index{Hopf algebra!antipode} is the antipode.

In the case of $U_h({\frak g})$ an explicit expression for $\mathcal R\in U_h({\frak g})\otimes U_h({\frak g})$ may be written by making use of the q--exponential
$$
{\exp}_q(x)={\exp}_q'(qx)=\sum_{k=0}^\infty q^{\frac{1}{2}k(k+1)}{x^k \over [k]_q!}
$$
in terms of which the element $\mathcal R$ takes the form
\begin{equation}\label{univr}
{\mathcal R}=\prod_{\beta\in \Delta_+}
{\exp}_{q_{\beta}}[(1-q_{\beta}^{-2})X_{\beta}^-\otimes X_{\beta}^+]{\exp}\left[h\sum_{i=1}^l(Y_i\otimes H_i)\right]=\prod_{\beta\in \Delta_+}(\theta_\beta)_{21}{\exp}\left[h\sum_{i=1}^l(Y_i\otimes H_i)\right],
\end{equation}
where $(\theta_\beta)_{21}=\sigma\theta_\beta$,
the product is over all the positive roots of $\frak g$, and the order of the terms is such that
the $\alpha$--term appears to the left of the $\beta$--term if $\alpha < \beta$ with respect to the normal
ordering 
$$
\beta_1=\alpha_{i_1},\beta_2=s_{i_1}\alpha_{i_2},\ldots,\beta_D=s_{i_1}\ldots s_{i_{D-1}}\alpha_{i_D}
$$
of $\Delta_+$ which is used in the definition of the quantum root vectors $X_\beta^\pm$.

One can calculate the action of the comultiplication on the root vectors $X_{\beta_k}^\pm$ in terms of the universal R--matrix. For instance for $\Delta_h(X_{\beta_k}^-)$ one has
\begin{equation}\label{comult}
\Delta_h(X_{\beta_k}^-)=\theta_{w_{k-1}}(X_{\beta_k}^-\otimes 1+e^{h\beta^\vee}\otimes X_{\beta_k}^-)\theta_{w_{k-1}}^{-1},
\end{equation}
where for $w_{k-1}:=s_{i_1}\ldots s_{i_{k-1}}$
$$
\theta_{w_{k-1}}=\theta_{\beta_{1}}\ldots \theta_{\beta_{k-1}}. \index[not]{t@$\theta_w$}
$$

The r--matrix $r_-:=-\frac 12 h^{-1}({\mathcal R}-1\otimes 1)~~(\mbox{mod }h)$, which is the classical limit of $\mathcal R$, coincides with the classical r--matrix (\ref{rcl}).


\section{Realizations of quantum groups associated to Weyl group elements}\label{wqreal}

\pagestyle{myheadings}
\markboth{CHAPTER~\thechapter.~QUANTUM GROUPS}{\thesection.~REALIZATIONS OF QUANTUM GROUPS ASSICIATED TO WEYL GROUP ELEMENTS}

\setcounter{equation}{0}
\setcounter{theorem}{0}

q-W--algebras will be defined in terms of certain integral forms of non--standard realizations of quantum groups associated to Weyl group elements.

Let $s$ be an element of the Weyl group $W$ of the pair $(\g,\h)$, and $\h'$ the orthogonal complement, with respect to the symmetric bilinear form, to the subspace of $\h$ fixed by the natural action of $s$ on $\h$. Let $\h'^*$ be the image of $\h'$ in $\h^*$ under the identification $\h^*\simeq \h$ induced by the symmetric invariant bilinear form on $\g$.
The restriction of the natural action of $s$ on $\h^*$ to the subspace $\h'^*$ has no fixed points. Therefore one can define the Cayley transform ${1+s \over 1-s }$ of the restriction of $s$ to $\h'^*$. Denote by $P_{\h'^*}$ \index[not]{P@$P_{\h'^*}$} the orthogonal projection operator onto ${\h'^*}$ in $\h^*$, with respect to the symmetric invariant bilinear form.

Let $\kappa \in \mathbb{Z}$ \index[not]{k@$\kappa$} be an integer number and 
$U_h^{s}({\frak g})$ \index[not]{U@$U_h^s(\g)$} the topological algebra over ${\Bbb C}[[h]]$ topologically generated by elements
$e_i , f_i , H_i,~i=1, \ldots l$ \index[not]{e@$e_i$} \index[not]{f@$f_i$} \index[not]{H@$H_i$} subject to the relations: \index{commutation relations!defining for $U_h^{s}({\frak g})$}
\begin{equation}\label{sqgr}
\begin{array}{l}
[H_i,H_j]=0,~~ [H_i,e_j]=a_{ij}e_j, ~~ [H_i,f_j]=-a_{ij}f_j,~~e_i f_j -q^{ c_{ij}} f_j e_i = \delta _{i,j}{K_i -K_i^{-1} \over q_i -q_i^{-1}} ,\\
\\
c_{ij}:=\kappa\left\langle  {1+s \over 1-s }P_{\h'^*}\alpha_i , \alpha_j \right\rangle,~~K_i=e^{d_ihH_i}, \\ \index[not]{c@$c_{ij}$}
\\
\sum_{r=0}^{1-a_{ij}}(-1)^r q^{r c_{ij}}
\left[ \begin{array}{c} 1-a_{ij} \\ r \end{array} \right]_{q_i}
(e_i )^{1-a_{ij}-r}e_j (e_i)^r =0 ,~ i \neq j , \\
\\
\sum_{r=0}^{1-a_{ij}}(-1)^r q^{r c_{ij}}
\left[ \begin{array}{c} 1-a_{ij} \\ r \end{array} \right]_{q_i}
(f_i )^{1-a_{ij}-r}f_j (f_i)^r =0 ,~ i \neq j .
\end{array}
\end{equation}

\begin{proposition}\label{newreal}
For every solution $n_{ij}\in {\Bbb C},~i,j=1,\ldots ,l$ \index[not]{n@$n_{ij}$} of equations
\begin{equation}\label{eqpi}
d_jn_{ij}-d_in_{ji}=c_{ij}
\end{equation}
there exists an algebra
isomorphism $\psi_{\{ n_{ij}\}} : U_h^{s}({\frak g}) \rightarrow
U_h({\frak g})$ \index[not]{p@$\psi_{\{ n_{ij}\}}$} defined  by the formulas:
$$
\psi_{\{ n_{ij}\}}(e_i)=X_i^+ \prod_{p=1}^lL_p^{n_{ip}},~~
\psi_{\{ n_{ij}\}}(f_i)=\prod_{p=1}^lL_p^{-n_{ip}}X_i^- ,~~
\psi_{\{ n_{ij}\}}(H_i)=H_i .
$$
\end{proposition}

\begin{proof} 
The proof of this proposition is by direct verification of defining relations (\ref{sqgr}). The most nontrivial part is to verify the deformed quantum Serre relations, i.e. the last two relations in (\ref{sqgr}). For instance, the defining relations of $U_h({\frak g})$  imply the following relations for $\psi_{\{ n_{ij}\}}(e_i)$,
$$
\sum_{k=0}^{1-a_{ij}}(-1)^k
\left[ \begin{array}{c} 1-a_{ij} \\ k \end{array} \right]_{q_i}
q^{k({d_j}n_{ij}-d_in_{ji})}\psi_{\{ n_{ij}\}}(e_i)^{1-a_{ij}-k}\psi_{\{ n_{ij}\}}(e_j)\psi_{\{ n_{ij}\}}(e_i)^k =0 ,
$$
for any $i\neq j$.
Now using equation (\ref{eqpi}) we arrive at the quantum Serre relations for $e_i$ in (\ref{sqgr}).

\end{proof}

The general solution of equation (\ref{eqpi})
is given by
\begin{equation}\label{eq3}
n_{ij}=\frac 1{2d_j} (c_{ij} + {s_{ij}}),
\end{equation}
where $s_{ij}=s_{ji}$.

We shall only use the solution for which $s_{ij}=0$ for all $i,j=1,\ldots l$. Then
\begin{equation}\label{nijfix}
n_{ij}=\frac 1{2d_j} c_{ij}	\index[not]{n@$n_{ij}$} \index[not]{n@$n_{ij}$}
\end{equation}
From now on we assume that solution (\ref{nijfix}) is used to identify $U_h^{s}({\frak g})$ and $U_h({\frak g})$.

The algebra $U_h^{s}({\frak g})$ is called {\it the realization of the quantum group $U_h({\frak g})$ associated to the Weyl group element $s\in W$}. \index{quantum group!realization associated to a Weyl group element}

Denote by $U_h^{s}({\frak n}_\pm)$ \index[not]{U@$U_h^s(\n_\pm)$} the subalgebra in $U_h^{s}({\frak g})$ generated by
$e_i$ (resp. $f_i$), $i=1, \ldots, l$.
Let $U_h^{s}({\frak h})$ be the subalgebra in $U_h^{s}({\frak g})$ \index[not]{U@$U_h^s(\h)$} generated by $H_i,~i=1,\ldots ,l$.

We shall construct analogues of quantum root vectors for $U_h^{s}({\frak g})$.
It is convenient to introduce an operator $K_s\in {\rm End}_{\mathbb{C}}{\frak h}$ \index[not]{K@$K_s$} defined by
\begin{equation}\label{Kdef}
K_sH_i=\sum_{j=1}^l{n_{ij} \over d_i}Y_j.
\end{equation}
From (\ref{nijfix}) we obtain that
\begin{equation}\label{Kdefs}
K_sh=\frac{\kappa}{2}{1+s \over 1-s}P_{\h'}h, h\in \h,
\end{equation}
where $P_{\h'}$ \index[not]{P@$P_{\h'}$} is the orthogonal projection operator onto ${\h'}$ in $\h$, with respect to the symmetric invariant bilinear form.

\begin{proposition}\label{rootss}
Let $s\in W$ be an element of the Weyl group $W$ of the pair $(\g,\h)$, $\Delta$ the root system of the pair $(\g,\h)$. Let $U_h^{s}({\frak g})$ be the realization of the quantum group $U_h({\frak g})$ associated to $s$.

For any normal ordering of the root system $\Delta_+$
the elements 
$$
e_{\beta}(s)=\psi_{\{ n_{ij}\}}^{-1}(X_{\beta}^+e^{hK_s\beta^\vee}) ~\text{and}
~f_{\beta}(s)=\psi_{\{ n_{ij}\}}^{-1}(e^{-hK_s\beta^\vee}X_{\beta}^-),~\beta \in \Delta_+ \index[not]{e@$e_{\beta}(s)$} \index[not]{f@$f_{\beta}(s)$}
$$
lie in the subalgebras $U_h^{s}({\frak n}_+)$ and $U_h^{s}({\frak n}_-)$, respectively.

The elements $f_\beta(s)\in U_h^{s}({\n_-})$, $\beta \in \Delta_{\m_+}$ defined with the help of the normal ordering (\ref{NO}) generate a subalgebra $U_h^{s}({\frak m}_-)\subset U_h^{s}({\frak g})$ \index[not]{U@$U_h^s(\m_-)$} such that 
$$
U_h^{s}({\frak m}_-)/hU_h^{s}({\frak m}_-)\simeq U({\m_-}), \index[not]{U@$U(\m_-)$}
$$ 
where ${\m_-}$ is the Lie subalgebra of $\g$ generated by the root vectors $X_{-\alpha}$, $\alpha\in \Delta_{\m_+}$. \index[not]{m@$\m_-$}
\end{proposition}

\begin{proof}
Fix a normal ordering of the root system $\Delta_+$.
Let $\beta=\sum_{i=1}^l m_i\alpha_i \in \Delta_+$ be a positive root,
$X_{\beta}^+\in U_h({\frak g})$ the corresponding quantum root vector constructed with the help of the fixed normal ordering of $\Delta_+$. Then $\beta^\vee=\sum_{i=1}^l m_id_iH_i$, and so
$K_s\beta^\vee=\sum_{i,j=1}^l m_in_{ij}Y_j$. Now the proof of the first statement follows immediately from
Proposition \ref{rootprop}, commutation relations (\ref{weight-root}) and the definition of the isomorphism $\psi_{\{ n_{ij}\}}$. 

The second assertion is a consequence of (\ref{qrootclass}).

\end{proof}

\begin{remark}
To simplify the notation we shall often write $e_{\beta}(s)=e_{\beta}$ and $f_{\beta}(s)=f_{\beta}$ if it does not cause any confusion. \index[not]{e@$e_{\beta}$} \index[not]{f@$f_{\beta}$}
\end{remark}

The realizations $U_h^{s}({\frak g})$ of the quantum group $U_h({\frak g})$
are related to quantizations of some nonstandard bialgebra structures on $\frak g$. At the quantum level changing bialgebra structure corresponds to the so--called Drinfeld twist. \index{Hopf algebra!Drinfeld twist}
The relevant class of such twists is described in the following proposition which is a combination of Propositions 4.2.13, 16.1.5, of formula (15) in \S 16.1 in \cite{ChP}, and of the results of \S 1 in \cite{Dri}. 
\begin{proposition}\label{twdef}
Let $(A,\mu , \imath , \Delta , \varepsilon , S)$ \index[not]{ZZZZ@$(A,\mu , \imath , \Delta , \varepsilon , S)$} be a (topological) Hopf algebra \index{Hopf algebra}  over a commutative ring with multiplication $\mu$, \index{Hopf algebra!multiplication} \index[not]{m@$\mu(~\cdot~,~\cdot~)$} unit $\imath$, \index{Hopf algebra!unit} \index[not]{i@$\imath$} comultiplication $\Delta$, \index[not]{D@$\Delta(~\cdot~)$} \index{Hopf algebra!comultiplication} counit $\varepsilon$ \index{Hopf algebra!counit} \index[not]{e@$\varepsilon(~\cdot~)$} and antipode $S$. \index{Hopf algebra!antipode} \index[not]{S@$S(~\cdot~)$}

(i) Let $\mathcal F$ \index[not]{F@$\mathcal F$} be an invertible element of the (completed) tensor product $A\otimes A$
such that
\begin{equation}\label{twist1}
(\varepsilon \otimes id)({\mathcal F})=(id \otimes \varepsilon )({\mathcal F})=1,
\end{equation}
and
\begin{equation}\label{twist2}
{\mathcal F}_{12}(\Delta \otimes id)({\mathcal F})={\mathcal F}_{23}(id \otimes \Delta)({\mathcal F}).
\end{equation}
Then 
\begin{equation}\label{defv}
v^{\mathcal F}:=\mu (id\otimes S)({\mathcal F}) \index[not]{v@$v^{\mathcal F}$}
\end{equation}
is an invertible element of $A$ with
$$
(v^{\mathcal F})^{-1}=\mu (S\otimes id)({\mathcal F}^{-1}).
$$

Moreover, if we define $\Delta^{\mathcal F}:A\rightarrow A\otimes A$ \index[not]{D@$\Delta^{\mathcal F}(~\cdot~)$} and $S^{\mathcal F}:A\rightarrow A$ \index[not]{S@$S^{\mathcal F}(~\cdot~)$} by
\begin{equation}\label{deltaF}
\Delta^{\mathcal F}(a)={\mathcal F}\Delta(a){\mathcal F}^{-1},~~S^{\mathcal F}(a)=v^{\mathcal F}S(a)(v^{\mathcal F})^{-1},
\end{equation}
then $(A,\mu , \imath , \Delta^{\mathcal F} , \varepsilon , S^{\mathcal F})$ \index[not]{ZZZZ@$(A,\mu , \imath , \Delta^{\mathcal F} , \varepsilon , S^{\mathcal F})$} is a (topological) Hopf algebra denoted by $A^{\mathcal F}$ \index[not]{A@$A^{\mathcal F}$}
and called the Drinfeld twist of $A$ by ${\mathcal F}$. \index{Hopf algebra!Drinfeld twist}

(ii) Suppose that $A$ and ${\mathcal F}$ are as in part (i), and assume in addition that $A$ is quasitriangular with universal R--matrix $\mathcal R$. Then $A^{\mathcal F}$ is quasitriangular with universal R--matrix
\begin{equation}\label{rf}
{\mathcal R}^{\mathcal F}={\mathcal F}_{21}{\mathcal R}{\mathcal F}^{-1}, \index[not]{R@${\mathcal R}^{\mathcal F}$}
\end{equation}
where ${\mathcal F}_{21}=\sigma {\mathcal F}$.

(iii) If ${\mathcal F}, {\mathcal G}\in A\otimes A$ are invertible elements, ${\mathcal G}$ satisfies (\ref{twist1}) and (\ref{twist2}), ${\mathcal F}$ satisfies (\ref{twist1}) and (\ref{twist2}) with $\Delta$ replaced by $\Delta^{\mathcal G}$ then ${\mathcal F}{\mathcal G}$ satisfies (\ref{twist1}) and (\ref{twist2}) and the twist of $A$ by ${\mathcal F}{\mathcal G}$ is the composition of the twists of $A$ by  ${\mathcal F}$ and by ${\mathcal G}$.
\end{proposition}

Let
\begin{equation}\label{Ftw}
\mathcal{F}_s={\exp}(-h\sum_{i,j=1}^l {n_{ij} \over d_i}Y_i\otimes Y_j)={\exp}(-h\sum_{i=1}^l Y_i\otimes K_sH_i) \in U_h({\frak h})\otimes U_h({\frak h}). \index[not]{F@$\mathcal{F}_s$}
\end{equation}
Then from the definition of the commutative Hopf subalgebra $U({\frak h})[[h]]\simeq U_h({\frak h})\subset U_h({\frak g})$ it immediately follows that $\mathcal{F}_s\in U_h({\frak g})\otimes U_h({\frak g})$ satisfies (\ref{twist1}) and (\ref{twist2}) (see also \cite{ChP}, Proposition 6.5.8), and according to Proposition \ref{twdef} (i)
\begin{equation}\label{defds}
\Delta_s(a)=(\psi_{\{ n_{ij}\}}^{-1}\otimes \psi_{\{ n_{ij}\}}^{-1})\mathcal{F}_s\Delta_h(\psi_{\{ n_{ij}\}}(a))\mathcal{F}_s^{-1}. \index[not]{D@$\Delta_s(~\cdot~)$}
\end{equation}
defines a comultiplication on $U_h^{s}({\frak g})$.

From the formulas in part (i) of Proposition \ref{twdef} it follows that on the generators the comultiplication $\Delta_{s}$ is explicitly given by the following formulas
$$
\begin{array}{c}
\Delta_{s}(H_i)=H_i\otimes 1+1\otimes H_i,\\
\\
\Delta_{s}(e_i)=e_i\otimes e^{-hd_iH_i}+e^{h\kappa d_i {1+s \over 1-s}P_{\h'}H_i}\otimes e_i,~~
\Delta_{s}(f_i)=f_i\otimes 1+e^{-h\kappa d_i {1+s \over 1-s}P_{\h'}H_i+hd_iH_i}\otimes f_i,
\end{array}
$$
the corresponding antipode $S_s$ \index[not]{S@$S_s(~\cdot~)$} is given by
\begin{equation}\label{S_s}
S_s(e_i)=-e^{-h\kappa d_i {1+s \over 1-s}P_{\h'}H_i}e_ie^{hd_iH_i},~S_s(f_i)=-e^{h\kappa d_i {1+s \over 1-s}P_{\h'}H_i-hd_iH_i}f_i,~S_s(H_i)=-H_i,
\end{equation}
and the corresponding counit $\varepsilon_{s}$ \index[not]{e@$\varepsilon_s(~\cdot~)$} is given by
$$
\varepsilon_{s}(H_i)=\varepsilon_s(e_i)=\varepsilon_s(f_i)=0.
$$
We shall always assume that the algebra $U_h^{s}({\frak g})$ is equipped with this Hopf algebra structure.

Note that the Hopf algebra $U_h^{s}({\frak g})$ is a quantization of the bialgebra structure on $\frak g$ defined by the cocycle $\delta^s$ \index[not]{d@$\delta^s$} associated to the Weyl group element $s$, \index{cocycle!associated to a Weyl group element} 
\begin{equation}\label{cocycles}
\delta^s (x)=({\rm ad}x\otimes 1+1\otimes {\rm ad}x)2r^{s}_\pm,~~ r^{s}_\pm\in {\frak g}\otimes {\frak g},
\end{equation}
where $r^{s}_\pm:=r_\pm + \frac 12 \sum_{i=1}^l \kappa {1+s \over 1-s }P_{\h'}H_i\otimes Y_i$, \index[not]{r@$r^{s}_\pm$} and $r_\pm$ are given by (\ref{rcl}). We call $r^s:=r^{s}_++r^{s}_-$ \index[not]{r@$r^{s}$} {\it the classical r--matrix associated to the Weyl group element $s$}. \index{r--matrix!classical!associated to a Weyl group element} 

By Proposition \ref{twdef} (ii) $U_h^{s}({\frak g})$ is a quasitriangular topological Hopf algebra with the universal R--matrix
$\mathcal{R}_s=(\psi_{\{ n_{ij}\}}^{-1}\otimes \psi_{\{ n_{ij}\}}^{-1})((\mathcal{F}_s)_{21}{\mathcal R}{\mathcal{F}_s}^{-1})$, \index[not]{R@$\mathcal{R}_s$}
\begin{equation}\label{rmatrspi}
\begin{array}{l}
\mathcal{R}_s=
\prod_{\beta\in \Delta_+}
{\exp}_{q_{\beta}}[(1-q_{\beta}^{-2})f_{\beta} \otimes
e_{\beta}e^{-h\kappa {1+s \over 1-s}P_{\h'} \beta^\vee}]\times \\
\times {\exp}\left[ h(\sum_{i=1}^l(Y_i\otimes H_i)-
\sum_{i=1}^l \kappa {1+s \over 1-s }P_{\h'}H_i\otimes Y_i) \right]=
\end{array}
\end{equation}
$$
\begin{array}{l}
={\exp}\left[ h(\sum_{i=1}^l(Y_i\otimes H_i)-
\sum_{i=1}^l \kappa {1+s \over 1-s }P_{\h'}H_i\otimes Y_i) \right]\times \\
\times
\prod_{\beta\in \Delta_+}
{\exp}_{q_{\beta}}[(1-q_{\beta}^{-2})e^{h\left(\kappa {1+s \over 1-s}P_{\h'} -id\right)\beta^\vee}f_{\beta} \otimes
e_{\beta}e^{h\beta^\vee}],
\end{array}
$$
where the order of the terms in the product over the positive roots is such that the $\alpha$--term appears to the left of the $\beta$--term if $\alpha < \beta$ in the normal ordering of $\Delta_+$ with the help of which the quantum root vectors $e_{\beta}, f_{\beta}$ are defined in Proposition \ref{rootss}.

Similarly to (\ref{dhtw}) one obtains that for a reduced decomposition $w=s_{i_1}\ldots s_{i_k}$ and $T_w=T_{i_1}\ldots T_{i_k}$ only depending on $w$ one has from (\ref{dhtw}) and (\ref{defds}) 
\begin{align}\label{DTs}
\Delta_s(T_w)=\prod^{k}_{p=1}\theta_{\beta_p}^s\mathcal{F}_s(T_w \otimes T_w)(\mathcal{F}_s^{-1})T_w\otimes T_w=
\\
=T_w\otimes T_w(T_{w^{-1}}\otimes T_{w^{-1}})(\mathcal{F}_s)\mathcal{F}_s^{-1}\prod_{p=1}^k\overline{\theta}_{\beta_p'}^s, \nonumber
\end{align}
where in the products $\theta_{\beta_p}^s$ \index[not]{t@$\theta_{\beta_p}^s$} (resp. $\overline{\theta}_{\beta_p'}^s$) \index[not]{t@$\overline{\theta}_{\beta_p'}^s$} appears on the left from $\theta_{\beta_q}^s$ (resp. $\overline{\theta}_{\beta_q'}^s$) if $p<q$, and for $p=1,\ldots, k$
$$
e_{\beta_p}=\psi_{\{ n_{ij}\}}^{-1}(X_{\beta_p}^+e^{hK_s\beta_p^\vee}),
f_{\beta}=\psi_{\{ n_{ij}\}}^{-1}(e^{-hK_s\beta_p^\vee}X_{\beta_p}^-), \beta_p=s_{i_1}\ldots s_{i_{p-1}}\alpha_{i_p},
$$
$$
X_{\beta_p}^\pm=T_{i_1}\ldots T_{i_{p-1}}X_{i_p}^\pm,~{\overline{X}_{\beta_p'}^\pm}=T_{i_k}^{-1}\ldots T_{i_{p+1}}^{-1}X_{i_p}^\pm,
$$
$$
\overline{e}_{\beta_p'}=\psi_{\{ n_{ij}\}}^{-1}(\overline{X}_{\beta_p'}^+e^{hK_s{\beta_p'}^\vee}),
\overline{f}_{\beta_p'}=\psi_{\{ n_{ij}\}}^{-1}(e^{-hK_s{\beta_p'}^\vee}\overline{X}_{\beta_p'}^-),
~ K_{\beta_p'}=T_{i_k}^{-1}\ldots T_{i_{p+1}}^{-1}K_{i_p}, \beta_p'=s_{i_k}\ldots s_{i_{p+1}}\alpha_{i_p}
$$
$$
\theta_{\beta_p}^s=\mathcal{F}_s\theta_{\beta_p}\mathcal{F}_s^{-1}={\exp}_{q_{\beta_p}}[(1-q_{\beta_p}^{-2})e_{\beta_p}e^{-h\kappa {1+s \over 1-s}P_{\h'} \beta_p^\vee}\otimes f_{\beta_p}],
$$
$$\overline{\theta}_{\beta_p'}^s=\mathcal{F}_s\overline{\theta}_{\beta_p'}\mathcal{F}_s^{-1}={\exp}_{q_{\beta_p'}}[(1-q_{\beta_p'}^{-2})K_{\beta_p'}^{-1}e^{h\kappa {1+s \over 1-s}P_{\h'} {\beta_p'}^\vee}{\overline{f}_{\beta_p'}}\otimes {\overline{e}_{\beta_p'}}K_{\beta_p'}]. 
$$

In the same way, for $\overline{T}_w=T_{i_1}^{-1}\ldots T_{i_k}^{-1}$ only depending on $w$ one has from (\ref{dhotw}) and (\ref{defds})
\begin{align}\label{DTw}
\Delta_s(\overline{T}_w)=\prod^{k}_{p=1}{\overline{\theta}_{\beta_p}^s}'\mathcal{F}_s(T_w\otimes T_w)(\mathcal{F}_s^{-1})\overline{T}_w\otimes \overline{T}_w=
\\
=\overline{T}_w\otimes \overline{T}_w(T_{w^{-1}}\otimes T_{w^{-1}})(\mathcal{F}_s)\mathcal{F}_s^{-1}\prod_{p=1}^k{\theta_{\beta_p'}^s}',\nonumber
\end{align}
where in the products ${\overline{\theta}_{\beta_p}^s}'$ \index[not]{t@${\overline{\theta}_{\beta_p}^s}'$} (resp. ${\theta_{\beta_p'}^s}'$) \index[not]{t@${\theta_{\beta_p'}^s}'$} appears on the left from ${\overline{\theta}_{\beta_q}^s}'$ (resp. ${\theta_{\beta_q'}^s}'$) if $p<q$, and for $p=1,\ldots, k$
\begin{align*}
\overline{e}_{\beta_p}=\psi_{\{ n_{ij}\}}^{-1}(\overline{X}_{\beta_p}^+e^{hK\beta_p^\vee}),
\overline{f}_{\beta_p}=\psi_{\{ n_{ij}\}}^{-1}(e^{-hK\beta_p^\vee}\overline{X}_{\beta_p}^-),
\overline{K}_{\beta_p}=T_{i_1}^{-1}\ldots T_{i_{p-1}}^{-1}K_{i_p},
\\
{\overline{X}_{\beta_p}^\pm}=T_{i_1}^{-1}\ldots T_{i_{p-1}}^{-1}X_{i_p}^\pm,
\beta_p=s_{i_1}\ldots s_{i_{p-1}}\alpha_{i_p},
\end{align*}
$$
e_{\beta_p'}'=\psi_{\{ n_{ij}\}}^{-1}({X_{\beta_p'}^+}'e^{hK{\beta_p'}^\vee}),
f_{\beta_p'}'=\psi_{\{ n_{ij}\}}^{-1}(e^{-hK{\beta_p'}^\vee}{X_{\beta_p'}^-}'),
\beta_p'=s_{i_k}\ldots s_{i_{p+1}}\alpha_{i_p}, {X_{\beta_p'}^\pm}'=T_{i_k}\ldots T_{i_{p+1}}X_{i_p}^\pm,
$$
\begin{equation}\label{ths}
{\theta_{\beta_p'}^s}'=\mathcal{F}_s\theta_{\beta_p'}'\mathcal{F}_s^{-1}={\exp}_{q_{\beta_p'}^{-1}}[(1-q_{\beta_p'}^{2})e_{\beta_p'}'e^{-h\kappa {1+s \over 1-s}P_{\h'} {\beta_p'}^\vee}\otimes f_{\beta_p'}'], 
\end{equation}
$$
{\overline{\theta}_{\beta_p}^s}'=\mathcal{F}_s\overline{\theta}_{\beta_p}'\mathcal{F}_s^{-1}={\exp}_{q_{\beta_p}^{-1}}[(1-q_{\beta_p}^{2})\overline{K}_{\beta_p}^{-1}e^{h\kappa {1+s \over 1-s}P_{\h'} \beta_p^\vee}\overline{f}_{\beta_p}\otimes \overline{e}_{\beta_p}\overline{K}_{\beta_p}]. 
$$


\section{The adjoint action}\label{Adj}

\pagestyle{myheadings}
\markboth{CHAPTER \thechapter.~QUANTUM GROUPS}{\thesection.~THE ADJOINT ACTION}

\setcounter{equation}{0}
\setcounter{theorem}{0}

Next we discuss the properties of the adjoint action of a Hopf algebra on itself with respect to Drinfeld twists.
Define the {\it right adjoint action} of a Hopf algebra $(A,\mu , \imath , \Delta , \varepsilon , S)$ \index{action!adjoint!right, of a Hopf algebra} on itself by the formula
\begin{equation}\label{ad}
{\rm Ad}x(z)=S(x^1)zx^2, \index[not]{A@${\rm Ad}$}
\end{equation}
and the {\it left adjoint action} by \index{action!adjoint!left, of a Hopf algebra}
\begin{equation}\label{ad'}
{\rm Ad}'x(z)=x^1zS(x^2),~x, z\in A, \index[not]{A@${\rm Ad}'$}
\end{equation}
where here and below we use the abbreviated Sweedler notation \index{Sweedler notation} for the comultiplication $\Delta(x)=x^1\otimes x^2$, $(\Delta\otimes id)\Delta(x)=(id\otimes\Delta)\Delta(x)=x^1\otimes x^2\otimes x^3$, etc.

Note that by Lemma 2.2 in \cite{JL}
\begin{equation}\label{adm}
{\rm Ad}x(wz)={\rm Ad}x^1(w){\rm Ad}x^2(z).
\end{equation}

\begin{proposition}
Let $(A,\mu , \imath , \Delta , \varepsilon , S)$ be a Hopf algebra, $\mathcal{F}\in A\otimes A$ an invertible element satisfying conditions (\ref{twist1}) and (\ref{twist2}), $(A,\mu , \imath , \Delta^{\mathcal F} , \varepsilon , S^{\mathcal F})$ the twist of $A$ by ${\mathcal F}$. Denote by ${\rm Ad}^{\mathcal F}$ \index[not]{A@${\rm Ad}^{\mathcal F}$} the right adjoint action of $(A,\mu , \imath , \Delta^{\mathcal F} , \varepsilon , S^{\mathcal F})$ on itself. Then for all $x,z\in A$
\begin{equation}\label{adF}
{\rm Ad}^{\mathcal F}x(z)=\psi^{\mathcal{F}}{\rm Ad}x(\psi^{\mathcal{F}})^{-1}(z),
\end{equation}
where $\psi^{\mathcal{F}}:A\to A$ \index[not]{p@$\psi^{\mathcal{F}}$} is an invertible morphism of $A$, regarded as a module over the commutative ground ring, defined by
\begin{equation}\label{ps}
\psi^{\mathcal{F}}(z)=v^{\mathcal{F}}\sum_{i}S(c_i)zd_i,
\end{equation}
where $v^{\mathcal{F}}$ is given by (\ref{defv}), and ${\mathcal F}^{-1}=\sum_{i}c_i\otimes d_i$.

The inverse map $(\psi^{\mathcal{F}})^{-1}$ is given by the formula
\begin{equation}\label{ps-1}
(\psi^{\mathcal{F}})^{-1}(z)=\sum_{i}S(a_i)(v^{\mathcal{F}})^{-1}zb_i,
\end{equation}
where ${\mathcal F}=\sum_{i}a_i\otimes b_i$.
\end{proposition}

\begin{proof}
Firstly we check that formula (\ref{ps-1}) defines the inverse to $\psi^{\mathcal{F}}$. Indeed,
$$
(\psi^{\mathcal{F}})^{-1}\psi^{\mathcal{F}} (z)=\sum_{j}S(a_j)(v^{\mathcal F})^{-1}v^{\mathcal F}\sum_{i}S(c_i)zd_ib_j=\sum_{i,j}S(c_ia_j)zd_ib_j=z,
$$
as $1\otimes 1={\mathcal F}^{-1}{\mathcal F}=\sum_{i,j}c_ia_j\otimes d_ib_j$. 

Similarly,
$$
\psi^{\mathcal{F}}(\psi^{\mathcal{F}})^{-1} (z)=v^{\mathcal F}\sum_{j}S(c_j)\sum_{i}S(a_i)(v^{\mathcal F})^{-1}zb_id_j=\sum_{i,j}v^{\mathcal F}S(a_ic_j)(v^{\mathcal F})^{-1}zb_id_j=v^{\mathcal F}(v^{\mathcal F})^{-1}z=z,
$$
as
$1\otimes 1={\mathcal F}{\mathcal F}^{-1}=\sum_{i,j}a_ic_j\otimes b_id_j$.

Now by (\ref{deltaF}), (\ref{ad}), (\ref{ps}) and (\ref{ps-1}) we have
$$
{\rm Ad}^{\mathcal F}x(z)=v^{\mathcal F}S(a_ix^1c_j)(v^{\mathcal F})^{-1}zb_ix^2d_j=v^{\mathcal F}S(c_j)S(x^1)S(a_i)(v^{\mathcal F})^{-1}zb_ix^2d_j=\psi^{\mathcal{F}}{\rm Ad}x(\psi^{\mathcal{F}})^{-1}(z).
$$
This justifies (\ref{adF}). 

\end{proof}

\begin{proposition}
Let $(A,\mu , \imath , \Delta , \varepsilon , S)$ be a Hopf algebra, $\mathcal{F}\in A\otimes A$ an invertible element satisfying conditions (\ref{twist1}) and (\ref{twist2}), $v^{\mathcal F}\in A$ given by (\ref{defv}), and ${\mathcal F}^{-1}=\sum_{i}c_i\otimes d_i$.

Then 
\begin{equation}\label{dv}
\Delta v^{\mathcal F}={\mathcal F}^{-1}(v^{\mathcal F}\otimes v^{\mathcal F})(S\otimes S)({\mathcal F}^{-1}_{21}),
\end{equation}
and for all $z\in A$ 
\begin{equation}\label{dps}
\Delta \psi^{\mathcal{F}}(z)={\mathcal F}^{-1}(1\otimes v^{\mathcal F})\sum_{i,j}(\psi^{\mathcal{F}}({\rm Ad}(d_i^1c_j)(z^1))\otimes S(c_i)z^2d_i^2d_j) {\mathcal F},
\end{equation}
where $\psi^{\mathcal{F}}:A\to A$ is defined by (\ref{ps}) and we use the Sweedler notation for the comultiplication.
\end{proposition}

\begin{proof}
First we prove identity (\ref{dv}).
Denote as before ${\mathcal F}=\sum_{i}a_i\otimes b_i$. Then condition (\ref{twist2}) can be written in the form
\begin{equation}\label{twist2'}
\sum_{i,j}a_ia_j^1\otimes b_ia_j^2\otimes b_j=\sum_{i,j}a_i\otimes a_jb_i^1\otimes b_jb_i^2.
\end{equation}
Applying the antipode and the comultiplication to the last factor in the tensor product we obtain
$$
\sum_{i,j}a_ia_j^1\otimes b_ia_j^2\otimes (Sb_j)^1\otimes (Sb_j)^2=\sum_{i,j}a_i\otimes a_jb_i^1\otimes (S(b_jb_i^2))^1\otimes (S(b_jb_i^2))^2.
$$
Multiplying the first and the third, and the second and the fourth tensor factors we deduce 
$$
\sum_{i,j}a_ia_j^1(Sb_j)^1\otimes b_ia_j^2(Sb_j)^2=\sum_{i,j}a_i(S(b_jb_i^2))^1\otimes a_jb_i^1 (S(b_jb_i^2))^2,
$$
or, recalling that $\Delta$ is an algebra homomorphism, and $S$ is an anticoautomorphism,
$$
\sum_{i,j}a_i(a_jSb_j)^1\otimes b_i(a_jSb_j)^2=\sum_{i,j}a_iS((b_jb_i^2)^2)\otimes a_jb_i^1 S((b_jb_i^2)^1).
$$

Now, using the definition of $v^{\mathcal F}$ and the identity ${\mathcal F}=\sum_{i}a_i\otimes b_i$ in the left hand side, and recalling that $S$ is an algebra antiautomorphism and $\Delta$ is coassociative in the right hand side we have
$$
{\mathcal F}\Delta v^{\mathcal F}=\sum_{i,j}a_iS(b_j^2b_i^3)\otimes a_jb_i^1 S(b_j^1b_i^2)=\sum_{i,j}a_iS(b_i^3)S(b_j^2)\otimes a_jb_i^1S(b_i^2) S(b_j^1).
$$

Applying the identity $S(b_i^3)\otimes b_i^1S(b_i^2)=S(b_i^2)\otimes \varepsilon(b_i^1)=S(b_i)$ following from the defining properties of the counit and of the antipode we get
\begin{equation}\label{Fv}
{\mathcal F}\Delta v^{\mathcal{F}}=\sum_{i,j}a_iS(b_i^2)S(b_j^2))\otimes a_j\varepsilon(b_i^1) S(b_j^1)=\sum_{i,j}a_iS(b_i)S(b_j^2)\otimes a_j S(b_j^1)=(v^{\mathcal{F}}\otimes 1)\sum_{j}S(b_j^2)\otimes a_j S(b_j^1),
\end{equation}
where at the last step we also used the definition of $v^{\mathcal F}$.

Now we simplify the sum in the right hand side of the last formula. Applying the antipode to the second and the third factor in the tensor product in (\ref{twist2'}) and multiplying the first and the second factors after that we have
$$
\sum_{i,j}a_ia_j^1S(a_j^2)S(b_i)\otimes S(b_j)=\sum_{i,j}a_i S(b_i^1) S(a_j)\otimes S(b_i^2)S(b_j).
$$
Using the defining property of the antipode in the left hand side, and the formula ${\mathcal F}=\sum_{i}a_i\otimes b_i$ in the right hand side we can rewrite the previous identity in the form
$$
\sum_{i,j}a_i\varepsilon(a_j)S(b_i)\otimes S(b_j)=\sum_{i}a_i S(b_i^1) \otimes S(b_i^2)(S\otimes S)({\mathcal F}),
$$
or by (\ref{twist1}) and by the definition of $v^{\mathcal F}$ applied in the left hand side
$$
v^{\mathcal F}\otimes 1=\sum_{i}a_i S(b_i^1) \otimes S(b_i^2)(S\otimes S)({\mathcal F}),
$$
which is, by swapping the tensor factors, equivalent to
$$
1\otimes v^{\mathcal F}=\sum_{i}S(b_i^2)\otimes a_i S(b_i^1)(S\otimes S)({\mathcal F}_{21}),
$$
or
$$
\sum_{i}S(b_i^2)\otimes a_i S(b_i^1)=(1\otimes v^{\mathcal F})(S\otimes S)({\mathcal F}_{21}^{-1}).
$$

Substituting this expression into the right hand side of (\ref{Fv}) and multiplying by ${\mathcal F}^{-1}$ from the left we obtain (\ref{dv}).

Formula (\ref{dps}) is established in a similar way. Firstly we rewrite (\ref{twist2}) in the following form
\begin{equation}\label{twist3}
(\Delta \otimes id)({\mathcal F}^{-1}){\mathcal F}_{12}^{-1}=(id \otimes \Delta)({\mathcal F}^{-1}){\mathcal F}_{23}^{-1},
\end{equation}
or explicitly, using the expression ${\mathcal F}^{-1}=\sum_{i}c_i\otimes d_i$,
\begin{equation}\label{twist3'}
\sum_{i,j}c_i^1c_j\otimes c_i^2d_j\otimes d_i=\sum_{i,j}c_i\otimes d_i^1c_j\otimes d_i^2d_j.
\end{equation}

Applying the antipode to the first and to the second factor of the tensor product and the comultiplication to the last factor we obtain from this identity
$$
\sum_{i,j}S(c_j)S(c_i^1)\otimes S(d_j)S(c_i^2)\otimes d_i^1\otimes d_i^2=\sum_{i,j}S(c_i)\otimes S(c_j)S(d_i^1)\otimes (d_i^2d_j)^1\otimes (d_i^2d_j)^2,
$$
or, recalling that ${\mathcal F}^{-1}=\sum_{j}c_j\otimes d_j$ and swapping the first and the second factors in the tensor product,
$$
((S\otimes S)({\mathcal F}^{-1}_{21})\otimes 1\otimes 1)\sum_{i,j}S(c_i^2)\otimes S(c_i^1)\otimes d_i^1\otimes d_i^2=\sum_{i,j}S(c_j)S(d_i^1)\otimes S(c_i)\otimes (d_i^2d_j)^1\otimes (d_i^2d_j)^2.
$$ 
Multiplying this identity by $z^1\otimes z^2\otimes 1\otimes 1=\Delta(z)\otimes 1\otimes 1$ from the right we also have
$$
((S\otimes S)({\mathcal F}^{-1}_{21})\otimes 1\otimes 1)\sum_{i}S(c_i^2)z^1\otimes S(c_i^1)z^2\otimes d_i^1\otimes d_i^2=\sum_{i,j}S(c_j)S(d_i^1)z^1\otimes S(c_i)z^2\otimes (d_i^2d_j)^1\otimes (d_i^2d_j)^2.
$$ 

Now we multiply the first and the third, and the second and the fourth factors in the tensor product in this identity. Together with the coassociativity of the comultiplication in the right hand side this yields
\begin{equation}\label{SSF}
(S\otimes S)({\mathcal F}^{-1}_{21})\sum_{i}S(c_i^2)z^1 d_i^1\otimes S(c_i^1)z^2d_i^2=\sum_{i,j}S(c_j)S(d_i^1)z^1d_i^2d_j^1\otimes S(c_i)z^2 d_i^3d_j^2.
\end{equation}

Applying the identity
$$
\Delta(\sum_{i}S(c_i)zd_i)=\sum_{i}S(c_i^2)z^1 d_i^1\otimes S(c_i^1)z^2d_i^2
$$
in the left hand side of (\ref{SSF}) we obtain
$$
(S\otimes S)({\mathcal F}^{-1}_{21})\Delta(\sum_{i}S(c_i)zd_i)=\sum_{i,j}S(c_j)S(d_i^1)z^1d_i^2d_j^1\otimes S(c_i)z^2 d_i^3d_j^2
$$
which implies together with (\ref{ps}) and (\ref{dv})
\begin{equation}\label{SSF1}
\Delta \psi^{\mathcal{F}}(z)=\Delta(v^{\mathcal F})\Delta(\sum_{i}S(c_i)zd_i)={\mathcal F}^{-1}(v^{\mathcal F}\otimes v^{\mathcal F})(S\otimes S)({\mathcal F}^{-1}_{21})\Delta(\sum_{i}S(c_i)zd_i)=
\end{equation}
$$
={\mathcal F}^{-1}(v^{\mathcal F}\otimes v^{\mathcal F})\sum_{i,j}S(c_j)S(d_i^1)z^1d_i^2d_j^1\otimes S(c_i)z^2 d_i^3d_j^2.
$$

Now we bring the right hand side of this identity to the form indicated in (\ref{dps}).
For this purpose we rewrite (\ref{twist3'}) multiplying it by $\mathcal{F}_{23}$ from the right, by applying the antipode to the first factor in the tensor product, and by swapping the right hand side and the left hand side,
\begin{equation}\label{twist3''}
\sum_{j}S(c_j)\otimes d_j^1\otimes d_j^2=\sum_{j,k}S(c_k)S(c_j^1)\otimes c_j^2d_k\otimes d_j\mathcal{F}_{23},
\end{equation}
where we also renamed some summation indexes.

Applying this identity to the expression in the right hand side of (\ref{SSF1}) and recalling (\ref{ps}) we finally obtain 
$$
\Delta \psi^{\mathcal{F}}(z)
={\mathcal F}^{-1}(v^{\mathcal F}\otimes v^{\mathcal F})\sum_{i,j,k}S(c_k)S(c_j^1)S(d_i^1)z^1d_i^2c_j^2d_k\otimes S(c_i)z^2 d_i^3d_j{\mathcal F}=
$$
$$
={\mathcal F}^{-1}(1\otimes v^{\mathcal F})\sum_{i,j}(\psi^{\mathcal{F}}({\rm Ad}(d_i^1c_j)(z^1))\otimes S(c_i)z^2d_i^2d_j) {\mathcal F}
$$
which completes the proof.

\end{proof}

Apart from Drinfeld twists, there is another natural way for obtaining new comultiplications on a given Hopf algebra described in the following obvious proposition (see e.g. \cite{KT1}, Proposition 5.2).
\begin{proposition}\label{KhT}
Let $(A,\mu , \imath , \Delta , \varepsilon , S)$ be a Hopf algebra, $\upsilon:A\to A$ \index[not]{y@$\upsilon$} an algebra automorphism. Then $(A,\mu , \imath , \Delta^\upsilon , \varepsilon , S^\upsilon)$ \index[not]{ZZZZ@$(A,\mu , \imath , \Delta^\upsilon , \varepsilon , S^\upsilon)$} is a Hopf algebra, where
\begin{equation}\label{Du}
\Delta^\upsilon(x)=(\upsilon\otimes \upsilon)\Delta(\upsilon^{-1}x),~S^\upsilon=\upsilon S(\upsilon^{-1}x),~x\in A. \index[not]{D@$\Delta^\upsilon(~\cdot~)$} \index[not]{S@$S^\upsilon(~\cdot~)$}
\end{equation}

If $(A,\mu , \imath , \Delta , \varepsilon , S)$ is quasitriangular with universal R--matrix $\mathcal R$ then $(A,\mu , \imath , \Delta^\upsilon , \varepsilon , S^\upsilon)$ is quasitriangular with universal R--matrix
$$
{\mathcal R}^\upsilon=(\upsilon\otimes \upsilon){\mathcal R}. \index[not]{R@${\mathcal R}^\upsilon$}
$$
\end{proposition}

We denote the Hopf algebra $(A,\mu , \imath , \Delta^\upsilon , \varepsilon , S^\upsilon)$ by $A^\upsilon$ \index[not]{A@$A^\upsilon$} and call it {\it the twist of $A$ by $\upsilon$}. Denote by ${\rm Ad}^\upsilon$ \index[not]{A@${\rm Ad}^\upsilon$} the right adjoint action of $A^\upsilon$. One immediately has the following relation between the adjoint actions of  $A$ and of $A^\upsilon$.
\begin{proposition} 
Let $(A,\mu , \imath , \Delta , \varepsilon , S)$ be a Hopf algebra, $\upsilon:A\to A$ an algebra automorphism. Then
\begin{equation}\label{Adups}
{\rm Ad}^\upsilon x(z)=\upsilon({\rm Ad}\upsilon^{-1}(x)(\upsilon^{-1}(z))), x,z\in A.
\end{equation}
\end{proposition}

\begin{proof}
By (\ref{Du}) we have
$$
{\rm Ad}^\upsilon x(z)=\upsilon (S(\upsilon^{-1}\upsilon(\upsilon^{-1}x)^1))z\upsilon(\upsilon^{-1}x)^2=\upsilon (S(\upsilon^{-1}x)^1\upsilon^{-1}(z)(\upsilon^{-1}x)^2)=
\upsilon({\rm Ad}\upsilon^{-1}(x)(\upsilon^{-1}(z))).
$$
This completes the proof.

\end{proof}

Now we shall relate the two types of twists in the case when $A=U_h(\g)$. The following statements can be found in \cite{KT1}, Theorem 5.1 and Proposition 5.4.
\begin{proposition}\label{DhFT}
Let $w=s_{i_1}\ldots s_{i_k}$ be any reduced decomposition of an element $w\in W$, $T_w=T_{i_1}\ldots T_{i_k}$. Denote $\theta_w=\prod^{k}_{p=1}\theta_{\beta_p}$, \index[not]{t@$\theta_w$} where $\theta_{\beta_p}$ are defined in (\ref{dhtw}), and in the product $\theta_{\beta_p}$ appears on the left from $\theta_{\beta_q}$ if $p<q$. Then the following statements are true.

(i) The element $\theta_w^{-1}\in U_h(\g)\otimes U_h(\g)$  satisfies (\ref{twist1}) and (\ref{twist2}), and $\Delta_h^{\theta_w^{-1}}$ \index[not]{D@$\Delta_h^{\theta_w^{-1}}(~\cdot~)$} is a comultiplication on $U_h(\g)$. 

(ii) One has $\Delta_h^{\theta_w^{-1}}=\Delta_h^{T_w}$. \index[not]{D@$\Delta_h^{T_w}(~\cdot~)$}
\end{proposition}

\begin{proof}
(i) The proof is by induction over the length of $w$. When $w=s_i$ is a simple reflection 
condition (\ref{twist1}) for $\theta_w^{-1}=\theta_i^{-1}$ is clear from its definition and from the definition of the counit $\varepsilon_h$, and condition (\ref{twist2}) is the statement of Proposition 4.2.4 in \cite{Lus}, up to some change of the notation and different conventions on the comultiplication. The second claim in (i) for $\theta_w^{-1}=\theta_i^{-1}$ follows from the first one and from  part (i) of Proposition \ref{twdef}.

Now assume that the statement is true for all elements $w$ of length $k$, and for some $w=s_{i_1}\ldots s_{i_k}$ and $1\leq i\leq l$ the element $w':=s_iw$ has length $k+1$. Then by the definitions of $\theta_w$ and $\theta_{s_iw}$ one has 
$$
\theta_{s_iw}=\theta_{i}(T_i\otimes T_i)(\theta_{w}),
$$ 
so
\begin{equation}\label{tiw-1}
\theta_{s_iw}^{-1}=(T_i\otimes T_i)(\theta_{w}^{-1})\theta_{i}^{-1}.
\end{equation}

Condition (\ref{twist1}) for $\theta_{s_iw}^{-1}$ is again clear from its definition and from the definition of the counit $\varepsilon_h$. 

To check condition (\ref{twist2}) for $\theta_{s_iw}^{-1}$ we write using (\ref{tiw-1}) 
$$
(\theta_{s_iw}^{-1})_{12}(\Delta_h\otimes id)(\theta_{s_iw}^{-1})=(T_i\otimes T_i\otimes id)(\theta_{w}^{-1})_{12}(\theta_{i}^{-1})_{12}(\Delta_h\otimes id)(T_i\otimes T_i\theta_{w}^{-1}T_i^{-1}\otimes T_i^{-1})(\Delta_h\otimes id)(\theta_{i}^{-1}).
$$

Observe that by (\ref{DT}) $\Delta_h(T_i)=\theta_iT_i\otimes T_i$, so $\Delta_h(T_i(x))=\theta_i(T_i\otimes T_i)(\Delta_h(x))\theta_i^{-1}$, $x\in U_h(\g)$, and that $T_i$ is an algebra automorphism. Hence the previous identity can be rewritten as follows
$$
(\theta_{s_iw}^{-1})_{12}(\Delta_h\otimes id)(\theta_{s_iw}^{-1})=
$$
$$
=(T_i\otimes T_i\otimes T_i)(\theta_{w}^{-1})_{12}(\theta_{i}^{-1})_{12}(\theta_{i})_{12}(T_i\otimes T_i\otimes T_i)((\Delta_h\otimes id)(\theta_{w}^{-1}))(\theta_{i}^{-1})_{12}(\Delta_h\otimes id)(\theta_{i}^{-1})=
$$
$$
=(T_i\otimes T_i\otimes T_i)(\theta_{w}^{-1})_{12}(T_i\otimes T_i\otimes T_i)((\Delta_h\otimes id)(\theta_{w}^{-1}))(\theta_{i}^{-1})_{12}(\Delta_h\otimes id)(\theta_{i}^{-1})=
$$
$$
=(T_i\otimes T_i\otimes T_i)(\theta_{w}^{-1})_{12}(T_i\otimes T_i\otimes T_i)((\Delta_h\otimes id)(\theta_{w}^{-1}))(\theta_{i}^{-1})_{12}(\Delta_h\otimes id)(\theta_{i}^{-1})=
$$
\begin{equation}\label{ttti}
=(T_i\otimes T_i\otimes T_i)((\theta_{w}^{-1})_{12}(\Delta_h\otimes id)(\theta_{w}^{-1}))(\theta_{i}^{-1})_{12}(\Delta_h\otimes id)(\theta_{i}^{-1}).
\end{equation}

Now recall that by the induction assumption $\theta_{i}$ and $\theta_{w}$ satisfy (\ref{twist2}), so that after applying (\ref{twist2}) in the right hand side of the previous identity we have
$$
(\theta_{s_iw}^{-1})_{12}(\Delta_h\otimes id)(\theta_{s_iw}^{-1})=(T_i\otimes T_i\otimes T_i)((\theta_{w}^{-1})_{23}(id \otimes \Delta_h)(\theta_{w}^{-1}))(\theta_{i}^{-1})_{23}(id \otimes \Delta_h)(\theta_{i}^{-1}).
$$

Finally repeating in the opposite order arguments similar to those which lead us to (\ref{ttti}) we arrive at
$$
(\theta_{s_iw}^{-1})_{12}(\Delta_h\otimes id)(\theta_{s_iw}^{-1})=
$$
$$
=(T_i\otimes T_i)((\theta_{w}^{-1}))_{23}(\theta_{i}^{-1})_{23}(\theta_{i})_{23}(T_i\otimes T_i\otimes T_i)((id \otimes \Delta_h)(\theta_{w}^{-1}))(\theta_{i}^{-1})_{23}(id \otimes \Delta_h)(\theta_{i}^{-1})=
$$
$$
=(\theta_{s_iw}^{-1})_{23}(id \otimes \Delta_h)((T_i\otimes T_i)(\theta_{w}^{-1}))(id \otimes \Delta_h)(\theta_{i}^{-1})=(\theta_{s_iw}^{-1})_{23}(id \otimes \Delta_h)(\theta_{s_iw}^{-1})
$$
which confirms (\ref{twist2}) for $\theta_{s_iw}^{-1}$ and establishes the induction step. Thus the first claim in (i) is proved. 

The second claim in (i) for $\theta_w^{-1}$ with arbitrary $w$ follows from the first one and from  part (i) of Proposition \ref{twdef}.

(ii) Observe that by (\ref{dhtw}) $\Delta_h(T_w^{-1})=T_w^{-1}\otimes T_w^{-1}\theta_w^{-1}$, and hence
\begin{equation}\label{DhTi}
\Delta_h(T_w^{-1}(x))=(T_w^{-1}\otimes T_w^{-1})(\theta_w^{-1}\Delta_h(x)\theta_w), x\in U_h(\g),
\end{equation}
or, by the definitions of $\Delta_h^{T_w}$ and of $\Delta_h^{\theta_w^{-1}}$,
$$
\Delta_h^{T_w}(x)=(\theta_w^{-1}\Delta_h(x)\theta_w)=\Delta_h^{\theta_w^{-1}}(x), x\in U_h(\g).
$$
This completes the proof.

\end{proof}

We shall need the following proposition which is a version of Proposition \ref{DhFT} for the algebras $U_h^{s}({\frak g})$. To simplify the notation we shall identify $U_h^{s}({\frak g})$ with $U_h({\frak g})$ as algebras using the isomorphisms from Proposition \ref{newreal} with $n_{ij}$ given by (\ref{nijfix}), and omit $\psi_{\{ n_{ij}\}}$ in all formulas. 
\begin{proposition}
Let $w=s_{i_1}\ldots s_{i_k}$ be any reduced decomposition of an element $w\in W$, $T_w=T_{i_1}\ldots T_{i_k}$. Denote $\theta_w^s=\prod^{k}_{p=1}\theta_{\beta_p}^s$, where $\theta_{\beta_p}^s$ \index[not]{t@$\theta_w^s$} are defined in (\ref{DTs}), and in the product $\theta_{\beta_p}^s$ appears on the left from $\theta_{\beta_q}^s$ if $p<q$. Then the following statements are true.

(i) The element ${\theta_w^s}^{-1}\in U_h^s(\g)\otimes U_h^s(\g)$  satisfies (\ref{twist1}) and (\ref{twist2}), for $\varepsilon=\varepsilon_s$ and $\Delta=\Delta_s$, and $\Delta_s^{{\theta_w^s}^{-1}}$ \index[not]{D@$\Delta_s^{{\theta_w^s}^{-1}}(~\cdot~)$} is a comultiplication on $U_h^s(\g)$. 

(ii) One has $\Delta_s^{{\theta_w^s}^{-1}}=\Delta_{w^{-1}sw}^{T_w}$. \index[not]{D@$\Delta_{w^{-1}sw}^{T_w}(~\cdot~)$}
\end{proposition}

\begin{proof}

(i) Firstly by (\ref{DTs})
\begin{equation}\label{thths}
\mathcal{F}_s\theta_w\mathcal{F}_s^{-1}=\theta_w^s.
\end{equation}

From the definition of the commutative Hopf subalgebra $U({\frak h})[[h]]\simeq U_h({\frak h})\subset U_h({\frak g})\simeq U_h^s({\frak g})$ it immediately follows that $\mathcal{F}_s^{-1}\in U_h^s({\frak g})\otimes U_h^s({\frak g})$ satisfies (\ref{twist1}) and (\ref{twist2}) for $\varepsilon=\varepsilon_s$ and $\Delta=\Delta_s$, so $\Delta_s^{\mathcal{F}_s^{-1}}$ is a comultiplication on $U_h^s({\frak g})$ according to Proposition \ref{twdef} (i). In fact by formula (\ref{defds})
\begin{equation}\label{defds1}
\Delta_s^{\mathcal{F}_s^{-1}}=\Delta_h.
\end{equation}

By Proposition \ref{DhFT} (i) and (ii) the element $\theta_w^{-1}\in U_h(\g)\otimes U_h(\g)$  satisfies (\ref{twist1}) and (\ref{twist2}) with $\varepsilon=\varepsilon_h$ and $\Delta=\Delta_h$, $\Delta_h^{\theta_w^{-1}}$ is a comultiplication on $U_h(\g)$, and $\Delta_h^{\theta_w^{-1}}=\Delta_h^{T_w}$. Using this identity and (\ref{defds1}) we obtain
\begin{equation}\label{defds2}
(\Delta_s^{\mathcal{F}_s^{-1}})^{\theta_w^{-1}}=\Delta_h^{\theta_w^{-1}}=\Delta_h^{T_w}.
\end{equation}

Finally note that on the subalgebra $U_h({\frak h})\subset U_h({\frak g})$ one has $\Delta_h^{T_w}=\Delta_h$ and from the definition of the commutative Hopf subalgebra $U({\frak h})[[h]]\simeq U_h({\frak h})\subset U_h({\frak g})\simeq U_h^s({\frak g})$ it immediately follows that $\mathcal{F}_s\in U_h({\frak g})\otimes U_h({\frak g})$ satisfies (\ref{twist1}) and (\ref{twist2}) for $\varepsilon=\varepsilon_h$ and $\Delta=\Delta_h^{T_w}$, so $(\Delta_h^{T_w})^{\mathcal{F}_s}$ is a comultiplication on $U_h({\frak g})$ according to Proposition \ref{twdef} (i). 

Using this observation together with (\ref{defds2}) and Proposition \ref{twdef} (iii) twice we deduce that $\mathcal{F}_s\theta_w^{-1}\mathcal{F}_s^{-1}={\theta_w^s}^{-1}$ satisfies (\ref{twist1}) and (\ref{twist2}) for $\varepsilon=\varepsilon_s$ and $\Delta=\Delta_s$, $\Delta_s^{{\theta_w^s}^{-1}}$ is a comultiplication on $U_h^s(\g)$, and
\begin{equation}\label{defds3}
\Delta_s^{{\theta_w^s}^{-1}}=((\Delta_s^{\mathcal{F}_s^{-1}})^{\theta_w^{-1}})^{\mathcal{F}_s}=(\Delta_h^{T_w})^{\mathcal{F}_s}.
\end{equation}

(ii) Observe that by (\ref{DTw}) and by the definition (\ref{Ftw}) of $\mathcal{F}_s$ one has $\Delta_s(T_w^{-1})=T_w^{-1}\otimes T_w^{-1}(T_w\otimes T_w)(\mathcal{F}_s)\mathcal{F}_s^{-1}{\theta_w^s}^{-1}$, and hence by the definition of $\Delta_s^{T_w}$ we obtain
$$
\Delta_s^{T_w}(x)=(T_w\otimes T_w)(\mathcal{F}_s)\mathcal{F}_s^{-1}{\theta_w^s}^{-1}\Delta_s(x)\theta_w^s\mathcal{F}_s(T_w\otimes T_w)(\mathcal{F}_s^{-1}), x\in U_h^s(\g).
$$

Using the definition of $\Delta_s^{T_w}$ and conjugating by $\mathcal{F}_s(T_w\otimes T_w)(\mathcal{F}_s^{-1})$ one can rewrite this identity as follows
$$
\mathcal{F}_s(T_w\otimes T_w)(\mathcal{F}_s^{-1})(T_w\otimes T_w)(\Delta_s(T_w^{-1}x))(T_w\otimes T_w)(\mathcal{F}_s)\mathcal{F}_s^{-1}={\theta_w^s}^{-1}\Delta_s(x)\theta_w^s.
$$
Using (\ref{defds}) and the definition of the braid group action on the generators of $U_h({\frak h})$ we can further transform the left hand side and obtain
$$
(T_w\otimes T_w)((T_w^{-1}\otimes T_w^{-1})(\mathcal{F}_s)\mathcal{F}_s^{-1}\mathcal{F}_s\Delta_h(T_w^{-1}x)\mathcal{F}_s^{-1}\mathcal{F}_s(T_w^{-1}\otimes T_w^{-1})(\mathcal{F}_s^{-1}))={\theta_w^s}^{-1}\Delta_s(x)\theta_w^s,
$$
or
\begin{equation}\label{ttdsf}
(T_w\otimes T_w)((T_w^{-1}\otimes T_w^{-1})(\mathcal{F}_s)\Delta_h(T_w^{-1}x)(T_w^{-1}\otimes T_w^{-1})(\mathcal{F}_s^{-1}))={\theta_w^s}^{-1}\Delta_s(x)\theta_w^s.
\end{equation}

By definition (\ref{Ftw}) $(T_w^{-1}\otimes T_w^{-1})(\mathcal{F}_s)=\mathcal{F}_{w^{-1}sw}$, and hence by (\ref{defds})
$$
(T_w^{-1}\otimes T_w^{-1})(\mathcal{F}_s)\Delta_h(T_w^{-1}x)(T_w^{-1}\otimes T_w^{-1})(\mathcal{F}_s^{-1})=\mathcal{F}_{w^{-1}sw}\Delta_h(T_w^{-1}x)\mathcal{F}_{w^{-1}sw}^{-1}=\Delta_{w^{-1}sw}(T_w^{-1}x).
$$ 
Using this expression in the left hand side of (\ref{ttdsf}) and recalling the definitions of $\Delta_s^{{\theta_w^s}^{-1}}$ and of $\Delta_{w^{-1}sw}^{T_w}$ we derive
$\Delta_{w^{-1}sw}^{T_w}=\Delta_s^{{\theta_w^s}^{-1}}$. This completes the proof.

\end{proof}

Rewriting the identity in part (ii) of the previous proposition in the form 
$$
\Delta_{w^{-1}sw}=(\Delta_s^{{\theta_w^s}^{-1}})^{T_w^{-1}}
$$
and applying formulas (\ref{adF}) and (\ref{Adups}) we obtain the following corollary. 
\begin{corollary}
Denote by ${\rm Ad}_s$ and ${\rm Ad}'_s$ \index[not]{A@${\rm Ad}_s$} \index[not]{A@${\rm Ad}'_s$} the right and the left adjoint action of the Hopf algebra $U_h^s(\g)$, respectively. Then
\begin{equation}\label{adstt}
{\rm Ad}_{w^{-1}sw}x=({\rm Ad}_s^{{\theta_w^s}^{-1}})^{T_w^{-1}}(x)=T_w^{-1}\psi^{{\theta_w^s}^{-1}}{\rm Ad}_s(T_wx)(\psi^{{\theta_w^s}^{-1}})^{-1}T_w,
\end{equation} 
where $\psi^{{\theta_w^s}^{-1}}$ \index[not]{p@$\psi^{{\theta_w^s}^{-1}}$} is defined by (\ref{ps}) with the help of the antipode $S_s$ and of the multiplication in the Hopf algebra $U_h^s(\g)$.
\end{corollary}

From (\ref{dps}) and (\ref{DTs}) we also obtain the following proposition which will be crucial for the definition of the Zhelobenko type operators.
\begin{proposition}
Let $w=s_{i_1}\ldots s_{i_k}$ be any reduced decomposition of an element $w\in W$, $T_w=T_{i_1}\ldots T_{i_k}$. Then
\begin{equation}\label{w-1sw}
\Delta_{w^{-1}sw} (T_w^{-1}\psi^{{\theta_w^s}^{-1}}(z)T_w)
=T_w^{-1}\otimes T_w^{-1}(1\otimes v^{{\theta_w^s}^{-1}})\sum_{i,j}(\psi^{{\theta_w^s}^{-1}}({\rm Ad}_s(d_i^1c_j)(z^1))\otimes S_s(c_i)z^2d_i^2d_j) T_w\otimes T_w,
\end{equation}
where ${\theta_w^s}^{-1}=\sum_ia_i\otimes b_i$, $v^{{\theta_w^s}^{-1}}=\sum_ia_iS_s(b_i)$, \index[not]{v@$v^{{\theta_w^s}^{-1}}$} $\theta_w^s=\sum_i c_i\otimes d_i$, $\Delta_s z=z^1\otimes z^2$, $\Delta_s d_i=d_i^1\otimes d_i^2$.
\end{proposition}

\begin{proof}
From (\ref{dps}) with $\Delta=\Delta_s$, $\mathcal{F}={\theta_w^s}^{-1}$, $S=S_s$ we obtain
\begin{equation}\label{dpst}
\Delta_s (\psi^{{\theta_w^s}^{-1}}(z))=\theta_w^s(1\otimes v^{{\theta_w^s}^{-1}})\sum_{i,j}(\psi^{{\theta_w^s}^{-1}}({\rm Ad}_s(d_i^1c_j)(z^1))\otimes S_s(c_i)z^2d_i^2d_j) {\theta_w^s}^{-1},
\end{equation}
where ${\theta_w^s}^{-1}=\sum_ia_i\otimes b_i$, $v^{{\theta_w^s}^{-1}}=\sum_ia_iS_s(b_i)$, $\theta_w^s=\sum_i c_i\otimes d_i$, $\Delta_s z=z^1\otimes z^2$.

By (\ref{DTs}) and (\ref{Kdefs})
\begin{align}\label{DTst}
\Delta_s(T_w)=\prod^{k}_{p=1}\theta_{\beta_p}^s\mathcal{F}_s(T_w \otimes T_w)(\mathcal{F}_s^{-1})T_w\otimes T_w=\theta_w^s\mathcal{F}_s(T_w \otimes T_w)(\mathcal{F}_s^{-1})T_w\otimes T_w= \\
=\theta_w^s T_w\otimes T_w (T_w^{-1}\otimes T_w^{-1})(\mathcal{F}_s)\mathcal{F}_s^{-1}=\theta_w^s T_w\otimes T_w \mathcal{F}_{w^{-1}sw}\mathcal{F}_s^{-1}. \nonumber
\end{align}

Now (\ref{dpst}) and (\ref{DTst}) imply
$$
\Delta_s (T_w^{-1}\psi^{{\theta_w^s}^{-1}}(z)T_w)=\Delta_s (T_w^{-1})\Delta_s(\psi^{{\theta_w^s}^{-1}}(z))\Delta_s(T_w)= 
$$
$$
=\mathcal{F}_s\mathcal{F}_{w^{-1}sw}^{-1}(T_w^{-1}\otimes T_w^{-1}){\theta_w^s}^{-1}\theta_w^s(1\otimes v^{{\theta_w^s}^{-1}})\sum_{i,j}(\psi^{{\theta_w^s}^{-1}}({\rm Ad}_s(d_i^1c_j)(z^1))\otimes S_s(c_i)z^2d_i^2d_j) {\theta_w^s}^{-1}\theta_w^s (T_w\otimes T_w) \mathcal{F}_{w^{-1}sw}\mathcal{F}_s^{-1}= 
$$
$$
=\mathcal{F}_s\mathcal{F}_{w^{-1}sw}^{-1}T_w^{-1}\otimes T_w^{-1}(1\otimes v^{{\theta_w^s}^{-1}})\sum_{i,j}(\psi^{{\theta_w^s}^{-1}}({\rm Ad}_s(d_i^1c_j)(z^1))\otimes S_s(c_i)z^2d_i^2d_j) T_w\otimes T_w \mathcal{F}_{w^{-1}sw}\mathcal{F}_s^{-1}.
$$

Conjugating this formula by $\mathcal{F}_{w^{-1}sw}\mathcal{F}_s^{-1}$ and observing that by the definition of $\Delta_s$ one has 
\begin{align*}
\mathcal{F}_{w^{-1}sw}\mathcal{F}_s^{-1}\Delta_s(~\cdot~)\mathcal{F}_s\mathcal{F}_{w^{-1}sw}^{-1}=
\mathcal{F}_{w^{-1}sw}\mathcal{F}_s^{-1}\mathcal{F}_s\Delta_h(~\cdot~)\mathcal{F}_s^{-1}\mathcal{F}_s\mathcal{F}_{w^{-1}sw}^{-1}= \\
=\mathcal{F}_{w^{-1}sw}\Delta_h(~\cdot~)\mathcal{F}_{w^{-1}sw}^{-1}=\Delta_{w^{-1}sw}(~\cdot~)
\end{align*}
we obtain
$$
\Delta_{w^{-1}sw} (T_w^{-1}\psi^{{\theta_w^s}^{-1}}(z)T_w)
=T_w^{-1}\otimes T_w^{-1}(1\otimes v^{{\theta_w^s}^{-1}})\sum_{i,j}(\psi^{{\theta_w^s}^{-1}}({\rm Ad}_s(d_i^1c_j)(z^1))\otimes S_s(c_i)z^2d_i^2d_j) T_w\otimes T_w.
$$
This completes the proof.

\end{proof}


\section{Some forms and specializations of quantum groups}\label{QGspec}

\pagestyle{myheadings}
\markboth{CHAPTER~\thechapter.~QUANTUM GROUPS}{\thesection.~FORMS AND SPECIALIZATIONS OF QUANTUM GROUPS}

\setcounter{equation}{0}
\setcounter{theorem}{0}

In order to define q-W--algebras we shall actually need not the algebras $U_h^{s}({\frak g})$ themselves but some their forms defined over certain rings. They are similar to the rational form and the restricted integral form for the standard quantum group $U_h(\g)$.  \index{form!of a quantum group} The motivations of the definitions given below will be clear in Section \ref{qplgroups}. The results below are slight modifications of similar statements for $U_h(\g)$.

We start with a very important technical lemma which will play the key role in the definition of q-W--algebras. Below we keep the notation introduced in Section \ref{background}. 

Let $s\in W$ be an element of the Weyl group. Recall that by formula (\ref{inv}) $s$ can be represented as a product of two involutions, $s=s^1s^2$, where $s^1=s_{\gamma_1}\ldots s_{\gamma_{\widetilde{l}}}$, $s^2=s_{\gamma_{\widetilde{l}+1}}\ldots s_{\gamma_{l'}}$, and the roots $\gamma_1, \ldots , \gamma_{l'}$ form a basis of the subspace ${\h'}^*\subset \h^*$ on which $s$ acts without fixed points. We shall study the matrix elements of the Cayley transform of the restriction of
$s$ to ${\h'}^*$ with respect to this basis.
\begin{lemma}\label{tmatrel}
Let $P_{{\h'}^*}$ be the orthogonal projection operator onto ${{\h'}^*}$ in $\h^*$, with respect to the symmetric invariant bilinear form.
Then the matrix elements of the operator ${1+s \over 1-s }P_{{\h'}^*}$ in the basis $\gamma_1, \ldots , \gamma_{l'}$ are of the form
\begin{equation}\label{matrel}
\left\langle  {1+s \over 1-s }P_{{\h'}^*}\gamma_i , \gamma_j \right\rangle=
\varepsilon_{ij}\left\langle \gamma_i,\gamma_j\right\rangle,
\end{equation}
where
$$
\varepsilon_{ij} =\left\{ \begin{array}{ll}
-1 & i <j \\
0 & i=j \\
1 & i >j
\end{array}
\right  . .
$$
\end{lemma}
\begin{proof} 
First we calculate the matrix of the element $s$ with respect to the
basis $\gamma_1, \ldots , \gamma_{l'}$. We obtain this matrix in the form of the Gauss
decomposition of the operator $s:{\h'}^*\rightarrow {\h'}^*$.

Let $z_i=s \gamma_{i}$. Recall that $s_{\gamma_i}(\gamma_j)=\gamma_j-A_{ij}\gamma_i$, $A_{ij}=\gamma_i^\vee(\gamma_j)$. \index[not]{A@$A_{ij}$}
Define
\begin{equation}\label{y}
y_{i}=s_{\gamma_1}\ldots s_{\gamma_{i-1}}\gamma_{i}.
\end{equation}
Using this definition the elements $z_i$ may be represented as 
$$
z_i=s_{\gamma_1}\ldots s_{\gamma_{l'}}\gamma_{i}=s_{\gamma_1}\ldots s_{\gamma_{l'-1}}(\gamma_{i}-A_{l'i}\gamma_{l'})=s_{\gamma_1}\ldots  s_{\gamma_{l'-1}}\gamma_{i}-A_{l'i}y_{l'}=\ldots 
=y_i -\sum_{k \geq i} A_{k i}y_{k}.
$$

Using the matrix notation we can rewrite the last formula as follows
\begin{equation}\label{2*}
z_{i}=
\sum_{k =1}^{l'}(I-V)_{k i}y_{k},~ \mbox{ where } V_{k i}=
\left\{ \begin{array}{ll}
A_{k i} & k\geq i \\
0 & k < i
\end{array}
\right  .
\end{equation}

To calculate the matrix of the operator $s:{\h'}^*\rightarrow {\h'}^*$ with respect to the basis $\gamma_1, \ldots , \gamma_{l'}$ we have to express
the elements $y_{i}$ via $\gamma_1, \ldots , \gamma_{l'}$.
Applying the definition of reflections to (\ref{y}) we obtain
\[
y_{i}=s_{\gamma_1}\ldots s_{\gamma_{i-1}}\gamma_{i}=s_{\gamma_1}\ldots s_{\gamma_{i-2}}\gamma_{i}-A_{i-1 i}y_{i-1}=\ldots=\gamma_{i}-\sum_{k<i}A_{ki}y_{k}.
\]
Therefore
\[
\gamma_{i}=\sum_{k =1}^{l'}(I+U)_{k i}y_{k} ~, \mbox{ where } U_{ki}=
\left\{ \begin{array}{ll}
A_{ki} & k<i \\
0 & k \geq i
\end{array}
\right .
\]
Thus
\begin{equation}\label{1*}
y_{k}=\sum_{k =1}^{l'}(I+U)^{-1}_{jk}\gamma_{j}.
\end{equation}

Summarizing (\ref{1*}) and (\ref{2*}) we obtain
\begin{equation}\label{**}
s \gamma_i=\sum_{k =1}^{l'}\left( (I+U)^{-1}(I-V) \right)_{ki}\gamma_k .
\end{equation}
This implies
\begin{equation}\label{3*}
{1+s \over 1-s}P_{{\h'}^*}\gamma_i=\sum_{k =1}^{l'}\left( {2I+U-V \over U+V}\right)_{ki}\gamma_k .
\end{equation}

Observe that $(U+V)_{ki}=A_{ki}$ and $(2I+U-V)_{ij}=-A_{ij}\varepsilon_{ij}$.
Substituting these expressions into (\ref{3*}) we get
\begin{eqnarray}
\left\langle  {1+s \over 1-s }P_{{\h'}^*}\gamma_i , \gamma_j \right\rangle =
-(A^{-1})_{kp}\varepsilon_{pi} A_{pi}\left\langle \gamma_j,\gamma_k\right\rangle=\varepsilon_{ij}\left\langle \gamma_i,\gamma_j\right\rangle.
\end{eqnarray}
This completes the proof of the lemma.

\end{proof}

Let $\gamma_i^*$, $i=1,\ldots, l'$ \index[not]{g@$\gamma_i^*$} be the basis of $\h'^*$ dual to $\gamma_i$, $i=1,\ldots, l'$ with respect to the restriction of the symmetric invariant bilinear form $\left\langle ~\cdot~,~\cdot~\right\rangle$ to $\h'^*$. Since the numbers $\left\langle \gamma_i,\gamma_j\right\rangle$ are integer, each element $\gamma_i^*$ has the form $\gamma_i^*=\sum_{j=1}^{l'}\bar{m}_{ij}\gamma_j$, where $\bar{m}_{ij}\in \mathbb{Q}$. Therefore by Lemma \ref{tmatrel} and using for simple roots $\alpha_i$ the decomposition of the form $P_{{\h'}^*}\alpha_i=\sum_{p=1}^{l'}\left\langle \alpha_i,\gamma_p\right\rangle\gamma_p^*=\sum_{p,q=1}^{l'}\left\langle \alpha_i,\gamma_p\right\rangle \bar{m}_{pq}\gamma_q$ we deduce that the numbers
\begin{eqnarray}
p_{ij}:=\frac{1}{2d_j}\left\langle  {1+s \over 1-s }P_{{\h'}^*}\alpha_i,\alpha_j\right\rangle= \qquad \qquad \qquad \qquad \qquad \qquad \qquad \label{rat} \\ \qquad \qquad =\frac{1}{2d_j}\sum_{k,l,p,q=1}^{l'}\left\langle \gamma_k,\alpha_i\right\rangle\left\langle \gamma_l,\alpha_j\right\rangle\left\langle  {1+s \over 1-s }P_{{\h'}^*}\gamma_p,\gamma_q\right\rangle \bar{m}_{kp}\bar{m}_{lq},~i,j=1,\ldots,l \nonumber \index[not]{p@$p_{ij}$}
\end{eqnarray}
are rational, $p_{ij}\in \mathbb{Q}$, as all factors in the products in the sum in the right hand side are rational. Denote by $d$ \index[not]{d@$d$} an integer number divisible by all the denominators of the rational numbers $p_{ij}$, $i,j=1,\ldots, l$.

Let $\bar{r}\in \mathbb{N}$ \index[not]{r@$\bar{r}$} be such that $a^{-1}_{ij}\in \frac{1}{\bar{r}}\mathbb{Z}$, $i,j=1,\ldots ,l$. Let  
$U_q^{s}({\frak g})$ \index[not]{U@$U_q^s(\g)$} be the $\mathbb{C}(q^{\frac{1}{d{\bar{r}}^2}})$-algebra generated by the elements $e_i , f_i , L_i^{\pm 1}, t_i^{\pm 1},~i=1, \ldots, l$ with the same relations as the relations in $U_h^{s}({\frak g})$ for the generators denoted by the same symbols, where we assume that $t_i^{\pm 1}=\exp(\pm\frac{h\kappa}{d}Y_i)$. \index[not]{t@$t_i$} The coefficients of these relations indeed belong to $\mathbb{C}(q^{\frac{1}{d{\bar{r}}^2}})$, where $q^{\frac{1}{d{\bar{r}}^2}}=e^{h\frac{1}{d{\bar{r}}^2}}$.

Let $U_q({\frak g})$ \index[not]{U@$U_q(\g)$} be the $\mathbb{C}(q^{\frac{1}{d{\bar{r}}^2}})$-algebra generated by the elements $X_i^\pm, L_i^{\pm 1}, t_i^{\pm 1},~i=1, \ldots, l$ subject to the same relations as the relations in $U_h({\frak g})$ for the generators denoted by the same symbols, where we assume that $t_i^{\pm 1}=\exp(\pm\frac{h\kappa}{d}Y_i)$. The coefficients of these relations indeed belong to $\mathbb{C}(q^{\frac{1}{d{\bar{r}}^2}})$, where $q^{\frac{1}{d{\bar{r}}^2}}=e^{h\frac{1}{d{\bar{r}}^2}}$.


Note that by the choice of $d$ we have $q^{c_{ij}}\in \mathbb{C}[q^{\frac{\kappa}{d}},q^{-\frac{\kappa}{d}}]$.

The second form \index{form!of a quantum group} of $U_h^{s}({\frak g})$  is a subalgebra $U_\mathcal{A}^{s}(\g)$ \index[not]{U@$U_\mathcal{A}^s(\g)$} in $U_q^s(\g)$ over the ring $\mathcal{A}:=\mathcal{P}[q^{\frac{1}{d{\bar{r}}^2}},q^{-\frac{1}{d{\bar{r}}^2}}]$. \index[not]{A@$\mathcal{A}$} $U_\mathcal{A}^{s}(\g)$ is the subalgebra in $U_q^s(\g)$ generated over $\mathcal{A}$ by the elements 
$$
L_i^{\pm 1}, t_i^{\pm 1},~{K_i -K_i^{-1} \over q_i -q_i^{-1}},~e_i,~f_i,~i=1,\ldots ,l.
$$ 
Denote also by $U_\mathcal{A}(\g)$ \index[not]{U@$U_\mathcal{A}(\g)$} the subalgebra in $U_q(\g)$ generated over $\mathcal{A}$ by the elements 
$$
L_i^{\pm 1}, t_i^{\pm 1},~{K_i -K_i^{-1} \over q_i -q_i^{-1}},~X_i^\pm,~i=1,\ldots ,l.
$$

For the solution $n_{ij}=\frac{1}{2d_j}c_{ij}$ to equations (\ref{eqpi}) the root vectors $e_{\beta},f_\beta$ belong to all the above introduced forms of $U_h({\frak g})$.

For any normal ordering of $\Delta_+$ we denote $e_{\beta}^{(k)}=\frac{e_{\beta}^k}{[k]_{q_\alpha}!},f_\beta^{(k)}=\frac{f_{\beta}^k}{[k]_{q_\alpha}!}$, \index[not]{e@$e_{\beta}^{(k)}$} \index[not]{f@$f_{\beta}^{(k)}$} where $e_{\beta}, f_{\beta}$ are the corresponding quantum root vectors from Proposition \ref{rootss}. Let 
$U_{\mathcal{B}}^{s, res}(\g)$ \index[not]{U@$U_\mathcal{B}^{s, res}(\g)$} be the subalgebra in $U_q^s(\g)$ generated over $\mathcal{B}:=\mathbb{C}[q^{\frac{1}{d{\bar{r}}^2}},q^{-\frac{1}{d{\bar{r}}^2}}]$ \index[not]{B@$\mathcal{B}$} by the elements 
$$
L_i^{\pm 1}, t_i^{\pm 1},~e_{i}^{(k)},~f_{i}^{(k)},~i=1,\ldots ,l,  k\geq 1.
$$ 
Denote also by $U_\mathcal{B}^{res}(\g)$ \index[not]{U@$U_\mathcal{B}^{res}(\g)$} the subalgebra in $U_q(\g)$ generated over $\mathcal{B}$ by the elements 
$$
L_i^{\pm 1}, t_i^{\pm 1},~(X_i^\pm)^{(k)},~i=1,\ldots ,l, k\geq 1.
$$

Let $\varepsilon\in \mathbb{C}^*$. Fix a root of $\varepsilon$ of degree $d\bar{r}^2$, $\varepsilon^{\frac{1}{d{\bar{r}}^2}}$ \index[not]{e@$\varepsilon^{\frac{1}{d{\bar{r}}^2}}$} and if $\varepsilon=1$ put $\varepsilon^{\frac{1}{r^2d}}=1$.
Then  we define the specialization $U_\varepsilon^{s, res}(\g)$ of $U_\mathcal{B}^{s,res}(\g)$ by
\index{specialization!of a quantum group}
$$
U_\varepsilon^{s, res}(\g)=U_\mathcal{B}^{s,res}(\g)/(q^{\frac{1}{d{\bar{r}}^2}}-\varepsilon^{\frac{1}{d{\bar{r}}^2}})U_\mathcal{B}^{s, res}(\g). \index[not]{U@$U_\varepsilon^{s, res}(\g)$}
$$

If $[2]_\varepsilon\neq 0$  for $\g$ of type $B_l$, $C_l$ or $F_4$, and $[2]_\varepsilon, [3]_\varepsilon\neq 0$ for $\g$ of type $G_2$, we also define the specialization $U_\varepsilon^s(\g)$ \index[not]{U@$U_\varepsilon^s(\g)$} of $U_\mathcal{A}^{s}(\g)$,
$$
U_\varepsilon^s(\g)=U_\mathcal{A}^{s}(\g)/(q^{\frac{1}{d{\bar{r}}^2}}-\varepsilon^{\frac{1}{d{\bar{r}}^2}})U_\mathcal{A}^{s}(\g).
$$

The algebras $U_q^s(\g)$, $U_\mathcal{A}^{s}(\g)$, $U_\mathcal{B}^{s, res}(\g)$, $U_\varepsilon^{s,res}(\g)$ and $U_\varepsilon^s(\g)$ are in fact Hopf algebras with the comultiplication given on generators by the same formulas as in $U_h^s(\g)$ with $q^{\frac{1}{d{\bar{r}}^2}}=e^{h\frac{1}{d{\bar{r}}^2}}$. 

The elements $t_i$ and $L_i$ are central in the algebra $U_1^s(\g)$, and the quotient of $U_1^s(\g)$ by the two--sided ideal generated by $t_i-1$ and $L_i-1$ is isomorphic to $U(\g)$. Note that none of the forms and specializations \index{specialization!of a quantum group} \index{form!of a quantum group} of $U_h^s(\g)$ introduced above is quasitriangular. Similarly, the quotient of $U_1^{s,res}(\g)$ by the two--sided ideal generated by $t_i-1$ and $L_i-1$ is isomorphic to $U(\g)$, and the quotient of $U_1^{res}(\g)$ by the two--sided ideal generated by and $K_i-1$ is isomorphic to $U(\g)$ (see \cite{ChP}, Proposition 9.3.10).

The algebra isomorphism $\psi_{\{n_{ij}\}}$ with $n_{ij}=\frac{1}{2d_j}c_{ij}$ induces isomorphisms of $U_h(\g)$ and $U_h^{s}(\g)$ and of the forms and specializations \index{specialization!of a quantum group} \index{form!of a quantum group} of $U_h(\g)$ and $U_h^{s}(\g)$ with the superscript $s$ defined above and of their counterparts with the superscript s dropped. We shall always identify them using these isomorphisms.
To simplify the notation we shall also write, if it does not cause confusion, $e_{\beta}=X_{\beta}^+e^{hK_s\beta^\vee}$, $f_{\beta}=e^{-hK_s\beta^\vee}X_{\beta}^-$, $\beta \in \Delta_+$.

Note that the structure constants in the commutation relations for $U_\mathcal{A}^{s}(\g)= U_\mathcal{A}(\g)$ actually belong to $\mathcal{P}[q^{\frac{\kappa}{d}},q^{-\frac{\kappa}{d}}]$, so if $\kappa$ is divisible by $d$ the specialization $U_\varepsilon^s(\g)$ actually depends on $\varepsilon\in \mathbb{C}^*$ but not on its root $\varepsilon^{\frac{1}{d{\bar{r}}^2}}$. As we shall see below  one can define an action of the universal R-matrix $\mathcal{R}_s$ on tensor products of finite rank $U_\mathcal{B}^{s, res}(\g)$--modules. This action will play a crucial role in subsequent considerations.

The algebras $U_\mathcal{A}(\g)$, $U_\mathcal{B}^{res}(\g)$, $U_\mathcal{A}^s(\g)$ and $U_\mathcal{B}^{s, res}(\g)$ can be also regarded as subalgebras of $U_h^{s}(\g)\simeq U_h(\g)$, and $U_q^{res}(\g)$ can be regarded as a $\mathbb{C}[q,q^{-1}]$--subalgebra in $U_\mathcal{B}^{s, res}(\g)\simeq U_\mathcal{B}^{res}(\g)$.

Denote by $U_q^s({\frak n}_+)$, $U_q^s({\frak n}_-)$ \index[not]{U@$U_q^s(\n_\pm)$} and $U_q^s({\frak h})$ \index[not]{U@$U_q^s(\h)$} the subalgebras of $U_q^s({\frak g})$ generated by the $e_i$, $f_i$ and by the $t_i, L_i$, respectively, and let $U_q^s({\frak b}_\pm)$ \index[not]{U@$U_q^s(\b_\pm)$} be  the subalgebra in $U_q^s({\frak g})$ generated by $U_q^s({\frak n}_\pm)$ and by $U_q^s({\frak h})$, $U_q^s({\frak b}_\pm)=U_q^s({\frak n}_\pm)U_q^s({\frak h})$.

Let $U_\mathcal{A}({\frak n}_+)$, $U_\mathcal{A}({\frak n}_-)$ \index[not]{U@$U_\mathcal{A}(\n_\pm)$} (resp. $U_\mathcal{B}^{res}({\frak n}_+)$, $U_\mathcal{B}^{res}({\frak n}_-)$) \index[not]{U@$U_\mathcal{B}^{res}(\n_\pm)$} be the subalgebras of $U_\mathcal{A}({\frak g})$ (resp. of $U_\mathcal{B}^{res}(\g)$) generated by the $X_i^+$ and by the $X_i^-$, $i=1,\ldots,l$ (resp. by the $(X_i^+)^{(r)}$ and by the $(X_i^-)^{(r)}$, $i=1,\ldots,l$, $r\geq 0$), respectively.

Let $U_\mathcal{A}^{s}({\frak n}_+)$, $U_\mathcal{A}^{s}({\frak n}_-)$ \index[not]{U@$U_\mathcal{A}^s(\n_\pm)$} (resp. $U_\mathcal{B}^{s, res}({\frak n}_+)$, $U_\mathcal{B}^{s, res}({\frak n}_-)$) \index[not]{U@$U_\mathcal{B}^{s, res}(\n_\pm)$} be the subalgebras of $U_\mathcal{A}^{s}({\frak g})$ (resp. of $U_\mathcal{B}^{s, res}(\g)$) generated by the $e_i$ and by the $f_i$, $i=1,\ldots,l$ (resp. by the $e_i^{(r)}$ and by the $f_i^{(r)}$, $i=1,\ldots,l$, $r\geq 0$), respectively, and $U_\mathcal{A}(\h)\simeq U_\mathcal{A}^{s}(\h)$ \index[not]{U@$U_\mathcal{A}(\h)$} \index[not]{U@$U_\mathcal{A}^s(\h)$} the subalgebra in $U_\mathcal{A}({\frak g})\simeq U_\mathcal{A}^{s}({\frak g})$ generated by $t_i,L_i$, $i=1,\ldots ,l$.

We denote by $U_\varepsilon^{s}({\frak n}_+)$, $U_\varepsilon^{s}({\frak n}_-)$, \index[not]{U@$U_\varepsilon^s(\n_\pm)$} $U_\varepsilon^{s, res}({\frak n}_+)$, $U_\varepsilon^{s, res}({\frak n}_-)$,  \index[not]{U@$U_\varepsilon^{s, res}(\n_\pm)$} $U_\varepsilon(\h)\simeq U_\varepsilon^{s}(\h)$,  \index[not]{U@$U_\varepsilon(\h)$} \index[not]{U@$U_\varepsilon^s(\h)$} $U_\varepsilon({\frak n}_+)$, $U_\varepsilon({\frak n}_-)$, \index[not]{U@$U_\varepsilon(\n_\pm)$} $U_\varepsilon^{res}({\frak n}_+)$, $U_\varepsilon^{res}({\frak n}_-)$ \index[not]{U@$U_\varepsilon^{res}(\n_\pm)$} their specializations at $q^{\frac{1}{d{\bar{r}}^2}}=\varepsilon^{\frac{1}{d{\bar{r}}^2}}$ whenever they are defined.

The elements
$$
\left[ \begin{array}{l}
K_i;c \\
r
\end{array} \right]_{q_i}=\prod_{s=1}^r \frac{K_i q_i^{c+1-s}-K_i^{-1}q_i^{s-1-c}}{q_i^s-q_i^{-s}}~,~i=1,\ldots,l,~c\in \mathbb{Z},~r\in \mathbb{N}
$$
belong to $U_\mathcal{B}^{s, res}(\g)$. Denote by $U_\mathcal{B}^{s, res}(\h)$ \index[not]{U@$U_\mathcal{B}^{s, res}(\h)$} the subalgebra of $U_\mathcal{B}^{s, res}(\g)$ generated by those elements and by $t_i^{\pm 1}, L_i^{\pm 1}$, $i=1,\ldots,l$ and by $U_\varepsilon^{s, res}(\h)$ \index[not]{U@$U_\varepsilon^{s, res}(\h)$} its specialization at $q^{\frac{1}{d{\bar{r}}^2}}=\varepsilon^{\frac{1}{d{\bar{r}}^2}}$.

Let $U_\mathcal{B}^{s, res}({\frak b}_\pm)$ \index[not]{U@$U_\mathcal{B}^{s, res}(\b_\pm)$} be the subalgebra in $U_\mathcal{B}^{s, res}(\g)$ generated by $U_\mathcal{B}^{s, res}({\frak n}_\pm)$ and by $U_\mathcal{B}^{s, res}(\h)$, $U_\mathcal{B}^{s, res}({\frak b}_\pm)=U_\mathcal{B}^{s, res}({\frak n}_\pm)U_\mathcal{B}^{s, res}(\h)$.
As above we denote by $U_\varepsilon^{s, res}({\frak b}_\pm)$ \index[not]{U@$U_\varepsilon^{s, res}(\b_\pm)$} its specialization at $q^{\frac{1}{d{\bar{r}}^2}}=\varepsilon^{\frac{1}{d{\bar{r}}^2}}$.

The algebra anti-involution $\omega$ \index[not]{o@$\omega$} of $U_h({\frak g})$ defined by (\ref{omega}) gives rise to an algebra anti-involution of $U_h^{s}({\frak g})\simeq U_h({\frak g})$. We denote it by the same letter. 
By (\ref{omegaX}) for any $\alpha\in \Delta_+$ it satisfies
\begin{equation}\label{omega1}
\omega(f_\alpha)=e_\alpha, \omega(e_\alpha)=f_\alpha.
\end{equation}

The algebra anti-involution $\omega_0'$ \index[not]{o@$\omega_0'$} of  $U_h({\frak g})$ also induces an algebra anti-involution of $U_h^{s}({\frak g})\simeq U_h({\frak g})$ which we denote by the same letter. From (\ref{omega0'}) we deduce that it satisfies
$$
\omega_0'(e_i)=e_i, \omega_0'(f_i)=f_i.
$$
Formula (\ref{Xo0}) implies
\begin{equation}\label{eo0}
\omega_0'(e_\alpha)=c_\alpha e_\alpha,	
\end{equation}
where $c_{\alpha}=\epsilon_\alpha p_\alpha$,  $\epsilon_\alpha= \pm 1$, $p_\alpha\in q^{\mathbb{Z}}$. 
We also have
\begin{equation}\label{fo0}
\omega_0'(f_\alpha)=\omega_0'\omega(e_\alpha)=\omega\omega_0'(e_\alpha)=\omega(c_\alpha e_\alpha)=c_\alpha^{-1}\omega(e_\alpha)=c_\alpha^{-1}f_\alpha	
\end{equation}

One can check straightforwardly that $\omega_0'$ is a coalgebra homomorphism of $U_h^{s}({\frak g})$.

The anti-involution $\omega_0$ (resp. the involution $\tau$) give rise to an anti-involution (resp. to an involution) of $U_q^{s}({\frak g})$, $U_{\mathcal{B}}^{s,res}(\g)$ and $U_{\mathcal{A}}^{s}(\g)$ which we denote by the same letters.

Using the root vectors $e_{\beta}$ and $f_{\beta}$ we can construct bases for the algebras introduced above.
\begin{lemma}\label{segmPBWs}
Fix a normal ordering of the system of positive roots $\Delta_+$ and let $e_\beta, f_\beta$ be the corresponding quantum root vectors defined in Proposition \ref{rootss}. Then the following statements are true.

(i) The elements $f_\beta$ satisfy the following commutation relations
\begin{equation}\label{cmrelf}
{f}_{\alpha}{f}_{\beta} - q^{\left\langle \alpha,\beta\right\rangle+\kappa \left\langle {1+s \over 1-s}P_{\h'^*}\alpha,\beta\right\rangle}{f}_{\beta}{f}_{\alpha}= \sum_{p_1,\ldots,p_k\in \mathbb{N}}C(p_1,\ldots,p_k)
{f}_{\zeta_1}^{(p_1)}{f}_{\zeta_2}^{(p_2)}\ldots {f}_{\zeta_k}^{(p_k)}=
\end{equation}
$$
=\sum_{p_1,\ldots,p_k\in \mathbb{N}}C'(p_1,\ldots,p_k)
{f}_{\zeta_1}^{p_1}{f}_{\zeta_2}^{p_2}\ldots {f}_{\zeta_k}^{p_k},~\alpha<\beta,
$$
where $\alpha<\zeta_1<\ldots<\zeta_k<\beta$, $[\alpha,\beta]=\{\alpha,\zeta_1,\ldots,\zeta_k,\beta\}$ as a set, $C(p_1,\ldots,p_k)\in \mathcal{B}$, $C'(p_1,\ldots,p_k)\in \mathcal{A}$, and only finitely many of these coefficients are non--zero.

(ii) The elements $e^{\bf r}:=e_{\beta_1}^{r_1}\ldots e_{\beta_D}^{r_D}$, \index[not]{e@$e^{\bf r}$} $f^{\bf r}:=f_{\beta_D}^{r_D}\ldots f_{\beta_1}^{r_1}$, \index[not]{f@$f^{\bf r}$} for ${\bf r}:=(r_1,\ldots, r_D)\in {\Bbb N}^D$, form bases of $U_q^s({\frak n}_+),U_q^s({\frak n}_-)$, respectively.

(iii) The multiplication defines an isomorphism of $\mathbb{C}(q^{\frac{1}{d{\bar{r}}^2}})$--modules:
$$
U_q^{s}({\frak n}_-)\otimes U_q^{s}({\frak h}) \otimes U_q^{s}({\frak n}_+)\rightarrow U_q^{s}({\frak g}).
$$

(iv) The elements $e^{\bf r}$, $f^{\bf r}$ (resp. $e^{\bf (r)}:=e_{\beta_1}^{(r_1)}\ldots e_{\beta_D}^{(r_D)}$, \index[not]{e@$e^{\bf (r)}$} $f^{\bf (r)}:=f_{\beta_D}^{(r_D)}\ldots f_{\beta_1}^{(r_1)}$) \index[not]{f@$f^{\bf (r)}$} for ${\bf r}\in {\Bbb N}^D$
form bases of $U_\mathcal{A}^{s}({\frak n}_+), U_\mathcal{A}^{s}({\frak n}_-)$ (resp. $U_\mathcal{B}^{s, res}({\frak n}_+),U_\mathcal{B}^{s, res}({\frak n}_-)$),
respectively.

(v) The multiplication defines an isomorphisms of $\mathcal{B}$--modules:
$$
U_\mathcal{B}^{s, res}({\frak n}_-)\otimes U_\mathcal{B}^{s, res}({\frak h}) \otimes U_\mathcal{B}^{s, res}({\frak n}_+)\rightarrow U_\mathcal{B}^{s, res}({\frak g}).
$$

(vi) Let $[\alpha,\beta]=\{\beta_p,\ldots , \beta_q\}$ be a minimal segment in $\Delta_+$,
$U_\mathcal{A}^s([\alpha,\beta])$, \index[not]{U@$U_\mathcal{A}^s([\pm\alpha,\pm\beta])$} $U_\mathcal{A}^s([-\alpha,-\beta])$ (resp. $U_{\mathcal{B}}^{s,res}([\alpha,\beta])$, \index[not]{U@$U_{\mathcal{B}}^{s,res}([\pm\alpha,\pm\beta])$} $U_{\mathcal{B}}^{s, res}([-\alpha,-\beta])$) the $\mathcal{A}$(resp. $\mathcal{B}$)--subalgebras of $U_h^s({\frak g})$ generated by the $e_\gamma$ and by the $f_\gamma$, $\gamma \in [\alpha, \beta]$ (resp. by the $(e_\gamma)^{(r)}$ and by the $(f_\gamma)^{(r)}$, $\gamma \in [\alpha, \beta]$, $r\in \mathbb{N}$), respectively. Then the elements $(e_{\beta_p})^{r_p}\ldots (e_{\beta_q})^{r_q}$, $(f_{\beta_q})^{r_q}\ldots (f_{\beta_p})^{r_p}$ (resp. $(e_{\beta_p})^{(r_p)}\ldots (e_{\beta_q})^{(r_q)}$, $(f_{\beta_q})^{(r_q)}\ldots (f_{\beta_p})^{(r_p)}$), $r_i\in {\Bbb N}$ form bases of $U_\mathcal{A}^s([\alpha,\beta])$, $U_\mathcal{A}^s([-\alpha,-\beta])$ (resp. $U_{\mathcal{B}}^{s, res}([\alpha,\beta])$, $U_{\mathcal{B}}^{s, res}([-\alpha,-\beta])$),
respectively.

(vii) The elements $f^{\bf r}$ (resp. $f^{\bf (r)}$) for ${\bf r}\in {\Bbb N}^D$ with $r_i>0$ for at least one $i\geq p$ form a basis in the right ideal $\bar{Y}_p$ (resp. $\bar{Y}_p^{res}$) \index[not]{Y@$\bar{Y}_p$} \index[not]{Y@$\bar{Y}_p^{res}$} of $U_\mathcal{A}^{s}({\frak n}_-)$ (resp. $U_\mathcal{B}^{s, res}({\frak n}_-)$) generated by $f_\gamma$, $\gamma \in [\beta_p, \beta_D]$ (resp. by $(f_\gamma)^{(r)}$, $\gamma \in [\beta_p, \beta_D]$, $r>0$).

(viii) For any $n> 0$, $p=1,\ldots, D-1$ one has
$$
f_{\beta_p}^{(n)}U_{\mathcal{B}}^{s, res}([-\beta_{p+1},-\beta_D])\subset \sum_{i=0}^{n-1}(U_{\mathcal{B}}^{s, res}([-\beta_{p+1},-\beta_D]))_0f_{\beta_p}^{(i)}+U_{\mathcal{B}}^{s, res}([-\beta_{p+1},-\beta_D])f_{\beta_p}^{(n)},
$$
where $(U_{\mathcal{B}}^{s, res}([-\beta_{p+1},-\beta_D]))_0=\bar{Y}_{p+1}^{res}\cap U_{\mathcal{B}}^{s, res}([-\beta_{p+1},-\beta_D])$. \index[not]{U@$(U_{\mathcal{B}}^{s, res}([-\beta_{p+1},-\beta_D]))_0$}

(ix) Let $[\alpha,\beta]\subset \Delta_+$ or $[\alpha,\beta]\subset \Delta_-$ be any minimal segment,
such that $[\alpha, \beta]=[\alpha,\gamma]\cup [\delta,\beta]$ (disjoint union of minimal segments). Then the multiplication in $U_{\mathcal{B}}^{s, res}(\g)$ defines isomorphisms of $\mathcal{B}$--modules
$$
U_{\mathcal{B}}^{s, res}([\alpha,\gamma])\otimes U_{\mathcal{B}}^{s, res}([\delta,\beta])\to U_{\mathcal{B}}^{s, res}([\alpha,\gamma])U_{\mathcal{B}}^{s, res}([\delta,\beta])=U_{\mathcal{B}}^{s, res}([\alpha,\beta]),
$$
$$
U_{\mathcal{B}}^{s, res}([\delta,\beta])\otimes U_{\mathcal{B}}^{s, res}([\alpha,\gamma])\to U_{\mathcal{B}}^{s, res}([\delta,\beta])U_{\mathcal{B}}^{s, res}([\alpha,\gamma])=U_{\mathcal{B}}^{s, res}([\alpha,\beta]).
$$
\end{lemma}

\begin{proof}
Commutation relations (\ref{cmrelf}) follow from commutation relations (\ref{qcom}), (\ref{weight-root}), (\ref{comml}), Proposition \ref{rootprop}, the definition of the elements $e_\beta, f_\beta$ and the definition of the isomorphism $\psi_{\{ n_{ij}\}}$. This proves (i). 

Statements (ii)-(vi) of this lemma follow straightforwardly from parts (i), (iii) and (iv) of Lemma \ref{segmPBW} and Propositions \ref{newreal} and \ref{rootss}.

For (vii), using commutation relations (\ref{cmrelf}) we can represent any element of the right ideal of $U_\mathcal{A}^{s}({\frak n}_-)$ generated by $f_\gamma$, $\gamma \in [\beta_p, \beta_D]$  as an $\mathcal{A}$--linear combination of the elements $f^{\bf r}$ for ${\bf r}\in {\Bbb N}^D$ with $r_i>0$ for at least one $i\geq p$. This presentation is unique by the Poincar\'{e}--Birkhoff--Witt decomposition for $U_\mathcal{A}^{s}({\frak n}_-)$ stated in (vi).

Note that a similar result holds for the algebra $U_q^{s}({\frak n}_-)=U_\mathcal{A}^{s}({\frak n}_-)\otimes_{\mathcal{A}}\mathbb{C}(q^{\frac{1}{d{\bar{r}}^2}})$ for the same reasons. 

We can apply it to represent any element of the right ideal $\bar{Y}_p^{res}$ of $U_\mathcal{B}^{s, res}({\frak n}_-)\subset U_q^{s}({\frak n}_-)$ generated by $(f_\gamma)^{(r)}$, $\gamma \in [\beta_p, \beta_D]$, $r>0$ as a $\mathbb{C}(q^{\frac{1}{d{\bar{r}}^2}})$--linear combination of the elements $f^{\bf (r)}$ for ${\bf r}\in {\Bbb N}^D$ with $r_i>0$ for at least one $i\geq p$. This presentation is unique and by the uniqueness of the Poincar\'{e}--Birkhoff--Witt decomposition for $U_\mathcal{B}^{s, res}({\frak n}_-)$ stated in (vi) the coefficients in this decomposition belong to $\mathcal{B}$. This completes the proof of part (vii). 

(viii) is justified using similar arguments, commutation relations (\ref{cmrelf}) and induction over $n$. 

(ix) follows from Corollary \ref{segmq} (iii) and Propositions \ref{newreal} and \ref{rootss}.

\end{proof}

A basis for $U_\mathcal{B}^{s, res}(\h)$ is a little bit more difficult to describe. We do not need its explicit description.

\begin{remark}\label{segmPBWsrev}
Applying the anti-involution $\omega_0$ to the elements of the bases constructed in Lemma \ref{segmPBWs} and using (\ref{eo0}) and (\ref{fo0}) we obtain other bases of similar types where the order of the quantum root vectors in the products defining the elements of the bases is reversed.

By specializing the above constructed bases at $q^{\frac{1}{rd^2}}=\varepsilon^{\frac{1}{rd^2}}$ one can obtain similar bases and similar subalgebras for $U_\varepsilon^s(\g)$ and $U_\varepsilon^{s, res}(\g)$.

Note that only the anti-involution $\omega_0$ gives rise to an anti-involution of $U_\varepsilon^s(\g)$ and $U_\varepsilon^{s, res}(\g)$ for arbitrary $\varepsilon^{\frac{1}{rd^2}}\neq 0$ as $\omega_0(h)=h$ and $\omega(h)=\omega_0'(h)=\tau(h)=-h$.
\end{remark}

Similarly to the case of $U_q^{res}(\g)$, using parts (ii), (iii), (iv) and (v) of Lemma \ref{segmPBWs} one can introduce $Q$--gradings on the algebras $U_q^s(\g)$, $U_\mathcal{A}^{s}(\g)$, $U_\mathcal{B}^{s, res}(\g)$, $U_\varepsilon^{s,res}(\g)$ and $U_\varepsilon^s(\g)$. For instance, we have 
$$
U_\mathcal{B}^{s, res}(\g)=\bigoplus_{\mu\in Q} (U_\mathcal{B}^{s, res}(\g))_\mu,
$$
where
$$
(U_\mathcal{B}^{s, res}(\g))_\mu=\{x\in U_\mathcal{B}^{s, res}(\g):[h,x]=\mu(h)x~\text{for all}~h\in \h\subset U_h(\g)\simeq U_h^s(\g)\}, \index[not]{U@$(U_\mathcal{B}^{s, res}(\g))_\mu$}
$$ 
and
$$
U_\mathcal{A}^{s}(\g)=\bigoplus_{\mu\in Q} (U_\mathcal{A}^{s}(\g))_\mu,
$$
where
$$
(U_\mathcal{A}^{s}(\g))_\mu=\{x\in U_\mathcal{A}^{s}(\g):[h,x]=\mu(h)x~\text{for all}~h\in \h\subset U_h(\g)\simeq U_h^s(\g)\}. \index[not]{U@$(U_\mathcal{A}^{s}(\g))_\mu$}
$$

The algebra $U_h^{s}(\g)$ is, of course, not $Q$--graded. But one can define its weight subspaces $(U_h^{s}(\g))_\mu$, $\mu\in Q$ in a similar way,
$$
(U_h^{s}(\g))_\mu=\{x\in U_h^{s}(\g):[h,x]=\mu(h)x~\text{for all}~h\in \h\subset U_h(\g)\simeq U_h^s(\g)\}. \index[not]{U@$(U_h^s(\g))_\mu$}
$$ 

The anti-involutions $\omega_0'$, $\omega$ and the involution $\tau$ preserve weights of elements while the anti-involution $\omega$ inverses the signs of weights.

Note that isomorphisms $\psi_{\{n_{ij}\}}$ also preserve weights of elements.

Using formulas (\ref{comult}) and (\ref{defds}) one can also find that
\begin{equation}\label{comults}
\Delta_s(f_{\beta_k})=\theta_{w_{k-1}}^s(e^{-h\kappa {1+s \over 1-s}P_{\h'}\beta_k^\vee+h\beta_k^\vee} \otimes f_{\beta_k} + f_{\beta_k}\otimes 1)(\theta_{w_{k-1}}^s)^{-1}=
\end{equation}
$$
=G_{\beta_k}^{-1} \otimes f_{\beta_k} + f_{\beta_k}\otimes 1+\sum_i y_i\otimes x_i,
$$
where
$$
G_\beta= e^{h\kappa {1+s \over 1-s}P_{\h'}\beta^\vee-h\beta^\vee}, y_i=e^{-h\kappa {1+s \over 1-s}P_{\h'}\gamma_{x_i}^\vee+h\gamma_{x_i}^\vee}\overline{y}_i, \index[not]{G@$G_\beta$} \index[not]{x@$x_i$} \index[not]{y@$y_i$}
$$
$$
\overline{y}_i\in U_\mathcal{A}^s([-\beta_{k+1},-\beta_{D}])\cap U_{\mathcal{B}}^{s, res}([-\beta_{k+1},-\beta_D]),
$$
$$x_i\in U_\mathcal{A}^s([-\beta_{1},-\beta_{k-1}])\cap U_{\mathcal{B}}^{s, res}([-\beta_{1},-\beta_{k-1}]),$$ $\overline{y}_i,x_i$ belong to weight subspaces and have non-zero weights, $\gamma_{x_i}$ is the weight of $x_i$ (see \cite{Dr}, Corollary 4.3.2), for $w_{k-1}=s_{i_1}\ldots s_{i_{k-1}}$
\begin{equation}\label{Rsinv1}
\theta_{w_{k-1}}^s=\theta^s_{\beta_{1}}\ldots \theta^s_{\beta_{k-1}},~\theta^s_{\beta_r}={\exp}_{q_{\beta_r}}[(1-q_{\beta_r}^{-2})e_{\beta_r}e^{-h\kappa {1+s \over 1-s}P_{\h'}\beta^\vee}\otimes f_{\beta_r}],
\end{equation}
and
\begin{equation}\label{Rsinv2}
(\theta_{w_{k-1}}^s)^{-1}=(\theta_{\beta_{k-1}}^s)^{-1}\ldots (\theta^s_{\beta_{1}})^{-1},~(\theta^s_{\beta_r})^{-1}={\exp}_{q_{\beta_r}^{-1}}[(1-q_{\beta_r}^{2})e_{\beta_r}e^{-h\kappa {1+s \over 1-s}P_{\h'}\beta^\vee}\otimes f_{\beta_r}].
\end{equation}

From (\ref{comults}) we also obtain
\begin{equation}\label{comultsn}
\Delta_s(f_{\beta_k}^{(n)})=\frac{1}{[n]_{q_{\beta_k}}!}\theta_{w_{k-1}}^s(G_{\beta_k}^{-1} \otimes f_{\beta_k} + f_{\beta_k}\otimes 1)^n(\theta_{w_{k-1}}^s)^{-1}=
\end{equation}
$$
=\theta_{w_{k-1}}^s(\sum_{k=0}^nq_{\beta_k}^{k(n-k)}G_{\beta_k}^{-k}f_{\beta_k}^{(n-k)}\otimes f_{\beta_k}^{(k)})(\theta_{w_{k-1}}^s)^{-1}=
$$
$$
=\sum_{k=0}^nq_{\beta_k}^{k(n-k)}G_{\beta_k}^{-k}f_{\beta_k}^{(n-k)}\otimes f_{\beta_k}^{(k)}+\sum_i y_i^{(n)}\otimes x_i^{(n)},
$$
where
$$
y_i^{(n)}=e^{-h\kappa {1+s \over 1-s}P_{\h'}\gamma_{x_i^{(n)}}^\vee+h\gamma_{x_i^{(n)}}^\vee}\overline{y}_i^{(n)},
$$
$\overline{y}_i^{(n)}\in I_k^>$, \index[not]{y@$\overline{y}_i^{(n)}$} \index[not]{I@$I_k^>$} $x_i^{(n)}\in I_k^<$ \index[not]{x@$x_i^{(n)}$} \index[not]{I@$I_k^<$} belong to weight subspaces and have non-zero weights, $\gamma_{x_i^{(n)}}$ is the weight of $x_i^{(n)}$, $I_k^>$ is the ideal in $U_{\mathcal{B}}^{s, res}([-\beta_{k},-\beta_D])$ generated by $f_{\beta_i}^{(p)}$, $i=k+1,\ldots ,D$, $p>0$, and $I_k^<$ is the ideal in $U_{\mathcal{B}}^{s, res}([-\beta_{1},-\beta_k])$ generated by $f_{\beta_i}^{(p)}$, $i=1,\ldots, k-1$, $p>0$.

Similarly
\begin{equation}\label{comultse}
\Delta_s(e_{\beta_k})=\theta_{w_{k-1}}^s(e_{\beta_k}\otimes e^{-h\beta_k^\vee}+e^{h\kappa {1+s \over 1-s}P_{\h'}\beta_k^\vee} \otimes e_{\beta_k})(\theta_{w_{k-1}}^s)^{-1}=
\end{equation}
$$
=e_{\beta_k}\otimes e^{-h\beta_k^\vee}+e^{h\kappa {1+s \over 1-s}P_{\h'}\beta_k^\vee} \otimes e_{\beta_k}+\sum_i x_i'\otimes y_i',
$$
where 
$$
y_i'\in U_\mathcal{A}^s([\beta_{k+1},\beta_{D}])U_\mathcal{A}^s(\h)\cap U_{\mathcal{B}}^{s, res}([\beta_{k+1},\beta_D])U_\mathcal{B}^{s, res}(\h), \index[not]{y@$y_i'$}
$$
$$
x_i'\in U_\mathcal{A}^s([\beta_{1},\beta_{k-1}])U_\mathcal{A}^s(\h)\cap U_{\mathcal{B}}^{s, res}([\beta_{1},\beta_{k-1}])U_\mathcal{B}^{s, res}(\h), \index[not]{x@$x_i'$}
$$
and ${y}_i',x_i'$ belong to weight subspaces and have non-zero weights.

From these formulas we deduce
\begin{equation}\label{Ss}
S_s(f_\beta)=-G_\beta f_\beta-\sum_i S_s(y_i)x_i=-G_\beta (f_\beta +\sum_i y_iS_s(x_i)),
\end{equation}
\begin{equation}\label{Sinv}
S_s^{-1}(f_\beta)=-(f_\beta +\sum_i S_s^{-1}(x_i)y_i)G_\beta=-f_\beta G_\beta-\sum_i x_iS_s^{-1}(y_i).
\end{equation}

We also have
$$
\omega_0S_s^{-1}(f_\beta)=-G_\beta^{-1}(\omega_0 f_\beta +\sum_i \omega_0(y_i)\omega_0 S_s^{-1}(x_i))=
$$
$$
=-G_\beta^{-1} \omega_0 f_\beta -\sum_i \omega_0S_s^{-1}(y_i)\omega_0 (x_i),
$$
$$
\omega_0 f_\beta=-G_\beta(\omega_0S_s^{-1}(f_\beta)+\sum_i \omega_0S_s^{-1}(y_i)\omega_0 (x_i)).
$$

As usual, one can define highest weight, Verma and finite-dimensional modules for all forms and specializations of the quantum group $U_h^{s}({\frak g})$ introduced above.
We recall that by Propositions 6.5.5 and 6.5.7 in \cite{ChP} $U_h(\g)\simeq U(\g)[[h]]$, and this isomorphism of algebras restricts to the identity map on $U(\h)$ and induces a canonical isomorphism of the center of $U_h(\g)$ and of $Z(U(\g))[[h]]$, where $Z(U(\g))$ is the center of $U(\g)$.
Therefore if $V$ is a $U_h(\g)$--module free and of finite rank \index{representation!of a quantum group!finite rank} over $\mathbb{C}[[h]]$ then $V_1:=V/hV$ \index[not]{V@$V_1$} is a finite-dimensional $U(\g)$--module which is completely reducible, and its irreducible components are highest weight irreducible finite-dimensional representations of $U(\g)$. By Corollary 6.5.6 in \cite{ChP}, one has an isomorphism of $U_h(\g)$--modules, $V\simeq V_1[[h]]$, where the action of $U_h(\g)$ on $V_1[[h]]$ is defined using the algebra isomorphism $U_h(\g)\simeq U(\g)[[h]]$.

Using the algebra isomorphism $U_h(\g)\simeq U(\g)[[h]]$, which restricts to the identity map on $U(\h)$, one can define weight vectors   in $V$ by requiring that $v\in V$ has weight $\lambda\in P$ if $hv=\lambda(h)v$ for any $h\in \h\subset U(\h)\subset U(\g)[[h]]\simeq U_h(\g)$. For $\lambda\in P$ one also defines the corresponding weight subspace $(V)_\lambda=\{v\in V:hv=\lambda(h)v \text{ for all } h\in \h\}$. \index{representation!of a quantum group!weight subspace} \index[not]{V@$(V)_\lambda$}

If $V_1$ is the highest weight irreducible representation of $U(\g)$ with highest weight $\lambda\in P_+$ we call the corresponding representation $V=V_\lambda$ \index[not]{V@$V_\lambda$}  the highest weight indecomposable representation of highest weight $\lambda$. \index{representation!of a quantum group!highest weight indecomposable} $V_\lambda$ is generated by a highest weight vector \index{representation!of a quantum group!highest weight vector} with respect to the action of $\h\subset U_h(\h)$. All indecomposable $U_h(\g)$--modules free and of  finite rank over $\mathbb{C}[[h]]$ can be obtained this way. 

Let $V$ be a topologically free finite rank $U_h(\g)$--module. Recall that there is a contravariant non--degenerate form $(~\cdot~,~\cdot~)$ \index[not]{ZZ@$(~\cdot~ ,~\cdot~ )$} \index{form!contravariant non--degenerate!on a representation of a quantum group} on $V$ such that $(u,xv)=(\omega(x)u,v)$ for any $u,v\in V$, $x\in U_h(\g)$. Different weight subspaces are orthogonal with respect to this form. If $V=V_\lambda$, $\lambda\in P_+$ is indecomposable and $v_\lambda$ is a highest weight vector generating it then the condition $(v_\lambda,v_\lambda)=1$ uniquely determines the contravariant non--degenerate form. We shall always consider the contravariant non--degenerate form on $V_\lambda$ fixed by this condition.

More generally, recall that for arbitrary topologically free finite rank $U_h(\g)$--module $V$ one has $V\simeq V_1[[h]]$, where $V_1$ is a finite-dimensional $U(\g)$--module, and the action of $U_h(\g)$ on $V_1[[h]]$ is defined using the algebra isomorphism $U_h(\g)\simeq U(\g)[[h]]$. Since the $U(\g)$--module $V_1$ is completely reducible, $V$ is isomorphic to a direct sum of highest weight indecomposable representations, and hence the sum of the contravariant non--degenerate forms on the indecomposable components $V_\lambda$ of $V$ normalized by the condition $(v_\lambda,v_\lambda)=1$ uniquely determines a contravariant non--degenerate form on $V$. We shall always consider the contravariant non--degenerate form on $V$ fixed by this condition.  

The definition (\ref{omega}) of the anti-involution $\omega$ implies that the $\omega$ maps elements of negative weights to elements of positive weights. From this remark we obtain the following obvious corollary which will be often used later.
\begin{lemma}\label{hwv}
Let $V$ be a topologically free finite rank $U_h(\g)$--module. Let $u, v\in V$, $x, y\in U_h(\g)$.

(i) If $u$ is a highest weight vector and $x\in U_h(\g)$ has weight which is not non--negative then $(u,xyv)=0$.

(ii) If $v$ is a highest weight vector and $y\in U_h(\g)$ has weight which is not non--positive then $(u,xyv)=0$.
\end{lemma}

Recall that $U_{\mathcal{B}}^{s,res}(\g)\simeq U_{\mathcal{B}}^{res}(\g)$ can be regarded as a subalgebra of $U_h^{s}({\frak g})\simeq U_h(\g)$. Let $V$ be a $U_h(\g)$--module topologically free and of  finite rank over $\mathbb{C}[[h]]$ (for brevity we shall call such modules finite rank $U_h(\g)$--modules). Then a $U_{\mathcal{B}}^{s,res}(\g)$--module $V^{res}$ \index[not]{V@$V^{res}$} is called a {\it $U_{\mathcal{B}}^{s,res}(\g)$--lattice} in $V$ if $V^{res}\otimes_{\mathcal{B}}\mathbb{C}[[h]]\simeq V$. \index{lattice!in a representation of a quantum group} 

For any $U_h(\g)$--module $V$ topologically free and of finite rank over $\mathbb{C}[[h]]$ a $U_{\mathcal{B}}^{s,res}(\g)$--lattice $V^{res}$ exists. Indeed, as we observed above $V$ is a direct sum of indecomposable highest weight representations. Therefore it suffices to consider the case of indecomposable $V=V_\lambda$ generated by a highest weight vector $v$. In this case, similarly to the proof of Proposition 4.2 in \cite{L3} (see also Proposition 10.1.4 in \cite{ChP}) one can show that $V^{res}_\lambda:=U_{\mathcal{B}}^{s,res}(\g)v$ \index[not]{V@$V^{res}_\lambda$} is a $U_{\mathcal{B}}^{s,res}(\g)$--lattice, $V^{res}_\lambda$  is the direct sum of its intersections with the weight spaces of $V_\lambda$ and each such intersection is a finitely generated free $\mathcal{B}$--module of finite rank. Moreover, using the arguments from the proof of Proposition 4.2 in Section 4.9 of \cite{L3} one can see that the last two properties hold for any $U_{\mathcal{B}}^{s,res}(\g)$--lattice in any topologically free $U_h(\g)$--module $V$ of finite rank over $\mathbb{C}[[h]]$. Together with the results of Section 10.1 in \cite{ChP} this yields the following statement.
\begin{proposition}\label{wresdec}
(i) For any $U_h(\g)$--module $V$ topologically free and of finite rank over $\mathbb{C}[[h]]$ $V^{res}$ is the direct sum of its intersections with the weight spaces of $V$ and each such intersection is a finitely generated free $\mathcal{B}$--module of finite rank
$$
V^{res}=\bigoplus_{\lambda\in P}(V^{res})_\lambda, (V^{res})_\lambda:=V^{res}\cap (V)_\lambda. \index[not]{V@$(V^{res})_\lambda$}
$$

(ii) The specialization of $V^{res}$ at $q^{\frac{1}{d{\bar{r}}^2}}=1$ is naturally a finite--dimensional $U(\g)$--module. \index{specialization!of a module over a quantum group}
\end{proposition}

For any $U_h(\g)$--module $V$ topologically free and of finite rank over $\mathbb{C}[[h]]$, one can equip the module $V^{res}$ with a natural action of $\h$, where $h\in \h$ acts on the weight subspace of $V^{res}$ of weight $\lambda$ by multiplication by $\lambda(h)$.

The contravariant non--degenerate form $(~\cdot~,~\cdot~)$ on $V$ restricts to a contravariant non--degenerate form $(~\cdot~,~\cdot~)$ on the $U_{\mathcal{B}}^{s,res}(\g)$--submodule $V^{res}\subset V$ such that $(u,xv)=(\omega(x)u,v)$ for any $u,v\in V^{res}$, $x\in U_{\mathcal{B}}^{s,res}(\g)$, and different weight subspaces $(V^{res})_\lambda$, $\lambda \in P$ are orthogonal with respect to this form.

The form $(~\cdot~,~\cdot~)$ on $V^{res}$ gives rise to a contravariant non--degenerate form $(~\cdot~,~\cdot~)$ on the $U_q^s(\g)\simeq U_{\mathcal{B}}^{s,res}(\g)\otimes_{\mathcal{B}}\mathbb{C}(q^{\frac{1}{d{\bar{r}}^2}})$--module $V_q:=V^{res}\otimes_{\mathcal{B}}\mathbb{C}(q^{\frac{1}{d{\bar{r}}^2}})$ \index[not]{V@$V_q$} such that $(u,xv)=(\omega(x)u,v)$ for any $u,v\in V_q$, $x\in U_q^s(\g)$, and different weight subspaces $(V_q)_\lambda:=V_q\cap (V)_\lambda$, \index[not]{V@$(V_q)_\lambda$} $\lambda \in P$ are orthogonal with respect to this form.

From the explicit formulas and the results in \cite{Lus}, Section 5.2 it follows that the braid group elements given by (\ref{T1}) and (\ref{T2}), or by applying $\tau$ to (\ref{T1}) and (\ref{T2}), act in $U_h(\g)$--modules topologically free and of finite rank over $\mathbb{C}[[h]]$ and leave invariant $U_{\mathcal{B}}^{s,res}(\g)$--lattices in them.

More generally, if $T_w$ (resp. $T_w^\tau$) is the braid group element corresponding to $w\in W$ then $T_w(V)_\lambda =(V)_{w\lambda}$ (resp. $T_w^\tau(V)_\lambda =(V)_{w\lambda}$). This action satisfies the property $T_wxv=T_w(x)T_wv$ (resp. $T_w^\tau xv=T_w^\tau (x)T_w^\tau v$) for any $x\in U_h(\g)$ and $v\in V$. The $U_{\mathcal{B}}^{s,res}(\g)$--submodule $V^{res}\subset V$ is invariant under the action of the elements $T_w$ and $T_w^\tau$.

The following lemma will be useful for some calculations later.
\begin{lemma}\label{wHact}
Let $V$ be a $U_h(\g)$--module $V$ topologically free and of finite rank over $\mathbb{C}[[h]]$. Then the elements of the subalgebra $U_q^{res}(H)$ act on weight vectors of $V$ by multiplication by elements of $\mathbb{C}[q,q^{-1}]$, and the elements of the subalgebra $U_\mathcal{B}^{s, res}(\h)\subset U_\mathcal{B}^{s, res}(\g)\subset U_h^s(\g)\simeq U_h(\g)$ act on weight vectors of $V$ by multiplication by elements of $\mathcal{B}$. Moreover,
\begin{equation}\label{wcond} 
(V)_\lambda=\{v\in V: K_iv=q_i^{\lambda(H_i)}v, i=1,\ldots, l \}, (V^{res})_\lambda=\{v\in V^{res}: K_iv=q_i^{\lambda(H_i)}v, i=1,\ldots, l \}.
\end{equation}
\end{lemma}

\begin{proof}
By the definition of the action of $\h$ on $V$ the elements
$$
\left[ \begin{array}{l}
K_i;c \\
r
\end{array} \right]_{q_i}=\prod_{s=1}^r \frac{K_i q_i^{c+1-s}-K_i^{-1}q_i^{s-1-c}}{q_i^s-q_i^{-s}}~,~i=1,\ldots,l,~c\in \mathbb{Z},~r\in \mathbb{N}
$$
act on a vector of weight $\lambda$ by multiplication by $\prod_{s=1}^r \frac{q_i^{\lambda(H_i)+c+1-s}-q_i^{-\lambda(H_i)+s-1-c}}{q_i^s-q_i^{-s}}=\prod_{s=1}^r(-[\lambda(H_i)+c+1]_{q_i^{-s}})$, and $K_i^{\pm1}$, $i=1,\ldots,l$,  act by multiplication by $q_i^{\pm\lambda(H_i)}$. Also for $v\in V$ the condition $K_iv=q_i^{\lambda(H_i)}v$, $i=1,\ldots, l$ is equivalent to $hv=\lambda(h)v$ for any $h\in \h\subset U_h(\g)$. This implies (\ref{wcond}). 

The case of the generators $t_i^{\pm 1}, L_i^{\pm 1}$, $i=1,\ldots,l$ of the algebra $U_\mathcal{B}^{s, res}(\h)$ is considered in a similar way.

\end{proof}

The specialization of $V^{res}_\lambda$ at $q^{\frac{1}{d{\bar{r}}^2}}=\varepsilon^{\frac{1}{d{\bar{r}}^2}}$ \index{specialization!of a module over a quantum group} is a highest weight $U_{\varepsilon}^{s,res}(\g)$--module. In general, for a $U_h(\g)$--module $V$ topologically free and of finite rank over $\mathbb{C}[[h]]$, the specialization of $V^{res}$ at $q^{\frac{1}{d{\bar{r}}^2}}=\varepsilon^{\frac{1}{d{\bar{r}}^2}}$ is not irreducible even if $V$ is indecomposable.

The R-matrix $\mathcal{R}_s$ acts in tensor products of $U_{\mathcal{B}}^{s,res}(\g)$--lattices in $U_h(\g)$--modules topologically free and of finite rank over $\mathbb{C}[[h]]$. Namely, recalling (\ref{rmatrspi}) we can represent $\mathcal{R}_s$ in the form $\mathcal{R}_s=\mathcal{E}_s\mathcal{R}_s^0=\mathcal{R}_s^0\mathcal{E}_s^1$, where \index[not]{R@$\mathcal{R}_s^0$} \index[not]{R@$\mathcal{R}_s^0$} \index[not]{E@$\mathcal{E}_s$} \index[not]{E@$\mathcal{E}_s^1$}
\begin{eqnarray}\label{rmatrspires}
\mathcal{R}_s^0={\exp}\left[ h\left(\sum_{i=1}^l(Y_i\otimes H_i)-
\sum_{i=1}^l \kappa {1+s \over 1-s }P_{\h'}H_i\otimes Y_i\right) \right], \\
\mathcal{E}_{s}=
\prod_{\beta\in \Delta_+}\left(
\sum_{k=0}^\infty q_{\beta}^{{k(k+1)} \over 2}(1-q_{\beta}^{-2})^kf_{\beta}^{(k)} \otimes
e_{\beta}^ke^{-kh\kappa {1+s \over 1-s}P_{\h'} \beta^\vee}\right) = \nonumber \\
=\prod_{\beta\in \Delta_+}\left(
\sum_{k=0}^\infty q_{\beta}^{{k(k+1)} \over 2}(1-q_{\beta}^{-2})^kf_{\beta}^{k} \otimes
e_{\beta}^{(k)}e^{-kh\kappa {1+s \over 1-s}P_{\h'} \beta^\vee}\right), \nonumber 
\end{eqnarray}
\begin{eqnarray*}
\mathcal{E}_s^1=
\prod_{\beta\in \Delta_+}\left(
\sum_{k=0}^\infty q_{\beta}^{{k(k+1)} \over 2}(1-q_{\beta}^{-2})^ke^{kh\left(\kappa {1+s \over 1-s}P_{\h'} -id\right)\beta^\vee} f_{\beta}^{(k)} \otimes
e_{\beta}^ke^{kh\beta^\vee}\right) = \nonumber \\
=\prod_{\beta\in \Delta_+}\left(
\sum_{k=0}^\infty q_{\beta}^{{k(k+1)} \over 2}(1-q_{\beta}^{-2})^ke^{kh\left(\kappa {1+s \over 1-s}P_{\h'} -id\right)\beta^\vee} f_{\beta}^{k} \otimes
e_{\beta}^{(k)}e^{kh\beta^\vee}\right), \nonumber 
\end{eqnarray*}
and for $\beta=\sum_{i=1}^l m_i\alpha_i$ 
$$
 e^{-kh\kappa {1+s \over 1-s}P_{\h'} \beta^\vee}=e^{-kh\sum_{i,j=1}^l\frac{c_{ij}}{d_j}m_iY_j}\in U_{\mathcal{B}}^{s,res}(\h)\cap U_{\mathcal{A}}^{s}(\h).
$$

We can define the action of $\mathcal{R}_s^0$ on tensor products of modules of the form $V^{res}$, where $V$ is a highest weight finite rank $U_h({\frak g})$--module as follows. If $V^{res}$, $W^{res}$ are two such modules and $v_\lambda\in V^{res}$, $w_\mu\in W^{res}$ are vectors of weights $\lambda$ and $\mu$ then we define
$$
\mathcal{R}_s^0v_\lambda\otimes w_\mu=q^{\left\langle \lambda,\mu\right\rangle-\left\langle \kappa {1+s \over 1-s}P_{\h'^*}\mu,\lambda\right\rangle}v_\lambda\otimes w_\mu,
$$
and $q^{\left\langle \lambda,\mu\right\rangle-\left\langle \kappa {1+s \over 1-s}P_{\h'^*}\mu,\lambda\right\rangle}\in \mathcal{B}$.

Since for any $\beta\in \Delta_+$ and any module of the form $V^{res}$, where $V$ is a finite rank $U_h({\frak g})$--module, $f_{\beta}^{(k)}$ and $e_{\beta}^k$ belong to $U_{\mathcal{B}}^{s,res}(\g)$ and act as zero operators on $V^{res}$ for $k$ large enough, 
and for any $k$
$$
e^{-kh\kappa {1+s \over 1-s}P_{\h'} \beta^\vee}\in U_{\mathcal{B}}^{s,res}(\h)\cap U_{\mathcal{A}}^{s}(\h),
$$
$\mathcal{E}_s$ naturally acts on tensor products of such modules.

For two modules of the form $V^{res}, W^{res}$, where $V$ and $W$ are highest weight finite rank $U_h({\frak g})$--modules, we denote by $R^{VW}$ \index[not]{R@$R^{VW}$} the operator corresponding to the action of $\mathcal{R}_s$ in $V^{res}\otimes W^{res}$.

With this definition the identity 
$$
\Delta_s^{opp}(x)=\mathcal{R}_s\Delta_s(x){\mathcal{R}_s}^{-1}, x\in U_{\mathcal{B}}^{s,res}(\g)
$$
still holds being evaluated in tensor products of modules of the form $V^{res}$, where $V$ is a finite rank $U_h({\frak g})$--module.

We shall also need the following technical lemma regarding the action of the elements $T_w$ and $\overline{T}_w$ on finite rank indecomposable modules.

\begin{lemma}\label{bga}
Let $V$ be a finite rank $U_h(\g)$--module. If $T, \tilde{T}$ are two elements of the braid group $\mathcal{B}_\g$ which act as the same transformation on $\h\subset U_h(\h)$  and $v$ is a highest weight vector in $V$ then $Tv=t\tilde{T}v$ and $T^\tau v=t'\tilde{T}^\tau v$, where $t,t'$ are non--zero multiples of powers of $q$.
\end{lemma}

\begin{proof}
Since $v$ generates a highest weight indecomposable submodule $V_\lambda\subset V$ of highest weight $\lambda$ equal to the weight of $v$, and $V_\lambda$ is invariant under the braid group action, we can assume without loss of generality that $V=V_\lambda$.
  
First observe that similarly to the proof of Proposition 4.2 in \cite{L3} (see also Proposition 10.1.4 in \cite{ChP}) one can show that $V':=U_q^{res}(\g)v$ is a $U_q^{res}(\g)$--lattice in $V_\lambda$ in the sense that $V'\otimes_{\mathbb{C}[q,q^{-1}]}\mathbb{C}[[h]]\simeq V_\lambda$. This module coincides with the one defined in Section 4.1 in \cite{L3} (see also Proposition 10.1.4 in \cite{ChP}).

Below we shall consider the proof in the case of $T$ and $\tilde{T}$. The case of $T^\tau$ and $\tilde{T}^\tau$ is treated in a similar way.

It is well known that any element $T$ of the braid group acts as an invertible linear automorphism of $V'$ which can be specialized to any non--zero numeric value of $q$ in the sense that for any $\varepsilon\in \mathbb{C}^*$ $T$ gives rise to a linear automorphism of $U_q^{res}(\g)/(q-\varepsilon)U_q^{res}(\g)$--module $V'_\varepsilon:=V'/(q-\varepsilon)V'$. It suffices to verify this statement when $T=T_i$ for $i=1,\ldots ,l$, and in this case it follows from the explicit formulas and the results in \cite{Lus}, Section 5.2 or from formula (\ref{T1}).

Recall that elements of the braid group act as Weyl group elements on $\h\subset U_h(\h)$. Assume that the action of $T$ and $T'$ on $\h\subset U_h(\h)$ coincides with the action of a Weyl group element $w$. 
Since the $\mathbb{C}[q,q^{-1}]$--submodule of $V'$ which consists of elements of weight $w\mu$ has rank one and $Tv$ and $T'v$ must belong to this submodule, the relation $Tv=t(q)T' v$ must hold for some rational function $t(q)$ of $q$ with poles or zeroes only at zero and infinity. Indeed, if $t(q_0)=0$, $q_0\neq 0, \infty$ then in $V'_{q_0}$ we have $Tv=0$, i.e. $T$ does not induce an automorphism of $V'_{q_0}$, and if $t^{-1}(q_0')=0$, $q_0'\neq 0, \infty$ then in $V'_{q_0'}$ we have $T'v=0$, i.e. $T'$ does not induce an automorphism of $V'_{q_0'}$. In both cases we arrive at a contradiction. Thus $t(q)$ must be a non--zero multiple of a power of $q$.

\end{proof}


\section{Bibliographic comments}

\pagestyle{myheadings}
\markboth{CHAPTER \thechapter.~QUANTUM GROUPS}{\thesection.~BIBLIOGRAPHIC COMMENTS}

\setcounter{equation}{0}
\setcounter{theorem}{0}

The material presented in Sections \ref{QGdef}, \ref{BGact}, \ref{QRvect}, \ref{PBWbases} and \ref{URmatr} is mostly standard and we refer to books \cite{ChP,Jan,Lus} for more details and omitted proofs.  Formula (\ref{T1}) can be found in \cite{Saito}.  

Realizations of quantum groups associated to Weyl group elements were introduced in \cite{S0} in the case of Coxeter elements and in \cite{S10} in general. 

The results of Section \ref{Adj} are new, except for Propositions \ref{KhT} and \ref{DhFT} the statements of which can be found in \cite{KT1} as Proposition 5.2, Theorem 5.1 and Proposition 5.4.
 
Lemma \ref{tmatrel} is a generalization of the result of Exercise 3 in Chapter V, \S 6, \cite{Bur}.

Specializations of quantum groups similar to those which appear in this book were considered in \cite{SDM}. In this book we introduce slightly different specializations of quantum groups in order to use restricted specializations as well.


\chapter{q-W--algebras}\label{part3}

\pagestyle{myheadings}
\markboth{CHAPTER~\thechapter.~Q-W--ALGEBRAS}{\thesection.~SOME FACTS ON POISSON--LIE GROUPS}

In this chapter we introduce q--W--algebras and study the structure of their the quasi--classical versions, Poisson q-W--algebras. In the next chapter similar results will be obtained for q-W--algebras.

As we briefly mentioned in the introduction the naive definition of q-W--algebras as Hecke type algebras $Hk(A,B, \chi)$ requires some modification. In fact the main ingredient of the definition of q-W--algebras is the adjoint action of the quantum group on itself, and they are defined using a $\mathcal{B}$--subalgebra $\mathbb{C}_\mathcal{B}^s[G_*]$ of the quantum group the restriction of the adjoint action to which is locally finite. When $q$ is specialized to $\varepsilon\in \mathbb{C}^*$ which is not a root of unity the algebra $\mathbb{C}_\mathcal{B}^s[G_*]$ becomes the locally finite part of the quantum group with respect to the adjoint action which was introduced and studied by Joseph. 

The algebra $\mathbb{C}_\mathcal{B}^s[G_*]$ is a quantization of the algebra of regular functions on a Poisson manifold $G_*$ which is isomorphic to $G$ as a manifold and the Poisson structure of which is closely related to that of the Poisson--Lie group $G^*$ dual to a quasitriangular Poisson--Lie group $G$.

After recalling basic facts on Poisson--Lie groups in Section \ref{plgroups} we introduce an algebra $\mathbb{C}^s[G^*]$ of functions on $G^*$ in Section \ref{qplgroups}, its quantization $\mathbb{C}_\mathcal{B}^s[G^*]\subset U_h^s(\g)$ and the subalgebra $\mathbb{C}_\mathcal{B}^s[G_*]\subset \mathbb{C}_\mathcal{B}^s[G^*]$. 

A special choice of the bialgebra structure entering the definitions of $\mathbb{C}^s[G^*]$, $\mathbb{C}_\mathcal{B}^s[G^*]\subset U_h^s(\g)$ and $\mathbb{C}_\mathcal{B}^s[G_*]$ is crucial for the definition of q-W--algebras. It depends on the choice of a Weyl group element $s\in W$ and ensures that one can define a subalgebra $\mathbb{C}_\mathcal{B}^s[M_+]\subset \mathbb{C}_\mathcal{B}^s[G^*]$ equipped with a non--trivial character, so that the q-W--algebra $W_{\mathcal{B}}^s(G)$ can be defined as the result of a quantum constrained reduction with respect to the subalgebra $\mathbb{C}_\mathcal{B}^s[M_+]$.

Next, in Section \ref{wpsred} we proceed with the study of the specialization $W^s(G)$ of the algebra $W_{\mathcal{B}}^s(G)$ at $q^{\frac{1}{rd^2}}=1$. The algebra $W^s(G)$ is naturally a Poisson algebra which can be regarded as the algebra of regular functions on a reduced Poisson manifold which is also an algebraic variety. Poisson reduction works well for differential Poisson manifolds. Therefore it is easier firstly to describe the reduced Poisson structure on the algebra of $C^\infty$--functions on the reduced Poisson manifold and then to recover the structure of the algebraic variety on it. This is done in Proposition \ref{constrt} and Theorem \ref{var}. 

In Section \ref{pZhel} we define a projection operator $\Pi$ onto the algebra $W^s(G)$. In Theorem \ref{proj}, which is central in this chapter, we obtain a formula for the operator $\Pi$ suitable for quantization. This formula plays the key role in the proof of Theorem \ref{Piqmain} describing a localization of the algebra $W_{\mathcal{B}}^s(G)$ in terms of a quantum counterpart of the operator $\Pi$. Miraculously the formula for $\Pi$ from Theorem \ref{proj} can be directly extrapolated to the quantum case.


\section{Some facts on Poisson--Lie groups}\label{plgroups}

\setcounter{equation}{0}
\setcounter{theorem}{0}

In this section, following \cite{ChP}, Ch. 1, we recall some notions related to Poisson--Lie groups. These facts will be needed for the study of Poisson q-W--algebras.

Let $G$ be a finite-dimensional Lie group equipped with a Poisson bracket, \index{Poisson!bracket}
$\frak g$ its Lie algebra. $G$ is called
a {\it Poisson--Lie group} \index{group!Poisson--Lie} if the multiplication $G\times G \rightarrow G$ is a
Poisson map.
A Poisson bracket satisfying this axiom is degenerate and, in particular, is
identically zero
at the unit element of the group. Linearizing this bracket at the unit element $1\in G$
defines the
structure of a Lie algebra in the space $T^*_1 G\simeq {\frak g}^*$.
The pair $({\frak g},{\frak g}^{*})$ \index[not]{g@$({\frak g},{\frak g}^{*})$} is called the {\it tangent Lie bialgebra} \index{Lie!bialgebra!tangent} of $G$.

Lie brackets in $\frak{g}$ and $\frak{g}^{*}$ satisfy the following
compatibility condition:

{\em Let }$\delta: {\frak g}\rightarrow {\frak g}\wedge {\frak g}$ {\em be
the dual  of the commutator map } $[,]_{*}: {\frak g}^{*}\wedge
{\frak g}^{*}\rightarrow {\frak g}^{*}$. {\em Then } $\delta$ {\em is a
1-cocycle on} $  {\frak g}$ {\em (with respect to the adjoint action
of } $\frak g$ {\em on} ${\frak g}\wedge{\frak g}$). \index{cocycle!associated to a bialgebra}

Let $c_{ij}^{k}, f^{ab}_{c}$ be the structure constants of
${\frak g}, {\frak g}^{*}$ with respect to dual bases $\{e_{i}\},
\{e^{i}\}$ in ${\frak g},{\frak g}^{*}$. The compatibility condition
means that

$$
c_{ab}^{s} f^{ik}_{s} ~-~ c_{as}^{i} f^{sk}_{b} ~+~ c_{as}^{k}
f^{si}_{b} ~-~ c_{bs}^{k} f^{si}_{a} ~+~ c_{bs}^{i} f^{sk}_{a} ~~=
~~0.
$$
This condition is symmetric with respect to exchange of $c$ and
$f$. Thus if $({\frak g},{\frak g}^{*})$ is a Lie bialgebra, then
$({\frak g}^{*}, {\frak g})$ is also a Lie bialgebra.

The following proposition shows that the category of finite-dimensional Lie
bialgebras is isomorphic to
the category of finite-dimensional connected simply connected \index{group!Lie!simply connected} Poisson--Lie
groups.
\begin{proposition}{\bf (\cite{ChP}, Theorem 1.3.2)}\label{PLiecorr}
If $G$ is a connected simply connected finite-dimensional Lie group, every
bialgebra structure on $\frak g$
is the tangent bialgebra of a unique Poisson structure on $G$ which makes $G$
into a Poisson--Lie group.
\end{proposition}

Let $G$ be a finite-dimensional Poisson--Lie group, $({\frak g},{\frak g}^{*})$
the tangent bialgebra of $G$.
The connected simply connected finite-dimensional
Poisson--Lie group corresponding to the Lie bialgebra $({\frak g}^{*}, {\frak
g})$ is called the dual
Poisson--Lie group and denoted by $G^*$. \index{group!Poisson--Lie!dual} \index[not]{G@$G^*$}

$({\frak g},{\frak g}^{*})$ is called a factorizable Lie
bialgebra \index{Lie!bialgebra!factorizable} if the following conditions are satisfied (see \cite{fact}):
\begin{enumerate}
\item
${\frak g}${\em \ is equipped with a non--degenerate symmetric invariant
bilinear form} $\left\langle ~\cdot~ ,~\cdot~ \right\rangle$. \index{form!non--degenerate symmetric invariant bilinear!associated to a factorizable Lie bialgebra} \index[not]{ZZ@$\left\langle ~\cdot~ ,~\cdot~ \right\rangle$}

We shall always identify ${\frak g}^{*}$ and ${\frak g}$ by means of this form.

\item  {\em The dual Lie bracket on }${\frak g}^{*}\simeq {\frak g}${\em \
is given by}
\begin{equation}
\left[ X,Y\right] _{*}=\frac 12\left( \left[ rX,Y\right] +\left[ X,rY\right]
\right) ,X,Y\in {\frak g},  \label{rbr}
\end{equation}
{\em where }$r\in {\rm End}_{\mathbb{C}} {\frak g}${\em \ is a skew symmetric, with respect to $\left\langle ~\cdot~ ,~\cdot~ \right\rangle$, linear
operator
(classical r-matrix).} \index{r--matrix!classical!associated to a Lie bialgebra}

\item  $r${\em \ satisfies} {\em the} {\em modified classical Yang-Baxter
equation:} \index{Yang-Baxter equation!modified classical}
\begin{equation}
\left[ rX,rY\right] -r\left( \left[ rX,Y\right] +\left[ X,rY\right] \right)
=-\left[ X,Y\right] ,\;X,Y\in {\frak g}{\bf .}  \label{cybe}
\end{equation}
\end{enumerate}

Define operators $r_\pm \in {\rm End}_{\mathbb{C}} {\frak g}$ \index[not]{r@$r_\pm$} by
\[
r_{\pm }=\frac 12\left( r\pm id\right) .
\]
We shall need some properties of the operators $r_{\pm }$.
Denote by ${\frak i}_\pm$ \index[not]{i@${\frak i}_\pm$} and ${\frak k}_\mp$ \index[not]{k@${\frak k}_\mp$} the image and the kernel of the
operator
$r_\pm $:
\begin{equation}\label{bnpm}
{\frak i}_\pm = Im~r_\pm,~~{\frak k}_\mp = Ker~r_\pm.
\end{equation}

\begin{proposition}{\bf (\cite{BD}, Lemma 6.6; \cite{rmatr}, Sect. 4)}\label{bpm}
Let $({\frak g}, {\frak g}^*)$ be a factorizable Lie bialgebra. Then

(i) ${\frak i}_\pm \subset {\frak g}$ is a Lie subalgebra, \index{Lie!subalgebra} the subspace ${\frak k}_\pm$ is a Lie ideal in \index{Lie!ideal}
${\frak i}_\pm,~{\frak i}_\pm^\perp ={\frak k}_\pm$, where ${\frak i}_\pm^\perp$ is the orthogonal complement to ${\frak i}_\pm$ with respect to the non--degenerate invariant
bilinear form $\left\langle  ~\cdot~ ,~\cdot~ \right\rangle$.

(ii) ${\frak k}_\pm$ is an ideal in ${\frak {g}}^{*}$.

(iii) ${\frak i}_\pm$ is a Lie subalgebra in ${\frak {g}}^{*}$. Moreover ${\frak
i}_\pm ={\frak {g}}^{*}/ {\frak k}_\pm$.

(iv) $({\frak i}_\pm,{\frak i}_\pm ^*)$ is a subbialgebra of $({\frak
{g}},{\frak {g}}^{*})$ \index{Lie!subbialgebra} and
$({\frak i}_\pm,{\frak i}_\pm ^*)\simeq ({\frak i}_\pm,{\frak i}_\mp)$. The
canonical paring $\left\langle ~\cdot~, ~\cdot~ \right\rangle_\pm$
between ${\frak i}_\mp$ and ${\frak i}_\pm$ is given by
\begin{equation}
\left\langle X_\mp ,Y_\pm \right\rangle_\pm=\left\langle X_\mp,r_\pm^{-1}Y_\pm \right\rangle ,~ X_\mp \in {\frak i}_\mp ;~ Y_\pm
\in {\frak i}_\pm .
\end{equation}
\end{proposition}

The classical Yang--Baxter equation implies that $r_{\pm }$ , regarded as a
mapping from
${\frak g}^{*}$ into ${\frak g}$, is a Lie algebra homomorphism (see \cite{rmatr}, Proposition 7). \index{Lie!algebra!homomorphism}
Moreover, $r_{+}^{*}=-r_{-},$\ and $r_{+}-r_{-}=id.$

By Proposition 9 in \cite{rmatr}, the mapping
\begin{eqnarray}\label{imbd}
{\frak {g}}^{*}\rightarrow \g \oplus \g~~~:X\mapsto (X_{+},~X_{-}),~~~X_{\pm
}~=~r_{\pm }X,
\end{eqnarray}
is a Lie algebra embedding \index{Lie!algebra!embedding} (Here $\g \oplus \g$ is the Lie algebra direct sum of two copies of $\g$). Thus we may identify ${\frak g^{*}}$ with the image of this embedding in $\g \oplus \g$. \index[not]{g@$\g \oplus \g$}

Naturally, embedding (\ref{imbd}) extends to a group homomorphism \index{group!homomorphism}
$$
G^*\rightarrow G\times G,~~L\mapsto (L_+,L_-). \index[not]{L@$(L_+,L_-)$}
$$
In the situations considered later in this book this homomorphism will be always an embedding. In such situations we shall identify $G^*$ with the corresponding subgroup in $G\times G$.


\section{Quantization of Poisson--Lie groups}\label{qplgroups}

\pagestyle{myheadings}
\markboth{CHAPTER~\thechapter.~Q-W--ALGEBRAS}{\thesection.~QUANTIZATION OF POISSON--LIE GROUPS}

\setcounter{equation}{0}
\setcounter{theorem}{0}

In this section we introduce Poisson--Lie groups and their quantizations required for defining q-W--algebras. These results are well known in the case of the standard bialgebra structure and we only briefly outline the proofs of the statements in this section. We consider algebras defined over the ring $\mathcal{B}$ since later in our construction the restricted specialization of the quantum group $U_h^s(\g)$ defined over $\mathcal{B}$ will play the key role. 

Let $\frak g$ be a finite-dimensional complex semisimple Lie algebra, $\h\subset \g$ its Cartan subalgebra. Let $s\in W$ be an element of the Weyl group $W$ of the pair $(\g,\h)$ and $\Delta_+$ a system of positive roots in $\Delta=\Delta(\g,\h)$.
Observe that cocycle (\ref{cocycles}) equips
$\frak g$ with
the structure of a factorizable Lie bialgebra, where the non--degenerate symmetric invariant bilinear form is the non--degenerate symmetric invariant bilinear form on $\g$ introduced in Section \ref{notation}.
Using the isomorphism ${\rm End}_{\mathbb{C}}{\frak g}\cong {\frak g}\otimes {\frak g}$ the corresponding r--matrix $r^{s}\in {\rm End}_{\mathbb{C}}{\frak g}\cong {\frak g}\otimes {\frak g}$ may be represented as
$$
r^{s}=P_+-P_-+\kappa{1+s \over 1-s}P_{\h'},
$$
where $P_+,P_-$ and $P_{\h'}$ are the orthogonal projection operators onto the nilpotent radical ${\frak n}_+$ \index{radical!nilpotent} corresponding to $\Delta_+$, the opposite nilpotent radical ${\frak
n}_-$, and ${\frak h}'$, respectively, in
the direct sum
$$
{\frak g}={\frak n}_+ +{\frak h}'+\h'^\perp + {\frak n}_-,
$$
and $\h'^\perp=\h^s$ is the orthogonal complement to $\h'$ in $\h$ with respect to the symmetric bilinear form.

Let $G$ be the connected simply connected semisimple Poisson--Lie group with the
tangent Lie bialgebra $({\frak g},{\frak g}^*)$,
$G^*$ the dual Poisson--Lie group.

Observe that $G$ is a connected simply connected semisimple complex algebraic group (see e.g. \S 104, Theorem 12 in \cite{Z}). Note also that 
\begin{equation}\label{rspm}
r^{s}_+=P_+ +\frac{\kappa}{2}{1+s \over 1-s}P_{\h'}+\frac{1}{2}P_{\h},~~r^{s}_-=-P_- +\frac{\kappa}{2}{1+s \over 1-s}P_{\h'}-\frac{1}{2}P_{\h},
\end{equation}
where $P_{{\h}}$ is the orthogonal projection operator onto ${\frak h}\subset \g$ with respect to the symmetric bilinear form, and hence the subspaces ${\frak i}_\pm$ and ${\frak k}_\pm$ defined by (\ref{bnpm}) coincide with
the Borel subalgebras \index{Lie!subalgebra!Borel} ${\frak b}_\pm$ in $\frak g$ corresponding to $\Delta_\pm$ and their nilpotent radicals ${\frak n}_\pm$, respectively. \index{radical!nilpotent!of a Borel subalgebra}
Therefore every element $(L_+,L_-)\in G^*\subset G\times G$ may be uniquely written as
\begin{equation}\label{fact}
(L_+,L_-)=(n_+,n_-)(h_+,h_-), \index[not]{L@$(L_+,L_-)$}
\end{equation}
where $n_\pm \in N_\pm$, $h_+={\exp}((\frac{\kappa}{2}{1+s \over 1-s}P_{{\h'}}+\frac{1}{2}id)x),~h_-={\exp}((\frac{\kappa}{2}{1+s \over 1-s}P_{{\h'}}-\frac{1}{2}id)x),~x\in
{\frak h}$.
In particular, $G^*$ is a solvable subgroup \index{subgroup!solvable} in $G\times G$. In general, the group $G^*$ is not algebraic.

Note that embedding (\ref{imbd}) restricts to embeddings
$$
\b_\pm\rightarrow \g \oplus \g,
$$
and that by Proposition \ref{bpm} (iv) $(\b_\pm,\b_\pm^*)\simeq (\b_\pm,\b_\mp)$ $(\b_\pm^*\simeq \b_\mp)$ is a subbialgebra of $(\g,\g^*)$. Therefore $B_\pm\subset G^*$ are Poisson--Lie subgroups.

Observe also that one can consider $\n_\pm$ as Lie subalgebras in $\g^*$ via embeddings
\begin{equation}\label{npmg*inc1}
\n_+\rightarrow \g^*\subset \g\oplus \g,~x\mapsto (x,0),
\end{equation}
\begin{equation}\label{npmg*inc2}
\n_-\rightarrow \g^*\subset \g\oplus \g,~x\mapsto (0,x),
\end{equation}
where $\g^*$ is regarded as a Lie subalgebra of $\g\oplus \g$ using embedding (\ref{imbd}).

Using these embeddings the algebraic subgroups $N_\pm\subset G$ corresponding to the algebraic Lie subalgebras $\n_\pm\subset \g$ can be regarded as Lie subgroups in $B_\pm\subset G^*$ corresponding to the Lie subalgebras $\n_\pm\subset \b_\pm\subset \g^*$.
However, $N_\pm$ are not Poisson--Lie subgroups of $G^*$.

Our main object will be a certain algebra of functions on $G^*$, $\mathbb{C}^s[G^*]$. \index[not]{C@$\mathbb{C}^s[G^*]$}
This algebra may be explicitly described as follows.
Let $\pi_V$ \index[not]{p@$\pi_V$} be a finite-dimensional representation of $G$ in a vector space $V$. Then matrix
elements of $\pi_V(L_\pm)$ \index[not]{p@$\pi_V(L_\pm)$} are well--defined functions on $G^*$, and $\mathbb{C}^s[G^*]$ is defined as the subalgebra in $C^\infty(G^*)$ \index[not]{C@$C^\infty(G^*)$} generated by the matrix elements of $\pi_V(L_\pm)$, \index{matrix element!of a representation of an algebraic group} where $V$ runs through all finite--dimensional representations of $G$.
The elements $L^{\pm,V}:=\pi_V(L_\pm)$ \index[not]{L@$L^{\pm,V}$} may be viewed as elements of the space $\mathbb{C}^s[G^*]\otimes {\rm End}_{\mathbb{C}}V$. If we fix a basis in $V$, $L^{\pm,V}$ can be regarded as matrices with entries $L^{\pm,V}_{ij}$ being elements of $\mathbb{C}^s[G^*]$.

For every two finite--dimensional ${\frak g}$--modules $V$ and $W$ we denote ${r^s_\pm}^{VW}=(\pi_V\otimes \pi_W)r^s_\pm\in {\rm End}_{\mathbb{C}}V\otimes {\rm End}_{\mathbb{C}}W$, \index[not]{r@${r^s_\pm}^{VW}$} where $r^s_\pm$ is regarded as an element of ${\frak g}\otimes {\frak g}$.
\begin{proposition}{\bf (\cite{dual}, Section 2)}\label{pbff}
$\mathbb{C}^s[G^*]$ is a Poisson subalgebra \index{Poisson!subalgebra} in the Poisson algebra $C^\infty(G^*)$, \index{Poisson!algebra}
the Poisson bracket \index{Poisson!bracket} on $\mathbb{C}^s[G^*]$ is given by
\begin{equation}\label{pbf}
\begin{array}{l}
\{L^{\pm,V}_{1},L^{\pm,W}_{2}\}=
-2[{r_\pm^s}^{VW},L^{\pm,V}_{1}L^{\pm,W}_{2}],\\
\\
\{L^{-,V}_{1},L^{+,W}_{2}\}~=-2[{r_{\pm}^s}^{VW},L^{-,V}_{1}L^{+,W}_{2}],
\end{array}
\end{equation}
where
$$
L^{\pm,V}_1=L^{\pm,V}\otimes I_W,~~L^{\pm,W}_2=I_V\otimes L^{\pm,W},L^{\pm,V}_1,L^{\pm,W}_2\in \mathbb{C}^s[G^*]\otimes {\rm End}_{\mathbb{C}}V\otimes {\rm End}_{\mathbb{C}}W, \index[not]{L@$L^{\pm,V}_{1,2}$}
$$
$I_X\in {\rm End}_{\mathbb{C}}X$ is the identity endomorphism of a vector space $X$,  ${r^s_\pm}^{VW}$ is regarded as an element of $\mathbb{C}^s[G^*]\otimes {\rm End}_{\mathbb{C}}V\otimes {\rm End}_{\mathbb{C}}W$ via the embedding
$$
{\rm End}_{\mathbb{C}}V\otimes {\rm End}_{\mathbb{C}}W \to \mathbb{C}^s[G^*]\otimes {\rm End}_{\mathbb{C}}V\otimes {\rm End}_{\mathbb{C}}W, a\otimes b\mapsto 1\otimes a\otimes b,
$$
in the left hand sides in (\ref{pbf}) the Poisson brackets of $L^{\pm,V}_1,L^{\pm,W}_{2}\in \mathbb{C}^s[G^*]\otimes {\rm End}_{\mathbb{C}}V\otimes {\rm End}_{\mathbb{C}}W$ are taken with respect to the Poisson structure on $\mathbb{C}^s[G^*]$, the composition of endomorphisms is taken in ${\rm End}_{\mathbb{C}}V\otimes {\rm End}_{\mathbb{C}}W$, and the commutators in the right hand sides of (\ref{pbf}) are taken in $\mathbb{C}^s[G^*]\otimes {\rm End}_{\mathbb{C}}V\otimes {\rm End}_{\mathbb{C}}W$. 

Moreover, the map $\Delta^s:\mathbb{C}^s[G^*]\rightarrow \mathbb{C}^s[G^*]\otimes \mathbb{C}^s[G^*]$ \index[not]{D@$\Delta^s(~\cdot~)$} dual to the multiplication in $G^*$,
\begin{equation}\label{comultcl}
\Delta^s(L^{\pm,V}_{ij})=\sum_k L^{\pm,V}_{ik}\otimes L^{\pm,V}_{kj},
\end{equation}
is a homomorphism of Poisson algebras, \index{Poisson!algebra!homomorphism} and the map $S^s:\mathbb{C}^s[G^*]\rightarrow \mathbb{C}^s[G^*]$ \index[not]{S@$S^s(~\cdot~)$} dual to taking inverse in $G^*$,
$$
S^s(L^{\pm,V}_{ij})=((L^{\pm,V})^{-1})_{ij},
$$
is an antihomomorphism of Poisson algebras. \index{Poisson!algebra!antihomomorphism}
\end{proposition}
\begin{remark}
Recall that a Poisson--Hopf algebra \index{Poisson--Hopf algebra} is a Poisson algebra which is also a Hopf algebra such that the comultiplication is a homomorphism of Poisson algebras and the antipode is an antihomomorphism of Poisson algebras. According to Proposition \ref{pbff} $\mathbb{C}^s[G^*]$ is a Poisson--Hopf algebra.
\end{remark}


Now we define a quantization of the Poisson--Hopf algebra $\mathbb{C}^s[G^*]$.
For any finite rank representation $\pi_V: U_\mathcal{B}^{s, res}({\frak
g})\rightarrow V^{res}$, \index[not]{p@$\pi_V$} where $V$ is a finite rank representation of $U_h({\frak
g})$, one can define an action of elements $H_i$, $i=1,\ldots ,l$ on $V^{res}$ by requiring that $H_i$ acts on weight vectors of weight $\lambda$ by multiplication by $\lambda(H_i)$. 
Then from the definition of the R--matrix $\mathcal{R}_s$ and from formula (\ref{rmatrspires}) it follows that the elements ${^q{L^{\pm,V}}}$ \index[not]{L@${^q{L^{\pm,V}}}$} given by
$$
{^q{L^{-,V}}}=(id\otimes \pi_V){(\mathcal{R}_s)_{21}}^{-1}=(id\otimes
\pi_VS_{s})(\mathcal{R}_s)_{21}
,~~ {^q{L^{+,V}}}=(id\otimes \pi_V)\mathcal{R}_s.
$$
are well--defined invertible elements of $U_h^{s}({\frak g})\otimes {\rm End}_\mathcal{B} (V^{res})$.

If we fix a basis in $V^{res}$, ${^q{L^{\pm,V}}}$ may be regarded as matrices
with matrix elements $({^q{L^{\pm,V}}})_{ij}$ being elements of $U_h^{s}({\frak g})$.

From the quantum Yang--Baxter equation (\ref{YB}) for $\mathcal{R}_s$ we get
relations between ${^q{L^{\pm,V}}}$,
\begin{equation}\label{ppcomm}
\begin{array}{l}
R^{VW}{^q{L^{\pm,W}_2}}{^q{L^{\pm,V}_1}}={^q{L^{\pm,V}_1}}{^q{L^{\pm,W}_2}}R^{VW
},
\end{array}
\end{equation}
\begin{equation}\label{pmcomm}
R^{VW}{^q{L^{+,W}_2}}{^q{L^{-,V}_1}}={^q{L^{-,V}_1}}{^q{L^{+,W}_2}}R^{VW},
\end{equation}
where by ${^q{L^{\pm,W}_1}},~{^q{L^{\pm,V}_2}}$ we understand the following elements of 
$U_h^{s}({\frak g})\otimes {\rm End}_\mathcal{B} (V^{res})\otimes {\rm End}_\mathcal{B} (W^{res})$,
$$
{^q{L^{\pm,V}_1}}={^q{L^{\pm,V}}}\otimes I_{W^{res}},~~{^q{L^{\pm,W}_2}}=I_{V^{res}}\otimes
{^q{L^{\pm,W}}}, \index[not]{L@${^q{L^{\pm,V}_{1,2}}}$}
$$
$I_X$ is identity endomorphism of $X$, and $R^{VW}\in {\rm End}_\mathcal{B} (V^{res})\otimes {\rm End}_\mathcal{B} (W^{res})$ is regarded as an element of $U_h^{s}({\frak g})\otimes {\rm End}_\mathcal{B} (V^{res})\otimes {\rm End}_\mathcal{B} (W^{res})$ via the embedding
$$
{\rm End}_\mathcal{B} (V^{res})\otimes {\rm End}_\mathcal{B} (W^{res})\to U_h^{s}({\frak g})\otimes {\rm End}_\mathcal{B} (V^{res})\otimes {\rm End}_\mathcal{B} (W^{res}), a\otimes b\mapsto 1\otimes a\otimes b.
$$

From (\ref{rmprop}) we can obtain the action of the comultiplication on ${^q{L^{\pm,V}}}$:
\begin{equation}\label{comultGq}
\Delta_s({^q{L^{\pm,V}_{ij}}})=\sum_k {^q{L^{\pm,V}_{ik}}}\otimes
{^q{L^{\pm,V}_{kj}}},
\end{equation}
and of the antipode,
\begin{equation}\label{ants}
S_s({^q {L^{\pm,V}_{ij}}})=(({^q{L^{\pm,V}}})^{-1})_{ij}.
\end{equation}

We denote by $\mathbb{C}_{\mathcal{B}}^s[G^*]$ \index[not]{C@$\mathbb{C}_\mathcal{B}^s[G^*]$} the $\mathcal{B}$--subalgebra in $U_h^{s}({\frak
g})$ generated by the matrix elements of $({^q{L^{+,V}}})^{\pm 1}$ and of $({^q{L^{-,V}}})^{\pm 1}$, \index{matrix element!of a representation of a quantum group} where $V$ runs through all finite rank representation of $U_h({\frak g})$. Formulas (\ref{comultGq}) and (\ref{ants}) imply that $\mathbb{C}_{\mathcal{B}}^s[G^*]$ is in fact a Hopf subalgebra \index{Hopf subalgebra} in $U_h^{s}({\frak g})$.

Since $\mathcal{R}_s=1\otimes 1$ (mod $h$) relations (\ref{ppcomm}) and (\ref{pmcomm}) imply that the quotient algebra $\mathbb{C}_{\mathcal{B}}^s[G^*]/(q^{\frac{1}{d{\bar{r}}^2}}-1)\mathbb{C}_{\mathcal{B}}^s[G^*]$ is commutative, and one can equip it with a Poisson structure given by \index{Poisson!bracket}
\begin{equation}\label{quasipb}
\{x_1,x_2\}=\frac{1}{d{\bar{r}}^2}{[a_1,a_2] \over q^{\frac{1}{d{\bar{r}}^2}}-1}~(\mbox{mod }(q^{\frac{1}{d{\bar{r}}^2}}-1)), \index[not]{ZZZ@$\{~\cdot~ ,~\cdot~ \}$}
\end{equation}
where $a_1,a_2\in \mathbb{C}_{\mathcal{B}}^s[G^*]$ reduce to
$x_1,x_2\in \mathbb{C}_{\mathcal{B}}^s[G^*]/(q^{\frac{1}{d{\bar{r}}^2}}-1)\mathbb{C}_{\mathcal{B}}^s[G^*]$ mod $(q^{\frac{1}{d{\bar{r}}^2}}-1)$.

Obviously, the comultiplication and the antipode in $\mathbb{C}_{\mathcal{B}}^s[G^*]$ induce a comultiplication and an antipode in $\mathbb{C}_{\mathcal{B}}^s[G^*]/(q^{\frac{1}{d{\bar{r}}^2}}-1)\mathbb{C}_{\mathcal{B}}^s[G^*]$ compatible with the introduced Poisson structure, and the quotient $\mathbb{C}_{\mathcal{B}}^s[G^*]/(q^{\frac{1}{d{\bar{r}}^2}}-1)\mathbb{C}_{\mathcal{B}}^s[G^*]$ becomes a Poisson--Hopf algebra.

Let $\mathbb{C}_{\mathcal{B}}^s[G]$ \index[not]{C@$\mathbb{C}_\mathcal{B}^s[G]$} be the restricted Hopf algebra dual \index{Hopf algebra!restricted dual} to $U_{\mathcal{B}}^{s, res}(\g)$ which is spanned by the matrix elements of finite rank representations of $U_{\mathcal{B}}^{s, res}(\g)$ of the form $V^{res}$, where $V$ is a finite rank representation of $U_h^{s}({\frak g})$. By Proposition \ref{wresdec} (i) $\mathbb{C}_{\mathcal{B}}^s[G]$ is naturally $P\times P$--graded. Note that for different $s$ the Hopf algebras $\mathbb{C}_{\mathcal{B}}^s[G]$ are isomorphic as coalgebras since $U_{\mathcal{B}}^{s, res}(\g)$ are isomorphic as algebras for different $s$ with the help of  the isomorphisms $\psi_{\{n_{ij}\}}$ with $n_{ij}=\frac{1}{2d_j}c_{ij}$ as explained in Section \ref{QGspec}.

We shall use the following notation for elements of $\mathbb{C}_{\mathcal{B}}^s[G]$. Let $V^{res}$ be a $U_{\mathcal{B}}^{s,res}(\g)$--lattice in a finite rank $U_h(\g)\simeq U_h^s(\g)$--module $V$. Recall that there is a contravariant non--degenerate form $(~\cdot~,~\cdot~)$ on $V$, such that $(u,xv)=(\omega(x)u,v)$ for any $u,v\in V$, $x\in U_h(\g)$, which gives rise to a contravariant non--degenerate form on the $U_q^s(\g)$--module $V_q=V^{res}\otimes_{\mathcal{B}}\mathbb{C}(q^{\frac{1}{d{\bar{r}}^2}})$, and different weight subspaces $(V_q)_\lambda$, $\lambda \in P$ are orthogonal with respect to this form. 

Assume that $u\in V_q$ is such that $(u,w)\in \mathcal{B}$ for any $w\in V^{res}$. Then $u^*(~\cdot~):=(u,~\cdot~)$ is an element of the dual module  ${V^{res}}^*$. Since $V_q$ and $V^{res}$ are of finite ranks, the weight subspaces of $V^{res}$ are free $\mathcal{B}$--modules, and the restriction of the contravariant form $(~\cdot~,~\cdot~)$ to each weight subspace of $V_q$ is non--degenerate, all elements of ${V^{res}}^*$ can be obtained this way. Clearly, for any $v\in V^{res}$ $u^*(~\cdot~ v)=(u,~\cdot~ v)\in \mathbb{C}_{\mathcal{B}}^s[G]$, and by the definition $\mathbb{C}_{\mathcal{B}}^s[G]$ is spanned by such elements. 

If $V$ is a complex algebraic variety, \index{variety!algebraic} we denote by ${\mathbb{C}}[V]$ \index[not]{C@$\mathbb{C}[V]$} the algebra of regular functions \index{function!regular} on $V$. By the definition of $\mathbb{C}_{\mathcal{B}}^s[G]$ we have the following proposition.
\begin{proposition}\label{CG1}
The quotient $\mathbb{C}_{\mathcal{B}}^s[G]/(q^{\frac{1}{d{\bar{r}}^2}}-1)\mathbb{C}_{\mathcal{B}}^s[G]$ is isomorphic to ${\mathbb{C}}[G]$ \index[not]{C@$\mathbb{C}[G]$} as a Hopf algebra,  $\mathbb{C}_{\mathcal{B}}^s[G]/(q^{\frac{1}{d{\bar{r}}^2}}-1)\mathbb{C}_{\mathcal{B}}^s[G] \simeq {\mathbb{C}}[G]$.
\end{proposition}

The following properties and the description of the Hopf algebra $\mathbb{C}_{\mathcal{B}}^s[G^*]$ are completely analogous to similar results in the case of the standard bialgebra structure which corresponds to the situation when $s=1$ in our setting.
\begin{proposition}\label{qG*}
Let $s\in W$ be an element of $W$, and $e_{\beta}$, $f_{\beta}$, $\beta\in \Delta_+$ the quantum root vectors in $U_h^s(\g)$ defined with the help of a normal ordering in $\Delta_+$. Then the following statements are true.
 
(i) $\mathbb{C}_{\mathcal{B}}^s[G^*]$ is the $\mathcal{B}$--subalgebra in $U_h^{s}({\frak g})$ generated by the elements 
$$
q^{\pm(Y_i-\kappa \frac{1+s}{1-s}P_{\h'}Y_i)},~q^{\pm(-Y_i-\kappa \frac{1+s}{1-s}P_{\h'}Y_i)},~i=1,\ldots, l,~
\tilde f_{\beta}:=(1-q_\beta^{-2})f_{\beta},~\index[not]{f@$\tilde f_{\beta}$} \tilde e_{\beta}:=(1-q_\beta^{2})e_{\beta}e^{h\beta^\vee},~ \index[not]{e@$\tilde e_{\beta}$} \beta\in \Delta_+.
$$

(ii) The matrix elements of ${^q{L^{+,V}}}^{\pm 1}$ (resp. ${^q{L^{-,V}}}^{\pm 1}$) form a Hopf subalgebra $\mathbb{C}_{\mathcal{B}}^s[B_+]\subset \mathbb{C}_{\mathcal{B}}^s[G^*]$ \index[not]{C@$\mathbb{C}_\mathcal{B}^s[B_\pm]$} (resp. $\mathbb{C}_{\mathcal{B}}^s[B_-]\subset \mathbb{C}_{\mathcal{B}}^s[G^*]$). 

(iii) Let $\mathbb{C}_{\mathcal{B}}^s[H]\subset \mathbb{C}_{\mathcal{B}}^s[G^*]$ \index[not]{C@$\mathbb{C}_\mathcal{B}^s[H]$} be the subalgebra generated by the elements $q^{\pm(Y_i-\kappa \frac{1+s}{1-s}P_{\h'}Y_i)}$, $q^{\pm(-Y_i-\kappa \frac{1+s}{1-s}P_{\h'}Y_i)}$, $i=1,\ldots, l$. Then $\mathbb{C}_{\mathcal{B}}^s[H]$ is $\mathcal{B}$--free with a countable basis  $\bar{V}_i$, $i\in \mathbb{N}$. \index[not]{V@$\bar{V}_i$}

(iv) The algebra $\mathbb{C}_{\mathcal{B}}^s[G^*]$ contains the subalgebra $\mathbb{C}_{\mathcal{B}}^s[N_+]$ \index[not]{C@$\mathbb{C}_\mathcal{B}^s[N_\pm]$} (resp. $\mathbb{C}_{\mathcal{B}}^s[N_-]$)  generated by the elements $\tilde f_{\beta}$ (resp. $\tilde e_{\beta}$), $\beta\in \Delta_+$, and the algebra $\mathbb{C}_{\mathcal{B}}^s[N_+]$ (resp. $\mathbb{C}_{\mathcal{B}}^s[N_-]$) is free over $\mathcal{B}$, and the elements $\tilde f_{\beta_D}^{m_D}\ldots \tilde f_{\beta_1}^{m_1}$ or ${\tilde{f}_{\beta_1}}^{m_1}{\tilde{f}_{\beta_2}}^{m_2}\ldots {\tilde{f}_{\beta_D}}^{m_D}$ (resp. $\tilde e_{\beta_D}^{m_D}\ldots \tilde e_{\beta_1}^{m_1}$ or ${\tilde{e}_{\beta_1}}^{m_1}{\tilde{e}_{\beta_2}}^{m_2}\ldots {\tilde{e}_{\beta_D}}^{m_D}$) with $m_j,\in \mathbb{N}$ is a $\mathcal{B}$--basis of $\mathbb{C}_{\mathcal{B}}^s[N_+]$ (resp. $\mathbb{C}_{\mathcal{B}}^s[N_-]$).

(v) ${\mathbb{C}}_{\mathcal{B}}[{G}^*]$ is free over $\mathcal{B}$, and the multiplication defines an isomorphism of $\mathcal{B}$--modules,
\begin{equation}\label{CGqiso}
\mathbb{C}_{\mathcal{B}}^s[N_-]\otimes \mathbb{C}_{\mathcal{B}}^s[H]\otimes \mathbb{C}_{\mathcal{B}}^s[N_+]\to {\mathbb{C}}_{\mathcal{B}}[{G}^*],
\end{equation}
so that the elements $$\tilde e_{\beta_1}^{r_1}\ldots \tilde e_{\beta_D}^{r_D}\bar{V}_i \tilde f_{\beta_D}^{m_D}\ldots \tilde f_{\beta_1}^{m_1}$$ with $r_j,m_j,i\in \mathbb{N}$, $j=1,\ldots ,D$ form a $\mathcal{B}$--basis in $\mathbb{C}_{\mathcal{B}}^s[G^*]$.

(vi) The Poisson--Hopf algebra $\mathbb{C}_{\mathcal{B}}^s[G^*]/(q^{\frac{1}{d{\bar{r}}^2}}-1)\mathbb{C}_{\mathcal{B}}^s[G^*]$ is isomorphic to $\mathbb{C}^s[G^*]$ as a Poisson--Hopf
algebra. 

(vii) The quotients $\mathbb{C}_{\mathcal{B}}^s[N_\pm]/(q^{\frac{1}{d{\bar{r}}^2}}-1)\mathbb{C}_{\mathcal{B}}^s[N_\pm]$ are naturally Poisson subalgebras of $\mathbb{C}^s[G^*]\simeq \mathbb{C}_{\mathcal{B}}^s[G^*]/(q^{\frac{1}{d{\bar{r}}^2}}-1)\mathbb{C}_{\mathcal{B}}^s[G^*]$, and we have algebra isomorphisms ${\mathbb{C}}[N_\pm]\simeq \mathbb{C}_{\mathcal{B}}^s[N_\pm]/(q^{\frac{1}{d{\bar{r}}^2}}-1)\mathbb{C}_{\mathcal{B}}^s[N_\pm]$. \index[not]{C@$\mathbb{C}[N_\pm]$} Thus ${\mathbb{C}}[N_\pm]$ are naturally Poisson subalgebras of $\mathbb{C}^s[G^*]$, and the algebras $\mathbb{C}_{\mathcal{B}}^s[N_\pm]$ are deformations of the algebras of regular functions on the subgroups $N_\pm\subset G^*$.
\end{proposition}

\begin{proof}
(i) This property is analogous to the results obtained in Section 1.4 and Theorem 3.2 of \cite{DKPschubert}, and  in Proposition 4.2 and Theorem 4.6 in \cite{DL}, in the case of the quantum group associated to the standard bialgebra structure. 

We only briefly outline the main steps in the proof as the arguments used in the case of the standard bialgebra structure can be applied verbatim in the setting of this lemma.

Let $\omega_i$, $i=1, \ldots , l$ be the fundamental weights of $\g$. Then from the definition of the elements ${^q{L^{\pm,V}}}$ and from formula (\ref{rmatrspires}) it follows that $q^{Y_i-\kappa \frac{1+s}{1-s}P_{\h'}Y_i}$ is equal to the matrix element $(id\otimes g_i)\mathcal{R}_s$, where $g_i(~\cdot~)=(v_{\omega_i},~\cdot~ v_{\omega_i})$, and $v_{\omega_i}$ is the highest weight vector of $V^{res}_{\omega_i}$ normalized in such a way that $(v_{\omega_i},v_{\omega_i})=1$. Other elements from the set $q^{\pm(Y_i-\kappa \frac{1+s}{1-s}P_{\h'}Y_i)}, q^{\pm(-Y_i-\kappa \frac{1+s}{1-s}P_{\h'}Y_i)},~i=1,\ldots, l$ are obtained in a similar way.

The fact that $\tilde f_{\beta}=(1-q_\beta^{-2})f_{\beta},~\tilde e_{\beta}=(1-q_\beta^{2})e_{\beta}e^{h\beta^\vee},~\beta\in \Delta_+$ belong to $\mathbb{C}_{\mathcal{B}}^s[G^*]$ follows, e.g,  from a slight modification of Theorem 3.2 in \cite{GY} (one should just replace the quantum R-matrix associated with the standard bialgebra structure with R-matrix (\ref{rmatrspires})). This theorem along with the description of the generators $q^{\pm(Y_i-\kappa \frac{1+s}{1-s}P_{\h'}Y_i)}, q^{\pm(-Y_i-\kappa \frac{1+s}{1-s}P_{\h'}Y_i)},~i=1,\ldots, l$ as matrix elements and the definition of the algebra $\mathbb{C}_{\mathcal{B}}^s[G^*]$ also imply that the elements listed in the statement of part (i) of this proposition generate $\mathbb{C}_{\mathcal{B}}^s[G^*]$.

(ii) This follows from the definition of the elements ${^q{L^{\pm,V}}}$ and from formula (\ref{comultGq}).

(iii) Recall that by part (i) of this proposition $\mathbb{C}_{\mathcal{B}}^s[G^*]$ is the $\mathcal{B}$--subalgebra in $U_h^{s}({\frak g})$ generated by the elements $q^{\pm(Y_i-\kappa \frac{1+s}{1-s}P_{\h'}Y_i)}$, $q^{\pm(-Y_i-\kappa \frac{1+s}{1-s}P_{\h'}Y_i)}$, $i=1,\ldots, l$, $\tilde f_{\beta}=(1-q_\beta^{-2})f_{\beta}$, $\tilde e_{\beta}=(1-q_\beta^{2})e_{\beta}e^{h\beta^\vee}$, $\beta\in \Delta_+$. The subalgebra $\mathbb{C}_{\mathcal{B}}^s[H]$ of  $\mathbb{C}_{\mathcal{B}}^s[G^*]$  generated by the elements $q^{\pm(Y_i-\kappa \frac{1+s}{1-s}P_{\h'}Y_i)}, q^{\pm(-Y_i-\kappa \frac{1+s}{1-s}P_{\h'}Y_i)}$, $i=1,\ldots, l$ is in turn a subalgebra of the $\mathcal{B}$--subalgebra $U_{\mathcal{B}}'(\h)\subset U_h^{s}({\frak g})$ \index[not]{U@$U_\mathcal{B}'(\h)$} generated by the elements $\bar{U}_i:=q^{\frac{1}{d{\bar{r}}^2}Y_i}, \bar{U}_i^{-1}$, $i=1,\ldots, l$. \index[not]{U@$\bar{U}_i$} The last algebra is obviously $\mathcal{B}$--free with the countable basis consisting of the products $\bar{U}_1^{m_1}\ldots \bar{U}_l^{m_l}$, $m_1,\ldots ,m_l\in \mathbb{Z}$.  Since $\mathcal{B}$ is a principal ideal domain \index{principal ideal domain} the subalgebra in $\mathbb{C}_{\mathcal{B}}^s[H]\subset U_{\mathcal{B}}'(\h)$ is also $\mathcal{B}$--free by Theorem 6.5 in \cite{R}. We denote by $\bar{V}_i$, $i\in \mathbb{N}$ the elements of some $\mathcal{B}$--basis of this subalgebra.  

(iv) From Step 1 of the proof of the Theorem in Section 12.1 in \cite{DP} and from the following commutation relations in $U^s_h(\g)$
\begin{equation}\label{Ccomms}
e^{hX}\tilde{e}_\alpha e^{-hX}=q^{\alpha(X)}\tilde{e}_\alpha,~e^{hX}\tilde{f}_\alpha e^{-hX}=q^{-\alpha(X)}\tilde{f}_\alpha,~X\in \h\subset U^s_h(\g),
\end{equation}
which are consequences of (\ref{Ccomm}), it follows that the products $\tilde f_{\beta_D}^{m_D}\ldots \tilde f_{\beta_1}^{m_1}$ (resp. ${\tilde{e}_{\beta_1}}^{m_1}{\tilde{e}_{\beta_2}}^{m_2}\ldots {\tilde{e}_{\beta_D}}^{m_D}$) with $m_j,\in \mathbb{N}$ is a $\mathcal{B}$--basis of $\mathbb{C}_{\mathcal{B}}^s[N_+]$ (resp. $\mathbb{C}_{\mathcal{B}}^s[N_-]$). Applying the algebra anti-involution $\omega_0'$ to these elements we obtain that ${\tilde{f}_{\beta_1}}^{m_1}{\tilde{f}_{\beta_2}}^{m_2}\ldots {\tilde{f}_{\beta_D}}^{m_D}$ (resp. $\tilde e_{\beta_D}^{m_D}\ldots \tilde e_{\beta_1}^{m_1}$) with $m_j,\in \mathbb{N}$ is a $\mathcal{B}$--basis of $\mathbb{C}_{\mathcal{B}}^s[N_+]$ (resp. $\mathbb{C}_{\mathcal{B}}^s[N_-]$).

(v) From Step 1 of the proof of the Theorem in Section 12.1 in \cite{DP}, from part (iii) of this proposition, and from commutation relations (\ref{Ccomms}) it follows that the elements $$\tilde e_{\beta_1}^{r_1}\ldots \tilde e_{\beta_D}^{r_D}\bar{V}_i \tilde f_{\beta_D}^{m_D}\ldots \tilde f_{\beta_1}^{m_1}$$ with $r_j,m_j,i\in \mathbb{N}$, $j=1,\ldots ,D$ form a $\mathcal{B}$--basis in $\mathbb{C}_{\mathcal{B}}^s[G^*]$. This implies isomorphism (\ref{CGqiso}).

(vi) In the case when the underlying bialgebra structure is the standard one used in the definition of the Drinfeld--Jimbo quantum group the statement of part (vi) is the Theorem in Section 12.1 in \cite{DP}. The arguments in the proof given there can be used verbatim to prove this proposition. Note, however, that our approach to the definition of the algebra $\mathbb{C}_{\mathcal{B}}^s[G^*]$, which follows the ideas of \cite{FRT}, \S 2, is different from the one adapted in \cite{DP}. The two approaches are equivalent. But in our framework one can give a more straightforward and shorter proof of the statement of part (vi) which is presented below.

Denote by $p:\mathbb{C}_{\mathcal{B}}^s[G^*] \rightarrow \mathbb{C}_{\mathcal{B}}^s[G^*]/(q^{\frac{1}{d{\bar{r}}^2}}-1)\mathbb{C}_{\mathcal{B}}^s[G^*]:=\mathbb{C}^s[G^*]'$ \index[not]{C@$\mathbb{C}^s[G^*]'$} \index[not]{p@$p(~\cdot~)$} the
canonical projection. We also denote by the same letter the
canonical projection $p: U_\mathcal{B}^{s, res}({\frak g})\to U_\mathcal{B}^{s, res}({\frak g})/(q^{\frac{1}{d{\bar{r}}^2}}-1)U_\mathcal{B}^{s, res}({\frak g})=U_1^{s,res}(\g)$.

If $V$ is a free over $\mathbb{C}[[h]]$ of finite rank $U_h(\g)$--module, $\pi_V$ the corresponding representation of $U_\mathcal{B}^{s, res}({\frak g})$ in the space $V^{res}$ equipped also with the natural action of $\h$ with respect to which elements $H_i$, $i=1,\ldots ,l$ act on weight vectors of weight $\lambda$ by multiplication by $\lambda(H_i)$, then $\overline{V}:=V^{res}/(q^{\frac{1}{d{\bar{r}}^2}}-1)V^{res}$ \index[not]{V@$\overline{V}$} is a finite--dimensional module over $U^{s,res}_1(\g)$ equipped with the natural action of $\h$ with respect to which elements $H_i$, $i=1,\ldots ,l$ act on weight vectors of weight $\lambda$ by multiplication by $\lambda(H_i)$. We denote by $\pi_{\overline{V}}$ the corresponding representation of $U(\g)$.
 
Let $p_V:V^{res}\rightarrow \overline{V}=V^{res}/(q^{\frac{1}{d{\bar{r}}^2}}-1)V^{res}$ \index[not]{p@$p_V$} be the projection of $V^{res}$ onto $\overline{V}$. Then any element $T\in {\rm End}_{\mathcal{B}}(V^{res})$ gives rise to an element $\nu(T)\in {\rm End}_{\mathbb{C}}(\overline{V})$ \index[not]{n@$\nu(~\cdot~)$} defined by $\nu(T)\bar{v}=p_VTv$, where $v\in V^{res}$ is any representative of $\bar{v}\in\overline{V}=V^{res}/(q^{\frac{1}{d{\bar{r}}^2}}-1)V^{res}$. If $x\in U_\mathcal{B}^{s, res}({\frak g})$ then according to this definition $\nu(\pi_V(x))=\pi_{\overline{V}}(px)$.  

Let 
\begin{equation}\label{imathdef}
{\tilde L^{\pm,V}}=(p\otimes \nu)({^q{L^{\pm,V}}})\in
\mathbb{C}^s[G^*]'\otimes {\rm End}_{\mathbb{C}}\overline{V}. \index[not]{L@${\tilde L^{\pm,V}}$}
\end{equation}
Then the map
$$
\imath :\mathbb{C}^s[G^*]'\rightarrow \mathbb{C}^s[G^*] \index[not]{i@$\imath$}
$$
defined by
\begin{equation}\label{imath}
(\imath \otimes id){\tilde L^{\pm,V}}={L^{\pm,\overline{V}}} \index[not]{L@${L^{\pm,\overline{V}}}$}
\end{equation}
is a well--defined algebra isomorphism.
Indeed, consider, for instance, the element ${\tilde L^{+,V}}$.
From (\ref{rmatrspi}) it follows that
\begin{equation}\label{lv1}
\begin{array}{l}
\tilde L^{+,V}=
\prod_{\beta\in \Delta_+}
{\exp}[p(\tilde{f}_{\beta}) \otimes
\pi_{\overline{V}}(X_{\beta})] \times \\
\times (p\otimes id)\exp\left[ \sum_{i=1}^lhH_i\otimes \pi_{\overline{V}}((\kappa {1+s \over 1-s}P_{{\h'}}+id)Y_i)\right].
\end{array}
\end{equation}

On the other hand (\ref{fact}) implies that $L_+$ may be uniquely
represented in the form
\begin{equation}\label{l+dec}
L_+ = 
\prod_{\beta\in \Delta_+}
\exp[b_{\beta}X_{\beta}]\exp\left[ \sum_{i=1}^lb_i(\kappa{1+s \over 1-s}P_{{\h'}}+id)Y_i\right],~b_i,b_\beta\in {\Bbb C},
\end{equation}
and hence
\begin{equation}\label{LVbar}
\begin{array}{l}
L^{+,\overline{V}}=
\prod_{\beta\in \Delta_+}
{\exp}[b_{\beta} \pi_{\overline{V}}(X_{\beta})] {\exp}\left[ \sum_{i=1}^lb_i \pi_{\overline{V}}((\kappa{1+s \over 1-s}P_{{\h'}}+id)Y_i)\right],
\end{array}
\end{equation}
where the order of the terms in the product over the positive roots is the same as in formula (\ref{rmatrspires}) for the R--matrix.

Note that by part (v) $p(\tilde{f}_{\beta})$, $p(\tilde{e}_{\beta})$, $\beta\in \Delta_+$ and $p(\bar{V}_i)$, $i\in \mathbb{N}$ are non--zero, and $\mathbb{C}^s[G^*]'$ is the complex commutative algebra with the linear basis 
$$
p(\tilde e_{\beta_1})^{r_1}\ldots  p(\tilde e_{\beta_D})^{r_D}p(\bar{V}_i) p(\tilde f_{\beta_D})^{m_D}\ldots p(\tilde f_{\beta_1})^{m_1},~r_j,m_j,i\in \mathbb{N},~j=1,\ldots ,D.
$$ 
Therefore comparing (\ref{LVbar}) with (\ref{lv1}), and similar formulas that can be obtained for $L^{-,\overline{V}}$ and $\tilde L^{-,V}$, and recalling the definition of $\imath$, we deduce that  $\imath$ is an algebra isomorphism. We have to prove that $\imath$ is an
isomorphism of Poisson--Hopf algebras.

Observe that $\mathcal{R}_s=1\otimes 1 -2hr_-^{s}$ (mod $h^2$). Therefore from
commutation relations (\ref{ppcomm}), (\ref{pmcomm}) it follows that $\mathbb{C}^s[G^*]'$ is a commutative
algebra, and the Poisson brackets of matrix elements ${\tilde L^{\pm,\overline{V}}}_{ij}$
(see (\ref{quasipb}))
are given by (\ref{pbf}), where $L^{\pm,\overline{V}}$ are replaced by ${\tilde
L^{\pm,\overline{V}}}$. The factor $\frac{1}{d{\bar{r}}^2}$ in formula (\ref{quasipb}) normalizes the Poisson bracket in such a way that bracket (\ref{quasipb}) is in agreement with (\ref{pbf}).

From (\ref{comultGq})
we also obtain that the action of the comultiplication on the matrices ${\tilde
L^{\pm,\overline{V}}}$ is given by
(\ref{comultcl}), where $L^{\pm,\overline{V}}$ are replaced by ${\tilde L^{\pm,\overline{V}}}$.

Part (vii) follows from (iv) and from formulas (\ref{imath}), (\ref{lv1}) and (\ref{LVbar}) in part (vi), and similar formulas relating $\tilde L^{-,V}$ and $L^{-,\overline{V}}$.
This completes the proof.

\end{proof} 


We shall call the map $p:\mathbb{C}_{\mathcal{B}}^s[G^*] \rightarrow \mathbb{C}_{\mathcal{B}}^s[G^*]/(q^{\frac{1}{d{\bar{r}}^2}}-1)\mathbb{C}_{\mathcal{B}}^s[G^*]\simeq\mathbb{C}^s[G^*]$ the {\it quasiclassical limit}.

Next we define the algebra $\mathbb{C}_{\mathcal{B}}^s[G_*]$.
For any finite rank representation $V$ of $U_h({\frak
g})$, let ${^q{L^{V}}}={^q{L^{-,V}}}^{-1}{^q{L^{+,V}}}=(id\otimes \pi_V)(\mathcal{R}_s)_{21}{\mathcal{R}_s}$.
Let $\mathbb{C}_{\mathcal{B}}^s[G_*]$ be the $\mathcal{B}$--subalgebra in $\mathbb{C}_{\mathcal{B}}^s[G^*]$ generated by the matrix entries of ${^q{{L}^{V}}}$, where $V$ runs over all finite rank representations of $U_h({\frak
g})$. From the definition of $\mathcal{R}_s$ we have
\begin{equation}\label{L}
\begin{array}{l}
(\mathcal{R}_s)_{21}\mathcal{R}_s=
\prod_{\beta\in \Delta_+}
\sum_{k=0}^\infty q_{\beta}^{{k(k+1)} \over 2}[(1-q_{\beta}^{-2})e_{\beta}^ke^{-hk\kappa {1+s \over 1-s}P_{{\h'}} \beta^\vee} \otimes
f_{\beta}^{(k)}]\times \\
\times \exp\left[ 2h\sum_{i=1}^lY_i\otimes H_i \right]
 \prod_{\beta\in \Delta_+}
\sum_{k=0}^\infty q_{\beta}^{{k(k+1)} \over 2}[(1-q_{\beta}^{-2})e^{hk\kappa {1+s \over 1-s}P_{{\h'}} \beta^\vee-hk\beta^\vee}f_{\beta}^k \otimes
e_{\beta}^{(k)}q^{k\beta^\vee}].
\end{array}
\end{equation}
Using this formula one immediately checks that actually 
\begin{equation}\label{G*AG*}
\mathbb{C}_{\mathcal{B}}^s[G_*]\subset U_{\mathcal{A}}^s(\g)\cap\mathbb{C}_{\mathcal{B}}^s[G^*].
\end{equation}

Define the right adjoint action ${\rm Ad}_s$ of $U_q^s(\g)$ on $U_{q}^s(\g)$ by formula (\ref{ad})
and the left adjoint action ${\rm Ad}'_s$ of $U_q^s(\g)$ on $U_q^s(\g)$ by formula (\ref{ad'}).

Recall that $\mathbb{C}_{\mathcal{B}}^s[G]$ is the restricted Hopf algebra dual to $U_{\mathcal{B}}^{s, res}(\g)$ which is generated by the matrix elements of finite rank representations of $U_{\mathcal{B}}^{s, res}(\g)$ of the form $V^{res}$, where $V$ is a finite rank representation of $U_h({\frak g})$. 
The action ${\rm Ad}'_s$ induces a right adjoint action ${\rm Ad}_s$ of $U_{\mathcal{B}}^{s,res}(\g)$ on $\mathbb{C}_{\mathcal{B}}^s[G]$  defined by
\begin{equation}\label{ad''}
({\rm Ad}_sxf)(w)=f({\rm Ad}'_sx(w)), f\in \mathbb{C}_{\mathcal{B}}^s[G], x, w \in U_{\mathcal{B}}^{s, res}(\g). \index[not]{A@${\rm Ad}_s$}
\end{equation}

One can also equip finite rank representations of $U_{\mathcal{B}}^{s, res}(\g)$ of the form $V^{res}$, where $V$ is a finite rank representation of $U_h({\frak
g})$, with a natural action of $\mathbb{C}_{\mathcal{B}}^s[G^*]$, for which the elements $$q^{\pm(Y_i-\kappa \frac{1+s}{1-s}P_{\h'}Y_i)}, q^{\pm(-Y_i-\kappa \frac{1+s}{1-s}P_{\h'}Y_i)},~i=1,\ldots l$$ act on a weight vector $v_\lambda$ of weight $\lambda$ by multiplication by the elements 
$$
q^{\pm((Y_i,\lambda^\vee)-\kappa (\frac{1+s}{1-s}P_{\h'}Y_i,\lambda^\vee))}, q^{\pm(-(Y_i,\lambda^\vee)-\kappa (\frac{1+s}{1-s}P_{\h'}Y_i,\lambda^\vee))}\in \mathcal{B},~i=1,\ldots l,
$$ 
respectively, and all the other generators of $\mathbb{C}_{\mathcal{B}}^s[G^*]$ belonging to $U_{\mathcal{B}}^{s, res}(\g)$ act in a natural way. Therefore adjoint action (\ref{ad''}) induces an action of $\mathbb{C}_{\mathcal{B}}^s[G^*]$, where elements $x\in \mathbb{C}_{\mathcal{B}}^s[G^*]$ act by the same formula (\ref{ad''}). 

For $\varepsilon \in \mathbb{C}^*$ we define  ${\mathbb{C}}_{\varepsilon}^s[G]:=\mathbb{C}_{\mathcal{B}}^s[G]/(q^{\frac{1}{d{\bar{r}}^2}}-\varepsilon^{\frac{1}{d{\bar{r}}^2}})\mathbb{C}_{\mathcal{B}}^s[G]$, \index[not]{C@$\mathbb{C}_\varepsilon^s[G]$} $\mathbb{C}_{\varepsilon}^s[G_*]:=\mathbb{C}_{\mathcal{B}}^s[G_*]/(q^{\frac{1}{d{\bar{r}}^2}}-\varepsilon^{\frac{1}{d{\bar{r}}^2}})\mathbb{C}_{\mathcal{B}}^s[G_*]$, \index[not]{C@$\mathbb{C}_\varepsilon^s[G_*]$} where $\varepsilon^{\frac{1}{d{\bar{r}}^2}}$ is a root of $\varepsilon$ of degree $\frac{1}{d{\bar{r}}^2}$.

Recall that by Proposition \ref{CG1} $\mathbb{C}_{\mathcal{B}}^s[G]/(q^{\frac{1}{d{\bar{r}}^2}}-1)\mathbb{C}_{\mathcal{B}}^s[G] \simeq {\mathbb{C}}[G]$, and
denote the canonical projections $\mathbb{C}_{\mathcal{B}}^s[G]\to \mathbb{C}_{\mathcal{B}}^s[G]/(q^{\frac{1}{d{\bar{r}}^2}}-1)\mathbb{C}_{\mathcal{B}}^s[G] \simeq {\mathbb{C}}[G]$, $\mathbb{C}_{\mathcal{B}}^s[G_*]\to \mathbb{C}_{\mathcal{B}}^s[G_*]/(q^{\frac{1}{d{\bar{r}}^2}}-1)\mathbb{C}_{\mathcal{B}}^s[G_*]:= {\mathbb{C}}^s[G_*]$ \index[not]{C@$\mathbb{C}^s[G_*]$} by the same symbol $p$. \index[not]{p@$p(~\cdot~)$}

\begin{proposition}\label{locfin}
(i) The map 
\begin{equation}\label{gg*}
\phi_{\mathcal{B}}:\mathbb{C}_{\mathcal{B}}^s[G]\rightarrow \mathbb{C}_{\mathcal{B}}^s[G_*], \phi_{\mathcal{B}}(f)= (id\otimes f) ((\mathcal{R}_s)_{21}\mathcal{R}_s) \index[not]{f@$\phi_{\mathcal{B}}(~\cdot~)$}
\end{equation}
is an isomorphism of $U_{\mathcal{B}}^{s,res}(\g)$ and $\mathbb{C}_{\mathcal{B}}^s[G^*]$--modules with respect to the adjoint actions ${\rm Ad}_s$ defined by (\ref{ad''}) and (\ref{ad}), respectively. In particular, $\mathbb{C}_{\mathcal{B}}^s[G_*]$ is stable under the adjoint action of $U_{\mathcal{B}}^{s,res}(\g)$ and $\mathbb{C}_{\mathcal{B}}^s[G^*]$.

(ii) Let $G_*\subset G$ \index[not]{G@$G_*$} be the image of $G^*\subset G\times G$ under the map $q:G^*\to G$, $q (L_+,L_-)= L_-^{-1}L_+$ (see (\ref{fact}) \index[not]{q@$q(~\cdot~)$} for the description of $G^*$ in terms of $G\times G$). Then $G_*=B_-B_+$ is the big Bruhat cell \index{Bruhat!cell!big} in $G$, the algebra ${\mathbb{C}}^s[G_*]$ is generated by the restrictions of elements of ${\mathbb{C}}[G]$ to $G_*=B_-B_+$, the restriction map induces an isomorphism of algebras $\overline{\phi}:{\mathbb{C}}[G]\to {\mathbb{C}}^s[G_*]$, \index[not]{f@$\overline{\phi}$} and $\overline{\phi} p=p\phi_{\mathcal{B}}$. 

(iii) If $\varepsilon$ is not a root of unity the algebra $\mathbb{C}_{\varepsilon}^s[G_*]$ can be identified with  the ${\rm Ad}_s$--locally finite part $U_{\varepsilon}^s(\g)^{fin}$ \index[not]{U@$U_\varepsilon^s(\g)^{fin}$} of $U_{\varepsilon}^s(\g)$, \index{quantum group!locally finite part}
$$
U_{\varepsilon}^s(\g)^{fin}=\{x\in U_{\varepsilon}^s(\g):{\rm dim}({\rm Ad}_s U_{\varepsilon}^s(\g)(x))<+\infty \},
$$
where the adjoint action of the algebra $U_{\varepsilon}^s(\g)$ on itself is defined by formula (\ref{ad}).

\end{proposition}

\begin{proof}
Firstly, we prove (i) and (ii). From the definitions of the algebras $\mathbb{C}_{\mathcal{B}}^s[G]$ and $\mathbb{C}_{\mathcal{B}}^s[G_*]$ it follows that (\ref{gg*}) is surjective. 

Recall that by Proposition \ref{CG1} we have the following algebra isomorphism ${\mathbb{C}}_{1}[G]:=\mathbb{C}_{\mathcal{B}}^s[G]/(q^{\frac{1}{d{\bar{r}}^2}}-1)\mathbb{C}_{\mathcal{B}}^s[G] \simeq {\mathbb{C}}[G]$. 

By (\ref{fact}) $G_*=B_-B_+$ is the big Bruhat cell in $G$, and by the definition the algebra ${\mathbb{C}}^s[G_*]={\mathbb{C}}_{1}^s[G_*]=\mathbb{C}_{\mathcal{B}}^s[G_*]/(q^{\frac{1}{d{\bar{r}}^2}}-1)\mathbb{C}_{\mathcal{B}}^s[G_*]$ is generated by the restrictions of elements of ${\mathbb{C}}[G]\simeq \mathbb{C}_{\mathcal{B}}^s[G]/(q^{\frac{1}{d{\bar{r}}^2}}-1)\mathbb{C}_{\mathcal{B}}^s[G]$ to $G_*=B_-B_+$. As the big Bruhat cell is dense in $G$, we have an isomorphism of algebras $\overline{\phi}:{\mathbb{C}}[G]\to {\mathbb{C}}^s[G_*]$ induced by the restriction map. 
 
Now observe that by the definition of map (\ref{gg*}) one has $\overline{\phi} p=p\phi_{\mathcal{B}}$. 

Let $h\in \mathbb{C}_{\mathcal{B}}^s[G]$, $h\neq 0$ be such that $\phi_{\mathcal{B}}(h)=0$. Dividing by an appropriate power of $(q^{\frac{1}{d{\bar{r}}^2}}-1)$ we can assume without loss of generality that $p(h)\neq 0$. Then $\overline{\phi} p(h)=p\phi_{\mathcal{B}}(h)=0$ which implies $p(h)=0$ as $\overline{\phi}$ is an algebra isomorphism. Thus we arrive at a contradiction, and hence $\phi_{\mathcal{B}}$ is injective.
We conclude that $\phi_{\mathcal{B}}$ an isomorphism of $\mathcal{B}$--modules.

By Proposition 11(ii) in \cite{KlSm}, Section 10.1.3, (\ref{gg*}) is a morphism of $U_{\mathcal{B}}^{s,res}(\g)$ and $\mathbb{C}_{\mathcal{B}}^s[G^*]$--modules with respect to the adjoint actions ${\rm Ad}_s$ defined by (\ref{ad}) and (\ref{ad''}), respectively. 


(iii) It remains to establish the isomorphism $\mathbb{C}_{\varepsilon}^s[G_*]=U_{\varepsilon}^s(\g)^{fin}$ which can be done similarly to the proof of Proposition 1.7 in \cite{KolSt}. Since this statement is not used later in this book, we only briefly outline its proof.

Let $V_{\omega_i}$, $i=1,\ldots, l$ be the finite rank representation of $U_h(\g)$ with the highest weight $\omega_i$, $i=1,\ldots, l$. Below we use the notation introduced in the proof of Proposition \ref{qG*} (i).  From formula (\ref{L}) and from the definition of ${^q{L^{V}}}=(id\otimes \pi_V)(\mathcal{R}_s)_{21}\mathcal{R}_s$ it follows that the matrix element $(id\otimes g_i)(\mathcal{R}_s)_{21}\mathcal{R}_s$ of $^q{L^{V_{\omega_i}}}$, where $g_i(~\cdot~)=(v_{\omega_i},~\cdot~ v_{\omega_i})$, and $v_{\omega_i}$ is the highest weight vector of $V_{\omega_i}^{res}$ normalized in such a way that $(v_{\omega_i},v_{\omega_i})=1$, coincides with $L_i^2$. This implies that $L_i^2$ are elements of the algebra ${\mathbb{C}}_\varepsilon[G_*]\subset U_{\varepsilon}^s(\g)$ as well. 

Denote by $\mathfrak{H}\subset \mathbb{C}_{\varepsilon}^s[G_*]\subset U_{\varepsilon}^s(\g)$ the subalgebra generated by the elements $L_i^2\in \mathbb{C}_{\varepsilon}^s[G_*]$, $i=1,\ldots, l$. Similarly to Theorem 7.1.6 and Lemma 7.1.16 in \cite{Jos} one can obtain that $U_{\varepsilon}^s(\g)^{fin}={\rm Ad}_s U_{\varepsilon}^s(\g)\mathfrak{H}$. Since $\mathbb{C}_{\varepsilon}^s[G_*]$ is stable under the adjoint action we have an inclusion, $U_{\varepsilon}^s(\g)^{fin}\subset \mathbb{C}_{\varepsilon}^s[G_*]$. 

On the other hand the adjoint action of $U_{\varepsilon}^s(\g)$ on $\mathbb{C}_{\varepsilon}^s[G]$ is locally finite by the definition of $\mathbb{C}_{\varepsilon}^s[G]$ and of the adjoint action. 
Using isomorphism (\ref{gg*}) we deduce that the adjoint action of $U_{\varepsilon}^s(\g)$ on $\mathbb{C}_{\varepsilon}^s[G_*]$ is locally finite as well. Hence $\mathbb{C}_{\varepsilon}^s[G_*]\subset U_{\varepsilon}^s(\g)^{fin}$, and we deduce that $\mathbb{C}_{\varepsilon}^s[G_*]=U_{\varepsilon}^s(\g)^{fin}$. This completes the proof.

\end{proof}

We shall need the following property of the algebras $\mathbb{C}_{\mathcal{B}}^s[G^*]$, $\mathbb{C}_{\mathcal{B}}^s[G_*]$ and $\mathbb{C}_{\mathcal{B}}^s[G]$.

\begin{proposition}\label{Afree}
The subalgebra $\mathbb{C}_{\mathcal{B}}^s[G_*]\subset \mathbb{C}_{\mathcal{B}}^s[G^*]$, and the algebra $\mathbb{C}_{\mathcal{B}}^s[G]$ are free $\mathcal{B}$--modules.
\end{proposition}

\begin{proof}
Since $\mathcal{B}$ is a principal ideal domain, and the algebra $\mathbb{C}_{\mathcal{B}}^s[G^*]$ is $\mathcal{B}$--free by part (v) of Proposition \ref{qG*}, the subalgebra $\mathbb{C}_{\mathcal{B}}^s[G_*]\subset \mathbb{C}_{\mathcal{B}}^s[G^*]$ is also $\mathcal{B}$--free by Theorem 6.5 in \cite{R}, and isomorphism (\ref{gg*}) implies that $\mathbb{C}_{\mathcal{B}}^s[G]$ is $\mathcal{B}$--free. 

\end{proof}



\section{The definition of q-W--algebras}\label{qWalgpois}

\pagestyle{myheadings}
\markboth{CHAPTER~\thechapter.~Q-W--ALGEBRAS}{\thesection.~THE DEFINITION OF Q-W--ALGEBRAS}

\setcounter{equation}{0}
\setcounter{theorem}{0}

In this section we introduce the main object of this book, q-W--algebras.

Suppose that a positive root system $\Delta_+$ and its ordering are associated to $s$ as in Definition \ref{circorddef}.
Denote by $\mathbb{C}_{\mathcal{B}}^s[M_+]$ \index[not]{C@$\mathbb{C}_\mathcal{B}^s[M_+]$} the subalgebra in $\mathbb{C}_{\mathcal{B}}^s[N_+]$
generated by the elements
$\tilde f_{\beta}$, $\beta\in {\Delta_{\m_+}}$.

By Lemma \ref{segmsub} the vector subspace of $\g$ generated by the root vectors $X_{\alpha}$ (resp. $X_{-\alpha}$), $\alpha\in \Delta_{\m_+}$ is in fact a Lie subalgebra ${\m_+}\subset \g$ (resp. ${\m_-}\subset \g$). By definition $\Delta_{\m_+} \subset \Delta_+$, and hence ${\m_\pm}\subset \n_\pm$. 

Note that one can consider $\m_\pm$ as Lie subalgebras in $\g^*$ via embeddings
$$
\m_+\rightarrow \g^*\subset \g\oplus \g,~x\mapsto (x,0),
$$
$$
\m_-\rightarrow \g^*\subset \g\oplus \g,~x\mapsto (0,x),
$$
where $\g^*$ is regarded as a Lie subalgebra of $\g\oplus \g$ using the embedding (\ref{imbd}).

Using these embeddings the algebraic subgroups $M_\pm\subset G$ \index[not]{M@$M_\pm$} corresponding to the algebraic Lie subalgebras $\m_\pm\subset \g$ can be regarded as Lie subgroups in $G^*$ corresponding to the Lie subalgebras $\m_\pm\subset \g^*$. Note that by the definition $M_\pm\subset N_\pm\subset G^*$, where the embeddings $N_\pm\subset G^*$ correspond to the Lie algebra embeddings (\ref{npmg*inc1}) and (\ref{npmg*inc2}).

The following proposition gives the most important property of the subalgebra $\mathbb{C}_{\mathcal{B}}^s[M_+]$ which plays the key role in the definition of q-W--algebras.
\begin{proposition}\label{Qdefr}
Denote the roots in $\Delta_+$ ordered as in (\ref{NO}) by $\beta_1,\ldots, \beta_D$, $\beta_1<\ldots< \beta_D$, and the roots in the segment $\Delta_{\m_+}$ by $\beta_1,\ldots ,\beta_c$, $\beta_1<\ldots <\beta_c$, so that $\Delta_{\m_+}=\{\beta_1,\ldots, \beta_c\}\subset \Delta_+=\{\beta_1,\ldots, \beta_D\}$. Then the following statements are true. \index[not]{b@$\beta_1,\ldots ,\beta_c$}

(i) The defining relations in the subalgebra $\mathbb{C}_{\mathcal{B}}^s[M_+]$ for the generators $\tilde{f}_\beta=(1-q_\beta^{-2})f_\beta$, $\beta\in \Delta_{\m_+}$  are of the form
\begin{equation}\label{cmrel}
\tilde{f}_{\alpha}\tilde{f}_{\beta} - q^{\left\langle \alpha,\beta\right\rangle+\kappa\left\langle {1+s \over 1-s}P_{\h'^*}\alpha,\beta\right\rangle}\tilde{f}_{\beta}\tilde{f}_{\alpha}= \sum_{r_1,\ldots,r_k\in \mathbb{N}}C(r_1,\ldots,r_k)
{\tilde{f}_{\zeta_1}}^{r_1}{\tilde{f}_{\zeta_2}}^{r_2}\ldots {\tilde{f}_{\zeta_k}}^{r_k},~\alpha<\beta,
\end{equation}
where $\alpha<\zeta_1<\ldots<\zeta_k<\beta$, $[\alpha,\beta]=\{\alpha,\zeta_1,\ldots,\zeta_k,\beta\}$ as a set, $C(r_1,\ldots,r_k)\in \mathcal{B}$.

(ii) The products ${\tilde{f}_{\beta_1}}^{r_1}{\tilde{f}_{\beta_2}}^{r_2}\ldots {\tilde{f}_{\beta_c}}^{r_c}$, $r_1,\ldots,r_c\in \mathbb{N}$ form a $\mathcal{B}$--basis of $\mathbb{C}_{\mathcal{B}}^s[M_+]$. 

(iii) If $\kappa=1$ then for any $\bar{k}_i\in \mathcal{B}$, $i=1,\ldots ,l'$ \index[not]{k@$\bar{k}_i$} the map $\chi_q^{s}:\mathbb{C}_{\mathcal{B}}^s[M_+]\rightarrow \mathcal{B}$ \index[not]{x@$\chi_q^s$} defined by
\begin{equation}\label{charq}
\chi_q^s(\tilde f_\beta)=\left\{ \begin{array}{ll} 0 & \beta \not \in \{\gamma_1, \ldots, \gamma_{l'}\} \\ \bar{k}_i & \beta=\gamma_i 
\end{array}
\right  .
\end{equation}
is a character of $\mathbb{C}_{\mathcal{B}}^s[M_+]$.

(iv) Assume that $\varepsilon\in \mathbb{C}^*$ is such that $\varepsilon^{2d_i}\neq 1$, and $\varepsilon^4\neq 1$ if $\g$ is of type $G_2$. Then the algebra $\mathbb{C}_{\varepsilon}^s[M_+]:=\mathbb{C}_{\mathcal{B}}^s[M_+]/(q^{\frac{1}{d{\bar{r}}^2}}-\varepsilon^{\frac{1}{d{\bar{r}}^2}})\mathbb{C}_{\mathcal{B}}^s[M_+]$, \index[not]{C@$\mathbb{C}_\varepsilon^s[M_+]$} where $\varepsilon^{\frac{1}{d{\bar{r}}^2}}$ is a root of $\varepsilon$ of degree $\frac{1}{d{\bar{r}}^2}$, is isomorphic to the subalgebra $U_\varepsilon^{s}({\frak m}_-)\subset U_\varepsilon^{s}(\g)$ generated by the elements $f_{\beta_1},\ldots, f_{\beta_c}$. \index[not]{U@$U_\varepsilon^s(\m_-)$}

(v) Under the assumptions of part (iv), the elements
$f^{\bf r}:=f_{\beta_1}^{r_1}\ldots f_{\beta_c}^{r_c}$, $r_i\in \mathbb{N}$, $i=1,\ldots d$, ${\bf r}:=(r_1,\ldots, r_c)$ form a linear basis of the algebra $U_\varepsilon^{s}({\frak m}_-)$.

(vi) In addition to the assumptions of part (iv), suppose that there exists $\bar{n}\in \mathbb{Z}$ such that $\varepsilon^{\bar{n}d-1}=1$. Let $\kappa=\bar{n}d$. Then for any $\bar{c}_i\in \mathbb{C}$, $i=1,\ldots ,l'$ \index[not]{c@$\bar{c}_i$} the map $\chi_\varepsilon^{s}:U_{\varepsilon}^{s}({\m_-})\rightarrow \mathbb{C}$ \index[not]{x@$\chi_\varepsilon^s$} defined by
\begin{equation}\label{charqfe}
\chi_\varepsilon^s( f_\beta)=\left\{ \begin{array}{ll} 0 & \beta \not \in \{\gamma_1, \ldots, \gamma_{l'}\} \\ \bar{c}_i & \beta=\gamma_i
\end{array}
\right  .
\end{equation}
is a character of $U_{\varepsilon}^{s}({\m_-})$.

(vii) The algebra $\mathbb{C}_{\mathcal{B}}^s[M_+]$ is a deformation of the algebra of regular functions on the subgroup $M_+\subset G^*$ in the sense that $p(\mathbb{C}_{\mathcal{B}}^s[M_+])\simeq {\mathbb{C}}[M_+]$. \index[not]{C@$\mathbb{C}[M_+]$} Thus ${\mathbb{C}}[M_+]\simeq \mathbb{C}_{\mathcal{B}}^s[M_+]/(q^{\frac{1}{d{\bar{r}}^2}}-1)\mathbb{C}_{\mathcal{B}}^s[M_+]$ is naturally a Poisson subalgebra of $\mathbb{C}^s[G^*]\simeq \mathbb{C}_{\mathcal{B}}^s[G^*]/(q^{\frac{1}{d{\bar{r}}^2}}-1)\mathbb{C}_{\mathcal{B}}^s[G^*]$.
\end{proposition}

\begin{proof}
(ii) By Lemma \ref{segmPBWs} (vi) and Remark \ref{segmPBWsrev} any element of $\mathbb{C}_{\mathcal{B}}^s[M_+]$ can be uniquely represented as a $\mathbb{C}(q^{\frac{1}{d{\bar{r}}^2}})$--linear combination of the elements ${\tilde{f}_{\beta_1}}^{r_1}{\tilde{f}_{\beta_2}}^{r_2}\ldots {\tilde{f}_{\beta_c}}^{r_c}$. By Proposition \ref{qG*} (iv) the coefficients of this decomposition must belong to $\mathcal{B}$. 

(i) From (\ref{cmrelf}) we also obtain commutation relations (\ref{cmrel}) with $C(r_1,\ldots,r_k)\in \mathbb{C}(q^{\frac{1}{d{\bar{r}}^2}})$. As we already proved the products ${\tilde{f}_{\beta_1}}^{r_1}{\tilde{f}_{\beta_2}}^{r_2}\ldots {\tilde{f}_{\beta_c}}^{r_c}$ form a $\mathcal{B}$--basis of $\mathbb{C}_{\mathcal{B}}^s[M_+]$. Therefore the coefficients $C(r_1,\ldots,r_k)$ in (\ref{cmrel}) belong to $\mathcal{B}$. 

(iii) Assume that $\kappa=1$. In order to prove that the map $\chi_q^{s}:\mathbb{C}_{\mathcal{B}}^s[M_+]\rightarrow \mathcal{B}$ defined by (\ref{charq}) is a character we show that all relations (\ref{cmrel}) for $\tilde{f}_\alpha, \tilde{f}_\beta$ with $\alpha,\beta \in \Delta_{\m_+}$, which are defining relations in the subalgebra $\mathbb{C}_{\mathcal{B}}^s[M_+]$ by part (ii),  belong to the kernel of $\chi_q^s$. By the definition the only generators of $\mathbb{C}_{\mathcal{B}}^s[M_+]$ on which $\chi_q^s$ may not vanish are $\tilde{f}_{\gamma_i}$, $i=1,\ldots,l'$. By part (vi) of Proposition \ref{pord} for any two roots $\alpha, \beta\in \Delta_{\m_+}$ such that $\alpha<\beta$ the sum $\alpha+\beta$ cannot be represented as a linear combination $\sum_{k=1}^jc_k\gamma_{i_k}$, where $c_k\in \mathbb{N}$ and $\alpha<\gamma_{i_1}<\ldots <\gamma_{i_j}<\beta$. Hence for any two roots $\alpha, \beta\in \Delta_{\m_+}$ such that $\alpha<\beta$ the value of the map $\chi_q^s$ on the right hand side of the corresponding commutation relation (\ref{cmrel}) is equal to zero.

Therefore it suffices to prove that
$$
\chi_q^s(\tilde{f}_{\gamma_i}\tilde{f}_{\gamma_j} - q^{\left\langle \gamma_i,\gamma_j\right\rangle+\left\langle {1+s \over 1-s}P_{{\h'}^*}\gamma_i,\gamma_j\right\rangle}\tilde{f}_{\gamma_j}\tilde{f}_{\gamma_j})=\bar{k}_i\bar{k}_j(1-q^{\left\langle \gamma_i,\gamma_j\right\rangle+\left\langle {1+s \over 1-s}P_{{\h'}^*}\gamma_i,\gamma_j\right\rangle})=0,~i<j.
$$
The last identity holds provided $\left\langle \gamma_i,\gamma_j\right\rangle+\left\langle {1+s \over 1-s}P_{{\h'}^*}\gamma_i, \gamma_j\right\rangle=0$ for $i<j$ which is indeed the case by Lemma \ref{tmatrel}.

(iv) Assume now that $\varepsilon^{2d_i}\neq 1$, $i=1,\ldots, l$ and $\varepsilon^4\neq 1$ if $\g$ is of type $G_2$. Under these conditions imposed on $\varepsilon$ the map $\mathbb{C}_{\varepsilon}^s[M_+]\to U_{\varepsilon}^{s}({\m_-})$, $\tilde{f}_\alpha\mapsto (1-\varepsilon_\alpha^{-2})f_\alpha$, $\alpha\in \Delta_{\m_+}$ is obviously an algebra isomorphism. 

(v) By Lemma \ref{segmPBWs} (vi) and Remark \ref{segmPBWsrev} the elements
$f^{\bf r}=f_{\beta_1}^{r_1}\ldots f_{\beta_c}^{r_c}$, $r_i\in \mathbb{N}$, $i=1,\ldots d$ form a linear basis of $U_\varepsilon^{s}({\frak m}_-)$. 

(vi) Assume now that $\varepsilon^{2d_i}\neq 1$, $i=1,\ldots, l$, and $\varepsilon^4\neq 1$ if $\g$ is of type $G_2$. Suppose also that there exists $\bar{n}\in \mathbb{Z}$ such that $\varepsilon^{\bar{n}d-1}=1$. Let $\kappa=\bar{n}d$.
From (\ref{cmrelf}) we obtain the following commutation relations
\begin{equation}\label{cmrele}
f_{\alpha}f_{\beta} - \varepsilon^{\left\langle \alpha,\beta\right\rangle+\bar{n}d\left\langle {1+s \over 1-s}P_{\h'^*}\alpha,\beta\right\rangle}f_{\beta}f_{\alpha}= \sum_{r_1,\ldots,r_k\in \mathbb{N}}D(r_1,\ldots,r_k)
f_{\zeta_1}^{r_1}f_{\zeta_2}^{r_2}\ldots f_{\zeta_k}^{r_k},~\alpha<\beta,
\end{equation}
where $\alpha,\beta\in \Delta_{\m_+}$, $\alpha<\zeta_1<\ldots<\zeta_k<\beta$, $[\alpha,\beta]=\{\alpha,\zeta_1,\ldots,\zeta_k,\beta\}$ as a set, $D(r_1,\ldots,r_k)\in \mathbb{C}$. Relations (\ref{cmrele}) are defining relations in $U_{\varepsilon}^{s}({\m_-})$ by part (v).

In order to show that the map $\chi_\varepsilon^{s}:U_{\varepsilon}^{s}({\m_-})\rightarrow \mathbb{C}$ defined by (\ref{charqfe}) is a character we verify that all relations (\ref{cmrele}) for $f_\alpha, f_\beta$ with $\alpha,\beta \in \Delta_{\m_+}$  belong to the kernel of $\chi_\varepsilon^s$. By the definition the only generators of $U_\varepsilon^{s}({\frak m}_-)$ on which $\chi_\varepsilon^s$ may not vanish are $f_{\gamma_i}$, $i=1,\ldots,l'$. By part (vi) of Proposition \ref{pord} for any two roots $\alpha, \beta\in \Delta_{\m_+}$ such that $\alpha<\beta$ the sum $\alpha+\beta$ cannot be represented as a linear combination $\sum_{k=1}^jc_k\gamma_{i_k}$, where $c_k\in \mathbb{N}$ and $\alpha<\gamma_{i_1}<\ldots <\gamma_{i_j}<\beta$. Hence for any two roots $\alpha, \beta\in \Delta_{\m_+}$ such that $\alpha<\beta$ the value of the map $\chi_\varepsilon^s$ on the right hand side of the corresponding commutation relation (\ref{cmrel}) is equal to zero.

Therefore it suffices to prove that
$$
\chi_\varepsilon^s(f_{\gamma_i}f_{\gamma_j} - \varepsilon^{\left\langle \gamma_i,\gamma_j\right\rangle+\bar{n}d\left\langle {1+s \over 1-s}P_{{\h'}^*}\gamma_i,\gamma_j\right\rangle}f_{\gamma_j}f_{\gamma_j})=\bar{c}_i\bar{c}_j(1-\varepsilon^{\left\langle \gamma_i,\gamma_j\right\rangle+\bar{n}d\left\langle {1+s \over 1-s}P_{{\h'}^*}\gamma_i,\gamma_j\right\rangle})=0,~i<j.
$$
By Lemma \ref{tmatrel} $\left\langle {1+s \over 1-s}P_{{\h'}^*}\gamma_i, \gamma_j\right\rangle=-\left\langle \gamma_i,\gamma_j\right\rangle$ for $i<j$, and hence
$$
\chi_\varepsilon^s(f_{\gamma_i}f_{\gamma_j} - \varepsilon^{\left\langle \gamma_i,\gamma_j\right\rangle+\bar{n}d\left\langle {1+s \over 1-s}P_{{\h'}^*}\gamma_i,\gamma_j\right\rangle}f_{\gamma_j}f_{\gamma_j})=\bar{c}_i\bar{c}_j(1-\varepsilon^{\left\langle \gamma_i,\gamma_j\right\rangle(1-\bar{n}d)})=0
$$
for $i<j$ as by the assumption $\varepsilon^{\bar{n}d-1}=1$. 

Part (vii) follows from part (ii) and from Proposition \ref{qG*} (vi), (vii) (see especially formulas (\ref{imath}), (\ref{lv1}) and (\ref{LVbar}) in the proof of part (vi)of Proposition \ref{qG*}). This completes the proof.

\end{proof} 

Now we are ready to define q-W--algebras. In the rest of this section we assume that $\kappa=1$. Denote by $I_{\mathcal{B}}$ \index[not]{I@$I_{\mathcal{B}}$} the left ideal in $\mathbb{C}_{\mathcal{B}}^s[G^*]$ generated by the kernel of $\chi_q^s$, and by $\rho_{\chi^{s}_q}$ \index[not]{r@$\rho_{\chi^{s}_q}$} the canonical projection $\mathbb{C}_{\mathcal{B}}^s[G^*]\rightarrow \mathbb{C}_{\mathcal{B}}^s[G^*]/I_{\mathcal{B}}:=Q_{\mathcal{B}}'$. \index[not]{Q@$Q_{\mathcal{B}}'$} Let $Q_{\mathcal{B}}$ \index[not]{Q@$Q_{\mathcal{B}}$} be the image of $\mathbb{C}_{\mathcal{B}}^s[G_*]\subset \mathbb{C}_{\mathcal{B}}^s[G^*]$ under the projection $\rho_{\chi^{s}_q}$, $Q_{\mathcal{B}}=\rho_{\chi^{s}_q}(\mathbb{C}_{\mathcal{B}}^s[G_*])$.

We shall need the following properties of $Q_{\mathcal{B}}'$ and $Q_{\mathcal{B}}$.

\begin{proposition}\label{Afree1}
$Q_{\mathcal{B}}'$ and $Q_{\mathcal{B}}$ are free $\mathcal{B}$--modules.
\end{proposition}

\begin{proof}

Using the $\mathcal{B}$--basis $$\tilde e_{\beta_1}^{m_1}\ldots \tilde e_{\beta_D}^{m_D}\bar{V}_i \tilde f_{\beta_D}^{r_D}\ldots \tilde f_{\beta_1}^{r_1}$$ with $m_j,r_j,i\in \mathbb{N}$, $j=1,\ldots ,D$ of $\mathbb{C}_{\mathcal{B}}^s[G^*]$ from Proposition \ref{qG*} (v) and the definition of $\mathbb{C}_{\mathcal{B}}^s[G^*]/I_{\mathcal{B}}$ one immediately sees that the classes of the elements $\tilde e_{\beta_1}^{m_1}\ldots \tilde e_{\beta_D}^{m_D}\bar{V}_i \tilde f_{\beta_D}^{r_D}\ldots \tilde f_{\beta_{c+1}}^{r_{c+1}}$ with $m_j,r_m,i\in \mathbb{N}$, $j=1,\ldots ,D$, $m=c+1,\ldots, D$ form a $\mathcal{B}$--basis in $\mathbb{C}_{\mathcal{B}}^s[G^*]/I_{\mathcal{B}}$. Here we use the notation introduced in Proposition \ref{Qdefr}, $\Delta_{\m_+}=\{\beta_1,\ldots, \beta_c\}\subset \Delta_+=\{\beta_1,\ldots, \beta_D\}$.

Since $\mathcal{B}$ is a principal ideal domain the $\mathcal{B}$--submodule $Q_{\mathcal{B}}\subset \mathbb{C}_{\mathcal{B}}^s[G^*]/I_{\mathcal{B}}$ is $\mathcal{B}$--free.

\end{proof}

\begin{lemma}\label{AdM}
The restriction of the adjoint action of $\mathbb{C}_{\mathcal{B}}^s[G^*]$ on $\mathbb{C}_{\mathcal{B}}^s[G^*]$ to $\mathbb{C}_{\mathcal{B}}^s[M_+]$ induces an action on $Q_{\mathcal{B}}'$ and on $Q_{\mathcal{B}}$.
\end{lemma}
\begin{proof}
Observe that since $\chi_q^{s}$ is a character we have an inclusion $[\mathbb{C}_{\mathcal{B}}^s[M_+],{\rm Ker } \chi_q^{s}]\subset {\rm Ker} \chi_q^{s}$. Note also that $\Delta_s(\mathbb{C}_{\mathcal{B}}^s[M_+])\subset  \mathbb{C}_{\mathcal{B}}^s[B_+]\otimes   \mathbb{C}_{\mathcal{B}}^s[M_+]$ by (\ref{comults}). Now from  formula (\ref{ad}) which defines the adjoint action ${\rm Ad}_s$ of $\mathbb{C}_{\mathcal{B}}^s[G^*]$ on itself it follows that the left ideal $I_{\mathcal{B}}$ is invariant under the restriction of the adjoint action of $\mathbb{C}_{\mathcal{B}}^s[G^*]$ on $\mathbb{C}_{\mathcal{B}}^s[G^*]$ to $\mathbb{C}_{\mathcal{B}}^s[M_+]$. This implies that the restriction of the adjoint action of $\mathbb{C}_{\mathcal{B}}^s[G^*]$ on itself to $\mathbb{C}_{\mathcal{B}}^s[M_+]$ naturally induces an action on $Q_{\mathcal{B}}'$, and hence on $Q_{\mathcal{B}}$ as $\mathbb{C}_{\mathcal{B}}^s[G_*]$ is invariant under the adjoint action of $\mathbb{C}_{\mathcal{B}}^s[G^*]$ by Proposition \ref{locfin} (i).

\end{proof}

We call the action of $\mathbb{C}_{\mathcal{B}}^s[M_+]$ on $Q_{\mathcal{B}}'$ and on $Q_{\mathcal{B}}$ the {\it adjoint action} as well and denote it by ${\rm Ad}_s$. \index[not]{A@${\rm Ad}_s$}

Let $\mathcal{B}_{\varepsilon_s}$ \index[not]{B@$\mathcal{B}_{\varepsilon_s}$} be the trivial representation of $\mathbb{C}_{\mathcal{B}}^s[M_+]$ given by the counit. Consider the $\mathcal{B}$--submodule $W_{\mathcal{B}}^s(G)$ of ${\rm Ad}_s$--invariants in $Q_{\mathcal{B}}$,
\begin{equation}\label{QW}
W_{\mathcal{B}}^s(G)={\rm Hom}_{\mathbb{C}_{\mathcal{B}}^s[M_+]}(\mathcal{B}_{\varepsilon_s},Q_{\mathcal{B}}). \index[not]{W@$W_{\mathcal{B}}^s(G)$}
\end{equation}

\begin{proposition}\label{9.2}
(i) $W_{\mathcal{B}}^s(G)$ is isomorphic to the ${\mathcal{B}}$--submodule of all $v+I_{\mathcal{B}}\in Q_{\mathcal{B}}$ such that $mv\in I_{\mathcal{B}}$ (or $[m,v]\in I_{\mathcal{B}}$) in $\mathbb{C}_{\mathcal{B}}^s[G^*]$ for any $m\in I_{\mathcal{B}}$, where $v\in \mathbb{C}_{\mathcal{B}}^s[G^*]$ is any representative of $v+I_{\mathcal{B}}\in Q_{\mathcal{B}}$.

(ii) Multiplication in $\mathbb{C}_{\mathcal{B}}^s[G^*]$ induces a multiplication on $W_{\mathcal{B}}^s(G)$.
\end{proposition}

\begin{proof}
(i) For the proof we shall firstly derive a formula for the adjoint action of the generators $\tilde{f}_{\beta}$. We use the notation introduced in Proposition \ref{Qdefr}, so that $\Delta_{\m_+}=\{\beta_1,\ldots, \beta_c\}\subset \Delta_+=\{\beta_1,\ldots, \beta_D\}$.

From (\ref{comults}), using linear independence of weight components and the fact that $\mathbb{C}_{\mathcal{B}}^s[B_+]$ is a Hopf algebra by Proposition \ref{qG*} (ii), we obtain 
\begin{equation}\label{comultsG}
\Delta_s(\tilde{f}_{\beta_k})=G_{\beta_k}^{-1} \otimes \tilde{f}_{\beta_k} + \tilde{f}_{\beta_k}\otimes 1+\sum_i \tilde{y}_i\otimes \tilde{x}_i,
\end{equation}
where
$$
G_\beta= e^{h\kappa {1+s \over 1-s}P_{\h'}\beta^\vee-h\beta^\vee}\in \mathbb{C}_{\mathcal{B}}^s[B_+], \index[not]{G@$G_\beta$} \tilde{y}_i=e^{-h\kappa {1+s \over 1-s}P_{\h'}\gamma_{x_i}^\vee+h\gamma_{x_i}^\vee}\overline{\tilde{y}}_i, \index[not]{y@$\tilde{y}_i$}
$$
$$
\overline{\tilde{y}}_i\in \mathbb{C}_\mathcal{B}([-\beta_{k+1},-\beta_{D}]), \index[not]{y@$\overline{\tilde{y}}_i$}
\tilde{x}_i\in \mathbb{C}_\mathcal{B}([-\beta_{1},-\beta_{k-1}]), \index[not]{x@$\tilde{x}_i$}
$$ 
$\overline{\tilde{y}}_i,\tilde{x}_i$ belong to weight subspaces of $U_h^{s}(\g)$  and have non-zero weights, $\gamma_{x_i}$ \index[not]{g@$\gamma_{x_i}$} is the weight of $\tilde{x}_i$, 
$$
\mathbb{C}_\mathcal{B}([-\beta_{k+1},-\beta_{D}]) ~(\text{resp.} ~\mathbb{C}_\mathcal{B}([-\beta_{1},-\beta_{k-1}])) \index[not]{C@$\mathbb{C}_\mathcal{B}([-\beta_{k+1},-\beta_{D}])$} \index[not]{C@$\mathbb{C}_\mathcal{B}([-\beta_{1},-\beta_{k-1}])$}
$$ 
is the subalgebra in $\mathbb{C}_\mathcal{B}^s[N_+]$ generated by $\tilde{f}_{\beta_{k+1}},\ldots, \tilde{f}_{\beta_{D}}$ (resp. $\tilde{f}_{\beta_{1}},\ldots, \tilde{f}_{\beta_{k-1}}$).

By (\ref{comultsG}) we also have
$$
S_s(\tilde{f}_\beta)=-G_\beta \tilde{f}_\beta-\sum_i S_s(\tilde{y}_i)\tilde{x}_i.
$$

Using this formula, (\ref{comultsG}) and the definition of the action ${\rm Ad}_s$ in (\ref{ad}) we deduce
$$
{\rm Ad}_s\tilde{f}_{\beta_k}w=-G_{\beta_k}[\tilde{f}_{\beta_k},w]-\sum_i S_s(\tilde{y}_i)[\tilde{x}_i,w]
$$

The induced action of the elements $\tilde{f}_{\beta_k}\in \mathbb{C}_{\mathcal{B}}^s[M_+]$, $\beta_k\in \Delta_{\m_+}$, on $Q_{\mathcal{B}}'$ takes the form
\begin{equation}\label{adf}
{\rm Ad}_s\tilde{f}_{\beta_k}v=-G_{\beta_k}(\tilde{f}_{\beta_k}-\chi_q^s(\tilde{f}_{\beta_k}))v-\sum_i S_s(\tilde{y}_i)(\tilde{x}_i-\chi_q^s(\tilde{x}_i))v, \beta_k\in \Delta_{\m_+}.
\end{equation}

We have to show that $W_{\mathcal{B}}^s(G)$ is isomorphic to the ${\mathcal{B}}$--submodule of all $v\in Q_{\mathcal{B}}\subset Q_{\mathcal{B}}'$ such that $mv=0$ in $Q_{\mathcal{B}}'$ for any $m\in I_{\mathcal{B}}$. 

The left ideal $I_{\mathcal{B}}$ is generated by the elements $x-\chi_q^s(x)$, $x\in \mathbb{C}_{\mathcal{B}}^s[M_+]$. Therefore by (\ref{adf}) if for some $v\in Q_{\mathcal{B}}$ $mv=0$ in $Q_{\mathcal{B}}'$ for any $m\in I_{\mathcal{B}}$ then $v$ is invariant with respect to the adjoint action of all generators $\tilde{f}_{\beta_k}$ of $\mathbb{C}_{\mathcal{B}}^s[M_+]$, and hence $v\in W_{\mathcal{B}}^s(G)$.

Now assume that $v\in W_{\mathcal{B}}^s(G)$. We shall prove that $mv=0$ in $Q_{\mathcal{B}}'$ for any $m\in I_{\mathcal{B}}$. 

Since the left ideal $I_{\mathcal{B}}$ is generated by elements $x-\chi_q^s(x)$, $x\in \mathbb{C}_{\mathcal{B}}^s[M_+]$ it suffices to show that $(x-\chi_q^s(x))v=0$ for any $x\in \mathbb{C}_{\mathcal{B}}^s[M_+]$.
We shall prove this statement by induction using the subalgebras $\mathbb{C}_\mathcal{B}([-\beta_{1},-\beta_{k}])$, $k=1,\ldots, c$, so that $\mathbb{C}_\mathcal{B}([-\beta_{1},-\beta_{c}])=\mathbb{C}_{\mathcal{B}}^s[M_+]$. 

Observe that $\beta_1$ is a simple root, and hence from (\ref{adf}) we obtain
$$
0={\rm Ad}_s\tilde{f}_{\beta_1}v=-G_{\beta_1}(\tilde{f}_{\beta_1}-\chi_q^s(\tilde{f}_{\beta_1}))v.
$$
Since the element $G_{\beta_1}\in \mathbb{C}_\mathcal{B}^s[B_+]$ is invertible this implies
$$
(\tilde{f}_{\beta_1}-\chi_q^s(\tilde{f}_{\beta_1}))v=0,
$$
i.e. $(x-\chi_q^s(x))v=0$ in $Q_{\mathcal{B}}'$ for any $x\in \mathbb{C}_\mathcal{B}([-\beta_{1},-\beta_1])$ as the subalgebra $\mathbb{C}_\mathcal{B}([-\beta_{1},-\beta_1])$ is generated by $\tilde{f}_{\beta_1}$.

Now suppose that for some $1<k\leq c$ $(x-\chi_q^s(x))v=0$ in $Q_{\mathcal{B}}'$ for any $x\in \mathbb{C}_\mathcal{B}([-\beta_{1},-\beta_{k-1}])$. Then from (\ref{adf}) we obtain
$$
0={\rm Ad}_s\tilde{f}_{\beta_k}v=-G_{\beta_k}(\tilde{f}_{\beta_k}-\chi_q^s(\tilde{f}_{\beta_k}))v-\sum_i S_s(\tilde{y}_i)(\tilde{x}_i-\chi_q^s(\tilde{x}_i))v=-G_{\beta_k}(\tilde{f}_{\beta_k}-\chi_q^s(\tilde{f}_{\beta_k}))v
$$
since $\tilde{x}_i\in \mathbb{C}_\mathcal{B}([-\beta_{1},-\beta_{k-1}])$. The previous identity and the fact that the element $G_{\beta_k}\in \mathbb{C}_\mathcal{B}^s[B_+]$ is invertible yield
$$
(\tilde{f}_{\beta_k}-\chi_q^s(\tilde{f}_{\beta_k}))v=0,
$$
and hence
\begin{equation}\label{fti}
\tilde{f}_{\beta_k}^iv=\chi_q^s(\tilde{f}_{\beta_k}^i)v, i\in \mathbb{N}.
\end{equation}

Now observe that by Lemma \ref{segmPBWs} (vi), Remark \ref{segmPBWsrev} and by Proposition \ref{Qdefr} (ii) any element $x$ of $\mathbb{C}_\mathcal{B}([-\beta_{1},-\beta_{k}])$ can be uniquely represented in the form $x=\sum_{i=0}^r\tilde{f}_{\beta_k}^iz_i$ for some $r\in \mathbb{N}$ and
$z_i\in \mathbb{C}_\mathcal{B}([-\beta_{1},-\beta_{k-1}])$ for $i=0,\ldots, r$. Therefore by (\ref{fti}) and by the induction assumption
$$
xv=\sum_{i=0}^r\tilde{f}_{\beta_k}^iz_iv=\sum_{i=0}^r\tilde{f}_{\beta_k}^i\chi_q^s(z_i)v=\sum_{i=0}^r\chi_q^s(\tilde{f}_{\beta_k}^i)\chi_q^s(z_i)v=\chi_q^s(\sum_{i=0}^r\tilde{f}_{\beta_k}^iz_i)v=\chi_q^s(x)v.
$$
This establishes the induction step and proves the first claim of this proposition.

(ii) From the description of $W_{\mathcal{B}}^s(G)$ obtained in part (i) it follows that if $v_1,v_2\in \mathbb{C}_{\mathcal{B}}^s[G^*]$ are any representatives of elements $v_1+I_{\mathcal{B}},v_2+I_{\mathcal{B}}\in W_{\mathcal{B}}^s(G)$ then the formula
$$
(v_1+I_{\mathcal{B}})(v_2+I_{\mathcal{B}})=v_1v_2+I_q
$$
defines a multiplication in $W_{\mathcal{B}}^s(G)$. This completes the proof.

\end{proof}

The ${\mathcal{B}}$--module $W_{\mathcal{B}}^s(G)$ equipped with the multiplication opposite to the one defined in the previous proposition is called {\it the q-W--algebra} \index{q-W--algebra} associated to (the conjugacy class of) the Weyl group element $s\in W$.

In conclusion we obtain some results on the structure of $Q_{\mathcal{B}}$.
Consider the Lie algebra $\mathfrak{L}_{\mathcal{B}}$ \index[not]{L@$\mathfrak{L}_{\mathcal{B}}$} associated to the associative algebra $\mathbb{C}_{\mathcal{B}}^s[M_+]$, i.e. $\mathfrak{L}_{\mathcal{B}}$ is the Lie algebra which is isomorphic to $\mathbb{C}_{\mathcal{B}}^s[M_+]$ as a ${\mathcal{B}}$--module, and the Lie bracket in $\mathfrak{L}_{\mathcal{B}}$ is given by the usual commutator of elements in $\mathbb{C}_{\mathcal{B}}^s[M_+]$.

Note that since ${\chi_q^s}$ is a character of $\mathbb{C}_{\mathcal{B}}^s[M_+]$ the ideal $I_{\mathcal{B}}$ is stable under the action of $\mathbb{C}_{\mathcal{B}}^s[M_+]$ on $\mathbb{C}_{\mathcal{B}}^s[G^*]$ by commutators.
Therefore one can define an action of the Lie algebra $\mathfrak{L}_{\mathcal{B}}$ on  $Q_{\mathcal{B}}'=\mathbb{C}_{\mathcal{B}}^s[G^*]/I_{\mathcal{B}}$:
\begin{equation}\label{qmainactcl}
m\cdot (x+I_{\mathcal{B}}) =\rho_{\chi^{s}_q}([m,x] ).
\end{equation}
where $x\in \mathbb{C}_{\mathcal{B}}^s[G^*]$ is any representative of $x+I_{\mathcal{B}}\in \mathbb{C}_{\mathcal{B}}^s[G^*]/I_{\mathcal{B}}$ and $m\in \mathbb{C}_{\mathcal{B}}^s[M_+]$.
The algebra $W_{\mathcal{B}}^s(G)$ can be described now as the intersection of the $\mathcal{B}$--module of invariants $(\mathbb{C}_{\mathcal{B}}^s[G^*]/I_{\mathcal{B}})^{\mathfrak{L}_{\mathcal{B}}}$ with respect to action (\ref{qmainactcl}) with $Q_{\mathcal{B}}\subset \mathbb{C}_{\mathcal{B}}^s[G^*]/I_{\mathcal{B}}$,
\begin{equation}\label{WqL}
W_{\mathcal{B}}^s(G)\simeq (\mathbb{C}_{\mathcal{B}}^s[G^*]/I_{\mathcal{B}})^{\mathfrak{L}_{\mathcal{B}}}\cap Q_{\mathcal{B}}.
\end{equation}

Denote by $\mathcal{B}_{\chi_q^{s}}$ \index[not]{B@$\mathcal{B}_{\chi_q^{s}}$} the rank one representation of the algebra $\mathbb{C}_{\mathcal{B}}^s[M_+]$ defined by the character $\chi_q^{s}$. Using the description of the algebra $W_{\mathcal{B}}^s(G)$ in terms of action (\ref{qmainactcl}) and the isomorphism $\mathbb{C}_{\mathcal{B}}^s[G^*]/I_{\mathcal{B}}\simeq
\mathbb{C}_{\mathcal{B}}^s[G^*]\otimes_{\mathbb{C}_{\mathcal{B}}^s[M_+]}\mathcal{B}_{\chi_q^{s}}$ one can also define
the algebra $W_{\mathcal{B}}^s(G)$ as the following intersection
$$
W_{\mathcal{B}}^s(G)={\rm Hom}_{\mathbb{C}_{\mathcal{B}}^s[M_+]}(\mathcal{B}_{\chi_q^{s}},
\mathbb{C}_{\mathcal{B}}^s[G^*]\otimes_{\mathbb{C}_{\mathcal{B}}^s[M_+]}\mathcal{B}_{\chi_q^{s}})\cap Q_{\mathcal{B}}.
$$
Using Frobenius reciprocity we also have
$$
{\rm Hom}_{\mathbb{C}_{\mathcal{B}}^s[M_+]}(\mathcal{B}_{\chi_q^{s}},
\mathbb{C}_{\mathcal{B}}^s[G^*]\otimes_{\mathbb{C}_{\mathcal{B}}^s[M_+]}\mathcal{B}_{\chi_q^{s}})={\rm End}_{\mathbb{C}_{\mathcal{B}}^s[G^*]}(\mathbb{C}_{\mathcal{B}}^s[G^*]\otimes_{\mathbb{C}_{\mathcal{B}}^s[M_+]}\mathcal{B}_{\chi_q^{s}}).
$$
Hence the algebra $W_{\mathcal{B}}^s(G)$ acts on $\mathbb{C}_{\mathcal{B}}^s[G^*]\otimes_{\mathbb{C}_{\mathcal{B}}^s[M_+]}\mathcal{B}_{\chi_q^{s}}$ from the right by operators commuting with the natural left $\mathbb{C}_{\mathcal{B}}^s[G^*]$--action on $\mathbb{C}_{\mathcal{B}}^s[G^*]\otimes_{\mathbb{C}_{\mathcal{B}}^s[M_+]}\mathcal{B}_{\chi_q^{s}}$. By the definition of $W_{\mathcal{B}}^s(G)$ this action preserves $Q_{\mathcal{B}}\subset Q_{\mathcal{B}}'\simeq \mathbb{C}_{\mathcal{B}}^s[G^*]\otimes_{\mathbb{C}_{\mathcal{B}}^s[M_+]}\mathcal{B}_{\chi_q^{s}}\simeq \mathbb{C}_{\mathcal{B}}^s[G^*]/I_{\mathcal{B}}$ and by the  arguments of this paragraph it commutes with the natural left $\mathbb{C}_{\mathcal{B}}^s[G_*]$--action on $Q_{\mathcal{B}}$.

Thus $Q_{\mathcal{B}}$ is a $\mathbb{C}_{\mathcal{B}}^s[G_*]$--$W_{\mathcal{B}}^s(G)$ bimodule equipped also with the adjoint action of $\mathbb{C}_{\mathcal{B}}^s[M_+]$.  By (\ref{adm}) the adjoint action satisfies
\begin{equation}\label{adm4}
{\rm Ad}_s x(yv)={\rm Ad}_sx^1(y){\rm Ad}_s x^2(v),~x\in \mathbb{C}_{\mathcal{B}}^s[M_+],~y\in \mathbb{C}_{\mathcal{B}}^s[G_*],v\in Q_{\mathcal{B}},
\end{equation}
where $\Delta_s(x)=x^1\otimes x^2$.

Denote by $1\in Q_{\mathcal{B}}$ the image of the element $1 \in \mathbb{C}_{\mathcal{B}}^s[G_*]$ in the quotient $Q_{\mathcal{B}}$ under the canonical projection $\mathbb{C}_{\mathcal{B}}^s[G_*]\rightarrow Q_{\mathcal{B}}$. Obviously $1$ is the generating vector for $Q_{\mathcal{B}}$ as a module over $\mathbb{C}_{\mathcal{B}}^s[G_*]$.
Using formula (\ref{adm4}) and recalling that $Q_{\mathcal{B}}$ is a $\mathbb{C}_{\mathcal{B}}^s[G_*]$--$W_{\mathcal{B}}^s(G)$ bimodule, for $x\in \mathbb{C}_{\mathcal{B}}^s[M_+], y\in \mathbb{C}_{\mathcal{B}}^s[G_*]$, and for a representative $w\in \mathbb{C}_{\varepsilon}^s[G_*]$ of an element $w+I_{\mathcal{B}}\in W_{\mathcal{B}}^s(G)$ we have
$$
{\rm Ad}_s x(wy1)={\rm Ad}_sx(yw1)={\rm Ad}_sx^1(y){\rm Ad}_sx^2(w1)=
$$
$$
={\rm Ad}_sx^1(y)\varepsilon_s(x^2)w1={\rm Ad}_sx(y)w1=w{\rm Ad}_sx(y1).
$$
Since $Q_{\mathcal{B}}$ is generated by the vector $1$ over $\mathbb{C}_{\mathcal{B}}^s[G_*]$ the last relation implies that the action of $W_{\mathcal{B}}^s(G)$ on $Q_{\mathcal{B}}$ commutes with the adjoint action.

We can summarize the results of the discussion above in the following proposition.
\begin{proposition}\label{Qpr}
The ${\mathcal{B}}$--module $Q_{\mathcal{B}}$ is naturally equipped with the structure of a left $\mathbb{C}_{\mathcal{B}}^s[G_*]$--module, a right $\mathbb{C}_{\mathcal{B}}^s[M_+]$--module via the adjoint action and a right $W_{\mathcal{B}}^s(G)$--module in such a way that the left $\mathbb{C}_{\mathcal{B}}^s[G_*]$--action and the right $\mathbb{C}_{\mathcal{B}}^s[M_+]$--action commute with the right $W_{\mathcal{B}}^s(G)$--action and compatibility condition (\ref{adm4}) is satisfied.
\end{proposition}

Finally we remark that by specializing $q$ to a particular value $\varepsilon\in \mathbb{C}$, $\varepsilon\neq 0$, one can define a complex associative algebra
${\mathbb{C}}_\varepsilon^s[G_*]=\mathbb{C}_{\mathcal{B}}^s[G_*]/(q^{\frac{1}{d{\bar{r}}^2}}-\varepsilon^{\frac{1}{d{\bar{r}}^2}})
\mathbb{C}_{\mathcal{B}}^s[G_*]$, its subalgebra ${\mathbb{C}}_\varepsilon^s[M_+]$, \index[not]{C@$\mathbb{C}_\varepsilon^s[M_+]$} which is a specialization of ${\mathbb{C}}_{\mathcal{B}}^s[M_+]$ at $q^{\frac{1}{d{\bar{r}}^2}}=\varepsilon^{\frac{1}{d{\bar{r}}^2}}$, with a character $\chi_\varepsilon^{s}$ \index[not]{x@$\chi_\varepsilon^s$} induced by $\chi_q^{s}$, and the corresponding q-W--algebra
\begin{equation}\label{eW}
W_\varepsilon^s(G):={\rm Hom}_{\mathbb{C}_{\varepsilon}^s[M_+]}(\mathbb{C}_{\varepsilon_s},Q_\varepsilon), \index[not]{W@$W_\varepsilon^s(G)$}
\end{equation}
where $\mathbb{C}_{\varepsilon_s}$ \index[not]{C@$\mathbb{C}_{\varepsilon_s}$} is the trivial representation of the algebra $\mathbb{C}_{\varepsilon}^s[M_+]$ induced by the counit, $Q_\varepsilon:=Q_{\mathcal{B}}/Q_{\mathcal{B}}(q^{\frac{1}{d{\bar{r}}^2}}-\varepsilon^{\frac{1}{d{\bar{r}}^2}})$. \index[not]{Q@$Q_\varepsilon$}

Obviously, for generic $\varepsilon$ we have $W_\varepsilon^s(G)\simeq W_{\mathcal{B}}^s(G)/(q^{\frac{1}{d{\bar{r}}^2}}-\varepsilon^{\frac{1}{d{\bar{r}}^2}})W_{\mathcal{B}}^s(G)$.


\section{Poisson reduction}\label{poisred}

\pagestyle{myheadings}
\markboth{CHAPTER~\thechapter.~Q-W--ALGEBRAS}{\thesection.~POISSON REDUCTION}

\setcounter{equation}{0}
\setcounter{theorem}{0}

In this section we recall basic facts on Poisson reduction. They will be used in the next section to describe Poisson q-W--algebras as reduced Poisson algebras.

Let $\mathfrak{M},~\mathfrak{B},~\mathfrak{B}'$ \index[not]{M@$\mathfrak{M}$} \index[not]{B@$\mathfrak{B}$} \index[not]{B@$\mathfrak{B}'$} be Poisson manifolds. \index{Poisson!manifold} Two Poisson submersions \index{Poisson!submersion}
$$
\begin{array}{ccccc}
&  & \mathfrak{M} &  &  \\
& \stackrel{\pi^{\prime } }{\swarrow } &  & \stackrel{\pi }{\searrow } &  \\
\mathfrak{B}' &  &  &  & \mathfrak{B}
\end{array}
$$
onto $\mathfrak{B}'$ and $\mathfrak{B}$ form a {\it dual pair} \index{dual pair} if the pullback ${\pi^\prime}^*C^\infty(\mathfrak{B}')$ is the centralizer of $\pi^* C^\infty (\mathfrak{B})$ in the Poisson algebra
$C^\infty (\mathfrak{M})$. \index[not]{C@$C^\infty (\mathfrak{M})$} In this case the sets $\mathfrak{B}'_b:=\pi^{\prime } \left( \pi
^{-1}(b) \right),~b\in
\mathfrak{B}$ \index[not]{B@$\mathfrak{B}'_b$} are Poisson submanifolds \index{Poisson!submanifold} in $\mathfrak{B}'$ called {\it reduced Poisson manifolds} (see \cite{W}, \S 8). \index{Poisson!manifold!reduced}

Fix an element $b\in \mathfrak{B}$. Then the algebra of functions $C^\infty (\mathfrak{B}'_b)$ may be described as follows. Let $\mathfrak{I}_b$ \index[not]{I@$\mathfrak{I}_b$} be the ideal in $C^\infty (\mathfrak{M})$ generated by the elements
${\pi}^*(f),~f\in C^\infty (\mathfrak{B}),~f(b)=0$. Denote $\mathfrak{M}_b=\pi^{-1}(b)$. Then the algebra $C^\infty (\mathfrak{M}_b)$
is simply the quotient of $C^\infty (\mathfrak{M})$ by $\mathfrak{I}_b$.  Denote by
$P_b:C^\infty (\mathfrak{M})\rightarrow C^\infty (\mathfrak{M})/\mathfrak{I}_b\simeq C^\infty (\mathfrak{M}_b)$ \index[not]{P@$P_b$} the canonical projection onto the quotient.
\begin{lemma}\label{redspace}
Suppose that the map $f\mapsto f(b)$ is
a character of the Poisson algebra $C^\infty (\mathfrak{B})$. Then one can define an action of the Poisson algebra $C^\infty (\mathfrak{B})$ on the space $C^\infty (\mathfrak{M}_b)$ by
\begin{equation}\label{redact}
f\cdot \varphi=P_b(\{ {\pi}^*(f), \tilde \varphi \}),
\end{equation}
where $f\in C^\infty (\mathfrak{B})$, $\varphi \in C^\infty (\mathfrak{M}_b)$ and $\tilde \varphi \in
C^\infty (\mathfrak{M})$ is a
representative of $\varphi$ in $C^\infty (\mathfrak{M})$ such that $P_b(\tilde
\varphi)=\varphi$.
Moreover, $C^\infty (\mathfrak{B}'_b)$ is isomorphic to the subspace of invariants in $C^\infty(\mathfrak{M}_b)$ with respect to this action.
\end{lemma}

\begin{proof}
Let $\varphi \in C^\infty (\mathfrak{M}_b)$. Choose a representative
$\tilde \varphi \in C^\infty (\mathfrak{M})$ such that $P_b(\tilde \varphi)=\varphi$.
Since the map $f\mapsto f(b)$ is
a character of the Poisson algebra $C^\infty (\mathfrak{B})$, $b\in \mathfrak{B}$ is a Poisson submanifold of $\mathfrak{B}$ with the zero Poisson structure, and hence the right hand side of (\ref{redact}) only depends on $\varphi$ but not on the representative $\tilde \varphi$. 

Indeed, if $f,g\in C^\infty (\mathfrak{B}),~g(b)=0$ and $h\in C^\infty (\mathfrak{M})$, so that $h{\pi}^*(g)\in \mathfrak{I}_b$, then we have
\begin{equation}\label{Ib}
\{{\pi}^*(f),h{\pi}^*(g)\}=h\{{\pi}^*(f),{\pi}^*(g)\}+{\pi}^*(g)\{{\pi}^*(f),h\}=
h{\pi}^*\{f,g\}+{\pi}^*(g)\{{\pi}^*(f),h\}.
\end{equation}
Next, $h{\pi}^*\{f,g\} \in \mathfrak{I}_b$ since the map $f\mapsto f(b)$ is a character of the Poisson algebra $C^\infty (\mathfrak{B})$, and hence $\{f,g\}(b)=0$. Clearly, by the definition of $\mathfrak{I}_b$, ${\pi}^*(g)\{{\pi}^*(f),h\}\in \mathfrak{I}_b$. We deduce that the right hand side of (\ref{Ib}) belongs to $\mathfrak{I}_b$. 

Therefore for any $f\in C^\infty (\mathfrak{B})$ one has $\{f, \mathfrak{I}_b\}\subset \mathfrak{I}_b$, and hence formula (\ref{redact}) defines an action of the Poisson algebra $C^\infty (\mathfrak{B})$ on the space $C^\infty (\mathfrak{M}_b)$.

Note that at the same time this proves that Hamiltonian vector fields of functions ${\pi}^*(f),~f\in C^\infty (\mathfrak{B})$ are tangent to the submanifold $\mathfrak{M}_b$. 

Using the definition of the dual
pair we obtain that
$\varphi={\pi^{\prime }}^*(\psi)$ for some $\psi \in C^\infty(\mathfrak{B}'_b)$ if
and only if
$P_b(\{ {\pi}^*(f), \tilde \varphi\})=0$ for every $f\in C^\infty (\mathfrak{B})$.
This implies that $C^\infty (\mathfrak{B}'_b)$ is isomorphic to the subspace of
invariants in $C^\infty (\mathfrak{M}_b)$ with respect to action (\ref{redact}).

\end{proof}

The algebra $C^\infty (\mathfrak{B}'_b)$ is called a {\it reduced Poisson algebra}. \index{Poisson!algebra!reduced}
We also denote it by $C^\infty (\mathfrak{M}_b)^{C^\infty (\mathfrak{B})}$. \index[not]{C@$C^\infty (\mathfrak{M}_b)^{C^\infty (\mathfrak{B})}$}

\begin{remark}\label{redpoisalg}
Note that the description of the algebra $C^\infty (\mathfrak{M}_b)^{C^\infty (\mathfrak{B})}$
obtained in Lemma \ref{redspace}
is independent
of both the manifold $\mathfrak{B}'$ and the projection $\pi^{\prime }$.
Observe also that the Hamiltonian vector fields of functions $\pi^*(f), f\in C^\infty (\mathfrak{B})$ are tangent to $\mathfrak{M}_b$, and hence the reduced space $\mathfrak{B}'_b$ may be identified with a cross--section of the action of the Poisson algebra $C^\infty (\mathfrak{B})$ on $\mathfrak{M}_b$ by Hamiltonian vector fields in the case when this action is free.
In particular, in this case $\mathfrak{B}'_b$ may be regarded as a submanifold in $\mathfrak{M}_b$.
\end{remark}

In the case when the map $f\mapsto f(b)$ is a character of the Poisson algebra $C^\infty (\mathfrak{B})$, the Poisson structure on the algebra $C^\infty (\mathfrak{B}'_b)$ can be explicitly described as follows. Let $\varphi_1, \varphi_2 \in C^\infty (\mathfrak{M}_b)^{C^\infty (\mathfrak{B})}$. Choose representatives
$\tilde{\varphi}_{1}, \tilde{\varphi}_{2} \in C^\infty (\mathfrak{M})$ such that $P_b(\tilde{\varphi}_{1})=\varphi_{1}$, $P_b(\tilde{\varphi}_{2})=\varphi_{2}$. Then
\begin{equation}\label{redPstr}
\{\varphi_1, \varphi_2\}=\{\tilde{\varphi}_1, \tilde{\varphi}_2\}\quad {\rm mod}~\mathfrak{I}_b.
\end{equation}
By Lemma \ref{redact} the class in $C^\infty (\mathfrak{M})/\mathfrak{I}_b\simeq C^\infty (\mathfrak{M}_b)$ of the function in right hand side of this formula is ${C^\infty (\mathfrak{B})}$--invariant and independent of the choice of the representatives $\tilde{\varphi}_{1}, \tilde{\varphi}_{2} \in C^\infty (\mathfrak{M})$.

An important example of dual pairs is provided by Poisson--Lie group actions. \index{action!Poisson--Lie group}
We say that a Lie group $A$ \index[not]{A@$A$} with Lie algebra $\frak a$ \index[not]{a@$\frak a$} {\it locally acts on a manifold} $\mathfrak{M}$ if there is a Lie algebra homomorphism $X \mapsto \widehat X$ \index[not]{X@$\widehat X$} from $\frak a$ to the Lie algebra of vector fields \index{vector field!algebra of} on $\mathfrak{M}$. In this case, by the existence and uniqueness theorem for solutions to ordinary differential equations, this homomorphism can be integrated to a {\it local group action} \index{action!Lie group!local} of $A$ on $\mathfrak{M}$ in the sense that for every point $m\in \mathfrak{M}$ there exists an open neighborhood $U$ of $m$ and an open neighborhood $V$ of the identity element in $A$ such that there is a smooth map $V\times U\to U$, $(a,u)\mapsto a\circ u$, and if $a_1, a_2\in V$ and $a_1a_2\in V$ then $(a_1a_2)\circ u=a_1\circ (a_2\circ m)$. Moreover, if we denote by $\Theta_m$ \index[not]{T@$\Theta_m$} the differential  of the map $V\to U$, $a\to a\circ m$ at the identity of $A$ then $\Theta_m(X)=\widehat X(m)$ for any $X\in {\frak a}$, $m\in \mathfrak{M}$.

Note that in this case the Lie algebra $\frak a$ acts on $C^\infty(\mathfrak{M})$ by
$$
X\circ \varphi=\widehat X\varphi,~X\in{\frak a},~\varphi\in C^\infty(\mathfrak{M}).
$$ 
We shall denote the space of invariants for this action by $C^\infty(\mathfrak{M})^A$. \index[not]{C@$C^\infty(\mathfrak{M})^A$}  

Right local Lie group actions are defined in a similar way, $X\mapsto \widehat X$ being a Lie algebra anti-homomorphism in this case. \index{Lie!algebra!anti-homomorphism}

If all vector fields $\widehat X$, $X\in {\frak a}$ are complete 
\index{vector field!complete} then by the existence and uniqueness theorem for solutions to ordinary differential equations $U=\mathfrak{M}$ and $V=A$, i.e. there is a Lie group action of $A$ on $\mathfrak{M}$. \index{action!Lie group}

Recall that a (local) left {\it Poisson--Lie group action} \index{action!Poisson--Lie group} of a Poisson--Lie group $A$ on a Poisson manifold $\mathfrak{M}$ is a (local) left Lie group action $A\times \mathfrak{M}\rightarrow \mathfrak{M}$ which is also a (locally defined) Poisson map \index{Poisson!map} (as usual, we suppose that $A\times \mathfrak{M}$ is equipped with the product Poisson structure). If in this situation the space $\mathfrak{M}/A$ is a smooth manifold, there exists a unique Poisson structure on $\mathfrak{M}/A$ such that the canonical projection $\mathfrak{M}\rightarrow \mathfrak{M}/A$ is a Poisson map (see Proposition \ref{admiss} below for a more general statement). 

Right Poisson--Lie group actions are defined in a similar way. 

Note that the property for a map to be Poisson is a local property, so it makes sense to consider locally defined Poisson maps. 

Denote by $\langle~\cdot~,~\cdot~\rangle$ \index[not]{ZZ@$\left\langle ~\cdot~ ,~\cdot~ \right\rangle$}
the canonical paring between  ${\frak a}^*$ \index[not]{a@${\frak a}^*$} and $\frak a$.
A map $\mu_A :\mathfrak{M}\rightarrow A^*$ \index[not]{m@$\mu_A$} is called a {\it moment map} \index{moment map} for a (local) left Poisson--Lie group
action $A\times \mathfrak{M}\rightarrow \mathfrak{M}$ if
\begin{equation}
L_{\widehat X} \varphi =\langle \mu_A^*(\theta_{A^*}) , X \rangle (\xi_\varphi ) ,
\end{equation}
where $\theta_{A^*}$ \index[not]{t@$\theta_{A^*}$}  is the universal left--invariant Maurer--Cartan form \index{Maurer--Cartan form, left--invariant} on
$A^*$, $X \in {\frak a}$,
$\widehat X$ is the corresponding vector field on $\mathfrak{M}$, $L_{\widehat X}$ the Lie derivative with respect to $\widehat X$, \index[not]{L@$L_{\widehat X}$} \index{Lie!derivative} and
$\xi_\varphi $ \index[not]{k@$\xi_\varphi $} is the Hamiltonian vector field \index{vector field!Hamiltonian} of $\varphi \in C^\infty (\mathfrak{M})$.

\begin{proposition}{\bf (\cite{Lu}, Theorem 4.9)} 
Let $A\times \mathfrak{M}\rightarrow \mathfrak{M}$ be a left (local) Poisson group action of a Poisson--Lie group $A$ on a Poisson manifold $\mathfrak{M}$ with moment map $\mu_A:\mathfrak{M}\to A^*$.
Denote by $\Pi_{A^*}$ \index[not]{P@$\Pi_{A^*}$} the Poisson tensor \index{Poisson!tensor} of $A^*$. Then there exists a right invariant bivector field \index{bivector field, right invariant} $\Lambda$ \index[not]{L@$\Lambda$} on $A^*$ such that $\Pi_{\mu_A}:=\Pi_{A^*}+\Lambda$ \index[not]{P@$\Pi_{\mu_A}$} is a Poisson tensor on $A^*$ and the map $\mu_A:\mathfrak{M}\to A^*_{\mu_A}$ is Poisson, where $A^*_{\mu_A}$ \index[not]{A@$A^*_{\mu_A}$} is the manifold $A^*$ equipped with the Poisson structure associated to $\Pi_{\mu_A}$. 
\end{proposition}

From the definition of the moment map it follows that if $\mathfrak{M}/A$ is a smooth
manifold then the canonical projection $\mathfrak{M}\rightarrow \mathfrak{M}/A$ and the moment map
$\mu_A:\mathfrak{M}\rightarrow A^*$
form a dual pair.

The main example of Poisson--Lie group actions is the so--called {\it dressing action}. \index{action!dressing}
The dressing action may be described as follows (see \cite{RIMS}, \S 3; \cite{LuW}, Theorem 2.4; Example 4.3 in \cite{Lu}; and formula (2.24) in \cite{dual}).
\begin{proposition}\label{dressingact}
Let $G$ be a connected simply connected Poisson--Lie group with factorizable
tangent Lie bialgebra,
$G^*$ the dual group. Then there exists a unique left local Poisson--Lie group action
$$
G\times G^*\rightarrow G^*,~~(g,(L_+,L_-))\mapsto g\circ (L_+,L_-),
$$
such that the identity mapping $\mu_G: G^* \rightarrow G^*$ is the moment map for
this action.

Moreover, let $q:G^* \rightarrow G$ \index[not]{q@$q(~\cdot~)$} be the map defined by
$$
q(L_+,L_-)=L_-^{-1}L_+.
$$
Then
\begin{equation}\label{qdress}
q(g\circ (L_+,L_-))=gL_-^{-1}L_+g^{-1}.
\end{equation}
\end{proposition}

The notion of Poisson--Lie groups may be generalized as follows.
Let $\left( {\frak a},{\frak a}^{*}\right)$ \index[not]{a@$\left( {\frak a},{\frak a}^{*}\right)$} be the tangent
Lie bialgebra of a Poisson--Lie group $A$. A connected Lie subgroup $K\subset A$ \index[not]{K@$K$} with Lie algebra ${\frak k}\subset {\frak a}$ \index[not]{k@${\frak k}\subset {\frak a}$} is called {\it admissible} \index{subgroup!Lie!admissible} if the annihilator ${\frak k}^{\perp }$ \index[not]{k@${\frak k}^{\perp }$} of $\mathfrak{k}$ in ${\frak a}^{*}$ is a Lie subalgebra ${\frak k}^{\perp }\subset {\frak a}^{*}$. \index{Lie!subalgebra}

\begin{proposition}\label{admiss}{\bf (\cite{RIMS}, Theorem 6; \cite{Lu}, \S 2)}
Let $A\times \mathfrak{M} \rightarrow \mathfrak{M}$ be a (local) Poisson--Lie group action of a Poisson--Lie group $A$ on a Poisson manifold $\mathfrak{M}$.
If $K\subset A$ is an admissible subgroup of $A$ then the space
$C^\infty(\mathfrak{M})^K$ of $K$-invariants in $C^\infty(\mathfrak{M})$ is a Poisson subalgebra in $C^\infty(\mathfrak{M})$. 

If $\mathfrak{M}/K$ is a smooth manifold, there exists a Poisson
structure on $\mathfrak{M}/K$ such that the canonical projection $\mathfrak{M}\rightarrow \mathfrak{M}/K$ is a Poisson map,
and $C^\infty(\mathfrak{M}/K)\simeq C^\infty(\mathfrak{M})^K$ as Poisson algebras.
\end{proposition}

We shall need the following particular example of dual pairs arising from
Poisson group actions.

Let $A\times \mathfrak{M} \rightarrow \mathfrak{M}$ be a left (local) Poisson group action of a Poisson--Lie
group $A$ on a manifold $\mathfrak{M}$.
Suppose that this action possesses a moment map $\mu_A : \mathfrak{M}\rightarrow A^*$.
Let $K$ be an admissible subgroup in $A$. Denote by $\frak k$ the Lie algebra of
$K$, so that ${\frak k}^\perp \subset {\frak a}^*$ is a Lie subalgebra in ${\frak
a}^*$.
Suppose also that there is a splitting ${\frak a}^*={\frak t}+ {\frak
k}^\perp$ \index[not]{t@${\frak t}$} (direct sum of vector spaces), and that
$\frak t$ is a Lie subalgebra in ${\frak a}^*$. Then the vector space ${\frak
k}^*$ is naturally identified with $\frak t$.

Assume that there is a unique factorization $A^*=TK^\perp$, where $K^\perp$, $T$ \index[not]{K@$K^\perp$} \index[not]{T@$T$} are the Lie subgroups \index{subgroup!Lie} of $A^*$ corresponding to the Lie subalgebras ${\frak k}^\perp$, ${\frak t}\subset {\frak a}^*$, respectively. For any $a^*=tk^\perp \in A^*$ with $k^\perp\in K^\perp$, $t\in T$ denote $\pi_{K^\perp}(a^*)=k^\perp$, $\pi_{T}(a^*)=t$. \index[not]{p@$\pi_{K^\perp}$} \index[not]{p@$\pi_{T}$} This defines maps $\pi_{K^\perp}:A^*\to K^\perp$, $\pi_{T}:A^*\to T$. 

\begin{proposition}\label{QPmoment}
Suppose that for any $k^\perp \in K^\perp$ the transformation
\begin{eqnarray}\label{tinv}
{\frak t}\rightarrow {\frak t},\\
t\mapsto ({\rm Ad}(k^\perp)t)_{{\frak t}}, \nonumber
\end{eqnarray}
where the subscript ${{\frak t}}$ stands for the ${{\frak t}}$--component with respect to the decomposition ${\frak a}^*={\frak t}+ {\frak
k}^\perp$, is invertible.

Define a map $\overline{\mu}_T:M\rightarrow T$ \index[not]{m@$\overline{\mu}_T$} by
$$
\overline{\mu}_T=\pi_{T}\mu_A.
$$
Then

(i)
$\overline{\mu}_T^*(C^\infty(T))$ is a Poisson subalgebra
in $C^\infty(\mathfrak{M})$,
and hence one can equip $T$ with a Poisson structure such that
$\overline{\mu}_T:\mathfrak{M}\rightarrow T$ is
a Poisson map;

(ii) The algebra $C^\infty(\mathfrak{M})^K$ is the centralizer of
$\overline{\mu}_T^*(C^\infty(T))$ in the Poisson algebra
$C^\infty(\mathfrak{M})$.
In particular, if $\mathfrak{M}/K$ is a smooth manifold then the maps
\begin{equation}\label{dp}
\begin{array}{ccccc}
&  & \mathfrak{M} &  &  \\
& \stackrel{\pi }{\swarrow } &  & \stackrel{\overline{\mu}_T}{\searrow } &  \\
\mathfrak{M}/K &  &  &  & T
\end{array}
\end{equation}
form a dual pair.
\end{proposition}

\begin{proof}
(i) We claim that multiplication in $A^*$ gives rise to a right Poisson--Lie group action $A^*_{\mu_A}\times A^*\to A^*_{\mu_A}$. Indeed, for $g\in A^*$ denote by $l_g, r_g$ \index[not]{l@$l_g$} \index[not]{r@$r_g$} the left (right) translation by $g$ on $A^*$. \index{action!Lie group!by left translations} \index{action!Lie group!by right translations} By the definition of $\Pi_{\mu_A}$
\begin{equation}\label{pmga}
\Pi_{\mu_A}(gh)=\Pi_{A^*}(gh)+\Lambda(gh)={l_g}_*\Pi_{A^*}(h)+{r_h}_*\Pi_{A^*}(g)+{r_h}_*\Lambda(g)={l_g}_*\Pi_{A^*}(h)+{r_h}_*\Pi_{\mu_A}(g),
\end{equation}
where we used the fact that $\Pi_{A^*}(gh)={l_g}_*\Pi_{A^*}(h)+{r_h}_*\Pi_{A^*}(g)$ as $A^*$ is a Poisson--Lie group and that $\Lambda$ is right invariant. By the definition of  Poisson--Lie group actions and Poisson maps, identity (\ref{pmga}) is equivalent to the fact that multiplication in $A^*$ gives rise to a right Poisson--Lie group action $A^*_{\mu_A}\times A^*\to A^*_{\mu_A}$.

Since ${\frak k}\subset  {\frak  a}$ is a Lie subalgebra and ${{\frak k}^\perp}^\perp\simeq {\frak k}$, the subgroup $K^\perp\subset A^*$ is admissible. Therefore restricting the action $A^*_{\mu_A}\times A^*\to A^*_{\mu_A}$ to $K^\perp$ we deduce that $C^\infty(A^*_{\mu_A})^{K^\perp}$ is a Poisson subalgebra in $C^\infty(A^*_{\mu_A})$, where the action of $K^\perp$ is induced by the action of $K^\perp\subset A^*$ on $A^*$ by right translations.  

Now recall that there is a unique factorization $A^*=TK^\perp$, and hence $C^\infty(A^*_{\mu_A})^{K^\perp}\simeq C^\infty(T)$. Thus $T$ naturally becomes a Poisson manifold and the map $\pi_T:A^*_{\mu_A}\to T$ becomes Poisson. We deduce that the map $\overline{\mu}_T=\pi_T{\mu_A}$ is Poisson as the composition of the Poisson maps $\mu_A:M\to A^*_{\mu_A}$ and $\pi_T:A^*_{\mu_A}\to T$.

(ii) By the definition of the moment map we have
\begin{equation}\label{X5}
L_{\widehat X} \varphi =\langle \mu_A^*(\theta_{A^*}) , X \rangle (\xi_\varphi ) ,
\end{equation}
where $X \in {\frak a}$, $\widehat X$ is the corresponding vector field on $M$ and
$\xi_\varphi $ is the Hamiltonian vector field of $\varphi \in C^\infty (\mathfrak{M})$. Since $A^*=TK^\perp$, the pullback of the left--invariant Maurer--Cartan form $\mu_A^*(\theta_{A^*})$ may be represented as follows
$$
\mu_A^*(\theta_{A^*})= {\rm Ad}(\pi_{K^\perp}\mu_A )^{-1}({\overline{\mu}_T}^*\theta_{T})+(\pi_{K^\perp}\mu_A )^*\theta_{K^\perp},
$$
where the form $(\pi_{K^\perp}\mu_A )^*\theta_{K^\perp}$ takes values in ${\frak k}^\perp$.

Now let $X \in {\frak k}$. Then $\langle (\pi_{K^\perp}\mu_A )^*\theta_{K^\perp}),X\rangle =0$ and formula (\ref{X5}) takes the form
\begin{equation}\label{+}
\begin{array}{l}
L_{\widehat X} \varphi =
\langle {\rm Ad}(\pi_{K^\perp}\mu_A )^{-1}({\overline{\mu}_T}^*\theta_{T}),X \rangle (\xi_\varphi )=\\
\\
=\langle {\rm Ad}(\pi_{K^\perp}\mu_A )^{-1}(\theta_{T}), X \rangle ({\overline{\mu}_T}_*(\xi_\varphi )) .
\end{array}
\end{equation}

Since by the assumption transformation (\ref{tinv}) is invertible, $L_{\widehat X} \varphi =\langle {\rm Ad}(\pi_{K^\perp}\mu_A )^{-1}(\theta_{T}), X \rangle ({\overline{\mu}_T}_*(\xi_\varphi ))=0$ for any $X\in {\frak k}$ if and only if ${\overline{\mu}_T}_*(\xi_\varphi )=0$.
Thus the function $\varphi \in C^\infty (\mathfrak{M})$ is $K$--invariant if and only if $\{ \varphi ,\overline{\mu}_T^*(\psi) \}=0$ for any $\psi \in C^\infty (T)$. This completes the proof.

\end{proof}

From the previous proposition, from Lemma \ref{redspace} and Remark \ref{redpoisalg} we immediately obtain the following corollary.
\begin{corollary}\label{remred}
Suppose that the conditions of Proposition \ref{QPmoment} are satisfied.
Let $t\in T$ be such that the map $f\mapsto f(t)$ is a character of the Poisson algebra $C^\infty(T)$. Then the action of $K$ on $\mathfrak{M}$ induces a (local) action on ${\overline{\mu}_T}^{-1}(t)$ and a (local) action on $C^\infty({\overline{\mu}_T}^{-1}(t))$ given by
$$
X\circ \varphi
=\langle {\rm Ad}(\pi_{K^\perp}\mu_A )^{-1}(\theta_{T}), X \rangle ({\overline{\mu}_T}_*(\xi_{\tilde\varphi} )), X\in {\frak k}, \varphi \in C^\infty({\overline{\mu}_T}^{-1}(t)),
$$
where $\tilde\varphi$ is any representative of $\varphi\in C^\infty({\overline{\mu}_T}^{-1}(t))\simeq C^\infty(\mathfrak{M})/\mathfrak{I}_t$ in $C^\infty(\mathfrak{M})$. The algebra $C^\infty({\overline{\mu}_T}^{-1}(t))^K$ of invariants with respect to this action is isomorphic to the reduced Poisson algebra $C^\infty({\overline{\mu}_T}^{-1}(t))^{C^\infty(T)}$.
\end{corollary}


\section{Poisson reduction and q-W--algebras}\label{wpsred}

\pagestyle{myheadings}
\markboth{CHAPTER~\thechapter.~Q-W--ALGEBRAS}{\thesection.~POISSON REDUCTION AND Q-W--ALGEBRAS}

\setcounter{equation}{0}
\setcounter{theorem}{0}

In this section we realize the quasiclassical limit of the algebra $W_{\mathcal{B}}^s(G)$ as the algebra of functions on a reduced
Poisson manifold. In this section we always assume that $\kappa=1$ and use the notation and conventions introduced in Sections \ref{qplgroups} and \ref{qWalgpois}. In particular, we always assume that the system of positive roots $\Delta_+$ is associated to the Weyl group element $s$ as in Definition \ref{circorddef}.

By Proposition \ref{Qdefr} (vii) one can define a character $\chi^{s}$ \index[not]{x@$\chi^{s}$} of the Poisson subalgebra ${\mathbb{C}}[M_+]\subset \mathbb{C}^s[G^*]$ such that $\chi^{s}(p(x))=\chi_q^{s}(x)$ {mod } $(q^{\frac{1}{d{\bar{r}}^2}}-1)$ for every $x\in \mathbb{C}_{\mathcal{B}}^s[M_+]$.

Recall that the image of the algebra $\mathbb{C}_{\mathcal{B}}^s[G_*]$ under the projection $p:\mathbb{C}_{\mathcal{B}}^s[G^*]\rightarrow \mathbb{C}_{\mathcal{B}}^s[G^*]/(1-q^{\frac{1}{d{\bar{r}}^2}})\mathbb{C}_{\mathcal{B}}^s[G^*]$ is a certain Poisson subalgebra of $\mathbb{C}^s[G^*]$ that we denoted by ${\mathbb{C}}^s[G_*]$ in Section \ref{qplgroups}. By part (ii) of Proposition \ref{locfin}
${\mathbb{C}}^s[G_*]\simeq {\mathbb{C}}[G]$ as algebras.
Let $I=p(I_{\mathcal{B}})$ \index[not]{I@$I$} be the ideal in $\mathbb{C}^s[G^*]$ generated by the kernel of $\chi^{s}$.
Then by formula (\ref{WqL}) the Poisson algebra $W^s(G):=W^s_q(G)/(q^{\frac{1}{d{\bar{r}}^2}}-1)W^s_q(G)$ \index[not]{W@$W^s(G)$} is isomorphic to the subspace of all $x+I\in Q_{1}$, $Q_{1}:=Q_{\mathcal{B}}/(1-q^{\frac{1}{d{\bar{r}}^2}})Q_{\mathcal{B}}\subset \mathbb{C}^s[G^*]/I$, \index[not]{Q@$Q_1$} such that $\{m,x\}\in I$ for any $m\in {\mathbb{C}}[M_-]$, and the Poisson bracket in $W^s(G)$ takes the form $\{(x+I),(y+I)\}=\{x,y\}+I$, $x+I,y+I\in W^s(G)$. Using the Poisson analogues of formulas (\ref{qmainactcl}) and (\ref{WqL}) we can also write $W^s(G)\simeq (\mathbb{C}^s[G^*]/I)^{{\mathbb{C}}[M_+]}\cap Q_{1}$, where the action of the Poisson algebra ${{\mathbb{C}}[M_+]}$ on the space $\mathbb{C}^s[G^*]/I$ is defined as follows
\begin{equation}\label{mainactcl}
x\cdot (v+I) =\rho_{\chi^{s}}(\{x,v\} ),
\end{equation}
$v\in \mathbb{C}^s[G^*]$ is any representative of $v+I\in \mathbb{C}^s[G^*]/I$, $\rho_{\chi^{s}}:\mathbb{C}^s[G^*]\to \mathbb{C}^s[G^*]/I$ the canonical projection, \index[not]{r@$\rho_{\chi^{s}}$} and $x\in {\mathbb{C}}[M_+]$.

One can describe the space of invariants $(\mathbb{C}^s[G^*]/I)^{{\mathbb{C}}[M_+]}$  with respect to this action by analyzing the related underlying manifolds and varieties. First observe that the algebra $(\mathbb{C}^s[G^*]/I)^{{\mathbb{C}}[M_+]}$ is a particular example of reduced Poisson algebras introduced in Lemma \ref{redspace} .

Indeed, recall that according to (\ref{fact}) any element $(L_+,L_-)\in G^*$ may be uniquely written as
\begin{equation}\label{fact1}
(L_+,L_-)=(n_+,n_-)(h_+,h_-),
\end{equation}
where $n_\pm \in N_\pm$, $h_+={\exp}(({1 \over 1-s}P_{{\h'}}+\frac{1}{2}P_{\h'^\perp})x),~h_-={\exp}(({s \over 1-s}P_{{\h'}}-\frac{1}{2}P_{\h'^\perp})x),~x\in
{\frak h}$.

Formula (\ref{fact}) and a decomposition of elements of $N_+$ into products of elements which belong to the one--dimensional subgroups corresponding to roots also imply that any element $L_+$ can be represented as in (\ref{l+dec})
\begin{equation}\label{lm}
L_+ = \prod_{\beta\in \Delta_+}
{\exp}[b_{\beta}X_{\beta}]
{\exp}\left[ \sum_{i=1}^lb_i({1 \over 1-s}P_{{\h'}}+\frac{1}{2}P_{\h'^\perp})H_i\right],~b_i,b_\beta\in {\Bbb C},
\end{equation}
where the product over roots is taken according to the normal ordering associated to $s$.

Now define a surjective submersion $\mu_{M_-}:G^* \rightarrow M_+$ \index[not]{m@$\mu_{M_-}$} by
\begin{equation}\label{mun}
\mu_{M_-}(L_+,L_-)=m_+,
\end{equation}
where for $L_+$ given by (\ref{lm}) $m_+$ is defined as follows
$$
m_+=\prod_{\beta\in \Delta_{\m_+}}
{\exp}[b_{\beta}X_{\beta}],
$$
and the product over roots is taken according to the normal order in the segment $\Delta_{\m_+}$.

Note that by definition ${\mathbb{C}}[M_+]=\{ \varphi\in \mathbb{C}^s[G^*]:\varphi=
\varphi(m_+)\}$. Therefore ${\mathbb{C}}[M_+]$ is generated by the pullbacks of
regular functions on $M_+$ with respect to the map $\mu_{M_-}$.
Since ${\mathbb{C}}[M_+]$ is a Poisson subalgebra in $\mathbb{C}^s[G^*]$, and  regular
functions
on $M_+$ are dense in $C^\infty(M_+)$ on every compact subset, we can equip the
manifold $M_+$ with
the Poisson structure in such a way that $\mu_{M_-}$ becomes a Poisson mapping.

Let $u\in M_+$ be the element defined by
\begin{equation}\label{defu}
\bar{u}=\prod_{i=1}^{l'}{\exp}[\bar{t}_i X_{\gamma_i}], \bar{t}_i\in \mathbb{C}, \bar{t}_i=\bar{k}_i~({\rm mod}~(q^{\frac{1}{d{\bar{r}}^2}}-1)), \index[not]{u@$\bar{u}$} \index[not]{t@$\bar{t}_i$}
\end{equation}
where the product over roots is taken according to the normal order in the segment $\Delta_{\m_+}$, and $\bar{k}_i$, $i=1, \ldots ,l'$ are defined in formula (\ref{charq}).

By Proposition \ref{qG*} (vi) (see formulas (\ref{imathdef}) and (\ref{imath})) the elements $
\tilde{L}^{\pm,{V}}=(p\otimes \nu)({^q{L^{\pm,V}}})$ belong to the space $
\mathbb{C}^s[G^*]\otimes {\rm End}\overline{V}$, and the map
$$
\imath: \mathbb{C}_{\mathcal{B}}^s[G^*]/(q^{\frac{1}{d{\bar{r}}^2}}-1)\mathbb{C}_{\mathcal{B}}^s[G^*] \rightarrow \mathbb{C}^s[G^*],~
(\imath\otimes id) \tilde{L}^{\pm,V}= {L^{\pm,\overline{V}}}
$$
is an isomorphism of Poisson--Hopf algebras.
In particular, from (\ref{rmatrspi}) it follows that
\begin{equation}\label{lv}
\begin{array}{l}
\tilde{L}^{+,\overline{V}}=\prod_{\beta\in \Delta_+}
{\exp}[p(\tilde{f}_{\beta}) \otimes
\pi_{\overline{V}}(X_{\beta})]\times \\
(p\otimes id){\exp}\left[ \sum_{i=1}^lhH_i\otimes \pi_{\overline{V}}(({2 \over 1-s}P_{{\h'}}+P_{\h'^\perp})Y_i)\right].
\end{array}
\end{equation}

From (\ref{lv}) and the definition of $\chi^s$ we obtain that $\chi^s(\varphi)=\varphi (\bar{u})$
for every $\varphi \in {\mathbb{C}}[M_+]$. $\chi^s$ naturally extends to a character
of the
Poisson algebra $C^\infty(M_+)$.

Now applying Lemma \ref{redspace} we can
define the reduced Poisson algebra $C^\infty(\mu_{M_-}^{-1}(\bar{u}))^{C^\infty(M_+)}$ as follows.
Denote by $I_{\bar{u}}$ \index[not]{I@$I_{\bar{u}}$} the ideal in $C^\infty(G^*)$ generated by elements
$\mu_{M_-}^*\psi,~\psi \in C^\infty(M_+),
~\psi({\bar{u}})=0$. Let $P_{\bar{u}}:C^\infty(G^*)\rightarrow
C^\infty(G^*)/I_{\bar{u}}\simeq C^\infty(\mu_{M_-}^{-1}({\bar{u}}))$ \index[not]{P@$P_{\bar{u}}$} be the
canonical projection. Define an action of $C^\infty(M_+)$ on
$C^\infty(\mu_{M_-}^{-1}({\bar{u}}))$ by
\begin{equation}\label{actred}
\psi\cdot \varphi=P_{\bar{u}}(\{ \mu_{M_-}^*\psi, \tilde \varphi\}),
\end{equation}
where $\psi \in C^\infty(M_-),~\varphi \in C^\infty(\mu_{M_-}^{-1}({\bar{u}}))$ and
$\tilde \varphi \in C^\infty(G^*)$
is a representative of $\varphi$ such that $P_{\bar{u}}\tilde \varphi=\varphi$.
The reduced Poisson algebra $C^\infty(\mu_{M_-}^{-1}({\bar{u}}))^{C^\infty(M_+)}$ is the algebra of $C^\infty(M_+)$--invariants in $C^\infty(\mu_{M_-}^{-1}({\bar{u}}))$ with respect to action (\ref{actred}). The reduced Poisson algebra is naturally equipped with a Poisson structure induced from $C^\infty(G^*)$ as described in (\ref{redPstr}).

\begin{lemma}\label{redreg}
Let $\overline{q(\mu_{M_-}^{-1}({\bar{u}}))}$ \index[not]{q@$\overline{q(\mu_{M_-}^{-1}({\bar{u}}))}$}  be the closure of $q(\mu_{M_-}^{-1}({\bar{u}}))$ in $G$  with respect to the Zariski topology. \index{topology!Zariski} Then  $Q_1\simeq {\mathbb{C}}[\overline{q(\mu_{M_-}^{-1}({\bar{u}}))}]$, and  the algebra
$W^s(G)$ is isomorphic to the algebra of regular functions on $\overline{q(\mu_{M_-}^{-1}({\bar{u}}))}$
 pullbacks of which under the map $q$ are invariant with respect to the action
(\ref{actred}) of
$C^\infty(M_+)$ on $C^\infty(\mu_{M_-}^{-1}({\bar{u}}))$, i.e.
$$
W^s(G)\simeq {\mathbb{C}}[\overline{q(\mu_{M_-}^{-1}({\bar{u}}))}]\cap
C^\infty(\mu_{M_-}^{-1}({\bar{u}}))^{C^\infty(M_+)},
$$
where ${\mathbb{C}}[\overline{q(\mu_{M_-}^{-1}({\bar{u}}))}]$ is regarded as a subalgebra in $C^\infty(\mu_{M_-}^{-1}({\bar{u}}))$ using the map $q^*:C^\infty(q(\mu_{M_-}^{-1}({\bar{u}})))\rightarrow C^\infty(\mu_{M_-}^{-1}({\bar{u}}))$ and the imbedding ${\mathbb{C}}[\overline{q(\mu_{M_-}^{-1}({\bar{u}}))}]\subset C^\infty(q(\mu_{M_-}^{-1}({\bar{u}})))$. 
\end{lemma}

\begin{proof}
First observe that by the definition $\mu_{M_-}^{-1}({\bar{u}})$ is a submanifold in $G^*$ and that $I=\mathbb{C}^s[G^*]\cap I_{\bar{u}}$. Therefore by the definition of the algebra $\mathbb{C}^s[G^*]$ and of the map $\mu_{M_-}$ the quotient $\mathbb{C}^s[G^*]/I$ is
identified with the algebra of functions on $\mu_{M_-}^{-1}({\bar{u}})$ generated by the restrictions of elements of $\mathbb{C}^s[G^*]$ to $\mu_{M_-}^{-1}({\bar{u}})$. 

Also by the definition $Q_1\subset \mathbb{C}^s[G^*]/I$ is the algebra generated by the restrictions to $\mu_{M_-}^{-1}({\bar{u}})$ of the pullbacks of elements from the algebra of regular functions $\mathbb{C}[G]$ under the map $q:G^*\to G$. Therefore $Q_{1}\simeq {\mathbb{C}}[\overline{q(\mu_{M_-}^{-1}({\bar{u}}))}]$.
From these observations we deduce that
$W^s(G)\simeq (\mathbb{C}^s[G^*]/I)^{{\mathbb{C}}[M_-]}\cap Q_{1}\simeq (\mathbb{C}^s[G^*]/I)^{{\mathbb{C}}[M_-]}\cap {\mathbb{C}}[\overline{q(\mu_{M_-}^{-1}({\bar{u}}))}]$.

Since ${\mathbb{C}}[M_-]$ is dense in $C^\infty(M_-)$ on every compact subset in
$M_-$ we have
$$
C^\infty(\mu_{M_+}^{-1}({\bar{u}}))^{C^\infty(M_-)}=
C^\infty(\mu_{M_+}^{-1}({\bar{u}}))^{{\mathbb{C}}[M_-]}.
$$

Now observe that action (\ref{actred}) of elements from ${\mathbb{C}}[M_-]$ coincides with action (\ref{mainactcl}) when restricted to elements from $\mathbb{C}^s[G^*]/I$, and hence
$W^s(G)\simeq (\mathbb{C}^s[G^*]/I)^{{\mathbb{C}}[M_-]}\cap {\mathbb{C}}[\overline{q(\mu_{M_-}^{-1}({\bar{u}}))}]\simeq C^\infty(\mu_{M_+}^{-1}({\bar{u}}))^{{\mathbb{C}}[M_-]}\cap {\mathbb{C}}[\overline{q(\mu_{M_-}^{-1}({\bar{u}}))}]=C^\infty(\mu_{M_-}^{-1}({\bar{u}}))^{C^\infty(M_-)}\cap {\mathbb{C}}[\overline{q(\mu_{M_+}^{-1}({\bar{u}}))}]$. This completes the proof.

\end{proof}

We shall realize the
algebra $C^\infty(\mu_{M_-}^{-1}({\bar{u}}))^{C^\infty(M_+)}$ as the algebra of
functions on a reduced
Poisson manifold. In this construction we use
the dressing
action of the Poisson--Lie group $G$ on $G^*$.

Consider the restriction of the (local) dressing action $G\times G^* \rightarrow G^*$ to
the subgroup $M_-\subset G$. We shall describe the reduced Poisson algebra $C^\infty(\mu_{M_-}^{-1}({\bar{u}}))^{C^\infty(M_-)}$ in terms of the dressing action.
\begin{lemma}\label{redpois}
The preimage $\mu_{M_-}^{-1}({\bar{u}})\subset G^*$ is locally invariant under the (locally defined) dressing action of $M_-$, and the algebra $C^\infty(\mu_{M_-}^{-1}({\bar{u}}))^{M_-}$ is isomorphic to $C^\infty(\mu_{M_-}^{-1}({\bar{u}}))^{C^\infty(M_+)}$.
\end{lemma}

\begin{proof}
The proof will be based on Corollary \ref{remred}. We shall verify that the conditions of the corollary are satisfied with $A=B_- , K=M_- , A^*=B_+,
T=M_+$, and with an appropriate choice of $K^\perp$ and $\mu_A$ which will be specified below.

Recall that by (\ref{rspm}) the subspaces ${\frak i}_\pm$ and ${\frak k}_\pm$ defined by (\ref{bnpm}) coincide with
the Borel subalgebras ${\frak b}_\pm$ and their nilpotent radicals ${\frak n}_\pm$, respectively. Therefore by Proposition \ref{bpm} (iv) $(\b_\pm,\b_\pm^*)\simeq (\b_\pm,\b_\mp)$ $(\b_\pm^*\simeq \b_\mp)$ is a subbialgebra of $(\g,\g^*)$, and hence $B_\pm\subset G^*$ are Poisson--Lie subgroups.  

Restrict the local dressiong action of $G$ on $G^*$ to the Poisson--Lie subgroup $B_-\subset G$.
By Proposition \ref{dressingact}, by part (iv) of Proposition \ref{bpm}, and by the definition of the moment map we have for any $X \in {\frak b}_-$, $\varphi \in C^\infty(G^*)$
\begin{equation}
L_{\widehat X} \varphi(L_+,L_-) =\left\langle \theta_{G^*}(L_+,L_-) , X \right\rangle (\xi_\varphi )=
\left\langle r_+^{-1}\mu_{B_-}^*(\theta_{B_+}), X\right\rangle (\xi_\varphi ),
\end{equation}
where $\widehat X$ is the corresponding dressing action vector field on $G^*$,
$\xi_\varphi $ is the Hamiltonian vector field of $\varphi \in C^\infty (G^*)$,
and the map $\mu_{B_-}:G^*\rightarrow B_+$ \index[not]{m@$\mu_{B_-}$} is defined by $\mu_{B_-}(L_+,L_-)=L_+$.

From Proposition \ref{bpm} (iv) and from the definition of the moment map it
follows that $\mu_{B_-}$ is a moment map for the restriction of the dressing action to the subgroup $B_-$.

Next, observe that the complementary subset to $\Delta_{\m_+}$ in $\Delta_+$ is a minimal segment $\Delta_{\m_+}^0=\Delta_+\setminus \Delta_{\m_+}$. \index[not]{m@$\Delta_{\m_+}^0$} Now using Proposition \ref{bpm} (iv) the subspace $\m_-^\perp$ in $\b_+$ can be identified with the vector subspace in $\b_+$ spanned by the Cartan subalgebra $\h$ and by the root subspaces corresponding to the roots from the minimal segment $\Delta_{\m_+}^0$. Using Lemmas \ref{segmsub} and \ref{orthogr}, and the fact that the adjoint action of $\h$ normalizes each root subspace, we deduce that $\m_-^\perp\subset \b_+$ \index[not]{m@$\m_-^\perp$} is the Lie subalgebra generated by the root vectors corresponding to the roots from $\Delta_{\m_+}^0$ and by $\h$, and hence $M_-\subset B_-$ is an admissible subgroup.

Note that 
\begin{equation}\label{m-pdec}
\m_-^\perp=\h+{\m_-^\perp}_0\subset \b_+ \text{ (direct sum of vector spaces)}, 
\end{equation}
where ${\m_-^\perp}_0$ is the Lie subalgebra generated by root vectors corresponding to the roots from the minimal segment $\Delta_{\m_+}^0$.

Since we also have the disjoint union $\Delta_+=\Delta_{\m_+}\cup \Delta_{\m_+}^0$,
one also obtains, using (\ref{m-pdec}), a direct vector space decomposition
\begin{equation}\label{m+mdec}
\b_+=\m_++\m_-^\perp.
\end{equation}
 
Moreover, since we have the unique factorization $B_+=N_+H$, Lemma \ref{decNr} applied to the disjoint union $\Delta_+=\Delta_{\m_+}\cup \Delta_{\m_+}^0$ and decomposition (\ref{m-pdec}) imply that the Poisson--Lie group $B_+$ dual to $B_-$ can be uniquely factorized as $B_+=M_+M_-^\perp$, where $M_-^\perp\subset B_+$ \index[not]{M@$M_-^\perp$} is the Lie subgroup corresponding to the Lie subalgebra ${\frak m}_-^\perp \subset\b_+$, and $M_+\subset B_+$ is the Lie subgroup corresponding to the Lie subalgebra ${\frak m}_+$.

Now observe that by (\ref{m-pdec}) any element $k^\perp \in M_-^\perp$ can be uniquely represented in the form $k^\perp=hk^\perp_0$, $h\in H$, $k^\perp_0=\exp(x)$, $x\in {\m_-^\perp}_0$. Recall also that the Lie subalgebra $\m_+$ is generated by the root subspaces corresponding to the roots from the minimal segment $\Delta_{\m_+}$. Since all root subspaces are invariant under the adjoint action of $\h$ and the restriction of the adjoint action of the root vectors corresponding to the roots from the minimal segment $\Delta_{\m_+}^0$ is nilpotent we deduce that for any $m_+\in \m_+$, $k^\perp \in M_-^\perp$, $k^\perp=hk^\perp_0$, $h\in H$, $k^\perp_0=\exp(x)$, $x\in {\m_-^\perp}_0$ one has
$$
({\rm Ad}(k^\perp)(m_+))_{\m_+}=({\rm Ad}(hk^\perp_0)(m_+))_{\m_+}={\rm Ad}h(({\rm Ad}k^\perp_0(m_+))_{\m_+})=
{\rm Ad}h((\exp({\rm ad}x)(m_+))_{\m_+})={\rm Ad}h((id+U(x))(m_+)),
$$
where the subscript $\m_+$ stands for the $\m_+$--component in direct vector space decomposition (\ref{m+mdec}), and $U(x)$ is a linear nilpotent transformation of $\m_+$. 

The maps ${\rm Ad}h$ and $id+U(x)$ are obviously invertible. Hence for any  $k^\perp \in M_-^\perp$ the map
$$
\m_+\to \m_+, m_+\mapsto ({\rm Ad}(k^\perp)(m_+))_{\m_+}
$$
is invertible as well.

We conclude that all the conditions of Corollary \ref{remred}
are satisfied with $A=B_- , K=M_- , A^*=B_+,
T=M_+ , K^\perp = M_-^\perp, \mu_A =\mu_{B_-}$.
It follows that the preimage $\mu_{M_-}^{-1}({\bar{u}})\subset G^*$ is locally stable under the (locally defined) dressing action of $M_-$, and the algebra $C^\infty(\mu_{M_-}^{-1}({\bar{u}}))^{M_-}$ is isomorphic to $C^\infty(\mu_{M_-}^{-1}({\bar{u}}))^{C^\infty(M_+)}$.
This completes the proof.

\end{proof}

Observe that by (\ref{qdress}) under the map $q:G^*\rightarrow G$, $q(L_+,L_-)=L_-^{-1}L_+$ the dressing action becomes the action of $G$ on itself by conjugations. Consider the restriction
of this action to the subgroup $M_-$. Denote by $\pi_q:G\rightarrow G/M_-$ \index[not]{p@$\pi_q$}
the canonical projection onto the quotient with respect to this action. We shall see that $\pi_q(\overline{q(\mu_{M_-}^{-1}({\bar{u}}))})$ \index[not]{p@$\pi_q(\overline{q(\mu_{M_-}^{-1}({\bar{u}}))})$} is an algebraic variety and $\mathbb{C}[\pi_q(\overline{q(\mu_{M_-}^{-1}({\bar{u}}))})]\simeq W^s(G)$. We shall also obtain an explicit description of the variety $\overline{q(\mu_{M_-}^{-1}({\bar{u}}))}$ and of the projection $\pi_q(\overline{q(\mu_{M_-}^{-1}({\bar{u}}))})$.

First we describe the image of the level surface \index{level surface}  $\mu_{M_-}^{-1}({\bar{u}})$ of the map $\mu_{M_-}$  under the map $q$.
Let $X_\alpha(t)=\exp(tX_\alpha)\in G$, $t\in \mathbb{C}$ \index[not]{X@$X_\alpha(t)$} be the one--parameter subgroup \index{subgroup!one--parameter corresponding to a root} in $G$ corresponding to root $\alpha\in \Delta$. Recall that for any $\alpha \in \Delta_+$ and any $t\neq 0$ the element 
\begin{equation}\label{sarep}
s_\alpha(t)=X_{\alpha}(t)X_{-\alpha}(-t^{-1})X_{\alpha}(t)\in N_G(H) \index[not]{s@$s_\alpha(t)$}
\end{equation} 
is a representative for the reflection $s_\alpha$ corresponding to the root $\alpha$. Denote by $s\in N_G(H)$ the following representative of the Weyl group element $s\in W$,
\begin{equation}\label{defrep}
s=s_{\gamma_1}(\bar{t}_1)\ldots s_{\gamma_{l'}}(\bar{t}_{l'}),
\end{equation}
where the numbers $\bar{t}_i$ are defined in (\ref{defu}), and we assume that $\bar{t}_i\neq 0$ for any $i$.

We shall also use the following representatives for $s^1$ and $s^2$
$$
s^1=s_{\gamma_{1}}(\bar{t}_1)\ldots s_{\gamma_{\widetilde{l}}}(\bar{t}_{\widetilde{l}}),~s^2=s_{\gamma_{\widetilde{l}+1}}(\bar{t}_{\widetilde{l}+1})\ldots s_{\gamma_{l'}}(\bar{t}_{l'}).
$$

The following Proposition is an improved version of Proposition 7.2 in \cite{S11} suitable for the purposes of quantization.
\begin{proposition}\label{constrt}
Let $q:G^*\rightarrow G$ be the map defined by
$$
q(L_+,L_-)=L_-^{-1}L_+.
$$
Suppose that the numbers $\bar{t}_i$ defined in (\ref{defu}) are not equal to zero for all $i$. Then 
\begin{equation}\label{preim}
q(\mu_{M_-}^{-1}({\bar{u}})) \subset N_-sH^0Z_+M_-=N_-sH^0M_-Z_+=(N_-\cap N^s) Z_-sH^0M_-Z_+ =
\end{equation}
$$
=(N_-\cap N^s)Z_-sH^0Z_+M_-\subset N^ssZ^sN^s,
$$
where $H^0$ \index[not]{H@$H^0$} is the connected subgroup of $H$ corresponding to the Lie subalgebra $\h^s=\h'^\perp\subset \h$, $Z_\pm=Z^s\cap N_\pm$. \index[not]{Z@$Z_\pm$}

The closure $\overline{q(\mu_{M_-}^{-1}({\bar{u}}))}$ of $q(\mu_{M_-}^{-1}({\bar{u}}))$ in $G$ with respect to the Zariski topology is also contained in $N^ssZ^sN^s$. 
\end{proposition}

\begin{proof}
Using definition (\ref{mun}) of the map $\mu_{M_-}$  we can
describe the preimage $\mu_{M_-}^{-1}({\bar{u}})$ as follows
\begin{equation}\label{mun1}
\mu_{M_-}^{-1}({\bar{u}})=\{({\bar{u}}yh_+,n_-h_-) | n_- \in N_- , h_\pm=e^{r_\pm^s x}, x \in \h, y\in N_{\Delta_+\setminus \Delta_{\m_+}} \},
\end{equation}
where, as in Section \ref{slodowy}, for any additively closed subset of roots $\Xi\subset \Delta$ which does not contain opposite roots we denote by $N_\Xi$ the subgroup in $G$ generated by the one--parameter subgroups corresponding to the roots from $\Xi$.
Therefore
\begin{equation}\label{dva}
q(\mu_{M_-}^{-1}({\bar{u}}))=
\{ h_-^{-1}n_-^{-1}\bar{u}yh_+| n_- \in N_- , h_\pm=e^{r_\pm^s x}, x \in \h, y\in N_{\Delta_+\setminus \Delta_{\m_+}} \}.
\end{equation}

First we show that for any $y\in N_{\Delta_+\setminus \Delta_{\m_+}}$ and $n_-\in N_-$ the element $n_-^{-1}uy$ belongs to $N_-sH^0Z_+M_-$. Fix the circular normal ordering on $\Delta$ corresponding to the normal ordering of $\Delta_+$ associated to $s$ as in Definition \ref{circorddef}.

In the proof we shall frequently use the following lemma.
\begin{lemma}\label{comm}
Let $[\alpha,\beta]\subset \Delta$ be a minimal segment and assume that $[\alpha,\beta]=[\alpha,\gamma]\cup [\delta,\beta]$, where the segments $[\alpha, \gamma]$ and $[\delta, \beta]$ are disjoint and minimal as well. Then any element $m\in N_{[\alpha,\beta]}$ can be uniquely factorized as $m=g_1g_2=g_2'g_1'$, $g_1,g_1'\in N_{[\alpha,\gamma]}$, $g_2,g_2'\in N_{[\delta,\beta]}$. Moreover, if $\delta=\beta$ then for any $m'\in N_{[\alpha,\gamma]}$ and any $t\in \mathbb{C}$ one has
$m'X_\beta(t)=X_\beta(t)m''$, where $m''\in N_{[\alpha,\gamma]}$.
\end{lemma}

\begin{proof}
The proof is obtained by straightforward application of Lemmas \ref{decNr} and \ref{segmsub}.

 \end{proof}

Since the roots $\gamma_1, \ldots , \gamma_{\widetilde{l}}$ are mutually orthogonal the adjoint action of $s_{\gamma_i}(\bar{t}_i)$, $i=1,\ldots, \widetilde{l}$ on each of the root subspaces $\g_{\gamma_j}$, $j=1, \ldots ,n, j\neq i$ is given by multiplication by a non--zero constant. Therefore there are non--zero constants $c_1,\ldots , c_{\widetilde{l}}$ such that $X_{\gamma_{k}}(c_{k})s_{\gamma_1}\ldots s_{\gamma_{k-1}}=s_{\gamma_1}\ldots s_{\gamma_{k-1}}X_{\gamma_{k}}(-\bar{t}_{k}^{-1})$, $k=2,\ldots, \widetilde{l}$, and we define $c_1=-\bar{t}_{1}^{-1}$.

Obviously we have
$$
X_{\gamma_{1}}(\bar{t}_{1})\ldots X_{\gamma_{\widetilde{l}}}(\bar{t}_{\widetilde{l}})= X_{-\gamma_1}(-c_1)\ldots X_{-\gamma_{\widetilde{l}}}(-c_{\widetilde{l}})X_{-\gamma_{\widetilde{l}}}(c_{\widetilde{l}})\ldots X_{-\gamma_{1}}(c_{1})X_{\gamma_{1}}(\bar{t}_{1})\ldots X_{\gamma_{\widetilde{l}}}(\bar{t}_{\widetilde{l}})=
$$
$$
=n_1X_{-\gamma_{\widetilde{l}}}(c_{\widetilde{l}})\ldots X_{-\gamma_{1}}(c_{1})X_{\gamma_{1}}(\bar{t}_{1})\ldots X_{\gamma_{\widetilde{l}}}(\bar{t}_{\widetilde{l}}),~~n_1=X_{-\gamma_1}(-c_1)\ldots X_{-\gamma_{\widetilde{l}}}(-c_{\widetilde{l}})\in N_{\Delta_{(1)}^-},
$$
where $\Delta_{(1)}^-=\{\alpha\in \Delta_-: -\gamma_1\leq \alpha\leq -\gamma_{\widetilde{l}}\}$. \index[not]{D@$\Delta_{(1)}^-$}

Using the relation $X_{-\gamma_{1}}(-\bar{t}_{1}^{-1})X_{\gamma_1}(\bar{t}_1)=X_{\gamma_1}(-\bar{t}_1){s}_{\gamma_1}$, which follows from (\ref{sarep}), one can rewrite the last identity as follows
\begin{equation}\label{1"}
X_{\gamma_{1}}(\bar{t}_{1})\ldots X_{\gamma_{\widetilde{l}}}(\bar{t}_{\widetilde{l}})=
n_1X_{-\gamma_{\widetilde{l}}}(c_{\widetilde{l}})\ldots X_{-\gamma_{2}}(c_{2})X_{\gamma_1}(-\bar{t}_1){s}_{\gamma_1}X_{\gamma_{2}}(\bar{t}_{2})\ldots X_{\gamma_{\widetilde{l}}}(\bar{t}_{\widetilde{l}}).
\end{equation}

Now we can write 
$$
X_{-\gamma_{\widetilde{l}}}(c_{\widetilde{l}})\ldots X_{-\gamma_{2}}(c_{2})X_{\gamma_1}(-\bar{t}_1)=
$$
$$
=X_{-\gamma_{\widetilde{l}}}(c_{\widetilde{l}})\ldots X_{-\gamma_{2}}(c_{2})X_{\gamma_1}(-\bar{t}_1)X_{-\gamma_{2}}(-c_{2})\ldots X_{-\gamma_{\widetilde{l}}}(-c_{\widetilde{l}})X_{-\gamma_{\widetilde{l}}}(c_{\widetilde{l}})\ldots X_{-\gamma_{2}}(c_{2}).
$$ 
The product $X_{-\gamma_{\widetilde{l}}}(c_{\widetilde{l}})\ldots X_{-\gamma_{2}}(c_{2})X_{\gamma_1}(-t_1)X_{-\gamma_{2}}(-c_{2})\ldots X_{-\gamma_{\widetilde{l}}}(-c_{\widetilde{l}})$ belongs to the subgroup of $G$ generated by the one--parameter subgroups corresponding to roots from the set $\Delta^{(1)}:=\{\alpha\in \Delta: -\gamma_{2}\leq \alpha \leq \gamma_{1},s^1\alpha=-\alpha\}$. \index[not]{D@$\Delta^{(1)}$} By Lemma \ref{minsegm} the minimal segment $\{\alpha\in \Delta: -\gamma_{2}\leq \alpha \leq \gamma_{1}\}$ is additively closed and the set of roots on which $s^1$ acts by multiplication by $-1$ is also additively closed. Hence $\Delta^{(1)}$ is additively closed, and $X_{-\gamma_{\widetilde{l}}}(c_{\widetilde{l}})\ldots X_{-\gamma_{2}}(c_{2})X_{\gamma_1}(-\bar{t}_1)X_{-\gamma_{2}}(-c_{2})\ldots X_{-\gamma_{\widetilde{l}}}(-c_{\widetilde{l}})\in N_{\Delta^{(1)}}$. 

Let $\Delta_+^{(1)}=\Delta^{(1)}\cap \Delta_+=\{\alpha\in \Delta_+: \alpha\leq \gamma_1,s^1\alpha=-\alpha\}$, and $\Delta^{(1)}_-=\Delta^{(1)}\cap\Delta_-=\{\alpha\in \Delta_-: -\gamma_{2}\leq \alpha,s^1\alpha=-\alpha\}$. \index[not]{D@$\Delta^{(1)}_\pm$} Then $\Delta^{(1)}=\Delta_+^{(1)}\cup \Delta_-^{(1)}$ (disjoint union), and, using Lemma \ref{decNr}, the element 
$$
X_{-\gamma_{\widetilde{l}}}(c_{\widetilde{l}})\ldots X_{-\gamma_{2}}(c_{2})X_{\gamma_1}(-\bar{t}_1)X_{-\gamma_{2}}(-c_{2})\ldots X_{-\gamma_{\widetilde{l}}}(-c_{\widetilde{l}})\in N_{\Delta^{(1)}}= N_{\Delta_-^{(1)}}N_{\Delta_+^{(1)}}
$$ 
can be uniquely factorized as follows
$$
X_{-\gamma_{\widetilde{l}}}(c_{\widetilde{l}})\ldots X_{-\gamma_{2}}(c_{2})X_{\gamma_1}(-\bar{t}_1)X_{-\gamma_{2}}(-c_{2})\ldots X_{-\gamma_{\widetilde{l}}}(-c_{\widetilde{l}})=n_2'x_1',
$$
where $n_2'\in N_{\Delta_-^{(1)}}, x_1'\in N_{\Delta_+^{(1)}}$.

Substituting the last relation into (\ref{1"}) and using the definition of $c_2$ and the orthogonality of roots $\gamma_1$ and $\gamma_2$ we obtain
$$
X_{\gamma_1}(\bar{t}_1)\ldots X_{\gamma_{\widetilde{l}}}(\bar{t}_{\widetilde{l}})=n_2x_1'X_{-\gamma_{\widetilde{l}}}(c_{\widetilde{l}})\ldots X_{-\gamma_{3}}(c_{3}){s}_{\gamma_1}X_{-\gamma_{2}}(-\bar{t}_{2}^{-1})X_{\gamma_{2}}(\bar{t}_2)\ldots X_{\gamma_{\widetilde{l}}}(\bar{t}_{\widetilde{l}}),
$$
where $n_2=n_1n_2'\in N_{\Delta_-^{s^1}}$, $\Delta_-^{s^1}:=\{\alpha \in \Delta_-: s^1\alpha=-\alpha\}$. \index[not]{D@$\Delta_-^{s^1}$}

Now we can use the relation $X_{-\gamma_{2}}(-\bar{t}_{2}^{-1})X_{\gamma_{2}}(\bar{t}_2)=X_{\gamma_2}(-\bar{t}_2){s}_{\gamma_2}$ following from (\ref{sarep}), the orthogonality of roots $\gamma_1$ and $\gamma_2$, and apply similar arguments to get
\begin{equation}\label{2"}
X_{\gamma_1}(\bar{t}_1)\ldots X_{\gamma_{\widetilde{l}}}(\bar{t}_{\widetilde{l}})=n_2x_1'X_{-\gamma_{\widetilde{l}}}(c_{\widetilde{l}})\ldots X_{-\gamma_{3}}(c_{3})X_{\gamma_2}(a_2){s}_{\gamma_1}{s}_{\gamma_2}X_{\gamma_{3}}(\bar{t}_3)\ldots X_{\gamma_{\widetilde{l}}}(\bar{t}_{\widetilde{l}}),a_2\neq 0.
\end{equation}

We can also write 
$$
X_{-\gamma_{\widetilde{l}}}(c_{\widetilde{l}})\ldots X_{-\gamma_{3}}(c_{3})X_{\gamma_2}(a_2)=
$$
$$
=X_{-\gamma_{\widetilde{l}}}(c_{\widetilde{l}})\ldots X_{-\gamma_{3}}(c_{3})X_{\gamma_2}(a_2)X_{-\gamma_{3}}(-c_{3})\ldots X_{-\gamma_{\widetilde{l}}}(-c_{\widetilde{l}})X_{-\gamma_{\widetilde{l}}}(c_{\widetilde{l}})\ldots X_{-\gamma_{3}}(c_{3}).
$$ 
The product $X_{-\gamma_{\widetilde{l}}}(c_{\widetilde{l}})\ldots X_{-\gamma_{3}}(c_{3})X_{\gamma_2}(a_2)X_{-\gamma_{3}}(-c_{3})\ldots X_{-\gamma_{\widetilde{l}}}(-c_{\widetilde{l}})$ belongs to the subgroup of $G$ generated by the one--parameter subgroups corresponding to the roots from the set $\Delta^{(2)}:=\{\alpha\in \Delta: -\gamma_{3}\leq \alpha \leq \gamma_{2},s^1\alpha=-\alpha\}$. \index[not]{D@$\Delta^{(2)}$} By Lemma \ref{minsegm} the minimal segment $\{\alpha\in \Delta: -\gamma_{3}\leq \alpha \leq \gamma_{2}\}$ is additively closed and the set of roots on which $s^1$ acts by multiplication by $-1$ is also additively closed. Hence $\Delta^{(2)}$ is additively closed, and 
$X_{-\gamma_{\widetilde{l}}}(c_{\widetilde{l}})\ldots X_{-\gamma_{3}}(c_{3})X_{\gamma_2}(a_2)X_{-\gamma_{3}}(-c_{3})\ldots X_{-\gamma_{\widetilde{l}}}(-c_{\widetilde{l}})\in N_{\Delta^{(2)}}$.

Let $\Delta_+^{(2)}=\Delta^{(2)}\cap \Delta_+=\{\alpha\in \Delta_+: \alpha\leq \gamma_2, s^1\alpha=-\alpha\}$ and $\Delta_-^{(2)}=\Delta^{(2)} \cap\Delta_-=\{ \alpha\in \Delta_-:s^1\alpha=-\alpha, -\gamma_3\leq \alpha\}$. \index[not]{D@$\Delta^{(2)}_\pm$}

Then $\Delta^{(2)}=\Delta_+^{(2)}\cup \Delta_-^{(2)}$ (disjoint union), and, using Lemma \ref{decNr}, the element 
$$
X_{-\gamma_{\widetilde{l}}}(c_{\widetilde{l}})\ldots X_{-\gamma_{3}}(c_{3})X_{\gamma_2}(a_2)X_{-\gamma_{3}}(-c_{3})\ldots X_{-\gamma_{\widetilde{l}}}(-c_{\widetilde{l}})\in N_{\Delta^{(2)}}= N_{\Delta_-^{(2)}}N_{\Delta_+^{(2)}}
$$ 
can be represented as follows    
\begin{equation}\label{aux"}
X_{-\gamma_{\widetilde{l}}}(c_{\widetilde{l}})\ldots X_{-\gamma_{3}}(c_{3})X_{\gamma_2}(a_2)X_{-\gamma_{3}}(-c_{3})\ldots X_{-\gamma_{\widetilde{l}}}(-c_{\widetilde{l}})=n_3'x_2'',
\end{equation}
where $n_3'\in N_{\Delta_-^{(2)}}, x_2''\in N_{\Delta_+^{(2)}}$.

Substituting the last relation into (\ref{2"}) and using the definition of $c_3$ and the orthogonality of roots $\gamma_1$, $\gamma_2$ and $\gamma_3$ we obtain
\begin{equation}\label{3"}
X_{\gamma_1}(\bar{t}_1)\ldots X_{\gamma_{\widetilde{l}}}(\bar{t}_{\widetilde{l}})=n_2x_1'n_3'x_2''X_{-\gamma_{\widetilde{l}}}(c_{\widetilde{l}})\ldots X_{-\gamma_{4}}(c_{4})
{s}_{\gamma_1}{s}_{\gamma_2}X_{\gamma_{3}}(-\bar{t}_{3}^{-1})X_{\gamma_{3}}(\bar{t}_3)\ldots X_{\gamma_{\widetilde{l}}}(\bar{t}_{\widetilde{l}}).
\end{equation}

Since $\Delta_+^{(1)}\subset \Delta_+^{(2)}$, $N_{\Delta_+^{(1)}}\subset N_{\Delta_+^{(2)}}$, and we deduce that $x_1'n_3'x_2''\in N_{\Delta^{(2)}}$. Therefore using Lemma \ref{decNr} and the factorization $N_{\Delta^{(2)}}=N_{\Delta_-^{(2)}}N_{\Delta_+^{(2)}}$ we get $x_1'n_3'x_2''=n_3''x_2'$, $x_2'\in N_{\Delta_+^{(2)}}$, $n_3''\in N_{\Delta_-^{(2)}}$. Now (\ref{3"}) takes the form
$$
X_{\gamma_1}(\bar{t}_1)\ldots X_{\gamma_{\widetilde{l}}}(\bar{t}_{\widetilde{l}})=n_3x_2'X_{-\gamma_{\widetilde{l}}}(c_{\widetilde{l}})\ldots X_{-\gamma_{4}}(c_{4})
{s}_{\gamma_1}{s}_{\gamma_2}X_{-\gamma_{3}}(-\bar{t}_{3}^{-1})X_{\gamma_{3}}(\bar{t}_3)\ldots X_{\gamma_{\widetilde{l}}}(\bar{t}_{\widetilde{l}}), 
$$
where $n_3=n_2n_3''\in N_{\Delta_-^{s^1}}$.

We can proceed in a similar way to obtain the following representation
\begin{equation}\label{1?}
X_{\gamma_1}(\bar{t}_1)\ldots X_{\gamma_{\widetilde{l}}}(\bar{t}_{\widetilde{l}})=n\widetilde{x}{s}_{\gamma_1}\ldots{s}_{\gamma_{\widetilde{l}}}=n\widetilde{x}s^1,~n\in N_{\Delta_-^{s^1}},~\widetilde{x}\in N_{\Delta_{(1)}^+},
\end{equation}
where $\Delta_{(1)}^+=\{\alpha\in \Delta_+: \alpha\leq \gamma_{\widetilde{l}}, s^1\alpha=-\alpha\}=\{\alpha\in \Delta_+: \gamma_1\leq \alpha\leq \gamma_{\widetilde{l}}\}$. \index[not]{D@$\Delta_{(1)}^+$} 

Note that $s^1$ acts by multiplication by $-1$ on the roots from $\Delta_{(1)}^+$, so $\Delta_{(1)}^-=-\Delta_{(1)}^+=s^1(\Delta_{(1)}^+)$. Therefore
$N_{\Delta_{(1)}^-}=(s^1)^{-1}N_{\Delta_{(1)}^+}s^1\subset N_{\Delta_-^{s^1}}$ and (\ref{1?}) can be rewritten in the following form
\begin{equation}\label{2?}
X_{\gamma_1}(\bar{t}_1)\ldots X_{\gamma_{\widetilde{l}}}(\bar{t}_{\widetilde{l}})=ns^1(s^1)^{-1}\widetilde{x}s^1=ns^1n',~n\in N_{\Delta_-^{s^1}},~n'=(s^1)^{-1}\widetilde{x}s^1\in N_{\Delta_{(1)}^-}\subset N_{\Delta_-^s\setminus \Delta_0}.
\end{equation}

Similarly one has
\begin{equation}\label{3?}
X_{\gamma_{\widetilde{l}+1}}(\bar{t}_{\widetilde{l}+1})\ldots X_{\gamma_{l'}}(\bar{t}_{l'})=n''{s}_{\gamma_{\widetilde{l}+1}}\ldots {s}_{\gamma_{l'}}n'''=n''s^2n''',~n''\in N_{\Delta_-^{s^2}}\subset N_{\Delta_-^s\setminus \Delta_0},~n'''\in N_{\Delta_{(2)}^-},
\end{equation}
where $\Delta_-^{s^2}=\{\alpha\in \Delta_-:s^2\alpha=-\alpha\}$, and $\Delta_{(2)}^-=\{\alpha \in \Delta_-:-\gamma_{\widetilde{l}+1}\leq \alpha \leq -\gamma_{l'}\}$. \index[not]{D@$\Delta_{(2)}^-$} \index[not]{D@$\Delta_-^{s^2}$}

Substituting (\ref{2?}) and (\ref{3?}) into the formula defining ${\bar{u}}=X_{\gamma_1}(\bar{t}_1)\ldots X_{\gamma_{\widetilde{l}}}(\bar{t}_{\widetilde{l}})X_{\gamma_{\widetilde{l}+1}}(\bar{t}_{\widetilde{l}+1})\ldots X_{\gamma_{l'}}(\bar{t}_{l'})$ one can obtain that for any $n_- \in N_- , y\in N_{\Delta_+\setminus \Delta_{\m_+}}$
\begin{equation}\label{6?}
n_-^{-1}{\bar{u}}y=n_-^{-1}ns^1n'n''s^2n'''y=ks^1gs^2n'''y,~k=n_-^{-1}n\in N_-,~g=n'n''\in N_{\Delta_-^s\setminus \Delta_0}.
\end{equation}


Now note that the minimal segment $\Delta_+^s\setminus \Delta_0$ can be represented as the following disjoint union
\begin{equation}\label{DDec}
\Delta_+^s\setminus \Delta_0=\left(\Delta_+^s\setminus (\Delta_{s^1}^s\cup \Delta_0)\right)\cup \Delta_{s^1}^s.
\end{equation}

Observe that the segments $\Delta_{s^1}^s$, $\Delta_+^s\setminus (\Delta_{s^1}^s\cup \Delta_0)$ are minimal with respect to the circular normal ordering on $\Delta$ associated to $s$. Thus from (\ref{DDec}) and Lemma \ref{comm} we deduce  that the element $g\in N_{\Delta_-^s\setminus \Delta_0}$ from formula (\ref{6?}) can be uniquely factorized as the product
$g=g'g''$, where $g'\in N_{-\left(\Delta_+^s\setminus (\Delta_{s^1}^s\cup \Delta_0)\right)}$, and $g''\in N_{-\Delta_{s^1}^s}$. 

Note that by the definition of $\Delta_{s^1}^s$ we have
\begin{equation}\label{ss2}
s^{1}(\Delta_+^s\setminus (\Delta_{s^1}^s\cup \Delta_0))\subset \Delta_+^s\setminus (\Delta_{s^1}^s\cup \Delta_0)\subset \Delta_+ \setminus \Delta_0,
\end{equation}
and by Proposition \ref{pord} (iii)
\begin{equation}\label{ss1}
s^2\Delta_{s^{1}}^s\subset \Delta_+^s\setminus (\Delta_{s^{1}}^s \cup \Delta_{s^{2}}^s\cup \Delta_0)\subset \Delta_{\m_+}. 
\end{equation}
Therefore by (\ref{ss1}), (\ref{ss2}) we have $(s^2)^{-1}g''s^2\in M_-$ and $s^1g'(s^1)^{-1}\in N_-$, and (\ref{6?}) takes the form
\begin{equation}\label{levdec}
n_-^{-1}{\bar{u}}y=ks^1g'g''s^2n'''y=ks^1g'(s^1)^{-1}s(s^2)^{-1}g''s^2n'''y=k's\widehat{n}y,
\end{equation}
$$
\widehat{n}=(s^2)^{-1}g''s^2n'''\in M_-,~k'=ks^1g'(s^1)^{-1}\in N_-.
$$

The element $\widehat{n}y$ belongs to the subgroup $N_{\Delta'}$, $\Delta':=\{\alpha\in \Delta: \gamma_{l'}<\alpha\leq -\gamma_{l'}\}$. \index[not]{D@$\Delta'$} We have the following disjoint union of minimal segments (see Figure 4)
$$
\Delta'=(\Delta'\cap \Delta_{s^2}^s)\cup (\Delta'\cap (-\Delta_{s^1}^s))\cup (\Delta_0\cap \Delta_+)\cup (-\Delta_{\m_+}). 
$$

Therefore by Lemmas \ref{decNr} and \ref{segmsub} we obtain a unique factorization, $N_{\Delta'}=N_{\Delta'\cap \Delta_{s^2}^s}N_{\Delta'\cap (-\Delta_{s^1}^s)}N_{\Delta_0\cap \Delta_+}M_-$. Note that by the definition of $Z^s$, $N_{\Delta_0\cap \Delta_+}= Z_+=Z^s\cap N_+$. Therefore we have a unique factorization $\widehat{n}y=y'z_+m$, where $z_+\in Z_+$, $y'\in N_{\Delta'\cap \Delta_{s^2}^s}N_{\Delta'\cap (-\Delta_{s^1}^s)}$, $m\in M_-$.
The images of all roots from the set $(\Delta'\cap \Delta_{s^2}^s)\cup (\Delta'\cap (-\Delta_{s^1}^s))=(\Delta'\setminus \Delta_0)\cap \Delta_+$ under the action of $s$ belong to $\Delta_-$ by Proposition \ref{pord} (vii), (iii) and (i). Hence $sN_{\Delta'\cap \Delta_{s^2}^s}N_{\Delta'\cap (-\Delta_{s^1}^s)}s^{-1}\subset N_-$, and $sy's^{-1}\in N_-$. Thus substituting the factorization $\widehat{n}y=y'z_+m$ into (\ref{levdec}) we obtain  
\begin{equation}\label{8?}
n_-^{-1}{\bar{u}}y=k''sz_+m=k''z_+'sm,m\in M_-,~k''=k'sy's^{-1}\in N_-,~z_+,z_+'=sz_+s^{-1}\in Z_+.
\end{equation}
Hence 
\begin{equation}\label{bel}
n_-^{-1}{\bar{u}}y\in N_-sZ_+M_-.
\end{equation}

Next we prove that for any $n_- \in N_-, x \in \h$ and $y\in N_{\Delta_+\setminus \Delta_{\m_+}}$ we have
$h_-^{-1}n_-^{-1}{\bar{u}}yh_+\in N_-sH^0Z_+M_-$, where $h_\pm=e^{r_\pm^s x}$, i.e. $q(\mu_{M_-}^{-1}({\bar{u}})) \subset N_-sH^0Z_+M_-$.

Let $H'\subset H$ \index[not]{H@$H'$} be the subgroup corresponding to the Lie subalgebra $\h'\subset \h$.  We obviously have $H=H'H^0$. From the definition of $r_\pm^s$ it follows that for any $h^0\in H^0$ and $h'\in H'$ the elements $h_+=h^0h'$ and $h_-={h^0}^{-1}s(h')$ are of the form  $h_\pm=e^{r_\pm^s x}$ for some $x\in \h$ and all elements  $h_\pm$ are obtained in this way.

Next observe that the set $N_-sH^0Z_+M_-$ is invariant with respect to
the following action of the subgroup of $H\times H$ formed by elements of the form $(h_+,h_-)=({h^0}h',{h^0}^{-1}s(h'))$:
\begin{equation}\label{tri}
(h_+,h_-)\circ L= h_-^{-1}Lh_+, h= h_+={h^0}h', h_-={h^0}^{-1}s(h').
\end{equation}

Indeed, let $L=vskz_+w,~ v\in N_-,w \in M_-,z_+\in Z_+, k\in H^0$ be an element of $N_-sH^0Z_+M_-$. Then
\begin{equation}\label{hact}
(h_+,h_-)\circ L=h_-^{-1} vh_-h_-^{-1} skh_+h_+^{-1}  z_+wh_+=h_-^{-1} vh_-sk{h^0}^{2}h_+^{-1} z_+ wh_+
\end{equation}
since $s^{-1}h_-^{-1}sh_+={h^0}h'^{-1}{h^0}h'={h^0}^{2}$. The right hand side of the last equality belongs to $N_-sH^0Z_+M_-$ because $H$ normalizes $N_-$, $M_-$ and $Z_+$.

Comparing action (\ref{tri}) with (\ref{dva}) and recalling that by (\ref{bel}) for any $n_- \in N_-$ and $y\in N_{\Delta_+\setminus \Delta_{\m_+}}$ one has $n_-^{-1}{\bar{u}}y\in N_-sZ_+M_-\subset N_-sH^0Z_+M_-$ we deduce 
\begin{equation}\label{bel1}
q(\mu_{M_-}^{-1}({\bar{u}})) \subset N_-sH^0Z_+M_-. 
\end{equation}

Now we show that 
\begin{equation}\label{NNM}
N_-sZ_+M_-=N_-sM_-Z_+=(N_-\cap N^s) Z_-sM_-Z_+ =(N_-\cap N^s)Z_-sZ_+M_-=(N_-\cap N^s)sZ_-Z_+M_-.
\end{equation}

First observe that we have the following disjoint union of minimal segments (see Figure 4)
$$
\Delta_-=((\Delta_-^s\setminus \Delta_0)\cap \Delta_-)\cup (\Delta_0\cap \Delta_-)\cup (\Delta_{s^1}^s\cap \Delta_-).
$$
Therefore by Lemmas \ref{decNr} and \ref{segmsub} we obtain a unique factorization
$$
N_-=N_{(\Delta_-^s\setminus \Delta_0)\cap \Delta_-}N_{\Delta_0\cap \Delta_-}N_{\Delta_{s^1}^s\cap \Delta_-}=(N\cap N_-)Z_-(N_-\cap \overline{N}^s)
$$
as $N_{(\Delta_-^s\setminus \Delta_0)\cap \Delta_-}=N\cap N_-$, $N_{\Delta_0\cap \Delta_-}=Z_-=Z\cap N_-$, $N_{\Delta_{s^1}^s\cap \Delta_-}=N_-\cap \overline{N}^s$.
 
Note that the elements of the subgroup $N_-\cap \overline{N}^s=N_{\Delta_{s^1}^s\cap \Delta_-}$ are transformed to $M_-$ by the conjugation by $s^{-1}$ as by part (iii) of Proposition \ref{pord}
$$
s^{-1}(\Delta_{s^1}^s)=s^2s^1 (\Delta_{s^1}^s)=-s^2(\Delta_{s^1}^s)\subset -(\Delta_+^s\setminus(\Delta_{s^1}^s \cup \Delta_{s^2}^s\cup \Delta_0))\subset -\Delta_{\m_+}.
$$

We deduce that any $k''\in N_-$ can be uniquely represented in the form $k''=\check{n}z_-k'''$, $\check{n}\in N^s\cap N_-$, $z_-\in Z_-=Z^s\cap N_-$, $k'''\in N_-\cap \overline{N}^s$, $s^{-1}k'''s\in M_-$. Therefore for any $m\in M_-, z_+\in Z_+$, we have
$$
k''smz_+=\check{n}z_-sm'z_+, \check{n}\in N^s\cap N_-, m'=s^{-1}k'''sm\in M_-, z_-\in Z_-, 
$$
and hence $N_-sM_-Z_+\subset (N^s \cap N_-)Z_-sM_-Z_+$. The opposite inclusion is obvious. Thus $N_-sM_-Z_+= (N^s \cap N_-)Z_-sM_-Z_+$. 

Similarly one obtains $N_-sZ_+M_-=N_-sM_-Z_+= (N^s \cap N_-)Z_-sZ_+M_-=(N^s \cap N_-)Z_-sZ_+M_-=(N^s \cap N_-)sZ_-Z_+M_-$. This proves (\ref{NNM}).

Since the conjugation action of $H^0$ normalizes $N_-$, $(N_-\cap N^s)$, $Z_\pm$, and $H^0$ is fixed by the conjugation by $s$, one immediately obtains from (\ref{NNM})
$$
H^0N_-sZ_+M_-=N_-sH^0Z_+M_-=N_-sH^0M_-Z_+=(N_-\cap N^s) Z_-sH^0M_-Z_+ 
=
$$
$$
=(N_-\cap N^s)Z_-sH^0Z_+M_-=(N^s\cap N_-)sZ_-H^0Z_+M_-\subset N^ssZ^sN^s,
$$
where the last inclusion follows from the inclusions $Z_-H^0Z_+\subset Z^s$, $M_-\subset N^s$.
This proves the identities in (\ref{preim}) and together with inclusion (\ref{bel1}) establishes (\ref{preim}) completely.

Since by Proposition \ref{crosssect} $N^ssZ^sN^s$ is Zariski closed in $G$, the Zariski closure $\overline{q(\mu_{M_-}^{-1}({\bar{u}}))}$ of $q(\mu_{M_-}^{-1}({\bar{u}}))\subset N^ssZ^sN^s$ in $G$ is contained in $N^ssZ^sN^s$. This completes the proof.

\end{proof}

Now we are in a position to describe the closure $\overline{q(\mu_{M_-}^{-1}({\bar{u}}))}$, the projection $\pi_q(\overline{q(\mu_{M_-}^{-1}({\bar{u}}))})$ and the algebra $W^s(G)$.

\begin{theorem}\label{var}
Suppose that the numbers $\bar{t}_i$ defined in (\ref{defu}) are not equal to zero for all $i$. Let $N_{-\Delta_s'}\subset N_-\cap N_s$ be the subgroup generated by one--parameter subgroups corresponding to the roots from the minimal segment $-\Delta_s'$, where $\Delta_s'=\Delta_s^s\setminus \Delta_{m_+}=\{\alpha \in \Delta_+^s\setminus \Delta_0: \gamma_{l'}<\alpha\}$, \index[not]{D@$\Delta_s'$} and $M_-^s=M_-\cap N_s$. \index[not]{M@$M_-^s$}
Then 

(i) the variety $\overline{q(\mu_{M_-}^{-1}({\bar{u}}))}$ is invariant under conjugations by elements of $M_-$, the conjugation action of $M_-$ on $\overline{q(\mu_{M_-}^{-1}({\bar{u}}))}$ is free, and the quotient $\pi_q(\overline{q(\mu_{M_-}^{-1}({\bar{u}}))})$ is a smooth variety;  

(ii) $\overline{q(\mu_{M_-}^{-1}({\bar{u}}))}= N_-{s}Z^sM_-=N_-{s}Z^sM_-^s$; 

(iii) $\pi_q(\overline{q(\mu_{M_-}^{-1}({\bar{u}}))})\simeq N_{-\Delta_s'}sZ^sM_-^s\simeq \Sigma_s=sZ^sN_s$,  the conjugation action
\begin{equation}\label{cross1}
M_-\times N_{-\Delta_s'}sZ^sM_-^s \rightarrow N_-{s}Z^sM_-
\end{equation}
is an isomorphism of varieties, and hence the algebra ${\mathbb{C}}[\overline{q(\mu_{M_-}^{-1}({\bar{u}}))}]={\mathbb{C}}[N_-{s}Z^sM_-]$ is isomorphic to ${\mathbb{C}}[M_-] \otimes {\mathbb{C}}[N_{-\Delta_s'}sZ^sM_-^s]$;

(iv) The algebra $W^s(G)$ is isomorphic to the algebra of regular
functions on $N_{-\Delta_s'}sZ^sM_-^s$, 
$$
W^s(G)\simeq \mathbb{C}[\pi_q(\overline{q(\mu_{M_-}^{-1}({\bar{u}}))})]\simeq \mathbb{C}[N_{-\Delta_s'}sZ^sM_-^s]\simeq\mathbb{C}[\Sigma_s]\simeq {\mathbb{C}}[N_-{s}Z^sM_-]^{M_-}.
$$ 
Thus the algebra $W_{\mathcal{B}}^s(G)$ is a non--commutative deformation of the algebra of regular functions on the transversal slice $\Sigma_s\simeq N_{-\Delta_s'}sZ^sM_-^s$.
\end{theorem}

\begin{proof}

(i) Firstly, as we observed in Lemma \ref{redpois} the preimage $\mu_{M_-}^{-1}({\bar{u}})$ is locally stable under the (locally defined) dressing action of $M_-$. On the other hand by Proposition \ref{constrt} $q(\mu_{M_-}^{-1}({\bar{u}}))\subset N^ssZ^sN^s$, so by Proposition \ref{crosssect} (i) and Proposition \ref{dressingact}, $q(\mu_{M_-}^{-1}({\bar{u}}))$  is (locally) stable under the action of $M_-\subset N^s$ on $N^ssZ^sN^s$ by conjugations. Since the conjugation action of $N^s$ on $N^ssZ^sN^s$ is free the (locally defined) conjugation action of $M_-$ on $q(\mu_{M_-}^{-1}({\bar{u}}))$ is (locally) free as well.

Now recall that by Proposition \ref{constrt} $\overline{q(\mu_{M_-}^{-1}({\bar{u}}))}\subset N^ssZ^sN^s$.  Since by Proposition \ref{crosssect} (i) the conjugation action of $N^s$ on $N^ssZ^sN^s$ is free and regular, $sZ^sN_s$ being a cross--section for this action, and  $\overline{q(\mu_{M_-}^{-1}({\bar{u}}))}$ is closed, the local action of $M_-$ on $q(\mu_{M_-}^{-1}({\bar{u}}))\subset N^ssZ^sN^s$ by conjugations extends by continuity to the genuine regular action of $M_-\subset N^s$ on $\overline{q(\mu_{M_-}^{-1}({\bar{u}}))}$ which is free as well.  Therefore the quotient $\pi_q(\overline{q(\mu_{M_-}^{-1}({\bar{u}}))})$ is a smooth variety.

(ii) We show now that the closure of $q(\mu_{M_-}^{-1}({\bar{u}}))$ contains $N_-sZ^sM_-$.
Recall that by (\ref{preim}) $q(\mu_{M_-}^{-1}({\bar{u}}))\subset N_-sH^0Z_+M_-=N_-H^0Z_+sM_-$.
From (\ref{dva}) it follows that $q(\mu_{M_-}^{-1}({\bar{u}}))$ is closed with respect to the right multiplication by arbitrary elements from $Z_+$ and with respect to the left multiplication by arbitrary elements from $N_-$, as $Z_+\subset N_{\Delta_+\setminus \Delta_{\m_+}}$, and $q(\mu_{M_-}^{-1}({\bar{u}}))$ is closed with respect to the right multiplication by arbitrary elements from $N_{\Delta_+\setminus \Delta_{\m_+}}$. $q(\mu_{M_-}^{-1}({\bar{u}}))$ is also closed with respect to the restriction of action (\ref{tri}) to the subgroup of $H\times H$ which consists of elements of the form $({h^0},{h^0}^{-1})$, ${h^0}\in H^0$. Thus by (\ref{hact}), and since $Z_+$ normalizes $M_-$ and $s$ centralizes $H^0$, $q(\mu_{M_-}^{-1}({\bar{u}}))$ contains elements of the form $ksn$ for some $n\in M_-\subset N^s$ and arbitrary $k\in N_-H^0Z_+$.  

Now recall that $\bar{h}_0(\alpha)>0$ for $\alpha\in \Delta_+^s\setminus \Delta_0$ and $\bar{h}_0(\alpha)=0$ for $\alpha\in \Delta_0$ (see formula (\ref{barh0}) and the discussion after it), and hence the $\mathbb{C}^*$--action on $G$ induced by conjugations by the elements $h(t)$ from the one--parameter subgroup generated by $\bar{h}_0\in \h'$ is contracting on $N^s$ and fixes all elements of $Z^s$. Applying action (\ref{tri}) with $h=h(t)$ to the elements $ksn$ with arbitrary $k\in N_-H^0Z_+$ we immediately deduce, with the help of (\ref{tri}), that the $M_-$--component $n$ can be contracted to the identity element using the above defined contracting action, and hence the Zariski closure $\overline{q(\mu_{M_-}^{-1}({\bar{u}}))}$ of $q(\mu_{M_-}^{-1}({\bar{u}}))$ contains the set $N_-H^0Z_+s=N_-Z_-H^0Z_+s$, where the last identity holds since $Z_-\subset N_-$. Hence $\overline{q(\mu_{M_-}^{-1}({\bar{u}}))}$ also contains $N_-Z^ss$ as $N_-Z^ss\simeq N_-Z^s$ is the orbit of the identity element under the regular action of $N_-\times Z^s$ on $G$ defined by $(n_-,z)\circ g=n_-gz^{-1}$, and the Zariski closure of $Z_-H^0Z_+\subset Z^s$ is obviously $Z^s$.

By (\ref{dva}) $q(\mu_{M_-}^{-1}({\bar{u}}))$ is closed with respect to the left multiplication by arbitrary elements from $N_-$. Recall also that $M_-\subset N^s$ freely acts on $\overline{q(\mu_{M_-}^{-1}({\bar{u}}))}$ by conjugations by part (i). Therefore $\overline{q(\mu_{M_-}^{-1}({\bar{u}}))}$ also contains the set $N_-Z^ssM_-=N_-sZ^sM_-$.  

Note that 
$$
N_-sH^0Z_+M_-\subset N_-Z^ssM_- \subset \overline{q(\mu_{M_-}^{-1}({\bar{u}}))} \subset \overline{N_-sH^0Z_+M_-},
$$
where the last inclusion follows from (\ref{preim}), and all Zariski closures are taken in $G$.
This implies, after taking Zariski closures in $G$, that
\begin{equation}\label{qclos}
\overline{N_-Z^ssM_-}= \overline{q(\mu_{M_-}^{-1}({\bar{u}}))} = \overline{N_-sH^0Z_+M_-}.
\end{equation}

\begin{lemma}\label{lclos}
$N_-Z^ssM_-$ is a closed irreducible subvariety of $G$ of dimension ${\rm dim}~\Sigma_s+{\rm dim}~M_-$.
\end{lemma} 

\begin{proof}
First note that $\Delta_-=((\Delta_-^s\setminus \Delta_0)\cap \Delta_-)\cup (\Delta_+^s\cap \Delta_-)\cup (\Delta_-\cap \Delta_0)$ (disjoint union of minimal segments).
Now by Lemma \ref{comm} there is a unique factorization 
$$
N_-=N_{\Delta_-}=N_{(\Delta_-^s\setminus \Delta_0)\cap \Delta_-}N_{\Delta_+^s\cap \Delta_-}N_{\Delta_-\cap \Delta_0}=(N_-\cap N^s) N_{\Delta_+^s\cap \Delta_-}Z_-
$$ 
as $N_{(\Delta_-^s\setminus \Delta_0)\cap \Delta_-}=N_-\cap N^s$ and $N_{\Delta_-\cap \Delta_0}=Z^s\cap N_-=Z_-$. We deduce
\begin{equation}\label{NZM}
N_-Z^ssM_-=(N_-\cap N^s) N_{\Delta_+^s\cap \Delta_-}Z_-Z^ssM_-=(N_-\cap N^s) N_{\Delta_+^s\cap \Delta_-}Z^ssM_-.
\end{equation}

From the definition of $\Delta_+$ and $\Delta_+^s$ (see formula (\ref{NO}) and Definition \ref{circorddef}) it follows that $\Delta_+^s\cap \Delta_-\subset \Delta_{s^1}^s$. Therefore 
\begin{equation}\label{Ns+}
N_{\Delta_+^s\cap \Delta_-}\subset N_{\Delta_{s^1}^s}.
\end{equation} 

Recall that $Z^s$ is generated by $H^0$ and by the one-parameter subgroups corresponding to the roots from $\Delta_0$. Note that $H^0$ normalizes $N_{\Delta_{s^1}^s}$, and by Proposition \ref{pord} (ii) the one-parameter subgroups corresponding to the roots from $\Delta_0$ also normalize $N_{\Delta_{s^1}^s}$. Therefore $Z^s$ normalizes $N_{\Delta_{s^1}^s}$, and we deduce using (\ref{Ns+})
\begin{equation}\label{ND}
N_{\Delta_+^s\cap \Delta_-}Z^s\subset N_{\Delta_{s^1}^s}Z^s=Z^sN_{\Delta_{s^1}^s}
\end{equation}

By part (iii) of Proposition \ref{pord}
$$
s^{-1}(\Delta_{s^1}^s)=s^2s^1 (\Delta_{s^1}^s)=-s^2(\Delta_{s^1}^s)\subset -(\Delta_+^s\setminus(\Delta_{s^1}^s \cup \Delta_{s^2}^s\cup \Delta_0))\subset -\Delta_{\m_+},
$$
and hence, recalling (\ref{ND}), we obtain 
$$
N_{\Delta_+^s\cap \Delta_-}Z^ss\subset Z^sss^{-1}N_{\Delta_{s^1}^s}s=Z^ssN_{s^{-1}(\Delta_{s^1}^s)}\subset Z^ssM_-.
$$
Together with (\ref{NZM}) this implies (compare with (\ref{preim}))
$$
N_-Z^ssM_-=(N_-\cap N) Z^ssM_-.
$$

Now consider the group $N_{(\Delta_-^s\setminus \Delta_0)\cap \Delta_-}=N_-\cap N^s$. By the definition of $\Delta_+$ and $\Delta_+^s$ (see formula \ref{NO} and Definition \ref{circorddef}) $(\Delta_-^s\setminus \Delta_0)\cap \Delta_-=((-\Delta_{s^{-1}}^s)\cap \Delta_-)\cup (\Delta_-^s\setminus ((-\Delta_{s^{-1}}^s)\cup \Delta_0))$ (disjoint union of additively closed subsets of roots).  

Now by Lemma \ref{decNr} there is a unique factorization 
$$
N_-\cap N^s=N_{(\Delta_-^s\setminus \Delta_0)\cap \Delta_-}=N_{(-\Delta_{s^{-1}}^s)\cap \Delta_-}N_{\Delta_-^s\setminus ((-\Delta_{s^{-1}}^s)\cup \Delta_0)}.
$$ 
We deduce
\begin{equation}\label{NZM1}
N_-Z^ssM_-=(N_-\cap N^s) N_{\Delta_+^s\cap \Delta_-}Z_-Z^ssM_-=N_{(-\Delta_{s^{-1}}^s)\cap \Delta_-}N_{\Delta_-^s\setminus ((-\Delta_{s^{-1}}^s)\cup \Delta_0)} Z^ssM_-.
\end{equation}

Recall that $Z^s$ is generated by $H^0$ and by the one-parameter subgroups corresponding to the roots from $\Delta_0$. Note that $H^0$ normalizes $N_{\Delta_-^s\setminus ((-\Delta_{s^{-1}}^s)\cup \Delta_0)}$, and by Proposition \ref{pord} (ii) the one-parameter subgroups corresponding to the roots from $\Delta_0$ also normalize $N_{\Delta_-^s\setminus ((-\Delta_{s^{-1}}^s)\cup \Delta_0)}$. Therefore $Z^s$ normalizes $N_{\Delta_-^s\setminus ((-\Delta_{s^{-1}}^s)\cup \Delta_0)}$, and we deduce 
\begin{equation}\label{ND1}
N_-Z^ssM_-=N_{(-\Delta_{s^{-1}}^s)\cap \Delta_-}Z^sN_{\Delta_-^s\setminus ((-\Delta_{s^{-1}}^s)\cup \Delta_0)} sM_-
\end{equation}

By the definition of the sets $\Delta_{s^{-1}}^s$, $\Delta_s^s$ and $\Delta_0$ we have 
$$
s^{-1}(\Delta_-^s\setminus ((-\Delta_{s^{-1}}^s)\cup \Delta_0))=\Delta_-^s\setminus ((-\Delta_s^s)\cup \Delta_0),
$$
and hence, recalling (\ref{ND1}), we obtain 
$$
N_-Z^ssM_-=N_{(-\Delta_{s^{-1}}^s)\cap \Delta_-}Z^sss^{-1}N_{\Delta_-^s\setminus ((-\Delta_{s^{-1}}^s)\cup \Delta_0)} sM_-=N_{(-\Delta_{s^{-1}}^s)\cap \Delta_-}Z^ssN_{s^{-1}(\Delta_-^s\setminus ((-\Delta_{s^{-1}}^s)\cup \Delta_0))} M_-=
$$
$$
=N_{(-\Delta_{s^{-1}}^s)\cap \Delta_-}Z^ssN_{\Delta_-^s\setminus ((-\Delta_s^s)\cup \Delta_0)}M_-=N_{(-\Delta_{s^{-1}}^s)\cap \Delta_-}sZ^sN_{-[\beta_1^1,\gamma_{l'}]},
$$
where $N_{\Delta_-^s\setminus ((-\Delta_s^s)\cup \Delta_0)}M_-=N_{-[\beta_1^1,\gamma_{l'}]}\subset N^s$, and $[\beta_1^1,\gamma_{l'}]=\{\alpha\in \Delta_+^s:\alpha\leq \gamma_{l'}\}$ (we use the notation of Proposition \ref{pord} (iv)).

Note that $N_{(-\Delta_{s^{-1}}^s)\cap \Delta_-}=N_{s^{-1}}\cap N_-$, where $N_{s^{-1}}=\{n\in N^s: s^{-1}ns\in \overline{N}^s\}$, so we have
$$
N_-Z^ssM_-=(N_{s^{-1}}\cap N_-)sZ^sN_{-[\beta_1^1,\gamma_{l'}]}.
$$

Finally observe that 
$$
s^{-1}N_-Z^ssM_-=s^{-1}(N_{s^{-1}}\cap N_-)sZ^sN_{-[\beta_1^1,\gamma_{l'}]}\subset \overline{N}^sZ^sN^s.
$$
Therefore by Lemma \ref{NNclosed} $s^{-1}N_-Z^ssM_-$ is a closed subvariety of $\overline{N}^sZ^sN^s$ as
$$
\overline{N}^sZ^sN^s\simeq \overline{N}^s\times Z^s\times N^s,
$$
$s^{-1}(N_{s^{-1}}\cap N_-)s\subset \overline{N}^s$, $N_{-[\beta_1^1,\gamma_{l'}]}\subset N^s$ are closed algebraic subgroups, and $s^{-1}N_-Z^ssM_-$ is the image in $\overline{N}^sZ^sN^s$ of 
$s^{-1}(N_{s^{-1}}\cap N_-)s\times Z^s \times N_{-[\beta_1^1,\gamma_{l'}]}$ under the product map. In particular,
\begin{equation}\label{NZsM}
N_-Z^ssM_-\simeq s^{-1}(N_{s^{-1}}\cap N_-)s\times Z^s \times N_{-[\beta_1^1,\gamma_{l'}]}.
\end{equation}

The variety $\overline{N}^sZ^sN^s$ is closed in $G$ by Lemma \ref{NNclosed}, and hence $s^{-1}N_-Z^ssM_-$ is a closed subvariety of $G$. Thus $N_-Z^ssM_-$ is a closed subvariety of $G$.

From the isomorphism of varieties (\ref{NZsM}) using Theorem 1.5.4 in \cite{Sp} it follows that $N_-Z^ssM_-$ is also irreducible as $s^{-1}(N_{s^{-1}}\cap N_-)s$, $Z^s$, and $N_{-[\beta_1^1,\gamma_{l'}]}$ are irreducible.

To find the dimension of $N_-Z^ssM_-$ we observe that  by formula \ref{NO}, by Definition \ref{circorddef} and by part (vii) of Proposition \ref{pord} $(-\Delta_{s^{-1}}^s)\cap \Delta_-=(-\Delta_{s^{-1}}^s)\setminus \{\alpha\in \Delta_-^s: \alpha <-\gamma_1\}$, and $-[\beta_1^1,\gamma_{l'}]=(-\Delta_{\m_+})\cup \{\alpha\in \Delta_-^s: \alpha <-\gamma_1\}$ (disjoint union). Note also that $\{\alpha\in \Delta_-^s: \alpha <-\gamma_1\}\subset -\Delta_{s^{-1}}^s$ 
Thus  
$$
{\rm dim}~N_{s^{-1}}\cap N_-={\rm dim}~N_{(-\Delta_{s^{-1}}^s)\cap \Delta_-}=|(-\Delta_{s^{-1}}^s)\setminus \{\alpha\in \Delta_-^s: \alpha <-\gamma_1\}|=|-\Delta_{s^{-1}}^s|-| \{\alpha\in \Delta_-^s: \alpha <-\gamma_1\}|,
$$
and
$$
{\rm dim}~N_{-[\beta_1^1,\gamma_{l'}]}=|\{\alpha\in \Delta_-^s: \alpha <-\gamma_1\}|+|-\Delta_{\m_+}|
$$
Now using isomorphism (\ref{NZsM}) we obtain
$$
{\rm dim}~N_-Z^ssM_-={\rm dim}~N_{s^{-1}}\cap N_-+{\rm dim}~Z^s+{\rm dim}~N_{-[\beta_1^1,\gamma_{l'}]}=
$$
$$
=|-\Delta_{s^{-1}}^s|-| \{\alpha\in \Delta_-^s: \alpha <-\gamma_1\}|+|\{\alpha\in \Delta_-^s: \alpha <-\gamma_1\}|+|-\Delta_{\m_+}|+{\rm dim}~Z^s=
$$
$$
=|\Delta_{s^{-1}}^s|+{\rm dim}~Z^s+{\rm dim}~M_-={\rm dim}~N_s+{\rm dim}~Z^s+{\rm dim}~M_-={\rm dim}~\Sigma_s+{\rm dim}~M_-.
$$
This completes the proof of the lemma.

\end{proof}

Now from (\ref{qclos}) and Lemma \ref{lclos} we infer
$$
N_-Z^ssM_-= \overline{q(\mu_{M_-}^{-1}({\bar{u}}))} = \overline{N_-sH^0Z_+M_-}.
$$

Next, we show that $N_-Z^ssM_-=N_-Z^ssM_-^s$. Observe that $M_-=N_{-\Delta_{\m_+}}$, and by the definition of $\Delta_{\m_+}$ (see part (v) of Proposition \ref{pord}), and by part (vii) of Proposition \ref{pord}, $\Delta_{\m_+}=\Delta_{\m_+}^s\cup (\Delta_{\m_+}\setminus \Delta_s^s)$ (disjoint union of minimal segments), where $\Delta_{\m_+}^s=\Delta_{\m_+}\cap \Delta_s^s$. \index[not]{D@$\Delta_{\m_+}^s$}

Using Lemma \ref{comm} we have a unique decomposition 
\begin{equation}\label{M-s}
M_-=N_{-(\Delta_{\m_+}\setminus \Delta_s^s)}N_{-\Delta_{\m_+}^s}=N_{-(\Delta_{\m_+}\setminus \Delta_s^s)}M_-^s
\end{equation} 
as by the definition $M_-^s=M_-\cap N_s=N_{-\Delta_{\m_+}^s}$. If $\alpha\in \Delta_{\m_+}\setminus \Delta_s^s$ then $s\alpha\in \Delta_+^s$ by the definition of $\Delta_s^s$, and by Proposition \ref{pord} (viii)  we have $s\alpha >\alpha$, and if $s\alpha +\alpha_0 \in \Delta$  for $\alpha_0\in \Delta_0$ then $s\alpha +\alpha_0 \in \Delta_+^s$ and $s\alpha+\alpha_0>\alpha$, so in both cases $s\alpha, s\alpha+\alpha_0 \in \Delta_+^s$ and $s\alpha, s\alpha+\alpha_0>\alpha$. This implies $s\alpha, s\alpha+\alpha_0 \not\in \Delta_s^s$. In particular,
from the definition (\ref{circorddef}) of the root system $\Delta_+$ it follows that in both cases $s\alpha, s\alpha+\alpha_0\in \Delta_+$. Observing also that $Z^s$ is generated by the one--parameter subgroups corresponding to roots from $\Delta_0$ and by the connected component of the centralizer of $s$ in $H$ which normalizes all one--parameter subgroups corresponding to roots, we deduce $sN_{-(\Delta_{\m_+}\setminus \Delta_s^s)}s^{-1}\subset N_-$ and $zsN_{-(\Delta_{\m_+}\setminus \Delta_s^s)}s^{-1}z^{-1}\subset N_-$ for any $z\in Z^s$. Thus by (\ref{M-s})
$$
N_-Z^ssM_-=N_-Z^ssN_{-(\Delta_{\m_+}\setminus \Delta_s^s)}M_-^s=N_-Z^ssN_{-(\Delta_{\m_+}\setminus \Delta_s^s)}s^{-1}sM_-^s\subset N_-sZ^sM_-^s.
$$
The opposite inclusion is obvious. So we obtain that $N_-Z^ssM_-= N_-Z^ssM_-^s$.

(iii) We show that $N^sZ^ssM_-^s$ is a cross--section for the free conjugation action of $M_-$ on $N_-Z^ssM_-= \overline{q(\mu_{M_-}^{-1}({\bar{u}}))}$.  Note that by the definition $N^ssZ^sM_-^s\subset N_-Z^ssM_-= \overline{q(\mu_{M_-}^{-1}({\bar{u}}))}$.

From Figure 4 and Proposition \ref{pord} (i), (iii), (vii) we obtain $-\Delta_s^s=(-\Delta_s')\cup -(\Delta_{\m_+}^s)$ (disjoint union of minimal segments). By the definition $N_s=N_{-\Delta_s^s}$, $M_-^s=N_{-\Delta_{\m_+}^s}$. Therefore Lemma \ref{comm} implies the following unique factorization $N_s=M_-^sN_{-\Delta_s'}$. Thus if $szn_s\in sZ^sN_s$, $z\in Z^s$, $n_s\in N_s$ then $n_s$ can be uniquely factorized as $n_s=m_sn_s'$, $m_s\in M_-^s$, $n_s'\in N_{-\Delta_s'}$ and we have 
$$
szn_s=szm_sn_s'.
$$
Conjugating this element by $n_s'$ we deduce that $szn_s$ is uniquely conjugated to the element
$$
n_s'szm_s\in N_{-\Delta_s'}sZ^sM_-^s,
$$
and hence $N_{-\Delta_s'}sZ^sM_-^s\simeq sZ^sN_s=\Sigma_s$ is a cross--section for the conjugation action of $N^s$ on $N^ssZ^sN^s$ as well. At the same time by construction the bijection $N_{-\Delta_s'}sZ^sM_-^s\simeq sZ^sN_s$ is an isomorphism of varieties, where the variety structure on $N_{-\Delta_s'}sZ^sM_-^s$ is induced from $N_{-\Delta_s'}\times Z^s\times M_-^s$ using the product map 
$$
N_{-\Delta_s'}\times Z^s\times M_-^s\to N_{-\Delta_s'}sZ^sM_-^s,~(n_s,z,m_s)\mapsto n_sszm_s.
$$

Observe that any two points of $N_{-\Delta_s'}Z^ssM_-^s$ are not $M_-$--conjugate. Indeed, we have an inclusion $N_{-\Delta_s'}Z^ssM_-^s \subset \overline{q(\mu_{M_-}^{-1}({\bar{u}}))}$, and two points of $\overline{q(\mu_{M_-}^{-1}({\bar{u}}))}$ cannot be $M_-$-conjugate if they are not $N^s$--conjugate in $N^s{s}Z^sN^s\supset \overline{q(\mu_{M_-}^{-1}({\bar{u}}))}$ as $M_-\subset N^s$.  But $N_{-\Delta_s'}Z^ssM_-^s\simeq \Sigma_s$ is a cross--section for the conjugation action of $N^s$ on $N^ssZ^sN^s$ by Proposition \ref{cross} (i). Thus any two points of $N_{-\Delta_s'}Z^ssM_-^s$ are not $N^s$--conjugate, and hence they are not $M_-$--conjugate.
Therefore the closed variety $\pi_q(\overline{q(\mu_{M_-}^{-1}({\bar{u}}))})$ must contain the closed variety $N_{-\Delta_s'}Z^ssM_-^s\simeq \Sigma_s$.  

From formula (\ref{dimm}) for the cardinality $|\Delta_{\m_+}|$ of the set $\Delta_{\m_+}$ and from the definitions of $\overline{q(\mu_{M_-}^{-1}({\bar{u}}))}$ and of $N_{-\Delta_s'}Z^ssM_-^s$ we deduce that the dimension of the quotient $\pi_q(\overline{q(\mu_{M_-}^{-1}({\bar{u}}))})$ is equal to the dimension of the variety $N_{-\Delta_s'}Z^ssM_-^s$,
\begin{eqnarray*}
{\rm dim}~\pi_q(\overline{q(\mu_{M_-}^{-1}({\bar{u}}))})={\rm dim}~G-2{\rm dim}~M_-=2D+l-2| \Delta_{\m_+}| =2D+l -2\left(D-\frac{l(s)-l'}{2}-D_0\right)= \\
=l(s)+2D_0+l-l'={\rm dim}~N_{s}+{\rm dim}~Z^s={\rm dim}~sZ^sN_s={\rm dim}~\Sigma_s={\rm dim}~N_{-\Delta_s'}Z^ssM_-^s.
\end{eqnarray*}

Since $\pi_q$ is a morphism of varieties and the conjugation action of $M_-$ on $\overline{q(\mu_{M_-}^{-1}({\bar{u}}))}=N_-Z^ssM_-$ is free by part (i), $\pi_q^{-1}(N_{-\Delta_s'}Z^ssM_-^s)$ is a closed smooth subvariety of the smooth variety $N_-Z^ssM_-= \overline{q(\mu_{M_-}^{-1}({\bar{u}}))}$, and   ${\rm dim}~\pi_q^{-1}(N_{-\Delta_s'}Z^ssM_-^s)={\rm dim}~N_{-\Delta_s'}Z^ssM_-^s+{\rm dim}~M_-={\rm dim}~\Sigma_s+{\rm dim}~M_-$.  

Now recall that by Lemma \ref{lclos} $N_-Z^ssM_-= \overline{q(\mu_{M_-}^{-1}({\bar{u}}))}$ is irreducible and has the same dimension, ${\rm dim}~N_-Z^ssM_-={\rm dim}~\Sigma_s+{\rm dim}~M_-={\rm dim}~\pi_q^{-1}(N_{-\Delta_s'}Z^ssM_-^s)$.  By \cite{Har}, Ex. 1.10(c) the identity for the dimensions and the closed embedding $\pi_q^{-1}(N_{-\Delta_s'}Z^ssM_-^s)\subset N_-Z^ssM_-$ imply $\pi_q^{-1}(N_{-\Delta_s'}Z^ssM_-^s)=N_-Z^ssM_-$.

Therefore $\pi_q(\overline{q(\mu_{M_-}^{-1}({\bar{u}}))})\simeq N_{-\Delta_s'}Z^ssM_-^s$, $N_{-\Delta_s'}Z^ssM_-^s$ is a cross--section for the action of $M_-$ on $\overline{q(\mu_{M_-}^{-1}({\bar{u}}))}$, and the conjugation action
$$
M_-\times N_{-\Delta_s'}Z^ssM_-^s \rightarrow N_-Z^ssM_-
$$
is an isomorphism of varieties. We conclude that the algebra ${\mathbb{C}}[\overline{q(\mu_{M_-}^{-1}({\bar{u}}))}]$ is isomorphic to ${\mathbb{C}}[M_-] \otimes {\mathbb{C}}[N^sZ^ssM_-^s]$, ${\mathbb{C}}[\overline{q(\mu_{M_-}^{-1}({\bar{u}}))}]\simeq {\mathbb{C}}[M_-] \otimes {\mathbb{C}}[N_{-\Delta_s'}Z^ssM_-^s]$.

(iv) Now recall that by Lemma \ref{redreg}
$$
W^s(G)={\mathbb{C}}[\overline{q(\mu_{M_-}^{-1}({\bar{u}}))}]\cap
C^\infty(\mu_{M_-}^{-1}({\bar{u}}))^{C^\infty(M_+)},
$$
where ${\mathbb{C}}[\overline{q(\mu_{M_-}^{-1}({\bar{u}}))}]$ is regarded as a subalgebra in $C^\infty(\mu_{M_-}^{-1}({\bar{u}}))$ using the map $q^*:C^\infty(q(\mu_{M_-}^{-1}({\bar{u}})))\rightarrow C^\infty(\mu_{M_-}^{-1}({\bar{u}}))$ and the imbedding ${\mathbb{C}}[\overline{q(\mu_{M_-}^{-1}({\bar{u}}))}]\subset C^\infty(q(\mu_{M_-}^{-1}({\bar{u}})))$.

By Lemma \ref{redpois}  the algebra $C^\infty(\mu_{M_-}^{-1}({\bar{u}}))^{M_-}$ is isomorphic to $C^\infty(\mu_{M_-}^{-1}({\bar{u}}))^{C^\infty(M_+)}$, and hence
\begin{equation}\label{Wad}
W^s(G)={\mathbb{C}}[\overline{q(\mu_{M_-}^{-1}({\bar{u}}))}]\cap
C^\infty(\mu_{M_-}^{-1}({\bar{u}}))^{C^\infty(M_+)}={\mathbb{C}}[\overline{q(\mu_{M_-}^{-1}({\bar{u}}))}]\cap
C^\infty(\mu_{M_-}^{-1}({\bar{u}}))^{M_-}.
\end{equation}

As we already proved the variety $\overline{q(\mu_{M_-}^{-1}({\bar{u}}))}$ is stable under the conjugation action of $M_-$, and the map $\pi_q: \overline{q(\mu_{M_+}^{-1}({\bar{u}}))}\rightarrow
\pi_q\overline{q(\mu_{M_+}^{-1}({\bar{u}}))}$
is a morphism of varieties. Moreover, under the map $q:G^*\to G$ the local dressing action of $M_-$ on $G^*$ becomes the conjugation action on $G$. Therefore the map
\begin{equation}\label{Wad1}
\mathbb{C}[\pi_q\overline{q(\mu_{M_-}^{-1}({\bar{u}}))}]\rightarrow {\mathbb{C}}[\overline{q(\mu_{M_-}^{-1}({\bar{u}}))}]\cap
C^\infty(\mu_{M_-}^{-1}({\bar{u}}))^{M_-},~~\psi \mapsto \pi_q^*\psi
\end{equation}
is an algebra isomorphism, where ${\mathbb{C}}[\overline{q(\mu_{M_-}^{-1}({\bar{u}}))}]$ is regarded as a subalgebra in $C^\infty(\mu_{M_-}^{-1}({\bar{u}}))$ using the map $$q^*:C^\infty(q(\mu_{M_+}^{-1}({\bar{u}})))\rightarrow C^\infty(\mu_{M_-}^{-1}({\bar{u}}))$$ and the imbedding ${\mathbb{C}}[\overline{q(\mu_{M_-}^{-1}({\bar{u}}))}]\subset C^\infty(q(\mu_{M_-}^{-1}({\bar{u}})))$.

From (\ref{Wad}) and (\ref{Wad1}) it follows that $W^s(G)\simeq \mathbb{C}[\pi_q\overline{q(\mu_{M_-}^{-1}({\bar{u}}))}]\simeq {\mathbb{C}}[N_-{s}Z^sM_-]^{M_-}$.
This completes the proof.

\end{proof}


\section{Zhelobenko type operators for Poisson q-W--algebras}\label{pZhel}

\pagestyle{myheadings}
\markboth{CHAPTER~\thechapter.~Q-W--ALGEBRAS}{\thesection.~ZHELOBENKO TYPE OPERATORS FOR POISSON Q-W--ALGEBRAS}

\setcounter{equation}{0}
\setcounter{theorem}{0}

In this section we present the main result of this chapter, a formula for a projection operator $\Pi: \mathbb{C}[N_-Z^ssM_-]\rightarrow \mathbb{C}[N_-Z^ssM_-]^{M_-}$ onto the subspace of invariants $\mathbb{C}[N_-Z^ssM_-]^{M_-}$ which is isomorphic to $W^s(G)$ as an algebra according to Theorem \ref{var} (iv). This formula has a direct quantum analogue which will be introduced in the next chapter. The definitions and the constructions presented in this section have purely quantum group motivations which will be clear in the next chapter.

The operator $\Pi$ can be defined following the philosophy of \cite{S13} where a similar projection operator onto the subspace $\mathbb{C}[N^sZ^ssN^s]^{N^s}\subset \mathbb{C}[N^sZ^ssN^s]$ was defined and studied.  
More precisely, according to Theorem \ref{var} (iii) any $g\in N_-Z^ssM_-$ can be uniquely represented in the form 
\begin{equation}\label{gslice}
g=nn_szsm_sn^{-1}, n\in M_-, n_s\in N_{-\Delta_s'}, m_s\in M_-^s, z\in Z^s.
\end{equation}
If for $f\in \mathbb{C}[N_-Z^ssM_-]$ we define $\Pi f\in \mathbb{C}[N_-Z^ssM_-]$ by
\begin{equation}\label{Pidef}
(\Pi f)(g)=f(n^{-1}gn)=f(n_szsm_s) \index[not]{P@$\Pi(~\cdot~)$}
\end{equation}
then $\Pi f$ is an $M_-$--invariant function, and any $M_-$--invariant regular function on $N_-Z^ssM_-$ can be obtained this way. Moreover, by the definition $\Pi^2=\Pi$, i.e. $\Pi$ is a projection onto $\mathbb{C}[N_-Z^ssM_-]^{M_-}$. \index{projection operator}

To obtain an explicit formula for the operator $\Pi$ suitable for quantization we firstly find an explicit formula for $n$ in terms of $g$ in (\ref{gslice}). Denote by $\omega$ the Chevalley anti-involution \index{Chevalley anti-involution} on $\g$ which is induced by the anti-involution $\omega$ of $U_h^s(\g)$ on $U(\g)\simeq U_h^s(\g)/hU_h^s(\g)$. We also denote the corresponding anti-involution of $G$ by the same letter. \index[not]{o@$\omega$}

An explicit formula for $\Pi$ suitable for quantization will be given in terms of matrix elements of finite--dimensional irreducible representations of $G$. As before, let $H\subset G$ be a maximal torus, $\Delta_+\subset \Delta=\Delta(G,H)$ a system of positive roots, $P_+$ the corresponding set of integral dominant weights, $B_+$ the corresponding Borel subgroup, $B_-$ the opposite Borel subgroup, $N_\pm$ the unipotent radicals of $B_\pm$, respectively. For any $\lambda\in P_+$ we denote by $v_\lambda$ \index[not]{v@$v_\lambda$} a non--zero highest weight vector \index{representation!of a semisimple algebraic group!highest weight vector} in the irreducible highest weight $G$--module $V_\lambda$ \index[not]{V@$V_\lambda$} \index{representation!of a semisimple algebraic group!irreducible highest weight} of highest weight $\lambda$, and by $(~\cdot~,~\cdot~)$ \index{form!contravariant non--degenerate!on a representation of a semisimple algebraic group} \index[not]{ZZ@$(~\cdot~ ,~\cdot~ )$} the contravariant non--degenerate form on $V_\lambda$ such that $(v,xw)=(\omega(x)v,w)$ for any $v,w\in V_\lambda$, $x\in \g$, and $(v_\lambda,v_\lambda)=1$.

Let $\beta_1,\ldots, \beta_D$ be a normal order on $\Delta_+$, $X_{\pm \beta_1},\ldots, X_{\pm \beta_D}\in \g$ the corresponding root vectors defined in (\ref{rootvectg}). We shall need some special matrix elements of finite-dimensional irreducible representations of $G$. These matrix elements can be defined by specializing the results of Lemma 2.3 in \cite{DKPschubert} at $q=1$. By this lemma there are integral dominant weights $\mu_p\in P_+$, $p=1,\ldots ,D$ \index[not]{m@$\mu_p$}  and elements $v_p\in V_{\mu_p}$ \index[not]{v@$v_p$} such that
\begin{equation}\label{vpkcl}
(v_p,X_{-\beta_{D}}^{(m_D)}\ldots X_{-\beta_1}^{(m_1)}v_{\mu_p})=\left\{ \begin{array}{ll} 1 & {\rm if}~X_{-\beta_{D}}^{(m_D)}\ldots X_{-\beta_1}^{(m_1)}= X_{-\beta_{p}}\\ 0 & {\rm otherwise}
\end{array}
\right. ,
\end{equation}
where for $\alpha\in \Delta$, $k\in \mathbb{N}$ we define $X_\alpha^{(k)}=\frac{X_\alpha^{k}}{k!}$. \index[not]{X@$X_\alpha^{(k)}$} 

Since by this definition $(v_p,X_{-\beta_p}v_{\mu_p})=(\omega(X_{-\beta_p})v_p,v_{\mu_p})=1$, $(v_{\mu_p},v_{\mu_p})=1$, different weight subspaces of $V_{\mu_p}$ are orthogonal with respect to the contravariant non--degenerate form, and the highest weight subspace in $V_{\mu_p}$ is one--dimensional, we deduce that 
\begin{equation}\label{x-bvm}
\omega(X_{-\beta_p})v_p=v_{\mu_p}.
\end{equation}

We shall need the following properties of the matrix elements $(v_p,\cdot~ v_{\mu_p})$.
\begin{lemma}\label{vancl}
(i) For any $1\leq q<p\leq D$, any $y=X_{-\beta_{D}}^{(m_D)}\ldots X_{-\beta_p}^{(m_p)}\in U(\g)$, where $m_i\in \mathbb{N}$, $y\neq 1$, one has $\omega(y)v_q=0$. In particular, for any $u\in U(\g)$
\begin{equation}\label{aq1cl}
(v_q, yuv_{\mu_q})=0.
\end{equation}

(ii) For any $1\leq p\leq D$, any $y=X_{-\beta_{D}}^{(m_D)}\ldots X_{-\beta_p}^{(m_p)}\in U(\g)$, where $m_i\in \mathbb{N}$, $y\neq 1, X_{-\beta_p}$, one has $\omega(y)v_q=0$. In particular, for any $u\in U(\g)$ 
\begin{equation}\label{aq1lcl}
(v_p, yuv_{\mu_p})=0.
\end{equation}
\end{lemma}

We shall later prove Lemma \ref{Apvan} which is a quantum group analogue of this lemma. Lemma \ref{vancl} follows from Lemma \ref{Apvan} by specializing at $q^{\frac{1}{d{\bar{r}}^2}}=1$.

\begin{corollary}\label{vanclcor}
(i) For any $1\leq q<p\leq D$ and any $g\in N_{[-\beta_p,-\beta_D]}$ one has $\omega(g)v_q=v_q$.

(ii) For any $1\leq p\leq D$ and $t\in \mathbb{C}$ one has $\omega(X_{-\beta_p}(t))v_p=v_p+tv_{\mu_p}$, where $X_{-\beta_p}(t)=\exp(tX_{-\beta_p})\in N_{\{-\beta_p\}}$.

(iii) For any $1\leq p\leq D$ and $n_1\in N_{[-\beta_{p+1},-\beta_D]}$, $n_2\in N_{[-\beta_{p-1},-\beta_1]}$, $t\in \mathbb{C}$ one has $(v_p,n_1X_{-\beta_p}(t)n_2v_{\mu_p})=t$, where $X_{-\beta_p}(t)=\exp(tX_{-\beta_p})\in N_{\{-\beta_p\}}$.
\end{corollary}

\begin{proof}
Using the unique factorization $N_-=N_{\{-\beta_D\}}\ldots N_{\{-\beta_1\}}$ introduced in Lemma \ref{decNr} and the exponential map in the one--parameter subgroups $N_{\{-\beta_D\}},\ldots, N_{\{-\beta_1\}}$ one reduces the proof of parts (i) and (iii) of this corollary to the statements of parts (i) and (ii) of the previous lemma, respectively, and to the definition of the matrix elements $(v_p,\cdot~ v_{\mu_p})$.

Part (ii) follows from (\ref{x-bvm}).
\end{proof}

Let $G^0:=N_-HN_+\subset G$ \index[not]{G@$G^0$} be the big Bruhat cell. \index{Bruhat!cell!big} Recall that by Lemma \ref{vart} (iii)  multiplication gives rise to an isomorphism of verieties 
\begin{equation}\label{g0iso}
N_-HN_+\simeq N_-\times H \times N_+
\end{equation} 
as $N_-HN_+$ is the orbit of the unit element in $G$ for the following action
$$
(N_-\times B_+) \times G\to G, (n_-, b_+)\circ g=n_-gb_+^{-1},
$$
and $B_+\simeq N_+\times H$ as a variety by Lemma \ref{vart} (i) (see also Lemma 8.3.6 in \cite{Sp}).

By Lemma \ref{decNr} group multiplication yields isomorphisms of varieties $N_\pm\simeq N_{\{\pm\beta_D\}}\times\ldots \times N_{\{\pm\beta_1\}}$. Composing them with the exponential maps in the one--parameter subgroups $N_{\{\pm\beta_i\}}$, $i=1,\ldots, D$ and with isomorphism (\ref{g0iso}) we obtain the following isomorphism of varieties 
\begin{equation}\label{decbr}
\Omega: N_-HN_+ \to \mathbb{C}^{D}\times H \times \mathbb{C}^D, 
\end{equation}
$$
\Omega(n_-hn_+)=(u_1^-,\ldots, u_D^-,h,u_1^+,\ldots, u_D^+),~n_\pm=\exp(u_D^\pm X_{\pm\beta_D})\ldots \exp(u_1^\pm X_{\pm\beta_1}),~u_i^\pm\in \mathbb{C},~i=1,\ldots, D,~h\in H. \index[not]{O@$\Omega$} 
$$

The big Bruhat cell can also be defined as the compliment in $G$ of the common zero locus of the regular functions $(v_\mu, \cdot~ v_\mu)$, $\mu\in P_{+}$, $\mu\neq 0$. These functions together with the constant function equal to $1=(v_0, \cdot~ v_0)$ form a multiplicatively closed set $\mathcal{S}$, \index[not]{S@$\mathcal{S}$} and the localization of $\mathbb{C}[G]$ by this set \index{localization by a multiplicatively closed set} is isomorphic to $\mathbb{C}[N_-HN_+]$ (see e.g. \cite{Z}, \S 9 and \S 100),
\begin{equation}\label{GlocS}
\mathbb{C}[N_-HN_+]\simeq \mathbb{C}[G][\mathcal{S}^{-1}]. \index[not]{C@$\mathbb{C}[G][\mathcal{S}^{-1}]$}
\end{equation}

For any complex algebraic variety $V$ and any subset $X\subset V$ we denote by $\mathcal{V}_V(X)\subset \mathbb{C}[V]$ \index[not]{V@$\mathcal{V}_V(X)$} the vanishing ideal of $X$ in $\mathbb{C}[V]$. \index{ideal!vanishing}

The formula for the operator $\Pi$ will be given in terms of functions introduced in the following lemma.
\begin{lemma}\label{pprop}
Let $\varphi_p\in \mathbb{C}(G)$, $p=1,\ldots, D$ be the rational function on $G$ defined by 
$$
\varphi_p(g)=\frac{(v_p, g v_{\mu_p})}{(v_{\mu_p},g v_{\mu_p})}. \index[not]{f@$\varphi_p(~\cdot~)$}
$$
Then the following statements are true. 

(i) $\varphi_p\in \mathbb{C}[N_-HN_+]$. Moreover,  
\begin{equation}\label{Omp}
\varphi_p={\bf x}_p\Omega,
\end{equation}
where the function ${\bf x}_p:\mathbb{C}^{D}\times H \times \mathbb{C}^D\to \mathbb{R}$ is defined by
$$
{\bf x}_p(u_1^-,\ldots, u_D^-,h,u_1^+,\ldots, u_D^+)=u_p^-, (u_1^-,\ldots, u_D^-,h,u_1^+,\ldots, u_D^+)\in \mathbb{C}^{D}\times H \times \mathbb{C}^D,~u_i^\pm\in \mathbb{C},~i=1,\ldots, D,~h\in H. \index[not]{x@${\bf x}_p$}
$$

(ii) $\varphi_p$ and $\mathcal{V}_{G^0}(N_{[-\beta_p,-\beta_D]}HN_+)$ generate $\mathcal{V}_{G^0}(N_{[-\beta_{p+1},-\beta_D]}HN_+)$.

(iii) $\varphi_q$, $q=1,\ldots, p$ generate $\mathcal{V}_{G^0}(N_{[-\beta_{p+1},-\beta_D]}HN_+)$.
\end{lemma}

\begin{proof}

(i) $\varphi_p\in \mathbb{C}[N_-HN_+]$ as $(v_{\mu_p},\cdot~ v_{\mu_p})\in \mathcal{S}$, $\varphi_p\in\mathbb{C}[G][\mathcal{S}^{-1}]$, and the localization of $\mathbb{C}[G]$ by $\mathcal{S}$ is isomorphic to $\mathbb{C}[N_-HN_+]$ by (\ref{GlocS}).


Next, if $n_-hn_+\in N_-HN_+$, $n_-\in N_-$, $h\in H$, $n_+\in N_+$ then, since $\Delta_-=[-\beta_{p+1},-\beta_D]\cup \{-\beta_p\}\cup [-\beta_1,-\beta_{p-1}]$ (disjoint union of minimal segments), by Lemma \ref{comm} $n_-=n_1X_{-\beta_p}(u_p^-)n_2$,  $n_1\in N_{[-\beta_{p+1},-\beta_D]}$, $n_2\in N_{[-\beta_1,-\beta_{p-1}]}$,  $X_{-\beta_p}(u_p^-)\in N_{\{-\beta_p\}}$, $u_p^-\in \mathbb{C}$, and we use the symbol $u_p^-$ for the parameter in the one--parameter subgroup $N_{\{-\beta_p\}}$ introduced in formula (\ref{decbr}). Since $v_{\mu_p}$ is a highest weight vector $(v_p,n_-hn_+v_{\mu_p})=\mu_p(h)(v_p,n_-v_{\mu_p})$ and $(v_{\mu_p},n_-hn_+v_{\mu_p})=\mu_p(h)\neq 0$. By Corollary \ref{vanclcor} (iii) one has $(v_p,n_-v_{\mu_p})=(v_p,n_1X_{-\beta_p}(u_p^-)n_2v_{\mu_p})=u_p^-$. Using these identities in the definition of $\varphi_p$ we obtain
$$
\varphi_p(n_-hn_+)=\frac{(v_p,n_-hn_+ v_{\mu_p})}{(v_{\mu_p},n_-hn_+ v_{\mu_p})}=u_p^-.
$$
This proves (\ref{Omp}).

(ii) Under isomorphism (\ref{decbr}) 
$$
N_{[-\beta_p,-\beta_D]}HN_+\simeq \underbrace{0\times\ldots\times 0}_{\text{$p-1$ times}}\times\mathbb{C}^{D-p+1}\times H \times \mathbb{C}^D\subset \mathbb{C}^D\times H \times \mathbb{C}^D,
$$
and $\varphi_p$ becomes the first coordinate function on $\mathbb{C}^{D-p+1}$. We deduce that the ideal generated by $\varphi_p$ and by $\mathcal{V}_{G^0}(N_{[-\beta_p,-\beta_D]}HN_+)$ is the vanishing ideal of $\mathcal{V}_{G^0}(N_{[-\beta_{p+1},-\beta_D]}HN_+)$ as under $\Omega$
$$
N_{[-\beta_{p+1},-\beta_D]}HN_+\simeq \underbrace{0\times\ldots\times 0}_{\text{$p$ times}}\times\mathbb{C}^{D-p}\times H \times \mathbb{C}^D,
$$ 
and the first coordinate function on $\mathbb{C}^{D-p+1}$ in 
$$
\underbrace{0\times\ldots\times 0}_{\text{$p-1$ times}}\times\mathbb{C}^{D-p+1}\times H \times \mathbb{C}^D
$$ 
generates the vanishing ideal of 
$$
\underbrace{0\times\ldots\times 0}_{\text{$p$ times}}\times\mathbb{C}^{D-p}\times H \times \mathbb{C}^D\subset \underbrace{0\times\ldots\times 0}_{\text{$p-1$ times}}\times\mathbb{C}^{D-p+1}\times H \times \mathbb{C}^D.
$$
This proves part (ii).

Part (iii) follows from part (ii) by induction over $p$.

\end{proof}

Next, to define the functions of type $\varphi_p$ which will appear in the formula for the operator $\Pi$ we shall need to consider some normal orderings on the set $\Delta_+$  associated to $s$ in Definition \ref{circorddef}. For convenience in the future, we label the roots in the initial segment $\Delta_{\m_+}\subset \Delta_+$ of the system $\Delta_+$ associated to $s$ and ordered as in Definition \ref{circorddef} as follows
$$
\beta_{1 1},\ldots, \beta_{1 n_1},\beta_{2 1},\ldots, \beta_{2 n_2}, \ldots, \beta_{R-1 1}\ldots \beta_{R-1 n_{R-1}}, 
$$
where $\{\beta_{1 1},\ldots, \beta_{1 n_1}\}=\Delta^1\cap \Delta_{\m_+}$, $\{\beta_{R-1 1}\ldots \beta_{R-1 n_{R-1}}\}=\Delta^{R-1}\cap \Delta_{\m_+}$, and for $1<j<R-1$ $\{\beta_{j 1},\ldots, \beta_{j n_j}\}=\Delta^j$. \index[not]{b@$\beta_{jk}$} \index[not]{n@$n_j$}

Let
$$
s_{i_{11}}\ldots s_{i_{1n_1}}s_{i_{21}}\ldots s_{i_{2n_2}}\ldots s_{i_{R-1 1}}\ldots s_{i_{R-1 n_{R-1}}} \index[not]{s@$s_{i_{jk}}$}
$$
be the corresponding initial part of the reduced decomposition of the longest element of the Weyl group associated to the normal ordering introduced in Definition \ref{circorddef}.

Let $w_j=s_{i_{j1}}\ldots s_{i_{jn_j}}$. \index[not]{w@$w_j$} Note that since 
\begin{equation}\label{wjmin}
\Delta_{(w_1\ldots w_{j-1})^{-1}}=(\Delta^1\cap \Delta_{\m_+})\cup \Delta^2\cup \ldots \cup \Delta^{j-1}
\end{equation} 
for $j=2,\ldots, R-1$ one has
\begin{equation}\label{j+iso}
(w_1\ldots w_{j-1})^{-1}\Delta_+^j=\Delta_+
\end{equation}
by the definition of $\Delta_+^j$.

\begin{remark}\label{Djord}
Note that each root system $\Delta_+^j$ inherits a normal ordering from the circular normal ordering associated to (\ref{NO}). Therefore for each $j=2,\ldots, R-1$ the corresponding isomorphism (\ref{j+iso}) induces a normal ordering on $\Delta_+$ with respect to which the images of the roots from the segment $\beta_{j 1},\ldots, \beta_{j n_j}$ form an initial segment. We denote the roots in this segment by $\delta_{j 1},\ldots, \delta_{j n_j}$, $\delta_{j k}= (w_1\ldots w_{j-1})^{-1}\beta_{j k}$, $k=1,\ldots, n_j$, and the remaining positive roots by $\delta_{jk}$, $k=n_j+1,\ldots, D$ \index[not]{d@$\delta_{jk}$} in the increasing order according to the values of $k$. To keep the notation uniform we also write $\delta_{1 k}= \beta_{1 k}$, $k=1,\ldots, D$.
\end{remark}

For $j=1,\ldots, R-1$ define $\widehat{s}_j=(w_1\ldots w_{j-1})^{-1}s(w_1\ldots w_{j-1})$, where we assume that $w_0=1$, so $\widehat{s}_0=s$. \index[not]{s@$\widehat{s}_j$}

For each $j=1,\ldots, R-1$, $k=1,\ldots, n_j$ we define matrix elements $(v_{jk},\cdot~ v_{\mu_{jk}})$ \index[not]{v@$v_{jk}$} \index[not]{m@$\mu_{jk}$} by condition (\ref{vpkcl}), where for $j=1$ the normal ordering on $\Delta_+$ associated to $s$ introduced in Definition \ref{circorddef} is used, and for $j>1$ the normal ordering on $\Delta_+$ induced by the normal ordering on $\Delta_+^j$ with the help of isomorphism (\ref{j+iso}) is used, and (\ref{vpkcl}) is used with $\beta_p=\delta_{jk}$. Fix arbitrary representatives $w_1,\ldots, w_{R-1}, s\in N_G(H)$ of the elements $w_1,\ldots, w_{R-1}, s\in W$, and define the representative of $\widehat{s}_j\in W$ by  $\widehat{s}_j=(w_1\ldots w_{j-1})^{-1}s(w_1\ldots w_{j-1})\in N_G(H)$. Let 
\begin{equation}\label{Apc}
\varphi_{jk}(g)=\frac{(v_{jk},g \widehat{s}_j^{-1}v_{\mu_{jk}})}{(v_{\mu_{jk}},g \widehat{s}_j^{-1}v_{\mu_{jk}})}. \index[not]{f@$\varphi_{jk}$}
\end{equation}
The formula for the operator $\Pi$ will be given in terms of the functions $\varphi_{jk}$. To obtain this formula we shall have to firstly obtain some properties of these functions.

Let $G^j=G^0\widehat{s}_j=N_-HN_+\widehat{s}_j$, $j=1,\ldots, R-1$. \index[not]{G@$G^j$} Multiplication by $\widehat{s}_j^{-1}$ in $G$ from the right yields the isomorphism of $G^j$ with the open dense Bruhat cell $G^0=N_-HN_+\subset G$,
\begin{equation}\label{gopiso}
G^j\to G^0, g\mapsto g\widehat{s}_j^{-1}.
\end{equation}

Recall that by \cite{Z}, \S 9 and \S 100 $G^0$ can be defined in $G$ as the compliment of the common zero locus of the regular functions $(v_\mu, \cdot~ v_\mu)$, $\mu\in P_{+}$, $\mu\neq 0$, and by (\ref{GlocS}) $\mathbb{C}[N_-HN_+]\simeq \mathbb{C}[G][\mathcal{S}^{-1}]$, where $\mathcal{S}$ is the multiplicatively closed set formed by these functions and by the constant function equal to $1=(v_0, \cdot~ v_0)$.  

Thus if we introduce the multiplicatively closed sets $\mathcal{S}_j=\{(v_\mu, \cdot~ \widehat{s}_j^{-1}v_\mu), \mu\in P_{+}\}$, $j=1,\ldots, R-1$ \index[not]{S@$\mathcal{S}_j$} then 
\begin{equation}\label{Gpdef}
G^j=\{g\in G : (v_\mu, g \widehat{s}_j^{-1}v_\mu)\neq 0~ \forall \mu\in P_{+}\},
\end{equation} 
and 
\begin{equation}\label{locdef}
\mathbb{C}[G^j]\simeq \mathbb{C}[G][\mathcal{S}_j^{-1}].
\end{equation}
This isomorphism and  (\ref{Apc}) imply that $\varphi_{jk}$ can be regarded as an element of $\mathbb{C}[G^j]$.

By Lemma \ref{kparab} the set of roots $(w_1\ldots w_{j-1})^{-1}(-\Delta_+^j\cup \Delta_0)$ is parabolic. Let $P^j\subset G$ \index[not]{P@$P^j$} be the corresponding parabolic subgroup, $N^j_-$ \index[not]{N@$N^j_-$} its unipotent radical, and $P_s\subset G$ \index[not]{P@$P_s$} the subgroup generated by the one--parameter subgroups corresponding to the roots from the additively closed set $(w_1\ldots w_{j-1})^{-1}(-\Delta_s^s\cup \Delta_0)$. The Levi factor $L^j$ \index[not]{L@$L^j$} of $P^j$ is $(w_1\ldots w_{j-1})^{-1}L^sw_1\ldots w_{j-1}$ as $-(-\Delta_+^j\cup \Delta_0)\cap (-\Delta_+^j\cup \Delta_0)=\Delta_0$.

The main properties of the functions $\varphi_{jk}$, which will be used to obtain the formula for the operator $\Pi$, are summarized in the following lemma.
\begin{lemma}\label{jklem}
(i) For any $\mu\in P_+$, $j=1,\ldots, R-1$, $k=1,\ldots, n_j$, and $g=n_-z\widehat{s}_jn\in N_{[-\delta_{jk}, -\delta_{jD}]}\widehat{s}_j G$, where $n_-\in N_{[-\delta_{jk}, -\delta_{jD}]}$, $z\in L^j$, and $n\in N^j_-$ or $n\in N_{(w_1\ldots w_{j-1})^{-1}(-\Delta_s^s)}$ one has 
\begin{equation}\label{uform}
(v_{jk},g \widehat{s}_j^{-1}v_{\mu_{jk}})=u(v_{\mu_{jk}}, z v_{\mu_{jk}}),
\end{equation}
and 
\begin{equation}\label{zeroform}
(v_\mu,g \widehat{s}_j^{-1}v_\mu)=(v_\mu, z v_\mu),
\end{equation}
where $u\in \mathbb{C}$ is defined from the unique factorizations $n_-=X_{-\delta_{jk}}(u)n'=n''X_{-\delta_{jk}}(u)$, $n',n''\in N_{[-\delta_{jk+1}, -\delta_{jD}]}$ with the help of Lemma \ref{comm}.
Thus 
\begin{equation}\label{phijku}
\varphi_{jk}(g)=u.
\end{equation}

(ii) Let $Z^j_\pm=(w_1\ldots w_{j-1})^{-1}Z_\pm w_1\ldots w_{j-1}$, \index[not]{Z@$Z^j_\pm$} $M_-^j=N_-^j\cap (w_1\ldots w_{j-1})^{-1}M_- w_1\ldots w_{j-1}$, $j=1,\ldots, R-1$. \index[not]{M@$M_-^j$} Then for any $j=1,\ldots, R-1$, $k=1,\ldots, n_j$ one has $N_{[-\delta_{jk}, -\delta_{jD}]} Z^j_-HZ^j_+\widehat{s}_j N^j_-\subset G^j$. Moreover, for $k=1,\ldots, n_j$ the function $\varphi_{jk}\in \mathbb{C}[G^j]$ and the ideal ${J^{jk}}^{loc}:=\mathcal{V}_{G^j}(N_{[-\delta_{jk}, -\delta_{jD}]} Z^j_-HZ^j_+\widehat{s}_j M^j_-)$ \index[not]{J@${J^{jk}}^{loc}$} generate ${J^{jk+1}}^{loc}:=\mathcal{V}_{G^j}(N_{[-\delta_{j k+1}, -\delta_{jD}]} Z^j_-HZ^j_+\widehat{s}_j M^j_-)$.

(iii) For any $j=1,\ldots, R-1$, $k=1,\ldots, n_j$ the functions $\varphi_{jm}$, $m=1,\ldots, k$ and the ideal ${J^{j1}}^{loc}$ generate ${J^{jk+1}}^{loc}$.
\end{lemma} 

\begin{proof}

(i) By Corollary \ref{vanclcor} (i) and (ii) we have
\begin{equation}\label{vjk}
(v_{jk},g \widehat{s}_j^{-1}v_{\mu_{jk}})=(v_{jk},n''X_{-\delta_{jk}}(u)z\widehat{s}_jn \widehat{s}_j^{-1}v_{\mu_{jk}})=
(v_{jk},  z\widehat{s}_j n \widehat{s}_j^{-1}v_{\mu_{jk}})+u(v_{\mu_{jk}},  z\widehat{s}_jn \widehat{s}_j^{-1}v_{\mu_{jk}}).
\end{equation}

If $n\in N^j_-$ the weight components of the vector $z\widehat{s}_jn \widehat{s}_j^{-1}v_{\mu_{jk}}$ belong to $-(w_1\ldots w_{j-1})^{-1}s(\Delta_0\cup\Delta_+^j)+\mu_{jk}$ by the definition of $P^j$. By Lemma (\ref{jsegm}) (i) this set has empty intersection with $-(w_1\ldots w_{j-1})^{-1}\Delta^j+\mu_{jk}$. On the other hand  the weight of $v_{jk}$ is equal to $\mu_{jk}-\delta_{jk}\in -(w_1\ldots w_{j-1})^{-1}\Delta^j+\mu_{jk}$. Since different weight subspaces are orthogonal with respect to the contravariant form, the first term in (\ref{vjk}) vanishes in this case.

The same conclusion can be obtained if $n\in N_{(w_1\ldots w_{j-1})^{-1}(-\Delta_s^s)}$ with the help of part (ii) of Lemma \ref{jsegm}. 

In both cases the $\mu_{jk}$--weight component of $z\widehat{s}_jn \widehat{s}_j^{-1}v_{\mu_{jk}}$ can give a nontrivial contribution to the right hand side of (\ref{vjk}). If $n\in N^j_-$ we note that the set of roots $-(w_1\ldots w_{j-1})^{-1}s(\Delta_0\cup\Delta_+^j)$ is parabolic and $(-(w_1\ldots w_{j-1})^{-1}s(\Delta_0\cup\Delta_+^j))\cap(w_1\ldots w_{j-1})^{-1}s(\Delta_0\cup\Delta_+^j)=(w_1\ldots w_{j-1})^{-1}\Delta_0$. Therefore the $\mu_{jk}$--weight component of $z\widehat{s}_jn \widehat{s}_j^{-1}v_{\mu_{jk}}$ is equal to the $\mu_{jk}$--weight component of $zv_{\mu_{jk}}$, and we obtain (\ref{uform}). 

The case when $n\in N_{(w_1\ldots w_{j-1})^{-1}(-\Delta_s^s)}$ can be treated in a similar way noting that $L^jN_{(w_1\ldots w_{j-1})^{-1}(-\Delta_s^s)}\subset (w_1\ldots w_{j-1})^{-1}P^s(w_1\ldots w_{j-1})$, and that the subgroup $(w_1\ldots w_{j-1})^{-1}P^s(w_1\ldots w_{j-1})\subset G$ is parabolic with the Levi factor $L^j$, as the additively closed subset of roots $(w_1\ldots w_{j-1})^{-1}(\Delta_-^s\cup \Delta_0)$ is parabolic by the definition with $(w_1\ldots w_{j-1})^{-1}(\Delta_-^s\cup \Delta_0))\cap (-(w_1\ldots w_{j-1})^{-1}(\Delta_-^s\cup \Delta_0)=(w_1\ldots w_{j-1})^{-1}(\Delta_0)$, so the Levi factor of $(w_1\ldots w_{j-1})^{-1}P^s(w_1\ldots w_{j-1})$ is $L^j$ which has root system $(w_1\ldots w_{j-1})^{-1}(\Delta_0)$.

Formula (\ref{zeroform}) is established using the same arguments.

Formula (\ref{phijku}) follows from the definition of $\varphi_{jk}$, formula (\ref{uform}) and (\ref{zeroform}) with $\mu=\mu_{jk}$.

(ii) Note that $Z^j_\pm=(w_1\ldots w_{j-1})^{-1}Z_\pm w_1\ldots w_{j-1}\subset N_\pm$ as $\Delta_{(w_1\ldots w_{j-1})^{-1}}\cap \Delta_0=\emptyset$ by (\ref{wjmin}). Thus for $g=n_-z_-hz_+\widehat{s}_jn\in N_{[-\delta_{jk}, -\delta_{jD}]} Z^j_-HZ^j_+\widehat{s}_j N^j_-$, where $n_-\in N_{[-\delta_{jk}, -\delta_{jD}]}$, $z_\pm\in Z^j_\pm$, $h\in H$ and $n\in N^j_-$ one has from (\ref{zeroform})
$$
(v_\mu,g \widehat{s}_j^{-1}v_\mu)=(v_\mu, z_-hz_+ v_\mu)=(v_\mu, h v_\mu)\neq 0.
$$
From (\ref{Gpdef}) it follows that $N_{[-\delta_{jk}, -\delta_{jD}]} Z^j_-HZ^j_+\widehat{s}_j N^j_-\subset G^j$.

By the definition $\varphi_{jk}\in \mathbb{C}[G^j]\simeq \mathbb{C}[G][\mathcal{S}_j^{-1}]$.

Now let $I'$ be the ideal in $\mathbb{C}[G^j]$ generated by the function $\varphi_{jk}$ and the ideal ${J^{jk}}^{loc}$ 

From formula (\ref{phijku}) it follows that an element $g\in N_{[-\delta_{jk}, -\delta_{jD}]} Z^j_-HZ^j_+\widehat{s}_j M^j_-$ belongs to $N_{[-\delta_{jk+1}, -\delta_{jD}]} Z^j_-HZ^j_+\widehat{s}_j M^j_-$ if and only if $\varphi_{jk}(g)=0$. Therefore the zero locus of $I'$ coincides with the closure of $N_{[-\delta_{jk+1}, -\delta_{jD}]} Z^j_-HZ^j_+\widehat{s}_j M^j_-$. It remains to show that $I'$ is radical.

Indeed, let $a\in \mathbb{C}[G^p]$ be such that $a^n\in I'$ for some $n\in \mathbb{N}$, $n>0$, i.e. $a^n=f+f'$, $f=\varphi_{jk}^m\rho$ for some $m\in \mathbb{N}$, $m>0$ and $\rho\in \mathbb{C}[G^p]$, $f'\in {J^{jk}}^{loc}$. Then the restriction of $a^n$ to $N_{[-\delta_{jk}, -\delta_{jD}]} Z^j_-HZ^j_+\widehat{s}_j M^j_-$ coincides with the restriction of $f$ to $N_{[-\delta_{jk}, -\delta_{jD}]} Z^j_-HZ^j_+\widehat{s}_j M^j_-$. In particular, the restriction of $a^n$ to $N_{[-\delta_{jk}, -\delta_{jD}]} Z^j_-HZ^j_+\widehat{s}_j M^j_-$ vanishes on $N_{[-\delta_{jk+1}, -\delta_{jD}]} Z^j_-HZ^j_+\widehat{s}_j M^j_-$. This implies that the restriction of $a$ to $N_{[-\delta_{jk}, -\delta_{jD}]} Z^j_-HZ^j_+\widehat{s}_j M^j_-$ must also vanish on $N_{[-\delta_{jk+1}, -\delta_{jD}]} Z^j_-HZ^j_+\widehat{s}_j M^j_-$.

By Lemma \ref{pprop} (i) the composition of the isomorphisms (\ref{decbr}) and (\ref{gopiso}) yields an isomorphism of varieties\begin{equation}\label{decbrj}
G^j\simeq \mathbb{C}^{D}\times H \times \mathbb{C}^D,
\end{equation}
and under this isomorphism the function $\varphi_{jk}$ becomes the $k$-th coordinate function on $\mathbb{C}^D$ on the first factor in (\ref{decbrj}).

Now formula (\ref{phijku}) and simple induction over $k$ imply that under isomorphism (\ref{decbrj}) $N_{[-\delta_{jk}, -\delta_{jD}]} Z^j_-HZ^j_+\widehat{s}_j M^j_-$ becomes a subset of $\mathbb{C}^{D}\times H \times \mathbb{C}^D$ of the form 
$$
\underbrace{0\times\ldots \times 0}_{\text{$k-1$ factors}}\times \mathbb{C}\times X_{j_k},
$$
where the factor $\mathbb{C}$ in the last product corresponds to the $k$-th component in $\mathbb{C}^D$ in the first factor in the right hand side of (\ref{decbrj}), $X_{j_k}$ is a subset of all the remaining factors in the right hand side of (\ref{decbrj}), and the function $\varphi_{jk}$ becomes the coordinate function on the factor $\mathbb{C}$.

From formula (\ref{phijku}) it follows that an element $g\in N_{[-\delta_{jk}, -\delta_{jD}]} Z^j_-HZ^j_+\widehat{s}_j M^j_-$ belongs to $N_{[-\delta_{jk+1}, -\delta_{jD}]} Z^j_-HZ^j_+\widehat{s}_j M^j_-$ if and only if $\varphi_{jk}(g)=0$. Hence $X_{j_k}\simeq N_{[-\delta_{jk+1}, -\delta_{jD}]} Z^j_-HZ^j_+\widehat{s}_j M^j_-$, and under isomorphism (\ref{decbrj}) we have
$N_{[-\delta_{jk+1}, -\delta_{jD}]} Z^j_-HZ^j_+\widehat{s}_j M^j_-\simeq 0\times\ldots \times 0\times 0 \times X_{j_k}$.

We conclude that if $a$ vanishes on $N_{[-\delta_{jk+1}, -\delta_{jD}]} Z^j_-HZ^j_+\widehat{s}_j M^j_-$, its restriction to $N_{[-\delta_{jk}, -\delta_{jD}]} Z^j_-HZ^j_+\widehat{s}_j M^j_-$ coincides with the restriction to this subset of the function $\varphi_{jk}^p\psi$ for some $p\in \mathbb{N}, p>0$ and $\psi\in \mathbb{C}[G^p]$. Thus $a=\varphi_{jk}^p(\psi+\psi')+\psi''$, where $\psi', \psi''\in {J^{jk}}^{loc}$, and hence $a\in I'$ by the definition of $I'$. 

We deduce that $I'$ is radical. This implies $I'={{J}^{jk+1}}^{loc}$ and completes the proof of part (ii).

Part (iii) follows from (ii) by induction over $k$. 

\end{proof}

Now we come back to the description of the projection operator $\Pi$. 

\begin{proposition}\label{tp}
Let $g=nn_szsm_sn^{-1}$ be the unique presentation (\ref{gslice}) for an element $g\in N_-Z^ssM_-$. Then $n$ can be uniquely factorized as $n= X_{-\beta_{11}}(u_{11})\ldots X_{-\beta_{R-1 n_{R-1}}}(u_{R-1 n_{R-1}})$ where we assume that the root vectors $X_{-\beta_{jk}}$ used in the definition of the one--parameter subgroups are related to the root vectors $X_{-\delta_{jk}}$ used in the definition of the functions $\varphi_{jk}$ as follows $X_{-\beta_{jk}}={\rm Ad}(w_1\ldots w_{j-1})X_{-\delta_{jk}}$, and the numbers $u_{jk}$ can be found inductively by the following formula
\begin{equation}\label{tind}
u_{jk}=\varphi_{jk}((w_1\ldots w_{j-1})^{-1}g_{jk}(w_1\ldots w_{j-1})),
\end{equation}
where $g_{jk}=n_{jk}^{-1}gn_{jk}$, $n_{jk}=X_{-\beta_{11}}(u_{11})\ldots X_{-\beta_{jk-1}}(u_{jk-1})$, $j=1,\ldots, R-1$, $k=1, \ldots, n_j$ and it is assumed that $n_{10}=1$ and $X_{-\beta_{j0}}(u_{j0})=X_{-\beta_{j-1n_{j-1}}}(u_{jn_{j-1}})$.
\end{proposition}

\begin{proof}
The numbers $u_{jk}$ can be found by induction starting with $u_{11}$. We shall establish the induction step. The case of $u_{11}$ corresponding to the base of the induction can be considered in a similar way.

Assume that $u_{11}, \ldots , u_{jk-1}$ have already been found. Then
\begin{eqnarray*}
g_{jk}=n_{jk}^{-1}gn_{jk}=X_{-\beta_{jk}}(u_{jk})\ldots X_{-\beta_{R-1n_{R-1}}}(u_{R-1n_{R-1}})n_szsm_sX_{-\beta_{R-1n_{R-1}}}(-u_{R-1n_{R-1}})\ldots X_{-\beta_{jk}}(-u_{jk}).
\end{eqnarray*}

Now $(w_1\ldots w_{j-1})^{-1}g_{jk}(w_1\ldots w_{j-1})$ has the form of $g$ from Lemma \ref{jklem} (i) with $u=u_{jk}$. Therefore we obtain (\ref{tind}) by formula (\ref{phijku}).

\end{proof}
  
For a representative $w\in N_G(H)$ of a Weyl group element $w\in W$ we denote the operator on $\mathbb{C}(G)$ \index[not]{C@$\mathbb{C}(G)$} induced by the conjugation by $w$ on $G$ by the same letter,
$$
(wf)(g)=f(wgw^{-1}), f\in \mathbb{C}(G). 
$$

Observing that in the notation of Proposition \ref{tp} for $g=nn_szsm_sn^{-1}\in N_-Z^ssM_-$ we have $n^{-1}gn=n_szsm_s$ and recalling the definition of the operator $\Pi$ in (\ref{Pidef}) we infer the following theorem from Proposition \ref{tp}. 
\begin{theorem}\label{proj}
Let $\Pi_{jk}$ \index[not]{P@$\Pi_{jk}$} be the operator on the space of rational functions $\mathbb{C}(G)$ \index{function!rational} on $G$ induced by conjugation by the element $\exp(-\varphi_{jk}X_{-\delta_{jk}})$,
\begin{equation}\label{Pip}
(\Pi_{jk} f)(g)=f(\exp(-\varphi_{jk}(g)X_{-\delta_{jk}})g\exp(\varphi_{jk}(g)X_{-\delta_{jk}})).
\end{equation}

Then the restriction of the composition 
$$
{\Pi}_{11}\ldots {\Pi}_{1n_1}\circ w_1^{-1}\circ {\Pi}_{21}\ldots {\Pi}_{2n_2}\circ w_2^{-1} \ldots \circ w_{R-2}^{-1} \circ {\Pi}_{R-11}\ldots {\Pi}_{R-1 n_{R-1}}\circ w_1 \ldots w_{R-2}
$$ 
to $\mathbb{C}[N_-Z^ssM_-]$ is equal to the projection operator $\Pi$ onto the subspace $\mathbb{C}[N_-Z^ssM_-]^{M_-}$ of $M_-$--invariant regular functions on $N_-Z^ssM_-$, $\Pi: \mathbb{C}[N_-Z^ssM_-]\rightarrow \mathbb{C}[N_-Z^ssM_-]^{M_-}$,
\begin{equation}\label{Pi1}
\Pi={\Pi}_{11}\ldots {\Pi}_{1n_1}\circ w_1^{-1}\circ {\Pi}_{21}\ldots {\Pi}_{2n_2}\circ w_2^{-1} \ldots \circ w_{R-2}^{-1} \circ {\Pi}_{R-11}\ldots {\Pi}_{R-1 n_{R-1}}\circ w_1 \ldots w_{R-2}.
\end{equation}
\end{theorem}

\begin{corollary}\label{proj1}
The operator $\Pi_c: \mathbb{C}[(w_1 \ldots w_{R-2})^{-1}N_-Z^ssM_-(w_1 \ldots w_{R-2})]\rightarrow \mathbb{C}[N_-Z^ssM_-]^{M_-}$  \index[not]{P@$\Pi_c$}
defined by
\begin{equation}\label{pic}
\Pi_c={\Pi}_{11}\ldots {\Pi}_{1n_1}\circ w_1^{-1}\circ {\Pi}_{21}\ldots {\Pi}_{2n_2}\circ w_2^{-1} \ldots \circ w_{R-2}^{-1} \circ {\Pi}_{R-11}\ldots {\Pi}_{R-1 n_{R-1}}=\Pi \circ (w_1 \ldots w_{R-2})^{-1}.
\end{equation}
is surjective.
\end{corollary}

The operator $\Pi_c$ and this corollary have direct quantum analogues.


\section{Vanishing ideals}\label{vanI}

\setcounter{equation}{0}
\setcounter{theorem}{0}

In this section we obtain some properties of the ideals ${J^{jk}}^{loc}$ which will be used in the next chapter to introduce and study quantum analogues of the operators $\Pi_{jk}$.
Let $\iota :\mathbb{C}[G]\to \mathbb{C}[G^j]\simeq \mathbb{C}[G][\mathcal{S}_j^{-1}]$ be the canonical ring homomorphism, $J^{jk}:=\iota^{-1}{J^{jk}}^{loc}$. \index[not]{J@$J^{jk}$}

By this definition  
\begin{equation}
J^{jk}=\mathcal{V}_{G}(N_{[-\delta_{jk}, -\delta_{jD}]} Z^j_-HZ^j_+\widehat{s}_j M^j_-).
\end{equation}

\begin{lemma}\label{jloc}
(i) ${J^{jk}}^{loc}=J^{jk}[\mathcal{S}_j^{-1}]$.

(ii) ${{J}^{jk}}^{loc}={{J}^{jk}}'[\mathcal{S}_j^{-1}]$, where ${{J}^{jk}}'\subset \mathbb{C}[G]$ \index[not]{J@${J^{jk}}'$} is the ideal generated by $J^{j1}$ and by $(v_{jm},\cdot~ \widehat{s}_j^{-1} v_{\mu_{jm}})$, $m=1,\ldots, k-1$, and for convenience we put ${{J}^{j1}}'=J^{j1}$.
\end{lemma}

\begin{proof}
(i) 
If $f\in {J^{jk}}^{loc}$ then $f=g/g'$, $g\in \mathbb{C}[G]$, $g'\in \mathcal{S}_j$ and $f$ vanishes on $N_{[-\delta_{jk}, -\delta_{jD}]} Z^j_-HZ^j_+\widehat{s}_j M^j_-$. Since by Lemma \ref{jklem} (ii)
$$
N_{[-\delta_{jk}, -\delta_{jD}]} Z^j_-HZ^j_+\widehat{s}_j M^j_-\subset G^j,
$$
and by (\ref{Gpdef}) $g'$ does not vanish on $G^j$ we obtain that $g$ vanishes on $N_{[-\delta_{jk}, -\delta_{jD}]} Z^j_-HZ^j_+\widehat{s}_j M^j_-$, i.e. $g\in J^{jk}$. We deduce that ${J^{jk}}^{loc}\subset J^{jk}[\mathcal{S}_j^{-1}]$. The opposite inclusion is obvious.

(ii) Note that by the definition of $\varphi_{jm}$ one has 
$$
(v_{jm},g \widehat{s}_j^{-1}v_{\mu_{jk}})=\varphi_{jm}(g)(v_{\mu_{jm}}, g \widehat{s}_j^{-1}v_{\mu_{jm}}),
$$
and $(v_{\mu_{jm}}, \cdot~ \widehat{s}_j^{-1}v_{\mu_{jm}})\in \mathcal{S}_j$. Therefore the statement in part (ii) follows from Lemma \ref{jklem} (iii).

\end{proof}

Next we give several descriptions of the ideals  $J^{jk}$. Lemma \ref{jklem} (ii) and (iii), Corollary \ref{proj1} and the statements of this section are the only results of this chapter which will be used in Chapter \ref{part4} for the purposes of quantization.

\begin{proposition}\label{Jpdescr}
(i) For $j=1,\ldots, R-1$, $k=1,\ldots, n_j$ let $\n_{jk}, \z^j_\pm, \m_-^j\subset \g$\index[not]{n@$\n_{jk}$} \index[not]{z@$\z^j_\pm$} \index[not]{m@$\m_-^j$} be the Lie subalgebras of the subgroups $N_{[-\delta_{jk}, -\delta_{jD}]}$, $Z^j_\pm$, and $M^j_-$ of $G$, respectively. Then the ideal $J^{jk}$ consists of the matrix elements of the form $(w, \cdot~ v)\in \mathbb{C}[G]$, where $w,v\in V$, $V$ is a finite-dimensional representation of $\g$, and $(w,y\widehat{s}_jhz_+xv)=0$ for any $y\in U(\n_{jk})$, $h\in U(\h)$, $z_+\in U(\z_+^j)$, $x\in U(\m_-^j)$.

(ii) The ideal $J^{j1}$ is generated by the matrix elements of the form $(u, \cdot~ v)\in \mathbb{C}[G]$, where $u$ is a highest weight vector in a finite-dimensional representation $V$ of $\g$, and $v\in V$ satisfies, and $(u,\widehat{s}_jz_+xv)=0$ for any $z_+\in U(\z_+^j)$, $x\in U(\m_-^j)$.
\end{proposition}

The proof of this proposition follows from the following lemma.

\begin{lemma}\label{Jl1}
Let $G_1,\ldots, G_k\subset G$ be the Lie subgroups corresponding to Lie subalgebras $\g_1,\ldots, \g_k\subset \g$, respectively, and $g\in G$. The following statements are true.

(i) The ideal $\mathcal{V}_G(G_1 g G_2\ldots G_k)\subset \mathbb{C}[G]$ consists of the matrix elements of the form $(w, \cdot~ v)\in \mathbb{C}[G]$, where $w,v\in V$, $V$ is a finite-dimensional representation $V$ of $\g$, and $(w,x_1gx_2\ldots x_kv)=0$ for any $x_i\in U(\g_i)$, $i=1,\ldots,k$.

(ii) If $G_1=B_-$ the ideal $\mathcal{V}_G(G_1 g G_2\ldots G_k)\subset \mathbb{C}[G]$ is generated by the matrix elements of the form $(u, \cdot~ v)\in \mathbb{C}[G]$, where $u$ is a highest weight vector in a finite-dimensional representation $V$ of $\g$, and $v\in V$ is such that $(u, gx_2\ldots x_kv)=0$ for any $x_i\in U(\g_i)$, $i=2,\ldots,k$.
\end{lemma}

\noindent
{\em Proof of Proposition \ref{Jpdescr}.}
Since $s$ fixes all roots from $\Delta_0$, $Z_+$ is generated by one-parameter subgroups corresponding to roots from $(\Delta_0)_+$, and any representative of any Weyl group element normalizes $H$, we have 
$$
HZ_+^j\widehat{s}_j=H(w_1\ldots w_{j-1})^{-1}Z_+sw_1\ldots w_{j-1}=H(w_1\ldots w_{j-1})^{-1}sZ_+w_1\ldots w_{j-1})=
$$
$$
=(w_1\ldots w_{j-1})^{-1}sHZ_+w_1\ldots w_{j-1}=H\widehat{s}_jZ_+^j=\widehat{s}_jHZ_+^j.
$$

Therefore we can write
\begin{equation}\label{Jp1}
J^{jk}=\mathcal{V}_G(N_{[-\delta_{jk}, -\delta_{jD}]} Z^j_-HZ^j_+\widehat{s}_j M^j_-)=\mathcal{V}_G(N_{[-\delta_{jk}, -\delta_{jD}]} Z^j_-\widehat{s}_j HZ^j_+M^j_-).
\end{equation}

Since for $j=1,\ldots, R-1$, $k=1,\ldots n_j$ $(\Delta_0)_-\subset [-\beta_{jk}, -\beta_D]$, where $\beta_D$ is the last root in the normal ordering of $\Delta_+$ associated to $s$ in Definition \ref{circorddef}, we deduce that $(w_1\ldots w_{j-1})^{-1}((\Delta_0)_-)\subset (w_1\ldots w_{j-1})^{-1}([-\beta_{jk},-\beta_D])\subset [-\delta_{jk}, -\delta_{jD}]$, and hence $Z_-^j\subset N_{[-\delta_{jk}, -\delta_{jD}]} \subset N_-$. Thus $N_{[-\delta_{jk}, -\delta_{jD}]}Z_-^j=N_{[-\delta_{jk}, -\delta_{jD}]}$, and (\ref{Jp1}) takes the form
\begin{equation}\label{Jp2}
J^{jk}=\mathcal{V}_G(N_{[-\delta_{jk}, -\delta_{jD}]}\widehat{s}_jH Z^j_+M^j_-)=\mathcal{V}_G(N_{[-\delta_{jk}, -\delta_{jD}]}H\widehat{s}_j Z^j_+M^j_-).
\end{equation}

Now part (i) of Proposition \ref{Jpdescr} follows from part (i) of Lemma \ref{Jl1}, with $k=4$, $G_1=N_{[-\delta_{jk}, -\delta_{jD}]}$, $G_2=H$, $G_3=Z_+^j$, $G_4=M^j_-$, and $g=\widehat{s}_j$.

Part (ii) of Proposition \ref{Jpdescr} follows from parts (ii) of Lemma \ref{Jl1} with $k=3$, $G_1=N_{[-\delta_{jk}, -\delta_{jD}]}H=B_-$, $G_2=Z_+^j$, $G_3=M^j_-$, and $g=\widehat{s}_j$.

\qed

\noindent
{\em Proof of part (i) of Lemma \ref{Jl1}.}   
 
Since for any $j=1,\ldots,k$, $i,i_j\in \{1,2,\ldots\}$, $r_j^{i}\in \mathbb{R}$, $x_j^i\in \g_j$ one has 
$$
e^{r_1^1x_1^1}\ldots e^{r_1^{i_1}x_1^{i_1}}ge^{r_2^1x_2^1}\ldots e^{r_2^{i_2}x_2^{i_2}}\ldots e^{r_k^1x_k^1}\ldots e^{r_k^{i_k}x_k^{i_k}}\in G_1 g G_2\ldots G_k,
$$
we obtain that for every element $(w,\cdot~ v)\in \mathcal{V}_G(G_1 g G_2\ldots G_k)$
$$
0=\frac{d^{n}}{\prod_{j=1}^k\prod_{i=1}^{i_j}dr_j^i}\left|_{r_i^j=0}(w,e^{r_1^1x_1^1}\ldots e^{r_1^{i_1}x_1^{i_1}}ge^{r_2^1x_2^1}\ldots e^{r_2^{i_2}x_2^{i_2}}\ldots e^{r_k^1x_k^1}\ldots e^{r_k^{i_k}x_k^{i_k}} v)= \right.
$$
$$
=(w,x_1^1\ldots x_1^{i_1}gx_2^1\ldots x_2^{i_2}\ldots x_k^1\ldots x_k^{i_k} v),
$$
where $n=\sum_{j=1}^k\sum_{i=1}^{i_j}1=\sum_{j=1}^k\frac{(i_j+1)i_j}{2}$.
Thus $(w,x_1gx_2\ldots x_kv)=0$ for any $x_j\in U(\g_j)$, $j=1,\ldots,k$ of the form $x_j=x_j^1\ldots x_j^{i_j}$ as in the previous formula. By linearity $(w,x_1gx_2\ldots x_kv)=0$ for any $x_j\in U(\g_j)$, $j=1,\ldots,k$. 

Conversely, let $(w,\cdot~ v)$ be any element of $\mathbb{C}[G]$ such that $(w,x_1gx_2\ldots x_kv)=0$ for any $x_j\in U(\g_j)$, $j=1,\ldots,k$. Observe that the Lie groups $G_j$, $j=1,\ldots k$ are connected, so that  each $G_j$ is generated by $\exp(\g_j)$. Thus any element of $G_1 g G_2\ldots G_k$ can be written in the form $g_1gg_2\ldots g_k$, where $g_j=e^{x_j^1}\ldots e^{x_j^{i_j}}$, $j=1,\ldots,k$, $i_j\in \{1,2,\ldots\}$, $x_j^i\in \g_j$, $i=1,\ldots, i_j$. Now we have
$$
(w,g_1gg_2\ldots g_k v)=(w,e^{x_1^1}\ldots e^{x_1^{i_1}}ge^{x_2^1}\ldots e^{x_2^{i_2}}\ldots e^{x_k^1}\ldots e^{x_k^{i_k}}v)=
$$
$$
=\sum_{n_j^i\in \mathbb{N}}\frac{1}{\prod_{j=1}^k\prod_{i=1}^{i_j} n_j^i!}(w,(x_1^1)^{n_1^1}\ldots (x_1^{i_1})^{n_1^{i_1}}g(x_2^1)^{n_2^1}\ldots (x_2^{i_2})^{n_2^{i_2}}\ldots (x_k^1)^{n_k^1}\ldots (x_k^{i_k})^{n_k^{i_k}}v)=0.
$$

Thus $(w,\cdot~ v)\in \mathcal{V}_G(G_1 g G_2\ldots G_k)$. This completes the proof.

\qed

The proof of part (ii) of Lemma \ref{Jl1} is based on the description of closed subvarieties in the flag variety \index{variety!flag} $B_-\setminus G$ \index[not]{B@$B_-\setminus G$} in terms of the so--called generalized Pl\"{u}cker coordinates.

Recall that matrix elements of the form $(u, \cdot~ v)$, where $u$ is a highest weight vector in a finite-dimensional representation $V$ of $\g$, and $v\in V$, can be viewed as sections of line bundles on $B_-\setminus G$ (see \cite{FZ1}, Section 3.1). They are also called {\it generalized Pl\"{u}cker coordinates} on $B_-\setminus G$. \index{Pl\"{u}cker coordinates, generalized} More precisely, if $V$ has highest weight $\lambda$ then $(u, \cdot~ v)$ is a section of the line bundle on $B_-\setminus G$ associated to the one--dimensional representation of $B_+=\omega(B_-)$ corresponding to $\lambda$.  

\begin{lemma}\label{zerolp}
Any closed subvariety in $B_-\setminus G$ is the zero locus of a finite set of generalized Pl\"{u}cker coordinates. 
\end{lemma}

\begin{proof}
In this proof, for any complex vector space $V$ we denote by $P(V)$ \index[not]{P@$P(V)$} its projectivisation, \index{projectivisation of a vector space} and by $[v]\in P(V)$ the line defined by a non--zero element $v\in V$.

The flag variety $B_-\setminus G$ can be realized as the $G$--orbit $\mathcal{O}$ of the line $[u^*]$ defined by a non--zero lowest weight vector $u^*$ in the projectivisation $P(V_\mu^*)$ of the finite-dimensional irreducible representation $V_\mu^*$ \index[not]{V@$V_\mu^*$} of $G$ dual \index{representation!of a semisimple algebraic group!dual} to a highest weight irreducible representation $V_\mu$ with a regular dominant highest weight \index{weight!regular dominant} $\mu$ (see e.g. \cite{Gel}, Section 2 or \cite{GS}, \S 4). Note that $B_-\setminus G$ is a projective complete variety \index{variety!projective}\index{variety!complete} (see \cite{Sp}, Section 6.2, in particular, Lemma 6.2.2), and hence by Proposition 6.1.2 (iv) in \cite{Sp} it is closed in $P(V_\mu^*)$.

Let $i:V_\mu^*\to V_\mu$ be the linear isomorphism defined with the help of the contravariant form on $V_\mu$, so that if $y^*\in V_\mu^*$ then $y^*(v)=(i(y^*),v)$ for all $v\in V_\mu$. It induces an isomorphism $P(V_\mu^*)\simeq P(V_\mu)$ which we denote by the same letter. Then $u:=i(u^*)$ is a non--zero highest weight vector in $V_\mu$, and the class $[g]\in B_-\setminus G$ of an element $g\in G$ corresponds to $\omega(g)[u]=[\omega(g)u]\in i(\mathcal{O})\simeq \mathcal{O}\subset P(V_\mu^*)$. 

Now let $v\in V_\mu$. Then 
\begin{equation}\label{ugv}
(u,gv)=(\omega(g)u,v)=(x,v), x=\omega(g)u\in V_\mu\simeq V_\mu^*,
\end{equation}
so that $[x]=[\omega(g)u]\in i(\mathcal{O})\simeq \mathcal{O}\subset P(V_\mu^*)$.

Any $y\in V_\mu$ can be uniquely written in the form $y=\sum_{n=1}^{d(V_\mu)}c_ne_n\in V_\mu$, where $e_n$, $n=1,\ldots , d(V_\mu)$, $d(V_\mu):={\rm dim}V_\mu$, is a linear  basis of $V_\mu$ consisting of weight vectors. The functions $\chi_n(y):=c_n$ are linear coordinates on $V_\mu$. Moreover, since the contravariant form gives rise to a vector space isomorphism $V_\mu\simeq V_\mu^*$, one can find elements $v_n\in V_\mu$ such that $\chi_n(y)=(y,v_n)$. In particular, $\chi_n$ generate the algebra of polynomial functions on $V_\mu$, and any closed subvariety in $P(V_\mu)$, and hence in $i(\mathcal{O})\simeq \mathcal{O}$, which is closed in $P(V_\mu)$, is the zero locus of a finite collection of some homogeneous polynomials in $\chi_n$ (see \cite{Har}, Ch. 1, \S2). 

If $f(\chi_1,\ldots \chi_{d(V_\mu)})$ is such a polynomial of degree $d$ then the equation
$f(\chi_1(y),\ldots \chi_{d(V_\mu)}(y))=0$ is well--defined in $P(V_\mu)$, i.e. if $y\in V_\mu$ is its solution then any element of the line $[y]\in P(V_\mu)$ is also its solution. Now using (\ref{ugv}), the definition of $\chi_n$ and the homogeneity of $f$ we deduce that for $[x]\in i(\mathcal{O})$ $f(\chi_1(x),\ldots \chi_{d(V_\mu)}(x))=0$ if and only if $f((u,gv_1),\ldots ,(u,gv_{d(V_\mu)}))=0$ for any $g\in G$, where $[g]\in B_-\setminus G$ corresponds to $[x]$ under the isomorphism $i(\mathcal{O})\simeq B_-\setminus G$. Note also that using algebraic rules for matrix elements we immediately obtain that $f((u,\cdot~ v_1),\ldots ,(u,\cdot~ v_{d(V_\mu)}))$ is a generalized Pl\"{u}cker coordinate defined using the representation $V_\mu^{\otimes d}$ and the highest weight vector 
$$
u^{\otimes d}=\underbrace{u\otimes \ldots \otimes u}_{d~{\rm times}}\in V_\mu^{\otimes d},
$$
i.e. $f((u,\cdot~ v_1),\ldots ,(u,\cdot~ v_{d(V_\mu)}))=(u^{\otimes d},\cdot~ v)$ for some $v\in V_\mu^{\otimes d}$.
The last two facts imply that any closed subvariety in $B_-\setminus G$ is the zero locus of a finite set of generalized Pl\"{u}cker coordinates. 

\end{proof}

\vskip 0.3cm
\noindent
{\em Proof of part (ii) of Lemma \ref{Jl1}.}

Observe that the projection $pr: G\to B_-\setminus G$ \index[not]{p@$pr$} is an open map \index{map, open} by Theorem 5.5.5 in \cite{Sp}, so that for any $U\subset B_-\setminus G$ one has $pr^{-1}(\overline{U})=\overline{pr^{-1}(U)}$, where the closures are taken with respect to the Zariski topology. In particular, 
\begin{equation}\label{pclo}
pr^{-1}(\overline{B_-\setminus B_-g G_2\ldots G_k})=\overline{pr^{-1}(B_-\setminus B_-g G_2\ldots G_k)}=\overline{B_-g G_2\ldots G_k}.
\end{equation}

By Lemma \ref{zerolp} the closure $\overline{B_-\setminus B_-g G_2\ldots G_k}$ is the zero locus of a finite set of generalized Pl\"{u}cker coordinates. Therefore by (\ref{pclo}) $h \in \overline{B_-g G_2\ldots G_k}$ if and only if all generalized Pl\"{u}cker coordinates from this set vanish on the class $[h]$ of $h$ in $B_-\setminus G$. But if $(u,\cdot~ v)$ is the matrix element corresponding to one of these Pl\"{u}cker coordinates then by the definition this Pl\"{u}cker coordinate vanishes on $[h]\in B_-\setminus G$ if and only if $(u,h v)=0$. Therefore the matrix elements $(u,\cdot~ v)$ corresponding to the Pl\"{u}cker coordinates the common zero locus of which is $\overline{B_-\setminus B_-g G_2\ldots G_k}$ generate $\mathcal{V}_G(\overline{B_-g G_2\ldots G_k})=\mathcal{V}_G(B_-g G_2\ldots G_k)$. 

If $(u, \cdot~ v)\in \mathcal{V}_G(B_-g G_2\ldots G_k)$, where $u$ is a highest weight vector in a finite-dimensional representation $V$ of $\g$, and $v\in V$ then by part (i) of Lemma \ref{Jl1}
$(u,x_1gx_2\ldots x_kv)=0$ for any $x_i\in U(\g_i)$, $i=1,\ldots,k$, where $\g_1=\b_-$. Since
$\omega(U(\b_-))u=U(\b_+)u=\mathbb{C}u$, the last condition is equivalent to $(u,gx_2\ldots x_kv)=0$ for any $x_i\in U(\g_i)$, $i=2,\ldots,k$. Thus the ideal $\mathcal{V}_G(B_- g G_2\ldots G_k)\subset \mathbb{C}[G]$ is generated by the matrix elements of the form $(u, \cdot~ v)\in \mathbb{C}[G]$, where $u$ is a highest weight vector in a finite-dimensional representation $V$ of $\g$, and $v\in V$ satisfies $(u, gx_2\ldots x_kv)=0$ for any $x_i\in U(\g_i)$, $i=2,\ldots,k$. This completes the proof.

\qed


\section{Bibliographic comments}

\pagestyle{myheadings}
\markboth{CHAPTER \thechapter.~Q-W--ALGEBRAS}{\thesection.~BIBLIOGRAPHIC COMMENTS}

\setcounter{equation}{0}
\setcounter{theorem}{0}

The results on Poisson--Lie groups used in this book can be found in \cite{ChP}, \cite{Dm}, \cite{fact}, \cite{dual}. 

Proposition \ref{PLiecorr} is stated in \cite{ChP} as Theorem 1.3.2 and Proposition \ref{bpm} and the relevant properties of classical r-matrices can be found in \cite{BD}, \cite{rmatr}.

The result stated in Proposition \ref{pbff} can be found in \cite{dual}, Section 2.

Q-W--algebras for realizations of quantum groups associated to Weyl group elements were introduced in \cite{SThes,S2} in the case of Coxeter elements and in \cite{S10} in the general situation. However, in the definitions given in those papers other forms of the quantum group are used. The definition of q-W--algebras in this book is more close to the one given in \cite{SDM}; it uses the Ad locally finite part of the quantum group (see \cite{Jos}, Chapter 7) which reduces to the algebra of regular functions on $G$ when $q=1$. However, in this book we define all algebras over slightly different rings.

The exposition in Sections \ref{qplgroups} and \ref{wpsred} follows \cite{S10,SDM} with some appropriate modifications. 

The presentation of the results on Poisson reduction in Section \ref{poisred} is close to \cite{S3}, Section 2.3. More details on the notion and the properties of dual pairs and Poisson reduction can be found in \cite{RIMS}, and for statements related to the moment map for Poisson--Lie group actions the reader is referred to \cite{Lu}. 

The original definition of the Poisson algebras $W^s(G)$ using the classical Poisson reduction only was given in \cite{S6}. 

The definition of the classical Zhelobenko type operator $\Pi$ in Section \ref{pZhel} is a modified version of the definition given in \cite{S13}.


\chapter{Zhelobenko type operators for q-W--algebras}\label{part4}

\pagestyle{myheadings}
\markboth{CHAPTER~\thechapter.~ZHELOBENKO TYPE OPERATORS FOR Q-W--ALGEBRAS}{\thesection.~A QUANTUM ANALOGUE OF THE LEVEL SURFACE OF THE MOMENT MAP}

In this chapter we define a quantum analogue $\Pi_c^q$ of the operator $\Pi_c$ and apply it to describe q-W--algebras. Observe that the operator $\Pi$ is defined using the conjugation action and operators of multiplication by the functions $\varphi_{jk}$. The conjugation action has a natural quantum group analogue, the adjoint action. But multiplication by functions in $\mathbb{C}[G]$ is quite far from the multiplication in the algebra $\mathbb{C}_{\mathcal{B}}^s[G_*]$ which is used in the definition of q-W--algebras. However, using isomorphism (\ref{gg*}) of ${\rm Ad}_s$--modules $\mathbb{C}_{\mathcal{B}}^s[G_*]$ and $\mathbb{C}_{\mathcal{B}}^s[G]$ we can try to describe q-W--algebras in terms of the algebra $\mathbb{C}_{\mathcal{B}}^s[G]$ multiplication in which is more closely related to that of $\mathbb{C}[G]$. Therefore it is natural to expect that a quantum analogue of the operator $\Pi_c$, if it exists at all, should be defined in terms of the adjoint action and of operators of multiplication in $\mathbb{C}_{\mathcal{B}}^s[G]$ using appropriate quantum analogues of formulas (\ref{Pip}) and (\ref{pic}). We shall see that this conjecture is almost correct. In fact, $\mathbb{C}_{\mathcal{B}}^s[G]$ should be replaced with a certain localization. More precisely, recall that the operator $\Pi_c$ is defined using the functions $\varphi_{jk}$ given by (\ref{Apc}). Natural analogues of matrix elements which appear in formula (\ref{Apc}) can be defined. But formula (\ref{Apc}) contains some artificial denominators zeroes of which do not correspond to any singularities of the functions $\varphi_{jk}$, which are in fact regular, in the formula for $\Pi_c$. It turns out that in the quantum case formulas similar to (\ref{Apc}) make sense but the denominators in them are not canceled in the formula for $\Pi_c^q$, and we are forced to use localizations containing all such denominators. This will also force us to replace the algebra $W_{\mathcal{B}}^s(G)$ with a certain localization $W_{\mathcal{B}}^{s,loc}(G)$ of it.

The main difficulty in defining the quantum analogue $\Pi_c^q$ of the operator $\Pi_c$ is that the proof of the fact that the operator defined by (\ref{Pi1}) is a projection operator onto $W^s(G)$ is based on isomorphism (\ref{cross1}) a quantum counterpart of which does not make sense. Recall that $W_{\mathcal{B}}^s(G)$ is the space of invariants with respect to the adjoint action of $\mathbb{C}_{\mathcal{B}}^s[M_+]$ on $Q_{\mathcal{B}}$. Although quantum analogues of operators (\ref{Pip}) can be defined the proof of the fact that their composition similar to (\ref{pic}) is an operator with the image being the localization $W_{\mathcal{B}}^{s,loc}(G)$ of $W_{\mathcal{B}}^s(G)$ should only use the algebra structure of $\mathbb{C}_{\mathcal{B}}^s[G]$, the properties of the adjoint action of $\mathbb{C}_{\mathcal{B}}^s[M_+]$ on $Q_{\mathcal{B}}$, and the structure of $Q_{\mathcal{B}}$. These are the only technical tools in our disposal.  

Thus our first task is to describe in terms of $\mathbb{C}_{\mathcal{B}}^s[G]$ the $\mathbb{C}_{\mathcal{B}}^s[M_+]$--module $Q_{\mathcal{B}}$ originally defined using $\mathbb{C}_{\mathcal{B}}^s[G_*]$. In the classical case this would correspond to describing the vanishing ideal of the closed subvariety $N_-Z^ssM_-= N_-Z^ssM_-^s\subset G$ as by Lemma \ref{redreg} and by Theorem \ref{var} $Q_{\mathcal{B}}/(q^{\frac{1}{d{\bar{r}}^2}}-1)Q_{\mathcal{B}}\simeq \mathbb{C}[N_-Z^ssM_-^s]$. It turns out that not all elements of $\mathbb{C}[G]$ generating the vanishing ideal of $N_-Z^ssM_-^s$ have nice quantum counterparts in $\mathbb{C}_{\mathcal{B}}^s[G]$. Recall that $\mathbb{C}[G]$ is $P\times P$--graded via the left and the right regular action of $H$ on $G$. The subvariety $N_-Z^ssM_-^s\subset G$ is closed and some generators of its vanishing ideal belong to the graded components and some do not. It turns out that at least some of the generators of the latter type have no nice quantum counterparts. But for our purposes it suffices to replace $N_-Z^ssM_-^s$ with a larger set $N_-L^ssM_-^s$ the vanishing ideal $J^{11}$ of which, considered in the last two sections of the previous chapter, has a nice quantum counterpart $I^{11}_{\mathcal{B}}$ in $\mathbb{C}_{\mathcal{B}}^s[G]$. This counterpart is described in Proposition \ref{kerphi} and its image under the natural map  $\mathbb{C}_{\mathcal{B}}^s[G]\simeq \mathbb{C}_{\mathcal{B}}^s[G_*]\rightarrow Q_{\mathcal{B}}$ is zero. 

After recollecting some facts on the algebra $\mathbb{C}_{\mathcal{B}}^s[G]$ and on the adjoint action in Section \ref{CG} we study properties of $I^{11}_{\mathcal{B}}$ in Section \ref{vanid}. 

In order to show that $\Pi_c^q$ is an operator with the image $W_{\mathcal{B}}^{s,loc}(G)$ we shall need some relations which resemble relations in the algebras $\mathbb{C}[G^j]/{J^{jk}}$. 

In Section \ref{localG} we introduce the localizations mentioned above and study their properties and the relevant properties of the adjoint action. The results obtained in Sections \ref{Qlev}, \ref{CG}, \ref{vanid}, \ref{hvan} and \ref{localG} are prerequisites for the study of the properties of the quantum analogues $P_{jk}$ of the operators $\Pi_{jk}$ and of their compositions in Section \ref{Zheldef}, the main properties being summarized in Proposition \ref{BasP}. In Proposition \ref{Bpbas1} we also define quantum analogues of monomials in variables $\varphi_{jk}$ which play a crucial role in the study of equivariant modules over a quantum group and in the proof of the De Concini--Kac--Procesi conjecture.

Finally in Section \ref{qWdescr} we prove that the image of the operator $\Pi_c^q$ almost coincides with the localization $W_{\mathcal{B}}^{s,loc}(G)$ of the algebra $W_{\mathcal{B}}^s(G)$.


\section{A quantum analogue of the level surface of the moment map for q-W--algebras}\label{Qlev}

\setcounter{equation}{0}
\setcounter{theorem}{0}

In this section we describe a quantum counterpart $I^{11}_{\mathcal{B}}\subset \mathbb{C}_{\mathcal{B}}^s[G]$ of the vanishing ideal $J^{11}$. As we mentioned in the introduction to this chapter isomorphism (\ref{gg*}) of ${\rm Ad}_s$--modules $\mathbb{C}_{\mathcal{B}}^s[G_*]$ and $\mathbb{C}_{\mathcal{B}}^s[G]$ plays a central role in the description of $I^{11}_{\mathcal{B}}$. Note that since both $\omega_0$ and the antipode $S_s$ are algebra antiautomorphisms, the compositions $\omega_0 S_s^{-1}$ and $S_s\omega_0$ are algebra automorphisms.
For technical reasons, to ensure that the specializations at $q^{\frac{1}{d{\bar{r}}^2}}=1$ of the objects defined below agree with their algebraic group counterparts defined in the previous chapter, we shall replace the adjoint action of $U_{\mathcal{B}}^{s,res}(\g)$ and of $\mathbb{C}_{\mathcal{B}}^s[G^*]$ on $\mathbb{C}_{\mathcal{B}}^s[G]$ with the twisted adjoint action defined by
\begin{equation}\label{ado}
({\rm Ad}_s^0 xf)(w)=f(\omega_0 S_s^{-1}({\rm Ad}'_sx(S_s\omega_0 w))), \index[not]{A@${\rm Ad}_s^0$}
\end{equation}
where $f\in \mathbb{C}_{\mathcal{B}}^s[G], w \in U_{\mathcal{B}}^{s,res}(\g)$ and $x\in U_{\mathcal{B}}^{s,res}(\g)$ or $x\in \mathbb{C}_{\mathcal{B}}^s[G^*]$. Since $\omega_0'$ is an algebra antiautomorphism and a coalgebra automorphism and $\tau$ is an algebra involution we can also write
\begin{equation}\label{ado'}
({\rm Ad}_s^0 xf)(w)=f((\omega_0 S_s^{-1})(x^1)w\omega_0 (x^2))=f(S_s^\tau(\omega_0 x^1)w\omega_0 (x^2))=f({\rm Ad}_s^\tau\omega_0(x)(w))=f(\tau({\rm Ad}_s\omega_0'(x)(\tau^{-1} (w))), 
\end{equation}
where at the last step we used (\ref{Adups}).

Consider isomorphism (\ref{gg*}) twisted by the $\mathcal{B}$--linear automorphism $\omega_0 S_s^{-1}$, 
\begin{equation}\label{gg**}
\varphi:\mathbb{C}_{\mathcal{B}}^s[G]\rightarrow \mathbb{C}_{\mathcal{B}}^s[G_*], f\mapsto (id\otimes f) (id \otimes \omega_0 S_s^{-1})((\mathcal{R}_s)_{21}\mathcal{R}_s). \index[not]{f@$\varphi$}
\end{equation}
If $\kappa=1$, by the definition of $Q_{\mathcal{B}}=\rho_{\chi^{s}_q}(\mathbb{C}_{\mathcal{B}}^s[G_*])$, where $\rho_{\chi^{s}_q}:\mathbb{C}_{\mathcal{B}}^s[G^*]\rightarrow \mathbb{C}_{\mathcal{B}}^s[G^*]/I_{\mathcal{B}}=Q_{\mathcal{B}}'$ is the canonical projection, and by Lemma \ref{AdM} $\varphi$ induces a homomorphism of $\mathbb{C}_{\mathcal{B}}^s[M_+]$--modules
\begin{equation}\label{gg***}
\phi: \mathbb{C}_{\mathcal{B}}^s[G]\rightarrow Q_{\mathcal{B}}, \phi(f)=\varphi(f)1, \index[not]{f@$\phi$}
\end{equation}
where $\mathbb{C}_{\mathcal{B}}^s[G]$ is equipped with the restriction of action (\ref{ado}) to $\mathbb{C}_{\mathcal{B}}^s[M_+]\subset \mathbb{C}_{\mathcal{B}}^s[G^*]$, $Q_{\mathcal{B}}$ with the action induced by the adjoint action ${\rm Ad}_s$ of $\mathbb{C}_{\mathcal{B}}^s[M_+]$ (see Lemma \ref{AdM}), and $1$ is the image of $1\in \mathbb{C}_{\mathcal{B}}^s[G_*]$ in $Q_{\mathcal{B}}$ under $\rho_{\chi^{s}_q}$.

Now we make some preparations to state a quantum counterpart of Proposition \ref{constrt}. The proof of Proposition \ref{constrt} was based on the Chevalley commutation relations between one--parameter subgroups in $G$ as described in Lemma \ref{comm} and on formula (\ref{sarep}) for representatives of Weyl group elements in $N_G(H)$. In the quantum case instead of the Chevalley commutation relations we have commutation relations between quantum root vectors, and the Weyl group is replaced with the corresponding braid group generators of which are also expressed in terms of generators of the quantum group by formula (\ref{T1}). However, in the quantum case the generators of the braid group do not square to identity automorphisms of the quantum group and we are only allowed to use the braid group relations. The action of the braid group on quantum root vectors is also very difficult to control. It is much more complicated than the adjoint action of representatives in $N_G(H)$ of  Weyl group elements on root vectors in $\g$. All this brings additional complications to the proof of Proposition \ref{kerphi} below which is a quantum counterpart of Proposition \ref{constrt}.

As before, we assume that for a Weyl group element $s\in W$  a normally ordered system of positive roots $\Delta_+$ associated to $s$ is fixed as in Definition \ref{circorddef}, and denote by $\beta_1,\ldots ,\beta_D$ the ordered roots in $\Delta_+$. So if $\alpha_1, \ldots, \alpha_l$ are the simple roots in $\Delta_+$ and $\overline{w}=s_{i_1}\ldots s_{i_D}$ the corresponding decomposition of the longest element $\overline{w}\in W$ then 
$$
\beta_1=\alpha_{i_1},\beta_2=s_{i_1}\alpha_{i_2},\ldots ,\beta_D=s_{i_1}\ldots s_{i_{D-1}}\alpha_{i_D}.
$$

Unless explicitly stated otherwise, we shall assume that all quantum root vectors are defined using this normal ordering of $\Delta_+$. 

Let $U_q^{res}(w'(\b_+))=U_{U_q^{res}(H)}^{res}([\beta_{k_{l'}+1},-\beta_{k_{l'}}])$  \index[not]{U@$U_q^{res}(w'(\b_+))$} be the subalgebra in $U_q^{res}(\g)$
generated by the elements 
$$
({X}^-_{\beta_1})^{(r_{1})}, \ldots, ({X}^-_{\beta_{k_{l'}}})^{(r_{k_{l'}})}, (X^+_{\beta_{k_{l'}+1}})^{(r_{k_{l'}+1})},\ldots ,(X^+_{\beta_D})^{(r_{D})}, r_i\in \mathbb{N}, i=1,\ldots, D,
$$
where $\beta_{k_{l'}}=\gamma_{l'}$, and by $U_q^{res}(H)$. The reason for the use of the symbol $U_q^{res}(w'(\b_+))$ for this subalgebra will be clear later from the definition of $w'\in W$ in formula (\ref{w'}) and Remark \ref{w'bc}. 

We equip the $\mathcal{B}$--module $\mathbb{C}_{\mathcal{B}}^s[G]$ with the multiplication induced by the comultiplication $\Delta_s^\tau$ on $U_{\mathcal{B}}^{s,res}(\g)$ and denote this multiplication by $\otimes$. Sometimes for brevity we shall omit the symbol $\otimes$ if it does not cause any confusion. We shall only use this multiplication on $\mathbb{C}_{\mathcal{B}}^s[G]$.

The twist of the comultiplication $\Delta_s$ by the involution $\tau$ in the definition of the multiplication $\otimes$ is needed since $\tau$ appears in definition (\ref{ado}) of the action ${\rm Ad}_s^0$ via the algebra anti-involution $\omega_0=\tau\omega_0'$, and $\tau$ is neither coalgebra homomorphism nor coalgebra anti-homomorphism.
 
We recall the notation introduced in Section \ref{qplgroups} for elements of $\mathbb{C}_{\mathcal{B}}^s[G]$. Let $V^{res}$ be a $U_{\mathcal{B}}^{s,res}(\g)$--lattice in a finite rank $U_h(\g)$--module $V$. Recall that there is a contravariant non--degenerate form $(~\cdot~,~\cdot~)$ on $V$, which gives rise to a contravariant non--degenerate form on the $U_q^s(\g)$--module $V_q=V^{res}\otimes_{\mathcal{B}}\mathbb{C}(q^{\frac{1}{d{\bar{r}}^2}})$, such that $(u,xv)=(\omega(x)u,v)$ for any $u,v\in V_q$, $x\in U_q^s(\g)$, and different weight subspaces $(V_q)_\lambda\subset V_q$, $\lambda \in P$ are orthogonal with respect to this form. 

Assume that $u\in V_q$ is such that $(u,w)\in \mathcal{B}$ for any $w\in V^{res}$. Then $u^*(~\cdot~):=(u,\cdot~)$ is an element of the dual module  ${V^{res}}^*$. Since $V_q$ and $V^{res}$ are of finite ranks, the weight subspaces of $V^{res}$ are free $\mathcal{B}$--modules, and the restriction of the contravariant form $(~\cdot~,~\cdot~)$ to each weight subspace of $V_q$ is non--degenerate, all elements of ${V^{res}}^*$ can be obtained this way. Clearly, for any $v\in V^{res}$ $u^*(~\cdot~ v)=(u,\cdot~ v)\in \mathbb{C}_{\mathcal{B}}^s[G]$, and by the definition $\mathbb{C}_{\mathcal{B}}^s[G]$ is spanned by such elements.  

Let $c_\beta\in \mathcal{B}$, $\beta\in \Delta_{\m_+}$ be elements such that 
$$
c_\beta=\left\{\begin{array}{ll} \bar{k}_i\in \mathcal{B} & {\rm if}~\beta=\gamma_i, i=1,\ldots, l'\\ 0 & {\rm otherwise} \end{array} \right. .
$$

As we observed in Proposition \ref{qG*} (v) the elements $\tilde e_{\beta_1}^{r_1}\ldots \tilde e_{\beta_D}^{r_D}\bar{V}_i \tilde f_{\beta_D}^{m_D}\ldots \tilde f_{\beta_1}^{m_1}$ with $r_j,m_j,i\in \mathbb{N}$, $j=1,\ldots , D$ form a $\mathcal{B}$--basis in $\mathbb{C}_{\mathcal{B}}^s[G^*]$, where we keep the notation introduced in Proposition \ref{Qdefr}, so that $\Delta_{\m_+}=\{\beta_1,\ldots, \beta_c\}\subset \Delta_+=\{\beta_1,\ldots, \beta_D\}$.  

Clearly, the elements $\tilde e_{\beta_1}^{r_1}\ldots \tilde e_{\beta_D}^{r_D}\bar{V}_i \tilde f_{\beta_D}^{m_D}\ldots \tilde f_{\beta_{c+1}}^{m_{c+1}}(\tilde f_{\beta_c}-c_{\beta_c})^{m_c}\ldots (\tilde f_{\beta_1}-c_{\beta_1})^{m_1}$ with $r_j,m_j,i\in \mathbb{N}$, $j=1,\ldots , D$ also form a $\mathcal{B}$--basis in $\mathbb{C}_{\mathcal{B}}^s[G^*]$.
Let $I_{\mathcal{B}}^{\bf k}$ \index[not]{I@$I_{\mathcal{B}}^{\bf k}$} be the $\mathcal{B}$--submodule in $\mathbb{C}_{\mathcal{B}}^s[G^*]$ generated by the elements $$\tilde e_{\beta_1}^{r_1}\ldots \tilde e_{\beta_D}^{r_D}\bar{V}_i \tilde f_{\beta_D}^{m_D}\ldots \tilde f_{\beta_{c+1}}^{m_{c+1}}(\tilde f_{\beta_c}-c_{\beta_c})^{m_c}\ldots (\tilde f_{\beta_1}-c_{\beta_1})^{m_1}$$ with $r_j,m_j,i\in \mathbb{N}$, $j=1,\ldots , D$, and where at least one $m_j>0$ for $j<c+1$. Since these elements are linearly independent they form a $\mathcal{B}$--basis in $I_{\mathcal{B}}^{\bf k}$. 

\begin{proposition}\label{Iq}
Let ${J^{11}_{\mathcal{B}}}'$ \index[not]{J@${J^{11}_\mathcal{B}}'$} be the left ideal in $\mathbb{C}_{\mathcal{B}}^s[G]$ generated by the elements $(u, \cdot~ v)\in \mathbb{C}_{\mathcal{B}}^s[G]$, where $u$ is a highest weight vector in a finite rank representation $V$ of $U_h(\g)$, and $v\in V^{res}$ is such that $(u,\tau(T_sx) v)=0$ for any $x\in U_q^{res}(w'(\b_+))$. Denote $I^{11}_{\mathcal{B}}:=({J^{11}_{\mathcal{B}}}'\otimes_{\mathcal{B}}\mathbb{C}(q^{\frac{1}{d{\bar{r}}^2}}))\cap \mathbb{C}_{\mathcal{B}}^s[G]$. \index[not]{I@$I^{11}_{\mathcal{B}}$} Let $Q_{\mathcal{B}}^{\bf k}$ \index[not]{Q@$Q_{\mathcal{B}}^{\bf k}$} be the image of $\mathbb{C}_{\mathcal{B}}^s[G_*]\subset \mathbb{C}_{\mathcal{B}}^s[G^*]$ under the canonical projection $\mathbb{C}_{\mathcal{B}}^s[G^*]\rightarrow \mathbb{C}_{\mathcal{B}}^s[G^*]/I_{\mathcal{B}}^{\bf k}$. Denote by $1\in Q_{\mathcal{B}}^{\bf k}$ the image of $1\in \mathbb{C}_{\mathcal{B}}^s[G^*]$ in $\mathbb{C}_{\mathcal{B}}^s[G^*]/I_{\mathcal{B}}^{\bf k}$. Then the following statements are true.

(i) $\varphi({{J^{11}_{\mathcal{B}}}'})\subset I_{\mathcal{B}}^{\bf k}\cap \mathbb{C}_{\mathcal{B}}^s[G_*]$ and $\varphi({I^{11}_{\mathcal{B}}})\subset I_{\mathcal{B}}^{\bf k}\cap \mathbb{C}_{\mathcal{B}}^s[G_*]$. 

(ii) If $u$ is a highest weight vector in a finite rank indecomposable representation $V_\lambda$ of $U_h(\g)$ of highest weight $\lambda$ such that $(u,u)=1$ then for any $f\in \mathbb{C}_{\mathcal{B}}^s[G]$ we have  
\begin{equation}\label{varphiIq}
\varphi (f\otimes (u,\cdot~ \tau(T_s^{-1})u))1=c_\lambda \varphi({\rm Ad}_s^0(q^{-(\kappa{1+s \over 1-s }P_{{\h'}}+id)\lambda^\vee})(f))q^{(s^{-1}+id)(id-\kappa P_{\h'})\lambda^\vee}1=
\end{equation}
$$
=c_\lambda q^{(s^{-1}+id)(id-\kappa P_{\h'})\lambda^\vee}\varphi({\rm Ad}_s^0(q^{(-\kappa{1+s \over 1-s }s^{-1}P_{{\h'}}+s^{-1})\lambda^\vee})(f))1\in Q_{\mathcal{B}}^{\bf k},
$$
where $c_\lambda=c_0\prod_{i=1}^{l'} \bar{k}_i^{n_i(\lambda)}$, \index[not]{c@$c_\lambda$} $c_0\in \mathcal{B}^*$ is an invertible element of $\mathcal{B}$ which only depends on $\lambda$, $\gamma_1, \ldots \gamma_{l'}$, and $n_i(\lambda)= \lambda^\vee(\gamma_i)\geq 0$ for $i=1, \ldots, \widetilde{l}$, $n_i(\lambda)= \lambda^\vee(s^1\gamma_i)\geq 0$ for $i=\widetilde{l}+1, \ldots, l'$. \index[not]{n@$n_i(\lambda)$}
The classes in the quotient $\mathbb{C}_{\mathcal{B}}^s[G^*]/I_{\mathcal{B}}^{\bf k}$ of the elements of $\mathbb{C}_{\mathcal{B}}^s[G^*]$ in the right hand side of (\ref{varphiIq}) belong to $Q_{\mathcal{B}}^{\bf k}\subset \mathbb{C}_{\mathcal{B}}^s[G^*]/I_{\mathcal{B}}^{\bf k}$.

In particular,
$$
\varphi ((u,\cdot~ \tau(T_s^{-1})u))1=c_\lambda q^{(s^{-1}+id)(id-\kappa P_{\h'})\lambda^\vee}1\in Q_{\mathcal{B}}^{\bf k}.
$$ 
$q^{(s^{-1}+id)(id-\kappa P_{\h'})\lambda^\vee}1$ should be understood as the class of the element $q^{(s^{-1}+id)(id-\kappa P_{\h'})\lambda^\vee}\in \mathbb{C}_{\mathcal{B}}^s[G^*]$ in the quotient $\mathbb{C}_{\mathcal{B}}^s[G^*]/I_{\mathcal{B}}^{\bf k}$. This class belongs to $Q_{\mathcal{B}}^{\bf k}$.           
\end{proposition}

\begin{proof}
The proof of this proposition is based on Lemma \ref{mainl} which will be proved in the end of this section.

(i) We start proving this proposition by obtaining a useful expression for $(id \otimes \omega_0 S_s^{-1})((\mathcal{R}_s)_{21}\mathcal{R}_s)$.
In order to do that we recall some properties of universal R-matrices (see (\ref{SR}), (\ref{S})),
\begin{equation}\label{SRs}
(S_s\otimes id)\mathcal{R}_s=(id \otimes S_s^{-1})\mathcal{R}_s= {\mathcal{R}_s}^{-1}, (S_s\otimes S_s)\mathcal{R}_s=\mathcal{R}_s,
\end{equation}
Using the first identity above we can write
$$
(\mathcal{R}_s)_{21}\mathcal{R}_s=(\mathcal{R}_s)_{21}(id \otimes S_s)({\mathcal{R}_s}^{-1})=(id \otimes S_s^{-1})({(\mathcal{R}_s)_{21}}^{-1})(id \otimes S_s)({\mathcal{R}_s}^{-1})=(id \otimes S_s)((id \otimes S_s^{-2})({(\mathcal{R}_s)_{21}}^{-1})\circ {\mathcal{R}_s}^{-1}),
$$
where
$$
(a\otimes b)\circ (c\otimes d)=ac\otimes db,~a,b,c,d\in U_h^s(\g).
$$

Now since $\omega_0$ an algebra antiautomorphism we have
\begin{equation}\label{SoR}
(id \otimes \omega_0 S_s^{-1})((\mathcal{R}_s)_{21}\mathcal{R}_s)=(id \otimes \omega_0)((id \otimes S_s^{-2})({(\mathcal{R}_s)_{21}}^{-1})\circ {\mathcal{R}_s}^{-1})= (id \otimes \omega_0 S_s^{-2})({(\mathcal{R}_s)_{21}}^{-1}) (id \otimes \omega_0) ({\mathcal{R}_s}^{-1}).
\end{equation}

Recalling the definition (\ref{rmatrspi}) of $\mathcal{R}_s$ and (\ref{Rsinv1}), (\ref{Rsinv2}) we obtain
\begin{eqnarray}\label{r-1}
{\mathcal{R}_s}^{-1}={\exp}\left[ -h(\sum_{i=1}^lY_i\otimes H_i-
\sum_{i=1}^l \kappa{1+s \over 1-s }P_{{\h'}}H_i\otimes Y_i) \right]\times \\
\times \prod_{\beta \in \Delta_+}
{\exp}_{q_{\beta}^{-1}}[(1-q_{\beta}^{2})f_{\beta} \otimes
e_{\beta}e^{-h\kappa{1+s \over 1-s}P_{{\h'}} \beta^\vee}]= \nonumber
\end{eqnarray}
\begin{eqnarray*}
=\prod_{\beta \in \Delta_+}
{\exp}_{q_{\beta}^{-1}}[(1-q_{\beta}^{2})e^{h(\kappa{1+s \over 1-s}P_{{\h'}} -id)\beta^\vee}f_{\beta} \otimes
e_{\beta}q^{\beta^\vee}]\times \\
\times {\exp}\left[ -h(\sum_{i=1}^lY_i\otimes H_i-
\sum_{i=1}^l \kappa{1+s \over 1-s }P_{{\h'}}H_i\otimes Y_i) \right],
\end{eqnarray*}
where the order of the terms in the products over the positive roots is opposite to the normal ordering in $\Delta_+$.

Let $\rho$ \index[not]{r@$\rho$} be a half of the sum of the positive roots. Using the fact that 
$$
S_s^{-2}={\rm Ad}_s~q^{2\rho^\vee},
$$
which can be derived straightforwardly from the definition (\ref{S_s}) of $S_s$, we also deduce
\begin{eqnarray}\label{rs-2}
(id \otimes S_s^{-2})({(\mathcal{R}_s)_{21}}^{-1})= 
\prod_{\beta \in \Delta_+}
{\exp}_{q_{\beta}^{-1}}[(1-q_{\beta}^{2})q^{-2\beta(\rho^\vee)}e_{\beta}q^{\beta^\vee}\otimes e^{h\left(\kappa{1+s \over 1-s}P_{{\h'}}-id\right) \beta^\vee}f_{\beta}] \times 
\\
\times {\exp}\left[ -h(\sum_{i=1}^l(Y_i\otimes H_i)+
\sum_{i=1}^l\kappa{1+s \over 1-s }P_{{\h'}}H_i\otimes Y_i) \right]. \nonumber
\end{eqnarray}
The order of the terms in the products in the formulas above is such that
the $\alpha$--term appears to the left of the $\beta$--term if $\alpha > \beta$ with respect to the normal
ordering of $\Delta_+$.

Substituting (\ref{r-1}) and (\ref{rs-2}) into (\ref{SoR}) we arrive at the following expression for $(id \otimes \omega_0 S_s^{-1})((\mathcal{R}_s)_{21}\mathcal{R}_s)$
\begin{eqnarray}
(id \otimes \omega_0 S_s^{-1})((\mathcal{R}_s)_{21}\mathcal{R}_s)
=\prod^{\leftarrow}_{\rightarrow}
{\exp}_{q_{\beta}^{-1}}[(1-q_{\beta}^{2})q^{-2\beta(\rho^\vee)}e_{\beta}q^{\beta^\vee}\otimes \omega_0(e^{h\left(\kappa{1+s \over 1-s}P_{{\h'}}-id\right) \beta^\vee}f_{\beta})] \times \nonumber
\\
\times {\exp}\left[ h(\sum_{i=1}^l(Y_i\otimes H_i^f)+
\sum_{i=1}^l \kappa{1+s \over 1-s }P_{{\h'}}H_i\otimes Y_i^f) \right] \times \label{RR}
\\
\times {\exp}\left[ h(\sum_{i=1}^l(Y_i\otimes H_i^r)-
\sum_{i=1}^l \kappa{1+s \over 1-s }P_{{\h'}}H_i\otimes Y_i^r) \right]\times \nonumber \\
\times \prod_{\rightarrow}^{\leftarrow}
{\exp}_{q_{\beta}^{-1}}[(1-q_{\beta}^{2})f_{\beta} \otimes
\omega_0(e_{\beta}e^{-h\kappa{1+s \over 1-s}P_{{\h'}} \beta^\vee})],\nonumber
\end{eqnarray}
where in the product $$\prod^{\leftarrow}_{\rightarrow}$$ the upper (the lower) arrow indicates the order of the terms in the first (the second) factor of the tensor product relative to the normal
ordering of $\Delta_+$, and superscripts $f$ (resp. $r$) indicate that the corresponding term appears in the front (resp. in the rear) of all the other terms in the product. 

Assume now that $u\in V$ has highest weight $\lambda$ and $v\in V^{res}$ is any vector of weight $\mu$ such that $g(~\cdot~):=(u,\cdot~ v)\in \mathbb{C}_{\mathcal{B}}^s[G]$. Observe that by (\ref{omega0}) and (\ref{fo0}) all elements $z_\beta:=\omega_0(e^{h\left(\kappa{1+s \over 1-s}P_{{\h'}}-id\right) \beta^\vee}f_{\beta})$ in the first product of the q-exponentials in (\ref{RR}) have strictly negative weights. As $u$ has the highest possible weight in $V$, we deduce by Lemma \ref{hwv} (i) that only constant terms in the expansions (\ref{expq}) of the q-exponentials in the first product in (\ref{RR}) will contribute to the formula for $\varphi (g)=(id\otimes g)(id \otimes \omega_0 S_s^{-1})((\mathcal{R}_s)_{21}\mathcal{R}_s)$. Using also the definition of the action of the generators $H_i^r$, $H_i^f$, $Y_i^r$, $Y_i^f$ on $u$ and $v$ we obtain from (\ref{RR})
\begin{eqnarray*}
\varphi ((u,\cdot~ v))=(id\otimes g)(id \otimes \omega_0 S_s^{-1})((\mathcal{R}_s)_{21}\mathcal{R}_s)
=q^{\lambda^\vee+\kappa{1+s \over 1-s }P_{{\h'}}\lambda^\vee+\mu^\vee-\kappa{1+s \over 1-s }P_{{\h'}}\mu^\vee}\times 
 \\
\times (id\otimes g)(\prod_{\rightarrow}^{\leftarrow}
{\exp}_{q_{\beta}^{-1}}[(1-q_{\beta}^{2})f_{\beta} \otimes
\omega_0(e_{\beta}e^{-h\kappa{1+s \over 1-s}P_{{\h'}} \beta^\vee})]). 
\end{eqnarray*}

Using the definition (\ref{expq}) of the q-exponential and formulas (\ref{omega0}) and (\ref{eo0}) this can be rewritten as follows
$$
\varphi ((u,\cdot~ v))=q^{\lambda^\vee+\kappa{1+s \over 1-s }P_{{\h'}}\lambda^\vee+\mu^\vee-\kappa{1+s \over 1-s }P_{{\h'}}\mu^\vee}\times 
$$
\begin{equation}\label{fpr1}
\times  \hspace{-1em}\sum_{\tiny\begin{array}{c} m_i\in \mathbb{N} \\ i=1,\ldots D \end{array}} \hspace{-1em} \bar{d}(m_1,\ldots,m_{D})  \tilde f_{\beta_D}^{m_D}\ldots \tilde f_{\beta_1}^{m_1} (u,
\tau(e^{-m_1 h\kappa{1+s \over 1-s}P_{{\h'}} \beta_1^\vee} e_{\beta_1}^{(m_1)}\ldots e^{-m_D h\kappa{1+s \over 1-s}P_{{\h'}} \beta_D^\vee}e_{\beta_D}^{(m_D)})v), 
\end{equation}
where $\bar{d}(m_1,\ldots,m_{D})\in \mathcal{B}^*$, and the sum is finite as the elements $e_{\beta}^{(n)}$ act by zero on $V$ for large enough $n$.

Recalling that according to the definition in Proposition \ref{rootss} $e_\beta^{(n)}=(X_\beta^+)^{(n)}q^{nK_s\beta^\vee}$, where $q^{nK_s\beta^\vee}=e^{n h\frac{\kappa}{2}{1+s \over 1-s}P_{{\h'}} \beta^\vee}$, and using commutation relations (\ref{Ccomm}) we can write 
$$
e^{-m_1 h\kappa{1+s \over 1-s}P_{{\h'}} \beta_1^\vee}e_{\beta_1}^{(m_1)}\ldots e^{-m_D h\kappa{1+s \over 1-s}P_{{\h'}} \beta_D^\vee}e_{\beta_D}^{(m_D)}=
$$
$$
=e^{-m_1 h\kappa{1+s \over 1-s}P_{{\h'}} \beta_1^\vee} (X_{\beta_1}^+)^{(m_1)}q^{m_1K_s\beta_1^\vee}
\ldots e^{-m_D h\kappa{1+s \over 1-s}P_{{\h'}} \beta_D^\vee}(X_{\beta_D}^+)^{(m_D)}q^{m_DK_s\beta_D^\vee}=
$$
\begin{equation}\label{eX}
=b(m_1,\ldots,m_D)(X_{\beta_1}^+)^{(m_1)}
\ldots (X_{\beta_D}^+)^{(m_D)}q^{-K_s(m_1\beta_1^\vee+\ldots+m_D\beta_D^\vee)},
\end{equation}
where $b(m_1,\ldots,m_D)\in \mathcal{B}^*$.

Next, since $v$ has weight $\mu$, we can rewrite (\ref{fpr1}) using (\ref{eX}) as follows
$$
\varphi ((u,\cdot~ v))=q^{\lambda^\vee+\kappa{1+s \over 1-s }P_{{\h'}}\lambda^\vee+\mu^\vee-\kappa{1+s \over 1-s }P_{{\h'}}\mu^\vee}\times 
$$
\begin{equation}\label{fpr2}
\times  \hspace{-1em} \sum_{\tiny\begin{array}{c} m_i\in \mathbb{N} \\ i=1,\ldots D \end{array}} \hspace{-1em} \bar{d}'(m_1,\ldots,m_{D})  \tilde f_{\beta_D}^{m_D}\ldots \tilde f_{\beta_1}^{m_1} (u,
\tau((X_{\beta_1}^+)^{(m_1)}\ldots (X_{\beta_D}^+)^{(m_D)})v), 
\end{equation}
where $\bar{d}'(m_1,\ldots,m_{D})=\bar{d}(m_1,\ldots,m_{D})\tau(b(m_1,\ldots,m_{D})q^{-\mu(K_s(m_1\beta_1^\vee+\ldots+m_D\beta_D^\vee))})\in \mathcal{B}^*$.

Now using  the definitions of $I_{\mathcal{B}}^{\bf k}$ and of $Q_{\mathcal{B}}^{\bf k}={\rm Im}(\mathbb{C}_{\mathcal{B}}^s[G_*]\rightarrow \mathbb{C}_{\mathcal{B}}^s[G^*]/I_{\mathcal{B}}^{\bf k})$ we have in $Q_{\mathcal{B}}^{\bf k}$
\begin{equation}\label{vp}
\varphi ((u,\cdot~ v))1=q^{\lambda^\vee+\kappa{1+s \over 1-s }P_{{\h'}}\lambda^\vee+\mu^\vee-\kappa{1+s \over 1-s }P_{{\h'}}\mu^\vee}\times 
\end{equation}
$$
\times  \hspace{-2em} \sum_{\tiny\begin{array}{c} r_i, m_j\in \mathbb{N} \\ i=1,\ldots,l' \\ j=c+1,\ldots, D \end{array}} \hspace{-2em} c(r_1,\ldots,r_{l'}, m_{c+1},\ldots, m_D)\prod_{i=1}^{l'}
 \bar{k}_i^{r_i}
 \tilde f_{\beta_D}^{m_D}\ldots \tilde f_{\beta_{c+1}}^{m_{c+1}} (u,
\tau((X_{\gamma_1}^+)^{(r_1)}\ldots (X_{\gamma_{l'}}^+)^{(r_{l'})}(X_{\beta_{c+1}}^+)^{(m_{c+1})}\ldots (X_{\beta_D}^+)^{(m_D)})v)1,
$$
where $c(r_1,\ldots,r_{l'},m_{c+1},\ldots, m_D)\in \mathcal{B}^*$, and the sum is finite.

Lemma \ref{mainl} applied to the products $\tau((X_{\gamma_1}^+)^{(r_1)}\ldots (X_{\gamma_{l'}}^+)^{(r_{l'})})$ in the previous formula and the fact that by the definition of the algebra $U_q^{res}(w'(\b_+))$ one has $(X_{\beta_{c+1}}^+)^{(m_{c+1})}\ldots (X_{\beta_D}^+)^{(m_D)}\in U_q^{res}(w'(\b_+))$, imply that in $Q_{\mathcal{B}}^{\bf k}$
\begin{equation}\label{uvel}
\varphi ((u,\cdot~ v))1= q^{\lambda^\vee+\kappa{1+s \over 1-s }P_{{\h'}}\lambda^\vee+\mu^\vee-\kappa{1+s \over 1-s }P_{{\h'}}\mu^\vee}\times
\end{equation}
$$
\times  \hspace{-2em} \sum_{\tiny\begin{array}{c} r_i, m_j\in \mathbb{N} \\ i=1,\ldots,l' \\ j=c+1,\ldots, D \end{array}} \hspace{-2em} c(r_1,\ldots,r_{l'}, m_{c+1},\ldots, m_D)\prod_{i=1}^{l'}
 \bar{k}_i^{r_i}
 \tilde f_{\beta_D}^{m_D}\ldots \tilde f_{\beta_{c+1}}^{m_{c+1}} (u,
\tau(T_s X(r_1,\ldots, r_{l'})(X_{\beta_{c+1}}^+)^{(m_{c+1})}\ldots (X_{\beta_D}^+)^{(m_D)})v)1=
$$
$$
=\sum_ix_i(u,\tau(T_sy_i)v)1,
$$
where $x_i\in \mathbb{C}_{\mathcal{B}}^s[G_*]$, $X(r_1,\ldots, r_{l'}), y_i\in U_q^{res}(w'(\b_+))$. 

Since every element of $V^{res}$ is the sum of its weight components a formula similar to that which appears in the last line of (\ref{uvel}) holds for $\varphi ((u,\cdot~ v))1$ with arbitrary $(u,\cdot~ v)\in \mathbb{C}_{\mathcal{B}}^s[G]$, where $u$ is a highest weight vector.

If $v$ is chosen in such a way that $(u,\tau(T_sx) v)=0$ for any $x\in U_q^{res}(w'(\b_+))$ we deduce from (\ref{uvel}) that $(u,\tau(T_sy_i)v)=0$. Thus $\varphi ((u,\cdot~ v))1=0$ in $Q_{\mathcal{B}}^{\bf k}$. This implies $\varphi((u,\cdot~ v))\in I_{\mathcal{B}}^{\bf k}\cap \mathbb{C}_{\mathcal{B}}^s[G_*]$, as $\varphi ((u,\cdot~ v)) \in \mathbb{C}_{\mathcal{B}}^s[G_*]$ by the definition of $\varphi$.    

Now using the properties (see formulas (\ref{rmprop}))
$$
(\Delta_s\otimes id)\mathcal{R}_s=(\mathcal{R}_s)_{13}(\mathcal{R}_s)_{23},(id \otimes \Delta_s)\mathcal{R}_s=(\mathcal{R}_s)_{13}(\mathcal{R}_s)_{12},
$$
and the fact that $\omega_0'$ is a coalgebra automorphism (see Section \ref{QGspec}) and $S_s$ is an anti-coautomorphism we get
$$
(id \otimes \Delta_s^\tau)(id \otimes \omega_0 S_s^{-1})((\mathcal{R}_s)_{21}\mathcal{R}_s)=
(id \otimes \omega_0 S_s^{-1}\otimes \omega_0 S_s^{-1})(id \otimes \Delta_s^{opp})((\mathcal{R}_s)_{21}\mathcal{R}_s)=
$$
$$
=(id \otimes \omega_0 S_s^{-2}\otimes \omega_0 S_s^{-2})({(\mathcal{R}_s)_{31}}^{-1}{(\mathcal{R}_s)_{21}}^{-1})
(id \otimes \omega_0 \otimes \omega_0)({(\mathcal{R}_s)_{12}}^{-1}{(\mathcal{R}_s)_{13}}^{-1}).
$$

In the case when $v$ has weight $\mu$, from this identity we obtain, similarly to (\ref{vp}), that for any $f\in \mathbb{C}_{\mathcal{B}}^s[G]$ in $Q_{\mathcal{B}}^{\bf k}$
\begin{equation}\label{phf}
\varphi (f\otimes (u,\cdot~ v))1=(id\otimes f\otimes (u,\cdot~ v))((id \otimes \Delta_s^\tau)(id \otimes \omega_0 S_s^{-1})((\mathcal{R}_s)_{21}\mathcal{R}_s))1 =  
\end{equation}
$$
=q^{(\kappa{1+s \over 1-s }P_{{\h'}}+id)\lambda^\vee}\varphi(f)q^{(-\kappa{1+s \over 1-s }P_{{\h'}}+id)\mu^\vee} \times 
$$
$$
\times  \hspace{-2em} \sum_{\tiny\begin{array}{c} r_i, m_j\in \mathbb{N} \\ i=1,\ldots,l' \\ j=c+1,\ldots, D \end{array}} \hspace{-2em} c(r_1,\ldots,r_{l'}, m_{c+1},\ldots, m_D)\prod_{i=1}^{l'}
 \bar{k}_i^{r_i}
 \tilde f_{\beta_D}^{m_D}\ldots \tilde f_{\beta_{c+1}}^{m_{c+1}} (u,
\tau((X_{\gamma_1}^+)^{(r_1)}\ldots (X_{\gamma_{l'}}^+)^{(r_{l'})}(X_{\beta_{c+1}}^+)^{(m_{c+1})}\ldots (X_{\beta_D}^+)^{(m_D)})v)1= 
$$
$$
=q^{(\kappa{1+s \over 1-s }P_{{\h'}}+id)\lambda^\vee}\varphi(f)q^{(-\kappa{1+s \over 1-s }P_{{\h'}}+id)\mu^\vee} \times 
$$
$$
\times  \hspace{-2em} \sum_{\tiny\begin{array}{c} r_i, m_j\in \mathbb{N} \\ i=1,\ldots,l' \\ j=c+1,\ldots, D \end{array}} \hspace{-2em} c(r_1,\ldots,r_{l'}, m_{c+1},\ldots, m_D)\prod_{i=1}^{l'}
 \bar{k}_i^{r_i}
 \tilde f_{\beta_D}^{m_D}\ldots \tilde f_{\beta_{c+1}}^{m_{c+1}} (u,
\tau(T_s X(r_1,\ldots, r_{l'})(X_{\beta_{c+1}}^+)^{(m_{c+1})}\ldots (X_{\beta_D}^+)^{(m_D)})v)1,
$$
which implies, similarly to (\ref{uvel}), that for arbitrary $v$
$$
\varphi (f\otimes (u,\cdot~ v))1=\sum_ix_i'\varphi(f)x_i''(u,\tau(T_sy_i')v)1,
$$
where $x_i', x_i''\in \mathbb{C}_{\mathcal{B}}^s[G_*]$, $y_i'\in U_q^{res}(w'(\b_+))$. Hence $\varphi (f\otimes (u,\cdot~ v))1=0$ in $Q_{\mathcal{B}}^{\bf k}$ by the choice of $v$, i.e. $\varphi({{J^{11}_{\mathcal{B}}}'})\subset I_{\mathcal{B}}^{\bf k}\cap \mathbb{C}_{\mathcal{B}}^s[G_*]$ as $\varphi (f\otimes (u,\cdot~ v)) \in \mathbb{C}_{\mathcal{B}}^s[G_*]$ by the definition of $\varphi$.

In order to show that $\varphi({I^{11}_{\mathcal{B}}})\subset I_{\mathcal{B}}^{\bf k}\cap \mathbb{C}_{\mathcal{B}}^s[G_*]$ we naturally extend $\varphi$ to and ${\rm Ad}_s$--module isomorphism $\varphi:\mathbb{C}_q^s[G]\rightarrow {\mathbb{C}}_q^s[{G}_*]$, where $\mathbb{C}_q^s[G]=\mathbb{C}_{\mathcal{B}}^s[G]\otimes_{\mathcal{B}}\mathbb{C}(q^{\frac{1}{d{\bar{r}}^2}})$, \index[not]{C@$\mathbb{C}_q^s[G]$} $\mathbb{C}_q^s[G_*]=\mathbb{C}_{\mathcal{B}}^s[G_*]\otimes_{\mathcal{B}}\mathbb{C}(q^{\frac{1}{d{\bar{r}}^2}})$. \index[not]{C@$\mathbb{C}_q^s[G_*]$} By the definition of ${I^{11}_{\mathcal{B}}}$ we have
$\varphi({I^{11}_{\mathcal{B}}})\subset (I_{\mathcal{B}}^{\bf k}\otimes_{\mathcal{B}}\mathbb{C}(q^{\frac{1}{d{\bar{r}}^2}}))\cap \mathbb{C}_{\mathcal{B}}^s[G_*]$ as obviously $\varphi({{J^{11}_{\mathcal{B}}}'}\otimes_{\mathcal{B}}\mathbb{C}(q^{\frac{1}{d{\bar{r}}^2}}))\subset I_{\mathcal{B}}^{\bf k}\otimes_{\mathcal{B}}\mathbb{C}(q^{\frac{1}{d{\bar{r}}^2}})$ since we already proved that $\varphi({{J^{11}_{\mathcal{B}}}'})\subset I_{\mathcal{B}}^{\bf k}\cap \mathbb{C}_{\mathcal{B}}^s[G_*]$, and $\varphi(\mathbb{C}_{\mathcal{B}}^s[G])\subset \mathbb{C}_{\mathcal{B}}^s[G_*]$.

We also have $(I_{\mathcal{B}}^{\bf k}\otimes_{\mathcal{B}}\mathbb{C}(q^{\frac{1}{d{\bar{r}}^2}}))\cap \mathbb{C}_{\mathcal{B}}^s[G_*]\subset (I_{\mathcal{B}}^{\bf k}\otimes_{\mathcal{B}}\mathbb{C}(q^{\frac{1}{d{\bar{r}}^2}}))\cap \mathbb{C}_{\mathcal{B}}^s[G^*]$ as $\mathbb{C}_{\mathcal{B}}^s[G_*]\subset \mathbb{C}_{\mathcal{B}}^s[G^*]$.

Recall that by the definition the elements $\tilde e_{\beta_1}^{r_1}\ldots \tilde e_{\beta_D}^{r_D}\bar{V}_i \tilde f_{\beta_D}^{m_D}\ldots \tilde f_{\beta_{c+1}}^{m_{c+1}}(\tilde f_{\beta_c}-c_{\beta_c})^{m_c}\ldots (\tilde f_{\beta_1}-c_{\beta_1})^{m_1}$ with $r_j,m_j,i\in \mathbb{N}$, $j=1,\ldots , D$, and where at least one $m_j>0$ for $j<c+1$ form a $\mathcal{B}$--basis in $I_{\mathcal{B}}^{\bf k}$, and this basis can be completed to a $\mathcal{B}$--basis of $\mathbb{C}_{\mathcal{B}}^s[G^*]$ which consists of the elements
$\tilde e_{\beta_1}^{r_1}\ldots \tilde e_{\beta_D}^{r_D}\bar{V}_i \tilde f_{\beta_D}^{m_D}\ldots \tilde f_{\beta_{c+1}}^{m_{c+1}} (\tilde f_{\beta_c}-c_{\beta_c})^{m_c}\ldots (\tilde f_{\beta_1}-c_{\beta_1})^{m_1}$ with $r_j,m_j,i\in \mathbb{N}$, $j=1,\ldots , D$.

This implies that $(I_{\mathcal{B}}^{\bf k}\otimes_{\mathcal{B}}\mathbb{C}(q^{\frac{1}{d{\bar{r}}^2}}))\cap \mathbb{C}_{\mathcal{B}}^s[G^*]=I_{\mathcal{B}}^{\bf k}$, and hence $\varphi({I^{11}_{\mathcal{B}}})\subset (I_{\mathcal{B}}^{\bf k}\otimes_{\mathcal{B}}\mathbb{C}(q^{\frac{1}{d{\bar{r}}^2}}))\cap \mathbb{C}_{\mathcal{B}}^s[G_*]=((I_{\mathcal{B}}^{\bf k}\otimes_{\mathcal{B}}\mathbb{C}(q^{\frac{1}{d{\bar{r}}^2}}))\cap \mathbb{C}_{\mathcal{B}}^s[G^*])\cap \mathbb{C}_{\mathcal{B}}^s[G_*]=I_{\mathcal{B}}^{\bf k}\cap \mathbb{C}_{\mathcal{B}}^s[G_*]$. This completes the proof of part (i).

(ii) Consider formula (\ref{phf}) with $v=\tau(T_s^{-1})u$,
\begin{equation}\label{phf1}
\varphi (f\otimes (u,\cdot~ \tau(T_s^{-1})u))1=q^{(\kappa{1+s \over 1-s }P_{{\h'}}+id)\lambda^\vee}\varphi(f)q^{s^{-1}(-\kappa{1+s \over 1-s }P_{{\h'}}+id)\lambda^\vee} \times
\end{equation}
$$
\times \hspace{-2em} \sum_{\tiny\begin{array}{c} r_i, m_j\in \mathbb{N} \\ i=1,\ldots,l' \\ j=c+1,\ldots, D \end{array}} \hspace{-2em} c(r_1,\ldots,r_{l'}, m_{c+1},\ldots, m_D)\prod_{i=1}^{l'}
 \bar{k}_i^{r_i}
 \tilde f_{\beta_D}^{m_D}\ldots \tilde f_{\beta_{c+1}}^{m_{c+1}} (u,
\tau((X_{\gamma_1}^+)^{(r_1)}\ldots (X_{\gamma_{l'}}^+)^{(r_{l'})}(X_{\beta_{c+1}}^+)^{(m_{c+1})}\ldots (X_{\beta_D}^+)^{(m_D)}T_s^{-1})u)1=
$$
$$
=q^{(\kappa{1+s \over 1-s }P_{{\h'}}+id)\lambda^\vee}\varphi(f)q^{s^{-1}(-\kappa{1+s \over 1-s }P_{{\h'}}+id)\lambda^\vee} \times 
$$
$$
\times  \hspace{-2em}\sum_{\tiny\begin{array}{c} r_i, m_j\in \mathbb{N} \\ i=1,\ldots,l' \\ j=c+1,\ldots, D \end{array}} \hspace{-2em} c(r_1,\ldots,r_{l'}, m_{c+1},\ldots, m_D)\prod_{i=1}^{l'}
 \bar{k}_i^{r_i}
 \tilde f_{\beta_D}^{m_D}\ldots \tilde f_{\beta_{c+1}}^{m_{c+1}} (u,
\tau(T_s X(r_1,\ldots, r_{l'})(X_{\beta_{c+1}}^+)^{(m_{c+1})}\ldots (X_{\beta_D}^+)^{(m_D)}T_s^{-1})u)1.
$$

Observe that the roots $-\beta_1, \ldots, -\beta_{k_{l'}}, \beta_{k_{l'}+1},\ldots ,\beta_D$ corresponding to the elements $$({X}^-_{\beta_1})^{(r_{1})}, \ldots, ({X}^-_{\beta_{k_{l'}}})^{(r_{k_{l'}})}, (X^+_{\beta_{k_{l'}+1}})^{(r_{k_{l'}+1})},\ldots ,(X^+_{\beta_D})^{(r_{D})},$$ which generate together with $U_q^{res}(H)$ the subalgebra $U_q^{res}(w'(\b_+))$, form a minimal segment $\{\alpha\in \Delta, \gamma_{l'}< \alpha\leq - \gamma_{l'}\}$ which is, in fact, a system of positive roots by Lemma \ref{circ+} (i) (see Figure 4). Therefore if an element of $U_q^{res}(w'(\b_+))$ has a non-zero weight, the product of this element and of any other element of $U_q^{res}(w'(\b_+))$ has also a non-zero weight. 

On the other hand, from the expression in the last line of (\ref{phf1}) and from the fact that $\tau$ preserves weights it follows that only terms with $$X(r_1,\ldots, r_{l'})(X_{\beta_{c+1}}^+)^{(m_{c+1})}\ldots (X_{\beta_D}^+)^{(m_D)}\in U_q^{res}(w'(\b_+))$$ of zero weight may give non-trivial contributions to the right hand side of (\ref{phf1}). Since the elements $(X_\beta^+)^{(n)}$, $\beta\in \{\alpha\in \Delta_+:\alpha>\gamma_{l'}\}\subset \{\alpha\in \Delta, \gamma_{l'}< \alpha\leq - \gamma_{l'}\}$ have non-zero weights if $n>0$, only the terms with $m_{c+1}=\ldots=m_D=0$ will give non-trivial contributions to the right hand side of (\ref{phf1}).

Now the second expression for $\varphi (f\otimes (u,\cdot~ \tau(T_s^{-1})u))1$ in (\ref{phf1}) yields
\begin{equation}\label{fo}
\varphi (f\otimes (u,\cdot~ \tau(T_s^{-1})u))1=
\end{equation}
$$
=q^{(\kappa{1+s \over 1-s }P_{{\h'}}+id)\lambda^\vee}\varphi(f)q^{s^{-1}(-\kappa{1+s \over 1-s }P_{{\h'}}+id)\lambda^\vee} \hspace{-1em} \sum_{\tiny\begin{array}{c} r_i\in \mathbb{N} \\ i=1,\ldots,l'\end{array}} \hspace{-1em} d(r_1,\ldots,r_{l'})
 \prod_{i=1}^{l'} \bar{k}_i^{r_i} (u,\tau((X_{\gamma_1}^+)^{(r_1)}\ldots (X_{\gamma_{l'}}^+)^{(r_{l'})}T_s^{-1})u)1,
$$
where $d(r_1,\ldots,r_{l'})\in \mathcal{B}^*$.

Next observe that since $u$ has weight $\lambda$, and different weight spaces in $V$ are orthogonal with respect to the contravariant form $(~\cdot~,~\cdot~)$, only elements $$\tau((X_{\gamma_1}^+)^{(r_1)} \ldots (X_{\gamma_{l'}}^+)^{(r_{l'})}T_s^{-1})u\in V^{res}$$ of weight $\lambda$ can contribute to the right hand side of formula (\ref{fo}).

Because the roots $\gamma_1, \ldots, \gamma_{l'}$ are linearly independent, for $r_1, \ldots, r_{l'}\in \mathbb{N}$ the element $$\tau((X_{\gamma_1}^+)^{(r_1)} \ldots (X_{\gamma_{l'}}^+)^{(r_{l'})}T_s^{-1})u\in V^{res}$$ has weight $\lambda$ if and only if $r_1\gamma_1+ \ldots +\gamma_{l'} r_{l'}=\lambda-s^{-1}\lambda=\sum_{i=1}^{\widetilde{l}}\lambda^\vee(\gamma_i)\gamma_i+\sum_{i=\widetilde{l}+1}^{l'}\lambda^\vee(s^1\gamma_i)\gamma_i$, where to obtain the last expression we also used the fact that the roots $\gamma_1, \ldots \gamma_{\widetilde{l}}$ are mutually orthogonal and the roots $\gamma_{\widetilde{l}+1}, \ldots \gamma_{l'}$ are also mutually orthogonal. This implies $r_i=n_i(\lambda)= \lambda^\vee(\gamma_i)$, $i=1, \ldots, \widetilde{l}$, and $r_i=n_i(\lambda)= \lambda^\vee(s^1\gamma_i)$, $i=\widetilde{l}+1, \ldots, l'$. Note that for $i=\widetilde{l}+1, \ldots, l'$, by Proposition \ref{pord} (i) one has $s^1\gamma_i\in s^1(\Delta_{s^2}^s)\subset \Delta_+^s\setminus \Delta_{s^1}^s\subset \Delta_+$, so $n_i(\lambda)= \lambda^\vee(s^1\gamma_i)\geq 0$ as $\lambda\in P_+$. Also, since $\lambda\in P_+$, we obtain $n_i(\lambda)= \lambda^\vee(\gamma_i)\geq 0$ for $i=1, \ldots, \widetilde{l}$.

By Proposition B.6. from \cite{AJS}, in this case we have $\tau((X_{\gamma_1}^+)^{(n_1(\lambda))} \ldots (X_{\gamma_{l'}}^+)^{(n_{l'}(\lambda))}T_s^{-1})u=\epsilon c_0'u$, where $\epsilon\in \{\pm 1\}$ and $c_0'\in q^{\mathbb{Z}}$ only depend on $\gamma_1, \ldots, \gamma_{l'}$ and $\lambda$. Therefore formula (\ref{fo}) takes the form
$$
\varphi (f\otimes (u,\cdot~ \tau(T_s^{-1})u))1=c_0\prod_{i=1}^{l'}
\bar{k}_i^{n_i(\lambda)}q^{(\kappa{1+s \over 1-s }P_{{\h'}}+id)\lambda^\vee}\varphi(f)q^{s^{-1}(-\kappa{1+s \over 1-s }P_{{\h'}}+id)\lambda^\vee} 1,
$$
where $c_0\in \mathcal{B}^*$ is an invertible element of $\mathcal{B}$ which only depends on $\lambda$, $\gamma_1, \ldots \gamma_{l'}$, and $n_i(\lambda)= \lambda^\vee(\gamma_i)$, $i=1, \ldots, \widetilde{l}$, $n_i(\lambda)= \lambda^\vee(s^1\gamma_i)$, $i=\widetilde{l}+1, \ldots, l'$.

The last formula can also be rewritten as follows 
$$
\varphi (f\otimes (u,\cdot~ \tau(T_s^{-1})u))1=c_0\prod_{i=1}^{l'}
\bar{k}_i^{n_i(\lambda)}q^{(\kappa{1+s \over 1-s }P_{{\h'}}+id)\lambda^\vee}\varphi(f)q^{s^{-1}(-\kappa{1+s \over 1-s }P_{{\h'}}+id)\lambda^\vee}1=
$$
$$
=c_0\prod_{i=1}^{l'}
\bar{k}_i^{n_i(\lambda)}q^{(s^{-1}+id)(id-\kappa P_{\h'})\lambda^\vee}{\rm Ad}_s(q^{(-\kappa{1+s \over 1-s }s^{-1}P_{{\h'}}+s^{-1})\lambda^\vee})(\varphi(f))1=
$$
$$
=c_0\prod_{i=1}^{l'}
\bar{k}_i^{n_i(\lambda)}q^{(s^{-1}+id)(id-\kappa P_{\h'})\lambda^\vee}\varphi({\rm Ad}_s^0(q^{(-\kappa{1+s \over 1-s }s^{-1}P_{{\h'}}+s^{-1})\lambda^\vee})(f))1=
$$
$$
=c_0\prod_{i=1}^{l'}
\bar{k}_i^{n_i(\lambda)}{\rm Ad}_s(q^{-(\kappa{1+s \over 1-s }P_{{\h'}}+id)\lambda^\vee})(\varphi(f))q^{(s^{-1}+id)(id-\kappa P_{\h'})\lambda^\vee}1=
$$
$$
=c_0\prod_{i=1}^{l'}
\bar{k}_i^{n_i(\lambda)}\varphi({\rm Ad}_s^0(q^{-(\kappa{1+s \over 1-s }P_{{\h'}}+id)\lambda^\vee})(f))q^{(s^{-1}+id)(id-\kappa P_{\h'})\lambda^\vee}1.
$$

In particular, in $Q_{\mathcal{B}}^{\bf k}$
$$
\varphi ((u,\cdot~ \tau(T_s^{-1})u))1=c_0\prod_{i=1}^{l'}
\bar{k}_i^{n_i(\lambda)}q^{(s^{-1}+id)(id-\kappa P_{\h'})\lambda^\vee}1.
$$

This completes the proof.

\end{proof}

Now consider the case when $\kappa=1$ and $\bar{k}_i\in \mathcal{B}$, $i=1,\ldots, l'$ are defined in (\ref{charq}). Recall that the elements $\tilde e_{\beta_1}^{r_1}\ldots \tilde e_{\beta_D}^{r_D}\bar{V}_i \tilde f_{\beta_D}^{m_D}\ldots \tilde f_{\beta_{c+1}}^{m_{c+1}}(\tilde f_{\beta_c}-c_{\beta_c})^{m_c}\ldots (\tilde f_{\beta_1}-c_{\beta_1})^{m_1}$ with $r_j,m_j,i\in \mathbb{N}$, $j=1,\ldots , D$ form a $\mathcal{B}$--basis in $\mathbb{C}_{\mathcal{B}}^s[G^*]$ and observe that when $\kappa=1$ the elements $(\tilde f_{\beta_c}-c_{\beta_c})^{m_c}\ldots (\tilde f_{\beta_1}-c_{\beta_1})^{m_1}$ with $m_j,i\in \mathbb{N}$, $j=1,\ldots , D$, and where at least one $m_j>0$ for $j<c+1$, form a $\mathcal{B}$--basis of ${\rm Ker}\chi_q^s$ by the definition of $\chi_q^s$. Therefore, if $\kappa=1$, we have $I_{\mathcal{B}}^{\bf k}=I_\mathcal{B}$, and hence $Q_{\mathcal{B}}^{\bf k}=Q_{\mathcal{B}}$, so we can apply the previous proposition to get the following statement. 

\begin{proposition}\label{kerphi}
Assume that $\kappa=1$ and $\bar{k}_i\in \mathcal{B}$ are defined in (\ref{charq}). Then the following statements are true.

(i) ${{J^{11}_{\mathcal{B}}}'}, {I^{11}_{\mathcal{B}}}\subset {\rm Ker}~\phi$.  

(ii) If $u$ is a highest weight vector in a finite rank indecomposable representation $V_\lambda$ of $U_h(\g)$ of highest weight $\lambda$ such that $(u,u)=1$ then for any $f\in \mathbb{C}_{\mathcal{B}}^s[G]$ 
\begin{equation}\label{phidef}
\phi (f\otimes (u,\cdot~ \tau(T_s^{-1})u))=c_\lambda \varphi({\rm Ad}_s^0(q^{-({1+s \over 1-s }P_{{\h'}}+id)\lambda^\vee})(f))q^{2P_{\h'^\perp}\lambda^\vee}1=
\end{equation}
$$
=c_\lambda q^{2P_{\h'^\perp}\lambda^\vee}\phi({\rm Ad}_s^0(q^{(-{1+s \over 1-s }s^{-1}P_{{\h'}}+s^{-1})\lambda^\vee})(f))\in Q_{\mathcal{B}},
$$
where $c_\lambda=c_0\prod_{i=1}^{l'} \bar{k}_i^{n_i(\lambda)}$, $c_0\in \mathcal{B}^*$ is an invertible element of $\mathcal{B}$ which only depends on $\lambda$, $\gamma_1, \ldots \gamma_{l'}$, and $n_i(\lambda)= \lambda^\vee(\gamma_i)\geq 0$, $i=1, \ldots, \widetilde{l}$, $n_i(\lambda)= \lambda^\vee(s^1\gamma_i)\geq 0$, $i=\widetilde{l}+1, \ldots, l'$.
Here the classes in the quotient $\mathbb{C}_{\mathcal{B}}^s[G^*]/I_{\mathcal{B}}=Q_{\mathcal{B}}'$ of the elements of $\mathbb{C}_{\mathcal{B}}^s[G^*]$ in the right hand side of (\ref{phidef}) belong to $Q_{\mathcal{B}}\subset \mathbb{C}_{\mathcal{B}}^s[G^*]/I_{\mathcal{B}}=Q_{\mathcal{B}}'$.

In particular,
$$
\phi ((u,\cdot~ \tau(T_s^{-1})u))=c_\lambda q^{2P_{\h'^\perp}\lambda^\vee}1\in Q_{\mathcal{B}}.
$$  
$q^{2P_{\h'^\perp}\lambda^\vee}1$ should be understood as the class of the element $q^{2P_{\h'^\perp}\lambda^\vee}\in \mathbb{C}_{\mathcal{B}}^s[G^*]$ in the quotient $\mathbb{C}_{\mathcal{B}}^s[G^*]/I_{\mathcal{B}}=Q_{\mathcal{B}}'$. This class belongs to $Q_{\mathcal{B}}$.           
\end{proposition}

The rest of this section will be devoted to the proof of Lemma \ref{mainl}. This proof is in turn split into several other lemmas. The appearance of the automorphism $\tau$ in all statements below is a bit annoying. An alternative would be to introduce the new quantum root vectors $\tau(X_\beta^\pm)$ and to state analogues of the results of Section \ref{PBWbases} for them. We prefer to keep $\tau$ in all formulas to avoid new symbols and statements which would double their number.

\begin{lemma}\label{l1}
Let $V$ be a finite rank representation of $U_h(\g)$, $u,v\in V$ weight vectors. Let $\overline{w}=s_{i_1}\ldots s_{i_D}$ be a reduced decomposition of the longest element of the Weyl group $W$. Then for any $\beta=s_{i_1}\ldots s_{i_{k-1}}\alpha_{i_k}\in \Delta_+$ and $n\in \mathbb{N}$
\begin{equation}\label{T4}
(u,\tau((X_\beta^+)^{(n)})v)=\sum_{p,p'\in \mathbb{N}}(u,K_{p,p'}\tau((X_\beta^-)^{(p)}(X_\beta^+)^{(p')}T_\beta) v),
\end{equation}
where the sum in the right hand side is finite, $X_\beta^\pm =T_{i_1}\ldots T_{i_{k-1}}X_{i_k}^\pm$, $K_{p,p'}\in \mathbb{C}[q,q^{-1}]$, and
\begin{equation}\label{tbeta}
T_\beta=T_{i_1}\ldots T_{i_{k-1}}T_{i_k}^{-1}T_{i_{k-1}}^{-1}\ldots T_{i_1}^{-1}. \index[not]{T@$T_\beta$}
\end{equation}
If $v$ has weight $\lambda$ then identity (\ref{T4}) may not be trivial only if $u$ has weight $\lambda+n\beta$. In this case the finite sum in (\ref{T4}) is taken over $p$ and $p'$ subject to the condition $p'-p-\beta^\vee(\lambda)=n$, so all terms $K_{p,p'}\tau((X_\beta^-)^{(p)}(X_\beta^+)^{(p')}T_\beta) v$ have also weight $\lambda+n\beta$.
\end{lemma}
\begin{proof}
Assume that $v$ is a weight vector of weight $\lambda$. Then the left hand side of (\ref{T4}) may be non-zero only if $u$ is of weight $\lambda+k\beta$ as different weight subspaces of $V$ are orthogonal with respect to the form $(~\cdot~,~\cdot~)$, so we can assume that $u$ has weight $\lambda+n\beta$. 

Conjugating (\ref{T3}) by $T_{i_1}\ldots T_{i_{k-1}}$ we get
\begin{equation}
{\exp}_{q_i}'(-X_\beta^+)=
{\exp}_{q_i^{-1}}'(-q_iX_\beta^-K_\beta^{-1})q_i^{\frac{H_\beta(H_\beta+1)}{2}}{\exp}_{q_i^{-1}}'(q_i^{-1}X_\beta^+)T_\beta,
\end{equation}
where $K_\beta=q^{\beta^\vee}=q_\beta^{H_\beta}$, $H_\beta= T_{i_1}\ldots T_{i_{k-1}}H_{i_k}$.

Applying the involution $\tau$ to this identity, evaluating it on the matrix element $(u,\cdot~ v)$ and using (\ref{qexp'}) we obtain
$$ 
-q_i^{-\frac{1}{2}n(n-1)}(u,\tau((X_\beta^+)^{(n)})v)=(u,\sum_{p=0}^\infty \tau(q_i^{-\frac{1}{2}p(p-1)}{(-q_iX_\beta^-K_\beta^{-1})^p \over [p]_{q_i}!}q_i^{\frac{H_\beta(H_\beta+1)}{2}}\sum_{p'=0}^\infty q_i^{-\frac{1}{2}p'(p'-1)}{(q_i^{-1}X_\beta^+)^{p'} \over [p']_{q_i}!}T_\beta) v)
$$

Observing that the elements $H_\beta$ act on weight spaces of $V$ by multiplication by integer numbers, that the elements $(X_\beta^-)^p$ and $(X_\beta^+)^{p'}$ map weight spaces of $V$ to weight spaces and for large enough $p$ and $p'$ they act on $V$ by zero endomorphisms, we obtain (\ref{T4}) for $v$ of weight $\lambda$, where the sum in the right hand side is such that all terms $K_{p,p'}(X_\beta^-)^{(p)}(X_\beta^+)^{(p')}T_\beta v$ have weight $\lambda+n\beta$, i.e. $p'-p-\beta^\vee(\lambda)=n$, and the number of these terms is finite. This completes the proof.

\end{proof}

Next, we obtain some useful relations in the Weyl group which lead to important formulas for the action of braid group elements on quantum root vectors.
Recall that according to Lemma \ref{invdec} (i) $s^1$ is the longest element in the Weyl group $W(\m_{s^1},\h_{s^1})$ \index[not]{W@$W(\m_{s^1},\h_{s^1})$} of the semisimple part $\m_{s^1}$ of a Levi subalgebra of $\g$, the Cartan subalgebra of $\m_{s^1}$ is denoted by $\h_{s^1}$. 

The system of positive roots $\Delta_+(\m_{s^1},\h_{s^1}):=\Delta_+\cap \Delta(\m_{s^1},\h_{s^1})$ \index[not]{D@$\Delta_+(\m_{s^1},\h_{s^1})$} of the root system $\Delta(\m_{s^1},\h_{s^1})=\Delta_{s^1}^{-1}\bigcup (-\Delta_{s^1}^{-1})$ is the set (we use the notation of (\ref{NO}))
\begin{equation}\label{m1}
\quad \Delta_+(\m_{s^1},\h_{s^1})=\{\gamma_1, \ldots, \gamma_2, 
\ldots, \gamma_3,\ldots, \gamma_{\widetilde{l}}, 
 -\beta_{t+1}^1, \ldots,-\beta_{t+\frac{p-\widetilde{l}}{2}}^1\}
\end{equation}
and $s^1$ acts on the elements of this set by multiplication by $-1$. According to (\ref{NO}) the number of roots in $\Delta_+(\m_{s^1},\h_{s^1})$ is equal to $p$. The roots in (\ref{m1}) are ordered as in the circular normal ordering of $\Delta$ associated to $s$. With respect to this normal ordering the set (\ref{m1}) is the disjoint union of an initial segment and a final segment in the normal ordering of $\Delta_+$ associated to $s$. Therefore $\m_{s^1}$ is in fact the semisimple part of a standard Levi subalgebra of $\g$ (relative to $\Delta_+$), for otherwise by Lemma \ref{invdec} (iii) there would be some roots preceding those from the set (\ref{m1}) in the normal ordering of $\Delta_+$. Thus by Lemma \ref{invdec} (iii) 
\begin{equation}\label{s1-}
\Delta_+(\m_{s^1},\h_{s^1})=\{\alpha\in \Delta_+:s^1\alpha \in \Delta_-\}=\{\alpha\in \Delta_+:s^1\alpha=-\alpha\}.
\end{equation}
This also implies that one can define the subalgebra $U_h(\m_{s^1})\subset U_h(\g)$ \index[not]{U@$U_h(\m_{s^1})$} generated by the elements $X_i^\pm$ and $H_i$ for $\alpha_i\in \Delta_+(\m_{s^1},\h_{s^1})$ and its restricted specialization $U_q^{res}(\m_{s^1})\subset U_q^{res}(\g)$. \index[not]{U@$U_q^{res}(\m_{s^1})$}

Let 
\begin{equation}\label{s1dec}
s^1=s_{i_1}\ldots s_{i_{p}}=s_{i_{k_1}}\ldots s_{i_{k_2}}\ldots  s_{i_{k_{\widetilde{l}}}}s_{i_{k_{\widetilde{l}}+1}}\ldots s_{i_{p}}, i_{k_1}=i_1,
\end{equation}
be the reduced decomposition of $s^1$ in $W(\m_{s^1},\h_{s^1})$ corresponding to (\ref{m1}) as described in Lemma \ref{invdec} (iii), where $s_i=s_{\alpha_i}$, $\alpha_i$ are simple roots in $\Delta_+(\m_{s^1},\h_{s^1})$, 
\begin{equation}\label{gamm}
\gamma_m=s_{i_{k_1}}\ldots s_{i_{k_2}}\ldots  s_{i_{k_m-1}}\alpha_{i_{k_m}}, ~m=1,\ldots , \widetilde{l}.
\end{equation} 

Since $s^1$ is an involution we also have the following reduced decomposition
\begin{equation}\label{s1}
s^1=s_{i_{p}}\ldots s_{i_1}=s_{i_{p}}\ldots s_{i_{k_{\widetilde{l}}+1}}s_{i_{k_{\widetilde{l}}}}\ldots s_{i_{k_2}}\ldots s_{i_{k_1}}.	
\end{equation}

Let $\gamma_1\leq \beta_q \leq \gamma_{\widetilde{l}}$, $\beta_q=s_{i_1}\ldots s_{i_{q-1}}\alpha_{i_q}$, $q=1,\ldots, k_{\widetilde{l}}$. 
Note that according to (\ref{m1}) (see also (\ref{pnum}))
\begin{equation}\label{kn}
k_{\widetilde{l}}=p-\frac{p-\widetilde{l}}{2}=\frac{p+\widetilde{l}}{2}.
\end{equation}

Since $s^1=-1$ in $W(\m_{s^1},\h_{s^1})$ 
$$
s^1=s_{i_q}\ldots s_{i_1}s^1s_{i_1}\ldots s_{i_q}=s_{i_q}\ldots s_{i_1}s_{i_1}\ldots s_{i_{p}}s_{i_1}\ldots s_{i_q}=s_{i_{q+1}}\ldots s_{i_{p}}s_{i_1}\ldots s_{i_q},
$$ 
and in the right hand side we obtain a reduced decomposition of $s^1$ as the number of simple reflections in it is equal to $p$, i.e. to the length of $s^1$. 

Now by the previous formula
$$
s^1\alpha_{i_q}=-\alpha_{i_q}=s_{i_{q+1}}\ldots s_{i_{p}}s_{i_1}\ldots s_{i_q}\alpha_{i_q}=-s_{i_{q+1}}\ldots s_{i_{p}}s_{i_1}\ldots s_{i_{q-1}}\alpha_{i_q},
$$
and hence
\begin{equation}\label{aq}
s_{i_{q+1}}\ldots s_{i_{p}}s_{i_1}\ldots s_{i_{q-1}}\alpha_{i_q}=\alpha_{i_q}.	
\end{equation}

Form the expressions $\gamma_m=s_{i_{k_1}}\ldots s_{i_{k_2}}\ldots  s_{i_{k_m-1}}\alpha_{i_{k_m}}$, $m=1,\ldots ,\widetilde{l}$ we deduce 
\begin{equation}\label{gqgq}
s_{\gamma_1}\ldots s_{\gamma_{q-1}}=
\end{equation}
$$
=s_{i_{k_1}}(s_{i_{k_1}}\ldots s_{i_{k_2-1}}s_{i_{k_2}}s_{i_{k_2-1}}^{-1}\ldots s_{i_{k_1}}^{-1})(s_{i_{k_1}}\ldots s_{i_{k_3-1}}s_{i_{k_3}}s_{i_{k_3-1}}^{-1}\ldots s_{i_{k_1}}^{-1})\ldots (s_{i_{k_1}}\ldots s_{i_{k_2}}\ldots  s_{i_{k_{q-1}-1}}s_{i_{k_{q-1}}}s_{i_{k_{q-1}-1}}^{-1}\ldots s_{i_{k_1}}^{-1})=
$$
$$
=s_{i_{k_1+1}}\ldots s_{i_{k_2-1}}s_{i_{k_2+1}}\ldots  s_{i_{k_{q-1}-1}}s_{i_{k_{q-1}}}s_{i_{k_{q-1}-1}}\ldots s_{i_{k_1}}.
$$

In particular, in the right hand sides of the following formulas we have reduced decompositions
\begin{equation}\label{s1red}
s^1=s_{\gamma_1}\ldots s_{\gamma_{\widetilde{l}}}=s_{i_{k_1+1}}\ldots s_{i_{k_2-1}}s_{i_{k_2+1}}\ldots  s_{i_{k_{\widetilde{l}}-1}}s_{i_{k_{\widetilde{l}}}}s_{i_{k_{\widetilde{l}}-1}}\ldots s_{i_{k_1}}	
\end{equation}
and
\begin{equation}\label{s2}
s^1=s_{i_{k_{\widetilde{l}}}}\ldots s_{i_{1}}s^1s_{i_{1}}\ldots s_{i_{k_{\widetilde{l}}}}=s_{i_{k_{\widetilde{l}}}}s_{i_{k_{\widetilde{l}}-1}}\ldots s_{i_{k_1}}s_{i_{k_1+1}}\ldots s_{i_{k_2-1}}s_{i_{k_2+1}}\ldots  s_{i_{k_{\widetilde{l}}-1}}
\end{equation}
as by (\ref{kn}) the number of simple reflections in them is equal to $k_{\widetilde{l}}+k_{\widetilde{l}}-\widetilde{l}=2k_{\widetilde{l}}-\widetilde{l}=p+\widetilde{l}-\widetilde{l}=p$, i.e. to the length of $s^1$. 

Multiplying (\ref{s1red}) and (\ref{s1}) by $(s_{i_{k_{\widetilde{l}}}}s_{i_{k_{\widetilde{l}}-1}}\ldots s_{i_{k_1}})^{-1}$ on the right we obtain the following identity for reduced decompositions
\begin{equation}\label{lsd}
s_{i_{k_1+1}}\ldots s_{i_{k_2-1}}s_{i_{k_2+1}}\ldots  s_{i_{k_{\widetilde{l}}-1}}=s_{i_p}\ldots s_{i_{k_{\widetilde{l}}+1}}.	
\end{equation}

As the roots $\gamma_m$, $m=1,\ldots ,\widetilde{l}$ are mutually orthogonal we deduce using (\ref{gqgq}) and (\ref{gamm})
\begin{align}\label{sss1}
s_{\gamma_1}\ldots s_{\gamma_{q-1}}\gamma_q=s_{i_{k_1+1}}\ldots s_{i_{k_2-1}}s_{i_{k_2+1}}\ldots  s_{i_{k_{q-1}-1}}s_{i_{k_{q-1}}}s_{i_{k_{q-1}-1}}\ldots s_{i_{k_1}}s_{i_{k_1}}\ldots s_{i_{k_q-1}}\alpha_{i_{k_q}}= \nonumber \\
=s_{i_{k_1+1}}\ldots s_{i_{k_2-1}}s_{i_{k_2+1}}\ldots  s_{i_{k_{q-1}-1}}s_{i_{k_{q-1}+1}}\ldots s_{i_{k_q-1}}\alpha_{i_{k_q}}=\gamma_q=s_{i_{k_1}}\ldots s_{i_{k_q-1}}\alpha_{i_{k_q}}.
\end{align}
Therefore 
$$
s_{i_{k_q-1}}\ldots s_{i_{k_1}}s_{i_{k_1+1}}\ldots s_{i_{k_2-1}}s_{i_{k_2+1}}\ldots  s_{i_{k_{q-1}-1}}s_{i_{k_{q-1}+1}}\ldots s_{i_{k_q-1}}\alpha_{i_{k_q}}=\alpha_{i_{k_q}},
$$
where
$$
s_{i_{k_q-1}}\ldots s_{i_{k_1}}s_{i_{k_1+1}}\ldots s_{i_{k_2-1}}s_{i_{k_2+1}}\ldots  s_{i_{k_{q-1}-1}}s_{i_{k_{q-1}+1}}\ldots s_{i_{k_q-1}}
$$
is a reduced decomposition since it is a part of reduced decomposition (\ref{s2}).

The last two properties and (\ref{treda}) imply
$$
T_{i_{k_q-1}}\ldots T_{i_{k_1}}T_{i_{k_1}+1}\ldots T_{i_{k_2-1}}T_{i_{k_2+1}}\ldots  T_{i_{k_{q-1}-1}}T_{i_{k_{q-1}+1}}\ldots T_{i_{k_q-1}}X_{i_{k_q}}^\pm=X_{i_{k_q}}^\pm,
$$
or
\begin{equation}\label{redX}
T_{i_{k_1}+1}\ldots T_{i_{k_2-1}}T_{i_{k_2+1}}\ldots  T_{i_{k_{q-1}-1}}T_{i_{k_{q-1}+1}}\ldots T_{i_{k_q-1}}X_{i_{k_q}}^\pm=T_{i_{k_1}}^{-1}\ldots T_{i_{k_q-1}}^{-1}X_{i_{k_q}}^\pm.
\end{equation}
Using the definition of $T_{\gamma_i}$ (see (\ref{tbeta})) we obtain
\begin{equation}\label{Tq-1}
T_{\gamma_1}\ldots T_{\gamma_{q-1}}=
\end{equation}
$$
T_{i_{k_1}}^{-1}(T_{i_{k_1}}\ldots T_{i_{k_2-1}}T_{i_{k_2}}^{-1}T_{i_{k_2-1}}^{-1}\ldots T_{i_{k_1}}^{-1})(T_{i_{k_1}}\ldots T_{i_{k_3-1}}T_{i_{k_3}}^{-1}T_{i_{k_3-1}}^{-1}\ldots T_{i_{k_1}}^{-1})\ldots (T_{i_{k_1}}\ldots T_{i_{k_{q-1}-1}}T_{i_{k_{q-1}}}^{-1}T_{i_{k_{q-1}-1}}^{-1}\ldots T_{i_{k_1}}^{-1})=
$$
$$
=T_{i_{k_1}+1}\ldots T_{i_{k_2-1}}T_{i_{k_2+1}}\ldots  T_{i_{k_{q-1}-1}}T_{i_{k_{q-1}}}^{-1}T_{i_{k_{q-1}-1}}^{-1}\ldots T_{i_{k_1}}^{-1}.
$$
Recalling also that $X_{\gamma_q}^\pm=T_{i_{k_1}}\ldots T_{i_{k_q-1}}X_{i_{k_q}}^\pm$ and applying formula (\ref{redX}) we arrive at
\begin{align}
T_{\gamma_1}\ldots T_{\gamma_{q-1}}X_{\gamma_q}^\pm=T_{i_{k_1}+1}\ldots T_{i_{k_2-1}}T_{i_{k_2+1}}\ldots  T_{i_{k_{q-1}-1}}T_{i_{k_{q-1}}}^{-1}T_{i_{k_{q-1}-1}}^{-1}\ldots T_{i_{k_1}}^{-1}T_{i_{k_1}}\ldots T_{i_{k_q-1}}X_{i_{k_q}}^\pm= \nonumber
\\
=T_{i_{k_1}+1}\ldots T_{i_{k_2-1}}T_{i_{k_2+1}}\ldots  T_{i_{k_{q-1}-1}}T_{i_{k_{q-1}+1}}\ldots T_{i_{k_q-1}}X_{i_{k_q}}^\pm=T_{i_{k_1}}^{-1}\ldots T_{i_{k_q-1}}^{-1}X_{i_{k_q}}^\pm=\overline{X}_{\gamma_q}^\pm. \label{gpm}
\end{align}

\begin{lemma}\label{l2}
Let $V$ be a finite rank representation of $U_h(\g)$, $u,v\in V$ weight vectors.
Then for any $m_1,\ldots, m_{\widetilde{l}}\in \mathbb{N}$ the following statements are true.

(i) 
\begin{equation}\label{T5}
(u,\tau((X_{\gamma_1}^+)^{(m_1)}\ldots (X_{\gamma_{\widetilde{l}}}^+)^{(m_{\widetilde{l}})})v)=(u,\tau(K(m_1,\ldots, m_{\widetilde{l}})T_{\gamma_1}\ldots T_{\gamma_{\widetilde{l}}})v),
\end{equation}
where $K(m_1,\ldots, m_{\widetilde{l}})\in U_q^{res}(\m_{s^1}^-)\overline{U}_q^{res}([\gamma_1,\gamma_{\widetilde{l}}])$,  
$U_q^{res}(\m_{s^1}^-)\subset U_q^{res}(\g)$ \index[not]{U@$U_q^{res}(\m_{s^1}^-)$} is the subalgebra generated by $({X}_{\beta}^-)^{(k)}$ for simple roots $\beta \in \Delta_+(\m_{s^1},\h_{s^1})$ and $k\geq 0$. Moreover, $K(m_1,\ldots, m_{\widetilde{l}})$ belongs to a weight subspace of $U_q^{res}(\g)$.

(ii) If $u\in V$ is a highest weight vector then
\begin{equation}\label{T50}
(u,\tau((X_{\gamma_1}^+)^{(m_1)}\ldots (X_{\gamma_{\widetilde{l}}}^+)^{(m_{\widetilde{l}})})v)=(u,\tau(K'(m_1,\ldots, m_{\widetilde{l}})T_{\gamma_1}\ldots T_{\gamma_{\widetilde{l}}})v),
\end{equation}
where $K'(m_1,\ldots, m_{\widetilde{l}})\in \overline{U}_q^{res}([\gamma_1,\gamma_{\widetilde{l}}])$. Moreover, $K'(m_1,\ldots, m_{\widetilde{l}})$ belongs to a weight subspace of $U_q^{res}(\g)$.
\end{lemma}

\begin{proof}
By Lemma \ref{l1}
$$
(u,\tau((X_{\gamma_1}^+)^{(m_1)}\ldots (X_{\gamma_{\widetilde{l}}}^+)^{(m_{\widetilde{l}})})v)= 
$$
$$
=\sum_{\tiny \begin{array}{l} c_1,\ldots ,c_{\widetilde{l}}, \\ c_1',\ldots ,c_{\widetilde{l}}'\in \mathbb{N} \end{array}} K^{c_1,\ldots ,c_{\widetilde{l}}}_{c_1',\ldots ,c_{\widetilde{l}}'}(u,\tau((X_{\gamma_1}^-)^{(c_1)}(X_{\gamma_1}^+)^{(c_1')}T_{\gamma_1}(X_{\gamma_2}^-)^{(c_2)}(X_{\gamma_2}^+)^{(c_2')}T_{\gamma_2}\ldots (X_{\gamma_{\widetilde{l}}}^-)^{(c_{\widetilde{l}})}(X_{\gamma_{\widetilde{l}}}^+)^{(c_{\widetilde{l}}')}T_{\gamma_{\widetilde{l}}})v),
$$
where the sum in the right hand side is finite and $K^{c_1,\ldots ,c_{\widetilde{l}}}_{c_1',\ldots ,c_{\widetilde{l}}'}\in \mathbb{C}[q,q^{-1}]$ .

Using (\ref{gpm}) and observing that $X_{\gamma_1}^\pm=\overline{X}_{\gamma_1}^\pm$, as $\gamma_1$ is the first simple root in the normal ordering of $\Delta_+$ associated to $s$, we obtain
\begin{equation}\label{st*}
(u,\tau((X_{\gamma_1}^+)^{(k_1)}\ldots (X_{\gamma_{\widetilde{l}}}^+)^{(k_{\widetilde{l}})})v)= 
\end{equation}
$$
=\sum_{\tiny \begin{array}{l} c_1,\ldots ,c_{\widetilde{l}}, \\ c_1',\ldots ,c_{\widetilde{l}}' \in \mathbb{N} \end{array}} K^{c_1,\ldots ,c_{\widetilde{l}}}_{c_1',\ldots ,c_{\widetilde{l}}'}(u,\tau((\overline{X}_{\gamma_1}^-)^{(c_1)}(\overline{X}_{\gamma_1}^+)^{(c_1')}(\overline{X}_{\gamma_2}^-)^{(c_2)}(\overline{X}_{\gamma_2}^+)^{(c_2')}\ldots (\overline{X}_{\gamma_{\widetilde{l}}}^-)^{(c_{\widetilde{l}})}(\overline{X}_{\gamma_{\widetilde{l}}}^+)^{(c_{\widetilde{l}}')}T_{\gamma_1}\ldots T_{\gamma_{\widetilde{l}}})v).
$$

Now to justify (\ref{T5}) we show that all elements 
\begin{equation}\label{rootprodX}
(\overline{X}_{\gamma_1}^-)^{(c_1)}(\overline{X}_{\gamma_1}^+)^{(c_1')}(\overline{X}_{\gamma_2}^-)^{(c_2)}(\overline{X}_{\gamma_2}^+)^{(c_2')}\ldots (\overline{X}_{\gamma_{\widetilde{l}}}^-)^{(c_{\widetilde{l}})}(\overline{X}_{\gamma_{\widetilde{l}}}^+)^{(c_{\widetilde{l}}')}
\end{equation}
belong to $U_q^{res}(\m_{s^1}^-)\overline{U}_q^{res}([\gamma_1,\gamma_{\widetilde{l}}])U_q^{res}(H)$.

Indeed, denote by $\beta_p$ the last root in normal ordering (\ref{m1}). Then $(\overline{X}_{\gamma_{\widetilde{l}-1}}^+)^{(c_{\widetilde{l}-1}')}(\overline{X}_{\gamma_{\widetilde{l}}}^-)^{(c_{\widetilde{l}})}\in \overline{U}_{U_q^{res}(H_{s^1})}^{res}([-\gamma_{\widetilde{l}},\gamma_{\widetilde{l}-1}])$, where here and in this proof below we consider only minimal segments $[\alpha,\beta]$ of the circular ordering of $\Delta(\m_{s^1},\h_{s^1})$ corresponding normal ordering (\ref{m1}) of $\Delta_+(\m_{s^1},\h_{s^1})$, and the corresponding subalgebras of $U_q^{res}(\m_{s^1})$,  so,  in particular, $[-\gamma_{\widetilde{l}},\gamma_{\widetilde{l}-1}]\subset \Delta(\m_{s^1},\h_{s^1})$, and $\overline{U}_{U_q^{res}(H_{s^1})}^{res}([-\gamma_{\widetilde{l}},\gamma_{\widetilde{l}-1}])\subset U_q^{res}(\m_{s^1})$, \index[not]{U@$\overline{U}_{U_q^{res}(H_{s^1})}^{res}([-\gamma_{\widetilde{l}},\gamma_{\widetilde{l}-1}])$} where $U_q^{res}(H_{s^1})\subset U_q^{res}(\m_{s^1})$ \index[not]{U@$U_q^{res}(H_{s^1})$} is the subalgebra generated by elements (\ref{Kic}) of $U_q^{res}(\g)$ and  by the $K_i^{\pm 1}$ corresponding to the simple roots $\alpha_i\in \Delta_+(\m_{s^1},\h_{s^1})$. 

Then by Corollary \ref{segmq} (i) we obtain
\begin{equation}\label{st1*}
(\overline{X}_{\gamma_{\widetilde{l}-1}}^+)^{(c_{\widetilde{l}-1}')}(\overline{X}_{\gamma_{\widetilde{l}}}^-)^{(c_{\widetilde{l}})}(\overline{X}_{\gamma_{\widetilde{l}}}^+)^{(c_{\widetilde{l}}')}\in\overline{U}_{U_q^{res}(H_{s^1})}^{res}([-\gamma_{\widetilde{l}},\gamma_{\widetilde{l}-1}])(\overline{X}_{\gamma_{\widetilde{l}}}^+)^{(c_{\widetilde{l}}')}=
\end{equation}
$$
=\overline{U}_{U_q^{res}(H_{s^1})}^{res}([-\gamma_{\widetilde{l}},-\beta_p])\overline{U}_{U_q^{res}(H_{s^1})}^{res}([\gamma_1,\gamma_{\widetilde{l}-1}])(\overline{X}_{\gamma_{\widetilde{l}}}^+)^{(c_{\widetilde{l}}')}\subset \overline{U}_{U_q^{res}(H_{s^1})}^{res}([-\gamma_{\widetilde{l}},-\beta_p])\overline{U}_{U_q^{res}(H_{s^1})}^{res}([\gamma_1,\gamma_{\widetilde{l}}]).
$$

Next, $(\overline{X}_{\gamma_{\widetilde{l}-2}}^+)^{(c_{\widetilde{l}-2}')}(\overline{X}_{\gamma_{\widetilde{l}-1}}^-)^{(c_{\widetilde{l}-1})}\in \overline{U}_{U_q^{res}(H_{s^1})}^{res}([-\gamma_{\widetilde{l}-1},\gamma_{\widetilde{l}-2}])$, and by (\ref{st1*}) and by Corollary \ref{segmq} (i) one has
$$
(\overline{X}_{\gamma_{\widetilde{l}-2}}^+)^{(c_{\widetilde{l}-2}')}(\overline{X}_{\gamma_{\widetilde{l}-1}}^-)^{(c_{\widetilde{l}-1})}(\overline{X}_{\gamma_{\widetilde{l}-1}}^+)^{(c_{\widetilde{l}-1}')}(\overline{X}_{\gamma_{\widetilde{l}}}^-)^{(c_{\widetilde{l}})}(\overline{X}_{\gamma_{\widetilde{l}}}^+)^{(c_{\widetilde{l}}')}\in 
$$
$$
\in\overline{U}_{U_q^{res}(H_{s^1})}^{res}([-\gamma_{\widetilde{l}-1},\gamma_{\widetilde{l}-2}])\overline{U}_{U_q^{res}(H_{s^1})}^{res}([-\gamma_{\widetilde{l}},-\beta_p])\overline{U}_{U_q^{res}(H_{s^1})}^{res}([\gamma_1,\gamma_{\widetilde{l}}])=
$$
$$
=\overline{U}_{U_q^{res}(H_{s^1})}^{res}([-\gamma_{\widetilde{l}-1},-\beta_p])\overline{U}_{U_q^{res}(H_{s^1})}^{res}([\gamma_1,\gamma_{\widetilde{l}-2}])\overline{U}_{U_q^{res}(H_{s^1})}^{res}([-\gamma_{\widetilde{l}},-\beta_p])\overline{U}_{U_q^{res}(H_{s^1})}^{res}([\gamma_1,\gamma_{\widetilde{l}}])=
$$
$$
=\overline{U}_{U_q^{res}(H_{s^1})}^{res}([-\gamma_{\widetilde{l}-1},-\beta_p])\overline{U}_{U_q^{res}(H_{s^1})}^{res}([-\gamma_{\widetilde{l}},\gamma_{\widetilde{l}-2}])\overline{U}_{U_q^{res}(H_{s^1})}^{res}([\gamma_1,\gamma_{\widetilde{l}}])=
$$
$$
=\overline{U}_{U_q^{res}(H_{s^1})}^{res}([-\gamma_{\widetilde{l}-1},-\beta_p])\overline{U}_{U_q^{res}(H_{s^1})}^{res}([-\gamma_{\widetilde{l}},-\beta_p])\overline{U}_{U_q^{res}(H_{s^1})}^{res}([\gamma_1,\gamma_{\widetilde{l}-2}])\overline{U}_{U_q^{res}(H_{s^1})}^{res}([\gamma_1,\gamma_{\widetilde{l}}])\subset 
$$
$$
\subset \overline{U}_{U_q^{res}(H_{s^1})}^{res}([-\gamma_{\widetilde{l}-1},-\beta_p])\overline{U}_{U_q^{res}(H_{s^1})}^{res}([\gamma_1,\gamma_{\widetilde{l}}]).
$$

One can proceed by induction in a similar way to get
\begin{equation}\label{st2*}
(\overline{X}_{\gamma_1}^-)^{(c_1)}(\overline{X}_{\gamma_1}^+)^{(c_1')}(\overline{X}_{\gamma_2}^-)^{(c_2)}(\overline{X}_{\gamma_2}^+)^{(c_2')}\ldots (\overline{X}_{\gamma_{\widetilde{l}}}^-)^{(c_{\widetilde{l}})}(\overline{X}_{\gamma_{\widetilde{l}}}^+)^{(c_{\widetilde{l}}')}\in \overline{U}_{U_q^{res}(H_{s^1})}^{res}([-\gamma_1,-\beta_p])\overline{U}_{U_q^{res}(H_{s^1})}^{res}([\gamma_1,\gamma_{\widetilde{l}}]).
\end{equation}

Next, Corollary \ref{segmq} (ii) implies
$$
\overline{U}_{U_q^{res}(H_{s^1})}^{res}([-\gamma_1,-\beta_p])\overline{U}_{U_q^{res}(H_{s^1})}^{res}([\gamma_1,\gamma_{\widetilde{l}}])=\overline{U}_q^{res}([-\gamma_1,-\beta_p])U_q^{res}(H_{s^1})\overline{U}_q^{res}([\gamma_1,\gamma_{\widetilde{l}}])U_q^{res}(H_{s^1})= 
$$
$$
=\overline{U}_q^{res}([-\gamma_1,-\beta_p])\overline{U}_q^{res}([\gamma_1,\gamma_{\widetilde{l}}])U_q^{res}(H_{s^1})=U_q^{res}(\m_{s^1}^-)\overline{U}_q^{res}([\gamma_1,\gamma_{\widetilde{l}}])U_q^{res}(H_{s^1}),
$$
where at the last step we also used the fact that by the definition $\overline{U}_q^{res}([-\gamma_1,-\beta_p])=U_q^{res}(\m_{s^1}^-)$.

Now (\ref{st2*}) takes the form
\begin{equation}\label{st3*}
(\overline{X}_{\gamma_1}^-)^{(c_1)}(\overline{X}_{\gamma_1}^+)^{(c_1')}(\overline{X}_{\gamma_2}^-)^{(c_2)}(\overline{X}_{\gamma_2}^+)^{(c_2')}\ldots (\overline{X}_{\gamma_{\widetilde{l}}}^-)^{(c_{\widetilde{l}})}(\overline{X}_{\gamma_{\widetilde{l}}}^+)^{(c_{\widetilde{l}}')}\in U_q^{res}(\m_{s^1}^-)\overline{U}_q^{res}([\gamma_1,\gamma_{\widetilde{l}}])U_q^{res}(H_{s^1}).
\end{equation}

Thus, recalling that $v$ is a weight vector we deduce (\ref{T5}) from (\ref{st*}) and (\ref{st3*}) using Lemma \ref{wHact}. 

Finally note that $u$ and $T_{\gamma_1}\ldots T_{\gamma_{\widetilde{l}}}v$ are weight vectors, and different weight subspaces of $V$ are orthogonal with respect to the contravariant form $(~\cdot~,~\cdot~)$. Therefore we can assume also that $K(m_1,\ldots, m_{\widetilde{l}})$ in (\ref{T5}) belongs to a weight subspace of $U_q^{res}(\g)$, so that the weight of $u$ is equal to that of $K(m_1,\ldots, m_{\widetilde{l}})T_{\gamma_1}\ldots T_{\gamma_{\widetilde{l}}}v$.

(\ref{T50}) follows from (\ref{T5}) by Lemma \ref{hwv} (i) because the only elements of $U_q^{res}(\m_{s^1}^-)$, whose weights are not strictly negative, belong to $\mathbb{C}[q,q^{-1}]$ and $\tau$ preserves weight subspaces.
This completes the proof.

\end{proof}

\begin{lemma}\label{l3}
Let $V$ be a finite rank representation of $U_h(\g)$, $u,v\in V$. Suppose that $u$ is a highest weight vector, and $v$ is a weight vector.
Then for any $m_1,\ldots, m_{\widetilde{l}}\in \mathbb{N}$
\begin{equation}\label{l3st}
(u,\tau((X_{\gamma_1}^+)^{(m_1)}\ldots (X_{\gamma_{\widetilde{l}}}^+)^{(m_{\widetilde{l}})})v)
=(u,\tau(T_{\gamma_1}\ldots T_{\gamma_{\widetilde{l}}}U(m_1,\ldots, m_{\widetilde{l}})) v), 
\end{equation}
where $U(m_1,\ldots, m_{\widetilde{l}})\in U_q^{res}([-\gamma_1,-\gamma_{\widetilde{l}}])$ belongs to a weight subspace of $U_q^{res}(\g)$.
\end{lemma}

\begin{proof}

Since $u$ is a highest weight vector, Lemma \ref{l2} (ii) implies
\begin{equation}\label{mmm*}
(u,\tau((X_{\gamma_1}^+)^{(m_1)}\ldots (X_{\gamma_{\widetilde{l}}}^+)^{(m_{\widetilde{l}})})v)=(u,\tau(K'(m_1,\ldots, m_{\widetilde{l}})T_{\gamma_1}\ldots T_{\gamma_{\widetilde{l}}})v),
\end{equation}
where $K'(m_1,\ldots, m_{\widetilde{l}})\in \overline{U}_q^{res}([\gamma_1,\gamma_{\widetilde{l}}])$.

Denote $T^1=T_{\gamma_1}\ldots T_{\gamma_{\widetilde{l}}}$. \index[not]{T@$T^1$} We find the action of $(T^1)^{-1}$ on the generators of the algebra $\overline{U}_q^{res}([\gamma_1,\gamma_{\widetilde{l}}])$. 

Consider reduced decomposition (\ref{s1dec}) of $s^1$,
\begin{equation}\label{s1dec'}
s^1=s_{i_1}\ldots s_{i_{p}}=s_{i_{k_1}}\ldots s_{i_{k_2}}\ldots  s_{i_{k_{\widetilde{l}}}}s_{i_{k_{\widetilde{l}}+1}}\ldots s_{i_{p}}, i_{k_1}=i_1,
\end{equation}
and the roots $\beta_q=s_{i_1}\ldots s_{i_{q-1}}\alpha_{i_q}$, $q=1,\ldots, k_{\widetilde{l}}$ forming the segment $[\gamma_1,\gamma_{\widetilde{l}}]$.

From (\ref{Tq-1}) with $q-1=\widetilde{l}$ and (\ref{lsd}) we obtain
\begin{equation}\label{T1red}
T^1=T_{\gamma_1}\ldots T_{\gamma_{\widetilde{l}}}=T_{i_{k_1}+1}\ldots T_{i_{k_2-1}}T_{i_{k_2+1}}\ldots  T_{i_{k_{\widetilde{l}}-1}}T_{i_{k_{\widetilde{l}}}}^{-1}T_{i_{k_{\widetilde{l}}-1}}^{-1}\ldots T_{i_{1}}^{-1}=
\end{equation}
$$
=T_{i_{p}}\ldots T_{i_{k_{\widetilde{l}}+1}}T_{i_{k_{\widetilde{l}}}}^{-1}T_{i_{k_{\widetilde{l}}-1}}^{-1}\ldots T_{i_{k_1}}^{-1}.
$$
Therefore for the generators $(\overline{X}_{\beta_q}^+)^{(k)}=T_{i_1}^{-1}\ldots T_{i_{q-1}}^{-1}(X_{i_q}^+)^{(k)}$, $q=1,\ldots, k_{\widetilde{l}}$, $k\in \mathbb{N}$, of the algebra $\overline{U}_q^{res}([\gamma_1,\gamma_{\widetilde{l}}])$ we obtain
\begin{equation}\label{oX}
(T^1)^{-1}(\overline{X}_{\beta_q}^+)^{(k)}=T_{i_{1}}\ldots T_{i_{k_{\widetilde{l}}}}T_{i_{k_{\widetilde{l}}+1}}^{-1} \ldots T_{i_{p}}^{-1}T_{i_1}^{-1}\ldots T_{i_{q-1}}^{-1}(X_{i_q}^+)^{(k)}=	
\end{equation}
$$
=T_{w_1}T_{i_{k_{\widetilde{l}}+1}}^{-1} \ldots T_{i_{p}}^{-1}T_{i_1}^{-1}\ldots T_{i_{q-1}}^{-1}(X_{i_q}^+)^{(k)},
$$
where $T_{\bar{w}_1}=T_{i_{1}}\ldots T_{i_{k_{\widetilde{l}}}}$ \index[not]{w@$\bar{w}_1$} for the reduced decomposition $\bar{w}_1=s_{i_{1}}\ldots s_{i_{k_{\widetilde{l}}}}$ which is an initial part of reduced decomposition (\ref{s1dec'}).

By (\ref{aq}), (\ref{treda})
$$
T_{i_{q+1}}^{-1}\ldots T_{i_{p}}^{-1}T_{i_1}^{-1}\ldots T_{i_{q-1}}^{-1}(X_{i_q}^+)^{(k)}=(X_{i_q}^+)^{(k)},
$$
and hence
$$
T_{i_{k_{\widetilde{l}}+1}}^{-1}\ldots T_{i_{p}}^{-1}T_{i_1}^{-1}\ldots T_{i_{q-1}}^{-1}(X_{i_q}^+)^{(k)}=T_{i_{k_{\widetilde{l}}}}\ldots T_{i_{q+1}}(X_{i_q}^+)^{(k)}.
$$

Substituting this into (\ref{oX}) we infer
\begin{equation}\label{t11}
(T^1)^{-1}(\overline{X}_{\beta_q}^+)^{(k)}=T_{\bar{w}_1}T_{i_{k_{\widetilde{l}}}}\ldots T_{i_{q+1}}(X_{i_q}^+)^{(k)}=T_{\bar{w}_1}(X_{\delta_q}^+)^{(k)},
\end{equation}
where for $q=1,\ldots, k_{\widetilde{l}}$ we denote $\delta_q=s_{i_{k_{\widetilde{l}}}}\ldots s_{i_{q+1}}\alpha_{i_q}$, \index[not]{d@$\delta_q$} and  $X_{\delta_q}^+= T_{i_{k_{\widetilde{l}}}}\ldots T_{i_{q+1}}X_{i_q}^+$.

Let $U_q^{res}(\m_{s^1}^+)$ \index[not]{U@$U_q^{res}(\m_{s^1}^+)$} be the subalgebra of $U_q^{res}(\m_{s^1})$ generated by $({X}_{\beta}^+)^{(k)}$ for simple roots $\beta \in \Delta_+(\m_{s^1},\h_{s^1})$ and $k\geq 0$. We show that for $q=1,\ldots, k_{\widetilde{l}}$ and $k\geq 0$ one has $({X}_{\delta_q}^+)^{(k)}\in U_q^{res}(\m_{s^1}^+)$.

Indeed, observe that $s_{i_{k_{\widetilde{l}}}}\ldots s_{i_{1}}s_{i_p}\ldots s_{i_{k_{\widetilde{l}}+1}}$ is the reduced decomposition of $s^1$, with respect to the system of simple roots in $\Delta_+(\m_{s^1},\h_{s^1})$, obtained according to (\ref{s2}) in the following way 
\begin{equation}\label{s2'}
s^1=s_{i_{k_{\widetilde{l}}}}\ldots s_{i_{1}}s^1s_{i_{1}}\ldots s_{i_{k_{\widetilde{l}}}}=s_{i_{k_{\widetilde{l}}}}\ldots s_{i_{1}}s_{i_p}\ldots s_{i_1}s_{i_{1}}\ldots s_{i_{k_{\widetilde{l}}}}=s_{i_{k_{\widetilde{l}}}}\ldots s_{i_{1}}s_{i_p}\ldots s_{i_{k_{\widetilde{l}}+1}}.
\end{equation}
$s_{i_{k_{\widetilde{l}}}}\ldots s_{i_1}$ is an initial part of this reduced decomposition.

Now by Lemma \ref{segmPBW} (i) we deduce that for $q=1,\ldots, k_{\widetilde{l}}$ one has  $X_{\delta_q}^+=T_{i_{k_{\widetilde{l}}}}\ldots T_{i_{q+1}}X_{i_q}^+\in U_q^{res}(\m_{s^1}^+)$, and also obviously 
\begin{equation}\label{Xd}
(X_{\delta_q}^+)^{(k)}\in U_q^{res}(\m_{s^1}^+)
\end{equation} 
for any $k\geq 0$. Since for $q=1,\ldots, k_{\widetilde{l}}$ and $k\geq 0$ the elements $(\overline{X}_{\beta_q}^+)^{(k)}$ generate $\overline{U}_q^{res}([\gamma_1,\gamma_{\widetilde{l}}])$, we deduce from (\ref{Xd}) and (\ref{t11}) that 
$$
(T^1)^{-1}(\overline{U}_q^{res}([\gamma_1,\gamma_{\widetilde{l}}]))\subset T_{\bar{w}_1}(U_q^{res}(\m_{s^1}^+)).
$$

By this inclusion $(T^1)^{-1}(K'(m_1,\ldots, m_{\widetilde{l}}))\in T_{\bar{w}_1}(U_q^{res}(\m_{s^1}^+))$, and (\ref{mmm*}) takes the form

\begin{equation}\label{G1}
(u,\tau((X_{\gamma_1}^+)^{(m_1)}\ldots (X_{\gamma_{\widetilde{l}}}^+)^{(m_{\widetilde{l}})})v)=(u,\tau(K'(m_1,\ldots, m_{\widetilde{l}})T_{\gamma_1}\ldots T_{\gamma_{\widetilde{l}}})v)=
\end{equation}
$$
=(u,\tau(T_{\gamma_1}\ldots T_{\gamma_{\widetilde{l}}} (T^1)^{-1}(K'(m_1,\ldots, m_{\widetilde{l}})))v)=(u,\tau(T_{\gamma_1}\ldots T_{\gamma_{\widetilde{l}}} K''(m_1,\ldots, m_{\widetilde{l}}))v),
$$
where $K''(m_1,\ldots, m_{\widetilde{l}})=(T^1)^{-1}(K'(m_1,\ldots, m_{\widetilde{l}}))\in T_{\bar{w}_1}(U_q^{res}(\m_{s^1}^+))$.

Now we factorize $U_q^{res}(\m_{s^1}^+)$ in an appropriate way to bring (\ref{G1}) to form (\ref{l3st}). Observe that using reduced decomposition (\ref{s2'}) one can define
the corresponding normal ordering of $\Delta_+(\m_{s^1},\h_{s^1})$,
$$
\phi_{k_{\widetilde{l}}}=\alpha_{i_{k_{\widetilde{l}}}},\ldots, \phi_1=s_{i_{k_{\widetilde{l}}}}\ldots s_{i_{2}}\alpha_{i_{1}},
\phi_p=s_{i_{k_{\widetilde{l}}}}\ldots s_{i_{1}}\alpha_{i_p}, \ldots,
\phi_{k_{\widetilde{l}}+1}=s_{i_{k_{\widetilde{l}}}}\ldots s_{i_{1}}s_{i_p}\ldots s_{i_{k_{\widetilde{l}}+2}}\alpha_{i_{k_{\widetilde{l}}+1}},
$$
\index[not]{f@$\phi_i$}
and the elements
$$
\begin{array}{l}
Y_{\phi_{k_{\widetilde{l}}}}=X_{i_{k_{\widetilde{l}}}}^+,\ldots,
\\
\\
Y_{\phi_{1}}=T_{i_{k_{\widetilde{l}}}}^{-1}\ldots T_{i_{2}}^{-1}X_{i_{1}}^+,
\\
\\
Y_{\phi_{p}}=T_{i_{k_{\widetilde{l}}}}^{-1}\ldots T_{i_{1}}^{-1}X_{i_p}^+, \ldots,
\\
\\
Y_{\phi_{k_{\widetilde{l}}+2}}=T_{i_{k_{\widetilde{l}}}}^{-1}\ldots T_{i_{1}}^{-1}T_{i_p}^{-1}\ldots T_{i_{k_{\widetilde{l}}+3}}^{-1}X_{i_{k_{\widetilde{l}}+2}}^+,
\\
\\
Y_{\phi_{k_{\widetilde{l}}+1}}=T_{i_{k_{\widetilde{l}}}}^{-1}\ldots T_{i_{1}}^{-1}T_{i_p}^{-1}\ldots T_{i_{k_{\widetilde{l}}+2}}^{-1}X_{i_{k_{\widetilde{l}}+1}}^+,
\end{array}
$$
\index[not]{f@$Y_{\phi_i}$}
which belong to $U_q^{res}(\m_{s^1}^+)$ by Lemma \ref{segmPBW} (v).

Let $\overline{U}_q^{res}([\phi_{k_{\widetilde{l}}},\phi_1])\subset U_q^{res}(\m_{s^1}^+)$ \index[not]{U@$\overline{U}_q^{res}([\phi_{k_{\widetilde{l}}},\phi_1])$} be the subalgebra generated by $(Y_{\phi_i})^{(k)}$, $i=1,\ldots k_{\widetilde{l}}$, $k\in \mathbb{N}$ , and $\overline{U}_q^{res}([\phi_p,\phi_{k_{\widetilde{l}}+1}])\subset U_q^{res}(\m_{s^1}^+)$ \index[not]{U@$\overline{U}_q^{res}([\phi_p,\phi_{k_{\widetilde{l}}+1}])$} the subalgebra generated by $(Y_{\phi_i})^{(k)}$, $i=k_{\widetilde{l}}+1,\ldots p$, $k\in \mathbb{N}$.

Now by Corollary \ref{segmq} (iii), 
\begin{equation}\label{ubdec}
U_q^{res}(\m_{s^1}^+)=\overline{U}_q^{res}([\phi_p,\phi_{k_{\widetilde{l}}+1}])\overline{U}_q^{res}([\phi_{k_{\widetilde{l}}},\phi_1]).
\end{equation}

We find the action of the automorphism $T_{\bar{w}_1}=T_{i_{1}}\ldots T_{i_{k_{\widetilde{l}}}}$ on $\overline{U}_q^{res}([\phi_p,\phi_{k_{\widetilde{l}}+1}])$ and on $\overline{U}_q^{res}([\phi_{k_{\widetilde{l}}},\phi_1])$. 

By the definition of the elements $Y_{\phi_{r}}$, for $r=k_{\widetilde{l}}+1,\ldots p$
\begin{equation}\label{tt}
T_{\bar{w}_1}(Y_{\phi_{r}})^{(k)}=(T_{i_{1}}\ldots T_{i_{k_{\widetilde{l}}}})(Y_{\phi_{r}})^{(k)}=(T_{i_p}^{-1}\ldots T_{i_{r+1}}^{-1})(X_{i_{r}}^+)^{(k)}.	
\end{equation}

Since $u$ is a highest weight vector, by Lemma \ref{bga} we have $\tau(T_{\gamma_{\widetilde{l}}} \ldots T_{\gamma_1})u=\omega(c)\tau(T_{i_p}\ldots T_{i_1})u$, where $c\in \mathbb{C}[q,q^{-1}]^*$. Therefore for any element $w\in V$ we can write
\begin{equation}\label{gT}
(u,\tau(T_{\gamma_1}\ldots T_{\gamma_{\widetilde{l}}})w)=c(u,\tau(T_{i_1}\ldots T_{i_p})w)=c(u,\tau(T_{s^1})w), 
\end{equation}
where $T_{s^1}=T_{i_1}\ldots T_{i_p}$ is the braid group element corresponding to reduced decomposition (\ref{s1dec'}).
Now by (\ref{tt}), using the definition (\ref{BGactqg}) of the braid group action and commutation relations (\ref{Ccomm1}), we obtain
for all generators $(Y_{\phi_r})^{(k)}$, $r=k_{\widetilde{l}}+1,\ldots p$, $k\in \mathbb{N}$ of $\overline{U}_q^{res}([\phi_p,\phi_{k_{\widetilde{l}}+1}])$
$$
T_{s^1}T_{\bar{w}_1}(Y_{\phi_{r}})^{(k)}=T_{i_1}\ldots T_{i_p}T_{i_{1}}\ldots T_{i_{k_{\widetilde{l}}}}(Y_{\phi_{r}})^{(k)}=
T_{i_1}\ldots T_{i_p}T_{i_p}^{-1}\ldots T_{i_{r+1}}^{-1}(X_{i_{r}}^+)^{(k)}=
$$
$$ 
=T_{i_1}\ldots T_{i_r}(X_{i_{r}}^+)^{(k)}=\frac{1}{[k]_{q_{i_r}}!}T_{i_1}\ldots T_{i_{r-1}}(-X_{i_{r}}^-K_{i_r})^k=
\frac{1}{[k]_{q_{i_r}}!}(-X_{\beta_r}^-K_{\beta_r})^k=
$$
$$
=(-1)^kq_{\beta_r}^{-k(k-1)} (X_{\beta_r}^-)^{(k)}K_{\beta_r}^k\in U_{U_q^{res}(H)}^{res}([-\beta_{k_{\widetilde{l}}+1},\ldots, -\beta_p]), 
$$
where $K_{\beta_r}=T_{i_1}\ldots T_{i_{r-1}}K_{i_r}$. Thus  
\begin{equation}\label{inc1*}
T_{s^1}T_{\bar{w}_1}(\overline{U}_q^{res}([\phi_p,\phi_{k_{\widetilde{l}}+1}]))\subset U_{U_q^{res}(H)}^{res}([-\beta_{k_{\widetilde{l}}+1},\ldots, -\beta_p]).
\end{equation}

On the other hand by the definition of the elements $Y_{\phi_{r}}$, for $r=1,\ldots k_{\widetilde{l}}$ we have for all generators $(Y_{\phi_r})^{(k)}$, $r=1,\ldots k_{\widetilde{l}}$, $k\in \mathbb{N}$ of $\overline{U}_q^{res}([\phi_{k_{\widetilde{l}}},\phi_1])$
\begin{equation}\label{tt1}
T_{\bar{w}_1}(Y_{\phi_{r}})^{(k)}=T_{i_{1}}\ldots T_{i_{k_{\widetilde{l}}}}(Y_{\phi_{r}})^{(k)}=\frac{1}{[k]_{q_{i_r}}!}T_{i_1}\ldots T_{i_{r}}(X_{i_{r}}^+)^k=\frac{1}{[k]_{q_{i_r}}!}T_{i_1}\ldots T_{i_{r-1}}(-X_{i_{r}}^-K_{i_r})^k= 
\end{equation}
$$
=\frac{1}{[k]_{q_{i_r}}!}(-X_{\beta_r}^-K_{\beta_r})^k=
$$
$$
=(-1)^kq_{\beta_r}^{-k(k-1)} (X_{\beta_r}^-)^{(k)}K_{\beta_r}^k\in U_{U_q^{res}(H)}^{res}([-\gamma_1,-\gamma_{\widetilde{l}}]),
$$
where $K_{\beta_r}=T_{i_1}\ldots T_{i_{r-1}}K_{i_r}$, and we also used the definition (\ref{BGactqg}) of the braid group action and commutation relations (\ref{Ccomm1}).

Thus 
\begin{equation}\label{inc2*}
T_{\bar{w}_1}(\overline{U}_q^{res}([\phi_{k_{\widetilde{l}}},\phi_1]))\subset U_{U_q^{res}(H)}^{res}([-\gamma_1,-\gamma_{\widetilde{l}}]).
\end{equation}

Now recall that in (\ref{G1}) $K''(m_1,\ldots, m_{\widetilde{l}})\in T_{\bar{w}_1}(U_q^{res}(\m_{s^1}^+))$, and hence we can write $K''(m_1,\ldots, m_{\widetilde{l}})=T_{\bar{w}_1}A,$ $A\in U_q^{res}(\m_{s^1}^+)$.

Factorizing $A$ according to (\ref{ubdec}),
$$
A=\sum_{i}a_ib_i,~a_i\in \overline{U}_q^{res}([\phi_p,\phi_{k_{\widetilde{l}}+1}]),~b_i\in \overline{U}_q^{res}([\phi_{k_{\widetilde{l}}},\phi_1]),
$$
we obtain 
$$
K''(m_1,\ldots, m_{\widetilde{l}})=T_{\bar{w}_1}(A)=\sum_{i}(T_{\bar{w}_1}(a_i))(T_{\bar{w}_1}(b_i))=\sum_{i}(T_{\bar{w}_1}(a_i))\bar{b}_i,
$$
where $\bar{b}_i=T_{\bar{w}_1}(b_i)\in U_{U_q^{res}(H)}^{res}([-\gamma_1,-\gamma_{\widetilde{l}}])$ by (\ref{inc2*}).

Now recalling (\ref{gT}), (\ref{G1}) can be rewritten as follows
$$
(u,\tau((X_{\gamma_1}^+)^{(m_1)}\ldots (X_{\gamma_{\widetilde{l}}}^+)^{(m_{\widetilde{l}})})v)=c(u,\tau(T_{i_1}\ldots T_{i_p} K''(m_1,\ldots, m_{\widetilde{l}}))v)=
$$
$$
=c\sum_{i}(u,\tau(T_{s^1}((T_{\bar{w}_1})(a_i))\bar{b}_i) v)=
c\sum_{i}(u,\tau(((T_{s^1}T_{\bar{w}_1})(a_i))T_{s^1} \bar{b}_i) v)=
$$
$$
=c\sum_{i}(u,\tau(g_iT_{s^1} \bar{b}_i) v),
$$
where $g_i=(T_{s^1}T_{\bar{w}_1})(a_i)\in U_{U_q^{res}(H)}^{res}([-\beta_{k_{\widetilde{l}}+1},\ldots, -\beta_p])$ by (\ref{inc1*}). 

Observe that the only elements of the algebra $U_{U_q^{res}(H)}^{res}([-\beta_{k_{\widetilde{l}}+1},\ldots, -\beta_p])$ whose weights are not strictly negative belong to $U_q^{res}(H)$, and hence by Lemma \ref{hwv} (i) the last formula takes the form
$$
(u,\tau((X_{\gamma_1}^+)^{(m_1)}\ldots (X_{\gamma_{\widetilde{l}}}^+)^{(m_{\widetilde{l}})})v)=c(u,\tau(T_{s^1} K''(m_1,\ldots, m_{\widetilde{l}}))v)=c\sum_{i}(u,\tau(g_i'T_{s^1} \bar{b}_i) v)=
$$
\begin{equation}\label{XXf}
=c(u,\tau(T_{s^1}U')v)=(u,\tau(T_{\gamma_1}\ldots T_{\gamma_{\widetilde{l}}}U')v),
\end{equation}
where $g_i'\in U_q^{res}(H)$, and $U'=\sum_{i}(T_{s^1})^{-1}(g_i')\bar{b}_i\in U_{U_q^{res}(H)}^{res}([-\gamma_1,-\gamma_{\widetilde{l}}])$.

Recall finally that by Corollary \ref{segmq} (ii) we have $U_{U_q^{res}(H)}^{res}([-\gamma_1,-\gamma_{\widetilde{l}}])=U_q^{res}([-\gamma_1,-\gamma_{\widetilde{l}}])U_q^{res}(H)$ and that according to Lemma \ref{wHact} elements of the subalgebra $U_q^{res}(H)$, which is preserved by the action  of $\tau$, act on weight vectors by multiplication by elements of $\mathbb{C}[q,q^{-1}]$. Therefore (\ref{XXf}) immediately implies (\ref{l3st}).

Finally note that $v$ and $\tau(T_{\gamma_{\widetilde{l}}}\ldots T_{\gamma_{1}})u$ are weight vectors, and different weight subspaces of $V$ are orthogonal with respect to the contravariant form $(~\cdot~,~\cdot~)$. Therefore we can assume also that $U(m_1,\ldots, m_{\widetilde{l}})$ in (\ref{l3st}) belongs to a weight subspace of $U_q^{res}(\g)$, so that the weight of $\tau(T_{\gamma_{\widetilde{l}}}\ldots T_{\gamma_{1}})u$ is equal to that of $\tau(U(m_1,\ldots, m_{\widetilde{l}}))v$. This completes the proof.

\end{proof}

Now we obtain an analogue of Lemma \ref{l2} for $s^2$. We argue in the way similar to the previous discussion in the case of $s^1$. According to Lemma \ref{invdec} (i) $s^2$ is the longest element in the Weyl group $W(\m_{s^2},\h_{s^2})$ \index[not]{W@$W(\m_{s^2},\h_{s^2})$} of the semisimple part $\m_{s^2}$ of a Levi subalgebra of $\g$, the Cartan subalgebra of $\m_{s^2}$ is denoted by $\h_{s^2}$. 

The system of positive roots $\Delta_+(\m_{s^2},\h_{s^2}):=\Delta_+\cap \Delta(\m_{s^2},\h_{s^2})$ \index[not]{D@$\Delta_+(\m_{s^2},\h_{s^2})$} of the root system $\Delta(\m_{s^2},\h_{s^2})=\Delta_{s^2}^{-1}\bigcup (-\Delta_{s^2}^{-1})$ is the set (we again use the notation of (\ref{NO}))
\begin{eqnarray}\label{m2}
\Delta_+(\m_{s^2},\h_{s^2})=\{\gamma_{\widetilde{l}+1}, \ldots, \gamma_{\widetilde{l}+2}, \ldots , \gamma_{\widetilde{l}+3},\ldots,
\gamma_{l'},\beta_{t'+\frac{p'+l'-\widetilde{l}}{2}+1}^2, \ldots,\beta_{t'+p'}^2\},
\end{eqnarray}
which is in fact a segment, and $s^2$ acts on the elements of this set by multiplication by $-1$. The roots in (\ref{m2}) are ordered as in the normal ordering of $\Delta_+$ associated to $s$. Note that, in fact, $\Delta_+(\m_{s^2},\h_{s^2})=\Delta_+^s(\m_{s^2},\h_{s^2})=\Delta(\m_{s^2},\h_{s^2})\cap \Delta_+^s=\Delta_{s^2}^{-1}$.

As before if $\overline{w}=s_{i_1}\ldots s_{i_D}$ is the reduced decomposition of the longest element $\overline{w}\in W$ corresponding to $\Delta_+$, we write $\gamma_m=s_{i_{1}}\ldots s_{i_{k_{\widetilde{l}+1}}}\ldots  s_{i_{k_m-1}}\alpha_{i_{k_m}}$ for $m=\widetilde{l}+1,\ldots ,l'$. 

Let $\hat{w}=s_{i_1}\ldots s_{i_{k_{\widetilde{l}+1}-1}}$, \index[not]{w@$\hat{w}$} $T_{\hat{w}}=T_{i_1}\ldots T_{i_{k_{\widetilde{l}+1}-1}}$. Then $\widetilde{\Delta}_+(\m_{s^2},\h_{s^2}):={\hat{w}}^{-1}(\Delta_+(\m_{s^2},\h_{s^2}))$ \index[not]{D@$\widetilde{\Delta}_+(\m_{s^2},\h_{s^2})$} is the set of positive roots of the semisimple part $\widetilde{\m}_{s^2}:={\rm Ad}{\hat{w}}^{-1}(\m_{s^2})$ \index[not]{m@$\widetilde{\m}_{s^2}$} of a standard Levi subalgebra of $\g$. Indeed, by part (iv) of Lemma \ref{wN} the reduced decomposition $s_{i_{k_{\widetilde{l}+1}}}\ldots s_{i_D}$, which is a part of the reduced decomposition $\overline{w}=s_{i_1}\ldots s_{i_D}$, can be completed to a reduced decomposition of $\overline{w}$, $\overline{w}=s_{i_{k_{\widetilde{l}+1}}}\ldots s_{i_D}s_{j_1}\ldots s_{j_{k_{\widetilde{l}+1}-1}}$, and by the definition the roots from the set $\widetilde{\Delta}_+(\m_{s^2},\h_{s^2})={\hat{w}}^{-1}(\Delta_+(\m_{s^2},\h_{s^2}))$ form an initial segment in the corresponding normal ordering of $\Delta_+$. Therefore, if $\widetilde{\m}_{s^2}$ was not the semisimple part of a standard Levi subalgebra of $\g$ (relative to $\Delta_+$), by Lemma \ref{invdec} (iii) there would be some roots preceding those from the set ${\hat{w}}^{-1}(\Delta_+(\m_{s^2},\h_{s^2}))$ in this normal ordering.

Therefore one can define the subalgebra $U_h(\widetilde{\m}_{s^2})\subset U_h(\g)$ \index[not]{U@$U_h(\widetilde{\m}_{s^2})$} generated by the elements $X_i^\pm$ and $H_i$ for $\alpha_i\in \widetilde{\Delta}_+(\m_{s^2},\h_{s^2})$ and its restricted specialization $U_q^{res}(\widetilde{\m}_{s^2})\subset U_q^{res}(\g)$. \index[not]{U@$U_q^{res}(\widetilde{\m}_{s^2})$}

Lemma \ref{invdec} (i) also implies that $\widetilde{s}^2:={\hat{w}}^{-1}s^2{\hat{w}}$ \index[not]{s@$\widetilde{s}^2$} is the longest element in the Weyl group $\widetilde{W}(\m_{s^2},\h_{s^2}):={\hat{w}}^{-1}W(\m_{s^2},\h_{s^2}){\hat{w}}$ \index[not]{W@$\widetilde{W}(\m_{s^2},\h_{s^2})$} with respect to the system of simple roots in ${\hat{w}}^{-1}(\Delta_+(\m_{s^2},\h_{s^2}))$.

For any root $\beta \in \Delta_+$ denote
\begin{equation}\label{Xtilde}
\widetilde{X}_{\beta}^\pm=T_{\hat{w}}^{-1}(X_{\beta}^\pm), ~~\widetilde{\overline{X}}_{\beta}^\pm=\overline{T}_{\hat{w}}^{-1}(\overline{X}_{\beta}^\pm), ~~\widetilde{T}_\beta=T_{\hat{w}}^{-1}T_\beta T_{\hat{w}}. \index[not]{X@$\widetilde{X}_{\beta}^\pm$} \index[not]{X@$\widetilde{\overline{X}}_{\beta}^\pm$} \index[not]{T@$\widetilde{T}_\beta$}
\end{equation}

Denote by $\beta_r$ \index[not]{b@$\beta_r$} the greatest root $\beta_{t'+p'}^2$ in the segment $\Delta_+(\m_{s^2},\h_{s^2})$ with respect to ordering (\ref{m2}). Then by the definition $U_q^{res}(\widetilde{\m}_{s^2}^-):=T_{\hat{w}}^{-1}(U_q^{res}([-\gamma_{\widetilde{l}+1},-\beta_r]))$ \index[not]{U@$U_q^{res}(\widetilde{\m}_{s^2}^-)$} is the subalgebra in the algebra $U_q^{res}(\widetilde{\m}_{s^2})$ generated by $({X}_{\beta}^-)^{(k)}$ for simple roots $\beta \in \widetilde{\Delta}_+(\m_{s^2},\h_{s^2})$ and $k\geq 0$.  

Similarly to Lemma \ref{l2} we have the following statement.
\begin{lemma}\label{l2'}
Let $V$ be a finite rank representation of $U_h(\g)$, $u,v\in V$ weight vectors.
Then for any $m_{\widetilde{l}+1},\ldots, m_{l'}\in \mathbb{N}$ one has
\begin{equation}\label{T5'}
(u,\tau((\widetilde{X}_{\gamma_{\widetilde{l}+1}}^+)^{(m_{\widetilde{l}+1})}\ldots (\widetilde{X}_{\gamma_{l'}}^+)^{(m_{l'})})v)=(u,\tau(\widetilde{K}(m_{\widetilde{l}+1},\ldots, m_{l'})\widetilde{T}_{\gamma_{\widetilde{l}+1}}\ldots \widetilde{T}_{\gamma_{l'}})v),
\end{equation}
where $\widetilde{K}(m_{\widetilde{l}+1},\ldots, m_{l'})\in  T_{\hat{w}}^{-1}(U_q^{res}([-\gamma_{\widetilde{l}+1},-\beta_r]))\overline{T}_{\hat{w}}^{-1}(\overline{U}_q^{res}([\gamma_{\widetilde{l}+1},\gamma_{l'}]))=U_q^{res}(\widetilde{\m}_{s^2}^-)\overline{T}_{\hat{w}}^{-1}(\overline{U}_q^{res}([\gamma_{\widetilde{l}+1},\gamma_{l'}]))$, and $\widetilde{K}(m_{\widetilde{l}+1},\ldots, m_{l'})$ belongs to a weight subspace of $U_q^{res}(\g)$. 
\end{lemma}

From Lemmas \ref{l3} and \ref{l2'} we obtain the following statement.
\begin{lemma}\label{mainlpr}
Let $V$ be a finite rank representation of $U_h(\g)$, $u,v\in V$ weight vectors. Suppose that $u$ is a highest weight vector.
Then for any $m_1,\ldots, m_{l'}\in \mathbb{N}$
$$
(u,\tau((X_{\gamma_1}^+)^{(m_1)}\ldots (X_{\gamma_{\widetilde{l}}}^+)^{(m_{\widetilde{l}})}({X}_{\gamma_{\widetilde{l}+1}}^+)^{(m_{\widetilde{l}+1})}\ldots ({X}_{\gamma_{l'}}^+)^{(m_{l'})})v)=
$$
\begin{equation}\label{prepm}
=(u,\tau(T_{\gamma_1}\ldots T_{\gamma_{\widetilde{l}}}
{T}_{\gamma_{\widetilde{l}+1}}\ldots {T}_{\gamma_{l'}}T_{w'}
X(m_1,\ldots, m_{l'})T_{w'}^{-1})v), 
\end{equation}
where $X(m_1,\ldots, m_{l'})\in U_q^{res}(\n_+)$, $T_{w'}=T_{i_1}\ldots T_{i_{k_{l'}}}$ for
$w'=s_{i_1}\ldots s_{i_{k_{l'}}}$, \index[not]{w@$w'$} and $X(m_1,\ldots, m_{l'})$ belongs to a weight subspace of $U_q^{res}(\g)$.
\end{lemma}

\begin{proof}
By Lemmas \ref{l3} and \ref{l2'} we have, using the notation introduced in these lemmas,
$$
(u,\tau((X_{\gamma_1}^+)^{(m_1)}\ldots (X_{\gamma_{\widetilde{l}}}^+)^{(m_{\widetilde{l}})}({X}_{\gamma_{\widetilde{l}+1}}^+)^{(m_{\widetilde{l}+1})}\ldots ({X}_{\gamma_{l'}}^+)^{(m_{l'})})v)=
$$
$$
=(u,\tau((X_{\gamma_1}^+)^{(m_1)}\ldots (X_{\gamma_{\widetilde{l}}}^+)^{(m_{\widetilde{l}})}T_{\hat{w}}(\widetilde{X}_{\gamma_{\widetilde{l}+1}}^+)^{(m_{\widetilde{l}+1})}\ldots (\widetilde{X}_{\gamma_{l'}}^+)^{(m_{l'})}T_{\hat{w}}^{-1})v)=
$$
\begin{equation}\label{a0}
=(u,\tau(T_{\gamma_1}\ldots T_{\gamma_{\widetilde{l}}}U(m_1,\ldots, m_{\widetilde{l}})T_{\hat{w}}\widetilde{K}(m_{\widetilde{l}+1},\ldots, m_{l'})\widetilde{T}_{\gamma_{\widetilde{l}+1}}\ldots \widetilde{T}_{\gamma_{l'}}T_{\hat{w}}^{-1})v),
\end{equation}
where $U(m_1,\ldots, m_{\widetilde{l}})\in U_q^{res}([-\gamma_1,-\gamma_{\widetilde{l}}])$, $\widetilde{K}(m_{\widetilde{l}+1},\ldots, m_{l'})\in  T_{\hat{w}}^{-1}(U_q^{res}([-\gamma_{\widetilde{l}+1},-\beta_r]))\overline{T}_{\hat{w}}^{-1}(\overline{U}_q^{res}([\gamma_{\widetilde{l}+1},\gamma_{l'}]))$, so we can write $\widetilde{K}(m_{\widetilde{l}+1},\ldots, m_{l'})=\sum_j T_{\hat{w}}^{-1}(a_j) b_j$, $a_j\in U_q^{res}([-\gamma_{\widetilde{l}+1},-\beta_r])$, $b_j\in \overline{T}_{\hat{w}}^{-1}(\overline{U}_q^{res}([\gamma_{\widetilde{l}+1},\gamma_{l'}]))$. Substituting this into (\ref{a0}) we obtain
\begin{equation}\label{a1}
(u,\tau((X_{\gamma_1}^+)^{(m_1)}\ldots (X_{\gamma_{\widetilde{l}}}^+)^{(m_{\widetilde{l}})}({X}_{\gamma_{\widetilde{l}+1}}^+)^{(m_{\widetilde{l}+1})}\ldots ({X}_{\gamma_{l'}}^+)^{(m_{l'})})v)=
\end{equation}
$$
=(u,\tau(T_{\gamma_1}\ldots T_{\gamma_{\widetilde{l}}}U(m_1,\ldots, m_{\widetilde{l}})T_{\hat{w}}\sum_j T_{\hat{w}}^{-1}(a_j) b_j\widetilde{T}_{\gamma_{\widetilde{l}+1}}\ldots \widetilde{T}_{\gamma_{l'}}T_{\hat{w}}^{-1})v)=
$$
$$
=\sum_j(u,\tau(T_{\gamma_1}\ldots T_{\gamma_{\widetilde{l}}}U(m_1,\ldots, m_{\widetilde{l}})a_j T_{\hat{w}} b_j\widetilde{T}_{\gamma_{\widetilde{l}+1}}\ldots \widetilde{T}_{\gamma_{l'}}T_{\hat{w}}^{-1})v),
$$
where $U(m_1,\ldots, m_{\widetilde{l}})a_j\in U_q^{res}([-\gamma_1,-\gamma_{\widetilde{l}}])U_q^{res}([-\gamma_{\widetilde{l}+1},-\beta_r])\subset U_q^{res}([-\gamma_1,-\beta_r])$.

Observe that by Corollary \ref{segmq} (iii) $U_q^{res}([-\gamma_1,-\beta_r])=U_q^{res}([-\beta_{k_{\widetilde{l}}+1},-\beta_r])U_q^{res}([-\gamma_1,-\gamma_{\widetilde{l}}])$, so $U(m_1,\ldots, m_{\widetilde{l}})a_j=\sum_ic_i^jo_i^j$ for some $c_i^j\in U_q^{res}([-\beta_{k_{\widetilde{l}}+1},-\beta_r])$, $o_i^j\in U_q^{res}([-\gamma_1,-\gamma_{\widetilde{l}}])$.


Now (\ref{a1}) takes the form
\begin{equation}\label{a2}
(u,\tau((X_{\gamma_1}^+)^{(m_1)}\ldots (X_{\gamma_{\widetilde{l}}}^+)^{(m_{\widetilde{l}})}({X}_{\gamma_{\widetilde{l}+1}}^+)^{(m_{\widetilde{l}+1})}\ldots ({X}_{\gamma_{l'}}^+)^{(m_{l'})})v)=
\end{equation}
$$
=\sum_{i,j}(u,\tau(T_{\gamma_1}\ldots T_{\gamma_{\widetilde{l}}}c_i^jo_i^j T_{\hat{w}} b_j\widetilde{T}_{\gamma_{\widetilde{l}+1}}\ldots \widetilde{T}_{\gamma_{l'}}T_{\hat{w}}^{-1})v)
=\sum_{i,j}(u,\tau((T_{\gamma_1}\ldots T_{\gamma_{\widetilde{l}}}(c_i^j))T_{\gamma_1}\ldots T_{\gamma_{\widetilde{l}}}o_i^j T_{\hat{w}} b_j\widetilde{T}_{\gamma_{\widetilde{l}+1}}\ldots \widetilde{T}_{\gamma_{l'}}T_{\hat{w}}^{-1})v).
$$

If $\gamma_{\widetilde{l}} <\beta \leq \beta_r$ then by (\ref{s1-}) we have $\beta \not \in \Delta_+(\m_{s^1},\h_{s^1})=\{\alpha\in \Delta_+:s^1\alpha \in \Delta_-\}$. This implies $s^1\beta \in \Delta_+$, and hence by the definition of the braid group action, for $\gamma_{\widetilde{l}} <\beta \leq \beta_r$, $k\in \mathbb{N}$ the element $T_{\gamma_1}\ldots T_{\gamma_{\widetilde{l}}}(X_{\beta}^-)^{(k)}$ has weight $-ks^1\beta<0$. Since the elements $(X_{\beta}^-)^{(k)}$ for $\gamma_{\widetilde{l}} <\beta \leq \beta_r$, $k\in \mathbb{N}$ generate $U_q^{res}([-\beta_{k_{\widetilde{l}}+1},-\beta_r])$, we deduce that the only elements of $T_{\gamma_1}\ldots T_{\gamma_{\widetilde{l}}}(U_q^{res}([-\beta_{k_{\widetilde{l}}+1},-\beta_r]))$ which do not have strictly negative weights belong to $\mathbb{C}[q,q^{-1}]$. Thus by Lemma \ref{hwv} (i) the right hand side of (\ref{a2}) takes the form
\begin{equation}\label{a3}
(u,\tau((X_{\gamma_1}^+)^{(m_1)}\ldots (X_{\gamma_{\widetilde{l}}}^+)^{(m_{\widetilde{l}})}({X}_{\gamma_{\widetilde{l}+1}}^+)^{(m_{\widetilde{l}+1})}\ldots ({X}_{\gamma_{l'}}^+)^{(m_{l'})})v)=
\end{equation}
$$
=\sum_{i,j}(u,\tau((T_{\gamma_1}\ldots T_{\gamma_{\widetilde{l}}}({c'}_i^j))T_{\gamma_1}\ldots T_{\gamma_{\widetilde{l}}}o_i^j T_{\hat{w}} b_j\widetilde{T}_{\gamma_{\widetilde{l}+1}}\ldots \widetilde{T}_{\gamma_{l'}}T_{\hat{w}}^{-1})v)
=\sum_{j}(u,\tau(T_{\gamma_1}\ldots T_{\gamma_{\widetilde{l}}}p_j T_{\hat{w}} b_j\widetilde{T}_{\gamma_{\widetilde{l}+1}}\ldots \widetilde{T}_{\gamma_{l'}}T_{\hat{w}}^{-1})v),
$$
where ${c'}_i^j\in \mathbb{C}[q,q^{-1}]$, $p_j=\sum_i {c'}_i^jo_i^j\in U_q^{res}([-\gamma_1,-\gamma_{\widetilde{l}}])$.

Now we proceed as in the proof of Lemma \ref{l3}. Namely, denote $\widetilde{T}^2=\widetilde{T}_{\gamma_{\widetilde{l}+1}}\ldots \widetilde{T}_{\gamma_{l'}}$. \index[not]{T@$\widetilde{T}^2$} We find the action of $(\widetilde{T}^2)^{-1}$ on the generators of the algebra $\overline{T}_{\hat{w}}^{-1}(\overline{U}_q^{res}([\gamma_{\widetilde{l}+1},\gamma_{l'}]))$.

For $\gamma_{\widetilde{l}+1}\leq \beta_{q} \leq \gamma_{l'}$ we have by the definition $\beta_q ={\hat{w}}s_{k_{\widetilde{l}+1}}\ldots s_{i_{q-1}}\alpha_{i_q}$, $q=k_{\widetilde{l}+1},\ldots, k_{l'}$, where we use the reduced decomposition $\overline{w}=s_{i_1}\ldots s_{i_D}$ corresponding to the normal ordering in $\Delta_+$ associated to $s$, and ${\hat{w}}=s_{i_1}\ldots s_{i_{k_{\widetilde{l}+1}-1}}$. 

Similarly to (\ref{t11}) we infer for $\gamma_{\widetilde{l}+1}\leq \beta_{q} \leq \gamma_{l'}$ that
\begin{equation}\label{t11'}
(\widetilde{T}^2)^{-1}\widetilde{\overline{X}}_{\beta_q}^+=T_{i_{k_{\widetilde{l}+1}}}\ldots T_{i_{k_{l'}}}T_{i_{k_{l'}}}\ldots T_{i_{q+1}}X_{i_q}^+.
\end{equation}

Denote for $q=k_{\widetilde{l}+1},\ldots, k_{l'}$, $\delta_q=s_{i_{k_{l'}}}\ldots s_{i_{q+1}}\alpha_{i_q}$, \index[not]{d@$\delta_q$} and
\begin{equation}\label{Xtilde1}
\widehat{X}_{\delta_q}^+=T_{i_{k_{l'}}}\ldots T_{i_{q+1}}X_{i_q}^+. \index[not]{X@$\widehat{X}_{\delta_q}^+$}
\end{equation}
We show that $\widehat{X}_{\delta_q}^+\in U_q^{res}(\n_+)$ for $q=k_{\widetilde{l}+1},\ldots, k_{l'}$.

Indeed, the element $\bar{w}_2^{-1}:=s_{i_{k_{l'}}}\ldots s_{k_{\widetilde{l}+1}}$ \index[not]{w@$\bar{w}_2$} is a part of the reduced decomposition
\begin{equation}\label{wbarinv}
\overline{w}=s_{i_{D}}\ldots s_{i_{k_{l'}+1}}s_{i_{k_{l'}}}\ldots  s_{i_{k_{\widetilde{l}+1}}}\ldots s_{i_1}   
\end{equation}
inverse to that corresponding to the normal ordering in $\Delta_+$ associated to $s$. 
Therefore by part (iv) of Lemma \ref{wN} we can consider it as an initial part of another reduced decomposition of the longest element $\overline{w}\in W$. Now by Lemma \ref{segmPBW} (i) one has $\widehat{X}_{\delta_q}^+\in U_q^{res}(\n_+)$ for $q=k_{\widetilde{l}+1},\ldots, k_{l'}$ which obviously implies that for the generators $(\widetilde{\overline{X}}_{\beta_q}^+)^{(k)}$, $q=k_{\widetilde{l}+1},\ldots, k_{l'}$, $k\in \mathbb{N}$ of the algebra $\overline{T}_{\hat{w}}^{-1}(\overline{U}_q^{res}([\gamma_{\widetilde{l}+1},\gamma_{l'}]))$ we have by (\ref{t11'})
$$
(\widetilde{T}^2)^{-1}(\widetilde{\overline{X}}_{\beta_q}^+)^{(k)}=T_{i_{k_{\widetilde{l}+1}}}\ldots T_{i_{k_{l'}}}(\widehat{X}_{\delta_q}^+)^{(k)}=T_{\bar{w}_2}(\widehat{X}_{\delta_q}^+)^{(k)}\in T_{i_{k_{\widetilde{l}+1}}}\ldots T_{i_{k_{l'}}}(U_q^{res}(\n_+))=T_{\bar{w}_2}(U_q^{res}(\n_+)),
$$
where $T_{\bar{w}_2}=T_{i_{k_{\widetilde{l}+1}}}\ldots T_{i_{k_{l'}}}$.
Thus
\begin{equation}\label{tild*1}
(\widetilde{T}^2)^{-1}(\overline{T}_{\hat{w}}^{-1}(\overline{U}_q^{res}([\gamma_{\widetilde{l}+1},\gamma_{l'}])))\subset T_{\bar{w}_2}(U_q^{res}(\n_+)).
\end{equation}   

Now, from (\ref{a3}) and (\ref{tild*1}) we obtain
\begin{equation}\label{a4}
(u,\tau((X_{\gamma_1}^+)^{(m_1)}\ldots (X_{\gamma_{\widetilde{l}}}^+)^{(m_{\widetilde{l}})}({X}_{\gamma_{\widetilde{l}+1}}^+)^{(m_{\widetilde{l}+1})}\ldots ({X}_{\gamma_{l'}}^+)^{(m_{l'})})v)=
\end{equation}
$$
=\sum_{j}(u,\tau(T_{\gamma_1}\ldots T_{\gamma_{\widetilde{l}}}p_j T_{\hat{w}} \widetilde{T}_{\gamma_{\widetilde{l}+1}}\ldots \widetilde{T}_{\gamma_{l'}}((\widetilde{T}^2)^{-1}(b_j))T_{\hat{w}}^{-1})v)
=\sum_{j}(u,\tau(T_{\gamma_1}\ldots T_{\gamma_{\widetilde{l}}}p_j T_{\hat{w}} \widetilde{T}_{\gamma_{\widetilde{l}+1}}\ldots \widetilde{T}_{\gamma_{l'}}(T_{\bar{w}_2}(b_j'))T_{\hat{w}}^{-1})v)=
$$
$$
=\sum_{j}(u,\tau(T_{\gamma_1}\ldots T_{\gamma_{\widetilde{l}}}p_j T_{\hat{w}} \widetilde{T}_{\gamma_{\widetilde{l}+1}}\ldots \widetilde{T}_{\gamma_{l'}}T_{\bar{w}_2}b_j'T_{\bar{w}_2}^{-1}T_{\hat{w}}^{-1})v)
=\sum_{j}(u,\tau(T_{\gamma_1}\ldots T_{\gamma_{\widetilde{l}}}p_j \widetilde{T} b_j'T_{w'}^{-1})v),
$$
where $(\widetilde{T}^2)^{-1}(b_j)=T_{\bar{w}_2}(b_j')$, $b_j'\in U_q^{res}(\n_+)$, 
$$
\widetilde{T}=T_{\hat{w}} \widetilde{T}_{\gamma_{\widetilde{l}+1}}\ldots \widetilde{T}_{\gamma_{l'}}T_{\bar{w}_2}, \index[not]{T@$\widetilde{T}$}
$$
and
$$
T_{w'}=T_{\hat{w}}T_{\bar{w}_2}=T_{i_1}\ldots T_{i_{k_{l'}}},
$$
for the reduced decomposition
\begin{equation}\label{w'}
w'={\hat{w}}\bar{w}_2=s_{i_1}\ldots s_{i_{k_{l'}}}. \index[not]{w@$w'$}
\end{equation}

We find the action of $\widetilde{T}^{-1}$ on the algebra $U_q^{res}([-\gamma_1,-\gamma_{\widetilde{l}}])$. First we obtain a convenient expression for $\widetilde{T}$.
If
\begin{equation}\label{wbar}
\overline{w}=s_{i_1}\ldots s_{i_{k_{\widetilde{l}+1}}}\ldots s_{i_{k_{l'}}}s_{i_{k_{l'}+1}}\ldots s_{i_{r}}\ldots s_{i_{D}}
\end{equation}
is the reduced decomposition of the longest element of the Weyl group corresponding to the normal ordering in $\Delta_+$ associated to $s$, so that $\widetilde{s}^2=s_{i_{k_{\widetilde{l}+1}}}\ldots s_{i_{k_{l'}}}s_{i_{k_{l'}+1}}\ldots s_{i_{r}}$ is the corresponding reduced decomposition of $\widetilde{s}^2$, then, similarly to (\ref{T1red}), we have
$$
\widetilde{T}_{\gamma_{\widetilde{l}+1}}\ldots \widetilde{T}_{\gamma_{l'}}=T_{i_{r}}\ldots T_{i_{k_{l'}+1}}T_{i_{k_{l'}}}^{-1}T_{i_{k_{l'}-1}}^{-1}\ldots T_{i_{k_{\widetilde{l}+1}}}^{-1},
$$
and hence
\begin{equation}\label{T2red}
\widetilde{T}_{\gamma_{\widetilde{l}+1}}\ldots \widetilde{T}_{\gamma_{l'}}T_{\bar{w}_2}=\widetilde{T}_{\gamma_{\widetilde{l}+1}}\ldots \widetilde{T}_{\gamma_{l'}}T_{i_{k_{\widetilde{l}+1}}}\ldots T_{i_{k_{l'}}}=T_{i_{r}}\ldots T_{i_{k_{l'}+1}}.
\end{equation}
This identity together with the definition of $T_{\hat{w}}$ imply 
\begin{equation}\label{Tdef}
\widetilde{T}=T_{\hat{w}}\widetilde{T}_{\gamma_{\widetilde{l}+1}}\ldots \widetilde{T}_{\gamma_{l'}}T_{\bar{w}_2}=T_{i_1}\ldots T_{i_{k_{\widetilde{l}+1}-1}}T_{i_{r}}\ldots T_{i_{k_{l'}+1}}.
\end{equation}

Observe that $\widetilde{T}$ is the braid group element corresponding to the reduced decomposition
\begin{equation}\label{Rd}
s_{i_1}\ldots s_{i_{k_{\widetilde{l}+1}-1}}s_{i_{r}}\ldots s_{i_{k_{l'}+1}}	
\end{equation}
which is a part of the reduced decomposition
$$
\overline{w}=s_{i_1}\ldots s_{i_{k_{\widetilde{l}+1}-1}}s_{i_{r}}\ldots s_{i_{k_{l'}+1}}s_{i_{k_{l'}}}\ldots
s_{i_{k_{\widetilde{l}+1}}} \ldots s_{i_D}
$$
obtained from reduced decomposition (\ref{wbar}),
$$
\overline{w}=s_{i_1}\ldots s_{i_{k_{\widetilde{l}+1}-1}}s_{i_{k_{\widetilde{l}+1}}}\ldots s_{i_{k_{l'}}}s_{i_{k_{l'}+1}}\ldots s_{i_{r}}\ldots s_{i_D}, 
$$
by inverting the part $\widetilde{s}^2=s_{i_{k_{\widetilde{l}+1}}}\ldots s_{i_{k_{l'}}}s_{i_{k_{l'}+1}}\ldots s_{i_{r}}$. This inversion gives a reduced decomposition again because $\widetilde{s}^2=-1$ is the longest element in $\widetilde{W}(\m_{s^2},\h_{s^2})$.

Now for $\beta_1\leq \beta_t\leq \beta_{k_{\widetilde{l}}}$, $\beta_t=s_{i_1}\ldots s_{i_{t-1}}\alpha_{i_{t}}$, $X_{\beta_{t}}^-=T_{i_1}\ldots T_{i_{t-1}}X_{i_{t}}^-$, $t=1,\ldots, k_{\widetilde{l}}$, we have by (\ref{Tdef}), by commutation relations (\ref{Ccomm1}), and by the definition of the braid group action (see (\ref{BGactg})), for all generators $(X_{\beta_{t}}^-)^{(k)}$, $t=1,\ldots, k_{\widetilde{l}}$, $k\in \mathbb{N}$, of the algebra $U_q^{res}([-\gamma_1,-\gamma_{\widetilde{l}}])$
\begin{equation}\label{dd}
\widetilde{T}^{-1}(X_{\beta_{t}}^-)^{(k)}=\frac{1}{[k]_{q_{i_t}}!}T_{i_{k_{l'}+1}}^{-1}\ldots T_{i_{r}}^{-1}T_{i_{k_{\widetilde{l}+1}-1}}^{-1}\ldots T_{i_{t+1}}^{-1}T_{i_t}^{-1}\ldots T_{i_1}^{-1}T_{i_1}\ldots T_{i_{t-1}}(X_{i_{t}}^-)^k=
\end{equation}
$$
=\frac{1}{[k]_{q_{i_t}}!}(T_{i_{k_{l'}+1}}^{-1}\ldots T_{i_{r}}^{-1}T_{i_{k_{\widetilde{l}+1}-1}}^{-1}\ldots T_{i_{t+1}}^{-1}T_{i_t}^{-1}X_{i_{t}}^-)^k=\frac{1}{[k]_{q_{i_t}}!}(-T_{i_{k_{l'}+1}}^{-1}\ldots T_{i_{r}}^{-1}T_{i_{k_{\widetilde{l}+1}-1}}^{-1}\ldots T_{i_{t+1}}^{-1}X_{i_t}^+K_{i_t})^k=
$$
$$
=\frac{1}{[k]_{q_{i_t}}!}(-T_{i_{k_{l'}+1}}^{-1}\ldots T_{i_{r}}^{-1}T_{i_{k_{\widetilde{l}+1}-1}}^{-1}\ldots T_{i_{t+1}}^{-1}(X_{i_{t}}^+)R_t)^k=(-1)^kq_{i_t}^{k(k-1)}T_{i_{k_{l'}+1}}^{-1}\ldots T_{i_{r}}^{-1}T_{i_{k_{\widetilde{l}+1}-1}}^{-1}\ldots T_{i_{t+1}}^{-1}((X_{i_{t}}^+)^{(k)})R_t^k, 
$$
where 
$$
R_t=T_{i_{k_{l'}+1}}^{-1}\ldots T_{i_{r}}^{-1}T_{i_{k_{\widetilde{l}+1}-1}}^{-1}\ldots T_{i_{t+1}}^{-1}(K_{i_{t}}). \index[not]{R@$R_t$}
$$

Since $s_{i_{k_{l'}+1}}\ldots s_{i_{r}}s_{i_{k_{\widetilde{l}+1}-1}}\ldots s_{i_1}$ is a reduced decomposition by (\ref{Rd}), by part (iv) of Lemma \ref{wN} we can consider it as an initial part of a reduced decomposition of the longest element $\overline{w}\in W$, and hence by Lemma \ref{segmPBW} (v)
\begin{equation}\label{Zqk}
Z_t^{(k)}:=T_{i_{k_{l'}+1}}^{-1}\ldots T_{i_{r}}^{-1}T_{i_{k_{\widetilde{l}+1}-1}}^{-1}\ldots T_{i_{t+1}}^{-1}(X_{i_{t}}^+)^{(k)}\in U_q^{res}(\n_+). \index[not]{Z@$Z_t^{(k)}$}
\end{equation}

Therefore for $\beta_1\leq \beta_t\leq \beta_{k_{\widetilde{l}}}$ using (\ref{Zqk}) one has from (\ref{dd}) for all generators $(X_{\beta_{t}}^-)^{(k)}$, $t=1,\ldots, k_{\widetilde{l}}$, $k\in \mathbb{N}$, of the algebra $U_q^{res}([-\gamma_1,-\gamma_{\widetilde{l}}])$
\begin{equation}\label{Xb+}
\widetilde{T}^{-1}(X_{\beta_{t}}^-)^{(k)}=(-1)^kq_{i_t}^{k(k-1)}Z_t^{(k)}R_t^k\in U_q^{res}({\frak b}_+).
\end{equation}
This implies 
$$
\widetilde{T}^{-1}(U_q^{res}([-\gamma_1,-\gamma_{\widetilde{l}}]))\subset U_q^{res}({\frak b}_+),
$$ 
and hence (\ref{a4}) takes the form
\begin{equation}\label{a5'}
(u,\tau((X_{\gamma_1}^+)^{(m_1)}\ldots (X_{\gamma_{\widetilde{l}}}^+)^{(m_{\widetilde{l}})}({X}_{\gamma_{\widetilde{l}+1}}^+)^{(m_{\widetilde{l}+1})}\ldots ({X}_{\gamma_{l'}}^+)^{(m_{l'})})v)=
\end{equation}
$$
=\sum_{j}(u,\tau(T_{\gamma_1}\ldots T_{\gamma_{\widetilde{l}}}\widetilde{T} p_j'  b_j'T_{w'}^{-1})v)=(u,\tau(T_{\gamma_1}\ldots T_{\gamma_{\widetilde{l}}}\widetilde{T} X'T_{w'}^{-1})v),
$$
where $p_j'=\widetilde{T}^{-1}(p_j)\in U_q^{res}({\frak b}_+)$, so that $X'=\sum_{j}p_j'  b_j'\in U_q^{res}({\frak b}_+)$.

We also have by the definitions of $\widetilde{T}_{\gamma_{\widetilde{l}+1}}\ldots \widetilde{T}_{\gamma_{l'}}=T_{\hat{w}}^{-1}{T}_{\gamma_{\widetilde{l}+1}}\ldots {T}_{\gamma_{l'}}T_{\hat{w}}$ and of $T_{w'}=T_{\hat{w}}T_{\bar{w}_2}$ that
$$
T_{\gamma_1}\ldots T_{\gamma_{\widetilde{l}}}\widetilde{T}=T_{\gamma_1}\ldots T_{\gamma_{\widetilde{l}}}T_{\hat{w}}\widetilde{T}_{\gamma_{\widetilde{l}+1}}\ldots \widetilde{T}_{\gamma_{l'}}T_{\bar{w}_2}=T_{\gamma_1}\ldots T_{\gamma_{\widetilde{l}}}
{T}_{\gamma_{\widetilde{l}+1}}\ldots {T}_{\gamma_{l'}}T_{\hat{w}}T_{\bar{w}_2}=
$$
$$
=T_{\gamma_1}\ldots T_{\gamma_{\widetilde{l}}}
{T}_{\gamma_{\widetilde{l}+1}}\ldots {T}_{\gamma_{l'}}T_{w'}.
$$
Therefore we can rewrite (\ref{a5'}) as follows
\begin{equation}\label{a6}
(u,\tau((X_{\gamma_1}^+)^{(m_1)}\ldots (X_{\gamma_{\widetilde{l}}}^+)^{(m_{\widetilde{l}})}({X}_{\gamma_{\widetilde{l}+1}}^+)^{(m_{\widetilde{l}+1})}\ldots ({X}_{\gamma_{l'}}^+)^{(m_{l'})})v)=(u,\tau(T_{\gamma_1}\ldots T_{\gamma_{\widetilde{l}}}{T}_{\gamma_{\widetilde{l}+1}}\ldots {T}_{\gamma_{l'}}T_{w'} X'T_{w'}^{-1})v).
\end{equation}

Now recall that by Corollary \ref{segmq} (ii) $X'\in U_q^{res}({\frak b}_+)=U_q^{res}(\n_+)U_q^{res}(H)$. Observe also that $\tau(T_{w'}^{-1})v\in V$ is a weight vector, and by Lemma \ref{wHact} elements of the subalgebra $U_q^{res}(H)$, which is preserved by the action of $\tau$,  act on it by multiplication by elements of $\mathbb{C}[q,q^{-1}]$. Therefore (\ref{a6}) immediately implies (\ref{prepm}).

Finally note that $\tau(T_{w'}^{-1})v$ and $\tau(T_{w'}^{-1}T_{\gamma_{l'}}\ldots T_{\gamma_{1}})u$ are weight vectors, and different weight subspaces of $V$ are orthogonal with respect to the contravariant form $(~\cdot~,~\cdot~)$. Therefore we can assume also that $X(m_1,\ldots, m_{l'})$ in (\ref{prepm}) belongs to a weight subspace of $U_q^{res}(\g)$, so that the weight of $\tau(T_{w'}^{-1}T_{\gamma_{l'}}\ldots T_{\gamma_{1}})u$ is equal to that of $\tau(X(m_1,\ldots, m_{l'})T_{w'}^{-1})v$.
This completes the proof.

\end{proof}

\begin{lemma}\label{mainl}
Let $V$ be a finite rank representation of $U_h(\g)$, $u,v\in V$ weight vectors. Suppose that $u$ is a highest weight vector.
Then for any $m_1,\ldots, m_{l'}\in \mathbb{N}$
$$
(u,\tau((X_{\gamma_1}^+)^{(m_1)}\ldots (X_{\gamma_{\widetilde{l}}}^+)^{(m_{\widetilde{l}})}({X}_{\gamma_{\widetilde{l}+1}}^+)^{(m_{\widetilde{l}+1})}\ldots ({X}_{\gamma_{l'}}^+)^{(m_{l'})})v)=
$$
$$
=(u,\tau(T_{\gamma_1}\ldots T_{\gamma_{\widetilde{l}}}
{T}_{\gamma_{\widetilde{l}+1}}\ldots {T}_{\gamma_{l'}}Y) v)
=\check{c}(u,\tau(T_s Y)v), 
$$
where 
$$
Y=\sum_{r_{1},\ldots , r_{D}\in \mathbb{N}}F_{r_{1},\ldots ,r_{D}}({X}^-_{\beta_{k_{l'}}})^{(r_{k_{l'}})}\ldots ({X}^-_{\beta_1})^{(r_{1})}(X^+_{\beta_{k_{l'}+1}})^{(r_{k_{l'}+1})}\ldots (X^+_{\beta_D})^{(r_{D})}\in U_q^{res}(w'(\b_+)),
$$
$F_{r_{1},\ldots ,r_{D}}\in \mathbb{C}[q,q^{-1}]$, the sum is finite, and $\check{c}\in \mathbb{C}[q,q^{-1}]^*$ is an integer power of $q$ up to a numeric factor.
\end{lemma}

\begin{proof}

We bring the right hand side of formula (\ref{prepm}) to the form stated in this lemma.
Firstly, using part (iv) of Lemma \ref{wN} we complete the final part $s_{i_{k_{l'}}}\ldots s_{i_1}$ of reduced decomposition (\ref{wbarinv}) to a reduced decomposition $\overline{w}=s_{i_{k_{l'}}}\ldots s_{i_1}s_{p_{k_{l'}+1}}\ldots s_{p_D}$, and consider the corresponding normal ordering $\beta_{k_{l'}}',\ldots , \beta_1',\beta_{k_{l'}+1}', \ldots , \beta _D'$ \index[not]{b@$\beta_i'$}   of $\Delta_+$. Define the corresponding quantum root vectors
$$
\overline{X}_{\beta_q'}^+=\left\{ \begin{array}{ll} T^{-1}_{i_{k_{l'}}}\ldots T^{-1}_{i_{q-1}}X_{i_q}^+ & 1\leq q \leq k_{l'} \\ T_{w'}^{-1}T_{p_{k_{l'}+1}}^{-1}\ldots T_{p_{q-1}}^{-1}X_{p_q}^+ & k_{l'}+1\leq q \leq D
\end{array}
\right. , \index[not]{X@$\overline{X}_{\beta_q'}^+$}
$$
and the basis of $U_q^{res}(\n_+)$ as in Lemma \ref{segmPBW} (v),
\begin{equation}\label{b1}
(\overline{X}^+_{\beta'_{k_{l'}}})^{(r_{k_{l'}})}\ldots (\overline{X}^+_{\beta'_1})^{(r_{1})}(\overline{X}^+_{\beta'_{k_{l'}+1}})^{(r_{k_{l'}+1})}\ldots (\overline{X}^+_{\beta'_D})^{(r_{D})}, r_i\in \mathbb{N},i=1,\ldots, D.	
\end{equation}

Similarly, we also complete the final segment $s_{i_{k_{l'}+1}}\ldots s_{i_D}$ of reduced decomposition (\ref{wbar}) to another reduced decomposition $\overline{w}=s_{i_{k_{l'}+1}}\ldots s_{i_D}s_{p_{1}}\ldots s_{p_{k_{l'}}}$, and consider the corresponding normal ordering $\beta_{k_{l'}+1}'',\ldots , \beta_D'',\beta_1'', \ldots , \beta _{k_{l'}}''$ of $\Delta_+$. \index[not]{b@$\beta_i''$} Define the corresponding quantum root vectors
$$
X_{\beta_q''}^+=\left\{ \begin{array}{ll} T_{w''}T_{p_1}\ldots T_{p_{q-1}}X_{p_q}^+ & 1\leq q \leq k_{l'} \\ T_{i_{k_{l'}+1}}\ldots T_{i_{q-1}}X_{i_q}^+ & k_{l'}+1\leq q \leq D
\end{array}
\right. , \index[not]{X@$X_{\beta_q''}^+$}
$$
where $T_{w''}=T_{i_{k_{l'}+1}}\ldots T_{i_D}$ for $w'':=s_{i_{k_{l'}+1}}\ldots s_{i_D}$, \index[not]{w@$w''$}
and the basis of $U_q^{res}(\n_+)$ as in Lemma \ref{segmPBW} (i),
\begin{equation}\label{b2}
(X^+_{\beta''_{k_{l'}+1}})^{(r_{k_{l'}+1})}\ldots (X^+_{\beta''_D})^{(r_{D})}(X^+_{\beta''_{1}})^{(r_{1})}\ldots (X^+_{\beta''_{k_{l'}}})^{(r_{k_{l'}})}, r_i\in \mathbb{N},i=1,\ldots, D.	
\end{equation}

We can represent the elements $(\overline{X}^+_{\beta'_{k_{l'}+1}})^{(r_{k_{l'}+1})}\ldots (\overline{X}^+_{\beta'_D})^{(r_{D})}$ using basis (\ref{b2}),
\begin{equation}\label{bdec}
(\overline{X}^+_{\beta'_{k_{l'}+1}})^{(r_{k_{l'}+1})}\ldots (\overline{X}^+_{\beta'_D})^{(r_{D})}=
\end{equation}
$$
=\sum_{j_{1}, \ldots, j_D\in \mathbb{N}}C_{j_{1}, \ldots j_D}(X^+_{\beta''_{k_{l'}+1}})^{(j_{k_{l'}+1})}\ldots (X^+_{\beta''_D})^{(j_{D})}(X^+_{\beta''_{1}})^{(j_{1})}\ldots (X^+_{\beta''_{k_{l'}}})^{(j_{k_{l'}})},	
$$
where $C_{j_{1}, \ldots j_D}\in \mathbb{C}[q,q^{-1}]$.
Applying $T_{w'}$ to this identity we get
\begin{equation}\label{bdec1}
(\widehat{\overline{X}}^+_{\beta'_{k_{l'}+1}})^{(r_{k_{l'}+1})}\ldots (\widehat{\overline{X}}^+_{\beta'_D})^{(r_{D})}=
\end{equation}
$$
=\sum_{j_{1}, \ldots, j_D\in \mathbb{N}}C_{j_{1}, \ldots j_D}(X^+_{\beta_{k_{l'}+1}})^{(j_{k_{l'}+1})}\ldots (X^+_{\beta_D})^{(j_{D})}(\check{X}^+_{\beta''_{1}})^{(j_{1})}\ldots (\check{X}^+_{\beta''_{k_{l'}}})^{(j_{k_{l'}})},	
$$
where by Lemma \ref{segmPBW} (v) for $k_{l'}+1\leq q \leq D$
\begin{equation}\label{bq'}
\widehat{\overline{X}}_{\beta_q'}^+:=T_{w'}\overline{X}_{\beta_q'}^+=T_{p_{k_{l'}+1}}^{-1}\ldots T_{p_{q-1}}^{-1}X_{p_q}^+\in U_q^{res}(\n_+), \index[not]{X@$\widehat{\overline{X}}_{\beta_q'}^+$}	
\end{equation}
and for $1\leq q \leq k_{l'}$
\begin{equation}\label{bq''}
\check{X}_{\beta_q''}^+:=T_{w'}X_{\beta_q''}^+=T_{w'}T_{w''}T_{p_1}\ldots T_{p_{q-1}}X_{p_q}^+=T_{\overline{w}}T_{p_1}\ldots T_{p_{q-1}}X_{p_q}^+\in U_q^{res}(\b_-), \index[not]{X@$\check{X}_{\beta_q''}^+$}	
\end{equation}
and we used the fact that $T_{p_1}\ldots T_{p_{q-1}}X_{p_q}^+\in U_q^{res}(\n_+)$ by Lemma \ref{segmPBW} (i) and $T_{\overline{w}}U_q^{res}(\n_+)\subset U_q^{res}(\b_-)$ by Lemma \ref{wlong}.
 
Observe that for $1\leq q \leq k_{l'}$, $s_{p_1}\ldots s_{p_{q-1}}\alpha_{p_q}\in \Delta_+$ by Lemma \ref{wN} (iv) as $s_{p_1}\ldots s_{p_{k_{l'}}}$ is a reduced decomposition being a part of the reduced decomposition $\overline{w}=s_{i_{k_{l'}+1}}\ldots s_{i_D}s_{p_1}\ldots s_{p_{k_{l'}}}$. Therefore $\overline{w}s_{p_1}\ldots s_{p_{q-1}}\alpha_{p_q}\in \Delta_-$. 

Thus the elements $\check{X}_{\beta_q''}^+\in U_q^{res}(\b_-)$, $1\leq q \leq k_{l'}$ have strictly negative weights, ${\rm wt}(\check{X}_{\beta_q''}^+)=\overline{w}s_{p_1}\ldots s_{p_{q-1}}\alpha_{p_q}<0$. Now by (\ref{bq'}) and (\ref{bq''}) we have that in (\ref{bdec1}) 
$$
(\widehat{\overline{X}}^+_{\beta'_{k_{l'}+1}})^{(r_{k_{l'}+1})}\ldots (\widehat{\overline{X}}^+_{\beta'_D})^{(r_{D})}\in U_q^{res}(\n_+),
$$
$$ 
(X^+_{\beta_{k_{l'}+1}})^{(j_{k_{l'}+1})}\ldots (X^+_{\beta_D})^{(j_{D})}\in U_q^{res}(\n_+),
$$
and the elements
$$
(\check{X}^+_{\beta''_{1}})^{(j_{1})}\ldots (\check{X}^+_{\beta''_{k_{l'}}})^{(j_{k_{l'}})}\in U_q^{res}(\b_-)=U_q^{res}(H)U_q^{res}(\n_-)
$$
have strictly negative weights unless $j_1=\ldots=j_{k_{l'}}=0$. 

By the uniqueness of the Poincar\'{e}--Birkhoff--Witt decomposition (\ref{resPBW}) in Lemma \ref{segmPBW} we deduce that (\ref{bdec1}) takes the form
$$
(\widehat{\overline{X}}^+_{\beta'_{k_{l'}+1}})^{(r_{k_{l'}+1})}\ldots (\widehat{\overline{X}}^+_{\beta'_D})^{(r_{D})}=
$$
$$
=\sum_{j_{k_{l'}+1}, \ldots, j_D\in \mathbb{N}}C_{j_{k_{l'}+1}, \ldots j_D}(X^+_{\beta_{k_{l'}+1}})^{(j_{k_{l'}+1})}\ldots (X^+_{\beta_D})^{(j_{D})},	
$$
where $C_{j_{k_{l'}+1}, \ldots, j_D}\in \mathbb{C}[q,q^{-1}]$.

Applying $T_{w'}^{-1}$ to this identity we obtain that (\ref{bdec}) takes the form
\begin{equation}\label{bdec2}
(\overline{X}^+_{\beta'_{k_{l'}+1}})^{(r_{k_{l'}+1})}\ldots (\overline{X}^+_{\beta'_D})^{(r_{D})}
=\sum_{j_{k_{l'}+1}, \ldots, j_D\in \mathbb{N}}C_{j_{k_{l'}+1}, \ldots j_D}(X^+_{\beta''_{k_{l'}+1}})^{(j_{k_{l'}+1})}\ldots (X^+_{\beta''_D})^{(j_{D})}.
\end{equation}

Recalling basis (\ref{b1}) we infer that every element of $U_q^{res}(\n_+)$ is a $\mathbb{C}[q,q^{-1}]$--linear combination of elements of the form
\begin{equation}\label{b3}
(\overline{X}^+_{\beta'_{k_{l'}}})^{(r_{k_{l'}})}\ldots (\overline{X}^+_{\beta'_1})^{(r_{1})}(X^+_{\beta''_{k_{l'}+1}})^{(r_{k_{l'}+1})}\ldots (X^+_{\beta''_D})^{(r_{D})}, r_i\in \mathbb{N},i=1,\ldots, D.	
\end{equation}
Kostant's formula \index{Kostant's formula} shows that they form a $\mathbb{C}[q,q^{-1}]$--basis of $U_q^{res}(\n_+)$.

Now for $1\leq q\leq k_{l'}$ we have
\begin{equation}\label{w'1}
T_{w'}\overline{X}_{\beta_q'}^+=T_{w'}T^{-1}_{i_{k_{l'}}}\ldots T^{-1}_{i_{q-1}}X_{i_q}^+=T_{i_1}\ldots T_{i_q}X_{i_q}^+=-T_{i_1}\ldots T_{i_{q-1}}X_{i_q}^-K_{i_q}=-X_{\beta_q}^-K_{\beta_q},
\end{equation}
where
$$
K_{\beta_q}=T_{i_1}\ldots T_{i_{q-1}}K_{i_q},
$$
and for $k_{l'}+1\leq q \leq D$
\begin{equation}\label{w'2}
T_{w'}X_{\beta_q''}^+=T_{w'}T_{i_{k_{l'}+1}}\ldots T_{i_{q-1}}X_{i_q}^+=X_{\beta_q}^+.
\end{equation}
The last two identities and commutation relations between elements $K_{\beta_q}$ and quantum root vectors following from relations (\ref{Ccomm}) imply
\begin{equation}\label{F'}
T_{w'}\left((\overline{X}^+_{\beta'_{k_{l'}}})^{(r_{k_{l'}})}\ldots (\overline{X}^+_{\beta'_1})^{(r_{1})}(X^+_{\beta''_{k_{l'}+1}})^{(r_{k_{l'}+1})}\ldots (X^+_{\beta''_D})^{(r_{D})}\right)=
\end{equation}
$$
=Q_{r_{1},\ldots ,r_{k_{l'}}}({X}^-_{\beta_{k_{l'}}})^{(r_{k_{l'}})}\ldots ({X}^-_{\beta_1})^{(r_{1})}(X^+_{\beta_{k_{l'}+1}})^{(r_{k_{l'}+1})}\ldots (X^+_{\beta_D})^{(r_{D})},
$$
where $Q_{r_{1},\ldots ,r_{k_{l'}}}$ is a monomial in $K_1^{\pm 1},\ldots, K_l^{\pm 1}$.

Recalling basis (\ref{b3}) we can express $X\in U_q^{res}(\n_+)$ in (\ref{prepm}) as follows
$$
X=\sum_{r_{1},\ldots ,r_{D}\in \mathbb{N}}\overline{F}_{r_{1},\ldots ,r_{D}}(\overline{X}^+_{\beta'_{k_{l'}}})^{(r_{k_{l'}})}\ldots (\overline{X}^+_{\beta'_1})^{(r_{1})}(X^+_{\beta''_{k_{l'}+1}})^{(r_{k_{l'}+1})}\ldots (X^+_{\beta''_D})^{(r_{D})},
$$
where $\overline{F}_{r_{1},\ldots ,r_{D}}\in \mathbb{C}[q,q^{-1}]$ and the sum is finite.

Now using (\ref{F'}) we deduce that in (\ref{prepm})
\begin{equation}\label{XTw}
T_{w'}XT_{w'}^{-1}=\sum_{r_{1},\ldots ,r_{D}\in\mathbb{N}}F'_{r_{1},\ldots ,r_{D}\in\mathbb{N}}({X}^-_{\beta_{k_{l'}}})^{(r_{k_{l'}})}\ldots ({X}^-_{\beta_1})^{(r_{1})}(X^+_{\beta_{k_{l'}+1}})^{(r_{k_{l'}+1})}\ldots (X^+_{\beta_D})^{(r_{D})},
\end{equation}
where $F'_{r_{1},\ldots ,r_{D}}=\overline{F}_{r_{1},\ldots ,r_{D}}Q_{r_{1},\ldots ,r_{k_{l'}}}\in U_q(H)\subset U_q^{res}(H)$.

Observing that  $\tau({T}_{\gamma_{l'}}\ldots {T}_{\gamma_{\widetilde{l}+1}}T_{\gamma_{\widetilde{l}}}\ldots T_{\gamma_1})u$ is a weight vector and using Lemma \ref{wHact}, formula (\ref{prepm}) can be rewritten as follows
\begin{equation}\label{a5}
(u,\tau((X_{\gamma_1}^+)^{(m_1)}\ldots (X_{\gamma_{\widetilde{l}}}^+)^{(m_{\widetilde{l}})}({X}_{\gamma_{\widetilde{l}+1}}^+)^{(m_{\widetilde{l}+1})}\ldots ({X}_{\gamma_{l'}}^+)^{(m_{l'})})v)=
\end{equation}
$$
=\sum_{r_{1},\ldots ,r_{D}\in\mathbb{N}}\tau(F_{r_{1},\ldots ,r_{D}})(u,\tau(T_{\gamma_1}\ldots T_{\gamma_{\widetilde{l}}}
{T}_{\gamma_{\widetilde{l}+1}}\ldots {T}_{\gamma_{l'}}({X}^-_{\beta_{k_{l'}}})^{(r_{k_{l'}})}\ldots ({X}^-_{\beta_1})^{(r_{1})}(X^+_{\beta_{k_{l'}+1}})^{(r_{k_{l'}+1})}\ldots (X^+_{\beta_D})^{(r_{D})})v)=
$$
$$
=(u,\tau(T_{\gamma_1}\ldots T_{\gamma_{\widetilde{l}}}
{T}_{\gamma_{\widetilde{l}+1}}\ldots {T}_{\gamma_{l'}} Y) v), 
$$
where 
$$
Y=\sum_{r_{1},\ldots ,r_{D}\in\mathbb{N}}F_{r_{1},\ldots ,r_{D}}({X}^-_{\beta_{k_{l'}}})^{(r_{k_{l'}})}\ldots ({X}^-_{\beta_1})^{(r_{1})}(X^+_{\beta_{k_{l'}+1}})^{(r_{k_{l'}+1})}\ldots (X^+_{\beta_D})^{(r_{D})}\in U_q^{res}(w'(\b_+)),
$$
$F_{r_{1},\ldots ,r_{D}}\in \mathbb{C}[q,q^{-1}]$ and the sum is finite. This proves the first formula in the statement of this lemma.

To justify the last formula in the statement of the lemma we observe that 
${T}_{\gamma_{l'}}\ldots {T}_{\gamma_{\widetilde{l}+1}}T_{\gamma_{\widetilde{l}}} \ldots T_{\gamma_1}$ and $T_s^{-1}$  
act as the same transformations of $\h\subset U_h(\h)$ and apply Lemma \ref{bga}. This completes the proof.


\end{proof}

In the course of the proof of the previous lemma we obtained the following result.
\begin{corollary}\label{w'b}
The products
\begin{equation}\label{basw'}
({X}^-_{\beta_{k_{l'}}})^{(r_{k_{l'}})}\ldots ({X}^-_{\beta_1})^{(r_{1})}(X^+_{\beta_{k_{l'}+1}})^{(r_{k_{l'}+1})}\ldots (X^+_{\beta_D})^{(r_{D})}, r_i\in \mathbb{N}, i=1,\ldots, D,
\end{equation}
or 
\begin{equation}\label{basw'1}
(X^+_{\beta_{k_{l'}+1}})^{(r_{k_{l'}+1})}\ldots (X^+_{\beta_D})^{(r_{D})}({X}^-_{\beta_{k_{l'}}})^{(r_{k_{l'}})}\ldots ({X}^-_{\beta_1})^{(r_{1})}, r_i\in \mathbb{N}, i=1,\ldots, D
\end{equation}
form $U_q^{res}(H)$--bases in the subalgebra $U_q^{res}(w'(\b_+))$ of $U_q^{res}(\g)$ generated over $U_q^{res}(H)$ by the elements 
\begin{equation}\label{w'bel}
({X}^-_{\beta_{k_{l'}}})^{(r_{k_{l'}})},\ldots ,({X}^-_{\beta_1})^{(r_{1})},(X^+_{\beta_{k_{l'}+1}})^{(r_{k_{l'}+1})},\ldots ,(X^+_{\beta_D})^{(r_{D})}, r_i\in \mathbb{N}, i=1,\ldots, D.
\end{equation}
\end{corollary}

\begin{proof}
If $Y\in U_q^{res}(w'(\b_+))$ then by (\ref{w'1}) and (\ref{w'2}) we have $T_{w'}^{-1}(Y)\in U_q^{res}(\b_+)$. By (\ref{XTw}) $Y=T_{w'}(T_{w'}^{-1}(Y))$ can be represented as a $U_q^{res}(H)$--linear combination of elements (\ref{basw'}). By Kostant's formula they form a $U_q^{res}(H)$--basis in $U_q^{res}(w'(\b_+))$.

The case of elements (\ref{basw'1}) is considered in a similar way. This completes the proof.

\end{proof}

\begin{remark}\label{w'bc}
The use of the symbol $U_q^{res}(w'(\b_+))$ is motivated by the fact that the algebra $U_q^{res}(w'(\b_+))$ is generated by elements (\ref{w'bel}) which are defined with the help of the quantum root vectors corresponding to the roots $-\beta_{k_{l'}},\ldots ,-\beta_1,\beta_{k_{l'}+1},\ldots ,\beta_D$. They form a system of positive roots $\Delta_+^{w'}$ \index[not]{D@$\Delta_+^{w'}$} in $\Delta$ such that ${w'}^{-1}\Delta_+^{w'}=\Delta_+$. Indeed, recalling (\ref{w'}) we obtain that $\Delta_{{w'}^{-1}}=\{\beta_1,\ldots, \beta_{k_{l'}}\}$ by Lemma \ref{wN} (iv), and hence ${w'}^{-1}\Delta_+^{w'}=\Delta_+$ as the number of the roots in $\Delta_+^{w'}$ is the same as in $\Delta_+$. So, we infer that $\Delta_+^{w'}={w'}\Delta_+$. Thus $w'(\b_+):={\rm Ad}w'(\b_+)$ \index[not]{w@$w'(\b_+)$} is the Borel subalgebra of $\g$ the nilradical of which is generated by the root vectors corresponding to the roots from $\Delta_+^{w'}$, and $U_q^{res}(w'(\b_+))$ is the restricted version of the quantum counterpart of the enveloping algebra $U(w'(\b_+))$.  
\end{remark}


\section{Some auxiliary results on the quantized algebra of regular functions on an algebraic Poisson--Lie group}
\label{CG}

\pagestyle{myheadings}
\markboth{CHAPTER~\thechapter.~ZHELOBENKO TYPE OPERATORS FOR Q-W--ALGEBRAS}{\thesection.~QUANTIZED ALGEBRA OF REGULAR FUNCTIONS ON A POISSON--LIE GROUP}

\setcounter{equation}{0}
\setcounter{theorem}{0}

In this section we give several formulas related to the adjoint action and commutation relations in the algebra $\mathbb{C}_{\mathcal{B}}^s[G]$.

Let $\beta_1,\ldots, \beta_D$ be an arbitrary normal ordering on a positive root system $\Delta_+\subset \Delta$, $f_{\beta_1},\ldots, f_{\beta_D}$, $e_{\beta_1},\ldots, e_{\beta_D}$$\in U_{\mathcal{B}}^{s,res}(\g)$ the corresponding quantum root vectors defined with the help of this normal ordering.

Firstly, following \cite{BG}, Theorem I.8.16 we recall the commutation relations in the algebra $\mathbb{C}_{\mathcal{B}}^s[G]$ which follow from the fact that $U_h^s(\g)$ is quasitriangular. 
Namely, if $V$, $V'$ are finite rank representations of $U_h(\g)$, $(V)_\eta, (V')_\rho$, $(V)_\beta, (V')_\gamma$ their weight subspaces of weights $\eta, \rho, \beta$ and $\gamma$, respectively, and $v\in (V)_\eta, v_1\in (V')_\rho$, $u\in (V)_\beta, u_1\in (V')_\gamma$ then evaluating the identity $(\Delta_s^\tau)^{opp}(x)\mathcal{R}_s^\tau=\mathcal{R}_s^\tau\Delta_s^\tau(x)$ on the matrix element $(u,\cdot ~ v)\otimes (u_1,\cdot ~ v_1)$ and recalling formula (\ref{rmatrspi}) we obtain
$$
q^{-\left\langle (\kappa{1+s \over 1-s }P_{{\h'}}+id)\eta^\vee, \rho^\vee\right\rangle}\left((u_1,\cdot ~ v_1)\otimes (u,\cdot ~ v)+\sum_{\nu\in Q_+,\nu\neq 0}\sum_{\tiny \begin{array}{c} m_1,\ldots, m_D \in \mathbb{N}, \\ m_1\beta_1+\ldots +m_D\beta_D=\nu\end{array}}\hspace{-2em} (u_1,\cdot ~ u_{m_1,\ldots, m_D}v_1)\otimes (u,\cdot ~ u_{-m_1,\ldots, -m_D}v)\right)=
$$
\begin{equation}\label{Gcomm}
=q^{-\left\langle (\kappa{1+s \over 1-s }P_{{\h'}}+id)\beta^\vee, \gamma^\vee\right\rangle}(u,\cdot ~ v)\otimes (u_1,\cdot ~ v_1)+
\end{equation}
$$
+ \sum_{\nu\in Q_+,\nu\neq 0}\sum_{\tiny \begin{array}{c} m_1,\ldots, m_D \in \mathbb{N}, \\ m_1\beta_1+\ldots +m_D\beta_D=\nu\end{array}}\hspace{-2em} q^{-\left\langle (\kappa{1+s \over 1-s }P_{{\h'}}+id)(\beta^\vee+\nu^\vee), \gamma^\vee-\nu^\vee\right\rangle}(\omega(u_{-m_1,\ldots, -m_D})u,\cdot ~ v)\otimes (\omega(u_{m_1,\ldots, m_D})u_1,\cdot ~ v_1),
$$
where
$$
u_{-m_1,\ldots, -m_D}=c_{-m_1,\ldots, -m_D}\tau(f_{\beta_1}^{(m_1)}\ldots f_{\beta_D}^{(m_D)}),
$$
$$
u_{m_1,\ldots, m_D}=c_{m_1,\ldots, m_D}\tau(e_{\beta_1}^{m_1}\ldots e_{\beta_D}^{m_D}),
$$
$c_{\pm m_1,\ldots, \pm m_D}\in \mathcal{B}$,
and similarly evaluating the identity $(\Delta_s^\tau)^{opp}(x){(\mathcal{R}_s^\tau)_{21}}^{-1}={(\mathcal{R}_s^\tau)_{21}}^{-1}\Delta_s^\tau(x)$ on the matrix element $(u,\cdot ~ v)\otimes (u_1,\cdot ~ v_1)$ we get
$$
q^{-\left\langle (\kappa{1+s \over 1-s }P_{{\h'}}-id)\eta^\vee, \rho^\vee\right\rangle}\left((u_1,\cdot ~ v_1)\otimes (u,\cdot ~ v)+\sum_{\nu\in Q_+,\nu\neq 0}\sum_{\tiny \begin{array}{c} m_1,\ldots, m_D \in \mathbb{N}, \\ m_1\beta_1+\ldots +m_D\beta_D=\nu\end{array}} \hspace{-2em} (u_1,\cdot ~ u_{-m_1,\ldots, -m_D}'v_1)\otimes (u,\cdot ~ u_{m_1,\ldots, m_D}'v)\right)=
$$
\begin{equation}\label{Melid}
=q^{-\left\langle (\kappa{1+s \over 1-s }P_{{\h'}}-id)\beta^\vee, \gamma^\vee\right\rangle}(u,\cdot ~ v)\otimes (u_1,\cdot ~ v_1)+
\end{equation}
$$
+\sum_{\nu\in Q_+,\nu\neq 0}\sum_{\tiny \begin{array}{c} m_1,\ldots, m_D \in \mathbb{N}, \\ m_1\beta_1+\ldots +m_D\beta_D=\nu\end{array}}\hspace{-2em} q^{-\left\langle (\kappa{1+s \over 1-s }P_{{\h'}}-id)(\beta^\vee+\nu^\vee), \gamma^\vee-\nu^\vee\right\rangle}(\omega(u_{m_1,\ldots, m_D}')u,\cdot ~ v)\otimes (\omega(u_{-m_1,\ldots, -m_D}')u_1,\cdot ~ v_1),
$$
where
\begin{equation}\label{defu-}
u_{-m_1,\ldots, -m_D}'=c_{-m_1,\ldots, -m_D}'\tau(f_{\beta_D}^{(m_D)}\ldots f_{\beta_1}^{(m_1)}),
\end{equation}
\begin{equation}\label{defu+}
u_{m_1,\ldots, m_D}'=c_{m_1,\ldots, m_D}'\tau(e_{\beta_D}^{m_D}\ldots e_{\beta_1}^{m_1}), 
\end{equation}
$c_{\pm m_1,\ldots, \pm m_D}'\in \mathcal{B}$.

If $v=\tau(T_s^{-1})v_\lambda \in (V)_{s^{-1}\lambda}, v_1=\tau(T_s^{-1})v_\mu\in (V')_{s^{-1}\mu}$, $u\in (V)_\beta, u_1=v_\mu\in (V')_\mu$, where $v_\lambda\in V$ and $v_\mu \in V'$ are highest weight vectors, then, since $v_\lambda$ and $v_\mu$ are highest weight vectors,  the only non--vanishing contributions in the sum over $\nu$  in the left hand side of (\ref{Melid}) may only appear if the weights of $\tau(T_s)(u_{-m_1,\ldots, -m_D}')$ and of $\tau(T_s)(u_{m_1,\ldots, m_D}')$ are both negative which by (\ref{defu-}) and (\ref{defu+}) amounts to  $s(\nu)\in Q_+$ and $s(\nu)\in -Q_+$, respectively. As $\nu\neq 0$ this is impossible, and hence all terms in the sum over $\nu$  in the left hand side of (\ref{Melid}) vanish.

Since $v_\mu$ is a highest weight vector, the definition of $\omega$ and (\ref{defu-}) also imply in this case that the second factor in all terms of the sum over $\nu$ in the right hand side of (\ref{Melid}) vanishes. 

Thus (\ref{Melid}) reduces in this case to
\begin{equation}\label{comm1}
q^{-\left\langle (\kappa{1+s \over 1-s }P_{{\h'}}-id)\lambda^\vee, \mu^\vee\right\rangle}(v_\mu,\cdot ~ \tau(T_s^{-1})v_\mu)\otimes (u,\cdot ~ \tau(T_s^{-1})v_\lambda)=
\end{equation}
$$
=q^{-\left\langle (\kappa{1+s \over 1-s }P_{{\h'}}-id)\beta^\vee, \mu^\vee\right\rangle}(u,\cdot ~ \tau(T_s^{-1})v_\lambda)\otimes (v_\mu,\cdot ~ \tau(T_s^{-1})v_\mu).	
$$

The next lemma shows how the adjoint action behaves with respect to the multiplication $\otimes$ in $\mathbb{C}_{\mathcal{B}}^s[G]$. 
\begin{lemma}
For any $f,g\in \mathbb{C}_{\mathcal{B}}^s[G], x \in U_{\mathcal{B}}^{s,res}(\g)$ we have
$$
{\rm Ad}_s^0x(f\otimes g)(\cdot ~, \cdot ~)=({\rm Ad}_s^0 x^2 f)(\cdot ~)\otimes g((\omega_0 S_s^{-1})(x^1)\cdot ~ \omega_0 x^3),
$$
where $\Delta_s^2x=(\Delta_s\otimes id)\Delta_s x=x^1\otimes x^2\otimes x^3$ in the Sweedler notation.

In particular,
\begin{equation}\label{ad0}
{\rm Ad}_s^0f_\beta (f\otimes g)(~\cdot ~, ~\cdot ~)= f(~\cdot ~)\otimes g((\omega_0 S_s^{-1})(f_\beta)~\cdot ~ )+
({\rm Ad}_s^0 f_\beta f)(~\cdot ~)\otimes g(G_\beta^{-1} ~\cdot ~)+\sum_i ({\rm Ad}_s^0 x_i f)(~\cdot ~)\otimes g((\omega_0 S_s^{-1})(y_i)~\cdot ~ )+	
\end{equation}
$$
+f(G_\beta^{-1} ~\cdot ~ G_\beta)\otimes g(G_\beta^{-1} ~\cdot ~ \omega_0(f_\beta))+\sum_i ({\rm Ad}_s^0y_i^2 f)(~\cdot ~)\otimes g((\omega_0 S_s^{-1})(y_i^1)~\cdot ~ \omega_0 x_i)=
$$
$$
=f(~\cdot ~)\otimes g((\omega_0 S_s^{-1})(f_\beta)~\cdot ~ )+
({\rm Ad}_s^0 f_\beta f)(~\cdot ~)\otimes g(G_\beta^{-1} ~\cdot ~)+f(G_\beta^{-1} ~\cdot ~ G_\beta)\otimes g(G_\beta^{-1} ~\cdot ~ \omega_0(f_\beta))+
$$
$$
+\sum_i ({\rm Ad}_s^0y_i f)(~\cdot ~)\otimes g(G_\beta^{-1} ~\cdot ~ \omega_0 x_i)+\sum_i ({\rm Ad}_s^0 x_i^1 f)(~\cdot ~)\otimes g((\omega_0 S_s^{-1})(y_i)~\cdot ~ \omega_0 x_i^2)
$$
and
$$
{\rm Ad}_s^0f_\beta^{(n)} (f\otimes g)(~\cdot ~, ~\cdot ~)=
\sum_{k=0}^{n}\sum_{p=0}^{n-k}q_{\beta}^{k(n-k)+p(n-k-p)}({\rm Ad}_s^0 (G_\beta^{-k}f_\beta^{(p)}) f)(~\cdot ~)\otimes g(\omega_0S_s^{-1}(G_\beta^{-k-p}f_\beta^{(n-k-p)}) ~\cdot ~ \omega_0(f_\beta^{(k)}))+	
$$
\begin{equation}\label{ad0n}
+\sum_{k=0}^{n-1}\sum_i q_\beta^{k(n-k)}({\rm Ad}_s^0 (G_\beta^{-k} x_i^{(n-k)}) f)(~\cdot ~)\otimes g((\omega_0 S_s^{-1})(G_\beta^{-k}y_i^{(n-k)})~\cdot ~ \omega_0(f_\beta^{(k)}))+
\end{equation}
$$
+\sum_i ({\rm Ad}_s^0({y_i^{(n)}}^2) f)(~\cdot ~)\otimes g((\omega_0 S_s^{-1})({y_i^{(n)}}^1)~\cdot ~ \omega_0 (x_i^{(n)})),
$$
where $G_\beta$, $x_i$, $y_i$, $x_i^{(p)}$, $y_i^{(p)}$ are defined in (\ref{comults}) and (\ref{comultsn}), and
$\Delta_s(x_i)=x_i^1\otimes x_i^2$, $\Delta_s(y_i)=y_i^1\otimes y_i^2$, $\Delta_s(y_i^{(p)})={y_i^{(p)}}^1\otimes {y_i^{(p)}}^2$ in the Sweedler notation.
\end{lemma}

\begin{proof}
Denote, using the Sweedler notation,
$$
\Delta_s^3x=(\Delta_s\otimes id\otimes id)(\Delta_s\otimes id)\Delta_s x=x^1\otimes x^2\otimes x^3\otimes x^4
$$
and observe that the definition of ${\rm Ad}_s'$ implies that for any $x,z \in U_{\mathcal{B}}^{s,res}(\g)$
$$
\Delta_s^{opp}({\rm Ad}'_sx z)=(x^2\otimes x^1)(z^2\otimes z^1)(S_sx^3\otimes S_s x^4)={\rm Ad}_s'x^2 z^2\otimes x^1z^1S_sx^3.
$$

Let $z=S_s\omega_0 y$, $y\in U_{\mathcal{B}}^{s,res}(\g)$. Then, since $\omega_0 S_s^{-1}$ is an algebra homomorphism and a coalgebra anti-homomorphism, we deduce
$$
\Delta_s^\tau\omega_0 S_s^{-1}({\rm Ad}'_sx S_s\omega_0 y)=(\omega_0 S_s^{-1}\otimes \omega_0 S_s^{-1})\Delta_s^{opp}({\rm Ad}_s'x z)=(\omega_0 S_s^{-1}\otimes \omega_0 S_s^{-1}){\rm Ad}'_sx^2 z^2\otimes x^1z^1S_sx^3=
$$
$$
=(\omega_0 S_s^{-1}\otimes \omega_0 S_s^{-1}){\rm Ad}'_sx^2 (S_s\omega_0) (y^1)\otimes x^1 (S_s\omega_0) (y^2)S_sx^3=\omega_0 S_s^{-1}{\rm Ad}'_sx^2 ((S_s\omega_0) (y^1))\otimes (\omega_0 S_s^{-1})(x^1)y^2\omega_0 x^3,
$$
where $\Delta_s(z)=z^1\otimes z^2$ and $\Delta_s^\tau(y)=y^1\otimes y^2$.
Evaluating the last identity on $f\otimes g$ we get the first formula in the statement of the lemma. (\ref{ad0}) and (\ref{ad0n}) are obtained from it using (\ref{comults}) and (\ref{comultsn}).

\end{proof}

It will be also convenient to consider a left action ${\rm Ad}^s$ of $U_{\mathcal{B}}^{s,res}(\g)$ on $\mathbb{C}_{\mathcal{B}}^s[G]$ defined as follows
\begin{equation}\label{ad1}
({\rm Ad}^s xf)(w)=f(\tau({\rm Ad}_sx\tau^{-1}(w))). \index[not]{A@${\rm Ad}^s$} 
\end{equation}
It is related to ${\rm Ad}_s^0$ by the formula 
\begin{equation}\label{Ads0}
{\rm Ad}^s xf={\rm Ad}_s^0 \omega_0'(x)f.
\end{equation}

To define quantum analogues of the ideals $J^{jk}$ we need to introduce quantum counterparts of matrix elements (\ref{vpkcl}). Proposition 8.3 in \cite{GY} implies that there exist integral dominant weights $\mu_p$, $p=1,\ldots ,D$, \index[not]{m@$\mu_p$} and elements $v_p\in (V_{\mu_p})_{\mu_p-\beta_p}$ \index[not]{v@$v_p$} such that $(v_p,\cdot~ v_{\mu_p})\in \mathbb{C}_{\mathcal{B}}^s[G]$, and
\begin{equation}\label{defvm}
(v_p,\omega_0 S_s^{-1}(f_{\beta_{D}}^{(m_D)}\ldots f_{\beta_1}^{(m_1)})v_{\mu_p})=\left\{ \begin{array}{ll} 1 & {\rm if}~f_{\beta_{D}}^{(m_D)}\ldots f_{\beta_1}^{(m_1)}= f_{\beta_{p}}\\ 0 & {\rm otherwise}
\end{array}
\right. ,
\end{equation}
where the highest weight vectors $v_{\mu_p}\in V_{\mu_p}$ are normalized by the condition $(v_{\mu_p},v_{\mu_p})=1$.

Since $\omega(\omega_0 S_s^{-1}(f_{\beta_{p}}))v_p\in (V_{\mu_p})_{\mu_p}$ and the subspace $(V_{\mu_p})_{\mu_p}$ has rank one this definition implies 
\begin{equation}\label{oof}
\omega(\omega_0 S_s^{-1}(f_{\beta_{p}}))v_p=v_{\mu_p}.
\end{equation} 
In particular, $(v_{\mu_p},\cdot~ v_{\mu_p})\in \mathbb{C}_{\mathcal{B}}^s[G]$, and 
\begin{equation}\label{Ap}
(v_p,\omega_0 S_s^{-1}(f_{\beta_{p}}) \cdot~ \tau(T_s^{-1})v_{\mu_p})=(v_{\mu_p}, \cdot~ \tau(T_s^{-1})v_{\mu_p}).	
\end{equation}

Let 
\begin{equation}\label{Apquant}
A_p=(v_p, \cdot~ \tau(T_s^{-1})v_{\mu_p}), A_p^0=(v_{\mu_p}, \cdot~ \tau(T_s^{-1})v_{\mu_p}). \index[not]{A@$A_p$} \index[not]{A@$A_p^0$}
\end{equation}

Denote $\Delta_\mu^s(~\cdot~)=(v_\mu, \cdot~ \tau(T_s^{-1})v_\mu)\in \mathbb{C}_{\mathcal{B}}^s[G]$, $\mu\in P_+$, where $v_\mu\in V_\mu$ are highest weight vectors normalized by $(v_\mu,v_\mu)=1$. 

If $v_\lambda\in V_\lambda$ and $v_\mu \in V_\mu$ are highest weight vectors, $u\in (V_\lambda)_\beta$, then by (\ref{comm1})
\begin{equation}\label{comm1'}
\Delta_\mu^s\otimes (u,\cdot ~ \tau(T_s^{-1})v_\lambda)=q^{-\left\langle (\kappa{1+s \over 1-s }P_{{\h'}}-id)(\beta^\vee-\lambda^\vee), \mu^\vee\right\rangle}(u,\cdot ~ \tau(T_s^{-1})v_\lambda)\otimes \Delta_\mu^s.
\end{equation}

For $\lambda=\mu_p$ and $u=v_p$ this yields
\begin{equation}\label{comm1''}
\Delta_\mu^s\otimes A_p=q^{\left\langle (\kappa{1+s \over 1-s }P_{{\h'}}-id)\beta_p^\vee, \mu^\vee\right\rangle}A_p\otimes \Delta_\mu^s=q^{-\left\langle (-2K_s+id)\beta_p^\vee, \mu^\vee\right\rangle}A_p\otimes \Delta_\mu^s.	
\end{equation}

From (\ref{comm1''}) with $\mu=\mu_p$ we also deduce
$$
(v_{\mu_p},\cdot ~\tau(T_s^{-1})v_{\mu_p})\otimes (v_p,\cdot~ \tau(T_s^{-1})v_{\mu_p})
=q^{-\left\langle (-2K_s+id)\beta_p^\vee, \mu_p^\vee\right\rangle}(v_p,\cdot~ \tau(T_s^{-1})v_{\mu_p})\otimes (v_{\mu_p},\cdot~ \tau(T_s^{-1})v_{\mu_p}),
$$
i.e.
\begin{equation}\label{comm2}
A_p^0\otimes A_p=q^{-\left\langle (-2K_s+id)\beta_p^\vee, \mu_p^\vee\right\rangle} A_p\otimes A_p^0.
\end{equation}

The following statement is a quantum counterpart of Lemma \ref{vancl}. In fact Lemma \ref{vancl} can be proved by specializing the statement of Lemma \ref{Apvan} at $q^{\frac{1}{d{\bar{r}}^2}}=1$.
\begin{lemma}\label{Apvan}
(i) For any $1\leq q<p\leq D$, $y\in \omega_0 S_s^{-1}(U_{\mathcal{B}}^{s,res}([-\beta_p,-\beta_D]))$ of the form $y=\omega_0 S_s^{-1}(f_{\beta_{D}}^{(m_D)}\ldots f_{\beta_p}^{(m_p)})$, $m_i\in \mathbb{N}$, $y\neq 1$, and any $u\in U_{\mathcal{B}}^{s,res}(\g)$ one has 
\begin{equation}\label{aq1ll}
(v_q, yuv_{\mu_q})=0.
\end{equation}

In particular, $\omega(y)v_q=0$ in this case.

(ii) For any $1\leq p\leq D$, $y\in \omega_0 S_s^{-1}(U_{\mathcal{B}}^{s,res}([-\beta_p,-\beta_D]))$ of the form $y=\omega_0 S_s^{-1}(f_{\beta_{D}}^{(m_D)}\ldots f_{\beta_p}^{(m_p)})$, $m_i\in \mathbb{N}$, $y\neq 1, \omega_0 S_s^{-1}(f_{\beta_p})$, and any $u\in U_{\mathcal{B}}^{s,res}(\g)$ one has 
\begin{equation}\label{aq1lll}
(v_p, yuv_{\mu_p})=0.
\end{equation}

In particular, $\omega(y)v_p=0$ in this case.

\end{lemma}

\begin{proof}
(i) Observe that, by the definitions of $\omega_0$ and of $S_s$, $\omega_0 S_s^{-1}$ is an algebra automorphism of $U_\mathcal{B}^{s,res}(\g)$ such that $\omega_0 S_s^{-1}(U_\mathcal{B}^{s,res}(\n_{\pm}))\subset U_\mathcal{B}^{s,res}(\b_{\pm})$ and $\omega_0 S_s^{-1}(U_\mathcal{B}^{s,res}(\h))= U_\mathcal{B}^{s,res}(\h)$. Therefore, since $v_{\mu_p}\in V_{\mu_p}$ is a highest weight vector generating $V_{\mu_p}$, by Lemmas \ref{segmPBWs} (iii), \ref{hwv} (ii) and \ref{wHact} we have 
\begin{equation}\label{osgen}
V_{\mu_p}^{res}=\omega_0 S_s^{-1}(U_\mathcal{B}^{s,res}(\g))v_{\mu_p}=\omega_0 S_s^{-1}(U_\mathcal{B}^{s,res}(\n_-))\omega_0 S_s^{-1}(U_\mathcal{B}^{s,res}(\h))\omega_0 S_s^{-1}(U_\mathcal{B}^{s,res}(\n_+))v_{\mu_p}=
\end{equation}
$$
=\omega_0 S_s^{-1}(U_\mathcal{B}^{s,res}(\n_-))U_\mathcal{B}^{s,res}(\h)U_\mathcal{B}^{s,res}(\b_+)v_{\mu_p}=\omega_0 S_s^{-1}(U_\mathcal{B}^{s,res}(\n_-))v_{\mu_p}.
$$
Thus we can write the left hand side of (\ref{aq1ll}) in the form
\begin{equation}\label{aq1}
(v_q, yuv_{\mu_q})=(v_q, y\omega_0 S_s^{-1}(y')v_{\mu_q})
\end{equation}
for some $y'\in U_\mathcal{B}^{s,res}(\n_-)$.  We claim that that the right hand side of (\ref{aq1})vanishes.  

Indeed, $y\omega_0 S_s^{-1}(y')=\omega_0 S_s^{-1}(f_{\beta_{D}}^{(m_D)}\ldots f_{\beta_p}^{(m_p)}y')$ and $f_{\beta_{D}}^{(m_D)}\ldots f_{\beta_p}^{(m_p)}y'$ belongs to the right ideal $\bar{Y}_p^{res}$ in $U_\mathcal{B}^{s,res}(\n_-)$ generated by $f_{\beta_u}^{(m_u)}$, $u\geq p>q$ and $m_u>0$. By Lemma \ref{segmPBWs} (vii) the elements $f_{\beta_{D}}^{(m_D)}\ldots f_{\beta_1}^{(m_1)}$ with at least one $m_u>0$ for $u\geq p>q$ form a linear basis of this ideal, and this basis does not contain multiples of $f_{\beta_q}$. So by (\ref {defvm}) we have $$(v_q, \omega_0 S_s^{-1}(f_{\beta_{D}}^{(m_D)}\ldots f_{\beta_p}^{(m_p)}y')v_{\mu_q})=0.$$ 

The result that we proved implies that $\omega(y)v_q$ is orthogonal, with respect to the contravariant form, to any element of the form $uv_{\mu_q}$, $u\in U_\mathcal{B}^{s,res}(\g)$. Since any element of $V_{\mu_q}^{res}$ is of this form and the contravariant form is non--degenerate we deduce that $\omega(y)v_q=0$.

(ii) By (\ref{osgen}) we can write, similarly to (\ref{aq1}),
\begin{equation}\label{ap2}
(v_p, yuv_{\mu_p})=(v_p, y\omega_0 S_s^{-1}(y')v_{\mu_p})=(v_p, \omega_0 S_s^{-1}(f_{\beta_{D}}^{(m_D)}\ldots f_{\beta_p}^{(m_p)}y')v_{\mu_p})
\end{equation}
for some $y'\in U_\mathcal{B}^{s,res}(\n_-)$. 

If at least one $m_i>0$ for some $i>p$, $f_{\beta_{D}}^{(m_D)}\ldots f_{\beta_p}^{(m_p)}y'$ belongs to the right ideal $\bar{Y}_{p+1}^{res}$ of $U_\mathcal{B}^{s,res}(\n_-)$ generated by $f_{\beta_u}^{(m_u)}$, $u>p$ and $m_u>0$. By Lemma \ref{segmPBWs} (vii) the elements $f_{\beta_{D}}^{(m_D)}\ldots f_{\beta_1}^{(m_1)}$ with at least one $m_u>0$ for $u>p$ form a linear basis of this ideal, and this basis does not contain multiples of $f_{\beta_p}$. Therefore according to (\ref{defvm}) 
\begin{equation}\label{Yp}
(v_p, \omega_0 S_s^{-1}(y'')v_{\mu_p})=0, \text{ for any } y''\in \bar{Y}_{p+1}^{res},
\end{equation}
and hence the right hand side of (\ref{ap2}) vanishes if at least one $m_i>0$ for some $i>p$. 

If $m_i=0$ for all $i>p$ we have $f_{\beta_{D}}^{(m_D)}\ldots f_{\beta_p}^{(m_p)}y'=f_{\beta_p}^{(m_p)}y'$, $m_p>1$. By Lemma \ref{segmPBWs} (ix) $y'\in U_\mathcal{B}^{s,res}(\n_-)=U_\mathcal{B}^{s,res}([-\beta_{p+1},-\beta_D])U_\mathcal{B}^{s,res}([-\beta_1,-\beta_p])$, and hence by Lemma \ref{segmPBWs} (viii) 
$$
f_{\beta_p}^{(m_p)}y'\in f_{\beta_p}^{(m_p)}U_\mathcal{B}^{s,res}([-\beta_{p+1},-\beta_D])U_\mathcal{B}^{s,res}([-\beta_1,-\beta_p])\subset \sum_{i=0}^{m_p-1}(U_{\mathcal{B}}^{s, res}([-\beta_{p+1},-\beta_D]))_0f_{\beta_p}^{(i)}U_\mathcal{B}^{s,res}([-\beta_1,-\beta_p])+
$$
$$
+U_{\mathcal{B}}^{s, res}([-\beta_{p+1},-\beta_D])f_{\beta_p}^{(m_p)}U_\mathcal{B}^{s,res}([-\beta_1,-\beta_p]),
$$
where $(U_{\mathcal{B}}^{s, res}([-\beta_{p+1},-\beta_D]))_0=\bar{Y}_{p+1}^{res}\cap U_{\mathcal{B}}^{s, res}([-\beta_{p+1},-\beta_D])$, \index[not]{U@$(U_\mathcal{B}^{s, res}([-\beta_{p+1},-\beta_D]))_0$} so for $i<m_p$ one has 
$$
(U_{\mathcal{B}}^{s, res}([-\beta_{p+1},-\beta_D]))_0f_{\beta_p}^{(i)}U_\mathcal{B}^{s,res}([-\beta_1,-\beta_p])\subset \bar{Y}_{p+1}^{res}.
$$
Thus by (\ref{Yp}) 
$$ 
(v_p, \omega_0 S_s^{-1}(f_{\beta_p}^{(m_p)}y')v_{\mu_p})=(v_p, \omega_0 S_s^{-1}(y''')v_{\mu_p}), y'''\in U_{\mathcal{B}}^{s, res}([-\beta_{p+1},-\beta_D])f_{\beta_p}^{(m_p)}U_\mathcal{B}^{s,res}([-\beta_1,-\beta_p]).
$$

By Lemma \ref{segmPBWs} (vi) and by Remark \ref{segmPBWsrev} the decompositions of elements of $U_{\mathcal{B}}^{s, res}([-\beta_{p+1},-\beta_D])f_{\beta_p}^{(m_p)}U_\mathcal{B}^{s,res}([-\beta_1,-\beta_p])$ with respect to the basis $f_{\beta_{D}}^{(r_D)}\ldots f_{\beta_1}^{(r_1)}$ only contain non-zero multiples of elements of the basis with $r_p\geq m_p>1$. None of these elements is equal to $f_{\beta_p}$.  Therefore by (\ref{defvm}) 
$$
(v_p, yuv_{\mu_p})=(v_p, y\omega_0 S_s^{-1}(y')v_{\mu_p})=(v_p, \omega_0 S_s^{-1}(f_{\beta_p}^{(m_p)}y')v_{\mu_p})=(v_p, \omega_0 S_s^{-1}(y''')v_{\mu_p})=0.
$$

Thus
$$
(v_p, yuv_{\mu_p})=(v_p,f_{\beta_{D}}^{(m_D)}\ldots f_{\beta_p}^{(m_p)}uv_{\mu_p})=0
$$
if $y=\omega_0 S_s^{-1}(f_{\beta_{D}}^{(m_D)}\ldots f_{\beta_p}^{(m_p)})\neq 1, \omega_0 S_s^{-1}(f_{\beta_p})$. 

The result that we proved implies that $\omega(y)v_p$ is orthogonal, with respect to the contravariant form, to any element of the form $uv_{\mu_p}$, $u\in U_\mathcal{B}^{s,res}(\g)$. Since any element of $V_{\mu_p}^{res}$ is of this form and the contravariant form is non--degenerate we deduce that $\omega(y)v_p=0$.
This completes the proof of Lemma \ref{Apvan}.

\end{proof}

In order co define quantum counterparts of the functions $\varphi_p$ we have to introduce a certain localization $\mathbb{C}_{\mathcal{B}}^{s,loc}[G]$ of $\mathbb{C}_{\mathcal{B}}^s[G]$. By Lemma 9.1.10 in \cite{Jos} it is possible to define a localization of $\mathbb{C}_{\mathcal{B}}^s[G]$, in the sense of localization for non--commutative algebras, which contains the quantum counterparts of the denominators in the definition of the functions $\varphi_p$ given by (\ref{Apquant}). We shall not need these results in full generality. In fact, we shall only need right denominators and $\mathbb{C}_{\mathcal{B}}^{s,loc}[G]$ will be defined as a ``right'' localization of $\mathbb{C}_{\mathcal{B}}^s[G]$. The exact meaning of this term will be explained below. 

We start by introducing a subalgebra $\mathbb{C}_{\mathcal{B}}^s[G]_0\subset \mathbb{C}_{\mathcal{B}}^s[G]$ a proper localization of which contains all the required denominators.  
\begin{lemma}\label{go}
The set of elements $(u,\cdot~ \tau(T_s^{-1})v)\in \mathbb{C}_{\mathcal{B}}^s[G]$, where $(u,\cdot~ v)\in \mathbb{C}_{\mathcal{B}}^s[G]$, $v$ is a highest weight vector in a finite rank representation $V$ of  $U_h(\g)$ and $u\in V$, form a subalgebra $\mathbb{C}_{\mathcal{B}}^s[G]_0$ \index[not]{C@$\mathbb{C}_\mathcal{B}^s[G]_0$} in $\mathbb{C}_{\mathcal{B}}^s[G]$.
\end{lemma}

\begin{proof}
It suffices to show that the product of two elements from $\mathbb{C}_{\mathcal{B}}^s[G]_0$ belongs to $\mathbb{C}_{\mathcal{B}}^s[G]_0$.
Let $s=s_{i_1}\ldots s_{i_k}$ be a reduced decomposition of $s$. Then by (\ref{DTs})
\begin{align}\label{DTs'}
\Delta_s^\tau(T_s^\tau)=(\tau\otimes\tau)(\prod^{k}_{p=1}\theta_{\beta_p}^s\mathcal{F}_s(T_s \otimes T_s)(\mathcal{F}_s^{-1})T_s\otimes T_s),
\end{align}
where in the products $\theta_{\beta_p}^s$ appears on the left from $\theta_{\beta_q}^s$ if $p<q$, and for $p=1,\ldots, k$
$$
e_{\beta_p}=\psi_{\{ n_{ij}\}}^{-1}(X_{\beta_p}^+e^{hK_s\beta_p^\vee}),
f_{\beta}=\psi_{\{ n_{ij}\}}^{-1}(e^{-hK_s\beta_p^\vee}X_{\beta_p}^-), \beta_p=s_{i_1}\ldots s_{i_{p-1}}\alpha_{i_p},
$$
$$
X_{\beta_p}^\pm=T_{i_1}\ldots T_{i_{p-1}}X_{i_p}^\pm,
$$
$$
\theta_{\beta_p}^s={\exp}_{q_{\beta_p}}[(1-q_{\beta_p}^{-2})e_{\beta_p}e^{-h\kappa {1+s \over 1-s}P_{\h'} \beta_p^\vee}\otimes f_{\beta_p}].
$$

Next, 
\begin{equation}\label{qt1}
\mathcal{F}_s(T_s \otimes T_s)(\mathcal{F}_s^{-1})=q^{\sum_{i=1}^l(T_sY_i\otimes T_sK_sH_i-Y_i\otimes K_sH_i)}=q^{\sum_{i=1}^l(T_sY_i\otimes K_sT_{\widehat{s}_p}H_i-Y_i\otimes K_sH_i)}=
\end{equation}
$$
=q^{\sum_{i=1}^l(Y_i\otimes K_sH_i-Y_i\otimes K_sH_i)}=1.
$$
Here we also used the identity $\sum_{i=1}^lT_sY_i\otimes K_sT_sH_i=\sum_{i=1}^lY_i\otimes K_sH_i$ which holds since $T_sY_i$, $T_sH_i$, $i=1,\ldots ,l$ is a pair of dual bases in $\h$ with respect to the bilinear form.

Now by (\ref{qt1}) formula (\ref{DTs'}) takes the form
\begin{equation}\label{DTss}
\Delta_s^\tau(T_s^\tau)
=(\tau\otimes \tau)(\prod^{k}_{p=1}\theta_{\beta_p}^s T_s\otimes T_s),
\end{equation}
and for two highest weight vectors $v\in V, v'\in V'$ we have by Lemma \ref{hwv} (ii)
\begin{equation}\label{ds}
\Delta_s^\tau(\tau(T_s^{-1}))v\otimes v'=\tau(T_s^{-1})v\otimes \tau(T_s^{-1})v'.	
\end{equation}

Therefore for any $u\in V, u'\in V'$ we have by (\ref{ds})
\begin{equation}\label{multsys}
(u,\cdot~ \tau(T_s^{-1})v)\otimes (u',\cdot~ \tau(T_s^{-1})v')=(u\otimes u',\cdot~ \tau(T_s^{-1})v\otimes \tau(T_s^{-1})v')=(u\otimes u',\cdot~ \tau(T_s^{-1})(v\otimes v')).
\end{equation}
Since $v\otimes v'$ is a highest weight vector in $V\otimes V'$, the last identity implies $(u,\cdot~ \tau(T_s^{-1})v)\otimes (u',\cdot~ \tau(T_s^{-1})v')\in \mathbb{C}_{\mathcal{B}}^s[G]_0$ if $(u,\cdot~ v)\in \mathbb{C}_{\mathcal{B}}^s[G]$ and $(u',\cdot~ v')\in \mathbb{C}_{\mathcal{B}}^s[G]$. This completes the proof.

\end{proof}

Note that from (\ref{comm1'}) for $u=v_\lambda$ we obtain using (\ref{multsys}) with $v=u=v_\lambda$, $v'=u'=v_\mu$ that
\begin{equation}\label{comm1'''}
\Delta_\mu^s\otimes \Delta_\lambda^s=\Delta_\lambda^s \otimes \Delta_\mu^s=\Delta_{\lambda+\mu}^s.
\end{equation}

From formulas (\ref{comm1}) and (\ref{multsys}) it follows that
the set $\mathfrak{S}_s=\{\hat{c}q^{n\frac{1}{d{\bar{r}}^2}}\Delta_\mu^s| \mu\in P_+, \hat{c}\in \mathbb{C}^*,n\in \mathbb{Z}\}$ \index[not]{S@$\mathfrak{S}_s$} is a multiplicative set of normal elements in $\mathbb{C}_{\mathcal{B}}^s[G]_0$.

Let $\mathbb{C}_{\mathcal{B}}^{s,loc}[G]_0$ \index[not]{C@$\mathbb{C}_\mathcal{B}^{s,loc}[G]_0$}  be the localization of $\mathbb{C}_{\mathcal{B}}^s[G]_0$ by $\mathfrak{S}_s$.
Denote $\mathfrak{S}_s^*=\{f\otimes g^{-1}\in \mathbb{C}_{\mathcal{B}}^{s,loc}[G]_0| f,g \in \mathfrak{S}_s\}$, \index[not]{S@$\mathfrak{S}_s^*$} $\mathfrak{S}_s^{-1}=\{f^{-1}\in \mathbb{C}_{\mathcal{B}}^{s,loc}[G]_0| f\in \mathfrak{S}_s\}$, where we denote the multiplication in $\mathbb{C}_{\mathcal{B}}^{s,loc}[G]_0$ by the same symbol $\otimes$ which will be sometimes omitted if it does not lead to a confusion. \index[not]{S@$\mathfrak{S}_s^{-1}$}

We shall need more information on the structure of subalgebras $U_k$ in $\mathbb{C}_{\mathcal{B}}^{s,loc}[G]_0$ generated by the elements  
\begin{equation}\label{Bpdef}
B_p:={A_p^0}^{-1}\otimes A_p\in \mathbb{C}_{\mathcal{B}}^{s,loc}[G]_0, \index[not]{B@$B_p$}
\end{equation}
where $p=1,\ldots, k$. The elements $B_p$ are quantum analogues of the functions $\varphi_p$. From (\ref{comm1''}) and (\ref{comm1'''}) we also obtain
\begin{equation}\label{comm1''''}
\Delta_\mu^s\otimes B_p=q^{\left\langle (\kappa{1+s \over 1-s }P_{{\h'}}-id)\beta_p^\vee, \mu^\vee\right\rangle}B_p\otimes \Delta_\mu^s.	
\end{equation}

The following Lemma is similar to Proposition 8.3 in \cite{GY}.
\begin{lemma}
Let $U_k$ \index[not]{U@$U_k$} be the non-unital subalgebra in $\mathbb{C}_{\mathcal{B}}^{s,loc}[G]_0$ generated by the elements  $B_p={A_p^0}^{-1}\otimes A_p\in \mathbb{C}_{\mathcal{B}}^{s,loc}[G]_0$, $p=1,\ldots, k$, and $U_{\mathcal{B}}^k$ \index[not]{U@$U_\mathcal{B}^k$} the non-unital subalgebra in $\mathbb{C}_{\mathcal{B}}^s[G_*]$ generated by the elements $\bar{e}_p:=(q_{\beta_p}^{-1}-q_{\beta_p})q^{\beta_p^\vee}e_{\beta_p}$, $p=1,\ldots, k$. \index[not]{e@$\bar{e}_p$} Then the map
$$
B_p\mapsto \bar{e}_p
$$
gives rise to an algebra isomorphism $\vartheta :U_k \rightarrow U_{\mathcal{B}}^k$. \index[not]{t@$\vartheta$}

In particular, the elements $B_p$ satisfy the following commutation relations
\begin{equation}\label{Bcmrel}
B_pB_r - q^{\left\langle \beta_p,\beta_r\right\rangle+\left\langle \kappa{1+s \over 1-s}P_{\h'^*}\beta_p,\beta_r\right\rangle}B_rB_p= \sum_{m_{p+1},\ldots,m_{r-1}\in \mathbb{N}}C(m_{p+1},\ldots,m_{r-1})
B_{p+1}^{m_{p+1}}\ldots B_{r-1}^{m_{r-1}},~p<r,
\end{equation}
where $C(m_{p+1},\ldots,m_{r-1})\in \mathcal{B}$.
\end{lemma}

\begin{proof}
First observe that the map
$$
\psi:\mathbb{C}_{\mathcal{B}}^s[G]\rightarrow \mathbb{C}_{\mathcal{B}}^s[B_-],
(u,\cdot~ v)\mapsto (u\otimes id,(\omega_0\otimes id)(\mathcal{R}_s)v\otimes id)
$$
is an algebra homomorphism, where $(\omega_0\otimes id)(\mathcal{R}_s)\in U_h^s(\g)\otimes U_h^s(\g)$, and the matrix element $(u,\cdot~ v)$ is evaluated on the first factor of the tensor product. Indeed, since $\omega_0'$ is a coautomorphism we have by the first property in (\ref{rmprop}) in the case of $\Delta_s$ and $\mathcal{R}_s$
$$
(\Delta_s^\tau\otimes id)((\omega_0\otimes id)(\mathcal{R}_s))=(\omega_0\otimes \omega_0\otimes id)(\Delta_s\otimes id)\mathcal{R}_s=
(\omega_0\otimes \omega_0\otimes id)((\mathcal{R}_s)_{13}(\mathcal{R}_s)_{23})=(\omega_0\otimes id)((\mathcal{R}_s)_{13})(\omega_0\otimes id)((\mathcal{R}_s)_{23}),
$$
and hence
$$
\psi((u,\cdot~ v)\otimes (u',\cdot~ v'))=(u\otimes id,(\omega_0\otimes id)((\mathcal{R}_s)_{13})v\otimes id)(u'\otimes id,(\omega_0\otimes id)((\mathcal{R}_s)_{23})v'\otimes id)=\psi((u,\cdot~ v))\psi((u',\cdot~ v')), \index[not]{p@$\psi$}
$$
i.e. $\psi$ is an algebra homomorphism.

Recall that by Lemma \ref{go} elements of the form $(u,\cdot~ \tau(T_s^{-1})v)\in \mathbb{C}_{\mathcal{B}}^s[G]$, where $v$ is a highest weight vector in a finite rank representation $V$ of  $U_h(\g)$ and $u\in V$, form a subalgebra $\mathbb{C}_{\mathcal{B}}^s[G]_0$ in $\mathbb{C}_{\mathcal{B}}^s[G]$.

From (\ref{ds}) we also deduce that the map
$$
\mathbb{C}_{\mathcal{B}}^s[G]_0\rightarrow \mathbb{C}_{\mathcal{B}}^s[G],
(u,\cdot~ \tau(T_s^{-1})v)\mapsto (u,\cdot~ v)
$$
is an algebra homomorphism.
Composing this map with $\psi$ we obtain another algebra homomorphism $\psi^0$. \index[not]{p@$\psi^0$}

Next using (\ref{SRs}), (\ref{r-1})  and the definition of $A_p$ we obtain 
$$
\psi^0(A_p)=(v_p\otimes id,(\omega_0S_s^{-1}S_s\otimes id)(\mathcal{R}_s)v_{\mu_p}\otimes id)=(v_p\otimes id,(\omega_0S_s^{-1}\otimes id)({\mathcal{R}_s}^{-1})v_{\mu_p}\otimes id)=
$$
$$
=(q_{\beta_p}^{-1}-q_{\beta_p})q^{-(\kappa{1+s \over 1-s }P_{{\h'}}+id)\mu_p^\vee}q^{\beta_p^\vee}e_{\beta_p}.
$$

Similarly
$$
\psi^0(A_p^0)=q^{-(\kappa{1+s \over 1-s }P_{{\h'}}+id)\mu_p^\vee}.
$$
From the last two formulas we deduce that $\psi^0$ gives rise to an algebra homomorphism $\vartheta :U_k \rightarrow U_{\mathcal{B}}^k$ such that
$$
\vartheta(B_p)=\vartheta({A_p^0}^{-1}\otimes A_p)=(q_{\beta_p}^{-1}-q_{\beta_p})q^{\beta_p^\vee}e_{\beta_p}=\bar{e}_p.
$$
This homomorphism is surjective by construction. $\vartheta$ is also injective as $\psi^0$ is injective (see Proposition 8.3 in \cite{GY} and \cite{Y}, Theorem 2.6 for the proof). 

Commutation relations (\ref{Bcmrel}) follow from (\ref{cmrel}) by applying $\omega$ and by multiplying by $q^{\beta_p^\vee}q^{\beta_r^\vee}$.

\end{proof}

Denote $\mathbb{C}_{\mathcal{B}}^{s,loc}[G]=\mathbb{C}_{\mathcal{B}}^s[G]\otimes_{\mathbb{C}_{\mathcal{B}}^s[G]_0}\mathbb{C}_{\mathcal{B}}^{s,loc}[G]_0$. \index[not]{C@$\mathbb{C}_\mathcal{B}^{s,loc}[G]$} $\mathbb{C}_{\mathcal{B}}^{s,loc}[G]$ is naturally a left $\mathbb{C}_{\mathcal{B}}^s[G]$--module and a right $\mathbb{C}_{\mathcal{B}}^{s,loc}[G]_0$--module. We denote by $\otimes$ 
both the left $\mathbb{C}_{\mathcal{B}}^s[G]$--action and the right $\mathbb{C}_{\mathcal{B}}^{s,loc}[G]_0$--action on $\mathbb{C}_{\mathcal{B}}^{s,loc}[G]$ and call these actions multiplications. We shall often omit the symbol $\otimes$ to shorten the notation if it does not lead to a confusion.

From (\ref{comm1''}) it follows that 
\begin{equation}\label{ABS}
\mathfrak{S}_s^*A_p\subset \mathcal{B}A_p\mathfrak{S}_s^*, \mathfrak{S}_s^*B_p\subset \mathcal{B}A_p\mathfrak{S}_s^*, U_r\Delta_\mu^s=\Delta_\mu^s U_r,
\end{equation}
where $\mathcal{B}A_p\mathfrak{S}_s^*$ (resp. $\mathcal{B}A_p\mathfrak{S}_s^*$) is the $\mathcal{B}$--submodule in $\mathbb{C}_{\mathcal{B}}^{s,loc}[G]$ generated by $A_p\mathfrak{S}_s^*$ (resp. by $A_p\mathfrak{S}_s^*$).  

By (\ref{Bcmrel}) we also have
\begin{equation}
U_pB_r\subset B_rU_{r-1}+U_{r-1}, r> p,
\end{equation}
and hence if we denote by $\overline{U}_p$ \index[not]{U@$\overline{U}_p$} the non-unital subalgebra in $\mathbb{C}_{\mathcal{B}}^s[G]$ generated by $A_1,\ldots, A_p$ then by (\ref{comm1''}) and by the relations $A_r=\Delta_{\mu_r}B_r$, $r=1,\ldots, D$ following from the definition of $B_r$ we have
\begin{equation}\label{AUp}
\overline{U}_pA_r\subset A_r\overline{U}_{r-1}\mathfrak{S}_s^*+\overline{U}_{r-1}\mathfrak{S}_s^*, r> p.
\end{equation}

Also by (\ref{comm1''})
\begin{equation}\label{OUm}
\overline{U}_p\Delta_\mu^s =\Delta_\mu^s \overline{U}_p.
\end{equation}

Let $J^p=\mathbb{C}_{\mathcal{B}}^s[G]\overline{U}_p\subset \mathbb{C}_{\mathcal{B}}^s[G]$ \index[not]{J@$J_p$} be the left ideal generated by $\overline{U}_p$, ${J^p}^{loc}\subset \mathbb{C}_{\mathcal{B}}^{s,loc}[G]$ \index[not]{J@${J^p}^{loc}$} the image of $J^p\otimes \mathfrak{S}_s^*\subset \mathbb{C}_{\mathcal{B}}^s[G]\otimes \mathbb{C}_{\mathcal{B}}^{s,loc}[G]_0$ under the projection $\mathbb{C}_{\mathcal{B}}^s[G]\otimes \mathbb{C}_{\mathcal{B}}^{s,loc}[G]_0 \rightarrow \mathbb{C}_{\mathcal{B}}^s[G]\otimes_{\mathbb{C}_{\mathcal{B}}^s[G]_0}\mathbb{C}_{\mathcal{B}}^{s,loc}[G]_0=\mathbb{C}_{\mathcal{B}}^{s,loc}[G]$. 
\begin{lemma}\label{rightloc}
${J^p}^{loc}A_r={J^p}^{loc}B_r\subset {J^{r-1}}^{loc}$ for $r>p$.
\end{lemma}
 
\begin{proof}
Firstly, the identity ${J^p}^{loc}A_r={J^p}^{loc}B_r$ follows from the definitions of $A_r, B_r$ and of ${J^p}^{loc}$.

Now we establish the inclusion.
Since by the definition ${J^p}^{loc}={J^p}\mathfrak{S}_s^*=\mathbb{C}_{\mathcal{B}}^s[G]\overline{U}_p\mathfrak{S}_s^*$, for any $r> p$ one has by the first formula in (\ref{ABS})
\begin{equation}\label{IpAr}
{J^p}^{loc}A_r={J^p}\mathfrak{S}_s^*A_r={J^p}A_r\mathfrak{S}_s^*.
\end{equation}

By (\ref{AUp}) we have
\begin{equation}\label{JpAr}
{J^p}A_r=\mathbb{C}_{\mathcal{B}}^s[G]\overline{U}_{p}A_r\subset \mathbb{C}_{\mathcal{B}}^s[G]A_r\overline{U}_{r-1}\mathfrak{S}_s^*+\overline{U}_{r-1}\mathfrak{S}_s^*)\subset {J^{r}_{\mathcal{B}}}\mathfrak{S}_s^*={J^r}^{loc}.
\end{equation}
This completes the proof.

\end{proof}


\section{Quantized vanishing ideals}\label{vanid}

\pagestyle{myheadings}
\markboth{CHAPTER~\thechapter.~ZHELOBENKO TYPE OPERATORS FOR Q-W--ALGEBRAS}{\thesection.~QUANTIZED VANISHING IDEAL OF THE LEVEL SURFACE OF THE MOMENT MAP}

\setcounter{equation}{0}
\setcounter{theorem}{0}

In this section we introduce and study quantum analogues of the ideals $J^{j1}$, $j=2,\ldots, R-1$.

Firstly we obtain an alternative description of ${J^{11}_{\mathcal{B}}}'$ which agrees at the classical level with the description of $J^{11}$ in Proposition \ref{Jpdescr} (ii). We start with a technical lemma. In the rest of this chapter we assume that for $s\in W$ a positive root system $\Delta_+^s$ is equipped with a normal ordering (\ref{NO}), $\Delta_+$ is the corresponding system of positive roots associated to $s$ in Definition \ref{circorddef}.
\begin{lemma}\label{ii10}
Let $\Delta_{\m_+}^s=\Delta_{\m_+}\cap \Delta_s^s$ and note that $(\Delta_0)_+=\Delta_0\cap\Delta_+^s=\Delta_0\cap \Delta_+$. Both $\Delta_{\m_+}^s, (\Delta_0)_+\subset \Delta_+$ are minimal segments. Denote by $U_q^{res}(-\Delta_{\m_+}^s)$ and $U_q^{res}((\Delta_0)_+)$ the subalgebras of $U_q^{res}(\g)$ corresponding to $-\Delta_{\m_+}^s$ and $(\Delta_0)_+$, respectively (see Proposition \ref{segmPBW} (iv) for their definition).
Let $b\in U_q^{res}(w'(\b_+))$, $(u,\cdot~ u')\in \mathbb{C}_{\mathcal{B}}^s[G]$ an element of $\mathbb{C}_{\mathcal{B}}^s[G]$ such that $u$ is a highest weight vector in a finite rank representation $V$ of $U_h(\g)$, $u'\in V$. Then
\begin{equation}\label{tsb}
(u,\tau(T_sb) u')=(u,\tau(T_s\sum_i z_+^ix_i') u'),
\end{equation}
where $z_+^i\in U_q^{res}((\Delta_0)_+)$, $x_i'\in U_q^{res}(-\Delta_{\m_+}^s)$.
\end{lemma}

\begin{proof}
First note that $U_q^{res}(w'(\b_+))=U_{U_q^{res}(H)}^{res}([\beta_{k_{l'}+1},-\beta_{k_{l'}}])$, and we have a disjoint union of minimal segments (see Figure 4)
$$
[\beta_{k_{l'}+1},-\beta_{k_{l'}}]=(-\Delta_{\m_+}^s)\cup(-(\Delta_+^s\setminus (\Delta_s^s \cup (\Delta_0)_+)))\cup (\Delta_0)_+\cup (\Delta_{s^2}^s\setminus \Delta_{\m_+}),
$$
where the order of the segments in the union agrees with the circular normal ordering of $\Delta$ corresponding to normal ordering (\ref{NO}) of $\Delta_+^s$.

Applying iteratively Corollary \ref{segmq} (i) we obtain using this union
$$
U_q^{res}(w'(\b_+))=U_{U_q^{res}(H)}^{res}([\beta_{k_{l'}+1},-\beta_{k_{l'}}])=
$$
$$
=U_{U_q^{res}(H)}^{res}((-(\Delta_+^s\setminus (\Delta_s^s \cup (\Delta_0)_+)))\cup(\Delta_0)_+\cup (\Delta_{s^2}^s\setminus \Delta_{\m_+}))U_{U_q^{res}(H)}^{res}(-\Delta_{\m_+}^s)=
$$
$$
=U_{U_q^{res}(H)}^{res}((-(\Delta_+^s\setminus (\Delta_s^s \cup (\Delta_0)_+)))U_{U_q^{res}(H)}^{res}((\Delta_0)_+\cup (\Delta_{s^2}^s\setminus \Delta_{\m_+}))U_{U_q^{res}(H)}^{res}(-\Delta_{\m_+}^s)=
$$
$$
=U_{U_q^{res}(H)}^{res}(-(\Delta_+^s\setminus (\Delta_s^s \cup (\Delta_0)_+)))U_{U_q^{res}(H)}^{res}(\Delta_{s^2}^s\setminus \Delta_{\m_+})U_{U_q^{res}(H)}^{res}((\Delta_0)_+)U_{U_q^{res}(H)}^{res}(-\Delta_{\m_+}^s).
$$

Note that $(\Delta_0)_+\subset \Delta_+$ and $-\Delta_{\m_+}^s\subset \Delta_-$. Therefore by Corollary \ref{segmq} (ii)
$$
U_q^{res}(w'(\b_+))=U_{U_q^{res}(H)}^{res}(-(\Delta_+^s\setminus (\Delta_s^s \cup (\Delta_0)_+)))U_{U_q^{res}(H)}^{res}(\Delta_{s^2}^s\setminus \Delta_{\m_+})U_{U_q^{res}(H)}^{res}((\Delta_0)_+)U_q^{res}(H)U_q^{res}(-\Delta_{\m_+}^s)=
$$
\begin{equation}\label{uqresb}
U_{U_q^{res}(H)}^{res}(-(\Delta_+^s\setminus (\Delta_s^s \cup (\Delta_0)_+)))U_{U_q^{res}(H)}^{res}(\Delta_{s^2}^s\setminus \Delta_{\m_+})U_q^{res}(H)U_q^{res}((\Delta_0)_+)U_q^{res}(-\Delta_{\m_+}^s)=
\end{equation}
$$
=U_{U_q^{res}(H)}^{res}(-(\Delta_+^s\setminus (\Delta_s^s \cup (\Delta_0)_+)))U_{U_q^{res}(H)}^{res}(\Delta_{s^2}^s\setminus \Delta_{\m_+})U_q^{res}((\Delta_0)_+)U_q^{res}(-\Delta_{\m_+}^s).
$$ 

Now observe that by the definition of $\Delta_s^s$ we have $s(-(\Delta_+^s\setminus (\Delta_s^s \cup (\Delta_0)_+)))\subset -\Delta_+^s$, and that $s^{-1}(s(-(\Delta_+^s\setminus (\Delta_s^s \cup (\Delta_0)_+))))=-(\Delta_+^s\setminus (\Delta_s^s \cup (\Delta_0)_+))\subset \Delta_-^s$, so in fact $s(-(\Delta_+^s\setminus (\Delta_s^s \cup (\Delta_0)_+)))\subset -(\Delta_+^s\setminus \Delta_{s^{-1}}^s)\subset \Delta_-$, where the last inclusion follows from the definition of $\Delta_-$ (see also Figure 4). Recalling that for $\mu\in Q$ one has $T_s(U_q^{res}(\g))_\mu=(U_q^{res}(\g))_{s\mu}$ we deduce  
\begin{equation}\label{ts1} 
T_s(U_{U_q^{res}(H)}^{res}(-(\Delta_+^s\setminus (\Delta_s^s \cup (\Delta_0)_+))))\subset \bigoplus_{\mu\leq 0}(U_q^{res}(\g))_\mu.
\end{equation}

Next, $\Delta_{s^2}^s\setminus \Delta_{\m_+}\subset \Delta_{s^2}^s$, so  
$$
s(\Delta_{s^2}^s\setminus \Delta_{\m_+})\subset s^1s^2(\Delta_{s^2}^s)=-s^1(\Delta_{s^2}^s),
$$
where at the last step we used the fact that $s^2$ is an involution, and hence $s^2(\Delta_{s^2}^s)=-\Delta_{s^2}^s$. Now by Proposition \ref{pord} (i)
$$
s(\Delta_{s^2}^s\setminus \Delta_{\m_+})\subset -s^1(\Delta_{s^2}^s)\subset -(\Delta_{s^{-1}}^s\setminus \Delta_{s^1}^s)\subset \Delta_-,
$$
where at the last step we used the definition of $\Delta_-$ (see also Figure 4).
Similarly to (\ref{ts1}) the last inclusion implies
\begin{equation}\label{ts2} 
T_s(U_{U_q^{res}(H)}^{res}(\Delta_{s^2}^s\setminus \Delta_{\m_+}))\subset \bigoplus_{\mu\leq 0}(U_q^{res}(\g))_\mu.
\end{equation}

Now we can express any $b\in U_q^{res}(w'(\b_+))$ using (\ref{uqresb}) as follows
\begin{equation}\label{xi1}
b=\sum_iy_iw_i'z_+^ix_i,
\end{equation}
where
$y_i\in U_{U_q^{res}(H)}^{res}(-(\Delta_+^s\setminus (\Delta_s^s \cup (\Delta_0)_+)))$, $w_i'\in U_{U_q^{res}(H)}^{res}(\Delta_{s^2}^s\setminus \Delta_{\m_+})$, $z_+^i\in U_q^{res}((\Delta_0)_+)$, $x_i\in U_q^{res}(-\Delta_{\m_+}^s)$.

Let $(u, \cdot ~ u')\in \mathbb{C}_{\mathcal{B}}^s[G]$ be such that $u$ is a highest weight vector in a finite rank representation $V$ of $U_h(\g)$, and $u'\in V$.
Using (\ref{xi1}) we obtain 
$$
(u,\tau(T_sb) u')=(u,\sum_i\tau(T_s(y_i)T_s(w_i')T_sz_+^ix_i) u').
$$
Note that by (\ref{ts1}) and (\ref{ts2}) the weights of the elements $\tau(T_s(y_i))$ and $\tau(T_s(w_i'))$ are non--positive, and hence Lemma \ref{hwv} yields
\begin{equation}\label{Tsx}
(u,\tau(T_sb) u')=(u,\sum_i\tau(T_s(y_i')T_s(w_i'')T_sz_+^ix_i) u'),
\end{equation}
where $y_i'\in U_{U_q^{res}(H)}^{res}(-(\Delta_+^s\setminus (\Delta_s^s \cup (\Delta_0)_+)))$, $w_i''\in U_{U_q^{res}(H)}^{res}(\Delta_{s^2}^s\setminus \Delta_{\m_+})$ have zero weights. Recalling that $-(\Delta_+^s\setminus (\Delta_s^s \cup (\Delta_0)_+))$ and $\Delta_{s^2}^s\setminus \Delta_{\m_+}$ are minimal segments we deduce that the only zero weight elements of the algebras $U_{U_q^{res}(H)}^{res}(-(\Delta_+^s\setminus (\Delta_s^s \cup (\Delta_0)_+)))$ and  $U_{U_q^{res}(H)}^{res}(\Delta_{s^2}^s\setminus \Delta_{\m_+})$ belong to $U_q^{res}(H)$. Therefore $y_i', w_i''\in U_q^{res}(H)$, and by Lemma \ref{wHact} (\ref{Tsx}) takes the form
$$
(u,T_sb u')=\sum_ik_i'(u,\tau(T_sz_+^ix_i) u')=(u,\sum_i\tau(T_sz_+^ix_i') u'),
$$
where $k_i'\in \mathbb{C}[q,q^{-1}]$, $z_+^i\in U_q^{res}((\Delta_0)_+)$, $x_i'=\tau(k_i')x_i\in U_q^{res}(-\Delta_{\m_+}^s)$. This completes the proof

\end{proof}

The following description of ${J^{11}_{\mathcal{B}}}'$ is analogous to the description of $J^{11}$ in Proposition \ref{Jpdescr} (ii).
\begin{lemma}\label{ii1}
${J^{11}_{\mathcal{B}}}'$ coincides with the left ideal generated by the elements $(u,\cdot~ v)\in \mathbb{C}_{\mathcal{B}}^s[G]$ such that $u$ is a highest weight vector in a finite rank representation $V$ of $U_h(\g)$, and $v\in V$ satisfies $(u,\tau(t_szx) v)=0$ for any  $z\in U_{\mathcal{B}}^{s,res}((\Delta_0)_+)$, $x\in U_{\mathcal{B}}^{s,res}(-\Delta_{\m_+})$ and any element $t_s$ \index[not]{t@$t_s$} of the braid group acting on $\h\subset U_h(\h)$ in the same way as $s$.
\end{lemma}

\begin{proof}
Firstly we prove the statement for $t_s=T_s$ as in the definition of ${J^{11}_{\mathcal{B}}}'$.

Assume that $(u,\tau(T_sb) v)=0$ for any $b \in U_q^{res}(w'(\b_+))$. Observe that $U_{\mathcal{B}}^{s,res}((\Delta_0)_+)\subset U_{\mathcal{B}}^{s,res}(\h)U_q^{res}((\Delta_0)_+)$ and $U_{\mathcal{B}}^{s,res}(-\Delta_{\m_+})\subset U_{\mathcal{B}}^{s,res}(\h)U_q^{res}(-\Delta_{\m_+})$. Thus 
\begin{equation}\label{pzm}
U_{\mathcal{B}}^{s,res}((\Delta_0)_+)U_{\mathcal{B}}^{s,res}(-\Delta_{\m_+})\subset U_{\mathcal{B}}^{s,res}(\h)U_q^{res}((\Delta_0)_+)U_{\mathcal{B}}^{s,res}(\h)U_q^{res}(-\Delta_{\m_+})=
\end{equation}
$$
=U_{\mathcal{B}}^{s,res}(\h)U_q^{res}((\Delta_0)_+)U_q^{res}(-\Delta_{\m_+})
$$
as $U_q^{res}((\Delta_0)_+)U_{\mathcal{B}}^{s,res}(\h)=U_{\mathcal{B}}^{s,res}(\h)U_q^{res}((\Delta_0)_+)$.

Since $(\Delta_0)_+,-\Delta_{\m_+}\subset [\beta_{k_{l'}+1},-\beta_{k_{l'}}]$, we have $U_q^{res}((\Delta_0)_+)U_q^{res}(-\Delta_{\m_+})\subset U_{U_q^{res}(H)}^{res}([\beta_{k_{l'}+1},-\beta_{k_{l'}}])=U_q^{res}(w'(\b_+))$. Thus by (\ref{pzm})
$$
U_{\mathcal{B}}^{s,res}((\Delta_0)_+)U_{\mathcal{B}}^{s,res}(-\Delta_{\m_+})\subset U_{\mathcal{B}}^{s,res}(\h)U_q^{res}(w'(\b_+)),
$$
and for any  $z\in U_{\mathcal{B}}^{s,res}((\Delta_0)_+)$, $x\in U_{\mathcal{B}}^{s,res}(-\Delta_{\m_+})$ one has $zx=\sum_ih_ib_i$, $h_i\in U_{\mathcal{B}}^{s,res}(\h)$, $b_i\in U_q^{res}(w'(\b_+))$. So by Lemma \ref{wHact}
$$
(u,\tau(T_szx) v)=\sum_i(u,\tau(T_sh_ib_i)v)=\sum_ic_i(u,\tau(T_sb_i)v)=0, 
$$
where $c_i\in \mathcal{B}$ are defined by $c_iu=\tau(T_s(h_i))u$.

Conversely, if $(u,\tau(T_szx) v)=0$ for any  $z\in U_{\mathcal{B}}^{s,res}((\Delta_0)_+)$, $x\in U_{\mathcal{B}}^{s,res}(-\Delta_{\m_+}^s)$
then by Lemma \ref{ii10} with $u'=v$ 
\begin{equation}\label{utb}
(u,\tau(T_sb) v)=(u,\tau(T_s\sum_i z_ix_i) v),
\end{equation}
where for all $i$ $z_i\in U_q^{res}((\Delta_0)_+)$, $x_i\in U_q^{res}(-\Delta_{\m_+}^s)\subset U_q^{res}(-\Delta_{\m_+})$. 

Observe that $U_{\mathcal{B}}^{s,res}(\h)U_{\mathcal{B}}^{s,res}((\Delta_0)_+)\supset U_q^{res}((\Delta_0)_+)$ and $U_{\mathcal{B}}^{s,res}(\h)U_{\mathcal{B}}^{s,res}(-\Delta_{\m_+})\supset U_q^{res}(-\Delta_{\m_+})$. Thus 
\begin{equation}\label{pzm1}
U_q^{res}((\Delta_0)_+)U_q^{res}(-\Delta_{\m_+})\subset U_{\mathcal{B}}^{s,res}(\h)U_{\mathcal{B}}^{s,res}((\Delta_0)_+)U_{\mathcal{B}}^{s,res}(\h)U_{\mathcal{B}}^{s,res}(-\Delta_{\m_+})=
\end{equation}
$$
=U_{\mathcal{B}}^{s,res}(\h)U_{\mathcal{B}}^{s,res}((\Delta_0)_+)U_{\mathcal{B}}^{s,res}(-\Delta_{\m_+})
$$
as $U_{\mathcal{B}}^{s,res}((\Delta_0)_+)U_{\mathcal{B}}^{s,res}(\h)=U_{\mathcal{B}}^{s,res}(\h)U_{\mathcal{B}}^{s,res}((\Delta_0)_+)$.

By (\ref{pzm1}) formula (\ref{utb}) takes the form
$$
(u,\tau(T_sb) v)=(u,\tau(T_s\sum_i h_i'z_i'x_i') v),
$$
where $z_i'\in U_{\mathcal{B}}^{s,res}((\Delta_0)_+)$, $x_i'\in U_{\mathcal{B}}^{s,res}(-\Delta_{\m_+})$, $h_i'\in U_{\mathcal{B}}^{s,res}(\h)$.

Thus by Lemma \ref{wHact}
$$
(u,\tau(T_sb) v)=\sum_i(u,\tau(T_s h_i'z_i'x_i')v)=\sum_ic_i'(u,\tau(T_s z_i'x_i')v)=0,
$$
where $c_i'\in \mathcal{B}$ are defined by $c_i'u=\tau(T_s(h_i'))u$.

Hence ${J^{11}_{\mathcal{B}}}'$ coincides with the left ideal generated by the elements $(u,\cdot~ v)\in \mathbb{C}_{\mathcal{B}}^s[G]$ such that $u$ is a highest weight vector in a finite rank representation $V$ of $U_h(\g)$, and $v\in V$ satisfies $(u,\tau(T_szx) v)=0$ for any  $z\in U_{\mathcal{B}}^{s,res}((\Delta_0)_+)$, $x\in U_{\mathcal{B}}^{s,res}(-\Delta_{\m_+})$. 

Finally by Lemma \ref{bga} $T_s$ used in the definition of ${J^{11}_{\mathcal{B}}}'$ can be replaced in the statement of this lemma with any element $t_s$ of the braid group acting on $\h\subset U_h(\h)$ in the same way as $s$. This completes the proof.

\end{proof}

We proceed with the definition of the quantum analogues of the ideals $J^{j1}$, $j=2,\ldots, R-1$.
Denote $\delta_{jo_j}= (w_1\ldots w_{j-1})^{-1}\beta_{R-1 n_{R-1}}=(w_1\ldots w_{j-1})^{-1}\beta_c$. \index[not]{d@$\delta_{jo_j}$}
Let $U_{\mathcal{B}}^{\widehat{s}_j,res}([-\delta_{jk},-\delta_{jo_j}])=T_{w_1\ldots w_{j-1}}^{-1}U_{\mathcal{B}}^{s,res}([-\beta_{jk},-\beta_c])$ \index[not]{U@$U_\mathcal{B}^{\widehat{s}_j,res}([-\delta_{jk},-\delta_{jo_j}])$}  be the subalgebra in $U_{\mathcal{B}}^{\widehat{s}_j,res}(\g)$ generated by the elements $T_{w_1\ldots w_{j-1}}^{-1}f_\beta^{(n)}(s)$, $\beta\in [\beta_{jk},\beta_{R-1 n_{R-1}}]$, $n\in \mathbb{N}$. Note that from the definition of the elements $f_\beta$ it follows that 
\begin{equation}\label{twjb}
T_{w_1\ldots w_{j-1}}^{-1}f_\beta^{(n)}(s)=f_{(w_1\ldots w_{j-1})^{-1}\beta}^{(n)}(\widehat{s}_j)
\end{equation}
for $\beta\in [\beta_{j1},\beta_D]$, where for $\alpha\in \Delta_+$ and $j>1$ the elements $f_{\alpha}(\widehat{s}_j)$ are defined using the element $\widehat{s}_j$ and the normal ordering on $\Delta_+$ introduced in Remark \ref{Djord}. To shorten the notation we write from now on $f_{\alpha}^{(n)}(\widehat{s}_j)=f_{\alpha}^{(n)}$ for $\alpha\in [\delta_{j1},\delta_{jD}]$ if it does not cause a confusion. We shall also need the subalgebra $U_{\mathcal{B}}^{\widehat{s}_j,res}([-\delta_{jk},-\delta_{jD}])\subset U_{\mathcal{B}}^{\widehat{s}_j,res}(\g)$. \index[not]{U@$U_\mathcal{B}^{\widehat{s}_j,res}([-\delta_{jk},-\delta_{jD}])$}

In complete analogy with the definition of the ideals $J^{j1}$ for $j=1,\ldots, R-1$ we define ${J^{j1}_{\mathcal{B}}}'$ \index[not]{J@${J^{j1}_\mathcal{B}}'$} as the left ideal in $\mathbb{C}_{\mathcal{B}}^{\widehat{s}_j}[G]$ \index[not]{C@$\mathbb{C}_\mathcal{B}^{\widehat{s}_j}[G]$}  generated by the elements $(u,\cdot~ v)\in \mathbb{C}_{\mathcal{B}}^{\widehat{s}_j}[G]$ such that $u$ is a highest weight vector in a finite rank representation $V$ of $U_h(\g)$, and $v\in V$ satisfies $(u,\tau(t_{\widehat{s}_j}z_+x) v)=0$ for any  $z_+\in U_{\mathcal{B}}^{\widehat{s}_j,res}(\z_+^j):=T_{w_1\ldots w_{j-1}}^{-1}U_{\mathcal{B}}^{s,res}((\Delta_0)_+)$, \index[not]{U@$U_\mathcal{B}^{\widehat{s}_j,res}(\z_+^j)$} $x\in U_{\mathcal{B}}^{\widehat{s}_j,res}([-\delta_{j1},-\delta_{jo_j}])$ and an arbitrary fixed element $t_{\widehat{s}_j}$ of the braid group acting on $\h\subset U_h(\h)$ in the same way as $\widehat{s}_j$. 

Note that by Lemma \ref{bga} this definition does not depend on the choice of $t_{\widehat{s}_j}$.
For $j=1$ this definition agrees with the previously given definition of ${J^{11}_{\mathcal{B}}}'$ due to Lemma \ref{ii1}. 

Also similarly to ${I^{11}_{\mathcal{B}}}$ defined in Proposition \ref{Iq} we introduce left ideals ${I^{j1}_{\mathcal{B}}}={J^{j1}_q}'\cap \mathbb{C}_{\mathcal{B}}^{\widehat{s}_j}[G]$ \index[not]{I@${I^{j1}_{\mathcal{B}}}$} in $\mathbb{C}_{\mathcal{B}}^{\widehat{s}_j}[G]$, where ${J^{j1}_q}'={{J^{j1}_{\mathcal{B}}}'}\otimes_{\mathcal{B}}\mathbb{C}(q^{\frac{1}{d{\bar{r}}^2}})$, $j=1,\ldots, R-1$. \index[not]{J@${J^{j1}_q}'$}

Next, for the left ideals ${I^{j1}_{\mathcal{B}}}$ we obtain a description similar to the one given in Proposition \ref{Jpdescr} (i).
We start with a technical lemma which will allow us to use properties of the vanishing ideals $J^{jk}$ to prove some properties of their quantum counterparts.
\begin{lemma}\label{incl}
Let $X$ be a free $\mathcal{B}$--module, $V\subset V'\subset X$ two its submodules such that $V=V'$ mod $(q^{\frac{1}{d{\bar{r}}^2}}-1)V'$. Let $V_q=V\otimes_{\mathcal{B}}{\mathbb{C}}(q^{\frac{1}{d{\bar{r}}^2}})$, $V_q'=V'\otimes_{\mathcal{B}}{\mathbb{C}}(q^{\frac{1}{d{\bar{r}}^2}})$. Then $V_q=V'_q$. 
\end{lemma} 

\begin{proof}
First observe that by Theorem 6.5 in \cite{R} $V'\subset X$ is $\mathcal{B}$--free, as $\mathcal{B}$ is a principal ideal domain, and $V\subset V'$ is a free submodule of $V'$. Let $e_a$, $a\in \mathbb{A}$ be a basis of $V'$. As noted in the proof of Theorem 6.5 in \cite{R}, $\mathbb{A}$ is a well ordered set, i.e. it is a totally ordered set in which
any nonempty subset has a smallest element. Then as shown in the proof of Theorem 6.5 in \cite{R} $V$ has a basis elements of which have the form $f_a=\sum_{b\leq a}c_a^be_b$, $a\in \mathbb{B}\subset \mathbb{A}$, where $c_a^a\neq 0$, $c_a^b\in \mathcal{B}$, and the sum is finite. Since $V=V'$ mod $(q^{\frac{1}{d{\bar{r}}^2}}-1)V'$ we must have $\mathbb{B}=\mathbb{A}$, and $c_a^a\neq 0$ mod $(q^{\frac{1}{d{\bar{r}}^2}}-1)$ for all $a\in \mathbb{A}$.


Now using a simple transfinite induction we can express the elements of the basis $e_a$ regarded as elements of $V'_q$ in terms of $f_a$. Indeed, if $a_0\in \mathbb{A}$ is the minimal element then $e_{a_0}={c_{a_0}^{a_0}}^{-1}f_{a_0}$, and assuming that $e_d=\sum_{b\leq d}g_d^bf_b$, $g_d^b\in {\mathbb{C}}(q^{\frac{1}{d{\bar{r}}^2}})$, $g_d^d\neq 0$ holds for all $d<a$ we get $e_a={c_a^a}^{-1}(f_a-\sum_{b<a}c_a^be_b)={c_a^a}^{-1}(f_a-\sum_{b<a}c_a^b\sum_{h\leq b}g_b^hf_h)$, where all sums are finite. This establishes the induction step and completes the proof.

\end{proof}

\begin{lemma}\label{barJ1}
(i) For $j=1,\ldots, R-1$, let $t_{\widehat{s}_j}$ be an arbitrary element of the braid group acting on $\h\subset U_h(\h)$ in the same way as $\widehat{s}_j$. Let $J^{j1}_{\mathcal{B}}$ \index[not]{J@$J^{j1}_\mathcal{B}$} be the left ideal in $\mathbb{C}_{\mathcal{B}}^{\widehat{s}_j}[G]$ generated by the elements $(w,\cdot~ v)$ such that $(w,\tau(yt_{\widehat{s}_j}z_+x)v)=0$ for any $y\in U_{\mathcal{B}}^{\widehat{s}_j,res}(\b_-), z_+\in U_{\mathcal{B}}^{\widehat{s}_j,res}(\z_+^j), x\in U_{\mathcal{B}}^{\widehat{s}_j,res}([-\delta_{j1},-\delta_{jo_j}])$. 

Then ${{J^{j1}_{\mathcal{B}}}'}\subset J^{j1}_{\mathcal{B}}$ and ${{J^{j1}_{\mathcal{B}}}'}=J^{j1}_{\mathcal{B}}=J^{j1}$  mod $(q^{\frac{1}{d{\bar{r}}^2}}-1)J^{j1}_{\mathcal{B}}$.

(ii) Denote ${J^{j1}_q}'={{J^{j1}_{\mathcal{B}}}'}\otimes_{\mathcal{B}}\mathbb{C}(q^{\frac{1}{d{\bar{r}}^2}})$, ${J^{j1}_q}=J^{j1}_{\mathcal{B}}\otimes_{\mathcal{B}}\mathbb{C}(q^{\frac{1}{d{\bar{r}}^2}})$. \index[not]{J@$J^{j1}_q$} Then ${J^{j1}_q}'={J^{j1}_q}$. Thus ${I^{j1}_{\mathcal{B}}}={J^{j1}_q}\cap \mathbb{C}_{\mathcal{B}}^{\widehat{s}_j}[G]={J^{j1}_q}'\cap \mathbb{C}_{\mathcal{B}}^{\widehat{s}_j}[G]$.

(iii) ${{J^{j1}_{\mathcal{B}}}'}$, $J^{j1}_{\mathcal{B}}$ and ${I^{j1}_{\mathcal{B}}}$ are stable under the ${\rm Ad}_{\widehat{s}_j}^0$-action of $U_{\mathcal{B}}^{\widehat{s}_j,res}([-\delta_{j1},-\delta_{jo_j}])$. Thus $\mathbb{C}_{\mathcal{B}}^{\widehat{s}_j}[G]/{{J^{j1}_{\mathcal{B}}}'}$ and $\mathbb{C}_{\mathcal{B}}^{\widehat{s}_j}[G]/{I^{j1}_{\mathcal{B}}}$ are naturally equipped with the $U_{\mathcal{B}}^{\widehat{s}_j,res}([-\delta_{j1},-\delta_{jo_j}])$--action induced by the ${\rm Ad}_{\widehat{s}_j}^0$-action of $U_{\mathcal{B}}^{\widehat{s}_j,res}([-\delta_{j1},-\delta_{jo_j}])$ on $\mathbb{C}_{\mathcal{B}}^{\widehat{s}_j}[G]$.
\end{lemma}

\begin{proof}

(i) If $(u,\cdot~ v)\in {J^{j1}_{\mathcal{B}}}'$ is one of the elements generating ${J^{j1}_{\mathcal{B}}}'$, i.e. $u$ is a highest weight vector in a finite rank representation $V$ of $U_h(\g)$, and $v\in V$ obeys $(u,\tau(t_{\widehat{s}_j}z_+x) v)=0$ for any  $z_+\in U_{\mathcal{B}}^{\widehat{s}_j,res}(\z_+^j)$, $x\in U_{\mathcal{B}}^{\widehat{s}_j,res}([-\delta_{j1},-\delta_{jo_j}])$, then for any $y\in U_{\mathcal{B}}^{\widehat{s}_j,res}(\b_-)$, $z_+\in U_{\mathcal{B}}^{\widehat{s}_j,res}(\z_+^j)$, $x\in U_{\mathcal{B}}^{\widehat{s}_j,res}([-\delta_{j1},-\delta_{jo_j}])$ one has by Lemma \ref{hwv} (i) and by the definition of ${J^{j1}_{\mathcal{B}}}'$ that
$$
(u,\tau(yt_{\widehat{s}_j}z_+x)v)=(u,\tau(y_0t_{\widehat{s}_j}z_+x)v)=\check{c}(u,\tau(t_{\widehat{s}_j}z_+x)v)=0,
$$
where $y_0\in U_{\mathcal{B}}^{\widehat{s}_j,res}(\h)$ is the zero weight component of $y$, and $\check{c}\in \mathcal{B}$ is defined by $\check{c}u=\tau(y_0)u$.

Thus $(u,\cdot~ v)$ is also one of the elements generating $J^{j1}_{\mathcal{B}}$, so that ${{J^{j1}_{\mathcal{B}}}'}\subset J^{j1}_{\mathcal{B}}$. This proves the first claim of part (i) of the lemma. 

Next, observe that the specializations of $\tau(U_{\mathcal{B}}^{\widehat{s}_j,res}(\z_+^j))$ and of $\tau(U_{\mathcal{B}}^{\widehat{s}_j,res}([-\delta_{j1},-\delta_{jo_j}]))$ at $q^{\frac{1}{d{\bar{r}}^2}}=1$ are isomorphic to $U(\z_+^j)$ and $U(\m_-^j)$, respectively, the quotient of the specialization of $\tau(U_{\mathcal{B}}^{\widehat{s}_j,res}(\b_-))$ at $q^{\frac{1}{d{\bar{r}}^2}}=1$ by the ideal generated by $L_i-1$, and by $t_i-1$, $i=1,\ldots, l$ is isomorphic to $U(\b_-)$, and the elements $L_i$ and $t_i$ act by identity transformations on the specialization of any module $V^{res}$ at $q^{\frac{1}{d{\bar{r}}^2}}=1$. Therefore recalling Proposition \ref{wresdec} (ii) we deduce that the specialization of $J^{j1}_{\mathcal{B}}$ at $q^{\frac{1}{d{\bar{r}}^2}}=1$ is generated by the matrix elements of the form $(w, \cdot~ v)\in \mathbb{C}[G]$, where $w,v\in V$, $V$ is a finite-dimensional representation of $\g$, and $(w,y\widehat{s}_jz_+xv)=0$ for any $y\in U(\b_-)$, $z_+\in U(\z_+^j)$, $x\in U(\m_-^j)$, and the specialization of ${{J^{j1}_{\mathcal{B}}}'}$ at $q^{\frac{1}{d{\bar{r}}^2}}=1$ is generated by the matrix elements of the form $(u, \cdot~ v)\in \mathbb{C}[G]$, where $u$ is a highest weight vector in a finite-dimensional representation $V$ of $\g$, and $v\in V$ satisfies, and $(v,\widehat{s}_jz_+xv)=0$ for any $z_+\in U(\z_+^j)$, $x\in U(\m_-^j)$.

Since the enveloping algebra of any Lie subalgebra of $\g$ is also a Hopf subalgebra of $U(\g)$, the product of any element of $\mathbb{C}[G]$ and of $(w, \cdot~ v)\in \mathbb{C}[G]$, where $w,v\in V$, $V$ is a finite-dimensional representation of $\g$, and $(w,y\widehat{s}_jz_+xv)=0$ for any $y\in U(\b_-)$, $z_+\in U(\z_+^j)$, $x\in U(\m_-^j)$ satisfies the same condition. So in fact the specialization of $J^{j1}_{\mathcal{B}}$ at $q^{\frac{1}{d{\bar{r}}^2}}=1$ consists of the matrix elements of the form $(w, \cdot v)\in \mathbb{C}[G]$, where $w,v\in V$, $V$ is a finite-dimensional representation of $\g$, and $(w,y\widehat{s}_jz_+xv)=0$ for any $y\in U(\b_-)$, $z_+\in U(\z_+^j)$, $x\in U(\m_-^j)$.

Now from parts (i) and (ii) of Proposition \ref{Jpdescr} it follows that ${{J^{j1}_{\mathcal{B}}}'}=J^{j1}_{\mathcal{B}}=J^{j1}$  mod $(q^{\frac{1}{d{\bar{r}}^2}}-1)J^{j1}_{\mathcal{B}}$.

(ii) Note that by part (i) ${{J^{j1}_{\mathcal{B}}}'}\subset J^{j1}_{\mathcal{B}}$ are submodules of the $\mathcal{B}$--module $\mathbb{C}_{\mathcal{B}}^{\widehat{s}_j}[G]$ which is free by Proposition \ref{Afree} (iv). Also, by part (i) $J^{j1}_{\mathcal{B}}={{J^{j1}_{\mathcal{B}}}'}=J^{j1}$ mod $(q^{\frac{1}{d{\bar{r}}^2}}-1)J^{j1}_{\mathcal{B}}$. Therefore by Lemma \ref{incl} ${J^{j1}_q}={J^{j1}_q}'$.

(iii) Firstly we prove the statement for ${{J^{j1}_{\mathcal{B}}}'}$. 

From formula (\ref{ad0n}) with $s=\widehat{s}_j$ we have for any $\beta\in [\delta_{j1},\delta_{jo_j}]$
$$
{\rm Ad}_{\widehat{s}_j}^0f_\beta^{(n)} (f\otimes g)(~\cdot ~, ~\cdot ~)=
\sum_{k=0}^{n}\sum_{p=0}^{n-k}q_{\beta}^{k(n-k)+p(n-k-p)}({\rm Ad}_{\widehat{s}_j}^0 (G_\beta^{-k}f_\beta^{(p)}) f)(\cdot ~)\otimes g(\omega_0S_{\widehat{s}_j}^{-1}(G_\beta^{-k-p}f_\beta^{(n-k-p)}) \cdot ~ \omega_0(f_\beta^{(k)}))+	
$$
\begin{equation}\label{AdbarJ10}
+\sum_{k=0}^{n-1}\sum_i q_\beta^{k(n-k)}({\rm Ad}_{\widehat{s}_j}^0 (G_\beta^{-k} x_i^{(n-k)}) f)(~\cdot ~)\otimes g((\omega_0 S_{\widehat{s}_j}^{-1})(G_\beta^{-k}y_i^{(n-k)})~\cdot ~ \omega_0(f_\beta^{(k)}))+
\end{equation}
$$
+\sum_i ({\rm Ad}_{\widehat{s}_j}^0({y_i^{(n)}}^2) f)(~\cdot ~)\otimes g((\omega_0 S_{\widehat{s}_j}^{-1})({y_i^{(n)}}^1)~\cdot ~ \omega_0 (x_i^{(n)})).
$$

Assume that $g(~\cdot ~)=(u, \cdot ~ v)\in \mathbb{C}_{\mathcal{B}}^{\widehat{s}_j}[G]$, where $u$ is a highest weight vector in a finite rank representation $V$ of $U_h(\g)$, and $v\in V$ is such that $(u,\tau(t_{\widehat{s}_j}z_+x) v)=0$ for any $z_+\in U_{\mathcal{B}}^{\widehat{s}_j,res}(\z_+^j)$, $x\in U_{\mathcal{B}}^{\widehat{s}_j,res}([-\delta_{j1},-\delta_{jo_j}])$, and $f\in \mathbb{C}_{\mathcal{B}}^{\widehat{s}_j}[G]$ is arbitrary.

We claim now that all terms in the right hand side of the last formula belong to ${{J^{j1}_{\mathcal{B}}}'}$. More generally, one can prove the following statement.
\begin{lemma}
Let $g(~\cdot ~)=(u, \cdot ~ v)\in \mathbb{C}_{\mathcal{B}}^{\widehat{s}_j}[G]$, where $u$ is a highest weight vector in a finite rank representation $V$ of $U_h(\g)$, and $v\in V$ is such that $(u,\tau(t_{\widehat{s}_j}z_+x) v)=0$ for any $z_+\in U_{\mathcal{B}}^{\widehat{s}_j,res}(\z_+^j)$, $x\in U_{\mathcal{B}}^{\widehat{s}_j,res}([-\delta_{j1},-\delta_{jo_j}])$. Then for any $y'\in U_{\mathcal{B}}^{\widehat{s}_j,res}(\b_-)$ and $b\in U_{\mathcal{B}}^{\widehat{s}_j,res}([-\delta_{j1},-\delta_{jo_j}])$ the element $g(\tau(y')~\cdot~ \tau(b))$ satisfies the same properties as $g$.
\end{lemma}
\begin{proof}
By Lemma \ref{wHact}
$$
g(\tau(y')~\cdot~ \tau(b))=(u,\tau(y)~\cdot~ \tau(b)v)=(u,\tau(y'_0)~\cdot~ \tau(b)v)=c_0(u, \cdot~ \tau(b)v),
$$ 
where $y'_0\in U_{\mathcal{B}}^{\widehat{s}_j,res}(\h)$ is the zero weight component of $y'$ and $c_0\in \mathcal{B}$ is defined by $\tau(y'_0)u=c_0u$. 

Also, for any $z_+\in U_{\mathcal{B}}^{\widehat{s}_j,res}(\z_+^j)$, $x\in U_{\mathcal{B}}^{\widehat{s}_j,res}([-\delta_{j1},-\delta_{jo_j}])$ we have $xb\in U_{\mathcal{B}}^{\widehat{s}_j,res}([-\delta_{j1},-\delta_{jo_j}])$, and hence by the definition of $g$
$$
(u,\tau(y't_{\widehat{s}_j}z_+x b)v)=c_0(u,\tau(t_{\widehat{s}_j}z_+x b)v)=0.
$$
 
\end{proof}

Note that by (\ref{fo0}) $\omega_0'(U_{\mathcal{B}}^{\widehat{s}_j,res}([-\delta_{j1},-\delta_{jo_j}]))=U_{\mathcal{B}}^{\widehat{s}_j,res}([-\delta_{j1},-\delta_{jo_j}])$, and hence 
\begin{equation}\label{om0t}
\omega_0(U_{\mathcal{B}}^{\widehat{s}_j,res}([-\delta_{j1},-\delta_{jo_j}]))=\tau(U_{\mathcal{B}}^{\widehat{s}_j,res}([-\delta_{j1},-\delta_{jo_j}])).
\end{equation}

Since $\omega_0S_{\widehat{s}_j}^{-1}(G_\beta^{-k-p}f_\beta^{(n-k-p)}),\omega_0S_{\widehat{s}_j}^{-1}(G_\beta^{-k}y_i^{(n-k)}),\omega_0S_{\widehat{s}_j}^{-1}({y_i^{(n)}}^1)\in \tau(U_{\mathcal{B}}^{\widehat{s}_j,res}(\b_-))=U_{\mathcal{B}}^{\widehat{s}_j,res}(\b_-)$ and by (\ref{om0t}) one has $\omega_0(f_\beta^{(k)}), \omega_0(x_i^{(n)})\in \tau(U_{\mathcal{B}}^{\widehat{s}_j,res}([-\delta_{j1},-\delta_{jo_j}]))$ in (\ref{AdbarJ10}), the previous lemma implies that all terms in the right hand side of (\ref{AdbarJ10}) belong to ${{J^{j1}_{\mathcal{B}}}'}$, i.e. ${\rm Ad}_{\widehat{s}_j}^0f_\beta^{(n)} (f\otimes g)(~\cdot~, ~\cdot~)\in {{J^{j1}_{\mathcal{B}}}'}$.
Since the elements $f\otimes g$, with $f, g\in \mathbb{C}_{\mathcal{B}}^{\widehat{s}_j}[G]$, $g(\cdot ~)=(u, \cdot ~ v)\in \mathbb{C}_{\mathcal{B}}^{\widehat{s}_j}[G]$, where $u$ is a highest weight vector in a finite rank representation $V$ of $U_h(\g)$, and $v\in V$ is such that $(u,\tau(t_{\widehat{s}_j}z_+x) v)=0$ for any $z_+\in U_{\mathcal{B}}^{\widehat{s}_j,res}(\z_+^j)$, $x\in U_{\mathcal{B}}^{\widehat{s}_j,res}([-\delta_{j1},-\delta_{jo_j}])$, span ${{J^{j1}_{\mathcal{B}}}'}$ we obtain ${\rm Ad}_{\widehat{s}_j}^0f_\beta^{(n)} ({{J^{j1}_{\mathcal{B}}}'})\subset {{J^{j1}_{\mathcal{B}}}'}$. The elements $f_\beta^{(n)}$, $\beta\in [\delta_{j1},\delta_{jo_j}]$ generate $U_{\mathcal{B}}^{\widehat{s}_j,res}([-\delta_{j1},-\delta_{jo_j}])$, and hence ${{J^{j1}_{\mathcal{B}}}'}$ is stable under the ${\rm Ad}_{\widehat{s}_j}^0$-action of $U_{\mathcal{B}}^{\widehat{s}_j,res}([-\delta_{j1},-\delta_{jo_j}])$ by definition.

If we naturally extend the ${\rm Ad}_{\widehat{s}_j}^0$-action of $U_{\mathcal{B}}^{\widehat{s}_j,res}([-\delta_{j1},-\delta_{jo_j}])$ to $\mathbb{C}_q^{\widehat{s}_j}[G]=\mathbb{C}_{\mathcal{B}}^{\widehat{s}_j}[G]\otimes_{\mathcal{B}}\mathbb{C}(q^{\frac{1}{d{\bar{r}}^2}})$ we immediately deduce that ${J^{j1}_q}={{J^{j1}_{\mathcal{B}}}'}\otimes_{\mathcal{B}}\mathbb{C}(q^{\frac{1}{d{\bar{r}}^2}})$ is stable under this action.

Hence ${I^{j1}_{\mathcal{B}}}={J^{j1}_q}\cap \mathbb{C}_{\mathcal{B}}^{\widehat{s}_j}[G]$ is also stable under this action as clearly $\mathbb{C}_{\mathcal{B}}^{\widehat{s}_j}[G]$ is ${\rm Ad}_{\widehat{s}_j}^0$--stable. 

Assume now that $g(~\cdot ~)=(w, ~\cdot ~ v)\in \mathbb{C}_{\mathcal{B}}^{\widehat{s}_j}[G]$ is such that $(w,\tau(yt_{\widehat{s}_j}z_+x)v)=0$ for any $y\in U_{\mathcal{B}}^{\widehat{s}_j,res}(\b_-), z_+\in U_{\mathcal{B}}^{\widehat{s}_j,res}(\z_+^j), x\in U_{\mathcal{B}}^{\widehat{s}_j,res}([-\delta_{j1},-\delta_{jo_j}])$, and $f\in \mathbb{C}_{\mathcal{B}}^{\widehat{s}_j}[G]$ is arbitrary.

We claim again that all terms in the right hand side of (\ref{AdbarJ10}) belong to $J^{j1}_{\mathcal{B}}$. More generally, one can prove the following statement.
\begin{lemma}
Let $g(~\cdot ~)=(w, \cdot ~ v)\in \mathbb{C}_{\mathcal{B}}^{\widehat{s}_j}[G]$ is such that $(w,\tau(yt_{\widehat{s}_j}z_+x)v)=0$ for any $y\in U_{\mathcal{B}}^{\widehat{s}_j,res}(\b_-), z_+\in U_{\mathcal{B}}^{\widehat{s}_j,res}(\z_+^j), x\in U_{\mathcal{B}}^{\widehat{s}_j,res}([-\delta_{j1},-\delta_{jo_j}])$. Then for any $y'\in U_{\mathcal{B}}^{\widehat{s}_j,res}(\b_-)$ and $b\in U_{\mathcal{B}}^{\widehat{s}_j,res}([-\delta_{j1},-\delta_{jo_j}])$ the element $g(\tau(y')\cdot~ \tau(b))$ satisfies the same properties as $g$.
\end{lemma}

\begin{proof}
The results is obvious as for any $y\in U_{\mathcal{B}}^{\widehat{s}_j,res}(\b_-), x\in U_{\mathcal{B}}^{\widehat{s}_j,res}([-\delta_{j1},-\delta_{jo_j}])$ one has $y'y\in U_{\mathcal{B}}^{\widehat{s}_j,res}(\b_-), xb\in U_{\mathcal{B}}^{\widehat{s}_j,res}([-\delta_{j1},-\delta_{jo_j}])$, and hence
$$
g(\tau(y' yt_{\widehat{s}_j}z_+xb))=0.
$$
\end{proof}

Since $\omega_0S_{\widehat{s}_j}^{-1}(G_\beta^{-k-p}f_\beta^{(n-k-p)}),\omega_0S_{\widehat{s}_j}^{-1}(G_\beta^{-k}y_i^{(n-k)}),\omega_0S_{\widehat{s}_j}^{-1}({y_i^{(n)}}^1)\in \tau(U_{\mathcal{B}}^{\widehat{s}_j,res}(\b_-))=U_{\mathcal{B}}^{\widehat{s}_j,res}(\b_-)$ and by (\ref{om0t}) one has $\omega_0(f_\beta^{(k)}), \omega_0(x_i^{(n)})\in \tau(U_{\mathcal{B}}^{\widehat{s}_j,res}([-\delta_{j1},-\delta_{jo_j}]))$ in (\ref{AdbarJ10}), the previous lemma implies that all terms in the right hand side of (\ref{AdbarJ10}) belong to $J^{j1}_{\mathcal{B}}$, i.e. ${\rm Ad}_{\widehat{s}_j}^0f_\beta^{(n)} (f\otimes g)(~\cdot~, ~\cdot~)\in J^{j1}_{\mathcal{B}}$.
Since the elements $f\otimes g$, with $f, g\in \mathbb{C}_{\mathcal{B}}^{\widehat{s}_j}[G]$ and $g$ as in the previous lemma, span ${{J^{j1}_{\mathcal{B}}}}$ we obtain ${\rm Ad}_{\widehat{s}_j}^0f_\beta^{(n)} (J^{j1}_{\mathcal{B}})\subset J^{j1}_{\mathcal{B}}$. The elements $f_\beta^{(n)}$, $\beta\in [\delta_{j1},\delta_{jo_j}]$ generate $U_{\mathcal{B}}^{\widehat{s}_j,res}([-\delta_{j1},-\delta_{jo_j}])$, and hence $J^{j1}_{\mathcal{B}}$ is stable under the ${\rm Ad}_{\widehat{s}_j}^0$-action of $U_{\mathcal{B}}^{\widehat{s}_j,res}([-\delta_{j1},-\delta_{jo_j}])$.
This completes the proof.

\end{proof}


\section{Higher quantized vanishing ideals}\label{hvan}

\pagestyle{myheadings}
\markboth{CHAPTER~\thechapter.~ZHELOBENKO TYPE OPERATORS FOR Q-W--ALGEBRAS}{\thesection.~HIGHER QUANTIZED VANISHING IDEALS}

\setcounter{equation}{0}
\setcounter{theorem}{0}

The next step in our approach is the definition and the study of quantum analogues ${{J}^{jk}_{\mathcal{B}}}$ and ${{J}^{jk}_{\mathcal{B}}}'$ of the ideals ${J}^{jk}$ and ${{J}^{jk}}'$ for $k>1$.

Firstly, for each $j=1,\ldots, R-1$, $k=1,\ldots, n_j$ we define the matrix elements $(v_{jk},\cdot~ v_{\mu_{jk}})\in \mathbb{C}_\mathcal{B}^{\widehat{s}_j}[G]$ \index[not]{v@$v_{jk}$} by condition (\ref{defvm}), \index[not]{m@$\mu_{jk}$} where for $j=1$ the normal ordering on $\Delta_+$ introduced in Definition \ref{circorddef} is used, and for $j>1$ the normal ordering on $\Delta_+$ defined in Remark \ref{Djord} is used, and $\beta_p=\delta_{jk}$. Let 
$$
A_{jk}(~\cdot~ )=(v_{jk}, ~\cdot~ \tau(T_{\widehat{s}_j}^{-1})v_{\mu_{jk}}). \index[not]{A@$A_{jk}$}
$$

The left ideals ${{J}^{jk}_{\mathcal{B}}}$ and ${{J}^{jk}_{\mathcal{B}}}'$ are defined in the following lemma.
\begin{lemma}\label{barJp}
For $j=1,\ldots, R-1$, $k=1,\ldots, n_j+1$, let ${{J}^{jk}_{\mathcal{B}}}'$ \index[not]{J@${J^{jk}_\mathcal{B}}'$} be the left ideal in $\mathbb{C}_{\mathcal{B}}^{\widehat{s}_j}[G]$ generated by ${J^{j1}_{\mathcal{B}}}'$ and by $A_{jp}$, $p=1,\ldots, k-1$.
 
Let $t_{\widehat{s}_j}$ be an arbitrary element of the braid group acting on $\h\subset U_h(\h)$ in the same way as $\widehat{s}_j$, ${{J}^{jk}_{\mathcal{B}}}$ \index[not]{J@$J^{jk}_\mathcal{B}$} the left ideal in $\mathbb{C}_{\mathcal{B}}^{\widehat{s}_j}[G]$ generated by the elements $(w,\cdot~ v)\in \mathbb{C}_{\mathcal{B}}^{\widehat{s}_j}[G]$ such that $(w,y\tau(ht_{\widehat{s}_j}z_+x)v)=0$ for any $y\in \omega_0S_{\widehat{s}_j}^{-1}(U_{\mathcal{B}}^{\widehat{s}_j,res}([-\delta_{jk},-\delta_{jD}]))$, $h\in U_{\mathcal{B}}^{\widehat{s}_j,res}(\h)$, $z_+\in U_{\mathcal{B}}^{\widehat{s}_j,res}(\z_+^j), x\in U_{\mathcal{B}}^{\widehat{s}_j,res}([-\delta_{j1},-\delta_{jo_j}])$. 

Then the following statements are true.

(i) If $p<k$  then for any $y'\in \omega_0 S_{\widehat{s}_j}^{-1}(U_{\mathcal{B}}^{\widehat{s}_j,res}([-\delta_{jk},-\delta_{jD}]))$ of the form $y'=\omega_0 S_{\widehat{s}_j}^{-1}(f_{\delta_{D}^j}^{(n_D^j)}\ldots f_{\delta_{jk}}^{(n^j_k)})$, $n_i^j\in \mathbb{N}$, $b\in U_{\mathcal{B}}^{\widehat{s}_j,res}([-\delta_{j1},-\delta_{jo_j}])$, $h'\in U_{\mathcal{B}}^{\widehat{s}_j,res}(\h)$ one has 
$$
A_{jp}(y'h'\cdot~ \tau(b))\in {{J}^{jk}_{\mathcal{B}}}.
$$

If in addition $y'\neq 1,\omega_0 S_{\widehat{s}_j}^{-1}(f_{\delta_{jk}})$ or the zero weight component of $b$ is zero then 
\begin{equation}\label{A1}
A_{jk}(y'h'~\cdot~ \tau(b))\in {{J}^{jk}_{\mathcal{B}}}.
\end{equation}

One also has
\begin{equation}\label{A2}
A_{jk}^0(~\cdot~):=\Delta_{\mu_{jk}}^{\widehat{s}_j}(~\cdot~)=A_{jk}(\omega_0S_{\widehat{s}_j}^{-1}(f_{\delta_{jk}})~\cdot~ ). \index[not]{A@$A_{jk}^0(~\cdot~)$}
\end{equation}

(ii) ${{J}^{j1}_{\mathcal{B}}}'\subset J^{j1}_{\mathcal{B}}\subset {{J}^{jk}_{\mathcal{B}}}$ and ${{J}^{j1}_{\mathcal{B}}}'\subset {{J}^{jk}_{\mathcal{B}}}'\subset {{J}^{jk}_{\mathcal{B}}}$.

(iii) For any $y'\in U_{\mathcal{B}}^{\widehat{s}_j,res}(\b_-)$, $b\in U_{\mathcal{B}}^{\widehat{s}_j,res}([-\delta_{j1},-\delta_{jo_j}])$, if the zero weight component of $y'$ is zero or the zero weight component of $b$ is zero then 
$$
{\Delta_\mu^{\widehat{s}_j}}(\tau(y')~\cdot~ \tau(b))\in {{J}^{j1}_{\mathcal{B}}}'\subset {{J}^{jk}_{\mathcal{B}}}'\subset {{J}^{jk}_{\mathcal{B}}}.
$$

(iv) ${\rm Ad}_{\widehat{s}_j}^0(U_{\mathcal{B}}^{\widehat{s}_j,res}([-\delta_{jk},-\delta_{jo_j}]){{J}^{jk}_{\mathcal{B}}}'\subset {{J}^{jk}_{\mathcal{B}}}$. 
\end{lemma}

\begin{remark}
Although the definition of the left ideals ${{J}^{jk}_{\mathcal{B}}}$ depends on the choice of $t_{\widehat{s}_j}$, many properties of the ideals ${{J}^{jk}_{\mathcal{B}}}$ are the same for different $t_{\widehat{s}_j}$. Therefore in statements and in the proofs below valid for arbitrary $t_{\widehat{s}_j}$ we shall not comment on the choice of $t_{\widehat{s}_j}$.
\end{remark}

\begin{proof}

(i) If $p<k$ and $y'\neq 1$  then $A_{jp}(y'h'~\cdot~ \tau(b))=0$ by Lemma \ref{Apvan} (i).

If $y'=1$ we have to show that $A_{jp}(h'y\tau(ht_{\widehat{s}_j}z_+xb))=0$ for any $y\in \omega_0S_{\widehat{s}_j}^{-1}(U_{\mathcal{B}}^{\widehat{s}_j,res}([-\delta_{jk},-\delta_{jD}]))$, $h\in U_{\mathcal{B}}^{\widehat{s}_j,res}(\h)$, $z_+\in U_{\mathcal{B}}^{\widehat{s}_j,res}(\z_+^j), x\in U_{\mathcal{B}}^{\widehat{s}_j,res}([-\delta_{j1},-\delta_{jo_j}])$.

By linearity it suffices to check this condition for $y$ of the form $y=\omega_0 S_{\widehat{s}_j}^{-1}(f_{\delta_{D}^j}^{(n_D^j)}\ldots f_{\delta_{jk}}^{(n^j_k)})$, $n_i^j\in \mathbb{N}$ as these elements form a basis of $\omega_0S_{\widehat{s}_j}^{-1}(U_{\mathcal{B}}^{\widehat{s}_j,res}([-\delta_{jk},-\delta_{jD}]))$ by Lemma \ref{segmPBWs} (vi).

Note that for any $y=\omega_0 S_{\widehat{s}_j}^{-1}(f_{\delta_{D}^j}^{(n_D^j)}\ldots f_{\delta_{jk}}^{(n^j_k)})$, $n_i^j\in \mathbb{N}$, and $h'\in U_{\mathcal{B}}^{\widehat{s}_j,res}(\h)$ commutation relations (\ref{Ccomm}) and the definition of the elements $f_\alpha$ imply  
\begin{equation}\label{yh'}
h'y=yh'' 
\end{equation}
for some $h''\in U_{\mathcal{B}}^{\widehat{s}_j,res}(\h)$.

Therefore if $y\neq 1$ then by (\ref{yh'}) we obtain $A_{jp}(h'y\tau(ht_{\widehat{s}_j}z_+xb))=A_{jp}(yh''\tau(ht_{\widehat{s}_j}z_+xb))=0$ by Lemma \ref{Apvan} (i). 

If $y=1$ 
\begin{equation}\label{Ajph'}
A_{jp}(h'y\tau(ht_{\widehat{s}_j}z_+xb))=(v_{jp},h'\tau(ht_{\widehat{s}_j}z_+xbT_{\widehat{s}_j}^{-1})v_{\mu_{jp}})=c_{jp}c_{jp}'(v_{jp},\tau(t_{\widehat{s}_j}z_+xbt_{\widehat{s}_j}^{-1})v_{\mu_{jp}})=
\end{equation}
$$
=c_{jp}c_{jp}'(v_{jp},\tau(t_{\widehat{s}_j}(z_+xb))v_{\mu_{jp}}),
$$
where $h'\tau(h)v_{jp}=c_{jp}v_{jp}$, $c_{jp}\in \mathcal{B}$ and $\tau(T_{\widehat{s}_j}^{-1})v_{\mu_{jp}}=c_{jp}'\tau(t_{\widehat{s}_j}^{-1})v_{\mu_{jp}}$, $c_{jp}'\in \mathcal{B}$ by Lemma \ref{bga}.

We proceed with the following lemma.
\begin{lemma}\label{noweight}
(i) For any $v\in U_{\mathcal{B}}^{\widehat{s}_j,res}(\z_+^j)U_{\mathcal{B}}^{\widehat{s}_j,res}([-\delta_{j1},-\delta_{jo_j}])$ the element $t_{\widehat{s}_j}(v)$ has no weight components of weights which belong to $-(w_1\ldots w_{j-1})^{-1}\Delta^j$.

(ii) For any $v\in U_{\mathcal{B}}^{\widehat{s}_j,res}(\z_+^j)U_{\mathcal{B}}^{\widehat{s}_j,res}([-\delta_{j1},-\delta_{jo_j}])$and any $b\in U_{\mathcal{B}}^{\widehat{s}_j,res}(\z_+^j)U_{\mathcal{B}}^{\widehat{s}_j,res}([-\delta_{j1},-\delta_{jo_j}])$ which has no non--trivial zero weight component the element $t_{\widehat{s}_j}(vb)$ has no non--trivial zero weight component.
\end{lemma}

\begin{proof}
(i) The weights of the weight components of $t_{\widehat{s}_j}(v)$ belong to $-\mathbb{N}(w_1\ldots w_{j-1})^{-1}s(\Delta_+^j\cup \Delta_0)$, and by Lemma \ref{jsegm} (ii) this set has empty intersection with $-(w_1\ldots w_{j-1})^{-1}\Delta^j$.

(ii) The weights of the weight components of $t_{\widehat{s}_j}(vb)$ belong to $\mathbb{N}(\widehat{s}_j([-\delta_{j1},-\delta_{jo_j}])\cup (w_1\ldots w_{j-1})^{-1}s(\Delta_0)_+)=\mathbb{N}(w_1\ldots w_{j-1})^{-1}s(([-\beta_{j1},-\beta_c])\cup (\Delta_0)_+)$, and the set $[-\beta_{j1},-\beta_c]\cup (\Delta_0)_+$ is contained in the minimal segment $[\beta_1^0,-\beta_c]$. Therefore if $b$ has no non--trivial zero weight component the element $vb$ has no non--trivial zero weight component as well. Hence $t_{\widehat{s}_j}(vb)$ has no non--trivial zero weight component.

\end{proof}

Now we complete the proof of part (i) of Lemma \ref{barJp}. Since $z_+xb\in U_{\mathcal{B}}^{\widehat{s}_j,res}(\z_+^j)U_{\mathcal{B}}^{\widehat{s}_j,res}([-\delta_{j1},-\delta_{jo_j}])$, by part (i) of the previous lemma the element $t_{\widehat{s}_j}(z_+xb)$ has no weight components of weights which belong to $-(w_1\ldots w_{j-1})^{-1}\Delta^j$, and $-\delta_{jp}\in -(w_1\ldots w_{j-1})^{-1}\Delta^j$. Since $v_{jp}$ has weight $\mu_{jp}-\delta_{jp}$, and $\tau$ preserves weights, we deduce from (\ref{Ajph'}) that
$$
A_{jp}(h'\tau(ht_{\widehat{s}_j}z_+xb))=c_{jp}c_{jp}'(v_{jp},\tau (t_{\widehat{s}_j}(z_+xb))v_{\mu_{jp}})=0
$$
by orthogonality of different weight subspaces with respect to the contravariant form. This proves part (i) for $p<k$.

If $p=k$ and $y'\neq 1, \omega_0 S_{\widehat{s}_j}^{-1}(f_{\delta_{jk}})$  then $A_{jk}(y'\cdot~ \tau(b))=0$ by Lemma \ref{Apvan} (ii).

If $y'=1$ we have to show that $A_{jk}(h'y\tau(ht_{\widehat{s}_j}z_+xb))=0$ for any $y=\omega_0 S_{\widehat{s}_j}^{-1}(f_{\delta_{D}^j}^{(n_D^j)}\ldots f_{\delta_{jk}}^{(n^j_k)})$, $n_i^j\in \mathbb{N}$, $h\in U_{\mathcal{B}}^{\widehat{s}_j,res}(\h)$, $z_+\in U_{\mathcal{B}}^{\widehat{s}_j,res}(\z_+^j), x\in U_{\mathcal{B}}^{\widehat{s}_j,res}([-\delta_{j1},-\delta_{jo_j}])$. 

Therefore if $y\neq 1, \omega_0 S_{\widehat{s}_j}^{-1}(f_{\delta_{jk}})$ we obtain using (\ref{yh'}) that $A_{jk}(h'y\tau(ht_{\widehat{s}_j}z_+xb))=A_{jk}(yh''\tau(ht_{\widehat{s}_j}z_+xb))=0$ by Lemma \ref{Apvan} (i).

If $y=1$ then
\begin{equation}\label{Ajkh'}
A_{jk}(h'y\tau(ht_{\widehat{s}_j}z_+xb))=(v_{jk},h'\tau(ht_{\widehat{s}_j}z_+xbT_{\widehat{s}_j}^{-1})v_{\mu_{jk}})=c_{jk}c_{jk}'(v_{jk},\tau(t_{\widehat{s}_j}z_+xbt_{\widehat{s}_j}^{-1})v_{\mu_{jk}})=
\end{equation}
$$
=c_{jk}c_{jk}'(v_{jk},\tau(t_{\widehat{s}_j}(z_+xb))v_{\mu_{jk}}),
$$
where $\tau(h)h'v_{jk}=c_{jk}v_{jk}$, $c_{jk}\in \mathcal{B}$ and $\tau(T_{\widehat{s}_j}^{-1})v_{\mu_{jk}}=c_{jk}'\tau(t_{\widehat{s}_j}^{-1})v_{\mu_{jk}}$, $c_{jk}'\in \mathcal{B}$ by Lemma \ref{bga}.

Since $z_+xb\in U_{\mathcal{B}}^{\widehat{s}_j,res}(\z_+^j)U_{\mathcal{B}}^{\widehat{s}_j,res}([-\delta_{j1},-\delta_{jo_j}])$, by Lemma \ref{noweight} (i) the element $t_{\widehat{s}_j}(z_+xb)$ has no weight components of weights which belong to $-(w_1\ldots w_{j-1})^{-1}\Delta^j$, and $-\delta_{jk}\in -(w_1\ldots w_{j-1})^{-1}\Delta^j$. Since $v_{jk}$ has weight $\mu_{jk}-\delta_{jk}$ and $\tau$ preserves weights, one has using (\ref{Ajkh'})
$$
A_{jk}(h'y\tau(ht_{\widehat{s}_j}z_+xb))=c_{jk}c_{jk}'(v_{jk},\tau(t_{\widehat{s}_j}(z_+xb))v_{\mu_{jk}})=0
$$
by orthogonality of different weight subspaces with respect to the contravariant form.

If $y=\omega_0 S_{\widehat{s}_j}^{-1}(f_{\delta_{jk}})$ then by (\ref{oof})
$$
A_{jk}(h'y\tau(ht_{\widehat{s}_j}z_+xb))=(v_{jk},h'\omega_0 S_{\widehat{s}_j}^{-1}(f_{\delta_{jk}})\tau(ht_{\widehat{s}_j}z_+xbT_{\widehat{s}_j}^{-1})v_{\mu_{jk}})=
$$
$$
=c_{jk}''c_{jk}'(v_{jk},\omega_0 S_{\widehat{s}_j}^{-1}(f_{\delta_{jk}})\tau(ht_{\widehat{s}_j}z_+xbt_{\widehat{s}_j}^{-1})v_{\mu_{jp}})=c_{jk}'''c_{jk}''c_{jk}'(v_{\mu_{jk}},\tau(t_{\widehat{s}_j}(z_+xb))v_{\mu_{jk}}),
$$
where $h'v_{jk}=c_{jk}''v_{jk}$, $\tau(h)v_{\mu_{jk}}=c_{jk}'''v_{\mu_{jk}}$ and $\tau(T_{\widehat{s}_j}^{-1})v_{\mu_{jk}}=c_{jk}'\tau(t_{\widehat{s}_j}^{-1})v_{\mu_{jk}}$, $c_{jk}'',c_{jk}'''\in \mathcal{B}$.

The element $t_{\widehat{s}_j}(z_+xb)$ has no non--trivial zero weight component by Lemma \ref{noweight} (ii) as by the assumption $b$ has no zero weight components in this case. So, since $\tau$ preserves weights, we deduce that 
$$
A_{jk}(y'h'y\tau(ht_{\widehat{s}_j}z_+xb))=c_{jk}'''c_{jk}''c_{jk}'(v_{\mu_{jk}},\tau(t_{\widehat{s}_j}(z_+xb))v_{\mu_{jk}})=0
$$
by orthogonality of different weight subspaces with respect to the contravariant form.

If $p=k$ and $y'=\omega_0 S_{\widehat{s}_j}^{-1}(f_{\delta_{jk}})$ we have to show that $A_{jk}(\omega_0S_{\widehat{s}_j}^{-1}(f_{\delta_{jk}})h'y\tau(ht_{\widehat{s}_j}z_+xb))=0$ for any $y=\omega_0 S_{\widehat{s}_j}^{-1}(f_{\delta_{D}^j}^{(n_D^j)}\ldots f_{\delta_{jk}}^{(n^j_k)})$, $n_i^j\in \mathbb{N}$, $h\in U_{\mathcal{B}}^{\widehat{s}_j,res}(\h)$, $z_+\in U_{\mathcal{B}}^{\widehat{s}_j,res}(\z_+^j), x\in U_{\mathcal{B}}^{\widehat{s}_j,res}([-\delta_{j1},-\delta_{jo_j}])$.

In this case by (\ref{oof})
$$ 
A_{jk}(\omega_0S_{\widehat{s}_j}^{-1}(f_{\delta_{jk}})h'y\tau(ht_{\widehat{s}_j}z_+xb))=(v_{\mu_{jk}},h'y\tau(ht_{\widehat{s}_j}z_+xbT_{\widehat{s}_j}^{-1})v_{\mu_{jk}})=d_{jk}c_{jk}'(v_{\mu_{jk}},y\tau(ht_{\widehat{s}_j}(z_+xb))v_{\mu_{jk}}),
$$
where $h'v_{\mu_{jk}}=d_{jk}v_{\mu_{jk}}$, $d_{jk}\in \mathcal{B}$.

Again if $y\neq 1$  then 
$$
A_{jk}(\omega_0S_{\widehat{s}_j}^{-1}(f_{\delta_{jk}})h'y\tau(ht_{\widehat{s}_j}z_+xb))=d_{jk}c_{jk}'(v_{\mu_{jk}},y\tau(ht_{\widehat{s}_j}(z_+xb))v_{\mu_{jk}})=0
$$ 
as $v_{\mu_{jk}}$ is a highest weight vector.

If $y=1$ 
$$
A_{jk}(\omega_0S_{\widehat{s}_j}^{-1}(f_{\delta_{jk}})h'y\tau(ht_{\widehat{s}_j}z_+xb))=d_{jk}c_{jk}'''c_{jk}'(v_{\mu_{jk}},\tau(t_{\widehat{s}_j}(z_+xb))v_{\mu_{jk}}).
$$
The element $t_{\widehat{s}_j}(z_+xb)$ has no non--trivial zero weight component by Lemma \ref{noweight} (ii) as by the assumption $b$ has no zero weight components in this case. Since $\tau$ preserves weights, we deduce that 
$$
A_{jk}(\omega_0S_{\widehat{s}_j}^{-1}(f_{\delta_{jk}})h'y\tau(ht_{\widehat{s}_j}z_+xb))=d_{jk}c_{jk}'''c_{jk}'(v_{\mu_{jk}},\tau(t_{\widehat{s}_j}(z_+xb))v_{\mu_{jk}})=0
$$
by orthogonality of different weight subspaces with respect to the contravariant form.

Formula (\ref{A2}) follows from (\ref{Ap}).

This completes the proof of part (i).

(ii) We show that ${{J}^{jk}_{\mathcal{B}}}'\subset {{J}^{jk}_{\mathcal{B}}}$.
Indeed, by Lemma \ref{barJ1} (i) any element from ${{J}^{j1}_{\mathcal{B}}}'$ belongs to $J^{j1}_{\mathcal{B}}$. We claim that $J^{j1}_{\mathcal{B}}\subset {{J}^{jk}_{\mathcal{B}}}$, where the same braid group element $t_{\widehat{s}_j}$ is used in the definitions of both algebras.  

Let $(w,\cdot~ v)\in J^{j1}_{\mathcal{B}}$ be one of the elements generating $J^{j1}_{\mathcal{B}}$, i.e. $(w,\tau(yt_{\widehat{s}_j}z_+x)v)=0$ for any for any $y\in U_{\mathcal{B}}^{\widehat{s}_j,res}(\b_-), z_+\in U_{\mathcal{B}}^{\widehat{s}_j,res}(\z_+^j), x\in U_{\mathcal{B}}^{\widehat{s}_j,res}([-\delta_{j1},-\delta_{jo_j}])$. Then $(w,y\tau(ht_{\widehat{s}_j}z_+x)v)=0$ for any $y\in \omega_0S_{\widehat{s}_j}^{-1}(U_{\mathcal{B}}^{\widehat{s}_j,res}([-\delta_{jk},-\delta_{jD}]))$, $h\in U_{\mathcal{B}}^{\widehat{s}_j,res}(\h)$, $z_+\in U_{\mathcal{B}}^{\widehat{s}_j,res}(\z_+^j), x\in U_{\mathcal{B}}^{\widehat{s}_j,res}([-\delta_{j1},-\delta_{jo_j}])$ as $\omega_0S_{\widehat{s}_j}^{-1}(U_{\mathcal{B}}^{\widehat{s}_j,res}([-\delta_{jk},-\delta_{jD}]))U_{\mathcal{B}}^{\widehat{s}_j,res}(\h)\subset \tau(U_{\mathcal{B}}^{\widehat{s}_j,res}(\b_-))=U_{\mathcal{B}}^{\widehat{s}_j,res}(\b_-)$. Thus $(w,\cdot~ v)\in {{J}^{jk}_{\mathcal{B}}}$ is one of the elements generating ${{J}^{jk}_{\mathcal{B}}}$.

We deduce using Lemma \ref{barJ1} (ii) that ${{J}^{j1}_{\mathcal{B}}}'\subset J^{j1}_{\mathcal{B}}\subset {{J}^{jk}_{\mathcal{B}}}$.

By part (i) $A_{jp}\in {{J}^{jk}_{\mathcal{B}}}$ for $p=1,\ldots, k-1$. Since ${{J}^{j1}_{\mathcal{B}}}'$ and $A_{jp}$ with $p<k$ generate ${{J}^{jk}_{\mathcal{B}}}'$ as a left ideal, we obtain that  ${{J}^{jk}_{\mathcal{B}}}'\subset {{J}^{jk}_{\mathcal{B}}}$.

(iii) In order to show that ${\Delta_\mu^{\widehat{s}_j}}(\tau(y')~\cdot~ \tau(b))\in {{J}^{j1}_{\mathcal{B}}}'$ we shall check that it is one of the generating elements of ${{J}^{j1}_{\mathcal{B}}}'$. More precisely, we verify that ${\Delta_\mu^{\widehat{s}_j}}(\tau(y')\cdot~ \tau(b))=\tilde{c}(v_\mu,\cdot~ \tau(bT_{\widehat{s}_p}^{-1})v_\mu)$ for some $\tilde{c}\in \mathcal{B}$, and
\begin{equation}\label{aq12}
{\Delta_\mu^{\widehat{s}_j}}(\tau(y't_{\widehat{s}_p}z_+xb))=(v_\mu, \tau(y't_{\widehat{s}_p}z_+xbT_{\widehat{s}_p}^{-1})v_\mu)=(\omega\tau(y')v_\mu, \tau(t_{\widehat{s}_p}z_+xbT_{\widehat{s}_p}^{-1})v_\mu)=0
\end{equation}
for any $z_+\in U_{\mathcal{B}}^{\widehat{s}_j,res}(\z_+^j), x\in U_{\mathcal{B}}^{\widehat{s}_j,res}([-\delta_{j1},-\delta_{jo_j}])$. 

Since $\omega$ changes the signs of the weights, and $\tau$ preserves weights, the weights of the weight components of $\omega\tau(y')$ belong to $\mathbb{N}(\Delta_+)$. Thus
$$
{\Delta_\mu^{\widehat{s}_j}}(\tau(y')~\cdot~ \tau(b))=(\omega\tau(y'_0)v_\mu, \cdot~ \tau(bT_{\widehat{s}_p}^{-1})v_\mu)
=\tilde{c}(v_\mu, \cdot~ \tau(bT_{\widehat{s}_p}^{-1})v_\mu),
$$
where $y'_0$ is the zero weight component of $y'$, $\tilde{c}\in \mathcal{B}$ is defined by $\tau(y'_0) v_\mu=\tilde{c}v_\mu$.

Next, for any $z_+\in U_{\mathcal{B}}^{\widehat{s}_j,res}(\z_+^j), x\in U_{\mathcal{B}}^{\widehat{s}_j,res}([-\delta_{j1},-\delta_{jo_j}])$ one has using the previous formula
\begin{equation}\label{aq13}
{\Delta_\mu^{\widehat{s}_j}}(\tau(y't_{\widehat{s}_p}z_+xb))=\tilde{c}(v_{\mu}, \tau(t_{\widehat{s}_j}z_+xbT_{\widehat{s}_j}^{-1})v_\mu)=\tilde{c}c''(v_\mu, \tau(t_{\widehat{s}_j}z_+xbt_{\widehat{s}_j}^{-1})v_\mu)
=\tilde{c}c''(v_\mu, \tau(t_{\widehat{s}_p}(z_+xb))v_\mu),
\end{equation}
where $c''\in \mathbb{C}q^{\mathbb{Z}}$, $c''\neq 0$ is given by Lemma \ref{bga} using the condition $\tau(T_{\widehat{s}_j}^{-1})v_\mu=c''\tau(t_{\widehat{s}_j}^{-1})v_\mu$.

If the zero weight component of $y'$ is zero then $\tilde{c}=0$. Thus (\ref{aq12}) holds. In fact ${\Delta_\mu^{\widehat{s}_j}}(\tau(y')~\cdot~ \tau(b))=0$ in this case.

By Lemma \ref{noweight} (ii), if the zero weight component of $b$ is zero then the zero weight component of $t_{\widehat{s}_p}(z_+xb)$ is zero as well, and hence the right hand side of (\ref{aq13}) vanishes as $\tau$ preserves weights, and different weight subspaces are orthogonal with respect to the contravariant form. Thus (\ref{aq12}) holds.

This completes the proof of part (iii).

(iv) Firstly, the ideal ${{J}^{j1}_{\mathcal{B}}}'$ is stable under the ${\rm Ad}_{\widehat{s}_j}^0$--action of $U_{\mathcal{B}}^{\widehat{s}_j,res}([-\delta_{j1},-\delta_{jo_j}])$ by Lemma \ref{barJ1} (iii), and hence under the ${\rm Ad}_{\widehat{s}_j}$--action of $U_{\mathcal{B}}^{\widehat{s}_j,res}([-\delta_{jk},-\delta_{jo_j}])\subset U_{\mathcal{B}}^{\widehat{s}_j,res}([-\delta_{j1},-\delta_{jo_j}])$. Therefore by part (ii) and by the definition of ${{J}^{jk}_{\mathcal{B}}}'$ we infer
$$
{\rm Ad}_{\widehat{s}_j}^0(U_{\mathcal{B}}^{\widehat{s}_j,res}([-\delta_{jk},-\delta_{jo_j}])){{J}^{j1}_{\mathcal{B}}}'\subset {{J}^{j1}_{\mathcal{B}}}'\subset {{J}^{jk}_{\mathcal{B}}}'\subset {{J}^{jk}_{\mathcal{B}}}.
$$

Since the elements $f_\beta^{(n)}$, $\beta \in [\delta_{jk},\delta_{jo_j}]$, $n\in \mathbb{N}$ generate $U_{\mathcal{B}}^{\widehat{s}_j,res}([-\delta_{jk},-\delta_{jo_j}])$, and ${{J}^{j1}_{\mathcal{B}}}'$ and $A_{jp}$ with $p<k$ generate ${{J}^{jk}_{\mathcal{B}}}'$, it remains to show that for arbitrary $f\in \mathbb{C}_\mathcal{B}^{\widehat{s}_j}[G]$, $\beta \in [\delta_{jk},\delta_{jo_j}]$, and $n\in \mathbb{N}$ one has ${\rm Ad}_{\widehat{s}_j}^0f_\beta^{(n)}(f\otimes A_{jp})\in {{J}^{jk}_{\mathcal{B}}}$ for $p<k$.

Indeed, consider formula (\ref{AdbarJ10}) with arbitrary $f\in \mathbb{C}_\mathcal{B}^{\widehat{s}_j}[G]$, $g=A_{jp}$, $p<k$, $n\in \mathbb{N}$, and $\beta \in [\delta_{jk},\delta_{jo_j}]$. By part (i), by (\ref{yh'}) and by (\ref{om0t}) the second factors in all $\otimes$--products in $\mathbb{C}_\mathcal{B}^{\widehat{s}_j}[G]$ in the terms in the right hand side of (\ref{AdbarJ10}) belong to ${{J}^{jk}_{\mathcal{B}}}$. Thus ${\rm Ad}_{\widehat{s}_j}^0f_\beta^{(n)}(f\otimes A_{jp})\in {{J}^{jk}_{\mathcal{B}}}$. This completes the proof of Lemma \ref{barJp}.

\end{proof}

We shall also need some properties of the left ideals ${{J}^{jk}_{\mathcal{B}}}$ and ${{J}^{jk}_{\mathcal{B}}}'$ with respect to the multiplication from the right.
\begin{lemma}\label{mainl2}
Assume that the element $t_{\widehat{s}_j}=\overline{T}_{\widehat{s}_j}$ is used in the definition of ${J^{jk}_{\mathcal{B}}}$.
Then for any $j\in \{1,\ldots, R-1\}$, $1\leq k\leq n_j+1$ the following statements are true.

(i) ${J^{j1}_{\mathcal{B}}}'A_{jk}\subset {J^{jk}_{\mathcal{B}}}$.

(ii) For any $\mu \in P_+$ ${{J^{j1}_{\mathcal{B}}}'}{\Delta_\mu^{\widehat{s}_j}}\subset J^{j1}_{\mathcal{B}}\subset {J^{jk}_{\mathcal{B}}}$.

(iii) For any $\mu \in P_+$ ${{J}^{jk}_{\mathcal{B}}}'{\Delta_\mu^{\widehat{s}_j}}\subset {{J}^{jk}_{\mathcal{B}}}$.
\end{lemma}
 
\begin{proof}
(i) To show that ${J^{j1}_{\mathcal{B}}}'A_{jk}\subset {J^{jk}_{\mathcal{B}}}$ it suffices to verify that for any $(u, \cdot~ v)\in \mathbb{C}_{\mathcal{B}}^{\widehat{s}_j}[G]$, where $u$ is a highest weight vector in a finite rank representation $V$ of $U_h(\g)$, and $v\in V$ is such that $(u,\tau(\overline{T}_{\widehat{s}_j}z_+x) v)=0$ for any $z_+\in U_{\mathcal{B}}^{\widehat{s}_j,res}(\z_+^j), x\in U_{\mathcal{B}}^{\widehat{s}_j,res}([-\delta_{j1},-\delta_{jo_j}])$ one has 
\begin{equation}\label{Bcond}
((u,\cdot~ v)\otimes A_{jk})(y\tau(h\overline{T}_{\widehat{s}_j}z_+x))=0
\end{equation}
for any $y\in \omega_0S_{\widehat{s}_j}^{-1}(U_{\mathcal{B}}^{\widehat{s}_j,res}([-\delta_{jk},-\delta_{jD}]))$, $h\in U_{\mathcal{B}}^{\widehat{s}_j,res}(\h)$, $z_+\in U_{\mathcal{B}}^{\widehat{s}_j,res}(\z_+^j), x\in U_{\mathcal{B}}^{\widehat{s}_j,res}([-\delta_{j1},-\delta_{jo_j}])$.

From (\ref{comults}) and the fact that $\omega_0' S_{\widehat{s}_j}^{-1}$ is an anti-coautomorphism preserving weights we obtain 
\begin{equation}\label{Dsy}
\Delta_{\widehat{s}_j}^\tau(y)=1\otimes y+\sum_i{v_i}\otimes v_i',
\end{equation} 
where the weights of the elements $v_i$ are strictly negative.

From formula (\ref{DTw}) with $k$ equal to the length of $\widehat{s}_j$, $k=l(\widehat{s}_j)$, we have
\begin{equation}\label{DTW1}
\Delta_{\widehat{s}_j}(\overline{T}_{\widehat{s}_j})=\prod^{l(\widehat{s}_j)}_{p=1}{\overline{\theta}_{\beta_p}^{\widehat{s}_j}}'q^{\sum_{i=1}^l(-Y_i\otimes K_{\widehat{s}_j}H_i+T_{\widehat{s}_j}Y_i\otimes T_{\widehat{s}_j}K_{\widehat{s}_j}H_i)}\overline{T}_{\widehat{s}_j}\otimes \overline{T}_{\widehat{s}_j}=\prod^{l(\widehat{s}_j)}_{p=1}{\overline{\theta}_{\beta_p}^{\widehat{s}_j}}'\overline{T}_{\widehat{s}_j}\otimes \overline{T}_{\widehat{s}_j}
\end{equation}
where we used the fact that $K_{\widehat{s}_j}h=\frac\kappa2{1+{\widehat{s}_j} \over 1-{\widehat{s}_j} }w_{j-1}^{-1}P_{{\h'}}w_{j-1}h$ for any $h\in \h$, so have $K_{\widehat{s}_j}s_{\widehat{s}_j}=s_{\widehat{s}_j}K_{\widehat{s}_j}$, and hence by (\ref{qt1})
\begin{equation}\label{qt1'}
q^{\sum_{i=1}^l(T_{\widehat{s}_j}Y_i\otimes T_{\widehat{s}_j}K_{\widehat{s}_j}H_i-Y_i\otimes K_{\widehat{s}_j}H_i)}=1.
\end{equation}

Since $u$ is a highest weight vector we obtain using Lemma \ref{hwv} (i), (\ref{DTW1}), the definition of ${\overline{\theta}_{\beta_p}^s}'$ after (\ref{DTw}), and (\ref{Dsy}) that the left hand side of (\ref{Bcond}) takes the form
$$
((u, \cdot~ v)\otimes A_{jk})(y\tau(h\overline{T}_{\widehat{s}_j}z_+x))=(u, \tau(h^1\overline{T}_{\widehat{s}_j}z_+^1x^1) v)(v_{jk},y\tau(h^2\overline{T}_{\widehat{s}_j}z_+^2x^2 T_{\widehat{s}_j}^{-1})v_{\mu_{jk}})=
$$
\begin{equation}\label{Bcond1}
=t(u, \tau(h^1\overline{T}_{\widehat{s}_j}z_+^1x^1) v)(v_{jk},y\tau(h^2\overline{T}_{\widehat{s}_j}z_+^2x^2 \overline{T}_{\widehat{s}_j}^{-1})v_{\mu_{jk}})=t(u, \tau(h^1\overline{T}_{\widehat{s}_j}z_+^1x^1) v)(v_{jk},y\tau(h^2\overline{T}_{\widehat{s}_j}(z_+^2x^2))v_{\mu_{jk}}),
\end{equation}
where $t$ is a non--zero multiple of a power of $q$ defined in Lemma \ref{bga} by the condition $\tau(T_{\widehat{s}_j}^{-1})v_{\mu_{jk}}=t\tau(\overline{T}_{\widehat{s}_j}^{-1})v_{\mu_{jk}}$ and we use the Sweedler notation for the coproducts $\Delta_s(z_+)=z_+^1\otimes z_+^2$, $\Delta_s(x)=x^1\otimes x^2$, $\Delta_s(h)=h^1\otimes h^2$.

By linearity and by Lemma \ref{segmPBWs} (vi) it suffices to check condition (\ref{Bcond}) for $y\in \omega_0S_{\widehat{s}_j}^{-1}(U_{\mathcal{B}}^{\widehat{s}_j,res}([-\delta_{jk},-\delta_{jD}]))$ of the form 
$$
y=\omega_0 S_s^{-1}(f_{\delta_{D}^j}^{(n_D^j)}\ldots f_{\delta_{jk}}^{(n_k^j)}), n_i^j\in \mathbb{N}.
$$ 
In this case by Lemma \ref{Apvan} (i) the second factor in the right hand side of (\ref{Bcond1}) vanishes if $y\neq 1, \omega_0 S_s^{-1}(f_{\delta_{jk}})$, and hence (\ref{Bcond}) holds.

If $y=1$ (\ref{Bcond1}) takes the form
\begin{equation}\label{Bcond2}
((u, \cdot~ v)\otimes A_{jk})(y\tau(h\overline{T}_{\widehat{s}_j}z_+x))=tt'(u, \tau(\overline{T}_{\widehat{s}_j}z_+^1x^1) v)(v_{jk},\tau(\overline{T}_{\widehat{s}_j}(z_+^2x^2))v_{\mu_{jk}}),
\end{equation}
where $t'\in \mathcal{B}$ is defined by $\tau(h^1)u\otimes \tau(h^2)v_{jk}=t'u\otimes v_{jk}$.
 
From (\ref{comults}) or (\ref{comultsn}) it follows that $x^2\in U_{\mathcal{B}}^{\widehat{s}_j,res}([-\delta_{1k},-\delta_{jo_j}])$, so the weights of the weight components of $x^2$ belong to $\mathbb{N}[-\delta_{j1},-\delta_{jo_j}]\subset \mathbb{N}(w_1\ldots w_{j-1})^{-1}(-\Delta_+^j\cup (\Delta_0)_+)$. The first line in (\ref{comultse}) implies that the weights of the weight components of $z_+^2$ belong to $\mathbb{N}(w_1\ldots w_{j-1})^{-1}((\Delta_0)_+ \cup[-\beta_{jk},-\beta_{D_0}^0])\subset \mathbb{N}(w_1\ldots w_{j-1})^{-1}(-\Delta_+^j\cup (\Delta_0)_+)$. We conclude that the weights of the weight components of $\overline{T}_{\widehat{s}_j}(z_+^2x^2)$ belong to $\mathbb{N}(w_1\ldots w_{j-1})^{-1}s(-\Delta_+^j\cup (\Delta_0)_+)$, and by Lemma \ref{jsegm} (ii) this set has empty intersection with the set $-(w_1\ldots w_{j-1})^{-1}\Delta^j$ which contains $-\delta_{jk}$. Since the weight of $v_{jk}$ is $\mu_{jk}-\delta_{jk}$ the right hand side of (\ref{Bcond2}) vanishes as different weight subspaces of $V_{\mu_{jk}}$ are orthogonal with respect to the contravariant form and $\tau$ preserves weights.

If $y=\omega_0 S_{\widehat{s}_j}^{-1}(f_{\delta_{jk}})$ then by (\ref{A2}) formula (\ref{Bcond1}) takes the form
\begin{equation}\label{Bcond3}
((u, \cdot~ v)\otimes A_{jk})(y\tau(h\overline{T}_{\widehat{s}_j}z_+x))=t(u, \tau(h^1\overline{T}_{\widehat{s}_j}z_+^1x^1) v)(v_{\mu_{jk}},\tau(h^2\overline{T}_{\widehat{s}_j}(z_+^2x^2))v_{\mu_{jk}})=
\end{equation}
$$
=tt''(u,\tau(\overline{T}_{\widehat{s}_j}z_+^1x^1) v)(v_{\mu_{jk}},\tau(\overline{T}_{\widehat{s}_j}(z_+^2x^2))v_{\mu_{jk}}),
$$
where $t''\in \mathcal{B}$ is defined by $\tau(h^1)u\otimes \tau(h^2)v_{\mu_{jk}}=t''u\otimes v_{\mu_{jk}}$.

Recall that from (\ref{comults}) or (\ref{comultsn}) it follows that $x^2\in U_{\mathcal{B}}^{\widehat{s}_j,res}([-\delta_{j1},-\delta_{jo_j}])=U_{\mathcal{B}}^{\widehat{s}_j,res}((w_1\ldots w_{j-1})^{-1}[-\beta_{j1},-\beta_c])$, and the second line in (\ref{comultse}) implies that $z_+^2\in U_{\mathcal{B}}^{\widehat{s}_j,res}((w_1\ldots w_{j-1})^{-1}[\beta_1^0,-\beta_{j-1 n_{j-1}}])U_{\mathcal{B}}^{\widehat{s}_j,res}(\h)$. Since $[-\beta_{j1},-\beta_c]$, $[\beta_1^0,-\beta_{j-1 n_{j-1}}]\subset [\beta_1^0,-\beta_c]$, and $[\beta_1^0,-\beta_c]$ is a minimal segment, the zero weight component $(z_+^2x^2)_0$ of $z_+^2x^2$ belongs to $U_{\mathcal{B}}^{\widehat{s}_j,res}(\h)$, and by (\ref{comults}) or (\ref{comultsn}), and by the second line in (\ref{comultse}) we have
$z_+^1x^1\otimes (z_+^2x^2)_0=\sum_n z_n x_n\otimes h_n$, $z_n\in U_{\mathcal{B}}^{\widehat{s}_p,res}((w_1\ldots w_{j-1})^{-1}(\Delta_0)_+), x_n\in U_{\mathcal{B}}^{\widehat{s}_p,res}((w_1\ldots w_{j-1})^{-1}[-\beta_{j1},-\beta_c])$, $h_n\in U_{\mathcal{B}}^{\widehat{s}_j,res}(\h)$. Thus (\ref{Bcond3}) takes the form
$$
((u, \cdot~ v)\otimes A_{jk})(y\tau(h\overline{T}_{\widehat{s}_j}z_+x))=tt'' (u,\tau(\overline{T}_{\widehat{s}_j}z_+^1x^1) v)(v_{\mu_{jk}}, \tau(\overline{T}_{\widehat{s}_j}(z_+^2x^2)_0)v_{\mu_{jk}})=
$$
$$
=tt''\sum_n (u,\tau(\overline{T}_{\widehat{s}_j}z_nx_n) v)(v_{\mu_{jk}}, \tau(\overline{T}_{\widehat{s}_j}h_n)v_{\mu_{jk}})=0
$$
since $(u,\tau(\overline{T}_{\widehat{s}_j}z_nx_n) v)=0$ for all $n$ by the choice of $u$ and $v$. Thus (\ref{Bcond}) holds, and the proof of part (i) is completed.

The proof of part (ii) is similar to that of part one. The same arguments are applied with $A_{jk}$ replaced by ${\Delta_\mu^{\widehat{s}_j}}$ and $k=1$. Formula (\ref{Bcond3}) will be replaced with 
$$
((u,\cdot~ v)\otimes {\Delta_\mu^{\widehat{s}_j}})(y\tau(h\overline{T}_{\widehat{s}_j}z_+x))=t\hat{t}(u,\tau(\overline{T}_{\widehat{s}_j}z_+^1x^1) v)(v_\mu, \tau(\overline{T}_{\widehat{s}_j}(z_+^2x^2))v_\mu),
$$
where $\hat{t}\in \mathcal{B}$ is defined by the condition $\Delta_{\widehat{s}_j}^\tau(y\tau(h))u\otimes v_\mu=\hat{t}u\otimes v_\mu$. The rest of the proof is repeated verbatim.

For part (iii) we recall that ${{J}^{jk}_{\mathcal{B}}}'$ is the left ideal in $\mathbb{C}_{\mathcal{B}}^{\widehat{s}_j}[G]$ generated by ${J^{j1}_{\mathcal{B}}}'$ and by $A_{jp}$, $p<k$. By part (ii) ${J^{j1}_{\mathcal{B}}}'{\Delta_\mu^{\widehat{s}_j}}\subset J^{j1}_{\mathcal{B}}\subset {J^{jk}_{\mathcal{B}}}$. 

Using commutation relations (\ref{comm1''}) we deduce that the left ideal in $\mathbb{C}_{\mathcal{B}}^{\widehat{s}_j}[G]$ generated by $A_{jp}$, $p<k$ is invariant with respect to multiplication by ${\Delta_\mu^{\widehat{s}_j}}$ from the right.
This observation together with the inclusion ${J^{j1}_{\mathcal{B}}}'{\Delta_\mu^{\widehat{s}_j}}\subset {J}^{jk}_{\mathcal{B}}$ imply the inclusion in part (iii).

\end{proof}


\section{Localizations}\label{localG}

\pagestyle{myheadings}
\markboth{CHAPTER~\thechapter.~ZHELOBENKO TYPE OPERATORS FOR Q-W--ALGEBRAS}{\thesection.~LOCALIZATIONS}

\setcounter{equation}{0}
\setcounter{theorem}{0}

In this section we introduce quantum group analogues of the ideals ${{J}^{jk}}^{loc}$ and study their properties.
Let ${{J}^{jk}_{\mathcal{B}}}^{loc}$ \index[not]{J@${J^{jk}_\mathcal{B}}^{loc}$} and ${{J^{jk}_\mathcal{B}}'}^{loc}$ \index[not]{J@${{J^{jk}_\mathcal{B}}'}^{loc}$} be the images of ${J}^{jk}_{\mathcal{B}}\otimes \mathfrak{S}_{\widehat{s}_j}^*\subset \mathbb{C}_{\mathcal{B}}^{{\widehat{s}_j}}[G]\otimes \mathbb{C}_{\mathcal{B}}^{{\widehat{s}_j},loc}[G]_0$ and of ${{J}^{jk}_{\mathcal{B}}}'\otimes \mathfrak{S}_{\widehat{s}_j}^*\subset \mathbb{C}_{\mathcal{B}}^{{\widehat{s}_j}}[G]\otimes \mathbb{C}_{\mathcal{B}}^{{\widehat{s}_j},loc}[G]_0$, respectively, in $\mathbb{C}_{\mathcal{B}}^{{\widehat{s}_j},loc}[G]$ under the projection $\mathbb{C}_{\mathcal{B}}^{{\widehat{s}_j}}[G]\otimes \mathbb{C}_{\mathcal{B}}^{{\widehat{s}_j},loc}[G]_0 \rightarrow \mathbb{C}_{\mathcal{B}}^{{\widehat{s}_j}}[G]\otimes_{\mathbb{C}_{\mathcal{B}}^{{\widehat{s}_j}}[G]_0}\mathbb{C}_{\mathcal{B}}^{{\widehat{s}_j},loc}[G]_0=\mathbb{C}_{\mathcal{B}}^{{\widehat{s}_j},loc}[G]$.

\begin{lemma}\label{locid}
Let ${{J}^{jk}_q}^{loc}={{J}^{jk}_{\mathcal{B}}}^{loc}\otimes_{\mathcal{B}} \mathbb{C}(q^{\frac{1}{d{\bar{r}}^2}})$, \index[not]{J@${J^{jk}_q}^{loc}$} ${{{J}^{jk}_q}'}^{loc}={{{J}^{jk}_{\mathcal{B}}}'}^{loc}\otimes_{\mathcal{B}} \mathbb{C}(q^{\frac{1}{d{\bar{r}}^2}})$. \index[not]{J@${{J^{jk}_q}'}^{loc}$}

Then for any choice of $t_{\widehat{s}_j}$ in the definition of ${J}^{jk}_{\mathcal{B}}$ one has ${{J}^{jk}_q}^{loc}={{{J}^{jk}_q}'}^{loc}$. Thus ${{I}^{jk}_{\mathcal{B}}}^{loc}:={{{J}^{jk}_q}'}^{loc}\cap \mathbb{C}_{\mathcal{B}}^{{\widehat{s}_j},loc}[G]={{{J}^{jk}_q}}^{loc}\cap \mathbb{C}_{\mathcal{B}}^{{\widehat{s}_j},loc}[G]$. \index[not]{I@${{I}^{jk}_{\mathcal{B}}}^{loc}$}
\end{lemma}

\begin{proof}
Let $\iota_j: \mathbb{C}_{\mathcal{B}}^{{\widehat{s}_j}}[G]\to \mathbb{C}_{\mathcal{B}}^{{\widehat{s}_j},loc}[G]$ \index[not]{i@$\iota_j$} be the canonical $\mathcal{B}$--module homomorphism. Since ${{J}^{jk}_{\mathcal{B}}}\subset {{J}^{jk}_{\mathcal{B}}}'$ by Lemma \ref{barJp} (ii), ${{{J}^{jk}_{\mathcal{B}}}}^{loc}\subset {{{J}^{jk}_{\mathcal{B}}}'}^{loc}$, and hence 
\begin{equation}\label{ioe}
\iota_j^{-1}({{{J}^{jk}_{\mathcal{B}}}}^{loc})\subset \iota_j^{-1}({{{J}^{jk}_{\mathcal{B}}}'}^{loc}).
\end{equation} 

When specializing $q^{\frac{1}{d{\bar{r}}^2}}$ to $1$ the matrix elements $A_{jk}$ become $(v_{jk}, \cdot~ \widehat{s}_j^{-1}v_{\mu_{jk}})$, and the set $\mathfrak{S}_{\widehat{s}_j}$ becomes the set $\mathcal{S}_j$ by their definitions. By Lemma \ref{barJ1} (i) ${{J^{j1}_{\mathcal{B}}}'}=J^{j1}_{\mathcal{B}}=J^{j1}$  mod $(q^{\frac{1}{d{\bar{r}}^2}}-1)J^{j1}_{\mathcal{B}}$. Therefore ${{J}^{jk}_{\mathcal{B}}}'={{J}^{jk}}'$ mod $(q^{\frac{1}{d{\bar{r}}^2}}-1)$ and ${{{J}^{jk}_{\mathcal{B}}}'}^{loc}={{{J}^{jk}}}^{loc}$ mod $(q^{\frac{1}{d{\bar{r}}^2}}-1)$ by Lemma \ref{jloc} (ii).

On the other hand ${{J}^{jk}_{\mathcal{B}}}={{J}^{jk}}$ mod $(q^{\frac{1}{d{\bar{r}}^2}}-1)$ by Lemma \ref{Jpdescr} (i), and hence ${{{J}^{jk}_{\mathcal{B}}}}^{loc}={{{J}^{jk}}}^{loc}$ mod $(q^{\frac{1}{d{\bar{r}}^2}}-1)$ by Lemma \ref{jloc} (i).

Thus ${{{J}^{jk}_{\mathcal{B}}}'}^{loc}={{{J}^{jk}_{\mathcal{B}}}}^{loc}={{{J}^{jk}}}^{loc}$ mod $(q^{\frac{1}{d{\bar{r}}^2}}-1)$, and hence $\iota_j^{-1}({{{J}^{jk}_{\mathcal{B}}}}^{loc})= \iota_j^{-1}({{{J}^{jk}_{\mathcal{B}}}'}^{loc})$ mod $(q^{\frac{1}{d{\bar{r}}^2}}-1)$. But $\iota_j^{-1}({{{J}^{jk}_{\mathcal{B}}}}^{loc})$, $\iota_j^{-1}({{{J}^{jk}_{\mathcal{B}}}'}^{loc})\subset \mathbb{C}_{\mathcal{B}}^{{\widehat{s}_j}}[G]$ which is $\mathcal{B}$--free. Thus by Lemma \ref{incl} $\iota_j^{-1}({{{J}^{jk}_{\mathcal{B}}}}^{loc})\otimes_{\mathcal{B}} \mathbb{C}(q^{\frac{1}{d{\bar{r}}^2}})=\iota_j^{-1}({{{J}^{jk}_{\mathcal{B}}}'}^{loc})\otimes_{\mathcal{B}} \mathbb{C}(q^{\frac{1}{d{\bar{r}}^2}})$, or 
\begin{equation}\label{incloc}
\iota_j^{-1}({{J}^{jk}_q}^{loc})=\iota_j^{-1}({{{J}^{jk}_q}'}^{loc}),
\end{equation}
where we denote by the same symbol the natural extension 
$$
\iota_j: \mathbb{C}_{\mathcal{B}}^{{\widehat{s}_j}}[G]\otimes_{\mathcal{B}} \mathbb{C}(q^{\frac{1}{d{\bar{r}}^2}})\to \mathbb{C}_{\mathcal{B}}^{{\widehat{s}_j},loc}[G]\otimes_{\mathcal{B}} \mathbb{C}(q^{\frac{1}{d{\bar{r}}^2}}).
$$ 

Formula (\ref{incloc}) implies that  ${{J}^{jk}_q}^{loc}={{{J}^{jk}_q}'}^{loc}$. This completes the proof.

\end{proof}

\begin{corollary}\label{corloc}
For any $j\in \{1,\ldots, R-1\}$, $1\leq k\leq n_j+1$ the following statements are true.

(i) For any $j\in \{1,\ldots, R-1\}$, $1\leq k\leq n_j+1$ one has ${I^{jm}_{\mathcal{B}}}^{loc}A_{jk}={I^{jm}_{\mathcal{B}}}^{loc}B_{jk}\subset {I^{jk}_{\mathcal{B}}}^{loc}$, $m\leq k$, where $B_{jk}={A_{jk}^0}^{-1}\otimes A_{jk}$, \index[not]{B@$B_{jk}$} $A_{jk}^0=\Delta_{\mu_{jk}}^{{\widehat{s}_j}}$. 

(ii) For any $j\in \{1,\ldots, R-1\}$, $1\leq k\leq n_j+1$, $\mu \in P_+$ one has ${{I}^{jk}_{\mathcal{B}}}^{loc}{\Delta_\mu^{\widehat{s}_j}}\subset {{I}^{jk}_{\mathcal{B}}}^{loc}$.

(iii) For any $j\in \{1,\ldots, R-1\}$, $1\leq k\leq n_j$ one has ${I^{jk}_{\mathcal{B}}}^{loc}\subset {I^{jk+1}_{\mathcal{B}}}^{loc}$.
\end{corollary}

\begin{proof}
Recall that ${{J}^{jm}_{\mathcal{B}}}'$ is the left ideal in $\mathbb{C}_{\mathcal{B}}^{\widehat{s}_j}[G]$ generated by ${J^{j1}_{\mathcal{B}}}'$ and by $A_{jp}$, $p=1,\ldots, m-1$. By Lemma \ref{mainl2} (i) for $j\in \{1,\ldots, R-1\}$, $1\leq k\leq n_j+1$ one has ${J^{j1}_{\mathcal{B}}}'A_{jk}\subset {J^{jk}_{\mathcal{B}}}$, where $t_{\widehat{s}_j}=\overline{T}_{\widehat{s}_j}$ is used in the definition of ${J^{jk}_{\mathcal{B}}}$. Now from the first two inclusions in (\ref{ABS}) with $s=\widehat{s}_j$, $A_p=A_{jk}$, $B_p=B_{jk}$ we deduce that ${{J^{j1}_{\mathcal{B}}}'}^{loc}A_{jk}\subset {J^{jk}_{\mathcal{B}}}^{loc}$ and  ${{J^{j1}_{\mathcal{B}}}'}^{loc}B_{jk}\subset {J^{jk}_{\mathcal{B}}}^{loc}$. By Lemma \ref{locid} the inclusion in part (i) follows from these inclusions and from Lemma \ref{rightloc} with $s=\widehat{s}_j$, $\overline{U}_p$ being the algebra generated by $A_{jn}$, $n<m$.

Part (ii) follows the definition of ${{I}^{jk}_{\mathcal{B}}}^{loc}$.

Pert (iii) follows from the obvious inclusion ${J^{jk}_{\mathcal{B}}}'\subset {J^{jk+1}_{\mathcal{B}}}'$ and from Lemma \ref{locid}

\end{proof}

For $j\in \{1,\ldots, R-1\}$, $1\leq k\leq n_j+1$ let
$$
\mathbb{C}_{jk}^{loc}[G]=\mathbb{C}_{\mathcal{B}}^{\widehat{s}_j,loc}[G]/{{I}^{jk}_{\mathcal{B}}}^{loc}. \index[not]{C@$\mathbb{C}_{jk}^{loc}[G]$}
$$
$\mathbb{C}_{jk}^{loc}[G]$ is naturally a left $\mathbb{C}_{\mathcal{B}}^{\widehat{s}_j}[G]$--module. 

From Corollary \ref{corloc} (i) and (ii) we deduce the following statement.
\begin{lemma}\label{quotact}
For $j\in \{1,\ldots, R-1\}$, $1\leq k\leq n_j$ multiplication from the right on $\mathbb{C}_{\mathcal{B}}^{\widehat{s}_j,loc}[G]$ induces a natural action of ${\Delta_\mu^{\widehat{s}_j}}\in \mathfrak{S}_{\widehat{s}_j}$, $\mu \in P_+$ on $\mathbb{C}_{jk}^{loc}[G]$ commuting with the left $\mathbb{C}_{\mathcal{B}}^{\widehat{s}_j}[G]$--action, and for $k\geq m$, right multiplication by $A_{jk}$ and by $B_{jk}$ gives rise to well--defined homomorphisms of left $\mathbb{C}_{\mathcal{B}}^{\widehat{s}_j}[G]$--modules $\mathbb{C}_{jm}^{loc}[G]\rightarrow \mathbb{C}_{jk}^{loc}[G]$.
\end{lemma}

In the next two lemmas we study the properties of the adjoint action. These properties will be needed to study properties of quantum analogues of the operators $\Pi_{jk}$ defined in (\ref{Pip}).  
\begin{lemma}\label{Adloc}
(i) For any $j=1,\ldots, R-1$, $n\in \mathbb{N}$, $\beta\in [\delta_{j1},\delta_{jo_j}]$ and $f\in \mathbb{C}_{\mathcal{B}}^{\widehat{s}_j}[G]$ we have  
\begin{equation}\label{admu}
{\rm Ad}_{\widehat{s}_j}^0 f_{\beta}^{(n)} (f\otimes {\Delta_\mu^{\widehat{s}_j}})=q^{-n\left\langle (2K_{\widehat{s}_j}-id)\beta^\vee, \mu^\vee\right\rangle}{\rm Ad}_{\widehat{s}_j}^0 f_{\beta}^{(n)} (f)\otimes {\Delta_\mu^{\widehat{s}_j}}	
\end{equation}
in $\mathbb{C}_{\mathcal{B}}^{\widehat{s}_j}[G]/{{J}^{j1}_{\mathcal{B}}}'$. 

(ii) The adjoint action of $U_{\mathcal{B}}^{{\widehat{s}_j},res}([-\delta_{j1},-\delta_{jo_j}])$ on $\mathbb{C}_{\mathcal{B}}^{\widehat{s}_j}[G]/{{J}^{j1}_{\mathcal{B}}}'$ defined in Lemma \ref{barJ1} (iii) induces an action on $\mathbb{C}_{j1}^{loc}[G]$ satisfying
\begin{equation}\label{Apo}
{\rm Ad}_{\widehat{s}_j}^0 f_{\beta}^{(n)} (f\otimes {\Delta_\mu^{\widehat{s}_j}}^{-1})=q^{n\left\langle (2K_{\widehat{s}_j}-id)\beta^\vee, \mu^\vee\right\rangle}{\rm Ad}_{\widehat{s}_j}^0 f_{\beta}^{(n)} (f)\otimes {\Delta_\mu^{\widehat{s}_j}}^{-1}, f\in \mathbb{C}_{j1}^{loc}[G].
\end{equation}
This action is locally finite.
\end{lemma}

\begin{proof}
The proof of part (i) follows from formula (\ref{ad0n}) applied with $g={\Delta_\mu^{\widehat{s}_j}}$. By Lemma \ref{barJp} (iii) all terms in the right hand side of (\ref{ad0n}) belong to ${{J}^{j1}_{\mathcal{B}}}'$, except for the term in the first sum which corresponds to $p=n$, $k=0$. It gives the right hand side of (\ref{admu}).

Part (ii) follows from Lemma \ref{locid} and from part (i).

Local finiteness of the action of $U_{\mathcal{B}}^{{\widehat{s}_j},res}([-\delta_{j1},-\delta_{jo_j}])$ on $\mathbb{C}_{j1}^{loc}[G]$ follows from the local finiteness of the action of $U_{\mathcal{B}}^{{\widehat{s}_j},res}([-\delta_{j1},-\delta_{jo_j}])$ on $\mathbb{C}_{\mathcal{B}}^{\widehat{s}_j}[G]$ and from formula (\ref{Apo}).

\end{proof}

\begin{lemma}\label{mainl3}
The ${\rm Ad}_{\widehat{s}_j}^0$--action of $U_{\mathcal{B}}^{\widehat{s}_j,res}([-\delta_{jk},-\delta_{jo_j}])$ on $\mathbb{C}_{\mathcal{B}}^{\widehat{s}_j}[G]$ induces a locally finite action on $\mathbb{C}_{jk}^{loc}[G]=\mathbb{C}_{\mathcal{B}}^{\widehat{s}_j,loc}[G]/{{I}^{jk}_{\mathcal{B}}}^{loc}$ satisfying (\ref{Apo}) for any $f\in \mathbb{C}_{jk}^{loc}[G]$.
\end{lemma}

\begin{proof}
By Lemma \ref{barJp} (iv) ${\rm Ad}_{\widehat{s}_j}^0(U_{\mathcal{B}}^{\widehat{s}_j,res}([-\delta_{jk},-\delta_{jo_j}]){{J}^{jk}_{\mathcal{B}}}'\subset {{J}^{jk}_{\mathcal{B}}}$ for any choice of $t_{\widehat{s}_j}$ in the definition of ${J}^{jk}_{\mathcal{B}}$. Since ${{J}^{j1}_{\mathcal{B}}}'\subset {{J}^{jk}_{\mathcal{B}}}'$ the adjoint action on $\mathbb{C}_{\mathcal{B}}^{\widehat{s}_j}[G]/{{J}^{j1}_{\mathcal{B}}}'$ defined in Lemma \ref{Adloc} (i) satisfies the property ${\rm Ad}_{\widehat{s}_j}^0(U_{\mathcal{B}}^{\widehat{s}_j,res}([-\delta_{jk},-\delta_{jo_j}])({{J}^{jk}_{\mathcal{B}}}'/{{J}^{j1}_{\mathcal{B}}}')\subset {{J}^{jk}_{\mathcal{B}}}/{{J}^{j1}_{\mathcal{B}}}'$. Hence by Lemma \ref{locid}  the action of $U_{\mathcal{B}}^{\widehat{s}_j,res}([-\delta_{jk},-\delta_{jo_j}])$ on $\mathbb{C}_{j1}^{loc}[G]$ defined in part (ii) of Lemma \ref{Adloc} satisfies 
$$
{\rm Ad}_{\widehat{s}_j}^0(U_{\mathcal{B}}^{\widehat{s}_j,res}([-\delta_{jk},-\delta_{jo_j}])({{I}^{jk}_{\mathcal{B}}}^{loc}/{{I}^{j1}_{\mathcal{B}}}^{loc})\subset {{I}^{jk}_{\mathcal{B}}}^{loc}/{{I}^{j1}_{\mathcal{B}}}^{loc}.
$$ 
Thus this action induces an action on $\mathbb{C}_{jk}^{loc}[G]=\mathbb{C}_{\mathcal{B}}^{\widehat{s}_j,loc}[G]/{{I}^{jk}_{\mathcal{B}}}^{loc}$. This action is locally finite as ${\rm Ad}_{\widehat{s}_j}^0$--action of $U_{\mathcal{B}}^{\widehat{s}_j,res}([-\delta_{jk},-\delta_{jo_j}]$ on $\mathbb{C}_{\mathcal{B}}^{\widehat{s}_j}[G]$ is locally finite.

\end{proof}

Now, in view of formula (\ref{Pi1}), we relate the adjoint action on $\mathbb{C}_{j+1 1}^{loc}[G]$ and on $\mathbb{C}_{jn_j+1}^{loc}[G]$. We start with the following proposition.
\begin{proposition}\label{phijprop}
For $j=1,\ldots, R-2$ the morphism of $\mathbb{C}[[h]]$--modules $\tau T_{w_j}^{-1}\psi^{{\theta_{w_j}^{\widehat{s}_j}}^{-1}}\tau : U_h^{\widehat{s}_j}(\g)\to U_h^{\widehat{s}_{j+1}}(\g)$ induces an invertible morphism of $\mathcal{B}$--modules $\varpi_j: \mathbb{C}_{\mathcal{B}}^{\widehat{s}_{j+1}}[G]\to\mathbb{C}_{\mathcal{B}}^{\widehat{s}_j}[G]$ \index[not]{p@$\varpi_j$} which is defined by 
\begin{equation}\label{phijdef}
(\varpi_j f)(x)=f(\tau T_{w_j}^{-1}\psi^{{\theta_{w_j}^{\widehat{s}_j}}^{-1}}\tau(x)), x\in U_{\mathcal{B}}^{\widehat{s}_j,res}(\g), f\in \mathbb{C}_{\mathcal{B}}^{\widehat{s}_{j+1}}[G]
\end{equation}
and satisfies
\begin{equation}\label{phiAd}
\varpi_j{\rm Ad}^{\widehat{s}_{j+1}}(T_{w_j}^{-1}x)={\rm Ad}^{\widehat{s}_j}(x)\varpi_j, x\in U_\mathcal{B}^{\widehat{s}_j,res}(\g),
\end{equation} 
\begin{equation}\label{phipr}
\varpi_j(f\otimes g)(~\cdot~,~\cdot~)
=\sum_{m,n}({\rm Ad}^{\widehat{s}_j}(d_m^1c_n)\varpi_j f)(~\cdot~ )\otimes g(\tau T_{w_j}^{-1}(v_{\widehat{s}_j}S_{\widehat{s}_j}(c_m)\tau(~\cdot~) d_m^2d_n)), f,g\in \mathbb{C}_{\mathcal{B}}^{\widehat{s}_{j+1}}[G],
\end{equation}
where ${\theta_{w_j}^{\widehat{s_j}}}^{-1}=\sum_ma_m\otimes b_m$, $v_{\widehat{s}_j}=\sum_ma_mS_{\widehat{s}_j}(b_m)$, \index[not]{v@$v_{\widehat{s}_j}$} $\theta_{w_j}^{\widehat{s}_j}=\sum_n c_n\otimes d_n$, $\Delta_{\widehat{s}_j} d_m=d_m^1\otimes d_m^2$.

The inverse to $\varpi_j$ is given by
\begin{equation}\label{phijdefinv}
(\varpi_j^{-1} f)(x)=f(\tau(\psi^{{\theta_{w_j}^{\widehat{s}_j}}^{-1}})^{-1}T_{w_j}\tau(x)), x\in U_{\mathcal{B}}^{\widehat{s}_{j+1},res}(\g), f\in \mathbb{C}_{\mathcal{B}}^{\widehat{s}_j}[G].
\end{equation}
\end{proposition}

\begin{proof}
Recall that by (\ref{DTs}) for the reduced decomposition $w_j=s_{n_{j1}}\ldots s_{n_{jn_{j}}}$ and $T_{w_j}=T_{n_{j1}}\ldots T_{n_{jn_j}}$ we defined 
\begin{equation}\label{thj}
\theta_{w_j}^{\widehat{s}_j}=\prod^{n_j}_{k=1}\theta_{\delta_{jk}}^{\widehat{s}_j},
\end{equation}
where in the products $\theta_{\delta_{jk}}^{\widehat{s}_j}$ appears on the left from $\theta_{\delta_{jm}}^{\widehat{s}_j}$ if $k<m$, and 
$$
\theta_{\delta_{jk}}^{\widehat{s}_j}={\exp}_{q_{\delta_{jk}}}[(1-q_{\delta_{jk}}^{-2})e_{\delta_{jk}}e^{-2K_{\widehat{s}_j} \delta_{jk}^\vee}\otimes f_{\delta_{jk}}], {\theta_{\delta_{jk}}^{\widehat{s}_j}}^{-1}={\exp}_{q_{\delta_{jk}}^{-1}}[(1-q_{\delta_{jk}}^{2})e_{\delta_{jk}}e^{-2K_{\widehat{s}_j} \delta_{jk}^\vee}\otimes f_{\delta_{jk}}],
$$
$$
e_{\delta_{jk}}=(X_{\delta_{jk}}^+e^{hK_{\widehat{s}_j}\delta_{jk}^\vee}),
f_{\delta_{jk}}=(e^{-hK_{\widehat{s}_j}\delta_{jk}^\vee}X_{\delta_{jk}}^-), \delta_{jk}=s_{n_{j1}}\ldots s_{n_{jk-1}}\alpha_{i_{jk}},
$$
$$
X_{\delta_{jk}}^\pm=T_{n_{j1}}\ldots T_{n_{jk-1}}X_{i_{jk}}^\pm.
$$

Observe that for $f\in \mathbb{C}_{\mathcal{B}}^{\widehat{s}_{j+1}}[G]$ one has $f(\tau T_{w_j}^{-1}\tau(~\cdot~))\in \mathbb{C}_{\mathcal{B}}^{\widehat{s}_{j+1}}[G]$ as the braid group acts on finite rank $U_{\mathcal{B}}^{\widehat{s}_{j+1},res}(\g)$--modules.  

Formula (\ref{ps}) for $\psi^{{\theta_{w_j}^{\widehat{s}_j}}^{-1}}$ is expressed in terms of $\theta_{w_j}^{\widehat{s}_j}$ and of $v_{\widehat{s}_j}$ the definitions of which contain infinite series in divided powers of quantum root vectors (see (\ref{thj})). Hence the same observation applies to $\tau\psi^{{\theta_{w_j}^{\widehat{s}_j}}^{-1}}\tau$. Only finitely many divided powers of quantum root vectors act on a given finite rank $U_{\mathcal{B}}^{\widehat{s}_{j+1},res}(\g)$--module in a non--trivial way. Therefore only finitely many terms in the infinite series will contribute to formula (\ref{phijdef}).

The elements $e^{2K_{\widehat{s}_j} \delta_{jk}^\vee}$ which appear in the formula for $\tau(\theta_{w_j}^{\widehat{s}_j})$ also act on finite rank $U_{\mathcal{B}}^{\widehat{s}_{j+1},res}(\g)$--modules as such modules are direct sums of their weight subspaces, and $e^{2K_{\widehat{s}_j} \delta_{jk}^\vee}$ acts on a subspace of weight $\lambda$ by multiplication by $e^{\lambda(2K_{\widehat{s}_j} \delta_{jk}^\vee)}\in \mathcal{B}$.

We conclude that (\ref{phijdef}) defines a morphism of $\mathcal{B}$--modules $\varpi_j: \mathbb{C}_{\mathcal{B}}^{\widehat{s}_{j+1}}[G]\to\mathbb{C}_{\mathcal{B}}^{\widehat{s}_j}[G]$.

Formula (\ref{phiAd}) follows from (\ref{adstt}) by applying conjugation by $\tau$, and (\ref{phipr}) immediately follows from (\ref{w-1sw}) with $s=\widehat{s}_j$, $w=w_j$ and from the relation
$$
(f\otimes g)(\tau(x))=f(\tau(x^1))g(\tau(x^2)),~f,g\in \mathbb{C}_{\mathcal{B}}^w[G],~w\in W,~ x\in U_{\mathcal{B}}^{w,res}(\g),~ \Delta_w(x)=x^1\otimes x^2,
$$
with $w=\widehat{s}_{j+1}$ and $w=\widehat{s}_j$, which is a consequence of the definition of the multiplication on $\mathbb{C}_{\mathcal{B}}^w[G]$ induced by $\Delta_w^\tau(~\cdot~)=(\tau\otimes \tau)\Delta_w(\tau^{-1}(~\cdot~))$.

\end{proof}

To show that $\varpi_j$ gives rise to a morphism $\mathbb{C}_{j+1 1}^{loc}[G]\to \mathbb{C}_{jn_j+1}^{loc}[G]$ we need the following lemma.
\begin{lemma}\label{phijlem}
(i) For any choice of the element $t_{\widehat{s}_j}$ in the definition of ${{J}^{jn_j+1}_{\mathcal{B}}}$ and for $j=1,\ldots, R-2$ one has $\varpi_j({{J}^{j+11}_{\mathcal{B}}}')\subset {{J}^{jn_j+1}_{\mathcal{B}}}$.

(ii) Assume that $t_{\widehat{s}_j}=T_{\widehat{s}_j}$ in the definition of ${{J}^{jn_j+1}_{\mathcal{B}}}$. Then for $j=1,\ldots, R-2$, $f\in \mathbb{C}_{\mathcal{B}}^{\widehat{s}_{j+1}}[G]$ one has $\varpi_j(f\otimes {\Delta_\mu^{\widehat{s}_{j+1}}})=(\varpi_j f)\otimes T_{w_j}^{-1}({\Delta_\mu^{\widehat{s}_{j+1}}})$ mod ${{J}^{jn_j+1}_{\mathcal{B}}}$, where $T_{w_j}^{-1}({\Delta_\mu^{\widehat{s}_{j+1}}})(~\cdot~ )={\Delta_\mu^{\widehat{s}_{j+1}}}(\tau T_{w_j}^{-1}\tau(~\cdot~) )$. \index[not]{D@$T_{w_j}^{-1}({\Delta_\mu^{\widehat{s}_{j+1}}})$}

(iii) Assume again that $t_{\widehat{s}_j}=T_{\widehat{s}_j}$ in the definition of ${{J}^{jn_j+1}_{\mathcal{B}}}$. If $\mu, \nu \in P_+$ are such that $w_j\mu+\nu\in P_+$ and $j=1,\ldots, R-2$ then $T_{w_j}^{-1}({\Delta_\mu^{\widehat{s}_{j+1}}})\otimes {\Delta_\nu^{\widehat{s}_{j}}}=\kappa_j^\mu{\Delta_{w_j\mu+\nu}^{\widehat{s}_{j}}}$ mod ${{J}^{jn_j+1}_{\mathcal{B}}}$, where $\kappa_j^\mu\in \mathcal{B}^*$ is defined by the condition $\tau(T_{w_j}T_{\widehat{s}_{j+1}}^{-1})v_\mu=\kappa_j^\mu \tau(T_{\widehat{s}_j}^{-1}T_{w_j})v_\mu$ using Lemma \ref{bga}. \index[not]{k@$\kappa_j^\mu$}
\end{lemma}

\begin{proof}
(i) By the definition of ${{J}^{jn_j}_{\mathcal{B}}}$ and of ${J^{j+11}_{\mathcal{B}}}'$ it suffices to show that for any $f\in \mathbb{C}_{\mathcal{B}}^{\widehat{s}_{j+1}}[G]$ and $g(~\cdot~)=(u,\cdot~ v)\in \mathbb{C}_{\mathcal{B}}^{\widehat{s}_{j+1}}[G]$ such that $u$ is a highest weight vector in a finite rank representation $V$ of $U_h(\g)$, and $v\in V$ satisfies $(u,\tau(t_{\widehat{s}_{j+1}}z_+'x') v)=0$ for any  $z_+'\in U_{\mathcal{B}}^{\widehat{s}_{j+1},res}(\z_+^{j+1})$, $x'\in U_{\mathcal{B}}^{\widehat{s}_{j+1},res}([-\delta_{j+11},-\delta_{j+1}])$ and some element $t_{\widehat{s}_{j+1}}$ of the braid group acting on $\h\subset U_h(\h)$ in the same way as $\widehat{s}_{j+1}$ one has 
$$
\varpi_j(f\otimes g)(w \otimes y\tau(ht_{\widehat{s}_j}z_+x))=0
$$
for any $w\in U_{\mathcal{B}}^{\widehat{s}_j,res}(\g)$, $y\in \omega_0S_{\widehat{s}_j}^{-1}(U_{\mathcal{B}}^{\widehat{s}_j,res}([-\delta_{jn_j+1},-\delta_{jD}]))$, $h\in U_{\mathcal{B}}^{\widehat{s}_j,res}(\h)$, $z_+\in U_{\mathcal{B}}^{\widehat{s}_j,res}(\z_+^j)$, $x\in U_{\mathcal{B}}^{\widehat{s}_j,res}([-\delta_{j1},-\delta_{jo_j}])$, and for the element $t_{\widehat{s}_j}$ of the braid group used in the definition of ${{J}^{jn_j+1}_{\mathcal{B}}}$ and acting on $\h\subset U_h(\h)$ in the same way as $\widehat{s}_j$.

Indeed, from the explicit formula for $v_{\widehat{s}_j}$ and for $\theta_{w_j}^{\widehat{s}_j}$ and from the fact that $S_{\widehat{s}_j}$ and $\tau$ preserve weights, and for $k=1,\ldots n_j$ $T_{w_j}^{-1}(e_{\delta_{jk}})$ has strictly negative weights it follows by Lemma \ref{hwv} (i) and by formula (\ref{phipr}) that
$$
\varpi_j(f\otimes g)(w \otimes y\tau(ht_{\widehat{s}_j}z_+x))=\sum_{m,n}({\rm Ad}^{\widehat{s}_j}(d_m^1c_n)\varpi_j f)(w)g(\tau T_{w_j}^{-1}(v_{\widehat{s}_j}S_{\widehat{s}_j}(c_m)\tau(y)ht_{\widehat{s}_j}z_+x d_m^2d_n))=
$$
$$
=\sum_{m,n}({\rm Ad}^{\widehat{s}_j}(d_m^1c_n)\varpi_j f)(w)(u,\tau T_{w_j}^{-1}(v_{\widehat{s}_j}S_{\widehat{s}_j}(c_m)\tau(y)ht_{\widehat{s}_j}z_+x d_m^2d_n)v)=\sum_n({\rm Ad}^{\widehat{s}_j}(c_n)\varpi_j f)(w)(u,\tau T_{w_j}^{-1}(\tau(y)ht_{\widehat{s}_j}z_+x d_n)v).
$$

Since ${w_j}^{-1}[-\delta_{jn_j+1},-\delta_{jD}]\subset \Delta_-$ Lemma \ref{hwv} (i) implies that the right hand side of the last formula takes the form
\begin{equation}\label{pfg}
\varpi_j(f\otimes g)(w\otimes y\tau(ht_{\widehat{s}_j}z_+x))=y_0\sum_n({\rm Ad}^{\widehat{s}_j}(c_n)\varpi_j f)(w)(u,\tau T_{w_j}^{-1}(ht_{\widehat{s}_j}z_+x d_n)v),
\end{equation}
where $y_0\in \mathcal{B}$ is the zero weight component of $y$. 

Next, $x d_n\in U_{\mathcal{B}}^{\widehat{s}_j,res}([-\delta_{j1},-\delta_{jo_j}])$, therefore by Lemma \ref{segmPBWs} (ix) we can write $xd_n=\sum_px_{np}x_{np}'+y_n$, where $x_{np}\in U_{\mathcal{B}}^{\widehat{s}_j,res}([-\delta_{j1},-\delta_{jn_j}])$ have strictly negative weights, $y_n,x_{np}'\in U_{\mathcal{B}}^{\widehat{s}_j,res}([-\delta_{jn_j+1},-\delta_{jo_j}])$,   and (\ref{pfg}) can be rewritten as follows
\begin{equation}\label{pfg1}
\varpi_j(f\otimes g)(w\otimes y\tau(ht_{\widehat{s}_j}z_+x))=y_0\sum_{p,n}({\rm Ad}^{\widehat{s}_j}(c_n)\varpi_j f)(w)(u,\tau T_{w_j}^{-1}(ht_{\widehat{s}_j}z_+x_{np}x_{np}')v)+
\end{equation}
$$
+y_0\sum_n({\rm Ad}^{\widehat{s}_j}(c_n)\varpi_j f)(w)(u,\tau T_{w_j}^{-1}(ht_{\widehat{s}_j}z_+y_n)v)=
$$
$$
=y_0\sum_{p,n}({\rm Ad}^{\widehat{s}_j}(c_n)\varpi_j f)(w)(u,\tau T_{w_j}^{-1}(ht_{\widehat{s}_j}(z_+x_{np}))\tau T_{w_j}^{-1}(t_{\widehat{s}_j}x_{np}')v)+y_0c\sum_n({\rm Ad}^{\widehat{s}_j}(c_n)\varpi_j f)(w)(u,\tau (t_{\widehat{s}_{j+1}})\tau T_{w_j}^{-1}(z_+y_n)v),
$$
where $c\in \mathcal{B}$ is defined by the condition $\tau T_{w_j}^{-1}(h)u=cu$ and we choose $t_{\widehat{s}_{j+1}}$ as follows $t_{\widehat{s}_{j+1}}=T_{w_j}^{-1}t_{\widehat{s}_{j}}T_{w_j}$.

The weights of the weight components of $t_{\widehat{s}_j}(z_+x_{np})$ belong to $\mathbb{N}\widehat{s}_j((w_1\ldots w_{j-1})^{-1}(\Delta_0)_+\cup [-\delta_{j1},-\delta_{jn_j}])=\mathbb{N}(w_1\ldots w_{j-1})^{-1}(s(-\Delta^j\cup (\Delta_0)_+))$ and by the choice of $x_{np}$ they do not belong to $\mathbb{N}(w_1\ldots w_{j-1})^{-1}(\Delta_0)_+$. Since by Lemma \ref{jsegm} (i) $s(\Delta^j)\subset \Delta_+^{j+1}$, and $(w_1\ldots w_{j})^{-1}\Delta_+^{j+1}=\Delta_+$ we deduce that the weights of the weight components of $T_{w_j}^{-1}(ht_{\widehat{s}_j}(z_+x_{np}))$ are non--zero; they belong to $\mathbb{N}(\Delta_-\cup (w_1\ldots w_{j})^{-1}(\Delta_0)_+)$ and do not belong to $\mathbb{N}(w_1\ldots w_{j})^{-1}(\Delta_0)_+\subset \Delta_+$. In particular, these weights are not non--negative, and hence by Lemma \ref{hwv} (i) the first term in the left hand side of (\ref{pfg1}) vanishes as $u$ is a highest weight vector and $\tau$ preserves weights.

For the second term we note that by the definition $$T_{w_j}^{-1}(U_{\mathcal{B}}^{\widehat{s}_{j+1},res}(\z_+^{j}))=U_{\mathcal{B}}^{\widehat{s}_{j+1},res}(\z_+^{j+1}),T_{w_j}^{-1}(U_{\mathcal{B}}^{\widehat{s}_j,res}([-\delta_{jn_j+1},-\delta_{jo_j}]))=U_{\mathcal{B}}^{\widehat{s}_{j+1},res}([-\delta_{j+11},-\delta_{j+1}]),$$ so that $T_{w_j}^{-1}(z_+y_n)\in U_{\mathcal{B}}^{\widehat{s}_{j+1},res}(\z_+^{j+1})U_{\mathcal{B}}^{\widehat{s}_{j+1},res}([-\delta_{j+11},-\delta_{j+1}])$.
Therefore the second term in the right hand side of (\ref{pfg1}) vanishes by the choice of $u$ and $v$. This completes the proof of part (i)

(ii)
By the definition of ${{J}^{jn_j+1}_{\mathcal{B}}}$ it suffices to show that for any $f\in \mathbb{C}_{\mathcal{B}}^{\widehat{s}_{j+1}}[G]$ one has 
$$
\varpi_j(f\otimes {\Delta_\mu^{\widehat{s}_{j+1}}})(w \otimes y\tau(hT_{\widehat{s}_j}z_+x))=(\varpi_j f)(w) {\Delta_\mu^{\widehat{s}_{j+1}}}
(\tau T_{w_j}^{-1}\tau (y\tau(hT_{\widehat{s}_j}z_+x)))=
$$
$$
=(\varpi_j f)(w) {\Delta_\mu^{\widehat{s}_{j+1}}}
(\tau T_{w_j}^{-1}(\tau (y)hT_{\widehat{s}_j}z_+x))
$$
for any $w\in U_{\mathcal{B}}^{\widehat{s}_j,res}(\g)$, $y\in \omega_0S_{\widehat{s}_j}^{-1}(U_{\mathcal{B}}^{\widehat{s}_j,res}([-\delta_{jn_j+1},-\delta_{jD}]))$, $h\in U_{\mathcal{B}}^{\widehat{s}_j,res}(\h)$, $z_+\in U_{\mathcal{B}}^{\widehat{s}_j,res}(\z_+^j)$, $x\in U_{\mathcal{B}}^{\widehat{s}_j,res}([-\delta_{j1},-\delta_{jo_j}])$.

Similarly to (\ref{pfg}) we obtain, recalling the definition of ${\Delta_\mu^{\widehat{s}_{j+1}}}$, that
\begin{equation}\label{pfg2}
\varpi_j(f\otimes {\Delta_\mu^{\widehat{s}_{j+1}}})(w\otimes y\tau(hT_{\widehat{s}_j}z_+x))=y_0\sum_n({\rm Ad}^{\widehat{s}_j}(c_n)\varpi_j f)(w)(v_\mu,\tau T_{w_j}^{-1}(hT_{\widehat{s}_j}z_+x d_n)\tau(T_{\widehat{s}_{j+1}}^{-1})v_\mu)=
\end{equation}
$$
=y_0\kappa_j^\mu\sum_n({\rm Ad}^{\widehat{s}_j}(c_n)\varpi_j f)(w)(v_\mu,\tau T_{w_j}^{-1}(hT_{\widehat{s}_j}(z_+x d_n))v_\mu),
$$
where $y_0\in \mathcal{B}$ is the zero weight component of $y$, and $\kappa_j^\mu\in \mathcal{B}^*$ is chosen in such a way that $\tau(T_{w_j}T_{\widehat{s}_{j+1}}^{-1})v_\mu=\kappa_j^\mu \tau(T_{\widehat{s}_{j}}^{-1}T_{w_j})v_\mu$ using Lemma \ref{bga}.

Note that only the zero weight component of $z_+x d_n$ can contribute to the right hand side of (\ref{pfg2}).
Since $z_+xd_n\in U_{\mathcal{B}}^{\widehat{s}_j,res}(\z_+^j)U_{\mathcal{B}}^{\widehat{s}_j,res}([-\delta_{j1},-\delta_{jo_j}])$, the weights of the weight components of $z_+xd_n$ belong to $$\mathbb{N}((w_1\ldots w_{j-1})^{-1}(\Delta_0)_+\cup [-\delta_{j1},-\delta_{jo_j}]=\mathbb{N}(w_1\ldots w_{j-1})^{-1}((\Delta_0)_+\cup [-\beta_{j1},-\beta_c]).$$ Observe that $(\Delta_0)_+\cup [-\beta_{j1},-\beta_c]\subset [\beta^0_1,-\beta_c]$ which is a minimal segment. Therefore only the product of the zero weight components of $z_+$, $x$ and $d_n$ can contribute to the right hand side of (\ref{pfg2}). From formula (\ref{thj}) for $\theta_{w_j}^{\widehat{s}_j}=\sum_n c_n\otimes d_n$ it follows that only one term with $d_n=1$ has this property, and the corresponding $c_n=1$ as well. We conclude that (\ref{pfg2}) takes the form
$$
\varpi_j(f\otimes {\Delta_\mu^{\widehat{s}_{j+1}}})(w\otimes y\tau (hT_{\widehat{s}_j}z_+x))=y_0\kappa_j^\mu(\varpi_j f)(w)(v_\mu,\tau T_{w_j}^{-1}(hT_{\widehat{s}_j}(z_+x))v_\mu).
$$

Finally since ${w_j}^{-1}[-\delta_{jn_j+1},-\delta_{jD}]\subset \Delta_-$ and $y\in \omega_0S_{\widehat{s}_j}^{-1}(U_{\mathcal{B}}^{\widehat{s}_j,res}([-\delta_{jn_j+1},-\delta_{jD}]))$,  Lemma \ref{hwv} (i) implies that the right hand side of the last formula takes the form
$$
\varpi_j(f\otimes {\Delta_\mu^{\widehat{s}_{j+1}}})(w\otimes y\tau(hT_{\widehat{s}_j}z_+x))=(\varpi_j f)(w)(v_\mu,\tau T_{w_j}^{-1}(\tau(y)hT_{\widehat{s}_j}z_+x)\tau(T_{\widehat{s}_{j+1}}^{-1})v_\mu)=
$$
$$
=(\varpi_j f)(w) {\Delta_\mu^{\widehat{s}_{j+1}}}
(\tau T_{w_j}^{-1}\tau (y\tau (hT_{\widehat{s}_j}z_+x))).
$$
This completes the proof of part (ii).

(iii) 
By the definition of ${{J}^{jn_j}_{\mathcal{B}}}$ it suffices to show that
\begin{equation}\label{dmus}
({\Delta_\mu^{\widehat{s}_{j+1}}}(\tau T_{w_j}^{-1}\tau(~\cdot~))\otimes {\Delta_\nu^{\widehat{s}_{j}}})(y\tau(hT_{\widehat{s}_j}z_+x))=\kappa_j^\mu {\Delta_{w_j\mu+\nu}^{\widehat{s}_{j}}}(y\tau(hT_{\widehat{s}_j}z_+x))
\end{equation}
for any $w\in U_{\mathcal{B}}^{\widehat{s}_j,res}(\g)$, $y\in \omega_0S_{\widehat{s}_j}^{-1}(U_{\mathcal{B}}^{\widehat{s}_j,res}([-\delta_{jn_j+1},-\delta_{jD}]))$, $h\in U_{\mathcal{B}}^{\widehat{s}_j,res}(\h)$, $z_+\in U_{\mathcal{B}}^{\widehat{s}_j,res}(\z_+^j), x\in U_{\mathcal{B}}^{\widehat{s}_j,res}([-\delta_{j1},-\delta_{jo_j}])$.

From (\ref{comults}) and the fact that $\omega_0' S_{\widehat{s}_j}^{-1}$ is an anti-coautomorphism preserving weights and $\omega=\tau\omega_0'$ we obtain 
\begin{equation}\label{Dsy1}
\Delta_{\widehat{s}_j}^\tau (y)=y_0+\sum_i y_i\otimes h_i'+\sum_k{v_k}\otimes v_k',
\end{equation} 
where the weights of the elements $v_k'$ are strictly negative, $y_i\in \omega_0S_{\widehat{s}_j}^{-1}(U_{\mathcal{B}}^{\widehat{s}_j,res}([-\delta_{jn_j+1},-\delta_{jD}]))$ have non--zero weights, $h_i'\in U_{\mathcal{B}}^{\widehat{s}_j,res}(\h)$, and $y_0\in \mathcal{B}$ is the zero weight component of $y$.

By (\ref{DTss}) 
$$
\Delta_{\widehat{s}_j}(T_{\widehat{s}_j})=\theta_{\widehat{s}_j}^{\widehat{s}_j}T_{\widehat{s}_j}\otimes T_{\widehat{s}_j}
$$
Recalling the definition of $\theta_{\widehat{s}_j}^{\widehat{s}_j}$, and observing that $v_\mu$ and $v_\nu$ are highest weight vectors, and that the weights of the non--trivial weights components of $T_{w_j}^{-1}(\omega_0S_{\widehat{s}_j}^{-1}(U_{\mathcal{B}}^{\widehat{s}_j,res}([-\delta_{jn_j+1},-\delta_{jD}])))$ are non--positive as ${w_j}^{-1}[-\delta_{jn_j+1},-\delta_{jD}]\subset \Delta_-$, we deduce, using also (\ref{Dsy1}), that the left hand side of (\ref{dmus}) takes the form
$$
({\Delta_\mu^{\widehat{s}_{j+1}}}(\tau T_{w_j}^{-1}\tau(~\cdot~))\otimes {\Delta_\nu^{\widehat{s}_{j}}})(y\tau(hT_{\widehat{s}_j}z_+x))=y_0(v_\mu, \tau(T_{w_j}^{-1}h^1T_{\widehat{s}_j}z_+^1x^1 T_{w_j}T_{\widehat{s}_{j+1}}^{-1})v_\mu) (v_\nu, \tau(h^2T_{\widehat{s}_j}z_+^2x^2T_{\widehat{s}_j})v_\nu)=
$$
\begin{equation}\label{dmus1}
=y_0\kappa_j^\mu(v_\mu, \tau T_{w_j}^{-1}(h^1T_{\widehat{s}_j}(z_+^1x^1))v_\mu) (v_\nu, \tau h^2T_{\widehat{s}_j}(z_+^2x^2)v_\nu),
\end{equation}
where $\kappa_j^\mu\in \mathcal{B}^*$ is defined by the condition $\tau(T_{w_j}T_{\widehat{s}_{j+1}}^{-1})v_\mu=\kappa_j^\mu \tau(T_{\widehat{s}_j}^{-1}T_{w_j})v_\mu$ using Lemma \ref{bga}, $\Delta_{\widehat{s}_j}h=h^1\otimes h^2$, $\Delta_{\widehat{s}_j}z_+=z_+^1\otimes z_+^2$, and $\Delta_{\widehat{s}_j}x=x^1\otimes x^2$.

We conclude that only the zero weight components of $z_+^1x^1$ and of $z_+^2x^2$ can contribute to the right hand side of the last formula. Since $z_+x\in U_{\mathcal{B}}^{\widehat{s}_j,res}(\z_+^j)U_{\mathcal{B}}^{\widehat{s}_j,res}([-\delta_{j1},-\delta_{jo_j}])$, the weights of the weight components of $z_+x$ belong to $$\mathbb{N}((w_1\ldots w_{j-1})^{-1}(\Delta_0)_+\cup [-\delta_{j1},-\delta_{jo_j}]=\mathbb{N}(w_1\ldots w_{j-1})^{-1}((\Delta_0)_+\cup [-\beta_{j1},-\beta_c]).$$ Observe that $(\Delta_0)_+\cup [-\beta_{j1},-\beta_c]\subset [\beta^0_1,-\beta_c]$ which is a minimal segment. Therefore, since the comultiplication preserves weights, only the product of the zero weight components $(z_+)_0$ and $x_0$ of $z_+$ and $x$, respectively, can contribute to the right hand side of (\ref{dmus1}). By the definition of the algebras $U_{\mathcal{B}}^{\widehat{s}_j,res}(\z_+^j)$ and $U_{\mathcal{B}}^{\widehat{s}_j,res}([-\delta_{j1},-\delta_{jo_j}])$ we have $(z_+)_0,x_0\in \mathcal{B}$. Therefore we can rewrite (\ref{dmus1}) as follows
\begin{equation}\label{dmus2}
({\Delta_\mu^{\widehat{s}_{j+1}}}(\tau T_{w_j}^{-1}(~\cdot~))\otimes {\Delta_\nu^{\widehat{s}_{j}}})(y\tau(hT_{\widehat{s}_j}z_+x))=\kappa_j^\mu y_0\tau((z_+)_0)\tau(x_0)(v_\mu, \tau(T_{w_j}^{-1}(h^1))v_\mu) (v_\nu, \tau(h^2)v_\nu)=
\end{equation}
$$
=\kappa_j^\mu y_0\tau((z_+)_0)\tau(x_0)(\tau(T_{w_j})(v_\mu)\otimes v_\nu, \tau(h) (\tau(T_{w_j})(v_\mu)\otimes v_\nu))=\kappa_j^\mu y_0\tau((z_+)_0)\tau(x_0)(v_{w_j\mu+\nu},\tau(h)v_{w_j\mu+\nu})=
$$
$$
=\kappa_j^\mu \tau(y_0)(v_{w_j\mu+\nu},\tau(h)\tau(T_{\widehat{s}_j})(\tau(z_+x))v_{w_j\mu+\nu})=\kappa_j^\mu(v_{w_j\mu+\nu},y\tau(hT_{\widehat{s}_j}z_+xT_{\widehat{s}_j}^{-1})v_{w_j\mu+\nu})=\kappa_j^\mu\Delta_{w_j\mu+\nu}^{\widehat{s}_{j}}(y\tau(hT_{\widehat{s}_j}z_+x)).
$$
This establishes (\ref{dmus}) and completes the proof.

\end{proof}

Note that by (\ref{fo0}) $\omega_0'$ gives rise to algebra antiautomorphisms of the algebras $U_\mathcal{B}^{\widehat{s}_j,res}([-\delta_{jk},-\delta_{jo_j}])$. Therefore by (\ref{Ads0}) and by Lemma \ref{mainl3} the ${\rm Ad}^{\widehat{s}_j}$--action of $U_{\mathcal{B}}^{\widehat{s}_j,res}([-\delta_{jk},-\delta_{jo_j}])$ on $\mathbb{C}_{\mathcal{B}}^{\widehat{s}_j}[G]$ induces a locally finite action on $\mathbb{C}_{jk}^{loc}[G]=\mathbb{C}_{\mathcal{B}}^{\widehat{s}_j,loc}[G]/{{I}^{jk}_{\mathcal{B}}}^{loc}$.
Now we can state the following proposition.
\begin{proposition}\label{phijloc}
(i) If $\mu, \nu \in P_+$ are such that $w_j\mu+\nu\in P_+$ then for $j=1,\ldots, R-2$ the following relation holds in $\mathbb{C}_{jn_j+1}^{loc}[G]$:
$$
T_{w_j}^{-1}({\Delta_\mu^{\widehat{s}_{j+1}}})=\kappa_j^\mu {\Delta_{w_j\mu+\nu}^{\widehat{s}_{j}}}\otimes ({\Delta_\nu^{\widehat{s}_{j}}})^{-1}.
$$

(ii) For $j=1,\ldots, R-2$ $\varpi_j$ induces a morphism of $\mathcal{B}$--modules 
\begin{equation}\label{phijhom}
\varpi_j: \mathbb{C}_{j+1 1}^{loc}[G]\to \mathbb{C}_{jn_j+1}^{loc}[G]
\end{equation}
which satisfies 
\begin{equation}\label{phiAdloc}
\varpi_j{\rm Ad}^{\widehat{s}_{j+1}}(T_{w_j}^{-1}x)={\rm Ad}^{\widehat{s}_j}(x)\varpi_j, x\in U_\mathcal{B}^{\widehat{s}_j,res}([-\delta_{jn_j+1},-\delta_{jo_j}]),
\end{equation} 
and
$$
\varpi_j(f\otimes {\Delta_\mu^{\widehat{s}_{j+1}}})=\kappa_j^\mu (\varpi_j f)\otimes {\Delta_{w_j\mu+\nu}^{\widehat{s}_{j}}}\otimes (\Delta_\nu^{\widehat{s}_{j}})^{-1}, f\in \mathbb{C}_{j+1 1}^{loc}[G],
$$
$$
\varpi_j(f\otimes {\Delta_\mu^{\widehat{s}_{j+1}}}^{-1})=(\kappa_j^\mu)^{-1} (\varpi_j f)\otimes {\Delta_\nu^{\widehat{s}_{j}}}\otimes (\Delta_{w_j\mu+\nu}^{\widehat{s}_{j}})^{-1}, f\in \mathbb{C}_{j+1 1}^{loc}[G],
$$
where $\mu, \nu \in P_+$ are such that $w_j\mu+\nu\in P_+$.
\end{proposition}

\begin{proof}
(i) follows from Lemma \ref{phijlem} (iii) and from the definition of $\mathbb{C}_{jn_j+1}^{loc}[G]$ together with Lemma \ref{locid}.

(ii) By Lemma \ref{phijlem} (ii) for $f\in \mathbb{C}_{\mathcal{B}}^{\widehat{s}_{j+1}}[G]$ one has $\varpi_j(f\otimes {\Delta_\mu^{\widehat{s}_{j+1}}})=(\varpi_j f)\otimes T_{w_j}^{-1}({\Delta_\mu^{\widehat{s}_{j+1}}})$ mod ${{J}^{jn_j+1}_{\mathcal{B}}}$ assuming that the element $t_{\widehat{s}_j}=T_{\widehat{s}_j}$ is used in the definition of ${{J}^{jn_j+1}_{\mathcal{B}}}$.

Therefore by Lemma \ref{phijlem} (iii) $\varpi_j(f\otimes {\Delta_\mu^{\widehat{s}_{j+1}}})=\kappa_j (\varpi_j f)\otimes {\Delta_{w_j\mu+\nu}^{\widehat{s}_{j}}}\otimes (\Delta_\nu^{\widehat{s}_{j}})^{-1}$ mod ${{J}^{jn_j+1}_{\mathcal{B}}}^{loc}$ in $\mathbb{C}_{\mathcal{B}}^{\widehat{s}_j,loc}[G]$.

Thus $\varpi_j$ gives rise to a morphism of $\mathcal{B}$--modules $\varpi_j:\mathbb{C}_{\mathcal{B}}^{\widehat{s}_{j+1},loc}[G]\to \mathbb{C}_{\mathcal{B}}^{\widehat{s}_j,loc}[G]/ {{J}^{jn_j+1}_{\mathcal{B}}}^{loc}$ defined by 
\begin{equation}\label{pjl}
\varpi_j(f\otimes {\Delta_\mu^{\widehat{s}_{j+1}}}\otimes (\Delta_{\mu'}^{\widehat{s}_{j+1}})^{-1})=\kappa_j^\mu(\kappa_j^{\mu'})^{-1}(\varpi_j f)\otimes {\Delta_{w_j\mu+\nu}^{\widehat{s}_{j}}}\otimes (\Delta_\nu^{\widehat{s}_{j}})^{-1}\otimes (\Delta_{w_j\mu'+\nu'}^{\widehat{s}_{j}})^{-1}\otimes {\Delta_{\nu'}^{\widehat{s}_{j}}}=
\end{equation}
$$
=\kappa_j^\mu(\kappa_j^{\mu'})^{-1}(\varpi_j f)\otimes {\Delta_{w_j\mu+\nu+\nu'}^{\widehat{s}_{j}}}\otimes (\Delta_{\nu+w_j\mu'+\nu'}^{\widehat{s}_{j}})^{-1},
$$
where $\mu', \nu' \in P_+$ are such that $w_j\mu'+\nu'\in P_+$.

By Lemma \ref{phijlem} (i) and by (\ref{pjl}) one has $\varpi_j({{{J}^{j+11}_{\mathcal{B}}}'}^{loc})=0$ in $\mathbb{C}_{\mathcal{B}}^{\widehat{s}_j,loc}[G]/ {{J}^{jn_j+1}_{\mathcal{B}}}^{loc}$. Now (\ref{phijhom}) follows from the definitions of $\mathbb{C}_{j+1 1}^{loc}[G]$ and $\mathbb{C}_{jn_j+1}^{loc}[G]$ and from Lemma \ref{locid}, and (\ref{phiAdloc}) follows from formula (\ref{phiAd}).

The last two formulas in part (ii) are consequences of (\ref{pjl}).

\end{proof}


\section{Zhelobenko type operators for q-W--algebras}\label{Zheldef}

\pagestyle{myheadings}
\markboth{CHAPTER~\thechapter.~ZHELOBENKO TYPE OPERATORS FOR Q-W--ALGEBRAS}{\thesection.~ZHELOBENKO TYPE OPERATORS FOR Q-W--ALGEBRAS}

\setcounter{equation}{0}
\setcounter{theorem}{0}

This section is central in this chapter and in the whole book. We are going to introduce and study some quantum analogues of the operators $\Pi_{jk}$ and $\Pi_c$ defined in (\ref{Pip}) and (\ref{pic}). It turns out that analogues $P_{jk}$ of $\Pi_{jk}$ can be obtained by proper extrapolation of the expansion of the conjugation operator in (\ref{Pip}) in terms of the adjoint action operator and by replacing the coefficients $\varphi_{jk}$ with their quantum counterparts $B_{jk}$ introduced in Corollary \ref{corloc} (i).  However, the proof of Proposition \ref{BasP} which asserts that the image of their composition $P_c$ consists of invariant elements with respect to the adjoint action of $U_{\mathcal{B}}^{s,res}(\m_-)$ is rather complicated. It entirely relies on the properties of ${{I}^{jk}_{\mathcal{B}}}^{loc}$, of the quotients $\mathbb{C}_{jk}^{loc}[G]$ and of the adjoint action, which were obtained in the previous sections of this chapter.

The point is that in the quantum case we do not have in our disposal isomorphism (\ref{cross1}) which plays a crucial role in the proof of a similar property for the operator $\Pi_c$, given by (\ref{Pip}), (\ref{pic}), as one can see from the proof of Proposition \ref{tp}.  

We start with the definition of quantum analogues $P_{jp}$ of the operators $\Pi_{jp}$. For technical reasons we shall also need more general operators $P_{jp}^k$. More precisely, by Lemmas \ref{quotact} and \ref{mainl3} for $f\in \mathbb{C}_{jp}^{loc}[G]$, $j=1,\ldots ,R-1$, $p=1,\ldots, n_j$, $n,k\in \mathbb{N}$, $n\geq k$, we have a well--defined element ${\rm Ad}_{\widehat{s}_j}^0 f_{\delta_{jp}}^{(n-k)}(f)\otimes B_{jp}^n\in \mathbb{C}_{jp}^{loc}[G]$. Note that by Lemma \ref{mainl3} for each $f\in \mathbb{C}_{jp}^{loc}[G]$  one has ${\rm Ad}_{\widehat{s}_j}^0 f_{\delta_{jp}}^{(n-k)}(f)=0$ for $n-k>N(f)$, where $N(f)\in \mathbb{N}$ depends on $f$. Therefore  we can define define an element $P_{jp}^k(f)\in \mathbb{C}_{jp}^{loc}[G]$ by 
\begin{equation}\label{Ppk}
P_{jp}^k(f)=\sum_{n=k}^\infty (-1)^nq_{\delta_{jp}}^{-\frac{(n-1)(n-2k)}{2}}{\rm Ad}_{\widehat{s}_j}^0 f_{\delta_{jp}}^{(n-k)}(f)\otimes B_{jp}^n.	\index[not]{P@$P_{jp}^k$}
\end{equation} 

In particular, formula (\ref{Ppk}) defines an operator
$$
P_{jp}^k:\mathbb{C}_{jp}^{loc}[G]\rightarrow \mathbb{C}_{jp}^{loc}[G].
$$

Denote $P_{jp}^0=P_{jp}$. \index[not]{P@$P_{jp}$}

The following proposition summarizes the main properties of the operators $P_{jp}$. 
\begin{proposition}\label{BasP}
(i) For any $j=1,\ldots, R-1$, $p=1,\ldots, n_j$ the composition 
\begin{equation}\label{comp}
P_{\leq jp}:=P_{11}\ldots P_{1n_1}\circ \varpi_1 \circ P_{21}\ldots P_{2n_2}\circ \varpi_2 \ldots \circ \varpi_{j-1} \circ P_{j1}\ldots P_{jp}:\mathbb{C}_{jp}^{loc}[G]\rightarrow \mathbb{C}_{11}^{loc}[G] \index[not]{P@$P_{\leq jp}$}
\end{equation}
is well--defined. In fact $P_{\leq jp}$ is well--defined as an operator with the domain $\mathbb{C}_{jp+1}^{loc}[G]$, i.e. ${I_{\mathcal{B}}^{j p+1}}^{loc}/{I_{\mathcal{B}}^{jp}}^{loc}\subset \mathbb{C}_{jp}^{loc}[G]$ belongs to the kernel of this operator.

(ii) For any $\beta\in [\beta_{11}, \beta_{jp}]$, $n>0$ and for any $f\in \mathbb{C}_{jp+1}^{loc}[G]$  we have  
\begin{equation}\label{projp}
{\rm Ad}_s^0 f_{\beta}^{(n)} (P_{\leq jp}(f))=0.	
\end{equation}

In particular, for any $x\in U_{\mathcal{B}}^{s,res}([-\beta_{11},-\beta_{jp}])$, $f\in \mathbb{C}_{jp+1}^{loc}[G]$
$$
{\rm Ad}_s^0 x (P_{\leq jp}(f))=\varepsilon_s(x)P_{\leq jp}(f),
$$
and if we denote $P_c=P_{\leq R-1 n_{R-1}}$ \index[not]{P@$P_c$} then for any $x\in U_{\mathcal{B}}^{s,res}(\m_-)$, $f\in \mathbb{C}_{R-1 n_{R-1}+1}^{loc}[G]$
$$
{\rm Ad}_s^0 x (P_c(f))=\varepsilon_s(x)P_c(f).
$$

Moreover, for $j=1,\ldots, R-1$, $p=1,\ldots, n_j$, $(j,p)\neq (R-1,n_{R-1})$ one has
\begin{equation}\label{Pcomm}
{\rm Ad}_s^0 f_{\beta_{jp+1}}^{(n)} (P_{\leq jp}(f))=c_{jp+1}P_{\leq jp}({\rm Ad}_{\widehat{s}_j}^0 f_{\delta_{jp+1}}^{(n)}f),	
\end{equation}
where we assume that $\beta_{jn_j+1}=\beta_{j+11}$ and $c_{jp+1}\in \{\pm q^{\mathbb{Z}}\}$.
\end{proposition}

The proof of Proposition \ref{BasP} is quite long. It will be split into several lemmas. Firstly we shall show that $P_{\leq jp}$ is well--defined as an operator with the domain of definition $\mathbb{C}_{jp+1}^{loc}[G]$, i.e. ${I_{\mathcal{B}}^{j p+1}}^{loc}/{I_{\mathcal{B}}^{jp}}^{loc}\subset \mathbb{C}_{jp}^{loc}[G]$ belongs to the kernel of this operator. 

Then we prove that for $n>0$ ${\rm Ad}_{\widehat{s}_j}^0 f_{\delta_{jp}}^{(n)}P_{jp}$ is the zero operator. This is done using an explicit calculation and observing that ${\rm Ad}_{\widehat{s}_j}^0 f_{\delta_{jp}}^{(n)}$ has some properties analogues to properties of derivations of order $n$. 

Finally to show that for $\beta\in [\beta_{11}, \beta_{jp}]$ and $n>0$ ${\rm Ad}_s^0 f_{\beta}^{(n)} (P_{\leq jp}(f))=0$ we shall use the property of the images of operators $P_{jp}$ mentioned in the previous paragraph and commutation relations between $P_{jp}$ and ${\rm Ad}_{\widehat{s}_j}^0 f_{\beta}^{(n)}$ for $\beta>\delta_{jp}$ which, in particular, lead to commutation relations (\ref{Pcomm}).

We start by showing using an explicit calculation that for $n>0$   
${\rm Ad}_{\widehat{s}_j}^0 f_{\delta_{jp}}^{(n)}P_{jp}$ is the zero operator.
Since the right hand side of formula (\ref{Ppk}) contains products of elements from $\mathbb{C}_{jp}^{loc}[G]$ and of $B_{jp}^{n}$, we have to study the adjoint action on such products. 
\begin{lemma}\label{mainl31}
(i) Let $B_{jp}={A_{jp}^0}^{-1}\otimes A_{jp}\in \mathbb{C}_{\mathcal{B}}^{\widehat{s}_j,loc}[G]_0$, $f\in \mathbb{C}_{jp}^{loc}[G]$, $j=1,\ldots, R-1$, $p=1,\ldots ,n_j$. Then for any $k,n\in \mathbb{N}$ the identity 
\begin{equation}\label{mainf}
{\rm Ad}_{\widehat{s}_j}^0 f_{\delta_{jp}}^{(n)} (f\otimes B_{jp}^k)=\sum_{r=0}^{{\rm min}(n,k)}q_{\delta_{jp}}^{-(2n-r)k-\frac{r(r-1)}2+rn}\left[ \begin{array}{c} k \\ r \end{array} \right]_{q_{\delta_{jp}}}{\rm Ad}_{\widehat{s}_j}^0 f_{\delta_{jp}}^{(n-r)}(f)\otimes B_{jp}^{k-r}
\end{equation}
holds in $\mathbb{C}_{jp}^{loc}[G]$.

(ii) For any $\alpha, \beta \in Q$ we denote
$$
c_{\alpha\beta}^s=q^{\left\langle (2K_s+id)\alpha, \beta\right\rangle}. \index[not]{c@$c_{\alpha\beta}^s$}
$$ 
Then for any $j=1,\ldots, R-1$, $p=1,\ldots ,n_j$, $m,n\in \mathbb{N}$, $f\in \mathbb{C}_{jp}^{loc}[G]$, $\beta \in [\delta_{jp+1},\delta_{jo_j}]$ the identity
\begin{equation}\label{mainfpq}
{\rm Ad}_{\widehat{s}_j}^0 f_{\beta}^{(m)} (f\otimes B_{jp}^n)=(c_{\delta_{jp}\beta}^{\widehat{s}_j})^{-mn}{\rm Ad}_{\widehat{s}_j}^0 f_{\beta}^{(m)}(f)\otimes B_{jp}^{n}	
\end{equation}
holds in $\mathbb{C}_{jp}^{loc}[G]$.
\end{lemma}

\begin{proof}
(i) We prove (\ref{mainf}) by induction over $k$. We start with the case when $f$ is the image of an element of $\mathbb{C}_{\mathcal{B}}^{\widehat{s}_j}[G]$ in $\mathbb{C}_{jp}^{loc}[G]$ under the canonical map $\mathbb{C}_{\mathcal{B}}^{\widehat{s}_j}[G]\to \mathbb{C}_{\mathcal{B}}^{\widehat{s}_j,loc}[G]\to \mathbb{C}_{jp}^{loc}[G]$. Firstly by (\ref{ad0n}), we have for any $n\in \mathbb{N}$
$$
{\rm Ad}_{\widehat{s}_j}^0f_{\delta_{jp}}^{(n)} (f\otimes A_{jp})=q^{-n\left\langle \kappa{1+s \over 1-s}P_{{\h'}}\delta_{jp}-\delta_{jp},\mu_{jp}-\delta_{jp}\right\rangle}{\rm Ad}_{\widehat{s}_j}^0 (f_{\delta_{jp}}^{(n)}) (f)\otimes A_{jp}+
$$
$$
+q^{-(n-1)\left\langle \kappa{1+s \over 1-s}P_{{\h'}}\delta_{jp}-\delta_{jp},\mu_{jp}-\delta_{jp}\right\rangle}q_{{\delta_{jp}}}^{n-1}{\rm Ad}_{\widehat{s}_j}^0 (f_{\delta_{jp}}^{(n-1)}) (f)\otimes A_{jp}(\omega_0S_{\widehat{s}_j}^{-1}(f_{\delta_{jp}}) ~\cdot~ )+
$$
$$
+\sum_{k=1}^{n}q_{{\delta_{jp}}}^{k(n-k)}q^{n\left\langle \kappa{1+s \over 1-s}P_{{\h'}}\delta_{jp}-\delta_{jp},\mu_{jp}-\delta_{jp}\right\rangle}{\rm Ad}_{\widehat{s}_j}^0 (G_{\delta_{jp}}^{-k}f_{\delta_{jp}}^{(n-k)}) (f)\otimes A_{jp}( ~\cdot~ \omega_0(f_{\delta_{jp}}^{(k)}))+
$$
\begin{equation}\label{Adfp}
+\sum_{k=1}^{n}q_{{\delta_{jp}}}^{k(n-k)+(n-k-1)}{\rm Ad}_{\widehat{s}_j}^0 (G_{\delta_{jp}}^{-k}f_{\delta_{jp}}^{(n-k-1)}) (f)\otimes A_{jp}(\omega_0S_{\widehat{s}_j}^{-1}(G_{\delta_{jp}}^{-n+1}f_{\delta_{jp}}) ~\cdot~ \omega_0(f_{\delta_{jp}}^{(k)}))+	
\end{equation}
$$
+\sum_{k=0}^{n}\sum_{p=0}^{n-k-2}q_{{\delta_{jp}}}^{k(n-k)+p(n-k-p)}{\rm Ad}_{\widehat{s}_j}^0 (G_{\delta_{jp}}^{-k}f_{\delta_{jp}}^{(p)}) (f)\otimes A_{jp}(\omega_0S_{\widehat{s}_j}^{-1}(G_{\delta_{jp}}^{-k-p}f_{\delta_{jp}}^{(n-k-p)}) ~\cdot~ \omega_0(f_{\delta_{jp}}^{(k)}))+	
$$
$$
+\sum_{k=0}^{n-1}\sum_i q_{\delta_{jp}}^{k(n-k)}{\rm Ad}_{\widehat{s}_j}^0 (G_{\delta_{jp}}^{-k} x_i^{(n-k)}) (f)\otimes A_{jp}((\omega_0 S_{\widehat{s}_j}^{-1})(G_{\delta_{jp}}^{-k}y_i^{(n-k)})~\cdot~ \omega_0(f_{\delta_{jp}}^{(k)}))+
$$
$$
+\sum_i {\rm Ad}_{\widehat{s}_j}^0({y_i^{(n)}}^2) (f)\otimes A_{jp}((\omega_0 S_{\widehat{s}_j}^{-1})({y_i^{(n)}}^1)~\cdot~ \omega_0 (x_i^{(n)})).
$$

By Lemma \ref{barJp} (i) (see, in particular, (\ref{A1})), by the definition of ${{J}^{jk}_{\mathcal{B}}}^{loc}$ in the beginning of Section \ref{localG} and by Lemma \ref{locid} all terms in the right hand side of (\ref{Adfp}), except for the first two, vanish in $\mathbb{C}_{jp}^{loc}[G]$. Also by (\ref{A2}) we have $A_{jp}((\omega_0 S_{\widehat{s}_j}^{-1})(f_{\delta_{jp}})~\cdot~ )=A_{jp}^0(~\cdot~)$, and hence (\ref{Adfp}) takes the form
$$
{\rm Ad}_{\widehat{s}_j}^0f_{\delta_{jp}}^{(n)} (f\otimes A_{jp})=q^{-n\left\langle \kappa{1+s \over 1-s}P_{{\h'}}\delta_{jp}-\delta_{jp},\mu_{jp}-\delta_{jp}\right\rangle}{\rm Ad}_{\widehat{s}_j}^0 (f_{\delta_{jp}}^{(n)}) (f)\otimes A_{jp}+
$$
$$
+q^{-(n-1)\left\langle \kappa{1+s \over 1-s}P_{{\h'}}\delta_{jp}-\delta_{jp},\mu_{jp}-\delta_{jp}\right\rangle}q_{{\delta_{jp}}}^{n-1}{\rm Ad}_{\widehat{s}_j}^0 (f_{\delta_{jp}}^{(n-1)}) (f)\otimes A_{jp}^0.
$$

Since by definition elements of $\mathfrak{S}_{\widehat{s}_j}^{-1}$ naturally act on $\mathbb{C}_{jp}^{loc}[G]$ by multiplication from the right one can multiply the last identity by $(A_{jp}^0)^{-1}$ from the right. Using commutation relations (\ref{comm1}), formula (\ref{Apo}) and recalling the definition of $B_{jp}$ we obtain 
\begin{equation}\label{basead}
{\rm Ad}_{\widehat{s}_j}^0 f_{\delta_{jp}}^{(n)} (f\otimes B_{jp})=q_{\delta_{jp}}^{-(n-1)}{\rm Ad}_{\widehat{s}_j}^0 f_{\delta_{jp}}^{(n-1)}(f)+q_{\delta_{jp}}^{-2n}{\rm Ad}_{\widehat{s}_j}^0 f_{\delta_{jp}}^{(n)}(f)\otimes B_{jp}.	
\end{equation}

Multiplying the last identity by $(\Delta_\mu^{\widehat{s}_j})^{-1}$ from the right and using the same arguments we get
$$
{\rm Ad}_{\widehat{s}_j}^0 f_{\delta_{jp}}^{(n)} ({f\otimes (\Delta_\mu^{\widehat{s}_j})^{-1}}\otimes B_{jp})=q_{\delta_{jp}}^{-(n-1)}{\rm Ad}_{\widehat{s}_j}^0 f_{\delta_{jp}}^{(n-1)}({f\otimes (\Delta_\mu^{\widehat{s}_j})^{-1}})+q_{\delta_{jp}}^{-2n}{\rm Ad}_{\widehat{s}_j}^0 f_{\delta_{jp}}^{(n)}({f\otimes (\Delta_\mu^{\widehat{s}_j})^{-1}})\otimes B_{jp},
$$
Since elements ${f\otimes (\Delta_\mu^{\widehat{s}_j})^{-1}}$ span $\mathbb{C}_{jp}^{loc}[G]$ formula (\ref{basead}) holds for any $n\in \mathbb{N}$, $f\in \mathbb{C}_{jp}^{loc}[G]$, i.e. (\ref{mainf}) holds for any $n\in \mathbb{N}$, $k=1$. 
This establishes the base of induction.

Now we assume that (\ref{mainf}) holds for some natural $k$ and for all natural $n$ and prove that it holds for $k+1$ and all natural $n$. 
Since by Lemma \ref{quotact} right multiplication gives rise to a well--defined action of $B_{jp}$ on $\mathbb{C}_{jp}^{loc}[G]$, we have for any $f\in \mathbb{C}_{jp}^{loc}[G]$ using the base of induction and the induction assumption
$$
{\rm Ad}_{\widehat{s}_j}^0 f_{\delta_{jp}}^{(n)} (f\otimes B_{jp}^{k+1})=q_{\delta_{jp}}^{-(n-1)}{\rm Ad}_{\widehat{s}_j}^0 f_{\delta_{jp}}^{(n-1)} (f\otimes B_{jp}^{k})+q_{\delta_{jp}}^{-2n}{\rm Ad}_{\widehat{s}_j}^0 f_{\delta_{jp}}^{(n)}(f\otimes B_{jp}^{k})\otimes B_{jp}=
$$
$$
=q_{\delta_{jp}}^{-(n-1)}\sum_{r=0}^{{\rm min}(n-1,k)}q_{\delta_{jp}}^{-(2(n-1)-r)k-\frac{r(r-1)}2+r(n-1)}\left[ \begin{array}{c} k \\ r \end{array} \right]_{q_{\delta_{jp}}}{\rm Ad}_{\widehat{s}_j}^0 f_{\delta_{jp}}^{(n-r-1)}(f)\otimes B_{jp}^{k-r}+
$$
$$
+q_{\delta_{jp}}^{-2n}\sum_{r=0}^{{\rm min}(n,k)}q_{\delta_{jp}}^{-(2n-r)k-\frac{r(r-1)}2+rn}\left[ \begin{array}{c} k \\ r \end{array} \right]_{q_{\delta_{jp}}}{\rm Ad}_{\widehat{s}_j}^0 f_{\delta_{jp}}^{(n-r)}(f)\otimes B_{jp}^{k-r+1}=
$$
$$
=\sum_{r=0}^{{\rm min}(n,k+1)}q_{\delta_{jp}}^{-(2n-r)(k+1)-\frac{r(r-1)}2+rn}\left(q_{\delta_{jp}}^{-r}\left[ \begin{array}{c} k \\ r \end{array} \right]_{q_{\delta_{jp}}}+q_{\delta_{jp}}^{-r+(k+1)}\left[ \begin{array}{c} k \\ r-1 \end{array} \right]_{q_{\delta_{jp}}}\right){\rm Ad}_{\widehat{s}_j}^0 f_{\delta_{jp}}^{(n-r)}(f)\otimes B_{jp}^{k-r+1}=
$$
$$
=\sum_{r=0}^{{\rm min}(n,k+1)}q_{\delta_{jp}}^{-(2n-r)(k+1)-\frac{r(r-1)}2+rn}\left[ \begin{array}{c} k+1 \\ r \end{array} \right]_{q_{\delta_{jp}}}{\rm Ad}_{\widehat{s}_j}^0 f_{\delta_{jp}}^{(n-r)}(f)\otimes B_{jp}^{k+1-r},
$$
where we used the identity
$$
q_{\delta_{jp}}^{-r}\left[ \begin{array}{c} k \\ r \end{array} \right]_{q_{\delta_{jp}}}+q_{\delta_{jp}}^{-r+(k+1)}\left[ \begin{array}{c} k \\ r-1 \end{array} \right]_{q_{\delta_{jp}}}=\left[ \begin{array}{c} k+1 \\ r \end{array} \right]_{q_{\delta_{jp}}}
$$
which can be found e.g. in \cite{KQ}, Proposition 6.1.
This establishes the induction step and completes the proof of (\ref{mainf}).

Formula (\ref{mainfpq}) is proved in a similar way by induction using Lemma \ref{barJp} (i), formulas (\ref{ad0n}), (\ref{A2}), (\ref{Apo}), and the definition of $B_{jp}$.

\end{proof}

The previous lemma shows that the operator ${\rm Ad}_{\widehat{s}_j}^0 f_{\delta_{jp}}^{(n)}$ acts on $f\otimes B_{jp}^{k}$ as a quantum analogue of a derivation of order $n$ such that the derivative of $B_{jp}$ by ${\rm Ad}_{\widehat{s}_j}^0 f_{\delta_{jp}}$ is equal to one. Recalling the definition of the classical counterpart $\varphi_{jp}$ of $B_{jp}$ one can observe that a similar formula can be obtained in the classical case as well. It will be given by specializing (\ref{mainf}) at $q^{\frac{1}{d{\bar{r}}^2}}=1$.  

In the next lemma we show that for $n>0$ the composition ${\rm Ad}_{\widehat{s}_j}^0 f_{\delta_{jp}}^{(n)}P_{jp}$ is the zero operator. At the same time this lemma lays a basis for the proof of the fact that the kernel of $P_{\leq jp}$ contains ${{I}_{\mathcal{B}}^{jp+1}}^{loc}/{{I}_{\mathcal{B}}^{jp}}^{loc}$.  
\begin{lemma}
For any $f\in \mathbb{C}_{jp}^{loc}[G]$, $j=1,\ldots, R-1$, $p=1,\ldots ,n_j$, $k\in \mathbb{N}$ we have
\begin{equation}\label{mainprop1}
P_{jp}^k(f\otimes B_{jp})=0,	
\end{equation}

\begin{equation}\label{mainprop2}
{\rm Ad}_{\widehat{s}_j}^0 f_{\delta_{jp}}^{(n)}(P_{jp}(f))=0,~n>0.	
\end{equation}
\end{lemma}

\begin{proof}
Since by Lemma \ref{quotact} right multiplication by $B_{jp}$ gives rise to an action of $B_{jp}$ on $\mathbb{C}_{jp}^{loc}[G]$ we have by (\ref{Ppk}) and  (\ref{mainf}) with $k=1$ 
$$
P_{jp}^k(f\otimes B_{jp})=\sum_{n=k+1}^\infty (-1)^nq_{\delta_{jp}}^{-\frac{(n-1)(n-2k)}{2}}q_{\delta_{jp}}^{-(n-k-1)}{\rm Ad}_{\widehat{s}_j}^0 f_{\delta_{jp}}^{(n-k-1)}(f)\otimes B_{jp}^n+
$$
$$
+\sum_{n=k}^\infty (-1)^nq_{\delta_{jp}}^{-\frac{(n-1)(n-2k)}{2}} q_{\delta_{jp}}^{-2(n-k)}{\rm Ad}_{\widehat{s}_j}^0 f_{\delta_{jp}}^{(n-k)}(f)\otimes B_{jp}^{n+1}=
$$
$$
=-\sum_{n=k}^\infty {(-1)^nq_{\delta_{jp}}^{-\frac{n(n+1-2k)}{2}-n+k}}{\rm Ad}_{\widehat{s}_j}^0 f_{\delta_{jp}}^{(n-k)}(f)\otimes B_{jp}^{n+1}+
$$
$$
+\sum_{n=k}^\infty {(-1)^nq_{\delta_{jp}}^{-\frac{(n-1)(n-2k)}{2}-2n+2k}}{\rm Ad}_{\widehat{s}_j}^0 f_{\delta_{jp}}^{(n-k)}(f)\otimes B_{jp}^{n+1}=0
$$
which proves (\ref{mainprop1}).

Similarly by (\ref{Ppk}) and (\ref{mainf}) we obtain 
$$
{\rm Ad}_{\widehat{s}_j}^0 f_{\delta_{jp}}^{(n)}(P_{jp}(f))={\rm Ad}_{\widehat{s}_j}^0 f_{\delta_{jp}}^{(n)}\left(\sum_{k=0}^\infty (-1)^kq_{\delta_{jp}}^{-\frac{(k-1)k}{2}}
{\rm Ad}_{\widehat{s}_j}^0 f_{\delta_{jp}}^{(k)}(f)\otimes B_{jp}^{k}\right)=
$$
$$
=\sum_{k=0}^\infty\sum_{t=0}^{{\rm min}(n,k)}(-1)^kq_{\delta_{jp}}^{-\frac{(k-1)k}{2}-(2n-t)k-\frac{t(t-1)}{2}+tn}\left[ \begin{array}{c} k \\ t \end{array} \right]_{q_{\delta_{jp}}}{\rm Ad}_{\widehat{s}_j}^0 (f_{\delta_{jp}}^{(k)}f_{\delta_{jp}}^{(n-t)})(f)\otimes B_{jp}^{k-t}.
$$

Introducing a new variable of summation $k-t=r$ and using the identity $$f_{\delta_{jp}}^{(k)}f_{\delta_{jp}}^{(n-t)}=\left[ \begin{array}{c} n+k-t \\ k \end{array} \right]_{q_{\delta_{jp}}}f_{\delta_{jp}}^{(n+k-t)}$$  we get
$$
{\rm Ad}_{\widehat{s}_j}^0 f_{\delta_{jp}}^{(n)}(P_{jp}(f))=
$$
$$
=\sum_{r=0}^\infty\sum_{t=0}^{n}(-1)^{r+t}q_{\delta_{jp}}^{-\frac{(r+t-1)(r+t)}{2}-(2n-t)(r+t)-\frac{t(t-1)}{2}+tn}\left[ \begin{array}{c} r+t \\ t \end{array} \right]_{q_{\delta_{jp}}}\left[ \begin{array}{c} r+n \\ r+t \end{array} \right]_{q_{\delta_{jp}}}
{\rm Ad}_{\widehat{s}_j}^0 (f_{\delta_{jp}}^{(n+r)})(f)\otimes B_{jp}^{r}. 
$$

Now recalling that $$\left[ \begin{array}{c} r+t \\ t \end{array} \right]_{q_{\delta_{jp}}}\left[ \begin{array}{c} r+n \\ r+t \end{array} \right]_{q_{\delta_{jp}}}=\left[ \begin{array}{c} n \\ t \end{array} \right]_{q_{\delta_{jp}}}\left[ \begin{array}{c} r+n \\ r \end{array} \right]_{q_{\delta_{jp}}}$$ we obtain
$$
{\rm Ad}_{\widehat{s}_j}^0 f_{\delta_{jp}}^{(n)}(P_{jp}(f))=
$$
$$
=\sum_{r=0}^\infty(-1)^{r}q_{\delta_{jp}}^{-\frac{(r-1)r}{2}-2nr}\left[ \begin{array}{c} r+n \\ r \end{array} \right]_{q_{\delta_{jp}}}\left(\sum_{t=0}^{n}(-1)^{t}q_{\delta_{jp}}^{t-tn}\left[ \begin{array}{c} n \\ t \end{array} \right]_{q_{\delta_{jp}}}\right){\rm Ad}_{\widehat{s}_j}^0 (f_{\delta_{jp}}^{(n+r)})(f)\otimes B_{jp}^{r}=0,
$$
where we used the identity 
$$
\sum_{t=0}^{n}(-1)^{t}q_{\delta_{jp}}^{t-tn}\left[ \begin{array}{c} n \\ t \end{array} \right]_{q_{\delta_{jp}}}=0
$$
which follows from the q--binomial theorem (see e.g. \cite{GR}, Ch. 1). \index{q--binomial theorem}

This proves (\ref{mainprop2}).

\end{proof}

\begin{lemma}\label{Jvan}
(i) For any $k\in \mathbb{N}$, $j=1,\ldots, R-1$, $p-1,\ldots, n_j$ one has $P_{jp}^k({I_{\mathcal{B}}^{jp+1}}^{loc}/{I_{\mathcal{B}}^{jp}}^{loc})=0$. Thus $P_{jp}^k$ is well--defined as an operator with the domain $\mathbb{C}_{jp+1}^{loc}[G]$.

(ii) Moreover, for any $1\leq i\leq j\leq R-1$ $1\leq q\leq n_i$, $1\leq p\leq n_j$ such that $q\leq p$ if $i=j$, and any $k_{ij}\in \mathbb{N}$ the composition 
$$
P_{iq}^{k_{iq}}\ldots P_{in_i}^{k_{in_i}}\circ \varpi_i\circ P_{i+1 1}^{k_{i+1 1}}\ldots P_{i+1 n_{i+1}}^{k_{i+1 n_{i+1}}}\circ \varpi_{i+1}  \ldots  \varpi_{j-1} \circ P_{j1}^{k_{j1}}\ldots P_{jp}^{k_{jp}}
$$ 
is well--defined and gives rise to an operator with the domain $\mathbb{C}_{jp+1}^{loc}[G]$ , and the target space being $\mathbb{C}_{iq}^{loc}[G]$.

In particular, the composition 
$$
P_{\leq jp}=P_{11}^{k_{11}}\ldots P_{1n_1}^{k_{1n_1}}\circ \varpi_1\circ P_{2 1}^{k_{2 1}}\ldots P_{2 n_{2}}^{k_{2 n_2}}\circ \varpi_2  \ldots  \varpi_{j-1} \circ P_{j1}^{k_{j1}}\ldots P_{jp}^{k_{jp}}
$$ 
is well--defined and gives rise to an operator  with the domain $\mathbb{C}_{jp+1}^{loc}[G]$, and the target space being $\mathbb{C}_{11}^{loc}[G]$.
\end{lemma}

\begin{proof}
Let us show that $P_{jp}^k({{I}_{\mathcal{B}}^{jp+1}}^{loc}/{{I}_{\mathcal{B}}^{jp}}^{loc})=0$. 
Indeed, by Lemma \ref{locid} ${{I}^{jp+1}_{\mathcal{B}}}^{loc}={{{J}^{jp+1}_q}'}^{loc}\cap \mathbb{C}_{\mathcal{B}}^{{\widehat{s}_j},loc}[G]$ and ${{I}^{jp}_{\mathcal{B}}}^{loc}={{{J}^{jp}_q}'}^{loc}\cap \mathbb{C}_{\mathcal{B}}^{{\widehat{s}_j},loc}[G]$. Therefore from the definitions of ${{{J}^{jp+1}_q}'}^{loc}$ and ${{{J}^{jp}_q}'}^{loc}$, and from commutation relations (\ref{comm1}), it follows that for any element $f\in {{I}_{\mathcal{B}}^{jp+1}}^{loc}$ there exist elements $u\in \mathcal{B}$, $u\neq 0$, $g\in {{I}_{\mathcal{B}}^{jp}}^{loc}$ and $h\in \mathbb{C}_{\mathcal{B}}^{{\widehat{s}_j},loc}[G]$ such that $uf=g+h\otimes B_{jp}$. Let $\bar{f},\bar{g}$, $\overline{h\otimes B_{jp}}$ and $\bar{h}$ be the corresponding classes of these elements in $\mathbb{C}_{jp}^{loc}[G]$. Note that by Corollary \ref{corloc} (i) and by Lemma \ref{quotact} one has $\overline{h\otimes B_{jp}}=\bar{h}\otimes B_{jp}$ and $\bar{g}=0$, and by the previous Lemma $P_{jp}^{k_{jp}}(\overline{h\otimes B_{jp}})=P_{jp}^{k_{jp}}(\bar{h}\otimes B_{jp})=0$. Thus $P_{jp}^{k_{jp}}(\overline{uf})=uP_{jp}^{k_{jp}}(\bar{f})=0$. Since $u\neq 0$ this implies $P_{jp}^{k_{jp}}(\bar{f})=0$.
This completes the proof of the first statement.

The remaining statements of the lemma are simple corollaries of the first assertion and of (\ref{phijhom}).

\end{proof}

Next, we are going to study how the adjoint action of quantum root vectors commutes with the operators $P_{jp}$. For this purpose we shall need some commutation relations between quantum root vectors stated in the following lemma.
\begin{lemma}
Let $f_\beta\in U_{\mathcal{B}}^{s,res}(\g)$ be quantum root vectors defined with the help of an arbitrary fixed normal ordering on $\Delta_+$. Then for any $\alpha<\beta$, $\alpha,\beta\in \Delta_+$ and any $m,n\in \mathbb{N}$ we have
\begin{equation}\label{fmnm}
f_\alpha^{(m)}f_\beta^{(n)}=(c_{\alpha \beta}^s)^{mn}\sum_{p=0}^{m}q_\alpha^{p(m-1)}\sum_{p_1,\ldots, p_k \in \mathbb{N}}d_{p_1,\ldots, p_k}^{np}f_{\zeta_1}^{(p_1)}\ldots f_{\zeta_k}^{(p_k)}f_\alpha^{(m-p)},	
\end{equation}
where $[\alpha,\beta]=\{\alpha,\zeta_1,\ldots,\zeta_k\}$ as a set, $\alpha<\zeta_1<\ldots<\zeta_k=\beta$ in the fixed normal ordering of $\Delta_+$, and the coefficients $d_{p_1,\ldots, p_k}^{np}\in \mathcal{B}$ do not depend on $m$ and satisfy the following properties
\begin{equation}\label{dnp1}
d_{p_1,\ldots, p_{k-1},n}^{np}=\left\{ \begin{array}{l} 1~~{\rm if}~~p=p_1=\ldots = p_{k-1}=0 \\ 0 ~~{\rm otherwise} \end{array}
\right.,
\end{equation}
and for $0\leq r<n$ and any $p$
\begin{equation}\label{dnp2}
d_{0,\ldots, 0,r}^{np}=0.
\end{equation}
\end{lemma}

\begin{proof}
To prove this lemma it suffices to show that
\begin{equation}\label{fmn}
f_\alpha^mf_\beta^n=(c_{\alpha \beta}^s)^{mn}\sum_{p=0}^{m}S_m^p\sum_{p_1,\ldots, p_k \in \mathbb{N}}c_{p_1,\ldots, p_k}^{np}f_{\zeta_1}^{p_1}\ldots f_{\zeta_k}^{p_k}f_\alpha^{m-p},	
\end{equation}
where
$$
S_m^p=q_\alpha^{p(m-1)}\left[\begin{array}{c} m \\ p \end{array}\right]_{q_\alpha}, \index[not]{S@$S_m^p$}
$$
and the coefficients $c_{p_1,\ldots, p_k}^{np}\in \mathcal{A}$ do not depend on $m$.

Indeed, dividing (\ref{fmn}) by $[n]_{q_\alpha}![m]_{q_\beta}!$ we arrive at an identity of the form (\ref{fmnm}) where the coefficients $d_{p_1,\ldots, p_k}^{np}$ a priori belong to $\mathbb{C}(q^{\frac{1}{d{\bar{r}}^2}})$. 
But by the uniqueness of the Poincar\'{e}-Birkhoff-Witt decomposition in $U_{\mathcal{B}}^{s,res}(\n_-)$ (see Lemma \ref{segmPBWs} (iv) and Remark \ref{segmPBWsrev}) we have $d_{p_1,\ldots, p_k}^{np}\in \mathcal{B}$.

Now we establish (\ref{fmn}). By commutation relations (\ref{cmrelf}) we have for any $\gamma>\alpha$
\begin{equation}\label{fm1}
{f}_{\alpha}{f}_{\gamma} = c_{\alpha \gamma}^s{f}_{\gamma}{f}_{\alpha}+ D_{\alpha}f_{\gamma},	
\end{equation}
where 
$$
D_{\alpha}f_{\gamma}=\sum_{p_1,\ldots, p_j \in \mathbb{N}}C(p_1,\ldots,p_j)
{f}_{\zeta'_1}^{p_1}{f}_{\zeta'_2}^{p_2}\ldots {f}_{\zeta'_j}^{p_j}, \index[not]{D@$D_\alpha$}
$$
where $\alpha<\zeta'_1<\ldots<\zeta'_j<\gamma$, $[\alpha,\gamma]=\{\alpha,\zeta'_1,\ldots,\zeta'_j,\gamma\}$ as a set, and $C(p_1,\ldots,p_j)\in \mathcal{A}$.

Let $\rho_\alpha$ be the automorphism of $U_{\mathcal{A}}^s(\n_-)$ defined by
\begin{equation}\label{thetaadef}
\rho_\alpha(x)=c_{\alpha \lambda}^sx, x\in (U_{\mathcal{A}}^s(\n_-))_\lambda. \index[not]{r@$\rho_\alpha$}
\end{equation} 
By the definition $\rho_\alpha$ restricts to an automorphism of the algebra $U_{\mathcal{A}}^s([-\gamma,-\beta])$ for any $\gamma<\beta$.

Assume that $\alpha$ is not the last root in the fixed normal ordering of $\Delta_+$. Let $\alpha'>\alpha$ \index[not]{a@$\alpha'$} be the root following $\alpha$ in this normal ordering.
By the results of \cite{GY}, Section 2.3 or \cite{Y}, Section 6.5 for any $\beta>\alpha$ $D_\alpha$ extends to a $\rho_\alpha$--derivation of $U_{\mathcal{A}}^s([-\alpha',-\beta])$ which satisfies
\begin{equation}\label{dera}
D_\alpha(xy)=D_\alpha(x)y+\rho_\alpha(x)D_\alpha(y), x,y\in U_{\mathcal{A}}^s([-\alpha',-\beta]),
\end{equation}
and $U_{\mathcal{A}}^s([-\alpha,-\beta])$ is the skew polynomial extension \index{skew polynomial extension} generated by $U_{\mathcal{A}}^s([-\alpha',-\beta])$ and by $f_\alpha$ satisfying 
\begin{equation}\label{fax}
f_\alpha x=\rho_\alpha(x)f_\alpha+D_\alpha(x), x\in U_{\mathcal{A}}^s([-\alpha',-\beta]).
\end{equation}

To establish (\ref{fmn}) we prove the following more general lemma.
\begin{lemma}
Let $f_\delta\in U_{\mathcal{A}}^{s}(\g)$, $\delta\in \Delta_+$ be the quantum root vectors defined with the help of an arbitrary fixed normal ordering on $\Delta_+$.
Then for any $m\in \mathbb{N}$, $\alpha,\beta\in \Delta_+$, $\alpha<\beta$ and $x\in (U_{\mathcal{A}}^s([-\alpha',-\beta]))_\lambda$, $\lambda\in Q$ one has
\begin{equation}\label{fmax}
f_\alpha^m x=(c_{\alpha \lambda}^s)^m\sum_{p=0}^{m}S_m^pq_\alpha^{-p(p-1)}(c_{\alpha \lambda}^s)^{-p}D_\alpha^p(x)f_{\alpha}^{m-p}.
\end{equation}
\end{lemma}

\begin{proof}
Let 
$$
\left(\begin{array}{c} m \\ p \end{array}\right)_{q^2}=q^{mp-p^2}\left[\begin{array}{c} m \\ p \end{array}\right]_{q} \index[not]{ZZZ@$\left(\begin{array}{c} m \\ p \end{array}\right)_{q^2}$}
$$
be the q-binomial coefficients normalized as in \cite{KQ}. In terms of these coefficients identity (\ref{fmax}) can be rewritten as follows
\begin{equation}\label{fmax1}
f_\alpha^m x=\sum_{p=0}^{m}(c_{\alpha \lambda}^s)^{m-p}\left(\begin{array}{c} m \\ p \end{array}\right)_{q_\alpha^2}D_\alpha^p(x)f_{\alpha}^{m-p}.
\end{equation}

We prove (\ref{fmax1}) by induction over $m$. For $m=1$ this identity is equivalent to (\ref{fax}) as $\rho_\alpha(x)=c_{\alpha \lambda}^s x$ by the definition of $\rho_\alpha$.

Now assume that (\ref{fmax1}) holds for some $m$. Then  
\begin{equation}\label{fmax2}
f_\alpha^{m+1} x=\sum_{p=0}^{m}(c_{\alpha \lambda}^s)^{m-p}\left(\begin{array}{c} m \\ p \end{array}\right)_{q_\alpha^2}f_\alpha D_\alpha^p(x)f_{\alpha}^{m-p}.
\end{equation}
Using (\ref{fax}) with $x$ replaced by $D_\alpha^p(x)$ in the terms of the sum in (\ref{fmax2}) and recalling that by the definition of $D_\alpha$ one has $D_\alpha^p(x)\in (U_{\mathcal{A}}^s([-\alpha',-\beta]))_{\lambda+p\alpha}$, we obtain with the help of (\ref{thetaadef}) and (\ref{fax})
\begin{equation}\label{fmax3}
f_\alpha^{m+1} x=\sum_{p=0}^{m}(c_{\alpha \lambda}^s)^{m-p}(\left(\begin{array}{c} m \\ p \end{array}\right)_{q_\alpha^2}c_{\alpha \lambda+p\alpha}^s D_\alpha^p(x)f_{\alpha}^{m+1-p}+D_\alpha^{p+1}(x)f_{\alpha}^{m-p}).
\end{equation} 

Recalling the definition of $c_{\delta \gamma}=q^{\left\langle (2K_s+id)\delta,\gamma\right\rangle}$, $\delta,\gamma\in Q$ we infer $c_{\alpha \lambda+p\alpha}^s=q^{\left\langle (2K_s+id)\alpha,\lambda+p\alpha\right\rangle}=c_{\alpha \lambda}^sq^{\left\langle \alpha,\alpha\right\rangle}=c_{\alpha \lambda}^sq_\alpha^2$ as $K_s$ is a skew-symmetric operator. Thus (\ref{fmax3}) takes the form
\begin{equation}\label{fmax4}
f_\alpha^{m+1} x=\sum_{p=0}^{m}(c_{\alpha \lambda}^s)^{m-p}(\left(\begin{array}{c} m \\ p \end{array}\right)_{q_\alpha^2}c_{\alpha \lambda}^sq_\alpha^{2p} D_\alpha^p(x)f_{\alpha}^{m+1-p}+D_\alpha^{p+1}(x)f_{\alpha}^{m-p})=
\end{equation} 
$$
=\sum_{p=0}^{m+1}(c_{\alpha \lambda}^s)^{m+1-p}(\left(\begin{array}{c} m \\ p \end{array}\right)_{q_\alpha^2}q_\alpha^{2p}+\left(\begin{array}{c} m \\ p-1 \end{array}\right)_{q_\alpha^2})D_\alpha^{p}(x)f_{\alpha}^{m-p}=
$$
$$
=\sum_{p=0}^{m+1}(c_{\alpha \lambda}^s)^{m+1-p}\left(\begin{array}{c} m+1 \\ p \end{array}\right)_{q_\alpha^2}D_\alpha^p(x)f_{\alpha}^{m+1-p},
$$
where at the last step we used the identity
$$
\left(\begin{array}{c} m \\ p \end{array}\right)_{q_\alpha^2}q_\alpha^{2p}+\left(\begin{array}{c} m \\ p-1 \end{array}\right)_{q_\alpha^2}=\left(\begin{array}{c} m+1 \\ p \end{array}\right)_{q_\alpha^2}
$$
which can be found in \cite{KQ}, Proposition 6.1. This establishes the induction step and completes the proof of the lemma.

\end{proof}

Now (\ref{fmn}) follows from (\ref{fmax}) with $x=f_\beta^n$ as for each $p$ one has  $D_\alpha^p(f_\beta^n)\in U_{\mathcal{A}}^s([-\alpha',-\beta])$, and hence by Lemma \ref{segmPBWs} (vi) and by Remark \ref{segmPBWsrev} it can be represented as an $\mathcal{A}$--linear combination of the elements $f_{\zeta_1}^{p_1}\ldots f_{\zeta_k}^{p_k}$, $p_i\in \mathbb{N}$, where $\zeta_1<\ldots<\zeta_k=\beta$, and $\{\alpha,\zeta_1,\ldots,\zeta_k\}= [\alpha,\beta]$ as a set. 

Properties (\ref{dnp1}) and (\ref{dnp2}) follow from Lemma \ref{circ+} (ii) and using the fact that the weights of the left hand side and and of the right hand side of (\ref{fmn}) must be the same. This completes the proof of Lemma \ref{fmnm}.

\end{proof}

The next lemma shows how the adjoint action of the quantum root vectors commutes with the operators $P_{jp}$. 
\begin{lemma}\label{mainl4}
For any $j=1,\ldots ,R-1$, $p=1,\ldots, n_j$, $p<q\leq m_j$, and any $f\in \mathbb{C}_{jp}^{loc}[G]$ we have 
\begin{equation}\label{Pdef}
{\rm Ad}_{\widehat{s}_j}^0 f_{\delta_{jq}}^{(m)} (P_{jp}(f))=\sum_{k=0}^\infty P_{jp}^k({\rm Ad}_{\widehat{s}_j}^0(\hspace{-1em}\sum_{p_{jp+1},\ldots ,p_{jq}\in \mathbb{N}}\hspace{-1em} d_{p_{jp+1},\ldots, p_{jq}}^{mk}f_{\delta_{jp+1}}^{(p_{jp+1})}\ldots f_{\delta_{jq}}^{(p_{jq})})(f)),
\end{equation}
where only a finite number of terms in the sum are non-zero, $d_{p_{jp+1},\ldots, p_{jq}}^{mk}\in \mathcal{B}$,
\begin{equation}\label{dmk1}
d_{p_{jp+1},\ldots, p_{jq-1},m}^{mk}=\left\{ \begin{array}{l} 1~~{\rm if}~~k=p_{jp+1}=\ldots = p_{jq-1}=0 \\ 0 ~~{\rm otherwise} \end{array}
\right.,
\end{equation}
and for $0\leq n<m$ and any $k$
\begin{equation}\label{dmk2}
d_{0,\ldots, 0,n}^{mk}=0.
\end{equation}

In particular, 
\begin{equation}\label{Pdef1}
{\rm Ad}_{\widehat{s}_j}^0 f_{\delta_{jp+1}}^{(m)} (P_{jp}(f))=P_{jp}({\rm Ad}_{\widehat{s}_j}^0(f_{\delta_{jp+1}}^{(m)})(f)).	
\end{equation}

\end{lemma}

\begin{proof}
By formula (\ref{mainfpq}) and by the definition of $P_{jp}$ we have 
$$
{\rm Ad}_{\widehat{s}_j}^0 f_{\delta_{jq}}^{(m)} (P_{jp}(f))=\sum_{n=0}^\infty (-1)^nq_{\delta_{jp}}^{-\frac{(n-1)n}{2}}(c_{\delta_{jp}\delta_{jq}}^{\widehat{s}_j})^{-mn}{\rm Ad}_{\widehat{s}_j}^0 (f_{\delta_{jp}}^{(n)}f_{\delta_{jq}}^{(m)})(f)\otimes B_{jp}^n.	
$$
Using (\ref{fmnm}) we can rewrite this formula as follows
$$
{\rm Ad}_{\widehat{s}_j}^0 f_{\delta_{jq}}^{(m)} (P_{jp}(f))=\sum_{n=0}^\infty \sum_{k=0}^n(-1)^nq_{\delta_{jp}}^{-\frac{(n-1)n}{2}}\omega_0(q_{\delta_{jp}}^{k(n-1)})(c_{\delta_{jp}\delta_{jq}}^{\widehat{s}_j})^{-mn}\omega_0((c_{\delta_{jp}\delta_{jq}}^{\widehat{s}_j})^{mn})\times
$$
$$
\times {\rm Ad}_{\widehat{s}_j}^0 f_{\delta_{jp}}^{(n-k)}({\rm Ad}_{\widehat{s}_j}^0(\hspace{-1em}\sum_{p_{jp+1},\ldots,p_{jq}\in \mathbb{N}}\hspace{-1em} d_{p_{jp},\ldots, p_{jq}}^{mk}f_{\delta_{jp+1}}^{(p_{jp+1})}\ldots f_{\delta_{jq}}^{(p_{jq})})(f))\otimes B_{jp}^n.	
$$ 
Now recalling that $\omega_0(q)=q$ and swapping the order of summation in the last formula we get
$$
{\rm Ad}_{\widehat{s}_j}^0 f_{\delta_{jq}}^{(m)} (P_{jp}(f))=
$$
$$
=\sum_{n=0}^\infty \sum_{k=0}^n(-1)^nq_{\delta_{jp}}^{-\frac{(n-1)n}{2}}q_{\delta_{jp}}^{k(n-1)}{\rm Ad}_{\widehat{s}_j}^0 f_{\delta_{jp}}^{(n-k)}({\rm Ad}_{\widehat{s}_j}^0(\hspace{-1em} \sum_{p_{jp+1},\ldots, p_{jq}\in \mathbb{N}}\hspace{-1em} d_{p_{jp+1},\ldots, p_{jq}}^{mk}f_{\delta_{jp+1}}^{(p_{j p+1})}\ldots f_{\delta_{jq}}^{(p_{jq})})(f))\otimes B_{jp}^n=
$$ 
$$
=\sum_{k=0}^\infty \sum_{n=k}^\infty (-1)^nq_{\delta_{jp}}^{-\frac{(n-1)(n-2k)}{2}}{\rm Ad}_{\widehat{s}_j}^0 f_{\delta_{jp}}^{(n-k)}({\rm Ad}_{\widehat{s}_j}^0(\hspace{-1em} \sum_{p_{jp+1},\ldots,p_{jq}\in \mathbb{N}}\hspace{-1em} d_{p_{jp+1},\ldots, p_{jq}}^{mk}f_{\delta_{j p+1}}^{(p_{j p+1})}\ldots f_{\delta_{jq}}^{(p_{jq})})(f))\otimes B_{jp}^n=
$$
$$
=\sum_{k=0}^\infty P_{jp}^k({\rm Ad}_{\widehat{s}_j}^0(\hspace{-1em} \sum_{p_{jp+1},\ldots,p_{jq}\in \mathbb{N}}\hspace{-1em} d_{p_{jp+1},\ldots, p_{jq}}^{mk}f_{\delta_{jp+1}}^{(p_{jp+1})}\ldots f_{\delta_{jq}}^{(p_{jq})})(f)).
$$

Properties (\ref{dmk1}) and (\ref{dmk2}) follow from (\ref{dnp1}) and (\ref{dnp2}), respectively.

Formula (\ref{Pdef1}) is obtained in a similar way using the relation $f_{\delta_{jp}}^{(n)}f_{\delta_{jp+1}}^{(m)}=(c_{\delta_{jp}\delta_{jp+1}}^{\widehat{s}_j})^{mn}f_{\delta_{jp+1}}^{(m)}f_{\delta_{jp}}^{(n)}$.

\end{proof}

Now we have all prerequisites to prove Proposition \ref{BasP}.

\vskip 0.3cm
\noindent
{\em Proof of Proposition \ref{BasP}}
Firstly, composition (\ref{comp}) is well defined by Lemma \ref{Jvan} (ii). This proves the first part of Proposition \ref{BasP}. 

Next, we shall prove a statement which implies (\ref{projp}). To state it we note that by Lemma \ref{Jvan} (ii) for any $1\leq i\leq j\leq R-1$, $1\leq q\leq n_i$, $1\leq p\leq n_j$ such that $q\leq p$ if $i=j$ the composition
$$
P_{iq,jp}:=P_{iq}\ldots P_{in_i}\circ \varpi_i\circ P_{i+1 1}\ldots P_{i+1 n_{i+1}}\circ \varpi_{i+1}  \ldots  \varpi_{j-1} \circ P_{j1}\ldots P_{jp}: \mathbb{C}_{jp+1}^{loc}[G]\to \mathbb{C}_{iq}^{loc}[G] \index[not]{P@$P_{iq,jp}$}
$$ 
is a well--defined operator. We shall show that 
\begin{equation}\label{adpij}
{\rm Ad}_{\widehat{s}_i}^0 x P_{iq,jp}(f)=0
\end{equation}
for all $x\in U_{\mathcal{B}}^{\widehat{s}_i,res}([-\delta_{iq},-\delta_{jp}^i])$ which have no weight components of zero weight and $f\in \mathbb{C}_{jp+1}^{loc}[G]$, where $\delta_{jp}^i=(w_1\ldots w_{i-1})^{-1}\beta_{jp}=w_i\ldots w_{j-1}\delta_{jp}\in [\delta_{iq},\delta_{i o_i}]$.

Since the algebra $U_{\mathcal{B}}^{\widehat{s}_i,res}([-\delta_{iq},-\delta_{jp}^i])$ is generated by the elements $f_{\delta_{ir}}^{(m)}=f_{\delta_{ir}}^{(m)}(\widehat{s}_i)$, $m\in \mathbb{N}$, $\delta_{ir}\in [\delta_{iq},\delta_{jp}^i]$ it suffices to verify (\ref{adpij}) for such $x$ with $m>0$.

Equip the set of pairs $(iq)$, $1\leq i\leq R-1$, $1\leq q\leq n_i$ with the lexicographic order \index{lexicographic order} such that $(iq)<(jp)$ if $i<j$ or $i=j$ and $q<p$.
Fix $1\leq j\leq R-1$, and $1\leq p\leq n_j$. To establish (\ref{adpij}) we argue by backwards induction over the pairs $(iq)\leq (jp)$ with respect to the lexicographic order starting from $i=j$ and $p=q$. In this case property (\ref{adpij}) follows from (\ref{mainprop2}).

Now assume that (\ref{adpij}) holds for all pairs $(iq)$ such that $(ht)<(iq)\leq (jp)$ for some pair $(ht)$. We shall prove that (\ref{adpij}) is valid for $(iq)=(ht)$.

We have to consider two cases. Firstly suppose that $1\leq t<n_h$. Then for $m>0$, $\delta_{hr}\in [\delta_{ht},\delta_{jp}^h]$ one has by (\ref{Pdef})
\begin{equation}\label{Pdef2}
{\rm Ad}_{\widehat{s}_h}^0 f_{\delta_{hr}}^{(m)} (P_{ht,jp}(f))={\rm Ad}_{\widehat{s}_h}^0 f_{\delta_{hr}}^{(m)} P_{ht}P_{ht+1,jp}(f)=
\end{equation}
$$
=\sum_{k=0}^\infty P_{ht}^k({\rm Ad}_{\widehat{s}_h}^0(\hspace{-1em}\sum_{p_{ht+1},\ldots ,p_{hr}\in \mathbb{N}}\hspace{-1em} d_{p_{ht+1},\ldots, p_{hr}}^{mk}f_{\delta_{ht+1}}^{(p_{ht+1})}\ldots f_{\delta_{hr}}^{(p_{hr})})P_{ht+1,jp}(f))=
$$
$$
=\sum_{k=0}^\infty P_{ht}^k\omega_0(d_{0,\ldots, 0}^{mk})P_{ht+1,jp}(f),
$$
where at the last step we used the induction assumption. 

Now recall that by (\ref{dmk2}) for $m>0$ and any $k$ one has
$$
d_{0,\ldots, 0,0}^{mk}=0.
$$
Therefore the right hand side of (\ref{Pdef2}) vanishes. This establishes the induction step if $1\leq t< n_h$.

Now suppose that $t=n_h$. Then for $m>0$, $\delta_{hr}\in [\delta_{ht},\delta_{jp}^i]$ one has by (\ref{Pdef})
\begin{equation}\label{Pdef3}
{\rm Ad}_{\widehat{s}_h}^0 f_{\delta_{hr}}^{(m)} (P_{hn_h,jp}(f))={\rm Ad}_{\widehat{s}_h}^0 f_{\delta_{hr}}^{(m)} P_{hn_h} \varpi_h P_{h+11,jp}(f)=
\end{equation}
$$
=\sum_{k=0}^\infty P_{hn_h}^k({\rm Ad}_{\widehat{s}_h}^0(\hspace{-1em} \sum_{p_{hn_{h+1}},\ldots ,p_{hr}\in \mathbb{N}}\hspace{-1em} d_{p_{hn_h+1},\ldots, p_{hr}}^{mk}f_{\delta_{hn_h+1}}^{(p_{hn_h+1})}(\widehat{s}_h)\ldots f_{\delta_{hr}}^{(p_{hr})}(\widehat{s}_h))\varpi_h P_{h+11,jp}(f)).
$$

By (\ref{Ads0}) and (\ref{phiAdloc}) we have
\begin{equation}\label{phiAdloc0}
{\rm Ad}_{\widehat{s}_h}^0(x)\varpi_h=\varpi_h{\rm Ad}_{\widehat{s}_{h+1}}^0(\omega_0'(T_{w_h}^{-1}\omega_0'(x))), x\in U_\mathcal{B}^{\widehat{s}_h,res}([-\delta_{hn_h+1},-\delta_{ho_h}]).
\end{equation} 
Therefore (\ref{Pdef3}) takes the form
\begin{equation}\label{Pdef4}
{\rm Ad}_{\widehat{s}_h}^0 f_{\delta_{hr}}^{(m)} (P_{hn_h,jp}(f))=
\end{equation}
$$
=\sum_{k=0}^\infty P_{hn_h}^k(\varpi_h ({\rm Ad}_{\widehat{s}_{h+1}}^0(\hspace{-1em} \sum_{p_{hn_{h+1}},\ldots ,p_{hr}\in \mathbb{N}}\hspace{-1em} d_{p_{hn_h+1},\ldots, p_{hr}}^{mk}\omega_0'(T_{w_h}^{-1}\omega_0'(f_{\delta_{hn_h+1}}^{(p_{hn_h+1})}(\widehat{s}_h)\ldots f_{\delta_{hr}}^{(p_{hr})}(\widehat{s}_h))))P_{h+11,jp}(f)))=
$$
$$
=\sum_{k=0}^\infty P_{hn_h}^k(\varpi_h ({\rm Ad}_{\widehat{s}_{h+1}}^0(\hspace{-1em} \sum_{p_{hn_{h+1}},\ldots ,p_{hr}\in \mathbb{N}}\hspace{-1em} d_{p_{hn_h+1},\ldots, p_{hr}}^{mk}c_{\delta_{hn_h+1}}^{-p_{hn_h+1}}\ldots c_{\delta_{hr}}^{-p_{hr}}\omega_0'(T_{w_h}^{-1} (f_{\delta_{hr}}^{(p_{hr})}(\widehat{s}_h)\ldots f_{\delta_{hn_h+1}}^{(p_{hn_h+1})}(\widehat{s}_h))))P_{h+11,jp}(f))),
$$
where at the last step we used (\ref{fo0}) and the fact that $\omega_0'$ is an  algebra antiautomorphism such that $\omega_0'(q)=q^{-1}$.

Next, recall that by (\ref{twjb})
$$
T_{w_1\ldots w_{h-1}}^{-1}f_\beta^{(n)}(s)=f_{(w_1\ldots w_{h-1})^{-1}\beta}^{(n)}(\widehat{s}_h)
$$
for $\beta\in [\beta_{h1},\beta_D]$, where for $\alpha\in \Delta_+$ and $h>1$ the elements $f_{\alpha}(\widehat{s}_h)$ are defined using the element $\widehat{s}_h$ and the normal ordering on $\Delta_+$ introduced in Remark \ref{Djord}.

Applying this identity to $\beta=w_1\ldots w_{h-1}\delta_{hu}$, $u=n_h+1,\ldots o_h$ we obtain
$$
T_{w_1\ldots w_{h-1}}^{-1}f_\beta^{(n)}(s)=f_{\delta_{hu}}^{(n)}(\widehat{s}_h),
$$
and hence
$$
T_{w_h}^{-1}f_{\delta_{hu}}^{(n)}(\widehat{s}_h)=T_{w_h}^{-1}T_{w_1\ldots w_{h-1}}^{-1}f_\beta^{(n)}(s)=T_{w_1\ldots w_h}^{-1}f_\beta^{(n)}(s)
$$
as the decomposition $w_1\ldots w_h$ is reduced.

Using (\ref{twjb}) with $j=h+1$ in the right hand side of the last identity we obtain
\begin{equation}\label{Tdhu}
T_{w_h}^{-1}f_{\delta_{hu}}^{(n)}(\widehat{s}_h)=f_{(w_1\ldots w_h)^{-1}(\beta)}^{(n)}(\widehat{s}_{h+1})=f_{\delta_{h+1 u-n_h}}^{(n)}(\widehat{s}_{h+1}),
\end{equation}
where at the last step we used the identity $(w_1\ldots w_h)^{-1}w_1\ldots w_{h-1}\delta_{hu}=w_h^{-1}\delta_{hu}=\delta_{h+1 u-n_h}$, $u=n_h+1,\ldots o_h$ following from the definition of the roots $\delta_{hu}$.

Now applying (\ref{Tdhu}) in the right hand side of (\ref{Pdef4}) we get
\begin{equation}\label{Pdef5}
{\rm Ad}_{\widehat{s}_h}^0 f_{\delta_{hr}}^{(m)} (P_{hn_h,jp}(f))=
\end{equation}
$$
=\sum_{k=0}^\infty P_{hn_h}^k(\varpi_h ({\rm Ad}_{\widehat{s}_{h+1}}^0(\hspace{-1em}\sum_{p_{hn_{h+1}},\ldots ,p_{hr}\in \mathbb{N}}\hspace{-1em} d_{p_{hn_h+1},\ldots, p_{hr}}^{mk}c_{\delta_{hn_h+1}}^{-p_{hn_h+1}}\ldots c_{\delta_{hr}}^{-p_{hr}}\times
$$
$$
\times\omega_0'((f_{\delta_{h+1 r-n_h}}^{(p_{hr})}(\widehat{s}_{h+1})\ldots f_{\delta_{h+1 1}}^{(p_{hn_h+1})}(\widehat{s}_{h+1}))))P_{h+11,jp}(f)))=
$$
$$
=\sum_{k=0}^\infty P_{hn_h}^k(\varpi_h ({\rm Ad}_{\widehat{s}_{h+1}}^0(\hspace{-1em} \sum_{p_{hn_{h+1}},\ldots ,p_{hr}\in \mathbb{N}}\hspace{-1em} d_{p_{hn_h+1},\ldots, p_{hr}}^{mk}c_{\delta_{hn_h+1}}^{-p_{hn_h+1}}\ldots c_{\delta_{hr}}^{-p_{hr}}c_{\delta_{h+1 r-n_h}}^{p_{hr}}\ldots c_{\delta_{h+1 1}}^{p_{hn_h+1}}\times
$$
$$
\times f_{\delta_{h+1 1}}^{(p_{hn_h+1})}(\widehat{s}_{h+1})\ldots f_{\delta_{h+1 r-n_h}}^{(p_{hr})}(\widehat{s}_{h+1})))P_{h+11,jp}(f)),
$$
where at the last step we used (\ref{fo0}) and the fact that $\omega_0'$ is an  algebra antiautomorphism such that $\omega_0'(q)=q^{-1}$.

By the induction assumption (\ref{Pdef5}) takes the form
\begin{equation}\label{Pdef6}
{\rm Ad}_{\widehat{s}_h}^0 f_{\delta_{hr}}^{(m)} (P_{hn_h,jp}(f))=\sum_{k=0}^\infty P_{hn_h}^k(\varpi_h (\omega_0(d_{0,\ldots, 0}^{mk})P_{h+11,jp}(f))).
\end{equation}

Finally recall that by (\ref{dmk2})
for $m>0$ and any $k$ one has
$$
d_{0,\ldots, 0,0}^{mk}=0.
$$
Therefore the right hand side of (\ref{Pdef6}) vanishes. This establishes the induction step if $t=n_h$ and completes the proof of (\ref{adpij}). Property (\ref{projp}) follows from (\ref{adpij}) with $i=1$ and $q=1$.

Property (\ref{Pcomm}) is established in a similar way with the help of (\ref{dmk1}), (\ref{Pdef1}) and (\ref{adpij}). This completes the proof of part (ii).

\qed

In conclusion, using the operators $P_{jp}$ and the elements $B_{jp}$ we define proper quantum analogues of monomials in functions $\varphi_{jp}$. They will play a crucial rule in establishing a quantum group version of the Skryabin equivalence in Section \ref{Skr} and in the proof of the De--Concini--Kac--Procesi conjecture.

\begin{proposition}\label{Bpbas1}
For $m_{ij}, k_{ij}\in \mathbb{N}$ the elements 
$$
B_{\{m_{ij}\}}:=P_{11}(\ldots P_{1 n_1-1}(P_{1n_1}(\varpi_1 (P_{21}(\ldots P_{2 n_2-1} (P_{2n_2} (\varpi_2( \ldots  \varpi_{R-2} ( P_{R-11}(\ldots 
$$
$$
\ldots P_{R-1 n_{R-1}-2}(P_{R-1n_{R-1}-1}(B_{R-1n_{R-1}}^{m_{R-1n_{R-1}}})B_{R-1n_{R-1}-1}^{m_{R-1 n_{R-1}-1}})\ldots B_{R-1 1}^{m_{R-1 1}}))\ldots ))B_{2n_2}^{m_{2n_2}})\ldots 
$$
$$
\ldots )B_{21}^{m_{21}}))B_{1 n_1}^{m_{1n_1}})\ldots)B_{11}^{m_{11}} \in \mathbb{C}_{11}^{loc}[G] \index[not]{B@$B_{\{m_{ij}\}}$}
$$
satisfy 
$$
{\rm Ad}_s^0(f_{\beta_{11}}^{(k_{11})}\ldots f_{\beta_{R-1n_{R-1}}}^{(k_{R-1 n_{R-1}})})(B_{\{m_{ij}\}})=\left\{\begin{array}{l} c_{\{m_{ij}\}}\prod_{i=1}^{R-1} \prod_{j=1}^{n_i}q_{\delta_{ij}}^{-\frac{m_{ij}(m_{ij}-1)}{2}} ~{\rm if}~m_{ij}=k_{ij}~{\rm for}~ i=1,\ldots , R-1, j=1,\ldots n_i \\ 
\\
0 ~{\rm if}~k_{ij}=m_{ij} ~{\rm if}~ (ij)<(pq)~{\rm and}~ k_{pq}>m_{pq}~{\rm for~some}~(pq) \end{array}\right. ,
$$
where the pairs $(ij)$ are ordered according to the lexicographic order \index{lexicographic order} and $c_{\{m_{ij}\}}\in \{\pm q^{\mathbb{Z}}\}$. \index[not]{c@$c_{\{m_{ij}\}}$}
Thus
$$
{\rm Ad}_s^0(f_{\beta_{11}}^{k_{11}}\ldots f_{\beta_{R-1n_{R-1}}}^{k_{R-1 n_{R-1}}})(B_{\{m_{ij}\}})=\left\{\begin{array}{l} c_{\{m_{ij}\}}\prod_{i=1}^{R-1} \prod_{j=1}^{n_i}q_{\delta_{ij}}^{-\frac{m_{ij}(m_{ij}-1)}{2}}[m_{ij}]_{q_{\delta_{ij}}}! ~{\rm if}~m_{ij}=k_{ij}~{\rm for}~ i=1,\ldots , R-1, j=1,\ldots n_i \\ 
\\
0 ~{\rm if}~k_{ij}=m_{ij} ~{\rm if}~ (ij)<(pq)~{\rm and}~ k_{pq}>n_{pq}~{\rm for~some}~(pq) \end{array}\right. .
$$

In particular, the elements $B_{\{m_{ij}\}}$ are linearly independent for different $\{m_{ij}\}$, $i=1,\ldots, R-1$, $j=1,\ldots n_i$.
\end{proposition}

\begin{proof}
The proof follows from Lemmas \ref{mainl2}, \ref{mainl3}, \ref{Jvan}, \ref{mainl4}, formula  (\ref{mainprop2}) and Proposition \ref{BasP}. We shall prove the statement by induction.

First observe that by (\ref{mainprop2}) 
$$
{\rm Ad}_s^0f_{\beta_{11}}^{(k)}P_{11}(\ldots P_{1 n_1-1}(P_{1n_1}(\varpi_1 (P_{21}(\ldots P_{2 n_2-1} (P_{2n_2} (\varpi_2( \ldots  \varpi_{R-2} ( P_{R-11}(\ldots 
$$
$$
\ldots P_{R-1 n_{R-1}-2}(P_{R-1n_{R-1}-1}(B_{R-1n_{R-1}}^{m_{R-1n_{R-1}}})B_{R-1n_{R-1}-1}^{m_{R-1 n_{R-1}-1}})\ldots B_{R-1 1}^{m_{R-1 1}}))\ldots ))B_{2n_2}^{m_{2n_2}})\ldots )B_{21}^{m_{21}}))B_{1 n_1}^{m_{1n_1}})\ldots B_{12}^{m_{12}})=0
$$
for any $k>0$ and recall that right multiplication by $B_{11}$ gives rise to an operator on $\mathbb{C}_{11}^{loc}[G]$ by Lemma \ref{quotact}. Therefore from (\ref{mainf}) for $(jp)=(11)$ and
$$
f=P_{11}(\ldots P_{1 n_1-1}(P_{1n_1}(\varpi_1 (P_{21}(\ldots P_{2 n_2-1} (P_{2n_2} (\varpi_2( \ldots  \varpi_{R-2} ( P_{R-11}(\ldots 
$$
$$
\ldots P_{R-1 n_{R-1}-2}(P_{R-1n_{R-1}-1}(B_{R-1n_{R-1}}^{m_{R-1n_{R-1}}})B_{R-1n_{R-1}-1}^{m_{R-1 n_{R-1}-1}})\ldots B_{R-1 1}^{m_{R-1 1}}))\ldots ))B_{2n_2}^{m_{2n_2}})\ldots )B_{21}^{m_{21}}))B_{1 n_1}^{m_{1n_1}})\ldots B_{12}^{m_{12}})
$$
we have 
$$
{\rm Ad}_s^0f_{\beta_{11}}^{(k_{11})}(B_{\{m_{ij}\}})=	
$$
$$
=q_{\delta_{11}}^{-\frac{k_{11}(k_{11}-1)}{2}} P_{11}(\ldots P_{1 n_1-1}(P_{1n_1}(\varpi_1 (P_{21}(\ldots P_{2 n_2-1} (P_{2n_2} (\varpi_2( \ldots  \varpi_{R-2} ( P_{R-11}(\ldots 
$$
$$
\ldots P_{R-1 n_{R-1}-2}(P_{R-1n_{R-1}-1}(B_{R-1n_{R-1}}^{m_{R-1n_{R-1}}})B_{R-1n_{R-1}-1}^{m_{R-1 n_{R-1}-1}})\ldots B_{R-1 1}^{m_{R-1 1}}))\ldots ))B_{2n_2}^{m_{2n_2}})\ldots )B_{21}^{m_{21}}))B_{1 n_1}^{m_{1n_1}})\ldots B_{12}^{m_{12}})
$$
if $k_{11}=m_{11}$ and
$$
{\rm Ad}_s^0f_{\beta_{11}}^{(k_{11})}(B_{\{m_{ij}\}})=0
$$
if $k_{11}>m_{11}$.

Now suppose that for some $1\leq p\leq R-1$, $1\leq q\leq n_i$
\begin{equation}\label{basind}
{\rm Ad}_s^0(f_{\beta_{11}}^{(k_{11})}\ldots f_{\beta_{pq}}^{(k_{pq})})(B_{\{k_{11},\ldots , k_{pq}, m_{p q+1},\ldots , m_{R-1 n_{R-1}}\}})=	
\end{equation}
$$
=d_{pq}\prod_{(ij)\leq (pq)} q_{\delta_{ij}}^{-\frac{k_{ij}(k_{ij}-1)}{2}} P_{\leq {pq}}(P_{p q+1}\ldots P_{p n_p-1}(P_{pn_p}(\varpi_p (P_{p+11}\ldots P_{p+1 n_{p+1}-1} (P_{p+1 n_{p+1}} (\varpi_{p+1}( \ldots  \varpi_{R-2} ( P_{R-11}\ldots 
$$
$$
\ldots P_{R-1 n_{R-1}-2}(P_{R-1n_{R-1}-1}(B_{R-1n_{R-1}}^{m_{R-1n_{R-1}}})B_{R-1n_{R-1}-1}^{m_{R-1 n_{R-1}-1}})\ldots B_{R-1 1}^{m_{R-1 1}}))\ldots )B_{p+2 n_{p+2}}^{m_{p+2 n_{p+2}}})\ldots B_{p+1 1}^{m_{p+11}}))B_{p n_p}^{m_{pn_p}})\ldots B_{pq+1}^{m_{p q+1}}),
$$
where we assume that $m_{p n_p+1}=m_{p+1 1}$, $(p~ n_p+1)=(p+1 ~1)$, $P_{p n_p+1}=\varpi_p P_{p+1 1}$ and $d_{pq}\in \{\pm q^{\mathbb{Z}}\}$.

Then by (\ref{Pcomm}) and by the induction assumption
$$
{\rm Ad}_s^0(f_{\beta_{11}}^{(k_{11})}\ldots f_{\beta_{p q+1}}^{(k_{p q+1})})(B_{\{k_{11},\ldots, k_{pq}, m_{p q+1},\ldots ,m_{R-1 n_{R-1}}\}})=
$$
$$
=d_{pq}c_{p q+1}\prod_{(ij)\leq (pq)} q_{\delta_{ij}}^{-\frac{k_{ij}(k_{ij}-1)}{2}} P_{\leq {pq}}{\rm Ad}_{\widehat{s}_p}^0(f_{\delta_{p q+1}}^{(k_{p q+1})})(P_{p q+1}\ldots P_{p n_p-1}(P_{pn_p}(\varpi_p (P_{p+11}\ldots P_{p+1 n_{p+1}-1} (P_{p+1 n_{p+1}} (\varpi_{p+1}( \ldots   
$$
$$
\ldots\varpi_{R-2} ( P_{R-11}\ldots P_{R-1 n_{R-1}-2}(P_{R-1n_{R-1}-1}(B_{R-1n_{R-1}}^{m_{R-1n_{R-1}}})B_{R-1n_{R-1}-1}^{m_{R-1 n_{R-1}-1}})\ldots B_{R-1 1}^{m_{R-1 1}}))\ldots )B_{p+2 n_{p+2}}^{m_{p+2 n_{p+2}}})\ldots 
$$
$$
\ldots B_{p+1 1}^{m_{p+11}}))B_{p n_p}^{m_{pn_p}})\ldots B_{pq+1}^{m_{p q+1}}),
$$
where $c_{p q+1}\in \{\pm q^{\mathbb{Z}}\}$ is defined in (\ref{Pcomm}). Denote $d_{p q+1}=d_{pq}c_{p q+1}\in \{\pm q^{\mathbb{Z}}\}$.

The last formula can be simplified using (\ref{mainf}) with $(jp)=(pq+1)$,
$$
f=P_{p q+1}(\ldots P_{p n_p-1}(P_{pn_p}(\varpi_p (P_{p+11}\ldots P_{p+1 n_{p+1}-1} (P_{p+1 n_{p+1}} (\varpi_{p+1}( \ldots  \varpi_{R-2} ( P_{R-11}\ldots 
$$
$$
\ldots P_{R-1 n_{R-1}-2}(P_{R-1n_{R-1}-1}(B_{R-1n_{R-1}}^{m_{R-1n_{R-1}}})B_{R-1n_{R-1}-1}^{m_{R-1 n_{R-1}-1}})\ldots B_{R-1 1}^{m_{R-1 1}}))\ldots )B_{p+2 n_{p+2}}^{m_{p+2 n_{p+2}}})\ldots
$$
$$
\ldots B_{p+1 1}^{m_{p+11}}))B_{p n_p}^{m_{pn_p}})\ldots B_{pq+2}^{m_{p q+2}}),
$$ 
the fact that 
$$
{\rm Ad}_{\widehat{s}_p}^0(f_{\delta_{p q+1}}^{(k)})(P_{p q+1}(\ldots P_{p n_p-1}(P_{pn_p}(\varpi_p (P_{p+11}\ldots P_{p+1 n_{p+1}-1} (P_{p+1 n_{p+1}} (\varpi_{p+1}( \ldots  \varpi_{R-2} ( P_{R-11}\ldots 
$$
$$
\ldots P_{R-1 n_{R-1}-2}(P_{R-1n_{R-1}-1}(B_{R-1n_{R-1}}^{m_{R-1n_{R-1}}})B_{R-1n_{R-1}-1}^{m_{R-1 n_{R-1}-1}})\ldots B_{R-1 1}^{m_{R-1 1}}))\ldots )B_{p+2 n_{p+2}}^{m_{p+2 n_{p+2}}})\ldots
$$
$$
\ldots B_{p+1 1}^{m_{p+11}}))B_{p n_p}^{m_{pn_p}})\ldots B_{pq+2}^{m_{p q+2}}))=0
$$
for any $k>0$ by (\ref{mainprop2}), Lemma \ref{Jvan} and recalling that right multiplication by $B_{p q+1}$ gives rise to an operator on $\mathbb{C}_{p q+1}^{loc}[G]$ by Lemma \ref{quotact}. This yields 
$$
{\rm Ad}_s^0(f_{\beta_{11}}^{(k_{11})}\ldots f_{\beta_{p q+1}}^{(k_{p q+1})})(B_{\{k_{11},\ldots, k_{pq}, m_{p q+1},\ldots ,m_{R-1 n_{R-1}}\}})=
$$
$$
=d_{p q+1}\prod_{(ij)\leq (pq+1)} q_{\delta_{ij}}^{-\frac{k_{ij}(k_{ij}-1)}{2}}(P_{\leq p q+1}(P_{p q+2}\ldots P_{p n_p-1}(P_{pn_p}(\varpi_p (P_{p+11}\ldots P_{p+1 n_{p+1}-1} (P_{p+1 n_{p+1}} (\varpi_{p+1}( \ldots   
$$
$$
\ldots \varpi_{R-2} ( P_{R-11}\ldots\ldots P_{R-1 n_{R-1}-2}(P_{R-1n_{R-1}-1}(B_{R-1n_{R-1}}^{m_{R-1n_{R-1}}})B_{R-1n_{R-1}-1}^{m_{R-1 n_{R-1}-1}})\ldots B_{R-1 1}^{m_{R-1 1}}))\ldots
$$
$$
\ldots )B_{p+2 n_{p+2}}^{m_{p+2 n_{p+2}}})\ldots B_{p+1 1}^{m_{p+11}}))B_{p n_p}^{m_{pn_p}})\ldots B_{pq+2}^{m_{p q+2}}))
$$
if $k_{p q+1}=m_{p q+1}$ and 
$$
{\rm Ad}_s^0(f_{\beta_{11}}^{(k_{11})}\ldots f_{\beta_{p q+1}}^{(k_{p q+1})})(B_{\{k_{11},\ldots, k_{pq}, m_{p q+1},\ldots ,m_{R-1 n_{R-1}}\}})=0
$$
if $k_{pq+1}>m_{pq+1}$. This establishes the induction step and completes the proof of the proposition.

\end{proof}

Since the ordering on the set of roots $\beta_{ij}$, $i=1,\ldots, R-1$, $j=1,\ldots, n_i$ induced by the lexicographic ordering of the pairs $(ij)$ coincides with the normal ordering on the segment $\Delta_{\m_+}$, we can rewrite the properties of the elements $B_{\{m_{ij}\}}$ as follows to simplify the notation for later use. 

\begin{corollary}\label{Bpbas}
For $m_1,\ldots ,m_c\in \mathbb{N}$ there are elements 
$$
B_{m_1\ldots m_c}\in \mathbb{C}_{11}^{loc}[G] \index[not]{B@$B_{m_1\ldots m_c}$}
$$
which satisfy for $r_1,\ldots, r_c\in \mathbb{N}$
$$
{\rm Ad}_s^0(f_{\beta_1}^{(r_1)}\ldots f_{\beta_c}^{(r_c)})(B_{m_1\ldots m_c})=\left\{\begin{array}{l} c_{m_1\ldots m_c}(q) ~~{\rm if}~~m_p=r_p~~{\rm for}~~ p=1,\ldots , c \\ 
\\
0 ~~{\rm if}~~r_i=m_i, i=1,\ldots, p-1~~{\rm and}~~ r_p>m_p~~{\rm for~~some}~~p\in \{1,\ldots , c\} \end{array}\right.,
$$
where $c_{m_1\ldots m_c}(q)\in \{\pm q^{\mathbb{Z}}\}$, \index[not]{c@$c_{m_1\ldots m_c}(q)$} and hence
$$
{\rm Ad}_s^0(f_{\beta_1}^{r_1}\ldots f_{\beta_c}^{r_c})(B_{m_1\ldots m_c})=\left\{\begin{array}{l}  c_{m_1\ldots m_c}'(q)~~{\rm if}~~m_p=r_p~~{\rm for}~~ p=1,\ldots , d \\ 
\\
0 ~~{\rm if}~~r_i=m_i, i=1,\ldots, p-1~~{\rm and}~~ r_p>m_p~~{\rm for~~some}~~p\in \{1,\ldots , c\} \end{array}\right. ,
$$
where $c_{m_1\ldots m_c}'(q)=c_{m_1\ldots m_c}(q)\prod_{p=1}^c [m_p]_{q_{\beta_p}}!$. \index[not]{c@$c_{m_1\ldots m_c}'(q)$}
\end{corollary}


\section{A description of q-W--algebras in terms of Zhelobenko type operators}\label{qWdescr}

\pagestyle{myheadings}
\markboth{CHAPTER~\thechapter.~ZHELOBENKO TYPE OPERATORS FOR Q-W--ALGEBRAS}{\thesection.~A DESCRIPTION OF Q-W--ALGEBRAS IN TERMS OF ZHELOBENKO TYPE OPERATORS}

\setcounter{equation}{0}
\setcounter{theorem}{0}

Now we are in a position to describe q-W--algebras in terms of the operator $P_c$ introduced in the previous section. Recall that the q-W--algebras $W_{\mathcal{B}}^s(G)$ are only defined when the value of the parameter $\kappa$ is equal to one. Therefore in this section we always assume that $\kappa=1$. As a $\mathcal{B}$--module the q-W--algebra $W_{\mathcal{B}}^s(G)$ is the space of $\mathbb{C}_{\mathcal{B}}^s[M_+]$--invariants in $Q_{\mathcal{B}}$ with respect to the adjoint action. In order to use the operator $P_c$ for the description of this space we shall transfer the results of Proposition \ref{BasP} from $\mathbb{C}_{11}^{loc}[G]$ to a localization $Q_{\mathcal{B}}^{loc}$ of $Q_{\mathcal{B}}$ using a natural extension of the  $\mathbb{C}_{\mathcal{B}}^s[M_+]$--module homomorphism $\phi:\mathbb{C}_{\mathcal{B}}^s[G] \rightarrow Q_{\mathcal{B}}$ to a homomorphism $\mathbb{C}_{\mathcal{B}}^{s,loc}[G] \rightarrow Q_{\mathcal{B}}^{loc}$. Recall that according to Proposition \ref{kerphi} ${{I}^{11}_{\mathcal{B}}}$ belongs to the kernel of the homomorphism $\phi$, and as we shall see ${{I}^{11}_{\mathcal{B}}}^{loc}$ belongs to the kernel of the extension of $\phi$ to $\mathbb{C}_{\mathcal{B}}^{s,loc}[G]$. Therefore one can compose this extension with the operator $P_c$, and by Proposition \ref{BasP} the image of this composition is invariant with respect to the natural extension of the $\mathbb{C}_{\mathcal{B}}^s[M_+]$--adjoint action to $Q_{\mathcal{B}}^{loc}$. The operator $\Pi_c$ is a classical counterpart of this composition and using the description of $\Pi_c$ given in Corollary \ref{proj1} we shall show that the image of the composition is a localization of the algebra $W_{\mathcal{B}}^s(G)$.

More precisely, formula (\ref{phidef}) and the surjectivity of the map $\phi$ imply that one can define a natural action of the algebra generated by the elements $q^{2P_{\h'^\perp}\lambda^\vee}\in \mathbb{C}_{\mathcal{B}}^s[G^*]$, $\lambda \in P_+$ on $Q_{\mathcal{B}}$ as follows
\begin{equation}\label{phidefh}
q^{2P_{\h'^\perp}\lambda^\vee}\phi(f)
=\varphi({\rm Ad}_s^0(q^{-2P_{\h'^\perp}\lambda^\vee})(f))q^{2P_{\h'^\perp}\lambda^\vee}1=
\end{equation}
$$
=c_\lambda^{-1}{\rm Ad}_s(q^{-2P_{\h'^\perp}\lambda^\vee})(\varphi(f))q^{2P_{\h'^\perp}\lambda^\vee}1,
$$
where the last identity follows from part (i) of Proposition \ref{locfin}.
Let $Q_{\mathcal{B}}^{loc}$ \index[not]{Q@$Q_{\mathcal{B}}^{loc}$} be the localization of $Q_{\mathcal{B}}$ by the multiplicative set of the elements $q^{2P_{\h'^\perp}\lambda^\vee}$, $\lambda \in P_+$.

Now consider the subalgebra  $\mathbb{C}_{\mathcal{B}}^{s,loc}[G_*]\subset \mathbb{C}_{\mathcal{B}}^s[G^*]$ \index[not]{C@$\mathbb{C}_\mathcal{B}^{s,loc}[G_*]$} generated by $\mathbb{C}_{\mathcal{B}}^s[G_*]$ and by the elements $q^{2P_{\h'^\perp}\lambda^\vee}$, $\lambda \in P$. Note that the ${\rm Ad}_s$--action of these elements normalizes $\mathbb{C}_{\mathcal{B}}^s[G_*]$ in $\mathbb{C}_{\mathcal{B}}^s[G^*]$ as $\mathbb{C}_{\mathcal{B}}^s[G_*]$ is the direct sum of its weight components.
Therefore $Q_{\mathcal{B}}^{loc}$ is the image of $\mathbb{C}_{\mathcal{B}}^{s,loc}[G_*]$ under the natural projection $\rho_{\chi_q^s}:\mathbb{C}_{\mathcal{B}}^s[G^*]\rightarrow \mathbb{C}_{\mathcal{B}}^s[G^*]/I_{\mathcal{B}}$, and hence the adjoint action of $\mathbb{C}_{\mathcal{B}}^s[M_+]$ on $Q_{\mathcal{B}}$ naturally extends to $Q_{\mathcal{B}}^{loc}$.

\begin{lemma}\label{phiext}
Let $\kappa=1$. Assume that for $i=1,\ldots l'$ $\bar{k}_i\in \mathcal{B}^*$, where $\bar{k}_i$ are defined in (\ref{charq}).
Then $c_\lambda \in \mathcal{B}^*$ for any $\lambda \in P_+$, and  
$\phi$ extends to a $\mathbb{C}_{\mathcal{B}}^s[M_+]$--module homomorphism $\phi: \mathbb{C}_{\mathcal{B}}^{s,loc}[G] \rightarrow Q_{\mathcal{B}}^{loc}$ which satisfies
\begin{equation}\label{philoc}
\phi (f\otimes {\Delta_\lambda^s}^{-1})=
c_\lambda^{-1}q^{2\left\langle P_{\h'}\lambda^\vee,\lambda^\vee\right\rangle-2P_{\h'^\perp}\lambda^\vee}\phi({\rm Ad}_s^0(q^{({1+s \over 1-s }s^{-1}P_{{\h'}}-s^{-1})\lambda^\vee})(f)), f\in \mathbb{C}_{\mathcal{B}}^s[G],
\end{equation}
and ${{I}^{11}_{\mathcal{B}}}^{loc}$ belongs to the kernel of this homomorphism, so $\phi$ naturally gives rise to a $\mathbb{C}_{\mathcal{B}}^s[M_+]$--module homomorphism
$$
\phi:\mathbb{C}_{11}^{loc}[G]\rightarrow Q_{\mathcal{B}}^{loc}.
$$
\end{lemma}

\begin{proof}
From formula (\ref{phidef}) it follows that $\phi$ extends to a $\mathbb{C}_{\mathcal{B}}^s[M_+]$--module homomorphism $\phi: \mathbb{C}_{\mathcal{B}}^{s,loc}[G] \rightarrow Q_{\mathcal{B}}^{loc}$ 
which is defined by (\ref{philoc}), and by Proposition \ref{kerphi} ${{I}^{11}_{\mathcal{B}}}^{loc}$ belongs to the kernel of this homomorphism.

\end{proof}

By (\ref{admu}) for $\lambda\in P_+$ the element ${\Delta_\lambda^s}$ is $\mathbb{C}_{\mathcal{B}}^s[M_+]$--invariant with respect to the ${\rm Ad}_s^0$--action on $\mathbb{C}_{11}^{loc}[G]$, and hence $\phi({\Delta_\lambda^s})=q^{2P_{\h'^\perp}\lambda^\vee}1$ is $\mathbb{C}_{\mathcal{B}}^s[M_+]$--invariant with respect to the ${\rm Ad}_s$--action of $\mathbb{C}_{\mathcal{B}}^s[M_+]$ on $Q_{\mathcal{B}}^{loc}$. Thus by (\ref{phidef}) and (\ref{adm4}) we have for $\beta\in \Delta_{\m_+}$
$$
\phi ({\rm Ad}_s^0(\tilde{f}_\beta)(f\otimes {\Delta_\lambda^s}))=c_\lambda {\rm Ad}_s(\tilde{f}_\beta)(\varphi({\rm Ad}_s^0(q^{-({1+s \over 1-s }P_{{\h'}}+id)\lambda^\vee})(f)))q^{2P_{\h'^\perp}\lambda^\vee}1,
$$
and
$$
\phi ({\rm Ad}_s^0(\tilde{f}_\beta)(f\otimes {\Delta_\lambda^s}^{-1}))=c_\lambda^{-1}q^{2(P_{\h'}\lambda^\vee,\lambda^\vee)} {\rm Ad}_s(\tilde{f}_\beta)(\varphi({\rm Ad}_s^0(q^{({1+s \over 1-s }P_{{\h'}}+id)\lambda^\vee})(f)))q^{-2P_{\h'^\perp}\lambda^\vee}1.
$$
The last formula completely determines the adjoint action of $\mathbb{C}_{\mathcal{B}}^s[M_+]$ on $Q_{\mathcal{B}}^{loc}$.

Now we can describe q-W--algebras in terms of the operator $P_c$.
\begin{theorem}\label{Piqmain}
Suppose that $\kappa=1$ and $\bar{k}_i\in \mathcal{B}^*$ for $i=1,\ldots l'$, where $\bar{k}_i$ are defined in (\ref{charq}). Then the composition $\phi P_c$ gives rise to a well--defined  operator 
\begin{equation}\label{Piqdef}
\Pi_c^q:=\phi P_c: \mathbb{C}_{R-1 n_{R-1}+1}^{loc}[G]\rightarrow W_{\mathcal{B}}^{s,loc}(G),	\index[not]{P@$\Pi_c^q$}
\end{equation}
where $W_{\mathcal{B}}^{s,loc}(G)={\rm Hom}_{\mathbb{C}_{\mathcal{B}}^s[M_+]}(\mathbb{C}_{\mathcal{B}},Q_{\mathcal{B}}^{loc})$. \index[not]{W@$W_{\mathcal{B}}^{s,loc}(G)$}

Moreover, we have an imbedding of $\mathbb{C}_{\mathcal{B}}^s[M_+]$--modules, $Q_{\mathcal{B}}^{loc}\subset \mathbb{C}_{\mathcal{B}}^s[G^*]/I_{\mathcal{B}}$, multiplication in $\mathbb{C}_{\mathcal{B}}^s[G^*]$ induces a multiplication on $W_{\mathcal{B}}^{s,loc}(G)$, and we have an imbedding of algebras $W_{\mathcal{B}}^s(G)\subset W_{\mathcal{B}}^{s,loc}(G)$. 

Let $W_q^{s,loc}(G)=W_{\mathcal{B}}^{s,loc}(G)\otimes_{\mathcal{B}}\mathbb{C}(q^{\frac{1}{d{\bar{r}}^2}})$, ${\rm Im}\Pi_c^q$ the image of $\Pi_c^q$, and denote ${\rm Im}_q\Pi_c^q={\rm Im}\Pi_c^q\otimes_{\mathcal{B}}\mathbb{C}(q^{\frac{1}{d{\bar{r}}^2}})$.

Then $W_q^{s,loc}(G)={\rm Im}_q\Pi_c^q$ and $W_{\mathcal{B}}^{s,loc}(G)=W_{\mathcal{B}}^{s,loc}(G)\cap {\rm Im}_q\Pi_c^q$.

\end{theorem}

\begin{proof}
By Lemma \ref{phiext} and by Proposition \ref{BasP} the composition $\Pi_c^q=\phi P_c: \mathbb{C}_{R-1 n_{R-1}+1}^{loc}[G]\rightarrow W_{\mathcal{B}}^{s,loc}(G)$ is well--defined. All other claims of this theorem, except for the isomorphisms $W_q^{s,loc}(G)={\rm Im}_q\Pi_c^q$ and $W_{\mathcal{B}}^{s,loc}(G)=W_{\mathcal{B}}^{s,loc}(G)\cap {\rm Im}_q\Pi_c^q$ follow from the definitions and are established similarly to Proposition \ref{9.2}.

Using Lemma \ref{incl} one sees that in order to establish these isomorphisms it suffices to verify that the specialization of $\Pi_c^q$ at $q^{\frac{1}{d{\bar{r}}^2}}=1$ is surjective. 

In order to do that we observe that by the definition and by Theorem \ref{var} 
the specialization of $Q_{\mathcal{B}}^{loc}$ at $q^{\frac{1}{d{\bar{r}}^2}}=1$ is isomorphic to the localization $\mathbb{C}^{loc}[N_-Z{s}M_-]$ \index[not]{C@$\mathbb{C}^{loc}[N_-Z{s}M_-]$} of the algebra $\mathbb{C}[N_-Z{s}M_-]$ by the multiplicative set of the classical counterparts of the elements ${\Delta_\lambda^s}$, which we denote by the same symbols, ${\Delta_\lambda^s}:=(v_\lambda,\cdot~ s^{-1}v_\lambda)$. \index[not]{D@${\Delta_\lambda^s}$} Denote by $\mathbb{C}^{loc}[(w_1\ldots w_{R-2})^{-1}N_-Z{s}M_-(w_1\ldots w_{R-2})]$ \index[not]{C@$\mathbb{C}^{loc}[(w_1\ldots w_{R-2})^{-1}N_-Z{s}M_-(w_1\ldots w_{R-2})]$} the image of $\mathbb{C}^{loc}[N_-Z{s}M_-]$ under the map $f\mapsto (w_1\ldots w_{R-2})(f)$. This map sends ${\Delta_\lambda^s}$ to a function which we denote $\Delta_{(w_1\ldots w_{R-2})^{-1}\lambda}^{\widehat{s}_{R-1}}$. \index[not]{D@$\Delta_{(w_1\ldots w_{R-2})^{-1}\lambda}^{\widehat{s}_{R-1}}$}

By (\ref{admu}) for $\lambda\in P_+$ ${\Delta_\lambda^s}$ regarded as elements of $\mathbb{C}_{11}^{loc}[G]$ are $\mathbb{C}_{\mathcal{B}}^s[M_+]$--invariant with respect to the ${\rm Ad}_s^0$--action on $\mathbb{C}_{11}^{loc}[G]$, and hence $\phi({\Delta_\lambda^s})=q^{2P_{\h'^\perp}\lambda^\vee}1$ is $\mathbb{C}_{\mathcal{B}}^s[M_+]$--invariant with respect to the ${\rm Ad}_s$--action on $Q_{\mathcal{B}}^{loc}$. Therefore their classical counterparts ${\Delta_\lambda^s}=(v_\lambda,\cdot~ s^{-1}v_\lambda)\in \mathbb{C}[N_-Z{s}M_-]$ are $M_-$--invariant, and hence $\mathbb{C}^{loc}[N_-Z{s}M_-]^{M_-}$ \index[not]{C@$\mathbb{C}^{loc}[N_-Z{s}M_-]^{M_-}$} is the localization of $\mathbb{C}[N_-Z{s}M_-]^{M_-}$ by the elements ${\Delta_\lambda^s}$, $\lambda \in P_+$. 

This result and explicit formulas (\ref{Pip}), (\ref{Pi1}), (\ref{tind}), (\ref{pic}) for the operator $\Pi_c$, formulas (\ref{Ppk}) with $k=0$ and the definition of the operator $P_c$ imply that the specialization of the operator $\Pi_c^q$ at $q^{\frac{1}{d{\bar{r}}^2}}=1$ gives rise to a natural extension of the operator $\Pi_c: \mathbb{C}[(w_1\ldots w_{R-2})N_-Z{s}M_-(w_1\ldots w_{R-2})^{-1}]\rightarrow \mathbb{C}[N_-Z{s}M_-]^{M_-}$ to an operator $\Pi_c^{loc}:\mathbb{C}^{loc}[(w_1\ldots w_{R-2})^{-1}N_-Z{s}M_-(w_1\ldots w_{R-2})]\rightarrow \mathbb{C}^{loc}[N_-Z{s}M_-]^{M_-}$ \index[not]{P@$\Pi_c^{loc}$} defined by
$$
\Pi_c^{loc}(f({\Delta_{(w_1\ldots w_{R-2})^{-1}\lambda}^{\widehat{s}_{R-1}}})^{-1})=\Pi_c(f){\Delta_\lambda^s}^{-1}, f\in \mathbb{C}[(w_1\ldots w_{R-2})N_-Z{s}M_-(w_1\ldots w_{R-2})^{-1}].
$$
Since the operator $\Pi_c$ is surjective, $\Pi_c^{loc}$ is also surjective. 

From the surjectivity of the specialization of the operator $\Pi_c^q$ at $q^{\frac{1}{d{\bar{r}}^2}}=1$ it follows that $W_{\mathcal{B}}^{s,loc}(G)={\rm Im}\Pi_c^q$  mod $(q^{\frac{1}{d{\bar{r}}^2}}-1)W_{\mathcal{B}}^{s,loc}(G)$.
Note also that ${\rm Im}\Pi_c^q\subset W_{\mathcal{B}}^{s,loc}(G)$ are submodules of the $\mathcal{B}$--module $Q_{\mathcal{B}}^{loc}$, and $Q_{\mathcal{B}}^{loc}$ is a $\mathcal{B}$--submodule of the $\mathcal{B}$--module $\mathbb{C}_{\mathcal{B}}^s[G^*]/I_{\mathcal{B}}$ which is free over $\mathcal{B}$ by Proposition \ref{Afree1}. 
Since $\mathcal{B}$ is a principal ideal domain $Q_{\mathcal{B}}^{loc}$ is $\mathcal{B}$--free by Theorem 6.5 in \cite{R}.

The properties mentioned in the previous paragraph and Lemma \ref{incl} imply that $W_q^{s,loc}(G)={\rm Im}_q\Pi_c^q$ and $W_{\mathcal{B}}^{s,loc}(G)=W_{\mathcal{B}}^{s,loc}(G)\cap {\rm Im}_q\Pi_c^q$. This completes the proof.

\end{proof}


\section{Bibliographic comments}

\pagestyle{myheadings}
\markboth{CHAPTER \thechapter.~~ZHELOBENKO TYPE OPERATORS FOR Q-W--ALGEBRAS}{\thesection.~BIBLIOGRAPHIC COMMENTS}

\setcounter{equation}{0}
\setcounter{theorem}{0}

The results presented in this chapter are entirely new. 

Commutation relations in the algebra $\mathbb{C}_{\mathcal{B}}^s[G]$ which appear in Section \ref{CG} can be found, for instance, in \cite{BG}, Theorem I.8.16.

The definition of the Zhelobenko type operators for q-W--algebras was inspired by the construction of extremal projection operators and of the Zhelobenko operators due to Zhelobenko. The definitions and the statements in this chapter are conceptually close to the definition and the properties of the Zhelobenko operators introduced and studied in \cite{Z1}--\cite{Z2} (see also \cite{KO}). Below, for the convenience of the reader who is familiar with these papers, we give references to similar statements from them. However, the results of \cite{Z1}--\cite{Z2} and \cite{KO} are not used in this book and not directly related to it.

For $k=0$ the operators $P_{jp}^k$ are counterparts of the Zhelobenko operators $q_\alpha$ introduced in \cite{Z1}, \S2 and \S5, in \cite{Z2}, Definition 5.2.1, and for $k>0$ the operators $P_{jp}^k$ are counterparts of the operators $q^{(k)}_{\alpha,\m}$ defined in \cite{KO}, formula (4.9).

Properties of the Zhelobenko operators similar to those of the operators $P_{jp}$ mentioned in Proposition \ref{BasP}  can be found in \cite{Z1}, \S5, Proposition 1 (iii) and (iv), \cite{Z2}, Proposition 5.2.4 (b) and (c), \cite{KO}, Lemma 4.5 (ii), (iv) and (v).

Zhelobenko operators $q_w$ conceptually analogous to $P_{\leq jp}$ were defined in \cite{Z1}, \S5, Definition 1, \cite{Z2}, Definition 5.2.4.

Properties (\ref{mainprop1}) and (\ref{mainprop2}) are counterparts of Proposition 3, parts (i) and (iii) in \cite{Z1}, \S2, and of properties ($\alpha$) and ($\beta$) in Section  5.2.4 in \cite{Z2}.

Lemma \ref{Jvan} is analogous to a similar property for the Zhelobenko operators stated in \cite{Z2}, Proposition and Corollary 5.2.3, and in \cite{KO}, Lemmas 4.3 and 4.5 (iii).

Formulas similar to (\ref{Pdef}) are used in the proof of Proposition 5.2.4 in \cite{Z2} and in Proposition 4.4 in \cite{KO}.

Theorem \ref{Piqmain} is analogous to Theorem 2 in
\cite{Z1}, \S6 and to Theorem 5.5.1 in \cite{Z2} for the Zhelobenko operators.


\chapter{Application of q-W--algebras to the description of the category of equivariant modules over a quantum group}\label{part5}

\pagestyle{myheadings}
\markboth{CHAPTER~\thechapter.~APPLICATION TO EQUIVARIANT MODULES OVER A QUANTUM GROUP}{\thesection.~A CATEGORY OF EQUIVARIANT MODULES OVER A QUANTUM GROUP}

In this rather short chapter we apply Corollary \ref{Bpbas} to establish an equivalence between the category of finitely generated representations of a q-W--algebra and a category ${\mathbb{C}}_{\varepsilon}^{s,loc}[G_*]-{\rm mod}_{\mathbb{C}_{\varepsilon}^s[M_+]}^{\chi_\varepsilon^s}$ of equivariant modules over a specialization of the algebra $\mathbb{C}_{\mathcal{B}}^{s,loc}[G_*]$. Categories of this kind were denoted $A-{\rm  mod}_B^\chi$ in the introduction.
The structure of modules from the category ${\mathbb{C}}_{\varepsilon}^{s,loc}[G_*]-{\rm mod}_{\mathbb{C}_{\varepsilon}^s[M_+]}^{\chi_\varepsilon^s}$ is similar to that of  $\g-K$--modules or of principal series representations over complex semisimple Lie algebras. 

The proof of the main theorem of this chapter, Theorem \ref{sqeq}, is based on Corollary \ref{Bpbas}. 
In this framework one can give precise values of $\varepsilon$ for which the categorical equivalence holds. Remarkably, with slight modifications this method is also applicable to the study of the structure of finite--dimensional representations over quantum groups at roots of unity. This will be done in the next chapter.


\section{A category of equivariant modules over a quantum group}\label{cateq}

\setcounter{equation}{0}
\setcounter{theorem}{0}

In this section we define a category of equivariant representations over a quantum group.  

Suppose that $\kappa=1$ and let $\varepsilon \in \mathbb{C}^*$ be such that $[n]_{\varepsilon_{\alpha_i}}\neq 0$ and $\varepsilon^{2d_i}\neq 1$ for $i=1,\ldots ,l$, $n\in \mathbb{N}\setminus\{0\}$. Assume that for $i=1,\ldots l'$ $\bar{k}_i\in \mathcal{B}^*$, where $\bar{k}_i$ are defined in (\ref{charq}). Fix a root $\varepsilon^{\frac{1}{d{\bar{r}}^2}}$ of $\varepsilon$ of order $\frac{1}{d{\bar{r}}^2}$. Let $U_{\varepsilon}^s(\g)$, \index[not]{U@$U_\varepsilon^s(\g)$} $U_{\varepsilon}^s(\m_-)$, \index[not]{U@$U_\varepsilon^s(\m_-)$} ${\mathbb{C}}_{\varepsilon}^{s,loc}[G_*]$, \index[not]{C@$\mathbb{C}_\varepsilon^{s,loc}[G_*]$} $\mathbb{C}_{\varepsilon}^s[M_+]$, \index[not]{C@$\mathbb{C}_\varepsilon^s[M_+]$} ${\mathbb{C}}_{\varepsilon}^s[B_+]$, \index[not]{C@$\mathbb{C}_\varepsilon^s[B_+]$} ${\mathbb{C}}_{\varepsilon}^s[G^*]$, \index[not]{C@$\mathbb{C}_\varepsilon^s[G^*]$} $Q_{\varepsilon}^{loc}$, \index[not]{Q@$Q_{\varepsilon}^{loc}$} $W_\varepsilon^{s,loc}(G)$, \index[not]{W@$W_\varepsilon^{s,loc}(G)$} $\chi_\varepsilon^s$, \index[not]{x@$\chi_\varepsilon^s$} $\mathbb{C}_{\varepsilon_s}$, \index[not]{C@$\mathbb{C}_{\varepsilon_s}$} $I_\varepsilon$, \index[not]{I@$I_\varepsilon$} $\phi_\varepsilon$, \index[not]{f@$\phi_\varepsilon$} $B_{m_1\ldots m_c}^\varepsilon$, \index[not]{B@$B_{m_1\ldots m_c}^\varepsilon$} ${\mathbb{C}}_{11}^{loc}[G]_\varepsilon$ \index[not]{C@$\mathbb{C}_{11}^{loc}[G]_\varepsilon$} be the natural specializations at $q^{\frac{1}{d{\bar{r}}^2}}=\varepsilon^{\frac{1}{d{\bar{r}}^2}}$ of $U_{\mathcal{A}}^s(\g)$, $U_{\mathcal{A}}^s(\m_-):=U_h^s(\m_-)\cap U_{\mathcal{A}}^s(\g)$, ${\mathbb{C}}_{\mathcal{B}}^{s,loc}[G_*]$, $\mathbb{C}_{\mathcal{B}}^s[M_+]$, $\mathbb{C}_{\mathcal{B}}^s[B_+]$, $\mathbb{C}_{\mathcal{B}}^s[G^*]$, $Q_{\mathcal{B}}^{loc}$, $W_{\mathcal{B}}^{s,loc}(G)$, $\chi_q^s$, $\mathcal{B}_{\varepsilon_s}$, $I_{\mathcal{B}}$, $\phi$, $B_{m_1\ldots m_c}$, ${\mathbb{C}}_{11}^{loc}[G]$, respectively. 

Note that under our assumptions on $\varepsilon$ one has $\mathbb{C}_{\varepsilon}^s[M_+]=U_{\varepsilon}^s(\m_-)=U_{\varepsilon}^{s,res}(\m_-)$, $U_{\varepsilon}^s(\g)=U_{\varepsilon}^{s,res}(\g)$, and ${\mathbb{C}}_{\varepsilon}^s[G^*]$ is a subalgebra in $U_{\varepsilon}^{s,res}(\g)$ as $\mathbb{C}_{\mathcal{B}}^s[G^*]$ is a subalgebra in $U_{\mathcal{B}}^{s,res}(\g)$.

Let $J={\rm Ker}~\varepsilon_s|_{\mathbb{C}_{\varepsilon}^s[M_+]}$ \index[not]{J@$J$} be the augmentation ideal \index{ideal!augmentation} of $\mathbb{C}_{\varepsilon}^s[M_+]$ related to the counit $\varepsilon_s$, and
$\mathbb{C}_{\varepsilon_s}$ \index[not]{C@$\mathbb{C}_{\varepsilon_s}$} the trivial representation of $\mathbb{C}_{\varepsilon}^s[M_+]$ given by the counit.
Let $V$ be a finitely generated ${\mathbb{C}}_{\varepsilon}^{s,loc}[G_*]$--module which satisfies the following conditions:
\begin{enumerate}
\item
$V$ is a right $\mathbb{C}_{\varepsilon}^s[M_+]$--module with respect to an action ${\rm Ad}_s$ \index[not]{A@${\rm Ad}_s$} such that the action of the augmentation ideal $J$ on $V$ is locally nilpotent. We shall also call this action {\it the adjoint action}. \index{action!adjoint!on a ${\mathbb{C}}_{\varepsilon}^{s,loc}[G_*]$--module} 

\item
The following compatibility condition holds for the two actions
\begin{equation}\label{compat}
{\rm Ad}_sx(yv)={\rm Ad}_sx^1(y){\rm Ad}_sx^2(v),~x\in \mathbb{C}_{\varepsilon}^s[M_+],~y\in {\mathbb{C}}_{\varepsilon}^{s,loc}[G_*],~v\in V,
\end{equation}
where $\Delta_s(x)=x^1\otimes x^2$, ${\rm Ad}_sx^1(y)$ is the adjoint action of $x^1\in {\mathbb{C}}_{\varepsilon}^s[B_+]$ on $y\in {\mathbb{C}}_{\varepsilon}^{s,loc}[G_*]$. 

An element $v\in V$ is called a {\it Whittaker vector} \index{Whittaker vector!for quantum groups at a generic deformation parameter} if
${\rm Ad}_sxv=\varepsilon_s(x)v$ for any $x\in \mathbb{C}_{\varepsilon}^s[M_+]$.
The space
\begin{equation}\label{winv}
{\rm Wh}(V):={\rm Hom}_{\mathbb{C}_{\varepsilon}^s[M_+]}(\mathbb{C}_{\varepsilon_s},V)\subset V \index[not]{W@${\rm Wh}(V)$}
\end{equation}
is called {\it the space of Whittaker vectors} of $V$. \index{space!of Whittaker vectors!for quantum groups at a generic deformation parameter}

Consider the induced ${\mathbb{C}}_{\varepsilon}^s[G^*]$--module $V'={\mathbb{C}}_{\varepsilon}^s[G^*]\otimes_{{\mathbb{C}}_{\varepsilon}^{s,loc}[G_*]}V$. Using the adjoint action of ${\mathbb{C}}_{\varepsilon}^s[G^*]$ on itself one can naturally extend the adjoint action of $\mathbb{C}_{\varepsilon}^s[M_+]$ from $V$ to $V'$ in such a way that compatibility condition $(\ref{compat})$ is satisfied for the natural action of ${\mathbb{C}}_{\varepsilon}^s[G^*]$ and the adjoint action ${\rm Ad}_s$ of $\mathbb{C}_{\varepsilon}^s[M_+]$ on $V'$. 

As we observed in Section \ref{qWalgpois} $\Delta_s(\mathbb{C}_{\varepsilon}^s[M_+])\subset {\mathbb{C}}_{\varepsilon}^s[B_+]\otimes \mathbb{C}_{\varepsilon}^s[M_+]$. 
We shall require that
\item
For any $x\in \mathbb{C}_{\varepsilon}^s[M_+]$ the natural action of the element $(S_s\otimes\chi_\varepsilon^s)\Delta_s(x)\in {\mathbb{C}}_{\varepsilon}^s[G^*]$ on $V'$ coincides with the adjoint action ${\rm Ad}_sx$ of $x$ on $V'$.

As in the second part of the proof of Proposition \ref{9.2} one can check that the last condition implies that
for any $z\in {\mathbb{C}}_{\varepsilon}^{s,loc}[G_*]\cap I_\varepsilon$ and $v\in {\rm Wh}(V)$ one has $zv=0$.
\end{enumerate}

Generalizing the definition of Whittaker vectors given in Condition 2 above, for any $\mathbb{C}_{\varepsilon}^s[M_+]$--module $X$ we denote by ${\rm Wh}(X)$ \index[not]{W@${\rm Wh}(X)$} the subspace of $X$ which consists of the elements $v$ such that $xv=\varepsilon_s(x)v$, $x\in \mathbb{C}_{\varepsilon}^s[M_+]$. We shall also call the space ${\rm Wh}(X)$ {\it the space of Whittaker vectors}, and the elements of this space will be called {\it Whittaker vectors}.

Denote by ${\mathbb{C}}_{\varepsilon}^{s,loc}[G_*]-{\rm mod}_{\mathbb{C}_{\varepsilon}^s[M_+]}^{\chi_\varepsilon^s}$ \index[not]{C@$\mathbb{C}_\varepsilon^{s,loc}[G_*]-{\rm mod}_{\mathbb{C}_\varepsilon^s[M_+]}^{\chi_\varepsilon^s}$} the category of finitely generated ${\mathbb{C}}_{\varepsilon}^{s,loc}[G_*]$--modules which satisfy Conditions 1--3 above. Morphisms in the category ${\mathbb{C}}_{\varepsilon}^{s,loc}[G_*]-{\rm mod}_{\mathbb{C}_{\varepsilon}^s[M_+]}^{\chi_\varepsilon^s}$ are ${\mathbb{C}}_{\varepsilon}^{s,loc}[G_*]$- and $\mathbb{C}_{\varepsilon}^s[M_+]$--module homomorphisms.
We call ${\mathbb{C}}_{\varepsilon}^{s,loc}[G_*]-{\rm mod}_{\mathbb{C}_{\varepsilon}^s[M_+]}^{\chi_\varepsilon^s}$ {\it the category of $(\mathbb{C}_{\varepsilon}^s[M_+],\chi_\varepsilon^s)$--equivariant modules over ${\mathbb{C}}_{\varepsilon}^{s,loc}[G_*]$}. \index{category of equivariant modules over ${\mathbb{C}}_{\varepsilon}^{s,loc}[G_*]$}

\begin{lemma}
The algebra $W_{\mathcal{B}}^{s,loc}(G)$ naturally acts on the space of Whittaker vectors of any object $V$ of the category ${\mathbb{C}}_{\varepsilon}^{s,loc}[G_*]-{\rm mod}_{\mathbb{C}_{\varepsilon}^s[M_+]}^{\chi_\varepsilon^s}$. 
\end{lemma}

\begin{proof}
Indeed, if $w,w'\in {\mathbb{C}}_{\varepsilon}^{s,loc}[G_*]$ are two representatives of an element from $W_{\mathcal{B}}^{s,loc}(G)$ then $w-w'\in {\mathbb{C}}_{\varepsilon}^{s,loc}[G_*]\cap I_\varepsilon$ by the definition of $Q_\varepsilon^{loc}$, and hence by Condition 3 above for any $v\in {\rm Wh}(V)$ one has $wv=w'v$. 

Moreover, by the definition of the algebra $W_{\mathcal{B}}^{s,loc}(G)$ and by condition (\ref{compat}) we have
$$
{\rm Ad}_sx(wv)={\rm Ad}_sx^1(w){\rm Ad}_sx^2(v)={\rm Ad}_sx^1(w)\varepsilon_s(x^2)v={\rm Ad}_sx(w)v=\varepsilon_s(x)wv.
$$
Therefore $wv$ is a Whittaker vector independent of the choice of the representative $w$.

\end{proof}

For any $\mathbb{C}$--module $R$ we denote by ${\rm hom}_{\mathbb{C}}(\mathbb{C}_{\varepsilon}^s[M_+]), R)$ \index[not]{h@${\rm hom}_{\mathbb{C}}(\mathbb{C}_{\varepsilon}^s[M_+]), R)$} the subspace in ${\rm Hom}_{\mathbb{C}}(\mathbb{C}_{\varepsilon}^s[M_+],R)$ which consists of the linear maps vanishing on some power of the augmentation ideal $J={\rm Ker }~\varepsilon_s$ of $\mathbb{C}_{\varepsilon}^s[M_+]$, ${\rm hom}_{\mathbb{C}}(\mathbb{C}_{\varepsilon}^s[M_+],R)=\{f\in {\rm Hom}_{\mathbb{C}}(\mathbb{C}_{\varepsilon}^s[M_+],R):f(J^n)=0~{\rm for~some}~n>0\}$. Note that for every element $f$ of ${\rm hom}_{\mathbb{C}}(\mathbb{C}_{\varepsilon}^s[M_+],R)$ one has $f(x)=0$ if $x$ does not belong to a finite--dimensional subspace of $\mathbb{C}_{\varepsilon}^s[M_+]$, and hence 
\begin{equation}\label{rt}
{\rm hom}_{\mathbb{C}}(\mathbb{C}_{\varepsilon}^s[M_+],R)\simeq {\rm hom}_{\mathbb{C}}(\mathbb{C}_{\varepsilon}^s[M_+],\mathbb{C})\otimes R.
\end{equation} 

Equip the space ${\rm hom}_{\mathbb{C}}(\mathbb{C}_{\varepsilon}^s[M_+],R)$ with the right action of $\mathbb{C}_{\varepsilon}^s[M_+]$ induced by the multiplication in $\mathbb{C}_{\varepsilon}^s[M_+]$ from the left. To study the properties of this module, which will be required to establish an equivalence between the category of finitely generated representations of $W_{\mathcal{B}}^{s,loc}(G)$ and the category ${\mathbb{C}}_{\varepsilon}^{s,loc}[G_*]-{\rm mod}_{\mathbb{C}_{\varepsilon}^s[M_+]}^{\chi_\varepsilon^s}$, we shall need a special filtration on the algebra $\mathbb{C}_{\varepsilon}^s[M_+]=U_{\varepsilon}^s(\m_-)$.

To define this filtration we recall that the algebra $U_\varepsilon^s(\g)$ can be equipped with the De Concini--Kac filtration \index{De Concini--Kac filtration} such that the associated graded algebra \index{algebra!associated graded} is almost commutative. Below we remind the definition of this filtration.

Assume, as before, that the system of positive roots $\Delta_+$ associated to $s$ is equipped with a normal ordering $\beta_1,\ldots,\beta_D$ associated to $s$ in Definition \ref{circorddef}. Let $e_\beta,f_\beta\in U_\varepsilon^s(\g)$ be the corresponding quantum root vectors, $e^{\bf r}=e_{\beta_1}^{r_1}\ldots e_{\beta_D}^{r_D}$, $f^{\bf m}=f_{\beta_D}^{m_D}\ldots f_{\beta_1}^{m_1}$ for ${\bf r},~{\bf m}\in {\Bbb N}^D$ the corresponding elements of the bases of $U_\varepsilon^s(\n_+)$ and $U_\varepsilon^s(\n_-)$, respectively (see Remark \ref{segmPBWsrev}). Define {\it the height} of the element \index{height!of an element of a quantum group} $u_{{\bf r},{\bf m},t}:=e^{\bf r}tf^{\bf m}$, \index[not]{u@$u_{{\bf r},{\bf m},t}$} $t\in U_\varepsilon^s({\frak h})$ as follows ${\rm ht}(u_{{\bf r},{\bf m},t})=\sum_{i=1}^D(m_i+r_i){\rm ht}~\beta_i\in \mathbb{N}$, \index[not]{h@${\rm ht}~(u_{{\bf r},{\bf m},t})$} where ${\rm ht}~\beta_i$ \index[not]{h@${\rm ht}~\beta$} is the height \index{height!of a root} of the root $\beta_i$. 

Introduce also the {\it degree} \index{degree!of an element of a quantum group} of $u_{{\bf r},{\bf m},t}$ by
 $$
 d(u_{{\bf r},{\bf m},t})=(r_1,\ldots,r_D,m_D,\ldots,m_1,{\rm ht}(u_{{\bf r},{\bf m},t}))\in \mathbb{N}^{2D+1}. \index[not]{d@$d(u_{{\bf r},{\bf m},t})$}
 $$
By Lemma \ref{segmPBWs} the elements $e^{\bf r}tf^{\bf m}$ span $U_\varepsilon^s(\g)$ as a vector space. 

 Equip $\mathbb{N}^{2D+1}$ with the total lexicographic order \index{order, lexicographic} and denote by $(U_\varepsilon^s(\g))_{\bf k}$ \index[not]{U@$(U_\varepsilon^s(\g))_{\bf k}$} the span of the elements $u_{{\bf r},{\bf m},t}$ with $d(u_{{\bf r},{\bf m},t})\leq {\bf k}\in \mathbb{N}^{2D+1}$ in $U_\varepsilon^s(\g)$. Then Proposition 1.7 in \cite{DK} implies that $(U_\varepsilon^s(\g))_{\bf k}$, ${\bf k}\in \mathbb{N}^{2D+1}$ is a filtration of $U_\varepsilon^s(\g)$ such that the associated graded algebra is the complex associative algebra with generators $e_\alpha,~f_\alpha$, $\alpha\in \Delta_+$, $t_i^{\pm 1}$, $i=1,\ldots l$ subject to the relations
 \begin{equation}\label{scommDK}
 \begin{array}{l}
 t_it_j=t_jt_i,~~t_it_i^{-1}=t_i^{-1}t_i=1,~~ t_ie_\alpha t_i^{-1}=\varepsilon^{\frac{Y_i(\alpha)}{d}}e_\alpha, ~~t_if_\alpha t_i^{-1}=\varepsilon^{-\frac{Y_i(\alpha)}{d}}f_j,\\
\\
e_\alpha f_\beta =\varepsilon^{\left\langle {1+s \over 1-s}P_{{\h'}^*}\alpha,\beta\right\rangle} f_\beta e_\alpha,\\
 \\
 e_{\alpha}e_{\beta} = \varepsilon^{\left\langle \alpha,\beta\right\rangle+ \left\langle {1+s \over 1-s}P_{{\h'}^*}\alpha,\beta\right\rangle}e_{\beta}e_{\alpha},~ \alpha<\beta, \\
 \\
 f_{\alpha}f_{\beta} = \varepsilon^{\left\langle \alpha,\beta\right\rangle+\left\langle {1+s \over 1-s}P_{{\h'}^*}\alpha,\beta\right\rangle}f_{\beta}f_{\alpha},~ \alpha<\beta.
 \end{array}
 \end{equation}
Such algebras are called {\it semi--commutative}. \index{algebra!semi--commutative}

The main property of the modules ${\rm hom}_{\mathbb{C}}(\mathbb{C}_{\varepsilon}^s[M_+],R)$ is stated in the following lemma.
\begin{lemma}\label{inj1}
Let $J={\rm Ker }~\varepsilon_s$ be the augmentation ideal of $\mathbb{C}_{\varepsilon}^s[M_+]$, $R$ a $\mathbb{C}$--module, ${\rm hom}_{\mathbb{C}}(\mathbb{C}_{\varepsilon}^s[M_+],R)=\{f\in {\rm Hom}_{\mathbb{C}}(\mathbb{C}_{\varepsilon}^s[M_+],R):f(J^n)=0~{\rm for~some}~n>0\}$. Equip ${\rm hom}_{\mathbb{C}}(\mathbb{C}_{\varepsilon}^s[M_+],R)$ with the right action of $\mathbb{C}_{\varepsilon}^s[M_+]$ induced by multiplication on $\mathbb{C}_{\varepsilon}^s[M_+]$ from the left. Then the $\mathbb{C}_{\varepsilon}^s[M_+]$--module ${\rm hom}_{\mathbb{C}}(\mathbb{C}_{\varepsilon}^s[M_+],R)$ is injective.
\end{lemma}

\begin{proof}
First observe that the algebra $\mathbb{C}_{\varepsilon}^s[M_+]= U_{\varepsilon}^s(\m_-)$ is Notherian \index{algebra!Notherian} and the ideal $J$ satisfies the so-called weak Artin--Rees property, \index{Artin--Rees property, weak} i.e. for every finitely generated left $\mathbb{C}_{\varepsilon}^s[M_+]$--module $N$ and any its submodule $N'$ there exists an integer $n>0$ such that $J^nN\cap N'\subset JN'$.

Indeed, observe that the algebra $\mathbb{C}_{\varepsilon}^s[M_+]=U_{\varepsilon}^s(\m_-)$ can be equipped with a filtration induced by the De Concini--Kac filtration on the algebra $U_\varepsilon^s(\g)$ in such a way that the associated graded algebra is finitely generated and semi--commutative (see (\ref{scommDK})).
The fact that $\mathbb{C}_{\varepsilon}^s[M_+]$ is Notherian follows from the existence of the filtration on it for which the associated graded algebra is  semi--commutative and from Theorem 4 in  Ch. 5, \S 3 in \cite{Jac} (compare also with Theorem 4.8 in \cite{Mc}). 

The ideal $J$ satisfies the weak Artin--Rees property because the subring $\mathbb{C}_{\varepsilon}^s[M_+] +Jp+J^2p^2+\ldots \subset \mathbb{C}_{\varepsilon}^s[M_+][p]$, where $p$ is a central indeterminate, is Notherian (see \cite{Pas}, Ch. 11, \S 2, Lemma 2.1). The last fact follows from the existence of the filtration on $\mathbb{C}_{\varepsilon}^s[M_+] +Jp+J^2p^2+\ldots$ induced by the filtration on $\mathbb{C}_{\varepsilon}^s[M_+]$ for which the associated graded algebra is semi--commutative and again from Theorem 4 in  Ch. 5, \S 3 in \cite{Jac}.

Finally, the module ${\rm Hom}_{\mathbb{C}}(\mathbb{C}_{\varepsilon}^s[M_+],R)$ is obviously injective. By Lemma 3.2 in Ch. 3, \cite{Har} the module ${\rm hom}_{\mathbb{C}}(\mathbb{C}_{\varepsilon}^s[M_+],R)=\{f\in {\rm Hom}_{\mathbb{C}}(\mathbb{C}_{\varepsilon}^s[M_+],R):f(J^n)=0~{\rm for~some}~n>0\}$ is also injective since the ideal $J$ satisfies the weak Artin--Rees property.

\end{proof}


\section{Skryabin equivalence for equivariant modules over a quantum group}\label{Skr}

\pagestyle{myheadings}
\markboth{CHAPTER~\thechapter.~APPLICATION TO EQUIVARIANT MODULES OVER A QUANTUM GROUP}{\thesection.~SKRYABIN EQUIVALENCE FOR EQUIVARIANT MODULES OVER A QUANTUM GROUP}

\setcounter{equation}{0}
\setcounter{theorem}{0}

Now we can formulate the main theorem on the structure of the category ${\mathbb{C}}_{\varepsilon}^{s,loc}[G_*]-{\rm mod}_{\mathbb{C}_{\varepsilon}^s[M_+]}^{\chi^{s}_\varepsilon}$ and on the properties of its objects.

\begin{theorem}\label{sqeq}
Suppose that $\kappa=1$ and let $\varepsilon \in \mathbb{C}^*$. Assume that for $i=1,\ldots l'$ $\bar{k}_i\in \mathcal{B}^*$, where $\bar{k}_i$ are defined in (\ref{charq}). Fix a root $\varepsilon^{\frac{1}{d{\bar{r}}^2}}$ of $\varepsilon$ of order $\frac{1}{d{\bar{r}}^2}$. Assume also that $[n]_{\varepsilon_{\alpha_i}}\neq 0$ and $\varepsilon^{2d_i}\neq 1$ for $i=1,\ldots ,l$, $n\in \mathbb{N}$, $n>0$. Then the following statements are true.

(i) The functor $E\mapsto Q_\varepsilon^{loc}\otimes_{W_\varepsilon^{s,loc}(G)}E$, is an equivalence of the category of finitely generated left $W_\varepsilon^{s,loc}(G)$--modules and of the category ${\mathbb{C}}_{\varepsilon}^{s,loc}[G_*]-{\rm mod}_{\mathbb{C}_{\varepsilon}^s[M_+]}^{\chi^{s}_\varepsilon}$. The inverse equivalence is given by the functor $V\mapsto {\rm Wh}(V)$. In particular, the latter functor is exact.

(ii) Every object $V$ of the category ${\mathbb{C}}_{\varepsilon}^{s,loc}[G_*]-{\rm mod}_{\mathbb{C}_{\varepsilon}^s[M_+]}^{\chi^{s}_\varepsilon}$ is isomorphic to ${\rm hom}_{\mathbb{C}}(\mathbb{C}_{\varepsilon}^s[M_+],\mathbb{C})\otimes {\rm Wh}(V)$ as a right $\mathbb{C}_{\varepsilon}^s[M_+]$--module, where ${\rm hom}_{\mathbb{C}}(\mathbb{C}_{\varepsilon}^s[M_+],\mathbb{C})$ is equipped with the right action of $\mathbb{C}_{\varepsilon}^s[M_+]$ induced by the multiplication in $\mathbb{C}_{\varepsilon}^s[M_+]$ from the left. In particular, $V$ is $\mathbb{C}_{\varepsilon}^s[M_+]$--injective, and ${\rm Ext}^\bullet_{\mathbb{C}_{\varepsilon}^s[M_+]}(\mathbb{C}_{\varepsilon_s},V)={\rm Wh}(V)$.

(iii) $Q_{\varepsilon}^{loc}$ is an object of ${\mathbb{C}}_{\varepsilon}^{s,loc}[G_*]-{\rm mod}_{\mathbb{C}_{\varepsilon}^s[M_+]}^{\chi_\varepsilon^s}$. Moreover, $Q_\varepsilon^{loc}$ is isomorphic to ${\rm hom}_{\mathbb{C}}(\mathbb{C}_{\varepsilon}^s[M_+],\mathbb{C})\otimes {W_\varepsilon^{s,loc}(G)}$ as a $\mathbb{C}_{\varepsilon}^s[M_+]$--$W_\varepsilon^{s,loc}(G)$--bimodule, where the right $W_\varepsilon^{s,loc}(G)$--action is induced by the multiplication in $W_\varepsilon^{s,loc}(G)$ from the right. 
\end{theorem}

\begin{proof}
(ii) Let $V$ be an object of the category ${\mathbb{C}}_{\varepsilon}^{s,loc}[G_*]-{\rm mod}_{\mathbb{C}_{\varepsilon}^s[M_+]}^{\chi_\varepsilon^s}$. 
Fix any linear map $\overline{\rho}: V\rightarrow {\rm Wh}(V)$ the restriction of which to ${\rm Wh}(V)$ is the identity map, and let for any $v\in V$ $\sigma_\varepsilon(v):\mathbb{C}_{\varepsilon}^s[M_+]\rightarrow {\rm Wh}(V)$ be the $\mathbb{C}$--linear homomorphism given by 
\begin{equation}\label{se}
\sigma_\varepsilon(v)(x)=\overline{\rho}({\rm Ad}_sx(v)), v\in V, x\in \mathbb{C}_{\varepsilon}^s[M_+]. 
\end{equation}

Since the adjoint action of $J$ on $V$ is locally nilpotent $\sigma_\varepsilon(v)\in {\rm hom}_{\mathbb{C}}(\mathbb{C}_{\varepsilon}^s[M_+], {\rm Wh}(V))$, and we have a map $\sigma_\varepsilon:V\rightarrow {\rm hom}_{\mathbb{C}}(\mathbb{C}_{\varepsilon}^s[M_+], {\rm Wh}(V))\simeq {\rm hom}_{\mathbb{C}}(\mathbb{C}_{\varepsilon}^s[M_+],\mathbb{C})\otimes {\rm Wh}(V)$, where the last isomorphism is due to (\ref{rt}).

By the definition $\sigma_\varepsilon$ is a homomorphism of right $\mathbb{C}_{\varepsilon}^s[M_+]$--modules, where the right action of $\mathbb{C}_{\varepsilon}^s[M_+]$ on $${\rm hom}_{\mathbb{C}}(\mathbb{C}_{\varepsilon}^s[M_+], {\rm Wh}(V))$$ is induced by the multiplication in $\mathbb{C}_{\varepsilon}^s[M_+]$ from the left. 
The first claim in part (ii) will follow from the following lemma
\begin{lemma}\label{selem}
For any object $V$ in the category ${\mathbb{C}}_{\varepsilon}^{s,loc}[G_*]-{\rm mod}_{\mathbb{C}_{\varepsilon}^s[M_+]}^{\chi_\varepsilon^s}$ and for any linear map $\overline{\rho}: V\rightarrow {\rm Wh}(V)$ the restriction of which to ${\rm Wh}(V)$ is the identity map, the map $\sigma_\varepsilon: V\rightarrow {\rm hom}_{\mathbb{C}}(\mathbb{C}_{\varepsilon}^s[M_+], {\rm Wh}(V))\simeq {\rm hom}_{\mathbb{C}}(\mathbb{C}_{\varepsilon}^s[M_+],\mathbb{C})\otimes {\rm Wh}(V)$ defined by (\ref{se}) is an isomorphism of right $\mathbb{C}_{\varepsilon}^s[M_+]$--modules. 
\end{lemma}

\begin{proof}
First we prove that $\sigma_\varepsilon$ is injective. The proof will be based on the following lemma that will be also used later.
\begin{lemma}\label{inj}
Let $g:X\rightarrow Y$ be a homomorphism of $\mathbb{C}_{\varepsilon}^s[M_+]$--modules.  Assume that the action of the augmentation ideal of $\mathbb{C}_{\varepsilon}^s[M_+]$ on $X$ is locally nilpotent and that the restriction of $g$ to ${\rm Wh}(X)$ is injective. Then $g$ is injective.
\end{lemma}

\begin{proof}
Let $Z \subset X$ be the kernel of $g$. Assume that $Z$ is not trivial. Observe that $Z$ is invariant with respect to the action induced by the action of $\mathbb{C}_{\varepsilon}^s[M_+]$ on $X$, and that the augmentation ideal of $\mathbb{C}_{\varepsilon}^s[M_+]$ acts on $X$ by locally nilpotent transformations. Therefore by the Engel theorem \index{Engel's theorem} $Z$ must contain a nonzero vector $v\in {\rm Wh}(X)$. But since the restriction of $g$ to the subspace ${\rm Wh}(X)$ is injective by the assumption, one has $g(v)\neq 0$. Thus we arrive at a contradiction, and hence $g$ is injective.

\end{proof}

Now recall that the action of $J$ on $V$ is locally nilpotent. All non--zero Whittaker vectors in $V$ belong to ${\rm Wh}(V)$ and by the definition of $\sigma_\varepsilon$ their images in ${\rm hom}_{\mathbb{C}}(\mathbb{C}_{\varepsilon}^s[M_+], {\rm Wh}(V))$ are non--zero homomorphisms non-vanishing at $1$. Therefore by Lemma \ref{inj} $\sigma_\varepsilon$ is injective.

Next we show that $\sigma_\varepsilon$ is also surjective. Denote $x_{m_1\ldots m_c}={c_{m_1\ldots m_c}'(\varepsilon)}^{-1}\phi_\varepsilon(B_{m_1\ldots m_c}^\varepsilon)\in Q_\varepsilon^{loc}$, \index[not]{x@$x_{m_1\ldots m_c}$} where $c_{m_1\ldots m_c}'(\varepsilon)$ are the values of $c_{m_1\ldots m_c}'(q)$ at $q^{\frac{1}{d{\bar{r}}^2}}=\varepsilon^{\frac{1}{d{\bar{r}}^2}}$ which are all non--zero by the choice of $\varepsilon^{\frac{1}{d{\bar{r}}^2}}$. By Corollary \ref{Bpbas}
\begin{equation}\label{xbas}
{\rm Ad}_s(f_{\beta_1}^{r_1}\ldots f_{\beta_c}^{r_c})(x_{m_1\ldots m_c})=\left\{\begin{array}{l} 1 ~~{\rm if}~~m_p=r_p~~{\rm for}~~ p=1,\ldots , c \\ 0 ~~{\rm if}~~r_i=m_i, i=1,\ldots, p-1~~{\rm and}~~ r_p>m_p~~{\rm for~~some}~~p\in \{1,\ldots , c\} \end{array}\right..
\end{equation}

Since for any $v\in {\rm Wh}(V)$ and $z\in {\mathbb{C}}_{\varepsilon}^{s,loc}[G_*]\cap I_\varepsilon$ we have $zv=0$ (see Condition 3 in the definition of the category ${\mathbb{C}}_{\varepsilon}^{s,loc}[G_*]-{\rm mod}_{\mathbb{C}_{\varepsilon}^s[M_+]}^{\chi^{s}_\varepsilon}$), and $x_{m_1\ldots m_c}\in Q_\varepsilon^{loc}$, the definition of $Q_\varepsilon^{loc}$ implies that the elements $x_{m_1\ldots m_c}v$ are well--defined and by Condition 2 in the definition of the category ${\mathbb{C}}_{\varepsilon}^{s,loc}[G_*]-{\rm mod}_{\mathbb{C}_{\varepsilon}^s[M_+]}^{\chi^{s}_\varepsilon}$ and by the definition of Whittaker vectors one has ${\rm Ad}_sx(x_{m_1\ldots m_c}v)={\rm Ad}_sx(x_{m_1\ldots m_c})v$, $x\in \mathbb{C}_{\varepsilon}^s[M_+]$. For this reason formula (\ref{xbas}) implies 
$$
{\rm Ad}_s(f_{\beta_1}^{r_1}\ldots f_{\beta_c}^{r_c})(x_{m_1\ldots m_c}v)=\left\{\begin{array}{l} v ~~{\rm if}~~m_p=r_p~~{\rm for}~~ p=1,\ldots , c \\ 0 ~~{\rm if}~~r_i=m_i, i=1,\ldots, p-1~~{\rm and}~~ r_p>m_p~~{\rm for~~some}~~p\in \{1,\ldots , c\} \end{array}\right.,
$$
and hence
\begin{equation}\label{basisgen}
\sigma_\varepsilon(x_{m_1\ldots m_c}v)(f_{\beta_1}^{r_1}\ldots f_{\beta_c}^{r_c})=\left\{\begin{array}{l} v ~~{\rm if}~~m_p=r_p~~{\rm for}~~ p=1,\ldots , c \\ 0 ~~{\rm if}~~r_i=m_i, i=1,\ldots, p-1~~{\rm and}~~ r_p>m_p~~{\rm for~~some}~~p\in \{1,\ldots , c\} \end{array}\right..
\end{equation}

Observe that by Proposition \ref{Qdefr} (v) the elements $f_{\beta_1}^{r_1}\ldots f_{\beta_c}^{r_c}$ form a linear basis of $\mathbb{C}_{\varepsilon}^s[M_+]$. 
Elements of this basis are labeled by the elements of the set $\mathbb{N}^c$. Introduce the lexicographic  order on this set, so that $(r_1,\ldots, r_c)>(m_1,\ldots ,m_c)$  if $r_i=m_i$ for $i=1,\ldots, p-1$ and $r_p>m_p$ for some $p\in \{1,\ldots , c\}$. 


Now let ${\bf m}=(m_1,\ldots, m_c)\in \mathbb{N}^c$, $v\in {\rm Wh}(V)$ and denote 
$$
f^{\bf m}_v=\sigma_\varepsilon(x_{m_1\ldots m_c}v). \index[not]{f@$f^{\bf m}_v$}
$$

Since by the definition $f^{\bf m}_v\in {\rm hom}_{\mathbb{C}}(\mathbb{C}_{\varepsilon}^s[M_+], {\rm Wh}(V))$, and by (\ref{basisgen}) $f^{\bf m}_v(f_{\beta_1}^{r_1}\ldots f_{\beta_c}^{r_c})=0$ for ${\bf m}<{\bf r}=(r_1,\ldots, r_c)$. Therefore the homomorphism $f^{\bf m}_v$ may not vanish only on a finite number of elements $f_{\beta_1}^{r_1}\ldots f_{\beta_c}^{r_c}$ with ${\bf m}=(m_1,\ldots, m_c)\geq (r_1,\ldots ,r_c)={\bf r}$.

Let ${\bf n}^1=(n_1^1,\ldots, n_c^1)\in \mathbb{N}^c$  be the largest element such that ${\bf m}>{\bf n}^1$ and 
\begin{equation}\label{f1}
v_1:=f^{\bf m}_v(f_{\beta_1}^{n_1^1}\ldots f_{\beta_c}^{n_c^1})\neq 0.
\end{equation} 

Denote
$$
f^{{\bf n}^1}_{v_1}=\sigma_\varepsilon(x_{n_1^1\ldots n_c^1}v_1). 
$$
Then (\ref{basisgen}) and (\ref{f1}) imply that for $g_1:=f^{\bf m}_v-f^{{\bf n}^1}_{v_1}\in {\rm hom}_{\mathbb{C}}(\mathbb{C}_{\varepsilon}^s[M_+], {\rm Wh}(V))$ one has 
\begin{equation}\label{g1}
g_1(f_{\beta_1}^{r_1}\ldots f_{\beta_c}^{r_c})=\left\{\begin{array}{l} 0~{\rm if}~{\bf r}=(r_1,\ldots, r_c)\geq {\bf n}^1, {\bf r}\neq {\bf m} \\ v ~{\rm if}~{\bf r}={\bf m}\end{array}\right..
\end{equation}

By the definition $g_1\in {\rm hom}_{\mathbb{C}}(\mathbb{C}_{\varepsilon}^s[M_+], {\rm Wh}(V))$, and by (\ref{g1}) the homomorphism $g_1$ may not vanish on a finite number of elements $f_{\beta_1}^{r_1}\ldots f_{\beta_c}^{r_c}$ with ${\bf n}^1>{\bf r}=(r_1,\ldots, r_c)$. Let ${\bf n}^2=(n_1^2,\ldots, n_c^2)\in \mathbb{N}^c$ be the largest element such that ${\bf n}^1>{\bf n}^2$ and 
\begin{equation}\label{f2}
v_2:=g_1(f_{\beta_1}^{n_1^2}\ldots f_{\beta_c}^{n_c^2})\neq 0.
\end{equation} 

Denote
$$
f^{{\bf n}^2}_{v_2}=\sigma_\varepsilon(x_{n_1^2\ldots n_c^2}v_2). 
$$
Then (\ref{basisgen}) and (\ref{f2}) imply that for $g_2:=f^{\bf r}_v-f^{{\bf n}^1}_{v_1}-f^{{\bf n}^2}_{v_2}=g_1-f^{{\bf n}^2}_{v_2}\in {\rm hom}_{\mathbb{C}}(\mathbb{C}_{\varepsilon}^s[M_+], {\rm Wh}(V))$ one has 
\begin{equation}
g_2(f_{\beta_1}^{r_1}\ldots f_{\beta_c}^{r_c})=\left\{\begin{array}{l} 0~{\rm if}~{\bf r}=(r_1,\ldots, r_c)\geq {\bf n}^2, {\bf r}\neq {\bf m} \\ v ~{\rm if}~{\bf r}={\bf m}\end{array}\right.. 
\end{equation}

Iterating this procedure we obtain a sequence of elements ${\bf m}>{\bf n}^1>{\bf n}^2>\ldots>{\bf n}^i>\ldots$, ${\bf n}^i\in \mathbb{N}^c$ \index[not]{n@${\bf n}^i$} and a sequence of elements $g_i\in {\rm hom}_{\mathbb{C}}(\mathbb{C}_{\varepsilon}^s[M_+], {\rm Wh}(V))$ such that 
\begin{equation}
g_i(f_{\beta_1}^{r_1}\ldots f_{\beta_c}^{r_c})=\left\{\begin{array}{l} 0~{\rm if}~{\bf r}=(r_1,\ldots, r_c)\geq {\bf n}^i, {\bf r}\neq {\bf m} \\ v ~{\rm if}~{\bf r}={\bf m}\end{array}\right..
\end{equation}

Since by Theorem 1.13 in \cite{Ha} (see also Theorem 2.4.2 in \cite{Ba}) the lexicographic order on $\mathbb{N}^c$ is a well--order the sequence ${\bf m}>{\bf n}^1>{\bf n}^2>\ldots>{\bf n}^i>\ldots$ must be finite, i.e. for some $j\in \mathbb{N}$ 
\begin{equation}
g_{j}(f_{\beta_1}^{r_1}\ldots f_{\beta_c}^{r_c})=\left\{\begin{array}{l} 0~{\rm if}~{\bf r}=(r_1,\ldots, r_c)\geq {\bf n}^{j}, {\bf r}\neq {\bf m} \\ v ~{\rm if}~{\bf r}={\bf m}\end{array}\right.,
\end{equation}  
and there is no element ${\bf r}=(r_1,\ldots, r_c)<{\bf n}^{j}$ such that 
\begin{equation}
g_j(f_{\beta_1}^{r_1}\ldots f_{\beta_c}^{r_c})\neq 0.
\end{equation} 

Therefore $h^{\bf m}_v:=g_j$ \index[not]{h@$h^{\bf m}_v$} satisfies
$$
h^{\bf m}_v(f_{\beta_1}^{r_1}\ldots f_{\beta_c}^{r_c})=\left\{\begin{array}{l} v ~~{\rm if}~~(r_1,\ldots, r_c)=(m_1,\ldots, m_c),\\ 0 ~~{\rm if}~~(r_1,\ldots, r_c)\neq (m_1,\ldots, m_c)  \end{array}\right..
$$

This implies that the elements $h^{\bf m}_v$ with $(m_1,\ldots, m_c)\in \mathbb{N}^c$, $v\in {\rm Wh}(V)$ span
${\rm hom}_{\mathbb{C}}(\mathbb{C}_{\varepsilon}^s[M_+], {\rm Wh}(V))$.
We deduce that the elements $\sigma_\varepsilon(x_{m_1\ldots m_c}v)$ with $(m_1,\ldots, m_c)\in \mathbb{N}^c$, $v\in {\rm Wh}(V)$ span ${\rm hom}_{\mathbb{C}}(\mathbb{C}_{\varepsilon}^s[M_+], {\rm Wh}(V))$ as well. Therefore $\sigma_\varepsilon$ is surjective.

So we deduce that $\sigma_\varepsilon$ is in fact an isomorphism which proves Lemma \ref{selem}.

\end{proof}

Now the first claim in part (ii) is a consequence of Lemma \ref{selem}. The remaining two properties states in part (ii) obviously follow from the first claim. This completes the proof of (ii).

(iii) First we prove that
$Q_{\varepsilon}^{loc}$ is an object of ${\mathbb{C}}_{\varepsilon}^{s,loc}[G_*]-{\rm mod}_{\mathbb{C}_{\varepsilon}^s[M_+]}^{\chi_\varepsilon^s}$.
We shall show that the  adjoint action of the augmentation ideal $J$ of $\mathbb{C}_{\varepsilon}^s[M_+]$ on $Q_{\varepsilon}^{loc}$ is locally nilpotent. All the other properties of objects of the category ${\mathbb{C}}_{\varepsilon}^{s,loc}[G_*]-{\rm mod}_{\mathbb{C}_{\varepsilon}^s[M_+]}^{\chi_\varepsilon^s}$ for $Q_{\varepsilon}^{loc}$ follow from Proposition \ref{Qpr}.

Indeed, recalling the $\mathbb{C}_{\varepsilon}^s[M_+]$--module homomorphism of $\phi_\varepsilon:{\mathbb{C}}_{11}^{loc}[G]_\varepsilon\rightarrow Q_{\varepsilon}^{loc}$ and the definition of the adjoint action on ${\mathbb{C}}_{11}^{loc}[G]_\varepsilon$ in formula (\ref{Apo}) we deduce that in order to show that the  adjoint action of the augmentation ideal $J$ of $\mathbb{C}_{\varepsilon}^s[M_+]$ on $Q_{\varepsilon}^{loc}$ is locally nilpotent it suffices to show that the  ${\rm Ad}_s^0$--action of the augmentation ideal $J$ of $\mathbb{C}_{\varepsilon}^s[M_+]$ on ${\mathbb{C}}_{\varepsilon}^s[G]$ is locally nilpotent. But the last fact is true as ${\mathbb{C}}_{\varepsilon}^s[G]=\bigoplus_{\lambda \in P_+}V_\lambda^* \otimes V_\lambda$, where $V_\lambda$ is the finite--dimensional irreducible representation of $U_{\varepsilon}^s(\g)$ of highest weight $\lambda$, and the action of $J$ on $V_\lambda^* \otimes V_\lambda$ induced by the adjoint action is locally nilpotent since the action of $J$ on finite--dimensional irreducible representations is locally nilpotent.

Now we show that $Q_\varepsilon^{loc}$ is isomorphic to ${\rm hom}_{\mathbb{C}}(\mathbb{C}_{\varepsilon}^s[M_+],\mathbb{C}))\otimes {W_\varepsilon^{s,loc}(G)}$ as a $\mathbb{C}_{\varepsilon}^s[M_+]$--$W_\varepsilon^{s,loc}(G)$--bimodule.
Indeed, by (\ref{rt}) with $R=W_\varepsilon^{s,loc}(G)$ we have an obvious right $\mathbb{C}_{\varepsilon}^s[M_+]$--module isomorphism 
$$
{\rm hom}_{\mathbb{C}}(\mathbb{C}_{\varepsilon}^s[M_+], W_\varepsilon^{s,loc}(G))\simeq {\rm hom}_{\mathbb{C}}(\mathbb{C}_{\varepsilon}^s[M_+],\mathbb{C})\otimes W_\varepsilon^{s,loc}(G).
$$

Now consider the $\mathbb{C}_{\varepsilon}^s[M_+]$--submodule $\sigma_{\varepsilon}^{-1}({\rm hom}_{\mathbb{C}}(\mathbb{C}_{\varepsilon}^s[M_+],\mathbb{C}))$ of $Q_{\varepsilon}^{loc}$, where $${\rm hom}_{\mathbb{C}}(\mathbb{C}_{\varepsilon}^s[M_+],\mathbb{C})\subset {\rm hom}_{\mathbb{C}}(\mathbb{C}_{\varepsilon}^s[M_+], W_\varepsilon^{s,loc}(G)).$$ Obviously $\sigma_{\varepsilon}^{-1}({\rm hom}_{\mathbb{C}}(\mathbb{C}_{\varepsilon}^s[M_+],\mathbb{C}))\simeq {\rm hom}_{\mathbb{C}}(\mathbb{C}_{\varepsilon}^s[M_+],\mathbb{C})$ as a right $\mathbb{C}_{\varepsilon}^s[M_+]$--module.

Let $f_\varepsilon:\sigma_{\varepsilon}^{-1}({\rm hom}_{\mathbb{C}}(\mathbb{C}_{\varepsilon}^s[M_+],\mathbb{C}))\otimes W_\varepsilon^{s,loc}(G)\rightarrow Q_{\varepsilon}^{loc}$  be the map induced  by right the action of $W_\varepsilon^{s,loc}(G)$ on $Q_{\varepsilon}^{loc}$. Since by Proposition \ref{Qpr} this action commutes with the adjoint action of $\mathbb{C}_{\varepsilon}^s[M_+]$ on $Q_{\varepsilon}^{loc}$, and $\sigma_{\varepsilon}^{-1}({\rm hom}_{\mathbb{C}}(\mathbb{C}_{\varepsilon}^s[M_+],\mathbb{C}))\otimes W_\varepsilon^{s,loc}(G)$ is naturally a $\mathbb{C}_{\varepsilon}^s[M_+]$--$W_\varepsilon^{s,loc}(G)$--bimodule with the right action of $\mathbb{C}_{\varepsilon}^s[M_+]$ induced by the multiplication on $\mathbb{C}_{\varepsilon}^s[M_+]$ from the left and the right action of $W_\varepsilon^{s,loc}(G)$ induced by the multiplication on itself from the right, we infer that $f_\varepsilon$ is a homomorphism of $\mathbb{C}_{\varepsilon}^s[M_+]$--$W_\varepsilon^{s,loc}(G)$--bimodules.

We claim that $f_\varepsilon$ is injective. This follows straightforwardly from Lemma \ref{inj} because by the definition all Whittaker vectors of $\sigma_{\varepsilon}^{-1}({\rm hom}_{\mathbb{C}}(\mathbb{C}_{\varepsilon}^s[M_+],\mathbb{C}))\otimes W_\varepsilon^{s,loc}(G)$ belong to the subspace $$1\otimes W_\varepsilon^{s,loc}(G)\subset \sigma_{\varepsilon}^{-1}({\rm hom}_{\mathbb{C}}(\mathbb{C}_{\varepsilon}^s[M_+],\mathbb{C}))\otimes W_\varepsilon^{s,loc}(G),$$ and the restriction of $f_\varepsilon$ to this subspace is injective by the definition.

Now we show that $f_\varepsilon$ is surjective. We shall need the following lemma.
\begin{lemma}\label{surj}
Let $g':X\rightarrow Y$ be an injective homomorphism of $\mathbb{C}_{\varepsilon}^s[M_+]$--modules. Assume that $g'$ induces an isomorphism of the spaces of Whittaker vectors of $X$  and of $Y$, and that ${\rm Ext}^{1}_{\mathbb{C}_{\varepsilon}^s[M_+]}(\mathbb{C}_{\varepsilon_s}, X)=0$, where $\mathbb{C}_{\varepsilon_s}$ is the trivial representation of $\mathbb{C}_{\varepsilon}^s[M_+]$ induced by the counit. Suppose also that the action of the augmentation ideal $J$ of $\mathbb{C}_{\varepsilon}^s[M_+]$ on the cokernel of $g'$ is locally nilpotent. Then $g'$ is surjective.
\end{lemma}

\begin{proof}
Consider the exact sequence
$$
0\rightarrow X \rightarrow Y \rightarrow Y' \rightarrow 0,
$$
where $Y'$ is the cokernel of $g'$, and the corresponding long exact sequence of cohomology,
$$
0\rightarrow {\rm Ext}^{0}_{\mathbb{C}_{\varepsilon}^s[M_+]}(\mathbb{C}_{\varepsilon_s},
X)\rightarrow {\rm Ext}^{0}_{\mathbb{C}_{\varepsilon}^s[M_+]}(\mathbb{C}_{\varepsilon_s},
Y)\rightarrow {\rm Ext}^{0}_{\mathbb{C}_{\varepsilon}^s[M_+]}(\mathbb{C}_{\varepsilon_s},
Y')\rightarrow
$$
$$
\rightarrow {\rm Ext}^{1}_{\mathbb{C}_{\varepsilon}^s[M_+]}(\mathbb{C}_{\varepsilon_s},
X)\rightarrow \ldots .
$$

Since $g'$ induces an isomorphism of the spaces of Whittaker vectors of $X$  and of $Y$, and \\ ${\rm Ext}^{1}_{\mathbb{C}_{\varepsilon}^s[M_+]}(\mathbb{C}_{\varepsilon_s},
X)=0$, the initial part of the long exact cohomology sequence takes the form
$$
0\rightarrow {\rm Wh}(X)\rightarrow {\rm Wh}(Y)\rightarrow {\rm Wh}(Y')\rightarrow 0,
$$
where the second map in the last sequence is an isomorphism. Using the last exact sequence we deduce that ${\rm Wh}(Y')=0$. But the augmentation ideal $J$ acts on $Y'$ by locally nilpotent transformations. Therefore, by the Engel theorem, if $Y'$ is not trivial there should exists a nonzero Whittaker vector in it. Thus we arrive at a contradiction, and $Y'=0$. Therefore $g'$ is surjective.

\end{proof} 

Now we apply the previous lemma with $X=\sigma_{\varepsilon}^{-1}({\rm hom}_{\mathbb{C}}(\mathbb{C}_{\varepsilon}^s[M_+],\mathbb{C}))\otimes W_\varepsilon^{s,loc}(G)$, $Y=Q_\varepsilon^{loc}$ and $g'=f_\varepsilon$.
By Lemma \ref{inj1} one can immediately deduce that the right $\mathbb{C}_{\varepsilon}^s[M_+]$--module $\sigma_{\varepsilon}^{-1}({\rm hom}_{\mathbb{C}}(\mathbb{C}_{\varepsilon}^s[M_+],\mathbb{C}))\otimes W_\varepsilon^{s,loc}(G)\simeq {\rm hom}_{\mathbb{C}}(\mathbb{C}_{\varepsilon}^s[M_+], W_\varepsilon^{s,loc}(G))$ is injective. In particular, $${\rm Ext}^{1}_{\mathbb{C}_{\varepsilon}^s[M_+]}(\mathbb{C}_{\varepsilon},\sigma_{\varepsilon}^{-1}({\rm hom}_{\mathbb{C}}(\mathbb{C}_{\varepsilon}^s[M_+],\mathbb{C}))\otimes W_\varepsilon^{s,loc}(G)
)=0.$$ 

In the beginning of the proof of this part it was already shown that the  adjoint action of the augmentation ideal $J$ of $\mathbb{C}_{\varepsilon}^s[M_+]=U_{\varepsilon}^s(\m_-)$ on $Q_{\varepsilon}^{loc}$ is locally nilpotent.

Also by the definition $f_\varepsilon$ gives rise to an isomorphism of the spaces of Whittaker vectors of $\sigma_{\varepsilon}^{-1}({\rm hom}_{\mathbb{C}}(\mathbb{C}_{\varepsilon}^s[M_+],\mathbb{C}))\otimes W_\varepsilon^{s,loc}(G)$ and of $Q_\varepsilon^{loc}$, which are both isomorphic to $W_\varepsilon^{s,loc}(G)$.

Therefore all conditions of Lemma \ref{surj} are satisfied for the map $f_\varepsilon$, and hence $f_\varepsilon$ is surjective.
Thus $Q_\varepsilon^{loc}$ is isomorphic to ${\rm hom}_{\mathbb{C}}(\mathbb{C}_{\varepsilon}^s[M_+],\mathbb{C}))\otimes {W_\varepsilon^{s,loc}(G)}$ as a $\mathbb{C}_{\varepsilon}^s[M_+]$--$W_\varepsilon^{s,loc}(G)$--bimodule. This completes the proof of part (iii).

(i)
Now we prove the main claim of this theorem.
Let $E$ be a finitely generated $W_\varepsilon^{s,loc}(G)$--module. Since by part (iii) $Q_\varepsilon^{loc}\simeq {\rm hom}_{\mathbb{C}}(\mathbb{C}_{\varepsilon}^s[M_+],\mathbb{C}))\otimes {W_\varepsilon^{s,loc}(G)}$ as a $\mathbb{C}_{\varepsilon}^s[M_+]$--$W_\varepsilon^{s,loc}(G)$--bimodule, we have 
$$
Q_{\varepsilon}^{loc}\otimes_{W_\varepsilon^{s,loc}(G)}E\simeq {\rm hom}_{\mathbb{C}}(\mathbb{C}_{\varepsilon}^s[M_+],\mathbb{C}))\otimes {W_\varepsilon^{s,loc}(G)}\otimes_{W_\varepsilon^{s,loc}(G)}E\simeq {\rm hom}_{\mathbb{C}}(\mathbb{C}_{\varepsilon}^s[M_+],\mathbb{C}))\otimes E
$$ 
as a $\mathbb{C}_{\varepsilon}^s[M_+]$--module. This implies
$$
{\rm Wh}(Q_{\varepsilon}^{loc}\otimes_{W_\varepsilon^{s,loc}(G)}E)\simeq {\rm Hom}_{\mathbb{C}_{\varepsilon}^s[M_+]}(\mathbb{C}_{\varepsilon_s},Q_{\varepsilon}^{loc}\otimes_{W_\varepsilon^{s,loc}(G)}E)\simeq {\rm Hom}_{\mathbb{C}_{\varepsilon}^s[M_+]}(\mathbb{C}_{\varepsilon_s}, {\rm hom}_{\mathbb{C}}(\mathbb{C}_{\varepsilon}^s[M_+],\mathbb{C})\otimes E)\simeq 
$$
\begin{equation}\label{whiso}
\simeq {\rm Hom}_{\mathbb{C}_{\varepsilon}^s[M_+]}(\mathbb{C}_{\varepsilon_s}, {\rm hom}_{\mathbb{C}}(\mathbb{C}_{\varepsilon}^s[M_+],E))\simeq {\rm hom}_{\mathbb{C}}(\mathbb{C}_{\varepsilon}^s[M_+]\otimes_{\mathbb{C}_{\varepsilon}^s[M_+]} \mathbb{C}_{\varepsilon_s},E)\simeq   E.
\end{equation}
Therefore to prove part (i) of the theorem it suffices to check that for any object $V$ of the category ${\mathbb{C}}_{\varepsilon}^{s,loc}[G_*]-{\rm mod}_{\mathbb{C}_{\varepsilon}^s[M_+]}^{\chi^{s}_\varepsilon}$ the canonical map $f:Q_{\varepsilon}^{loc}\otimes_{W_\varepsilon^{s,loc}(G)}{\rm Wh}(V)\rightarrow V$ is an isomorphism. 

In order to do this we observe that $Q_{\varepsilon}^{loc}\otimes_{W_\varepsilon^{s,loc}(G)}{\rm Wh}(V)$ is an object of the category ${\mathbb{C}}_{\varepsilon}^{s,loc}[G_*]-{\rm mod}_{\mathbb{C}_{\varepsilon}^s[M_+]}^{\chi^{s}_\varepsilon}$ since $Q_{\varepsilon}^{loc}$ is an object of this category by part (iii).
In particular, the action of the augmentation ideal $J$ on $Q_{\varepsilon}^{loc}\otimes_{W_\varepsilon^{s,loc}(G)}{\rm Wh}(V)$ is locally nilpotent. By (\ref{whiso}) the space of Whittaker vectors of $Q_{\varepsilon}^{loc}\otimes_{W^{s,loc}_{\varepsilon}(G)}{\rm Wh}(V)$ is $1\otimes {\rm Wh}(V)$, and the restriction of $f$ to $1\otimes {\rm Wh}(V)$ induces an isomorphism of the spaces of Whittaker vectors of $Q_{\varepsilon}^{loc}\otimes_{W^{s,loc}_{\varepsilon}(G)}{\rm Wh}(V)$ and of $V$. Therefore $f$ is injective by Lemma \ref{inj}.

As we showed in Lemma \ref{selem}, $\sigma_\varepsilon:V\rightarrow {\rm hom}_{\mathbb{C}}(\mathbb{C}_{\varepsilon}^s[M_+], {\rm Wh}(V))\simeq {\rm hom}_{\mathbb{C}}(\mathbb{C}_{\varepsilon}^s[M_+],\mathbb{C})\otimes {\rm Wh}(V)$ is an isomorphism of $\mathbb{C}_{\varepsilon}^s[M_+]$--modules for any object $V$ of the category ${\mathbb{C}}_{\varepsilon}^{s,loc}[G_*]-{\rm mod}_{\mathbb{C}_{\varepsilon}^s[M_+]}^{\chi^{s}_\varepsilon}$. By Lemma \ref{inj1} ${\rm hom}_{\mathbb{C}}(\mathbb{C}_{\varepsilon}^s[M_+],\mathbb{C})\otimes {\rm Wh}(V)$ is injective over $\mathbb{C}_{\varepsilon}^s[M_+]$. Therefore $V$ is injective as an $\mathbb{C}_{\varepsilon}^s[M_+]$--module with respect to the adjoint action. In particular, $Q_{\varepsilon}^{loc}\otimes_{W^{s,loc}_{\varepsilon,\xi}(G)}{\rm Wh}(V)$ is injective over $\mathbb{C}_{\varepsilon}^s[M_+]$, and hence ${\rm Ext}^{1}_{\mathbb{C}_{\varepsilon}^s[M_+]}(\mathbb{C}_{\varepsilon_s},
Q_{\varepsilon}^{loc}\otimes_{W^{s,loc}_{\varepsilon}(G)}{\rm Wh}(V))=0$. 

Recall also that $f$ induces an isomorphism of the spaces of Whittaker vectors of $Q_{\varepsilon}^{loc}\otimes_{W^{s,loc}_{\varepsilon}(G)}{\rm Wh}(V)$ and of $V$ and that the adjoint action of $J$ on $V$ is locally nilpotent. Therefore $f$ is surjective by Lemma \ref{surj} with $X=Q_{\varepsilon}^{loc}\otimes_{W^{s,loc}_{\varepsilon}(G)}{\rm Wh}(V)$, $Y=V$, $g'=f$.
This completes the proof of the theorem.

\end{proof}


\section{Bibliographic comments}

\pagestyle{myheadings}
\markboth{CHAPTER \thechapter.~APPLICATION TO EQUIVARIANT MODULES OVER A QUANTUM GROUP}{\thesection.~BIBLIOGRAPHIC COMMENTS}

\setcounter{equation}{0}
\setcounter{theorem}{0}

A categorical equivalence for Lie algebras, called the Skryabin equivalence, similar to that considered in this chapter was established in the Appendix to \cite{Pr}.

The main theorem of this chapter is an improvement of Theorem 7.7 in \cite{SDM} where a similar equivalence was established in a quantum group case for $q$ specialized to generic values $q=\varepsilon\in \mathbb{C}$. The proof of Theorem 7.7 in \cite{SDM} relies on homological methods and arguments related to the properties of the quasiclassical limits $W^s(G)$ of q-W--algebras. In this book we use the approach  similar to the original Skryabin's idea. 

The definition of the category of equivariant representations over a quantum group given in Section \ref{cateq} is a slight modification of a similar definition given in Section 7 in \cite{SDM}, minor changes being related to the fact that in Corollary \ref{Bpbas}, and more generally in the previous chapter, we dealt not with quantum groups themselves but with their localizations.


\chapter{Application of q-W--algebras to quantum groups at roots of unity and the proof of De Concini--Kac--Procesi conjecture}\label{part6}

In this chapter we are going to use the elements $B_{m_1\ldots m_c}$ introduced in Corollary \ref{Bpbas} to study the structure of representations of quantum groups at roots of unity. De Concini and Kac observed that every irreducible representation of a quantum group $U_\varepsilon(\g)$ at a root of unity $\varepsilon$ is in fact a representation of a finite--dimensional quotient $U_\eta(\g)$ of the quantum group, and hence every such representation is finite--dimensional itself. The quotient $U_\eta(\g)$ here depends on the representation. Later De Concini, Kac and Procesi also conjectured that the dimension of every such representation is divisible by a number $b_\eta$ which depends on (an isomorphism class of) $U_\eta(\g)$; a precise definition of $b_\eta$ will be given in Theorem \ref{fdfree}. 

Our main goal in this chapter is to prove this conjecture. We shall also obtain other related results on the structure of finite--dimensional representations of $U_\eta(\g)$. Firstly we are going to use an observation that every finite--dimensional representation of $U_\eta(\g)$ can be equipped with a second right action of a finite--dimensional subalgebra $U_{\eta_1}(\m_-)$ of the so--called small quantum group, and the dimension of this subalgebra is equal to $b_\eta$. The choice of the subalgebra $U_{\eta_1}(\m_-)$ depends on $U_\eta(\g)$ and the action of $U_{\eta_1}(\m_-)$ satisfies a compatibility condition similar to condition (\ref{compat}) for equivariant modules over quantum groups at generic $\varepsilon$. Thus every finite--dimensional representation of $U_\eta(\g)$ is in fact an equivariant $U_\eta(\g)-U_{\eta_1}(\m_-)$--bimodule. Next we prove that every finite--dimensional representation of $U_\eta(\g)$ is cofree over the corresponding subalgebra $U_{\eta_1}(\m_-)$ which confirms, in particular, the De Concini--Kac--Procesi conjecture. Remarkably, to prove this statement one can apply almost verbatim the arguments from the proof of Theorem \ref{sqeq} on the quantum group version of the Skryabin equivalence for generic $\varepsilon$ which overemphasizes again a striking similarity between the categories of finite--dimensional representations of algebras $U_\eta(\g)$ and the categories of equivariant modules over quantum groups introduced in Section \ref{cateq}.
 
The peculiarity of the root of unity case is that one can explicitly construct $U_{\eta_1}(\m_-)$--cofree bases of finite--dimensional $U_\eta(\g)$--modules using the elements $B_{m_1\ldots m_c}$ from Corollary \ref{Bpbas}.


\section{Quantum groups at roots of unity}\label{1root}

\pagestyle{myheadings}
\markboth{CHAPTER~\thechapter. APPLICATION TO QUANTUM GROUPS AT ROOTS OF UNITY}{\thesection.~QUANTUM GROUPS AT ROOTS OF UNITY}

\setcounter{equation}{0}
\setcounter{theorem}{0}

In this section we recall some results on representation theory of quantum groups at roots of unity. 

Let $\bar{m}$ \index[not]{m@$\bar{m}$} be an odd positive integer number such that $\bar{m}>d_i$ is coprime to all $d_i$ for all $i=1,\ldots, l$, $\varepsilon$ a primitive $\bar{m}$-th root of unity. An appropriate number $d$, which appears in the definition of the algebras $U_\varepsilon({\frak g})$ and $U_\varepsilon^s({\frak g})$, can be found from the following proposition.
\begin{proposition}\label{d'}
Let $\Delta$ be an irreducible root system, $\Delta_+^s$ the system of positive roots associated to the conjugacy class of a Weyl group element $s\in W$ in Theorem \ref{mainth}, $s=s_{\gamma_1}\ldots s_{\gamma_{l'}}$ representation (\ref{inv}) for $s$, $\alpha_1,\ldots,\alpha_l$ the set of simple roots in $\Delta_+^s$. Then

(i) if $\Delta$ is of exceptional type the lowest common multiple $d'$ \index[not]{d@$d'$} of the denominators of the numbers $\frac{1}{d_j}\left\langle  {1+s \over 1-s }P_{{\h'}^*}\alpha_i,\alpha_j\right\rangle$, where $i,j=1,\ldots,l$ is given in the tables in Appendix 2;

(ii) if $\Delta$ is of classical type then the conjugacy class of $s$ corresponds to the sum of a number of blocks as in (\ref{san}), (\ref{scn}), (\ref{sbn}) or (\ref{sdn}). To each block of type $X$ we associate an integer $d_{ij}(X)$, $i,j=1,\ldots , l$ as follows:

if $\Delta$ is not of type $A_l,D_l$, an orbit with the smallest number of elements for the action of the group $\left< s \right>$ on $E$ corresponds to a block of type $A_n$ and $s$ does not fix any root from $\Delta$ then

for $\Delta=B_l$
\begin{equation}\label{abn}
d_{ij}(A_n)=\left\{
         \begin{array}{ll}
           2p+1 & \hbox{if $n=2p$ is even;} \\
           p+1 & \hbox{if $n=2p+1$, $n\neq 4p-1$ is odd;} \\
           p  &  \hbox{if $n=4p-1$ is odd and $i<j$;} \\
           2p  &  \hbox{if $n=4p-1$ is odd and $i>j$;}
         \end{array}
       \right.
\end{equation}

for $\Delta=C_l$
\begin{equation}\label{acn}
d_{ij}(A_n)=\left\{
         \begin{array}{ll}
           2p+1 & \hbox{if $n=2p$ is even;} \\
           p+1 & \hbox{if $n=2p+1$, $n\neq 4p-1$ is odd;} \\
           2p  &  \hbox{if $n=4p-1$ is odd and $i<j$;} \\
           p  &  \hbox{if $n=4p-1$ is odd and $i>j$;}
         \end{array}
       \right.
\end{equation}

for $\Delta=D_l$ if $A_{l-1}\subset D_l$ is the only nontrivial block of the conjugacy class of $s$ then
\begin{equation}\label{adn}
d_{ij}(A_{l-1})=\left\{
         \begin{array}{ll}
           2p+1 & \hbox{if $l=2p+1$ is odd;} \\
           p+1 & \hbox{if $l=2p+2$, $l\neq 4p$ is even;} \\
           p  &  \hbox{if $l=4p$ is even;}
         \end{array}
       \right.
\end{equation}

for $\Delta=A_l$ if $s$ is a representative in the Coxeter conjugacy class, i.e. the conjugacy class of $s$ corresponds to the block of type $A_l$, then
\begin{equation}\label{acox}
d_{ij}(A_l)=1;
\end{equation}

in all other cases

\begin{equation}\label{asimplylaced}
d_{ij}(A_k)=\left\{
         \begin{array}{ll}
           k+1 & \hbox{if $k$ is even;} \\
           \frac{k-1}{2}+1 & \hbox{if $k$ is odd;}
         \end{array}
       \right.
\end{equation}

in all cases
$$
d_{ij}(C_n)=d_{ij}(B_n)=d_{ij}(D_{v+w}(a_{w-1}))=1,
$$
where, as before, we use the notation of \cite{C}, Section 7 for (blocks of) Weyl group conjugacy classes.

Then a common multiple $d'$ \index[not]{d@$d'$} of the denominators of the numbers $\frac{1}{d_j}\left\langle {1+s \over 1-s }P_{{\h'}^*}\alpha_i,\alpha_j\right\rangle$, where $i,j=1,\ldots,l$ is the lowest common multiple of the numbers $d_{ij}(X)$, \index[not]{d@$d_{ij}(X)$} $i,j=1,\ldots,l$ for all blocks $X$ of (the conjugacy class of) $s$.

(iii) If $\alpha_1',\ldots \alpha_l'$ is another system of simple roots in $\Delta$ then a common multiple of the denominators of the numbers $\frac{1}{d_j}\left\langle {1+s \over 1-s }P_{{\h'}^*}\alpha_i,\alpha_j\right\rangle$ will be also a common multiple of the denominators of the numbers $\frac{1}{d_j}\left\langle {1+s \over 1-s }P_{{\h'}^*}\alpha_i',\alpha_j'\right\rangle$ and vice versa.

\end{proposition}

\begin{proof}

(iii) First observe that
if $\Delta'_+$ is another system of positive roots in $\Delta$ with the simple roots $\alpha_1',\ldots ,\alpha_l'$ then $\alpha_i=\sum_{k=1}^lc_i^k\alpha_k'$, $\alpha_j^\vee=\sum_{k=1}^lb_i^k\alpha_k'^\vee$, where $c_i^k,b_i^k$ are integer coefficients. Hence
$$\frac{1}{d_j}\left\langle {1+s \over 1-s }P_{{\h'}^*}\alpha_i,\alpha_j\right\rangle=\left\langle {1+s \over 1-s }P_{{\h'}^*}\alpha_i,\alpha_j^\vee\right\rangle=\sum_{k,p=1}^lc_i^kb_j^p\left\langle {1+s \over 1-s }P_{{\h'}^*}\alpha_k',\alpha_p'^\vee\right\rangle=\sum_{k,p=1}^lc_i^kb_j^p\frac{1}{d_p}\left\langle {1+s \over 1-s }P_{{\h'}^*}\alpha_k',\alpha_p'\right\rangle,$$
and a common multiple of the denominators of the numbers $\frac{1}{d_j}\left\langle {1+s \over 1-s }P_{{\h'}^*}\alpha_i',\alpha_j'\right\rangle$ will be also a common multiple of the denominators of the numbers $\frac{1}{d_j}\left\langle {1+s \over 1-s }P_{{\h'}^*}\alpha_i,\alpha_j\right\rangle$ and vice versa.

(ii) In the case of classical irreducible root systems we shall compute a common multiple $d'$ of the denominators of the numbers $\frac{1}{d_j}\left\langle {1+s \over 1-s }P_{{\h'}^*}\alpha_i',\alpha_j'\right\rangle$, where $\Delta'_+$ is chosen in such a way that $s$ is elliptic in a parabolic Weyl subgroup $W'\subset W$ generated by the simple reflections corresponding to the roots from a subset of $\alpha_1',\ldots, \alpha_l'$ (for instance, one can take $\Delta_+'=\Delta_+^{s'}$ from the proof of Theorem \ref{mainth}).

Since different blocks of the conjugacy class of $s$ correspond to different disjoint mutually orthogonal subsets of simple roots in $\alpha_1',\ldots, \alpha_l'$
it suffices to consider the case when the conjugacy class of $s$ corresponds to a diagram with a single nontrivial block.
We shall compute $d'$ in the case when this block is of type $A_{k}$, $k>1$. Other cases can be considered in a similar way.
Assume that the root system $\Delta$ is realized as in Section \ref{stt}, where $V$ is a real Euclidean $n$--dimensional vector space equipped with the standard scalar product, with an orthonormal basis $\varepsilon_1,\ldots ,\varepsilon_n$. In that case simple roots are

\subsection*{$\bf A_n$} $\alpha'_i=\varepsilon_i-\varepsilon_{i+1}$, $1\leq i\leq n$;

\subsection*{$\bf B_n$} $\alpha'_i=\varepsilon_i-\varepsilon_{i+1}$, $1\leq i< n$, $\alpha'_n=\varepsilon_n$;

\subsection*{$\bf C_n$} $\alpha'_i=\varepsilon_i-\varepsilon_{i+1}$, $1\leq i< n$, $\alpha'_n=2\varepsilon_n$;

\subsection*{$\bf D_n$} $\alpha'_i=\varepsilon_i-\varepsilon_{i+1}$, $1\leq i< n$, $\alpha'_n=\varepsilon_{n-1}+\varepsilon_n$;

Then $s$ is of the form
$$
s=s^1s^2,~s^1=s_{\alpha'_{p+1}}s_{\alpha'_{p+3}}\ldots,
~s^2=s_{\alpha'_{p+2}}s_{\alpha'_{p+4}}\ldots,
$$
for some $p\in \mathbb{N}$, where in the formulas for $s^{i}, ~i=1,2$ the products are taken over mutually orthogonal simple roots labeled by indexes of the same parity; the last simple root which appears in those products is $\alpha'_{p+k}=\varepsilon_{p+k}-\varepsilon_{p+k+1}$, so $l'=k$, and $\{\gamma_1,\ldots,\gamma_k\}=\{\alpha'_{p+1},\ldots,\alpha'_{p+k}\}$.

Note that in the case of the root system $B_n$ the canonical bilinear form introduced in Section \ref{notation} is different from the standard scalar product on $V$ used in the description of the root system $B_n$ used in this proof. But the numbers $\left\langle {1+s \over 1-s }P_{{\h'}^*}\alpha_i',\alpha_j'^\vee\right\rangle$, that we have to compute, do not depend on the normalization of the scalar product, so we shall use the standard scalar product on $V$ in all cases in this proof. 

We consider the case when $i<j$. The case when $i>j$ can be obtained from it by observing that
\begin{equation}\label{ijji}
\left\langle {1+s \over 1-s }P_{{\h'}^*}\alpha_i',\alpha_j'^\vee\right\rangle=-\left\langle {1+s \over 1-s }P_{{\h'}^*}\alpha_j',\alpha_i'^\vee\right\rangle\frac{\left\langle \alpha_i',\alpha_i'\right\rangle}{\left\langle \alpha_j',\alpha_j'\right\rangle}.
\end{equation}

First recall that
by Lemma \ref{tmatrel}
\begin{equation}\label{matrel'}
\left\langle {1+s \over 1-s }P_{{\h'}^*}\gamma_i , \gamma_j \right\rangle=
\varepsilon_{ij}\left\langle \gamma_i,\gamma_j\right\rangle,
\end{equation}
where
$$
\varepsilon_{ij} =\left\{ \begin{array}{ll}
-1 & i <j \\
0 & i=j \\
1 & i >j
\end{array}
\right. .
$$

Let $\omega_t'$ be the fundamental weights of the root subsystem $A_{k}\subset \Delta$ with respect to the basis of simple roots $\alpha'_i$, $i=p+1,\ldots,p+k$,
$$
\omega_t'=\varepsilon_{p+1}+\ldots+\varepsilon_{p+t}-\frac{t}{k+1}\sum_{j=1}^{k+1}\varepsilon_{p+j}, t=1,\ldots, k.
$$
Since $\alpha'_{p+t}$, $t=1,\ldots, k$ form a linear basis of ${\h'}^*$, and $\omega_t'$, $t=1,\ldots, k$ form the dual basis we have
$$
\left\langle {1+s \over 1-s }P_{{\h'}^*}\alpha_i',\alpha_j'^\vee\right\rangle=\sum_{t,u=1}^{k}\left\langle \omega_t'^\vee,\alpha_i'\right\rangle\left\langle {1+s \over 1-s }P_{{\h'}^*}\alpha_{p+t}',\alpha_{p+u}'^\vee\right\rangle\left\langle \omega_u',\alpha_j'^\vee\right\rangle.
$$
Since the scalar product in $V$ is normalized in such a way that $\alpha_{p+u}'^\vee=\alpha_{p+u}'$, $u=1,\ldots,k$ we obtain using (\ref{matrel'})
\begin{eqnarray}\label{saa}
\left\langle {1+s \over 1-s }P_{{\h'}^*}\alpha_i',\alpha_j'^\vee\right\rangle=\sum_{t,u=1}^{k}\left\langle \omega_t'^\vee,\alpha_i'\right\rangle\left\langle {1+s \over 1-s }P_{{\h'}^*}\alpha_{p+t}',\alpha_{p+u}'\right\rangle\left\langle \omega_u',\alpha_j'^\vee\right\rangle= \\
\qquad \qquad \qquad \qquad \qquad \qquad =\sum_{t=1}^{k}(-1)^t\left\langle \omega_t'^\vee,\alpha_i'\right\rangle
\left\langle \omega_{t-1}'+\omega_{t+1}',\alpha_j'^\vee\right\rangle, \nonumber
\end{eqnarray}
where we assume that $\omega_0'=\omega_{k+1}'=0$.

Now one has to consider several cases.

If one of the roots $\alpha_i',\alpha_j'$ is orthogonal to ${\h'}^*$ then the left hand side of the last equality is zero.

If $\alpha_i',\alpha_j'\in \{\gamma_1,\ldots,\gamma_k\}$ then by (\ref{matrel'}) the left hand side of (\ref{saa}) is equal to $\pm 1$.

If $\alpha_i'=\alpha_{p+t}'$, $1<t<k$, $\alpha_j'=\alpha_{p+k+1}'$ then
\begin{equation}\label{dl}
\left\langle {1+s \over 1-s }P_{{\h'}^*}\alpha_i',\alpha_j'^\vee\right\rangle=(-1)^t\left\langle \omega_{t-1}'+\omega_{t+1}',\alpha_{p+k+1}'^\vee\right\rangle=
(-1)^t(\vartheta_1-\vartheta_2\frac{2t}{k+1}),
\end{equation}
where $\vartheta_2=2$ if $\alpha_j'^\vee=2\varepsilon_{p+k+1}$ or $\alpha_j'^\vee=\varepsilon_{p+k}+\varepsilon_{p+k+1}$, $\vartheta_1=0$ in the former case, and $\vartheta_1=1$ in the latter case. In all other cases $\vartheta_1=0$ and $\vartheta_2=1$. 

Note that $\vartheta_2\neq 1$ only in the case when $\Delta$ is of type $B_n$ or $D_n$; for arbitrary $s$ this situation can only be realized if an orbit with the smallest number of elements for the action of the group $\left< s \right>$ on $E$ corresponds to a block of type $A_k$ and $s$ does not fix any root from $\Delta$. The denominator $r$ of the number in the right hand side of (\ref{dl}) is given by
\begin{equation}\label{d1'}
r=\left\{
         \begin{array}{ll}
           2p+1 & \hbox{if $k=2p$ is even;} \\
           p+1 & \hbox{if $k=2p+1$, $k\neq 4p-1$ is odd;} \\
           \frac{2p}{\vartheta_2}  &  \hbox{if $k=4p-1$ is odd.}
         \end{array}
       \right.
\end{equation}

If $\alpha_i'=\alpha_{p+1}'$, $\alpha_j'=\alpha_{p+k+1}'$ then
$$
\left\langle {1+s \over 1-s }P_{{\h'}^*}\alpha_i',\alpha_j'^\vee\right\rangle=-\left\langle \omega_{2}',\alpha_{p+k+1}'^\vee\right\rangle=
\vartheta_2\frac{2}{k+1}-\vartheta_1',
$$
where as before $\vartheta_2=2$ if $\alpha_j'^\vee=2\varepsilon_{p+k+1}$ or $\alpha_j'^\vee=\varepsilon_{p+k}+\varepsilon_{p+k+1}$, $\vartheta_1'=0$ in the former case, and $\vartheta_1'=1$ in the latter case if in addition $k=2$. In all other cases $\vartheta_1'=0$ and $\vartheta_2=1$. We again obtain (\ref{d1'}).

If $\alpha_i'=\alpha_{p+k}'$, $\alpha_j'=\alpha_{p+k+1}'$ then
$$
\left\langle {1+s \over 1-s }P_{{\h'}^*}\alpha_i',\alpha_j'^\vee\right\rangle=(-1)^{k}\left\langle \omega_{k-1}',\alpha_{p+k+1}'^\vee\right\rangle=
-(-1)^k\vartheta_2\frac{k-1}{k+1},
$$
and we obtain (\ref{d1'}).

If $\alpha_i'=\alpha_{p}'$, $\alpha_j'=\alpha_{p+k+1}'$ then
\begin{eqnarray*}
\left\langle {1+s \over 1-s }P_{{\h'}^*}\alpha_i',\alpha_j'^\vee\right\rangle=\sum_{t=1}^{k}(-1)^t\left\langle \omega_t',\alpha_p'^\vee\right\rangle
\left\langle \omega_{t-1}'+\omega_{t+1}',\alpha_{p+k+1}'^\vee\right\rangle= \qquad \qquad \qquad \qquad \qquad \qquad \qquad\\ =-\sum_{t=1}^{k-1}(-1)^t\left(-1+\frac{t}{k+1}\right)
\frac{2t\vartheta_2}{k+1}-
(-1)^k\left(-1+\frac{k}{k+1}\right)\frac{k-1}{k+1}\vartheta_2+(-1)^{k-1}\vartheta_1\left(-1+\frac{k-1}{k+1}\right).
\end{eqnarray*}

Using the fact that
$$
\sum_{r=1}^n(-1)^{r+1}r^2=(-1)^{n+1}\frac{n(n+1)}{2}~{\rm and}~\sum_{r=1}^n(-1)^{r+1}r=
\left\{
\begin{array}{ll}
\frac{n+1}{2} & \hbox{ if $n$ is odd;} \\
-\frac{n}{2} & \hbox{if $n$ is even}
\end{array}
\right.
$$
we obtain
\begin{eqnarray*}
\left\langle {1+s \over 1-s }P_{{\h'}^*}\alpha_i',\alpha_j'^\vee\right\rangle=\left\{
\begin{array}{ll}
-\frac{\vartheta_2}{k+1}+\vartheta_1\frac{2}{k+1} & \hbox{if $k$ is even;} \\
-\vartheta_1\frac{2}{k+1}  & \hbox{if $k$ is odd.}
\end{array}
\right.
\end{eqnarray*}

The denominator $r$ of the number in the right hand side of the last equality is given by
$$
r=\left\{
         \begin{array}{ll}
           2p+1 & \hbox{if $k=2p$ is even;} \\
           1 & \hbox{if $k=2p+1$, $n$ is odd and $\vartheta_1=0$;} \\
           p+1  &  \hbox{if $k=2p+1$ is odd and $\vartheta_1=1$.}
         \end{array}
       \right.
$$

Summarizing all cases considered above and adding the case $i>j$ (see (\ref{ijji})) we arrive at (\ref{abn}), (\ref{acn}), (\ref{adn}), (\ref{acox}) and (\ref{asimplylaced}).

Other cases can be treated in a similar way.

Part (i) is proved by a direct verification using a computer code.

\end{proof}

Now one can choose $d=2d'$, where $d'$ is defined in the previous proposition. However, we shall not always assume that $d$ is chosen in this way, and the system of positive roots $\Delta_+^s$ will not be chosen as in Proposition \ref{d'} unless it is explicitly specified.

We shall always assume that $d$ and $\bar{m}$ are coprime. This condition is equivalent to the existence of an integer $\bar{n}$ \index[not]{n@$\bar{n}$} such that $\varepsilon^{\bar{n}d-1}=1$. From now on we shall also assume that $\kappa=\bar{n}d$. With this choice of $\kappa$ we have the following relation between the generators $t_i$ and $L_i$ of the quantum group $U_{\mathcal{A}}({\frak g})$, $t_i=L_i^{\bar{n}}$. In particular, the specialization $U_\varepsilon({\frak g})$ of $U_{\mathcal{A}}({\frak g})$ coincides with the specialization of the simply connected form of the standard Drinfeld--Jimbo quantum group without generators $t_i$ at $q=\varepsilon$.

Let $Z_\varepsilon$ \index[not]{Z@$Z_\varepsilon$} be the center of $U_\varepsilon({\frak g})$. In the following proposition we summarize the results on the structure of $Z_\varepsilon$. In particular, we recall that in the case when $\varepsilon$ is a root of unity $Z_\varepsilon$ is much larger than in the case of a generic $\varepsilon$. In fact, in the former case $Z_\varepsilon$ contains a remarkable subalgebra $Z_0$ the properties of which impose  very strong restrictions on the structure of irreducible representations of $U_\varepsilon({\frak g})$.

\begin{proposition}{\bf (\cite{DK}, Corollary 3.3; \cite{DKP1}, Theorems 3.5, 7.6 and Proposition 4.5)}\label{ue}
Fix the normal ordering in the positive root system $\Delta_+$ corresponding a reduced decomposition $\overline{w}=s_{i_1}\ldots s_{i_D}$ of the longest element $\overline{w}$ of the Weyl group $W$ of $\g$ and let $X_{\alpha}^\pm$ be the corresponding quantum root vectors in $U_\varepsilon({\frak g})$, and $X_{\alpha}$ the corresponding root vectors in $\g$. Let $x_{\alpha}^-=(\varepsilon_{\alpha}-\varepsilon_{\alpha}^{-1})^{\bar{m}}(X_{\alpha}^-)^{\bar{m}}$, \index[not]{x@$x_\alpha^-$} $x_{\alpha}^+=(\varepsilon_{\alpha}-\varepsilon_{\alpha}^{-1})^{\bar{m}}T_{\overline{w}}(X_{\alpha}^-)^{\bar{m}}$, \index[not]{x@$x_\alpha^+$} where $T_{\overline{w}}=T_{i_1}\ldots T_{i_D}$, $\alpha\in \Delta_+$, and $l_i=L_i^{\bar{m}}$, \index[not]{l@$l_i$} $i=1,\ldots ,l$.

Then the following statements are true.

(i) The elements $x_{\alpha}^\pm$, $\alpha\in \Delta_+$, $l_i$, $i=1,\ldots ,l$ lie in $Z_\varepsilon$.

(ii) Let $Z_0$ \index[not]{Z@$Z_0$} (resp. $Z_0^\pm$ \index[not]{Z@$Z_0^\pm$} and $Z_0^0$) \index[not]{Z@$Z_0^0$} be the subalgebras of $Z_\varepsilon$ generated by the $x_{\alpha}^\pm$ and the $l_i^{\pm 1}$ (respectively by the $x_{\alpha}^\pm$ and by the $l_i^{\pm 1}$). Then $Z_0^\pm\subset U_\varepsilon({\frak n}_\pm)$, $Z_0^0\subset U_\varepsilon({\frak h})$, $Z_0^\pm$ is the polynomial algebra with generators $x_{\alpha}^\pm$,
$Z_0^0$ is the algebra of Laurent polynomials in the $l_i$, $Z_0^\pm=U_\varepsilon({\frak n}_\pm)\cap Z_0$, and multiplication defines an isomorphism of algebras
$$
Z_0^-\otimes Z_0^0 \otimes Z_0^+ \rightarrow Z_0.
$$
The subalgebra $Z_0$ is independent of the choice of the reduced decomposition $\overline{w}=s_{i_1}\ldots s_{i_D}$.

(iii) $U_\varepsilon({\frak g})$ is a free $Z_0$--module with basis the set of monomials $(X^+)^{\bf r}L^{\bf p}(X^-)^{\bf m}$ for which $0\leq r_k,m_k,p_i<\bar{m}$ for $i=1,\ldots ,l$, $k=1,\ldots ,D$, where for ${\bf p}=(p_1,\ldots p_l)\in {\Bbb Z}^{l}$,
$$
L^{{\bf p}}=L_1^{{p_1}}\ldots L_l^{{p_l}}. \index[not]{L@$L^{{\bf p}}$}
$$

(iv) ${\rm Spec}(Z_0)\simeq \mathbb{C}^{2D}\times(\mathbb{C}^*)^l$; this variety has dimension equal to ${\rm dim}~\g$.

(v) The subalgebra $Z_0$ is preserved by the action of the braid group automorphisms $T_i$.

(vi) Let $G$ be the connected simply connected Lie group corresponding to the Lie algebra $\g$, $B_+\subset G$ the Borel subgroup corresponding to $\Delta_+$, $H\subset B_+$ the maximal torus contained in it, $N_+\subset B_+$ its unipotent radical, $N_-$ the opposite unipotent radical, and $G^*_0$ \index[not]{G@$G^*_0$} the solvable algebraic subgroup in $G\times G$ which consists of elements of the form $(L_+,L_-)\in G\times G$,
$$
(L_+,L_-)=(t,t^{-1})(n_+,n_-),~n_\pm \in N_\pm,~t\in H.
$$

Then ${\rm Spec}(Z_0^0)$ can be naturally identified with the maximal torus $H$ in $G$, and the map
$$
\widetilde{\pi}: {\rm Spec}(Z_0)={\rm Spec}(Z_0^+)\times {\rm Spec}(Z_0^0) \times {\rm Spec}(Z_0^-)\rightarrow G^*_0, \index[not]{p@$\widetilde{\pi}$}
$$
$$
\widetilde{\pi}(u_+,t,u_-)=(t{\bf X}^+(u_+),t^{-1}{\bf X}^-(u_-)^{-1}),~u_\pm\in {\rm Spec}(Z_0^\pm),~t\in {\rm Spec}(Z_0^0),
$$
$$
{\bf X}^\pm: {\rm Spec}(Z_0^\pm) \rightarrow N_\pm, \index[not]{X@${\bf X}^\pm$}
$$
$$
{\bf X}^-=\exp(x_{\beta_D}^-X_{-\beta_D})\exp(x_{\beta_{D-1}}^-X_{-\beta_{D-1}})\ldots \exp(x_{\beta_1}^-X_{-\beta_1}),
$$
$$
{\bf X}^+=\exp(x_{\beta_D}^+T_{\overline{w}}(X_{-\beta_D}))\exp(x_{\beta_{D-1}}^+T_{\overline{w}}(X_{-\beta_{D-1}}))\ldots \exp(x_{\beta_1}^+T_{\overline{w}}(X_{-\beta_1})),
$$
where $x_{\beta_i}^\pm$ should be regarded as complex-valued functions on ${\rm Spec}(Z_0)$, is an isomorphism of varieties independent of the choice of reduced decomposition of $\overline{w}$.
\end{proposition}

Parts (ii) and (iii) of Proposition \ref{ue} can also be reformulated in terms of the quantum root vectors $e_{\alpha}$ and $f_{\alpha}$.
\begin{proposition}\label{uq1z} 
Let $s\in W$ be a Weyl group element, and $e_{\alpha}$, $f_{\alpha}$ the quantum root vectors defined in Proposition \ref{rootss}. Then the following statements are true.

(i) The subalgebra $Z_0$ is the tensor product of the polynomial algebra with generators $e_{\alpha}^{\bar{m}}$, $f_{\alpha}^{\bar{m}}$, $\alpha \in \Delta_+$ and of the algebra of Laurent polynomials in $l_i$, $i=1,\ldots,l$.

(ii) $U_\varepsilon({\frak g})$ is a free $Z_0$--module with basis the set of monomials $f^{\bf r}L^{\bf p}e^{\bf m}$ for which $0\leq r_k,m_k,p_i<\bar{m}$ for $i=1,\ldots ,l$, $k=1,\ldots ,D$.
\end{proposition}

Let ${\bf K}:{\rm Spec}(Z_0^0)\rightarrow H$ \index[not]{K@${\bf K}$} be the map defined by ${\bf K}(h)=h^2$, $h\in H$.

\begin{proposition}{\bf (\cite{DKP1}, Corollary 4.7)}\label{LO}
Let $G^0=N_-HN_+$ be the big Bruhat cell in $G$. Then the map
$$
\pi:={\bf X}^-{\bf K}{\bf X}^+:{\rm Spec}(Z_0)\rightarrow G^0 \index[not]{p@$\pi$}
$$
is independent of the choice of reduced decomposition of $\overline{w}$, and is an unramified covering of degree $2^l$.
\end{proposition}

Note that $G_0^*$ is an example of dual Poisson--Lie groups $G^*$, considered as groups, and introduced in Section \ref{qplgroups}, $G_0^*$ being associated to the bialgebra structure defined by the r--matrix $r^s$ with $s=1$. 

As in Propositions \ref{locfin} and \ref{dressingact} we denote by $q:G^*_0\rightarrow G^0$ the map defined by $q(L_+,L_-)=L_-^{-1}L_+$. Then obviously
$\pi=q\circ \widetilde{\pi}$.

Another important property of quantum groups at roots of unity, which distinguishes the root of unity case,  is the existence of the so--called quantum coadjoint action which is an automorphism group action on an extension of ${U}_\varepsilon(\g)$. It is defined with the help of derivations $\underline{x}_i^\pm$ of $U_{\mathcal{A}}(\g)$ given by
\begin{equation}\label{qder}
\underline{x}_i^+(u)=\left[ \frac{(X_i^+)^{\bar{m}}}{[{\bar{m}}]_{q_i}!},u \right],~\underline{x}_i^-(u)=T_{\overline{w}}\underline{x}_i^+T_{\overline{w}}^{-1}(u),~i=1,\ldots ,l,~u\in U_{\mathcal{A}}(\g). \index[not]{x@$\underline{x}_i^+$}
\end{equation}

Let $\widehat{Z}_0$ \index[not]{Z@$\widehat{Z}_0$} be the algebra of formal power series in the $x_{\alpha}^\pm$, $\alpha\in \Delta_+$, and the $l_i^{\pm1}$, $i=1,\ldots ,l$, which define holomorphic functions on ${\rm Spec}(Z_0)\simeq \mathbb{C}^{2D}\times(\mathbb{C}^*)^l$. Let
$$
\widehat{U}_\varepsilon(\g)=U_\varepsilon(\g)\otimes_{Z_0}\widehat{Z}_0,~ \index[not]{U@$\widehat{U}_\varepsilon(\g)$}
\widehat{Z}_\varepsilon=Z_\varepsilon\otimes_{Z_0}\widehat{Z}_0. \index[not]{Z@$\widehat{Z}_\varepsilon$}
$$

\begin{proposition}{\bf (\cite{DK}, Propositions 3.4, 3.5; \cite{DKP1} Proposition 6.1 and Theorem 6.6)}\label{qcoadj}

(i) On specializing to $q^{\frac{1}{d\bar{r}^2}}=\varepsilon^{\frac{1}{d\bar{r}^2}}$, (\ref{qder}) induces a well--defined derivation $\underline{x}_i^\pm$ of ${U}_\varepsilon(\g)$.

(ii) The series
$$
\exp(t\underline{x}_i^\pm)=\sum_{k=0}^\infty \frac{t^k}{k!}(\underline{x}_i^\pm)^k
$$
converges for all $t\in \mathbb{C}$ to a well--defined automorphism of the algebra $\widehat{U}_\varepsilon(\g)$.

(iii) Let $\mathcal{G}$ \index[not]{G@$\mathcal{G}$} be the group of automorphisms generated by the one--parameter groups $\exp(t\underline{x}_i^\pm)$, $i=1,\ldots, l$. The action of $\mathcal{G}$ on $\widehat{U}_\varepsilon(\g)$ preserves the subalgebras $\widehat{Z}_\varepsilon$ and $\widehat{Z}_0$, and hence $\mathcal{G}$ acts by holomorphic automorphisms on the complex algebraic varieties ${\rm Spec}(Z_\varepsilon)$ and ${\rm Spec}(Z_0)$.

(iv) Let $\mathcal{O}$ be a conjugacy class in $G$. The intersection $\mathcal{O}^0:=\mathcal{O}\cap G^0$ is a smooth connected variety, and the variety $\pi^{-1}(\mathcal{O}^0)$ is a $\mathcal{G}$--orbit in ${\rm Spec}(Z_0)$.
\end{proposition}

Given a homomorphism $\eta:{Z}_0 \rightarrow \mathbb{C}$, let
$$
U_\eta(\g)={U}_\varepsilon(\g)/I_\eta, \index[not]{U@$U_\eta(\g)$}
$$
where $I_\eta$ \index[not]{I@$I_\eta$} is the ideal in ${U}_\varepsilon(\g)$ generated by elements $z-\eta(z)$, $z\in Z_0$. By part (iii) of Proposition \ref{ue} $U_\eta(\g)$ is an algebra of dimension $\bar{m}^{{\rm dim}~\g}$ with linear basis the set of monomials $(X^+)^{\bf r}L^{\bf p}(X^-)^{\bf m}$ for which $0\leq r_k,m_k,p_i<\bar{m}$ for $i=1,\ldots ,l$, $k=1,\ldots ,D$.

If $\widetilde{g}\in \mathcal{G}$ then for any $\eta\in {\rm Spec}(Z_0)$ we have $\widetilde{g}\eta\in {\rm Spec}(Z_0)$ by part (iii) of Proposition \ref{qcoadj}, and by part (ii) of the same proposition $\widetilde{g}$ induces an isomorphism of algebras,
\begin{equation}\label{qcoadje}
\widetilde{g}:U_\eta(\g) \rightarrow U_{\widetilde{g}\eta}(\g).
\end{equation}

Since on every irreducible representation of ${U}_\varepsilon(\g)$ the subalgebra $Z_0$ of the center $Z_\varepsilon$ acts by a character $\eta:{Z}_0 \rightarrow \mathbb{C}$, every irreducible representation of ${U}_\varepsilon(\g)$ is a representation of some algebra $U_\eta(\g)$ for a unique $\eta$. This reduces the study of irreducible representations of ${U}_\varepsilon(\g)$ to the study of representations of finite--dimensional algebras $U_\eta(\g)$. Moreover, taking into account isomorphisms (\ref{qcoadje}) it suffices to consider a representative in each isomorphism class of these algebras under the isomorphisms induced by the action of the elements of the group $\mathcal{G}$ on $\widehat{U}_\varepsilon(\g)$.


\section{Whittaker vectors in modules over quantum groups at roots of unity}\label{Whittr1}

\pagestyle{myheadings}
\markboth{CHAPTER~\thechapter. APPLICATION TO QUANTUM GROUPS AT ROOTS OF UNITY}{\thesection.~WHITTAKER VECTORS AT ROOTS OF UNITY}

\setcounter{equation}{0}
\setcounter{theorem}{0}

It turns out that any finite--dimensional representation $V$ of ${U}_\eta(\g)$ can be equipped with another action of a subalgebra $U_{\eta_1}(\m_-)$ of a small quantum group which is a root of unity ``truncated'' version of the algebra $U_{\mathcal{A}}^s(\m_-)$ for an appropriate $s$ depending on $\eta$. The new action is compatible with the original action of ${U}_\eta(\g)$ in a certain equivariant way, and the dimension of $U_{\eta_1}(\m_-)$ is equal to $b_\eta:=m^{\frac{1}{2}{\rm dim}~\mathcal{O}_{\pi\eta}}=m^{{\rm dim}~\m_-}$, \index[not]{b@$b_\eta$} where $\mathcal{O}_{\pi\eta}$ is the conjugacy class of $\pi\eta \in G$. The existence of the second action is crucial for the proof of the De Concini--Kac--Procesi conjecture and for the study of other properties of finite--dimensional representations of ${U}_\eta(\g)$. 

In this section we define the algebras $U_{\eta_1}(\m_-)$ and their actions on finite--dimensional representations of ${U}_\eta(\g)$. These definitions are related to the notion of Whittaker vectors for finite--dimensional ${U}_\eta(\g)$--modules which are defined using root of unity versions of the characters $\chi_q^s$. We start by reminding the definitions of these characters.  

From now on we assume that the system of positive roots $\Delta_+$ and its normal ordering are associated to $s\in W$ as in Definition \ref{circorddef}. 
Firstly, let us observe that $U_\varepsilon^{s}({\frak m}_-)$ can be regarded as a subalgebra in $U_\varepsilon({\frak g})$. Therefore for every  character $\eta: Z_0 \rightarrow \mathbb{C}$ one can define the corresponding subalgebra in $U_\eta({\frak g})$ generated by $f_\alpha$, $\alpha\in \Delta_{\m_+}$. We denote this subalgebra by $U_\eta({\frak m}_-)$. \index[not]{U@$U_\eta(\m_-)$} By part (ii) of Proposition \ref{uq1z} we have ${\rm dim}~U_\eta({\frak m}_-)=\bar{m}^{{\rm dim}~\m_-}$.

In order to define analogues of characters $\chi_q^s$ for quantum groups at roots of unity we shall need some properties of the finite dimensional algebras $U_\eta({\frak g})$ and $U_\eta({\frak m}_-)$ and auxiliary results on non--zero irreducible representations of the algebra $U_\eta({\frak m}_-)$.

Observe that by Proposition \ref{pord} for any two roots $\alpha, \beta\in \Delta_{\m_+}$ such that $\alpha<\beta$ the sum $\alpha+\beta$ can not be represented as a linear combination $\sum_{k=1}^jm_k\gamma_{i_k}$, where $m_k\in \mathbb{N}$ and $\alpha<\gamma_{i_1}<\ldots <\gamma_{i_j}<\beta$, and hence from commutation relations (\ref{cmrelf}) one can deduce that
\begin{equation}\label{eqJ}
f_{\alpha}f_{\beta}-\varepsilon^{\left\langle \alpha,\beta\right\rangle+\bar{n}d\left\langle \frac{1+s}{1-s}P_{{\h'}^*}\alpha,\beta\right\rangle}f_{\beta}f_{\alpha}=\sum_{m_1,\ldots, m_k\in\mathbb{N}}C'(m_1,\ldots,m_k)
f_{\zeta_k}^{m_k}f_{\zeta_{k-1}}^{m_{k-1}}\ldots f_{\zeta_1}^{m_1},
\end{equation}
where $\alpha<\zeta_1<\ldots<\zeta_k<\beta$, $[\alpha,\beta]=\{\alpha,\zeta_1,\ldots,\zeta_k,\beta\}$ as a set, and $m_i>0$ for at least one root $\zeta_i\in \Theta:=\{\alpha \in \Delta_{\m_+}:\alpha \not\in \{\gamma_1, \ldots ,\gamma_{l'}\}\}$. \index[not]{T@$\Theta$} Therefore the elements $f_\beta\in U_{\eta}({\frak m}_-)$, $\beta\in \Theta$ generate an ideal $\mathcal{J}$ \index[not]{J@$\mathcal{J}$} in $U_{\eta}({\frak m}_-)$. Thus from Lemma \ref{segmPBWs} (vi), Remark \ref{segmPBWsrev} together with part (ii) of Proposition \ref{uq1z} and commutation relations (\ref{eqJ}) we obtain the following lemma.

\begin{lemma}\label{pbwr1}
Let $\beta_1<\beta_2<\ldots<\beta_c$ be the roots in the segment $\Delta_{\m_+}$. Then the elements
\begin{equation}\label{basJJ}
f_{\beta_c}^{m_c}f_{\beta_{c-1}}^{m_{c-1}}\ldots f_{\beta_1}^{m_1}
\end{equation}
or the elements
$$
f_{\beta_1}^{m_1}f_{\beta_1}^{m_1}\ldots f_{\beta_c}^{m_c}
$$
with $m_i\in \Bbb{N}$, $m_i<\bar{m}$, $i=1,\ldots, c$ form a linear basis of $U_{\eta}({\frak m}_-)$, and elements (\ref{basJJ}) with $m_i\in \Bbb{N}$, $m_i<\bar{m}$, $i=1,\ldots, c$ and $m_j>0$ for at least one $\beta_j\in \Theta$ form a linear basis of $\mathcal{J}$.
\end{lemma}

The following Lemma describes the Jacobson radical of $U_{\eta}({\frak m}_-)$.
\begin{lemma}\label{Jacob}
Let $\eta$ be an element of ${\rm Spec}(Z_0)$. Assume that $a_i:=\eta(f_{\gamma_i}^{\bar{m}})\neq 0$ for $i=1,\ldots ,l'$ and that $\eta(f_\beta^{\bar{m}})=0$ for $\beta\in \Delta_{\m_+}$, $\beta\not\in \{\gamma_1, \ldots ,\gamma_{l'}\}$, and hence $f_{\gamma_i}^{\bar{m}}=\eta(f_{\gamma_i}^{\bar{m}})=a_i\neq 0$ in $U_{\eta}({\frak m}_-)$ for $i=1,\ldots ,l'$ and $f_\beta^{\bar{m}}=0$ in $U_{\eta}({\frak m}_-)$ for $\beta\in \Delta_{\m_+}$, $\beta\not\in \{\gamma_1, \ldots ,\gamma_{l'}\}$.
Then the ideal $\mathcal{J}$ is the Jacobson radical \index{radical!Jacobson} of $U_{\eta}({\frak m}_-)$ and $U_{\eta}({\frak m}_-)/\mathcal{J}$ is isomorphic to the truncated polynomial algebra \index{algebra!truncated polynomial}
$$
\mathbb{C}[f_{\gamma_1},\ldots,f_{\gamma_{l'}}]/\{f_{\gamma_i}^{\bar{m}}=a_i\}_{i=1,\ldots ,l'}.
$$
\end{lemma}

\begin{proof}
First we show that $\mathcal{J}$ is nilpotent.


Let $i$ be the largest number such that $m_j=0$ for $j>i$ in (\ref{basJJ}) and $m_i\neq 0$. Then we define the {\it truncated degree} \index{degree!truncated} of $f_{\beta_c}^{m_c}f_{\beta_{c-1}}^{m_{c-1}}\ldots f_{\beta_1}^{m_1}$ by
$$
deg'(f_{\beta_c}^{m_c}f_{\beta_{c-1}}^{m_{c-1}}\ldots f_{\beta_1}^{m_1})=(m_i,i)\in \{1,\ldots,\bar{m}-1\}\times \{1,\ldots,c\}. \index[not]{d@$deg'$}
$$

Equip $\{1,\ldots,\bar{m}-1\}\times \{1,\ldots,c\}$ with the order such that $(k,i)<(k',j)$ if $j>i$ or $j=i$ and $k'>k$.

For any given $(k,i)\in \{1,\ldots,\bar{m}-1\}\times \{1,\ldots,c\}$ denote by $(U_{\eta}({\frak m}_-))_{(k,i)}$ \index[not]{U@$(U_{\eta}(\m_-))_{(k,i)}$}
 the linear span of the elements $f_{\beta_c}^{m_c}f_{\beta_{c-1}}^{m_{c-1}}\ldots f_{\beta_1}^{m_1}$ with $deg'(f_{\beta_c}^{m_c}f_{\beta_{c-1}}^{m_{c-1}}\ldots f_{\beta_1}^{m_1})\leq (k,i)$ and define $\mathcal{J}_{(k,i)}=\mathcal{J}\cap (U_{\eta}({\frak m}_-))_{(k,i)}$. \index[not]{J@$\mathcal{J}_{(k,i)}$}
We also have $(U_{\eta}({\frak m}_-))_{(k,i)}\subset (U_{\eta}({\frak m}_-))_{(k',j)}$ and
$\mathcal{J}_{(k,i)}\subset \mathcal{J}_{(k',j)}$ if $(k,i)<(k',j)$, and $\mathcal{J}_{(\bar{m}-1,c)}=\mathcal{J}$. Note that for the first few $i$ spaces $\mathcal{J}_{(k,i)}$ may be trivial, and these are all possibilities when those spaces can be trivial.

We shall prove that $\mathcal{J}$ is nilpotent by induction over the order in $\{1,\ldots,\bar{m}-1\}\times \{1,\ldots, c\}$.
Let $(k,i)$ be minimal possible such that $\mathcal{J}_{(k,i)}$ is not trivial. Then we must have $k=1$. If $y\in \mathcal{J}_{(1,i)}$ then $y$ must be of the form
\begin{equation}\label{yf}
y=f_{\beta}v,
\end{equation}
where $v$ is a linear combination of elements of the form $f_{\beta_{i_1}}^{m_{1}}\ldots f_{\beta_{i_r}}^{m_r}$ for ${\beta_{i_1}},\ldots ,{\beta_{i_r}}\in \{\gamma_{1},\ldots , \gamma_{\widetilde{l}}\}$, $\beta>\beta_{i_1}$, and $\beta$ is the first root from the set $\Theta$ greater than $\gamma_1$ in the normal ordering of $\Delta_+$  associated to $s$. Here it is assumed that $f_{\beta_{i_1}}^{m_{1}}\ldots f_{\beta_{i_r}}^{m_r}=1$ if the set $\{\gamma_{1},\ldots , \gamma_{\widetilde{l}}\}$ is empty.

Now equation (\ref{eqJ}) implies that for any $f_{\beta_{i_j}}$ which appears in the expression for $v$ one has
\begin{equation}\label{rr1}
f_{\beta}f_{\beta_{i_j}}-\varepsilon^{\left\langle \beta,\beta_{i_j}\right\rangle+\bar{n}d\left\langle \frac{1+s}{1-s}P_{{\h'}^*}\beta,\beta_{i_j}\right\rangle}f_{\beta_{i_j}}f_{\beta}\in \mathcal{J}_{(\bar{m}-1,i-1)}=0
\end{equation}
as by our choice of $i$ $\mathcal{J}_{(\bar{m}-1,i-1)}=0$.

Formula (\ref{rr1}) implies that the product of $\bar{m}$ elements of type (\ref{yf}) can be represented in the form
$$
f_{\beta}^{\bar{m}}v',
$$
where $v'$ is of the same form as $v$. Since $f_{\beta}^{\bar{m}}=0$ we deduce that $\mathcal{J}_{(1,i)}^{\bar{m}}=0$.

Now assume that $\mathcal{J}_{(k,i)}^m=0$ for some $m>0$. Let $(k',i')$ be the smallest element of $\{1,\ldots,\bar{m}-1\}\times \{1,\ldots,c\}$ which satisfies $(k,i)<(k',i')$. Then by Lemma \ref{segmPBWs} (vi), by Remark \ref{segmPBWsrev}, by Proposition \ref{uq1z}, and by (\ref{eqJ}), any element of  $\mathcal{J}_{(k',i')}$ is of the form 
\begin{equation}\label{fut}
f_{\beta_{i'}}u+u', 
\end{equation}
where $u'\in \mathcal{J}_{(k,i)}$ and if $\beta_{i'}\in \Theta$ then $u\in (U_{\eta}({\frak m}_-))_{(k,i)}$; if $\beta_{i'}\not\in \Theta$ then $u\in \mathcal{J}_{(k,i)}$.

Now equation (\ref{eqJ}) implies that for any $u\in (U_{\eta}({\frak m}_-))_{(k,i)}$ of the form (\ref{basJJ}) one has
\begin{equation}\label{ufc}
uf_{\beta_{i'}}=c_{u,i'}f_{\beta_{i'}}u+w,
\end{equation}
where $c_{u,i'}$ is a non--zero constant depending on $u$ and $i'$, and $w\in \mathcal{J}_{(k,i)}$.
By formula (\ref{ufc}) the product of $\bar{m}$ elements $f_{\beta_{i'}}u_p+u'_p$, $p=1,\ldots,\bar{m}$ of type (\ref{fut}) can be represented in the form
\begin{equation}\label{sj}
\sum_{j=0}^{\bar{m}}f_{\beta_{i'}}^jc_j,
\end{equation}
where $c_j\in \mathcal{J}_{(k,i)}$ for $j=0,\ldots ,\bar{m}-1$ and if $\beta_{i'}\in \Theta$ then $c_{\bar{m}}\in (U_{\eta}({\frak m}_-))_{(k,i)}$; if $\beta_{i'}\not\in \Theta$ then $c_{\bar{m}}\in \mathcal{J}_{(k,i)}$. In the former case $f_{\beta_{i'}}^{\bar{m}}=0$, and the last term in sum (\ref{sj}) is zero; in the latter case $f_{\beta_{i'}}^{\bar{m}}=\eta(f_{\beta_{i'}}^{\bar{m}})\neq 0$, and the last term in sum (\ref{sj}) is from $\mathcal{J}_{(k,i)}$. So we can combine it with the term corresponding to $j=0$. In both cases sum (\ref{sj}) takes the form
\begin{equation}\label{st}
\sum_{j=0}^{\bar{m}-1}f_{\beta_{i'}}^jc_j',
\end{equation}
where $c_j'\in \mathcal{J}_{(k,i)}$. By (\ref{ufc}) the product of $m$ sums of type (\ref{st}) is of the form
$$
\sum_{j=0}^{(\bar{m}-1)m}f_{\beta_{i'}}^jc_j'',
$$
where each $c_j''$ is a linear combination of elements from $\mathcal{J}_{(k,i)}^m$.
By our assumption $\mathcal{J}_{(k,i)}^m=0$, and hence the product of any $\bar{m}m$ elements of $\mathcal{J}_{(k',i')}$ is zero. This justifies the induction step and proves that $\mathcal{J}_{(\bar{m}-1,c)}=\mathcal{J}$ is nilpotent.
Hence $\mathcal{J}$ is contained in the Jacobson radical of $U_{\eta}({\frak m}_-)$.

Using commutation relations (\ref{eqJ}) we also have (see the proof of Proposition \ref{Qdefr})
$$
f_{\gamma_i}f_{\gamma_j} - f_{\gamma_j}f_{\gamma_i}\in \mathcal{J}.
$$
Therefore the quotient algebra $U_{\eta}({\frak m}_-)/\mathcal{J}$ is isomorphic to the truncated polynomial algebra $$\mathbb{C}[f_{\gamma_1},\ldots,f_{\gamma_{l'}}]/\{f_{\gamma_i}^{\bar{m}}=a_i\}_{ i=1,\ldots ,l'}$$ which is semisimple. \index{algebra!semisimple} Therefore $\mathcal{J}$ coincides with the Jacobson radical of $U_{\eta}({\frak m}_-)$.

\end{proof}

In Proposition \ref{Qdefr} we constructed some characters of the algebra $U_\varepsilon^{s}({\frak m}_-)$. Now we show that the algebra $U_\eta({\frak m}_-)$ has a finite number of irreducible representations which are one--dimensional, and all these representations can be obtained from each other by twisting with the help of automorphisms of $U_\eta({\frak m}_-)$.

\begin{proposition}\label{irrepod}
Let $\eta$ be an element of ${\rm Spec}(Z_0)$. Assume that $a_i:=\eta(f_{\gamma_i}^{\bar{m}})\neq 0$ for $i=1,\ldots ,l'$ and that $\eta(f_\beta^{\bar{m}})=0$ for $\beta\in \Delta_{\m_+}$, $\beta\not\in \{\gamma_1, \ldots ,\gamma_{l'}\}$, and hence $f_{\gamma_i}^{\bar{m}}=\eta(f_{\gamma_i}^{\bar{m}})=a_i\neq 0$ in $U_{\eta}({\frak m}_-)$ for $i=1,\ldots ,l'$ and $f_\beta^{\bar{m}}=0$ in $U_{\eta}({\frak m}_-)$ for $\beta\in \Delta_{\m_+}$, $\beta\not\in \{\gamma_1, \ldots ,\gamma_{l'}\}$. Then all non--zero irreducible representations of the algebra $U_\eta({\frak m}_-)$ are one--dimensional and have the form
\begin{equation}\label{chichar}
\chi(f_\beta)=\left\{ \begin{array}{ll} 0 & \beta \not \in \{\gamma_1, \ldots, \gamma_{l'}\} \\ \bar{c}_i & \beta=\gamma_i,~i=1,\ldots ,l' \index[not]{c@$\bar{c}_i$}
\end{array}
\right  .,
\end{equation}
where complex numbers $\bar{c}_i$ satisfy the conditions $\bar{c}_i^{\bar{m}}=a_i$, $i=1,\ldots ,l'$.
Moreover, all non--zero irreducible representations of $U_\eta({\frak m}_-)$ can be obtained from each other by twisting with the help of automorphisms of $U_\eta({\frak m}_-)$.
\end{proposition}

\begin{proof}
Let $V$ be a non--zero finite--dimensional irreducible $U_{\eta}({\frak m}_-)$--module.
By Corollary 54.13 in \cite{CR} elements of the ideal $\mathcal{J}\subset U_{\eta}({\frak m}_-)$ act by zero transformations on $V$. Hence $V$ is in fact an irreducible representation of the algebra $U_{\eta}({\frak m}_-)/\mathcal{J}$ which is isomorphic to the truncated polynomial algebra $$\mathbb{C}[f_{\gamma_1},\ldots,f_{\gamma_{l'}}]/\{f_{\gamma_i}^{\bar{m}}=a_i\}_{ i=1,\ldots ,l'}.$$ The last algebra is commutative and all its complex irreducible representations are
one--dimensional. Therefore $V$ is one--dimensional, and if $v$ is a non--zero element of $V$ then $f_{\gamma_i}v=\bar{c}_iv$, for some $\bar{c}_i\in \Bbb{C}$, $i=1,\ldots ,l'$.

Note that $f_{\gamma_i}^{\bar{m}}=\eta(f_{\gamma_i}^{\bar{m}})=a_i\neq 0$, $i=1,\ldots ,l'$ in $U_{\eta}({\frak m}_-)$ and hence $\bar{c}_i^{\bar{m}}=a_i\neq 0$, $i=1,\ldots ,l'$. In particular, the elements $f_{\gamma_i}$ act on $V$ by semisimple automorphisms.

If we denote by $\chi:U_{\eta}({\frak m}_-) \rightarrow \mathbb{C}$ the character of $U_{\eta}({\frak m}_-)$ such that
$$
\chi(f_\beta)=\left\{ \begin{array}{ll} 0 & \beta \not \in \{\gamma_1, \ldots, \gamma_{l'}\} \\ \bar{c}_i & \beta=\gamma_i, ~i=1,\ldots ,l' \index[not]{x@$\chi$}
\end{array}
\right  .
$$
and by $\mathbb{C}_{\chi}$ \index[not]{C@$\mathbb{C}_\chi$} the corresponding one--dimensional representation of $U_{\eta}({\frak m}_-)$ then we have $V=\mathbb{C}_{\chi}$. This proves the first claim of the proposition.

Now we have to show that the representations $\mathbb{C}_{\chi}$ for different characters $\chi$ are obtained from each other by twisting with the help of automorphisms of $U_{\eta}({\frak m}_-)$.

Since $\bar{c}_i^{\bar{m}}=a_i$, $i=1,\ldots ,l'$ there are only finitely many possible characters $\chi$ corresponding to the given $\eta$ in the statement of this proposition.
If $\chi$ and $\chi'$ are two such characters, $\chi(f_{\gamma_i})=\bar{c}_i$, $i=1,\ldots ,l'$ and $\chi'(f_{\gamma_i})=\bar{c}_i'$, $i=1,\ldots ,l'$ then the relations $\bar{c}_i^{\bar{m}}=(\bar{c}_i')^{\bar{m}}=a_i$, $i=1,\ldots ,l'$ imply that
$\bar{c}_i'=\varepsilon^{m_i}\bar{c}_i$, $0\leq m_i \leq \bar{m}-1$, $m_i\in \mathbb{N}$, $i=1,\ldots ,l'$.

Now observe that for any $h\in \h$ the map defined by $f_\alpha \mapsto \varepsilon^{\alpha(h)}f_\alpha$, $\alpha \in \Delta_{\m_+}$ is an automorphism of the algebra $U_\varepsilon^s({\frak m}_-)$ generated by the elements $f_\alpha$, $\alpha \in \Delta_{\m_+}$ with defining relations (\ref{eqJ}). Here the principal branch of the analytic function $\varepsilon^z$, $z\in \mathbb{C}$ is used to define $\varepsilon^{\alpha(h)}$, so that $\varepsilon^{\alpha(h)}\varepsilon^{\beta(h)}=\varepsilon^{(\alpha+\beta)(h)}$ for any $\alpha,\beta\in \Delta_{\m_+}$. 

If in addition $\varepsilon^{\bar{m}\gamma_i(h)}=1$, $i=1,\ldots ,l'$ the above defined automorphism gives rise to an automorphism $\varsigma$ of $U_{\eta}({\frak m}_-)$. Indeed, in this case $(\varepsilon^{\gamma_i(h)}f_{\gamma_i})^{\bar{m}}=f_{\gamma_i}^{\bar{m}}$, $i=1,\ldots ,l'$ and all the remaining defining relations $f_{\gamma_i}^{\bar{m}}=\eta(f_{\gamma_i}^{\bar{m}})=a_i\neq 0$, $i=1,\ldots ,l'$, $f_{\beta}^{\bar{m}}=\eta(f_{\beta}^{\bar{m}})=0$, $\beta \in \Delta_{\m_+}$, $\beta \not \in \{\gamma_1, \ldots, \gamma_{l'}\}$ of the algebra $U_{\eta}({\frak m}_-)$ are preserved by the action of $\varsigma$.

Now fix $h\in \h$ such that $\gamma_i(h)=m_i$, $i=1,\ldots ,l'$. Obviously we have $\varepsilon^{\bar{m}m_i}=1$, $i=1,\ldots ,l'$. We claim that the representation $\mathbb{C}_{\chi}$ twisted by the corresponding automorphism  $\varsigma$ coincides with $\mathbb{C}_{\chi'}$. Indeed, we obtain
$$
\chi(\varsigma f_{\gamma_i})=\chi(\varepsilon^{m_i}f_{\gamma_i})=\varepsilon^{m_i}\bar{c}_i=\bar{c}_i',~i=1,\ldots ,l'.
$$
This completes the proof of the proposition.

\end{proof}

Now we can define the notion of Whittaker vectors for quantum groups at roots of unity.
Let $V$ be a $U_{\eta}({\frak g})$--module, where $\eta$ is an element of ${\rm Spec}(Z_0)$ such that $a_i=\eta(f_{\gamma_i}^{\bar{m}})\neq 0$ for $i=1,\ldots ,l'$ and $\eta(f_\beta^{\bar{m}})=0$ for $\beta\in \Delta_{\m_+}$, $\beta\not\in \{\gamma_1, \ldots ,\gamma_{l'}\}$. Let $\chi:U_\eta({\frak m}_-)\rightarrow \mathbb{C}$ be a character defined in the  Proposition \ref{irrepod}, $\mathbb{C}_{\chi}$ the corresponding one--dimensional $U_\eta({\frak m}_-)$--module. Then the space $V_\chi={\rm Hom}_{U_\eta({\frak m}_-)}(\mathbb{C}_{\chi},V)$ \index[not]{V@$V_\chi$} is called {\it the space of Whittaker vectors} \index{space!of Whittaker vectors!for quantum groups at roots of unity} of $V$. Elements of $V_\chi$ are called {\it Whittaker vectors}. \index{Whittaker vector!for quantum groups at roots of unity}

Now we describe the space of Whittaker vectors in terms of a nilpotent action of the unital subalgebra $U_{\eta_1}({\frak m}_-)$ generated by $f_\alpha$, $\alpha\in \Delta_{\m_+}$ in the small quantum group \index{quantum group!small} $U_{\eta_1}({\frak g}):=U_\varepsilon^s(\g)/I_{\eta_1}$ \index[not]{U@$U_{\eta_1}({\frak g})$} corresponding to the trivial central character $\eta_1$ \index[not]{i@$\eta_1$} such that $\widetilde{\pi}(\eta_1)=1\in G_0^*$. so that $\eta_1(x_\alpha^\pm)=0$, $\alpha\in \Delta_+$, $\eta_1(l_i)=1$, $i=1,\ldots, l$. 

One can also equip the algebra $U_\varepsilon^{s}({\frak m}_-)$ with a character given by formula (\ref{charqfe}), where the numbers $\bar{c}_i$ are the same as in the definition of the character $\chi$. We denote this character by the same letter, $\chi: U_\varepsilon^{s}({\frak m}_-)\rightarrow \mathbb{C}$.

Note that any $U_{\eta}({\frak g})$--module $V$ can be regarded as a $U_\varepsilon(\g)$--module and a $U_\varepsilon^{s}(\g)$--module assuming that the ideal $I_\eta$ acts on $V$ in the trivial way.

Recall that by (\ref{comults}) $U_\varepsilon^{s}({\frak m}_-)\subset U_\varepsilon^{s}(\g)$ is a right coideal in $U_\varepsilon^{s}(\g)$ and observe that $\Delta_{s}:U_\varepsilon^{s}({\frak m}_-)\rightarrow U_\varepsilon^{s}(\g)\otimes U_\varepsilon^{s}({\frak m}_-)$ is a homomorphism of algebras. Composing it with the tensor product $S_{s}\otimes \chi$ of the anti--homomorphism $S_{s}$ and of the character $\chi$, which can be regarded as an anti--homomorphism as well, one can define an anti--homomorphism, $U_\varepsilon^{s}({\frak m}_-)\rightarrow U_\varepsilon^{s}(\g)$, $x\mapsto S_{s}(x^1)\chi(x^2)$, $\Delta_{s}(x)=x^1\otimes x^2$, $x\in U_\varepsilon^{s}({\frak m}_-)$.

Using this anti--homomorphism one can introduce a right $U_\varepsilon^{s}({\frak m}_-)$--action on $V$ which we call {\it the adjoint action} \index{action!adjoint!on a $U_{\eta}({\frak g})$--module} and denote it by ${\rm Ad}_s$. \index[not]{A@${\rm Ad}_s$} It is given by the formula
\begin{equation}\label{AD}
{\rm Ad}_sx v=S_{s}(x^1)\chi(x^2)v, x\in U_\varepsilon^{s}({\frak m}_-), v\in V,
\end{equation}
where $\Delta_{s}(x)=x^1\otimes x^2$. 

Note that using the Sweedler notation for the comultiplication, $(\Delta_{s}\otimes id \otimes id) (\Delta_{s}\otimes id) \Delta_{s}(x)=x^1\otimes x^2\otimes x^3 \otimes x^4$, the coassociativity of the comultiplication and the definition of the antipode we have for any $x \in U_\varepsilon^{s}({\frak m}_-)$, $y\in U_\varepsilon^{s}(\g)$, $v\in V$ (compare with the proof of Lemma 2.2 in \cite{JL})
\begin{equation}\label{adme}
{\rm Ad}_sx(yv)=S_{s}(x^1)\chi(x^2)yv=S_{s}(x^1)yx^2S_{s}(x^3)\chi(x^4)v={\rm Ad}_sx^1(y){\rm Ad}_sx^2(v).
\end{equation}
 
We shall need the following formula for the action of the comultiplication on the quantum root vectors $f_\beta$, $\beta=\sum_{i=1}^lm_i\alpha_i\in \Delta_+$, $m_i\in \mathbb{N}$,
\begin{eqnarray}\label{cm1}
\Delta_{s}(f_{\beta})=\prod_{i=1}^lK_i^{m_i}\prod_{i,j=1}^lL_j^{-\frac{\bar{n}d}{d_j}\left\langle \frac{1+s}{1-s}P_{{\h'}^*}\alpha_i,\alpha_j\right\rangle m_i}\otimes f_\beta+f_\beta \otimes 1+ \\
+\sum_i y_i\otimes x_i,~x_i\in U_{<\beta}, y_i\in U_{>\beta}U_\varepsilon^s(\h), \nonumber
\end{eqnarray}
where $U_{<\beta}$ \index[not]{U@$U_{<\beta}$} is the subalgebra (without unit) in $U_\varepsilon^s(\n_-)$ generated by ${f}_{\alpha}$, $\alpha<\beta$ and $U_{>\beta}$ \index[not]{U@$U_{>\beta}$} is the subalgebra (without unit) in $U_\varepsilon^s(\n_-)$ generated by ${f}_{\alpha}$, $\alpha>\beta$. Formula (\ref{cm1}) is a straightforward consequence of (\ref{comults}). 

Similarly to the Proposition in Section 5.6 in \cite{DKP1} we infer that $Z_0$ is a Hopf subalgebra in $U_\varepsilon^{s}(\g)$. Namely,
\begin{eqnarray}
\Delta_{s}(f_i^{\bar{m}})=K_i^{\bar{m}}\prod_{j=1}^lL_j^{-\bar{m}\frac{\bar{n}d}{d_j}\left\langle \frac{1+s}{1-s}P_{{\h'}^*}\alpha_i,\alpha_j\right\rangle}\otimes f_i^{\bar{m}}+f_i^{\bar{m}} \otimes 1, \nonumber \\
\Delta_{s}(e_i^{\bar{m}})=e_i^{\bar{m}}\otimes K_i^{-\bar{m}}+\prod_{j=1}^lL_j^{\bar{m}\frac{\bar{n}d}{d_j}\left\langle \frac{1+s}{1-s}P_{{\h'}^*}\alpha_i,\alpha_j\right\rangle}\otimes e_i^{\bar{m}}, \label{Z0H} \\
\Delta_{s}(L_i^{\bar{m}})=L_i^{\bar{m}}\otimes L_i^{\bar{m}}. \nonumber
\end{eqnarray}
Therefore recalling that by the definition of $\chi$
for $x\in U_\varepsilon^{s}({\frak m}_-)\cap Z_0$ one has $\chi(x)=\eta(x)$ we deduce
\begin{equation}\label{adchi}
{\rm Ad}_sx v=S_{s}(x_1)\chi(x_2)v=\eta(S_{s}(x_1)x_2)v=\varepsilon_{s}(x)v, v\in V,
\end{equation}
where $\varepsilon_{s}$ is the counit of $U_\varepsilon^{s}(\g)$. 
Note that by the definition of the ideal $I_{\eta_1}$ the ideal $U_\varepsilon^{s}({\frak m}_-)\cap I_{\eta_1}\subset U_\varepsilon^{s}({\frak m}_-)$  is generated by the elements $f_\alpha^{\bar{m}}$, $\alpha\in \Delta_{\m_+}$ and  $\varepsilon_{s}(f_\alpha^{\bar{m}})=0$ for $\alpha\in \Delta_{\m_+}$ by the definition of $\varepsilon_{s}$.
Hence by (\ref{adchi}) the adjoint action of $U_\varepsilon^{s}({\frak m}_-)$ on $V$ induces an action of the subalgebra $U_{\eta_1}({\frak m}_-)$ of the small quantum group $U_{\eta_1}({\frak g})$. We call this action {\it the adjoint action} as well. 

Note that the small quantum group $U_{\eta_1}({\frak g})$ is a Hopf algebra with the comultiplication inherited from $U_\varepsilon^{s}({\frak g})$.

The space of Whittaker vectors $V_\chi$ can be characterized in terms of the adjoint action as follows.
\begin{lemma}\label{LW}
The space of Whittaker vectors $V_\chi$ coincides with the space of $U_{\eta_1}({\frak m}_-)$--invariants for the adjoint action on $V$,
\begin{equation}\label{ADW}
V_\chi=\{ v\in V: {\rm Ad}_sx(v)=\varepsilon_{s}(x)v ~\forall x\in U_{\eta_1}({\frak m}_-)\}.
\end{equation}
\end{lemma}

\begin{proof}
Indeed, as in (\ref{comults}), denote by $G_\beta^{-1}$ the factor $\prod_{i=1}^lK_i^{m_i}\prod_{i,j=1}^lL_j^{-\frac{\bar{n}d}{d_j}\left\langle \frac{1+s}{1-s}P_{{\h'}^*}\alpha_i,\alpha_j \right\rangle m_i}$ which appears in (\ref{cm1}), $G_\beta^{-1}=\prod_{i=1}^lK_i^{m_i}\prod_{i,j=1}^lL_j^{-\frac{\bar{n}d}{d_j}\left\langle \frac{1+s}{1-s}P_{{\h'}^*}\alpha_i,\alpha_j \right\rangle m_i}$. Then similarly to (\ref{Ss}) we obtain
\begin{equation}\label{cm2}
S_{s}(f_\beta)=-S_{s}({G_\beta}^{-1})f_\beta
-\sum_i S_{s}(y_i)x_i=-{G_\beta}f_\beta
-\sum_i S_{s}(y_i)x_i.
\end{equation}

Now for $\beta\in \Delta_{\m_+}$, formulas (\ref{cm1}), (\ref{cm2}) and definition (\ref{AD}) of the adjoint action imply
\begin{eqnarray}\label{cm3}
{\rm Ad}_sf_{\beta}v= G_\beta \chi(f_\beta)v-G_\beta f_\beta v - 
\sum_i S_{s}(y_i) x_i v +\sum_i S_{s}(y_i) \chi(x_i)v= \nonumber \\
=G_\beta (\chi(f_\beta)-f_\beta) v + 
\sum_i S_{s}(y_i) (\chi(x_i)-x_i)v, ~x_i\in U_{<\beta}, y_i\in U_{>\beta}U_\varepsilon^{s}(\h).
\end{eqnarray}

If $v\in V_\chi$ we immediately obtain from (\ref{cm3}) that ${\rm Ad}_sf_{\beta}v=0$ for any $\beta\in \Delta_{\m_+}$, i.e. $v$ belongs to the right hand side of (\ref{ADW}).

Conversely, suppose that $v$ belongs to the right hand side of (\ref{ADW}). We shall show that $xv=\chi(x)v$ for any $x\in U_\varepsilon^{s}({\frak m}_-)$. 

The subalgebra of $U_\varepsilon^{s}(\n_-)$ with unit generated by $U_{<\beta_k}$, $\beta_k\in \Delta_+=\{\beta_1,\ldots,\beta_D\}$, where $\beta_1<\beta_2<\ldots<\beta_D$ are the roots of $\Delta_+$ ordered as in Definition \ref{circorddef}, is the subalgebra $U_\varepsilon^s([-\beta_1,-\beta_k])\subset U_\varepsilon^{s}(\n_-)$ \index[not]{U@$U_\varepsilon^s([-\beta_1,-\beta_k])$} generated by $f_\alpha$, $\alpha\in [\beta_1,\beta_k]$. 

We proceed by induction over the subalgebras $U_\varepsilon^s([-\beta_1,-\beta_k])\subset U_\varepsilon^s(\g)$, $k=1,\ldots c$, where as before $\beta_1<\ldots <\beta_c$ are the roots in the normally ordered segment $\Delta_{\m_+}$. 

Observe that $\beta_1$ is a simple root. Therefore we deduce from (\ref{cm3}) for $\beta=\beta_1$
$$
{\rm Ad}_sf_{\beta_1}v= 
G_{\beta_1}(\chi(f_{\beta_1})-f_{\beta_1}) v=0.
$$
Since $G_{\beta_1}$ acts on $V$ by an invertible transformation this implies $(\chi(f_{\beta_1})-f_{\beta_1}) v=0$, and hence $xv=\chi(x)v$ for any $x\in U_\varepsilon^s([-\beta_1,-\beta_1])$ as $U_\varepsilon^s([-\beta_1,-\beta_1])$ is generated by $f_{\beta_1}$.

Now assume that for some $1<k\leq c$ one has $xv=\chi(x)v$ for any $x\in U_\varepsilon^s([-\beta_1,-\beta_{k-1}])$. Then by (\ref{cm3})
$$
{\rm Ad}_sf_{\beta_k}v= 
G_{\beta_k}(\chi(f_{\beta_k})-f_{\beta_k}) v=0.
$$
Since $G_{\beta_k}$ acts as an invertible endomorphism of $V$, this implies 
\begin{equation}\label{cm4}
(\chi(f_{\beta_k})-f_{\beta_k}) v=0.
\end{equation}

Now observe that by Lemma \ref{segmPBWs} (vi) and by Remark \ref{segmPBWsrev} any element $y\in U_\varepsilon^s([-\beta_1,-\beta_k])$ can be uniquely represented in the form $y=\sum_{i=0}^r f_{\beta_k}^iz_i$ for some $r\in \mathbb{N}$ and
$z_i\in U_\varepsilon^s([-\beta_1,-\beta_{k-1}])$ for $i=0,\ldots, r$. Therefore by (\ref{cm4}) and by the induction assumption
$$
yv=\sum_{i=0}^r f_{\beta_k}^iz_iv=\sum_{i=0}^r f_{\beta_k}^i\chi(z_i)v=\sum_{i=0}^r\chi(f_{\beta_k}^i)\chi(z_i)v=\chi(\sum_{i=0}^rf_{\beta_k}^iz_i)v=\chi(y)v.
$$
i.e. $yv=\chi(y)v$ for any $y\in U_\varepsilon^s([-\beta_1,-\beta_k])$. This establishes the induction step and completes the proof.

\end{proof}

The following proposition is an analogue of the Engel theorem for quantum groups at roots of unity.

\begin{proposition}\label{Whitt}
Let $\eta$ be an element of ${\rm Spec}(Z_0)$. Assume that $a_i:=\eta(f_{\gamma_i}^{\bar{m}})\neq 0$ for $i=1,\ldots ,l'$ and that  $\eta(f_\beta^{\bar{m}})=0$ for $\beta\in \Delta_{\m_+}$, $\beta\not\in \{\gamma_1, \ldots ,\gamma_{l'}\}$, and hence $f_{\gamma_i}^{\bar{m}}=\eta(f_{\gamma_i}^{\bar{m}})=a_i\neq 0$ in $U_{\eta}({\frak m}_-)$ for $i=1,\ldots ,l'$ and $f_\beta^{\bar{m}}=0$ in $U_{\eta}({\frak m}_-)$ for $\beta\in \Delta_{\m_+}$, $\beta\not\in \{\gamma_1, \ldots ,\gamma_{l'}\}$. Let $\chi:U_\eta({\frak m}_-)\rightarrow \mathbb{C}$ be any character defined in Proposition \ref{irrepod}.
Then any non--zero finite--dimensional $U_{\eta}({\frak g})$--module contains a non--zero Whittaker vector.
\end{proposition}

\begin{proof}

We begin the proof with the following lemma.
\begin{lemma}\label{J1Jac}
The augmentation ideal $\mathcal{J}^1$ \index[not]{J@$\mathcal{J}^1$} of $U_{\eta_1}({\frak m}_-)$ coincides with its Jacobson radical  which is nilpotent. 
\end{lemma}

\begin{proof}
The proof of this fact is similar to that of Lemma \ref{Jacob}, and we shall keep the notation used in that proof. 

We define $\mathcal{J}_{(k,i)}^1=\mathcal{J}^1\cap (U_{\eta_1}({\frak m}_-))_{(k,i)}$, \index[not]{J@$\mathcal{J}_{(k,i)}^1$} so that
$\mathcal{J}_{(k,i)}^1\subset \mathcal{J}_{(k',j)}^1$ if $(k,i)<(k',j)$, and $\mathcal{J}_{(\bar{m}-1,c)}^1=\mathcal{J}^1$. 

We shall prove that $\mathcal{J}^1$ is nilpotent by induction over the order in $\{1,\ldots,\bar{m}-1\}\times \{1,\ldots,c\}$.
Note that $(k,i)=(1,1)$ is minimal possible such that $\mathcal{J}_{(k,i)}$ is not trivial. If $y\in \mathcal{J}_{(1,1)}$ then $y$ must be of the form
\begin{equation}\label{yf1}
y=af_{\beta_1}, a\in \mathbb{C}.
\end{equation}

The product of ${\bar{m}}$ elements of type (\ref{yf1}) is equal to zero,
$$
f_{\beta_1}^{\bar{m}}a^{\bar{m}}=0,
$$
as $f_{\beta_1}^{\bar{m}}=0$ in $U_{\eta_1}({\frak m}_-)$. We deduce that $(\mathcal{J}_{(1,1)}^1)^{\bar{m}}=0$.

Now assume that $(\mathcal{J}_{(k,i)}^1)^m=0$ for some $m>0$. Let $(k',i')$ be the smallest element of $\{1,\ldots,\bar{m}-1\}\times \{1,\ldots,c\}$ which satisfies $(k,i)<(k',i')$. Then by Lemma \ref{segmPBWs} (vi), by Remark \ref{segmPBWsrev} and by Proposition \ref{uq1z} any element of  $\mathcal{J}_{(k',i')}^1$ is of the form 
\begin{equation}\label{Jki1}
f_{\beta_{i'}}u+u',
\end{equation} 
where $u'\in \mathcal{J}_{(k,i)}^1$ and $u\in (U_{\eta_1}({\frak m}_-))_{(k,i)}$.

Now equation (\ref{eqJ}) implies that for any $u\in (U_{\eta_1}({\frak m}_-))_{(k,i)}$ of the form (\ref{basJJ}) one has
\begin{equation}\label{ufc1}
uf_{\beta_{i'}}=c_{u,i'}f_{\beta_{i'}}u+w,
\end{equation}
where $c_{u,i'}$ is a non--zero constant depending on $u$ and $i'$, and $w\in \mathcal{J}_{(k,i)}^1$.

By formula (\ref{ufc1}) the product of ${\bar{m}}$ elements $f_{\beta_{i'}}u_p+u'_p$, $p=1,\ldots,m$ of type (\ref{Jki1}) can be represented in the form
\begin{equation}\label{sj1}
\sum_{j=0}^{\bar{m}}f_{\beta_{i'}}^jc_j,
\end{equation}
where $c_j\in \mathcal{J}_{(k,i)}^1$ for $j=0,\ldots ,{\bar{m}}-1$ and $c_{\bar{m}}\in (U_{\eta_1}({\frak m}_-))_{(k,i)}$. Since $f_{\beta_{i'}}^{\bar{m}}=0$ the last term in sum (\ref{sj1}) is zero. So sum (\ref{sj1}) takes the form
\begin{equation}\label{st1}
\sum_{j=0}^{m-1}f_{\beta_{i'}}^jc_j',
\end{equation}
where $c_j'\in \mathcal{J}_{(k,i)}^1$. By (\ref{ufc1}) the product of $m$ sums of type (\ref{st1}) is of the form
$$
\sum_{j=0}^{(\bar{m}-1)m}f_{\beta_{i'}}^jc_j'',
$$
where each $c_j''$ is a linear combination of elements from $(\mathcal{J}_{(k,i)}^1)^m$.
By our assumption $(\mathcal{J}_{(k,i)}^1)^m=0$, and hence the product of any $\bar{m}m$ elements of $\mathcal{J}_{(k',i')}^1$ is zero. This justifies the induction step and proves that $\mathcal{J}_{(\bar{m}-1,c)}^1=\mathcal{J}^1$ is nilpotent.
Hence $\mathcal{J}^1$ is contained in the Jacobson radical of $U_{\eta_1}({\frak m}_-)$.

The quotient algebra $U_{\eta_1}({\frak m}_-)/\mathcal{J}^1$ is isomorphic to $\mathbb{C}$. Therefore $\mathcal{J}^1$ coincides with the Jacobson radical of $U_{\eta_1}({\frak m}_-)$.

\end{proof}

Now let $V$ be a finite--dimensional $U_{\eta}({\frak g})$--module. Then $V$ is also a finite--dimensional  $U_{\eta_1}({\frak m}_-)$--module with respect to the adjoint action. Thus $V$ must contain a non--trivial irreducible $U_{\eta_1}({\frak m}_-)$--submodule with respect to the adjoint action on which the Jacobson radical $\mathcal{J}^1$ must act trivially. From (\ref{ADW}) it follows that this non--trivial irreducible submodule consists of Whittaker vectors. This completes the proof of the proposition.

\end{proof}

Now we show that for any $\eta \in {\rm Spec}(Z_0)$ subalgebras and characters which appear in Propositions \ref{irrepod}, \ref{Whitt} and in Lemma \ref{LW} indeed exist. Moreover, we shall see that to each $\eta \in {\rm Spec}(Z_0)$ one can associate a subalgebra $U_{\eta}^{\widetilde{g}}({\frak m}_-)$ of the same type as $U_\eta$ the dimension of which is equal to $\bar{m}^{\frac{1}{2}{\rm dim}~\mathcal{O}_{\pi\eta}}$, where $\mathcal{O}_{\pi\eta}$ \index[not]{O@$\mathcal{O}_{\pi\eta}$} is the conjugacy class of $\pi\eta \in G$.

\begin{proposition}\label{DKPconj}
Let
$$
G=\bigcup_{\mathcal{C}\in C(W)}G_\mathcal{C}.
$$
be the Lusztig partition of $G$, $\eta \in {\rm Spec}(Z_0)$ be an element  such that $\pi\eta\in G_\mathcal{C}$, $\mathcal{C}\in C(W)$ and $s^{-1}\in \mathcal{C}$.
Let $\Delta_+^{s}$ be a system of positive roots defined for $s$ in Theorem \ref{mainth}, $\Delta_+$ the corresponding system of positive roots associated to $s$ in Definition \ref{circorddef}, $d=2d'$, where $d'$ is defined in Proposition \ref{d'}. Assume that $\bar{m}$ and $d$ are coprime. Then the following statements are true.

(i) There exists a quantum coadjoint transformation $\widetilde{g}$ such that $\xi :=\widetilde{g}\eta$ satisfies $\xi(f_{\gamma_i}^{\bar{m}})\neq 0$ for $i=1,\ldots ,l'$ and $\xi(f_\beta^{\bar{m}})=0$  for $\beta\in \Delta_{\m_+}$, $\beta\not\in \{\gamma_1, \ldots ,\gamma_{l'}\}$, where $f_\alpha \in U_{\xi}({\frak m}_-)$ are the generators of the corresponding algebra $U_{\xi}({\frak m}_-)\subset U_{\xi}(\g)$. 

(ii) Let $\chi:U_{\xi}({\frak m}_-)\rightarrow \mathbb{C}$ be any character defined in Proposition \ref{irrepod}.
Then any finite--dimensional $U_{\xi}({\frak g})$--module contains a non--zero Whittaker vector with respect to the subalgebra $U_{\xi}({\frak m}_-)$ and the character $\chi$, and any $U_{\eta}({\frak g})$--module contains a non--zero Whittaker vector with respect to the subalgebra $U_{\eta}^{\widetilde{g}}({\frak m}_-):=\tilde{g}^{-1}(U_{\xi}({\frak m}_-))$ \index[not]{U@$U_\eta^{\widetilde{g}}(\m_-)$} and the character $\chi^{\widetilde{g}}$ \index[not]{x@$\chi^{\widetilde{g}}$} given by the composition of $\chi$ and $\widetilde{g}$, $\chi^{\widetilde{g}}:=\chi \circ \widetilde{g}:U_{\eta}^{\widetilde{g}}({\frak m}_-) \rightarrow \mathbb{C}$.

(iii) ${\rm dim}~U_{\xi}({\frak m}_-)={\rm dim}~U_{\eta}^{\widetilde{g}}({\frak m}_-)=\bar{m}^{\frac{1}{2}{\rm dim}~\mathcal{O}_{\pi\eta}}=\bar{m}^{\frac{1}{2}{\rm dim}~\mathcal{O}_{\pi\xi}}=\bar{m}^{{\rm dim}~\m_-}$, where $\mathcal{O}_{\pi\eta}$ is the conjugacy class of $\pi\eta \in G$, and $\mathcal{O}_{\pi\xi}$ is the conjugacy class of $\pi\xi \in G$.
\end{proposition}

\begin{proof}

(i) First observe that the system of positive roots $\Delta_+^s$ satisfies the conditions of Theorem \ref{mainth} when $s$ is replaced with $s^{-1}$. Indeed, in the case of classical root systems its definition only depends on the spectral decomposition of $\h$ under the action of $s$ which is the same as the spectral decomposition of $\h$ under the action of $s^{-1}$. In the case of exceptional root systems one has to note in addition that obviously ${\rm dim}~\Sigma_s={\rm dim}~\Sigma_{s^{-1}}$, and hence all properties of $\Delta_+^s$ used in the proof of Theorem \ref{mainth} are satisfied if $s$ is replaced with $s^{-1}$ in the proof. 

Let $\overline{N}_{s^{-1}}=\{n\in \overline{N}^s:s^{-1}ns\in N^s\}$. \index[not]{N@$\overline{N}_{s^{-1}}$}
Applying Theorem \ref{mainth} to $s^{-1}$ and to the system of positive rootrs $\Delta_+^s$ and swapping the roles of $N^s$ and of $\overline{N}^s$ we deduce  that all conjugacy classes in the stratum $G_\mathcal{C}$ intersect the set $s^{-1}H^0\overline{N}_{s^{-1}}$ which is a subset of the transversal slice $\overline{\Sigma}_{s^{-1}}:=s^{-1}Z^s\overline{N}_{s^{-1}}$ \index[not]{S@$\overline{\Sigma}_{s^{-1}}$} to the set of conjugacy classes in $G$. Note that $\overline{N}_{s^{-1}}$ is not a subgroup in $N_+$. But by Lemma \ref{decNr} every element $n\in \overline{N}_{s^{-1}}$ can be uniquely  factorized as follows $n=n^+n^-$, where $n^\pm\in \overline{N}_{s^{-1}}\cap N_\pm\subset N_\pm$. Therefore every element $s^{-1}{h^0}n\in s^{-1}H^0\overline{N}_{s^{-1}}$ can be represented as follows $s^{-1}{h^0}n=s^{-1}{h^0}n^+n^-$, and conjugating by $n^-$ we obtain that $s^{-1}{h^0}n$ is conjugate to $n^-s^{-1}{h^0}n^+$.

Note that $\overline{N}_{s^{-1}}\cap N_-$ is generated by the one--parameter subgroups corresponding to the roots from the set $\Delta_{s^{-1}}^s\cap \Delta_-\subset \Delta_{s^1}^s$. Since $s^1$ is an involution, from Proposition \ref{pord} (iii) it follows that (see also Figure 4)
$$
s(\Delta_{s^1}^s)=s^1(\Delta_+^s\setminus(\Delta_{s^1}^s\cup \Delta_{s^2}^s\cup \Delta_0))\subset \Delta_+^s\setminus\Delta_{s^1}^s \subset \Delta_+,
$$
and hence
$$
s(\overline{N}_{s^{-1}}\cap N_-)s^{-1}\subset N_+.
$$
Taking into account that $H^0$ normalizes $N_+$ we have $n^-s^{-1}{h^0}n^+=s^{-1}{h^0}{h^0}^{-1}sn^-s^{-1}{h^0}n^+=s^{-1}{h^0}m'$, $m'={h^0}^{-1}sn^-s^{-1}{h^0}n^+\in N_+$. We deduce that all conjugacy classes in the stratum $G_\mathcal{C}$ intersect the set $s^{-1}H^0N_+$.
 
Recall that by (\ref{levdec}) with the roles of $N_+$ and $N_-$ swapped we have $\bar{v}=n'sm$ for some $n'\in N_+, m\in M_+\subset N_+$, where
$\bar{v}\in G$ is an element of the form
$$
\bar{v}=\prod_{i=1}^{l'}{\exp}[\bar{t}_{i} X_{-\gamma_i}], \index[not]{v@$\bar{v}$}
$$
$\bar{t}_i\in \mathbb{C}$ are non--zero constants depending on the choice of the representative $s\in N_G(H)$ and the product over roots is taken in the normal order (\ref{NO}) associated to $s$. We deduce that 
$s^{-1}=m\bar{v}^{-1}n'$. 

Now let $s^{-1}{h^0}m',{h^0}\in H^0,m'\in N_+$ be an element of $s^{-1}H^0N_+$. Using the  expression $s^{-1}=m\bar{v}^{-1}n'$ for $s^{-1}$ we can write $s^{-1}{h^0}m'=m\bar{v}^{-1}n'{h^0}m'$. Conjugating this element by $m^{-1}$ and recalling that $H^0$ normalizes $N_+$ we infer that $s^{-1}{h^0}m'$ is conjugate to 
$$
\bar{v}^{-1}n'{h^0}m'm=\bar{v}^{-1}{h^0}n=q({h^0}^{\frac12}n,{h^0}^{-\frac12}\bar{v}),
$$
where $n={h^0}^{-1}n'{h^0}m'm\in N_+$, ${h^0}^{\frac12}\in H^0$ is any element such that ${h^0}^{\frac12}{h^0}^{\frac12}={h^0}$, $q:G_0^*\to G^0$ is defined immediately after Proposition \ref{LO}. We deduce that all conjugacy classes in the stratum $G_\mathcal{C}$ intersect the set $\bar{v}^{-1}H^0N_+$.

By part (iv) of Proposition \ref{qcoadj} we conclude that if $\eta \in {\rm Spec}(Z_0)$ satisfies $\pi\eta\in G_\mathcal{C}$ then there is a quantum coadjoint transformation $\widetilde{g}$ such that $\widetilde{\pi}(\widetilde{g}\eta)=({h^0}^{\frac12}n,{h^0}^{-\frac12}\bar{v})$ for some $n\in N_+$, ${h^0}^{\frac12}\in H^0$.

Denote $\xi=\widetilde{g}\eta$. From the definition of the map $\widetilde{\pi}$ and of the element $\bar{v}$ it follows that 
$$
\exp(\xi(x_{\beta_D}^-)X_{-\beta_D})\exp(\xi(x_{\beta_{D-1}}^-)X_{-\beta_{D-1}})\ldots \exp(\xi(x_{\beta_1}^-)X_{-\beta_1})=\bar{v}^{-1}
$$
which implies $\xi((X_{\gamma_i}^-)^{\bar{m}})=-\frac{\bar{t}_i}{(\varepsilon_{\gamma_i}-\varepsilon_{\gamma_i}^{-1})^{\bar{m}}}\neq 0$ for $i=1,\ldots ,l'$ and that $\xi((X_\beta^-)^{\bar{m}})=0$ for $\beta\in \Delta_{\m_+}$, $\beta\not\in \{\gamma_1, \ldots ,\gamma_{l'}\}$. 

By the definition of the elements $f_\beta$ with $\beta=\sum_{i=1}^lm_i\alpha_i$ we have $f_\beta=\prod_{i,j=1}^lL_j^{m_in_{ij}}X_\beta^-$. Therefore the commutation relations between elements $L_j$ and $X_\beta^-$ imply that $f_\beta^{\bar{m}}=c_\beta'\prod_{i,j=1}^lL_j^{\bar{m}m_in_{ij}}(X_\beta^-)^{\bar{m}}$, where $c_\beta'$ are non--zero constants, and hence $\xi(f_\beta^{\bar{m}})=c_\beta'\prod_{i,j=1}^l\xi(L_j^{\bar{m}})^{m_in_{ij}}\xi((X_\beta^-)^{\bar{m}})$. 

Since $\xi(L_j^{\bar{m}})\neq 0$ for $j=1,\ldots ,l$, $\xi((X_{\gamma_i}^-)^{\bar{m}})=-\frac{\bar{t}_i}{(\varepsilon_{\gamma_i}-\varepsilon_{\gamma_i}^{-1})^{\bar{m}}}\neq 0$ for $i=1,\ldots ,l'$ and $\xi((X_\beta^-)^{\bar{m}})=0$ for $\beta\in \Delta_{\m_+}$, $\beta\not\in \{\gamma_1, \ldots ,\gamma_{l'}\}$ we deduce 
$\xi(f_{\gamma_i}^{\bar{m}})\neq 0$ for $i=1,\ldots ,l'$ and $\xi(f_\beta^{\bar{m}})=0$  for $\beta\in \Delta_{\m_+}$, $\beta\not\in \{\gamma_1, \ldots ,\gamma_{l'}\}$. This proves part (i).

(ii) Note that $\xi$ satisfies the condition of Propositions \ref{irrepod} and \ref{Whitt}. Let $U_{\widetilde{g}\eta}({\frak m}_-)=U_{\xi}({\frak m}_-)$ be the corresponding subalgebra in $U_{\xi}(\g)$. Then part (ii) follows from Propositions \ref{irrepod} and \ref{Whitt}.

(iii) Note that by Theorem \ref{mainth} for any $g\in G_\mathcal{C}$ we have
$$
{\rm dim}~Z_G(g)={\rm dim}~\overline{\Sigma}_{s^{-1}},
$$
where $Z_G(g)$ is the centralizer of $g$ in $G$. 

By the definition of $\overline{\Sigma}_{s^{-1}}$ we also have ${\rm dim}~\overline{\Sigma}_{s^{-1}}=l(s)+2D_0+{\rm dim}~{\h'}^\perp$. 
Observe also that ${\rm dim}~G=2D+{\rm dim}~\h$ and ${\rm dim}~\h-{\rm dim}~{\h'}^\perp={\rm dim}~{\h'}=l'$,
and hence from (\ref{dimm}) we deduce that ${\rm dim}~\m_-=D-D_0-\frac12(l(s)-l')=\frac{1}{2}({\rm dim}~G-{\rm dim}~\overline{\Sigma}_{s^{-1}})=\frac{1}{2}{\rm dim}~\mathcal{O}_g$, and by Lemma \ref{pbwr1} ${\rm dim}~U_{\eta}^{\widetilde{g}}({\frak m}_-)={\rm dim}~U_{\widetilde{g}\eta}({\frak m}_-)=\bar{m}^{{\rm dim}~\m_-}=\bar{m}^{\frac{1}{2}{\rm dim}~\mathcal{O}_g}$, where $\mathcal{O}_g$ is the conjugacy class of $g \in G_\mathcal{C}$.

In particular, ${\rm dim}~U_{\eta}^{\widetilde{g}}({\frak m}_-)=\bar{m}^{\frac{1}{2}{\rm dim}~\mathcal{O}_{\pi\eta}}$, where $\mathcal{O}_{\pi\eta}$ is the conjugacy class of $\pi\eta \in G_\mathcal{C}$.

\end{proof}


\section{Skryabin equivalence for quantum groups at roots of unity and the proof of De Concini--Kac--Procesi conjecture}

\pagestyle{myheadings}
\markboth{CHAPTER~\thechapter. APPLICATION TO QUANTUM GROUPS AT ROOTS OF UNITY}{\thesection.~THE PROOF OF DE CONCINI--KAC--PROCESI CONJECTURE}

\setcounter{equation}{0}
\setcounter{theorem}{0}

In this section we shall study the $U_{\eta_1}(\m_-)$--action on finite--dimensional $U_\eta(\g)$--modules introduced in the previous section. We shall show that each such module is $U_{\eta_1}(\m_-)$--cofree. Taking into account that ${\rm dim}~U_{\eta_1}({\frak m}_-)=\bar{m}^{\frac{1}{2}{\rm dim}~\mathcal{O}_{\pi\eta}}$ this will imply the De Concini--Kac--Procesi conjecture. 

The main observation with the help of which we shall prove these statements is that the structure of finite--dimensional $U_\eta(\g)$--modules is similar to that of a root of unity analogue $Q_\chi$ of the module $Q_{\mathcal{B}}^{loc}$, and the results of Propositions \ref{BasP} and \ref{Bpbas} can be specialized to $q^{\frac{1}{d{\bar{r}}^2}}=\varepsilon^{\frac{1}{d{\bar{r}}^2}}$ and transferred to $Q_\chi$ and, more generally, to any finite--dimensional $U_\eta(\g)$--module using  root of unity analogues of the homomorphism $\phi$. In particular, using  the specializations $B^\varepsilon_{m_1\ldots m_c}$ of the elements $B_{m_1\ldots m_c}$ one can construct $U_{\eta_1}(\m_-)$--cofree bases in finite--dimensional $U_\eta(\g)$--modules.
 
To realize this program for any given $\eta\in {\rm Spec}(Z_0)$ we assume that $\Delta_+^{s}$ is a system of positive roots defined for $s$ in Theorem \ref{mainth}, $\Delta_+$ the corresponding system of positive roots associated to $s$, $d=2d'$, where $d'$ is defined in Proposition \ref{d'}. Assume also that $\bar{m}$ and $d$ are coprime.

Fix a quantum coadjoint transformation $\widetilde{g}\in \mathcal{G}$  as in Theorem \ref{DKPconj} and denote $\xi=\widetilde{g}\eta\in {\rm Spec}(Z_0)$. Since according to (\ref{qcoadje}) $\widetilde{g}$ gives rise to an isomorphism of the algebras $U_\eta(\g)$ and $U_\xi(\g)$ it suffices to consider the case of the algebra $U_\xi(\g)$. 

Our first objective is to obtain root of unity analogues of Proposition \ref{kerphi} and Lemma \ref{phiext}. We start by introducing the notions required for the formulations of these statements. For indeterminate $q$ these notions were introduced in Section \ref{Qlev}. 

Let $\chi$ be a character of $U_{\xi}({\frak m}_-)$ defined in Proposition \ref{irrepod}, $\mathbb{C}_\chi$ the corresponding representation of $U_{\xi}({\frak m}_-)$. Denote by $Q_\chi$ the induced left $U_{\xi}(\g)$--module, $Q_\chi=U_{\xi}(\g)\otimes_{U_{\xi}({\frak m}_-)}\mathbb{C}_\chi$. \index[not]{Q@$Q_\chi$} $Q_\chi$ can also be naturally regarded as a $U_{\varepsilon}^{s}(\g)$--module via the natural projection $U_{\varepsilon}^{s}(\g)\simeq U_{\varepsilon}(\g)\rightarrow U_{\xi}(\g)$. 

As in the previous chapter, let $\mathbb{C}_{\varepsilon}^s[G]$, $\mathbb{C}_{\varepsilon}^s[G_*]$, ${\mathbb{C}}_{\varepsilon}^{s,loc}[G]$, \index[not]{C@$\mathbb{C}_\varepsilon^{s,loc}[G]$} ${\mathbb{C}}_{\varepsilon}^{s,loc}[G_*]$, \index[not]{C@$\mathbb{C}_\varepsilon^{s,loc}[G_*]$} ${\mathbb{C}}_{\varepsilon}^s[G^*]$, \index[not]{C@$\mathbb{C}_\varepsilon^s[G^*]$} $B_{m_1\ldots m_c}^\varepsilon$, \index[not]{B@$B_{m_1\ldots m_c}^\varepsilon$} ${\mathbb{C}}_{11}^{loc}[G]_\varepsilon$ \index[not]{C@$\mathbb{C}_{11}^{loc}[G]_\varepsilon$} be the natural specializations at $q^{\frac{1}{d\bar{r}^2}}=\varepsilon^{\frac{1}{d\bar{r}^2}}$ of $\mathbb{C}_{\mathcal{B}}^s[G]$, $\mathbb{C}_{\mathcal{B}}^s[G_*]$, $\mathbb{C}_{\mathcal{B}}^{s,loc}[G]$, $\mathbb{C}_{\mathcal{B}}^{s,loc}[G_*]$, $\mathbb{C}_{\mathcal{B}}^s[G^*]$, $B_{m_1\ldots m_c}$, ${\mathbb{C}}_{11}^{loc}[G]$, respectively. 

We define the twisted adjoint action of $U_{\varepsilon}^{s,res}(\g)$ on $\mathbb{C}_{\varepsilon}^s[G]$ by specializing formula (\ref{ado}) at $q^{\frac{1}{d\bar{r}^2}}=\varepsilon^{\frac{1}{d\bar{r}^2}}$,
\begin{equation}\label{ado1}
({\rm Ad}_s^0 xf)(w)=f(\omega_0 S_s^{-1}({\rm Ad}_s'x(S_s\omega_0 w)))=f((\omega_0 S_s^{-1})(x^1)w\omega_0 (x^2)), f\in \mathbb{C}_{\varepsilon}^s[G], x,w \in U_{\varepsilon}^{s,res}(\g). \index[not]{A@${\rm Ad}_s^0$}
\end{equation}
Specializing isomorphism (\ref{gg**}) at $q^{\frac{1}{d{\bar{r}}^2}}=\varepsilon^{\frac{1}{d{\bar{r}}^2}}$ we obtain by Proposition \ref{locfin} (i) a $U_{\varepsilon}^{s,res}(\g)$--module homomorphism 
\begin{equation}\label{gg**1}
\varphi_\varepsilon:\mathbb{C}_{\varepsilon}^s[G]\rightarrow \mathbb{C}_{\varepsilon}^s[G_*], \index[not]{f@$\varphi_\varepsilon$}
\end{equation}
where $\mathbb{C}_{\varepsilon}^s[G_*]$ is equipped with the adjoint action ${\rm Ad}_s$ of $U_{\varepsilon}^{s,res}(\g)$.

Note that by (\ref{G*AG*}) $\mathbb{C}_{\varepsilon}^s[G_*]\subset U_{\varepsilon}^s(\g)$, and hence one can consider $\varphi_\varepsilon$ as a map $\varphi_\varepsilon:\mathbb{C}_{\varepsilon}^s[G]\rightarrow U_{\varepsilon}^s(\g)$. 

Observe now that the subalgebra in $U_{\varepsilon}^{s,res}(\g)$ generated by $f_\beta$, $\beta\in \Delta_{\m_+}$ is isomorphic to $U_{\eta_1}(\m_-)$. Therefore restricting the adjoint action ${\rm Ad}_s$ and action (\ref{ado1}) of $U_{\varepsilon}^{s,res}(\g)$ to $U_{\eta_1}(\m_-)$ we deduce that $\varphi_\varepsilon:\mathbb{C}_{\varepsilon}^s[G]\rightarrow U_{\varepsilon}^s(\g)$ is a homomorphism of $U_{\eta_1}(\m_-)$--modules, where the action of $U_{\eta_1}(\m_-)$ on $U_{\varepsilon}^s(\g)$ is defined by the same formula as ${\rm Ad}_s$. Note that that this action is well defined as by (\ref{Z0H}) $Z_0\subset U_{\varepsilon}^s(\g)$ is a Hopf subalgebra, and hence the ideal $I_{\eta_1}$ generated by the augmentation ideal of $Z_0$ belongs to the kernel of the action ${\rm Ad}_s$ of $U_{\varepsilon}^s(\g)$ on itself, so that this action induces an action of $U_{\varepsilon}^s(\g)/I_{\eta_1}= U_{\eta_1}(\g)$.

Since the ideal $I_\xi$ is generated by central elements, composing the homomorphism $\varphi_\varepsilon:\mathbb{C}_{\varepsilon}^s[G]\rightarrow U_{\varepsilon}^s(\g)$ with the natural projection $\mathbb{C}_{\varepsilon}^s[G_*]\subset U_\varepsilon^s(\g) \rightarrow U_\xi(\g)$ we obtain a homomorphism of $U_{\eta_1}(\m_-)$--modules
\begin{equation}\label{gg***1'}
\mathbb{C}_{\varepsilon}^s[G]\rightarrow U_\xi(\g)=U_\varepsilon(\g)/I_\xi\simeq U_{\varepsilon}^s(\g)/I_\xi, 
\end{equation}
where $\mathbb{C}_{\varepsilon}^s[G]$ is equipped with the action (\ref{ado1}) of $U_{\varepsilon}^{s,res}(\g)$ and $U_\xi(\g)$ with the action induced by the adjoint action ${\rm Ad}_s$.

Next, because $\chi$ is a character of $U_{\xi}({\frak m}_-)$, one has $[U_{\xi}({\frak m}_-),{\rm Ker}\chi ]\subset {\rm Ker}\chi$. By (\ref{comults}) $\Delta_s(U_\varepsilon^s(\m_-))\subset  U_\varepsilon^s(\g)\otimes   U_\varepsilon^s(\m_-)$, and hence formula (\ref{ad}) for the adjoint action ${\rm Ad}_s$ implies that the left ideal $U_{\xi}(\g){\rm Ker}\chi$ is invariant under the action of $U_{\eta_1}(\m_-)$ on $U_\xi(\g)$ induced by the adjoint action ${\rm Ad}_s$.
Therefore this action gives rise to an action on $Q_\chi=U_{\xi}(\g)\otimes_{U_{\xi}({\frak m}_-)}\mathbb{C}_\chi\simeq U_{\xi}(\g)/U_{\xi}(\g){\rm Ker}\chi$.

Thus composing homomorphism (\ref{gg**1}) with the natural projection $\mathbb{C}_{\varepsilon}^s[G_*]\subset U_\varepsilon^s(\g) \rightarrow U_\xi(\g) \rightarrow U_{\xi}(\g)/U_{\xi}(\g){\rm Ker}\chi\simeq Q_\chi$ we obtain 
 a homomorphism of $U_{\eta_1}(\m_-)$--modules
\begin{equation}\label{gg***1}
\phi_\xi: \mathbb{C}_{\varepsilon}^s[G]\rightarrow Q_{\chi}, \phi_\xi(f)=\varphi_\varepsilon(f)1, \index[not]{f@$\phi_\xi$}
\end{equation}
where $\mathbb{C}_{\varepsilon}^s[G]$ is equipped with the restriction of action (\ref{ado1}) to $U_{\eta_1}(\m_-)$ and $Q_{\chi}$ with the action induced by the adjoint action ${\rm Ad}_s$ of $U_{\eta_1}(\m_-)$, and $1$ is the image of $1\in U_{\xi}(\g)$ in $Q_{\xi}$ under the projection $U_{\xi}(\g)\rightarrow U_{\xi}(\g)/U_{\xi}(\g){\rm Ker}\chi\simeq Q_{\xi}$.

Similarly, using the fact that any Whittaker vector $w\in V_\chi$ in any finite--dimensional $U_{\xi}(\g)$--module $V$ is annihilated by the left ideal $U_{\xi}(\g){\rm Ker}\chi$, one can define a $U_{\eta_1}(\m_-)$--module homomorphism $\phi_\xi^w: \mathbb{C}_{\varepsilon}^s[G] \rightarrow V$ by 
$$
\phi_\xi^w: \mathbb{C}_{\varepsilon}^s[G]\rightarrow V, \phi_\xi(f)=\varphi_\varepsilon(f)w, \index[not]{f@$\phi_\xi^w$}
$$
where $\mathbb{C}_{\varepsilon}^s[G]$ is equipped with the restriction of action (\ref{ado1}) to $U_{\eta_1}(\m_-)$ and $V$ with the adjoint action ${\rm Ad}_s$ of $U_{\eta_1}(\m_-)$.

Now we can state a root of unity version of Proposition \ref{Iq}.
\begin{proposition}\label{phixidef}
For any $\lambda \in P_+$, any finite--dimensional $U_{\xi}(\g)$--module $V$ and any Whittaker vector $w\in V_\chi$ the following statements are true.

(i) The specialization ${{I}^{11}_{\varepsilon}}\subset \mathbb{C}_{\varepsilon}^s[G]$ \index[not]{I@${I}^{11}_{\varepsilon}$} of the left ideal ${{I}^{11}_{\mathcal{B}}}\subset \mathbb{C}_{\mathcal{B}}^s[G]$ at $q^{\frac{1}{d{\bar{r}}^2}}=\varepsilon^{\frac{1}{d{\bar{r}}^2}}$ lies in the kernel of $\phi_\xi^w$.

(ii) Denote the specialization at $q^{\frac{1}{d{\bar{r}}^2}}=\varepsilon^{\frac{1}{d{\bar{r}}^2}}$ of the element $\Delta_\lambda^s\in \mathbb{C}_{\mathcal{B}}^s[G]$ by the same symbol. Then for any $f\in \mathbb{C}_{\varepsilon}^s[G]$ 
\begin{equation}\label{phixi}
\phi_\xi^w (f\otimes \Delta_\lambda^s)=\varphi_\varepsilon({\rm Ad}_s^0(\varepsilon^{-(\bar{n}d{1+s \over 1-s }P_{{\h'}}+id)\lambda^\vee})(f))\phi_\xi^w (\Delta_\lambda^s),
\end{equation}
where ${\rm Ad}_s^0(\varepsilon^{-(\bar{n}d{1+s \over 1-s }P_{{\h'}}+id)\lambda^\vee})(f)$ is the adjoint action of the element 
$$
\varepsilon^{-(\bar{n}d{1+s \over 1-s }P_{{\h'}}+id)\lambda^\vee}:=q^{-(\bar{n}d{1+s \over 1-s }P_{{\h'}}+id)\lambda^\vee}~~({\rm mod}~(q^{\frac{1}{d{\bar{r}}^2}}-\varepsilon^{\frac{1}{d{\bar{r}}^2}})), 
\varepsilon^{-(\bar{n}d{1+s \over 1-s }P_{{\h'}}+id)\lambda^\vee}\in \mathbb{C}_{\varepsilon}^s[G^*]
$$ 
on $f\in \mathbb{C}_{\varepsilon}^s[G]$.

(iii) One can define an action of an operator $\varepsilon^{(id+s^{-1})(id-\bar{n}dP_{\h'})\lambda^\vee}$ on the image of $\mathbb{C}_{\varepsilon}^s[G]$ in $V$ by the formula
\begin{equation}\label{opdef}
\varepsilon^{(id+s^{-1})(id-\bar{n}dP_{\h'})\lambda^\vee}\phi_\xi^w(f)=\varphi_\varepsilon({\rm Ad}_s^0(\varepsilon^{-(id+s^{-1})(id-\bar{n}dP_{\h'})\lambda^\vee})(f))\phi_\xi^w (\Delta_\lambda^s),
\end{equation}
where ${\rm Ad}_s^0(\varepsilon^{-(id+s^{-1})(id-\bar{n}dP_{\h'})\lambda^\vee})(f)$ is the adjoint action of the element 
$$
\varepsilon^{-(id+s^{-1})(id-\bar{n}dP_{\h'})\lambda^\vee}:=q^{-(id+s^{-1})(id-\bar{n}dP_{\h'})\lambda^\vee}~~({\rm mod}~(q^{\frac{1}{d{\bar{r}}^2}}-\varepsilon^{\frac{1}{d{\bar{r}}^2}})), 
\varepsilon^{-(id+s^{-1})(id-\bar{n}dP_{\h'})\lambda^\vee}\in \mathbb{C}_{\varepsilon}^s[G^*]
$$ 
on $f\in \mathbb{C}_{\varepsilon}^s[G]$.

Using this operator formula (\ref{phixi}) can be rewritten as follows
\begin{equation}\label{phixi1}
\phi_\xi^w(f\otimes \Delta_\lambda^s)=\varepsilon^{(id+s^{-1})(id-\bar{n}dP_{\h'})\lambda^\vee}\phi_\xi^w({\rm Ad}_s^0(\varepsilon^{(-\bar{n}d{1+s \over 1-s }s^{-1}P_{{\h'}}+s^{-1})\lambda^\vee})(f)).
\end{equation}

(iv) The operator $\varepsilon^{(id+s^{-1})(id-\bar{n}dP_{\h'})\lambda^\vee}$ is invertible. More precisely, for some $n_0\in \mathbb{N}$, $n_0>0$ the action of $(\varepsilon^{(id+s^{-1})(id-\bar{n}dP_{\h'})\lambda^\vee})^{n_0}$ coincides with the action of an invertible element $\varepsilon^{n_0(id+s^{-1})(id-\bar{n}dP_{\h'})\lambda^\vee}\in U_{\varepsilon}^s(\h)$ satisfying
$$
\varepsilon^{n_0(id+s^{-1})(id-\bar{n}dP_{\h'})\lambda^\vee}:=q^{n_0(id+s^{-1})(id-\bar{n}dP_{\h'})\lambda^\vee}=\prod_{i=1}^lL_i^{r_i}~~({\rm mod}~(q^{\frac{1}{d{\bar{r}}^2}}-\varepsilon^{\frac{1}{d{\bar{r}}^2}}))
$$ 
for some $r_i\in \mathbb{Z}$, $i=1,\ldots, l$, and hence in $U_\xi(\g)$ one has $c_\lambda^\varepsilon:=(\varepsilon^{n_0(id+s^{-1})(id-\bar{n}dP_{\h'})\lambda^\vee})^{\bar{m}}=\xi(\varepsilon^{n_0\bar{m}(id+s^{-1})(id-\bar{n}dP_{\h'})\lambda^\vee})\in \mathbb{C}^*$, and $(\varepsilon^{(id+s^{-1})(id-\bar{n}dP_{\h'})\lambda^\vee})^{n_0\bar{m}}=c_\lambda^\varepsilon I_V$, where $I_V$ is the identity operator on $V$.

(v) $\phi_\xi^w$ extends to a $U_{\eta_1}(\m_-)$--module homomorphism $\phi_\xi^w: \mathbb{C}_{\varepsilon}^{s,loc}[G] \rightarrow V$ such that
\begin{equation}\label{phixiloc}
\phi_\xi^w (f\otimes {\Delta_\lambda^s}^{-1})=
\varepsilon^{({(1-\bar{n}d)s^{-1}+(1+\bar{n}d)s-2 \over 1-s}P_{\h'}\lambda^\vee,\lambda^\vee)}\left(\varepsilon^{(id+s^{-1})(id-\bar{n}dP_{\h'})\lambda^\vee}\right)^{-1}\phi_\xi^w({\rm Ad}_s^0(\varepsilon^{(\bar{n}d{1+s \over 1-s }s^{-1}P_{{\h'}}-s^{-1})\lambda^\vee})(f)),
\end{equation}
and the specialization ${I^{11}_{\varepsilon}}^{loc}\subset \mathbb{C}_{\varepsilon}^{s,loc}[G]$ of ${I^{11}_{\mathcal{B}}}^{loc}\subset \mathbb{C}_{\mathcal{B}}^{s,loc}[G]$ \index[not]{I@${I^{11}_{\varepsilon}}^{loc}$} at $q^{\frac{1}{d{\bar{r}}^2}}=\varepsilon^{\frac{1}{d{\bar{r}}^2}}$ belongs to the kernel of this homomorphism, so 
$$
\phi_\xi^w:\mathbb{C}_{11}^{loc}[G]_\varepsilon \rightarrow V.
$$            
\end{proposition}

\begin{proof}
In order to prove this proposition one can apply Proposition \ref{Iq} for $\kappa=\bar{n}d$ and $q^{\frac{1}{d{\bar{r}}^2}}$ specialized to $\varepsilon^{\frac{1}{d{\bar{r}}^2}}$, and the appropriately modified arguments before Proposition \ref{kerphi} and from the proof of Lemma \ref{phiext}. 

Indeed, by the definition of $I_{\mathcal{B}}^{\bf k}$ with $\kappa=\bar{n}d$, $\bar{k}_i=\bar{c}_i$, $i=1,\ldots ,l'$ for $\bar{c}_i$, $i=1,\ldots ,l'$ used in the definition of $\chi$, the specialization of $I_{\mathcal{B}}^{\bf k}$ at $q^{\frac{1}{d{\bar{r}}^2}}=\varepsilon^{\frac{1}{d{\bar{r}}^2}}$ belongs to the annihilator of $w$, and hence by Proposition \ref{Iq} the specialization $I^{11}_{\varepsilon}\subset \mathbb{C}_{\varepsilon}^s[G]$ of the left ideal $I^{11}_{\mathcal{B}}\subset \mathbb{C}_{\mathcal{B}}^s[G]$ at $q^{\frac{1}{d{\bar{r}}^2}}=\varepsilon^{\frac{1}{d{\bar{r}}^2}}$ lies in the kernel of $\phi_\xi^w$. This proves part (i) 

Formulas (\ref{phixi}), (\ref{opdef}) and (\ref{phixi1}) in parts (ii) and (iii) are obtained by specializing the formulas in Proposition \ref{Iq} (ii) at $q^{\frac{1}{d{\bar{r}}^2}}=\varepsilon^{\frac{1}{d{\bar{r}}^2}}$ with $\kappa=\bar{n}d$.

The only essential difference is that the operator $\varepsilon^{(id+s^{-1})(id-\bar{n}dP_{\h'})\lambda^\vee}$ is invertible and for some $n_0\in \mathbb{N}$, $n_0>0$ the action of its $n_0$--power coincides with the action of the element 
$\varepsilon^{n_0(id+s^{-1})(id-\bar{n}dP_{\h'})\lambda^\vee}\in U_{\varepsilon}^s(\h)$.
This can be justified as follows. 

Recall that the elements $\gamma_i^\vee$, $i=1,\ldots ,l'$ form a linear basis of $\h'$. Let $\gamma_i^*$, $i=1,\ldots, l'$ be the basis of $\h'$ dual to $\gamma_i^\vee$, $i=1,\ldots, l'$ with respect to the restriction of the bilinear form $\left\langle ~\cdot~,~\cdot~\right\rangle$ to $\h'$. Since the numbers $\left\langle \gamma_i^\vee,\gamma_j^\vee\right\rangle$ are integer each element $\gamma_i^*$ has the form $\gamma_i^*=\sum_{j=1}^{l'}m_{ij}\gamma_j^\vee$, where $m_{ij}\in \mathbb{Q}$. Therefore $P_{\h'}\lambda^\vee=\sum_{p=1}^{l'}\left\langle \lambda^\vee,\gamma_p^\vee\right\rangle\gamma_p^*=\sum_{p,q=1}^{l'}\left\langle \lambda^\vee,\gamma_p^\vee\right\rangle m_{pq}\gamma_q^\vee$ belongs to the rational span of the set of the simple coweights $Y_i$, $i=1,\ldots ,l$, and hence $(id+s^{-1})(id-\bar{n}dP_{\h'})\lambda^\vee$ belongs to the rational span of the set of the simple coweights $Y_i$, $i=1,\ldots ,l$ as well. We conclude that there exists an integer $n_0\in \mathbb{N}$, $n_0>0$ such that $n_0(id+s^{-1})(id-\bar{n}dP_{\h'})\lambda^\vee$ belongs to the integer span of the set of the simple coweights $Y_i$, $i=1,\ldots ,l$, and hence $q^{n_0(id+s^{-1})(id-\bar{n}dP_{\h'})\lambda^\vee}\in U_{\mathcal{A}}^s(\h)$ has the form $q^{n_0(id+s^{-1})(id-\bar{n}dP_{\h'})\lambda^\vee}=\prod_{i=1}^lL_i^{r_i}$, $r_i\in \mathbb{Z}$, $i=1,\ldots, l$. 
So if we define $\varepsilon^{n_0(id+s^{-1})(id-\bar{n}dP_{\h'})\lambda^\vee}=q^{n_0(id+s^{-1})(id-\bar{n}dP_{\h'})\lambda^\vee}=\prod_{i=1}^lL_i^{r_i}$ mod $(q^{\frac{1}{d{\bar{r}}^2}}-\varepsilon^{\frac{1}{d{\bar{r}}^2}})$ then $\varepsilon^{n_0(id+s^{-1})(id-\bar{n}dP_{\h'})\lambda^\vee}=\prod_{i=1}^lL_i^{r_i}\in U_{\varepsilon}^s(\h)$.

Now using (\ref{phixi}) and (\ref{opdef}) one immediately verifies that 
$$
(\varepsilon^{(id+s^{-1})(id-\bar{n}dP_{\h'})\lambda^\vee})^{n_0}\phi_\xi^w(f)=\varphi_\varepsilon({\rm Ad}_s^0(\varepsilon^{-n_0(id+s^{-1})(id-\bar{n}dP_{\h'})\lambda^\vee})(f))\phi_\xi^w ((\Delta_\lambda^s)^{n_0})
$$
and that
$$
\phi_\xi^w ((\Delta_\lambda^s)^{n_0})=\varepsilon^{n_0(id+s^{-1})(id-\bar{n}dP_{\h'})\lambda^\vee}w.
$$

Recalling the equivariance of $\varphi_\varepsilon$ with respect to the action of $U_{\varepsilon}^{s,res}(\g)\supset U_{\varepsilon}^s(\h)$ we obtain from the last two identities that
$$
(\varepsilon^{(id+s^{-1})(id-\bar{n}dP_{\h'})\lambda^\vee})^{n_0}\phi_\xi^w(f)=
$$
$$
=\varepsilon^{n_0(id+s^{-1})(id-\bar{n}dP_{\h'})\lambda^\vee}\varphi_\varepsilon (f)  \varepsilon^{-n_0(id+s^{-1})(id-\bar{n}dP_{\h'})\lambda^\vee}\varepsilon^{n_0(id+s^{-1})(id-\bar{n}dP_{\h'})\lambda^\vee}w=
$$
$$
=\varepsilon^{n_0(id+s^{-1})(id-\bar{n}dP_{\h'})\lambda^\vee}\varphi_\varepsilon (f)w=\varepsilon^{n_0(id+s^{-1})(id-\bar{n}dP_{\h'})\lambda^\vee}\phi_\xi^w (f),
$$
i.e. the action of $(\varepsilon^{(id+s^{-1})(id-\bar{n}dP_{\h'})\lambda^\vee})^{n_0}$ on the image of $\phi_\xi^w$ coincides with the action of $\varepsilon^{n_0(id+s^{-1})(id-\bar{n}dP_{\h'})\lambda^\vee}$.

Using this result, observing that $L_i^{\bar{m}r_i}\in Z_0$, $i=1,\ldots, l$, and $\xi(L_i^{\bar{m}})\neq 0$, $i=1,\ldots, l$  by the definition of $\xi$, one has the following identity for operators defined on the image of $\phi_\xi^w$ in $V$
$$
(\varepsilon^{(id+s^{-1})(id-\bar{n}dP_{\h'})\lambda^\vee})^{\bar{m}n_0}=\varepsilon^{\bar{m}n_0(id+s^{-1})(id-\bar{n}dP_{\h'})\lambda^\vee}=\prod_{i=1}^lL_i^{\bar{m}r_i}=\xi(\prod_{i=1}^lL_i^{\bar{m}r_i})I_V=\xi(\varepsilon^{\bar{m}n_0(id+s^{-1})(id-\bar{n}dP_{\h'})\lambda^\vee})I_V,
$$
where
$$
\xi(\prod_{i=1}^lL_i^{\bar{m}r_i})=\xi(\varepsilon^{\bar{m}n_0(id+s^{-1})(id-\bar{n}dP_{\h'})\lambda^\vee})\neq 0,
$$
which implies that the operator $\varepsilon^{(id+s^{-1})(id-\bar{n}dP_{\h'})\lambda^\vee}$ is invertible. This proves part (iv).

Finally applying verbatim the arguments from the proof of Lemma \ref{phiext} we deduce that $\phi_\xi^w$ extends to a $U_{\eta_1}(\m_-)$--module homomorphism $\phi_\xi^w: \mathbb{C}_{\varepsilon}^{s,loc}[G] \rightarrow V$ in such a way that (\ref{phixiloc}) holds, and the specialization ${I^{11}_{\varepsilon}}^{loc}\subset \mathbb{C}_{\varepsilon}^{s,loc}[G]$ of ${I^{11}_{\mathcal{B}}}^{loc}\subset \mathbb{C}_{\mathcal{B}}^{s,loc}[G]$ at $q^{\frac{1}{d{\bar{r}}^2}}=\varepsilon^{\frac{1}{d{\bar{r}}^2}}$ belongs to the kernel of this homomorphism, so 
$$
\phi_\xi^w:\mathbb{C}_{11}^{loc}[G]_\varepsilon \rightarrow V.
$$ 
Note that due to invertibility of the operator $\varepsilon^{(id+s^{-1})(id-\bar{n}dP_{\h'})\lambda^\vee}$ no localization of $V$, which appears in Lemma \ref{phiext} for $Q_{\mathcal{B}}$ in the case of indeterminate $q$, is required. This completes the proof of part (v).

\end{proof}

Now we define root of unity counterparts of q-W--algebras.
Let $W^s_{\varepsilon,\xi}(G)={\rm End}_{U_{\xi}(\g)}(Q_\chi)^{opp}$ \index[not]{W@$W^s_{\varepsilon,\xi}(G)$} be the algebra of $U_{\xi}(\g)$--endomorphisms of $Q_\chi$ with the opposite multiplication. The algebra $W^s_{\varepsilon,\xi}(G)$ is also called {\it a q-W--algebra} associated to $s\in W$ and to $\xi\in {\rm Spec}(Z_0)$. Denote by $U_{\xi}(\g)-{\rm mod}$ \index[not]{U@$U_\xi(\g)-{\rm mod}$} the category of finite--dimensional left $U_{\xi}(\g)$--modules and by $W^s_{\varepsilon,\xi}(G)-{\rm mod}$ \index[not]{W@$W^s_{\varepsilon,\xi}(G)-{\rm mod}$} the category of finite--dimensional left $W^s_{\varepsilon,\xi}(G)$--modules. Observe that if $V\in U_{\xi}(\g)-{\rm mod}$ then the algebra $W^s_{\varepsilon,\xi}(G)$ naturally acts on the finite--dimensional space $V_\chi={\rm Hom}_{U_{\xi}({\frak m}_-)}(\mathbb{C}_\chi,V)\simeq {\rm Hom}_{U_{\xi}(\g)}(Q_\chi,V)$ by compositions of homomorphisms.

The following theorem is a root of unity analogue of the Skryabin equivalence for equivariant modules over quantum groups.
This theorem uncovers some striking similarity between the structure of the category of finite--dimensional representations of $U_\xi(\g)$ and of the category of equivariant modules over a quantum group for generic $\varepsilon$.

\begin{theorem}\label{fdfree}
(i) Every module $V\in U_{\xi}(\g)-{\rm mod}$ is isomorphic to ${\rm Hom}_{\mathbb{C}}(U_{\eta_1}(\m_-),V_\chi)$ as a right $U_{\eta_1}(\m_-)$--module, where the right action of $U_{\eta_1}(\m_-)$ on ${\rm Hom}_{\mathbb{C}}(U_{\eta_1}(\m_-),V_\chi)$ is induced by the multiplication in $U_{\eta_1}(\m_-)$ from the left. In particular, $V$ is $U_{\eta_1}(\m_-)$--injective, ${\rm Ext}^\bullet_{U_{\eta_1}(\m_-)}(\mathbb{C}_{\varepsilon},V)=V_\chi$ and the dimension of $V$ is divisible by ${\rm dim}~U_{\xi}({\frak m}_-)=\bar{m}^{\frac{1}{2}{\rm dim}~\mathcal{O}_{\pi\xi}}$.

(ii) $Q_\chi$ is isomorphic to ${\rm Hom}_{\mathbb{C}}(U_{\eta_1}(\m_-))\otimes W^s_{\varepsilon,\xi}(G)$ as a $U_{\eta_1}(\m_-)$--$W^s_{\varepsilon,\xi}(G)$--bimodule. 

(iii) The functor $E\mapsto Q_\chi\otimes_{W^s_{\varepsilon,\xi}(G)}E$ establishes an equivalence of the category of finite--dimensional left $W^s_{\varepsilon,\xi}(G)$--modules and the category $U_{\xi}(\g)-{\rm mod}$. The inverse equivalence is given by the functor $V\mapsto V_\chi$. In particular, the latter functor is exact, and every finite--dimensional $U_{\xi}(\g)$--module is generated by Whittaker vectors.
\end{theorem}

\begin{proof}
(i) Let $V$ be an object in the category $U_{\xi}(\g)-{\rm mod}$. 
Fix any linear map $\rho: V\rightarrow V_\chi$ the restriction of which to $V_\chi$ is the identity map, and let for any $v\in V$ $\sigma_\varepsilon(v):U_{\eta_1}(\m_-)\rightarrow V_\chi$ be the $\mathbb{C}$--linear homomorphism given by $\sigma_\varepsilon(v)(x)=\rho({\rm Ad}_sx(v))$, so that we have a map $\sigma_\varepsilon:V\rightarrow {\rm Hom}_{\mathbb{C}}(U_{\eta_1}(\m_-), V_\chi)$.

By the definition $\sigma_\varepsilon$ is a homomorphism of right $U_{\eta_1}(\m_-)$--modules, where the right action of $U_{\eta_1}(\m_-)$ on $${\rm Hom}_{\mathbb{C}}(U_{\eta_1}(\m_-), V_\chi)$$ is induced by multiplication in $U_{\eta_1}(\m_-)$ from the left. 

We claim that $\sigma_\varepsilon$ is an isomorphism. Firstly, $\sigma_\varepsilon$ is injective for otherwise its kernel would contain a non--zero Whittaker vector. Indeed by Lemma \ref{J1Jac} the augmentation ideal of $U_{\eta_1}(\m_-)$ coincides with its Jacobson radical which is nilpotent. Therefore its action on the kernel of $\sigma_\varepsilon$ is nilpotent, and hence the kernel, if it is non--trivial, must contain a non--zero Whittaker vector annihilated by the augmentation ideal of $U_{\eta_1}(\m_-)$. But all non--zero Whittaker vectors in $V$ belong to $V_\chi$ and by the definition of $\sigma_\varepsilon$ their images in ${\rm Hom}_{\mathbb{C}}(U_{\eta_1}(\m_-), V_\chi)$ are non--zero homomorphisms non-vanishing at $1$. 

Next we show that $\sigma_\varepsilon$ is also surjective. One could apply modified arguments in the proof of a similar statement in the case of generic $\varepsilon$ from Theorem \ref{sqeq} (ii). In the case of roots of unity one can give a more explicit proof which we present below.

For $m_1,\ldots , m_c\in \{0,\ldots, \bar{m}-1\}$ and any $v\in V_\chi$ we introduce the elements $v^{m_1\ldots m_c}:={c_{m_1\ldots m_c}'(\varepsilon)}^{-1}\phi_\xi^v(B_{m_1\ldots m_d}^\varepsilon)$  \index[not]{v@$v^{m_1\ldots m_c}$} which are are well--defined, as all $c_{m_1\ldots m_c}'(\varepsilon)$ are non--zero by the choice of $\varepsilon$, and hence Corollary \ref{Bpbas} for $r_1,\ldots, r_c\in \{0,\ldots, \bar{m}-1\}$ implies 
$$
{\rm Ad}_s(f_{\beta_1}^{r_1}\ldots f_{\beta_c}^{r_c})v^{m_1\ldots m_c}=\left\{\begin{array}{l} v ~~{\rm if}~~m_p=r_p~~{\rm for}~~ p=1,\ldots , c \\ 0 ~~{\rm if}~~r_i=m_i, i=1,\ldots, p-1~~{\rm and}~~ r_p>m_p~~{\rm for~~some}~~p\in \{1,\ldots , c\} \end{array}\right.,
$$
so that
\begin{equation}\label{basis}
\sigma_\varepsilon(v^{m_1\ldots m_c})(f_{\beta_1}^{r_1}\ldots f_{\beta_c}^{r_c})=\left\{\begin{array}{l} v ~~{\rm if}~~m_p=r_p~~{\rm for}~~ p=1,\ldots , c \\ 0 ~~{\rm if}~~r_i=m_i, i=1,\ldots, p-1~~{\rm and}~~ r_p>m_p~~{\rm for~~some}~~p\in \{1,\ldots , c\} \end{array}\right..
\end{equation}

Observe that by Lemma \ref{pbwr1} the elements $f_{\beta_1}^{r_1}\ldots f_{\beta_c}^{r_c}$, $r_1,\ldots, r_c\in \{0,\ldots, \bar{m}-1\}$ form a linear basis of $U_{\eta_1}(\m_-)$. 
The elements of this basis are labeled by the elements of the set $\mathbb{N}_{\bar{m}}^c$, where $\mathbb{N}_{\bar{m}}=\{0,1,\ldots ,\bar{m}-1\}$. \index[not]{N@$\mathbb{N}_{\bar{m}}$} Introduce the lexicographic order on this set, so that $(r_1,\ldots, r_c)>(m_1,\ldots, m_c)$ if $r_i=m_i$ for $i=1,\ldots, p-1$ and $r_p>m_p$ for some $p\in \{1,\ldots, c\}$. 

Note that for any $(r_1,\ldots, r_c)\in \mathbb{N}_{\bar{m}}^c$ the number of elements $(m_1,\ldots, m_c)\in \mathbb{N}_{\bar{m}}^c$ such that $(r_1,\ldots, r_c)>(m_1,\ldots, m_c)$ is finite.

Now let $(r_1,\ldots, r_c)\in \mathbb{N}_{\bar{m}}^c$, $v\in V_\chi$. If for $(m_1,\ldots, m_c)\in \mathbb{N}_{\bar{m}}^c$  such that $(r_1,\ldots, r_c)\geq (m_1,\ldots, m_c)$ we denote
$$
\sigma_\varepsilon(v_{m_1\ldots m_c}^{m_1\ldots m_c})=f_{m_1\ldots m_c},
$$
where for $(r_1,\ldots, r_c)\geq (m_1,\ldots, m_c)$ the elements $v_{m_1\ldots m_c}\in V_\chi$ are defined by induction starting from $v_{r_1\ldots r_c}=v$ as follows
\begin{equation}\label{vnn}
v_{m_1\ldots m_c}=-\sum_{(r_1,\ldots, r_c)\geq (m_1',\ldots, m_c')>(m_1,\ldots, m_c)}f_{m_1'\ldots m_c'}(f_{\beta_1}^{m_1}\ldots f_{\beta_c}^{m_c}),
\end{equation}
then using (\ref{basis}) one obtains
\begin{equation}\label{mm'}
f_{m_1\ldots m_c}(f_{\beta_1}^{m_1'}\ldots f_{\beta_c}^{m_c'})=\left\{\begin{array}{l} v_{m_1\ldots m_c} ~~{\rm if}~~(m_1',\ldots, m_c')=(m_1,\ldots, m_c)\\ 0 ~~{\rm if}~~(m_1',\ldots, m_c')>(m_1,\ldots, m_c)  \end{array}\right..
\end{equation}

From this property and from (\ref{vnn}) one immediately verifies that if we define
$$
f^{r_1\ldots r_c}_v=\sum_{(r_1,\ldots, r_c)\geq (m_1,\ldots, m_c)}f_{m_1\ldots m_c} \index[not]{f@$f^{r_1\ldots r_c}_v$}
$$
then for any $(m_1,\ldots, m_c)\in \mathbb{N}_{\bar{m}}^c$
\begin{equation}\label{base}
f^{r_1\ldots r_c}_v(f_{\beta_1}^{m_1}\ldots f_{\beta_c}^{m_c})=\left\{\begin{array}{l} v ~~{\rm if}~~(m_1,\ldots, m_c)=(r_1,\ldots, r_c),\\ 0 ~~{\rm if}~~(m_1,\ldots, m_c)\neq (r_1,\ldots, r_c)  \end{array}\right..
\end{equation}

Indeed, if $(m_1,\ldots, m_c)\neq (r_1,\ldots, r_c)$ then by the definition of $f^{r_1\ldots r_c}_v$ and by (\ref{mm'}) 
$$
f^{r_1\ldots r_c}_v(f_{\beta_1}^{m_1}\ldots f_{\beta_c}^{m_c})=\sum_{(r_1,\ldots, r_c)\geq (m_1',\ldots, m_c')}f_{m_1'\ldots m_c'}(f_{\beta_1}^{m_1}\ldots f_{\beta_c}^{m_c})=
$$
$$
=\sum_{(r_1,\ldots, r_c)\geq (m_1',\ldots, m_c')\geq (m_1,\ldots, m_c)}f_{m_1'\ldots m_c'}(f_{\beta_1}^{m_1}\ldots f_{\beta_c}^{m_c})=
$$
$$
=f_{m_1\ldots m_c}(f_{\beta_1}^{m_1}\ldots f_{\beta_c}^{m_c})+\sum_{(r_1,\ldots, r_c)\geq (m_1',\ldots, m_c')> (m_1,\ldots, m_c)}f_{m_1'\ldots m_c'}(f_{\beta_1}^{m_1}\ldots f_{\beta_c}^{m_c})=
$$
$$
=v_{m_1\ldots m_c}-v_{m_1\ldots m_c}=0,
$$
where at the last step we also used (\ref{vnn}).

If $(m_1,\ldots, m_c)= (r_1,\ldots, r_c)$ then similar arguments yield 
$$
f^{r_1\ldots r_c}_v(f_{\beta_1}^{r_1}\ldots f_{\beta_c}^{r_c})=\sum_{(r_1,\ldots, r_c)\geq (m_1',\ldots, m_c')}f_{m_1'\ldots m_c'}(f_{\beta_1}^{r_1}\ldots f_{\beta_c}^{r_c})=
$$
$$
=\sum_{(r_1,\ldots, r_c)\geq (m_1',\ldots, m_c')\geq (r_1,\ldots, r_c)}f_{m_1'\ldots m_c'}(f_{\beta_1}^{r_1}\ldots f_{\beta_c}^{r_c})=f_{r_1\ldots r_c}(f_{\beta_1}^{r_1}\ldots f_{\beta_c}^{r_c})=v_{r_1\ldots r_c}=v,
$$
where at thew last step we used the definition of $v_{r_1\ldots r_c}=v$. This confirms (\ref{base}).

Since by Lemma \ref{pbwr1} the elements $f_{\beta_1}^{r_1}\ldots f_{\beta_c}^{r_c}$, $r_1,\ldots, r_c\in \{0,\ldots, \bar{m}-1\}$ form a linear basis of $U_{\eta_1}(\m_-)$, formula (\ref{base}) implies that the elements $f^{r_1\ldots r_c}_v$ with $(r_1,\ldots, r_c)\in \mathbb{N}_{\bar{m}}^c$, $v\in V_\chi$ span ${\rm Hom}_{\mathbb{C}}(U_{\eta_1}(\m_-), V_\chi)$, and hence the elements $\sigma_\varepsilon(v^{m_1\ldots m_c})$ with $(m_1,\ldots, m_c)\in \mathbb{N}_m^c$, $v\in V_\chi$ span ${\rm Hom}_{\mathbb{C}}(U_{\eta_1}(\m_-), V_\chi)$ as well. Therefore $\sigma_\varepsilon$ is surjective. We conclude that $\sigma_\varepsilon$ is an isomorphism of right $U_{\eta_1}(\m_-)$--modules. This implies all statements in part (i).

(ii) Similarly to the case of generic $\varepsilon$ (see last part of the proof of Theorem \ref{sqeq}) one shows that $Q_\chi$ is isomorphic to ${\rm Hom}_{\mathbb{C}}(U_{\eta_1}(\m_-))\otimes W^s_{\varepsilon,\xi}(G)$ as a $U_{\eta_1}(\m_-)$--$W^s_{\varepsilon,\xi}(G)$--bimodule.

(iii) Let $E$ be a finite--dimensional $W^s_{\varepsilon,\xi}(G)$--module. Using the isomorphism $Q_\chi\simeq {\rm Hom}_{\mathbb{C}}(U_{\eta_1}(\m_-))\otimes W^s_{\varepsilon,\xi}(G)$ of $U_{\eta_1}(\m_-)$--$W^s_{\varepsilon,\xi}(G)$--bimodules and the vector space isomorphism $W^s_{\varepsilon,\xi}(G)={\rm End}_{U_{\xi}(\g)}(Q_\chi)^{opp}\simeq {\rm Hom}_{U_{\xi}({\frak m}_-)}(\mathbb{C}_\chi,Q_\chi)=(Q_\chi)_\chi$ implied by the Frobenius reciprocity, one immediately deduces similarly to the case of generic $\varepsilon$ in the proof of Theorem \ref{sqeq} that
$(Q_{\chi}\otimes_{W^s_{\varepsilon,\xi}(G)}E)_\chi=E$. Therefore to complete the proof of part (iii) it suffices to check that for any $V\in U_{\xi}(\g)-{\rm mod}$ the canonical map $f:Q_{\chi}\otimes_{W^s_{\varepsilon,\xi}(G)}V_\chi \rightarrow V$ is an isomorphism.

Indeed, $f$ is injective because otherwise by Proposition \ref{Whitt} its kernel would contain a non--zero Whittaker vector with respect to $\chi$. But all Whittaker vectors of $Q_{\chi}\otimes_{W^s_{\varepsilon,\xi}(G)}V_\chi$ belong to the subspace $1\otimes V_\chi$, and the restriction of $f$ to $1\otimes V_\chi$ induces an isomorphism of the spaces of Whittaker vectors of $Q_{\chi}\otimes_{W^s_{\varepsilon,\xi}(G)}V_\chi$ and of $V$.

In order to prove that $f$ is surjective we consider the exact sequence
$$
0\rightarrow Q_{\chi}\otimes_{W^s_{\varepsilon,\xi}(G)}V_\chi \rightarrow V \rightarrow V'\rightarrow 0,
$$
where $V'$ is the cokernel of $f$, and the corresponding long exact sequence of cohomology,
$$
0\rightarrow {\rm Ext}^{0}_{U_{\eta_1}({\frak m}_-)}(\mathbb{C}_{\varepsilon},
Q_{\chi}\otimes_{W^s_{\varepsilon,\xi}(G)}V_\chi)\rightarrow {\rm Ext}^{0}_{U_{\eta_1}({\frak m}_-)}(\mathbb{C}_{\varepsilon},
V)\rightarrow {\rm Ext}^{0}_{U_{\eta_1}({\frak m}_-)}(\mathbb{C}_{\varepsilon},
V')\rightarrow
$$
$$
\rightarrow {\rm Ext}^{1}_{U_{\eta_1}({\frak m}_-)}(\mathbb{C}_{\varepsilon},
Q_{\chi}\otimes_{W^s_{\varepsilon,\xi}(G)}V_\chi)\rightarrow \ldots .
$$

Now recall that $f$ induces an isomorphism of the spaces of Whittaker vectors of $Q_{\chi}\otimes_{W^s_{\varepsilon,\xi}(G)}V_\chi$ and of $V$. As we proved in part (i) the finite--dimensional $U_{\xi}(\g)$--module $Q_{\chi}\otimes_{W^s_{\varepsilon,\xi}(G)}V_\chi$ is injective over $U_{\eta_1}({\frak m}_-)$, and hence ${\rm Ext}^{1}_{U_{\eta_1}({\frak m}_-)}(\mathbb{C}_{\varepsilon},
Q_{\chi}\otimes_{W^s_{\varepsilon,\xi}(G)}V_\chi)=0$. Therefore the initial part of the long exact cohomology sequence takes the form
$$
0\rightarrow V_\chi \rightarrow V_\chi \rightarrow V'_\chi \rightarrow 0,
$$
where the second map in the last sequence is an isomorphism. Using the last exact sequence we deduce that $V'_\chi=0$. But if $V'$ were non--trivial it would contain a non--zero Whittaker vector by Proposition \ref{Whitt}. Thus $V'=0$, and $f$ is surjective. This completes the proof of the theorem.

\end{proof}

By the previous theorem every module $V\in U_{\xi}(\g)-{\rm mod}$ is isomorphic to $${\rm Hom}_{\mathbb{C}}(U_{\eta_1}(\m_-),V_\chi)\simeq {\rm Hom}_{\mathbb{C}}(U_{\eta_1}(\m_-),\mathbb{C})\otimes V_\chi$$ as a right $U_{\eta_1}(\m_-)$--module. In fact, one can show that the algebra $U_{\eta_1}(\m_-)$ is Frobenius, \index{algebra!Frobenius} i.e. its left regular representation is isomorphic to the dual of the right regular representation and its right regular representation is isomorphic to the dual of the left regular representation. Thus as a right $U_{\eta_1}(\m_-)$--module $V$ is isomorphic to $U_{\eta_1}(\m_-)\otimes V_\chi$, where the right action of $U_{\eta_1}(\m_-)$ on $U_{\eta_1}(\m_-)\otimes V_\chi$ is induced by the multiplication in $U_{\eta_1}(\m_-)$ from the right. In particular, $V$ is $U_{\eta_1}(\m_-)$--free. 

More generally, we have the following proposition.
\begin{proposition}\label{frob}
For any character $\eta: Z_0 \rightarrow \mathbb{C}$ the algebra $U_\eta({\frak g})$ and its subalgebra $U_\eta({\frak m}_-)$ are Frobenius algebras.
\end{proposition}

\begin{proof}
The proof of this proposition is parallel to the proof of a similar statement for Lie algebras over fields of prime characteristic (see Proposition 1.2 in \cite{FP}) and for the restricted form of the quantum group in \cite{Kum}. We shall only briefly outline the main steps of the proof for $U_\eta({\frak g})$. The proof for $U_\eta({\frak m}_-)$ is similar.

The key ingredient of the proof is the De Concini-Kac filtration on $U_\varepsilon(\g)\simeq U_\varepsilon^{s}(\g)$ defined as in Section \ref{cateq}.
Similarly to the case of generic $\varepsilon$ in Section \ref{cateq}, for ${\bf r},~{\bf m}\in {\Bbb N}^D$, $t\in U_\varepsilon({\frak h})$ we introduce the element
 $u_{{\bf r},{\bf m},t}=e^{\bf r}tf^{\bf m}$, where we use the notation of Lemma \ref{segmPBWs}. Here the generators $f_\alpha,e_\alpha,\alpha \in  \Delta_+$ and the ordered products of them are defined with the help of the normal ordering of $\Delta_+$ associated to $s$. Define also the height of the element
 $u_{{\bf r},{\bf m},t}$ as follows ${\rm ht}~(u_{{\bf r},{\bf m},t})=\sum_{i=1}^D(m_i+r_i){\rm ht}~\beta_i\in \mathbb{N}$, where ${\rm ht}~\beta_i$ is the height of the root $\beta_i$. Introduce also the degree of $u_{{\bf r},{\bf m},t}$ by
 $$
 d(u_{{\bf r},{\bf m},t})=(r_1,\ldots,r_D,m_D,\ldots,m_1,{\rm ht}~(u_{{\bf r},{\bf m},t}))\in \mathbb{N}^{2D+1}.
 $$
 Equip $\mathbb{N}^{2D+1}$ with the total lexicographic order and for ${\bf k}\in \mathbb{N}^{2D+1}$ denote by $(U_\varepsilon(\g))_{\bf k}$ the span of elements $u_{{\bf r},{\bf m},t}$ with $d(u_{{\bf r},{\bf m},t})\leq {\bf k}$ in $U_\varepsilon(\g)$. Then Proposition 1.7 in \cite{DK} implies that $(U_\varepsilon(\g))_{\bf k}$ is a filtration of $U_\varepsilon(\g)$ such that the associated graded algebra is the associative semi--commutative algebra over $\mathbb{C}$ with generators $e_\alpha,~f_\alpha$, $\alpha\in \Delta_+$, $L_i^{\pm 1}$, $i=1,\ldots l$ subject to the relations
 \begin{equation}\label{scomm}
 \begin{array}{l}
 L_iL_j=L_jL_i,~~L_iL_i^{-1}=L_i^{-1}L_i=1,~~ L_ie_\alpha L_i^{-1}=\varepsilon^{\alpha(Y_i)}e_\alpha, ~~L_if_\alpha L_i^{-1}=\varepsilon^{-\alpha(Y_i)}f_\alpha,\\
\\
e_\alpha f_\beta =\varepsilon^{\bar{n}d({1+s \over 1-s}P_{\h'^*}\alpha,\beta)} f_\beta e_\alpha,\\
 \\
 e_{\alpha}e_{\beta} = \varepsilon^{(\alpha,\beta)+\bar{n}d({1+s \over 1-s}P_{\h'^*}\alpha,\beta)}e_{\beta}e_{\alpha},~ \alpha<\beta, \\
 \\
 f_{\alpha}f_{\beta} = \varepsilon^{(\alpha,\beta)+\bar{n}d({1+s \over 1-s}P_{\h'^*}\alpha,\beta)}f_{\beta}f_{\alpha},~ \alpha<\beta.
 \end{array}
 \end{equation}

By Theorem 61.3 in \cite{CR} it suffices to show that there is a non-degenerate bilinear form $B_\eta:U_\eta({\frak g}) \times U_\eta({\frak g})\rightarrow \mathbb{C}$ \index[not]{B@$B_\eta$} which is associative \index{form!non-degenerate bilinear associative} in the sense that
$$
B_\eta(ab,c)=B_\eta(a,bc),~~a,b,c\in U_\eta({\frak g}).
$$

Consider the free $Z_0$--basis of $U_\varepsilon({\frak g})$ introduced in part (ii) of Proposition \ref{uq1z}. This basis consists of the monomials $x_I:=f^{\bf r}L^{\bf s}e^{\bf m}$, $I=(r_1,\ldots,r_D,s_1,\ldots,s_l,m_1,\ldots,m_D)$ for which $0\leq r_k,m_k,s_i<m$ for $i=1,\ldots ,l$, $k=1,\ldots ,D$. Set $I'=(m-1-r_1,\ldots,m-1-r_D,m-1-s_1,\ldots,m-1-s_l,m-1-m_1,\ldots,m-1-m_D)$ and $P=(m-1,\ldots,m-1)$.

Let $\Phi: U_\varepsilon({\frak g})\rightarrow Z_0$ be the $Z_0$--linear map defined on the basis $x_I$ of monomials by
$$
\Phi(x_I)=\left\{ \begin{array}{ll}
1 & I = P \\
0 & {\rm otherwise}
\end{array}
\right  . .
$$

Let $x=\sum_I c_Ix_I, c_I\in Z_0$ be an element of $U_\varepsilon({\frak g})$, and $c_K\neq 0$ the coefficient such that $d(x_K)$ is maximal possible with $c_K\neq 0$ in the sum defining $x$.
 
Using the definition of the De Concini--Kac filtration and commutation relations (\ref{scomm}) one can check that $\Phi(xx_{K'})=a_xc_K$, where $a_x$ is a nonzero complex number (see \cite{Kum}, proof of Theorem 2.2, Assertion I for details). 

Therefore the bilinear form $B_\eta:U_\eta({\frak g}) \times U_\eta({\frak g})\rightarrow \mathbb{C}$ obtained by composing the associative $Z_0$--bilinear pairing $B:U_\varepsilon({\frak g})\otimes_{Z_0} U_\varepsilon({\frak g}) \rightarrow Z_0$, $B(x,y)=\Phi(xy)$ with the character $\eta$ of $Z_0$ is non--degenerate and associative. This completes the proof.

\end{proof}

We restate the results of the discussion before the previous proposition as its corollary.
\begin{corollary}
As a right $U_{\eta_1}(\m_-)$--module, every module $V\in U_{\xi}(\g)-{\rm mod}$ is isomorphic to $U_{\eta_1}(\m_-)\otimes V_\chi$, where the right action of $U_{\eta_1}(\m_-)$ on $U_{\eta_1}(\m_-)\otimes V_\chi$ is induced by the multiplication in $U_{\eta_1}(\m_-)$ from the right. In particular, $V$ is $U_{\eta_1}(\m_-)$--free. 
\end{corollary}


\section{Properties of q-W--algebras associated to quantum groups at roots of unity}

\pagestyle{myheadings}
\markboth{CHAPTER~\thechapter. APPLICATION TO QUANTUM GROUPS AT ROOTS OF UNITY}{\thesection.~Q-W--ALGEBRAS AT ROOTS OF UNITY}

\setcounter{equation}{0}
\setcounter{theorem}{0}

In conclusion we study some further properties of q-W--algebras at roots of unity and of the module $Q_\chi$. We keep the notation introduced in the previous section. First we prove the following lemma.

\begin{lemma}
The left $U_{\xi}(\g)$--module $Q_\chi$ is projective in the category $U_{\xi}(\g)-{\rm mod}$.
\end{lemma}

\begin{proof}
We have to show that the functor ${\rm Hom}_{U_{\xi}(\g)}(Q_\chi,\cdot~)$ is exact. Let $V^\bullet$ be an exact complex of finite--dimensional $U_{\xi}(\g)$--modules. Since by Theorem \ref{fdfree} (i) any object  $V$ of $U_{\xi}(\g)-{\rm mod}$ is isomorphic to ${\rm Hom}_{\mathbb{C}}(U_{\eta_1}(\m_-), V_\chi)$ as a right $U_{\eta_1}(\m_-)$--module we have
$$
V^\bullet\simeq {\rm Hom}_{\mathbb{C}}(U_{\eta_1}(\m_-), \overline{V}^\bullet),
$$
where $\overline{V}^\bullet$ is an exact complex of vector spaces and the action of $U_{\eta_1}({\frak m}_-)$ on ${\rm Hom}_{\mathbb{C}}(U_{\eta_1}(\m_-), \overline{V}^\bullet)$ is induced by multiplication from the left on $U_{\eta_1}({\frak m}_-)$.

Now by the Frobenius reciprocity we have obvious isomorphisms of complexes,
\begin{eqnarray*}
{\rm Hom}_{U_{\xi}(\g)}(Q_\chi,V^\bullet)\simeq {\rm Hom}_{U_{\xi}({\frak m}_-)}(\mathbb{C}_\chi,V^\bullet)\simeq {\rm Hom}_{U_{\eta_1}({\frak m}_-)}(\mathbb{C}_\varepsilon,{\rm Hom}_{\mathbb{C}}(U_{\eta_1}(\m_-), \overline{V}^\bullet))\simeq  
\\ 
\simeq{\rm Hom}_\mathbb{C}(U_{\eta_1}({\frak m}_-)\otimes_{U_{\eta_1}({\frak m}_-)}\mathbb{C}_\varepsilon,\overline{V}^\bullet)\simeq \overline{V}^\bullet,
\end{eqnarray*}
where the last complex is exact, and we used the fact that by Lemma \ref{LW} for any finite--dimensional $U_{\xi}(\g)$--module $V$ one has ${\rm Hom}_{U_{\xi}({\frak m}_-)}(\mathbb{C}_\chi,V)\simeq  {\rm Hom}_{U_{\eta_1}({\frak m}_-)}(\mathbb{C}_\varepsilon,V)$. We conclude that the functor ${\rm Hom}_{U_{\xi}({\frak m}_-)}(Q_\chi,\cdot~)$ is exact.

\end{proof}

The properties of q-W--algebras at roots of unity are summarized in the following proposition.
\begin{proposition}\label{Umatr}
Denote $b_\xi=\bar{m}^{{\rm dim}~\mathfrak{m}_-}=\bar{m}^{\frac{1}{2}{\rm dim}~\mathcal{O}_{\pi\xi}}$. \index[not]{b@$b_\xi$} Then the following statements are true.

(i) $U_{\xi}(\g)\simeq {\rm Mat}_{b_\xi}(W^s_{\varepsilon,\xi}(G))$ as algebras, where ${\rm Mat}_{b_\xi}(W^s_{\varepsilon,\xi}(G))$ \index[not]{M@${\rm Mat}_{b_\xi}(W^s_{\varepsilon,\xi}(G))$} is the algebra of square matrices of size $b_\xi$ with entries from $W^s_{\varepsilon,\xi}(G)$.

(ii) $Q_\chi^{b_\xi}\simeq U_{\xi}(\g)$ as left $U_{\xi}(\g)$--modules.

(iii) $Q_\chi\simeq (W^s_{\varepsilon,\xi}(G)^{opp})^{b_\xi}$ as right $W^s_{\varepsilon,\xi}(G)$--modules.
\end{proposition}

\begin{proof}
(i) Let $E_i$, $i=1,\ldots ,C$ be the simple finite--dimensional modules over the finite--dimensional algebra $U_{\xi}(\g)$. Denote by $P_i$ the projective cover of $E_i$. Since by Theorem \ref{fdfree} the dimension of $E_i$ is divisible by $b_\xi$ we have ${\rm dim}~E_i=b_\xi r_i$, $r_i\in \mathbb{N}$, where $r_i$ is the rank of $E_i$ over $U_{\eta_1}({\frak m}_-)$ equal to the dimension of the space of Whittaker vectors in $E_i$. By Proposition 2.1 in \cite{Pr}
$$
U_{\xi}(\g)={\rm Mat}_{b_\xi}({\rm End}_{U_{\xi}(\g)}(P)^{opp}),
$$
where $P=\bigoplus_{i=1}^CP_i^{r_i}$. Therefore to establish the first isomorphism in this proposition it suffices to show that $P\simeq Q_\chi$. Since by the previous lemma $Q_\chi$ is projective in the category $U_{\xi}(\g)-{\rm mod}$ we only need to verify that $$r_i={\rm dim}~{\rm Hom}_{U_{\xi}(\g)}(P,E_i)={\rm dim}~{\rm Hom}_{U_{\xi}(\g)}(Q_\chi,E_i).$$

Indeed, by the Frobenius reciprocity we have
$$
{\rm dim}~{\rm Hom}_{U_{\xi}(\g)}(Q_\chi,E_i)={\rm dim}~{\rm Hom}_{U_{\xi}({\frak m}_-)}(\mathbb{C}_\chi,E_i)=r_i.
$$

This proves the first statement of the proposition. 

(ii) From Proposition 2.1 in \cite{Pr} we also deduce that $P^{b_\xi}\simeq U_{\xi}(\g)$ as left $U_{\xi}(\g)$--modules. Together with the isomorphism $P\simeq Q_\chi$ this gives the second isomorphism in the proposition.

(iii) Using results of Section 6.4 in \cite{Pie} and the fact that $Q_\chi$ is projective one can find an idempotent $e\in  U_{\xi}(\g)$ such that $Q_\chi\simeq U_{\xi}(\g)e$ as modules and $(W^s_{\varepsilon,\xi}(G))^{opp}\simeq eU_{\xi}(\g)e$ as algebras.

By the first two statements of this proposition one can also find idempotents $e=e_1,e_2,\ldots ,e_{b_\xi}\in  U_{\xi}(\g)$ such that $e_1+\ldots +e_{b_\xi}=1$, $e_ie_j=0$ if $i\neq j$, $i,j=1,\ldots, b_\xi$ and $e_iU_{\xi}(\g)\simeq eU_{\xi}(\g)$, $i=1,\ldots, b_\xi$ as right $U_{\xi}(\g)$--modules. Therefore $e_iU_{\xi}(\g)e\simeq eU_{\xi}(\g)e$ as right $eU_{\xi}(\g)e$--modules, and
$$
Q_\chi\simeq U_{\xi}(\g)e=\bigoplus_{i=1}^{b_\xi}e_iU_{\xi}(\g)e\simeq (eU_{\xi}(\g)e)^{b_\xi}\simeq (W^s_{\varepsilon,\xi}(G)^{opp})^{b_\xi}
$$
as right $W^s_{\varepsilon,\xi}(G)$--modules.
This establishes the last isomorphism in the statement of this proposition and completes its proof. 

\end{proof}

\begin{corollary}
The algebra $W^s_{\varepsilon,\xi}(G)$ is finite--dimensional, and ${\rm dim}~W^s_{\varepsilon,\xi}(G)=\bar{m}^{{\rm dim}~\Sigma_s}$.
\end{corollary}

\begin{proof}
By Proposition \ref{var} $2{\rm dim}~\mathfrak{m}_-+{\rm dim}~\Sigma_s={\rm dim}~G$. Therefore by the definition of $Q_\chi$ we have ${\rm dim}~Q_\chi=\bar{m}^{{\rm dim}~G-{\rm dim}~\mathfrak{m}_-}=\bar{m}^{{\rm dim}~\mathfrak{m}_-+{\rm dim}~\Sigma_s}$. Finally from the last statement of the previous theorem one obtains that ${\rm dim}~W^s_{\varepsilon,\xi}(G)={\rm dim}~Q_\chi/\bar{m}^{{\rm dim}~\mathfrak{m}_-}=\bar{m}^{{\rm dim}~\Sigma_s}$.

\end{proof}

Using Proposition \ref{d'} we deduce from Proposition \ref{Umatr} the following statement on the structure of the algebra $U_{\eta}(\g)$.
\begin{corollary}
Let $\eta \in {\rm Spec}(Z_0)$ be an element  such that $\pi\eta\in G_\mathcal{C}$, $\mathcal{C}\in C(W)$ and $s^{-1}\in \mathcal{C}$, $d=2d'$, where $d'$ is defined in Proposition \ref{d'}. Assume that $\bar{m}$ and $d$ are coprime.

Then $U_{\eta}(\g)\simeq {\rm Mat}_{b_\eta}(W^s_{\varepsilon,\xi}(G))$, where $\xi\in  {\rm Spec}Z_0$ is chosen as in Proposition \ref{Umatr}, and $b_\eta=\bar{m}^{\frac{1}{2}{\rm dim}~\mathcal{O}_{\pi\eta}}$.

Let $\mathcal{L}$ be a sheaf of algebras over $ {\rm Spec}Z_0$ the stalk of which over $\eta\in  {\rm Spec}Z_0$ is $U_{\eta}(\g)$. Assume that the conditions imposed on $\bar{m}$ are satisfied for all Weyl group conjugacy classes in $C(W)$.
Then the sheaf $\mathcal{L}$ is isomorphic to a sheaf the stalk of which over any $\eta\in  {\rm Spec}Z_0$ with $\pi\eta \in G^0\cap G_\mathcal{C}$, $\mathcal{C}\in C(W)$ is ${\rm Mat}_{b_\eta}(W^s_{\varepsilon,\xi}(G))$, where $\xi\in {\rm Spec}Z_0$ is chosen as in Proposition \ref{Umatr}, $s^{-1}\in \mathcal{C}$, $b_\eta =\bar{m}^{\frac{1}{2}{\rm dim}~\mathcal{O}_{\pi \eta}}$.
\end{corollary}


\section{Bibliographic comments}

\pagestyle{myheadings}
\markboth{CHAPTER \thechapter.~APPLICATION TO QUANTUM GROUPS AT ROOTS OF UNITY}{\thesection.~BIBLIOGRAPHIC COMMENTS}

\setcounter{equation}{0}
\setcounter{theorem}{0}

The study of representations of quantum groups at roots of unity was initiated in \cite{DK}, where the quantum coadjoint action was defined as well. This action was studied in detail in \cite{DKP1} where the De Concini--Kac--Procesi conjecture on the dimensions of irreducible representations of quantum groups at roots of unity was formulated.

The results on quantum groups at roots of unity stated in Section \ref{1root} can be found in \cite{DK} and \cite{DKP1}.
Proposition \ref{d'} first appeared in Appendix A to \cite{S12}.
The statements of Proposition \ref{ue} can be found in \cite{DK}, Corollary 3.3, \cite{DKP1}, Theorems 3.5, 7.6 and Proposition 4.5.
Proposition \ref{LO} is Corollary 4.7 in \cite{DKP1}, and the statements of Proposition \ref{qcoadj} appear in Propositions 3.4, 3.5, \cite{DK}, and in Proposition 6.1 and Theorem 6.6 in \cite{DKP1}.
Finite--dimensional quotients $U_\eta(\g)$ were introduced in \cite{DK1}.

The notions of Whittaker vectors for representations of quantum groups at roots of unity, of the algebras $U_{\eta_1}(\m_-)$, and of their actions on finite--dimensional representations of ${U}_\eta(\g)$ were introduced in \cite{S11}, and the exposition in Section \ref{Whittr1} follows \cite{S11} as well.

In the representation theory of Lie algebras in prime characteristic there is a conjecture similar to the De Concini--Kac--Procesi conjecture. It is called the Kac--Weisfeiler conjecture. Our proof of the De Concini--Kac--Procesi conjecture is conceptually similar to the proof of the Kac--Weisfeiler conjecture given in \cite{Sk} which is in turn a straightforward prime characteristic generalization of the proof of the Skryabin equivalence for reductive Lie algebras over algebraically closed fields of zero characteristic suggested in the Appendix to \cite{PS}. All these proofs go back to the original Kostant's idea on the proof of the classification theorem for Whittaker representations of complex semisimple Lie algebras in \cite{K}, the proof of the Skryabin equivalence in \cite{PS} being a significantly refined and simplified version of the proof of the main Theorem 3.3 in \cite{K}.

The properties of q-W--algeras at roots of unity are similar to those of W--algebras associated to semisimple Lie algebras in prime characteristic proved in \cite{Pr},
Proposition \ref{Umatr} being an analogue of Theorem 2.3 in \cite{Pr}.


\chapter*{Appendix}
\addcontentsline{toc}{chapter}{Appendix}

\section*{Appendix 1. Normal orderings of root systems compatible with involutions in Weyl groups}
\addcontentsline{toc}{section}{Appendix 1. Normal orderings of root systems compatible with involutions in Weyl groups}

\renewcommand{\thetheorem}{A1.\arabic{theorem}}

\renewcommand{\thelemma}{A1.\arabic{lemma}}

\renewcommand{\theproposition}{A1.\arabic{proposition}}

\renewcommand{\thecorollary}{A1.\arabic{corollary}}

\renewcommand{\theremark}{A1.\arabic{remark}}

\renewcommand{\thedefinition}{A1.\arabic{definition}}

\renewcommand{\theequation}{A1.\arabic{equation}}

\setcounter{equation}{0}
\setcounter{theorem}{0}

\pagestyle{myheadings}
\markboth{APPENDIX}{APPENDIX 1. NORMAL ORDERINGS COMPATIBLE WITH INVOLUTIONS IN WEYL GROUPS}

By Theorem A in \cite{Ric} every involution $w$ in the Weyl group $W$ of the pair $(\g,\h)$ is the longest element of the Weyl group of a Levi subalgebra in $\g$ with respect to some system of positive roots, and $w$ acts by multiplication by $-1$ in the Cartan subalgebra $\h_w\subset \h$ of the semisimple part $\m_w$ of that Levi subalgebra. By Lemma 5 in \cite{C} the involution $w$ can also be expressed as a product of ${\rm dim}~\h_w$ reflections from the Weyl group of the pair $(\m_w,\h_w)$, with respect to mutually orthogonal roots, $w=s_{\gamma_1}\ldots s_{\gamma_{\widetilde{l}}}$.


If $w$ is the longest element in the Weyl group of the pair $(\m_w,\h_w)$ with respect to some system of positive roots, where $\m_w$ is a simple Lie algebra and $\h_w$ is a Cartan subalgebra of $\m_w$, then $w$ is an involution acting by multiplication by $-1$ in $\h_w$ if and only if $\m_w$ is of one of the following types: $A_1,B_l,C_l,D_{2n},E_7,E_8,F_4,G_2$.

Fix a system of positive roots $\Delta_+(\m_w,\h_w)$ of the pair $(\m_w,\h_w)$. Let $w=s_{\gamma_1}\ldots s_{\gamma_{\widetilde{l}}}$ be a representation of $w$ as a product of ${\rm dim}~\h_w$ reflections from the Weyl group of the pair $(\m_w,\h_w)$, with respect to mutually orthogonal positive roots. A normal ordering of $\Delta_+(\m_w,\h_w)$ is called {\it compatible} with the decomposition $w=s_{\gamma_1}\ldots s_{\gamma_{\widetilde{l}}}$ if it is of the following form
$$
\beta_1, \ldots, \beta_{\frac{p-\widetilde{l}}{2}}, \gamma_1, \ldots, \gamma_2, \ldots, \gamma_3,\ldots, \gamma_{\widetilde{l}},
$$
where $p$ is the number of positive roots, and for any two positive roots $\alpha, \beta\in \Delta_+(\m_w,\h_w)$ such that $\gamma_1\leq \alpha<\beta$ the sum $\alpha+\beta$ cannot be represented as a linear combination $\sum_{k=1}^qc_k\gamma_{i_k}$, where $c_k\in \mathbb{N}$ and $\alpha<\gamma_{i_1}<\ldots <\gamma_{i_k}<\beta$.
Note that from the definition it also follows that
\begin{equation}\label{pnuma}
|[\beta_1,\beta_{\frac{p-\widetilde{l}}{2}}]|=\frac{p-\widetilde{l}}{2}, ~|[\gamma_1,\gamma_{\widetilde{l}}]|=\frac{p+\widetilde{l}}{2}.
\end{equation}

Existence of such compatible normal orderings is checked straightforwardly for all simple Lie algebras of types $A_1,B_l,C_l,D_{2n},E_7,E_8,F_4$ and $G_2$. In the case of $A_1$ the existence of such an ordering is obvious since there is only one positive root. In the other cases normal orderings defined by the properties described below for each of the types $B_l,C_l,D_{2n},E_7,E_8,F_4,G_2$ exist and are compatible with decompositions of nontrivial involutions in the Weyl group.
We use the Bourbaki notation for the systems of positive and simple roots (see \cite{Bur}).

\begin{itemize}

\vskip 0.5cm

\item

$B_l$

\vskip 0.2cm

Dynkin diagram:

$$\xymatrix@R=.25cm{
\alpha_1&\alpha_2&&\alpha_{l-2}&\alpha_{l-1}&\alpha_l\\
{\bullet}\ar@{-}[r]&{\bullet}\ar@{-}[r]&\cdots
\ar@{-}[r]& {\bullet}\ar@{-}[r]&{\bullet}\ar@2{-}[r] &{\bullet}
}$$

Simple roots: $\alpha_1=\varepsilon_1-\varepsilon_2,\alpha_2=\varepsilon_2-\varepsilon_3,\ldots, \alpha_{l-1}=\varepsilon_{l-1}-\varepsilon_l,\alpha_l=\varepsilon_l$.

Positive roots: $\varepsilon_i$ $(1\leq i\leq l)$, $\varepsilon_i-\varepsilon_j,\varepsilon_i+\varepsilon_j$ $(1\leq i<j\leq l)$.

The longest element of the Weyl group expressed as a product of ${\rm dim}~\h_w$ reflections with respect to mutually orthogonal roots: $w=s_{\varepsilon_1}\ldots s_{\varepsilon_l}$.

Normal ordering of $\Delta_+(\m_w,\h_w)$ compatible with expression $w=s_{\varepsilon_1}\ldots s_{\varepsilon_l}$:

$$
\varepsilon_1-\varepsilon_2,\ldots,\varepsilon_{l-1}-\varepsilon_l,\varepsilon_1,\ldots,\varepsilon_2,\ldots,
\varepsilon_l,
$$
where the roots $\varepsilon_i-\varepsilon_j$ $(1\leq i<j\leq l)$ forming the subsystem $\Delta_+(A_{l-1})\subset \Delta_+(B_l)$ are situated to the left from $\varepsilon_1$, and the roots $\varepsilon_i+\varepsilon_j$ $(1\leq i<j\leq l)$ are situated to the right from $\varepsilon_1$.

\vskip 0.5cm

\item

$C_l$

\vskip 0.2cm

Dynkin diagram:

$$
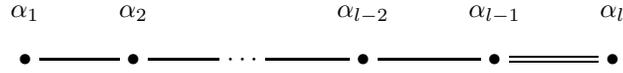
\xymatrix@R=.25cm{
\alpha_1&\alpha_2&&\alpha_{l-2}&\alpha_{l-1}&\alpha_l\\
{\bullet}\ar@{-}[r]&{\bullet}\ar@{-}[r]&\cdots
\ar@{-}[r]& {\bullet}\ar@{-}[r]&{\bullet}\ar@2{-}[r] &{\bullet}
}$$

Simple roots: $\alpha_1=\varepsilon_1-\varepsilon_2,\alpha_2=\varepsilon_2-\varepsilon_3,\ldots, \alpha_{l-1}=\varepsilon_{l-1}-\varepsilon_l,\alpha_l=2\varepsilon_l$.

Positive roots: $2\varepsilon_i$ $(1\leq i\leq l)$, $\varepsilon_i-\varepsilon_j,\varepsilon_i+\varepsilon_j$ $(1\leq i<j\leq l)$.

The longest element of the Weyl group expressed as a product of ${\rm dim}~\h_w$ reflections with respect to mutually orthogonal roots: $w=s_{2\varepsilon_1}\ldots s_{2\varepsilon_l}$.

Normal ordering of $\Delta_+(\m_w,\h_w)$ compatible with expression $w=s_{2\varepsilon_1}\ldots s_{2\varepsilon_l}$:

$$
\varepsilon_1-\varepsilon_2,\ldots,\varepsilon_{l-1}-\varepsilon_l,2\varepsilon_1,\ldots,2\varepsilon_2,\ldots,
2\varepsilon_l,
$$
where the roots $\varepsilon_i-\varepsilon_j$ $(1\leq i<j\leq l)$ forming the subsystem $\Delta_+(A_{l-1})\subset \Delta_+(C_l)$ are situated to the left from $2\varepsilon_1$, and the roots $\varepsilon_i+\varepsilon_j$ $(1\leq i<j\leq l)$ are situated to the right from $2\varepsilon_1$.

\vskip 0.5cm

\item

$D_{2n}$

\vskip 0.2cm

Dynkin diagram:

$$
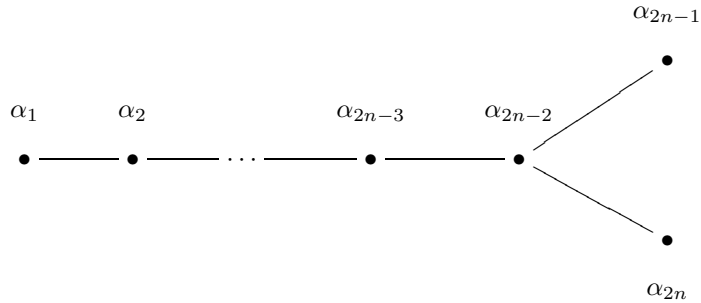
\xymatrix@R=.25cm{
&&&&&\alpha_{2n-1}\\
&&&&&{\bullet}\\
\alpha_1&\alpha_2&&\alpha_{2n-3}& \alpha_{2n-2}&\\
{\bullet}\ar@{-}[r]&{\bullet}\ar@{-}[r]&\cdots
\ar@{-}[r]&{\bullet}\ar@{-}[r]& {\bullet}\ar@{-}[uur]\ar@{-}[ddr]& \\
&&&&&\\
&&&&&{\bullet}\\
&&&&&\alpha_{2n}
}$$

Simple roots: $\alpha_1=\varepsilon_1-\varepsilon_2,\alpha_2=\varepsilon_2-\varepsilon_3,\ldots, \alpha_{2n-1}=\varepsilon_{2n-1}-\varepsilon_{2n},\alpha_{2n}=\varepsilon_{2n-1}+\varepsilon_{2n}$.

Positive roots: $\varepsilon_i-\varepsilon_j,\varepsilon_i+\varepsilon_j$ $(1\leq i<j\leq 2n)$.

The longest element of the Weyl group expressed as a product of ${\rm dim}~\h_w$ reflections with respect to mutually orthogonal roots: $$w=s_{\varepsilon_1-\varepsilon_2}s_{\varepsilon_1+\varepsilon_2}\ldots s_{\varepsilon_{2n-1}-\varepsilon_{2n}}s_{\varepsilon_{2n-1}+\varepsilon_{2n}}.$$

Normal ordering of $\Delta_+(\m_w,\h_w)$ compatible with expression $$w=s_{\varepsilon_1-\varepsilon_2}s_{\varepsilon_1+\varepsilon_2}\ldots s_{\varepsilon_{2n-1}-\varepsilon_{2n}}s_{\varepsilon_{2n-1}+\varepsilon_{2n}}:$$

\begin{eqnarray*}
\varepsilon_2-\varepsilon_3,\varepsilon_4-\varepsilon_5,\ldots,\varepsilon_{2n-2}-\varepsilon_{2n-1},\ldots,
\varepsilon_1-\varepsilon_2,\varepsilon_3-\varepsilon_4,\ldots,\varepsilon_{2n-1}-\varepsilon_{2n-2}, \\
\varepsilon_1+\varepsilon_2,\ldots,\varepsilon_3+\varepsilon_4,\ldots,
\varepsilon_{2n-1}+\varepsilon_{2n},
\end{eqnarray*}
where the roots $\varepsilon_i-\varepsilon_j$ $(1\leq i<j\leq l)$ forming the subsystem $\Delta_+(A_{l-1})\subset \Delta_+(C_l)$ are situated to the left from $\varepsilon_1+\varepsilon_2$, and the roots $\varepsilon_i+\varepsilon_j$ $(1\leq i<j\leq l)$ are situated to the right from $\varepsilon_1+\varepsilon_2$.

\vskip 0.5cm

\item

\vskip 0.2cm

$E_7$

Dynkin diagram:

$$
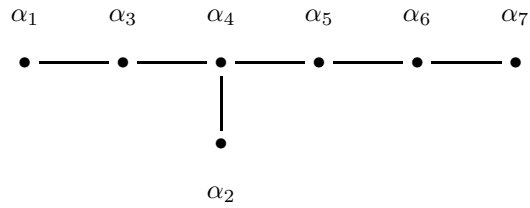
\xymatrix@R=.25cm@C=.25cm{
&\alpha_1&&\alpha_3&&\alpha_4&&\alpha_5&&\alpha_6&&\alpha_7\\
&{\bullet}\ar@{-}[rr]&&{\bullet}\ar@{-}[rr]&&{\bullet}
\ar@{-}[dd]\ar@{-}[rr]&& {\bullet}\ar@{-}[rr] &&{\bullet}\ar@{-}[rr]
&&{\bullet}\\ &&&&&&&&&&&\\
&&&&&{\bullet}&&&&&&\\
&&&&&\alpha_2&&&&&&}$$

Simple roots: $\alpha_1=\frac{1}{2}(\varepsilon_1+\varepsilon_8)-
\frac{1}{2}(\varepsilon_2+\varepsilon_3+\varepsilon_4+\varepsilon_5+\varepsilon_6+\varepsilon_7)$,
$\alpha_2=\varepsilon_1+\varepsilon_2,\alpha_3=\varepsilon_2-\varepsilon_1,$ $ \alpha_4=\varepsilon_3-\varepsilon_2,\alpha_5=\varepsilon_4-\varepsilon_3,
\alpha_6=\varepsilon_5-\varepsilon_4, \alpha_7=\varepsilon_6-\varepsilon_5$.

Positive roots: $\pm\varepsilon_i+\varepsilon_j$ $(1\leq i<j\leq 6)$, $\varepsilon_8-\varepsilon_7$,
$\frac{1}{2}(\varepsilon_8-\varepsilon_7+\sum_{i=1}^6(-1)^{\nu(i)}\varepsilon_i)$ with $\sum_{i=1}^6{\nu(i)}$ odd.

The longest element of the Weyl group expressed as a product of ${\rm dim}~\h_w$ reflections with respect to mutually orthogonal roots: $$w=s_{\varepsilon_2-\varepsilon_1}s_{\varepsilon_2+\varepsilon_1}
s_{\varepsilon_4-\varepsilon_3}s_{\varepsilon_4+\varepsilon_3} s_{\varepsilon_{6}-\varepsilon_{5}}s_{\varepsilon_{6}+\varepsilon_{5}}s_{\varepsilon_8-\varepsilon_7}.$$

Normal ordering of $\Delta_+(\m_w,\h_w)$ compatible with expression $$w=s_{\varepsilon_2-\varepsilon_1}s_{\varepsilon_2+\varepsilon_1}
s_{\varepsilon_4-\varepsilon_3}s_{\varepsilon_4+\varepsilon_3} s_{\varepsilon_{6}-\varepsilon_{5}}s_{\varepsilon_{6}+\varepsilon_{5}}s_{\varepsilon_8-\varepsilon_7}:$$

\begin{eqnarray*}
\alpha_1,\varepsilon_3-\varepsilon_2,\varepsilon_5-\varepsilon_4,\ldots,\varepsilon_{8}-\varepsilon_{7},\ldots,
\varepsilon_2-\varepsilon_1,\varepsilon_4-\varepsilon_3,\varepsilon_6-\varepsilon_5,\ldots, \\
\varepsilon_{6}+\varepsilon_{5},\ldots,\varepsilon_4+\varepsilon_3,\ldots,\varepsilon_2+\varepsilon_1,
\end{eqnarray*}
where the roots $\pm\varepsilon_i+\varepsilon_j$ $(1\leq i<j\leq 6)$ forming the subsystem $\Delta_+(D_6)\subset \Delta_+(E_7)$ are placed as in case of the compatible normal ordering of the system $\Delta_+(D_6)$, the only roots from the subsystem $\Delta_+(A_5)\subset \Delta_+(D_6)$ situated to the right from the maximal root $\varepsilon_{8}-\varepsilon_{7}$ are $\varepsilon_2-\varepsilon_1,\varepsilon_4-\varepsilon_3,\varepsilon_6-\varepsilon_5$, the roots $\varepsilon_i+\varepsilon_j$ $(1\leq i<j\leq 6)$ are situated to the right from $\varepsilon_6+\varepsilon_5$, and a half of the positive roots which do not belong to the subsystem $\Delta_+(D_6)\subset \Delta_+(E_7)$ are situated to the left from $\varepsilon_{8}-\varepsilon_{7}$ and the other half of those roots are situated to the right from $\varepsilon_{8}-\varepsilon_{7}$.

\vskip 0.5cm

\item

$E_8$

\vskip 0.2cm

Dynkin diagram:

$$
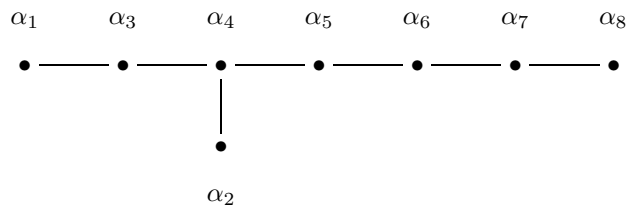
\xymatrix@R=.25cm@C=.25cm{
&\alpha_1&&\alpha_3&&\alpha_4&&\alpha_5&&\alpha_6&&\alpha_7&&\alpha_8\\
&{\bullet}\ar@{-}[rr]&&{\bullet}\ar@{-}[rr]&&{\bullet}
\ar@{-}[dd]\ar@{-}[rr]&& {\bullet}\ar@{-}[rr] &&{\bullet}\ar@{-}[rr]
&&{\bullet}\ar@{-}[rr]&& {\bullet}\\ &&&&&&&&&&&&&\\
&&&&&{\bullet}&&&&&&&& \\
&&&&&\alpha_2&&&&&&&&}$$

Simple roots: $\alpha_1=\frac{1}{2}(\varepsilon_1+\varepsilon_8)-
\frac{1}{2}(\varepsilon_2+\varepsilon_3+\varepsilon_4+\varepsilon_5+\varepsilon_6+\varepsilon_7)$,
$\alpha_2=\varepsilon_1+\varepsilon_2,\alpha_3=\varepsilon_2-\varepsilon_1,$ $ \alpha_4=\varepsilon_3-\varepsilon_2,\alpha_5=\varepsilon_4-\varepsilon_3,
\alpha_6=\varepsilon_5-\varepsilon_4, \alpha_7=\varepsilon_6-\varepsilon_5, \alpha_8=\varepsilon_7-\varepsilon_6$.

Positive roots: $\pm\varepsilon_i+\varepsilon_j$ $(1\leq i<j\leq 8)$,
$\frac{1}{2}(\varepsilon_8+\sum_{i=1}^7(-1)^{\nu(i)}\varepsilon_i)$ with $\sum_{i=1}^7{\nu(i)}$ even.

The longest element of the Weyl group expressed as a product of ${\rm dim}~\h_w$ reflections with respect to mutually orthogonal roots: $$w=s_{\varepsilon_2-\varepsilon_1}s_{\varepsilon_2+\varepsilon_1}
s_{\varepsilon_4-\varepsilon_3}s_{\varepsilon_4+\varepsilon_3} s_{\varepsilon_{6}-\varepsilon_{5}}s_{\varepsilon_{6}+\varepsilon_{5}}
s_{\varepsilon_8-\varepsilon_7}s_{\varepsilon_8+\varepsilon_7}.$$

Normal ordering of $\Delta_+(\m_w,\h_w)$ compatible with expression $$w=s_{\varepsilon_2-\varepsilon_1}s_{\varepsilon_2+\varepsilon_1}
s_{\varepsilon_4-\varepsilon_3}s_{\varepsilon_4+\varepsilon_3} s_{\varepsilon_{6}-\varepsilon_{5}}s_{\varepsilon_{6}+\varepsilon_{5}}
s_{\varepsilon_8-\varepsilon_7}s_{\varepsilon_8-\varepsilon_7}:$$

\begin{eqnarray*}
\alpha_1,\varepsilon_3-\varepsilon_2,\varepsilon_5-\varepsilon_4,\varepsilon_{7}-\varepsilon_{6},\ldots,
\varepsilon_2-\varepsilon_1,\varepsilon_4-\varepsilon_3,\varepsilon_6-\varepsilon_5,
\\
\varepsilon_{8}-\varepsilon_{7},\ldots,\varepsilon_{8}+\varepsilon_{7},\ldots,
\varepsilon_{6}+\varepsilon_{5},\ldots,\varepsilon_4+\varepsilon_3,\ldots,\varepsilon_2+\varepsilon_1,
\end{eqnarray*}
where the roots $\pm\varepsilon_i+\varepsilon_j$ $(1\leq i<j\leq 8)$ forming the subsystem $\Delta_+(D_8)\subset \Delta_+(E_8)$ are placed as in case of the compatible normal ordering of the system $\Delta_+(D_8)$, the roots $\varepsilon_i+\varepsilon_j$ $(1\leq i<j\leq 8)$ are situated to the right from $\varepsilon_8+\varepsilon_7$; the positive roots which do not belong to the subsystem $\Delta_+(D_8)\subset \Delta_+(E_8)$ can be split into two groups: the roots from the first group contain $\frac{1}{2}(\varepsilon_8+\varepsilon_7)$ in their decompositions with respect to the basis $\varepsilon_i,i=1,\ldots,8$, and the roots from the second group contain $\frac{1}{2}(\varepsilon_8-\varepsilon_7)$ in their decompositions with respect to the basis $\varepsilon_i,i=1,\ldots,8$; a half of the roots from the first group are situated to the left from $\varepsilon_{2}-\varepsilon_{1}$ and the other half of those roots are situated to the right from $\varepsilon_{8}+\varepsilon_{7}$; a half of the roots from the second group are situated to the left from $\varepsilon_{2}-\varepsilon_{1}$ and the other half of those roots are situated to the right from $\varepsilon_{8}-\varepsilon_{7}$.

\vskip 0.5cm

\item

$F_4$

\vskip 0.2cm

Dynkin diagram:

$$\xymatrix@R=.25cm{
\alpha_1&\alpha_2&\alpha_3&\alpha_4\\
{\bullet}\ar@{-}[r]&{\bullet}\ar@2{-}[r]
&{\bullet}\ar@{-}[r]&{\bullet} }$$

Simple roots: $\alpha_1=\varepsilon_2-\varepsilon_3,\alpha_2=\varepsilon_3-\varepsilon_4, \alpha_{3}=\varepsilon_{4},\alpha_4=\frac{1}{2}(\varepsilon_1-\varepsilon_2-\varepsilon_3-\varepsilon_4)$.

Positive roots: $\varepsilon_i$ $(1\leq i\leq 4)$, $\varepsilon_i-\varepsilon_j,\varepsilon_i+\varepsilon_j$ $(1\leq i<j\leq 4)$, $\frac{1}{2}(\varepsilon_1\pm\varepsilon_2\pm\varepsilon_3\pm\varepsilon_4)$.

The longest element of the Weyl group expressed as a product of ${\rm dim}~\h_w$ reflections with respect to mutually orthogonal roots: $w=s_{\varepsilon_1}s_{\varepsilon_2}s_{\varepsilon_3} s_{\varepsilon_4}$.

Normal ordering of $\Delta_+(\m_w,\h_w)$ compatible with expression $w=s_{\varepsilon_1}s_{\varepsilon_2}s_{\varepsilon_3} s_{\varepsilon_4}$:

$$
\alpha_4,\varepsilon_1-\varepsilon_2,\ldots,\varepsilon_{3}-\varepsilon_4,\ldots,\varepsilon_1,\ldots,\varepsilon_2,\ldots,
\varepsilon_4,
$$
where the roots $\varepsilon_i\pm \varepsilon_j$ $(1\leq i<j\leq l)$ forming the subsystem $\Delta_+(B_{4})\subset \Delta_+(F_4)$ are situated as in case of $B_4$, and a half of the positive roots which do not belong to the subsystem $\Delta_+(B_4)\subset \Delta_+(F_4)$ are situated to the left from $\varepsilon_{1}$ and the other half of those roots are situated to the right from $\varepsilon_{1}$.

\vskip 0.5cm

\item

$G_2$

\vskip 0.2cm

Dynkin diagram:

$$\xymatrix@R=.25cm{
\alpha_1&\alpha_2\\
{\bullet}\ar@3{-}[r] &{\bullet}
}$$

Simple roots: $\alpha_1=\varepsilon_1-\varepsilon_2,\alpha_2=-2\varepsilon_1+\varepsilon_2+\varepsilon_3$.

Positive roots: $\alpha_1,\alpha_1+\alpha_2, 2\alpha_1+\alpha_2, 3\alpha_1+\alpha_2,3\alpha_1+2\alpha_2,\alpha_2$.

The longest element of the Weyl group expressed as a product of ${\rm dim}~\h_w$ reflections with respect to mutually orthogonal roots: $w=s_{\alpha_1}s_{3\alpha_1+2\alpha_2}$.

Normal ordering of $\Delta_+(\m_w,\h_w)$ compatible with expression $w=s_{\alpha_1}s_{3\alpha_1+2\alpha_2}$:

$$
\alpha_2,\alpha_1+\alpha_2,3\alpha_1+2\alpha_2,2\alpha_1+\alpha_2,3\alpha_1+\alpha_2,\alpha_1.
$$

\end{itemize}


\section*{Appendix 2. Transversal slices for simple exceptional algebraic groups.}

\addcontentsline{toc}{section}{Appendix 2. Transversal slices for simple exceptional algebraic groups}

\renewcommand{\thetheorem}{A2.\arabic{theorem}}

\renewcommand{\thelemma}{A2.\arabic{lemma}}

\renewcommand{\theproposition}{A2.\arabic{proposition}}

\renewcommand{\thecorollary}{A2.\arabic{corollary}}

\renewcommand{\theremark}{A2.\arabic{remark}}

\renewcommand{\thedefinition}{A2.\arabic{definition}}

\renewcommand{\theequation}{A2.\arabic{equation}}

\setcounter{equation}{0}
\setcounter{theorem}{0}

\pagestyle{myheadings}
\markboth{APPENDIX}{APPENDIX 2. TRANSVERSAL SLICES FOR EXCEPTIONAL ALGEBRAIC GROUPS}

In this appendix, for simple exceptional algebraic groups we present the data related to the varieties $\Sigma_{{\bf k}, s}$ defined in Theorem \ref{mainth}.
Let $G_{\bf k}$ be a connected simple algebraic group of an exceptional type over an algebraically closed field ${\bf k}$, and $\mathcal{O}\in \widehat{\underline{\mathcal{N}}}(W)$. Let $H_{\bf k}$ be a maximal torus of $G_{\bf k}$, $W$ the Weyl group of the pair $(G_{\bf k},H_{\bf k})$, and $s\in W$ an element from the conjugacy class $\Psi^W(\mathcal{O})$. 

Let $\Delta$ be the root system of the pair $(G_{\bf k},H_{\bf k})$ and $\Delta_+^s$ the system of positive roots in $\Delta$ associated to $s$ and defined in Section \ref{background} with the help of decomposition (\ref{hdec}) where the subspaces $\h_i$ are ordered in such a way that in sum (\ref{hdec}) $\h_0$ is the vector subspace of $\h_{\mathbb{R}}$ fixed by the action of $s$, the one--dimensional subspaces $\h_i$ on which $s^1$ acts by multiplication by $-1$ follow immediately after $\h_0$ in (\ref{hdec}), and if $\h_i=\h_\lambda^k$, $\h_j=\h_\mu^l$ and $0\leq \lambda <\mu< 1$ then $i<j$. We also use a decomposition $s=s^1s^2$ for which the direct sum $\bigoplus_{k=0,i_k> 0}^r\h_{i_k}$ of  the one--dimensional subspaces $\h_{i_k}$ on which $s^1$ acts by multiplication by $-1$ is trivial. Such decomposition always exists. As a consequence, condition (\ref{cond2}) is satisfied. 
Let $\Sigma_{{\bf k}, s}$ be the corresponding  variety defined in Proposition \ref{crosssect} (ii).

Then straightforward calculation shows that
$$
{\rm dim}~Z_{G_p}(n)={\rm dim}~\Sigma_{{\bf k}, s}
$$
for any $n\in \mathcal{O}\in \underline{\mathcal{N}}(G_p)\subset \widehat{\underline{\mathcal{N}}}(W)$. The numbers ${\rm dim}~Z_{G_p}(n)$ can be found in \cite{Li}, Chapter 22 (note, however, that the notation in \cite{Li} for some classes is different from ours; we follow \cite{L3',Spal1}). The numbers ${\rm dim}~\Sigma_{{\bf k}, s}$ are contained in the tables below. These two numbers coincide in all cases. The tables below contain also the following information for each $\mathcal{O}\in \widehat{\underline{\mathcal{N}}}(W)$:
\\
-- The Weyl group conjugacy class $\Psi^W(\mathcal{O})$ which can be found in \cite{L3'};
\\
-- The two involutions $s^1$ and $s^2$ in the decomposition $s=s^1s^2\in \Psi^W(\mathcal{O})$; they are represented by sets of natural numbers which are the numbers of the roots appearing in the decompositions $s^1=s_{\gamma_1}\ldots s_{\gamma_{\widetilde{l}}}$, $s^2=s_{\gamma_{\widetilde{l}+1}}\ldots s_{\gamma_{l'}}$, where the system of positive roots $\Delta_+^s$ is chosen as in Theorem \ref{mainth}, and the numeration of positive roots is given in Appendix 3;
\\
-- The dimension of the fixed point space $\h_0$ for the action of $s$ on $\h$;
\\
-- The number $|\Delta_0|$ of roots fixed by $s$;
\\
-- The type of the root system $\Delta_0$ fixed by $s$;
\\
-- The Dynkin diagram $\Gamma_0^s$ of $\Delta_0$, where the numbers at the vertices of $\Gamma_0^s$ are the numbers of simple roots in $\Delta_+^s$ which appear in $\Gamma_0^s$; the numeration of simple roots is given in Appendix 3;
\\
-- The length $l(s)$ of $s$ with respect to the system of simple roots in $\Delta_+^s$;
\\
-- ${\rm dim}~\Sigma_{{\bf k}, s}={\rm dim}~\h_0+|\Delta_0|+l(s)$;
\\
-- The lowest common multiple $d'$ of the denominators of the numbers $\frac{1}{d_j}\left\langle {1+s \over 1-s }P_{{\h'}^*}\alpha_i,\alpha_j\right\rangle$, where $i,j=1,\ldots,l$;

\vskip 0.4cm

$\bf G_2.$

\vskip 0.1cm

\begin{tabular*}{1.03\textwidth}{|@{\extracolsep{\fill} }c|c|c|c|c|c|c|c|c|c|c|}
  \hline
 \rule[-.3cm]{0cm}{1cm}{}$\mathcal{O}$ & $\Psi^W(\mathcal{O})$ & $s^1$ & $s^2$ & ${\rm dim}~\h_0$ & $|\Delta_0|$ & $\Delta_0$ & $\Gamma_0^s$ & $l(s)$ & ${\rm dim}~\Sigma_{{\bf k}, s}$ & $d'$ \\ \hline \hline
$A_1$ & $A_1$ & -- & 6 & 1 & 2 & $A_1$ & $\xymatrix@R=.1cm@C=.1cm{
{\scriptstyle 1} \\
{\bullet}
}$
& 5 & 8 & 1  \\
  \hline
$(\tilde{A}_1)_3$ \rule[-.3cm]{0cm}{1cm}{} & $\tilde{A}_1$ & -- & 4 & 1 & 2 & $A_1$ &
$\xymatrix@R=.1cm@C=.1cm{
{\scriptstyle 2} \\
{\bullet}
}$
& 5 & 8 & 1  \\
  \hline
$\tilde{A}_1$ \rule[-.3cm]{0cm}{1cm}{} & $A_1+\tilde{A}_1$ & -- & $\begin{array}{c}
      1 \\
      6
    \end{array}$ & 0 & 0 & -- & -- & 6 & 6 & 1 \\
  \hline
$G_2(a_1) \rule[-.2cm]{0cm}{0.6cm}{} $ & $A_2$ & 5 & 2 & 0 & 0 & -- & -- & 4 & 4 & 3  \\
  \hline
$G_2$ \rule[-.2cm]{0cm}{0.6cm}{} & $G_2$ & 1 & 2 & 0 & 0 & -- & -- & 2 & 2 & 1  \\
  \hline
\end{tabular*}

\newpage

$\bf F_4.$

\vskip 0.1cm

\begin{tabular*}{1.03\textwidth}{|@{\extracolsep{\fill} }c|c|c|c|c|c|c|c|c|c|c|}
  \hline
 \rule[-.3cm]{0cm}{1cm}{} $\mathcal{O}$ & $\Psi^W(\mathcal{O})$ & $s^1$ & $s^2$ & ${\rm dim}~\h_0$ & $|\Delta_0|$ & $\Delta_0$ & $\Gamma_0^s$  & $l(s)$ & ${\rm dim}~\Sigma_{{\bf k}, s}$ & $d'$ \\ \hline \hline
 $A_1$ & $A_1$ & -- & 24 & 3 & 18 & $C_3$ &

$\xymatrix@R=.1cm@C=.1cm{
{\scriptstyle 4}&&{\scriptstyle 3}&&{\scriptstyle 2}\\
{\bullet}\ar@{-}[rr]&&{\bullet}\ar@2{-}[rr]&&{\bullet}
}$
 & 15 & 36 & 1 \\
  \hline
$(\tilde{A}_1)_2$ \rule[-.3cm]{0cm}{1cm}{} & $\tilde{A}_1$ & -- & 21 & 3 & 18 & $B_3$ &

$\xymatrix@R=.1cm@C=.1cm{
{\scriptstyle 1}&&{\scriptstyle 2}&&{\scriptstyle 3}\\
{\bullet}\ar@{-}[rr]&&{\bullet}\ar@2{-}[rr]&&{\bullet}
}$
& 15 & 36 & 1 \\
\hline
 $\tilde{A}_1$ \rule[-.3cm]{0cm}{1cm}{} & $2A_1$ & -- & $\begin{array}{c}
     16 \\
     24
   \end{array}
$ & 2 & 8 & $B_2$ &
$\xymatrix@R=.1cm@C=.1cm{
{\scriptstyle 2}&&{\scriptstyle 3}\\
{\bullet}\ar@2{-}[rr]&&{\bullet}
}$
& 20 & 30 & 1 \\
  \hline
$A_1+\tilde{A}_1$ & $4A_1$ & -- & $\begin{array}{c}
      5 \\
      11 \\
			18 \\
			23
    \end{array}$
 & 0 & 0 & -- & -- & 24 & 24 & 1 \\
  \hline
$A_2$ & $A_2$ & 23 & 1 & 2 & 6 & $A_2$ &
$\xymatrix@R=.1cm@C=.1cm{
{\scriptstyle 3}&&{\scriptstyle 4}\\
{\bullet}\ar@{-}[rr]&&{\bullet}
}$
& 14 & 22 & 3 \\
  \hline
$\tilde{A}_2$ \rule[-.3cm]{0cm}{1cm}{} & $\tilde{A}_2$ & 19 & 4 & 2 & 6 & $A_2$ &
$\xymatrix@R=.1cm@C=.1cm{
{\scriptstyle 1}&&{\scriptstyle 2}\\
{\bullet}\ar@{-}[rr]&&{\bullet}
}$
& 14 & 22 & 3 \\
\hline
 $(B_2)_2$ & $B_2$ & 16 & 8 & 2 & 8 & $B_2$ &
$\xymatrix@R=.1cm@C=.1cm{
{\scriptstyle 2}&&{\scriptstyle 3}\\
{\bullet}\ar@2{-}[rr]&&{\bullet}
}$
& 10 & 20 & 2 \\
\hline
$A_2+\tilde{A}_1$ & $A_2+\tilde{A}_1$ & 23 & $\begin{array}{c}
    1  \\
     7
   \end{array}
$ & 1 & 0 & -- & -- & 17 & 18 & 3 \\
\hline
$(\tilde{A}_2+A_1)_2$ & $\tilde{A}_2+A_1$ & 19 & $\begin{array}{c}
       4 \\
       5
     \end{array}
$ & 1 & 0 & -- & -- & 17 & 18 & 3 \\
\hline
$\tilde{A}_2+A_1$ & $A_2+\tilde{A}_2$ & $\begin{array}{c}
    23  \\
     3
   \end{array}
$ & $\begin{array}{c}
    1  \\
     4
   \end{array}
$ & 0 & 0 & -- & -- & 16 & 16 & 3 \\
\hline
$B_2$ & $A_3$ & 16 & $\begin{array}{c}
    1  \\
     14
   \end{array}
$ & 1 & 2 & $A_1$ &
$\xymatrix@R=.1cm@C=.1cm{
{\scriptstyle 3}\\
{\bullet}
}$
& 13 & 16 & 2 \\
\hline
$(C_3(a_1))_2$ & $B_2+A_1$ & 16 & $\begin{array}{c}
    8  \\
     9
   \end{array}
$ & 1 & 2 & $A_1$ &
$\xymatrix@R=.1cm@C=.1cm{
{\scriptstyle 2}\\
{\bullet}
}$
& 13 & 16 & 2 \\
\hline
$C_3(a_1)$ & $A_3+\tilde{A}_1$ & 16 & $\begin{array}{c}
    5  \\
     6 \\
		11
   \end{array}
$ & 0 & 0 & -- & -- & 14 & 14 & 2 \\
\hline
$F_4(a_3)$ & $D_4(a_1)$ & $\begin{array}{c}
    16  \\
     2
   \end{array}
$ & $\begin{array}{c}
    5  \\
     11
   \end{array}
$ & 0 & 0 & -- & -- & 12 & 12 & 1 \\
\hline
$B_3$ & $D_4$ & 1 & $\begin{array}{c}
    16  \\
     9 \\
     2
   \end{array}
$ & 0 & 0 & -- & -- & 10 & 10 & 1 \\
\hline
$C_3$ & $C_3+A_1$ & 4 & $\begin{array}{c}
    1 \\
    3 \\
		14
   \end{array}
$ & 0 & 0 & -- & -- & 10 & 10 & 1 \\
\hline
$F_4(a_2)$ & $F_4(a_1)$ & $\begin{array}{c}
    1  \\
     3
   \end{array}
$ & $\begin{array}{c}
    9  \\
     10
   \end{array}
$ & 0 & 0 & -- & -- & 8 & 8 & 1 \\
\hline
$F_4(a_1)$ & $B_4$ & $\begin{array}{c}
    9  \\
     2
   \end{array}
$ & $\begin{array}{c}
    1  \\
     4
   \end{array}
$ & 0 & 0 & -- & -- & 6 & 6 & 1 \\
\hline
$F_4$ & $F_4$ & $\begin{array}{c}
    1  \\
     3
   \end{array}
$ & $\begin{array}{c}
    2  \\
     4
   \end{array}
$ & 0 & 0 & -- & -- & 4 & 4 & 1 \\
\hline
\end{tabular*}

\newpage

$\bf E_6.$

\vskip 0.1cm

\begin{tabular*}{1.03\textwidth}{|@{\extracolsep{\fill} }c|c|c|c|c|c|c|c|c|c|c|}
  \hline
\rule[-.3cm]{0cm}{1cm}{}  $\mathcal{O}$ & $\Psi^W(\mathcal{O})$ & $s^1$ & $s^2$ & ${\rm dim}~\h_0$ & $|\Delta_0|$ & $\Delta_0$ & $\Gamma_0^s$ & $l(s)$ & ${\rm dim}~\Sigma_{{\bf k}, s}$ & $d'$  \\ \hline \hline
$A_1$ & $A_1$ & -- & 36 & 5 & 30 & $A_5$ &
$\xymatrix@R=.1cm@C=.1cm{
{\scriptstyle 1}&&{\scriptstyle 3}&&{\scriptstyle 4}&&{\scriptstyle 5}&&{\scriptstyle 6}\\
{\bullet}\ar@{-}[rr]&&{\bullet}\ar@{-}[rr]&&{\bullet}\ar@{-}[rr]&& {\bullet}\ar@{-}[rr] &&{\bullet}
}$
& 21 & 56 & 1 \\
\hline
$2A_1$ & $2A_1$ & -- & $\begin{array}{c}
    23 \\ 36
   \end{array}
$ & 4 & 12 & $A_3$ &

$\xymatrix@R=.1cm@C=.1cm{
{\scriptstyle 3}&&{\scriptstyle 4}&&{\scriptstyle 5}\\
{\bullet}\ar@{-}[rr]&&{\bullet}\ar@{-}[rr]&&{\bullet}
}$
& 30 & 46 & 1 \\
\hline
$3A_1$ & $4A_1$ & -- & $\begin{array}{c}
    8 \\ 19 \\ 27 \\ 35
   \end{array}
$ & 2 & 0 & -- & -- & 36 & 38 & 1 \\
\hline
$A_2$ & $A_2$ & 35 & 2 & 4 & 12 & $2A_2$ &
$\xymatrix@R=.1cm@C=.1cm{
{\scriptstyle 1}&&{\scriptstyle 3}\\
{\bullet}\ar@{-}[rr]&&{\bullet}
}$

$\xymatrix@R=.1cm@C=.1cm{
{\scriptstyle 5}&&{\scriptstyle 6}\\
{\bullet}\ar@{-}[rr]&&{\bullet}
}$
& 20 & 36 & 3 \\
\hline
$A_2+A_1$ & $A_2+A_1$ & 35 & $\begin{array}{c}
    2 \\ 7
   \end{array}
$ & 3 & 6 & $A_2$ &
$\xymatrix@R=.1cm@C=.1cm{
{\scriptstyle 5}&&{\scriptstyle 6}\\
{\bullet}\ar@{-}[rr]&&{\bullet}
}$
& 23 & 32 & 3 \\
\hline
$2A_2$ & $2A_2$ & $\begin{array}{c}
    1 \\ 35
   \end{array}
$ & $\begin{array}{c}
    2 \\ 3
   \end{array}
$ & 2 & 6 & $A_2$ &
$\xymatrix@R=.1cm@C=.1cm{
{\scriptstyle 5}&&{\scriptstyle 6}\\
{\bullet}\ar@{-}[rr]&&{\bullet}
}$
& 22 & 30 & 3 \\
\hline
$A_2+2A_1$ & $A_2+2A_1$ & 35 & $\begin{array}{c}
    2 \\ 7 \\ 11
   \end{array}
$ & 2 & 0 & -- & -- & 26 & 28 & 3 \\
\hline
$A_3$ & $A_3$ & 23 & $\begin{array}{c}
    2 \\ 24
   \end{array}
$ & 3 & 4 & $2A_1$ &
$\xymatrix@R=.1cm@C=.1cm{
{\scriptstyle 3}&&{\scriptstyle 5}\\
{\bullet}&&{\bullet}
}$
& 19 & 26 & 2 \\
\hline
$2A_2+A_1$ & $3A_2$ & $\begin{array}{c}
     1 \\ 6 \\ 35
   \end{array}
$ & $\begin{array}{c}
    2 \\ 3 \\ 5
   \end{array}
$ & 0 & 0 & -- & -- & 24 & 24 & 3 \\
\hline
$A_3+A_1$ & $A_3+2A_1$ & 23 & $\begin{array}{c}
  2 \\  3 \\ 5 \\ 24
   \end{array}
$ & 1 & 0 & -- & -- & 21 & 22 & 2 \\
\hline
$D_4(a_1)$ & $D_4(a_1)$ & $\begin{array}{c}
    4 \\ 23
   \end{array}
$ & $\begin{array}{c}
    8 \\ 19
   \end{array}
$ & 2 & 0 & -- & -- & 18 & 20 & 2 \\
\hline
$A_4$ & $A_4$ & $\begin{array}{c}
     21 \\ 24
   \end{array}
$ & $\begin{array}{c}
    1 \\ 2
   \end{array}
$ & 2 & 2 & $A_1$ &
$\xymatrix@R=.1cm@C=.1cm{
{\scriptstyle 5}\\
{\bullet}
}$
& 14 & 18 & 5 \\
\hline
$D_4$ & $D_4$ & 2 & $\begin{array}{c}
    23 \\ 4 \\ 15
   \end{array}
$ & 2 & 0 & -- & -- & 16 & 18 & 1  \\
\hline
$A_4+A_1$ & $A_4+A_1$ & $\begin{array}{c}
    21 \\ 24
   \end{array}
$ & $\begin{array}{c}
    1 \\ 2 \\ 5
   \end{array}
$ & 1 & 0 & -- & -- & 15 & 16 & 5 \\
\hline
$A_5$ & $A_5+A_1$ & $\begin{array}{c}
    1 \\ 6
   \end{array}
$ & $\begin{array}{c}
    8 \\ 9 \\ 10 \\ 19
   \end{array}
$ & 0 & 0 & -- & -- & 14 & 14 & 1 \\
\hline
$D_5(a_1)$ & $D_5(a_1)$ & $\begin{array}{c}
    2 \\ 7
   \end{array}
$ & $\begin{array}{c}
    15 \\ 12 \\ 16
   \end{array}
$ & 1 & 0 & -- & -- & 13 & 14 & 2 \\
\hline
$A_5+A_1$ & $E_6(a_2)$ & $\begin{array}{c}
    1 \\ 2 \\ 6
   \end{array}
$ & $\begin{array}{c}
    9 \\ 10 \\ 19
   \end{array}
$ & 0 & 0 & -- & -- & 12 & 12 & 1 \\
\hline
\end{tabular*}

\newpage

\begin{tabular*}{1.03\textwidth}{|@{\extracolsep{\fill} }c|c|c|c|c|c|c|c|c|c|c|}
  \hline
\rule[-.3cm]{0cm}{1cm}{}  $\mathcal{O}$ & $\Psi^W(\mathcal{O})$ & $s^1$ & $s^2$ & ${\rm dim}~\h_0$ & $|\Delta_0|$ & $\Delta_0$ & $\Gamma_0^s$ & $l(s)$ & ${\rm dim}~\Sigma_{{\bf k}, s}$ & $d'$  \\ \hline \hline
$D_5$ & $D_5$ & $\begin{array}{c}
    4 \\ 15
   \end{array}
$ & $\begin{array}{c}
    1 \\ 2 \\ 6
   \end{array}
$ & 1 & 0 & -- & -- & 9 & 10 & 2 \\
\hline
$E_6(a_1)$ & $E_6(a_1)$ & $\begin{array}{c}
    6 \\ 8 \\ 9
   \end{array}
$ & $\begin{array}{c}
    1 \\ 2 \\ 5
   \end{array}
$ & 0 & 0 & -- & -- & 8 & 8 & 1 \\
\hline
$E_6$ & $E_6$ & $\begin{array}{c}
    1 \\ 4 \\ 6
   \end{array}
$ & $\begin{array}{c}
    2 \\ 3 \\ 5
   \end{array}
$ & 0 & 0 & -- & -- & 6 & 6 & 1 \\
\hline
\end{tabular*}

\vskip 0.4cm

$\bf E_7.$

\vskip 0.1cm

\begin{tabular*}{1.03\textwidth}{|@{\extracolsep{\fill} }c|c|c|c|c|c|c|c|c|c|c|}
  \hline
\rule[-.3cm]{0cm}{1cm}{}  $\mathcal{O}$ & $\Psi^W(\mathcal{O})$ & $s^1$ & $s^2$ & ${\rm dim}~\h_0$ & $|\Delta_0|$ & $\Delta_0$ & $\Gamma_0^s$ & $l(s)$ & ${\rm dim}~\Sigma_{{\bf k}, s}$ & $d'$  \\ \hline \hline
$A_1$ & $A_1$ & -- & 63 & 6 & 60 & $D_6$ &
$\xymatrix@R=.01cm@C=.1cm{
&&&&&&&{\scriptstyle 3}\\
&&&&&&&{\bullet}\\
{\scriptstyle 7}&&{\scriptstyle 6}&&{\scriptstyle 5}&&{\scriptstyle 4}&\\
{\bullet}\ar@{-}[rr]&&{\bullet}\ar@{-}[rr]&&{\bullet}\ar@{-}[rr]&& {\bullet}\ar@{-}[uur]\ar@{-}[ddr]& \\
&&&&&&&\\
&&&&&&&{\bullet}\\
&&&&&&&{\scriptstyle 2}
}$
& 33 & 99 & 1 \\
\hline
$2A_1$ & $2A_1$ & -- & $\begin{array}{c}
    49 \\ 63
   \end{array}
$ & 5 & 26 & $A_1+D_4$ &
$\xymatrix@R=.01cm@C=.1cm{
&&&&&{\scriptstyle 3}\\
&&&&&{\bullet}\\
{\scriptstyle 7}&&{\scriptstyle 2}&&{\scriptstyle 4}&\\
{\bullet}&&{\bullet}\ar@{-}[rr]&& {\bullet}\ar@{-}[uur]\ar@{-}[ddr]& \\
&&&&&\\
&&&&&{\bullet}\\
&&&&&{\scriptstyle 5}
}$ & 50 & 81 & 1  \\
\hline
$(3A_1)''$ & $(3A_1)'$ & -- & $\begin{array}{c}
    7 \\ 49 \\ 63
   \end{array}
$ & 4 & 24 & $D_4$ & $\xymatrix@R=.01cm@C=.1cm{
&&&{\scriptstyle 3}\\
&&&{\bullet}\\
{\scriptstyle 2}&&{\scriptstyle 4}&\\
{\bullet}\ar@{-}[rr]&& {\bullet}\ar@{-}[uur]\ar@{-}[ddr]& \\
&&&\\
&&&{\bullet}\\
&&&{\scriptstyle 5}
}$ & 51 & 79 & 1  \\
\hline
$(3A_1)'$ & $(4A_1)''$ & -- & $\begin{array}{c}
    19 \\ 40 \\ 41 \\ 63
   \end{array}
$ & 3 & 6 & $3A_1$ & $\xymatrix@R=.1cm@C=.1cm{
{\scriptstyle 2}&&{\scriptstyle 3}&&{\scriptstyle 6}\\
{\bullet}&&{\bullet}&&{\bullet}
}$ & 60 & 69 & 1  \\
\hline
$A_2$ & $A_2$ & 62 & 1 & 5 & 30 & $A_5$ &
$\xymatrix@R=.1cm@C=.1cm{
{\scriptstyle 2}&&{\scriptstyle 4}&&{\scriptstyle 5}&&{\scriptstyle 6}&&{\scriptstyle 7}\\
{\bullet}\ar@{-}[rr]&&{\bullet}\ar@{-}[rr]&&{\bullet}\ar@{-}[rr]&& {\bullet}\ar@{-}[rr] &&{\bullet}
}$
 & 32 & 67 & 3  \\
\hline
$4A_1$ & $7A_1$ & -- & $\begin{array}{c}
    21 \\ 62 \\ 33 \\ 44 \\ 18 \\ 19 \\ 16
   \end{array}
$ & 0 & 0 & -- & -- & 63 & 63 & 1  \\
\hline
$A_2+A_1$ & $A_2+A_1$ & 62 & $\begin{array}{c}
    1 \\ 30
   \end{array}
$ & 4 & 12 & $A_3$ & $\xymatrix@R=.1cm@C=.1cm{
{\scriptstyle 4}&&{\scriptstyle 5}&&{\scriptstyle 6}\\
{\bullet}\ar@{-}[rr]&& {\bullet}\ar@{-}[rr] &&{\bullet}
}$ & 41 & 57 & 3 \\
\hline
$A_2+2A_1$ & $A_2+2A_1$ & 62 & $\begin{array}{c}
    1 \\ 18 \\ 30
   \end{array}
$ & 3 & 2 & $A_1$ & $\xymatrix@R=.1cm@C=.1cm{
{\scriptstyle 5}\\
{\bullet}
}$ & 46 & 51 & 3  \\
\hline
\end{tabular*}

\newpage

\begin{tabular*}{1.04\textwidth}{|@{\extracolsep{\fill} }c|c|c|c|c|c|c|c|c|c|c|}
  \hline
\rule[-.3cm]{0cm}{1cm}{}  $\mathcal{O}$ & $\Psi^W(\mathcal{O})$ & $s^1$ & $s^2$ & ${\rm dim}~\h_0$ & $|\Delta_0|$ & $\Delta_0$ & $\Gamma_0^s$ & $l(s)$ & ${\rm dim}~\Sigma_{{\bf k}, s}$ & $d'$  \\ \hline \hline
$A_2+3A_1$ & $A_2+3A_1$ & 62 & $\begin{array}{c}
    1 \\ 5 \\ 18 \\ 30
   \end{array}
$ & 2 & 0 & -- & -- & 47 & 49 & 3  \\
\hline
$2A_2$ & $2A_2$ & $\begin{array}{c}
    23 \\62
   \end{array}
$ & $\begin{array}{c}
    1 \\ 7
   \end{array}
$ & 3 & 6 & $A_2$ & $\xymatrix@R=.1cm@C=.1cm{
{\scriptstyle 4}&&{\scriptstyle 5}\\
{\bullet}\ar@{-}[rr] &&{\bullet}
}$ & 40 & 49 & 3   \\
\hline
$A_3$ & $A_3$ & 49 & $\begin{array}{c}
    1 \\ 37
   \end{array}
$ & 4 & 14 & $A_1+A_3$ & $\xymatrix@R=.1cm@C=.1cm{
{\scriptstyle 7}&&{\scriptstyle 2}&&{\scriptstyle 4}&&{\scriptstyle 5}\\
{\bullet}&&{\bullet}\ar@{-}[rr]&& {\bullet}\ar@{-}[rr] &&{\bullet}
}$ & 31 & 49 & 2  \\
\hline
$(A_3+A_1)''$ & $(A_3+A_1)'$ & 49 & $\begin{array}{c}
    1 \\ 7 \\ 37
   \end{array}
$ & 3 & 12 & $A_3$ & $\xymatrix@R=.1cm@C=.1cm{
{\scriptstyle 2}&&{\scriptstyle 4}&&{\scriptstyle 5}\\
{\bullet}\ar@{-}[rr]&& {\bullet}\ar@{-}[rr] &&{\bullet}
}$ & 32 & 47 & 2  \\
\hline
$2A_2+A_1$ & $3A_2$ & $\begin{array}{c}
    5 \\ 25 \\62
   \end{array}
$ & $\begin{array}{c}
    1 \\ 2 \\ 6
   \end{array}
$ & 1 & 0 & -- & -- & 42 & 43 & 3  \\
\hline
$(A_3+A_1)'$ & $(A_3+2A_1)''$ & 49 & $\begin{array}{c}
    7 \\ 14 \\ 26 \\ 28
   \end{array}
$ & 2 & 2 & $A_1$ & $\xymatrix@R=.1cm@C=.1cm{
{\scriptstyle 3}\\
{\bullet}
}$ & 37 & 41 & 2  \\
\hline
$A_3+2A_1$ & $A_3+3A_1$ & 49 & $\begin{array}{c}
    1 \\ 4 \\ 7 \\ 16 \\37
   \end{array}
$ & 1 & 0 & -- & -- & 38 & 39 & 2  \\
\hline
$D_4(a_1)$ & $D_4(a_1)$ & $\begin{array}{c}
    3 \\ 49
   \end{array}
$ & $\begin{array}{c}
    8 \\ 32
   \end{array}
$ & 3 & 6 & $3A_1$ & $\xymatrix@R=.1cm@C=.1cm{
{\scriptstyle 2} && {\scriptstyle 5} && {\scriptstyle 7}\\
{\bullet} && {\bullet} && {\bullet}
}$ & 30 & 39 & 2  \\
\hline
$D_4(a_1)+A_1$ & $D_4(a_1)+A_1$ & $\begin{array}{c}
    3 \\ 49
   \end{array}
$ & $\begin{array}{c}
   7 \\ 8 \\ 32
   \end{array}
$ & 2 & 4 & $2A_1$ & $\xymatrix@R=.1cm@C=.1cm{
{\scriptstyle 2} && {\scriptstyle 5} \\
{\bullet} && {\bullet}
}$ & 31 & 37 & 2  \\
\hline
$D_4$ & $D_4$ & 1 & $\begin{array}{c}
    3 \\ 28 \\ 49
   \end{array}
$ & 3 & 6 & $3A_1$ & $\xymatrix@R=.1cm@C=.1cm{
{\scriptstyle 2} && {\scriptstyle 5} && {\scriptstyle 7}\\
{\bullet} && {\bullet} && {\bullet}
}$ & 28 & 37 & 1  \\
\hline
$(A_3+A_2)_2$ & $A_3+A_2$ & $\begin{array}{c}
    22 \\ 49
   \end{array}
$ & $\begin{array}{c}
    4 \\ 20 \\ 21
   \end{array}
$ & 2 & 2 & $A_1$ & $\xymatrix@R=.1cm@C=.1cm{
{\scriptstyle 7}\\
{\bullet}
}$ & 33 & 37 & 6  \\
\hline
$A_3+A_2$ & $D_4(a_1)+2A_1$ & $\begin{array}{c}
    3 \\ 49
   \end{array}
$ & $\begin{array}{c}
    2 \\ 7 \\ 8 \\ 32
   \end{array}
$ & 1 & 2 & $A_1$ & $\xymatrix@R=.1cm@C=.1cm{
{\scriptstyle 5} \\
{\bullet}
}$ & 32 & 35 & 2  \\
\hline
$A_3+A_2+A_1$ & $2A_3+A_1$ & $\begin{array}{c}
    3 \\ 49
   \end{array}
$ & $\begin{array}{c}
    7 \\ 9 \\ 11 \\ 14 \\ 26
   \end{array}
$ & 0 & 0 & -- & -- & 33 & 33 & 2 \\
\hline
$A_4$ & $A_4$ & $\begin{array}{c}
    37 \\ 45
   \end{array}
$ & $\begin{array}{c}
    1 \\ 6
   \end{array}
$ & 3 & 6 & $A_2$ & $\xymatrix@R=.1cm@C=.1cm{
{\scriptstyle 2}&&{\scriptstyle 4}\\
{\bullet}\ar@{-}[rr]&& {\bullet}
}$ & 24 & 33 & 5  \\
\hline
$A_5''$ & $A_5'$ & $\begin{array}{c}
    6 \\ 40
   \end{array}
$ & $\begin{array}{c}
    7 \\ 20 \\ 21
   \end{array}
$ & 2 & 6 & $A_2$ & $\xymatrix@R=.1cm@C=.1cm{
{\scriptstyle 3}&&{\scriptstyle 4}\\
{\bullet}\ar@{-}[rr]&& {\bullet}
}$ & 23 & 31 & 3  \\
\hline
\end{tabular*}

\begin{tabular*}{1.03\textwidth}{|@{\extracolsep{\fill} }c|c|c|c|c|c|c|c|c|c|c|}
  \hline
\rule[-.3cm]{0cm}{1cm}{}  $\mathcal{O}$ & $\Psi^W(\mathcal{O})$ & $s^1$ & $s^2$ & ${\rm dim}~\h_0$ & $|\Delta_0|$ & $\Delta_0$ & $\Gamma_0^s$ & $l(s)$ & ${\rm dim}~\Sigma_{{\bf k}, s}$ & $d'$  \\ \hline \hline
$D_4+A_1$ & $D_4+3A_1$ & 1 & $\begin{array}{c}
   5\\18\\28\\29\\30\\31
   \end{array}
$ & 0 & 0 & -- & -- & 31 & 31 & 1  \\
\hline
$A_4+A_1$ & $A_4+A_1$ & $\begin{array}{c}
    37 \\ 45
   \end{array}
$ & $\begin{array}{c}
    1 \\ 6 \\ 9
   \end{array}
$ & 2 & 0 & -- & -- & 27 & 29 & 5  \\
\hline
$A_4+A_2$ & $A_4+A_2$ & $\begin{array}{c}
    4 \\ 37 \\ 45
   \end{array}
$ & $\begin{array}{c}
    1 \\ 2 \\ 6
   \end{array}
$ & 1 & 0 & -- & -- & 26 & 27 & 15  \\
\hline
$D_5(a_1)$ & $D_5(a_1)$ & $\begin{array}{c}
    1 \\ 18 
   \end{array}
$ & $\begin{array}{c}
    24 \\ 28 \\ 36
   \end{array}
$ & 2 & 2 & $A_1$ & $\xymatrix@R=.1cm@C=.1cm{
{\scriptstyle 5}\\
{\bullet}
}$ & 23 & 27 & 2  \\
\hline
$(A_5+A_1)''$ & $A_5+A_2$ & $\begin{array}{c}
    22 \\ 23 \\ 24
   \end{array}
$ & $\begin{array}{c}
    4 \\ 7 \\ 20 \\ 21
   \end{array}
$ & 0 & 0 & -- & -- & 25 & 25 & 3  \\
\hline
$A_5'$ & $(A_5+A_1)''$ & $\begin{array}{c}
    6 \\ 19
   \end{array}
$ & $\begin{array}{c}
    8 \\ 15 \\ 17 \\ 32
   \end{array}
$ & 1 & 0 & -- & -- & 24 & 25 & 3  \\
\hline
$D_5(a_1)+A_1$ & $D_5(a_1)+A_1$ & $\begin{array}{c}
    1 \\ 18 
   \end{array}
$ & $\begin{array}{c}
    5 \\ 24 \\ 28 \\ 36
   \end{array}
$ & 1 & 0 & -- & -- & 24 & 25 & 2  \\
\hline
$D_6(a_2)$ & $D_6(a_2)+A_1$ & $\begin{array}{c}
   1 \\ 2 
   \end{array}
$ & $\begin{array}{c}
    4 \\ 6 \\ 15 \\ 31 \\ 40
   \end{array}
$ & 0 & 0 & -- & -- & 23 & 23 & 1  \\
\hline
$(A_5+A_1)'$ & $E_6(a_2)$ & $\begin{array}{c}
   1 \\ 4 \\ 16
   \end{array}
$ & $\begin{array}{c}
    28 \\ 29 \\ 31
   \end{array}
$ & 1 & 0 & -- & -- & 22 & 23 & 3  \\
\hline
$D_5$ & $D_5$ & $\begin{array}{c}
    3 \\ 28
   \end{array}
$ & $\begin{array}{c}
    1 \\ 6 \\ 19
   \end{array}
$ & 2 & 2 & $A_1$ & $\xymatrix@R=.1cm@C=.1cm{
{\scriptstyle 2}\\
{\bullet}
}$ & 17 & 21 & 2  \\
\hline
$D_6(a_2)+A_1$ & $E_7(a_4)$ & $\begin{array}{c}
   1 \\ 4 \\ 7
   \end{array}
$ & $\begin{array}{c}
    12 \\ 22 \\ 31 \\ 35
   \end{array}
$ & 0 & 0 & -- & -- & 21 & 21 & 1  \\
\hline
$D_5+A_1$ & $D_5+A_1$ & $\begin{array}{c}
   3 \\ 28
   \end{array}
$ & $\begin{array}{c}
    1 \\ 2 \\ 6 \\ 19
   \end{array}
$ & 1 & 0 & -- & -- & 18 & 19 & 2  \\
\hline
$A_6$ & $A_6$ & $\begin{array}{c}
   11 \\ 19 \\ 26
   \end{array}
$ & $\begin{array}{c}
    6 \\ 9 \\ 10
   \end{array}
$ & 1 & 0 & -- & -- & 18 & 19 & 7  \\
\hline
$D_6(a_1)$ & $D_6(a_1)$ & $\begin{array}{c}
    3 \\ 5 \\ 28
   \end{array}
$ & $\begin{array}{c}
    1 \\ 12 \\ 13
   \end{array}
$ & 1 & 2 & $A_1$ & $\xymatrix@R=.1cm@C=.1cm{
{\scriptstyle 2}\\
{\bullet}
}$ & 16 & 19 & 1  \\
\hline
\end{tabular*}

\newpage

\begin{tabular*}{1.03\textwidth}{|@{\extracolsep{\fill} }c|c|c|c|c|c|c|c|c|c|c|}
  \hline
 \rule[-.3cm]{0cm}{1cm}{} $\mathcal{O}$ & $\Psi^W(\mathcal{O})$ & $s^1$ & $s^2$ & ${\rm dim}~\h_0$ & $|\Delta_0|$ & $\Delta_0$ & $\Gamma_0^s$ & $l(s)$ & ${\rm dim}~\Sigma_{{\bf k}, s}$ & $d'$  \\ \hline \hline
$D_6(a_1)+A_1$ & $A_7$ & $\begin{array}{c}
  1 \\ 12 \\ 13
   \end{array}
$ & $\begin{array}{c}
    9 \\ 10 \\ 11 \\ 22
   \end{array}
$ & 0 & 0 & -- & -- & 17 & 17 & 1  \\
\hline
$D_6$ & $D_6+A_1$ & $\begin{array}{c}
   1 \\ 6 
   \end{array}
$ & $\begin{array}{c}
    2 \\ 3 \\ 5 \\ 7 \\ 28
   \end{array}
$ & 0 & 0 & -- & -- & 15 & 15 & 1  \\
\hline
$E_6(a_1)$ & $E_6(a_1)$ & $\begin{array}{c}
   8 \\ 19 \\ 22
   \end{array}
$ & $\begin{array}{c}
    1 \\ 4 \\ 6
   \end{array}
$ & 1 & 0 & -- & -- & 14 & 15 & 3  \\
\hline
$E_6$ & $E_6$ & $\begin{array}{c}
   3  \\ 6  \\ 19
   \end{array}
$ & $\begin{array}{c}
    1 \\ 9 \\ 11
   \end{array}
$ & 1 & 0 & -- & -- & 12 & 13 & 3  \\
\hline
$D_6+A_1$ & $E_7(a_3)$ & $\begin{array}{c}
    1\\2\\6
   \end{array}
$ & $\begin{array}{c}
   7\\10\\11\\22
   \end{array}
$ & 0 & 0 & -- & -- & 13 & 13 & 1 \\
\hline
$E_7(a_2)$ & $E_7(a_2)$ & $\begin{array}{c}
   1\\4\\16
   \end{array}
$ & $\begin{array}{c}
    2\\3\\12\\13
   \end{array}
$ & 0 & 0 & -- & -- & 11 & 11 & 1  \\
\hline
$E_7(a_1)$ & $E_7(a_1)$ & $\begin{array}{c}
   6\\9\\10
   \end{array}
$ & $\begin{array}{c}
    1\\2\\5\\7
   \end{array}
$ & 0 & 0 & -- & -- & 9 & 9 & 1  \\
\hline
$E_7$ & $E_7$ & $\begin{array}{c}
   1\\4\\6
   \end{array}
$ & $\begin{array}{c}
    2\\3\\5\\7
   \end{array}
$ & 0 & 0 & -- & -- & 7 & 7 & 1  \\
\hline
\end{tabular*}

\vskip 0.4cm

$\bf E_8.$

\vskip 0.1cm

\begin{tabular*}{1.03\textwidth}{|@{\extracolsep{\fill} }c|c|c|c|c|c|c|c|c|c|c|}
  \hline
\rule[-.3cm]{0cm}{1cm}{}  $\mathcal{O}$ & $\Psi^W(\mathcal{O})$ & $s^1$ & $s^2$ & ${\rm dim}~\h_0$ & $|\Delta_0|$ & $\Delta_0$ & $\Gamma_0^s$ & $l(s)$ & ${\rm dim}~\Sigma_{{\bf k}, s}$ & $d'$  \\ \hline \hline
$A_1$ & $A_1$ & -- & 120 & 7 & 126 & $E_7$ &
$\xymatrix@R=.1cm@C=.1cm{
{\scriptstyle 1}&&{\scriptstyle 3}&&{\scriptstyle 4}&&{\scriptstyle 5}&&{\scriptstyle 6}&&{\scriptstyle 7}\\
{\bullet}\ar@{-}[rr]&&{\bullet}\ar@{-}[rr]&&{\bullet}\ar@{-}[dd]\ar@{-}[rr]&& {\bullet}\ar@{-}[rr] &&{\bullet}\ar@{-}[rr]&&{\bullet}\\
&&&&&&&&&&\\
&&&&{\bullet}&&&&&&\\
&&&&{\scriptstyle 2}&&&&&&}$
& 57 & 190 & 1 \\
\hline
$2A_1$ & $2A_1$ & -- & $\begin{array}{c}
    97 \\ 120
   \end{array}
$ & 6 & 60 & $D_6$ &
$\xymatrix@R=.01cm@C=.1cm{
&&&&&&&{\scriptstyle 2}\\
&&&&&&&{\bullet}\\
{\scriptstyle 7}&&{\scriptstyle 6}&&{\scriptstyle 5}&&{\scriptstyle 4}&\\
{\bullet}\ar@{-}[rr]&&{\bullet}\ar@{-}[rr]&&{\bullet}\ar@{-}[rr]&& {\bullet}\ar@{-}[uur]\ar@{-}[ddr]& \\
&&&&&&&\\
&&&&&&&{\bullet}\\
&&&&&&&{\scriptstyle 3}
}$ & 90 & 156 & 1  \\
\hline
\end{tabular*}

\newpage

\begin{tabular*}{1.03\textwidth}{|@{\extracolsep{\fill} }c|c|c|c|c|c|c|c|c|c|c|}
  \hline
\rule[-.3cm]{0cm}{1cm}{}  $\mathcal{O}$ & $\Psi^W(\mathcal{O})$ & $s^1$ & $s^2$ & ${\rm dim}~\h_0$ & $|\Delta_0|$ & $\Delta_0$ & $\Gamma_0^s$ & $l(s)$ & ${\rm dim}~\Sigma_{{\bf k}, s}$ & $d'$  \\ \hline \hline
$3A_1$ & $(4A_1)'$ & -- & $\begin{array}{c}
    7 \\ 61 \\ 97 \\ 120
   \end{array}
$ & 4 & 24 & $D_4$ & $\xymatrix@R=.01cm@C=.1cm{
&&&{\scriptstyle 3}\\
&&&{\bullet}\\
{\scriptstyle 2}&&{\scriptstyle 4}&\\
{\bullet}\ar@{-}[rr]&& {\bullet}\ar@{-}[uur]\ar@{-}[ddr]& \\
&&&\\
&&&{\bullet}\\
&&&{\scriptstyle 5}
}$ & 108 & 136 & 1  \\
\hline
$A_2$ & $A_2$ & 119 & 8 & 6 & 72 & $E_6$ &
$\xymatrix@R=.1cm@C=.1cm{
{\scriptstyle 1}&&{\scriptstyle 3}&&{\scriptstyle 4}&&{\scriptstyle 5}&&{\scriptstyle 6}\\
{\bullet}\ar@{-}[rr]&&{\bullet}\ar@{-}[rr]&&{\bullet}\ar@{-}[dd]\ar@{-}[rr]&& {\bullet}\ar@{-}[rr] &&{\bullet}\\
&&&&&&&&\\
&&&&{\bullet}&&&&\\
&&&&{\scriptstyle 2}&&&&}$
 & 56 & 134 & 3  \\
\hline
$4A_1$ & $8A_1$ & -- & $\begin{array}{c}
    9 \\ 13 \\ 19 \\ 50 \\ 67 \\ 69 \\ 83 \\ 119
   \end{array}
$ & 0 & 0 & -- & -- & 120 & 120 & 1  \\
\hline
$A_2+A_1$ & $A_2+A_1$ & 119 & $\begin{array}{c}
    8 \\ 69 
   \end{array}
$ & 5 & 30 & $A_5$ & $\xymatrix@R=.1cm@C=.1cm{
{\scriptstyle 1}&&{\scriptstyle 3}&&{\scriptstyle 4}&&{\scriptstyle 5}&&{\scriptstyle 6}\\
{\bullet}\ar@{-}[rr]&& {\bullet}\ar@{-}[rr] &&{\bullet}\ar@{-}[rr]&& {\bullet}\ar@{-}[rr] &&{\bullet}
}$ & 77 & 112 & 3  \\
\hline
$A_2+2A_1$ & $A_2+2A_1$ & 119 & $\begin{array}{c}
    8 \\ 31 \\ 69
   \end{array}
$ & 4 & 12 & $A_3$ & $\xymatrix@R=.1cm@C=.1cm{
{\scriptstyle 3}&&{\scriptstyle 4}&&{\scriptstyle 5}\\
{\bullet}\ar@{-}[rr] &&{\bullet}\ar@{-}[rr] &&{\bullet}
}$ & 86 & 102 & 3 \\
\hline
$A_3$ & $A_3$ & 97 & $\begin{array}{c}
    8 \\ 74
   \end{array}
$ & 5 & 40 & $D_5$ &
$\xymatrix@R=.01cm@C=.1cm{
&&&&&{\scriptstyle 2}\\
&&&&&{\bullet}\\
{\scriptstyle 6}&&{\scriptstyle 5}&&{\scriptstyle 4}&\\
{\bullet}\ar@{-}[rr]&&{\bullet}\ar@{-}[rr]&& {\bullet}\ar@{-}[uur]\ar@{-}[ddr]& \\
&&&&&\\
&&&&&{\bullet}\\
&&&&&{\scriptstyle 3}
}$
& 55 & 100 & 2  \\
\hline
$A_2+3A_1$ & $A_2+4A_1$ & 119 & $\begin{array}{c}
    2 \\ 8 \\ 32 \\ 45 \\ 57
   \end{array}
$ & 2 & 0 & -- & -- & 92 & 94 & 3  \\
\hline
$2A_2$ & $2A_2$ & $\begin{array}{c}
    63 \\119
   \end{array}
$ & $\begin{array}{c}
    2 \\ 8
   \end{array}
$ & 4 & 12 & $2A_2$ & $\xymatrix@R=.1cm@C=.1cm{
{\scriptstyle 1}&&{\scriptstyle 3}&&{\scriptstyle 5}&&{\scriptstyle 6}\\
{\bullet}\ar@{-}[rr] &&{\bullet}&&{\bullet}\ar@{-}[rr] &&{\bullet}
}$ & 76 & 92 & 3   \\
\hline
$2A_2+A_1$ & $3A_2$ & $\begin{array}{c}
    6 \\ 63 \\119
   \end{array}
$ & $\begin{array}{c}
    2 \\ 5 \\ 8
   \end{array}
$ & 2 & 6 & $A_2$ & $\xymatrix@R=.1cm@C=.1cm{
{\scriptstyle 1}&&{\scriptstyle 3}\\
{\bullet}\ar@{-}[rr] &&{\bullet}
}$ & 78 & 86 & 3  \\
\hline
$A_3+A_1$ & $(A_3+2A_1)'$ & 97 & $\begin{array}{c}
    7 \\ 22 \\ 61 \\ 62
   \end{array}
$ & 3 & 12 & $A_3$ & $\xymatrix@R=.1cm@C=.1cm{
{\scriptstyle 2}&&{\scriptstyle 4}&&{\scriptstyle 3}\\
{\bullet}\ar@{-}[rr]&& {\bullet}\ar@{-}[rr] &&{\bullet}
}$ & 69 & 84 & 2  \\
\hline
$D_4(a_1)$ & $D_4(a_1)$ & $\begin{array}{c}
    7 \\ 97
   \end{array}
$ & $\begin{array}{c}
    15 \\ 68
   \end{array}
$ & 4 & 24 & $D_4$ & $\xymatrix@R=.01cm@C=.1cm{
&&&{\scriptstyle 3}\\
&&&{\bullet}\\
{\scriptstyle 2}&&{\scriptstyle 4}&\\
{\bullet}\ar@{-}[rr]&& {\bullet}\ar@{-}[uur]\ar@{-}[ddr]& \\
&&&\\
&&&{\bullet}\\
&&&{\scriptstyle 5}
}$
& 54 & 82 & 2  \\
\hline
\end{tabular*}

\newpage

\begin{tabular*}{1.03\textwidth}{|@{\extracolsep{\fill} }c|c|c|c|c|c|c|c|c|c|c|}
  \hline
 \rule[-.3cm]{0cm}{1cm}{} $\mathcal{O}$ & $\Psi^W(\mathcal{O})$ & $s^1$ & $s^2$ & ${\rm dim}~\h_0$ & $|\Delta_0|$ & $\Delta_0$ & $\Gamma_0^s$ & $l(s)$ & ${\rm dim}~\Sigma_{{\bf k}, s}$ & $d'$  \\ \hline \hline
$D_4$ & $D_4$ & 8 & $\begin{array}{c}
    7 \\ 61 \\ 97
   \end{array}
$ & 4 & 24 & $D_4$ & $\xymatrix@R=.01cm@C=.1cm{
&&&{\scriptstyle 3}\\
&&&{\bullet}\\
{\scriptstyle 2}&&{\scriptstyle 4}&\\
{\bullet}\ar@{-}[rr]&& {\bullet}\ar@{-}[uur]\ar@{-}[ddr]& \\
&&&\\
&&&{\bullet}\\
&&&{\scriptstyle 5}
}$
& 52 & 80 & 1  \\
\hline
$2A_2+2A_1$ & $4A_2$ & $\begin{array}{c}
  1 \\  6 \\ 63 \\119
   \end{array}
$ & $\begin{array}{c}
   2 \\ 3 \\ 5 \\ 8
   \end{array}
$ & 0 & 0 & -- & -- & 80 & 80 & 3  \\
\hline
$A_3+2A_1$ & $A_3+4A_1$ & 97 & $\begin{array}{c}
    5 \\ 7 \\ 32 \\ 36 \\ 50 \\ 61
   \end{array}
$ & 1 & 0 & -- & -- & 75 & 76 & 2  \\
\hline
$D_4(a_1)+A_1$ & $D_4(a_1)+A_1$ & $\begin{array}{c}
    7 \\ 97
   \end{array}
$ & $\begin{array}{c}
    15 \\ 32 \\ 68
   \end{array}
$ & 3 & 6 & $3A_1$ & $\xymatrix@R=.1cm@C=.1cm{
{\scriptstyle 2} &&{\scriptstyle 3} && {\scriptstyle 5} \\
{\bullet} &&{\bullet} && {\bullet}
}$ & 63 & 72 & 2  \\
\hline
$(A_3+A_2)_2$ & $A_3+A_2$ & $\begin{array}{c}
    55 \\ 97
   \end{array}
$ & $\begin{array}{c}
    6 \\ 29 \\ 56
   \end{array}
$ & 3 & 4 & $2A_1$ & $\xymatrix@R=.1cm@C=.1cm{
{\scriptstyle 2}&&{\scriptstyle 3}\\
{\bullet}&&{\bullet}
}$ & 65 & 72 & 6  \\
\hline
$A_3+A_2$ & $(2A_3)'$ & $\begin{array}{c}
    7 \\ 97
   \end{array}
$ & $\begin{array}{c}
    13 \\ 22 \\ 40 \\ 62
   \end{array}
$ & 2 & 4 & $2A_1$ & $\xymatrix@R=.1cm@C=.1cm{
{\scriptstyle 2} && {\scriptstyle 3}\\
{\bullet} && {\bullet}
}$ & 64 & 70 & 2  \\
\hline
$A_4$ & $A_4$ & $\begin{array}{c}
    74 \\ 93
   \end{array}
$ & $\begin{array}{c}
    1 \\ 8
   \end{array}
$ & 4 & 20 & $A_4$ & $\xymatrix@R=.1cm@C=.1cm{
{\scriptstyle 2}&&{\scriptstyle 4}&&{\scriptstyle 5}&&{\scriptstyle 6}\\
{\bullet}\ar@{-}[rr]&&{\bullet}\ar@{-}[rr]&&{\bullet}\ar@{-}[rr]&& {\bullet}
}$ & 44 & 68 & 5  \\
\hline
$A_3+A_2+A_1$ & $2A_3+2A_1$ & $\begin{array}{c}
    7 \\ 97
   \end{array}
$ & $\begin{array}{c}
    5 \\ 26 \\ 27 \\ 32 \\ 36 \\ 50
   \end{array}
$ & 0 & 0 & -- & -- & 66 & 66 & 2  \\
\hline
$D_4(a_1)+A_2$ & $D_4(a_1)+A_2$ & $\begin{array}{c}
    7 \\ 25 \\97
   \end{array}
$ & $\begin{array}{c}
    4 \\ 15 \\ 68
   \end{array}
$ & 2 & 0 & -- & -- & 62 & 64 & 6  \\
\hline
$D_4+A_1$ & $D_4+4A_1$ & 8 & $\begin{array}{c}
    9 \\ 13 \\ 25 \\ 35 \\ 59 \\ 63 \\ 80
   \end{array}
$ & 0 & 0 & -- & --
& 64 & 64 & 1  \\
\hline
$2A_3$ & $2D_4(a_1)$ & $\begin{array}{c}
    2 \\ 3 \\ 7 \\97
   \end{array}
$ & $\begin{array}{c}
    11 \\ 12 \\ 15 \\ 68
   \end{array}
$ & 0 & 0 & -- & -- & 60 & 60 & 1  \\
\hline
$A_4+A_1$ & $A_4+A_1$ & $\begin{array}{c}
    74 \\ 93
   \end{array}
$ & $\begin{array}{c}
    1 \\ 8 \\ 26
   \end{array}
$ & 3 & 6 & $A_2$ & $\xymatrix@R=.1cm@C=.1cm{
{\scriptstyle 4}&&{\scriptstyle 5}\\
{\bullet}\ar@{-}[rr]&& {\bullet}
}$ & 51 & 60 & 5  \\
\hline
\end{tabular*}

\begin{tabular*}{1.03\textwidth}{|@{\extracolsep{\fill} }c|c|c|c|c|c|c|c|c|c|c|}
  \hline
 \rule[-.3cm]{0cm}{1cm}{} $\mathcal{O}$ & $\Psi^W(\mathcal{O})$ & $s^1$ & $s^2$ & ${\rm dim}~\h_0$ & $|\Delta_0|$ & $\Delta_0$ & $\Gamma_0^s$ & $l(s)$ & ${\rm dim}~\Sigma_{{\bf k}, s}$ & $d'$ \\ \hline \hline
$D_5(a_1)$ & $D_5(a_1)$ & $\begin{array}{c}
    8\\31
   \end{array}
$ & $\begin{array}{c}
    39\\61\\75
   \end{array}
$ & 3 & 12 & $A_3$ & $\xymatrix@R=.1cm@C=.1cm{
{\scriptstyle 3}&&{\scriptstyle 4}&&{\scriptstyle 5}\\
{\bullet}\ar@{-}[rr]&&{\bullet}\ar@{-}[rr]&& {\bullet}
}$ & 43 & 58 & 2  \\
\hline
$(D_4+A_2)_2$ & $D_4+A_2$ & $\begin{array}{c}
    8 \\ 63
   \end{array}
$ & $\begin{array}{c}
    2 \\ 49 \\ 59 \\ 71
   \end{array}
$ & 2 & 0 & -- & -- & 54 & 56 & 3  \\
\hline
$A_4+2A_1$ & $A_4+2A_1$ & $\begin{array}{c}
    74 \\ 93
   \end{array}
$ & $\begin{array}{c}
    1 \\ 8 \\ 12 \\ 26
   \end{array}
$ & 2 & 0 & -- & -- & 54 & 56 & 5  \\
\hline
$A_4+A_2$ & $A_4+A_2$ & $\begin{array}{c}
    18 \\ 74 \\ 93
   \end{array}
$ & $\begin{array}{c}
    1 \\ 6 \\ 8
   \end{array}
$ & 2 & 2 & $A_1$ & $\xymatrix@R=.1cm@C=.1cm{
{\scriptstyle 4}\\
{\bullet}
}$ & 50 & 54 & 15  \\
\hline
$A_4+A_2+A_1$ & $A_4+A_2+A_1$ & $\begin{array}{c}
    18 \\ 74 \\ 93
   \end{array}
$ & $\begin{array}{c}
    1 \\ 4 \\ 6 \\ 8
   \end{array}
$ & 1 & 0 & -- & -- & 51 & 52 & 15  \\
\hline
$A_5$ & $(A_5+A_1)'$ & $\begin{array}{c}
    1 \\ 44
   \end{array}
$ & $\begin{array}{c}
    15 \\ 34 \\ 35 \\ 68
   \end{array}
$ & 2 & 6 & $A_2$ & $\xymatrix@R=.1cm@C=.1cm{
{\scriptstyle 4}&&{\scriptstyle 5}\\
{\bullet}\ar@{-}[rr]&& {\bullet}
}$ & 44 & 52 & 3 \\
\hline
$D_5(a_1)+A_1$ & $D_5(a_1)+A_1$ & $\begin{array}{c}
    8 \\ 31
   \end{array}
$ & $\begin{array}{c}
    19 \\ 39 \\ 61 \\ 75
   \end{array}
$ & 2 & 2 & $A_1$ & $\xymatrix@R=.1cm@C=.1cm{
{\scriptstyle 4}\\
{\bullet}
}$ & 48 & 52 & 2  \\
\hline
$D_4+A_2$ & $D_4+A_3$ & $\begin{array}{c}
    8 \\ 31
   \end{array}
$ & $\begin{array}{c}
    2 \\ 32 \\ 53 \\ 61 \\ 64
   \end{array}
$ & 1 & 0 & -- & --
& 49 & 50 & 2  \\
\hline
$(A_5+A_1)''$ & $E_6(a_2)$ & $\begin{array}{c}
   1 \\ 8\\44
   \end{array}
$ & $\begin{array}{c}
    34\\35\\68
   \end{array}
$ & 2 & 6 & $A_2$ & $\xymatrix@R=.1cm@C=.1cm{
{\scriptstyle 4}&&{\scriptstyle 5}\\
{\bullet}\ar@{-}[rr]&& {\bullet}
}$ & 42 & 50 & 3  \\
\hline
$D_5$ & $D_5$ & $\begin{array}{c}
    7 \\ 61
   \end{array}
$ & $\begin{array}{c}
    1 \\ 8 \\ 44
   \end{array}
$ & 3 & 12 & $A_3$ & $\xymatrix@R=.1cm@C=.1cm{
{\scriptstyle 2}&&{\scriptstyle 4}&&{\scriptstyle 5}\\
{\bullet}\ar@{-}[rr]&&{\bullet}\ar@{-}[rr]&& {\bullet}
}$
& 33 & 48 & 2  \\
\hline
$A_4+A_3$ & $2A_4$ & $\begin{array}{c}
   2\\5\\74\\93
   \end{array}
$ & $\begin{array}{c}
    1\\4\\6\\8
   \end{array}
$ & 0 & 0 & -- & -- & 48 & 48 & 5 \\
\hline
$D_5(a_1)+A_2$ & $D_5(a_1)+A_3$ & $\begin{array}{c}
    4 \\ 8 \\ 31
   \end{array}
$ & $\begin{array}{c}
   3 \\ 5 \\ 39 \\ 61 \\ 75
   \end{array}
$ & 0 & 0 & -- & -- & 46 & 46 & 2  \\
\hline
$(A_5+A_1)'$ & $A_5+A_2+A_1$ & $\begin{array}{c}
    23 \\ 24 \\ 25
   \end{array}
$ & $\begin{array}{c}
    1 \\ 15 \\ 34 \\ 35 \\ 68
   \end{array}
$ & 0 & 0 & -- & -- & 46 & 46 & 3 \\
\hline
\end{tabular*}

\newpage

\begin{tabular*}{1.03\textwidth}{|@{\extracolsep{\fill} }c|c|c|c|c|c|c|c|c|c|c|}
  \hline
 \rule[-.3cm]{0cm}{1cm}{} $\mathcal{O}$ & $\Psi^W(\mathcal{O})$ & $s^1$ & $s^2$ & ${\rm dim}~\h_0$ & $|\Delta_0|$ & $\Delta_0$ & $\Gamma_0^s$ & $l(s)$ & ${\rm dim}~\Sigma_{{\bf k}, s}$ & $d'$ \\ \hline \hline
$D_6(a_2)$ & $2D_4$ & $\begin{array}{c}
    2 \\ 8
   \end{array}
$ & $\begin{array}{c}
    11 \\ 12 \\ 31 \\ 53 \\ 61 \\ 64
   \end{array}
$ & 0 & 0 & -- & --
& 44 & 44 & 1  \\
\hline
$A_5+2A_1$ & $E_6(a_2)+A_2$ & $\begin{array}{c}
   1 \\ 5 \\ 8 \\44
   \end{array}
$ & $\begin{array}{c}
    4 \\ 34 \\ 35 \\ 68
   \end{array}
$ & 0 & 0 & -- & -- & 44 & 44 & 3  \\
\hline
$A_5+A_2$ & $E_7(a_4)+A_1$ & $\begin{array}{c}
   2 \\ 5 \\ 8
   \end{array}
$ & $\begin{array}{c}
    9 \\ 19 \\ 41 \\ 59 \\ 76
   \end{array}
$ & 0 & 0 & -- & -- & 42 & 42 & 1  \\
\hline
$D_5+A_1$ & $D_5+2A_1$ & $\begin{array}{c}
   7 \\ 61
   \end{array}
$ & $\begin{array}{c}
    3 \\ 8 \\ 16 \\ 30 \\ 32
   \end{array}
$ & 1 & 0 & -- & -- & 39 & 40 & 2  \\
\hline
$2A_4$ & $E_8(a_8)$ & $\begin{array}{c}
   1 \\ 2 \\ 6 \\ 8
   \end{array}
$ & $\begin{array}{c}
    19 \\ 41 \\ 59 \\ 76
   \end{array}
$ & 0 & 0 & -- & -- & 40 & 40 & 1  \\
\hline
$D_6(a_1)$ & $D_6(a_1)$ & $\begin{array}{c}
    3\\7\\61
   \end{array}
$ & $\begin{array}{c}
    8\\9\\37
   \end{array}
$ & 2 & 4 & $2A_1$ & $\xymatrix@R=.1cm@C=.1cm{
{\scriptstyle 2}&&{\scriptstyle 5}\\
{\bullet}&&{\bullet}
}$ & 32 & 38 & 2  \\
\hline
$A_6$ & $A_6$ & $\begin{array}{c}
   40 \\ 44 \\ 62
   \end{array}
$ & $\begin{array}{c}
    1 \\ 13 \\ 14
   \end{array}
$ & 2 & 2 & $A_1$ & $\xymatrix@R=.1cm@C=.1cm{
{\scriptstyle 2}\\
{\bullet}
}$ & 34 & 38 & 7  \\
\hline
$A_6+A_1$ & $A_6+A_1$ & $\begin{array}{c}
   40 \\ 44 \\ 62
   \end{array}
$ & $\begin{array}{c}
    1 \\ 2 \\ 13 \\ 14
   \end{array}
$ & 1 & 0 & -- & -- & 35 & 36 & 7  \\
\hline
$D_6(a_1)+A_1$ & $A_7'$ & $\begin{array}{c}
  8 \\ 9 \\ 37
   \end{array}
$ & $\begin{array}{c}
   20 \\ 33 \\ 34 \\ 35
   \end{array}
$ & 1 & 2 & $A_1$ & $\xymatrix@R=.1cm@C=.1cm{
{\scriptstyle 5}\\
{\bullet}
}$ & 33 & 36 & 2 \\
\hline
$(D_5+A_2)_2$ & $D_5+A_2$ & $\begin{array}{c}
    7 \\ 25 \\ 61
   \end{array}
$ & $\begin{array}{c}
    4 \\ 8 \\ 23 \\ 24
   \end{array}
$ & 1 & 0 & -- & -- & 35 & 36 & 6  \\
\hline
$D_5+A_2$ & $A_7+A_1$ & $\begin{array}{c}
   8 \\ 9 \\ 37 
   \end{array}
$ & $\begin{array}{c}
   5 \\ 20 \\ 33 \\ 34 \\ 35 
   \end{array}
$ & 0 & 0 & -- & -- & 34 & 34 & 2 \\
\hline
$E_6(a_1)$ & $E_6(a_1)$ & $\begin{array}{c}
   15\\44\\55
   \end{array}
$ & $\begin{array}{c}
    1 \\ 6\\8
   \end{array}
$ & 2 & 6 & $A_2$ & $\xymatrix@R=.1cm@C=.1cm{
{\scriptstyle 2}&&{\scriptstyle 4}\\
{\bullet}\ar@{-}[rr]&& {\bullet}
}$ & 26 & 34 & 3 \\
\hline
\end{tabular*}

\newpage

\begin{tabular*}{1.03\textwidth}{|@{\extracolsep{\fill} }c|c|c|c|c|c|c|c|c|c|c|}
  \hline
 \rule[-.3cm]{0cm}{1cm}{} $\mathcal{O}$ & $\Psi^W(\mathcal{O})$ & $s^1$ & $s^2$ & ${\rm dim}~\h_0$ & $|\Delta_0|$ & $\Delta_0$ & $\Gamma_0^s$ & $l(s)$ & ${\rm dim}~\Sigma_{{\bf k}, s}$ & $d'$  \\ \hline \hline
$D_6$ & $D_6+2A_1$ & $\begin{array}{c}
   1 \\ 8
   \end{array}
$ & $\begin{array}{c}
    4 \\ 14 \\ 17 \\ 26 \\ 27 \\ 55
   \end{array}
$ & 0 & 0 & -- & -- & 32 & 32 & 1  \\
\hline
$D_7(a_2)$ & $D_7(a_2)$ & $\begin{array}{c}
   4\\7\\61
   \end{array}
$ & $\begin{array}{c}
    3\\8\\16\\30
   \end{array}
$ & 1 & 0 & -- & -- & 31 & 32 & 2 \\
\hline
$E_6$ & $E_6$ & $\begin{array}{c}
   1\\7\\44
   \end{array}
$ & $\begin{array}{c}
    8\\26\\27
   \end{array}
$ & 2 & 6 & $A_2$ & $\xymatrix@R=.1cm@C=.1cm{
{\scriptstyle 4}&&{\scriptstyle 5}\\
{\bullet}\ar@{-}[rr]&& {\bullet}
}$ & 24 & 32 & 3  \\
\hline
$(A_7)_3$ & $A_7''$ & $\begin{array}{c}
    24 \\ 38 \\ 48
   \end{array}
$ & $\begin{array}{c}
    10 \\ 11 \\ 14 \\ 15
   \end{array}
$ & 1 & 0 & -- & -- & 31 & 32 & 4  \\
\hline
$A_7$ & $D_8(a_3)$ & $\begin{array}{c}
   3\\5\\7\\61
   \end{array}
$ & $\begin{array}{c}
    4\\8\\23\\24
   \end{array}
$ & 0 & 0 & -- & -- & 30 & 30 & 1  \\
\hline
$E_6(a_1)+A_1$ & $E_6(a_1)+A_1$ & $\begin{array}{c}
   15\\44\\55
   \end{array}
$ & $\begin{array}{c}
    1 \\ 6\\8\\10
   \end{array}
$ & 1 & 0 & -- & -- & 29 & 30 & 3  \\
\hline
$D_8(a_3)$ & $A_8$ & $\begin{array}{c}
   22\\23\\26\\31
   \end{array}
$ & $\begin{array}{c}
    7\\11\\12\\25
   \end{array}
$ & 0 & 0 & -- & -- & 28 & 28 & 3  \\
\hline
$D_6+A_1$ & $E_7(a_3)$ & $\begin{array}{c}
    1\\2\\8
   \end{array}
$ & $\begin{array}{c}
   27\\28\\32\\41
   \end{array}
$ & 1 & 2 & $A_1$ & $\xymatrix@R=.1cm@C=.1cm{
{\scriptstyle 5}\\
{\bullet}
}$
& 25 & 28 & 1  \\
\hline
$(D_7(a_1))_2$ & $D_7(a_1)$ & $\begin{array}{c}
   1 \\ 8 \\ 12
   \end{array}
$ & $\begin{array}{c}
    19 \\ 21 \\ 33 \\ 49
   \end{array}
$ & 1 & 0 & -- & -- & 27 & 28 & 2  \\
\hline
$D_7(a_1)$ & $D_8(a_2)$ & $\begin{array}{c}
   1\\6\\8
   \end{array}
$ & $\begin{array}{c}
    12\\21\\25\\27\\49
   \end{array}
$ & 0 & 0 & -- & -- & 26 & 26 & 1  \\
\hline
$E_6+A_1$ & $E_6+A_2$ & $\begin{array}{c}
   1\\4\\7\\44
   \end{array}
$ & $\begin{array}{c}
    5\\8\\26\\27
   \end{array}
$ & 0 & 0 & -- & -- & 26 & 26 & 3  \\
\hline
$E_7(a_2)$ & $E_7(a_2)+A_1$ & $\begin{array}{c}
   6 \\ 8 \\ 48
   \end{array}
$ & $\begin{array}{c}
    7 \\ 10 \\ 12 \\ 16 \\ 30
   \end{array}
$ & 0 & 0 & -- & -- & 24 & 24 & 1  \\
\hline
\end{tabular*}

\newpage

\begin{tabular*}{1.03\textwidth}{|@{\extracolsep{\fill} }c|c|c|c|c|c|c|c|c|c|c|}
  \hline
\rule[-.3cm]{0cm}{1cm}{}  $\mathcal{O}$ & $\Psi^W(\mathcal{O})$ & $s^1$ & $s^2$ & ${\rm dim}~\h_0$ & $|\Delta_0|$ & $\Delta_0$ & $\Gamma_0^s$ & $l(s)$ & ${\rm dim}~\Sigma_{{\bf k}, s}$ & $d'$  \\ \hline \hline
$A_8$ & $E_8(a_6)$ & $\begin{array}{c}
   1\\2\\5\\8
   \end{array}
$ & $\begin{array}{c}
    21\\25\\27\\49
   \end{array}
$ & 0 & 0 & -- & -- & 24 & 24 & 1  \\
\hline
$D_7$ & $D_8(a_1)$ & $\begin{array}{c}
    10\\11\\15\\25
   \end{array}
$ & $\begin{array}{c}
    1\\20\\21\\22
   \end{array}
$ & 0 & 0 & -- & -- & 22 & 22 & 1  \\
\hline
$E_7(a_2)+A_1$ & $E_8(a_7)$ & $\begin{array}{c}
   4\\7\\23\\24
   \end{array}
$ & $\begin{array}{c}
    5\\8\\20\\33
   \end{array}
$ & 0 & 0 & -- & -- & 22 & 22 & 1  \\
\hline
$E_7(a_1)$ & $E_7(a_1)$ & $\begin{array}{c}
   1\\13\\14
   \end{array}
$ & $\begin{array}{c}
    3\\5\\8\\32
   \end{array}
$ & 1 & 2 & $A_1$ & $\xymatrix@R=.1cm@C=.1cm{
{\scriptstyle 2}\\
{\bullet}
}$ & 17 & 20 & 1  \\
\hline
$D_8(a_1)$ & $E_8(a_3)$ & $\begin{array}{c}
   3\\7\\23\\24
   \end{array}
$ & $\begin{array}{c}
    4\\8\\26\\27
   \end{array}
$ & 0 & 0 & -- & -- & 20 & 20 & 1 \\
\hline
$E_7(a_1)+A_1$ & $D_8$ & $\begin{array}{c}
   1 \\ 13 \\ 14
   \end{array}
$ & $\begin{array}{c}
    4 \\ 8 \\ 17 \\ 18 \\ 19
   \end{array}
$ & 0 & 0 & -- & -- & 18 & 18 & 1 \\
\hline
$D_8$ & $E_8(a_5)$ & $\begin{array}{c}
   10\\16\\20\\22
   \end{array}
$ & $\begin{array}{c}
    2\\3\\5\\7
   \end{array}
$ & 0 & 0 & -- & -- & 16 & 16 & 1 \\
\hline
$E_7$ & $E_7+A_1$ & $\begin{array}{c}
   1\\6\\8
   \end{array}
$ & $\begin{array}{c}
    2\\3\\5\\7 \\ 32
   \end{array}
$ & 0 & 0 & -- & -- & 16 & 16 & 1 \\
\hline
$E_7+A_1$ & $E_8(a_4)$ & $\begin{array}{c}
   7\\11\\12\\25
   \end{array}
$ & $\begin{array}{c}
    1\\2\\6\\8
   \end{array}
$ & 0 & 0 & -- & -- & 14 & 14 & 1 \\
\hline
$E_8(a_2)$ & $E_8(a_2)$ & $\begin{array}{c}
   2\\3\\5\\7
   \end{array}
$ & $\begin{array}{c}
    1\\8\\10\\20
   \end{array}
$ & 0 & 0 & -- & -- & 12 & 12 & 1 \\
\hline
$E_8(a_1)$ & $E_8(a_1)$ & $\begin{array}{c}
   6\\8\\10\\11
   \end{array}
$ & $\begin{array}{c}
    1\\2\\5\\7
   \end{array}
$ & 0 & 0 & -- & -- & 10 & 10 & 1  \\
\hline
$E_8$ & $E_8$ & $\begin{array}{c}
   1\\4\\6\\8
   \end{array}
$ & $\begin{array}{c}
    2\\3\\5\\7
   \end{array}
$ & 0 & 0 & -- & -- & 8 & 8 & 1 \\
\hline
\end{tabular*}


\section*{Appendix 3. Irreducible root systems of exceptional types.}

\addcontentsline{toc}{section}{Appendix 3. Irreducible root systems of exceptional types}

\renewcommand{\thetheorem}{A3.\arabic{theorem}}

\renewcommand{\thelemma}{A3.\arabic{lemma}}

\renewcommand{\theproposition}{A3.\arabic{proposition}}

\renewcommand{\thecorollary}{A3.\arabic{corollary}}

\renewcommand{\theremark}{A3.\arabic{remark}}

\renewcommand{\thedefinition}{A3.\arabic{definition}}

\renewcommand{\theequation}{A3.\arabic{equation}}

\setcounter{equation}{0}
\setcounter{theorem}{0}

\pagestyle{myheadings}
\markboth{APPENDIX}{APPENDIX 3. EXCEPTIONAL ROOT SYSTEMS}

In this Appendix we give the lists of positive roots in irreducible root systems of exceptional types. All simple roots are numbered as shown at the Dynkin diagrams. The other roots in each list are given in terms of their coordinates with respect to the basis of simple roots. The coordinates are indicated in the brackets ( ). Each set of coordinates is preceded by the number of the corresponding root. These numbers are used to indicate roots which appear in the columns $s^1$, $s^2$ and $\Gamma_0^s$ in the tables in Appendix 2.

\vskip 0.3cm

$\bf G_2.$

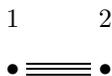
$$\xymatrix@R=.25cm{
1&2\\
{\bullet}\ar@3{-}[r] &{\bullet}
}$$

\begin{enumerate}[1~~]
\item
    (1 0)
\item
    (0 1)
\item
    (1 1)
\item
    (2 1)
\item
    (3 1)
\item
    (3 2)
\end{enumerate}

\vskip 1cm

$\bf F4.$

$$\xymatrix@R=.25cm{
1&2&3&4\\
{\bullet}\ar@{-}[r]&{\bullet}\ar@2{-}[r]
&{\bullet}\ar@{-}[r]&{\bullet} }$$

\begin{multicols}{3}
\begin{enumerate}[1~~]
\item
    (1 0 0 0)
\item
    (0 1 0 0)
\item
    (0 0 1 0)
\item
    (0 0 0 1)
\item
    (1 1 0 0)
\item
    (0 1 1 0)
\item
    (0 0 1 1)
\item
    (1 1 1 0)
\item
    (0 1 2 0)
\item
    (0 1 1 1)
\item
    (1 1 2 0)
\item
    (1 1 1 1)
\item
    (0 1 2 1)
\item
    (1 2 2 0)
\item
    (1 1 2 1)
\item
    (0 1 2 2)
\item
    (1 2 2 1)
\item
    (1 1 2 2)
\item
    (1 2 3 1)
\item
    (1 2 2 2)
\item
    (1 2 3 2)
\item
    (1 2 4 2)
\item
    (1 3 4 2)
\item
    (2 3 4 2)
\end{enumerate}
\end{multicols}

\vskip 1cm

$\bf E_6.$

$$\xymatrix@R=.25cm@C=.25cm{
1&&3&&4&&5&&6\\
{\bullet}\ar@{-}[rr]&&{\bullet}\ar@{-}[rr]&&{\bullet}
\ar@{-}[dd]\ar@{-}[rr]&& {\bullet}\ar@{-}[rr] &&{\bullet}
\\ &&&&&&&&\\
&&&&{\bullet}&&&& \\
&&&&2&&&&}$$

\begin{multicols}{3}
\begin{enumerate}[1~~]
\item
   (1 0 0 0 0 0)
\item
   (0 1 0 0 0 0)
\item
   (0 0 1 0 0 0)
\item
   (0 0 0 1 0 0)
\item
   (0 0 0 0 1 0)
\item
   (0 0 0 0 0 1)
\item
   (1 0 1 0 0 0)
\item
   (0 1 0 1 0 0)
\item
   (0 0 1 1 0 0)
\item
   (0 0 0 1 1 0)
\item
   (0 0 0 0 1 1)
\item
   (1 0 1 1 0 0)
\item
   (0 1 1 1 0 0)
\item
   (0 1 0 1 1 0)
\item
   (0 0 1 1 1 0)
\item
   (0 0 0 1 1 1)
\item
   (1 1 1 1 0 0)
\item
   (1 0 1 1 1 0)
\item
   (0 1 1 1 1 0)
\item
   (0 1 0 1 1 1)
\item
   (0 0 1 1 1 1)
\item
   (1 1 1 1 1 0)
\item
   (1 0 1 1 1 1)
\item
   (0 1 1 2 1 0)
\item
   (0 1 1 1 1 1)
\item
   (1 1 1 2 1 0)
\item
   (1 1 1 1 1 1)
\item
   (0 1 1 2 1 1)
\item
   (1 1 2 2 1 0)
\item
   (1 1 1 2 1 1)
\item
   (0 1 1 2 2 1)
\item
   (1 1 2 2 1 1)
\item
   (1 1 1 2 2 1)
\item
   (1 1 2 2 2 1)
\item
   (1 1 2 3 2 1)
\item
   (1 2 2 3 2 1)
\end{enumerate}
\end{multicols}

\vskip 1cm

$\bf E_7.$

$$\xymatrix@R=.25cm@C=.25cm{
1&&3&&4&&5&&6&&7\\
{\bullet}\ar@{-}[rr]&&{\bullet}\ar@{-}[rr]&&{\bullet}
\ar@{-}[dd]\ar@{-}[rr]&& {\bullet}\ar@{-}[rr] &&{\bullet}\ar@{-}[rr]&&{\bullet}
\\ &&&&&&&&&&\\
&&&&{\bullet}&&&&&& \\
&&&&2&&&&&&}$$

\begin{multicols}{3}
\begin{enumerate}[1~~]
\item
    (1 0 0 0 0 0 0)
\item
    (0 1 0 0 0 0 0)
\item
    (0 0 1 0 0 0 0)
\item
    (0 0 0 1 0 0 0)
\item
    (0 0 0 0 1 0 0)
\item
    (0 0 0 0 0 1 0)
\item
    (0 0 0 0 0 0 1)
\item
    (1 0 1 0 0 0 0)
\item
    (0 1 0 1 0 0 0)
\item
    (0 0 1 1 0 0 0)
\item
    (0 0 0 1 1 0 0)
\item
    (0 0 0 0 1 1 0)
\item
    (0 0 0 0 0 1 1)
\item
    (1 0 1 1 0 0 0)
\item
    (0 1 1 1 0 0 0)
\item
    (0 1 0 1 1 0 0)
\item
    (0 0 1 1 1 0 0)
\item
    (0 0 0 1 1 1 0)
\item
    (0 0 0 0 1 1 1)
\item
    (1 1 1 1 0 0 0)
\item
    (1 0 1 1 1 0 0)
\item
    (0 1 1 1 1 0 0)
\item
    (0 1 0 1 1 1 0)
\item
    (0 0 1 1 1 1 0)
\item
    (0 0 0 1 1 1 1)
\item
    (1 1 1 1 1 0 0)
\item
    (1 0 1 1 1 1 0)
\item
    (0 1 1 2 1 0 0)
\item
    (0 1 1 1 1 1 0)
\item
    (0 1 0 1 1 1 1)
\item
    (0 0 1 1 1 1 1)
\item
    (1 1 1 2 1 0 0)
\item
    (1 1 1 1 1 1 0)
\item
    (1 0 1 1 1 1 1)
\item
    (0 1 1 2 1 1 0)
\item
    (0 1 1 1 1 1 1)
\item
    (1 1 2 2 1 0 0)
\item
    (1 1 1 2 1 1 0)
\item
    (1 1 1 1 1 1 1)
\item
    (0 1 1 2 2 1 0)
\item
    (0 1 1 2 1 1 1)
\item
    (1 1 2 2 1 1 0)
\item
    (1 1 1 2 2 1 0)
\item
    (1 1 1 2 1 1 1)
\item
    (0 1 1 2 2 1 1)
\item
    (1 1 2 2 2 1 0)
\item
    (1 1 2 2 1 1 1)
\item
    (1 1 1 2 2 1 1)
\item
    (0 1 1 2 2 2 1)
\item
    (1 1 2 3 2 1 0)
\item
    (1 1 2 2 2 1 1)
\item
    (1 1 1 2 2 2 1)
\item
    (1 2 2 3 2 1 0)
\item
    (1 1 2 3 2 1 1)
\item
    (1 1 2 2 2 2 1)
\item
    (1 2 2 3 2 1 1)
\item
    (1 1 2 3 2 2 1)
\item
    (1 2 2 3 2 2 1)
\item
    (1 1 2 3 3 2 1)
\item
    (1 2 2 3 3 2 1)
\item
    (1 2 2 4 3 2 1)
\item
    (1 2 3 4 3 2 1)
\item
    (2 2 3 4 3 2 1)
\end{enumerate}
\end{multicols}

\vskip 1cm

$\bf E_8.$

$$\xymatrix@R=.25cm@C=.25cm{
1&&3&&4&&5&&6&&7&&8\\
{\bullet}\ar@{-}[rr]&&{\bullet}\ar@{-}[rr]&&{\bullet}
\ar@{-}[dd]\ar@{-}[rr]&& {\bullet}\ar@{-}[rr] &&{\bullet}\ar@{-}[rr]&&{\bullet}\ar@{-}[rr]&&{\bullet}
\\ &&&&&&&&&&&&\\
&&&&{\bullet}&&&&&&&& \\
&&&&2&&&&&&&&}$$

\begin{multicols}{3}
\begin{enumerate}[1~~]
\item
    (1 0 0 0 0 0 0 0)
\item
    (0 1 0 0 0 0 0 0)
\item
    (0 0 1 0 0 0 0 0)
\item
    (0 0 0 1 0 0 0 0)
\item
    (0 0 0 0 1 0 0 0)
\item
    (0 0 0 0 0 1 0 0)
\item
    (0 0 0 0 0 0 1 0)
\item
    (0 0 0 0 0 0 0 1)
\item
    (1 0 1 0 0 0 0 0)
\item
    (0 1 0 1 0 0 0 0)
\item
    (0 0 1 1 0 0 0 0)
\item
    (0 0 0 1 1 0 0 0)
\item
    (0 0 0 0 1 1 0 0)
\item
    (0 0 0 0 0 1 1 0)
\item
    (0 0 0 0 0 0 1 1)
\item
    (1 0 1 1 0 0 0 0)
\item
    (0 1 1 1 0 0 0 0)
\item
    (0 1 0 1 1 0 0 0)
\item
    (0 0 1 1 1 0 0 0)
\item
    (0 0 0 1 1 1 0 0)
\item
    (0 0 0 0 1 1 1 0)
\item
    (0 0 0 0 0 1 1 1)
\item
    (1 1 1 1 0 0 0 0)
\item
    (1 0 1 1 1 0 0 0)
\item
    (0 1 1 1 1 0 0 0)
\item
    (0 1 0 1 1 1 0 0)
\item
    (0 0 1 1 1 1 0 0)
\item
    (0 0 0 1 1 1 1 0)
\item
    (0 0 0 0 1 1 1 1)
\item
    (1 1 1 1 1 0 0 0)
\item
    (1 0 1 1 1 1 0 0)
\item
    (0 1 1 2 1 0 0 0)
\item
    (0 1 1 1 1 1 0 0)
\item
    (0 1 0 1 1 1 1 0)
\item
    (0 0 1 1 1 1 1 0)
\item
    (0 0 0 1 1 1 1 1)
\item
    (1 1 1 2 1 0 0 0)
\item
    (1 1 1 1 1 1 0 0)
\item
    (1 0 1 1 1 1 1 0)
\item
    (0 1 1 2 1 1 0 0)
\item
    (0 1 1 1 1 1 1 0)
\item
    (0 1 0 1 1 1 1 1)
\item
    (0 0 1 1 1 1 1 1)
\item
    (1 1 2 2 1 0 0 0)
\item
    (1 1 1 2 1 1 0 0)
\item
    (1 1 1 1 1 1 1 0)
\item
    (1 0 1 1 1 1 1 1)
\item
    (0 1 1 2 2 1 0 0)
\item
    (0 1 1 2 1 1 1 0)
\item
    (0 1 1 1 1 1 1 1)
\item
    (1 1 2 2 1 1 0 0)
\item
    (1 1 1 2 2 1 0 0)
\item
    (1 1 1 2 1 1 1 0)
\item
    (1 1 1 1 1 1 1 1)
\item
    (0 1 1 2 2 1 1 0)
\item
    (0 1 1 2 1 1 1 1)
\item
    (1 1 2 2 2 1 0 0)
\item
    (1 1 2 2 1 1 1 0)
\item
    (1 1 1 2 2 1 1 0)
\item
    (1 1 1 2 1 1 1 1)
\item
    (0 1 1 2 2 2 1 0)
\item
    (0 1 1 2 2 1 1 1)
\item
    (1 1 2 3 2 1 0 0)
\item
    (1 1 2 2 2 1 1 0)
\item
    (1 1 2 2 1 1 1 1)
\item
    (1 1 1 2 2 2 1 0)
\item
    (1 1 1 2 2 1 1 1)
\item
    (0 1 1 2 2 2 1 1)
\item
    (1 2 2 3 2 1 0 0)
\item
    (1 1 2 3 2 1 1 0)
\item
    (1 1 2 2 2 2 1 0)
\item
    (1 1 2 2 2 1 1 1)
\item
    (1 1 1 2 2 2 1 1)
\item
    (0 1 1 2 2 2 2 1)
\item
    (1 2 2 3 2 1 1 0)
\item
    (1 1 2 3 2 2 1 0)
\item
    (1 1 2 3 2 1 1 1)
\item
    (1 1 2 2 2 2 1 1)
\item
    (1 1 1 2 2 2 2 1)
\item
    (1 2 2 3 2 2 1 0)
\item
    (1 2 2 3 2 1 1 1)
\item
    (1 1 2 3 3 2 1 0)
\item
    (1 1 2 3 2 2 1 1)
\item
    (1 1 2 2 2 2 2 1)
\item
    (1 2 2 3 3 2 1 0)
\item
    (1 2 2 3 2 2 1 1)
\item
    (1 1 2 3 3 2 1 1)
\item
    (1 1 2 3 2 2 2 1)
\item
    (1 2 2 4 3 2 1 0)
\item
    (1 2 2 3 3 2 1 1)
\item
    (1 2 2 3 2 2 2 1)
\item
    (1 1 2 3 3 2 2 1)
\item
    (1 2 3 4 3 2 1 0)
\item
    (1 2 2 4 3 2 1 1)
\item
    (1 2 2 3 3 2 2 1)
\item
    (1 1 2 3 3 3 2 1)
\item
    (2 2 3 4 3 2 1 0)
\item
    (1 2 3 4 3 2 1 1)
\item
    (1 2 2 4 3 2 2 1)
\item
    (1 2 2 3 3 3 2 1)
\item
    (2 2 3 4 3 2 1 1)
\item
    (1 2 3 4 3 2 2 1)
\item
    (1 2 2 4 3 3 2 1)
\item
    (2 2 3 4 3 2 2 1)
\item
    (1 2 3 4 3 3 2 1)
\item
    (1 2 2 4 4 3 2 1)
\item
    (2 2 3 4 3 3 2 1)
\item
    (1 2 3 4 4 3 2 1)
\item
    (2 2 3 4 4 3 2 1)
\item
    (1 2 3 5 4 3 2 1)
\item
    (2 2 3 5 4 3 2 1)
\item
    (1 3 3 5 4 3 2 1)
\item
    (2 3 3 5 4 3 2 1)
\item
    (2 2 4 5 4 3 2 1)
\item
    (2 3 4 5 4 3 2 1)
\item
    (2 3 4 6 4 3 2 1)
\item
    (2 3 4 6 5 3 2 1)
\item
    (2 3 4 6 5 4 2 1)
\item
    (2 3 4 6 5 4 3 1)
\item
    (2 3 4 6 5 4 3 2)
\end{enumerate}
\end{multicols}


\begin{thebibliography}{99}
\addcontentsline{toc}{chapter}{Bibliography}

\pagestyle{myheadings}
\markboth{BIBLIOGRAPHY}{BIBLIOGRAPHY}

\bibitem{APK} Andersen, H. H., Polo, P., Kexin, W., Representations of quantum algebras, {\em Invent. Math.} {\bf 104} (1991), 1--60.

\bibitem{AJS} Andersen, H. H., Jantzen, J. C., Soergel, W., Representations of quantum groups at a $p$--th root of unity and of semisimple groups in
characteristic $p$: independence of $p$, {\em Ast\'{e}risque} {\bf 220} (1994).

\bibitem{AST1} Asherova, R. M., Smirnov, Yu. F., Tolstoy, V. N., Projection operators for the
simple Lie groups, {\em Theor. Math. Phys.} {\bf 8}(1971), 813--825.

\bibitem{AST2} Asherova, R. M., Smirnov, Yu. F., Tolstoy, V. N., Projection operators for the
simple Lie groups. II. General scheme for construction of lowering operators.
The case of the group SU(n), {\em Theor. Math. Phys.} {\bf 15} (1973), 392--393.

\bibitem{AST3} Asherova, R. M., Smirnov, Yu. F., Tolstoy, V. N., A description of some class of
projection operators for semisimple complex Lie algebras, {\em Matem. Zametki} {\bf 26}
(1979), 499--504.

\bibitem{Ba} Baader, F., Nipkow, T., Term Rewriting and All That, Cambridge University Press, Cambridge (1999).

\bibitem{BMR} Bezrukavnikov, R., Mirkovic, I., Rumynin, D., Localization of modules for a semisimple Lie algebra in prime characteristic. With an appendix by Bezrukavnikov and Simon Riche, {\em  Ann. of Math.} {\bf 167}  (2008), 945--991.

\bibitem{BD}  Belavin, A. A., Drinfeld, V. G., Solutions of the classical Yang-Baxter equation for simple Lie algebras, {\em Funct. Anal. Appl.} {\bf 16} (1981), 159--180.

\bibitem{Gel} Bernstein, J. H., Gelfand, I. M., Gelfand, S. I., Schubert cells and cohomology of the spaces G/P, {\em Uspekhi Mat. Nauk} {\bf 28}, no. 3 (171) (1973), 3--26; {\em Russian Math. Surveys} {\bf 28}, no. 3 (1973), 1--26.

\bibitem{Bur} Bourbaki, N., Groupes et algebras de Lie, Chap. 4,5,6, Paris, Hermann (1968).

\bibitem{Br} Brieskorn, E., Singular elements of semisimple algebraic groups, {\em Actes Congr\`{e}s Intern. Math.} {\bf 2} (1970), 279--284.

\bibitem{BG} Brown, K. A., Goodearl, K. R., Lectures on Algebraic Quantum Groups, Birkh\"{a}user (2002).

\bibitem{CCC} Cantarini, N., Carnovale, G., Costantini, M., Spherical orbits and representations of ${U}_\varepsilon({\g})$, {\em  Transf. Groups} {\bf  10}  (2005), 29–-62.

\bibitem{10} Cantarini, N., The quantized enveloping algebra $U_q(sl(n))$ at the roots of unity,
{\em Comm. Math. Phys.} {\bf 211} (2000), 207–-230.

\bibitem{9} Cantarini, N., Spherical orbits and quantized enveloping algebras, {\em Comm. Algebra}
{\bf 27} (1999), 3439–3458.

\bibitem{8} Cantarini, N., Mod-p reduction for quantum groups, {\em J. of Alg.} {\bf 202} (1998),
357–-366.

\bibitem{Car}  Carter, R. W., Simple groups of Lie type, John Wiley \& Sons, Inc., New York (1989).

\bibitem{C} Carter, R. W., Conjugacy classes in the Weyl group, {\em Compositio Math.} {\bf 25} (1972), 1--59.

\bibitem{Car1}  Carter, R. W., Finite groups of Lie type. Conjugacy classes and complex characters, John Wiley \& Sons, Ltd., Chichester (1993).

\bibitem{ChP} Chari, V., Pressley, A., A guide to quantum groups, Cambridge Univ. Press (1994).

\bibitem{Cox} Coxeter, H. S. M., Regular Polytopes, Methuen \& Co., London (1948).

\bibitem{DB} De Boer, J., Tjin, T., Quantization and representation theory of finite W--algebras, {\em Comm. math. Phys.}, {\bf 158} (1993), 485--516.

\bibitem{DK} De Concini, C., Kac, V. G., Representations of quantum groups at roots of $1$.  Operator algebras, unitary representations, enveloping algebras, and invariant theory (Paris, 1989), {\em Progr. Math.} {\bf 92}, Birkhäuser Boston, Boston, MA (1990), 471--506.

\bibitem{DK1} De Concini, C., Kac, V. G., Representations of quantum groups at roots of 1: reduction to the exceptional case, {\em Int. J. Mod. Phys. A} {\bf 7} {\em Suppl. 1A} (1992), 141--149.

\bibitem{DKP1} De Concini, C., Kac, V. G., Procesi, C., Quantum coadjoint action,  {\em J. Amer. Math. Soc.}  {\bf 5}  (1992), 151--189.

\bibitem{DKP2} De Concini, C., Kac, V. G., Procesi, C., Some remarkable degenerations of quantum groups,  {\em  Comm. Math. Phys.} {\bf  157}  (1993), 405--427.

\bibitem{DKP3} De Concini, C., Kac, V. G., Procesi, C., Some quantum analogues of solvable Lie groups. In: Geometry and analysis, Proceedings of the International Colloquium on Geometry and Analysis, Bombay, 1992, Oxford University Press (1995), 41--66.

\bibitem{DL} De Concini, C., Lyubashenko, V., Quantum function algebra at roots of 1, {\em Adv. in Math.} {\bf 108} (1994), 205--262.

\bibitem{DKPschubert} De Concini, C., Procesi, C., Quantum Schubert cells and representations at roots of 1. In: Algebraic groups and Lie groups (G. I. Lehrer ed.), {\em Austral. Math. Soc. Lect. Ser.} {\bf 9}, Cambridge Univ. Press, Cambridge (1997), 127--160.

\bibitem{DP} De Concini, C., Procesi, C., Quantum Groups. In: D-modules, Representation Theory, and Quantum Groups (L. Boutet de Monvel et al.), Proc. Venezia 1992, {\em Lecture
Notes in Mathematics} {\bf 1565}, Springer, Berlin (1993), 31--140. 

\bibitem{CR} Curtis, C. W., Reiner, I., Representation theory of finite groups and associative algebras, Interscience Publishers (1962).

\bibitem{Dm} Drinfeld, V.G., Quantum groups, Proc. Int. Congr. Math. Berkley, California, 1986, Amer. Math. Soc., Providence (1987), 718--820.

\bibitem{Dri} Drinfeld, V.G., Quasi--Hopf algebras, {\em Leningrad Math. Journal} {\bf 1} (1990), 1419--1457.

\bibitem{Dr} Drupieski, C. M., On injective modules and support varieties for the small quantum group, {\em  Int. Math. Res. Not.} {\bf  2011} (2011), 2263--2294.

\bibitem{EG} Ellers, E. W., Gordeev, N., Intersection of conjugacy classes with Bruhat cells in Chevalley groups, {\em Pacific J. Math.} {\bf 214} (2004), 245--260.

\bibitem{FRT}  Faddeev, L. D., Reshetikhin, N. Yu., Takhtajan, L. A., Quantization of Lie
groups and Lie algebras, {\em Leningrad Math. J.} {\bf 1} (1989), 178--206.


\bibitem{FZ} Fomin, S., Zelevinsky, A., Double Bruhat cells and total positivity, {\em J. Amer. Math. Soc.} {\bf 12} (1999), 335--380.

\bibitem{FZ1} Fomin, S., Zelevinsky, A., Recognizing Schubert cells, {\em Journal of Algebraic Combinatorics} {\bf 12} (2000), 37--57.


\bibitem{FP1} Friedlander, E. M., Parshall, B. J., Support varieties for restricted Lie algebras {\em  Invent. Math.} {\bf  86}  (1986),  553--562.

\bibitem{FP2} Friedlander, E. M., Parshall, B. J., Geometry of $p$-unipotent Lie algebras, {\em  J. of Alg.}  {\bf 109}  (1987), 25--45.

\bibitem{FP} Friedlander, E. M., Parshall, B. J., Modular representation theory of Lie algebras, {\em Amer. J. Math.} {\bf 110} (1988), 1055--1093.

\bibitem{GG} Gan, W. L., Ginzburg, V., Quantization of Slodowy
slices, {\em Int. Math. Res. Not.} {\bf 5} (2002), 243--255.

\bibitem{GR} Gasper, G., Rahman, M., Basic hypergeometric series, Cambridge Univ. Press (1990).

\bibitem{Gck} Geck, M., Pfeiffer, G., Characters of finite Coxeter groups and Iwahori-Hecke algebras, {\em London Mathematical Society Monographs. New Series} {\bf 21}, The Clarendon Press, Oxford University Press, New York (2000).

\bibitem{Gk1}  Geck, M., Malle, G., On the existence of a unipotent support for the irreducible characters of a finite group of Lie type, {\em Trans. Amer. Math. Soc.} {\bf 352} (2000), 429–-456.

\bibitem{GS} Gelfand, I. M., Serganova, V. V., Combinatorial geometries and the strata of a torus on homogeneous compact manifolds, {\em Russian Math. Surveys} {\bf 42} (1987), 133--168.

\bibitem{GK} Ginzburg, V., Kumar, S., Cohomology of quantum groups at roots of unity, {\em  Duke Math. J. } {\bf  69}  (1993),  179--198.

\bibitem{GY} Goodearl, K. R., Yakimov, M. T., The Berenstein-Zelevinsky quantum cluster algebra conjecture, {\em J. Eur. Math. Soc.}, {\bf 22} (2020), 2453--2509.

\bibitem{GG1} Goto, M., and Grosshans, F. D., Semisimple Lie algebras, Marcel Dekker, Inc., New York and Basel (1978).


\bibitem{Har} Hartshorne, R., Algebraic geometry, Springer, New York (2006).

\bibitem{Ha} Harzheim, E., Ordered Sets, Springer, New York (2006).

\bibitem{XL} He, X., Lusztig, G., A generalization of Steinberg's cross-section,  {\em J. Amer. Math. Soc.}  {\bf 25}  (2012), 739--757.

\bibitem{XN} He, X., Nie, S., Minimal length elements of finite Coxeter groups, {\em 	Duke Math. J.} {\bf 161} (2012), 2945--2967.

\bibitem{Hess} Hesselink, W. H., Nilpotency in classical groups over a field of characteristic 2, {\em Math. Z.} {\bf 166} (1979), 165–-181.

\bibitem{Hu} Humphreys, J., E., Linear algebraic groups, Springer, New York (1981). 

\bibitem{I} Iwahori, N., On the structure of a Hecke ring of a Chevalley group over a finite field, {\em J. Fac, Sci. Uni. Tokyo}, Sect. {\bf IA 10} (1964), 215--236.

\bibitem{Jac} Jacobson, N., Lie algebras, Interscience Publishers, New York (1962).

\bibitem{Jan} Jantzen, J. C., Lectures on Quantum Groups, {\em Graduate Studies in Mathematics} {\bf 6}, AMS (1996).

\bibitem{Jos} Joseph, A., Quantum groups and their primitive ideals. Ergebnisse der Mathematik und ihrer Grenzgebiete (3), {\bf 29}, Springer-Verlag, Berlin (1995).

\bibitem{JL} Joseph, A., Letzter, G., Local finiteness of the adjoint action for quantized enveloping algebras, {\em  J. of Alg.}  {\bf 153}  (1992), 289--318.

\bibitem{KW}  Kac, V. G., Weisfeiler, B., The irreducible representations of Lie $p$-algebras, {\em Funk. Anal. i ego Pril.}  {\bf 5}  (1971), no. 2, 28–-36.

\bibitem{KQ} Kac, V., Cheung, P., Quantum Calculus, Springer, New York (2002).

\bibitem{Ka1} Kawanaka, N., Generalized Gelfand–Graev representations and Ennola duality, in: Algebraic Groups and Related Topics, {\em Advanced Studies in Pure Mathematics} {\bf 6}, North-Holland, Amsterdam/New York/Oxford (1985), 175--206.

\bibitem{Ka2} Kawanaka, N., Generalized Gelfand–Graev representations of exceptional simple groups over a finite field I, {\em Invent. Math.} {\bf 84} (1986), 575--616.

\bibitem{KL}  Kazhdan, D., Lusztig, G., Fixed point varieties on affine flag manifolds, {\em Israel J. Math.} {\bf 62} (1988), 129--168.

\bibitem{KO} Khoroshkin, S. M., Ogievetsky O., Mickelsson algebras and Zhelobenko operators, {\em J. of Alg.} {\bf 319} (2008),  2113--2165.

\bibitem{KT1} Khoroshkin, S. M., Tolstoy, V. N., Twisting of quantized Lie (super)algebras. Quantum groups (Karpacz, 1994), 63--84, Wydawnictwo Naukowe PWN, Warsaw (1995).

\bibitem{KT3} Khoroshkin, S. M., Tolstoy, V. N., The Cartan-Weyl basis and the universal R-matrix for quantum Kac-Moody algebras and superalgebras.  Quantum symmetries (Clausthal, 1991),  336--351, World Sci. Publ., River Edge, NJ (1993).


\bibitem{kh-t}  Khoroshkin, S. M., Tolstoy, V. N., Universal R--matrix for quantized
(super)algebras, {\em Comm. Math. Phys.} {\bf 141} (1991), 599--617.

\bibitem{KlSm} Klimyk, A., Schm\"{u}dgen, K., Quantum groups and their representations, Springer-Verlag, Heidelberg (1997).

\bibitem{KolSt} Kolb, S., Stokman, J. V., Reflection equation algebras, coideal subalgebras, and their centres, {\em Selecta Math. (New Ser.)} {\bf 15} (2009), 621--664.

\bibitem{K} Kostant, B., On Whittaker vectors and representation theory, {\em Invent. Math.} {\bf 48} (1978), 101--184.

\bibitem{KR} Kremnizer, K., Proof of the De Concini-Kac-Procesi conjecture, arXiv:math/0611236.

\bibitem{Kum} Kumar, S., Representations of quantum groups at roots of unity,
{\em In: Quantum topology} (ed. by D.N. Yetter), World Scientific, Singapore (1994) 187--224.


\bibitem{Li}  Liebeck, M. W., Seitz, G. M., Unipotent and nilpotent classes in simple algebraic groups and Lie algebras, {\em Mathematical Surveys and Monographs} {\bf 180}, American Mathematical Society, Providence, RI (2012).

\bibitem{Los1} Losev, I., Finite dimensional representations of W-algebras, {\em Duke Math. J.} {\bf 159} (2011), 99--143.

\bibitem{Los2} Losev, I., Ostrik, V., Classification of finite dimensional irreducible modules over W-algebras, {\em Compositio Math.} {\bf 150} (2014), 1024--1076.

\bibitem{L} L\"{o}wdin, P.-O., Angular momentum wave functions constructed by projector operators, {\em Rev. Mod. Phys.} {\bf 36} (1964), 966--976.

\bibitem{Lu}  Lu J. H., Momentum mapping and reduction of Poisson actions. 
{\em In}: Symplectic geometry, groupoids and integrable systems, Berkeley,
1989. P.Dazord and A.Weinstein (eds), pp.209-226. Springer-Verlag.

\bibitem{LuW} Lu, J. H., Weinstein, A., Poisson Lie groups, dressing transformations, and
Bruhat decompositions, {\em J. Diff. Geom} {\bf 31} (1990), 501--526.

\bibitem{Lus} Lusztig, G., Introduction to quantum groups, Birkh\"{a}user (1994).


\bibitem{L1}  Lusztig, G., Quantum groups at roots of 1, {\em Geom. Dedicata} {\bf 35} (1990), no. 1--3, 89--113.

\bibitem{L2} Lusztig, G., On conjugacy classes in a reductive group. In: Representations of Reductive Groups (M. Nevins, P. E. Trapa eds.), {\em  Progress in Mathematics} {\bf 312}, Birkh\"{a}user (2015), 333--363. 

\bibitem{L2'} Lusztig, G., Spaltenstein, N., On the generalized Springer correspondence for classical groups. Algebraic groups and related topics (Kyoto/Nagoya, 1983), 289--316, {\em Adv. Stud. Pure Math.} {\bf 6}, North-Holland, Amsterdam, (1985).

\bibitem{L3'} Lusztig, G., From conjugacy classes in the Weyl group to unipotent classes II. {\em Represent. Theory} {\bf 16} (2012), 189--211.

\bibitem{L4'} Lusztig, G., From conjugacy classes in the Weyl group to unipotent classes, {\em Represent. Theory} {\bf 15} (2011), 494–-530.

\bibitem{L5'} Lusztig, G., Unipotent elements in small characteristic, {\em Transform. Groups} {\bf 10} (2005), 449--487.

\bibitem{L3} Lusztig, G., Quantum Deformations of Certain Simple Modules over Enveloping Algebras, {\em Adv. Math.} {\bf 70} (1988), 237--249.

\bibitem{Ly} Lynch, T. E.,  Generalized Whittaker vectors and representation theory, Thesis, M.I.T. (1979).

\bibitem{Malt1} Malten, W., Minimally dominant elements of finite Coxeter groups, arXiv:2110.09266.

\bibitem{Malt2} Malten, W., From braids to transverse slices in reductive groups, arXiv:2111.01313.

\bibitem{Mc} McConnell, J.C., Quantum groups, filtered rings and Gelfand--KIrillov dimension, In {\em Lecture Notes in Math} {\bf 1448}, Springer, Berlin (1991), 139--149.

\bibitem{O} Ostrik, V. V., Cohomological supports for quantum groups, {\em Funct. Anal. Appl.} {\bf  32}  (1999),  237--246.

\bibitem{Pas} Passman, D. S., The algebraic structure of group rings, Jonh Wiley \& Sons, New York (1977).

\bibitem{Pie} Pierce, R.S., Associative algebras. {\em Graduate Texts in Mathematics} {\bf 88}, Springer-Verlag, New York-Berlin (1982).

\bibitem{PS} Premet, A., Skryabin, S., Representations of restricted Lie algebras and families of associative $\mathcal{L}$--algebras, {\em J. Reine Angew. Math.} {\bf 507} (1999), 189--218.


\bibitem{Pr} Premet, A., Special transverse slices and their enveloping
algebras. With an appendix by Serge Skryabin,  {\em Adv. Math.}
{\bf 170} (2002), 1--55.

\bibitem{Pr1} Premet, A., Irreducible representations of Lie algebras of reductive groups and the Kac-Weisfeiler conjecture, {\em  Invent. Math.} {\bf  121}  (1995), 79--117.

\bibitem{Pr3} Premet, A., Support varieties of non-restricted modules over Lie algebras of reductive groups, {\em J. London Math. Soc.} (2) {\bf 55} (1997), 236--250.

\bibitem{fact}  Reshetikhin, N.Yu., Semenov-Tian-Shansky, M.A., Factorization problems for quantum groups, in: {\em  Geometry and Physics, essays in honour of I. M. Gelfand (S. Gindikin and I. M. Singer eds.)}, North Holland, Amsterdam - London - New York (1991), 533--550.

\bibitem{Ric} Richardson, R.W., Conjugacy classes of involutions in Coxeter groups, {\em Bull. Austral. Math. Soc.} {\bf 26} (1982), 1--15.

\bibitem{R} Rotman, S., Advanced Linear Algebra, Springer, New York (2005).

\bibitem{Saito} Saito, Y., PBW Basis of Quantized Universal Enveloping Algebras, {\em Publ. Math. RIMS}
{\bf 30} (1994), 209--232. 


\bibitem{rmatr}  Semenov-Tian-Shansky, M.A., What is a classical $r$-matrix, {\em Funct. Anal. Appl.} {\bf 17} (1983), 17--33.

\bibitem{RIMS}  Semenov-Tian-Shansky, M.A., Dressing transformations and Poisson--Lie group actions, {\em Publ. Math. RIMS} {\bf 21} (1985), 1237--1260.

\bibitem{dual}  Semenov-Tian-Shansky, M. A., Poisson Lie groups, quantum
duality principle and the quantum double, {\em Contemp. Math.} {\bf 175} (1994)
, 219--248.

\bibitem{SemSev} Semenov-Tian-Shansky, M. A., Sevostyanov, A.,
Drinfeld-Sokolov reduction for difference operators and
deformations of W-algebras. II. General Semisimple Case, {\em Comm.
Math. Phys.} {\bf 192} (1998), 631--647.


\bibitem{S0} Sevostyanov, A., Regular nilpotent elements and quantum groups, {\em Comm. Math. Phys.} {\bf 204} (1999) 1--16.

\bibitem{S1} Sevostyanov, A., Reduction of quantum systems with arbitrary first-class constraints and Hecke algebras, {\em Comm. Math. Phys.} {\bf 204} (1999), 137--146.

\bibitem{SThes} Sevostyanov, A., The Whittaker model of the center
of the quantum group and Hecke algebras, Ph.D. thesis, Uppsala
(1999).

\bibitem{S2} Sevostyanov, A., Quantum deformation of Whittaker
modules and the Toda lattice, {\em Duke Math. J.} {\bf 105} (2000),
211--238.

\bibitem{S3} Sevostyanov, A., Towards Drinfeld-Sokolov reduction for quantum groups,  {\em J. Geom. Phys.} {\bf 33}  (2000), 235--256.

\bibitem{S6} Sevostyanov, A., Algebraic group analogues of Slodowy slices and deformations of Poisson W--algebras, {\em Int. Math. Res. Not.} {\bf 2011} (2011), 1880--1925.

\bibitem{S10} Sevostyanov, A., Conjugacy classes in Weyl groups and q-W algebras, {\em Adv. Math.} {\bf 228} (2011), 1315--1376.

\bibitem{SZ} Sevostyanov, A., The geometric meaning of Zhelobenko operators, {\em Transf. Groups.} {\bf 18} (2013), 865–-875.

\bibitem{SDM} Sevostyanov, A., Localization of quantum biequivariant $\mathcal{D}$-modules and
q–W algebras, {\em Proc. London Math. Soc.} {\bf 113} (2016), 583–-626. 

\bibitem{S11} Sevostyanov A., Representations of quantum groups at roots of unity, Whittaker vectors and q-W algebras, {\em J. of Alg.} {\bf 511} (2018), 63--94.

\bibitem{S12} Sevostyanov A., Strictly transversal slices to conjugacy classes in algebraic groups, {\em J. of Alg.} {\bf 520} (2019), 309--366.

\bibitem{S13} Sevostyanov A., The structure of q-W algebras, {\em Transf. Groups} {\bf 25}, no. 1 (2020), 279--304. 

\bibitem{Sk} Skryabin, S., Representations of the Poisson algebra in prime characteristic, {\em Math. Z.} {\bf 243} (2003), 563--597.  

\bibitem{SL} Slodowy, P., {\em Simple Singularities and Simple Algebraic
Groups}, {\em Lecture Notes in Mathematics} {\bf 815}, Springer-Verlag, New York-Berlin
(1980).

\bibitem{Spal} Spaltenstein, N., Classes unipotentes et sous-groupes de Borel,
{\em Lecture Notes in Mathematics} {\bf 946}, Springer-Verlag, Berlin-New York (1982).

\bibitem{Spal1} Spaltenstein, N., On the generalized Springer correspondence for exceptional groups. Algebraic groups and related topics (Kyoto/Nagoya, 1983), {\em Adv. Stud. Pure Math.} {\bf 6}, North-Holland, Amsterdam (1985), 317--338 .

\bibitem{Sp} Springer, T. A., Linear Algebraic Groups, Birkh\"{a}user, Boston (2009).

\bibitem{Sp1} Springer, T. A., Some Results on Algebraic Groups with Involutions. Algebraic Groups and Related Topics (Kyoto/Nagoya, 1983), {\em Adv. Stud. Pure Math.} {\bf 6}, North-Holland, Amsterdam (1985), 525--543.

\bibitem{St1} Steinberg, R., Finite Reflection Groups, {\em Trans. Amer. Math. Soc.} {\bf 91} (1959), 493--504.

\bibitem{St2} Steinberg, R., Regular elements of semisimple
algebraic groups, {\em Publ. Math. I.H.E.S.} {\bf 25} (1965), 49--80.


\bibitem{T} Tolstoy, V. N., Fortieth anniversary of extremal projector method for Lie symmetries.  Noncommutative geometry and representation theory in mathematical physics, in: {\em Contemp. Math.} {\bf 391}, Amer. Math. Soc., Providence, RI (2005), 371--384. 


\bibitem{W}  Weinstein, A., Local structure of Poisson manifolds, {\em J.
Different. Geom.} {\bf 18} (1983), 523--558.

\bibitem{Y} Yakimov, M., On the spectra of quantum groups, {\em Mem. Amer. Math. Soc.} {\bf 229} (2014).

\bibitem{YU} Yun, Z., Minimal reduction type and the Kazhdan-Lusztig map, {\em Indag. Math.}
{\bf 32}, (2021), 1240--1274.

\bibitem{Z} Zelobenko, D. P., Compact Lie groups and their representations, Amer. Math. Soc., Providence (1973).

\bibitem{Z1} Zelobenko, D. P., Extremal cocycles on Weyl groups, {\em Funk. Anal. i ego Pril.}
{\bf 21}, no. 3 (1987), 11--21.

\bibitem{Z3} Zelobenko, D. P., Constructive Modules and Extremal Projectors over Chevalley Algebras, {\em Funk. Anal. i ego Pril.} {\bf 27}, no. 3 (1993),  5--14; {\em Funct. Anal. Appl.} {\bf 27}, no. 3 (1993), 158--165.

\bibitem{Z4} Zelobenko, D. P., S--algebras and Harish--Chandra modules over symmetric Lie algebras, {\em Izv. Akad. Nauk SSSR Ser. Mat.} {\bf 54}, no. 4 (1990),  659--675; {\em Math. USSR-Izv.} {\bf 37}, no 1 (1991), 1--17.

\bibitem{Z5} Zelobenko, D. P., An introduction to the theory of S--algebras over reductive Lie algebras. In: Representation of Lie groups and related topics, {\em Adv. Stud. Contemp. Math.} {\bf 7}, Gordon and Breach, New York (1990), 155--221.

\bibitem{Z6} Zelobenko, D. P., Extremal projectors and generalized Mickelsson algebras over reductive Lie algebras, {\em Izv. Akad. Nauk SSSR Ser. Mat.} {\bf 52}, no. 4 (1988),  758--773; {\em Math. USSR-Izv.} {\bf 33}, no. 1 (1989), 85--100.

\bibitem{Z7} Zelobenko, D. P., S--algebras and Harish--Chandra modules over reductive Lie algebras, {\em Dokl. Akad. Nauk SSSR} {\bf 283}, no 6 (1985),  1306--1308.

\bibitem{Z8} Zelobenko, D. P., Z-algebras over reductive Lie algebras, {\em Dokl. Akad. Nauk SSSR} {\bf 273}, no. 6 (1983), 1301--1304.

\bibitem{Z9} Zelobenko, D. P., S-algebras and Verma modules over reductive Lie algebras, {\em Dokl. Akad. Nauk SSSR} {\bf 273}, no. 4 (1983), 785--788.

\bibitem{Z2} Zelobenko, D. P., Representation Theory of Reductive Lie Algebras (in Russian), Nauka, Moscow (1994).


\end{thebibliography}
\end{document}